\documentclass[12 pt,a4paper,reqno]{amsart}

\makeatletter
\let\@font@info\@gobble
\let\@font@warning\@gobble
\makeatother
\usepackage{amsmath,amsthm,amssymb,mathtools}
\usepackage{bbm}
\usepackage{cancel}							
\usepackage{empheq}							
\usepackage{shuffle}						

\allowdisplaybreaks[1]			    		
\numberwithin{equation}{subsection} 		

\usepackage{fullpage}
\usepackage{float}
\usepackage{adjustbox}						
\usepackage{microtype}						
\usepackage[shortlabels, inline]{enumitem}	
\usepackage{caption}						
\usepackage{subcaption}

\usepackage{array}      					
\newcolumntype{L}{>{$}l<{$}} 				

\usepackage{tikz-cd}
\usepackage{quiver}

\usepackage[backend=biber,style=numeric]{biblatex}

\usepackage{xargs}				
\usepackage{xifthen}			

\usepackage{dsfont}				

\usepackage{bm}					
\usepackage{stmaryrd}			

\usepackage{stackrel}    		

\usepackage{scalerel}[2014/03/10]
\usepackage{stackengine}

\let\originalleft\left
\let\originalright\right
\renewcommand{\left}{\mathopen{}\mathclose\bgroup\originalleft}
\renewcommand{\right}{\aftergroup\egroup\originalright}

\newcommand\reallywidetilde[1]{\ThisStyle{%
		\setbox0=\hbox{$\SavedStyle#1$}%
		\stackengine{-.1\LMpt}{$\SavedStyle#1$}{%
			\stretchto{\scaleto{\SavedStyle\mkern.2mu\sim}{.5467\wd0}}{.7\ht0}%
		}{O}{c}{F}{T}{S}%
	}}

\newcommand{\wt}[1]{\reallywidetilde{#1}}

\newcommand{\extrawidehat}[1]{%
	\savestack{\tmpbox}{\stretchto{%
			\scaleto{%
				\scalerel*[\widthof{\ensuremath{#1}}]                         {\kern-.6pt\bigwedge\kern-.6pt}%
				{\rule[-\textheight/2]{1ex}{\textheight}}
			}{\textheight}%
		}{0.5ex}}%
	\ensurestackMath{\stackon[1pt]{#1}{\tmpbox}}%
}

\newcommand{\id}{\mathrm{id}}
\newcommand{\idd}{\mathds{1}}

\let\Im\relax
\DeclareMathOperator{\Im}{im}

\newcommand{\defeq}{\coloneqq{}}
\newcommand{\eqdef}{=\vcentcolon}
\newcommand{\rest}[2]{\left. #1 \right|_{#2}}

\newcommandx{\Set}[2][2=]{
	\ifthenelse{\isempty{#2}}
	{\left\{ {#1} \right\}}
	{\left\{ {#1}  \, \middle| \, {#2} \right\}}
}

\newcommandx{\fspace}[7][3=, 4=, 5=, 6=, 7={,}]{
	\ifthenelse{\isempty{#5}}
	{\ifthenelse{\isempty{#6}}
		{\ifthenelse{\isempty{#3}}                  
			{{#1}_{#4}^{#2}}
			{{#1}_{#4}^{#2 #7 #3}}}
		{\fspace{#1}{#2}[#3][#4][#6][]}}            
	{\ifthenelse{\isempty{#6}}
		{\ifthenelse{\isempty{#3}}                  
			{{#1}_{#4}^{#2} \left( #5 \right)}               
			{{#1}_{#4}^{#2, #3} \left( #5 \right)}}          
		{\ifthenelse{\isempty{#3}}                  
			{{#1}_{#4}^{#2} \left( #5 #7 #6 \right)}           
			{{#1}_{#4}^{#2, #3} \left( #5 #7 #6 \right)}}}     
}

\newcommandx{\slist}[4][2=, 3=, 4={,}]{
	\ifthenelse{\isempty{#2} \AND \isempty{#3}}
	{#1}
	{\ifthenelse{\isempty{#3}}
		{#1#4 #2}
		{#1#4 #2#4 #3}}}

\newcommandx{\fspacei}[9][2=,3=,4=,5=,6=,7=,8=,9=]{
	\ifthenelse{\isempty{#7}}
	{{#1}^{\slist{#2}[#3][][,]}_{\slist{#4}[#5][#6][,]}}
	{{#1}^{\slist{#2}[#3][][,]}_{\slist{#4}[#5][#6][,]}
		\left( \slist{{\slist{{#7}}[{#8}][][,]}}[{#9}][][;] \right)}
}

\newcommandx{\df}[4][2=,3=,4=]{\fspacei{\Omega}[#4][][][][][{#1}][#3][{#2}]}
\newcommandx{\cdiff}[4][2=,3=,4=,]{\fspacei{A}[#2][][][][][#1][#3][#4]}

\newcommandx{\homo}[4][2=*, 3=, 4=]{\fspacei{H}[][][#2][][][#1][#3][#4]}
\newcommandx{\cohom}[5][2=*,3=, 4=, 5=]{\fspacei{H}[#2][][#5][][][#1][#3][#4]}

\newcommand{\dr}{\text{DR}}

\newcommand{\Cone}[1]{\operatorname{Cone} \left( #1 \right)}

\newcommandx{\inncur}[2][1={\cdot},2={\cdot}]{{\prec {#1}, {#2} \succ}}

\DeclarePairedDelimiter\abs{\lvert}{\rvert}
\newcommandx{\nnorm}[1][1=\cdot]{\left\lVert{#1}\right\rVert}

\DeclareMathOperator{\HomOp}{Hom}
\newcommandx{\Hom}[4][3=,4=]{\fspace{\HomOp}{#3}[][{#4}][{#1}][{#2}]}
\DeclareMathOperator{\EndOp}{End}
\newcommandx{\End}[3][2=,3=]{\fspace{\EndOp}{#2}[][{#3}][{#1}][]}
\DeclareMathOperator{\InnerHomOp}{\ul{Hom}}
\newcommandx{\InnHom}[4][3=,4=]{\fspace{{\InnerHomOp}}{#3}[][{#4}][{#1}][{#2}]}
\DeclareMathOperator{\InnerEndOp}{\ul{End}}
\newcommandx{\InnEnd}[3][2=,3=]{\fspace{{\InnerEndOp}}{#2}[][{#3}][{#1}][]}

\newcommandx{\Mult}[4][3=,4=]{\fspacei{{\ul{\operatorname{Mult}}}}[{#4}][][{#3}][][][{#1}][][{#2}]}

\newcommandx{\Ainf}{A_{\infty}}

\newcommand{\RR}{\mathbb{R}}      
\newcommand{\ZZ}{\mathbb{Z}}      
\newcommand{\QQ}{\mathbb{Q}}	  
\newcommand{\NN}{\mathbb{N}}      
\newcommand{\NZ}{\mathbb{N}_{0}}  
\newcommand{\RPL}{\mathbb{R}_{\geq 0}} 

\newcommandx{\cp}[1][1=n]{\mathbb{CP}^{#1}}	
\newcommandx{\rp}[1][1=n]{\mathbb{RP}^{#1}}     

\newcommandx{\coker}[1]{\operatorname{coker} \left( #1 \right)}

\DeclareMathOperator{\monunit}{\mathbbm{1}}

\DeclareMathOperator{\CoderOp}{\ul{Coder}}
\newcommandx{\Coder}[3][2=,3=]{\fspacei{{\CoderOp}}[#3][][#2][][][#1][][]}
\DeclareMathOperator{\CoDerOp}{\ul{CoDer}}
\newcommandx{\CoDer}[3][2=,3=]{\fspacei{{\CoDerOp}}[#3][][#2][][][#1][][]}

\DeclareFontFamily{U}{mathx}{}
\DeclareFontShape{U}{mathx}{m}{n}{<-> mathx10}{}
\DeclareSymbolFont{mathx}{U}{mathx}{m}{n}
\DeclareMathAccent{\widecheck}{0}{mathx}{"71}

\newcommandx{\tc}[1]{\widecheck{#1}}

\newcommandx{\corest}[1]{#1_{\bullet}}

\newcommandx{\tens}[4][2=,3=,4=]{\fspace{T}{#3}[][#2][{#1}][#4]}
\newcommandx{\tenscyc}[2][2=]{\cycl{T}^{#2} \left( #1 \right)}

\newcommandx{\tensr}[2][2=]{\overline{T}_{#2} \left( #1 \right)}
\newcommandx{\tensrcyc}[2][2=]{\overline{T}_{\textrm{cyc}}^{#2} \left( #1 \right)}

\newcommandx{\tensf}[2][2=]{\widehat{T}_{#2} \left( #1 \right)}

\newcommand{\ootimesbase}{\bm{\otimes}}
\newcommand{\ootimes}{\mathbin{\ootimesbase}}

\newcommandx{\gen}[1]{\left< {#1} \right>}

\newcommandx{\eqcl}[1]{\left[ {#1} \right]}		

\newcommandx{\choch}[2][2=*]{C^{#2}_{\textrm{hoch}} \left( #1 \right)}
\newcommandx{\hhoch}[2][2=*]{HH^{#2} \left( #1 \right)}

\AtBeginDocument{%
	\renewcommand{\t}{\tau} 
}
\newcommand{\N}{N}

\DeclareMathOperator{\TrOp}{Tr}
\newcommandx{\tr}[1][1=]{
	\ifthenelse{\isempty{#1}}
	{\TrOp}
	{\TrOp \left( #1 \right)}
}
\DeclareMathOperator{\TrOpBold}{\mathbf{Tr}}
\newcommandx{\trb}[1][1=]{
	\ifthenelse{\isempty{#1}}
	{\TrOpBold}
	{\TrOpBold \left( #1 \right)}
}

\newcommandx{\cyccl}[1]{\eqcl{#1}_{\operatorname{cyc}}}

\newcommandx{\cconnes}[2][2=*]{C^{#2}_{\lambda} \left( #1 \right)}
\newcommandx{\hcyc}[2][2=*]{HC^{#2} \left( #1 \right)}

\newcommandx{\hcyce}[2][2=*]{
	\ifthenelse{\isempty{#1}}
	{HC_{+}^{#2}}
	{HC_{+}^{#2} \left( #1 \right)}
}
\newcommandx{\hcycred}[2][2=*]{HC^{#2}_{\textrm{red}} \left( #1 \right)}
\newcommandx{\hcycered}[2][2=*]{
	\ifthenelse{\isempty{#1}}
	{HC_{+,\textrm{red}}^{#2}}
	{HC_{+,\textrm{red}}^{#2} \left( #1 \right)}
}

\newcommandx{\degen}[3][1=,2=,3=]{
	\ifthenelse{\isempty{#1}}
	{{\mathcal{D}}}
	{\prescript{}{#2}{{\mathcal{D}}_{#3} \left( \mathcal{#1} \right)}}
}
\newcommandx{\degenu}[3][1=,2=,3=]{
	\ifthenelse{\isempty{#1}}
	{{\ul{\mathcal{D}}}}
	{\prescript{}{#2}{\ul{\mathcal{D}}_{{#3}} \left( \mathcal{#1} \right)}}
}

\newcommandx{\cycbi}[1][1=]{
	\ifthenelse{\isempty{#1}}
	{\mathcal{CC}}
	{\mathcal{CC} \left( \mathcal{#1} \right)}
}

\newcommandx{\cycbiaug}[1][1=]{
	\ifthenelse{\isempty{#1}}
	{\widetilde{\mathcal{CC}}}
	{\widetilde{\mathcal{CC}} \left( \mathcal{#1} \right)}
}

\newcommandx{\cycl}[1]{#1_{\operatorname{cyc}}}

\newcommandx{\degr}[1]{\left| {#1} \right|}		
\newcommandx{\degb}[1]{| {#1} |}			

\newcommandx{\braidop}[2][1=\cdot,2=\cdot]{\left< {#1}, {#2} \right>}
\newcommandx{\braidd}[2]{\braidop[\degb{#1}][\degb{#2}]}
\newcommandx{\braid}[2]{\braidop[#1][#2]}

\newcommandx{\cexp}[1]{\cycl{\exp} \left( #1 \right)}
\newcommandx{\Exp}[1]{e^{#1}}

\newcommand{\s}{s}

\newcommandx{\weight}[1]{\operatorname{w} \left( #1 \right)}

\newcommand{\ul}[1]{\underline{#1}}

\newcommandx{\ndf}[3][2=*,3=*]{\prescript{}{#2}{{\mho}^{#3} \left({#1}\right)}}
\newcommandx{\ndfr}[3][2=*,3=*]{\prescript{}{#2}{{\overline{\mho}}^{#3} \left({#1}\right)}}
\newcommandx{\ncdf}[3][2=*,3=*]{\prescript{}{#2}{{\mho}^{#3}_{\textrm{cyc}} \left({#1}\right)}}
\newcommandx{\ncdfr}[3][2=*,3=*]{\prescript{}{#2}{{\overline{\mho}}^{#3}_{\textrm{cyc}}
		\left({#1}\right)}}
\newcommandx{\ncdfred}[3][2=*,3=*]{\prescript{}{#2}{{\mho}^{#3}_{\textrm{cyc},\textrm{red}}
		\left({#1}\right)}}
\newcommandx{\ncdfrred}[3][2=*,3=*]{\prescript{}{#2}{{\overline{\mho}}^{#3}_{\textrm{cyc},\textrm{red}}
		\left({#1}\right)}}

\DeclareMathOperator{\qdr}{q}

\newcommandx{\cont}[1]{i_{#1}}
\newcommandx{\ccont}[1]{\iota_{#1}}

\newcommandx{\lie}[1]{L_{#1}}
\newcommandx{\clie}[1]{\mathcal{L}_{#1}}

\newcommandx{\indmap}[1]{\mathfrak{#1}}
\newcommandx{\cindmap}[2][2=]{
	\ifthenelse{\isempty{#2}}
	{\cycl{\indmap{#1}}}
	{\cycl{\indmap{#1}}^{\geq #2}}
}
\newcommandx{\cindmape}[2][2=]{
	\ifthenelse{\isempty{#2}}
	{\cycl{\indmap{#1}}^{+}}
	{\cycl{\indmap{#1}}^{+, \geq #2}}
}

\newcommandx{\totcomp}[3][2=,3=*]{
	\ifthenelse{\isempty{#2}}
	{\mathcal{E}^{#3} \left( #1 \right)}
	{\mathcal{E}_{\geq #2}^{#3} \left( #1 \right)}
}
\newcommandx{\totcompe}[3][2=,3=*]{
	\ifthenelse{\isempty{#1}}
	{
		\ifthenelse{\isempty{#2}}
		{\mathcal{E}^{+}}
		{\mathcal{E}_{\geq {#2}}^{+}}
	}
	{
		\ifthenelse{\isempty{#2}}
		{\mathcal{E}^{+} \left( {#1} \right)^{#3}}
		{\mathcal{E}_{\geq {#2}}^{+} \left( {#1} \right)^{#3}}
	}
}

\newcommandx{\totcompred}[3][2=,3=*]{
	\ifthenelse{\isempty{#1}}
	{
		\ifthenelse{\isempty{#2}}
		{\mathcal{E}_{\textrm{red}}}
		{\mathcal{E}_{\geq {#2}, \textrm{red}}}
	}
	{
		\ifthenelse{\isempty{#2}}
		{\mathcal{E}_{\textrm{red}} \left( {#1} \right)^{#3}}
		{\mathcal{E}_{\geq {#2}, \textrm{red}} \left( {#1} \right)^{#3}}
	}
}

\newcommandx{\totcompsred}[3][2=,3=*]{
	\ifthenelse{\isempty{#2}}
	{\mathcal{E}_{\textrm{sred}} \left( {#1} \right)^{#3}}
	{\mathcal{E}_{\geq {#2}, \textrm{sred}} \left( {#1} \right)^{#3}}
}

\newcommandx{\totcompered}[3][2=,3=*]{
	\ifthenelse{\isempty{#1}}
	{
		\ifthenelse{\isempty{#2}}
		{\mathcal{E}^{+}_{\textrm{red}}}
		{\mathcal{E}_{\geq {#2}, \textrm{red}}^{+}}
	}
	{
		\ifthenelse{\isempty{#2}}
		{\mathcal{E}^{+}_{\textrm{red}} \left( {#1} \right)^{#3}}
		{\mathcal{E}_{\geq {#2}, \textrm{red}}^{+} \left( {#1} \right)^{#3}}
	}
}

\newcommandx{\totcompesred}[3][2=,3=*]{
	\ifthenelse{\isempty{#2}}
	{\mathcal{E}^{+}_{\textrm{sred}} \left( {#1} \right)^{#3}}
	{\mathcal{E}_{\geq {#2}, \textrm{sred}}^{+} \left( {#1} \right)^{#3}}
}

\newcommandx{\be}{\mathbf{e}}

\DeclareMathOperator{\evalm}{eval}
\DeclareMathOperator{\evalmf}{\mathfrak{eval}}

\newcommandx{\base}[1]{{#1}_{\operatorname{base}}}
\newcommandx{\mcfunc}[1]{#1_{\star}}

\newcommandx{\mc}[2][2=]{{\fspacei{MC}[][][][][][{#1}][{#2}]}}

\DeclareMathOperator{\tot}{Tot}
\DeclareMathOperator{\totl}{tot}

\newcommandx{\totc}[4][2=*,3=, 4=]{
	\ifthenelse{\isempty{#4}}
	{\tot_{#3} \left( #1 \right)^{#2}}
	{\tot^{#4}_{#3} \left( #1 \right)^{#2}}
}

\newcommandx{\Der}[3][2=,3=]{\ul{\operatorname{Der}}^{#2}_{#3} \left( #1 \right)}

\newcommandx{\resover}[1]{\overline{#1}}
\newcommandx{\rescoho}[1]{\overline{#1}}

\newcommandx{\resunder}[1]{\ul{#1}}

\AtBeginDocument{%
	\renewcommandx{\H}[3][2=,3=]{H_{#2}^{#3} \left( #1 \right)}%
}

\newcommandx{\G}[3][2=,3=]{G_{#2}^{#3} \left( #1 \right)}
\newcommandx{\Go}[3][2=,3=]{\overline{G}_{#2}^{#3} \left( #1 \right)}
\newcommandx{\Ho}[3][2=,3=]{\overline{H}_{#2}^{#3} \left( #1 \right)}

\newcommandx{\SP}[2][1=,2=]{
	\ifthenelse{\isempty{#1}}
	{\Omega_{#2}}
	{\Omega_{{#2}} \left( {#1} \right)}
}

\newcommandx{\TM}[2][2=]{
	\ifthenelse{\isempty{#1}}
	{\Theta_{#2}}
	{\Theta_{{#2}} \left( {#1} \right)}
}

\newcommandx{\ulz}[1]{\ul{ {#1}_0 \left( 1 \right) }}

\newcommandx{\clsub}[1]{{#1}_{\textrm{cl}}}
\newcommandx{\clsup}[1]{#1^{\textrm{cl}}}

\newcommandx{\imsub}[1]{{#1}_{\textrm{im}}}
\newcommandx{\quotsub}[1]{{#1}_{\textrm{quot}}}

\newcommandx{\Alg}[1][1=]{
	\ifthenelse{\isempty{#1}}
	{\mathbf{Alg}}
	{\mathbf{Alg} \left( #1 \right)}
}

\newcommandx{\CAlg}[1][1=]{
	\ifthenelse{\isempty{#1}}
	{\mathbf{CAlg}}
	{\mathbf{CAlg} \left( #1 \right)}
}

\newcommandx{\CoAlg}[1][1=]{
	\ifthenelse{\isempty{#1}}
	{\mathbf{CoAlg}}
	{\mathbf{CoAlg} \left( #1 \right)}
}

\newcommand{\B}{\mathbf{B}} 

\newcommand{\Ab}{\mathbf{Ab}}
\newcommand{\SNAb}{\mathbf{SNAb}}
\newcommand{\BAb}{\mathbf{BAb}}

\newcommand{\Ring}{\mathbf{Ring}}
\newcommand{\SNRing}{\mathbf{SNRing}}
\newcommand{\BRing}{\mathbf{BRing}}

\newcommand{\GRing}{\mathbf{GRing}}
\newcommand{\GSNRing}{\mathbf{GSNRing}}
\newcommand{\GBRing}{\mathbf{GBRing}}

\newcommandx{\PDGSNAlg}[1][1=]{
	\ifthenelse{\isempty{#1}}
	{\mathbf{PDGSNAlg}}
	{\mathbf{PDGSNAlg} \left( #1 \right)}
}

\newcommandx{\GSNAlg}[1][1=]{
	\ifthenelse{\isempty{#1}}
	{\mathbf{GSNAlg}}
	{\mathbf{GSNAlg} \left( #1 \right)}
}

\newcommandx{\GBAlg}[1][1=]{
	\ifthenelse{\isempty{#1}}
	{\mathbf{GBAlg}}
	{\mathbf{GBAlg} \left( #1 \right)}
}

\newcommandx{\PDGBAlg}[1][1=]{
	\ifthenelse{\isempty{#1}}
	{\mathbf{PDGBAlg}}
	{\mathbf{PDGBAlg} \left( #1 \right)}
}

\newcommandx{\Mod}[2][1=,2=]{
	\ifthenelse{\isempty{#1}}
	{\mathbf{Mod}}
	{\mathbf{Mod}_{#2} \left( #1 \right)}
}

\newcommandx{\SNMod}[1][1=]{
	\ifthenelse{\isempty{#1}}
	{\mathbf{SNMod}}
	{\mathbf{SNMod} \left( #1 \right)}
}

\newcommandx{\BMod}[2][1=,2=]{
	\ifthenelse{\isempty{#1}}
	{\mathbf{BMod}}
	{\ifthenelse{\isempty{#2}}
		{\mathbf{BMod} \left( #1 \right)}
		{\mathbf{BMod} ( #1 )}
	}
}

\newcommandx{\PDGMod}[1][1=]{
	\ifthenelse{\isempty{#1}}
	{\mathbf{PDGMod}}
	{\mathbf{PDGMod} \left( #1 \right)}
}

\newcommandx{\PDGSNMod}[1][1=]{
	\ifthenelse{\isempty{#1}}
	{\mathbf{PDGSNMod}}
	{\mathbf{PDGSNMod} \left( #1 \right)}
}

\newcommandx{\PDGBMod}[1][1=]{
	\ifthenelse{\isempty{#1}}
	{\mathbf{PDGBMod}}
	{\mathbf{PDGBMod} \left( #1 \right)}
}

\newcommandx{\DGMod}[1][1=]{
	\ifthenelse{\isempty{#1}}
	{\mathbf{DGMod}}
	{\mathbf{DGMod} \left( #1 \right)}
}

\newcommandx{\DGSNMod}[1][1=]{
	\ifthenelse{\isempty{#1}}
	{\mathbf{DGSNMod}}
	{\mathbf{DGSNMod} \left( #1 \right)}
}

\newcommandx{\DGBMod}[1][1=]{
	\ifthenelse{\isempty{#1}}
	{\mathbf{DGBMod}}
	{\mathbf{DGBMod} \left( #1 \right)}
}

\newcommandx{\GMod}[2][1=,2=]{
	\ifthenelse{\isempty{#1}}
	{\mathbf{GMod}}
	{
		\ifthenelse{\isempty{#2}}
		{\mathbf{GMod} \left( #1 \right)}
		{\mathbf{GMod}_{{#2}} \left( #1 \right)}
	}
}
\newcommandx{\GSNMod}[1][1=]{
	\ifthenelse{\isempty{#1}}
	{\mathbf{GSNMod}}
	{\mathbf{GSNMod} \left( #1 \right)}
}

\newcommandx{\GBMod}[2][1=,2=]{
	\ifthenelse{\isempty{#1}}
	{\mathbf{GBMod}}
	{\ifthenelse{\isempty{#2}}
		{\mathbf{GBMod} \left( #1 \right)}
		{\mathbf{GBMod} ( #1 )}
	}
}

\newcommandx{\clball}[2]{{#1}^{\bullet} \left( #2 \right)}

\DeclareMathOperator{\ord}{ord}

\DeclareMathOperator{\triv}{triv}
\newcommandx{\trivnorm}[1][1=\cdot]{\nnorm[#1]_{\triv}}

\newcommand{\cotimes}{\mathbin{\widehat{\otimes}}} 

\newcommand{\cootimes}{\mathbin{\widehat{\ootimesbase}}} 

\DeclareMathOperator{\coplus}{\widehat{\oplus}}
\DeclareMathOperator*{\cbigoplus}{\widehat{\bigoplus}}

\newcommand{\GG}{\mathbb{G}}

\newcommandx{\pows}[2][2=z]{{{#1} \llbracket #2 \rrbracket}}

\DeclareMathOperator{\hor}{hor}
\DeclareMathOperator{\ver}{ver}

\newcommand{\go}{\mathbbm{1}}

\makeatletter
\newcommand{\subalign}[1]{%
	\vcenter{%
		\Let@ \restore@math@cr \default@tag
		\baselineskip\fontdimen10 \scriptfont\tw@
		\advance\baselineskip\fontdimen12 \scriptfont\tw@
		\lineskip\thr@@\fontdimen8 \scriptfont\thr@@
		\lineskiplimit\lineskip
		\ialign{\hfil$\m@th\scriptstyle##$&$\m@th\scriptstyle{}##$\hfil\crcr
			#1\crcr
		}%
	}%
}
\makeatother

\newcommandx{\eqwithref}[2][1=,2=]{
	\ifthenelse{\isempty{#1}}
	{\hspace{0.5cm} = \hspace{0.5cm}}
	{
		\ifthenelse{\isempty{#2}}
		{\hspace{0.5cm} \stackrel{\mathclap{\eqref{#1}}}{=} \hspace{0.5cm}}
		{\hspace{0.5cm} \stackrel[\mathclap{\eqref{#2}}]{\mathclap{\eqref{#1}}}{=} \hspace{0.5cm}}
	}
}

\newcommandx{\eqwithtext}[2][1=,2=]{
	\ifthenelse{\isempty{#1}}
	{\hspace{0.5cm} = \hspace{0.5cm}}
	{
		\ifthenelse{\isempty{#2}}
		{\hspace{0.5cm} \stackrel{\mathclap{\text{#1}}}{=} \hspace{0.5cm}}
		{\hspace{0.5cm} \stackrel[\mathclap{\text{#2}}]{\mathclap{\text{#1}}}{=} \hspace{0.5cm}}
	}
}

\newcommand{\pstar}{\mathop{*}} 
\newcommand{\hpstar}{\mathop{\hat{\pstar}}}

\ExplSyntaxOn
\NewDocumentCommand{\mur}{ O{\mu} O{b} m }
{
	#1\sb{#3} \left(
	\int_compare:nNnTF { #3 } = { 0 }
	{ 1 }
	{
		#2
		\prg_replicate:nn { #3 - 1 } { , #2 }
	}
	\right)
}
\ExplSyntaxOff

\PassOptionsToPackage{hyphens}{url}
\usepackage[bookmarksopen,bookmarksdepth=3]{hyperref}
\usepackage{pdflscape}						
\usepackage{cleveref}

\AddToHook{cmd/appendix/before}{%
  \crefalias{section}{appendix}%
  \crefalias{subsection}{appendix}
  \crefalias{subsubsection}{appendix}
}

\crefname{section}{Section}{Sections}
\crefname{figure}{Figure}{Figures}
\crefname{appendix}{Appendix}{Appendices}
\crefname{table}{Table}{Tables}

\crefname{footnote}{Footnote}{Footnotes}
\crefformat{footnote}{#2Footnote~#1#3}
\Crefformat{footnote}{#2Footnote~#1#3}

\newtheorem{thmx}{Theorem}

\crefname{thmx}{Theorem}{Theorems}

\newtheorem{thm}{Theorem}[subsection]
\crefname{thm}{Theorem}{Theorems}

\newtheorem{prop}[thm]{Proposition}
\crefname{prop}{Proposition}{Propositions}

\newtheorem{lm}[thm]{Lemma}
\crefname{lm}{Lemma}{Lemmas}

\newtheorem{cor}[thm]{Corollary}
\crefname{cor}{Corollary}{Corollaries}

\crefname{conj}{Conjecture}{Conjectures}

\theoremstyle{definition}
\newtheorem{dfn}[thm]{Definition}
\crefname{dfn}{Definition}{Definitions}

\newtheorem{ex}[thm]{Example}
\crefname{ex}{Example}{Examples}

\theoremstyle{remark}
\newtheorem{rem}[thm]{Remark}
\crefname{rem}{Remark}{Remarks}

\newif\ifjake
\jaketrue 

\ifjake
  \AddToHook{env/prop/begin}{\crefalias{thm}{prop}}
  \AddToHook{env/cor/begin}{\crefalias{thm}{cor}}
  \AddToHook{env/lm/begin}{\crefalias{thm}{lm}}
  \AddToHook{env/dfn/begin}{\crefalias{thm}{dfn}}
  \AddToHook{env/ex/begin}{\crefalias{thm}{ex}}
  \AddToHook{env/rem/begin}{\crefalias{thm}{rem}}
\fi

\usepackage{xr} 							
\usepackage{subfiles}						

\makeatletter
\ifSubfilesClassLoaded{%
	\addbibresource{../superpotential.bib}
  \externaldocument{../superpotential}
  
  \global\expandafter\let\csname r@tocindent-1\endcsname\relax
  \global\expandafter\let\csname r@tocindent0\endcsname\relax
  \global\expandafter\let\csname r@tocindent1\endcsname\relax
  \global\expandafter\let\csname r@tocindent2\endcsname\relax
  \global\expandafter\let\csname r@tocindent3\endcsname\relax
  \global\expandafter\let\csname r@tocindent4\endcsname\relax

  \AtBeginDocument{%
    \expandafter\def\csname r@tocindent-1\endcsname{0pt}%
    \expandafter\def\csname r@tocindent0\endcsname{0pt}%
    \expandafter\def\csname r@tocindent1\endcsname{0pt}%
    \expandafter\def\csname r@tocindent2\endcsname{0pt}%
    \expandafter\def\csname r@tocindent3\endcsname{0pt}%
    \expandafter\def\csname r@tocindent4\endcsname{0pt}%
    
    \setbox0=\vbox{%
      \title{}%
      \maketitle
    }%
  }%
  \AtEndDocument{\printbibliography}
}{%
	\addbibresource{superpotential.bib}
  
  \AtBeginDocument{%
    \@ifundefined{r@tocindent-1}{\expandafter\def\csname r@tocindent-1\endcsname{0pt}}{}%
    \@ifundefined{r@tocindent0}{\expandafter\def\csname r@tocindent0\endcsname{0pt}}{}%
    \@ifundefined{r@tocindent1}{\expandafter\def\csname r@tocindent1\endcsname{0pt}}{}%
    \@ifundefined{r@tocindent2}{\expandafter\def\csname r@tocindent2\endcsname{0pt}}{}%
    \@ifundefined{r@tocindent3}{\expandafter\def\csname r@tocindent3\endcsname{0pt}}{}%
    \@ifundefined{r@tocindent4}{\expandafter\def\csname r@tocindent4\endcsname{0pt}}{}%
  }%
}
\makeatother

\begin{document}

\title{Total Inner Products and Cyclic Chern--Simons Forms}
\keywords{$\Ainf$-algebra, curved $\Ainf$-algebra, Banach $\Ainf$-algebra, cyclic $\Ainf$-algebra,
	cyclic homology, non-commutative cyclic differential forms, homotopy inner products,
	superpotential, Chern--Simons theory, Gromov--Witten theory, bounding cochain, pseudoisotopy}
\subjclass[2020]{18G70, 16E40 (Primary), 53D45, 53D37, 16W80 (Secondary)}

\date{August 2026}

\author[P. Giterman]{Pavel Giterman}
\address{Institute of Mathematics\\ Hebrew University, Givat Ram\\Jerusalem, 91904, Israel}
\email{pavel.giterman@mail.huji.ac.il}
\author[J. Solomon]{Jake P. Solomon}
\address{Institute of Mathematics\\ Hebrew University, Givat Ram\\Jerusalem, 91904, Israel}
\email{jake@math.huji.ac.il}

\begin{abstract}
	We introduce the notion of a total inner product on an $\Ainf$-algebra, and a cyclic Chern--Simons
	form associated to a topologically nilpotent element. The total inner product applied to the cyclic
	Chern--Simons form gives a superpotential function that is gauge invariant on solutions of the
	Maurer--Cartan equation known as bounding cochains. The derivative of the superpotential is computed.
	The definition of total inner product replaces strict symmetries that appear in previous notions of
	homotopy inner products with symmetries up to an infinite family of coherent homotopies.
	We show how to recover previous notions of homotopy inner products as special cases of total inner
	products. Total inner products are designed to facilitate the definition of descendent open
	Gromov--Witten invariants.

	The definitions of the total inner product and the cyclic Chern--Simons form use the total complex
	of cyclic codifferential forms, which gives a chain model for cyclic homology. We give explicit formulas
	for homotopy equivalences with other known models. We discuss also the notion of an $\infty$-trace,
	which arises from Connes' cyclic complex, and the associated $\infty$-modulus, which gives another
	gauge-invariant function on bounding cochains. In open Gromov--Witten theory, the $\infty$-modulus is
	used to normalize the bounding cochains to which the superpotential is applied.

	We develop our definitions for general curved Banach $\Ainf$-algebras, and give both unital and non-unital
	versions of the main results. In the unital setting, we work with weak bounding cochains, solutions of an
	inhomogeneous Maurer--Cartan equation involving the unit. Weak bounding cochains arise in the open Gromov--Witten theory of
	Lagrangian submanifolds with non-vanishing Maslov class. The derivative of the superpotential at a weak bounding
	cochain is related to the derivative of the $\infty$-modulus.
\end{abstract}

\maketitle

\makeatletter
\def\l@paragraph{\@tocline{4}{0pt}{1pc}{7pc}{}}
\def\l@subparagraph{\@tocline{5}{0pt}{1pc}{7pc}{}}
\makeatother

\setcounter{tocdepth}{3}

\tableofcontents


\section{Introduction} \label{sec:introduction}

\setcounter{equation}{0}
\renewcommand{\theequation}{\arabic{equation}}

\subsection{Overview}
The concept of a cyclic structure on an $\Ainf$-algebra, i.e., a cyclically invariant
inner product, emerged in the early 1990s at the intersection of symplectic geometry,
string theory, and operad theory \cite{Getzler1995,Gaberdiel1997,Kontsevich1994}.
In symplectic geometry, the de Rham model of the Fukaya $\Ainf$-algebra~\cite{Fukaya2009} for a Lagrangian submanifold is equipped
with a cyclic structure given by integration~\cite{Fukaya2010}.
The cyclic structure is used to define the superpotential function on topologically nilpotent elements of the $\Ainf$-algebra.
The superpotential restricts to a gauge-invariant function on solutions to the Maurer--Cartan equation, known as bounding cochains.
Evaluating the superpotential on a canonically chosen gauge equivalence class of bounding cochains gives a generating
series for open Gromov--Witten invariants~\cite{Fukaya2011,Joyce2008,Solomon2016,Solomon2016a}.

The current work is part of a broader program to define gravitational descendent open Gromov--Witten invariants~\cite{Solomon2026}.
Once gravitational descendents are incorporated in the Fukaya $\Ainf$-algebra, strict cyclic invariance of the inner product no longer holds.
Homotopy versions of cyclic invariance for inner products have been studied widely~\cite{Cho2008,cho-homotopy-superpotential,Kontsevich:vh,Tradler2008}.
The gauge invariance of the superpotential in the homotopy cyclic setting depends on a new strict symmetry, called closedness~\cite{Cho2008,cho-homotopy-superpotential,Kontsevich:vh}.  It is not clear how to construct geometrically an inner product that satisfies this strict symmetry.

We formulate a homotopy version of closedness together with cyclic invariance
giving the notion of a \textbf{total inner product}, which admits a geometric construction
on the gravitational descendent Fukaya $\Ainf$-algebra~\cite{Solomon2026}.
We define a \textbf{superpotential} function associated to a total inner product and prove its naturality
and gauge invariance.
Furthermore, we compute the formal derivative of the superpotential and show that it vanishes on strong bounding cochains.
The formula for the derivative of the superpotential in the strict cyclic setting plays a crucial role in the derivation of the open WDVV equations~\cite{Solomon2024}. It is expected to play a similar role in the derivation of topological recursion relations for open descendent Gromov--Witten invariants, where strict cyclic invariance is not available. Topological recursion relations for open descendent integrals were proved in~\cite{Pandharipande2024}.

Total inner products are built on the framework of cyclic codifferential forms, a variant of the construction in~\cite{Herscovich2023}, which is dual to the cyclic differential forms introduced by Kontsevich and Soibelman~\cite{Kontsevich:vh}.
For unital $\Ainf$-algebras, the total complex of cyclic codifferential forms of degree greater than or equal to two
gives a chain model for cyclic homology. A total inner product is a linear functional on an extension of the total complex that is a chain map.
The superpotential is obtained by applying the total inner product to a canonical \textbf{cyclic Chern--Simons form} associated to a topologically nilpotent element
of the $\Ainf$-algebra. When the element is a bounding cochain, the cyclic Chern--Simons form is closed.
The extension of the total complex is introduced to accommodate the gauge invariance of the superpotential for curved $\Ainf$-algebras.

A proper Calabi--Yau structure on an $\Ainf$-algebra without curvature is a linear functional
on the cyclic homology satisfying a certain non-degeneracy condition~\cite{Kontsevich:vh}.
In the curved setting, it is natural to consider appropriate extensions of cyclic homology.
Thus, a Calabi--Yau structure can be realized at chain-level by a total inner product.
Another chain-level realization, previously considered in the literature~\cite{Shklyarov2017}, is obtained by using Connes' complex.
We define an $\infty$-\textbf{trace} to be a linear functional on an extension of Connes' cyclic complex that is a chain map.
By analogy with the definition of the superpotential, we define the $\infty$-\textbf{modulus} to be
the value of the $\infty$-trace on a canonical cyclic chain, the \textbf{cyclic exponential}, associated
with a topologically nilpotent element of the $\Ainf$-algebra.
The cyclic exponential is closed when the element is a bounding cochain.
We prove that the $\infty$-modulus is natural and gauge invariant.
The $\infty$-modulus is used to choose the canonical bounding cochain upon which the superpotential is evaluated in the definition of open Gromov--Witten descendent invariants~\cite{Solomon2026} and also in the context of matrix factorizations~\cite{Sela2024}.

We prove versions of the above results both for unital and non-unital $\Ainf$-algebras.
In the non-unital case, we consider strong bounding cochains, which are solutions to the
homogeneous Maurer--Cartan equation. In the unital case, we consider weak bounding cochains,
which are solutions to an inhomogeneous Maurer--Cartan equation involving the unit. Weak bounding
cochains play a crucial role in the definition of open Gromov--Witten invariants for Lagrangians
with non-vanishing Maslov class~\cite{Solomon2016a}. The cyclic Chern--Simons form and
cyclic exponential associated to a weak bounding cochain are closed in appropriate reduced
versions of the total complex and Connes' complex respectively.
Unitality conditions for total inner products and $\infty$-traces are formulated by requiring
that they descend to the reduced complexes. The unitality conditions are used in the proof of
gauge invariance for the superpotential and the $\infty$-modulus.

Using the total complex instead of its extension gives a notion of pre-total inner product.
The pre-total inner product contains slightly less information. In the case of
Fukaya $\Ainf$-algebras, the missing information is a version of the
$\mathfrak{m}_{-1}$ term that appears in~\cite{Fukaya2011,Joyce2008,Solomon2016,Solomon2016a}.
Similarly, using Connes' complex instead of its extension gives a notion of pre-$\infty$ trace.
We give a formula for a pre-$\infty$-trace in terms of a pre-total inner product and vice versa.
When the pre-total inner product satisfies the strict closedness symmetry of~\cite{Cho2008,cho-homotopy-superpotential,Kontsevich:vh}, these formulas respect unitality conditions.
Furthermore, the derivative of the superpotential on weak bounding cochains coincides
with a multiple of the derivative of the modulus function.

We give a detailed treatment of curved Banach $\Ainf$-algebras over a differential graded-commutative ground algebra,
used to define our notion of pseudoisotopy and gauge equivalence of bounding cochains over a pseudoisotopy.
We also study various complexes computing the cyclic homology and give
explicit chain-level homotopy equivalences between them, which may be of independent interest.

\subsection{Statement of Main Results} \label{sec:statement-results}
Let $\mathbbm{k}$ be a field of characteristic zero which is fixed for the rest of the section.
In what follows, we freely use the notions of graded Banach objects, which are graded objects (algebras, modules, coalgebras, etc.)
equipped with a non-Archimedean norm $\nnorm$ making them complete. Categorical constructions are done in the appropriate
categories of Banach objects and notions such as direct sum $\oplus$ (resp.\ tensor product $\otimes$) are to be interpreted
as the complete direct sum (resp.\ complete tensor product). For precise definitions, we refer to \cref{sec:non-archimedean-graded-setting}.
We use the notation $f \colon M \rightharpoonup N$ to denote graded maps of arbitrary degree,
reserving the notation $f \colon M \rightarrow N$ for degree zero maps. Objects equipped
with differentials are denoted with calligraphic letters, as in $\mathcal{R} = \left( R, d \right)$,
reserving the corresponding letter $R$ for the underlying graded object.

\subsubsection{\texorpdfstring{$\Ainf$-Algebra}{A-infinity Algebra} Preliminaries}

Let $R$ be a graded-commutative Banach $\mathbbm{k}$-algebra. Given a graded Banach $R$-module $A$, we denote by
$\tens{A} = \oplus_{i=0}^{\infty} A^{\otimes i}$ the tensor module, where the tensor product is taken over $R$.
Given a map $f \colon \tens{A} \rightharpoonup \tens{B}$, the corestriction $\corest{f} \colon \tens{A} \rightharpoonup B$
of $f$ is defined to be the composition of $f$ with the projection $\tens{B} \twoheadrightarrow B$.
When equipped with the deconcatenation coproduct, $\tens{A}$ has the structure of a graded Banach coalgebra, called the tensor coalgebra.
An element $a \in A^0$ is called topologically nilpotent if it satisfies
$\nnorm[a^{\otimes n}] \to 0$, and the set of topologically nilpotent elements of $A$ is denoted by $\tc{A}$.
Given $a \in \tc{A}$, we denote by $\Exp{a} = \sum_{n = 0}^{\infty} a^{\otimes n}$ the exponential
of $a$ in the sense of \cite{Fukaya2009}, which converges in $\tens{A}$.

Let $\mathcal{R} = \left( R, d \right)$ be a differential graded-commutative Banach $\mathbbm{k}$-algebra.
A Banach $\Ainf$-algebra $\mathcal{A} = \left( A, \mu \right)$ over $\mathcal{R}$ is a $\ZZ$-graded Banach
$R$-module $A$, together with a degree one coderivation $\mu \colon \tens{A} \rightharpoonup \tens{A}$
of the tensor coalgebra which satisfies $\mu^2 = 0$ and is compatible with $d$ in the sense that
$\mu \left( r \cdot a \right) = dr \cdot a + (-1)^{\degb{r}} r \cdot \mu \left( a \right)$
for $r \in R$ and $a \in A$. We also require that $\nnorm[\mu] \leq 1$ and $\nnorm[\mu_0 \left( 1 \right)] < 1$.
The coderivation $\mu$ is completely determined by its corestriction $\corest{\mu} \colon \tens{A} \rightharpoonup A$, which
encodes a sequence of degree one operations $\mu_k \colon A^{\otimes k} \rightharpoonup A$ for $k \geq 0$. The relation $\mu^2 = 0$
is equivalent to the $\Ainf$-identities
\begin{align*}
	\sum_{k_1 + k_2 + k_3 = k} \mu_{k_1 + 1 + k_3} \circ \left( \idd^{\otimes k_1} \otimes \mu_{k_2} \otimes \idd^{\otimes k_3} \right) = 0
\end{align*}
for $k \geq 0$.

A morphism $f \colon \mathcal{A} \rightarrow \mathcal{B}$ between two Banach $\Ainf$-algebras $\mathcal{A} = \left( A, \mu \right)$
and $\mathcal{B} = \left( B, \nu \right)$ is by definition a morphism $f \colon \tens{A} \rightarrow \tens{B}$ of Banach coalgebras
which satisfies $f \circ \mu = \nu \circ f$. We note that in the Banach setting, it is possible for a morphism $f$ to have a
non-zero topologically nilpotent ``change of connection'' component $f_0 \left( 1 \right) \in \tc{B}$
(see \cref{subsubsec:morphisms-formal-tensor-coalgebra}) and we require that $\nnorm[f_0 \left( 1 \right)] < 1$.
The morphism $f$ is determined uniquely by its corestriction $\corest{f} \colon \tens{A} \rightarrow B$, encoding a sequence of
operations $f_k \colon A^{\otimes k} \rightarrow B$ for $k \geq 0$ and the condition for $f$ to be a morphism can be written
explicitly in terms of the corestrictions of $f, \mu$ and $\nu$ (see \cref{eq:ainf_morphism_explicit}).

The notion of a morphism can be extended naturally to encompass morphisms $f \colon \mathcal{A} \rightarrow \mathcal{B}$
of  Banach $\Ainf$-algebras defined over different ground differential graded-commutative Banach $\mathbbm{k}$-algebras
(see \cref{sub:graded-seminormed-banach-coalgebras}).
In this case, the morphism $f$ comes equipped with an underlying morphism $\base{f} \colon \mathcal{R} \rightarrow \mathcal{S}$
of differential graded-commutative Banach $\mathbbm{k}$-algebras, determined uniquely by $f$.

Given a Banach $\Ainf$-algebra $\mathcal{A} = \left( A, \mu \right)$, an element $e \in A^{-1}$ is called a unit
for $\mathcal{A}$ if $\nnorm[e] \leq 1$, $\mu_2 \left( e,a \right) = (-1)^{\degb{a} + 1} \mu_2 \left( a, e \right) = a$ for all $a \in A$,
and, in addition, $\mu_k \left( a_1, \dots, a_k \right) = 0$ whenever $k \neq 2$ and $a_i = e$ for some $1 \leq i \leq k$.
When $e$ is a unit for $\mathcal{A}$, the triple $\mathcal{A} = \left( A,\mu, e \right)$ is called a
unital Banach $\Ainf$-algebra.

\subsubsection{Bounding Cochains and Gauge Equivalence}
Let $\mathcal{A} = \left( A, \mu, e \right)$ be a unital Banach $\Ainf$-algebra over $\mathcal{R} = \left( R, d \right)$.
Given an element $c \in R^2$ with $dc = 0$ and $\nnorm[c] < 1$, an element $b \in A^0$ with $\nnorm[b] < 1$
is called a \textbf{weak bounding cochain} if it satisfies the inhomogeneous Maurer--Cartan equation
\begin{equation*}
	\corest{\mu} \left( \Exp{b} \right) = \sum_{i=0}^{\infty} \mu_i \left( b^{\otimes i} \right) = c \cdot e.
\end{equation*}
We will denote the set of all weak bounding cochains in $A$ with a fixed $c$ by  $\mc{\mathcal{A}}[c]$.
When $c = 0$, $b$ is called a
\textbf{strong bounding cochain}, and the set of all strong bounding cochains in $A$ is denoted by $\mc{\mathcal{A}}$.
Strong bounding cochains are also defined when $\mathcal{A}$ is non-unital.

Given a morphism $f \colon \mathcal{A} \rightarrow \mathcal{B}$ of Banach $\Ainf$-algebras, possibly over different bases,
there is a natural way to pushforward bounding cochains, or more generally, topologically nilpotent elements.
The pushforward map $\mcfunc{f} \colon \tc{A} \rightarrow \tc{B}$ is given by
\begin{equation*}
	\mcfunc{f} \left( b \right) = \corest{f} \left( \Exp{b} \right) = \sum_{i=0}^{\infty} f_i \left( b^{\otimes i} \right).
\end{equation*}
The pushforward of a strong bounding cochain is a strong bounding cochain, and the pushforward of a weak
bounding cochain $b \in \mc{\mathcal{A}}[c]$ along a unital morphism is a weak bounding cochain
$\mcfunc{f} \left( b \right) \in \mc{\mathcal{B}}[\base{f}(c)]$.

Given two Banach $\Ainf$-algebras $\mathcal{A}_0, \mathcal{A}_1$ over the same differential graded-commutative
$\mathbbm{k}$-algebra $\mathcal{S}$, a \textbf{pseudoisotopy} between $\mathcal{A}_0$ and $\mathcal{A}_1$
is given by the following data:
\begin{enumerate}
	\item A differential graded-commutative Banach $\mathbbm{k}$-algebra $\mathfrak{R}$.
	\item A Banach $\Ainf$-algebra $\mathfrak{A}$ over $\mathfrak{R}$ together with two morphisms
	      $\evalmf^i \colon \mathfrak{A} \rightarrow \mathcal{A}_i$
	      such that the underlying DGA morphisms
	      $\evalm^i \defeq \base{\evalmf}^i \colon \mathfrak{R} \rightarrow \mathcal{S}$
	      are homotopic as maps of differential graded $\mathbbm{k}$-modules.
\end{enumerate}
We will often succinctly denote the data of a pseudoisotopy between
$\mathcal{A}_0$ and $\mathcal{A}_1$ by $\mathfrak{A}$ and say that $\mathfrak{A}$
is a pseudoisotopy between
$\mathcal{A}_0$ and $\mathcal{A}_1$, leaving the base differential graded algebras and the morphisms implicit.
When the $\Ainf$-algebras $\mathcal{A}_0, \mathcal{A}_1$ and $\mathfrak{A}$ are unital,
and the morphisms $\evalmf^i$ are also unital, we say that the pseudoisotopy $\mathfrak{A}$ is \textbf{unital}.

Let $\mathcal{A}_0$ and $\mathcal{A}_1$ be two unital Banach
$\Ainf$-algebras over $\mathcal{S}$ and assume we have a fixed unital pseudoisotopy $\mathfrak{A}$
over $\mathfrak{R} = \left( R, d_R \right)$ between $\mathcal{A}_0$ and $\mathcal{A}_1$.
Let $b_0 \in \mc{\mathcal{A}_0}[c_0]$ and $b_1 \in \mc{\mathcal{A}_1}[c_1]$
be two weak bounding cochains. We will say that $b_0$ and $b_1$ are
$\mathfrak{A}$-\textbf{gauge-equivalent} if there exists an element $c \in R^2$ with
$d_R \left( c \right) = 0$ and $\nnorm[c] < 1$ such that
\begin{equation*}
	\evalm^0 \left( c \right) = c_0, \,\,\,
	\evalm^1 \left( c \right) = c_1
\end{equation*}
and an element $b \in \mc{\mathfrak{A}}[c]$ such that
\begin{equation} \label{eq:intro-gauge-equivalence-bounding-chains}
	\mcfunc{\evalmf}^0 \left( b \right) = b_0, \,\,\,
	\mcfunc{\evalmf}^1 \left( b \right) = b_1.
\end{equation}
When $b_0 \in \mc{\mathcal{A}_0}$ and $b_1 \in \mc{\mathcal{A}_1}$ are strong bounding cochains, we say that $b_0$
and $b_1$ are $\mathfrak{A}$-\textbf{gauge-equivalent} if there exists $b \in \mc{\mathfrak{A}}$ such that
\cref{eq:intro-gauge-equivalence-bounding-chains} holds.

\subsubsection{Noncommutative Cyclic Codifferential Forms}
Let $\mathcal{A} = \left( A, \mu \right)$ be a Banach $\Ainf$-algebra over $\mathcal{R} = \left( R, d \right)$.
In what follows, we will work with $\ZZ^2$-graded objects, identifying $\ZZ$-graded objects such
as $A^{*}$ and $R^{*}$ with the corresponding $\ZZ^2$-graded objects concentrated in bidegree $(0, *)$.
When working with bigraded objects, we adhere in this section to the Koszul sign conventions using the parity form
$\braidop \colon \ZZ^2 \times \ZZ^2 \rightarrow \ZZ_2$ given by
\begin{equation*}
	\braidop[\left( a_1, a_2 \right)][\left( b_1, b_2 \right)] = a_1 \cdot b_1 + a_2 \cdot b_2 \mod 2,
\end{equation*}
i.e., whenever exchanging two elements of bidegrees $a = (a_1,a_2)$ and $b = (b_1,b_2)$,
we introduce the sign $(-1)^{\braidop[a][b]}$.

Denote by $\ul{A}$ the shift $\ul{A} = A[(-1,0)]$ and let
$\ndf{A}[*][*] \defeq \tens{A \oplus \ul{A}}^{(*,*)}$ be the
tensor module on $A \oplus \ul{A}$.
An elementary tensor $x \in \ndf{A}[k][n]$ has the form
\begin{equation*}
	x = l^0 \otimes \ul{a_1} \otimes l^1 \otimes \dots \otimes \ul{a_k} \otimes l^{k},
\end{equation*}
where $l^0, \dots, l^{k} \in \tens{A}$ are elementary tensors, $a_1, \dots, a_k \in A$, and
$\sum_{i=0}^{k} \degr{l^i} + \sum_{i=1}^k \degr{a_i} = n$.
When no confusion is possible, we omit the
tensor product symbol for brevity and write $x = l^0 \, \ul{a_1} \, l^1 \, \cdots \, \ul{a_k} \, l^{k}$.
Elements of $\ndf{A}[k][n]$ will be called
\textbf{codifferential forms} on $A$ of \textbf{line degree} $k$ and \textbf{cohomological degree} $n$.
We have the following operations on $\ndf{A}[][]$:

\begin{enumerate}
	\item A coderivation $\qdr \colon \ndf{A}[][] \rightharpoonup \ndf{A}[][]$
	      of bidegree $(-1,0)$, called the \textbf{de Rham differential}, given by
	      \begin{equation*}
		      \qdr \left( l^0 \, \ul{a_1} \, l^1 \, \dots \, \ul{a_k} \, l^{k} \right) =
		      \sum_{i=1}^{k} (-1)^{i - 1} \,
		      l^0 \, \ul{a_1} \, l^1 \, \dots \, \ul{a_{i-1}} \, l^{i-1} \, a_i \, l^{i} \, \ul{a_{i+1}} \, \dots \ul{a_k} \, l^{k},
	      \end{equation*}
	      which satisfies $\qdr^2 = 0$.
	\item A coderivation $\lie{\mu} \colon \ndf{A}[][] \rightharpoonup \ndf{A}[][]$ of bidegree $(0,1)$,
	      called the \textbf{Lie derivative} along $\mu$, which is determined uniquely by the requirement
	      that $\lie{\mu}$ extends $\mu$ on $\ndf{A}[0][] = \tens{A}$ and satisfies
	      $\left[ \qdr, \lie{\mu} \right] = \qdr \circ \lie{\mu} - \lie{\mu} \circ \qdr = 0$.
	      An explicit formula for $\lie{\mu}$ is given in \cref{eq:lie-mu-full-formula}.
	      Just like $\mu$, the coderivation $\lie{\mu}$ satisfies $\lie{\mu}^2 = 0$.
\end{enumerate}

The construction of $\ndf{A}[][]$ appears in \cite{Herscovich2023} and is dual to a construction introduced by Kontsevich and Soibelman
in \cite{Kontsevich:vh}. The assignment $A \mapsto \ndf{A}[][]$ is functorial, i.e., given an $\Ainf$-morphism $f \colon \mathcal{A} \rightarrow \mathcal{B}$,
there is an induced morphism $\indmap{f} \colon \ndf{A}[][] \rightarrow \ndf{B}[][]$ which commutes with $\qdr$ and the Lie derivative
(see \cref{sec:functoriality-ndf-and-ncdf}).
Endowing $\ndf{A}[][]$ with the differentials $\qdr, \lie{\mu}$, we obtain a bicomplex  $\left( \ndf{A}[*][*], \qdr, \lie{\mu} \right)$,
which is not very interesting.
The rows of the bicomplex  $\left( \ndf{A}[*][*] / R, \qdr, \lie{\mu} \right)$ are contractible by the
formal Poincar\'{e} \cref{lm:formal-poincare-ncdfr}, and, when $\mathcal{A}$ is unital, the columns are also contractible.
However, if we replace the tensor module $\tens{A \oplus \ul{A}}$ with the cyclic tensor module
$\tenscyc{A \oplus \ul{A}}$, which is a quotient of $\tens{A \oplus \ul{A}}$ in which we identify
elements related by the rotation operator
\begin{equation*}
	\t \left( x_1 \otimes \dots \otimes x_k \right) =
	(-1)^{\braid{\deg{x_k}}{\deg{x_1} + \dots + \deg{x_{k-1}}}} x_k \otimes x_1 \otimes \dots \otimes x_{k-1},
\end{equation*}
we obtain non-trivial complexes.

Define the \textbf{cyclic codifferential forms} on $A$ to be $\ncdf{A}[*][*] = \tenscyc{A \oplus \ul{A}}[(*,*)]$.
There is a general functorial construction in which a coderivation $\eta$ (resp.\ a morphism $f$) on the tensor coalgebra
induces an operator $\cycl{\eta}$ (resp.\ $\cycl{f}$) on cyclic codifferential forms called its \textbf{cyclization}, explained in
\cref{sec:cyc-tensor-coalgebra}.
Applying this construction to $\lie{\mu}$, one obtains a bidegree $(0,1)$ differential $\clie{\mu} = \cycl{\left( \lie{\mu} \right)}$
on $\ncdf{A}[][]$, called the \textbf{cyclic Lie derivative}.
The differential $\qdr$ descends to $\ncdf{A}[][]$, and we still have $\left[ \qdr, \clie{\mu} \right] = 0$,
so we obtain a bicomplex
\begin{equation*}
	\ncdf{\mathcal{A}}[][] = \left( \ncdf{A}[*][*], \qdr, \clie{\mu} \right),
\end{equation*}
called the bicomplex of cyclic codifferential forms (see \cref{fig:bicomplex-cyc-codiff-forms}).
The construction of $\ncdf{\mathcal{A}}[][]$ is also functorial, i.e., given an $\Ainf$-morphism $f \colon \mathcal{A} \rightarrow \mathcal{B}$,
there is an induced morphism $\cindmap{f} \colon \ncdf{A}[][] \rightarrow \ncdf{B}[][]$ which commutes with $\qdr$ and the
cyclic Lie derivative.
\begin{figure}[htb]
	\centering
	\begin{tikzcd}
		& \vdots & {\vdots} & {\vdots } & {\vdots } \\
		{\cdots} & 0 & \ncdf{A}[0][1] & \ncdf{A}[1][1] & \ncdf{A}[2][1] & {\cdots } \\
		{\cdots} & 0 & \ncdf{A}[0][0] & \ncdf{A}[1][0] & \ncdf{A}[2][0] & {\cdots } \\
		{\cdots} & 0 & \ncdf{A}[0][-1] & \ncdf{A}[1][-1] & \ncdf{A}[2][-1] & {\cdots } \\
		& {\vdots} & {\vdots } & {\vdots } & {\vdots } \\
		\arrow[from=2-2, to=2-1]
		\arrow[from=2-3, to=2-2]
		\arrow["{\qdr}"', from=2-4, to=2-3]
		\arrow["{\qdr}"', from=2-5, to=2-4]
		\arrow[from=2-6, to=2-5]
		\arrow[from=3-2, to=3-1]
		\arrow[from=3-3, to=3-2]
		\arrow["{\qdr}"', from=3-4, to=3-3]
		\arrow["{\qdr}"', from=3-5, to=3-4]
		\arrow[from=3-6, to=3-5]
		\arrow[from=4-2, to=4-1]
		\arrow[from=4-3, to=4-2]
		\arrow["{\qdr}"', from=4-4, to=4-3]
		\arrow["{\qdr}"', from=4-5, to=4-4]
		\arrow[from=4-6, to=4-5]
		\arrow[from=2-2, to=1-2]
		\arrow[from=3-2, to=2-2]
		\arrow[from=4-2, to=3-2]
		\arrow[from=5-2, to=4-2]
		\arrow["{\clie{\mu}}", from=2-3, to=1-3]
		\arrow["{\clie{\mu}}", from=3-3, to=2-3]
		\arrow["{\clie{\mu}}", from=4-3, to=3-3]
		\arrow["{\clie{\mu}}", from=5-3, to=4-3]
		\arrow["{\clie{\mu}}", from=2-4, to=1-4]
		\arrow["{\clie{\mu}}", from=3-4, to=2-4]
		\arrow["{\clie{\mu}}", from=4-4, to=3-4]
		\arrow["{\clie{\mu}}", from=5-4, to=4-4]
		\arrow["{\clie{\mu}}", from=2-5, to=1-5]
		\arrow["{\clie{\mu}}", from=3-5, to=2-5]
		\arrow["{\clie{\mu}}", from=4-5, to=3-5]
		\arrow["{\clie{\mu}}", from=5-5, to=4-5]
	\end{tikzcd}
	\caption{The Bicomplex $\ncdf{\mathcal{A}}[][]$ of Cyclic Codifferential Forms.}
	\label{fig:bicomplex-cyc-codiff-forms}
\end{figure}

\subsubsection{Total Inner Products} \label{sec:intro-total-inner-products}
Given a Banach $\Ainf$-algebra $\mathcal{A} = \left( A, \mu \right)$, we denote by
$\totcomp{\mathcal{A}}[2][]$ the total complex of the bicomplex
obtained from $\ncdf{\mathcal{A}}[][]$
by removing the first two columns.
Elements $x \in \totcomp{\mathcal{A}}[2][]$ of degree $d$ are written as
$x = \sum_{k = 2}^{\infty} \s_k x_k$, where $x_k \in \ncdf{A}[k][d+k]$ and $\s_k \colon \ncdf{A}[k][] \rightharpoonup \ncdf{A}[k][][k]$
is the suspension map. The total differential on $\totcomp{\mathcal{A}}[2][]$ is given by
$D(x) = \sum_{k = 2}^\infty \s_k ( \qdr(x_{k+1}) + (-1)^k \clie{\mu}x_k )$.
By adjoining a formal symbol $\ul{1}$ of degree $-1$ to $\totcomp{\mathcal{A}}[2][]$ and setting
\begin{equation*}
	D \left( \ul{1} \right) = \sum_{k = 2}^{\infty}
	\frac{(-1)^{\frac{k \left( k + 1 \right)}{2}}}{k} \s_k \left( {\ulz{\mu}}^{\otimes k} \right),
\end{equation*}
we obtain the \textbf{extended total complex} $\totcompe{\mathcal{A}}[2][]$.

An $n$-\textbf{dimensional total inner product Banach} $\Ainf$-\textbf{algebra over} $\mathcal{R}$
is a triple $\mathcal{A} = \left( A, \mu, \phi \right)$, where $\left( A, \mu \right)$ is a
Banach $\Ainf$-algebra over $\mathcal{R}$, and
$\phi \colon \totcompe{\mathcal{A}}[2][] \rightarrow \mathcal{R}[4-n]$
is a morphism of differential graded Banach $\mathcal{R}$-modules, called
an $n$-\textbf{dimensional total inner product on} $\mathcal{A}$.
A total inner product gives a notion of pairing on an $A_\infty$-algebra that relaxes strict symmetries
in previous definitions to symmetries up to an infinite family of coherent homotopies.
To illustrate this, we assume for simplicity that $d = 0$, $\mu_0 \left( 1 \right) = 0,$ and $\phi \left( \ul{1} \right) = 0$.
Denote by $\phi_k \colon \ncdf{A}[k][] \rightharpoonup R$ the induced map of $\phi$ on cyclic codifferential forms of
line degree $k$.

A \textbf{cyclic structure}~\cite{Fukaya2010,Solomon2016} arises from a total inner product for which $\phi_k = 0$ when $k \geq 3$,
and such that
\begin{equation} \label{eq:phi_2-strict}
	\phi_2 \left( \ul{a}, b_1, \dots, b_i, \ul{c}, d_1, \dots, d_j \right) = 0
\end{equation}
whenever $i \geq 1$ or $j \geq 1$.
Given~\eqref{eq:phi_2-strict}, the relation $\phi_2 \circ \qdr = 0$ is satisfied tautologically, while
the relation $\phi_2 \circ \clie{\mu} = 0$ is equivalent to the
relations
\begin{equation}\label{eq:cycsym}
	\phi_2 \left( \ul{\mu_k \left( a_1, \dots, a_k \right)}, \ul{a_{k+1}} \right) =
	(-1)^{\degb{a_{k+1}} \cdot \left( \degb{a_1} + \dots + \degb{a_k} \right)}
	\phi_2 \left( \ul{\mu_k \left( a_{k+1}, a_1, \dots, a_{k-1} \right)}, \ul{a_k} \right)
\end{equation}
for $k \geq 1$ (see \cref{appendix:cyclic-structures}).

A \textbf{strong homotopy inner product}~\cite{Cho2008,cho-homotopy-superpotential} or
\textbf{non-constant symplectic structure}~\cite{Kontsevich:vh} arises from a total inner product
for which $\phi_k = 0$ when $k \geq 3$. In this case, the strict symmetry~\eqref{eq:cycsym} is relaxed
to the condition  $\phi_2 \circ \clie{\mu} = 0$, which can be understood as~\eqref{eq:cycsym} up to
an infinite family of coherent homotopies. The condition $\phi_2 \circ \qdr = 0$ is no longer tautological,
but rather reads,
\begin{equation}\label{eq:closed}
	(-1)^{\left( \degb{y} + \degb{b} + \degb{z} + \degb{c} \right) \left( \degb{x} + \degb{a} \right)}
	\phi_2 \left(
	\ul{y}, b, \ul{z}, c, x, a
	\right)
	-
	\phi_2 \left(
	\ul{x}, a, y, b, \ul{z}, c
	\right)
	+
	\phi_2 \left(
	\ul{x}, a, \ul{y}, b, z, c
	\right) = 0,
\end{equation}
for $x,y,z \in A$ and $a,b,c \in \tens{A}$. This condition is called closedness in~\cite{Cho2008,cho-homotopy-superpotential,Kontsevich:vh}.

A general total inner product on $\mathcal{A}$ relaxes the strict closedness symmetry~\eqref{eq:closed} to the infinite family of coherent homotopy conditions
\begin{equation*}
	\phi_2 \circ \clie{\mu} = 0, \qquad
	(-1)^k \phi_k \circ \clie{\mu} + \phi_{k-1} \circ \qdr = 0 \quad (k \geq 3).
\end{equation*}
When the differential $d$ on $R$ is non-zero, one adds
$(-1)^n d \circ \phi_k$ to the right-hand side of the preceding relations. When $\mu_0 \left( 1 \right)$ and $\phi \left( \ul{1} \right)$ are non-zero, one obtains the extra identity
\begin{equation} \label{eq:extra-identity-phi}
	(-1)^{n} d \left( \phi \left( \ul{1} \right) \right) =
	\sum_{k=2}^{\infty} \frac{(-1)^{\frac{k \left( k + 1 \right)}{2}}}{k} \phi_k \left( {\ulz{\mu}}^{\otimes k} \right).
\end{equation}

When $\phi_k = 0$ for $k \geq 3$, we call $\phi$ a \textbf{homotopy inner product}. In this case, \eqref{eq:extra-identity-phi} simplifies to
\begin{equation} \label{eq:d-phi-ul-1-simp}
	(-1)^{n} d \left( \phi \left( \ul{1} \right) \right) = - \frac{1}{2}\phi_2 \left( \ulz{\mu}, \ulz{\mu} \right).
\end{equation}
\Cref{eq:d-phi-ul-1-simp} appears in Proposition~4.20 of \cite{Solomon2016},
with $\mathfrak{m}_{-1}$ replacing $(-1)^n \phi \left( \ul{1} \right)$ and an extra term coming from $J$-holomorphic spheres.
There, it plays a role in the reformulation of~\cite[Definition~3.1]{Fukaya2011}. It is needed in the definition of open
Gromov--Witten invariants~\cite{Fukaya2011,Solomon2016a}. General total inner products appear naturally
in the context of descendent open Gromov--Witten invariants~\cite{Solomon2026} and
are used in this work to define the superpotential, as explained in \cref{sec:intro-superpotential}.

When $\mathcal{A}$ is unital, instead of using the total complex $\totcomp{\mathcal{A}}[2][]$, we work with reduced versions, which are quotients of $\totcomp{\mathcal{A}}[2][]$ in which
elements involving the unit in a specific form are annihilated. We introduce two versions,
$\totcompred{\mathcal{A}}[2][]$ and $\totcompsred{\mathcal{A}}[2][]$, called
reduced and strongly reduced respectively, and their extended versions
$\totcompered{\mathcal{A}}[2][], \totcompesred{\mathcal{A}}[2][]$
(see \cref{sec:extended-reduced-total-complexes}). A total inner product
which descends to $\totcompered{\mathcal{A}}[2][]$ (resp.\ $\totcompesred{\mathcal{A}}[2][]$)
is called \textbf{unital} (resp.\ \textbf{strongly unital}), and a unital (resp.\ strongly unital)
total inner product $\Ainf$-algebra is a unital $\Ainf$-algebra equipped with a unital
(resp.\ strongly unital) total inner product.

\subsubsection{The Superpotential and its Properties} \label{sec:intro-superpotential}

Given an $n$-dimensional total inner product Banach $\Ainf$-algebra $\mathcal{A} = \left( A, \mu, \phi \right)$ over $\mathcal{R}$,
the \textbf{superpotential} associated to $\mathcal{A}$ is the function
\begin{equation*}
	\SP[] = \SP[][\mathcal{A}] \colon \tc{A} \rightarrow R^{3-n}
\end{equation*}
defined by
\begin{equation*}
	\SP[b] =
	\phi \left( \ul{1} \right) +
	\sum_{\substack{k = 2 \\ i_1, \dots, i_{k-1} = 0 \\ j_1, \dots, j_k = 0}}^{\infty}
	\frac{(-1)^{\frac{(k - 2) \cdot (k - 1)}{2}}}{1 + \sum_{r=1}^{k-1} i_r + \sum_{r = 1}^k j_r}
	\phi_k \left(
	\ul{ \mu_{i_1} \left( b^{i_1} \right)} \, b^{j_1} \, \cdots \, \ul{ \mu_{i_{k-1}} \left( b^{i_{k-1}} \right) } \,
	b^{j_{k-1}} \, \ul{ \vphantom{ \mu_{i_1} \left( b^{i_1} \right) } b} \, b^{j_k}
	\right).
\end{equation*}
When $\phi_k = 0$ for $k \geq 3$, i.e., $\phi$ corresponds to a homotopy inner product, the superpotential takes the form
\begin{equation*}
	\SP[b] = \phi \left( \ul{1} \right) +
	\sum_{i,j,k=0}^{\infty} \frac{1}{i+j+1+k} \phi_2 \left(
	\ul{ \mu_i \left( b^i \right) } \raisebox{-2.5pt}{\,,\,} b^j \raisebox{-2.5pt}{\,,\,}
	\ul{ \vphantom{\mu_i \left( b^i \right)} b} \raisebox{-2.5pt}{\,,\,} b^k
	\right).
\end{equation*}
When $\phi \left( \ul{1} \right) = 0$ and $\mathcal{A}$ is uncurved, this superpotential was introduced
and studied by Cho and Lee \cite{cho-homotopy-superpotential}.
Furthermore, when $\phi$ corresponds to a cyclic structure, $\SP[b]$ reduces to the familiar form
\begin{equation*}
	\SP[b] =
	\sum_{i=0}^{\infty} \frac{1}{i+1} \phi_2 \left(
	\ul{ \mu_i \left( b^i \right) } \raisebox{-2.5pt}{\,,\,}  \ul{ \vphantom{ \mu_i \left( b^i \right) } b}
	\right)
\end{equation*}
appearing in \cite{Fukaya2011,Solomon2016a}.

To state the properties of superpotential, we need to discuss several additional preliminaries.
Given a morphism $f \colon \mathcal{A} \rightarrow \mathcal{B}$ of Banach $\Ainf$-algebras, there is
an induced chain map
$\cindmape{f} \colon \totcompe{\mathcal{A}}[2][] \rightarrow \totcompe{\mathcal{B}}[2][]$,
which is functorial on the level of cohomology. When $\mathcal{A}, \mathcal{B}$ and $f$ are unital,
the map descends to the reduced versions. See \cref{sec:extended-tot-comp-geq2-braidop-2}.
This leads to a notion
of a morphism between (unital) total inner product Banach $\Ainf$-algebras, which is a morphism
of the underlying $\Ainf$-algebras respecting the total inner products (see \cref{dfn:morphism-total-inner-products}).
A pseudoisotopy $\mathfrak{A}$ between two (unital) total inner product $\Ainf$-algebras $\mathcal{A}_0$ and $\mathcal{A}_1$ is then
a (unital) pseudoisotopy $\mathfrak{A}$ between the underlying $\Ainf$-algebras equipped with a (unital) total inner product such that the pseudoisotopy morphisms $\evalmf^i \colon \mathfrak{A} \rightarrow \mathcal{A}_i$ become morphisms of total inner product Banach
$\Ainf$-algebras.
\begin{thmx}[Properties of the Superpotential Function] \hfill \ \label{thm:superpotential-properties}
	\begin{enumerate}[label=(\arabic*), ref=(\arabic*)]
		\item
		      Let $\mathcal{A}$ and $\mathcal{B}$ be two $n$-dimensional total inner product
		      Banach $\Ainf$-algebras over $\mathcal{R} = \left( R, d_R \right)$ and $\mathcal{S} = \left( S, d_S \right)$ respectively.
		      Given a morphism $f \colon \mathcal{A} \rightarrow \mathcal{B}$ of total inner product Banach $\Ainf$-algebras
		      and $b \in \tc{A}$, we have
		      \begin{equation*}
			      \base{f} \left( \SP[b][\mathcal{A}] \right) =
			      \SP[\mcfunc{f} \left( b \right)][\mathcal{B}] + d_S \left( x \right)
		      \end{equation*}
		      for some $x \in S^{2-n}$.
		      \label{item:superpotential-1}
		\item
		      Let $\mathcal{A}$ be an $n$-dimensional non-unital (resp.\ unital)
		      total inner product Banach $\Ainf$-algebra over $\mathcal{R}$ and let $b \in \tc{A}$ be a strong (resp.\ weak)
		      bounding cochain. Then $\SP[b][\mathcal{A}]$ is a cocycle of degree $3 - n$.
		      \label{item:superpotential-2}
		\item
		      Let $\mathcal{A}_0$ and $\mathcal{A}_1$ be two non-unital (resp.\ unital) total inner product
		      Banach $\Ainf$-algebras over $\mathcal{S}$, and let $\mathfrak{A}$ be a non-unital (resp.\ unital)
		      pseudoisotopy of total inner product Banach $\Ainf$-algebras between $\mathcal{A}_0$ and $\mathcal{A}_1$.
		      Let $b_0 \in A_0$ and $b_1 \in A_1$ be two $\mathfrak{A}$-gauge-equivalent strong (resp.\ weak) bounding cochains.
		      Then
		      \begin{equation*}
			      \eqcl{ \SP[b_0][\mathcal{A}_0] } = \eqcl{ \SP[b_1][\mathcal{A}_1] }
			      \in \cohom{\mathcal{S}}[3-n].
		      \end{equation*}
		      \label{item:superpotential-3}
	\end{enumerate}
\end{thmx}
\cref{thm:superpotential-properties} is deduced from the properties of the cyclic Chern--Simons form,
which is defined in \cref{sec:Chern-Simons-intro}. \Cref{thm:Gb-geq-2-properties} gives
an analog of~\cref{thm:superpotential-properties} for cyclic Chern--Simons forms and the associated
cohomology classes.

In applications to open Gromov--Witten theory, one considers bounding cochains that depend
on a formal variable, which is adjoined by scalar extension~\cite{Solomon2016a}.
The superpotential of such a bounding cochain, which is a generating function for open
Gromov--Witten invariants, satisfies a system of PDE known as the open WDVV
equations~\cite{Solomon2024}. The proof of the open WDVV equations involves a computation
of the derivative of the superpotential with respect to the formal variable.
The formalism of the present paper is designed for the definition of descendent open
Gromov--Witten invariants~\cite{Solomon2026}. The derivative of the descendent
superpotential plays a role in proving the open topological recursion relations.
See~\cite{Pandharipande2024} for the case of open descendent integrals.

With the preceding paragraph in mind, we give a formula for the derivative of the superpotential of a total inner product Banach $\Ainf$-algebra.
Since all the ingredients involved in defining the extended total complex commute with scalar extension, one has
a natural notion of scalar extension for total inner products Banach $\Ainf$-algebras
(see \cref{sec:scalar-extension-total-inner-products}).
\begin{thmx}[Formal Derivative] \label{thm:superpotential-derivative}
	Let $R$ be a graded-commutative Banach $\mathbbm{k}$-algebra and let
	$\mathcal{A}$ be an $n$-dimensional total inner product Banach $\Ainf$-algebra
	over $R$. Assume that $\mathcal{B} = \left( B, \mu, \phi \right)$ is obtained from $\mathcal{A}$ by
	scalar extension along $R \rightarrow \pows{R}[t]$, where $t$ is an even formal variable
	with $\nnorm[t] < 1$.
	Then, given $b = b(t) \in \tc{B}$, we have
	\begin{equation} \label{eq:derivative-sp-braid-1}
		\partial_t \left( \SP[b][\mathcal{B}] \right) =
		\sum_{k=2}^{\infty} (-1)^{\frac{(k-2)(k-1)}{2}}
		\phi_k \left(
		\underbrace{\ul{ \corest{\mu} \left( \Exp{b} \right) } \, \Exp{b} \, \dots \, \ul{ \corest{\mu} \left( \Exp{b} \right) } \, \Exp{b}}_{k-1} \,
		\ul{ \partial_t \left( b \right) } \, \Exp{b}
		\right).
	\end{equation}
	In particular, if $b$ is a strong bounding cochain, then $\partial_t \left( \SP[b][\mathcal{B}] \right) = 0$.
\end{thmx}

\subsubsection{The Cyclic Chern--Simons Form}\label{sec:Chern-Simons-intro}
Given a Banach $\Ainf$-algebra $\mathcal{A}$ and $b \in \tc{A}$,
the superpotential $\SP[b][\mathcal{A}]$ is obtained by applying the total inner product $\phi$
to the cyclic Chern--Simons form $\G{b}[\geq 2][\mathcal{A}]$
given by
\begin{equation*}
	\G{b}[\geq 2][\mathcal{A}] = \ul{1} +
	\sum_{\substack{k = 2 \\ i_1, \dots, i_{k-1} = 0 \\ j_1, \dots, j_k = 0}}^{\infty}
	\frac{(-1)^{\frac{(k - 2) \cdot (k - 1)}{2}}}{1 + \sum_{r=1}^{k-1} i_r + \sum_{r = 1}^k j_r}
	\s_k \left(
	\ul{ \mu_{i_1} \left( b^{i_1} \right)} \, b^{j_1} \, \cdots \, \ul{ \mu_{i_{k-1}} \left( b^{i_{k-1}} \right) } \,
	b^{j_{k-1}} \, \ul{ \vphantom{ \mu_{i_1} \left( b^{i_1} \right) } b} \, b^{j_k}
	\right).
\end{equation*}
\cref{sec:recovering-chern-simons} shows how to recover the classical Chern--Simons form and action
from $\G{b}[\geq 2][\mathcal{A}]$ and $\SP[b][\mathcal{A}]$ respectively.

Let $\mathcal{A}_0$ and $\mathcal{A}_1$ be two non-unital Banach $\Ainf$-algebras over
the same differential graded-commutative Banach $\mathbbm{k}$-algebra $\mathcal{S}$.
A pseudoisotopy $\mathfrak{A}$ over $\mathfrak{R}$ between
$\mathcal{A}_0$ and $\mathcal{A}_1$ is called
$\cohom{}[] \totcompe{}[2]$-\textbf{strong}
if the underlying homotopic morphisms $\evalm^i = \base{\evalmf}^i \colon \mathfrak{R} \rightarrow \mathcal{S}$
of differential graded algebras
induce an isomorphism
$\cohom{\mathfrak{R}}[] \overset{\sim}{\rightarrow} \cohom{\mathcal{S}}[]$
on cohomology,
and the $\Ainf$-morphisms $\evalmf^i \colon \mathfrak{A} \rightarrow \mathcal{A}_i$
induce isomorphisms
$\cohom{ {\totcompe{{\mathfrak{A}}}[2][]}}[] \overset{\sim}{\rightarrow}
	\cohom{ {\totcompe{\mathcal{A}_i}[2][]}}[]$,
for $i = 0, 1$.
A $\cohom{}[] \totcompe{}[2]$-strong
pseudoisotopy yields a canonical isomorphism
\begin{equation*}
	\mathfrak{a} \colon \cohom{{\totcompe{\mathcal{A}_0 / \mathcal{S}}[2][]}}[] \overset{\sim}{\rightarrow} \cohom{{\totcompe{\mathcal{A}_1 / \mathcal{S}}[2][]}}[].
\end{equation*}
When $\mathcal{A}_0$ and $\mathcal{A}_1$ are unital, one similarly defines a
$\cohom{}[] \totcompered{}[2]$-\textbf{strong} pseudoisotopy by requiring all maps to be unital
and replacing $\totcompe{}[2]$ with $\totcompered{}[2]$.

The properties of the cyclic Chern--Simons form $\G{b}[\geq 2]$ are given in the following theorem,
which can be seen as a categorification of~\cref{thm:superpotential-properties}:
\begin{thmx}[Properties of the Cyclic Chern--Simons Form] \hfill \ \label{thm:Gb-geq-2-properties}
	\begin{enumerate}[label=(\arabic*), ref=(\arabic*)]
		\item
		      Let $\mathcal{A}$ and $\mathcal{B}$ be two Banach $\Ainf$-algebras over $\mathcal{R}$
		      and $\mathcal{S}$ respectively.
		      Given a morphism $f \colon \mathcal{A} \rightarrow \mathcal{B}$ of Banach $\Ainf$-algebras
		      and $b \in \tc{A}$, we have
		      \begin{equation*}
			      \cindmape{f} \left( \G{b}[\geq 2][\mathcal{A}] \right) =
			      \G{\mcfunc{f} \left( b \right)}[\geq 2][\mathcal{B}] +
			      D_{\mathcal{B}} \left( R \left( b; f \right) \right)
		      \end{equation*}
		      for some $R \left( b; f \right) \in \totcomp{B}[2][-2]$.
		      \label{item:Gb-geq-2-properties-1}
		\item
		      Let $\mathcal{A}$ be a non-unital (resp.\ unital) Banach $\Ainf$-algebra and
		      let $b \in \tc{A}$ be a strong (resp.\ weak) bounding cochain. Then $\G{b}[\geq 2][\mathcal{A}]$ is closed
		      in $\totcompe{\mathcal{A}}[2][]$ (resp.\ $\totcompered{\mathcal{A}}[2][]$).
		      \label{item:Gb-geq-2-properties-2}
		\item
		      Let $\mathcal{A}_0$ and $\mathcal{A}_1$ be two non-unital (resp.\ unital) Banach $\Ainf$-algebras
		      over $\mathcal{S}$ and let $\mathfrak{A}$ be a $\cohom{}[] \totcompe{}[2]$-strong
		      (resp.\ $\cohom{}[] \totcompered{}[2]$-strong) pseudoisotopy between $\mathcal{A}_0$
		      and $\mathcal{A}_1$.
		      Let $b_0 \in \tc{A_0}$ and $b_1 \in \tc{A_1}$ be
		      $\mathfrak{A}$-gauge-equivalent strong (resp.\ weak) bounding cochains.
		      Then
		      \begin{equation*}
			      \qquad\qquad
			      \mathfrak{a}(\eqcl{ \G{b_0}[\geq 2][\mathcal{A}_0] }) = \eqcl{ \G{b_1}[\geq 2][\mathcal{A}_1] }
			      \textrm{ in } \cohom{{\totcompe{\mathcal{A}_1}[2][]}}[-1]
			      \,\, \left(
			      \textrm{resp.\ in } \cohom{{\totcompered{\mathcal{A}_1}[2][]}}[-1]
			      \right).
		      \end{equation*}
		      \label{item:Gb-geq-2-properties-3}
	\end{enumerate}
\end{thmx}

\subsubsection{The \texorpdfstring{$\infty$}{Infinity}-Modulus Function and its Properties}
Let $\mathcal{A} = \left( A, \mu \right)$ be a Banach $\Ainf$-algebra over $\mathcal{R} = \left( R, d \right)$.
The zeroth column $\ncdf{\mathcal{A}}[0][] = \left( \ncdf{A}[0][], \clie{\mu} \right)$ of the bicomplex
of cyclic codifferential forms is called the \textbf{extended Connes complex}. As a graded Banach $R$-module,
we have $\ncdf{A}[0][] = R \oplus \ncdfr{A}[0][]$, where $\ncdfr{\mathcal{A}}[0][]$ is
\textbf{Connes' cyclic complex} computing the cyclic homology of $\mathcal{A}$, and the action of $\clie{\mu}$
on the $R$ summand is determined by $\clie{\mu} \left( 1 \right) = \mu_0 \left( 1 \right)$.

A morphism $\theta \colon \ncdf{\mathcal{A}}[0][] \rightarrow \mathcal{R}[1-n]$
of differential graded Banach $\mathcal{R}$-modules is called an
$n$-\textbf{dimensional} $\infty$-\textbf{trace} on $\mathcal{A}$, and $\mathcal{A}$ together
with $\theta$ is called an $n$-\textbf{dimensional} $\infty$-\textbf{trace Banach}
$\Ainf$-\textbf{algebra}.
An $\infty$-trace is completely determined by the associated sequence
$\left( \theta_k \colon A^{\otimes k} \rightharpoonup R \right)_{k \geq 0}$ of cyclically invariant operations of
degree $1 - n$, in terms of which the condition
$d_{\mathcal{R}[1-n]} \circ \theta = \theta \circ \clie{\mu}$ can be written explicitly (see \cref{eq:infty-trace-k-geq-1-rel,eq:infty-trace-0-rel}).

Given an $n$-dimensional $\infty$-trace Banach $\Ainf$-algebra $\mathcal{A} = \left( A, \mu, \theta \right)$,
and a topologically nilpotent element $b \in \tc{A}$, the $\infty$-\textbf{modulus}
function $\TM{} = \TM{}[\mathcal{A}] \colon \tc{A} \rightarrow R^{1-n}$ is defined by
\begin{equation} \label{eq:infty-modulus-intro}
	\TM{b} = \theta_0 \left( 1 \right) + \sum_{k=1}^{\infty} \frac{1}{k} \theta_k \left( b^k \right).
\end{equation}
When $\mathcal{A}$ is unital, instead of using the complex $\ncdf{\mathcal{A}}[0][]$, one
can work with a reduced version $\ncdfred{\mathcal{A}}[0][]$, which is the quotient of $\ncdf{\mathcal{A}}[0][]$
by the module generated by elementary tensors which contain the unit $e$. An $\infty$-trace which
descends to $\ncdfred{\mathcal{A}}[0][]$, i.e., satisfies $\theta_k \left( a_1, \dots, a_k \right) = 0$
whenever $a_i = e$ for some $1 \leq i \leq k$, is called \textbf{unital}, and a
\textbf{unital} $\infty$\textbf{-trace Banach} $\Ainf$\textbf{-algebra} is a unital $\Ainf$-algebra equipped with a unital $\infty$-trace.

\Cref{dfn:morphism-inf-trace} gives a notion of a morphism between (unital) $\infty$-trace Banach $\Ainf$-algebras.
A pseudoisotopy $\mathfrak{A}$ between two (unital) $\infty$-trace Banach $\Ainf$-algebras $\mathcal{A}_0$ and $\mathcal{A}_1$ is
a (unital) pseudoisotopy $\mathfrak{A}$ between the underlying $\Ainf$-algebras equipped with a (unital) $\infty$-trace, such
that the pseudoisotopy morphisms $\evalmf^i \colon \mathfrak{A} \rightarrow \mathcal{A}_i$ are morphisms of $\infty$-trace Banach
$\Ainf$-algebras.
Then we have the following:
\begin{thmx}[Properties of the $\infty$-modulus Function] \hfill \ \label{thm:infinity-modulus-properties}
	\begin{enumerate}[label=(\arabic*), ref=(\arabic*)]
		\item
		      Let $\mathcal{A}$ and $\mathcal{B}$ be $n$-dimensional
		      $\infty$-trace Banach $\Ainf$-algebras,
		      possibly over different bases. Given a morphism $f \colon \mathcal{A} \rightarrow \mathcal{B}$ of $\infty$-trace
		      Banach $\Ainf$-algebras and $b \in \tc{A}$, we have
		      \begin{equation*}
			      \base{f} \left( \TM{b}[\mathcal{A}] \right) =  \TM{\mcfunc{f} \left( b \right)}[\mathcal{B}].
		      \end{equation*}
		      \label{item:infinity-modulus-1}
		\item
		      Let $\mathcal{A}$ be a non-unital (resp.\ unital) $n$-dimensional $\infty$-trace
		      Banach $\Ainf$-algebra over $\mathcal{R}$ and let $b \in \tc{A}$ be a strong (resp.\ weak)
		      bounding cochain. Then $\TM{b}[\mathcal{A}]$ is a cocycle of degree $1 - n$.
		      \label{item:infinity-modulus-2}
		\item 
		      Let $\mathcal{A}_0$ and $\mathcal{A}_1$ be non-unital (resp.\ unital) $n$-dimensional
		      $\infty$-trace Banach $\Ainf$-algebras
		      over $\mathcal{S}$ and let $\mathfrak{A}$ be a non-unital (resp.\ unital) pseudoisotopy
		      of $\infty$-trace Banach $\Ainf$-algebras between $\mathcal{A}_0$ and $\mathcal{A}_1$.
		      Let $b_0 \in \tc{A_0}$ and $b_1 \in \tc{A_1}$ be
		      $\mathfrak{A}$-gauge-equivalent strong (resp.\ weak) bounding cochains.
		      Then
		      \begin{equation*}
			      \eqcl{ \TM{b_0}[\mathcal{A}_0] } = \eqcl{ \TM{b_1}[\mathcal{A}_1] }
			      \in \cohom{\mathcal{S}}[1-n].
		      \end{equation*}
		      \label{item:infinity-modulus-3}
	\end{enumerate}
\end{thmx}

\Cref{thm:infinity-modulus-properties} is deduced from the properties of
the \textbf{cyclic exponential}
\begin{equation*}
	\G{b}[0] = 1 + \sum_{k=1}^{\infty} \frac{1}{k} b^k
\end{equation*}
proven in \cref{subsec:canonical-elements-ncdf-0}. \Cref{lm:invariance-Gb0-pseudoisotopy}
gives an analog for the cyclic exponential of \cref{thm:Gb-geq-2-properties}.
The cyclic exponential, together with a variant adapted to work for a different reduced complex,
has been used in~\cite{Sela2024} in the context of matrix factorizations.

\subsubsection{\texorpdfstring{Pre-total Inner Products, Pre-$\infty$-traces and Pre-homotopy Inner Products}{Pre-total Inner Products, Pre-infinity Traces and Pre-homotopy Inner Products}}

Both total inner products and $\infty$-traces were defined in terms of \textit{extended} complexes to allow
total inner products (resp.\ $\infty$-traces) to carry along a constant term $\phi \left( \ul{1} \right) \in R^{3-n}$
(resp.\ $\theta_0 \left( 1 \right) \in R^{1-n}$), independent of the input $b$, to the superpotential (resp.\ $\infty$-modulus).
It turns out that if we ignore the constant terms, the notions are equivalent
when $\mathcal{A}$ is unital. More precisely,
let us call a morphism $\phi \colon \totcomp{\mathcal{A}}[2][] \rightarrow \mathcal{R}[4-n]$
of differential graded Banach $\mathcal{R}$-modules an $n$-\textbf{dimensional pre-total inner product} on $\mathcal{A}$.
Similarly, a morphism $\theta \colon \ncdfr{\mathcal{A}}[0][] \rightarrow \mathcal{R}[1-n]$
of differential graded Banach $\mathcal{R}$-modules is called an $n$-\textbf{dimensional pre}-$\infty$ \textbf{trace} on $\mathcal{A}$.

Consider the projection
$p_2 \colon \totcomp{\mathcal{A}}[2][] \rightarrow
	( {{\ncdf{\mathcal{A}}[2][]}/ \Im \left( \qdr \right)} )[ 2 ]$.
As the rows of $\totcomp{\mathcal{A}}[2][]$ are contractible,
$p_2$ is a homotopy equivalence with an explicit homotopy inverse $i_2$
(see \cref{lm:p_2-homotopy-equivalence}).
We call a morphism $\phi_{2} \colon \ncdf{\mathcal{A}}[2][] / \Im \left( \qdr \right) \rightarrow \mathcal{R}[2-n]$
an $n$-\textbf{dimensional pre}-\textbf{homotopy inner product}, and using $p_2$ and $i_2$, we have a way of converting pre-total inner products to pre-homotopy inner products and vice versa.

We also show that when $\mathcal{A}$ is unital, the complexes $\ncdf{\mathcal{A}}[2][] / \Im \left( \qdr \right)$ and ${\ncdfr{\mathcal{A}}[0][]}[1]$
are homotopy equivalent by constructing a chain map $\psi \colon {\ncdf{\mathcal{A}}[2][]} / \Im \left( \qdr \right) \rightarrow
	{\ncdfr{\mathcal{A}}[0][]}[1]$, together with an explicit homotopy inverse $\psi'$ (see \cref{lm:ncdf-0-ncdf-2-mod-q3-equiv}).
The maps $\psi$ and $\psi'$ give a correspondence between pre-homotopy inner products and pre-$\infty$ traces, which is compatible with the unitality conditions for both notions.

\subsubsection{Derivative of the Superpotential and the \texorpdfstring{$\infty$}{Infinity}-Modulus}
Recall from \cref{sec:intro-total-inner-products} that a homotopy inner product is
a total inner product $\phi$ whose components $\phi_k$ vanish for $k \geq 3$.
This is equivalent to the data of a pre-homotopy inner product $\phi_2$, together with the extra element $\phi \left( \ul{1} \right)$, related via the identity
\eqref{eq:d-phi-ul-1-simp}.
Given a pre-$\infty$-trace $\theta$, we define the pre-$\infty$-modulus function $\TM{}[>0]$
as in \cref{eq:infty-modulus-intro}, without the constant term. We show the following:
\begin{thmx} \label{thm:derivative-sp-multiple-modulus}
	Let $R$ be a graded-commutative Banach $\mathbbm{k}$-algebra and let
	$\mathcal{A}$ be a unital Banach $\Ainf$-algebra over $R$ equipped with a strongly unital homotopy inner product.
	Assume that $\mathcal{B} = \left( B, \mu, e, \phi \right)$ is obtained from $\mathcal{A}$ by
	scalar extension along $R \rightarrow \pows{R}[t]$, where $t$ is an even formal variable
	with $\nnorm[t] < 1$.

	Let $\theta = \phi_2 \circ \psi'$ be the pre-$\infty$-trace corresponding to the pre-homotopy inner product $\phi_2$,
	and let $b = b(t) \in \tc{B}$. Then
	\begin{enumerate}
		\item The pre-$\infty$-modulus function $\TM{}[> 0]$ is expressed in terms of $\phi_2$ by
		      \begin{equation*}
			      \TM{b}[> 0] = \sum_{j,k=0}^{\infty} \frac{1}{j+1+k} \phi_2 \left( \ul{e}, b^j, \ul{b}, b^k \right).
		      \end{equation*}
		\item We have
		      \begin{equation*}
			      \partial_t \, \TM{b}[> 0] = \sum_{j,k=0}^{\infty} \phi_2 \left( \ul{e}, b^j, \ul{\partial_t \left( b \right)}, b^k \right).
		      \end{equation*}
		\item \label{it:partial-t-sp-c} When $b \in \mc{\mathcal{B}}[c]$ is a weak bounding cochain, we have
		      \begin{equation*}
			      \partial_t \, \SP[b] = c \cdot \partial_t \, \TM{b}[>0].
		      \end{equation*}
	\end{enumerate}
\end{thmx}

In the case that $\phi$ corresponds to a cyclic structure,
\cref{thm:derivative-sp-multiple-modulus} \eqref{it:partial-t-sp-c} plays an important role in the proof of
the open WDVV equation in~\cite{Solomon2024}. We expect a generalization of this result for
total inner-products to play an analogous role in the proof of the topological recursion relations for
open Gromov--Witten descendent invariants.

\subsection{Organization}
The remainder of this work is organized as follows. Each section begins with a more detailed overview of its contents.

\Cref{sec:prelim} establishes the notation, algebraic framework, and sign conventions
for working with graded objects, possibly equipped with (pre)-differentials.

\Cref{sec:non-archimedean-graded-setting} defines non-Archimedean graded seminormed and Banach objects,
their categories, and properties, providing the setting for the rest of the work.

\Cref{sec:formal-tensor-coalgebra} develops the completed tensor coalgebra, describing its grouplike elements,
morphisms, coderivations, and base change.

\Cref{sec:a-inf-algebras} introduces curved Banach $\Ainf$-algebras over differential
graded-commutative Banach ground algebras, allowing morphisms with change of connection elements.
It also establishes working definitions of pseudoisotopy and gauge equivalence for our framework.

\Cref{sec:cyc-tensor-coalgebra} constructs cyclic versions of coderivations and coalgebra morphisms
in the presence of curvature and change of connection elements. It establishes
their functoriality and compatibility with base change, and extends the constructions from the
reduced to the full tensor module and its cyclic quotient.

\Cref{sec:noncomm-diff-calc} constructs noncommutative codifferential forms and their cyclic counterparts.
It develops a Cartan calculus for the ordinary forms and deduces the corresponding calculus for the cyclic ones.
It shows that both versions are functorial, compatible with base change, and satisfy a formal Poincar\'e lemma.

\Cref{sec:cyclic-homology-models} develops several different models for the cyclic homology of curved Banach
$\Ainf$-algebras, some of which are built from cyclic codifferential forms.
Under appropriate hypotheses it establishes explicit homotopy equivalences between the different models.
It also provides a contraction-based chain-level map on Connes' cyclic complex that induces the periodicity operator on
homology for curved Banach $\Ainf$-algebras; finally it defines extended and reduced cyclic complexes used later on.

\Cref{sec:generalized-trace} uses the cyclic exponential to associate an extended cyclic homology class
to a bounding cochain, and shows it is natural and invariant under gauge equivalence in an appropriate sense.
It then introduces $\infty$-traces and pairs them with the cyclic exponential to define the
$\infty$-modulus function. Finally, it shows that the $\infty$-modulus of a
bounding cochain defines a natural gauge-invariant cohomology class.

\Cref{sec:generalized-superpotential} introduces total inner products and the cyclic Chern--Simons form, which combine to give the superpotential. It establishes the naturality, closedness, gauge invariance, and derivative properties of these
constructions in both the unital and non-unital settings, and describes the relations
satisfied by total inner products in terms of their components.

\Cref{sec:homotopy-inner-products} discusses (pre)-homotopy inner products
and gives explicit translations between pre-homotopy inner products and pre-$\infty$-traces.
It shows that the derivative of the superpotential on weak bounding cochains is a
multiple of derivative of the pre-$\infty$-modulus under appropriate unitality conditions.
Finally, it gives a relation between the pre-superpotential and the periodicity operator
on cyclic homology.

\Crefrange{appendix:cat-theory-background}{appendix:bicomplexes} provide foundational background on category theory, homological algebra, and non-Archimedean algebraic structures,
as well as results on Banach bicomplexes used in the comparison arguments.
\Cref{appendix:bar-complex-not-necessarily-contractible} discusses the role of completion, \cref{appendix:sign-conversions} compares sign conventions, and \cref{appendix:parity-forms-equiv} explains parity conversions.
\Cref{appendix:cyclic-structures} identifies the classical notion of a cyclic structure with a strict chain map on cyclic codifferential forms.
\Cref{appendix:sign-conversions-jake} clarifies the relation between the constructions developed here and the corresponding notions appearing in earlier works on
open Gromov--Witten theory.

\subsection{Acknowledgements}
The authors would like to thank Or Kedar, Paul Seidel and Sara Tukachinsky for helpful discussions.
The authors were partially supported by ERC starting grant 337560 and ISF grant 569/18.
P.~G.\ was partially supported by the Hoffman program for the doctoral students at the Hebrew University of Jerusalem.
J.~S.\ was partially supported by ISF grant 1127/22 and the Miriam and Julius Vinik Chair in Mathematics.
The authors would like to thank the Institute for Advanced Study for its hospitality in the initial stages of the
research that led to this paper, funded by the Erik Ellentuck Fellowship and the IAS Fund for Math.

\subsection{Conventions and Notation} \label{sec:conv-and-not}
For the convenience of the reader, we collect here some of the conventions and notation used repeatedly in this work.

\paragraph{\textit{Algebraic Setting}}

We denote by $\mathbbm{k}$ a fixed ungraded commutative ground ring.
Most constructions and definitions are done relative to a fixed graded-commutative
$\mathbbm{k}$-algebra, often denoted by $R$, working with left $R$-modules.
Unless otherwise stated, all tensor products are taken over $R$ and not over $\mathbbm{k}$.
All (co)algebras and (co)modules, as well as their morphisms, are (co)unital.
However, we do not assume that morphisms of coaugmented coalgebras are compatible with the
coaugmentations.

\paragraph{\textit{Objects and Maps}}
We reserve $1$ for the unit of $\mathbbm{k}$ or $R$, while denoting
the identity map by $\idd$ or $\id$.
A calligraphic letter $\mathcal{M}$ denotes an object equipped with extra data,
reserving the plain letter $M$ for the underlying object,
and a morphism $f \colon \mathcal{M} \rightarrow \mathcal{N}$ is
assumed to preserve all the extra data.
For example, we write $\mathcal{M} = \left( M, d \right)$
when $d$ is a differential on $M$ and $f \colon \mathcal{M} \rightarrow \mathcal{M}$
when $f \colon M \rightarrow M$ is a morphism which commutes with $d$.

\paragraph{\textit{Grading and Signs}}
Objects $M$ are graded with an upper index by a fixed abelian group $\GG$ equipped with a parity form
$\braidop$. Elements $m \in M$ are always homogeneous of degree $\degb{m} \in \GG$.
Exchanging homogeneous data of degrees $a,b \in \GG$ contributes the Koszul sign $(-1)^{\braid{a}{b}}$.
We reserve the notation $f \colon M \rightarrow N$ for actual (degree zero) morphisms,
and use the notation $f \colon M \rightharpoonup N$ to denote graded maps of arbitrary degree.
We use the shift convention $M[h]^g = M^{g+h}$ and for the purpose of degree tracking,
write elements in $M[h]$ as $\s_h m$ for $m \in M$.

\paragraph{\textit{Banach Setting}}
From \cref{sec:cyc-tensor-coalgebra} onward, we work exclusively in the non-Archimedean graded
Banach setting. Objects and morphisms are taken from, and categorical constructions are done in,
one of the categories of graded Banach objects described in detail in \Cref{sec:non-archimedean-graded-setting}.
In this context, morphisms $f \colon M \rightarrow N$ are contractive maps of degree zero,
and we use the notation $f \colon M \rightharpoonup N$ to denote graded maps of arbitrary
degree with uniformly bounded components.
Differentials are assumed to be bounded but not necessarily contractive.
To streamline notation, we suppress completion marks from direct sums and tensor products
(e.g., writing $\otimes$ instead of $\cotimes$), with the implicit understanding
that all constructions are completed.

\paragraph{$\Ainf$-\textit{algebras}}
Given a map $g \colon C \rightharpoonup \tens{W}$ into a tensor coalgebra,
its corestriction to $W$ is denoted by $\corest{g} \colon C \rightharpoonup W$.
When $C = \tens{V}$ is itself a tensor coalgebra, the restrictions of $\corest{g}$ to
$V^{\otimes k}$ are denoted by $g_k \colon V^{\otimes k} \rightharpoonup W$.
Our $\Ainf$-algebras $\mathcal{A} = (A, \mu)$ are cohomological and shifted unless
otherwise stated, i.e., defined as a degree one coderivation on $\tens{A}$
and not on $\tens{A[1]}$.
We allow curvature and morphisms $f \colon \mathcal{A} \rightarrow \mathcal{B}$ between
$\Ainf$-algebras may have a non-zero change-of-connection term $f_0(1) \in B$.
Although our grading is cohomological, we retain standard homological terminology
(e.g., chain complexes, chains, homology) when discussing classical constructions
such as Hochschild and cyclic homology and their variants.

\section{Preliminaries} \label{sec:prelim}

In this work, we work primarily with graded Banach modules, algebras and coalgebras.
As these are graded objects endowed with additional structure (a non-Archimedean norm, or more precisely, a family of non-Archimedean norms),
we begin by establishing the necessary algebraic foundations and setting up our notation for graded structures.

Let $\mathbbm{k}$ be a commutative ground ring which will be fixed for the duration of this section.
We start in \cref{subsec:grading-data} by discussing the grading datum that governs our sign conventions.
We then recall the standard properties of graded $\mathbbm{k}$-modules and algebras in
\cref{subsec:graded-k-modules,subsec:graded-k-algebras}.
In \cref{sub:graded-R-modules,sec:scalar-extension-restriction-graded-modules}, we generalize these notions to the relative setting,
discussing graded modules over graded algebras and the functorial properties of scalar extension and restriction.
In \cref{sec:graded-coalgebras-over-graded-algebras}, we discuss graded coalgebras over graded algebras.
\Cref{subsec:pre-differential-graded-algebras-modules} introduces the pre-differential structure,
where we define pre-differential graded objects and the notion of derivations and $d$-operators.
We extend the pre-differential graded framework to coalgebras and coderivations in \cref{sec:pre-differential-graded-coalgebras}.
Finally, in \cref{subsec:differential-graded-algebras-modules}, we discuss differential graded objects and their cohomology.

Our notation and discussion are mostly standard and can be skipped, referring to specific points if necessary.
There are however a few points in which we possibly differ from standard treatments:
\begin{enumerate}
	\item We work with an arbitrary grading datum and use it consistently for grading and inserting signs determining
	      the Koszul sign convention.
	\item We work with the notion of maps between modules over different ground algebras,
	      which is equivalent to the standard notion involving pullbacks.
	\item We discuss morphisms between coalgebras defined over different ground algebras.
	\item We work with pre-differential graded algebras equipped with a derivation $d$, not necessarily
	      satisfying $d^2 = 0$.
	\item We introduce the notion of a $d$-operator over an algebra derivation $d$ and use it to
	      define pre-differential graded objects over pre-differential graded algebras, and we explore their structure.
	\item We introduce and discuss the notion of a generalized coderivation on a coalgebra,
	      which is simultaneously a coderivation and a $d$-operator.
\end{enumerate}

\subsection{Grading Datum} \label{subsec:grading-data}
As we will often encounter objects graded by different groups and work with
several different symmetries, it will prove useful to isolate and make precise the data
needed to work with graded objects in our setting.

A \textbf{grading datum} is a pair $\left( \GG, \braidop \right)$ where $\GG$ is an abelian group
and $\braidop \colon \GG \times \GG \rightarrow \ZZ_2$ is a \textbf{parity form}, i.e., a $\ZZ$-bilinear
symmetric map. We will use the group
$\GG$ to grade objects and the parity form $\braidop$ to endow the category of $\GG$-graded objects with
symmetry maps
\begin{equation*}
	m \otimes n \mapsto (-1)^{\braidd{m}{n}} n \otimes m
\end{equation*}
which will play a role in various sign rules. Some examples of grading data we will work with include:

\begin{enumerate}
	\item The group $\GG = \ZZ$ with the parity form $\braid{a}{b} = a \cdot b \mod 2$. In this case
	      we will call the grading datum $\left( \ZZ, \braidop \right)$ the \textbf{standard Koszul
		      grading datum}.
	\item The group $\GG = \ZZ^2$ with the parity form given by
	      \begin{equation*}
		      \braid{\left( a_1, a_2 \right)}{\left( b_1, b_2 \right)}_1 \defeq a_1 \cdot b_1 + a_2 \cdot b_2
		      \mod 2.
	      \end{equation*}
	\item The group $\GG = \ZZ^2$ with the parity form given by
	      \begin{equation*}
		      \braid{\left( a_1, a_2 \right)}{\left( b_1, b_2 \right)}_2 \defeq
		      \left( a_2 - a_1 \right) \cdot \left( b_2 - b_1 \right) \mod 2.
	      \end{equation*}
\end{enumerate}

When working with pre-differential graded objects, we will assume in addition that we have
a fixed element $\go \in \GG$ which is considered part of the grading datum. Pre-differentials
will raise the degree of elements by $\go$. Finally, when working
with differential graded objects, we will assume that the fixed element $\go$
is \textbf{odd} in the sense that
\begin{equation}
	\braid{\go}{\go} \equiv 1 \mod 2.
\end{equation}

Although much of the following will depend on the grading data $\GG, \braidop$ and $\go$,
we will often suppress them from our notation, as long as they are fixed.

\subsection{Graded \texorpdfstring{$\mathbbm{k}$}{k}-modules} \label{subsec:graded-k-modules}
Recall that a graded $\mathbbm{k}$-module can be defined in two equivalent ways:
externally, as an indexed family of $\mathbbm{k}$-modules, or internally, as a $\mathbbm{k}$-module together
with a direct sum decomposition. In our work, we adopt the external point of view.

\begin{dfn} \label{dfn:graded-k-module}
	A \textbf{graded} $\mathbbm{k}$\textbf{-module}
	$M = \left( M^g \right)_{g \in \GG}$ is an indexed family of $\mathbbm{k}$-modules, called
	the \textbf{components} of $M$.
\end{dfn}
An element $m \in M^d$ is said to be a \textbf{homogeneous element of degree}
$d \in \GG$. By an \textbf{element} $m \in M$ we will always mean a homogeneous
element $m \in M^d$ and use the notation $\degb{m} = d$ for the degree of $m$. We will work
only with homogeneous elements and never add elements of different degrees.

Given two graded $\mathbbm{k}$-modules $M$ and $N$, a \textbf{graded map} $f \colon M \rightharpoonup N$
\textbf{of degree} $d \in \GG$ is a family $f = \left( f^g \colon M^g \rightarrow N^{g+d} \right)_{g \in \GG}$
of homomorphisms of $\mathbbm{k}$-modules, called the \textbf{components} of $f$. We will denote the degree $d$
of $f$ by $\degb{f} = d$. Strictly speaking, a graded map $f$ is a family of maps. However, given
an element $m \in M$ in the sense defined above, we will write $f \left( m \right)$
for the element $f^{\degb{m}} \left( m \right) \in N^{\degb{m} + \degb{f}}$ of $N$. Composition of
graded maps is defined component-wise in the natural way.

A \textbf{morphism} $f \colon M \rightarrow N$ \textbf{of graded} $\mathbbm{k}$-\textbf{modules} is a graded map of degree zero. We denote the set of morphisms between $M$ and $N$ by
$\Hom{M}{N}[][\mathbbm{k}]$ and by $\GMod[\mathbbm{k}]$ the category of graded $\mathbbm{k}$-modules with morphisms
as defined above. Since we work with both degree zero and arbitrary degree graded maps, we
use the notation $f \colon M \rightarrow N$ for degree zero maps (i.e., the morphisms in the category
$\GMod[\mathbbm{k}]$) while reserving the notation $f \colon M \rightharpoonup N$ for general graded maps.

Note that the category $\GMod[\mathbbm{k}]$ is precisely the category of graded objects of the category
$\Mod[\mathbbm{k}]$ of ungraded $\mathbbm{k}$-modules.
Since $\Mod[\mathbbm{k}]$ is bicomplete, the category $\GMod[\mathbbm{k}]$ is also
bicomplete and both limits and colimits are computed component-wise in $\Mod[\mathbbm{k}]$.
The category $\GMod[\mathbbm{k}]$ is $\mathbbm{k}$-linear and abelian. An isomorphism
in the category $\GMod[\mathbbm{k}]$ is a morphism $f \colon M \rightarrow N$ for which each component
$f^g \colon M^g \rightarrow N^g$ is an isomorphism of $\mathbbm{k}$-modules.

The category $\GMod[\mathbbm{k}]$ has a standard closed symmetric monoidal structure which we now recall.
Given two graded $\mathbbm{k}$-modules $M$ and $N$,
their (graded) tensor product $M \otimes_{\mathbbm{k}} N$ is the graded $\mathbbm{k}$-module
whose components are given by
\begin{equation}
	\left( M \otimes_{\mathbbm{k}} N \right)^g =
	\bigoplus_{g_1 + g_2 = g} M^{g_1} \otimes_{\mathbbm{k}} N^{g_2}.
	\label{eq:tensor-product-graded-k-modules}
\end{equation}

Similar to the ungraded case, the tensor product of graded $\mathbbm{k}$-modules can be
characterized by a universal property involving graded $\mathbbm{k}$-bilinear maps. The notion
of $\mathbbm{k}$-bilinear and $\mathbbm{k}$-multilinear maps extends naturally to the graded setting as follows:
Given graded $\mathbbm{k}$-modules $M_1,\dots,M_n$ and $N$, a \textbf{graded}
$\mathbbm{k}$-\textbf{multilinear map}
$B \colon M_1 \times \dots \times M_n \rightharpoonup N$ \textbf{of degree} $d \in \GG$ is a family
\begin{equation*}
	B = \left( B^{g_1, \dots, g_n} \colon M_{1}^{g_1} \times \dots \times M_{n}^{g_n} \rightarrow
	N^{g_1 + \dots + g_n + d} \right)_{g_1,\dots,g_n \in \GG}
\end{equation*}
of $\mathbbm{k}$-multilinear maps, called the \textbf{components} of $B$.
We will denote the degree $d$ of $B$ by $\degb{B} = d$. Again, although strictly speaking a graded
$\mathbbm{k}$-multilinear map $B$ is not really a map, we will treat it as a map and
given elements $m_1 \in M_1, \dots, m_n \in M_n$, we write
$B \left( m_1, \dots, m_n \right)$ for the element
$B^{\degb{m_1} + \dots + \degb{m_n}} \left( m_1, \dots, m_n \right) \in N^{\degb{m_1} + \dots + \degb{m_n} + \degb{B}}$.
Similar to our convention with graded maps, we will use
the notation $B \colon M_1 \times \dots \times M_n \rightharpoonup N$ to denote graded $\mathbbm{k}$-multilinear
maps of arbitrary degree while reserving the notation $B \colon M_1 \times \dots \times M_n \rightarrow N$
to denote graded $\mathbbm{k}$-multilinear maps of degree zero.

The tensor product $M \otimes_{\mathbbm{k}} N$ comes equipped with a canonical graded $\mathbbm{k}$-bilinear map
$\otimes_{\mathbbm{k}} \colon M \times N \rightarrow M \otimes_{\mathbbm{k}} N$ of degree zero
characterized by the following universal property: Given a graded $\mathbbm{k}$-module $L$ and
a graded $\mathbbm{k}$-bilinear map $B \colon M \times N \rightharpoonup L$, there exists a unique
graded map $\varphi_B \colon M \otimes_{\mathbbm{k}} N \rightharpoonup L$
with $\degb{\varphi_B} = \degb{B}$ such that
$\varphi_B \left( m \otimes_{\mathbbm{k}} n \right) = B(m,n)$ for all $m \in M$ and $n \in N$ (see
\cref{fig:tensor-product-graded-k-modules-universal-property}).
\begin{figure}[htb]
	\centering
	\begin{tikzcd}
		{M \times N} && {M \otimes_{\mathbbm{k}} N} \\
		&& L
		\arrow["\otimes_{\mathbbm{k}}", from=1-1, to=1-3]
		\arrow["\substack{\exists! \, \varphi_B \\ \textrm{ graded map}}", dashed, harpoon, from=1-3, to=2-3]
		\arrow["\substack{B \textrm{ graded} \\ \mathbbm{k}\textrm{-bilinear}}"', harpoon, from=1-1, to=2-3]
	\end{tikzcd}
	\caption{Universal property of the tensor product of graded $\mathbbm{k}$-modules.}
	\label{fig:tensor-product-graded-k-modules-universal-property}
\end{figure}

The tensor product of two graded maps $f \colon M \rightharpoonup M'$ and
$g \colon N \rightharpoonup N'$ between graded $\mathbbm{k}$-modules is
the unique graded map
$f \otimes_{\mathbbm{k}} g \colon M \otimes_{\mathbbm{k}} N \rightharpoonup M' \otimes_{\mathbbm{k}} N'$
of degree $\degb{f} + \degb{g}$ which satisfies
\begin{equation}
	\left( f \otimes_{\mathbbm{k}} g \right) \left( m \otimes_{\mathbbm{k}} n \right) =
	(-1)^{\braidd{g}{m}} f \left( m \right) \otimes_{\mathbbm{k}} g \left( n \right)
	\label{eq:tensor-product-graded-maps}
\end{equation}
for all $m \in M$ and $n \in N$. Given graded maps $f_i \colon M_i \rightharpoonup N_i$ and
$g_i \colon L_i \rightharpoonup M_i$ for $i=1,2$, we have the formula
\begin{equation}
	\left( f_1 \otimes_{\mathbbm{k}} f_2 \right) \circ \left( g_1 \otimes_{\mathbbm{k}} g_2 \right) =
	(-1)^{\braidd{g_1}{f_2}} \left( f_1 \circ g_1 \right) \otimes_{\mathbbm{k}} \left( f_2 \circ g_2 \right)
	\label{eq:interaction-composition-tensor-product}
\end{equation}
which shows that the interaction between the composition and tensor product of graded maps is consistent
with the Koszul sign rule. In particular, when all the maps involved have degree zero, the tensor product
is functorial on the nose without a sign factor, and we have a bifunctor
$\otimes_{\mathbbm{k}} \colon \GMod[\mathbbm{k}] \times \GMod[\mathbbm{k}] \rightarrow \GMod[\mathbbm{k}]$.
The bifunctor $\otimes_{\mathbbm{k}}$, the symmetry maps
$M \otimes_{\mathbbm{k}} N \rightarrow N \otimes_{\mathbbm{k}} M$ given by
\begin{equation}
	m \otimes_{\mathbbm{k}} n \mapsto (-1)^{\braidd{m}{n}} n \otimes_{\mathbbm{k}} m,
	\label{eq:symmetry-graded-k-modules}
\end{equation}
and the standard associators and unitors,
endow the category $\GMod[\mathbbm{k}]$ with the structure of a symmetric monoidal category whose unit is
the ground ring $\mathbbm{k}$ (considered as graded $\mathbbm{k}$-module concentrated in degree zero).
It follows from the universal property of the graded tensor product that the symmetric monoidal category
$\GMod[\mathbbm{k}]$ is closed with the internal hom object
given by the graded $\mathbbm{k}$-module $\InnHom{M}{N}[][\mathbbm{k}]$ of all graded
maps where
\begin{equation}
	\begin{aligned}
		\InnHom{M}{N}[][\mathbbm{k}]^d & \defeq \Set{f \colon M \rightharpoonup N}
		                                        [f \textrm{ is a graded map of degree } d].
	\end{aligned} \label{eq:inner-hom-graded-k-module}
\end{equation}

\phantomsection
\label{sec:suspension-graded-k-module}
Given a graded $\mathbbm{k}$-module $M$ and $h \in \GG$, the
$h$-\textbf{suspension} or $h$-\textbf{shifted} module $M[h]$ is defined by setting
$M[h]^g \defeq M^{g+h}$ for all $g \in \GG$. Given an element $m \in M^g$ of degree $g$,
we will denote the same element, considered as an element of $M[h]$ of degree $g - h$ by
$\s_h \left( m \right)$. We also think of $\s_h \colon M \rightharpoonup M[h]$ as a graded map
of degree $-h$ whose underlying components $\s_h^g \colon M^g \rightarrow M[h]^{g-h} = M^g$ are the identity maps.

\subsection{Graded \texorpdfstring{$\mathbbm{k}$}{k}-algebras} \label{subsec:graded-k-algebras}
Since the category $\GMod[\mathbbm{k}]$ is monoidal, one can talk about algebra objects in $\GMod[\mathbbm{k}]$
(see \cref{subsec:alg-in-monoidal-cat}).
A graded $\mathbbm{k}$-algebra is an algebra object of $\GMod[\mathbbm{k}]$. Unwinding the definition, we have:

\begin{dfn} \label{dfn:graded-k-algebra}
	A \textbf{graded} $\mathbbm{k}$-\textbf{algebra} is a graded $\mathbbm{k}$-module
	$R = \left( R^g \right)_{g \in \GG}$
	together with a degree zero $\mathbbm{k}$-bilinear multiplication $\cdot \colon R \times R \rightarrow R$
	and a unit element $1_R \in R^0$ such that $a \cdot \left( b \cdot c \right) = \left( a \cdot b \right) \cdot c$
	and $a \cdot 1_R = 1_R \cdot a = a$ for all $a,b,c \in R$.\footnote{An algebra object $R$ of $\GMod[\mathbbm{k}]$
		is an object of $\GMod[\mathbbm{k}]$ equipped with a multiplication morphism of the form
		$m \colon R \otimes R \rightarrow R$ but by the universal property of the graded tensor product,
		such a multiplication corresponds bijectively to a degree zero $\mathbbm{k}$-bilinear multiplication
		$\cdot \colon R \times R \rightarrow R$. Similarly, we identify the unit morphism
		$u \colon \mathbbm{k} \rightarrow R$ with the element $u \left( 1_{\mathbbm{k}} \right) \defeq 1_R$
		which belongs to $R^0$ since $u$ has degree zero.}

	When $\mathbbm{k} = \ZZ$, a graded $\mathbbm{k}$-algebra
	is called a \textbf{graded ring}.
\end{dfn}

A \textbf{morphism of graded} $\mathbbm{k}$-\textbf{algebras} $f \colon R \rightarrow S$
is a morphism of graded $\mathbbm{k}$-modules (i.e., a degree zero graded map) which satisfies
$f \left( a \cdot b \right) = f \left( a \right) \cdot f \left( b \right)$ for all $a,b \in R$
and $f \left( 1_R \right) = 1_S$.

\begin{ex} \label{ex:inner-end-as-a-graded-algebra}
	Let $M$ be a graded $\mathbbm{k}$-module. Then the graded $\mathbbm{k}$-module
	\begin{equation*}
		\InnEnd{M}[][\mathbbm{k}] \defeq \InnHom{M}{M}[][\mathbbm{k}],
	\end{equation*}
	consisting of all graded maps $f \colon M \rightharpoonup M$, is a graded $\mathbbm{k}$-algebra,
	with composition of graded maps as multiplication. The graded $\mathbbm{k}$-module
	$\InnEnd{M}[][\mathbbm{k}]$ also has the structure of a graded Lie algebra over $\mathbbm{k}$, where the
	Lie bracket is given by the \textbf{graded commutator}
	\begin{equation}
		\left[ f, g \right] \defeq f \circ g - (-1)^{\braidd{f}{g}} g \circ f.
		\label{eq:graded-commutator}
	\end{equation}
\end{ex}

A graded $\mathbbm{k}$-algebra $R$ is called \textbf{graded-commutative} if
$r \cdot s = (-1)^{\braidd{r}{s}} s \cdot r$
for all $r,s \in R$. Equivalently, a graded-commutative $\mathbbm{k}$-algebra is a commutative algebra object
of $\GMod[\mathbbm{k}]$ with respect to the symmetry maps given by \cref{eq:symmetry-graded-k-modules}.

\subsection{Graded Modules over Graded Algebras}
\label{sub:graded-R-modules}
Generalizing \cref{subsec:graded-k-modules}, we discuss graded modules over graded $\mathbbm{k}$-algebras.
Since the category $\GMod[\mathbbm{k}]$ is monoidal, one can talk about module objects over algebra
objects of $\GMod[\mathbbm{k}]$ (see \cref{subsec:modules-in-monoidal-cat}).
Let $R = \left( R^g \right)_{g \in \GG}$ be a graded $\mathbbm{k}$-algebra, i.e., an algebra object
of $\GMod[\mathbbm{k}]$. Unwinding the definition, we see that a graded left $R$-module is given by:

\begin{dfn} \label{dfn:graded-R-module}
	A \textbf{graded left} $R$-\textbf{module} is a graded $\mathbbm{k}$-module
	$M = \left( M^g \right)_{g \in \GG}$ together with a degree zero $\mathbbm{k}$-bilinear action map
	$\cdot \colon R \times M \rightarrow M$ such that
	$\left( r \cdot s \right) \cdot m = r \cdot \left( s \cdot m \right)$ and $1_R \cdot m = m$
	for all $r, s \in R$ and $m \in M$.
\end{dfn}
Similarly, one can define graded right modules and bimodules. In what follows,
unless explicitly mentioned otherwise, the term ``module'' will always mean left module.
We note that every graded $R$-module $M$ has an \textbf{underlying graded} $\mathbbm{k}$-\textbf{module}
obtained by forgetting the $R$-action. The $R$-action on $M$ is implicit in our notation and, when necessary,
we forget the $R$-action and think of $M$ as a graded $\mathbbm{k}$-module without further mention.
We also note that if we consider $\mathbbm{k}$ as a graded $\mathbbm{k}$-algebra concentrated in degree zero,
then a graded $\mathbbm{k}$-module in the sense of \cref{dfn:graded-R-module} is the same as a
graded $\mathbbm{k}$-module in the sense of \cref{dfn:graded-k-module}
(i.e., an indexed family of $\mathbbm{k}$-modules).

Given two graded $R$-modules $M$ and $N$, we can endow the graded $\mathbbm{k}$-module
$\InnHom{M}{N}[][\mathbbm{k}]$ with a natural structure of
an $(R,R)$-bimodule via the \textbf{outer action}
\begin{equation}
	\left( r \cdot f \right) \left( m \right) \defeq r \cdot f \left( m \right) \label{eq:left-R-action-on-hom}
\end{equation}
and \textbf{inner action}
\begin{equation}
	\left( f \cdot r \right) \left( m \right) \defeq f \left( r \cdot m \right). \label{eq:right-R-action-on-hom}
\end{equation}
A graded map $f \colon M \rightharpoonup N$ is called (left) $R$-\textbf{linear} if
\begin{equation*}
	f \left( r \cdot m \right) = (-1)^{\braidd{f}{r}} r \cdot f \left( m \right)
\end{equation*}
for all $r \in R$ and $m \in M$, or, equivalently, in terms of the $R$-actions on $\InnHom{M}{N}[][\mathbbm{k}]$,
if $f \cdot r = (-1)^{\braidd{r}{f}} r \cdot f$.
The collection of all graded $R$-linear maps is denoted by
\begin{equation}
	\InnHom{M}{N}[][R] \defeq \Set{f \in \InnHom{M}{N}[][\mathbbm{k}]}
	[f \textrm{ is } R\textrm{-linear}]
	\label{eq:inner-hom-graded-R-module}
\end{equation}
and forms a graded $\mathbbm{k}$-submodule of $\InnHom{M}{N}[][\mathbbm{k}]$.

A \textbf{morphism} $f \colon M \rightarrow N$ \textbf{of graded}
$R$\textbf{-modules} is a graded $R$-linear map of degree zero. We denote
the set of morphisms between $M$ and $N$ by $\Hom{M}{N}[][R]$
and by $\GMod[R]$ the category of graded $R$-modules with morphisms
as defined above. The category $\GMod[R]$ is bicomplete
where both limits and colimits are computed in the category of graded $\mathbbm{k}$-modules and
endowed with the natural $R$-module structures. In addition, the category $\GMod[R]$
is $R^0$-linear and abelian. An isomorphism
in the category $\GMod[R]$ is a degree zero $R$-linear map $f \colon M \rightarrow N$ for which each component
$f^g \colon M^g \rightarrow N^g$ is an isomorphism of $\mathbbm{k}$-modules.

\phantomsection
\label{sec:suspension-graded-R-module}
Given a graded $R$-module $M$ and $h \in \GG$, the
$h$-\textbf{suspension} or $h$-\textbf{shifted} module $M[h]$ is defined
to be the suspension of the underlying graded $\mathbbm{k}$-module (see \cref{sec:suspension-graded-k-module})
together with the $R$-action given by
\begin{equation}
	r \cdot \s_h \left( m \right) \defeq (-1)^{\braid{\degb{r}}{-h}} \s_h \left( r \cdot m \right)
	= (-1)^{\braid{\degb{r}}{h}} \s_h \left( r \cdot m \right).
	\label{eq:R-action-on-suspension}
\end{equation}
With the definition given by \cref{eq:R-action-on-suspension}, the suspension map
$\s_h \colon M \rightharpoonup M[h]$ becomes a graded $R$-linear map of degree $-h$.

Given a graded right $R$-module $M$ and a graded left $R$-module $N$,
the graded tensor product $M \otimes_R N$ is the graded $\mathbbm{k}$-module defined by
\begin{align}
	\MoveEqLeft
	\left( M \otimes_R N \right)^g \defeq{}
	\left( \bigoplus_{a + b = g} M^{a} \otimes_{\mathbbm{k}} N^{b} \right)
	\Big/ \label{eq:algebraic-graded-tensor-product}
	\\
	 &
	\left< mr \otimes_{\mathbbm{k}} n - m \otimes_{\mathbbm{k}} rn \, \middle| \,
	r \in R^{g_1}, m \in M^{g_2}, n \in N^{g_3}, \, g_1 + g_2 + g_3 = g \right>.	\nonumber
\end{align}
Given elements $m \in M$ and $n \in N$, we will denote by $m \otimes_R n \in M \otimes_R N$
the equivalence class of $m \otimes_{\mathbbm{k}} n$ and call such an element an \textbf{elementary tensor}.
Note that each graded component $\left( M \otimes_R N \right)^d$ of $M \otimes_R N$ is generated as a
$\mathbbm{k}$-module by elementary tensors of the form $m \otimes_R n$ where $\degb{m} + \degb{n} = d$.

Similar to the ungraded case, the tensor product $M \otimes_R N$ of
a graded right $R$-module $M$ and a graded left $R$-module $N$ can be characterized by a universal property
involving graded $R$-balanced maps. Given a graded $\mathbbm{k}$-module $L$, a graded $\mathbbm{k}$-bilinear map
$B \colon M \times N \rightharpoonup L$ is called $R$-\textbf{balanced} if it satisfies
\begin{equation}
	B \left( mr, n \right) = B \left( m, rn \right) \label{eq:R-balanced-map}
\end{equation}
for all $r \in R, m \in M, n \in N$.\footnote{Recall that
	we always work with homogeneous elements and suppress components. Written explicitly, the $R$-balanced condition means that $B^{g_1 + h, g_2} \left( m^{g_1} \cdot r^{h}, n^{g_2} \right)
		= B^{g_1, h + g_2} \left( m^{g_1}, r^{h} \cdot n^{g_2} \right)$ for all
	elements $m^{g_1} \in M^{g_1}, n^{g_2} \in N^{g_2}, r^{h} \in R^{h}$ and all
	$g_1,g_2,h \in \GG$.}
The tensor product $M \otimes_{R} N$ comes equipped with a canonical graded $R$-balanced map
$\otimes_{R} \colon M \times N \rightarrow M \otimes_{R} N$ of degree zero
characterized by the following universal property: Given a graded $\mathbbm{k}$-module $L$ and
a graded $R$-balanced map $B \colon M \times N \rightharpoonup L$, there exists a unique
graded map $\varphi_B \colon M \otimes_{R} N \rightharpoonup L$ of graded $\mathbbm{k}$-modules
with $\degb{\varphi_B} = \degb{B}$ such that
$\varphi_B \left( m \otimes_{R} n \right) = B(m,n)$ for all $m \in M$ and $n \in N$ (see
\cref{fig:tensor-product-graded-R-modules-universal-property-balanaced}).
\begin{figure}[htb]
	\centering
	\begin{subfigure}{0.45\textwidth}
		\centering
		\begin{tikzcd}
			{M \times N} && {M \otimes_{R} N} \\
			&& L \\
			\arrow["\otimes_{R}", from=1-1, to=1-3]
			\arrow["{\substack{\exists! \, \varphi_B \\ \mathbbm{k}\textrm{-linear}}}", dashed, harpoon, from=1-3, to=2-3]
			\arrow["\substack{B \\ R\textrm{-balanced}}"', harpoon, from=1-1, to=2-3]
		\end{tikzcd}
		\caption{$R$ is a graded $\mathbbm{k}$-algebra, \\ $M$ is a graded right $R$-module, \\
			$N$ is a graded left $R$-module, \\ $L$ is a graded $\mathbbm{k}$-module.}
		\label{fig:tensor-product-graded-R-modules-universal-property-balanaced}
	\end{subfigure}
	\begin{subfigure}{0.45\textwidth}
		\centering
		\begin{tikzcd}
			{M \times N} && {M \otimes_{R} N} \\
			&& L \\
			\arrow["\otimes_{R}", from=1-1, to=1-3]
			\arrow["{\substack{\exists! \, \varphi_B \\ R\textrm{-linear}}}", dashed, harpoon, from=1-3, to=2-3]
			\arrow["\substack{B \\ R\textrm{-bilinear}}"', harpoon, from=1-1, to=2-3]
		\end{tikzcd}
		\caption{$R$ is a graded-commutative $\mathbbm{k}$-algebra, \\ $M,N,L$ are graded left $R$-modules. \\ \, \\ \,}
		\label{fig:tensor-product-graded-R-modules-universal-property}
	\end{subfigure}
	\caption{Universal properties of the tensor product of graded $R$-modules.}
\end{figure}

In what follows, we will assume that $R$ is graded-commutative.
When working over a graded-commutative $\mathbbm{k}$-algebra, any left $R$-module $M$ has a natural structure
of a right $R$-module via the action
\begin{equation}
	m \cdot r \defeq (-1)^{\braidd{r}{m}} r \cdot m. \label{eq:convert-left-to-right-module}
\end{equation}
Endowing $M$ with both actions, one obtains a \textbf{symmetric bimodule}. When $R$ is
a graded-commutative $\mathbbm{k}$-algebra, we identify left, right and symmetric bimodules as necessary.
In this case, given two graded left $R$-modules $M$ and $N$, we can convert $M$ to a graded right $R$-module
and form the tensor product $M \otimes_R N$. The resulting object $M \otimes_R N$ is not only a graded
$\mathbbm{k}$-module but also has a natural structure of a graded left $R$-module where the
$R$-action is given on elementary tensors by
\begin{equation}
	r \cdot \left( m \otimes_R n \right) \defeq \left( r \cdot m \right) \otimes_R n =
		(-1)^{\braidd{r}{m}} m \otimes_R \left( r \cdot n \right). \label{eq:R-action-R-tensor-product}
\end{equation}

The tensor product $M \otimes_R N$ of two graded $R$-modules over a graded-commutative ground algebra
$R$ can be
characterized by a universal property involving graded $R$-bilinear maps.
Let $M_1,\dots,M_n$ and $N$ be graded $R$-modules. A \textbf{graded} $R$-\textbf{multilinear map}
$B \colon M_1 \times \dots \times M_n \rightharpoonup N$ \textbf{of degree}
$d \in \GG$ is a degree $d$ $\mathbbm{k}$-multilinear map of graded $\mathbbm{k}$-modules
which is $R$-\textbf{multilinear} in the sense that
\begin{equation*}
	B \left( m_1, \dots, m_{i-1}, r \cdot m_i, m_{i+1}, \dots, m_n \right) =
	(-1)^{\braid{\degb{r}}{d + \degb{m_1} + \dots + \degb{m_{i-1}}}} r \cdot B \left( m_1, \dots, m_n \right)
\end{equation*}
for all $1 \leq i \leq n$ and $r \in R, m_1 \in M_1, \dots, m_n \in M_n$.

The canonical map $\otimes_R \colon M \times N \rightarrow M \otimes_R N$ is not only $R$-balanced
but also $R$-bilinear and is characterized by the following universal property:
Given a graded $R$-module $L$ and a graded $R$-bilinear map $B \colon M \times N \rightharpoonup L$, there exists a
unique graded $R$-linear map $\varphi_B \colon M \otimes_{R} N \rightharpoonup L$ of graded
$R$-modules with $\degb{\varphi_B} = \degb{B}$ such that
$\varphi_B \left( m \otimes_{R} n \right) = B(m,n)$ for all $m \in M$ and $n \in N$ (see
\cref{fig:tensor-product-graded-R-modules-universal-property}).

The universal properties of the tensor product allow us to abuse notation and define maps
on $M \otimes_R N$ by specifying their action on elementary tensors. For such a definition to be
well-defined, one must check that the resulting map is $R$-balanced (or $R$-bilinear) and then invoke
the universal property.

The tensor product of two graded $R$-linear maps $f \colon M \rightharpoonup M'$ and
$g \colon N \rightharpoonup N'$ is the unique graded $R$-linear map
$f \otimes_R g \colon M \otimes_R N \rightharpoonup M' \otimes_R N'$ of degree $\degb{f} + \degb{g}$
which satisfies
\begin{equation}
	\left( f \otimes_R g \right) \left( m \otimes_R n \right) =
	(-1)^{\braidd{g}{m}} f \left( m \right) \otimes_R g \left( n \right)
	\label{eq:tensor-product-graded-R-linear-maps}
\end{equation}
for all $m \in M$ and $n \in N$.\footnote{Note that the sign in \cref{eq:tensor-product-graded-R-linear-maps}
	is necessary and not merely a convention consistent with the Koszul sign rule, as the ``map''
	$m \otimes_R n \mapsto f \left( m \right) \otimes_R g \left( n \right)$ is not necessarily well-defined when
	$f$ and $g$ are $R$-linear of non-zero degree.\label[footnote]{foot:sign-necessary-tensor-product}}
The interaction between composition and tensor product of graded maps is the same as in
\cref{eq:interaction-composition-tensor-product}, with $\otimes_{\mathbbm{k}}$ replaced by $\otimes_R$.
In particular, when all the graded $R$-linear maps have degree zero, the tensor product is functorial
on the nose without a sign factor, and we have a bifunctor
$\otimes_R \colon \GMod[R] \times \GMod[R] \rightarrow \GMod[R]$. The bifunctor $\otimes_R$,
symmetry maps $M \otimes_R N \rightarrow N \otimes_R M$ given by
\begin{equation}
	m \otimes_R n \mapsto (-1)^{\braidd{m}{n}} n \otimes_R m, \label{eq:symmetry-graded-R-modules}
\end{equation}
and the standard associators and unitors, endow the category $\GMod[R]$ with the structure of a symmetric monoidal category whose unit is the ground algebra $R$, considered as a graded $R$-module over itself.

Since $R$ is graded-commutative, given a graded $R$-linear map $f \colon M \rightharpoonup N$ and $r \in R$,
the map $r \cdot f$ given by \cref{eq:left-R-action-on-hom} is also $R$-linear and hence
$\InnHom{M}{N}[][R]$ has the structure of a graded $R$-module, i.e., it is an object of $\GMod[R]$.
It follows from the universal property of the graded tensor product that the symmetric monoidal
category $\GMod[R]$ is closed with the internal hom object given by $\InnHom{M}{N}[][R]$.

More generally, let us denote by $\Mult{M_1,\dots,M_n}{N}[R]$ the graded $\mathbbm{k}$-module
of all graded $R$-multilinear maps $B \colon M_1 \times \dots \times M_n \rightharpoonup N$. Since
$R$ is graded-commutative, $\Mult{M_1,\dots,M_n}{N}[R]$ has a natural $R$-action given by
\begin{equation}
	\left( r \cdot B \right) \left( m_1, \dots, m_n \right) \defeq
	r \cdot B \left( m_1, \dots, m_n \right) \label{eq:R-action-on-multilinear-maps}
\end{equation}
with respect to which $\Mult{M_1,\dots,M_n}{N}[R]$ becomes a graded $R$-module,\footnote{In particular,
	the map $r \cdot B$ given by \cref{eq:R-action-on-multilinear-maps}
	is also $R$-multilinear of degree $\degb{r} + \degb{B}$.}
and we have natural isomorphisms of graded $R$-modules
\begin{align*}
	 & \Mult{M_1,\dots,M_n}{N}[R] \cong \InnHom{M_1 \otimes_R \dots \otimes_R M_n}{N}[][R],
	\\
	 & \Mult{M_1,M_2}{N}[R] \cong \InnHom{M_1}{\InnHom{M_2}{N}[][R]}[][R] \cong \InnHom{M_1 \otimes_R M_2}{N}[][R],
\end{align*}
given by the standard formulas generalizing the ungraded case.

Since $\GMod[R]$ is closed symmetric monoidal,
the tensor product of graded $R$-modules commutes with colimits in both variables, and we have natural
isomorphisms
\begin{equation*}
	\left( \bigoplus_{i \in I} M_i \right) \otimes_R N \cong
	\bigoplus_{i \in I} M_i \otimes_R N, \quad
	M \otimes_R \left( \bigoplus_{i \in I} N_i \right) \cong
	\bigoplus_{i \in I} M \otimes_R N_i.
\end{equation*}

\phantomsection
\label{sec:scalar-extension-restriction-graded-modules}

The scalar extension and restriction functors extend naturally to the graded setting.
Let $R$ and $S$ be two graded-commutative $\mathbbm{k}$-algebras and let
$\varphi \colon R \rightarrow S$ be a morphism of graded $\mathbbm{k}$-algebras. Given a graded
$S$-module $N$, denote by $\varphi^{*} \left( N \right)$ the graded $R$-module whose
underlying graded $\mathbbm{k}$-module is the same as $N$, equipped with the $R$-action given by
$r \cdot n \defeq \varphi \left( r \right) \cdot n$. The $R$-module $\varphi^{*} \left( N \right)$
is called the \textbf{restriction of scalars}, or \textbf{pullback}, of $N$ along $\varphi$.
Given a graded $S$-linear map
$f \colon N_1 \rightharpoonup N_2$ between two graded $S$-modules, we denote by
$\varphi^{*} \left( f \right) \colon \varphi^{*} \left( N_1 \right)
	\rightharpoonup \varphi^{*} \left( N_2 \right)$ the same graded map $f$, considered as
a graded map between two graded $R$-modules. Note that $\varphi^{*} \left( f \right)$ is indeed
$R$-linear since $\varphi$ has degree zero. This gives us a functor
$\varphi^{*} \colon \GMod[S] \rightarrow \GMod[R]$ whose action on morphisms is the trivial one.

The functor $\varphi^{*}$ has a left adjoint, denoted by $\varphi_{!}$, which we now describe.
Note that the graded-commutative $\mathbbm{k}$-algebra $S$ can be considered as a graded $(S,R)$-bimodule
where the left action is given by the multiplication on $S$ and the right action is given by $\left( s, r \right) \mapsto s \cdot \varphi(r)$. Given a graded $R$-module $M$, we can form the graded tensor product
$\varphi_{!} \left( M \right) \defeq S \otimes_R M$
which has a canonical structure of a graded $S$-module via the left action on $S$.
The $S$-module $\varphi_{!} \left( M \right)$ is called the \textbf{extension of scalars} of $M$ along $\varphi$.
Given a graded
$R$-linear map $f \colon M_1 \rightharpoonup M_2$ between two graded $R$-modules, define
a graded $S$-linear map $\varphi_{!} \left( f \right) \colon \varphi_{!} \left( M_1 \right) \rightharpoonup
	\varphi_{!} \left( M_2 \right)$ by
\begin{equation*}
	\varphi_{!} \left( f \right) \left( s \otimes_R m_1 \right) \defeq
	\left( \id_S \otimes_R f \right) \left( s \otimes_R m_1 \right) =
	(-1)^{\braidd{f}{s}} s \otimes_R f \left( m_1 \right).
\end{equation*}
This gives us a functor $\varphi_{!} \colon \GMod[R] \rightarrow \GMod[S]$ such that we have
an adjunction
\begin{equation}
	\varphi_{!} \colon \GMod[R] \stackrel[]{\dashv}{\rightleftarrows} \GMod[S] \colon \varphi^{*}
	\label{eq:restriction-extension-adjunction-gmod-R}
\end{equation}
with unit and counit maps given by
\begin{align*}
	 & M \rightarrow \varphi^{*} \left( \varphi_{!} \left( M \right) \right) \colon &  & m \mapsto 1_S \otimes_R m,
	\\
	 & \varphi_{!} \left( \varphi^{*} \left( N \right) \right) \rightarrow N \colon &  & s \otimes_R n \mapsto s \cdot n.
\end{align*}

The adjunction is in fact a monoidal adjunction.
The left adjoint functor $\varphi_{!}$ is strong symmetric
monoidal via the tensor constraints
\begin{equation*}
	\varphi_{!} \left( M_1 \right) \otimes_S \varphi_{!} \left( M_2 \right) =
	\left( S \otimes_R M_1 \right) \otimes_S \left( S \otimes_R M_2 \right)
	\xrightarrow[]{\cong}
	S \otimes_R \left( M_1 \otimes_R M_2 \right) = \varphi_{!} \left( M_1 \otimes_R M_2 \right)
\end{equation*}
given by
\begin{equation*}
	\left( s_1 \otimes_R m_1 \right) \otimes_S \left( s_2 \otimes_R m_2 \right) \mapsto
	(-1)^{\braidd{m_1}{s_2}} \left( s_1 \cdot s_2 \right) \otimes_R \left( m_1 \otimes_R m_2 \right),
\end{equation*}
and the unit constraint $S \rightarrow \varphi_{!} \left( R \right) = S \otimes_R R \cong S$.
The right adjoint functor $\varphi^{*}$ is lax symmetric monoidal via the tensor constraints
\begin{equation*}
	\varphi^{*} \left( N_1 \right) \otimes_R \varphi^{*} \left( N_2 \right) \rightarrow
	\varphi^{*} \left( N_1 \otimes_S N_2 \right)
\end{equation*}
given by
\begin{equation}
	n_1 \otimes_R n_2 \mapsto n_1 \otimes_S n_2 \label{eq:pullback-tensor-constraint-graded-modules}
\end{equation}
and the unit constraint $\varphi$, considered as an $R$-linear map
$\varphi \colon R \rightarrow \varphi^{*} \left( S \right)$.

\begin{rem}
	We have defined the action of the functors $\varphi^{*}$ and $\varphi_{!}$ not only on morphisms but also on graded maps
	of arbitrary degree, i.e., elements of $\InnHom{M}{N}[][R]$. Since $\GMod[R]$ is a closed symmetric monoidal category,
	it is enriched over itself. Since the functors $\varphi^*$ and $\varphi_!$ act on the internal hom objects, they
	can be enhanced to \textit{enriched} functors. The adjunction \eqref{eq:restriction-extension-adjunction-gmod-R} then
	becomes an \textit{enriched} adjunction.
\end{rem}

Given a graded $R$-module $M$ and a graded $S$-module $N$, a graded $\mathbbm{k}$-linear
map $f \colon M \rightarrow N$ of degree zero which satisfies
$f \left( r \cdot m \right) = \varphi(r) \cdot f(m)$
for all $r \in R$ and $m \in M$ will be called a \textbf{morphism of graded modules over}
$\varphi$. Equivalently, a morphism of graded modules over $\varphi$ is a morphism
$f \colon M \rightarrow \varphi^{*} \left( N \right)$
of graded $R$-modules between $M$ and the restriction of scalars of $N$ along $\varphi$.
Morphisms of graded modules over algebra morphisms can be composed so that if $g \colon L \rightarrow M$
is a morphism of graded modules over
$\psi \colon Q \rightarrow R$ and $f \colon M \rightarrow N$ is a morphism of graded modules over
$\varphi \colon R \rightarrow S$ then $f \circ g \colon L \rightarrow N$ is a morphism of graded modules
over $\varphi \circ \psi \colon Q \rightarrow S$.

Given two morphisms $f_1 \colon M_1 \rightarrow N_1$ and $f_2 \colon M_2 \rightarrow N_2$
of graded modules over $\varphi$, we will denote by
$f_1 \otimes_{\varphi} f_2 \colon M_1 \otimes_R M_2 \rightarrow N_1 \otimes_S N_2$ the morphism of
graded modules over $\varphi$ given by
\begin{equation}
	\left( f_1 \otimes_{\varphi} f_2 \right) \left( m_1 \otimes_R m_2 \right) \defeq
	f_1 \left( m_1 \right) \otimes_S f_2 \left( m_2 \right). \label{eq:tensor-product-morphisms-over-phi}
\end{equation}
Thinking of $f_1, f_2$ and $f_1 \otimes_{\varphi} f_2$ as $R$-linear morphisms whose codomains are
the appropriate pullbacks, the morphism $f_1 \otimes_{\varphi} f_2$ is given by the composition
\begin{equation*}
	M_1 \otimes_R M_2 \xrightarrow{f_1 \otimes_R f_2}
	\varphi^{*} \left( N_1 \right) \otimes_R \varphi^{*} \left( N_2 \right) \rightarrow
	\varphi^{*} \left( N_1 \otimes_S N_2 \right)
\end{equation*}
where the last map is the canonical map given by \cref{eq:pullback-tensor-constraint-graded-modules}.
Instead of working over a fixed base and composing everything with various canonical maps, we prefer to adopt
the notation above in order to reduce clutter.
When $R = S$ and $\varphi = \id_R$, the definition of $f_1 \otimes_{\varphi} f_2$ reduces
to the standard tensor product $f_1 \otimes_R f_2$ of graded $R$-linear maps.

\subsection{Graded Coalgebras over Graded Algebras} \label{sec:graded-coalgebras-over-graded-algebras}

Let $R$ be a graded-commutative $\mathbbm{k}$-algebra.
Since the category $\GMod[R]$ is monoidal, one can talk about coalgebra objects in $\GMod[R]$
(see \cref{subsec:coalg-in-monoidal-cat}). A graded $R$-coalgebra is a coalgebra object of $\GMod[R]$.
More explicitly, we have:

\begin{dfn}
	A \textbf{graded} $R$\textbf{-coalgebra} is a graded $R$-module $C$ together with a degree zero
	$R$-linear comultiplication $\Delta = \Delta_C \colon C \rightarrow C \otimes_R C$ and a degree zero
	$R$-linear counit map $\varepsilon = \varepsilon_C \colon C \rightarrow R$ such that
	$\left( \Delta \otimes_R \id \right) \circ \Delta = \left( \id \otimes_R \Delta \right) \circ \Delta$
	and $\left( \varepsilon \otimes_R \id \right) \circ \Delta = \left( \id \otimes_R \varepsilon \right) \circ \Delta = \id$.\footnote{
		As usual, we identify both sides using the associator and unitor isomorphisms coming
		from the monoidal structure.}
	A \textbf{morphism} $f \colon C \rightarrow D$ \textbf{of graded} $R$\textbf{-coalgebras} is a degree zero
	$R$-linear map which satisfies $\left( f \otimes_R f \right) \circ \Delta_C = \Delta_D \circ f$ and
	$\varepsilon_D \circ f = \varepsilon_C$.
\end{dfn}

Note that the graded-commutative ground $\mathbbm{k}$-algebra $R$ has a canonical structure
of a graded $R$-coalgebra where $\Delta_R \colon R \rightarrow R \otimes_R R$ is the
canonical isomorphism and $\varepsilon_R \colon R \rightarrow R$ is the identity map.
Given a graded $R$-coalgebra, a morphism $u \colon R \rightarrow C$ of graded $R$-coalgebras
which satisfies $\varepsilon_C \circ u = \id_R$ is called a \textbf{coaugmentation} on $C$.
A graded $R$-coalgebra $C$ together with a choice of a coaugmentation $u \colon R \rightarrow C$ is called
a \textbf{graded coaugmented} $R$\textbf{-coalgebra}.

Let $R$ and $S$ be two graded-commutative $\mathbbm{k}$-algebras and let
$\varphi \colon R \rightarrow S$ be a morphism of graded $\mathbbm{k}$-algebras. Let
$C$ be a graded $R$-coalgebra. Since the extension of scalars functor
$\varphi_{!} \colon \GMod[R] \rightarrow \GMod[S]$ is strong, and in particular oplax, monoidal,
the $S$-module $\varphi_{!} \left( C \right) = S \otimes_R C$
has a natural structure of a graded $S$-coalgebra (see \cref{subsec:coalg-in-monoidal-cat}). The
coproduct on $\varphi_{!} \left( C \right)$ is given by
\begin{equation*}
	\varphi_{!} \left( C \right) = S \otimes_R C
	\xrightarrow{\id_S \otimes_R \Delta_C} S \otimes_R \left( C \otimes_R C \right)
	\xrightarrow{\cong} \left( S \otimes_R C \right) \otimes_S \left( S \otimes_R C \right) =
	\varphi_{!} \left( C \right) \otimes_S \varphi_{!} \left( C \right)
\end{equation*}
and the counit on $\varphi_{!} \left( C \right)$ is given by
\begin{equation*}
	\varphi_{!} \left( C \right) = S \otimes_R C
	\xrightarrow{\id \otimes_R \varepsilon_C} S \otimes_R R \cong S.
\end{equation*}
The resulting coalgebra $\varphi_! \left( C \right)$ is called the \textbf{scalar extension of} $C$ \textbf{along}
$\varphi$.

\phantomsection
\label{sec:coalgebra-morphism-over-different-ground-algebras}

We will also need the notion of a coalgebra morphism between two coalgebras over different ground algebras.
Let $R$ and $S$ be two graded-commutative $\mathbbm{k}$-algebras and let
$\varphi \colon R \rightarrow S$ be a morphism of graded $\mathbbm{k}$-algebras.
Given a graded $R$-coalgebra $(C,\Delta_C,\varepsilon_C)$ and a graded $S$-coalgebra
$(D,\Delta_D,\varepsilon_D)$ a \textbf{coalgebra morphism} $f \colon C \rightarrow D$ \textbf{over an
	algebra morphism} $\varphi \colon R \rightarrow S$ is a morphism of graded modules over $\varphi$
such that the following diagrams commute:

\begin{figure}[H]
	\centering
	\begin{subfigure}{0.45\textwidth}
		\centering
		\begin{tikzcd}
			C & D \\
			{C \otimes_R C} & {D \otimes_S D} \\
			\arrow["f", from=1-1, to=1-2]
			\arrow["{\Delta_C}"', from=1-1, to=2-1]
			\arrow["{f \otimes_{\varphi} f}", from=2-1, to=2-2]
			\arrow["{\Delta_D}", from=1-2, to=2-2]
		\end{tikzcd}
		\caption{Compatibility with the coproducts.}
		\label{fig:morphisms-coalgebra-different-ground-algebra-succinct-coproduct}
	\end{subfigure}
	\begin{subfigure}{0.45\textwidth}
		\centering
		\begin{tikzcd}
			C & D \\
			R & S \\
			\arrow["f", from=1-1, to=1-2]
			\arrow["{\varepsilon_D}", from=1-2, to=2-2]
			\arrow["{\varepsilon_C}"', from=1-1, to=2-1]
			\arrow["\varphi", from=2-1, to=2-2]
		\end{tikzcd}
		\caption{Compatibility with the counits.}
		\label{fig:morphisms-coalgebra-different-ground-algebra-succinct-counit}
	\end{subfigure}
	\caption{Morphism between two graded coalgebras over different ground algebras.}
	\label{fig:morphisms-coalgebra-different-ground-algebra-succinct}
\end{figure}

When $R = S$ and $\varphi = \id_R$, we recover the usual notion of a coalgebra morphism.

\begin{rem} \label{rem:coalgebra-morphism-over-algebra-morphism-adjunction}
	We can use restriction of scalars and, thinking of a morphism $f \colon C \rightarrow D$
	over $\varphi$ as a morphism $f \colon C \rightarrow \varphi^{*} \left( D \right)$
	of graded $R$-modules, rewrite
	\cref{fig:morphisms-coalgebra-different-ground-algebra-succinct-coproduct}
	so that all the objects involved are graded $R$-modules and all maps are $R$-linear. This
	gives us \cref{fig:morphisms-coalgebra-different-ground-rings-restriction},
	in which the map $\eta$ is the natural map
	$d_1 \otimes_R d_2 \mapsto d_1 \otimes_S d_2$ coming from the tensor constraints of $\varphi^{*}$
	(see \cref{eq:pullback-tensor-constraint-graded-modules}).

	Note that \cref{fig:morphisms-coalgebra-different-ground-rings-restriction} is
	not a diagram of graded $R$-coalgebras. Since the functor $\varphi^{*}$ is only lax monoidal,
	the pullback $\varphi^{*} \left( D \right)$ has no natural structure of a graded $R$-coalgebra,
	and we can't think of a coalgebra morphism along $\varphi$ as a morphism of graded $R$-coalgebras
	$f \colon C \rightarrow \varphi^{*} \left( D \right)$.
	Instead, one can use the adjunction \eqref{eq:restriction-extension-adjunction-gmod-R} and
	take the adjunct $\tilde{f} \colon \varphi_{!} \left( C \right) \rightarrow D$
	of $f \colon C \rightarrow \varphi^{*} \left( D \right)$.
	The scalar extension $\varphi_{!} \left( C \right)$	has a natural structure of a graded
	$S$-coalgebra and $f$ is a morphism of coalgebras over $\varphi$
	if and only if the adjunct
	$\tilde{f} \colon \varphi_{!} \left( C \right) \rightarrow D$ is a morphism
	of graded $S$-coalgebras.
\end{rem}

\begin{figure}[htb]
	\begin{tikzcd}
		C && {\varphi^{*} \left( D \right)} \\
		{C \otimes_R C} && {\varphi^{*} \left( D \otimes_S D \right)} \\
		& {\varphi^{*} \left( D \right) \otimes_R \varphi^{*} \left( D \right)}
		\arrow["f", from=1-1, to=1-3]
		\arrow["{\Delta_C}"', from=1-1, to=2-1]
		\arrow["{f \otimes_R f}"', from=2-1, to=3-2]
		\arrow["{\varphi^{*} \left( \Delta_D \right)}", from=1-3, to=2-3]
		\arrow["\eta"', from=3-2, to=2-3]
	\end{tikzcd}
	\caption{Morphism of graded coalgebras over different ground algebras using restriction of scalars.}
	\label{fig:morphisms-coalgebra-different-ground-rings-restriction}
\end{figure}

Clearly coalgebra morphisms over algebra morphisms can be composed so that if $g \colon B \rightarrow C$
is a coalgebra morphism over
$\psi \colon Q \rightarrow R$ and $f \colon C \rightarrow D$ is a coalgebra morphism over
$\varphi \colon R \rightarrow S$ then $f \circ g \colon B \rightarrow D$ is a coalgebra morphism over
$\varphi \circ \psi \colon Q \rightarrow S$.

\subsection{Pre-Differential Graded Algebras and Modules} \label{subsec:pre-differential-graded-algebras-modules}
Starting with this subsection and for the remainder of the section,
we fix an element $\go \in \GG$
which is considered as part of the grading datum (see \cref{subsec:grading-data}).
Pre-differentials on objects will raise the degree of elements by $\go \in \GG$.

\phantomsection
\label{sec:pre-differential-graded-k-modules}

\begin{dfn}
	A \textbf{pre-differential graded} $\mathbbm{k}$\textbf{-module} is a pair $\mathcal{M} = \left( M, d \right)$
	where $M$ is a graded $\mathbbm{k}$-module, called the \textbf{underlying graded} $\mathbbm{k}$-\textbf{module}
	and $d \colon M \rightharpoonup M$ is a graded map of degree
	$\go$ called the \textbf{pre-differential}.
\end{dfn}

A \textbf{morphism} of pre-differential graded modules
$f \colon \left( M, d_M \right) \rightarrow \left( N, d_N \right)$ is a morphism $f \colon M \rightarrow N$
of the underlying graded $\mathbbm{k}$-modules compatible with the pre-differentials in the sense
that $f \circ d_M = d_N \circ f$.
Let us denote by $\PDGMod[\mathbbm{k}]$ the category of pre-differential graded $\mathbbm{k}$-modules
with morphisms as defined above. The category $\PDGMod[\mathbbm{k}]$ is bicomplete where
both limits and colimits are given by limits and colimits
of the underlying graded $\mathbbm{k}$-modules endowed with the natural pre-differentials.

Let $\mathcal{M} = \left( M, d_M \right)$ and $\mathcal{N} = \left( N, d_N \right)$
be two pre-differential graded $\mathbbm{k}$-modules. A \textbf{graded map} $f \colon
	\mathcal{M} \rightharpoonup \mathcal{N}$ between pre-differential graded $\mathbbm{k}$-modules
is by definition a graded map $f \colon M \rightharpoonup N$ of the underlying graded $\mathbbm{k}$-modules.
We use the notation $f \colon \mathcal{M} \rightharpoonup \mathcal{N}$ to emphasize
that $f$ is a graded map of arbitrary degree which is not required to be compatible with the pre-differentials
in any way while reserving the notation $f \colon \mathcal{M} \rightarrow \mathcal{N}$ for morphisms
(i.e., graded maps of degree zero compatible with the pre-differentials).
The \textbf{internal hom object}
$\InnHom{\mathcal{M}}{\mathcal{N}}[][\mathbbm{k}]$ is defined to be the graded $\mathbbm{k}$-module
$\InnHom{M}{N}[][\mathbbm{k}]$ endowed with the pre-differential $\partial$ given by
\begin{equation}
	\partial \left( f \right) \defeq d_N \circ f - (-1)^{\braid{\go}{\degb{f}}} f \circ d_M.
	\label{eq:hom-differential}
\end{equation}
The pre-differential $\partial$ satisfies a graded Leibniz rule with respect to composition, i.e.,
\begin{equation}
	\partial \left( f \circ g \right) = \partial f \circ g +
		(-1)^{\braid{\degb{\partial}}{\degb{f}}} f \circ \partial g.
	\label{eq:interaction-differential-composition}
\end{equation}
for	graded maps $f \colon \mathcal{M} \rightharpoonup \mathcal{N}$ and
$g \colon \mathcal{L} \rightharpoonup \mathcal{M}$.
A graded map $f \colon \mathcal{M} \rightharpoonup \mathcal{N}$ is called \textbf{closed}
or a \textbf{chain map} if $\partial \left( f \right) = 0$. In particular, a morphism
$f \colon \mathcal{M} \rightarrow \mathcal{N}$ is a degree zero chain map. As a direct consequence
of \cref{eq:interaction-differential-composition}, the composition of closed maps is also closed.

The \textbf{tensor product} $\mathcal{M} \otimes_{\mathbbm{k}} \mathcal{N}$ of two
pre-differential graded $\mathbbm{k}$-modules is defined to be the
graded $\mathbbm{k}$-module $M \otimes_{\mathbbm{k}} N$ endowed with the pre-differential
\begin{equation}
	d_{M \otimes_{\mathbbm{k}} N} \defeq d_M \otimes_{\mathbbm{k}} \id + \id \otimes_{\mathbbm{k}} d_N.
	\label{eq:k-tensor-product-differential}
\end{equation}
The pre-differential $\partial$ interacts with the tensor product of graded maps via a graded Leibniz formula, i.e.,
\begin{equation}
	\partial \left( f \otimes_{\mathbbm{k}} g \right) = \partial f \otimes_{\mathbbm{k}} g +
		(-1)^{\braid{\degb{\partial}}{\degb{f}}} f \otimes_{\mathbbm{k}} \partial g
	\label{eq:interaction-differential-tensor-product-of-maps}
\end{equation}
for graded maps $f \colon M_1 \rightharpoonup N_1$ and
$g \colon M_2 \rightharpoonup N_2$. In particular,
the tensor product of two chain maps is a chain map.

With the definitions above, the closed symmetric monoidal structure on graded $\mathbbm{k}$-modules extends
to pre-differential graded $\mathbbm{k}$-modules and the category $\PDGMod[\mathbbm{k}]$ becomes closed symmetric monoidal.

Given a pre-differential graded $\mathbbm{k}$-module $\mathcal{M} = \left( M, d \right)$
and $h \in \GG$, the $h$-\textbf{suspension} or $h$-\textbf{shifted} module $\mathcal{M}[h]$ is
defined to be the graded $\mathbbm{k}$-module $M[h]$ (see \cref{sec:suspension-graded-k-module})
endowed with the pre-differential
$d_{\mathcal{M}[h]}$ given by
\begin{equation}
	d_{\mathcal{M}[h]} \left( \s_h \left( m \right) \right) \defeq
	(-1)^{\braid{\go}{-h}} \s_h \left( d \left( m \right) \right) =
	(-1)^{\braid{\go}{h}} \s_h \left( d \left( m \right) \right). \label{eq:differential-on-suspension}
\end{equation}
With the definition given by \cref{eq:differential-on-suspension}, the suspension map
$\s_h \colon \mathcal{M} \rightharpoonup \mathcal{M}[h]$ becomes a \textbf{chain map} of degree $-h$.

\begin{dfn} \label{def-graded-algebra-derivation}
	Let $R$ be a graded $\mathbbm{k}$-algebra. A \textbf{graded algebra derivation} is a graded map
	$d \colon R \rightharpoonup R$ which satisfies the graded Leibniz rule
	\begin{equation*}
		d \left( r \cdot s \right) = d \left( r \right) \cdot s + (-1)^{\braidd{d}{r}} r \cdot d \left( s \right).
	\end{equation*}
\end{dfn}

Given a graded $\mathbbm{k}$-algebra $R$, we denote by $\Der{R}$ the graded $\mathbbm{k}$-module
of all graded algebra derivations on $R$. The $\mathbbm{k}$-module $\Der{R}$ is a graded Lie subalgebra of
$\InnEnd{R}[][\mathbbm{k}]$ where the
Lie bracket is given by the graded commutator
\begin{equation*}
	\left[ d_1, d_2 \right] = d_1 \circ d_2 - (-1)^{\braidd{d_1}{d_2}} d_2 \circ d_1.
\end{equation*}
When $R$ is graded-commutative, $\Der{R}$ has a natural outer $R$-action given by \cref{eq:left-R-action-on-hom}
which turns it into a graded $R$-module.

\phantomsection
\label{sec:pre-differential-graded-k-algebras}

\begin{dfn}
	A \textbf{pre-differential graded} $\mathbbm{k}$\textbf{-algebra} is a pair $\mathcal{R} = \left( R, d \right)$
	where $R$ is a graded $\mathbbm{k}$-algebra, called the \textbf{underlying graded}
	$\mathbbm{k}$-\textbf{algebra}, and $d \colon R \rightharpoonup R$ is a graded algebra derivation of
	degree $\go$.
\end{dfn}

We note that a pre-differential graded $\mathbbm{k}$-algebra is the same thing as
an algebra object $\mathcal{R}$ of the monoidal category $\PDGMod[\mathbbm{k}]$.
From this point of view, the graded Leibniz rule is a consequence of the fact that the algebra multiplication
$\mathcal{R} \otimes_{\mathbbm{k}} \mathcal{R} \rightarrow \mathcal{R}$ is a morphism of
pre-differential graded $\mathbbm{k}$-modules, hence compatible with the differentials on both sides.
A \textbf{morphism of pre-differential graded} $\mathbbm{k}$\textbf{-algebras}
$f \colon (R,d_R) \rightarrow (S,d_S)$ is a $\mathbbm{k}$-algebra morphism $f \colon R \rightarrow S$ compatible
with the pre-differentials in the sense that $f \circ d_R = d_S \circ f$.

\begin{ex}
	Given a pre-differential graded $\mathbbm{k}$-module $\mathcal{M}$, the formula in
	\cref{eq:interaction-differential-composition} shows that $\partial$ is a graded algebra
	derivation on $\InnEnd{M}[][\mathbbm{k}]$ with composition as multiplication and so the pair
	$\left( \InnEnd{M}[][\mathbbm{k}], \partial \right)$ is a pre-differential graded
	$\mathbbm{k}$-algebra.
\end{ex}

\begin{dfn} \label{dfn:d-operator}
	Let $R$ be a graded $\mathbbm{k}$-algebra and $M$ a graded left $R$-module.
	Given a graded algebra derivation $d \colon R \rightharpoonup R$, a \textbf{(left module) derivation on} $M$
	\textbf{over} $d$,
	or, more succinctly, a \textbf{left} $d$\textbf{-operator}, is a graded map
	$\mu \colon M \rightharpoonup M$ of the same degree as $d$ which satisfies a
	graded Leibniz rule with respect to $d$, i.e.,
	\begin{equation}
		\mu \left( r \cdot m \right) = dr \cdot m + (-1)^{\braidd{\mu}{r}} r \cdot \mu(m)
		\label{eq:d-operator-Leibniz-rule}
	\end{equation}
	for all $r \in R$ and $m \in M$. Similarly, when $M$ is a graded right $R$-module, a graded map
	$\mu \colon M \rightharpoonup M$ is called a \textbf{right} $d$\textbf{-operator} if it satisfies
	\begin{equation}
		\mu \left( m \cdot r \right) = \mu \left( m \right) \cdot r + (-1)^{\braidd{\mu}{m}} m \cdot dr
		\label{eq:d-operator-Leibniz-rule-right}
	\end{equation}
	for all $r \in R$ and $m \in M$.
\end{dfn}

We note that the notion of $d$-operator generalizes graded $R$-linear maps $\mu \colon M \rightharpoonup M$.
By considering $d = 0$ as a graded algebra derivation on $R$ having an arbitrary fixed degree $g \in \GG$, we see
that a $d$-operator $\mu \colon M \rightharpoonup M$ in this case is simply a graded $R$-linear map of degree $g$.
Note also that a graded map $\mu \colon M \rightharpoonup M$ might be a derivation over two different graded
algebra derivations $d_1,d_2 \colon R \rightharpoonup R$ so the ``underlying'' derivation cannot be recovered from
$\mu$.\footnote{When $M$ is a graded-faithful $R$-module, the underlying algebra derivation of a module derivation is unique.}
Let us state explicitly some useful elementary properties of $d$-operators:
\begin{enumerate}
	\item Given a $d$-operator $\mu$ and a graded $R$-linear map $f \colon M \rightharpoonup M$ with
	      $\degb{f} = \degb{d}$, the map $\mu + f$ is also a $d$-operator. The collection
	      of all $d$-operators over a fixed derivation $d$ is an affine space modelled on
	      $\InnEnd{M}[\degb{d}][R]$.
	\item Given a $d_1$-operator $\mu_1$ and a $d_2$-operator $\mu_2$, the sum $\mu_1 + \mu_2$ is
	      a $(d_1 + d_2)$-operator and the graded commutator $[\mu_1, \mu_2]$ is a $[d_1,d_2]$-operator.
	\item Given a $d$-operator $\mu$ and $\lambda \in \mathbbm{k}$, the map $\lambda \mu$ is
	      a $\left( \lambda d \right)$-operator.
	\item When $R$ is graded-commutative, $r \in R$ and $\mu$ is a $d$-operator then the map
	      $r \mu$ is a $\left( rd \right)$-operator.
\end{enumerate}
Let us denote by $\Der{M}$ the graded $\mathbbm{k}$-module whose graded components are given by
\begin{align*}
	\Der{M}[g] & \defeq \Set{\mu \in \InnEnd{M}[g][\mathbbm{k}]}[\exists d \in {\Der{R}[g]} \textrm{ such that }
		                    \mu \textrm{ is a module derivation over } d].
\end{align*}
It follows immediately that $\Der{M}$ is a graded Lie subalgebra of $\InnEnd{M}[][\mathbbm{k}]$ with respect
to the graded commutator which contains $\InnEnd{M}[][R]$. When $R$ is graded-commutative, $\Der{M}$
has a natural outer $R$-action given by \cref{eq:left-R-action-on-hom} which turns it into a graded $R$-module.

\begin{rem}
	The terminology in \cref{dfn:d-operator} is not consistent across the literature.
	We decided to adopt the terminology from \cite{SuarezAlvarez2017} where the notion appears in the
	ungraded case. In \cite[][Section 1.2.3]{DeParis2009}, $d$-operators are called ``der-operators''.
	In the special case where $R = \df{X}$ is the algebra of differential forms on a smooth manifold $X$
	and $M = \df{X}[E]$ is the module of vector-bundle valued differential forms,
	the notion appears in \cite[][Section 3.1]{Michor1987} simply under the name \textit{derivation}, where
	a full classification of such derivations is obtained.
\end{rem}

\phantomsection
\label{sec:pre-differential-graded-R-modules}

Let us fix a pre-differential graded $\mathbbm{k}$-algebra $\mathcal{R} = (R,d_R)$.
\begin{dfn}
	A \textbf{pre-differential graded left} (resp.\ \textbf{right}) $\mathcal{R}$\textbf{-module}
	is a pair $\mathcal{M} = \left( M, d_M \right)$ where $M$ is a graded left (resp.\ right) $R$-module,
	called the \textbf{underlying graded} $R$-\textbf{module}, and
	$d_M \colon M \rightharpoonup M$ is a left (resp.\ right)
	$d_R$-operator on $M$ called the \textbf{pre-differential}.
\end{dfn}

Similarly, one can define pre-differential graded bimodules. In what follows, unless explicitly stated otherwise,
the term ``module'' will always mean a left module. To lessen the burden of notation, when no confusion is possible,
we call a pre-differential graded $\mathcal{R}$-module $\mathcal{M}$ simply an $\mathcal{R}$-module, dropping
the adjectives ``pre-differential'' and ``graded'' and relying on context and the calligraphic font to remind
us that the module is graded and equipped with a pre-differential.

We note that a pre-differential graded $\mathcal{R}$-module is the same thing as
a module object $\mathcal{M}$ over an algebra object $\mathcal{R}$ in the monoidal category $\PDGMod[\mathbbm{k}]$.
From this point of view, the graded Leibniz rule of \cref{eq:d-operator-Leibniz-rule}
is a consequence of the fact that the module action $\mathcal{R} \otimes_{\mathbbm{k}} \mathcal{M} \rightarrow \mathcal{M}$ is a morphism of pre-differential graded $\mathbbm{k}$-modules, hence compatible with the
pre-differentials on both sides.
Note also that a pre-differential graded $\mathcal{R}$-module $\mathcal{M} = \left( M, d_M \right)$
is a pre-differential graded $\mathbbm{k}$-module, together with an $R$-action on $M$ such that the
pre-differential $d_M$ is compatible with the $R$-action and $d_R$ via \cref{eq:d-operator-Leibniz-rule}.
By forgetting the $R$-action and the compatibility condition of $d_M$, we obtain the
\textbf{underlying pre-differential graded} $\mathbbm{k}$-\textbf{module} of $\mathcal{M}$. Hence, we
can think of an $\mathcal{R}$-module $\mathcal{M}$ as a pre-differential graded $\mathbbm{k}$-module,
and we do so when necessary without further mention.

A \textbf{morphism of pre-differential graded} $\mathcal{R}$\textbf{-modules}
$f \colon (M,d_M) \rightarrow (N, d_N)$ is a
morphism $f \colon M \rightarrow N$ of graded $R$-modules (i.e., a degree zero $R$-linear map)
compatible with the pre-differentials in the sense that
$f \circ d_M = d_N \circ f$.\footnote{In other words, a morphism $f \colon \mathcal{M} \rightarrow \mathcal{N}$
	of pre-differential graded $\mathcal{R}$-modules is both a morphism of the underlying
	graded $R$-modules and the underlying pre-differential graded $\mathbbm{k}$-modules.}
Let us denote by $\PDGMod[\mathcal{R}]$ the category of pre-differential graded
$\mathcal{R}$-modules. The category $\PDGMod[\mathcal{R}]$
is bicomplete where both limits and colimits are given by limits and colimits of the underlying graded $R$-modules endowed with natural pre-differentials which become $d_R$-operators.

Let $\mathcal{M} = \left( M, d_M \right)$ and $\mathcal{N} = \left( N, d_N \right)$
be two pre-differential graded $\mathbbm{k}$-modules.
The pre-differential graded $\mathbbm{k}$-module
$\InnHom{\mathcal{M}}{\mathcal{N}}[][\mathbbm{k}]$ has left and right $R$-actions given by
\cref{eq:left-R-action-on-hom} and \cref{eq:right-R-action-on-hom}. The pre-differential
$\partial \colon \InnHom{M}{N}[][\mathbbm{k}] \rightharpoonup \InnHom{M}{N}[][\mathbbm{k}]$
given by \cref{eq:hom-differential} is actually a $d_R$-operator with respect to both actions and hence
$\InnHom{\mathcal{M}}{\mathcal{N}}[][\mathbbm{k}]$ is actually an $\mathcal{R}$-bimodule. This implies that
the pre-differential $\partial$ descends to the graded $\mathbbm{k}$-module $\InnHom{M}{N}[][R]$
of $R$-linear maps (i.e., the pre-differential of a graded $R$-linear map is $R$-linear).
We denote the resulting pre-differential graded $\mathbbm{k}$-module
by $\InnHom{\mathcal{M}}{\mathcal{N}}[][\mathcal{R}]$.

A \textbf{graded} $\mathcal{R}$-\textbf{linear} map $f \colon \mathcal{M} \rightharpoonup \mathcal{N}$
between $\mathcal{R}$-modules is by definition a graded $R$-linear map $f \colon M \rightharpoonup N$ of the underlying graded $R$-modules, i.e., an element of $\InnHom{\mathcal{M}}{\mathcal{N}}[][\mathcal{R}]$.
We use the notation $f \colon \mathcal{M} \rightharpoonup \mathcal{N}$ to emphasize
that $f$ is a graded $R$-linear map of arbitrary degree which is not required to be compatible with
the pre-differentials in any way while reserving the notation $f \colon \mathcal{M} \rightarrow \mathcal{N}$
for morphisms (i.e., graded $R$-linear maps of degree zero compatible with the pre-differentials).

Given a right $\mathcal{R}$-module $\mathcal{M} = \left( M, d_M \right)$ and a
left $\mathcal{R}$-module $\mathcal{N} = \left( N, d_N \right)$, the
tensor product $\mathcal{M} \otimes_{\mathcal{R}} \mathcal{N}$ is a pre-differential graded $\mathbbm{k}$-module
whose underlying graded $\mathbbm{k}$-module is $M \otimes_R N$, endowed with the pre-differential given by
\begin{equation}
	d_{M \otimes_R N} \defeq d_M \otimes_R \id + \id \otimes_R d_N. \label{eq:R-tensor-product-differential}
\end{equation}

\begin{rem} \label{rem:tensor-product-d-operators-abuse}
	Note that there is some slight abuse of notation involved in \cref{eq:R-tensor-product-differential}.
	Even though neither $d_M$ nor $d_N$ are $R$-linear, the graded Leibniz rules
	(\cref{eq:d-operator-Leibniz-rule,eq:d-operator-Leibniz-rule-right})
	guarantee that the map
	\begin{equation*}
		\left( m, n \right) \mapsto
		d_M \left( m \right) \otimes_R n + (-1)^{\braidd{d_N}{m}} m \otimes_R d_N \left( n \right)
	\end{equation*}
	is a graded $R$-balanced map (in the sense of \cref{eq:R-balanced-map}),
	hence induces a unique graded map $M \otimes_R N \rightharpoonup M \otimes_R N$
	acting on elementary tensors by
	\begin{equation*}
		m \otimes_R n \mapsto
		d_M \left( m \right) \otimes_R n + (-1)^{\braidd{d_N}{m}} m \otimes_R d_N \left( n \right)
	\end{equation*}
	which we denote by $d_M \otimes_R \id + \id \otimes_R d_N$.
\end{rem}

In what follows, we will assume that $\mathcal{R} = \left( R, d_R \right)$ is
a pre-differential graded-commutative $\mathbbm{k}$-algebra. In this case, a left $\mathcal{R}$-module
$\mathcal{M} = \left( M, d_M \right)$ has a natural structure of a
right $\mathcal{R}$-module where the right $R$-action on $M$ is given by
\cref{eq:convert-left-to-right-module} and the pre-differential is the same.
The pre-differential $d_M$ becomes a right $d_R$-operator so that $\mathcal{M}$ indeed becomes a right $\mathcal{R}$-module. Endowing $\mathcal{M}$ with both structures, $\mathcal{M}$ becomes a \textbf{symmetric} $\mathcal{R}$-bimodule. As usual in the graded-commutative setting, we will identify left $\mathcal{R}$-modules, right $\mathcal{R}$-modules
and symmetric $\mathcal{R}$-bimodules as necessary.

In particular, given two (left) $\mathcal{R}$-modules $\mathcal{M}$ and $\mathcal{N}$, we can think of
$\mathcal{M}$ and $\mathcal{N}$ as symmetric $\mathcal{R}$-bimodules and form the tensor product $\mathcal{M}
	\otimes_{\mathcal{R}} \mathcal{N}$ which has a natural structure of an $\mathcal{R}$-module
with the $R$-action given by \cref{eq:R-action-R-tensor-product}. The pre-differential $d_{M \otimes_R N}$
given by \cref{eq:R-tensor-product-differential} is readily verified to be a $d_R$-operator,
hence $\mathcal{M} \otimes_{\mathcal{R}} \mathcal{N}$ is not only a pre-differential graded
$\mathbbm{k}$-module but also an $\mathcal{R}$-module.

In addition, the graded $\mathbbm{k}$-module $\InnHom{M}{N}[][R]$ becomes a graded $R$-module
with respect to the $R$-action given by \cref{eq:left-R-action-on-hom} and the pre-differential $\partial$ on
$\InnHom{M}{N}[][R]$ becomes a $d_R$-operator. Hence, $\InnHom{\mathcal{M}}{\mathcal{N}}[][\mathcal{R}]$
is not only a pre-differential graded $\mathbbm{k}$-module but also an $\mathcal{R}$-module,
hence an object of $\PDGMod[\mathcal{R}]$.

With the definitions above, the closed symmetric monoidal structure on graded $R$-modules extends
to pre-differential graded $\mathcal{R}$-modules and the category $\PDGMod[\mathcal{R}]$
becomes closed symmetric monoidal.

Given an $\mathcal{R}$-module $\mathcal{M}$ and $h \in \GG$, the
$h$-\textbf{suspension} or $h$-\textbf{shifted} module $\mathcal{M}[h]$ is defined
to be the suspension $M[h]$ of the underlying graded $R$-module (see \cref{sec:suspension-graded-R-module})
endowed with the (pre)-differential given by \cref{eq:differential-on-suspension}. The (pre)-differential
is readily verified to be not only $\mathbbm{k}$-linear but a $d_R$-operator, hence
$\mathcal{M}[h]$ is indeed an $\mathcal{R}$-module.
With the definition given above, the suspension map
$\s_h \colon M \rightharpoonup M[h]$ becomes a graded $R$-linear \textit{chain map} of degree $-h$.

\subsection{Pre-Differential Graded Coalgebras} \label{sec:pre-differential-graded-coalgebras}

\begin{dfn} \label{def:coderivation-graded-coalgebra}
	Let $C$ be a graded $R$-coalgebra. A \textbf{(graded) coderivation} on $C$ is a graded $R$-linear map
	$\mu \colon C \rightharpoonup C$ which satisfies the graded co-Leibniz rule
	\begin{equation}
		\Delta \circ \mu =  \left( \mu \otimes_R \id + \id \otimes_R \mu \right) \circ \Delta.
		\label{eq:coderivation-equation}
	\end{equation}
\end{dfn}

More generally, we can talk about graded coderivations on $C$ over graded algebra derivations of
the ground $\mathbbm{k}$-algebra $R$:

\begin{dfn} \label{def:coderivation-over-d}
	Let $C$ be a graded $R$-coalgebra and let $d \colon R \rightharpoonup R$
	be a graded algebra derivation. A graded $\mathbbm{k}$-linear map
	$\mu \colon C \rightharpoonup C$ of the same degree as $d$ will be called a
	\textbf{(graded) coderivation over} $d$ (or a \textbf{generalized (graded) coderivation})
	if $\mu$ is a graded module derivation over $d$ in the sense of
	\cref{dfn:d-operator} and, in addition, it satisfies the graded co-Leibniz rule \eqref{eq:coderivation-equation}.
\end{dfn}

A few elementary remarks about the definition are in order:
\begin{enumerate}[label=\textbf{(R.\arabic*)},ref=R.\arabic*]
	\item While coderivations $\mu$ on $C$ are $R$-linear, i.e., they satisfy
	      \begin{equation*}
		      \mu \left( r \cdot c \right) = (-1)^{\braidd{\mu}{r}} r \cdot \mu \left( c \right),
	      \end{equation*}
	      coderivations over $d$ satisfy instead the identity
	      \begin{equation*}
		      \mu \left( r \cdot c \right) = dr \cdot c + (-1)^{\braidd{\mu}{r}} r \cdot \mu \left( c \right).
	      \end{equation*}
	      Even though generalized coderivations are not $R$-linear, the coderivation equation
	      \eqref{eq:coderivation-equation} still makes sense for them by
	      \cref{rem:tensor-product-d-operators-abuse}.
	\item Assume that $d = 0$, interpreted as a derivation on $R$ of an arbitrary fixed degree $g \in \GG$.
	      Then a coderivation $\mu \colon C \rightharpoonup C$ over $d$ in the sense of
	      \cref{def:coderivation-over-d} is a coderivation of degree $g$ in the sense of
	      \cref{def:coderivation-graded-coalgebra}.
	\item Given two coderivations $\mu_1,\mu_2$ over the same derivation $d$, their difference
	      $\mu_2 - \mu_1$ is an $R$-linear coderivation. The set of all coderivations over a fixed derivation $d$
	      is an affine space modelled on
	      the $\mathbbm{k}$-module of all $R$-linear coderivations of degree $\degb{d}$.
	\item Given coderivations $\mu_i \colon C \rightharpoonup C$ over $d_i \colon R \rightharpoonup R$ for
	      $i = 1,2$, their graded commutator $\left[ \mu_1, \mu_2 \right]$ is a coderivation over
	      $\left[ d_1, d_2 \right]$. In particular, the graded commutator of a generalized coderivation and an
	      $R$-linear coderivation is an $R$-linear coderivation.
	\item Using the counitality of $C$, one can show that any coderivation
	      $\mu \colon C \rightharpoonup C$ over $d$ satisfies the identity
	      \begin{equation}
		      d \circ \varepsilon = \varepsilon \circ \mu. \label{eq:compatibility-coderivation-counit}
	      \end{equation}
	      When $d = 0$, we have the identity $\varepsilon \circ \mu = 0$ which holds for $R$-linear
	      coderivations and is the dual of the identity $d \left( 1 \right) = 0$ which holds for algebra
	      derivations.
	\item \label{item:coderivation-underlying-derivation}
	      Assume that the coalgebra $C$ is coaugmented with a coaugmentation $u \colon R \rightarrow C$
	      which is a morphism of (counital) coalgebras. Given a coderivation
	      $\mu \colon C \rightharpoonup C$ over $d$, we have
	      $\varepsilon \circ \mu \circ u = d \circ \varepsilon \circ u = d \circ \id = d$. Hence,
	      in the presence of coaugmentation, the
	      coderivation map $\mu$ determines the underlying derivation $d$ uniquely, and we don't need to specify
	      explicitly over which $d$ the coderivation $\mu$ lies.
\end{enumerate}

\begin{dfn}
	Let $\mathcal{R} = \left( R, d \right)$ be a pre-differential graded-commutative $\mathbbm{k}$-algebra.
	A \textbf{pre-differential graded} $\mathcal{R}$\textbf{-coalgebra} is a pair
	$\mathcal{C} = \left( C, \mu \right)$, where $C$ is a graded $R$-coalgebra called
	the \textbf{underlying graded} $R$\textbf{-coalgebra}, and $\mu \colon C \rightharpoonup C$
	is a coderivation over $d$.
\end{dfn}

To lessen the burden of notation, when no confusion is possible,
we call a pre-differential graded $\mathcal{R}$-coalgebra $\mathcal{C}$ simply an $\mathcal{R}$-coalgebra, dropping
the adjectives ``pre-differential'' and ``graded'' and relying on context and the calligraphic font to remind
us that the coalgebra is graded and equipped with a coderivation compatible with the derivation on $R$.

We note that a pre-differential graded $\mathcal{R}$-coalgebra is the same thing as
a coalgebra object $\mathcal{C}$ of the monoidal category $\PDGMod[\mathcal{R}]$.
From this point of view, the graded co-Leibniz rule \eqref{eq:coderivation-equation}
stems from the requirement that the comultiplication
$\Delta \colon \mathcal{C} \rightarrow \mathcal{C} \otimes_{\mathcal{R}} \mathcal{C}$
is a morphism of pre-differential graded $\mathcal{R}$-modules, hence compatible with the
pre-differentials on both sides. Similarly, the identity \eqref{eq:compatibility-coderivation-counit}
stems from the requirement that the counit $\varepsilon \colon \mathcal{C} \rightarrow \mathcal{R}$ is compatible
with the pre-differentials on both sides.

A \textbf{morphism of pre-differential graded} $\mathcal{R}$\textbf{-coalgebras}
$f \colon (C,\mu) \rightarrow (D, \nu)$ is a
morphism $f \colon C \rightarrow D$ of the underlying graded $R$-coalgebras
compatible with the coderivations in the sense that
$f \circ \mu = \nu \circ f$.

\subsection{Differential Graded Algebras and Modules}  \label{subsec:differential-graded-algebras-modules}
For the remainder of this subsection, we fix an \textbf{odd} element $\go \in \GG$
which is considered as part of the grading datum (see \cref{subsec:grading-data}).
Differentials on objects will raise the degree of elements by $\go \in \GG$.

\begin{dfn}
	A \textbf{differential graded} $\mathbbm{k}$\textbf{-module} $\mathcal{M} = \left( M, d_M \right)$
	is a pre-differential graded $\mathbbm{k}$-module such that $d_M^2 = 0$ (i.e., $d_M$ is a \textbf{differential}).
\end{dfn}

Let us denote by $\DGMod[\mathbbm{k}]$ the full subcategory of $\PDGMod[\mathbbm{k}]$ consisting of differential graded
$\mathbbm{k}$-modules. Since a differential graded $\mathbbm{k}$-module is in particular a pre-differential graded
$\mathbbm{k}$-module, one can apply all the constructions described in
\cref{subsec:pre-differential-graded-algebras-modules} to differential graded $\mathbbm{k}$-modules
and the resulting object will also be a differential graded $\mathbbm{k}$-module.
This amounts to verifying that the pre-differential on the resulting object is actually a differential.

For example, for the tensor product of two differential graded $\mathbbm{k}$-modules
$\mathcal{M} = \left( M, d_M \right)$ and $\mathcal{N} = \left( N, d_N \right)$, we have
\begin{equation*}
	\begin{aligned}
		d_{M \otimes_{\mathbbm{k}} N}^2 & =
		\left( d_M \otimes_{\mathbbm{k}} \id + \id \otimes_{\mathbbm{k}} d_N \right) \circ
		\left( d_M \otimes_{\mathbbm{k}} \id + \id \otimes_{\mathbbm{k}} d_N \right) \\
		                                & =
		d_M^2 \otimes_{\mathbbm{k}} \id + d_M \otimes_{\mathbbm{k}} d_N +
			                                                              (-1)^{\braidd{d_M}{d_N}} d_M \otimes_{\mathbbm{k}} d_N + \id \otimes_{\mathbbm{k}} d_N^2
		\\
		                                & = 0
	\end{aligned}
\end{equation*}
by our assumption that $\degb{d_M} = \degb{d_N} = \go$ is odd. Similar verification is possible for all
other constructions such as the internal hom, limits, colimits, suspension, etc.

In particular, this means that the category $\DGMod[\mathbbm{k}]$ is also bicomplete, with limits and colimits
given by the limits and colimits of the underlying pre-differential graded $\mathbbm{k}$-modules,
the result being a differential graded $\mathbbm{k}$-module.
The category $\DGMod[\mathbbm{k}]$ is also closed symmetric
monoidal, the monoidal structure being the one restricted from $\PDGMod[\mathbbm{k}]$.

The \textbf{cohomology} $\cohom{\mathcal{M}}[]$
of a differential graded $\mathbbm{k}$-module $\mathcal{M} = \left( M, d_M \right)$
is a graded $\mathbbm{k}$-module whose components are given by
\begin{equation*}
	\cohom{\mathcal{M}}[g] \defeq \ker \left( d_M \colon M^g \rightarrow M^{g+\go} \right) /
	\Im \left( d_M \colon M^{g-\go} \rightarrow M^g \right)
\end{equation*}
for $g \in \GG$. We also write more succinctly
$\cohom{\mathcal{M}}[] = \ker \left( d_M \right) / \Im \left( d_M \right)$
with the understanding that kernels and cokernels of graded maps are taken component-wise.
A graded map $f \colon \mathcal{M} \rightharpoonup \mathcal{N}$ between two differential graded
$\mathbbm{k}$-modules which is closed (i.e., $\partial \left( f \right) = 0$) induces naturally a graded map
$\cohom{f}[] \colon \cohom{\mathcal{M}}[] \rightharpoonup \cohom{\mathcal{N}}[]$ by the formula
\begin{equation} \label{eq:induced-map-cohomology}
	\cohom{f}[] \left( \eqcl{m} \right) \defeq \eqcl{f \left( m \right)}.
\end{equation}
This gives us a cohomology functor $\cohom{}[] \colon \DGMod[\mathbbm{k}] \rightarrow \GMod[\mathbbm{k}]$ which is lax monoidal
via the coherence isomorphisms
$\cohom{\mathcal{M}}[] \otimes_{\mathbbm{k}} \cohom{\mathcal{N}}[] \rightarrow
	\cohom{\mathcal{M} \otimes_{\mathbbm{k}} \mathcal{N}}[]$ given by
\begin{equation} \label{eq:cohomology-lax-monoidal}
	\eqcl{m} \otimes_{\mathbbm{k}} \eqcl{n} \mapsto \eqcl{m \otimes_{\mathbbm{k}} n}.
\end{equation}

\begin{dfn}
	A \textbf{differential graded} $\mathbbm{k}$\textbf{-algebra} $\mathcal{R} = \left( R, d \right)$
	is a pre-differential graded $\mathbbm{k}$-algebra such that $d^2 = 0$ (i.e., the derivation $d$ is
	a \textbf{differential}).
\end{dfn}

We note that a differential graded $\mathbbm{k}$-algebra is the same thing as
an algebra object $\mathcal{R}$ of the monoidal category $\DGMod[\mathbbm{k}]$.
A \textbf{morphism of differential graded} $\mathbbm{k}$\textbf{-algebras}
$f \colon (R,d_R) \rightarrow (S,d_S)$ is the same thing as a morphism of pre-differential graded $\mathbbm{k}$-algebras,
i.e., a $\mathbbm{k}$-algebra morphism $f \colon R \rightarrow S$ compatible
with the differentials in the sense that $f \circ d_R = d_S \circ f$.
Since the cohomology functor $\cohom{}[] \colon \DGMod[\mathbbm{k}] \rightarrow \GMod[\mathbbm{k}]$ is
lax monoidal, the cohomology $\cohom{\mathcal{R}}[]$ of a differential graded $\mathbbm{k}$-algebra
has a natural structure of a graded $\mathbbm{k}$-algebra.

Let us fix a differential graded $\mathbbm{k}$-algebra $\mathcal{R} = (R,d_R)$.

\begin{dfn}
	A \textbf{differential graded
		left} (resp.\ \textbf{right}) 	$\mathcal{R}$\textbf{-module} $\mathcal{M} = \left( M, d_M \right)$
	is a pre-differential graded left (resp.\ right) $\mathcal{R}$-module such that $d_M^2 = 0$ (i.e.,
	$d_M$ is a \textbf{differential}).
\end{dfn}

Let us denote by $\DGMod[\mathcal{R}]$ the full subcategory of $\PDGMod[\mathcal{R}]$ consisting of differential graded
$\mathcal{R}$-modules. Since a differential graded $\mathcal{R}$-module is in particular a pre-differential graded
$\mathcal{R}$-module, one can apply all the constructions described in
\cref{subsec:pre-differential-graded-algebras-modules} to differential graded $\mathcal{R}$-modules
and verify that the resulting object will also be a differential graded $\mathcal{R}$-module.
In particular, this means that the category $\DGMod[\mathcal{R}]$ is also bicomplete, with limits and colimits
given by limits and colimits of pre-differential graded $\mathcal{R}$-modules, the result being
a differential graded $\mathcal{R}$-module. When $\mathcal{R}$ is a differential \textit{graded-commutative}
$\mathbbm{k}$-algebra, the category $\DGMod[\mathcal{R}]$, just like $\PDGMod[\mathcal{R}]$,
is also closed symmetric monoidal with the monoidal and closed structures being the ones restricted
from $\PDGMod[\mathcal{R}]$.

The \textbf{cohomology} $\cohom{\mathcal{M}}[]$
of a differential graded $\mathcal{R}$-module $\mathcal{M} = \left( M, d_M \right)$
is the cohomology of the underlying differential graded $\mathbbm{k}$-module and has
a natural structure of a graded $\cohom{\mathcal{R}}[]$-module over the
graded $\mathbbm{k}$-algebra $\cohom{\mathcal{R}}[]$.
A graded $\mathcal{R}$-linear map $f \colon \mathcal{M} \rightharpoonup \mathcal{N}$ between two
differential graded $\mathcal{R}$-modules which is closed (i.e., $\partial \left( f \right) = 0$)
induces naturally a graded $\cohom{\mathcal{R}}[]$-linear map
$\cohom{f}[] \colon \cohom{\mathcal{M}}[] \rightharpoonup \cohom{\mathcal{N}}[]$ by the same formula
as in \cref{eq:induced-map-cohomology}.
This gives us a cohomology functor
$\cohom{}[] \colon \DGMod[\mathcal{R}] \rightarrow \GMod[{\cohom{\mathcal{R}}[]}]$
which is lax monoidal via the coherence isomorphisms given by the same formula as in
\cref{eq:cohomology-lax-monoidal}, with $\otimes_{\mathbbm{k}}$ replaced by $\otimes_{R}$.

\section{The Non-Archimedean Graded Setting} \label{sec:non-archimedean-graded-setting}

In this section, we discuss the generalization of the notions of graded and (pre)-differential graded modules
and algebras to the non-Archimedean setting. The definitions are quite natural from a categorical
point of view but since we haven't found much written about the non-Archimedean graded case, we decided
to include the definitions and basic properties here in an organized manner.
For background on the non-Archimedean non-graded case, we refer to \cite{Bosch1984} and \cref{appendix:non-archimedean-groups-rings-modules}.
All seminorms in this work are non-Archimedean; hereafter, we often drop the adjective ``non-Archimedean''.

The section is structured similarly to \cref{sec:prelim}, and we present the basic definitions and properties,
including all the modifications required in the non-Archimedean setting.
To motivate our exposition, recall
that the realm of graded algebras and graded modules over graded algebras can be constructed and analyzed as follows:
\begin{enumerate}
	\item Fixing a commutative ground ring $\mathbbm{k}$, one constructs the category
	      $\Mod[\mathbbm{k}]$ of ungraded $\mathbbm{k}$-modules.
	\item Fixing a grading group $\GG$ and a symmetry, i.e., a grading datum, one constructs
	      the category $\mathbf{\GMod[\mathbbm{k}]}$ of $\GG$-graded $\mathbbm{k}$-modules
	      as the category of graded objects of $\Mod[\mathbbm{k}]$ and endows it with the structure
	      of a symmetric monoidal closed category.
	\item The category of $\GG$-graded $\mathbbm{k}$-algebras is then constructed as the category of
	      algebra objects of $\mathbf{\GMod[\mathbbm{k}]}$. A $\GG$-graded $\mathbbm{k}$-algebra $R$ is
	      graded-commutative if it is a commutative algebra object, with respect to the fixed symmetry.
	\item Finally, the category $\GMod[R]$ of $\GG$-graded modules over a $\GG$-graded $\mathbbm{k}$-algebra
	      $R$ is constructed as the category of module objects over an algebra object $R$ of $\mathbf{\GMod[\mathbbm{k}]}$.
\end{enumerate}

When $\GG = \ZZ$, by replacing $\mathbf{\GMod[\mathbbm{k}]}$ with the larger category $\mathbf{Ch(\mathbbm{k})}$ of
$\ZZ$-graded chain complexes of $\mathbbm{k}$-modules, one similarly obtains the categories of DG-algebras and
DG-modules. We will follow the same recipe, replacing the base category $\Mod[\mathbbm{k}]$ with the
non-Archimedean version $\SNMod[\mathbbm{k}]$ (resp.\ $\BMod[\mathbbm{k}]$) consisting of seminormed
(resp.\ Banach) $\mathbbm{k}$-modules.

In what follows, we fix a grading datum $\left( \GG, \braidop \right)$
which will be used to grade objects and endow them with symmetries (see \cref{subsec:grading-data}).
Fix also a commutative ungraded seminormed (possibly Banach) ground ring $\mathbbm{k}$.
We allow $\mathbbm{k}$ to be trivially normed and suppress the seminorm on $\mathbbm{k}$
from our notation.
Although much of the following will depend on $\GG, \braidop$ and $\mathbbm{k}$,
we will often suppress them from our notation as long as they are fixed.

Recall that in the algebraic setting, a graded $\mathbbm{k}$-module can be defined in two equivalent ways:
externally, as an indexed family of $\mathbbm{k}$-modules or, internally, as a $\mathbbm{k}$-module together
with a direct sum decomposition. Both definitions extend to the seminormed and Banach setting as long
as we use the appropriate notion of internal direct sum for each of the categories. In our work,
we adopt the external point of view. For a discussion of the internal point of view,
see \cref{sub:internal-point-of-view}.

\subsection{Graded Seminormed and Banach \texorpdfstring{$\mathbbm{k}$}{k}-Modules}
\label{sec:graded-seminormed-banach-k-modules}

\begin{dfn} \label{def:graded-seminormed-Banach-k-module}
	A \textbf{graded seminormed} $\mathbbm{k}$-\textbf{module} $M = \left( M^g, \nnorm^g \right)_{g \in \GG}$
	is an indexed family of (non-Archimedean) seminormed
	$\mathbbm{k}$-modules called the \textbf{components} of $M$.
	A \textbf{graded Banach} $\mathbbm{k}$-\textbf{module} is a graded seminormed $\mathbbm{k}$-module $M$
	for which each component $\left( M^g, \nnorm^g \right)$ is a Banach $\mathbbm{k}$-module (i.e., normed and complete).
	When $\mathbbm{k} = \left( \ZZ, \trivnorm \right)$, a graded seminormed (resp.\ Banach) $\mathbbm{k}$-module
	is called a \textbf{graded seminormed} (resp.\ \textbf{Banach}) \textbf{group}.
\end{dfn}

Equivalently, we can think of a graded seminormed $\mathbbm{k}$-module $M$ as a graded $\mathbbm{k}$-module
$M = \left( M^g \right)_{g \in \GG}$, called the \textbf{underlying graded} $\mathbbm{k}$-\textbf{module},
equipped with the extra structure of a \textbf{graded} $\mathbbm{k}$-\textbf{module seminorm},
i.e., a family $\nnorm = \left( \nnorm^g \right)_{g \in \GG}$ of $\mathbbm{k}$-module seminorms where
each $\nnorm^g$ is a seminorm on the component $M^g$. Hence, by forgetting the graded seminorm, we can apply notions on graded $\mathbbm{k}$-modules from \cref{subsec:graded-k-modules} to graded seminormed $\mathbbm{k}$-modules
and will do so without further mention. In what follows, we often suppress the seminorms from our notation when no confusion is possible.

Given an element $m \in M$, by which we always mean a homogeneous element $m \in M^d$ for some $d \in \GG$,
we use the notation $\nnorm[m]$ to denote the seminorm $\nnorm[m]^d$.
Given graded seminormed $\mathbbm{k}$-modules $M$ and $N$, we say that a graded map $f \colon M \rightharpoonup N$
of degree $d \in \GG$ with components $f^g \colon M^g \rightarrow N^{g+d}$ is \textbf{bounded} if each
component $f^g$ is bounded, and, furthermore, the components are \textit{uniformly bounded} in the sense
that the $\sup_{g \in \GG} \, \nnorm[f^g] < \infty$. When $f$ is bounded, we set
\begin{equation}
	\nnorm[f] \defeq \sup_{g \in \GG} \, \nnorm[f^g] < \infty \label{def:norm-of-graded-morphism2}
\end{equation}
and call $\nnorm[f]$ the
\textbf{operator seminorm} of $f$. \footnote{Thinking of a graded map $f$ of degree $d$ as an element of
	$\prod_{g \in \GG}^{\Mod[\mathbbm{k}]} \InnHom{M^g}{N^{g+d}}[][\mathbbm{k}]$, we see that $f$ is bounded
	if it belongs to the bounded product $\prod_{g \in \GG}^{\B} \Hom{M^g}{N^{g+d}}[\textrm{bounded}][\mathbbm{k}]$,
	i.e., the direct product in \textit{the category of seminormed modules} of the inner hom objects.
	The operator seminorm of $f$ coincides with the seminorm of $f$ as an element of the direct product.}
A graded bounded map $f \colon M \rightharpoonup N$ between graded
seminormed $\mathbbm{k}$-modules is \textbf{contractive} if $\nnorm[f] \leq 1$ and an \textbf{isometry}
if $\nnorm[f(m)] = \nnorm[m]$ for all $m \in M$.

A \textbf{morphism} $f \colon M \rightarrow N$ \textbf{of graded seminormed} $\mathbbm{k}$-\textbf{modules}
is a graded contractive map of degree zero. We denote by $\GSNMod[\mathbbm{k}]$ the category of graded seminormed
$\mathbbm{k}$-modules with morphisms as above and by $\GBMod[\mathbbm{k}]$ the full subcategory of
$\GSNMod[\mathbbm{k}]$ whose objects are graded Banach $\mathbbm{k}$-modules. In the graded seminormed
and Banach context, unless otherwise specified, we will assume that graded maps are always bounded.
An \textbf{isomorphism} $f \colon M \rightarrow N$ \textbf{of graded seminormed} (or \textbf{Banach})
$\mathbbm{k}$-\textbf{modules} is an isomorphism of the underlying graded $\mathbbm{k}$-modules which
is an isometry.

Even though the morphisms in our categories are graded contractive maps of degree zero, we will
see that graded bounded maps of arbitrary degree appear naturally as the
``internal hom'' object in the categories. Since we often work both with morphisms and
graded bounded maps of arbitrary degree, to avoid confusion,
we use the notation $f \colon M \rightharpoonup N$ to denote graded bounded
maps of arbitrary degree while reserving the notation $f \colon M \rightarrow N$ for
graded contractive maps of degree zero (the actual morphisms in our categories).

Given two graded seminormed $\mathbbm{k}$-modules $M$ and $N$,
we denote by
\begin{equation*}
	\Hom{M}{N}[][\mathbbm{k}] \defeq
	\Set{f \colon M \rightarrow N}[f \textrm{ is a degree zero map with } {\nnorm[f]} \leq 1 ]
\end{equation*}
the set of morphisms between $M$ and $N$ and by $\InnHom{M}{N}[][\mathbbm{k}]$ the graded seminormed
$\mathbbm{k}$-module of all graded bounded maps with components given by
\begin{equation}
	\begin{aligned}
		\InnHom{M}{N}[][\mathbbm{k}]^d & \defeq \Set{f \colon M \rightharpoonup N}
		                                        [f \textrm{ is a graded bounded map of degree } d],
	\end{aligned} \label{eq:inner-hom-graded-seminormed-k-module}
\end{equation}
endowed with the operator seminorm.

\begin{rem}
	Even though we have chosen to use the same notation as in the
	algebraic case (see \cref{eq:inner-hom-graded-k-module}), whenever working with graded
	seminormed $\mathbbm{k}$-modules, $\Hom{M}{N}[][\mathbbm{k}]$ will always mean the set of graded \textit{contractive} maps of degree zero and $\InnHom{M}{N}[][\mathbbm{k}]$ will always mean the space of graded \textit{bounded} maps.
	Note that in the graded seminormed setting, we have
	$\Hom{M}{N}[][\mathbbm{k}] \neq \InnHom{M}{N}[][\mathbbm{k}]^0$
	since maps in $\InnHom{M}{N}[][\mathbbm{k}]^0$ are bounded while morphisms are taken to be contractive.
\end{rem}

The definition of graded bounded maps and the operator seminorm extends naturally to
graded $\mathbbm{k}$-multilinear maps. Given graded seminormed $\mathbbm{k}$-modules
$M_1, \dots, M_n$ and $N$, a graded $\mathbbm{k}$-multilinear map
$B \colon M_1 \times \dots \times M_n \rightharpoonup N$ of degree $d \in \GG$ with components
\begin{equation*}
	B^{g_1,\dots,g_n} \colon M_1^{g_1} \times \dots \times M_n^{g_n} \rightarrow N^{g_1 + \dots + g_n + d}
\end{equation*}
is called \textbf{bounded} if each component
$B^{g_1,\dots,g_n}$ is bounded, and, furthermore, the components are \textit{uniformly bounded} in the sense that
$\sup_{g_1,\dots,g_n \in \GG} \nnorm[B^{g_1,\dots,g_n}] < \infty$.
When $B$ is bounded, we set
\begin{equation}
	\nnorm[B] \defeq \sup_{g_1,\dots,g_n \in \GG} \nnorm[B^{g_1,\dots,g_n}] < \infty
	\label{eq:norm-of-graded-multilinear-map}
\end{equation}
and call $\nnorm[B]$
the \textbf{operator seminorm} of $B$. A graded bounded $\mathbbm{k}$-multilinear map $B$ is called \textbf{contractive} if $\nnorm[B] \leq 1$.

Unless otherwise specified, a graded $\mathbbm{k}$-multilinear map between graded seminormed
$\mathbbm{k}$-modules will always be assumed to be bounded.
Similar to our convention with graded bounded maps, we will use the notation
$B \colon M_1 \times \dots \times M_n \rightharpoonup N$ to denote graded bounded $\mathbbm{k}$-multilinear
maps of arbitrary degree while reserving the notation $B \colon M_1 \times \dots \times M_n \rightarrow N$ to
denote graded contractive $\mathbbm{k}$-multilinear maps of degree zero.

\phantomsection
\label{sec:suspension-graded-seminormed-Banach-k-module}

Given a graded seminormed $\mathbbm{k}$-module $M$ and $h \in \GG$, the
$h$-\textbf{suspension} or $h$-\textbf{shifted} module $M[h]$ is defined
to be the $h$-suspension of the underlying graded $\mathbbm{k}$-module (see \cref{sec:suspension-graded-k-module})
together with the graded seminorm such that for all $m \in M$ we have
$\nnorm[\s_h \left( m \right)] \defeq \nnorm[m]$.
With the definition above, $M[h]$ becomes a graded seminormed $\mathbbm{k}$-module
and the suspension map $\s_h \colon M \rightharpoonup M[h]$
becomes a graded bounded map of degree $-h$ which is an isometry.
When $M$ is a graded Banach $\mathbbm{k}$-module, the suspension $M[h]$ is also a graded Banach
$\mathbbm{k}$-module.

\subsubsection{The Category \texorpdfstring{$\GSNMod[\mathbbm{k}]$}{GSNMod(k)}} \label{sub:category-GSNMod-k}
Note that since we are in the non-Archimedean setting, given two graded seminormed $\mathbbm{k}$-modules
$M$ and $N$, the hom set $\Hom{M}{N}[][\mathbbm{k}]$ is an abelian group (and even a seminormed
$\mathbbm{k}^{\bullet}(1)$-module) and composition of morphisms is $\ZZ$-bilinear so $\GSNMod[\mathbbm{k}]$
is preadditive.

Since the category $\SNMod[\mathbbm{k}]$ is bicomplete, the category $\GSNMod[\mathbbm{k}]$, being the category of
graded objects of $\SNMod[\mathbbm{k}]$, is also bicomplete and both limits and colimits are computed component-wise in $\SNMod[\mathbbm{k}]$. In particular, just like $\SNMod[\mathbbm{k}]$,
the category $\GSNMod[\mathbbm{k}]$ has finite biproducts, kernels and cokernels and hence is pre-abelian but
not abelian. We note that colimits of graded seminormed $\mathbbm{k}$-modules coincide
with colimits of the underlying graded $\mathbbm{k}$-modules, but this is not true for limits, i.e., the forgetful
functor $\GSNMod[\mathbbm{k}] \rightarrow \GMod[\mathbbm{k}]$ is cocontinuous but not
continuous.

We will use the same notation for limits and colimits of graded seminormed $\mathbbm{k}$-modules as in the
(ungraded) seminormed case. For example, the coproduct
$\oplus_{i \in I} M_i$ of a family $\left( M_i \right)_{i \in I}$ of
graded seminormed $\mathbbm{k}$-modules is the graded seminormed $\mathbbm{k}$-module with components
\begin{equation*}
	\left( \oplus_{i \in I} M_i \right)^g \defeq \oplus_{i \in I}^{\SNMod} M_i^g.
\end{equation*}
More explicitly, each component $\left( \oplus_{i \in I} M_i \right)^g$
is the algebraic direct sum $\oplus_{i \in I}^{\Mod} M_i^g$ endowed with the seminorm given by
\cref{eq:direct-sum-seminorm}.

Given two graded seminormed $\mathbbm{k}$-modules $M$ and $N$,
their (graded) tensor product $M \otimes_{\mathbbm{k}} N$ is the graded seminormed $\mathbbm{k}$-module
whose components are given by
\begin{equation}
	\left( M \otimes_{\mathbbm{k}} N \right)^g \defeq
	\bigoplus_{g_1 + g_2 = g} M^{g_1} \otimes_{\mathbbm{k}} N^{g_2}
	\label{eq:tensor-product-seminormed-k-modules}
\end{equation}
where we use the direct sum and tensor product of seminormed $\mathbbm{k}$-modules.\footnote{This is the natural
	way in which a category of graded objects inherits a monoidal structure from the base category.}
Since the direct sum and tensor product of seminormed $\mathbbm{k}$-modules coincides with the
direct sum and tensor product of the underlying $\mathbbm{k}$-modules,
the tensor product $M \otimes_{\mathbbm{k}} N$ coincides
with the algebraic tensor product of graded $\mathbbm{k}$-modules and the seminorm on each component
of $M \otimes_{\mathbbm{k}} N$ is the natural one
induced by the direct sum and tensor product seminorms of ungraded seminormed $\mathbbm{k}$-modules.

Similar to the graded and seminormed cases, the tensor product of graded seminormed $\mathbbm{k}$-modules can be
characterized by a universal property involving graded bounded $\mathbbm{k}$-bilinear maps.
The tensor product $M \otimes_{\mathbbm{k}} N$ of graded seminormed $\mathbbm{k}$-modules
comes equipped with a canonical graded contractive $\mathbbm{k}$-bilinear map
$\otimes_{\mathbbm{k}} \colon M \times N \rightarrow M \otimes_{\mathbbm{k}} N$ of degree zero
characterized by the following universal property: Given a graded seminormed $\mathbbm{k}$-module $L$ and
a graded bounded $\mathbbm{k}$-bilinear map $B \colon M \times N \rightharpoonup L$ there exists a unique
graded bounded map $\varphi_B \colon M \otimes_{\mathbbm{k}} N \rightharpoonup L$
with $\nnorm[\varphi_B] = \nnorm[B]$ and $\degb{\varphi_B} = \degb{B}$
such that $\varphi_B \left( m \otimes_{\mathbbm{k}} n \right) = B(m,n)$ for all $m \in M$ and $n \in N$ (see
\cref{fig:seminormed-tensor-product-k-modules-universal-property} and compare to
\cref{fig:tensor-product-graded-k-modules-universal-property} and
\cref{fig:seminormed-tensor-product-modules-universal-property}).
\begin{figure}[htb]
	\centering
	\begin{tikzcd}
		{M \times N} && {M \otimes_{\mathbbm{k}} N} \\
		&& L
		\arrow["\otimes_{\mathbbm{k}}", from=1-1, to=1-3]
		\arrow["\substack{{\exists! \, \varphi_B} \\ \textrm{graded bounded map}}",
			dashed, harpoon, from=1-3, to=2-3]
		\arrow["\substack{B \textrm{ graded} \\ \textrm{bounded } \mathbbm{k}\textrm{-bilinear}}"',
			harpoon, from=1-1, to=2-3]
	\end{tikzcd}
	\caption{Universal property of the tensor product of graded seminormed $\mathbbm{k}$-modules.}
	\label{fig:seminormed-tensor-product-k-modules-universal-property}
\end{figure}

Given two graded bounded maps $f \colon M \rightharpoonup M'$ and $g \colon N \rightharpoonup N'$ between
graded seminormed $\mathbbm{k}$-modules, their tensor product $f \otimes_{\mathbbm{k}} g$
(see \cref{eq:tensor-product-graded-maps})
is also bounded with $\nnorm[f \otimes_{\mathbbm{k}} g] \leq \nnorm[f] \cdot \nnorm[g]$.
Restricting our attention to graded contractive maps of degree zero, we obtain
a bifunctor $\otimes_{\mathbbm{k}} \colon \GSNMod[\mathbbm{k}] \times \GSNMod[\mathbbm{k}]
	\rightarrow \GSNMod[\mathbbm{k}]$. The bifunctor $\otimes_{\mathbbm{k}}$, the symmetry maps
$M \otimes_{\mathbbm{k}} N \rightarrow N \otimes_{\mathbbm{k}} M$ given by
\begin{equation}
	m \otimes_{\mathbbm{k}} n \mapsto (-1)^{\braidd{m}{n}} n \otimes_{\mathbbm{k}} m,
	\label{eq:symmetries-graded-seminormed-k-modules}
\end{equation}
and the standard associators and unitors\footnote{Note that the symmetry maps, the standard associators and unitors
	are all isomorphisms in $\GSNMod[\mathbbm{k}]$, i.e., \textit{isometric} isomorphisms of graded $\mathbbm{k}$-modules.},
endow the category $\GSNMod[\mathbbm{k}]$ with the structure of a symmetric monoidal category whose unit is the ground ring $\mathbbm{k}$,
considered as graded seminormed $\mathbbm{k}$-module concentrated in degree zero.
It follows from the universal property of the graded seminormed tensor product
that the symmetric monoidal category $\GSNMod[\mathbbm{k}]$ is closed with the internal hom object
$\InnHom{M}{N}[][\mathbbm{k}]$ given by the graded seminormed $\mathbbm{k}$-module of all graded bounded
maps \eqref{eq:inner-hom-graded-seminormed-k-module}.

\subsubsection{The Category \texorpdfstring{$\GBMod[\mathbbm{k}]$}{GBMod(k)}} \label{sub:category-GBMod-k}
Even though the notion of a graded Banach $\mathbbm{k}$-module
makes sense when $\mathbbm{k}$ is merely seminormed, when working with Banach $\mathbbm{k}$-modules we will always
assume that the ground ring $\mathbbm{k}$ is also Banach.\footnote{This is convenient and does not affect much since
	the categories $\GBMod[\mathbbm{k}]$ and $\GBMod (\widehat{\mathbbm{k}})$ are isomorphic.
	See \cref{sec:banach-modules-over-banach-rings} for the ungraded version.}
Before discussing the properties of $\GBMod[\mathbbm{k}]$, we describe the notion of a (separated) completion
of a graded seminormed $\mathbbm{k}$-module.

Given a graded seminormed $\mathbbm{k}$-module $M$, the (separated, graded) \textbf{completion} of $M$
is the graded Banach $\mathbbm{k}$-module $\widehat{M}$ whose components
are given by $\widehat{M}^g = \widehat{M^g}$ where each $\widehat{M^g}$ is the Banach $\mathbbm{k}$-module
obtained by completing the ungraded seminormed $\mathbbm{k}$-module $M^g$
(see \cref{sec:banach-modules-over-banach-rings}).
Similarly to the ungraded case, the completion $\widehat{M}$ comes equipped with a canonical \textbf{completion morphism},
i.e., a degree zero isometry $\eta_M \colon M \rightarrow \widehat{M}$ whose image is dense\footnote{By that, we mean
	that the image of each component $\eta_M^g \colon M^g \rightarrow \widehat{M}^g = \widehat{M^g}$ is dense in
	the codomain.}
in $\widehat{M}$ which is characterized by the following universal property:
Given a graded Banach $\mathbbm{k}$-module $N$ and a graded bounded map $f \colon M \rightharpoonup N$ there
exists a unique bounded ``extension'' $\tilde{f} \colon \widehat{M} \rightharpoonup N$ (the \textbf{adjunct} of $f$)
with $\degb{f} = \degb{\tilde{f}}$ and $\lVert \tilde{f} \rVert = \nnorm[f]$
such that $\tilde{f} \circ \eta_M = f$ (see \cref{fig:graded-completion-k-modules-adjunction}).
When $M$ is Banach, we have $\widehat{M} = M$ and $\eta_M = \id_M$ (see \cref{par:completion-of-Banach-identity}).

\begin{figure}[htb]
	\centering
	\begin{minipage}[b]{.5\textwidth}
		\centering
		\begin{tikzcd}
			M && {\widehat{M}} \\
			&& N
			\arrow["{\eta_M}", from=1-1, to=1-3]
			\arrow["f"', harpoon, from=1-1, to=2-3]
			\arrow["{\exists ! \tilde{f}}", harpoon, dotted, from=1-3, to=2-3]
		\end{tikzcd}
		\captionof{figure}{Universal property of the graded completion.}
		\label{fig:graded-completion-k-modules-adjunction}
	\end{minipage}%
	\begin{minipage}[b]{.5\textwidth}
		\centering
		\begin{tikzcd}
			M && {\widehat{M}} \\
			N && {\widehat{N}}
			\arrow["{\eta_M}", from=1-1, to=1-3]
			\arrow["f"', harpoon, from=1-1, to=2-1]
			\arrow["{\exists ! \widehat{f}}", harpoon, dotted, from=1-3, to=2-3]
			\arrow["{\eta_N}", from=2-1, to=2-3]
		\end{tikzcd}
		\captionof{figure}{\raggedright Completion of a graded bounded map.}
		\label{fig:graded-completion-k-modules-functoriality}
	\end{minipage}
\end{figure}

Given a graded bounded map $f \colon M \rightharpoonup N$, we will denote by
$\widehat{f} \colon \widehat{M} \rightharpoonup \widehat{N}$ the unique graded bounded map such that
$\widehat{f} \circ \eta_M = \eta_N \circ f$ (see \cref{fig:graded-completion-k-modules-functoriality}).
We have $\degb{f} = \degb{\widehat{f}}, \lVert \widehat{f} \rVert = \nnorm[f]$ and if $f$ is an isometry then so is
$\widehat{f}$. In addition, the
completion is functorial in the sense that
$\widehat{f \circ g} = \widehat{f} \circ \widehat{g}$ for all graded bounded maps $f \colon M \rightharpoonup N$ and
$g \colon L \rightharpoonup M$ and $\widehat{\id_M} = \id_{\widehat{M}}$.

Since the completion preserves the norm of morphisms, we obtain a completion functor
$\wedge \colon \GSNMod[\mathbbm{k}] \rightarrow \GBMod[\mathbbm{k}]$ which is left adjoint to the forgetful functor
$U \colon \GBMod[\mathbbm{k}] \rightarrow \GSNMod[\mathbbm{k}]$ with the map
$\eta_M \colon M \rightarrow \widehat{M}$ being the unit of the adjunction. Hence, we see that
$\GBMod[\mathbbm{k}]$ is a reflective replete full subcategory of $\GSNMod[\mathbbm{k}]$ with
reflector $\wedge$.

Since the category $\BMod[\mathbbm{k}]$ is bicomplete, the category $\GBMod[\mathbbm{k}]$, being
the category of graded objects of $\BMod[\mathbbm{k}]$, is also bicomplete and
both limits and colimits are computed component-wise in $\BMod[\mathbbm{k}]$. We will use the
same notation for limits and colimits of graded Banach $\mathbbm{k}$-modules as in the ungraded case.
For example, the complete direct sum (i.e., coproduct) $\coplus_{i \in I} M_i$ of a family
$\left( M_i \right)_{i \in I}$ of graded Banach $\mathbbm{k}$-modules is the graded Banach $\mathbbm{k}$-module
with components
\begin{equation*}
	\left( \coplus_{i \in I} M_i \right)^g \defeq \coplus_{i \in I} M_i^g.
\end{equation*}
More explicitly, each component $\left( \coplus_{i \in I} M_i \right)^g$
of $\coplus_{i \in I} M_i$ is given by the complete direct sum $\coplus_{i \in I} M_i^g$
with the natural norm (see \cref{eq:completed-direct-sum}).
Note that since $\GBMod[\mathbbm{k}]$ is a reflective replete full subcategory
of $\GSNMod[\mathbbm{k}]$, limits in $\GBMod[\mathbbm{k}]$ are the same as limits in $\GSNMod[\mathbbm{k}]$
(and are automatically Banach) while colimits are obtained by applying the completion functor
to colimits in $\GSNMod[\mathbbm{k}]$. Note also that in general,
neither limits nor colimits of Banach $\mathbbm{k}$-modules coincide with limits and colimits
of the underlying graded $\mathbbm{k}$-modules.

The notion of complete tensor product extends naturally to the graded setting.
Given two graded seminormed $\mathbbm{k}$-modules $M$ and $N$,
their \textbf{graded complete tensor product} $M \cotimes_{\mathbbm{k}} N$ is defined as the
completion $\extrawidehat{M \otimes_{\mathbbm{k}} N}$ of the graded seminormed tensor product
$M \otimes_{\mathbbm{k}} N$. Equivalently, we can define $M \cotimes_{\mathbbm{k}} N$ the same
way as in \cref{eq:tensor-product-seminormed-k-modules},
replacing the direct sum and tensor product with their complete versions as we have
\begin{equation*}
	\extrawidehat{M \otimes_{\mathbbm{k}} N}^g =
	\extrawidehat{\bigoplus_{g_1 + g_2 = g} M^{g_1} \otimes_{\mathbbm{k}} N^{g_2}} \cong
	\cbigoplus_{g_1 + g_2 = g} M^{g_1} \cotimes_{\mathbbm{k}} N^{g_2}.
\end{equation*}
Given $m \in M$ and $n \in N$, we will denote by $m \cotimes_{\mathbbm{k}} n$ the image of
$m \otimes_{\mathbbm{k}} n$ under the completion map
$\eta \colon M \otimes_{\mathbbm{k}} N \rightarrow M \cotimes_{\mathbbm{k}} N$
and call such elements \textbf{elementary tensors}. Note that unlike the algebraic or seminormed case,
the complete tensor product $M \cotimes_{\mathbbm{k}} N$ is not generated as a graded $\mathbbm{k}$-module by
elementary tensors. Instead, the graded $\mathbbm{k}$-module generated by elementary tensors is
dense in $M \cotimes_{\mathbbm{k}} N$. In other words, $M \cotimes_{\mathbbm{k}} N$ is generated by elementary tensors
as a graded \textit{Banach} $\mathbbm{k}$-module.

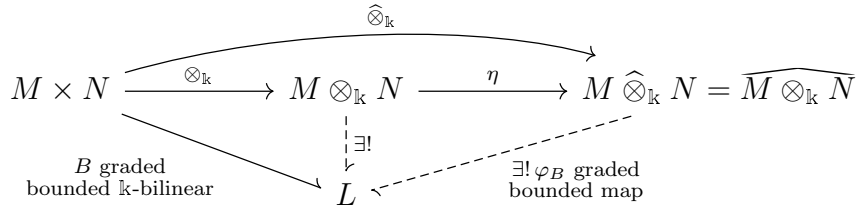
\begin{figure}[htb]
	\centering
	\begin{tikzcd}
		{M \times N} && {M \otimes_{\mathbbm{k}} N} &&
		{M \cotimes_{\mathbbm{k}} N = \extrawidehat{M \otimes_{\mathbbm{k}} N}} \\
		&& L
		\arrow["\cotimes_{\mathbbm{k}}", curve={height=-25pt}, from=1-1, to=1-5]
		\arrow["\otimes_{\mathbbm{k}}", from=1-1, to=1-3]
		\arrow["\eta", from=1-3, to=1-5]
		\arrow["{\exists!}", dashed, harpoon, from=1-3, to=2-3]
		\arrow["\substack{\exists! \, \varphi_B \textrm{ graded} \\ \textrm{bounded map}}",
			dashed, harpoon, from=1-5, to=2-3]
		\arrow["\substack{B \textrm{ graded} \\ \textrm{bounded } \mathbbm{k}\textrm{-bilinear}}"',
			harpoon, from=1-1, to=2-3]
	\end{tikzcd}
	\caption{Universal property of the graded complete tensor product.}
	\label{fig:Banach-tensor-product-k-modules-universal-property}
\end{figure}

Similar to the graded and seminormed cases, the graded complete tensor product
$M \cotimes_{\mathbbm{k}} N$ can be characterized by a universal property involving graded bounded
$\mathbbm{k}$-bilinear maps into graded \textit{Banach} $\mathbbm{k}$-modules.
The complete tensor product $M \cotimes_{\mathbbm{k}} N$ comes equipped with a canonical
graded contractive $\mathbbm{k}$-bilinear map
$\cotimes_{\mathbbm{k}} \colon M \times N \rightarrow M \cotimes_{\mathbbm{k}} N$ of degree zero
characterized by the following universal property: Given a graded Banach $\mathbbm{k}$-module $L$ and
a graded bounded $\mathbbm{k}$-bilinear map $B \colon M \times N \rightharpoonup L$ there exists a unique
graded bounded map $\varphi_B \colon M \cotimes_{\mathbbm{k}} N \rightharpoonup L$ of graded Banach
$\mathbbm{k}$-modules with $\nnorm[\varphi_B] = \nnorm[B]$ and $\degb{\varphi_B} = \degb{B}$
such that $\varphi_B \left( m \cotimes_{\mathbbm{k}} n \right) = B(m,n)$ for all $m \in M$ and $n \in N$ (see
\cref{fig:Banach-tensor-product-k-modules-universal-property} and compare to
\cref{fig:complete-tensor-product-modules-universal-property} and
\cref{fig:tensor-product-graded-k-modules-universal-property}).

Given two graded bounded maps $f \colon M \rightharpoonup M'$ and
$g \colon N \rightharpoonup N'$ between graded seminormed $\mathbbm{k}$-modules, their complete tensor product
$f \cotimes_{\mathbbm{k}} g \colon M \cotimes_{\mathbbm{k}} N \rightharpoonup M' \cotimes_{\mathbbm{k}} N'$ is
defined by $f \cotimes_{\mathbbm{k}} g \defeq \widehat{f \otimes_{\mathbbm{k}} g}$ and
is the unique graded bounded map of degree $\degb{f} + \degb{g}$ which satisfies
\begin{equation*}
	\left( f \cotimes_{\mathbbm{k}} g \right) \left( m \cotimes_{\mathbbm{k}} n \right) =
	(-1)^{\braidd{g}{m}} f \left( m \right) \cotimes_{\mathbbm{k}} g \left( n \right)
\end{equation*}
for all $m \in M$ and $n \in N$. In terms of the bounds on $f,g$ we have
$\nnorm[f \cotimes_{\mathbbm{k}} g] \leq \nnorm[f] \cdot \nnorm[g]$.
The interaction between composition and complete tensor product of graded bounded maps is the same as in
\cref{eq:interaction-composition-tensor-product},
with $\otimes_{\mathbbm{k}}$ replaced by $\cotimes_{\mathbbm{k}}$.
In particular, we see that the graded complete tensor product construction gives us a bifunctor
\begin{equation*}
	\cotimes_{\mathbbm{k}} \colon \GSNMod[\mathbbm{k}] \times \GSNMod[\mathbbm{k}] \rightarrow \GBMod[\mathbbm{k}].
\end{equation*}

Even though we have defined the complete tensor product for graded seminormed $\mathbbm{k}$-modules,
we now restrict our attention to the complete tensor product of graded Banach $\mathbbm{k}$-modules
and obtain a bifunctor
$\cotimes_{\mathbbm{k}} \colon \GBMod[\mathbbm{k}] \times \GBMod[\mathbbm{k}] \rightarrow \GBMod[\mathbbm{k}]$.
Given three graded Banach $\mathbbm{k}$-modules $M,N,L$,
there are natural associativity isomorphisms
$\left( M \cotimes_{\mathbbm{k}} N \right) \cotimes_{\mathbbm{k}} L \cong
	M \cotimes_{\mathbbm{k}} \left( N \cotimes_{\mathbbm{k}} L \right)$
which act on elementary tensors by the expected formula
$\left( m \cotimes_{\mathbbm{k}} n \right) \cotimes_{\mathbbm{k}} l \mapsto
	m \cotimes_{\mathbbm{k}} \left( n \cotimes_{\mathbbm{k}} l \right)$.
In addition, we also have natural isometric isomorphisms
\begin{align*}
	       & \mathbbm{k} \cotimes_{\mathbbm{k}} M \cong M
	\qquad & \lambda \cotimes_{\mathbbm{k}} m                             & \mapsto \lambda \cdot m,                                 \\
	       & M \cotimes_{\mathbbm{k}} \mathbbm{k} \cong M
	\qquad & m \cotimes_{\mathbbm{k}} \lambda                             & \mapsto \lambda \cdot m,                                 \\
	       & M \cotimes_{\mathbbm{k}} N \cong N \cotimes_{\mathbbm{k}} M,
	\qquad & m \cotimes_{\mathbbm{k}} n                                   & \mapsto (-1)^{\braidd{m}{n}} n \cotimes_{\mathbbm{k}} m.
\end{align*}
The bifunctor
$\cotimes_{\mathbbm{k}} \colon \GBMod[\mathbbm{k}] \times \GBMod[\mathbbm{k}] \rightarrow \GBMod[\mathbbm{k}]$, the associativity isomorphisms,
the unitors and the symmetry maps, described above, endow the category $\GBMod[\mathbbm{k}]$
with the structure of a symmetric monoidal category whose unit is the ground ring $\mathbbm{k}$,
considered as a graded Banach $\mathbbm{k}$-module concentrated in degree zero.\footnote{Recall
	that we assume that $\mathbbm{k}$ is Banach so the unit of $\GBMod[\mathbbm{k}]$ is the
	same as the unit of $\GSNMod[\mathbbm{k}]$.}

Given two graded seminormed $\mathbbm{k}$-modules $M$ and $N$, the graded seminormed $\mathbbm{k}$-module
$\InnHom{M}{N}[][\mathbbm{k}]$ of all graded bounded maps of arbitrary degree is Banach once
$N$ is Banach. In particular, $\InnHom{M}{N}[][\mathbbm{k}]$ is an object of $\GBMod[\mathbbm{k}]$ when
both $M$ and $N$ are objects of $\GBMod[\mathbbm{k}]$.
It follows from the universal property of the graded complete tensor product
that the symmetric monoidal category $\GBMod[\mathbbm{k}]$ is closed with internal hom object given by
$\InnHom{M}{N}[][\mathbbm{k}]$ (the same internal hom object as in $\SNMod[\mathbbm{k}]$).

With respect to the monoidal structures defined on $\GSNMod[\mathbbm{k}]$ and $\GBMod[\mathbbm{k}]$, the adjunction
\begin{equation}
	\wedge \colon \GSNMod[\mathbbm{k}] \stackrel[]{\dashv}{\rightleftarrows} \GBMod[\mathbbm{k}] \colon U
	\label{eq:adjunction-gsnmod-gbmod-k}
\end{equation}
is naturally enhanced into a monoidal adjunction in which the left adjoint functor $\wedge$ is strong symmetric
monoidal via the coherence isomorphism
\begin{equation}
	M \cotimes_{\mathbbm{k}} N = \extrawidehat{M \otimes_{\mathbbm{k}} N}
	\xrightarrow[\cong]{\extrawidehat{\eta_M \otimes_{\mathbbm{k}} \eta_N}}
	\extrawidehat{\widehat{M} \otimes_{\mathbbm{k}} \widehat{N}} = \widehat{M} \cotimes_{\mathbbm{k}} \widehat{N}
	\label{eq:graded-complete-tensor-product-k-modules-strong-monoidal}
\end{equation}
and the right adjoint forgetful functor $U$ is lax symmetric monoidal via the canonical map
\begin{equation*}
	M \otimes_{\mathbbm{k}} N \xrightarrow{\eta_{M \otimes_{\mathbbm{k}} N}} M \cotimes_{\mathbbm{k}} N.
\end{equation*}

\subsection{Graded Seminormed and Banach \texorpdfstring{$\mathbbm{k}$}{k}-Algebras}

\begin{dfn} \label{dfn:graded-seminormed-banach-k-algebra}
	A \textbf{graded seminormed} $\mathbbm{k}$-\textbf{algebra} is a graded seminormed $\mathbbm{k}$-module
	$R = \left( R^g, \nnorm^g \right)_{g \in \GG}$
	together with a degree zero $\mathbbm{k}$-bilinear multiplication $\cdot \colon R \times R \rightarrow R$
	and a unit element $1_R \in R^0$ which satisfies the usual axioms of a graded $\mathbbm{k}$-algebra.
	In addition, we require the following compatibility conditions between the unit, multiplication and
	norm:\footnote{Recall that we always assume elements are homogeneous. Explicitly,
		\cref{eq:multiplication-compatibility-seminorm} means that
		$\nnorm[r^g \cdot r^h]^{g+h} \leq \nnorm[r^g]^g \cdot \nnorm[r^h]^h$ for all $r^g \in R^g, r^h \in R^h$
		and \cref{eq:unit-compatibility-seminorm} means that $\nnorm[1_R]^0 \leq 1$.
		In other words, we require that both the multiplication $\cdot \colon R \times R \rightarrow R$
		and the unit $u \colon \mathbbm{k} \rightarrow R$ are contractive.}
	\begin{align}
		\nnorm[r_1 \cdot r_2] & \leq \nnorm[r_1] \cdot \nnorm[r_2] \textrm{ for all } r_1,r_2 \in R,
		\label{eq:multiplication-compatibility-seminorm}
		\\
		\nnorm[1_R]           & \leq 1.
		\label{eq:unit-compatibility-seminorm}
	\end{align}
	A \textbf{graded Banach} $\mathbbm{k}$-\textbf{algebra} is a graded seminormed $\mathbbm{k}$-algebra $R$
	for which each component $\left( R^g, \nnorm^g \right)$ is a Banach $\mathbbm{k}$-module (i.e., normed
	and complete).
	When $\mathbbm{k} = \left( \ZZ, \trivnorm \right)$, a graded seminormed (resp.\ Banach) $\mathbbm{k}$-algebra
	is called a \textbf{graded seminormed} (resp.\ \textbf{Banach}) \textbf{ring}.
\end{dfn}

Equivalently, we can think of a graded seminormed $\mathbbm{k}$-algebra $R$ as a graded $\mathbbm{k}$-algebra
$R = \left( R^g \right)_{g \in \GG}$ equipped with the extra structure of a \textbf{graded}
$\mathbbm{k}$-\textbf{algebra seminorm}, i.e., a family $\nnorm = \left( \nnorm^g \right)_{g \in \GG}$ of
$\mathbbm{k}$-module seminorms on each component
$R^g$ which are compatible with the multiplication and unit in the sense of
\cref{eq:multiplication-compatibility-seminorm,eq:unit-compatibility-seminorm}.

Every graded seminormed $\mathbbm{k}$-algebra $R$ has an \textbf{underlying graded} $\mathbbm{k}$-\textbf{algebra}
obtained by forgetting the graded seminorm on $R$ while retaining the algebra structure. In addition, $R$ has an
\textbf{underlying graded seminormed} $\mathbbm{k}$-\textbf{module} obtained by forgetting the algebra structure
while retaining the graded seminorm. A graded seminormed (or Banach) $\mathbbm{k}$-algebra
is \textbf{graded-commutative} if the underlying graded $\mathbbm{k}$-algebra is graded-commutative.

A \textbf{morphism of graded seminormed} $\mathbbm{k}$-\textbf{algebras} $f \colon R \rightarrow S$
is a morphism of graded seminormed $\mathbbm{k}$-modules which is compatible with the multiplication and unit in the usual way. In other words, a morphism $f \colon R \rightarrow S$ is both a morphism of the underlying
graded $\mathbbm{k}$-algebras and a morphism of the underlying seminormed graded $\mathbbm{k}$-modules
so that $f$ is contractive, of degree zero, and compatible with the multiplication and unit.

\begin{rem} \label{rem:graded-seminormed-banach-algebra-as-algebra-object}
	Since the category $\GSNMod[\mathbbm{k}]$ is monoidal, one can talk about algebra objects in
	$\GSNMod[\mathbbm{k}]$ and an algebra object of $\GSNMod[\mathbbm{k}]$ is precisely a graded seminormed
	$\mathbbm{k}$-algebra in the sense
	of \cref{dfn:graded-seminormed-banach-k-algebra}. The $\mathbbm{k}$-bilinear multiplication
	$\cdot \colon R \times R \rightarrow R$ is contractive, hence induces a multiplication morphism
	$m \colon R \otimes_{\mathbbm{k}} R \rightarrow R$ of graded seminormed $\mathbbm{k}$-modules.
	Taking into account the symmetries of $\GSNMod[\mathbbm{k}]$
	(see \cref{eq:symmetries-graded-seminormed-k-modules}),
	a graded-commutative seminormed $\mathbbm{k}$-algebra is the same thing as a commutative
	algebra object of the symmetric monoidal category $\GSNMod[\mathbbm{k}]$.

	Similarly, a graded Banach $\mathbbm{k}$-algebra $R$
	is the same thing as an algebra object of the category $\GBMod[\mathbbm{k}]$ and $R$ is graded-commutative if
	and only if it is a commutative algebra object of $\GBMod[\mathbbm{k}]$.
	A priori an algebra object of $\GBMod[\mathbbm{k}]$ is a graded Banach $\mathbbm{k}$-module $R$ equipped
	with a ``multiplication'' of the form $R \cotimes_{\mathbbm{k}} R \rightarrow R$ which satisfies
	associativity conditions stated using $\cotimes_{\mathbbm{k}}$. By the universal property of the graded
	complete tensor product, such a multiplication corresponds bijectively to a
	degree zero contractive bilinear map
	$\cdot \colon R \times R \rightarrow R$ which endows $R$ with an associative multiplication.

	Equivalently, but phrased differently, the forgetful functor
	$\GBMod[\mathbbm{k}] \rightarrow \GSNMod[\mathbbm{k}]$ is lax monoidal and hence sends an algebra
	object of $\GBMod[\mathbbm{k}]$ to an algebra object of $\GSNMod[\mathbbm{k}]$. The multiplication map
	$R \otimes_{\mathbbm{k}} R \rightarrow R$ is obtained by precomposing the multiplication map
	$R \cotimes_{\mathbbm{k}} R \rightarrow R$ with the canonical map
	$R \otimes_{\mathbbm{k}} R \rightarrow R \cotimes_{\mathbbm{k}} R$.
\end{rem}

Let us denote by $\GSNAlg[\mathbbm{k}]$ the category of graded seminormed $\mathbbm{k}$-algebras and
by $\GBAlg[\mathbbm{k}]$ the full subcategory of $\GSNAlg[\mathbbm{k}]$ whose objects are graded Banach
$\mathbbm{k}$-algebras.
Since the completion functor $\wedge \colon \GSNMod[\mathbbm{k}] \rightarrow \GBMod[\mathbbm{k}]$ is
lax (even strong) monoidal, it extends naturally from modules
to algebras and gives us a completion functor
$\wedge \colon \GSNAlg[\mathbbm{k}] \rightarrow \GBAlg[\mathbbm{k}]$
which is left adjoint to the forgetful functor $\GBAlg[\mathbbm{k}] \rightarrow \GSNAlg[\mathbbm{k}]$.
The algebra structure on the completion $\widehat{R}$ is the unique structure with respect to
which the canonical completion morphism $\eta_R \colon R \rightarrow \widehat{R}$ becomes an isometric morphism of
graded $\mathbbm{k}$-algebras.

\subsection{Graded Seminormed and Banach Modules over Graded Algebras}

Generalizing \cref{sec:graded-seminormed-banach-k-modules}, we discuss
graded seminormed (and Banach) modules over graded seminormed $\mathbbm{k}$-algebras.
Let $R = \left( R^g, \nnorm^g_{R} \right)_{g \in \GG}$ be a graded seminormed $\mathbbm{k}$-algebra.

\begin{dfn} \label{def:graded-seminormed-Banach-R-module}
	A \textbf{graded seminormed left} $R$-\textbf{module} is a graded seminormed $\mathbbm{k}$-module
	$M = \left( M^g, \nnorm^g_{M} \right)_{g \in \GG}$
	together with a $\mathbbm{k}$-bilinear action map $\cdot \colon R \times M \rightarrow M$
	of degree zero which satisfies the usual axioms of a graded unital left module. In addition,
	we require that the action is compatible with the norms in the sense that
	\begin{equation}
		\nnorm[r \cdot m]_{M} \leq \nnorm[r]_{R} \cdot \nnorm[m]_{M}
		\label{eq:graded-R-module-seminorm-compatibility}
	\end{equation}
	for all $r \in R$ and $m \in M$.\footnote{Recall that we always assume elements are homogeneous. Explicitly,
		this conditions means that $\nnorm[r^g \cdot m^h]^{g+h}_{M} \leq \nnorm[r^g]^g_{R} \cdot \nnorm[m^h]^h_{M}$
		for all $r^g \in R^g$ and $m^h \in M^h$. In other words, we require that the action map
		$\cdot \colon R \times M \rightarrow M$ is contractive.}
	A \textbf{left graded Banach} $R$-\textbf{module} is a left graded seminormed $R$-module for which
	each component $\left( M^g, \nnorm^g_M \right)$ is a Banach $\mathbbm{k}$-module.
\end{dfn}

Equivalently, we can think of a graded seminormed left $R$-module $M$ as a graded left $R$-module
$M = \left( M^g \right)_{g \in \GG}$ equipped with the extra structure of a \textbf{graded}
$R$-\textbf{module seminorm}, i.e., a family $\nnorm = \left( \nnorm^g \right)_{g \in \GG}$ of
$\mathbbm{k}$-module seminorms on each component $M^g$
which are compatible with the $R$-action in the sense of \cref{eq:graded-R-module-seminorm-compatibility}.

Similarly, one can define graded seminormed (resp.\ Banach) right modules and bimodules. In what follows,
unless explicitly mentioned otherwise, the term ``module'' will always mean left module.

\begin{rem}
	Since the category $\GSNMod[\mathbbm{k}]$ (resp.\ $\GBMod[\mathbbm{k}]$) is monoidal,
	one can talk about module objects over algebra objects of $\GSNMod[\mathbbm{k}]$
	(resp.\ $\GBMod[\mathbbm{k}]$).
	A module object over an algebra object $R$ of $\GSNMod[\mathbbm{k}]$ (resp.\ $\GBMod[\mathbbm{k}]$)
	is precisely a graded seminormed (resp.\ Banach) $R$-module in the sense of
	\cref{def:graded-seminormed-Banach-R-module}.
\end{rem}

Every graded seminormed $R$-module $M$ has an \textbf{underlying graded} $R$-\textbf{module} obtained by forgetting
the graded seminorms on $M$ and $R$. Similarly, $M$ has an \textbf{underlying graded seminormed}
$\mathbbm{k}$-\textbf{module} obtained by forgetting the $R$-action while retaining the seminorm on $M$.
Hence, we can apply notions on graded $R$-modules from \cref{sub:graded-R-modules} and notions on
seminormed graded $\mathbbm{k}$-modules
from \cref{sec:graded-seminormed-banach-k-modules} to graded seminormed $R$-modules and will do
so without further mention.

A \textbf{morphism of graded seminormed} $R$-\textbf{modules} $f \colon M \rightarrow N$
is a graded contractive $R$-linear map of degree zero.\footnote{This guarantees that
	$f$ is both a morphism of the underlying graded $R$-modules and a morphism of the underlying
	graded seminormed $\mathbbm{k}$-modules.}
We denote by $\GSNMod[R]$ the category of graded seminormed $R$-modules with morphisms as above
and by $\GBMod[R]$ the full subcategory of $\GSNMod[R]$ whose objects are graded Banach $R$-modules.
An \textbf{isomorphism} $f \colon M \rightarrow N$ \textbf{of graded seminormed} (or \textbf{Banach})
$R$-\textbf{modules} is an isomorphism of the underlying graded $R$-modules which is an isometry.

Given two graded seminormed $R$-modules $M$ and $N$,
we denote by
\begin{equation*}
	\Hom{M}{N}[][R] \defeq
	\Set{f \colon M \rightarrow N}[f \textrm{ is a degree zero } R\textrm{-linear map with }
		{\nnorm[f]} \leq 1 ]
\end{equation*}
the set of morphisms between $M$ and $N$ and by $\InnHom{M}{N}[][R]$ the graded seminormed
$\mathbbm{k}$-module of all graded bounded $R$-linear maps with components given by
\begin{equation}
	\begin{aligned}
		\InnHom{M}{N}[][R]^d \defeq \{ & f \colon M \rightharpoonup N \, \rvert
		\\
		                               & f \textrm{ is a graded bounded } R\textrm{-linear map of degree } d \},
	\end{aligned} 	\label{eq:inner-hom-graded-seminormed-R-module}
\end{equation}
endowed with the operator seminorm.

\begin{rem}
	Even though we have chosen to use the same notation as in the
	algebraic case (see \cref{eq:inner-hom-graded-R-module}), whenever working with graded
	seminormed $R$-modules, $\Hom{M}{N}[][R]$ will always mean the set of graded \textit{contractive}
	$R$-linear maps of degree zero and $\InnHom{M}{N}[][R]$ will always mean the space
	of graded \textit{bounded} $R$-linear maps.
	Note that in the graded seminormed and Banach setting, we have
	$\Hom{M}{N}[][R] \neq \InnHom{M}{N}[][R]^0$
	since maps in $\InnHom{M}{N}[][R]^0$ are bounded while morphisms are taken to be contractive.
	We continue to use our convention to denote morphisms by $f \colon M \rightarrow N$ while reserving
	the notation $f \colon M \rightharpoonup N$ for graded bounded $R$-linear maps of arbitrary degree.
\end{rem}

\phantomsection
\label{sec:suspension-graded-seminormed-R-module}
Given a graded seminormed $R$-module $M$ and $h \in \GG$, the
$h$-\textbf{suspension} or $h$-\textbf{shifted} module $M[h]$ is defined
to be the $h$-suspension of the underlying graded $R$-module (see \cref{sec:suspension-graded-R-module})
together with the family of seminorms such that for all $m \in M$ we have
$\nnorm[\s_h \left( m \right)] \defeq \nnorm[m]$.
With the definition above, $M[h]$ becomes a graded seminormed $R$-module and the suspension map
$\s_h \colon M \rightharpoonup M[h]$ becomes a graded bounded $R$-linear map of degree $-h$ which is an isometry.
When $M$ is a graded Banach $R$-module, the suspension $M[h]$ is also a graded Banach
$R$-module.

When working over a graded-commutative seminormed $\mathbbm{k}$-algebra $R$,
any graded seminormed left $R$-module $M$ has a natural structure of a graded seminormed right $R$-module
via the action given by \cref{eq:convert-left-to-right-module}.
Endowing $M$ with both actions, one obtains a graded seminormed symmetric $R$-bimodule. In what follows, when
working over a graded-commutative seminormed $\mathbbm{k}$-algebra $R$,
we will identify left, right and symmetric bimodules as necessary.

Given a graded seminormed (resp.\ Banach) module $M$, there is a natural notion
of a graded seminormed (resp.\ Banach) submodule $N$ of $M$, such that the inclusion
morphism $N \hookrightarrow M$ is an isometry.
Given graded seminormed (resp.\ Banach) $R$-submodules
$N_1, \dots, N_k$ of $M$, we denote by $\left< N_1, \dots, N_k \right>$ the
graded \textit{seminormed} (resp.\ \textit{Banach}) $R$-submodule of $M$ generated by $N_1, \dots, N_k$.

Note that by taking $R = \mathbbm{k}$, considered as a graded seminormed $\mathbbm{k}$-algebra concentrated
in degree zero, we recover the categories $\GSNMod[\mathbbm{k}]$ and $\GBMod[\mathbbm{k}]$
described in \cref{sec:graded-seminormed-banach-k-modules}.
More generally, when $R$ is graded-commutative, the properties of the categories $\GSNMod[R]$ and $\GBMod[R]$ are analogous to the properties of $\GSNMod[\mathbbm{k}]$ and $\GBMod[\mathbbm{k}]$.
Although our interest lies mainly in the category $\GBMod[R]$ of graded Banach modules over a
graded Banach $\mathbbm{k}$-algebra $R$, we also discuss the seminormed case since, similar
to what happens in the ungraded case (see \cref{appendix:non-archimedean-groups-rings-modules}),
some constructions in $\GBMod[R]$ are performed in $\GSNMod[R]$ and then reflected (i.e., completed)
to obtain an object of $\GBMod[R]$. More details are given below.

\subsubsection{The Category \texorpdfstring{$\GSNMod[R]$}{GSNMod(R)}} \label{sub:category-GSNMod-R}
The category $\GSNMod[R]$ is preadditive and bicomplete. Limits (resp.\ colimits) are given by limits
(resp.\ colimits) of the underlying graded seminormed $\mathbbm{k}$-modules endowed with the natural $R$-action,
and we use the same notation for them as in \cref{sub:category-GSNMod-k}. In particular,
the category $\GSNMod[R]$ has finite biproducts, kernels and cokernels and hence is pre-abelian but not
abelian.

Given a graded seminormed right $R$-module $M$ and a graded seminormed left $R$-module $N$,
we can endow the graded tensor product $M \otimes_R N$ (see \eqref{eq:algebraic-graded-tensor-product})
with the structure of a graded seminormed $\mathbbm{k}$-module by setting
\begin{equation}
	\nnorm[x]_{M \otimes_R N} \defeq \inf \Set{\max_{i \in I} \, \nnorm[m_i]_M \cdot \nnorm[n_i]_N}
	[x = \sum_{i \in I} m_i \otimes_R n_i, \, |I| < \infty] \label{eq:graded-projective-tensor-seminorm}
\end{equation}
for $x \in M \otimes_R N$. More precisely,
we endow each graded piece $\left( M \otimes_R N \right)^d$ with a
$\mathbbm{k}$-module seminorm $\nnorm^d_{M \otimes_R N}$
by setting
\begin{align}
	\nnorm[x]^d_{M \otimes_R N} \defeq \inf
	\bigg\{ & \max_{i \in I} \, \nnorm[m_i]^{\degb{m_i}}_M \cdot \nnorm[n_i]^{\degb{n_i}}_N \, \biggr\rvert
	\label{eq:graded-projective-tensor-seminorm-explicit}
	\\
	        & x =
	\sum_{i \in I} m_i \otimes_R n_i, \, |I| < \infty, \, \degb{m_i} + \degb{n_i} = d, \, m_i \in M^{\degb{m_i}}, \, n_i \in N^{\degb{n_i}}
	\bigg\}
	\nonumber
\end{align}
for $x \in \left( M \otimes_R N \right)^d$.
The graded $\mathbbm{k}$-module seminorm
$\nnorm_{M \otimes_R N}$ is called the \textbf{graded non-Archimedean projective tensor seminorm}.

Similar to the graded and seminormed cases, the tensor product $M \otimes_R N$ of
a graded seminormed right $R$-module $M$ and a graded seminormed left $R$-module $N$
can be characterized by a universal property involving graded bounded $R$-balanced maps.
The tensor product $M \otimes_{R} N$ comes equipped with a canonical graded contractive
$R$-balanced map $\otimes_{R} \colon M \times N \rightarrow M \otimes_{R} N$ of degree zero
characterized by the following universal property: Given a graded seminormed $\mathbbm{k}$-module $L$ and
a graded bounded $R$-balanced map $B \colon M \times N \rightharpoonup L$ there exists a unique
graded bounded map $\varphi_B \colon M \otimes_{R} N \rightharpoonup L$
of graded seminormed $\mathbbm{k}$-modules
with $\nnorm[\varphi_B] = \nnorm[B]$ and $\degb{\varphi_B} = \degb{B}$ such that
\begin{equation*}
	\varphi_B \left( m \otimes_{R} n \right) = B(m,n)
\end{equation*}
for all $m \in M$ and $n \in N$ (see \cref{fig:tensor-product-graded-seminormed-R-modules-universal-property-R-balanced}).

\begin{figure}[tb]
	\centering
	\subcaptionbox{$R$ is a graded seminormed $\mathbbm{k}$-algebra, \\ $M$ is a graded seminormed right $R$-module, \\
		$N$ is a graded seminormed left $R$-module, \\ $L$ is a graded seminormed $\mathbbm{k}$-module.
		\label{fig:tensor-product-graded-seminormed-R-modules-universal-property-R-balanced}}
	[.48\linewidth]{
		\begin{tikzcd}[ampersand replacement=\&]
			{M \times N} \& \& {M \otimes_{R} N} \\
			\& \& L \\
			\arrow["\otimes_{R}", from=1-1, to=1-3]
			\arrow["{\substack{\exists! \, \varphi_B \\ \textrm{graded bounded} \\ \mathbbm{k}\textrm{-linear}}}",
				dashed, harpoon, from=1-3, to=2-3]
			\arrow["\substack{B \\ \textrm{graded bounded} \\ R\textrm{-balanced}}"',
				harpoon, from=1-1, to=2-3]
		\end{tikzcd}
	}
	\subcaptionbox{$R$ is a graded-commutative seminormed $\mathbbm{k}$-algebra, \\
		$M,N,L$ are graded seminormed $R$-modules.
		\label{fig:tensor-product-graded-seminormed-R-modules-universal-property-R-linear}}[.48\linewidth]{
		\begin{tikzcd}[ampersand replacement=\&]
			{M \times N} \& \& {M \otimes_{R} N} \\
			\& \& L \\
			\arrow["\otimes_{R}", from=1-1, to=1-3]
			\arrow["{\substack{\exists! \, \varphi_B \\ \textrm{graded bounded} \\ R\textrm{-linear}}}",
				dashed, harpoon, from=1-3, to=2-3]
			\arrow["\substack{B \\ \textrm{graded bounded} \\ R\textrm{-bilinear}}"',
				harpoon, from=1-1, to=2-3]
		\end{tikzcd}
	}
	\caption{Universal properties of the tensor product of graded seminormed $R$-modules.}
	\label{fig:tensor-product-graded-seminormed-R-modules-universal-properties}
\end{figure}

In what follows, we will assume that $R$ is graded-commutative and endow $\GSNMod[R]$ with a
monoidal structure. Given two graded seminormed left $R$-modules $M$ and $N$, we can convert $M$ to a
graded seminormed right $R$-module using \cref{eq:convert-left-to-right-module}
and form the tensor product $M \otimes_R N$. The tensor
product has a natural structure of a graded left $R$-module with the action given by
\cref{eq:R-action-R-tensor-product}
and the seminorm given by \cref{eq:graded-projective-tensor-seminorm} on $M \otimes_R N$
is compatible with the $R$-action. Hence, $M \otimes_R N$ has the structure of a graded seminormed $R$-module.

Similar to the graded and seminormed case, the tensor product $M \otimes_R N$
of two graded seminormed $R$-modules over a graded-commutative seminormed $\mathbbm{k}$-algebra $R$
can be characterized by a universal property involving graded bounded $R$-bilinear maps.
The canonical map $\otimes_R \colon M \times N \rightarrow M \otimes_R N$ is
not only $R$-balanced but also $R$-bilinear and is characterized by the following universal property:
Given a graded seminormed $R$-module $L$ and a graded bounded $R$-bilinear map
$B \colon M \times N \rightharpoonup L$ there exists a
unique graded bounded $R$-linear map $\varphi_B \colon M \otimes_{R} N \rightharpoonup L$
with $\nnorm[\varphi_B] = \nnorm[B]$ and $\degb{\varphi_B} = \degb{B}$
such that $\varphi_B \left( m \otimes_{R} n \right) = B(m,n)$ for all $m \in M$ and $n \in N$ (see
\cref{fig:tensor-product-graded-seminormed-R-modules-universal-property-R-linear}).

Given two graded bounded $R$-linear maps $f \colon M \rightharpoonup M'$ and
$g \colon N \rightharpoonup N'$ between graded seminormed $R$-modules, their
tensor product $f \otimes_R g$ (see \eqref{eq:tensor-product-graded-R-linear-maps}) is also
bounded with $\nnorm[f \otimes_R g] \leq \nnorm[f] \cdot \nnorm[g]$.
Restricting our attention to graded contractive $R$-linear maps of degree zero, we obtain
a bifunctor $\otimes_{R} \colon \GSNMod[R] \times \GSNMod[R] \rightarrow \GSNMod[R]$. The bifunctor
$\otimes_{R}$, the symmetry maps $M \otimes_{R} N \rightarrow N \otimes_{R} M$ given by
\begin{equation}
	m \otimes_{R} n \mapsto (-1)^{\braidd{m}{n}} n \otimes_{R} m,
	\label{eq:symmetries-graded-seminormed-R-modules}
\end{equation}
and the standard associators and unitors, endow the category $\GSNMod[R]$ with the structure of a symmetric
monoidal category whose unit is the ground $\mathbbm{k}$-algebra $R$ (considered as graded
seminormed $R$-module over itself).

Since we work over a graded-commutative ground algebra $R$, the graded seminormed $\mathbbm{k}$-module
$\InnHom{M}{N}[][R]$ given by \cref{eq:inner-hom-graded-seminormed-R-module} has an
$R$-action given by \cref{eq:left-R-action-on-hom}.
The operator seminorm on $\InnHom{M}{N}[][R]$ is compatible
with the $R$-action and hence $\InnHom{M}{N}[][R]$ is not only a graded seminormed $\mathbbm{k}$-module but
also a graded seminormed $R$-module, i.e., an object of $\GSNMod[R]$.
It follows from the universal property of the graded seminormed tensor product $\otimes_R$
that the symmetric monoidal category $\GSNMod[R]$ is closed, with the internal hom object
given by $\InnHom{M}{N}[][R]$ which consists of all graded \textit{bounded} $R$-linear maps.

Given graded seminormed $R$-modules $M_1, \dots, M_n$ and $N$, we will use the notation
$\Mult{M_1,\dots,M_n}{N}[R]$ to denote the graded seminormed $R$-module of all graded \textit{bounded}
$R$-multilinear maps $B \colon M_1 \times \dots \times M_n \rightharpoonup N$, endowed with the operator seminorm
\eqref{eq:norm-of-graded-multilinear-map} and the natural $R$-module structure given by
\eqref{eq:R-action-on-multilinear-maps}. We have natural isomorphisms of graded seminormed
$R$-modules (i.e., \textit{isometric} bijections)
\begin{equation*}
	\Mult{M_1,\dots,M_n}{N}[R] \cong \InnHom{M_1 \otimes_R \dots \otimes_R M_n}{N}[][R],
\end{equation*}
given by the same formulas as in the graded case.
That is, we have bijective isometric correspondence
between graded bounded $R$-multilinear maps
$B \colon M_1 \times \dots \times M_n \rightharpoonup N$
and graded bounded $R$-linear maps $\varphi_B \colon M_1 \otimes_R \dots \otimes_R M_n \rightharpoonup N$
generalizing the correspondence shown in
\cref{fig:tensor-product-graded-seminormed-R-modules-universal-property-R-linear}.

Since $\GSNMod[R]$ is closed symmetric monoidal,
the tensor product of graded seminormed $R$-modules commutes with colimits in each variable,
and we have natural isomorphisms of graded seminormed $R$-modules (i.e., \textit{isometric} bijections)
\begin{equation*}
	\left( \bigoplus_{i \in I} M_i \right) \otimes_{R} N \cong
	\bigoplus_{i \in I} M_i \otimes_{R} N, \quad
	M \otimes_{R} \left( \bigoplus_{i \in I} N_i \right) \cong
	\bigoplus_{i \in I} M \otimes_{R} N_i
\end{equation*}
given by the same formulas as in the graded case.

The scalar extension and restriction constructions for graded $R$-modules
(see \cref{sec:scalar-extension-restriction-graded-modules}) extend verbatim to graded seminormed $R$-modules.
Given two graded seminormed $\mathbbm{k}$-algebras $R$ and $S$ and a morphism
$\varphi \colon R \rightarrow S$ of graded seminormed $\mathbbm{k}$-algebras, we get a pair
of adjoint functors
\begin{equation}
	\varphi_{!} \colon \GSNMod[R] \stackrel[]{\dashv}{\rightleftarrows} \GSNMod[S] \colon \varphi^{*},
	\label{eq:restriction-extension-adjunction-GSNMod-R}
\end{equation}
where all the formulas which play a role in the adjunction are the same as in the
graded case (see \cref{sec:scalar-extension-restriction-graded-modules}).

Given a graded seminormed $R$-module $M$ and a graded seminormed $S$-module $N$,
a \textbf{morphism of graded seminormed modules over} $\varphi$ is a graded \textit{contractive}
$\mathbbm{k}$-linear map $f \colon M \rightarrow N$ of degree zero which satisfies
$f \left( r \cdot m \right) = \varphi(r) \cdot f(m)$
for all $r \in R$ and $m \in M$. Equivalently, a morphism of graded seminormed modules over $\varphi$
is a morphism $f \colon M \rightarrow \varphi^{*} \left( N \right)$
of graded seminormed $R$-modules between $M$ and the restriction of scalars of $N$ along $\varphi$.
Given two morphisms $f_1 \colon M_1 \rightarrow N_1$ and $f_2 \colon M_2 \rightarrow N_2$
of graded seminormed modules over $\varphi$, the morphism
$f_1 \otimes_{\varphi} f_2 \colon M_1 \otimes_R M_2 \rightarrow N_1 \otimes_S N_2$ over $\varphi$
given by \cref{eq:tensor-product-morphisms-over-phi} is also contractive,
hence a morphism of graded seminormed modules over $\varphi$.

\subsubsection{The Category \texorpdfstring{$\GBMod[R]$}{GBMod(R)}}
\label{sub:category-GBMod-R}

Even though the notion of a graded Banach $R$-module
makes sense when $R$ is merely seminormed, when working with graded Banach $R$-modules we will always
assume that the graded ground $\mathbbm{k}$-algebra $R$ is also Banach.\footnote{This is convenient and does not affect much since the categories $\GBMod[R]$ and $\GBMod (\widehat{R})$ are isomorphic.
	See \cref{sec:banach-modules-over-banach-rings} for the ungraded version.}
Before discussing the properties of $\GBMod[R]$, we discuss the completion procedure
for graded seminormed $R$-modules.

\phantomsection
\label{sec:completion-graded-seminormed-R-module}
Given a graded seminormed $R$-module $M$, the completion $\widehat{M}$, which is a priori
a graded seminormed $\mathbbm{k}$-module, has a unique structure of a graded Banach $R$-module
with respect to which the canonical morphism $\eta_M \colon M \rightarrow \widehat{M}$ becomes
an isometric $R$-linear map of graded seminormed $R$-modules. The action
$R \times \widehat{M} \rightarrow \widehat{M}$ of $R$ on $\widehat{M}$ is obtained by extending
the action $R \times M \rightarrow M$ by density.
Given a graded bounded $R$-linear map
$f \colon M \rightharpoonup N$ from a graded seminormed $R$-module $M$ into a graded Banach
$R$-module $N$, the adjunct $\tilde{f} \colon \widehat{M} \rightharpoonup N$ becomes a graded
bounded $R$-linear map. In particular, this implies that given a graded bounded $R$-linear map
$f \colon M_1 \rightharpoonup M_2$ between graded seminormed $R$-modules, the completion
$\widehat{f} \colon \widehat{M_1} \rightharpoonup \widehat{M_2}$ also becomes a graded
bounded $R$-linear map. Hence, we obtain a completion functor
$\wedge \colon \GSNMod[R] \rightarrow \GBMod[R]$ which is left adjoint to the natural forgetful
functor $U \colon \GBMod[R] \rightarrow \GSNMod[R]$ and so
$\GBMod[R]$ is a reflective replete full subcategory of $\GSNMod[R]$.

The category $\GBMod[R]$ is preadditive and bicomplete. Limits (resp.\ colimits) are given by limits
(resp.\ colimits) of the underlying graded Banach $\mathbbm{k}$-modules endowed with the natural
$R$-action. We will use the same notation for limits and colimits of graded Banach $R$-modules
as for limits and colimits of graded Banach $\mathbbm{k}$-modules.
Note that since $\GBMod[R]$ is reflective replete full subcategory of
$\GSNMod[R]$, limits in $\GBMod[R]$ are the same as limits in $\GSNMod[R]$ (and are automatically
Banach) while colimits are obtained by computing colimits in $\GSNMod[R]$ and completing them.
In particular, the category $\GBMod[R]$ has finite biproducts, kernels and cokernels and hence
is pre-abelian but not abelian.

Given a graded seminormed right $R$-module $M$ and a graded seminormed left $R$-module $N$,
their \textbf{graded complete tensor product} is the graded Banach $\mathbbm{k}$-module
$\extrawidehat{M \otimes_R N}$ obtained by completing the graded seminormed tensor product $M \otimes_R N$.
Given $m \in M$ and $n \in N$, we will denote by $m \cotimes_{R} n$ the image of
$m \otimes_{R} n$ under the completion map
$\eta \colon M \otimes_{R} N \rightarrow M \cotimes_{R} N$
and call such elements \textbf{elementary tensors}. Note that unlike the algebraic or seminormed case,
the complete tensor product $M \cotimes_{R} N$ is not generated as a graded $\mathbbm{k}$-module by
elementary tensors. Instead, the graded $\mathbbm{k}$-module generated by elementary tensors is
dense in $M \cotimes_{R} N$. In other words, $M \cotimes_{R} N$ is generated by elementary tensors as a
graded \textit{Banach} $\mathbbm{k}$-module.

Similar to the graded case, the complete tensor product $M \cotimes_R N$ of a graded right seminormed
$R$-module $M$ and a graded left seminormed $R$-module $N$ can be characterized by a universal property involving
graded bounded $R$-balanced maps into graded \textit{Banach} $\mathbbm{k}$-modules. The complete
tensor product $M \cotimes_R N$ comes equipped with a canonical graded contractive $R$-balanced
map $\cotimes_R \colon M \times N \rightarrow M \cotimes_R N$ characterized by the following universal property:
Given a graded Banach $\mathbbm{k}$-module $L$ and a graded bounded $R$-balanced map
$B \colon M \times N \rightharpoonup L$, there exists a unique graded bounded map
$\varphi_B \colon M \cotimes_R N \rightharpoonup L$ of graded Banach $\mathbbm{k}$-modules with
$\nnorm[\varphi_B] = \nnorm[B]$ and $\degb{\varphi_B} = \degb{B}$ such that
\begin{equation*}
	\varphi_B \left( m \cotimes_R n \right) = B(m,n)
\end{equation*}
for all $m \in M$ and $n \in N$
(see \cref{fig:complete-tensor-product-graded-seminormed-R-modules-universal-property-R-balanced}
and compare to \cref{fig:tensor-product-graded-seminormed-R-modules-universal-property-R-balanced}).

In what follows, we will assume that $R$ is graded-commutative and endow $\GBMod[R]$ with a monoidal structure.
Given two graded seminormed left $R$-modules $M$ and $N$, the graded seminormed tensor product
$M \otimes_R N$ is not only a graded seminormed $\mathbbm{k}$-module but also a graded
seminormed $R$-module. Hence, the completion $M \cotimes_R N = \extrawidehat{M \otimes_R N}$ has
the structure of a graded Banach $R$-module.

Similar to the graded and seminormed case, the complete tensor product $M \cotimes_R N$
of two graded seminormed $R$-modules over a graded-commutative Banach $\mathbbm{k}$-algebra $R$
can be characterized by a universal property involving graded bounded $R$-bilinear maps
into graded \textit{Banach} $R$-modules.
The canonical map $\cotimes_R \colon M \times N \rightarrow M \cotimes_R N$ is
not only $R$-balanced but also $R$-bilinear and is characterized by the following universal property:
Given a graded Banach $R$-module $L$ and a graded bounded $R$-bilinear map
$B \colon M \times N \rightharpoonup L$ there exists a
unique graded bounded $R$-linear map $\varphi_B \colon M \cotimes_{R} N \rightharpoonup L$
with $\nnorm[\varphi_B] = \nnorm[B]$ and $\degb{\varphi_B} = \degb{B}$
such that
\begin{equation*}
	\varphi_B \left( m \cotimes_{R} n \right) = B(m,n)
\end{equation*}
for all $m \in M$ and $n \in N$
(see \cref{fig:complete-tensor-product-graded-seminormed-R-modules-universal-property-R-linear}
and compare to \cref{fig:tensor-product-graded-seminormed-R-modules-universal-property-R-linear}).

\begin{figure}
	\centering
	\subcaptionbox{$R$ is a graded Banach $\mathbbm{k}$-algebra, \\ $M$ is a graded seminormed right $R$-module, \\
		$N$ is a graded seminormed left $R$-module, \\ $L$ is a graded Banach $\mathbbm{k}$-module.
		\label{fig:complete-tensor-product-graded-seminormed-R-modules-universal-property-R-balanced}}
	[.48\linewidth]{
		\begin{tikzcd}[ampersand replacement=\&]
			{M \times N} \& \& {M \cotimes_{R} N} \\
			\& \& L \\
			\arrow["\cotimes_{R}", from=1-1, to=1-3]
			\arrow["{\substack{\exists! \, \varphi_B \\ \textrm{graded bounded} \\ \mathbbm{k}\textrm{-linear}}}",
				dashed, harpoon, from=1-3, to=2-3]
			\arrow["\substack{B \\ \textrm{graded bounded} \\ R\textrm{-balanced}}"',
				harpoon, from=1-1, to=2-3]
		\end{tikzcd}
	}
	\subcaptionbox{$R$ is a graded-commutative Banach $\mathbbm{k}$-algebra, \\
		$M,N$ are graded seminormed $R$-modules, \\ $L$ is a graded Banach $R$-module.
		\label{fig:complete-tensor-product-graded-seminormed-R-modules-universal-property-R-linear}}[.48\linewidth]{
		\begin{tikzcd}[ampersand replacement=\&]
			{M \times N} \& \& {M \cotimes_{R} N} \\
			\& \& L \\
			\arrow["\cotimes_{R}", from=1-1, to=1-3]
			\arrow["{\substack{\exists! \, \varphi_B \\ \textrm{graded bounded} \\ R\textrm{-linear}}}",
				dashed, harpoon, from=1-3, to=2-3]
			\arrow["\substack{B \\ \textrm{graded bounded} \\ R\textrm{-bilinear}}"',
				harpoon, from=1-1, to=2-3]
		\end{tikzcd}
	}
	\caption{Universal properties of the complete tensor product of graded seminormed $R$-modules.}
	\label{fig:complete-tensor-product-graded-seminormed-R-modules-universal-properties}
\end{figure}

Given two graded bounded maps $f \colon M \rightharpoonup M'$ and
$g \colon N \rightharpoonup N'$ between graded seminormed $R$-modules, their complete tensor product
$f \cotimes_{R} g \colon M \cotimes_{R} N \rightharpoonup M' \cotimes_{R} N'$ is
defined by $f \cotimes_{R} g \defeq \widehat{f \otimes_{R} g}$ and is the unique
graded bounded $R$-linear map of degree $\degb{f} + \degb{g}$ which satisfies
\begin{equation}
	\left( f \cotimes_{R} g \right) \left( m \cotimes_{R} n \right) =
	(-1)^{\braidd{g}{m}} f \left( m \right) \cotimes_{R} g \left( n \right)
	\label{eq:complete-tensor-product-graded-R-linear-maps}
\end{equation}
for all $m \in M$ and $n \in N$. In terms of the bounds on $f,g$ we have
$\nnorm[f \cotimes_{R} g] \leq \nnorm[f] \cdot \nnorm[g]$.
The interaction between composition and complete tensor product of graded bounded $R$-linear maps is the same as
in \cref{eq:interaction-composition-tensor-product},
with $\otimes_{\mathbbm{k}}$ replaced by $\cotimes_{R}$.
In particular, we see that the graded complete tensor product construction gives us a bifunctor
$\cotimes_{R} \colon \GSNMod[R] \times \GSNMod[R] \rightarrow \GBMod[R]$.

Even though we have defined the complete tensor product for graded seminormed $R$-modules,
we now restrict our attention to the complete tensor product of graded Banach $R$-modules
and obtain a bifunctor
$\cotimes_{R} \colon \GBMod[R] \times \GBMod[R] \rightarrow \GBMod[R]$.
Given three graded Banach $R$-modules $M,N,L$, there exist
natural associativity isomorphisms
$\left( M \cotimes_{R} N \right) \cotimes_{R} L \cong
	M \cotimes_{R} \left( N \cotimes_{R} L \right)$
which act on elementary tensors by the expected formula
$\left( m \cotimes_{R} n \right) \cotimes_{R} l \mapsto
	m \cotimes_{R} \left( n \cotimes_{R} l \right)$.
In addition, we also have natural isometric isomorphisms
\begin{align*}
	       & R \cotimes_{R} M \cong M \qquad          & r \cotimes_{R} m                               & \mapsto r \cdot m,                      \\
	       & M \cotimes_{R} R \cong M \qquad          & m \cotimes_{R} r                               & \mapsto (-1)^{\braidd{r}{m}} r \cdot m, \\
	       & M \cotimes_{R} N \cong N \cotimes_{R} M,
	\qquad & m \cotimes_{R} n                         & \mapsto (-1)^{\braidd{m}{n}} n \cotimes_{R} m.
\end{align*}
The bifunctor
$\cotimes_{R} \colon \GBMod[R] \times \GBMod[R] \rightarrow \GBMod[R]$,
the associativity isomorphisms, the unitors, and the symmetry maps described above,
endow the category $\GBMod[R]$
with the structure of a symmetric monoidal category whose unit is the ground
$\mathbbm{k}$-algebra $R$, considered as a graded Banach $R$-module over itself.

Given two graded seminormed $R$-modules $M$ and $N$, the graded seminormed $R$-module
$\InnHom{M}{N}[][R]$ of all graded bounded $R$-linear maps of arbitrary degree is Banach once
$N$ is Banach. In particular, $\InnHom{M}{N}[][R]$ is an object of $\GBMod[R]$ when
both $M$ and $N$ are objects of $\GBMod[R]$.
It follows from the universal property of the graded complete tensor product
that the symmetric monoidal category $\GBMod[R]$ is closed with internal hom object
$\InnHom{M}{N}[][R]$ of all graded bounded $R$-linear maps (the same internal hom object as in $\GSNMod[R]$).

Given graded Banach $R$-modules $M_1, \dots, M_n$ and $N$, the graded seminormed $R$-module
$\Mult{M_1,\dots,M_n}{N}[R]$ of all graded \textit{bounded} $R$-multilinear maps
$B \colon M_1 \times \dots \times M_n \rightharpoonup N$ is also Banach,
and we have natural isomorphisms of graded Banach $R$-modules
(i.e., \textit{isometric} isomorphisms of graded $R$-modules)
\begin{equation*}
	\Mult{M_1,\dots,M_n}{N}[R] \cong \InnHom{M_1 \cotimes_R \dots \cotimes_R M_n}{N}[][R]
\end{equation*}
given by the same formulas as in the graded case, using $\cotimes_R$ instead of $\otimes_R$.
That is, we have bijective isometric correspondence
between graded bounded $R$-multilinear maps
$B \colon M_1 \times \dots \times M_n \rightharpoonup N$
and graded bounded $R$-linear maps $\varphi_B \colon M_1 \cotimes_R \dots \cotimes_R M_n \rightharpoonup N$
generalizing the correspondence shown in
\cref{fig:complete-tensor-product-graded-seminormed-R-modules-universal-property-R-linear}.

Since $\GBMod[R]$ is closed symmetric monoidal,
the tensor product of graded Banach $R$-modules commutes with colimits in each variable,
and we have natural isomorphisms of graded Banach $R$-modules
\begin{equation*}
	\left( \cbigoplus_{i \in I} M_i \right) \cotimes_{R} N \cong
	\cbigoplus_{i \in I} M_i \cotimes_{R} N, \quad
	M \cotimes_{R} \left( \cbigoplus_{i \in I} N_i \right) \cong
	\cbigoplus_{i \in I} M \cotimes_{R} N_i.
\end{equation*}
In addition, the internal hom bifunctor
$\InnHom{-}{-}[][R]$ preserves limits in the second variable, and sends colimits in the
first variable to limits, so we have natural isomorphism of graded Banach $R$-modules
\begin{equation*}
	\InnHom{\cbigoplus_{i \in I} M_i}{N}[][R] \cong \prod_{i \in I}^{\B} \InnHom{M_i}{N}[][R], \quad
	\InnHom{M}{\prod_{i \in I}^{\B} N_i}[][R] \cong \prod_{i \in I}^{\B} \InnHom{M}{N_i}[][R].
\end{equation*}
The isomorphisms are given by standard formulas, adapted to the graded Banach case.

With respect to the monoidal structures defined on $\GSNMod[R]$ and $\GBMod[R]$, the adjunction
\begin{equation}
	\wedge \colon \GSNMod[R] \stackrel[]{\dashv}{\rightleftarrows} \GBMod[R] \colon U
	\label{eq:adjunction-gsnmod-gbmod-R}
\end{equation}
is naturally enhanced into a monoidal adjunction in which the left adjoint functor $\wedge$ is strong symmetric
monoidal via the coherence isomorphism
\begin{equation}
	M \cotimes_{R} N = \extrawidehat{M \otimes_{R} N}
	\xrightarrow[\cong]{\extrawidehat{\eta_M \otimes_{R} \eta_N}}
	\extrawidehat{\widehat{M} \otimes_{R} \widehat{N}} = \widehat{M} \cotimes_{R} \widehat{N}
	\label{eq:graded-complete-tensor-product-R-modules-strong-monoidal}
\end{equation}
and the right adjoint forgetful functor $U$ is lax symmetric monoidal via the canonical map
\begin{equation*}
	M \otimes_{R} N \xrightarrow{\eta_{M \otimes_R N}} M \cotimes_{R} N.
\end{equation*}

The scalar extension and restriction procedures for graded $R$-modules
(see \cref{sec:scalar-extension-restriction-graded-modules}) can be adapted naturally
to work for graded Banach $R$-modules. Given two graded-commutative Banach $\mathbbm{k}$-algebras $R$ and $S$
and a morphism $\varphi \colon R \rightarrow S$ of graded Banach $\mathbbm{k}$-algebras, we have a monoidal adjunction
\begin{equation}
	\varphi_{!} \colon \GBMod[R] \stackrel[]{\dashv}{\rightleftarrows} \GBMod[S] \colon \varphi^{*}
	\label{eq:restriction-extension-adjunction-GBMod-R}
\end{equation}
where all the formulas which play a role in the adjunction are the same as in the
graded case, with $\otimes_R$ replaced by $\cotimes_R$.

Given a graded Banach $R$-module $M$ and a graded Banach $S$-module $N$, a
\textbf{morphism of graded Banach modules over} $\varphi$ is a graded \textit{contractive}
$\mathbbm{k}$-linear map $f \colon M \rightarrow N$ of degree zero which satisfies
$f \left( r \cdot m \right) = \varphi(r) \cdot f(m)$ for all $r \in R$ and $m \in M$.
Equivalently, a morphism of graded Banach modules over $\varphi$ is a morphism
$f \colon M \rightarrow \varphi^{*} \left( N \right)$
of graded Banach $R$-modules between $M$ and the restriction of scalars of $N$ along $\varphi$.

Given two morphisms $f_1 \colon M_1 \rightarrow N_1$ and $f_2 \colon M_2 \rightarrow N_2$
of graded Banach modules over $\varphi$, we will denote by
$f_1 \cotimes_{\varphi} f_2 \colon M_1 \cotimes_R M_2 \rightarrow N_1 \cotimes_S N_2$ the morphism of
graded Banach modules over $\varphi$ given by
\begin{equation}
	\left( f_1 \cotimes_{\varphi} f_2 \right) \left( m_1 \cotimes_R m_2 \right) \defeq
	f_1 \left( m_1 \right) \cotimes_S f_2 \left( m_2 \right).
	\label{eq:complete-tensor-product-morphisms-over-phi}
\end{equation}
This is the same definition as for graded $R$-modules, with $\otimes_R$ replaced by $\cotimes_R$ (see
\cref{eq:tensor-product-morphisms-over-phi}).
When $R = S$ and $\varphi = \id_R$, we recover the definition of $f_1 \cotimes_R f_2$ of the
tensor product of morphisms of graded Banach $R$-modules given by
\cref{eq:complete-tensor-product-graded-R-linear-maps}.

We end this section with a simple but useful criterion for defining an inverse of an operator using a pointwise geometric series.

\begin{lm} \label{lm:inverse-geometric-series}
	Let $R$ be a graded-commutative Banach $\mathbbm{k}$-algebra and let $M$ be a graded Banach $R$-module.
	Let $f \colon M \rightarrow M$ be a morphism of graded Banach $R$-modules
	and assume we have $f^n(m) \to 0$ for all $m \in M$.
	Then $\idd - f$ is invertible with an inverse $g$ which satisfies $\nnorm[g] \leq 1$.
	The inverse $g$ is given by the pointwise converging geometric series $g \defeq \sum_{n=0}^{\infty} f^n$.
\end{lm}
\begin{proof}
	Since we are in the non-Archimedean setting, the condition $f^n(m) \to 0$ guarantees
	that the infinite sum $\sum_{n=0}^{\infty} f^n \left( m \right)$ converges for all $m \in M$.
	Define a map $g \colon M \rightarrow M$ by $g(m) \defeq \lim_{N \to \infty} \sum_{n=0}^N f^n \left( m \right)$.
	Then $g$ is a bounded $R$-linear map with $\nnorm[g] \leq 1$. We have
	\begin{equation*}
		\left( g \circ \left( \idd - f \right) \right) \left( m \right) =
		g \left( m - f \left( m \right) \right) =
		\lim_{N \to \infty} \sum_{n=0}^N f^n \left( m - f \left( m \right) \right) =
		\lim_{N \to \infty} \left( m - f^{N+1} \left( m \right) \right) = m
	\end{equation*}
	and since $f$ is continuous, we also have
	\begin{equation*}
		\left( \left( \idd - f \right) \circ g \right) \left( m \right) =
		\lim_{N \to \infty} \sum_{n=0}^N \left( \idd - f \right) \left( f^n \left( m \right) \right) =
		\lim_{N \to \infty} \left( m - f^{N+1} \left( m \right) \right) = m.
	\end{equation*}
\end{proof}

\subsection{Graded Seminormed and Banach Coalgebras}
\label{sub:graded-seminormed-banach-coalgebras}

Let $R$ be a graded-commutative Banach $\mathbbm{k}$-algebra.
We can consider the notion of a coalgebra object (see \cref{subsec:coalg-in-monoidal-cat})
in two different monoidal categories: $\GSNMod[R]$ and $\GBMod[R]$. This will result in two different notions,
described below.

A coalgebra object in $\GSNMod[R]$ is given by a graded seminormed $R$-module
$\left( C, \nnorm_C \right)$ together with a coproduct $\Delta \colon C \rightarrow C \otimes_R C$ and
a counit map $\varepsilon \colon C \rightarrow R$ satisfying the usual axioms of a counital coalgebra.
We will call $\left( C, \nnorm_C \right)$ a \textbf{graded seminormed} $R$-\textbf{coalgebra}.
Since we work in $\GSNMod[R]$, both the coproduct and the counit are required to satisfy
$\nnorm[\Delta], \nnorm[\varepsilon] \leq 1$ and are in particular continuous. Note that by applying
the forgetful functor $\GSNMod[R] \rightarrow \GMod[R]$ which forgets the seminorm
$\nnorm_C$ (and the seminorm on $R$ which is implicit in our notation),
we obtain a graded $R$-coalgebra $C$, called the \textbf{underlying graded} $R$-\textbf{coalgebra}.
Conversely, given a coalgebra object $C$ in $\GMod[R]$, any seminorm $\nnorm_C$ we put on $C$,
which is compatible with the norm on $R$ and with respect to
which the comultiplication and counit maps satisfy $\nnorm[\Delta], \nnorm[\varepsilon] \leq 1$, turns
$\left( C, \nnorm_C \right)$ into a graded seminormed $R$-coalgebra.
A \textbf{morphism} $f \colon C \rightarrow D$ \textbf{of graded seminormed} $R$-\textbf{coalgebras} is a degree zero
contractive $R$-linear map which satisfies $\left( f \otimes_R f \right) \circ \Delta_C = \Delta_D \circ f$ and
$\varepsilon_D \circ f = \varepsilon_C$.

A coalgebra object in $\GBMod[R]$ is given by a graded Banach $R$-module $\left( C, \nnorm_C \right)$
together with a coproduct $\Delta \colon C \rightarrow C \cotimes_R C$ and a counit map
$\varepsilon \colon C \rightarrow R$ satisfying
$\left( \Delta \cotimes_R \id \right) \circ \Delta = \left( \id \cotimes_R \Delta \right) \circ \Delta$
and $\left( \varepsilon \cotimes_R \id \right) \circ \Delta = \left( \id \cotimes_R \varepsilon \right) \circ \Delta = \id$.\footnote{
	As usual, we identify both sides using the associator and unitor isomorphisms coming
	from the monoidal structure on $\GBMod[R]$.}
Note that in the definition of $\Delta$ we use the complete tensor product $\cotimes_R$
(the monoidal structure on $\GBMod[R]$) and not the usual (seminormed) tensor product $\otimes_R$
(the monoidal structure on $\GSNMod[R]$).
Again, both the coproduct and the counit are required to be contractive.
We will call $\left( C, \nnorm_C \right)$ a \textbf{graded Banach} $R$-\textbf{coalgebra}.
Since the forgetful functor $\GBMod[R] \rightarrow \GSNMod[R]$ is only lax monoidal and not
strong, in general a graded Banach $R$-coalgebra is not a
graded seminormed $R$-coalgebra (or a graded $R$-coalgebra). While we have a natural map
$C \otimes_R C \rightarrow C \cotimes_R C$, this map is not an isomorphism (or even injective)
and in general we cannot factor
the coproduct $\Delta \colon C \rightarrow C \cotimes_R C$ through a map $C \rightarrow C \otimes_R C$.
A \textbf{morphism} $f \colon C \rightarrow D$ \textbf{of graded Banach} $R$-\textbf{coalgebras} is a degree zero
contractive $R$-linear map which satisfies $\left( f \cotimes_R f \right) \circ \Delta_C = \Delta_D \circ f$ and
$\varepsilon_D \circ f = \varepsilon_C$.

Since the completion functor $\GSNMod[R] \rightarrow \GBMod[R]$ is strong
monoidal, the completion $\widehat{C}$ of a graded seminormed $R$-coalgebra is naturally a graded Banach $R$-coalgebra.
Moreover, if $f \colon C \rightarrow D$
is a morphism of graded seminormed $R$-coalgebras then
$\widehat{f} \colon \widehat{C} \rightarrow \widehat{D}$ is a morphism
of graded Banach $R$-coalgebras.

We will also need the notion of a coalgebra morphism between two graded Banach coalgebras over
different ground Banach algebras. This is a straightforward adaptation of the corresponding
notation for graded coalgebras (see \cref{sec:coalgebra-morphism-over-different-ground-algebras}).
Let $R$ and $S$ be two graded-commutative Banach $\mathbbm{k}$-algebras and let
$\varphi \colon R \rightarrow S$ be a morphism of graded Banach $\mathbbm{k}$-algebras.
Given a graded Banach $R$-coalgebra $(C,\Delta_C,\varepsilon_C)$ and a graded Banach $S$-coalgebra
$(D,\Delta_D,\varepsilon_D)$ a \textbf{Banach coalgebra morphism} $f \colon C \rightarrow D$ \textbf{over an
	algebra morphism} $\varphi \colon R \rightarrow S$ is a morphism of graded Banach modules over $\varphi$
such that the following diagrams commute:

\begin{figure}[H]
	\centering
	\begin{subfigure}{0.45\textwidth}
		\centering
		\begin{tikzcd}
			C & D \\
			{C \cotimes_R C} & {D \cotimes_S D} \\
			\arrow["f", from=1-1, to=1-2]
			\arrow["{\Delta_C}"', from=1-1, to=2-1]
			\arrow["{f \cotimes_{\varphi} f}", from=2-1, to=2-2]
			\arrow["{\Delta_D}", from=1-2, to=2-2]
		\end{tikzcd}
		\caption{Compatibility with the coproducts.}
		\label{fig:morphisms-Banach-coalgebra-different-ground-algebra-succinct-coproduct}
	\end{subfigure}
	\begin{subfigure}{0.45\textwidth}
		\centering
		\begin{tikzcd}
			C & D \\
			R & S \\
			\arrow["f", from=1-1, to=1-2]
			\arrow["{\varepsilon_D}", from=1-2, to=2-2]
			\arrow["{\varepsilon_C}"', from=1-1, to=2-1]
			\arrow["\varphi", from=2-1, to=2-2]
		\end{tikzcd}
		\caption{Compatibility with the counits.}
		\label{fig:morphisms-Banach-coalgebra-different-ground-algebra-succinct-counit}
	\end{subfigure}
	\caption{Morphism between two graded Banach coalgebras over different ground Banach algebras.}
	\label{fig:morphisms-Banach-coalgebra-different-ground-algebra-succinct}
\end{figure}

Equivalently, a Banach coalgebra morphism over $\varphi$ is given by a Banach $S$-coalgebra morphism
$\tilde{f} \colon \varphi_{!} \left( C \right) \rightarrow D$ where
$\varphi_{!} \left( C \right) = S \cotimes_R C$ is the Banach $S$-coalgebra obtained by scalar extension
(see \cref{rem:coalgebra-morphism-over-algebra-morphism-adjunction} for a discussion of the graded case,
which translates directly to the Banach setting, with the obvious modifications).
The relation between $f$ and $\tilde{f}$ is given by
\begin{equation*}
	\tilde{f} \left( s \cotimes_R c \right) = s \cdot f \left( c \right).
\end{equation*}
When $R = S$ and $\varphi = \id_R$, we recover the notion of a Banach coalgebra morphism defined above.

Clearly Banach coalgebra morphisms over Banach algebra morphisms can be composed so that if
$g \colon B \rightarrow C$ is a Banach coalgebra morphism over
$\psi \colon Q \rightarrow R$ and $f \colon C \rightarrow D$ is a Banach coalgebra morphism over
$\varphi \colon R \rightarrow S$ then $f \circ g \colon B \rightarrow D$ is a Banach coalgebra morphism over
$\varphi \circ \psi \colon Q \rightarrow S$.

\subsection{Pre-Differential Graded Seminormed and Banach \texorpdfstring{$\mathbbm{k}$}{k}-Modules}
\label{sub:pre-differential-graded-seminormed-banach-k-modules}
Starting with this subsection and for the remainder of the section,
we fix an element $\go \in \GG$
which is considered as part of the grading datum (see \cref{subsec:grading-data}).
Pre-differentials on objects will raise the degree of elements by $\go \in \GG$.
We extend the notion of pre-differential graded $\mathbbm{k}$-module from
\cref{sec:pre-differential-graded-k-modules} to the seminormed
and Banach settings.

\begin{dfn}
	A \textbf{pre-differential graded seminormed} (resp.\ \textbf{Banach}) $\mathbbm{k}$-\textbf{module} is a pair
	$\mathcal{M} = \left( M, d \right)$ where $M$ is a graded seminormed (resp.\ Banach) $\mathbbm{k}$-module,
	called the \textbf{underlying graded seminormed} (resp.\ \textbf{Banach}) $\mathbbm{k}$-\textbf{module}
	and $d \colon M \rightharpoonup M$ is a graded bounded map of degree
	$\go$ called the \textbf{pre-differential}.
\end{dfn}

Every pre-differential graded seminormed $\mathbbm{k}$-module $\mathcal{M}$ has an
\textbf{underlying pre-differential graded} $\mathbbm{k}$-\textbf{module} obtained by forgetting
the graded seminorm on $M$. Hence, we can apply notions on pre-differential graded $\mathbbm{k}$-modules
from \cref{sec:pre-differential-graded-k-modules} to pre-differential graded seminormed
$\mathbbm{k}$-modules and will do so without further mention.

A \textbf{morphism} of pre-differential graded seminormed $\mathbbm{k}$-modules
$f \colon \left( M, d_M \right) \rightarrow \left( N, d_N \right)$ is a morphism $f \colon M \rightarrow N$
of the underlying graded seminormed $\mathbbm{k}$-modules (i.e., a graded contractive map of degree zero)
which is also compatible with the pre-differentials in the sense that $f \circ d_M = d_N \circ f$. In other words,
a morphism of pre-differential graded seminormed $\mathbbm{k}$-modules is both a morphism of the underlying
graded seminormed $\mathbbm{k}$-modules and a morphism of the underlying pre-differential graded $\mathbbm{k}$-modules. Let us denote by $\PDGSNMod[\mathbbm{k}]$ the category of pre-differential graded seminormed $\mathbbm{k}$-modules
with morphisms as defined above and by $\PDGBMod[\mathbbm{k}]$ the full subcategory of $\PDGSNMod[\mathbbm{k}]$
consisting of pre-differential graded Banach $\mathbbm{k}$-modules.

\begin{rem}
	Let $M$ be a graded seminormed $\mathbbm{k}$-module. While from a categorical perspective,
	it is more natural to consider pre-differentials $d \colon M \rightharpoonup M$ which are contractive,\footnote{
		A ``degree $\go$'' map in the graded category should be a family
		$d^g \colon M^g \rightarrow M^{g + \go}$ of morphisms in the ungraded category
		$\SNMod[\mathbbm{k}]$ between the
		components. Since the morphisms $\SNMod[\mathbbm{k}]$ are contractive, this
		means that we should require that $\nnorm[d^g] \leq 1$ for all $g \in \GG$.}
	we have chosen to allow the pre-differentials to be bounded. This has the unfortunate
	consequence that the categories $\PDGSNMod[\mathbbm{k}]$ and $\PDGBMod[\mathbbm{k}]$ do not
	admit arbitrary infinite limits and colimits, but they are still pre-abelian. This won't cause
	us any problems but if one wants bicomplete categories then one can work with contractive pre-differentials.
\end{rem}

Given a pre-differential graded seminormed $\mathbbm{k}$-module $\mathcal{M} = \left( M, d_M \right)$
and $h \in \GG$, the $h$-\textbf{suspension} or $h$-\textbf{shifted} module $\mathcal{M}[h]$
is the pre-differential graded seminormed $\mathbbm{k}$-module whose underlying
graded seminormed $\mathbbm{k}$-module is $M[h]$
(see \cref{sec:suspension-graded-seminormed-Banach-k-module}), endowed with the pre-differential
given by \cref{eq:differential-on-suspension} (the same pre-differential as in the graded case).
When $\mathcal{M}$ is a pre-differential graded Banach $\mathbbm{k}$-module, the suspension
$\mathcal{M}[h]$ is also a pre-differential graded Banach $\mathbbm{k}$-module.

The closed symmetric monoidal structure on graded seminormed (resp.\ Banach) $\mathbbm{k}$-modules
extends naturally to the pre-differential setting by endowing the (resp.\ complete) tensor product
and the internal hom with pre-differentials as in the graded case.
More precisely, given two pre-differential graded seminormed $\mathbbm{k}$-modules
$\mathcal{M} = \left( M, d_M \right)$ and $\mathcal{N} = \left( N, d_N \right)$,
the \textbf{internal hom object} $\InnHom{\mathcal{M}}{\mathcal{N}}[][\mathbbm{k}]$
is defined to be the graded seminormed $\mathbbm{k}$-module
$\InnHom{M}{N}[][\mathbbm{k}]$ of all graded \textit{bounded} maps (see
\cref{eq:inner-hom-graded-seminormed-k-module}) endowed with the pre-differential
$\partial$ given by \cref{eq:hom-differential}.\footnote{Note that we use the
	internal hom object in the category of graded seminormed $\mathbbm{k}$-modules and not the internal
	hom object in the category of graded $\mathbbm{k}$-modules which consists of all graded, not necessarily bounded,
	maps.}
The pre-differential $\partial$ satisfies
\begin{equation*}
	\nnorm[\partial \left( f \right)] \leq \max \Set{\nnorm[d_M], \nnorm[d_N]} \nnorm[f]
\end{equation*}
and hence $\partial$ maps graded bounded maps to graded bounded maps and $\partial$ is also bounded.
Hence, $\InnHom{\mathcal{M}}{\mathcal{N}}[][\mathbbm{k}]$ is indeed a pre-differential graded seminormed
$\mathbbm{k}$-module.
The \textbf{tensor product} $\mathcal{M} \otimes_{\mathbbm{k}} \mathcal{N}$ of two
pre-differential graded seminormed $\mathbbm{k}$-modules is defined to be the
graded seminormed $\mathbbm{k}$-module $M \otimes_{\mathbbm{k}} N$ endowed with the pre-differential
$d_{M \otimes_{\mathbbm{k}} N}$ given by \cref{eq:k-tensor-product-differential}.
The pre-differential $d_{M \otimes_{\mathbbm{k}} N}$ satisfies
$\nnorm[d_{M \otimes_{\mathbbm{k}} N}] \leq \max \Set{\nnorm[d_M], \nnorm[d_N]}$ and hence it is bounded
and $\mathcal{M} \otimes_{\mathbbm{k}} \mathcal{N}$ is indeed a pre-differential graded seminormed
$\mathbbm{k}$-module.

With the definitions above, the closed symmetric monoidal structure on graded seminormed $\mathbbm{k}$-modules
extends to pre-differential graded seminormed $\mathbbm{k}$-modules and the category $\PDGSNMod[\mathbbm{k}]$ becomes closed symmetric monoidal.

The situation for $\PDGBMod[\mathbbm{k}]$ is similar. The internal hom object of two pre-differential graded
Banach $\mathbbm{k}$-modules $\mathcal{M} = \left( M, d_M \right)$ and $\mathcal{N} = \left( N, d_N \right)$
is the same as in the seminormed case, while their \textbf{complete tensor product} is given by
the graded Banach $\mathbbm{k}$-module $M \cotimes_{\mathbbm{k}} N$ endowed with the pre-differential
$d_{M \cotimes_{\mathbbm{k}} N}$ given by
\begin{equation}
	d_{M \cotimes_{\mathbbm{k}} N} \defeq d_M \cotimes_{\mathbbm{k}} \id + \id \cotimes_{\mathbbm{k}} d_N,
	\label{eq:k-complete-tensor-product-differential}
\end{equation}
the same formula as in the graded and seminormed case, with $\otimes_{\mathbbm{k}}$ replaced by
$\cotimes_{\mathbbm{k}}$. With the definitions above, the closed symmetric monoidal structure
on graded Banach $\mathbbm{k}$-modules extends to pre-differential graded Banach $\mathbbm{k}$-modules
and the category $\PDGBMod[\mathbbm{k}]$ becomes closed symmetric monoidal.

Given a pre-differential graded seminormed $\mathbbm{k}$-module $\mathcal{M} = \left( M, d_M \right)$,
the \textbf{completion} of $\mathcal{M}$ is the pre-differential graded Banach $\mathbbm{k}$-module given by
$\widehat{\mathcal{M}} \defeq ( \widehat{M}, \widehat{d_M} )$.
Given a graded bounded chain map $f \colon M \rightharpoonup N$, i.e.,
a map which satisfies $\partial \left( f \right) = 0$, the completion
$\widehat{f} \colon \widehat{M} \rightharpoonup \widehat{N}$ also satisfies
$\partial ( \widehat{f} ) = 0$. In particular, the completion of a morphism
$f \colon \left( M, d_M \right) \rightarrow \left( N, d_N \right)$ of pre-differential graded
seminormed $\mathbbm{k}$-modules is also a morphism,\footnote{That is, we have
	$\widehat{f} \circ \widehat{d_M} = \widehat{d_N} \circ \widehat{f}$.}
and we obtain a completion functor $\wedge \colon \PDGSNMod[\mathbbm{k}] \rightarrow \PDGBMod[\mathbbm{k}]$
which is left adjoint to the forgetful
functor $U \colon \PDGBMod[\mathbbm{k}] \rightarrow \PDGSNMod[\mathbbm{k}]$.
Hence, we see that the category $\PDGBMod[\mathbbm{k}]$ is a reflective replete full subcategory of
$\PDGSNMod[\mathbbm{k}]$ with reflector $\wedge$.
The completion morphism $\eta_M \colon M \rightarrow \widehat{M}$ is a chain map and serves as the
unit $\mathcal{M} \rightarrow \widehat{\mathcal{M}}$ of the adjunction, just like it does for
the adjunction \eqref{eq:adjunction-gsnmod-gbmod-k}.

With respect to the monoidal structures defined on the categories $\PDGSNMod[\mathbbm{k}]$
and $\PDGBMod[\mathbbm{k}]$, the adjunction
\begin{equation*}
	\wedge \colon \PDGSNMod[\mathbbm{k}] \stackrel[]{\dashv}{\rightleftarrows} \PDGBMod[\mathbbm{k}] \colon U
\end{equation*}
is naturally enhanced into a monoidal adjunction in which the left adjoint functor $\wedge$ is strong symmetric
monoidal and the right adjoint forgetful functor $U$ is lax symmetric monoidal with
the same coherence maps as for the adjunction \eqref{eq:adjunction-gsnmod-gbmod-k} (i.e., the same
maps are also chain maps).

\subsection{Pre-Differential Graded Seminormed and Banach \texorpdfstring{$\mathbbm{k}$}{k}-Algebras}
Next, we extend the notion of pre-differential graded $\mathbbm{k}$-algebra
from \cref{sec:pre-differential-graded-k-algebras} to the seminormed and Banach settings.

\begin{dfn}
	A \textbf{pre-differential graded seminormed} (resp.\ \textbf{Banach}) $\mathbbm{k}$-\textbf{algebra}
	is a pair $\mathcal{R} = \left( R, d \right)$ where $R$ is a graded seminormed
	(resp.\ Banach) $\mathbbm{k}$-algebra,
	called the \textbf{underlying graded seminormed} (resp.\ \textbf{Banach}) $\mathbbm{k}$-\textbf{algebra},
	and $d \colon R \rightharpoonup R$ is a graded bounded algebra derivation of
	degree $\go$.
\end{dfn}

We note that a pre-differential graded seminormed (resp.\ Banach) $\mathbbm{k}$-algebra is the same thing as
an algebra object $\mathcal{R}$ of the monoidal category $\PDGSNMod[\mathbbm{k}]$ (resp.\ $\PDGBMod[\mathbbm{k}]$).
A \textbf{morphism of pre-differential graded seminormed} $\mathbbm{k}$-\textbf{algebras}
$f \colon \left( R, d_R \right) \rightarrow \left( S, d_S \right)$ is a morphism $f \colon R \rightarrow S$
of the underlying graded seminormed $\mathbbm{k}$-algebras which is compatible with the pre-differentials
in the sense that $f \circ d_R = d_S \circ f$. More explicitly,
a morphism $f \colon R \rightarrow S$ is a graded contractive map of degree zero which is compatible with the multiplication and unit in the usual way and also compatible with the pre-differentials.
A pre-differential graded seminormed (or Banach) $\mathbbm{k}$-algebra
is \textbf{graded-commutative} if the underlying graded $\mathbbm{k}$-algebra is graded-commutative.

Let us denote by $\PDGSNAlg[\mathbbm{k}]$ the category of pre-differential graded seminormed $\mathbbm{k}$-algebras
with morphisms as defined above and by $\PDGBAlg[\mathbbm{k}]$ the full subcategory of $\PDGSNAlg[\mathbbm{k}]$
consisting of pre-differential graded Banach $\mathbbm{k}$-algebras.
The completion functor $\wedge \colon \GSNAlg[\mathbbm{k}] \rightarrow \GBAlg[\mathbbm{k}]$ extends
to the pre-differential setting by completing the pre-differentials.
This gives us a
functor $\wedge \colon \PDGSNAlg[\mathbbm{k}] \rightarrow \PDGBAlg[\mathbbm{k}]$
which is left adjoint to the natural forgetful functor $\PDGBAlg[\mathbbm{k}] \rightarrow \PDGSNAlg[\mathbbm{k}]$.
The canonical unit morphism $\eta_{\mathcal{R}} \colon \mathcal{R} \rightarrow \widehat{\mathcal{R}}$ then
becomes a morphism of pre-differential graded seminormed $\mathbbm{k}$-algebras.

\subsection{Pre-Differential Graded Seminormed and Banach Modules over Algebras}
\label{sub:pre-differential-graded-seminormed-banach-R-modules}

Finally, we extend the notion of pre-differential graded module over a pre-differential graded algebra
from \cref{sec:pre-differential-graded-R-modules} to the seminormed and Banach setting.

Let us fix a pre-differential graded seminormed $\mathbbm{k}$-algebra $\mathcal{R} = (R,d_R)$.
\begin{dfn}
	A \textbf{pre-differential graded seminormed
		left} (resp.\ \textbf{right}) 	$\mathcal{R}$-\textbf{module} is a pair $\mathcal{M} = \left( M, d_M \right)$
	where $M$ is a graded seminormed left (resp.\ right) $R$-module,
	called the \textbf{underlying graded seminormed} $R$-\textbf{module}, and
	$d_M \colon M \rightharpoonup M$ is a bounded left (resp.\ right)
	$d_R$-operator on $M$ called the \textbf{pre-differential}. When $M$ is Banach, we call
	$\mathcal{M}$ a \textbf{pre-differential graded Banach left} (resp.\ \textbf{right}) $\mathcal{R}$-\textbf{module}.
\end{dfn}

In what follows, unless explicitly stated otherwise, the term ``module'' will always mean a left module.
To lessen the burden of notation, when no confusion is possible,
we call a pre-differential graded seminormed (resp.\ Banach) $\mathcal{R}$-module $\mathcal{M}$ simply a
seminormed (resp.\ Banach) $\mathcal{R}$-module, dropping
the adjectives ``pre-differential'' and ``graded'' and relying on context and the calligraphic font to remind
us that the module is graded and equipped with a pre-differential.

We note that a seminormed (resp.\ Banach) $\mathcal{R}$-module is the same thing as
a module object $\mathcal{M}$ over an algebra object $\mathcal{R}$ in the monoidal category $\PDGSNMod[\mathbbm{k}]$
(resp.\ $\PDGBMod[\mathbbm{k}]$).
Note also that a pre-differential graded seminormed (resp.\ Banach) $\mathcal{R}$-module
$\mathcal{M} = \left( M, d_M \right)$ is a pre-differential graded seminormed (resp.\ Banach) $\mathbbm{k}$-module, together with an $R$-action on $M$ such that the
pre-differential $d_M$ is compatible with the $R$-action and $d_R$ via \cref{eq:d-operator-Leibniz-rule}.
By forgetting the $R$-action and the compatibility condition of $d_M$, we obtain the
\textbf{underlying pre-differential graded seminormed} (resp.\ \textbf{Banach})
$\mathbbm{k}$-\textbf{module} of $\mathcal{M}$.

A \textbf{morphism} $f \colon (M,d_M) \rightarrow (N, d_N)$
\textbf{of pre-differential graded seminormed} $\mathcal{R}$-\textbf{modules}
is a
morphism $f \colon M \rightarrow N$ of the underlying graded seminormed $R$-modules (i.e., a degree zero
contractive $R$-linear map) compatible with the pre-differentials in the sense that
$f \circ d_M = d_N \circ f$.\footnote{In other words, a morphism $f \colon \mathcal{M} \rightarrow \mathcal{N}$
	of pre-differential graded seminormed $\mathcal{R}$-modules is both a morphism of the underlying
	graded seminormed $R$-modules and a morphism of the underlying pre-differential graded seminormed $\mathbbm{k}$-modules.}
Let us denote by $\PDGSNMod[\mathcal{R}]$ the category of pre-differential graded seminormed
$\mathcal{R}$-modules with morphisms as defined above and by $\PDGBMod[\mathcal{R}]$ the full subcategory of
$\PDGSNMod[\mathcal{R}]$ whose objects are pre-differential graded Banach $\mathcal{R}$-modules.

Given a seminormed $\mathcal{R}$-module $\mathcal{M} = \left( M, d_M \right)$ and $h \in \GG$, the
$h$-\textbf{suspension} or $h$-\textbf{shifted} module $\mathcal{M}[h]$ is defined
to be the suspension of the underlying graded seminormed $R$-module $M[h]$
(see \cref{sec:suspension-graded-seminormed-R-module})
endowed with the pre-differential given by \cref{eq:differential-on-suspension} (the same pre-differential
as in the graded case). When $\mathcal{M}$ is a Banach $\mathcal{R}$-module, the suspension $\mathcal{M}[h]$
is also a Banach $\mathcal{R}$-module.

In what follows, we will assume that $\mathcal{R} = \left( R, d_R \right)$ is
a pre-differential graded-commutative seminormed $\mathbbm{k}$-algebra.
The closed symmetric monoidal structure on graded seminormed $R$-modules described in \cref{sub:category-GSNMod-R}
extends naturally to the pre-differential setting by endowing the tensor product
and the internal hom with pre-differentials as in the graded case.
More precisely, given two seminormed $\mathcal{R}$-modules
$\mathcal{M} = \left( M, d_M \right)$ and $\mathcal{N} = \left( N, d_N \right)$,
the \textbf{internal hom object} $\InnHom{\mathcal{M}}{\mathcal{N}}[][\mathcal{R}]$ is defined to be the
graded seminormed $R$-module $\InnHom{M}{N}[][R]$ of all graded \textit{bounded} $R$-linear maps (see
\cref{eq:inner-hom-graded-seminormed-R-module}) endowed with the pre-differential
$\partial$ given by \cref{eq:hom-differential}.\footnote{Note that we use the
	internal hom object in the category of graded seminormed $R$-modules and not the internal
	hom object in the category of graded $R$-modules which consists of all graded, not necessarily bounded,
	maps.} The pre-differential $\partial$ is a bounded $d_R$-operator on $\InnHom{M}{N}[][R]$ and hence
$\InnHom{\mathcal{M}}{\mathcal{N}}[][\mathcal{R}]$ is indeed a pre-differential graded seminormed $\mathcal{R}$-module.

The \textbf{tensor product} $\mathcal{M} \otimes_{\mathcal{R}} \mathcal{N}$ of two
seminormed $\mathcal{R}$-modules is defined to be the
graded seminormed $R$-module $M \otimes_{R} N$ endowed with the pre-differential
$d_{M \otimes_{R} N}$ given by \cref{eq:R-tensor-product-differential}.
The pre-differential $d_{M \otimes_{R} N}$ is a bounded $d_R$-operator on $M \otimes_R N$ and
hence $\mathcal{M} \otimes_{\mathcal{R}} \mathcal{N}$ is indeed a pre-differential graded seminormed
$\mathcal{R}$-module.

With the definitions above, the closed symmetric monoidal structure on graded seminormed $R$-modules
extends to pre-differential graded seminormed $\mathcal{R}$-modules and the category $\PDGSNMod[\mathcal{R}]$
becomes closed symmetric monoidal.

The situation for $\PDGBMod[\mathcal{R}]$ is similar. Even though the notion of a
Banach $\mathcal{R}$-module makes sense when $\mathcal{R}$ is a pre-differential
graded seminormed $\mathbbm{k}$-algebra,
when working with Banach $\mathcal{R}$-modules we will always assume that $\mathcal{R}$ is
a pre-differential graded \textit{Banach} $\mathbbm{k}$-algebra. This is the same assumption
we made when discussing graded Banach $R$-modules without pre-differentials (see beginning of \cref{sub:category-GBMod-R}).

The internal hom object of two Banach
$\mathcal{R}$-modules $\mathcal{M} = \left( M, d_M \right)$ and $\mathcal{N} = \left( N, d_N \right)$
is the same as in the seminormed case, while their \textbf{complete tensor product} is given by
the graded Banach $R$-module $M \cotimes_{R} N$ endowed with the pre-differential
$d_{M \cotimes_{R} N}$ given by
\begin{equation}
	d_{M \cotimes_{R} N} \defeq d_M \cotimes_{R} \id + \id \cotimes_{R} d_N,
	\label{eq:R-complete-tensor-product-differential}
\end{equation}
the same formula as in the graded and seminormed case, with $\otimes_{R}$ replaced by
$\cotimes_{R}$. With the definitions above, the closed symmetric monoidal structure
on graded Banach $R$-modules described in \cref{sub:category-GBMod-R} extends to
pre-differential graded Banach $\mathcal{R}$-modules and the category $\PDGBMod[\mathcal{R}]$
becomes closed symmetric monoidal.

Given a seminormed $\mathcal{R}$-module $\mathcal{M} = \left( M, d_M \right)$,
the \textbf{completion} of $\mathcal{M}$ is the Banach $\mathcal{R}$-module given by
$\widehat{\mathcal{M}} \defeq ( \widehat{M}, \widehat{d_M} )$ (see \cref{sec:completion-graded-seminormed-R-module}).
The completion of a graded bounded $R$-linear map $f \colon M \rightharpoonup N$
is also bounded and $R$-linear, and we have
$\widehat{\partial \left( f \right)} = \partial ( \widehat{f} )$.
In particular, given a morphism
$f \colon \left( M, d_M \right) \rightarrow \left( N, d_N \right)$ of seminormed
$\mathcal{R}$-modules, the completion
$\widehat{f} \colon ( \widehat{M}, \widehat{d_M} ) \rightarrow ( \widehat{N}, \widehat{d_N} )$
is a morphism of Banach $\mathcal{R}$-modules, and we obtain a completion functor
$\wedge \colon \PDGSNMod[\mathcal{R}] \rightarrow \PDGBMod[\mathcal{R}]$ which is left adjoint to the forgetful
functor $U \colon \PDGBMod[\mathcal{R}] \rightarrow \PDGSNMod[\mathcal{R}]$.
Thus $\PDGBMod[\mathcal{R}]$ is a reflective replete full subcategory of $\PDGSNMod[\mathcal{R}]$.
The completion map $\eta_M \colon M \rightarrow \widehat{M}$ is an $R$-linear chain map and serves as the
unit $\mathcal{M} \rightarrow \widehat{\mathcal{M}}$ of the adjunction, just like it does for
the adjunction \eqref{eq:adjunction-gsnmod-gbmod-R}.

With respect to the monoidal structures defined on the categories $\PDGSNMod[\mathcal{R}]$
and $\PDGBMod[\mathcal{R}]$, the adjunction
\begin{equation*}
	\wedge \colon \PDGSNMod[\mathcal{R}] \stackrel[]{\dashv}{\rightleftarrows} \PDGBMod[\mathcal{R}] \colon U
\end{equation*}
is naturally enhanced into a monoidal adjunction in which the left adjoint functor $\wedge$ is strong symmetric
monoidal and the right adjoint forgetful functor $U$ is lax symmetric monoidal with
the same coherence maps as for the adjunction \eqref{eq:adjunction-gsnmod-gbmod-R} (i.e., the same
maps are also chain maps).

\subsection{Pre-Differential Graded Seminormed and Banach Coalgebras}
We extend the notion of pre-differential graded $\mathcal{R}$-coalgebras from \cref{sec:pre-differential-graded-coalgebras} to the seminormed and Banach setting.

\begin{dfn} \label{def:coderivation-R-linear-seminormed-Banach-coalgebra} \leavevmode
	\begin{enumerate}
		\item Let $R$ be a graded-commutative seminormed $\mathbbm{k}$-algebra and let $C$ be
		      a graded seminormed $R$-coalgebra. A (graded) \textbf{coderivation}
		      is a graded \textit{bounded} $R$-linear map $\mu \colon C \rightharpoonup C$
		      which satisfies the graded co-Leibniz rule
		      \begin{equation}
			      \Delta \circ \mu =
			      \left( \mu \otimes_R \id + \id \otimes_R \mu \right) \circ \Delta.
			      \label{eq:coderivation-equation-seminormed}
		      \end{equation}
		\item Let $R$ be a graded-commutative Banach $\mathbbm{k}$-algebra and let $C$ be
		      a graded Banach $R$-coalgebra. A (graded) \textbf{coderivation}
		      is a graded \textit{bounded} $R$-linear map $\mu \colon C \rightharpoonup C$
		      which satisfies the graded co-Leibniz rule with respect to $\cotimes$:
		      \begin{equation}
			      \Delta \circ \mu =
			      \left( \mu \cotimes_R \id + \id \cotimes_R \mu \right) \circ \Delta.
			      \label{eq:coderivation-equation-Banach}
		      \end{equation}
	\end{enumerate}
\end{dfn}

Note that a coderivation $\mu$ on a graded seminormed $R$-coalgebra $C$ is the same thing
as a coderivation on the underlying graded $R$-coalgebra which is bounded. Even though
one can talk about not necessarily bounded coderivations on $C$,
we will always assume that coderivations on a graded seminormed $R$-coalgebra are bounded.
In contrast, a graded Banach $R$-coalgebra $C$ has no underlying graded $R$-coalgebra and a coderivation
$\mu$ is defined using \cref{eq:coderivation-equation-Banach} instead of
\cref{eq:coderivation-equation-seminormed}.
Note also that it doesn't make sense to talk about coderivations on a graded Banach $R$-coalgebra
which are not bounded since \cref{eq:coderivation-equation-Banach} does not make sense if
$\mu$ is not bounded.
We have chosen to use the same terminology of ``coderivation''
both for the seminormed and Banach case even though they mean different things in each setting
since we believe confusion is highly unlikely.

More generally, we can talk about coderivations on a seminormed (resp.\ Banach) $R$-coalgebra $C$ over
graded algebra derivations of the ground seminormed (resp.\ Banach) $\mathbbm{k}$-algebra $R$:

\begin{dfn} \label{def:coderivation-graded-seminormed-Banach-coalgebra} \leavevmode
	\begin{enumerate}
		\item Let $R$ be a graded-commutative seminormed $\mathbbm{k}$-algebra and let
		      $d \colon R \rightharpoonup R$ be a graded bounded algebra derivation. Let $C$ be
		      a graded seminormed $R$-coalgebra. A \textbf{coderivation over} $d$ (or a
		      \textbf{generalized coderivation}) is a graded \textit{bounded} $\mathbbm{k}$-linear map
		      $\mu \colon C \rightharpoonup C$ of the same degree as $d$
		      such that $\mu$ is a graded module derivation
		      over $d$ in the sense of \cref{dfn:d-operator} and, in addition, it satisfies the graded
		      co-Leibniz rule \cref{eq:coderivation-equation-seminormed}.
		\item Let $R$ be a graded-commutative Banach $\mathbbm{k}$-algebra and let
		      $d \colon R \rightharpoonup R$ be a graded bounded algebra derivation. Let $C$ be
		      a graded Banach $R$-coalgebra. A \textbf{coderivation over} $d$ (or a
		      \textbf{generalized coderivation}) is a graded \textit{bounded} $\mathbbm{k}$-linear map
		      $\mu \colon C \rightharpoonup C$ of the same degree as $d$
		      such that $\mu$ is a graded module derivation
		      over $d$ in the sense of \cref{dfn:d-operator} which, in addition, satisfies the graded
		      co-Leibniz rule \eqref{eq:coderivation-equation-Banach} with respect to $\cotimes$.
	\end{enumerate}
\end{dfn}

Let $R$ be a graded-commutative Banach $\mathbbm{k}$-algebra and let $d \colon R \rightharpoonup R$ be
a graded bounded algebra derivation. Given a seminormed $R$-coalgebra $C$ and a coderivation
$\mu \colon C \rightharpoonup C$ over $d$, its completion
$\widehat{\mu} \colon \widehat{C} \rightharpoonup \widehat{C}$ is a coderivation over $d$
on the graded Banach $R$-coalgebra $\widehat{C}$.

\begin{dfn}
	Let $\mathcal{R} = \left( R, d \right)$ be a pre-differential graded-commutative
	seminormed (resp.\ Banach) $\mathbbm{k}$-algebra.
	A \textbf{pre-differential graded seminormed} (resp.\ \textbf{Banach}) $\mathcal{R}$-\textbf{coalgebra} is a
	pair $\mathcal{C} = \left( C, \mu \right)$, where $C$ is a graded seminormed (resp.\ Banach) $R$-coalgebra called
	the \textbf{underlying graded seminormed} (resp.\ \textbf{Banach}) $R$-\textbf{coalgebra},
	and $\mu \colon C \rightharpoonup C$ is a coderivation over $d$.
\end{dfn}

We note that a pre-differential graded seminormed (resp.\ Banach) $\mathcal{R}$-coalgebra is
the same thing as a coalgebra object $\mathcal{C}$ of the monoidal category $\PDGSNMod[\mathcal{R}]$
(resp.\ $\PDGBMod[\mathcal{R}]$). To lessen the burden of notation, when no confusion is possible,
we call a pre-differential graded seminormed (resp.\ Banach)
$\mathcal{R}$-coalgebra $\mathcal{C}$ simply a seminormed (resp.\ Banach) $\mathcal{R}$-coalgebra, dropping
the adjectives ``pre-differential'' and ``graded'' and relying on context and the calligraphic font to remind
us that the coalgebra is graded and equipped with a coderivation compatible with the derivation on $R$.
A \textbf{morphism of pre-differential graded seminormed} (resp.\ \textbf{Banach})
$\mathcal{R}-$\textbf{coalgebras} $f \colon (C,\mu) \rightarrow (D, \nu)$ is
a morphism $f \colon C \rightarrow D$ of the underlying graded seminormed (resp.\ Banach) $R$-coalgebras
which is also compatible with the coderivations in the sense that
$f \circ \mu = \nu \circ f$.

\subsection{Differential Graded Seminormed and Banach Algebras and Modules}
Starting with this subsection and for the remainder of the section,
we fix an odd element $\go \in \GG$ which is considered as part of the grading datum
(see \cref{subsec:grading-data}).
Differentials on objects will raise the degree of elements by $\go \in \GG$.
We extend the notion of differential graded modules and algebras from
\cref{subsec:differential-graded-algebras-modules} to the seminormed
and Banach settings.

\begin{dfn}
	A \textbf{differential graded seminormed} (resp.\ \textbf{Banach})
	$\mathbbm{k}$-\textbf{module} $\mathcal{M} = \left( M, d_M \right)$
	is a pre-differential graded seminormed (resp.\ Banach) $\mathbbm{k}$-module such that $d_M^2 = 0$
	(i.e., $d_M$ is a \textit{differential}).
\end{dfn}

Let us denote by $\DGSNMod[\mathbbm{k}]$ the full subcategory of $\PDGSNMod[\mathbbm{k}]$ consisting of differential graded seminormed $\mathbbm{k}$-modules and by $\DGBMod[\mathbbm{k}]$ the full subcategory of $\DGSNMod[\mathbbm{k}]$
consisting of differential graded Banach $\mathbbm{k}$-modules.
Since a differential graded seminormed (resp.\ Banach) $\mathbbm{k}$-module is in particular a
pre-differential graded seminormed (resp.\ Banach) $\mathbbm{k}$-module, one can apply all the constructions
described in \cref{sub:pre-differential-graded-seminormed-banach-k-modules}
to differential graded seminormed (resp.\ Banach) $\mathbbm{k}$-modules
and the resulting object will also be a differential graded seminormed (resp.\ Banach) $\mathbbm{k}$-module.

The category $\DGSNMod[\mathbbm{k}]$ (resp.\ $\DGBMod[\mathbbm{k}]$) is closed symmetric
monoidal, the symmetric monoidal structure being the one restricted from $\PDGSNMod[\mathbbm{k}]$
(resp.\ $\PDGBMod[\mathbbm{k}]$).

The \textbf{cohomology} $\cohom{\mathcal{M}}[]$
of a differential graded seminormed (or Banach) $\mathbbm{k}$-module $\mathcal{M} = \left( M, d_M \right)$
is defined to be the cohomology of the underlying differential graded $\mathbbm{k}$-module obtained by forgetting
the norm on $M$.

\begin{rem} \label{rem:cohomology-versions-seminormed-banach}
	There are finer versions of cohomology for differential-graded seminormed and Banach modules. Given a differential
	graded seminormed $\mathbbm{k}$-module $\mathcal{M} = \left( M, d_M \right)$,
	one can endow the cohomology $\cohom{\mathcal{M}}[]$ with a graded $\mathbbm{k}$-module seminorm
	by setting
	\begin{equation*}
		\nnorm[\eqcl{m}] \defeq \inf_{m' \in M^{g - \go}} \nnorm[m + dm']
	\end{equation*}
	for $m \in M^g$ such that $dm = 0$. This amounts to using the kernels and cokernels internal
	to the category $\SNMod[\mathbbm{k}]$. We will denote the resulting graded seminormed $\mathbbm{k}$-module
	by $\cohom{\mathcal{M}}[][][][\operatorname{SN}]$ and call $\cohom{\mathcal{M}}[][][][\operatorname{SN}]$
	the \textbf{seminormed cohomology}. This gives us a cohomology functor
	$\cohom{}[][][][\operatorname{SN}] \colon \DGSNMod[\mathbbm{k}] \rightarrow \GSNMod[\mathbbm{k}]$.

	Similarly, given a differential graded Banach $\mathbbm{k}$-module $\mathcal{M} = \left( M, d_M \right)$,
	one can use the kernels and cokernels internal to the category $\BMod[\mathbbm{k}]$ and define the
	\textbf{Banach cohomology} by
	\begin{equation*}
		\cohom{\mathcal{M}}[g][][][\operatorname{B}] \defeq
		\ker \left( d_M \colon M^g \rightarrow M^{g+\go} \right) /
		\overline{\Im \left( d_M \colon M^{g-\go} \rightarrow M^g \right)}
	\end{equation*}
	endowed with the natural norm. The resulting object $\cohom{\mathcal{M}}[][][][\operatorname{B}]$ is a
	graded \textit{Banach} $\mathbbm{k}$-module which is the graded separated completion of
	$\cohom{\mathcal{M}}[][][][\operatorname{SN}]$. This gives us a cohomology functor
	$\cohom{}[][][][\operatorname{B}] \colon \DGBMod[\mathbbm{k}] \rightarrow \BMod[\mathbbm{k}]$.

	We emphasize that in our work, the cohomology of a differential graded Banach $\mathbbm{k}$-module
	will always be taken to be the ``standard'' cohomology $\cohom{\mathcal{M}}[]$ obtained by ignoring
	the fact that $\mathcal{M}$ is Banach and taking the cohomology of the underlying differential
	graded $\mathbbm{k}$-module. We will never work with $\cohom{\mathcal{M}}[][][][\operatorname{SN}]$
	or $\cohom{\mathcal{M}}[][][][\operatorname{B}]$.
\end{rem}

\begin{dfn}
	A \textbf{differential graded seminormed} (resp.\ \textbf{Banach})
	$\mathbbm{k}$-\textbf{algebra} $\mathcal{R} = \left( R, d \right)$
	is a pre-differential graded seminormed (resp.\ Banach) $\mathbbm{k}$-algebra
	such that $d^2 = 0$ (i.e., the derivation $d$ is a bounded \textit{differential}).
\end{dfn}

We note that a differential graded seminormed (resp.\ Banach) $\mathbbm{k}$-algebra is the same thing as
an algebra object $\mathcal{R}$ of the monoidal category $\DGSNMod[\mathbbm{k}]$ (resp.\ $\DGBMod[\mathbbm{k}]$).
A \textbf{morphism} $f \colon (R,d_R) \rightarrow (S,d_S)$ \textbf{of differential graded seminormed} (resp.\ \textbf{Banach}) $\mathbbm{k}$-\textbf{algebras}
is the same thing as a morphism of pre-differential graded
seminormed (resp.\ Banach) $\mathbbm{k}$-algebras,
i.e., a degree zero contractive $\mathbbm{k}$-algebra morphism $f \colon R \rightarrow S$ compatible
with the differentials in the sense that $f \circ d_R = d_S \circ f$.

Let us fix a differential graded seminormed $\mathbbm{k}$-algebra $\mathcal{R} = (R,d_R)$.

\begin{dfn}
	A \textbf{differential graded seminormed} 	$\mathcal{R}$-\textbf{module} $\mathcal{M} = \left( M, d_M \right)$
	is a pre-differential graded seminormed $\mathcal{R}$-module such that $d_M^2 = 0$ (i.e.,
	$d_M$ is a \textbf{differential}). When $M$ is Banach, we call $\mathcal{M}$
	a \textbf{differential graded Banach} $\mathcal{R}$-\textbf{module}.
\end{dfn}

Let us denote by $\DGSNMod[\mathcal{R}]$ the full subcategory of $\PDGSNMod[\mathcal{R}]$ consisting of differential graded seminormed $\mathcal{R}$-modules and by $\DGBMod[\mathcal{R}]$ the full subcategory of $\DGSNMod[\mathcal{R}]$
consisting of differential graded Banach $\mathcal{R}$-modules.
Since a differential graded seminormed (resp.\ Banach) $\mathcal{R}$-module is in particular a pre-differential graded
seminormed (resp.\ Banach) $\mathcal{R}$-module, one can apply all the constructions described in
\cref{sub:pre-differential-graded-seminormed-banach-R-modules} to differential graded seminormed
(resp.\ Banach) $\mathcal{R}$-modules and the resulting object will also be differential graded
seminormed (resp.\ Banach) $\mathcal{R}$-module.

When $\mathcal{R}$ is a differential \textit{graded-commutative}
seminormed (resp.\ Banach) $\mathbbm{k}$-algebra, the category $\DGSNMod[\mathcal{R}]$ (resp.
$\DGBMod[\mathcal{R}]$) is closed symmetric monoidal with the symmetric monoidal and closed structures
being the ones restricted from $\PDGSNMod[\mathcal{R}]$ (resp.\ $\PDGBMod[\mathcal{R}]$).

The \textbf{cohomology} $\cohom{\mathcal{M}}[]$
of a differential graded seminormed (or Banach) $\mathcal{R}$-module
$\mathcal{M} = \left( M, d_M \right)$
is the cohomology of the underlying differential graded $\mathcal{R}$-module
obtained by forgetting the norms on $M$ and $R$. While it is possible to
work with alternative definitions (see \cref{rem:cohomology-versions-seminormed-banach}),
we emphasize that when taking the cohomology of a differential graded seminormed or Banach
$\mathcal{R}$-module, we ignore the norms on $M$ and $R$ and take the cohomology of the
underlying differential graded $R$-module. The result is a graded $\cohom{\mathcal{R}}[]$-module over the
graded $\mathbbm{k}$-algebra $\cohom{\mathcal{R}}[]$.

\subsection{Internal Point of View} \label{sub:internal-point-of-view}
\begin{figure}[htb]
	\centering
	\begin{tikzcd}
		\GBRing && \GSNRing && \GRing \\
		& {\cancel{\circlearrowleft}} && \circlearrowleft \\
		\BRing && \SNRing && \Ring
		\arrow[from=1-1, to=1-3]
		\arrow[from=1-3, to=1-5]
		\arrow["{R \mapsto \coplus_{g \in \GG} R^g}"', from=1-1, to=3-1]
		\arrow["{R \mapsto \oplus_{g \in \GG}^{\SNAb} R^g}"', from=1-3, to=3-3]
		\arrow["{R \mapsto \oplus_{g \in \GG} R^g}", from=1-5, to=3-5]
		\arrow[from=3-1, to=3-3]
		\arrow[from=3-3, to=3-5]
	\end{tikzcd}
	\caption{Underlying ring functors.}
	\label{fig:underlying-ring-func}
\end{figure}
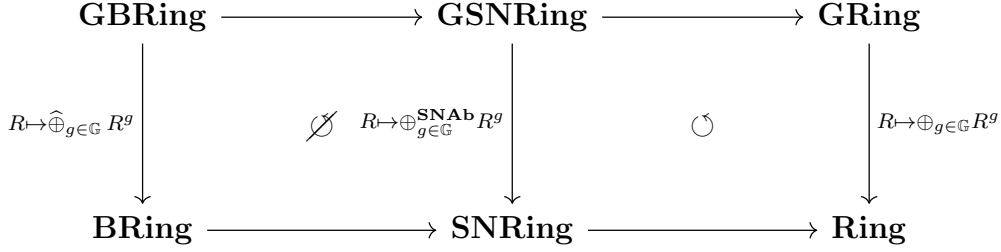

From the point of view we adopted in \cref{dfn:graded-k-algebra},
a graded ring $R = \left( R^g \right)_{g \in \GG}$ is not really a ring.
We can only talk about homogeneous elements and cannot add elements of different
degrees. However, it is often convenient to associate to a graded ring $R$
an honest ring $R_{\oplus} \defeq \oplus_{g \in \GG} R^g$ which we call the
\textbf{underlying ring}. Similarly, we can associate to each
graded seminormed ring $R$ an \textbf{underlying seminormed ring} $R_{\oplus} \defeq \oplus_{g \in \GG}^{\SNAb} R^g$
and to each graded Banach ring $R$ an \textbf{underlying Banach ring} $R_{\coplus} \defeq \coplus_{g \in \GG} R^g$
and obtain the vertical functors depicted in \cref{fig:underlying-ring-func}.

\begin{ex}[\textbf{Warning}] \label{ex:graded-polynomial-power-series-rings}
	Given a graded Banach ring $R$, the underlying (ungraded) Banach ring $R_{\coplus}$ will in general
	be different from the underlying seminormed ring $R_{\oplus}$ since the left square in
	\cref{fig:underlying-ring-func} does not commute. This has some possibly
	confusing consequences.

	Consider for example the polynomial ring $\ZZ[x]$ and the power series ring $\pows{\ZZ}[x]$
	endowed with the norms $\nnorm[p] = 2^{-\ord p}$
	where the order of a power series $p = \sum_{n \geq 0} a_n x^n$ is given by
	\begin{equation*}
		\ord p \defeq
		\begin{cases}
			\min \Set{n \geq 0}[a_n \neq 0] & p \neq 0, \\
			+\infty                         & p = 0.
		\end{cases}
	\end{equation*}
	The ring $\ZZ[x]$ is normed but not complete and has $\pows{\ZZ}[x]$ as its completion.

	Now, let's say we want to consider $\ZZ[x]$ and $\pows{\ZZ}[x]$ as graded seminormed/Banach rings
	with $\degb{x} = 1$ and $\nnorm[x] = \frac{1}{2}$.
	To do that, we let $R = \left( R^n \right)_{n \in \ZZ}$ be given by
	\begin{equation*}
		R^n \defeq
		\begin{cases}
			\ZZ \cdot x^n & n \geq 0, \\
			0             & n < 0.
		\end{cases}
	\end{equation*}
	Each $R^n$ is endowed with the norm $\nnorm[k \cdot x^n] = 2^{-n}$ where $0 \neq k \in \ZZ$.
	Since each $R^n$ is complete, $R$ is both a graded seminormed ring and
	a graded Banach ring. The underlying seminormed ring $R_{\oplus}$ is isometrically
	isomorphic to the polynomial ring $\ZZ[x]$ while the underlying Banach ring $R_{\coplus}$
	is isometrically isomorphic to the power series ring $\pows{\ZZ}[x]$.

	Hence, from the point of view of \cref{dfn:graded-seminormed-banach-k-algebra},
	there is no difference between the ``graded polynomial ring'' $\ZZ[x]$ and the ``graded power series ring''
	$\pows{\ZZ}[x]$. In fact, both are not really rings but defined by the same object
	(an indexed family of Banach groups). By applying different
	functors to the same object, we obtain $\ZZ[x]$ or $\pows{\ZZ}[x]$.
\end{ex}

The vertical functors depicted in \cref{fig:underlying-ring-func} give us an underlying ring but
forget the grading. To remember the grading, we introduce the following definition:

\begin{dfn} \label{def:semi-norm-complete-grading}
	\begin{enumerate}
		\item[]
		\item Let $\left( A, \nnorm \right)$ be a seminormed group. A $\GG$-\textbf{grading} on
		      $A$ is a family $\left( A^g \right)_{g \in \GG}$ of subgroups of $A$ such that
		      $A = \oplus_{g \in \GG}^{\SNAb} A^g$ in $\SNAb$.
		\item Let $\left( A, \nnorm \right)$ be a Banach group. A $\widehat{\GG}$-\textbf{grading} on $A$
		      is a family $\left( A^g \right)_{g \in \GG}$ of \textit{closed} subgroups of $A$ such that
		      $A = \coplus_{g \in \GG} A^g$ in $\BAb$.\footnote{The $\widehat{\quad}$ in $\widehat{\GG}$-grading
			      does not denote a completion of $\GG$. It was chosen to remind the reader that we work in
			      $\BAb$ and not in $\SNAb$ so that we can differentiate it from the notion of $\GG$-grading.}
	\end{enumerate}
\end{dfn}
Let us emphasize that in \cref{def:semi-norm-complete-grading} we require that $A$ is the
\textit{internal} direct sum in the appropriate category.
For example, in $\SNAb$ we require that the canonical map $\oplus_{g \in \GG} A^g \rightarrow A$ constructed from the
inclusions $\left( A^g, \nnorm \! \restriction_{A^g} \right) \hookrightarrow \left( A, \nnorm \right)$
is an isometric isomorphism and similarly in $\BAb$ (with $\oplus$ replaced by $\coplus$).
It is worthwhile to be explicit about the meaning of $\GG$-grading and $\widehat{\GG}$-grading:
\begin{enumerate}
	\item Let $\left( A, \nnorm \right)$ be a seminormed group. The group $A$ is $\GG$-graded
	      by a family $\left( A^g \right)_{g \in \GG}$ of subgroups if and only if:
	      \begin{enumerate}
		      \item Every element $a \in A$ can be represented uniquely as a \textit{finite sum}
		            $a = \sum_{g \in \GG} a^g$ where each $a^g \in A^g$.
		            That is, $A$ is the algebraic direct sum of $\left( A^g \right)_{g \in \GG}$.
		      \item We have $\nnorm[a] = \max_{g \in \GG} \nnorm[a^g]$.
	      \end{enumerate}
	\item Let $\left( A, \nnorm \right)$ be a Banach group. The group $A$ is $\widehat{\GG}$-graded
	      by a family $\left( A^g \right)_{g \in \GG}$ of closed subgroups if and only if:
	      \begin{enumerate}
		      \item Every element $a \in A$ can be represented uniquely as an unordered,
		            \textit{possibly infinite, convergent} sum $a = \sum_{g \in \GG} a^g$ where each
		            $a^g \in A^g$ and $a^g \to 0$. The sum converges since $A$ is Banach.
		      \item We have $\nnorm[a] = \max_{g \in \GG} \nnorm[a^g]$.
	      \end{enumerate}
\end{enumerate}

Given a graded seminormed ring $R = \left( R^g \right)_{g \in \GG}$ in the sense
of \cref{dfn:graded-seminormed-banach-k-algebra}, one can assign to it the underlying seminormed ring
$R_{\oplus} = \oplus_{g \in \GG}^{\SNAb} R^g$, given by the external direct sum, together with a canonical
$\GG$-grading. Hence, we can think of a graded seminormed ring in two equivalent ways:
\begin{enumerate}
	\item Externally, as an indexed family $\left( R^g \right)_{g \in \GG}$ of seminormed groups
	      together with multiplication maps as in \cref{dfn:graded-seminormed-banach-k-algebra}.
	\item Internally, as a pair $( R, \left( R^g \right)_{g \in \GG} )$ where $R$ is a seminormed ring
	      and $\left( R^g \right)_{g \in \GG}$ is a fixed $\GG$-grading on $R = \oplus_{g \in \GG}^{\SNAb} R^g$
	      which is compatible with the multiplication and unit in the sense that
	      $R^g \cdot R^h \subseteq R^{g+h}$ and $1_R \in R^0$.
\end{enumerate}
Similarly, if $R$ is a graded Banach ring, one can assign to it the underlying
Banach ring $R_{\coplus} = \coplus_{g \in \GG} R^g$ together with a canonical $\widehat{\GG}$-grading.
Hence, we can think of a graded Banach ring in two equivalent ways:
\begin{enumerate}
	\item Externally, as an indexed family $\left( R^g \right)_{g \in \GG}$ of Banach groups
	      together with multiplication maps as in \cref{dfn:graded-seminormed-banach-k-algebra}.
	\item Internally, as a pair $( R, \left( R^g \right)_{g \in \GG} )$ where $R$ is a Banach ring
	      and $\left( R^g \right)_{g \in \GG}$ is a fixed $\widehat{\GG}$-grading on
	      $R = \coplus_{g \in \GG} R^g$ compatible with the multiplication and unit.
\end{enumerate}

Analogously, one can think of graded seminormed (resp.\ Banach) modules as ungraded seminormed (resp.\ Banach) modules
equipped with a $\GG$-grading (resp.\ $\widehat{\GG}$-grading) compatible with the module action.

\begin{rem}
	Note that if one works with \cref{dfn:graded-seminormed-banach-k-algebra}, then a graded Banach ring
	is automatically a graded seminormed ring (and also a graded ring)
	and the category $\GBRing$ forms a full subcategory of $\GSNRing$.
	However, if one adopts the internal point of view discussed above then a graded Banach ring
	$( R, \left( R^g \right)_{g \in \GG} )$ is not a graded seminormed ring nor a graded ring in the algebraic sense
	since we don't have $R = \oplus_{g \in \GG} R^g$. In this case, the natural ``forgetful'' functor
	from $\GBRing$ to $\GSNRing$ sends a graded Banach ring $( R, \left( R^g \right)_{g \in \GG} )$
	to the graded seminormed ring $( \oplus_{g \in \GG}^{\SNAb} R^g, \left( R^g \right)_{g \in \GG} )$.
	For example, if we think of the power series ring $\pows{\ZZ}[x]$ as a graded Banach ring with $\degb{x} = 1$ and
	$\nnorm[x] = \frac{1}{2}$ as in \cref{ex:graded-polynomial-power-series-rings},
	then the ``forgetful'' functor sends $\pows{\ZZ}[x]$ to the polynomial ring $\ZZ[x]$.
\end{rem}

\section{The Formal Tensor Coalgebra} \label{sec:formal-tensor-coalgebra}

In this section, we introduce the formal tensor coalgebra and study its properties. The formal
tensor coalgebra is constructed analogously to the standard tensor coalgebra,
working instead in the category of non-Archimedean graded Banach modules,
using complete direct sums and tensor products.

We start in \cref{subsec:tensor-coalgebra} by discussing the standard tensor module and (co)algebra,
reviewing its properties and setting up notation used in the rest of this work.
In \cref{subsec:formal-tensor-coalgebra}, we define the formal tensor coalgebra $\tensf{V}$ of a graded Banach
$R$-module $V$ over a graded Banach $\mathbbm{k}$-algebra $R$ and study its properties.
We classify grouplike elements of $\tensf{V}$ as those which can be expressed as the exponential of topologically nilpotent
elements of $V^0$.
We classify morphisms between formal tensor coalgebras in terms of their corestrictions and show
that such morphisms can have a ``change of connection'' component, as long as it is topologically nilpotent.
We also classify coderivations $\mu$ on $\tensf{V}$ in terms of their corestrictions, both
when $\mu$ is $R$-linear and when $\mu$ is a coderivation over a derivation $d$ of $R$.
Finally, in \cref{sec:formal-tensor-coalgebra-different-ground}, we discuss scalar restriction
and extension for formal tensor coalgebras equipped with a generalized coderivation, used
to facilitate working with morphisms of formal tensor coalgebras over different ground algebras,
and to later define scalar restriction and extension of graded Banach $\Ainf$-algebras.

For the rest of this section, we fix a commutative ground ring $\mathbbm{k}$ and a grading datum
$\left( \GG, \braidop \right)$ (see \cref{subsec:grading-data}).

\subsection{The Tensor Coalgebra} \label{subsec:tensor-coalgebra}
Let $R$ be a graded-commutative $\mathbbm{k}$-algebra and let $V$ be a graded $R$-module. In what follows,
we write $\otimes$ for $\otimes_R$. The \textbf{tensor module on} $V$ is defined to be the graded
$R$-module
\begin{equation*}
	\tens{V} = \tens{V}[R] \defeq \bigoplus_{i=0}^{\infty} V^{\otimes i}.
\end{equation*}
The terminology tensor \textit{module} will prove useful as the tensor module $\tens{V}$ can be equipped
with several different algebraic structures. Elements of $\tens{V}$
carry two natural gradings,
the \textbf{weight grading} and the \textbf{degree grading}.
An elementary tensor $l = v_1 \otimes \dots \otimes v_k$ in $V^{\otimes k}$ has weight
$\weight{l} = k$ and degree $\degb{v_1} + \dots + \degb{v_k}$, where $\degb{v_i} \in \GG$ is the degree of $v_i$ in $V$.
We will often use both gradings but whenever the grading is not specified explicitly, we will assume that we work with the degree grading.

The \textbf{tensor coalgebra} on $V$ is given
by the tensor module $\tens{V}$ together with the coproduct
$\Delta \colon \tens{V} \rightarrow \tens{V} \ootimes \tens{V}$, called the \textbf{deconcatenation coproduct},
which acts on elementary tensors by
\begin{equation}
	\begin{aligned}
		\Delta \left( v_1 \otimes \dots \otimes v_n \right) \defeq{} &
		1 \ootimes \left( v_1 \otimes \dots \otimes v_n \right) +
		v_1 \ootimes \left( v_2 \otimes \dots \otimes v_n \right) + \dots
		\\
		                                                             & + \left( v_1 \otimes \dots \otimes v_n \right) \ootimes 1.
	\end{aligned}
	\label{eq:comult-tensor-coalgebra}
\end{equation}
The symbol $\ootimes$ used above is a synonym for the tensor product $\otimes$ and will be used
whenever we need to distinguish the ``external'' tensor product used in the definition of a coalgebra from
the ``internal'' tensor product used in describing elements $v_1 \otimes \dots \otimes v_n$ of $\tens{V}$.
The tensor coalgebra $\tens{V}$ is counital with the counit $\varepsilon \colon \tens{V} \twoheadrightarrow R$
given by the natural projection onto $V^{\otimes 0} = R$
and has a coaugmentation given by the inclusion $R \hookrightarrow \tens{V}$.
The kernel of the counit $\varepsilon$ is called the \textbf{reduced tensor module} and is denoted by
\begin{equation*}
	\tensr{V} = \tensr{V}[R] \defeq \bigoplus_{i=1}^{\infty} V^{\otimes i}.
\end{equation*}
As in any coalgebra, by using the coproduct and the counit one can define a sequence of iterated coproducts
$\Delta^n \colon \tens{V} \rightarrow \tens{V}^{\ootimes n}$ for $n \geq 0$ with $\Delta^2 = \Delta$
(see \cref{subsec:coalg-in-monoidal-cat} for the definitions).

Given an elementary tensor $l = v_1 \otimes \dots \otimes v_k \in \tens{V}$, we will often think informally of $l$ as
representing a list $(v_1,\dots,v_k)$ of $k$ elements from $V$.
When $k = 0$, the empty tensor product $v_1 \otimes \dots \otimes v_k$ is interpreted as
$1 = 1_R \in \tens{V}$ and thought of as the empty list. With this interpretation in mind, the action
of the iterated coproduct $\Delta^n \left( l \right)$ is given by the \textit{sum} over all possible splittings
of the list $l$ into $n$ consecutive, possibly empty, lists. For example, if $l = v_1 \otimes v_2$ we have
\begin{align*}
	\Delta^2 \left( l \right) ={} & 1 \ootimes \left( v_1 \otimes v_2 \right) + v_1 \ootimes v_2 +
	\left( v_1 \otimes v_2 \right) \ootimes 1,
	\\
	\Delta^3 \left( l \right) ={} & 1 \ootimes 1 \ootimes \left( v_1 \otimes v_2 \right)
	+ 1 \ootimes v_1 \ootimes v_2 + 1 \ootimes \left( v_1 \otimes v_2 \right) \ootimes 1
	\\
	                              & +
	v_1 \ootimes 1 \ootimes v_2 + v_1 \ootimes v_2 \ootimes 1 +
	\left( v_1 \otimes v_2 \right) \ootimes 1 \ootimes 1.
\end{align*}

The iterated coproducts are related to the natural projections on $\tens{V}$ as follows.
Let us denote by $\pi_k \colon \tens{V} \twoheadrightarrow V^{\otimes k}$ the canonical projections
and denote by
\begin{equation*}
	D_{V}^{k_1, \dots, k_n} \colon V^{\otimes \left( k_1 + \dots + k_n \right)} \rightarrow
	V^{\otimes k_1} \ootimes \dots \ootimes V^{\otimes k_n}
\end{equation*}
the natural associativity isomorphisms. Then for all $n \geq 0$ and $k_1,\dots,k_n \geq 0$, we have the identity
\begin{equation*}
	\left( \pi_{k_1} \ootimes \dots \ootimes \pi_{k_n} \right) \circ \Delta^n =
	D^{k_1,\dots,k_n}_V \circ \pi_{k_1 + \dots + k_n}.
\end{equation*}
When $n = 2$, the identities reduce to the definition \eqref{eq:comult-tensor-coalgebra} of the coproduct.

In addition to the coalgebra structure, the tensor module $\tens{V}$ also carries an algebra
structure which we use. The product $m \colon \tens{V} \ootimes \tens{V} \rightarrow \tens{V}$ is given on elementary tensors by
\begin{equation}
	m \left( \left( v_1 \otimes \dots \otimes v_k \right) \ootimes \left( u_1 \otimes \dots \otimes u_l \right) \right)
	\defeq
	v_1 \otimes \dots \otimes v_k \otimes u_1 \otimes \dots \otimes u_l.
	\label{eq:mult-tensor-algebra}
\end{equation}
As usual, we follow the convention that the empty tensor product $v_1 \otimes \dots \otimes v_k$ of elements
(i.e., when $k = 0$) is given by $1_R$ and always identify $1_R \otimes v$ and $v \otimes 1_R$ with $v$.
The tensor algebra is unital with unit map given by the natural inclusion $R \hookrightarrow \tens{V}$ and augmented
with augmentation given by the natural projection $\tens{V} \twoheadrightarrow R$.
With respect to the multiplication defined above, the tensor algebra $\tens{V}$ becomes
the free graded $R$-algebra on $V$. Both the product and coproduct on $\tens{V}$ are homogeneous of degree zero with respect to both
the weight grading and the degree grading.\footnote{Note that the multiplication and comultiplication on $\tens{V}$
	are not compatible in the sense that endowing $\tens{V}$ with both structures does not make $\tens{V}$ into a bialgebra.}

\subsubsection{Notation} \label{sec:splitting-tensor-coalgebra-notation}

In our work, given an elementary tensor $l = v_1 \otimes \dots \otimes v_k$ which we think of as a list,
we will often need to split $l$ into several consecutive lists, rotate them and apply maps to certain
parts of the splitting. To do that, it will be convenient to introduce the following pieces of notation:
\begin{enumerate}
	\item \textbf{Splitting}. Given $n \in \NZ$, we will denote by
	      $l_{(1)} \ootimes \dots \ootimes l_{(n)}$ the element of $\tens{V}^{\ootimes n}$ given by the
	      iterated coproduct $\Delta^n \left( l \right)$.
	      This notation is often called \textbf{Sweedler's notation} and allows us to write succinctly
	      the coproduct $\Delta$ as $l \mapsto l_{(1)} \ootimes l_{(2)}$.
	      After splitting a list once, we can split each of the factors and write expressions
	      such as
	      \begin{equation*}
		      l_{(1)} \ootimes l_{(2)} \ootimes l_{(3)} = l_{(11)} \ootimes l_{(12)} \ootimes l_{(2)} =
		      l_{(1)} \ootimes l_{(21)} \ootimes l_{(22)}
	      \end{equation*}
	      which are equal. More rigorously, the expression $l_{(11)} \ootimes l_{(12)} \ootimes l_{(2)}$ denotes
	      $\left( \left( \Delta^2 \ootimes \id \right) \circ \Delta^2 \right) \left( l \right)$, the expression
	      $l_{(1)} \ootimes l_{(21)} \ootimes l_{(22)}$ denotes
	      $\left( \left( \id \ootimes \Delta^2 \right) \circ \Delta^2 \right) \left( l \right)$
	      and since $\Delta$ is coassociative we identify both expressions with
	      $\Delta^3 \left( l \right) = l_{(1)} \ootimes l_{(2)} \ootimes l_{(3)}$.
	\item \textbf{Application}. Given a graded map $\psi \colon \tens{V} \rightharpoonup V$
	      (or $\psi \colon \tens{V} \rightharpoonup \tens{V}$), we will denote by
	      \begin{equation}
		      l_{(1)} \otimes \dots \otimes l_{(i-1)} \otimes \psi \left( l_{(i)} \right)
		      \otimes l_{(i+1)} \otimes \dots \otimes l_{(n)} \label{eq:application-general}
	      \end{equation}
	      the element of $\tens{V}$ which is the sum of all elements which are
	      obtained by splitting $l$ into $n$ consecutive lists, applying $\psi$ to the $i$-th list
	      and taking the multiplication \textbf{inside} $\tens{V}$ (hence, we use $\otimes$ and not $\ootimes$).
	      For example, if $\psi \colon \tens{V} \rightharpoonup V$ and $l = v_1 \otimes v_2$ we have
	      \begin{align*}
		      l_{(1)} \otimes \psi \left( l_{(2)} \right) & =
		      1 \otimes \psi \left( v_1 \otimes v_2 \right) + v_1 \otimes \psi \left( v_2 \right) +
		      \left( v_1 \otimes v_2 \right) \otimes \psi \left( 1 \right) \\
		                                                  & =
		      \underbrace{\psi(v_1 \otimes v_2)}_{\in V} +
		      \underbrace{v_1 \otimes \psi \left( v_2 \right)}_{\in V^{\otimes 2}} +
		      \underbrace{v_1 \otimes v_2 \otimes \psi(1)}_{\in V^{\otimes 3}}
		      \in \tens{V}.
	      \end{align*}
	      Note that we use only the ``internal'' tensor product and identify
	      $1 \otimes \psi \left( v_1 \otimes v_2 \right)$ with $\psi \left( v_1 \otimes v_2 \right)$.
	      We will also allow our notation to include ``sign factors'' as in the following expression:
	      \begin{equation}
		      \begin{aligned}
			      (-1)^{\braidd{\psi}{l_{(1)}}} l_{(1)} \otimes \psi \left( l_{(2)} \right) ={} &
			      \psi \left( v_1 \otimes v_2 \right) +
			           (-1)^{\braidd{\psi}{v_1}} v_1 \otimes \psi \left( v_2 \right)
			      \\
			                                                                                    & + (-1)^{\braid{\degb{\psi}}{\degb{v_1} + \degb{v_2}}}
			      v_1 \otimes v_2 \otimes \psi \left( 1 \right).
		      \end{aligned} \label{eq:application-with-sign}
	      \end{equation}
	\item \textbf{Rotation}. In addition to the expressions described above, we will often need to
	      rotate some parts of the splitting before applying maps. For example, given $x \in V, l \in \tens{V}$
	      and $\psi \colon \tens{V} \rightharpoonup V$, we will write expressions such as
	      \begin{equation}
		      (-1)^{\braid{\degb{l_{(3)}}}{\degb{x} + \degb{l_{(1)}} + \degb{l_{(2)}}}}
		      \psi \left( l_{(3)} \otimes x \otimes l_{(1)} \right) \otimes l_{(2)}
		      \label{eq:rotation-with-signs}
	      \end{equation}
	      In the expression above, we start with $x \otimes l$, split $l$ into three consecutive sublists
	      $l_{(1)},l_{(2)},l_{(3)}$, rotate $l_{(3)}$ across $x,l_{(1)},l_{(2)}$ to the beginning of the expression,
	      apply $\psi$ to $l_{(3)} \otimes x \otimes l_{(1)}$ (considered as an element of $\tens{V}$) and multiply
	      the result with $l_{(2)}$ \textbf{inside} $\tens{V}$ (hence, we use $\otimes$ and not $\ootimes$ both
	      inside $\psi$ and outside). The resulting element is an element of $\tens{V}$ (even $\tensr{V}$).
	      For example, if $l = v_1 \otimes v_2$ we have
	      \begin{equation*}
		      \begin{aligned}
			      (-1)^{\braid{\degb{l_{(3)}}}{\degb{x} + \degb{l_{(1)}} + \degb{l_{(2)}}}}
			      \psi \left( l_{(3)} \otimes x \otimes l_{(1)} \right) \otimes l_{(2)} ={} &
			      (-1)^{\braid{\degb{v_1} + \degb{v_2}}{\degb{x}}}
			      \psi \left( v_1 \otimes v_2 \otimes x \right)
			      \\
			                                                                                & +
			                                                                                  (-1)^{\braid{\degb{v_2}}{\degb{x} + \degb{v_1}}}
			      \psi \left( v_2 \otimes x \right) \otimes v_1
			      \\
			                                                                                & +
			      \psi \left( x \right) \otimes v_1 \otimes v_2
			      \\
			                                                                                & +
			                                                                                  (-1)^{\braid{\degb{v_2}}{\degb{x} + \degb{v_1}}}
			      \psi \left( v_2 \otimes x \otimes v_1 \right)
			      \\
			                                                                                & +
			      \psi \left( x \otimes v_1 \right) \otimes v_2
			      \\
			                                                                                & +
			      \psi \left( x \otimes v_1 \otimes v_2 \right).
		      \end{aligned}
	      \end{equation*}
\end{enumerate}

\begin{rem}
	Even though we haven't rigorously covered all possible options, we believe it is clear at this point how
	to interpret all the expressions written in the work. One warning is in order: Since we work
	over a graded-commutative ground algebra $R$, we should be careful with signs in order to make sure
	the expressions we write are in fact well-defined. For example, if $\psi \colon \tens{V} \rightharpoonup V$
	is a graded $R$-linear map which is not of degree zero then the expression in
	\cref{eq:application-general} is not well-defined without a sign factor
	(see \cref{foot:sign-necessary-tensor-product}). However, the expression
	in \cref{eq:application-with-sign} is well-defined and can be written as
	\begin{equation*}
		(-1)^{\braidd{\psi}{l_{(1)}}} l_{(1)} \otimes \psi \left( l_{(2)} \right) =
		\left( m_{\bullet,1} \circ \left( \id \ootimes \psi \right) \circ \Delta^2 \right)
		\left( l \right)
	\end{equation*}
	where $m_{\bullet,1} \colon \tens{V} \ootimes V \rightarrow \tens{V}$ is the natural
	multiplication map and the sign factor comes from the definition of the tensor product of
	graded maps (see \cref{eq:tensor-product-graded-R-linear-maps}).
	Similarly, the expression in \cref{eq:rotation-with-signs} is well-defined and can be written as
	\begin{equation*}
		\left(
		m_{1,\bullet} \circ \left( \psi \ootimes \id \right) \circ \left( m_{\bullet,1,\bullet} \ootimes \id \right)
		\circ \sigma_{1,\bullet,\bullet,\bullet} \circ \left( \id \ootimes \Delta^3 \right)
		\right) \left( x \otimes l \right)
	\end{equation*}
	where the map $\sigma_{1,\bullet,\bullet,\bullet} \colon V \ootimes \tens{V} \ootimes \tens{V} \ootimes \tens{V}
		\rightarrow \tens{V} \ootimes V \ootimes \tens{V} \ootimes \tens{V}$
	is the natural symmetry isomorphism induced by the permutation which rotates the last factor to the
	first place and involves signs (see \cref{eq:symmetry-graded-R-modules}), and the maps
	$m_{1,\bullet} \colon V \ootimes \tens{V} \rightarrow \tens{V}$ and
	$m_{\bullet,1,\bullet} \colon \tens{V} \ootimes V \ootimes \tens{V} \rightarrow \tens{V}$
	are the natural multiplication maps.

	All the expressions we write in this work using the notation above can always be expressed
	in some complicated way using a sequence of operations involving the coproducts, tensor product
	of maps, symmetry isomorphisms, and various multiplications, and
	hence will always be well-defined.
\end{rem}

\subsection{The Formal Tensor Coalgebra} \label{subsec:formal-tensor-coalgebra}
In what follows, let $\mathbbm{k}$ be a commutative non-Archimedean Banach ring (for example, a field
with trivial norm).
Let $R$ be a graded-commutative Banach $\mathbbm{k}$-algebra and let $V$ be a graded Banach $R$-module.
The construction of the tensor module and (co)algebra described in the previous section for
graded $R$-modules can be generalized and done internally in any nice enough category.\footnote{To construct
	$\tens{V}$ and endow it with an associative product and coassociative coproduct, it is enough to work
	in an additive monoidal category with countable coproducts which commute with the tensor product in each variable.
	This holds in particular in any cocomplete closed symmetric additive monoidal category.}
By thinking of $V$ as an object in $\GSNMod[R]$ (ignoring the fact that $V$ is Banach) or in $\GBMod[R]$,
we can construct the tensor (co)algebra in each of the categories. We describe explicitly the resulting
objects.

Let us denote temporarily by $\tens{V, \nnorm_{V}}$ the result of performing
the construction of the tensor module in the category $\GSNMod[R]$.
Since the direct sum, tensor product, and unit in $\GSNMod[R]$ coincide on the level of graded $R$-modules
with the direct sum, tensor product, and unit in $\GMod[R]$, we see that the resulting
seminormed tensor module coincides with the tensor module $\tens{V}$ as a
graded $R$-module. Hence, $\tens{V, \nnorm_V}$ is given by $( \tens{V}, \nnorm_{\tens{V}} )$ where
the seminorm $\nnorm_{\tens{V}}$ on $\tens{V}$ is induced by the norm $\nnorm_V$ on $V$.
Similarly, the seminormed tensor (co)algebra object on $V$ in $\GSNMod[R]$ coincides
with the tensor (co)algebra on $V$ in $\GMod[R]$ as a graded $R$-(co)algebra
and, with respect to the seminorm $\nnorm_{\tens{V}}$, all the structure maps of $\tens{V}$ (the (co)product, (co)unit and (co)augmentation) become contractive.

Next, consider the construction of the tensor module in the category $\GBMod[R]$. The resulting
object is the graded Banach $R$-module given by
\begin{equation*}
	\tensf{V} \defeq \extrawidehat{\bigoplus}_{i=0}^{\infty} V^{\cotimes i}
	= \extrawidehat{\bigoplus}_{i=0}^{\infty} V^{\cotimes_R i}
\end{equation*}
and will be called the \textbf{formal tensor module on} $V$.
Elements $x \in \tensf{V}$ have a unique representation $x = \sum_{i=0}^{\infty} x_i$ as an infinite convergent
sum where each $x_i \in V^{\cotimes i}$ and $\nnorm[x_i] \to 0$.
Even though general elements of $V^{\cotimes i}$ cannot be described explicitly as
finite sums of elementary tensors, the algebraic tensor product $V^{\otimes i}$
is dense in $V^{\cotimes i}$. Consequently, any bounded $R$-linear map from
$\tensf{V}$ to a graded Banach $R$-module is uniquely determined by its action on
elementary tensors $v_1 \cotimes \dots \cotimes v_i$
(see the discussion in \cref{sub:category-GBMod-R}).
For this reason, we will frequently define such maps, and write their
explicit formulas, solely on elementary tensors.

The formal tensor module $\tensf{V}$ carries also the structure of a graded Banach $R$-(co)algebra.
The multiplication $m$ and comultiplication $\Delta$ on $\tensf{V}$ are given on elementary tensors by the same
formulas as for $\tens{V}$ (\cref{eq:comult-tensor-coalgebra,eq:mult-tensor-algebra}),
with $\otimes$ replaced by $\cotimes$. When considering $\tensf{V}$ as a graded Banach $R$-(co)algebra,
we call $\tensf{V}$ the \textbf{formal tensor (co)algebra on} $V$.\footnote{A more natural name
	would be the \textbf{Banach tensor (co)algebra}, i.e., a tensor (co)algebra in the category of graded
	Banach $R$-modules, but we have chosen the terminology appearing in \cite{Cho2012,Fukaya2009}.}
We note that all the properties and notation conventions
for $\tens{V}$ appearing in \cref{subsec:tensor-coalgebra} work just as well for $\tensf{V}$ with
the internal tensor product $\otimes$ replaced by the internal complete tensor product $\cotimes$,
and the external tensor product $\ootimes$ replaced by the external complete tensor product $\cootimes$.
In particular, if we denote by
$\pi_k \colon \tensf{V} \twoheadrightarrow V^{\cotimes k}$ the canonical projections
and by
\begin{equation*}
	D = D_{V}^{k_1, \dots, k_n} \colon V^{\cotimes \left( k_1 + \dots + k_n \right)} \rightarrow
	V^{\cotimes k_1} \cootimes \dots \cootimes V^{\cotimes k_n}
\end{equation*}
the natural associativity isomorphisms, then for all $n \geq 0$ and $k_1,\dots,k_n \geq 0$ we have the identity
\begin{equation} \label{eq:projection-related-to-coproduct-general}
	\left( \pi_{k_1} \cootimes \dots \cootimes \pi_{k_n} \right) \circ \Delta^n =
	D^{k_1,\dots,k_n}_V \circ \pi_{k_1 + \dots + k_n}.
\end{equation}
When $k_1 = \dots = k_n = 1$, the isomorphism $D$ can be chosen to be the identity\footnote{The
	isomorphism $D = D_{V}^{1,\dots,1} \colon V^{\cotimes n} \rightarrow V^{\cootimes n}$ comes
	from identifying the internal tensor product $V^{\cotimes n} \hookrightarrow \tensf{V}$
	with the external tensor product $V^{\cootimes n} \hookrightarrow \tensf{V}^{\cootimes n}$
	in which the output of $\pi_1^{\cootimes n} \circ \Delta^n$ lands. For $n \geq 0$, one
	can choose to bracket the codomain $\tensf{V}^{\cootimes n}$ of $\Delta^n$ in exactly the same
	way as the internal component $V^{\cotimes n}$ is bracketed and take $D$ to be the identity.
},
yielding the relations
\begin{equation} \label{eq:projection-related-to-coproduct}
	\left( \pi_1 \right)^{\cotimes n} \circ \Delta^n = \pi_n
\end{equation}
which hold for $n \geq 0$.

An alternative and useful way of thinking about $\tensf{V}$ is as follows: Since
the completion functor $\wedge \colon \GSNMod[R] \rightarrow \GBMod[R]$ is strong monoidal, the
completion of a graded seminormed $R$-(co)algebra has a natural structure of a graded Banach
$R$-(co)algebra. In particular, we can start with $\tens{V}$ endowed with the seminorm
$\nnorm_{\tens{V}}$ and complete it to obtain the graded Banach $R$-module
\begin{equation*}
	\widehat{\tens{V}} = \extrawidehat{\bigoplus_{i=0}^{\infty} V^{\otimes i}} \cong
	\extrawidehat{\bigoplus}_{i=0}^{\infty} \widehat{V^{\otimes i}} \cong
	\extrawidehat{\bigoplus}_{i=0}^{\infty} V^{\cotimes i} = \tensf{V}.
\end{equation*}
Under the natural identifications above, the resulting object $\widehat{\tens{V}}$ is
naturally isomorphic as a graded Banach $R$-module and a graded Banach $R$-(co)algebra to $\tensf{V}$.

\begin{ex} \label{ex:tensor-coalgebra-as-a-formal-tensor-coalgebra}
	Assume that $R$ is a graded $\mathbbm{k}$-algebra endowed with the trivial norm
	and let $V$ be a graded $R$-module also endowed with the trivial norm. In this case,
	the induced norm on $\tens{V}$ is also trivial and so $\tens{V}$ is Banach and
	$\tensf{V} = \widehat{\tens{V}} = \tens{V}$ as a graded Banach $R$-coalgebra. Note that morphisms
	of $\tens{V}$ as a graded Banach $R$-coalgebra (see \cref{sub:graded-seminormed-banach-coalgebras} for
	the definition) coincide with morphisms of $\tens{V}$ as a
	graded $R$-coalgebra and the same holds for coderivations.
\end{ex}

\begin{ex}
	Assume that $R$ is a graded $\mathbbm{k}$-algebra endowed with the trivial norm. Let $0 < c < 1$,
	and let $V$ be a graded $R$-module equipped with a discrete norm such that each nonzero $v \in V$
	has norm $\nnorm[v] = c$. Then each nonzero $x \in V^{\otimes n}$ has norm $\nnorm[x] = c^n$,
	and so $V^{\cotimes n} = V^{\otimes n}$ and $\tensf{V} = \prod_{i=0}^{\infty} V^{\otimes i}$
	as a graded $R$-module. The deconcatenation coproduct is indeed well-defined
	on the infinite product $\tensf{V} = \prod_{i=0}^{\infty} V^{\otimes i}$ as a map
	\begin{equation*}
		\Delta \colon \tensf{V} \rightarrow \tensf{V} \cootimes \tensf{V}
		\cong \prod_{j,k = 0}^{\infty} V^{\otimes j} \otimes V^{\otimes k}.
	\end{equation*}
\end{ex}

In what follows, we will classify morphisms $f \colon \tensf{V} \rightarrow \tensf{W}$
between formal tensor coalgebras and coderivations $\mu \colon \tensf{V} \rightharpoonup \tensf{V}$
on the formal tensor coalgebra.
To set up notation, given a graded bounded map $f \colon U \rightharpoonup \tensf{W}$ from a graded
Banach $R$-module $U$ into the formal tensor module $\tensf{W}$,
the \textbf{corestriction} $\corest{f} \colon U \rightharpoonup W$ of $f$ is defined to be the map
$\corest{f} \defeq \pi_1 \circ f$ where $\pi_1 \colon \tensf{W} \twoheadrightarrow W$ is the natural projection.
When $U = \tensf{V}$ is another formal tensor module, we will denote by
$f_k \colon V^{\cotimes k} \rightharpoonup W$ the composition of the natural inclusion
$i_k \colon V^{\cotimes k} \rightarrow \tensf{V}$ with the corestriction
$\corest{f} \colon \tensf{V} \rightharpoonup W$ and call $f_k$ the \textbf{components} of $f$.\footnote{The
	components $f_k \colon V^{\cotimes k} \rightharpoonup W$ of $f$ are graded maps of the same degree
	as $f$ and are not to be confused
	with the \textit{graded components} $f^g \colon \tensf{V}^g \rightarrow \tensf{W}^{g + \degb{f}}$ of
	$f$ as a graded map, with $g \in \GG$, which are denoted differently.}
The corestriction $\corest{f} \colon \tensf{V} \rightharpoonup W$ is
determined uniquely by the sequence of components $\left( f_k \right)_{k \geq 0}$.

\subsubsection{Grouplike Elements and the Exponential Map} \label{sub:grouplike-exp}

\begin{dfn} \label{dfn:grouplike}
	Let $C$ be a graded Banach $R$-coalgebra. An element $g \in C^0$ is called \textbf{grouplike} if
	$\Delta g = g \cootimes g$ and $\varepsilon \left( g \right) = 1$. The set of all grouplike elements
	of $C$ is denoted by $\mathcal{G} \left( C \right)$.
\end{dfn}

In order to classify the grouplike elements of $\tensf{V}$, we will need the following notion:

\begin{dfn} \label{def:top-nil}
	Let $V$ be a graded Banach $R$-module. An element $v \in V^0$ is called \textbf{topologically nilpotent} if
	$\nnorm[v^{\otimes n}]_{V^{\otimes n}} \to 0$.
	The set of all topologically nilpotent elements of $V$ is denoted by $\tc{V}$ and forms an additive
	Banach subgroup of $V^0$.
\end{dfn}

Note that if $v \in V^0$ and $\nnorm[v] < 1$ then $v$ is topologically nilpotent since we have
$\nnorm[v^{\otimes n}]_{V^{\otimes n}} \leq \nnorm[v]_{V}^n \to 0$. Note also that
$\nnorm[v^{\otimes n}]_{V^{\otimes n}} = \nnorm[v^{\cotimes n}]_{V^{\cotimes n}}$.
Hence, thinking of $v \in \tc{V}$ as an element of the formal tensor \textit{algebra} $\tensf{V}$, the element
$v^{\cotimes n}$ is precisely the product of $v \in \tensf{V}$ with itself $n$ times and hence
$v$ is topologically nilpotent in the sense of \cite[Definition 1, Section 1.2.4]{Bosch1984}.
Since we are in the non-Archimedean setting, $v \in V^0$ is topologically nilpotent if and only if the
\textbf{exponential}
\begin{equation} \label{eq:exp-tensor-coalgebra}
	\exp \left( v \right) = \Exp{v} \defeq \sum_{n = 0}^{\infty} v^{\cotimes n} = 1 + v + v \cotimes v + \dots
\end{equation}
is well-defined and converges in $\tensf{V}$.\footnote{The exponential notation, without any factorials,
	is taken from \cite[page 108]{Fukaya2009}. It can be motivated by noting that the formal tensor coalgebra,
	endowed with the \textit{shuffle multiplication} $\shuffle$, is a Banach Hopf algebra.
	Assuming the positive integers are invertible in the ground ring, the standard exponential map
	of this Hopf algebra takes the form
	$v \mapsto \sum_{n=0}^{\infty} \frac{v^{\shuffle n}}{n!} = \sum_{n=0}^{\infty} v^{\cotimes n}$
	as $v^{\shuffle n} = n! \cdot v^{\cotimes n}$, resulting in no factorials.}
The following lemma gives a classification of the grouplike elements of the formal tensor algebra:

\begin{lm} \label{lm:exp-map-bijection}
	The exponential map $\exp \colon \tc{V} \rightarrow \mathcal{G} \left( \tensf{V} \right)$ is a
	bijection between the topologically nilpotent elements of $V$ and the grouplike elements of
	$\tensf{V}$. The inverse map $\mathcal{G} \left( \tensf{V} \right) \rightarrow \tc{V}$ is given
	by the restriction of the projection $\pi_1 \colon \tensf{V} \twoheadrightarrow V$.
\end{lm}
\begin{proof}
	Given $v \in \tc{V}$, we have
	\begin{equation*}
		\Delta \left( \Exp{v} \right) = \sum_{n=0}^{\infty} \Delta \left( v^{\cotimes n} \right) =
		\sum_{n = 0}^{\infty} \sum_{n_1 + n_2 = n} v^{\cotimes n_1} \cootimes v^{\cotimes n_2} =
		\Exp{v} \cootimes \Exp{v}
	\end{equation*}
	and hence $\Exp{v}$ is a grouplike element of $\tensf{V}$.
	Conversely, let $g \in \tensf{V}$ be grouplike.
	Then $\Delta^n \left( g \right) = g^{\cootimes n}$ for all $n \geq 0$
	and hence
	\begin{equation*}
		g = \sum_{n \geq 0} \pi_n \left( g \right) \stackrel{\eqref{eq:projection-related-to-coproduct}}{=}
		\sum_{n \geq 0} \left( \pi_1^{\cotimes n} \circ \Delta^n \right) \left( g \right) =
		\sum_{n \geq 0} \pi_1^{\cotimes n} \left( g^{\cotimes n} \right) =
		\sum_{n \geq 0} \pi_1 \left( g \right)^{\cotimes n} =
		\Exp{\pi_1 \left( g \right)}.
	\end{equation*}
	In particular, we also see that if $g$ is grouplike then $\pi_1 \left( g \right)$ is topologically nilpotent,
	since $\pi_1 \left( g \right)^{\cotimes n} = \pi_n \left( g \right) \to 0$ in $\tensf{V}$.
\end{proof}

Now, let $C$ be a graded Banach $R$-coalgebra and let $D$ be a graded Banach $S$-coalgebra.
Let $\varphi \colon R \rightarrow S$ be a morphism of graded Banach $\mathbbm{k}$-algebras
and let $f \colon C \rightarrow D$ be a morphism of graded Banach coalgebras over $\varphi$.
Then, given a grouplike element $g \in C$, the element $f \left( g \right)$
is a grouplike element of $D$, hence we have an induced map
$f \colon \mathcal{G} \left( C \right) \rightarrow \mathcal{G} \left( D \right)$.
When $C = \tensf{V}[R]$ and $D = \tensf{W}[S]$ are formal tensor coalgebras, we can use the
exponential maps and define a non-linear \textbf{pushforward map}
$\mcfunc{f} \colon \tc{V} \rightarrow \tc{W}$ by setting
$\mcfunc{f} = \exp^{-1} \circ f \circ \exp$ (see \cref{fig:def-pushforward-mc}).
By definition, the map $\mcfunc{f}$ is the unique map which satisfies the identity
\begin{equation} \label{eq:exp-mcfunc}
	f \left( \Exp{v} \right) = \Exp{\mcfunc{f} \left( v \right)}
\end{equation}
for all $v \in \tc{V}$ and, by \cref{lm:exp-map-bijection}, is given explicitly by the formula
\begin{equation}
	\mcfunc{f} \left( v \right) = \pi_1 \left( f \left( \Exp{v} \right) \right) =
	\corest{f} \left( \Exp{v} \right) = f_0(1) + f_1(v) + f_2 \left( v \cotimes v \right) + \dots
	\label{eq:mcfunc-explicit}
\end{equation}
It is clear from the description above that the pushforward map is functorial with respect to morphisms,
i.e., we have
\begin{equation} \label{eq:mcfunc-functoriality}
	\mcfunc{\left( f \circ g \right)} = \mcfunc{f} \circ \mcfunc{g}, \qquad
	\mcfunc{\id} = \id.
\end{equation}

\begin{figure}
	\begin{tikzcd}
		\tc{V} && {\mathcal{G} \left( \tensf{V}[R] \right)} \\
		\tc{W} && {\mathcal{G} \left( \tensf{W}[S] \right)}
		\arrow["\exp", from=1-1, to=1-3]
		\arrow["{\mcfunc{f}}"', from=1-1, to=2-1]
		\arrow["f", from=1-3, to=2-3]
		\arrow["\exp", from=2-1, to=2-3]
	\end{tikzcd}
	\caption{Definition of the pushforward map $\mcfunc{f} \colon \tc{V} \rightarrow \tc{W}$.}
	\label{fig:def-pushforward-mc}
\end{figure}

\subsubsection{Morphisms of the Formal Tensor Coalgebra} \label{subsubsec:morphisms-formal-tensor-coalgebra}
Let $f \colon \tensf{V} \rightarrow \tensf{W}$ be a morphism of Banach $R$-coalgebras.
We start by showing that $f$ is determined uniquely by its corestriction $\corest{f} \colon \tensf{V} \rightarrow W$
and that the corestriction must satisfy the condition that $\corest{f} \left( 1 \right) = f_0(1) \in W^0$ is topologically nilpotent.

First, note that since the element $1 \in \tensf{V}$ is grouplike,
the element $f(1)$ is also grouplike, and by \cref{lm:exp-map-bijection}, must be of the form
\begin{equation}
	f(1) = \Exp{\pi_1 \left( f \left( 1 \right) \right)} = \Exp{f_0(1)}
\end{equation}
where $f_0(1) \in W$ is topologically nilpotent. Next, let us give an explicit formula for $f$ in terms of its corestriction:

\begin{lm} \label{lm:map-coalg-into-tensor-coalg}
	Let $C$ be a Banach $R$-coalgebra and let $W$ be a Banach $R$-module.
	Let $f \colon C \rightarrow \tensf{W}$ be a morphism of Banach $R$-coalgebras.
	Then for each $c \in C$ we have
	\begin{equation}
		f(c) = \lim_{N \to \infty} \sum_{n=0}^N \corest{f}^{\cotimes n} \left( \Delta^n c \right)
		= \sum_{n = 0}^{\infty} \corest{f}^{\cotimes n} \left( \Delta^n c \right)
		\label{eq:morphism-in-terms-of-corestriction}
	\end{equation}
	where the limit is taken inside $\tensf{W}$. In particular,
	$\nnorm[\corest{f}^{\cotimes n} \left( \Delta^n c \right)]_{W^{\cotimes n}} \to 0$ for all $c \in C$.
\end{lm}
\begin{proof}
	Since $f$ is a \textit{counital} morphism of coalgebras, we have the sequence of identities
	$f^{\cootimes n} \circ \Delta^n_C = \Delta^n_{\tensf{W}} \circ f$ for all $n \geq 0$
	(see \cref{eq:morphism-commutes-iterated-coproducts}). Hence, given $c \in C$, we have
	\begin{equation*}
		\begin{aligned}
			f(c) & = \sum_{n \geq 0} \pi_n \left( f(c) \right) \stackrel{\eqref{eq:projection-related-to-coproduct}}{=}
			\sum_{n \geq 0} \left( \pi_1^{\cotimes n} \circ \Delta^n_{\tensf{W}} \circ f \right) \left( c \right)
			= \sum_{n \geq 0} \left( \pi_1^{\cotimes n} \circ f^{\cotimes n} \circ \Delta^n_C \right)
			\left( c \right)
			\\
			     & =
			\sum_{n \geq 0} \corest{f}^{\cotimes n} \left( \Delta^n c \right).
		\end{aligned}
	\end{equation*}
\end{proof}

\begin{rem}
	Note that the proof of \cref{lm:map-coalg-into-tensor-coalg} is similar to the proof
	of \cref{lm:exp-map-bijection}. In fact, one can deduce \cref{lm:exp-map-bijection}
	from \cref{lm:map-coalg-into-tensor-coalg} by noting that a grouplike element $g \in C$
	of a Banach $R$-coalgebra is the same thing as a counital Banach $R$-coalgebra morphism
	$g \colon R \rightarrow C$, where $R$ is endowed with the standard $R$-coalgebra structure.
\end{rem}

We can also show the converse:
\begin{lm}
	Let $V, W$ be graded Banach $R$-modules. Given a morphism $\varphi \colon \tensf{V} \rightarrow W$
	of graded Banach $R$-modules such that $\varphi(1)$ is topologically nilpotent, there exists a unique
	morphism $f \colon \tensf{V} \rightarrow \tensf{W}$ of graded Banach $R$-coalgebras,
	called the \textbf{(coalgebra) coextension} of $\varphi$, such that $\corest{f} = \varphi$.
\end{lm}
\begin{proof}
	Uniqueness follows from \cref{lm:map-coalg-into-tensor-coalg} which also gives us the
	formula \eqref{eq:morphism-in-terms-of-corestriction} 	for $f$ in terms of its corestriction.
	Hence, we can try to define $f$ using \cref{eq:morphism-in-terms-of-corestriction},
	verify that the infinite sum actually converges and check that we get a coalgebra morphism.

	Recall first that by morphism of graded Banach $R$-modules, we mean that $\varphi$ is $R$-linear, of degree zero
	and $\nnorm[\varphi] \leq 1$. In particular, $\nnorm[\varphi(1)] \leq 1$.
	Consider the maps $f_k^n \colon V^{\cotimes k} \rightarrow W^{\cotimes n}$ defined by
	$f_k^n \defeq \rest{\varphi^{\cotimes n} \circ \Delta^n_{\tensf{V}}}{V^{\cotimes k}}$.
	A crude estimate on the norm of $f_k^n$ is obtained by noting that
	\begin{equation} \label{eq:f_k^n-crude-estimate}
		\nnorm[f_k^n] \leq \nnorm[\varphi^{\cotimes n}] \cdot \nnorm[\Delta^n]
		\leq \nnorm[\varphi]^n \cdot \nnorm[\Delta^n] \leq 1.
	\end{equation}
	When $n > k$, we can improve this estimate. Let $x = v_1 \cotimes \dots \cotimes v_k \in V^{\cotimes k}$
	and write\footnote{To avoid confusion with Sweedler's notation $x_{(1)} \cootimes \dots \cootimes x_{(n)}$ which omits the summation sign and represents the entire coproduct, we use angular brackets $x_{\gen{1}} \cootimes \dots \cootimes x_{\gen{n}}$ to denote a specific fixed summand when the sum over the coproduct is written explicitly.}
	\begin{equation*}
		\begin{aligned}
			\Delta^n x                                     & = \sum x_{\gen{1}} \cootimes \dots \cootimes x_{\gen{n}},
			\\
			\varphi^{\cotimes n} \left( \Delta^n x \right) & =
			\sum \varphi \left( x_{\gen{1}} \right) \cotimes \dots \cotimes \varphi \left( x_{\gen{n}} \right),
		\end{aligned}
	\end{equation*}
	where $x_{\gen{1}}, \dots, x_{\gen{n}} \in \tensf{V}$ with
	$\nnorm[x_{\gen{1}} \cootimes \dots \cootimes x_{\gen{n}}]_{\tensf{V}^{\cootimes n}} =
		\nnorm[v_1 \cotimes \dots \cotimes v_k]_{V^{\cotimes k}}$.

	Since $n > k$, each fixed summand $x_{\gen{1}} \cootimes \dots \cootimes x_{\gen{n}}$ of $\Delta^n x$
	contains at least $n - k$ copies of $1$. More precisely, we can find a permutation $\sigma$
	of $\Set{1,\dots,n}$ (which depends on the specific summand) such that
	$x_{\gen{\sigma(k+1)}} = \dots = x_{\gen{\sigma(n)}} = 1$.
	The tensor norm is invariant under permutations,
	as well as natural unit and associativity isomorphisms, and so we have
	\begin{equation*}
		\begin{aligned}
			\nnorm[x_{\gen{\sigma(1)}} \cootimes \dots \cootimes x_{\gen{\sigma(k)}}]_{\tensf{V}^{\cootimes k}}
			 & =
			\nnorm[x_{\gen{\sigma(1)}} \cootimes \dots \cootimes x_{\gen{\sigma(k)}}
				\cootimes 1 \cootimes \dots \cootimes 1]_{\tensf{V}^{\cootimes n}}
			\\
			 & =
			\nnorm[x_{\gen{\sigma(1)}} \cootimes \dots \cootimes x_{\gen{\sigma(n)}}]_{\tensf{V}^{\cootimes n}}
			\\
			 & =
			\nnorm[x_{\gen{1}} \cootimes \dots \cootimes x_{\gen{n}}]_{\tensf{V}^{\cootimes n}}
			\\
			 & =
			\nnorm[v_1 \cotimes \dots \cotimes v_k]_{V^{\cotimes k}}.
		\end{aligned}
	\end{equation*}
	Hence,
		{ \small
			\begin{equation*}
				\begin{aligned}
					\nnorm[\varphi \left( x_{\gen{1}} \right) \cotimes \dots
						\cotimes \varphi \left( x_{\gen{n}} \right)]_{W^{\cotimes n}}
					 & = \nnorm[\varphi \left( x_{\gen{\sigma(1)}} \right) \cotimes \dots \cotimes
						     \varphi \left( x_{\gen{\sigma(n)}} \right)]_{W^{\cotimes n}}
					\\
					 & =
					\nnorm[\varphi \left( x_{\gen{\sigma(1)}} \right) \cotimes \dots \cotimes
						\varphi \left( x_{\gen{\sigma(k)}} \right) \cotimes
						\varphi(1)^{\cotimes (n-k)}]_{W^{\cotimes k} \cotimes W^{\cotimes (n-k)}}
					\\
					 & =
					\nnorm[
						\varphi^{\cotimes k} \left( x_{\gen{\sigma(1)}} \cootimes \dots \cootimes x_{\gen{\sigma(k)}} \right)
						\cotimes \varphi(1)^{\cotimes (n-k)}]_{W^{\cotimes k} \cotimes W^{\cotimes (n-k)}}
					\\
					 & \leq
					\nnorm[\varphi^{\cotimes k} \left( x_{\gen{\sigma(1)}} \cootimes \dots
						\cootimes x_{\gen{\sigma(k)}} \right)]_{W^{\cotimes k}} \cdot
					\nnorm[\varphi(1)^{\cotimes (n-k)}]_{W^{\cotimes (n-k)}}
					\\
					 & \leq
					\nnorm[\varphi^{\cotimes k}] \cdot
					\nnorm[x_{\gen{\sigma(1)}} \cootimes \dots \cootimes x_{\gen{\sigma(k)}}]_{\tensf{V}^{\cootimes k}} \cdot
					\nnorm[\varphi(1)^{\cotimes (n-k)}]_{W^{\cotimes (n-k)}}
					\\
					 & =
					\nnorm[\varphi^{\cotimes k}] \cdot \nnorm[v_1 \cotimes \dots \cotimes v_k]_{V^{\cotimes k}} \cdot
					\nnorm[\varphi(1)^{\cotimes (n-k)}]_{W^{\cotimes (n-k)}}
					\\
					 & \leq
					\nnorm[\varphi]^k \cdot \nnorm[\varphi(1)^{\cotimes (n-k)}] \cdot
					\nnorm[v_1] \cdots \nnorm[v_k]
					\\
					 & \leq
					\nnorm[\varphi(1)^{\cotimes (n-k)}] \cdot
					\nnorm[v_1] \cdots \nnorm[v_k].
				\end{aligned}
			\end{equation*}
		}
	Since the estimate above holds for each summand, we immediately get the estimate
	\begin{equation*}
		\nnorm[f_k^n \left( v_1 \cotimes \dots \cotimes v_k \right)]
		\leq
		\nnorm[\varphi(1)^{\cotimes (n-k)}] \cdot \nnorm[v_1] \cdots \nnorm[v_k].
	\end{equation*}
	The estimate above together with the definition of the tensor product norm imply that we have
	\begin{equation} \label{eq:f_k^n-estimate}
		\nnorm[f_k^n \left( x \right)] \leq \nnorm[\varphi(1)^{\cotimes (n-k)}] \cdot \nnorm[x]
	\end{equation}
	for all $x \in V^{\cotimes k}$, even if $x$ is not an elementary tensor.

	Now, fix $k \geq 0$ and define a map $f_k^{\bullet} \colon V^{\cotimes k} \rightarrow \tensf{W}$
	by the formula $f_k^{\bullet} \left( x \right) \defeq \sum_{n \geq 0} f_k^n \left( x \right)$,
	where we consider each $f_k^n \left( x \right)$ as an element of $\tensf{W}$ and take the limit
	in $\tensf{W}$. For a fixed $x \in V^{\cotimes k}$, \cref{eq:f_k^n-estimate} together
	with the assumption that $\varphi(1)$ is topologically nilpotent imply that
	$\nnorm[f_k^n \left( x \right)] \xrightarrow[n \to \infty]{} 0$ so that the infinite sum
	indeed converges. The maps $f_k^{\bullet}$ are $R$-linear and bounded with
	$\nnorm[f_k^{\bullet}] \leq 1$. Hence, by the universal property of the direct sum,
	there exists a unique contractive $R$-linear map $f \colon \tensf{V} \rightarrow \tensf{W}$
	such that
	$\rest{f}{V^{\cotimes k}} = f_k^{\bullet}$ for all $k \geq 0$.

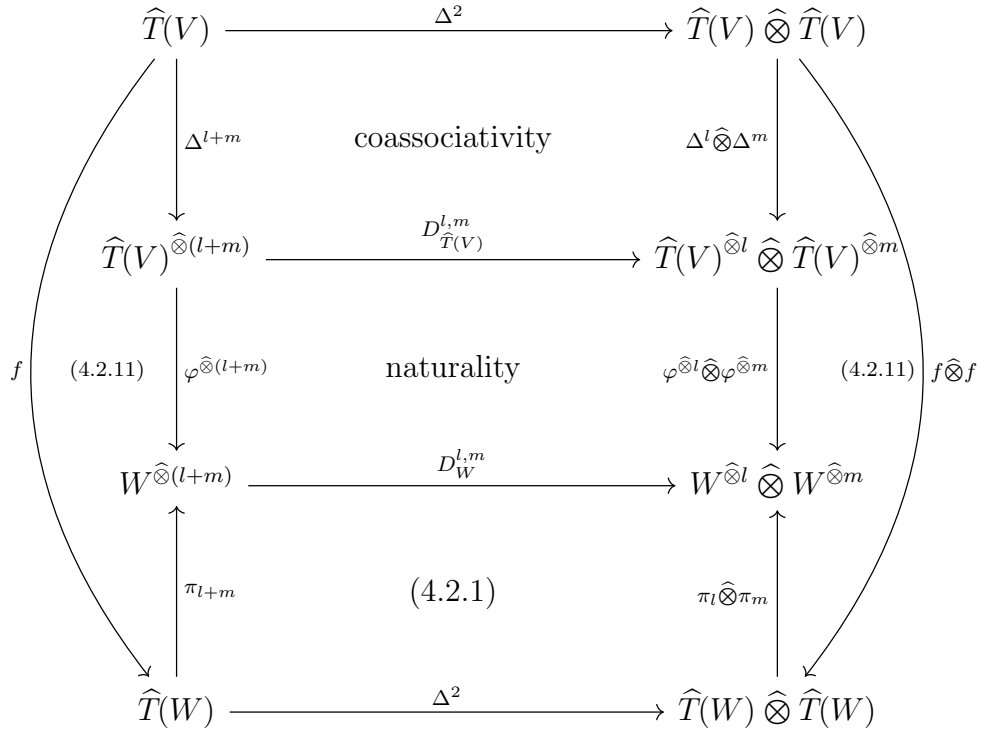
\begin{figure}[H]
	\begin{tikzcd}
		{\tensf{V}} && {\tensf{V} \cootimes \tensf{V}} \\
		& {\textrm{coassociativity}} \\
		{\tensf{V}^{\cotimes (l + m)}} && {\tensf{V}^{\cotimes l} \cootimes \tensf{V}^{\cotimes m}} \\
		& {\textrm{naturality}} \\
		{W^{\cotimes (l + m)}} && {W^{\cotimes l} \cootimes W^{\cotimes m}} \\
		& {\eqref{eq:projection-related-to-coproduct-general}} \\
		{\tensf{W}} && {\tensf{W} \cootimes \tensf{W}}
		\arrow["{\Delta^2}", from=1-1, to=1-3]
		\arrow["{\Delta^{l+m}}", from=1-1, to=3-1]
		\arrow["f"', "{\quad\eqref{eq:pi-l-f-varphi}}", curve={height=70pt}, from=1-1, to=7-1]
		\arrow["{\Delta^l \cootimes \Delta^m}"', from=1-3, to=3-3]
		\arrow["{f \cootimes f}", "{\qquad\eqref{eq:pi-l-f-varphi}}"', curve={height=-70pt}, from=1-3, to=7-3]
		\arrow["{D_{\tensf{V}}^{l,m}}", from=3-1, to=3-3]
		\arrow["{\varphi^{\cotimes (l+m)}}", from=3-1, to=5-1]
		\arrow["{\varphi^{\cotimes l} \cootimes \varphi^{\cotimes m}}"', from=3-3, to=5-3]
		\arrow["{D_W^{l,m}}", from=5-1, to=5-3]
		\arrow["{\pi_{l+m}}"', from=7-1, to=5-1]
		\arrow["{\Delta^2}", from=7-1, to=7-3]
		\arrow["{\pi_l \cootimes \pi_m}", from=7-3, to=5-3]
	\end{tikzcd}
	\caption{The commutative diagram behind the proof of \cref{eq:defined-f-is-a-coalgebra-morphism}.}
	\label{fig:proof-f-morphism}
\end{figure}

	It is clear that $\corest{f} = \varphi$. It remains to verify that $f$ is indeed a counital
	morphism of coalgebras. To show that $f$ is a morphism of coalgebras, we need to verify that
	\begin{equation*}
		\left( f \cootimes f \right) \circ \Delta_{\tensf{V}} = \Delta_{\tensf{W}} \circ f \colon \tensf{V}
		\rightarrow \tensf{W} \cootimes \tensf{W}.
	\end{equation*}
	It is enough to verify that both sides agree when composed
	with the maps $\pi_l \cootimes \pi_m$ where $\pi_l \colon \tensf{W} \rightarrow W^{\cotimes l}$ are the
	canonical projections. Note that we have
	\begin{equation} \label{eq:pi-l-f-varphi}
		\pi_l \circ f = \varphi^{\cotimes l} \circ \Delta^l_{\tensf{V}}
	\end{equation}
	by the definition of $f$ and hence we have
	\begin{equation} \label{eq:defined-f-is-a-coalgebra-morphism}
		\begin{aligned}
			\left( \pi_l \cootimes \pi_m \right) \circ \left( f \cootimes f \right) \circ \Delta_{\tensf{V}}
			\eqwithref                                             &
			\left( \pi_l \circ f \right) \cootimes \left( \pi_m \circ f \right) \circ \Delta_{\tensf{V}}
			\\
			\eqwithref[eq:pi-l-f-varphi]                           &
			\left( \varphi^{\cotimes l} \circ \Delta^l_{\tensf{V}} \right) \cootimes
			\left( \varphi^{\cotimes m} \circ \Delta^m_{\tensf{V}} \right) \circ \Delta_{\tensf{V}}
			\\
			\eqwithref                                             &
			\left( \varphi^{\cotimes l} \cootimes \varphi^{\cotimes m} \right) \circ
			\left( \Delta^l_{\tensf{V}} \cootimes \Delta^m_{\tensf{V}} \right) \circ \Delta_{\tensf{V}}
			\\
			\eqwithref                                             &
			\left( \varphi^{\cotimes l} \cootimes \varphi^{\cotimes m} \right) \circ D^{l,m}_{\tensf{V}}
			\circ \Delta^{l+m}_{\tensf{V}}
			\\
			\eqwithref                                             &
			D^{l,m}_W \circ \varphi^{\cotimes (l+m)} \circ \Delta^{l+m}_{\tensf{V}}
			\\
			\eqwithref[eq:pi-l-f-varphi]                           &
			D^{l,m}_W \circ \pi_{l+m} \circ f
			\\
			\eqwithref[eq:projection-related-to-coproduct-general] &
			\left( \pi_l \cootimes \pi_m \right) \circ \Delta_{\tensf{W}} \circ f,
		\end{aligned}
	\end{equation}
	where we used the coassociativity of the coproduct and the naturality of the associativity isomorphisms
	(see \cref{fig:proof-f-morphism}).
	Finally, we have
	\begin{equation*}
		\varepsilon_{\tensf{W}} \circ f = \pi_0 \circ f = \varphi^{\cotimes 0} \circ \Delta_{\tensf{V}}^0
		= \varepsilon_{\tensf{V}}
	\end{equation*}
	so $f$ is indeed counital.
\end{proof}

Hence, we have shown:
\begin{prop} \label{prop:classification-morphisms-tensor-coalgebra}
	Let $R$ be a graded-commutative Banach $\mathbbm{k}$-algebra and let $V$ and $W$ be two graded
	Banach $R$-modules.
	We have a bijective correspondence between
	\begin{enumerate}
		\item Morphisms $f \colon \tensf{V} \rightarrow \tensf{W}$ of graded Banach $R$-coalgebras.
		\item Morphisms $\corest{f} \colon \tensf{V} \rightarrow W$ of graded Banach $R$-modules
		      such that $f_0(1)$ is topologically nilpotent.
	\end{enumerate}
	Given a morphism $f \colon \tensf{V} \rightarrow \tensf{W}$ of graded Banach $R$-coalgebras,
	the morphism $\corest{f} \colon \tensf{V} \rightarrow W$ is the corestriction
	$\corest{f} \defeq \pi_1 \circ f$
	of $f$. Conversely, given a morphism $\corest{f} \colon \tensf{V} \rightarrow W$ of graded Banach $R$-modules
	satisfying $(2)$, there exists a unique morphism $f \colon \tensf{V} \rightarrow \tensf{W}$
	of graded Banach $R$-coalgebras, given by the formula
	\begin{align*}
		 & f(l) = \sum_{k = 0}^{\infty} \corest{f}(l_{(1)}) \cotimes \dots \cotimes \corest{f}(l_{(k)})
		 & \left( l \in \tensf{V} \right),
	\end{align*}
	whose corestriction is $\corest{f}$.\footnote{The $k = 0$ summand is to be interpreted as
		$\left( \corest{f}^{\cotimes 0} \circ \Delta^0 \right) \left( l \right) = \Delta^0 \left( l \right)$
		where $\Delta^0 \colon \tensf{V} \rightarrow R$ is the counit.}

	The morphism $f$ is coaugmented
	(i.e., $f(1) = 1$) if and only if $\corest{f} \left( 1 \right) = 0$. \qed
\end{prop}

\begin{ex}
	Assume that $R$ is a graded $\mathbbm{k}$-algebra endowed with the trivial norm and let $V,W$
	be graded $R$-modules endowed with the trivial norm. In this case, we have $\tensf{V} = \tens{V}$,
	$\tensf{W} = \tens{W}$ (see \cref{ex:tensor-coalgebra-as-a-formal-tensor-coalgebra}), and
	\cref{prop:classification-morphisms-tensor-coalgebra} recovers the usual correspondence between
	coaugmented coalgebra morphisms $f \colon \tens{V} \rightarrow \tens{W}$ and
	maps $\corest{f} \colon \tens{V} \rightarrow W$ such that $\corest{f}(1) = 0$ (see
	\cite[Section 1.2.6]{Loday2012}).
	When $f$ is not necessarily coaugmented, $f_0(1)$ must be ``topologically nilpotent'',
	which in this setting means that $f_0(1)^{\otimes n} = 0$ for some $n \geq 1$. When
	$R = \mathbbm{k}$ is a field, this is possible only if $f_0(1) = 0$ and so
	morphisms $f \colon \tens{V}[\mathbbm{k}] \rightarrow \tens{W}[\mathbbm{k}]$ between tensor coalgebras
	over a field must satisfy $f(1) = 1$, i.e., they are automatically coaugmented.
\end{ex}

\begin{rem} \label{rem:not-every-coalgebra-morphism-formal-is-completion}
	Let $V$ and $W$ be graded Banach $R$-modules. Any $R$-coalgebra morphism
	$f \colon \tens{V} \rightarrow \tens{W}$ with $\nnorm[f] \leq 1$
	(i.e., a morphism of graded seminormed $R$-coalgebras) can be completed
	to a morphism $\widehat{f} \colon \tensf{V} \rightarrow \tensf{W}$ of graded Banach $R$-coalgebras.
	However, not every morphism $g \colon \tensf{V} \rightarrow \tensf{V}$ between formal
	tensor coalgebras can be realized as a completion $g = \widehat{f}$. In fact,
	$g$ can be realized as a completion precisely when $g_0(1)^{\otimes n} = 0$ for some
	$n \geq 1$.
\end{rem}

The classification of morphisms $f \colon \tensf{V} \rightarrow \tensf{W}$ readily extends
to coalgebra morphisms between formal tensor coalgebra over different ground algebras:

\begin{prop} \label{prop:classification-morphisms-tensor-coalgebra-different-ground-algebras}
	Let $R,S$ be two graded-commutative Banach $\mathbbm{k}$-algebras and let $\varphi \colon R \rightarrow S$
	be a morphism of graded Banach $\mathbbm{k}$-algebras. Let $V$ be a graded Banach $R$-module and let $W$
	be a graded Banach $S$-module. We have a bijective correspondence between
	\begin{enumerate}
		\item Banach coalgebra morphisms $f \colon \tensf{V}[R] \rightarrow \tensf{W}[S]$ over $\varphi$.
		\item Morphisms $\corest{f} \colon \tensf{V}[R] \rightarrow W$ of graded Banach modules over $\varphi$
		      such that $f_0(1) \in W$ is topologically nilpotent.
	\end{enumerate}
	Given a Banach coalgebra morphism $f \colon \tensf{V}[R] \rightarrow \tensf{W}[S]$, the morphism
	$\corest{f} \colon \tensf{V}[R] \rightarrow W$ is the corestriction $\corest{f} \defeq \pi_1 \circ f$
	of $f$. Conversely, given a morphism $\corest{f} \colon \tensf{V}[R] \rightarrow W$ satisfying $(2)$,
	there exists a unique Banach coalgebra morphism $f \colon \tensf{V}[R] \rightarrow \tensf{W}[S]$
	over $\varphi$, given by the formula
	\begin{align*}
		 & f(l) = \sum_{k = 0}^{\infty}
		\corest{f}(l_{(1)}) \cotimes_S \dots \cotimes_S \corest{f}(l_{(k)})
		 & \left( l \in \tensf{V}[R] \right),
	\end{align*}
	whose corestriction is $\corest{f}$.\footnote{More precisely, we have the pointwise formula
		$f = \sum_{n = 0}^{\infty} \corest{f}^{\cotimes_{\varphi} n} \circ \Delta_{\tensf{V}[R]}^n$
		and the $n = 0$ summand in the sum $\corest{f}^{\cotimes_{\varphi} 0} \circ \Delta_{\tensf{V}[R]}^0$
		is given by the composition
		$\tensf{V}[R] \xrightarrow{\varepsilon_{\tensf{V}[R]}} R \xrightarrow{\corest{f}^{\cotimes_{\varphi} 0} = \varphi}
			S$.}
\end{prop}
\begin{proof}
	Given a morphism $f \colon \tensf{V}[R] \rightarrow \tensf{W}[S]$ of Banach coalgebras over
	$\varphi$, its corestriction $\corest{f} \colon \tensf{V}[R] \rightarrow W$ is the composition of
	$f$ with the projection $\tensf{W}[S] \rightarrow W$ and is also a morphism over $\varphi$.
	A map $f \colon \tensf{V}[R] \rightarrow \tensf{W}[S]$ over $\varphi$, i.e.,
	an $R$-linear map $f \colon \tensf{V}[R] \rightarrow \varphi^{*} ( \tensf{W}[S] )$,
	is a morphism of Banach coalgebras if and only if the adjunct
	$\tilde{f} \colon \varphi_{!} ( \tensf{V}[R] ) \rightarrow \tensf{W}[S]$
	under the adjunction \eqref{eq:restriction-extension-adjunction-GBMod-R}
	is a morphism of Banach $S$-coalgebras.
	Note that since $\varphi_{!}$ commutes with colimits and is strong monoidal, we have
	\begin{equation*}
		\varphi_{!} \left( \tensf{V}[R] \right) =
		\varphi_{!} \left( \extrawidehat{\bigoplus}_{i=0}^{\infty} V^{\cotimes_R i} \right) \cong
		\extrawidehat{\bigoplus}_{i=0}^{\infty} \varphi_{!} \left( V^{\cotimes_R i} \right) \cong
		\extrawidehat{\bigoplus}_{i=0}^{\infty} \varphi_{!} \left( V \right)^{\cotimes_S i} \cong
		\tensf{\varphi_{!} \left( V \right)}[S],
	\end{equation*}
	not only as graded Banach $S$-modules but also as graded Banach $S$-coalgebras. Hence, the result follows
	from \cref{prop:classification-morphisms-tensor-coalgebra}.
\end{proof}

\begin{rem} \label{rem:morphism-tensor-coalgebras-recover-base-map}
	Note that since the formal tensor coalgebras are coaugmented, given a morphism
	$f \colon \tensf{V}[R] \rightarrow \tensf{W}[S]$ of Banach coalgebras over $\varphi$,
	the underlying morphism $\varphi$ can be recovered from the
	map $f$ as $\pi_0 \circ f \circ i_0$ where $i_0 \colon R \rightarrow \tensf{V}[R]$ is the
	natural inclusion (the coaugmentation of $\tensf{V}[R]$) and
	$\pi_0 \colon \tensf{W}[S] \rightarrow S$ is the natural projection (the counit of $\tensf{W}[S]$).
	This means that the data of a Banach coalgebra morphism $f \colon \tensf{V}[R] \rightarrow \tensf{W}[S]$
	over $\varphi$ is encoded entirely in the map $f$.

	Given a map $f \colon \tensf{V}[R] \rightarrow \tensf{W}[S]$ between formal tensor coalgebras
	over different ground algebras, we will set $\base{f} \defeq \pi_0 \circ f \circ i_0$ and say
	that $f$ is a Banach coalgebra morphism if $\base{f} \colon R \rightarrow S$ is a morphism
	of graded Banach $\mathbbm{k}$-algebras and $f$ is a Banach coalgebra morphism over $\base{f}$.
	This way we won't have to specify in advance the underlying morphism of $f$.
	When $R = S$ and $\base{f} = \id$, we recover the standard notion of a morphism
	$f \colon \tensf{V}[R] \rightarrow \tensf{W}[R]$ of graded Banach $R$-coalgebras.
\end{rem}

\subsubsection{Coderivations on the Tensor Coalgebra} \label{subsubsec:coder-tensor-coalgebra}
In what follows, we classify coderivations $\mu \colon \tensf{V} \rightharpoonup \tensf{V}$
on the formal tensor coalgebra, viewed as a graded Banach $R$-coalgebra.
Unlike the case of morphisms, the classification of
coderivations on the formal tensor coalgebra $\tensf{V}$ is entirely parallel to the classification of
coderivations on the standard tensor coalgebra $\tens{V}$.
For completeness, we recall the argument. It will be beneficial to be slightly more general and
classify coderivations whose domain is a graded Banach bicomodule over $\tensf{V}$ (i.e.,
a bicomodule object over a coalgebra object in the category $\GBMod[R]$).
In what follows, we use the notation from \cref{subsec:comodules-in-monoidal-cat}.

\begin{dfn} \label{dfn:coderivation-domain-bicomodule}
	Let $C$ be a graded Banach $R$-coalgebra and let $M$ be a graded Banach $C$-bicomodule.
	A graded bounded $R$-linear map $\mu \colon M \rightharpoonup C$ is called a (\textbf{bicomodule})
	\textbf{coderivation} if it satisfies the graded co-Leibniz rule
	\begin{equation*}
		\Delta_C \circ \mu = \left( \mu \cotimes \id_C \right) \circ \Delta_{M}^{0|1} +
		\left( \id_C \cotimes \mu \right) \circ \Delta_{M}^{1|0}.
	\end{equation*}
\end{dfn}

When $M = C$ with the canonical $C$-bicomodule structure, in which
$\Delta_M^{1|0} = \Delta_M^{0|1} = \Delta_C$, one recovers the definition
of an $R$-linear coderivation $\mu \colon C \rightharpoonup C$ on a graded Banach $R$-coalgebra
(see \cref{def:coderivation-R-linear-seminormed-Banach-coalgebra}).
Using coassociativity and counitality, one can show that any coderivation satisfies
the generalized co-Leibniz rule\footnote{We suppress the canonical associativity isomorphisms
	required to identify the codomain $C^{\cotimes n_1} \cotimes C \cotimes C^{\cotimes n_2}$ of
	the right-hand side with the codomain $C^{\cotimes n}$ of the left-hand side.}
\begin{equation}
	\Delta_C^n \circ \mu =
	\sum_{n_1 + 1 + n_2 = n} \left( \id^{\cotimes n_1} \cotimes \mu \cotimes \id^{\cotimes n_2} \right) \circ
	\Delta_{M}^{n_1|n_2}.
	\label{eq:generalized-coLeibniz-comodules}
\end{equation}
for all $n \geq 0$.

\begin{lm} \label{lm:coderivation-determined-by-corestriction}
	Let $V$ be a graded Banach $R$-module and let $M$ be a graded Banach $\tensf{V}$-bicomodule.
	Let $\mu \colon M \rightharpoonup \tensf{V}$ be a coderivation. Then $\mu$ is determined uniquely by its
	corestriction $\corest{\mu} \colon M \rightharpoonup V$.
\end{lm}
\begin{proof}
	The proof is analogous to the proof of \cref{lm:map-coalg-into-tensor-coalg},
	the only difference being using \cref{eq:generalized-coLeibniz-comodules}
	instead of \cref{eq:morphism-commutes-iterated-coproducts}.
	For any $n \geq 0$, we have
	\begin{equation}
		\begin{aligned}
			\pi_n \circ \mu & = \left( \pi_1 \right)^{\cotimes n} \circ \Delta_{\tensf{V}}^n \circ \mu
			\\
			                & =
			\left( \pi_1 \right)^{\cotimes n}
			\circ \left( \sum_{n_1 + 1 + n_2 = n}
			\left( \id^{\cotimes n_1} \cotimes \mu \cotimes \id^{\cotimes n_2} \right)
			\circ \Delta_{M}^{n_1|n_2} \right)
			\\
			                & = \sum_{n_1 + 1 + n_2 = n} \left(
			\pi_1^{\cotimes n_1} \cotimes \corest{\mu} \cotimes \pi_1^{\cotimes n_2} \right)
			\circ \Delta_{M}^{n_1|n_2}.
		\end{aligned} \label{eq:weight-n-component-coderivation}
	\end{equation}
	\Cref{eq:weight-n-component-coderivation} gives us a formula for the weight $n$ component of $\mu$
	in terms of the corestriction $\corest{\mu}$ and shows that if we have two coderivations $\mu,\nu$
	with $\corest{\mu} = \corest{\nu}$ then $\mu = \nu$.
\end{proof}

Next, we show that any coderivation $\mu \colon M \rightharpoonup \tensf{V}$ factors uniquely via
a universal coderivation and obtain an alternative formula for $\mu$ in terms of its corestriction.
We start by recalling the notion of a cofree bicomodule.
Let $C$ be a graded Banach $R$-coalgebra, let $V$ be a graded Banach $R$-module and set
$M = C \cotimes V \cotimes C$. The graded Banach $R$-module $M$ has a natural structure of
a $C$-bicomodule with the coaction maps
\begin{align*}
	\Delta_{M}^{1|0} & \defeq  \Delta_C \cotimes \id_{V} \cotimes \id_{C}
	\colon C \cotimes V \cotimes C = M \rightarrow C \cotimes M = C \cotimes C \cotimes V \cotimes C,
	\\
	\Delta_{M}^{0|1} & \defeq \id_{C} \cotimes \id_{V} \cotimes \Delta_C
	\colon C \cotimes V \cotimes C = M \rightarrow M \cotimes C = C \cotimes V \cotimes C \cotimes C.
\end{align*}
The $C$-bicomodule $M$ is called the \textbf{cofree} $C$-bicomodule on $V$ and comes equipped with a
canonical morphism $p \colon C \cotimes V \cotimes C \rightarrow V$ of graded Banach $R$-modules
given by $p = m_{V}^{1|1} \circ \left( \varepsilon_C \cotimes \id_{V} \cotimes \varepsilon_C \right)$
where $m_{V}^{1|1} \colon R \cotimes V \cotimes R \rightarrow V$ is the natural isomorphism coming
from the $R$-module structure on $V$. The morphism $p$ is sometimes called the \textbf{cogenerator} of
the cofree $C$-bicomodule. The pair $\left( C \cotimes V \cotimes C, p \right)$ satisfies the following
universal property: Given a graded Banach $C$-bicomodule $N$ and a graded bounded $R$-linear map
$\varphi \colon N \rightharpoonup V$, there exists a unique graded bounded $R$-linear map
$\tilde{\varphi} \colon N \rightharpoonup C \cotimes V \cotimes C$ of graded Banach $C$-bicomodules
such that $p \circ \tilde{\varphi} = \varphi$ (see \cref{fig:cofree-bicomodule}).
The map $\tilde{\varphi}$
is given by $\tilde{\varphi} = \left( \id_{C} \cotimes \varphi \cotimes \id_{C} \right) \circ \Delta_M^{1|1}$
and is called the \textbf{(bicomodule) coextension} of $\varphi$.

\begin{figure}
	\begin{tikzcd}
		{C \cotimes V \cotimes C} && V \\
		M
		\arrow["p", from=1-1, to=1-3]
		\arrow["\varphi"', harpoon, from=2-1, to=1-3]
		\arrow["{\exists! \, \tilde{\varphi}}", dashed, harpoon, from=2-1, to=1-1]
	\end{tikzcd}
	\caption{The universal property of the cofree $C$-bicomodule.}
	\label{fig:cofree-bicomodule}
\end{figure}

\begin{dfn} \label{dfn:universal-coderivation}
	Let $V$ be a graded Banach $R$-module, and let
	$\tensf{V} \cotimes V \cotimes \tensf{V}$ be the cofree $\tensf{V}$-bicomodule on $V$.
	Let $m = m_{\bullet,1,\bullet} \colon \tensf{V} \cotimes V \cotimes \tensf{V} \rightarrow \tensf{V}$
	be the natural multiplication map. Then $m$ is a degree zero contractive coderivation called the
	\textbf{universal coderivation}.
\end{dfn}

One can verify directly that the corestriction $\corest{m} = \pi_1 \circ m$ of the universal coderivation
coincides with the cogenerator $p \colon \tensf{V} \cotimes V \cotimes \tensf{V} \rightarrow V$ of the cofree
$\tensf{V}$-bicomodule on $V$.

\begin{lm} \label{lm:coderivation-from-corestriction}
	Let $V$ be a graded Banach $R$-module and let $M$ be a graded Banach $\tensf{V}$-bicomodule. Given
	a graded bounded $R$-linear map $\varphi \colon M \rightharpoonup V$ there exists a unique
	coderivation $\mu \colon M \rightharpoonup \tensf{V}$,
	called the \textbf{(coderivation) coextension} of $\varphi$, with $\corest{\mu} = \varphi$.
	The coderivation $\mu$ is given in terms of $\varphi$ by the formula
	$\mu = m \circ \left( \id \cotimes \varphi \cotimes \id \right) \circ \Delta_{M}^{1|1}$.
\end{lm}
\begin{proof}
	Let $\tilde{\varphi} \colon M \rightharpoonup \tensf{V} \cotimes V \cotimes \tensf{V}$ be the bicomodule
	coextension of $\varphi$ and consider $\mu = m \circ \tilde{\varphi}$. Then $\mu$ is a coderivation
	(being the composition of a bicomodule map and a coderivation) and we have
	\begin{equation*}
		\corest{\mu} = \pi_1 \circ \mu = \pi_1 \circ m \circ \tilde{\varphi} = p \circ \tilde{\varphi} = \varphi
	\end{equation*}
	by the identity $\pi_1 \circ m = p$. Uniqueness follows from \cref{lm:coderivation-determined-by-corestriction}.
\end{proof}

The usage of the terminology ``universal coderivation'' for $m$ is explained by the following proposition:
\begin{prop} \label{prop:m-universal-coderivation}
	Let $V$ be a graded Banach $R$-module.
	The coderivation $m$	 is universal in the sense that
	given any $\tensf{V}$-bicomodule $M$ and any coderivation $\mu \colon M \rightharpoonup \tensf{V}$
	there exists a unique graded bounded $R$-linear
	map $f \colon M \rightharpoonup \tensf{V} \cotimes V \cotimes \tensf{V}$
	of graded Banach $\tensf{V}$-bicomodules such that $\mu = m \circ f$ (see \cref{fig:universal-coderivation}).
	The map $f$ is given by the formula
	$f \defeq \left( \id \cotimes \corest{\mu} \cotimes \id \right) \circ \Delta_M^{1|1}$.
	\begin{figure}[h]
		\begin{tikzcd}
			{\tensf{V} \cotimes V \cotimes \tensf{V}} && {\tensf{V}} \\
			M
			\arrow["m", from=1-1, to=1-3]
			\arrow["\mu"', harpoon, from=2-1, to=1-3]
			\arrow["{\exists! \, f}", dashed, harpoon, from=2-1, to=1-1]
		\end{tikzcd}
		\caption{Multiplication as a universal coderivation.}
		\label{fig:universal-coderivation}
	\end{figure}
\end{prop}
\begin{proof}
	Given a graded bounded map $f \colon M \rightharpoonup \tensf{V} \cotimes V \cotimes \tensf{V}$
	of $\tensf{V}$-bicomodules which satisfies $m \circ f = \mu$, we must have
	$\pi_1 \circ m \circ f = p \circ f = \corest{\mu}$ and hence $f$ is determined uniquely by
	$\mu$ (or, more precisely, $\corest{\mu}$) via the universal property of the cofree bicomodule.

	Conversely, given $\mu$, let $f \colon M \rightharpoonup \tensf{V} \cotimes V \cotimes \tensf{V}$
	be the coextension of $\corest{\mu}$ given by $f = \left( \id \cotimes \corest{\mu} \cotimes \id \right)
		\circ \Delta_{M}^{1|1}$. By \cref{lm:coderivation-from-corestriction},
	the map $m \circ f$ is a coderivation whose corestriction is
	$\corest{\mu}$ and hence, by uniqueness, $m \circ f = \mu$.
\end{proof}

Applying \cref{lm:coderivation-from-corestriction} to the case where $M = \tensf{V}$
with the canonical $\tensf{V}$-bicomodule structure, we obtain the following result:
\begin{prop} \label{prop:classification-coderivations-formal-tensor-coalgebra}
	Let $R$ be a graded-commutative Banach $\mathbbm{k}$-algebra and let $V$ be a graded Banach $R$-module.
	We have a bijective correspondence between
	\begin{enumerate}
		\item Banach coalgebra coderivations $\mu \colon \tensf{V} \rightharpoonup \tensf{V}$.
		\item Graded bounded $R$-linear maps $\corest{\mu} \colon \tensf{V} \rightharpoonup V$.
		\item Sequences $( \mu_k \colon V^{\cotimes k} \rightharpoonup V )_{k \geq 0}$ of graded bounded
		      $R$-linear maps having the same degree which satisfy $\sup_{k \geq 0} \nnorm[\mu_k] < \infty$.
		\item Sequences $\left( U_k \colon V^{\times k} \rightharpoonup V \right)_{k \geq 0}$ of graded
		      bounded $R$-multilinear maps having the same degree which satisfy $\sup_{k \geq 0} \nnorm[U_k] < \infty$.
	\end{enumerate}
	Given a Banach coalgebra coderivation $\mu \colon \tensf{V} \rightharpoonup \tensf{V}$,
	the map	$\corest{\mu} \colon \tensf{V} \rightharpoonup V$ is the corestriction
	$\corest{\mu} \defeq \pi_1 \circ \mu$ of $\mu$. Conversely, given a graded bounded $R$-linear map
	$\corest{\mu} \colon \tensf{V} \rightharpoonup V$,
	there exists a unique Banach coalgebra coderivation $\mu \colon \tensf{V} \rightharpoonup \tensf{V}$,
	given by the formula
	\begin{align*}
		 & \mu(l) = (-1)^{\braidd{\corest{\mu}}{l_{(1)}}}
		l_{(1)} \cotimes \corest{\mu} \left( l_{(2)} \right) \cotimes l_{(3)}
		 & \left( l \in \tensf{V} \right),
	\end{align*}
	whose corestriction is $\corest{\mu}$. \qed
\end{prop}

\begin{ex}
	Assume that $R$ is a graded $\mathbbm{k}$-algebra endowed with the trivial norm and let $V$
	be a graded $R$-module endowed with the trivial norm.
	In this case, we have $\tensf{V} = \tens{V}$ (see \cref{ex:tensor-coalgebra-as-a-formal-tensor-coalgebra})
	and 	\cref{prop:classification-coderivations-formal-tensor-coalgebra} recovers the usual correspondence
	between graded coderivations $\mu \colon \tens{V} \rightharpoonup \tens{V}$ and
	graded $R$-linear maps $\corest{\mu} \colon \tens{V} \rightharpoonup V$ (see for example
	\cite[Propositions 1.3.59, 1.3.63]{Anel2013}).
\end{ex}

\begin{rem}
	Let $V$ be a graded Banach $R$-module. Given a bounded coderivation
	$\nu \colon \tens{V} \rightharpoonup \tens{V}$ of the tensor coalgebra, we can complete it
	and obtain a coderivation $\widehat{\nu} \colon \tensf{V} \rightharpoonup \tensf{V}$
	of the formal tensor coalgebra.
	Conversely, the classification of coderivations of $\tensf{V}$
	shows in particular that any coderivation $\mu \colon \tensf{V} \rightharpoonup \tensf{V}$
	is the completion of a unique bounded coderivation $\nu \colon \tens{V} \rightharpoonup \tens{V}$
	(i.e., $\widehat{\nu} = \mu$).
	This is in contrast to the situation with coalgebra morphisms
	(see \cref{rem:not-every-coalgebra-morphism-formal-is-completion}).
\end{rem}

The classification of $R$-linear coderivations on $\tensf{V}$ readily extends to a classification
of generalized coderivations on $\tensf{V}$. Before stating the result,
let us characterize explicitly what the corestriction of
a generalized coderivation looks like. Given a coderivation $\mu \colon \tensf{V} \rightharpoonup \tensf{V}$ over an algebra derivation $d \colon R \rightharpoonup R$, the corestriction $\corest{\mu} = \pi_1 \circ \mu$ is a graded $\mathbbm{k}$-linear map which satisfies
\begin{equation*}
	\corest{\mu} \left( r \cdot l \right) = dr \cdot \pi_1 \left( l \right) +
	(-1)^{\braidd{\corest{\mu}}{r}} r \cdot \corest{\mu} \left( l \right)
\end{equation*}
for all $l \in \tensf{V}$ and $r \in R$. The map $\corest{\mu}$ is determined uniquely by the sequence of maps
$\mu_k \colon V^{\cotimes k} \rightharpoonup V$ and in terms of the maps $\mu_k$ we must have
\begin{equation*}
	\mu_k \left( r \cdot l \right) =
	\delta_{k,1} \left( dr \cdot l \right) + (-1)^{\braidd{\mu_k}{r}} r \cdot \mu_k \left( l \right)
\end{equation*}
for all $l \in V^{\cotimes k}$ and $r \in R$.
From this we see that the maps $\mu_k$ for $k \neq 1$ are actually $R$-linear and the map
$\mu_1 \colon V \rightharpoonup V$ is a module derivation on $V$ over $d$. In particular, the map
$\mu_0 \colon R \rightarrow V$ is also $R$-linear and determined uniquely by the element
$\mu_0(1) \in V$, which is often identified with the map itself.

\begin{prop} \label{prop:classification-generalized-coderivations-formal-tensor-coalgebra}
	Let $R$ be a graded-commutative Banach $\mathbbm{k}$-algebra and let $V$ be a graded Banach $R$-module.
	Let $d \colon R \rightharpoonup R$ be a bounded algebra derivation.
	We have a bijective correspondence between
	\begin{enumerate}
		\item Banach coalgebra coderivations $\mu \colon \tensf{V} \rightharpoonup \tensf{V}$
		      over $d$.
		\item Graded bounded $\mathbbm{k}$-linear maps $\corest{\mu} \colon \tensf{V} \rightharpoonup V$
		      of degree $\degb{d}$ which satisfy
		      \begin{equation} \label{eq:corest-over-derivation-projection}
			      \corest{\mu} \left( r \cdot l \right) = dr \cdot \pi_1 \left( l \right)
			      + (-1)^{\braidd{\corest{\mu}}{r}} r \cdot \corest{\mu} \left( l \right)
		      \end{equation}
		      for all $l \in \tensf{V}$ and $r \in R$, where $\pi_1 \colon \tensf{V} \rightarrow V$
		      is the canonical projection.
		\item Sequences $( \mu_k \colon V^{\cotimes k} \rightharpoonup V )_{k \geq 0}$ of graded
		      bounded $\mathbbm{k}$-linear maps of degree $\degb{d}$
		      which satisfy $\sup_{k \geq 0} \nnorm[\mu_k] < \infty$
		      and such that $\mu_1$ is a module derivation over $d$ while $\mu_k$ for $k \neq 1$ are
		      $R$-linear.
		\item Sequences $( U_k \colon V^{\times k} \rightharpoonup V )_{k \geq 0}$ of graded bounded
		      $\mathbbm{k}$-multilinear maps of degree $\degb{d}$
		      which satisfy $\sup_{k \geq 0} \nnorm[U_k] < \infty$
		      and such that $U_1$ is a module derivation over $d$ while $U_k$ for $k \neq 1$
		      are $R$-multilinear.
	\end{enumerate}
	Given a coderivation $\mu \colon \tensf{V} \rightharpoonup \tensf{V}$ over $d$,
	the map	$\corest{\mu} \colon \tensf{V} \rightharpoonup V$ is the corestriction
	$\corest{\mu} \defeq \pi_1 \circ \mu$ of $\mu$. Conversely, given a map
	$\corest{\mu} \colon \tensf{V} \rightharpoonup V$ satisfying $(2)$,
	there exists a unique Banach coalgebra coderivation $\mu \colon \tensf{V} \rightharpoonup \tensf{V}$
	over $d$, given by the formula
	\begin{equation} \label{eq:generalized-coder-coextension}
		\begin{aligned}
			\mu(l) & = (-1)^{\braidd{\corest{\mu}}{l_{(1)}}} l_{(1)} \cotimes
			\corest{\mu} \left( l_{(2)} \right) \cotimes l_{(3)},
			       & l \in V^{\cotimes k}, k \geq 1,
			\\
			\mu(r) & = dr + (-1)^{\braidd{\mu}{r}} r \cdot \mu_0(1) = dr + \mu_0(r),
			       & r \in R,
		\end{aligned}
	\end{equation}
	whose corestriction is $\corest{\mu}$.
\end{prop}
\begin{proof}
	Let $\mu,\nu \colon \tensf{V} \rightharpoonup \tensf{V}$ be two coderivations over the same derivation $d$
	such that $\corest{\mu} = \corest{\nu}$. Then $\mu - \nu$ is an $R$-linear coderivation
	which satisfies
	\begin{equation*}
		\corest{\left(\mu - \nu \right)} = \pi_1 \circ \left( \mu - \nu \right) =
		\pi_1 \circ \mu - \pi_1 \circ \nu = \corest{\mu} - \corest{\nu} = 0
	\end{equation*}
	and hence by \cref{lm:coderivation-determined-by-corestriction} we have $\mu - \nu = 0$, i.e., $\mu = \nu$.

	Conversely, let $\varphi \colon \tensf{V} \rightharpoonup V$ be a graded bounded $\mathbbm{k}$-linear map
	satisfying
	\begin{equation*}
		\varphi \left( r \cdot l \right) = dr \cdot \pi_1 \left( l \right) +
		(-1)^{\braidd{\varphi}{r}} r \cdot \varphi \left( l \right)
	\end{equation*}
	for all $l \in \tensf{V}$ and $r \in R$. We can construct a coderivation
	$\mu \colon \tensf{V} \rightharpoonup \tensf{V}$ over $d$ such that $\corest{\mu} = \varphi$ in two steps
	by treating the $R$-linear part and the $d$-operator part separately:
	\begin{enumerate}
		\item Define a map $\mu^d \colon \tensf{V} \rightharpoonup \tensf{V}$ over $d$
		      by the formula
		      \begin{equation*}
			      \mu^d \left( l \right) =
			      \begin{cases}
				      \sum_{i=0}^k \id^{\cotimes i} \cotimes \varphi_1 \cotimes \id^{\cotimes (k - i)}
				                         & l \in V^{\cotimes (k+1)}, \, k \geq 0, \\
				      d \left( r \right) & l = r \in V^{\cotimes 0} = R.
			      \end{cases}
		      \end{equation*}
		      Then $\mu^d$ is a coderivation over $d$ with $\mu^d_1 = \varphi_1$ while $\mu^d_k = 0$
		      for $k \neq 1$ (since $\mu^d$ is weight-preserving).
		\item Using \cref{prop:classification-coderivations-formal-tensor-coalgebra}, construct
		      an $R$-linear coderivation $\mu^R \colon \tensf{V} \rightharpoonup \tensf{V}$ such that
		      $\mu^R_k = \varphi_k$ for $k \neq 1$ while $\mu^R_1 = 0$.
	\end{enumerate}
	Then $\mu = \mu^R + \mu^d$ is a coderivation over $d$ with $\corest{\mu} = \varphi$.

	Hence, we have shown the equivalence between $(1)$ and $(2)$. The equivalence between $(2),(3)$ and $(4)$
	follows from the universal properties of the complete direct sum and tensor product.
\end{proof}

\begin{rem} \label{rem:coderivation-formula-d-operator-abuse}
	We note that the formula \eqref{eq:generalized-coder-coextension} for $\mu$ in terms
	of the corestriction $\corest{\mu}$, involves some abuse of notation which we now explain.
	Since $\corest{\mu}$ is not $R$-linear, the expression
	\begin{equation*}
		l \mapsto (-1)^{\braidd{\corest{\mu}}{l_{(1)}}} l_{(1)} \cotimes_R
		\corest{\mu} \left( l_{(2)} \right) \cotimes_R l_{(3)}
	\end{equation*}
	is a priori ill-defined. However, if we write it out explicitly for an elementary tensor
	$l = v_1 \cotimes_R \dots \cotimes_R v_k \in V^{\cotimes k}$, we see that it is the sum of
	\begin{equation*}
		\sum_{i=1}^k (-1)^{\braid{\degb{\mu_1}}{\degb{v_1} + \dots + \degb{v_{i-1}}}}
		v_1 \cotimes_R \dots \cotimes_R v_{i-1} \cotimes_R \mu_1 \left( v_i \right)
		\cotimes_R v_{i+1} \cotimes_R \dots \cotimes_R v_k
	\end{equation*}
	and
	\begin{equation*}
		\begin{aligned}
			\sum_{\substack{k_1 + k_2 + k_3 = k                             \\ k_1, k_2, k_3 \geq 0 \\ k_2 \neq 1}}
			 & (-1)^{\braid{\degb{\mu_{k_2}}}{\sum_{i=1}^{k_1} \degb{v_i}}}
			\\
			 & \quad
			v_1 \cotimes_R \dots \cotimes_R v_{k_1} \cotimes_{R}
			\mu_{k_2} \left( v_{k_1 + 1} \cotimes_R \dots \cotimes_R v_{k_1 + k_2} \right) \cotimes_{R}
			v_{k_1 + k_2 + 1} \cotimes_R \dots \cotimes_R v_k.
		\end{aligned}
	\end{equation*}
	The first term involves only $\mu_1$ and is well-defined since $\mu_1$ is a
	$d$-operator. In fact, if we think of $\mathcal{V} \defeq \left( V, \mu_1 \right)$
	as a pre-differential graded Banach module over $\mathcal{R} \defeq \left( R, d \right)$, then
	the first term is the pre-differential of the tensor product $\mathcal{V}^{\cotimes_{\mathcal{R}} k}$
	(see also \cref{rem:tensor-product-d-operators-abuse}).
	The second term involves only the operators $\mu_k$ for $k \neq 1$ which are $R$-linear and so it is
	clearly well-defined.
\end{rem}

We end this section with the notion of an $(f,g)$-coderivation and its classification.

\begin{dfn} \label{dfn:f-g-coderivation}
	Let $R$ be a graded-commutative Banach $\mathbbm{k}$-algebra and let $C,D$ be two graded Banach $R$-coalgebras.
	Let $f, g \colon C \rightarrow D$ be two morphisms of Banach $R$-coalgebras.
	A graded bounded $R$-linear map $\eta \colon C \rightharpoonup D$ is called an
	$(f, g)$-\textbf{coderivation} if it satisfies the generalized co-Leibniz rule
	\begin{equation*}
		\Delta_D \circ \eta = \left( f \cotimes \eta + \eta \cotimes g \right) \circ \Delta_C.
	\end{equation*}
\end{dfn}

\begin{prop} \label{prop:classification-f-g-coderivations}
	Let $R$ be a graded-commutative Banach $\mathbbm{k}$-algebra and let $W$ be a graded Banach $R$-module.
	Let $C$ be a graded Banach $R$-coalgebra and let $f,g \colon C \rightarrow \tensf{W}$ be two morphisms of
	Banach $R$-coalgebras. Then any $(f, g)$-coderivation $\eta \colon C \rightharpoonup \tensf{W}$
	is determined uniquely by its corestriction $\corest{\eta} \colon C \rightharpoonup W$ via the formula
	\begin{equation*}
		\eta = m \circ \left( f \cotimes \corest{\eta} \cotimes g \right) \circ \Delta^3_C.
	\end{equation*}
\end{prop}
\begin{proof}
	We can consider $C$ as a Banach $\tensf{W}$-bicomodule via the coactions
	\begin{equation*}
		\Delta^{1|0} = \left( f \cotimes \id \right) \circ \Delta_C, \qquad
		\Delta^{0|1} = \left( \id \cotimes g \right) \circ \Delta_C.
	\end{equation*}
	Then $\eta \colon C \rightharpoonup \tensf{W}$ becomes a coderivation in the sense of
	\cref{dfn:coderivation-domain-bicomodule} and hence by \cref{lm:coderivation-from-corestriction},
	the coderivation $\eta$ is determined uniquely by $\corest{\eta}$ via the formula
	\begin{equation*}
		\eta = m \circ \left( \id \cotimes \corest{\eta} \cotimes \id \right) \circ \Delta^{1|1}
		= m \circ \left( f \cotimes \corest{\eta} \cotimes g \right) \circ \Delta^3_C
	\end{equation*}
	since $\Delta^{1|1} = \left( f \cotimes \id_C \cotimes g \right) \circ \Delta^3_C$.
\end{proof}

\begin{cor} \label{cor:f-circ-mu-nu-circ-f-corest}
	Let $R$ be a graded-commutative Banach $\mathbbm{k}$-algebra and let $V, W$ be two graded Banach $R$-modules.
	Let $f \colon \tensf{V} \rightarrow \tensf{W}$ be a morphism of graded Banach $R$-coalgebras and let
	$\mu \colon \tensf{V} \rightharpoonup \tensf{V}$ and $\nu \colon \tensf{W} \rightharpoonup \tensf{W}$
	be two generalized coderivations over the same derivation $d \colon R \rightharpoonup R$.
	Then $f \circ \mu = \nu \circ f$ if and only if $\corest{f} \circ \mu = \corest{\nu} \circ f$.
\end{cor}
\begin{proof}
	The map $\eta \defeq f \circ \mu - \nu \circ f$
	is $R$-linear and an $(f,f)$-coderivation. Hence by \cref{prop:classification-f-g-coderivations},
	it is determined uniquely by the corestriction $\corest{\eta} = \corest{f} \circ \mu - \corest{\nu} \circ f$.
\end{proof}

\begin{rem}
	We note that the condition $f \circ \mu = \nu \circ f$ appearing in \cref{cor:f-circ-mu-nu-circ-f-corest}
	can hold only if $\mu$ and $\nu$ are coderivations over the \textit{same} derivation. To see this,
	we can compose both sides of the equality with the counit
	$\varepsilon_{\tensf{W}} \colon \tensf{W} \rightarrow R$ and the inclusion $i_R \colon R \rightarrow \tensf{V}$
	to obtain
	\begin{align*}
		\varepsilon_{\tensf{W}} \circ f \circ \mu \circ i_R & = \varepsilon_{\tensf{V}} \circ \mu \circ i_R =
		d_{\mu} \circ \varepsilon_{\tensf{V}} \circ i_R = d_{\mu},
		\\
		\varepsilon_{\tensf{W}} \circ \nu \circ f \circ i_R & = d_{\nu} \circ \varepsilon_{\tensf{W}} \circ f \circ i_R =
		d_{\nu} \circ \varepsilon_{\tensf{V}} \circ i_R = d_{\nu}.
	\end{align*}
\end{rem}

\subsection{Working over Different Ground Algebras} \label{sec:formal-tensor-coalgebra-different-ground}

Let $R,S$ be two graded-commutative Banach $\mathbbm{k}$-algebras and let
$\varphi \colon R \rightarrow S$ be a morphism of graded Banach $\mathbbm{k}$-algebras.
Consider the monoidal adjunction
\begin{equation*}
	\varphi_{!} \colon \GBMod[R] \stackrel[]{\dashv}{\rightleftarrows} \GBMod[S] \colon \varphi^{*}.
\end{equation*}
Since the scalar extension functor
$\varphi_{!}$ is strong monoidal, the scalar extension of
a graded Banach $R$-algebra (resp.\ coalgebra) has a natural structure of a graded Banach $S$-algebra
(resp.\ coalgebra). In contrast, the scalar restriction functor $\varphi^{*}$
is only lax monoidal and hence, the scalar restriction of a graded $S$-algebra has a natural structure
of a graded Banach $R$-algebra, but the same is not true for coalgebras.

In this work, we will work only with formal tensor coalgebras and not general Banach coalgebras.
Even though scalar restriction does not make sense for general Banach coalgebras, there is a reasonable
notion of scalar restriction for formal tensor coalgebras, possibly equipped with a generalized coderivation,
which we describe in \cref{sec:scalar-restriction-formal-tensor-coalgebras}.
In \cref{sec:scalar-extension-formal-tensor-coalgebras}, we discuss the notion of scalar extension
for tensor coalgebras, which coincides up to an isomorphism with the standard
scalar extension of coalgebras, and the relation between both constructions.

\subsubsection{Scalar Restriction for Formal Tensor Coalgebras}
\label{sec:scalar-restriction-formal-tensor-coalgebras}

Let $V$ be a graded Banach $R$-module and let $W$ be a graded Banach $S$-module. Given
the morphism $\varphi \colon R \rightarrow S$, there exists a universal Banach coalgebra morphism
$\resover{\varphi} = \resover{\varphi}^{W} \colon
	\tensf{\varphi^{*} \left( W \right)}[R] \rightarrow \tensf{W}[S]$ over $\varphi$ given by
\begin{align}\label{eq:canonical-pullback-map}
	 & \resover{\varphi} \left( r \right) \defeq \varphi(r),
	 & \resover{\varphi} \left( w_1 \cotimes_R \dots \cotimes_R w_n \right) \defeq
	w_1 \cotimes_S \dots \cotimes_S w_n.
\end{align}
Note that the only non-zero component of the corestriction $\corest{\resover{\varphi}}$ is given by
$\resover{\varphi}_1 \colon \varphi^{*} \left( W \right) \rightarrow W$ which is the canonical (identity) map
over $\varphi$.

\begin{dfn}
	The formal tensor coalgebra $\tensf{\varphi^{*} \left( W \right)}[R]$, together with the map
	$\resover{\varphi}$, is called the \textbf{pullback} (or \textbf{scalar restriction})
	of the formal tensor coalgebra $\tensf{W}[S]$ along $\varphi$.
\end{dfn}

The pullback construction satisfies the following universal property: Any Banach coalgebra morphism
$f \colon \tensf{V}[R] \rightarrow \tensf{W}[S]$ over $\varphi$ factors uniquely as
$f = \resover{\varphi} \circ \rescoho{f}$ where
$\rescoho{f} \colon \tensf{V}[R] \rightarrow \tensf{\varphi^{*} \left( W \right)}[R]$ is a morphism
of graded Banach $R$-coalgebras (see \cref{fig:scalar-rest-formal-tensor-coalg}).
\begin{figure}
	\begin{tikzcd}
		{\tensf{V}[R]} & {\tensf{\varphi^{*} \left( W \right)}[R]} & {\tensf{W}[S]} \\
		R & R & S
		\arrow[from=1-1, to=2-1]
		\arrow[from=1-2, to=2-2]
		\arrow[from=1-3, to=2-3]
		\arrow["{\exists ! \rescoho{f}}", dashed, from=1-1, to=1-2]
		\arrow["{\resover{\varphi}}", from=1-2, to=1-3]
		\arrow["{\id}", from=2-1, to=2-2]
		\arrow["\varphi", from=2-2, to=2-3]
		\arrow["f", bend left, from=1-1, to=1-3]
	\end{tikzcd}
	\caption{Universal property of the pullback of formal tensor coalgebras.}
	\small\textsuperscript{The unnamed vertical maps are the counit maps.}
	\label{fig:scalar-rest-formal-tensor-coalg}
\end{figure}

The corestriction $\corest{\rescoho{f}}$ of $\rescoho{f}$, and hence $\rescoho{f}$ itself,
is determined uniquely by the equation $\resover{\varphi}_1 \circ \corest{\rescoho{f}} = \corest{f}$.

\begin{rem}
	Since $\varphi^{*}$ is only lax monoidal, in general, the $R$-module $\varphi^{*} \left( \tensf{W}[S] \right)$
	has no natural structure of an $R$-coalgebra. This is the reason why we pull
	back the underlying module $W$ and work with $\tensf{\varphi^{*} \left( W \right)}[R]$
	instead of pulling back the whole tensor coalgebra and working with $\varphi^{*} \left( \tensf{W}[S] \right)$.
	See also \cref{rem:pullback-terminology-tensor-coalgebras}.
\end{rem}

The pullback construction for formal tensor coalgebras extends also to morphisms as follows.
Given a morphism
$f \colon \tensf{W_1}[S] \rightarrow \tensf{W_2}[S]$ of graded Banach $S$-coalgebras, there exists
a unique morphism $\varphi^{*} \left( f \right) \colon
	\tensf{\varphi^{*} \left( W_1 \right)}[R] \rightarrow \tensf{\varphi^{*} \left( W_2 \right)}[R]$
of graded Banach $R$-coalgebras such that
$\resover{\varphi}^{W_2} \circ \varphi^{*} \left( f \right) = f \circ \resover{\varphi}^{W_1}$
(see \cref{fig:scalar-rest-morphism-formal-tensor-coalg}). The morphism $\varphi^{*} \left( f \right)$
is defined via its corestriction
$\corest{\varphi^{*} \left( f \right)} \colon
	\tensf{\varphi^{*} \left( W_1 \right)}[R] \rightarrow \varphi^{*} \left( W_2 \right)$ by setting
\begin{align*}
	 & \varphi^{*} \left( f \right)_{0} \left( r \right) \defeq f_0 \left( \varphi \left( r \right) \right),
	 & \varphi^{*} \left( f \right)_{k} \left( x_1 \cotimes_R \dots \cotimes_R x_k \right) \defeq
	f_k \left( x_1 \cotimes_S \dots \cotimes_S x_k \right),
\end{align*}
for $r \in R$ and $x_1, \dots, x_k \in W_1$. Note that the corestriction $\corest{\varphi^{*} \left( f \right)}$
of $\varphi^{*} \left( f \right)$, and hence $\varphi^{*} \left( f \right)$ itself, is determined uniquely
by the equation
$\resover{\varphi}_1 \circ \corest{\varphi^{*} \left( f \right)} = \corest{f} \circ \resover{\varphi}$.

\begin{figure}
	\begin{tikzcd}
		{\tensf{W_1}[S]} && {\tensf{W_2}[S]} \\
		& {S} \\
		& {R} \\
		{\tensf{\varphi^{*} \left( W_1 \right)}[R]} &&
		{\tensf{\varphi^{*} \left( W_2 \right)}[R]}
		\arrow["f", from=1-1, to=1-3]
		\arrow["\varphi", from=3-2, to=2-2]
		\arrow[from=1-1, to=2-2] 
		\arrow[from=1-3, to=2-2] 
		\arrow["{\resover{\varphi}^{W_1}}", from=4-1, to=1-1]
		\arrow["{\resover{\varphi}^{W_2}}"', from=4-3, to=1-3]
		\arrow[from=4-1, to=3-2] 
		\arrow[from=4-3, to=3-2] 
		\arrow["{\exists ! \varphi^{*} \left( f \right)}"', dashed, from=4-1, to=4-3]
	\end{tikzcd}
	\caption{Pullback of a morphism between formal tensor $S$-coalgebras.}
	\small\textsuperscript{The unnamed maps are the counit maps.}
	\label{fig:scalar-rest-morphism-formal-tensor-coalg}
\end{figure}

\begin{dfn} \label{dfn:pullback-morphism-tensor-coalgebras}
	The morphism $\varphi^{*} \left( f \right) \colon \tensf{\varphi^{*} \left( W_1 \right)}[R]
		\rightarrow \tensf{\varphi^{*} \left( W_2 \right)}[R]$ is called the \textbf{pullback} (or
	\textbf{scalar restriction}) of $f \colon \tensf{W_1}[S] \rightarrow \tensf{W_2}[S]$ along $\varphi$.
\end{dfn}

The pullback construction $\varphi^{*}$ defines a covariant functor from the category of
formal tensor $S$-coalgebras to the category of formal tensor $R$-coalgebras. In addition,
the pullback is functorial with respect to morphisms of the underlying graded Banach $\mathbbm{k}$-algebras.
More precisely, given a morphism $\psi \colon P \rightarrow R$ of graded Banach $\mathbbm{k}$-algebras,
we have the identity $\psi^{*} \circ \varphi^{*} = \left( \varphi \circ \psi \right)^{*}$ as functors from
the category of formal tensor $S$-coalgebras to the category of formal tensor $P$-coalgebras.

Next, we can naturally extend the pullback construction to formal tensor coalgebras equipped with a generalized
coderivation. Before describing the extension, let us set up some notation. Given
a pre-differential graded Banach $\mathbbm{k}$-algebra $\mathcal{R} = \left( R, d_R \right)$,
a \textbf{formal tensor} $\mathcal{R}$\textbf{-coalgebra} is a pair $\big( \tensf{V}[R], \mu \big)$,
where $\tensf{V}[R]$ is the formal tensor coalgebra on a graded Banach $R$-module $V$ and $\mu$ is a coderivation on $\tensf{V}[R]$ over $d_R$. Note that a formal tensor $\mathcal{R}$-coalgebra is
a pre-differential graded Banach $\mathcal{R}$-coalgebra, i.e., a coalgebra object
of the monoidal category $\PDGBMod[\mathcal{R}]$.
A morphism $f \colon \big( \tensf{V}[R], \mu \big) \rightarrow \big( \tensf{W}[R], \nu \big)$ of formal tensor
$\mathcal{R}$-coalgebras is a morphism of the pre-differential graded Banach $\mathcal{R}$-coalgebras, i.e.,
a morphism $f \colon \tensf{V}[R] \rightarrow \tensf{W}[R]$ of graded Banach $R$-coalgebras such
that $f \circ \mu = \nu \circ f$.

Let $\mathcal{R} = \left( R, d_R \right)$ and $\mathcal{S} = \left( S, d_S \right)$ be
two pre-differential graded-commutative Banach $\mathbbm{k}$-algebras
and let $\varphi \colon \mathcal{R} \rightarrow \mathcal{S}$ be a morphism of pre-differential graded Banach
$\mathbbm{k}$-algebras. Let $\big( \tensf{W}[S], \nu \big)$ be a formal tensor $\mathcal{S}$-coalgebra.
Then there exists a unique coderivation
$\varphi^{*} \left( \nu \right) \colon \tensf{\varphi^{*} \left( W \right)}[R] \rightharpoonup
	\tensf{\varphi^{*} \left( W \right)}[R]$
over $d_R$ with $\degb{\varphi^{*} \left( \nu \right)} = \degb{\nu}$ such that
$\resover{\varphi} \circ \varphi^{*} \left( \nu \right) = \nu \circ \resover{\varphi}$
(see \cref{fig:pullback-of-coderivation}). In other words, the coderivation
$\varphi^{*} \left( \nu \right)$ is the unique coderivation with respect to which the morphism
$\resover{\varphi}$ becomes a morphism
\begin{equation*}
	\resover{\varphi} \colon
	\left( \tensf{\varphi^{*} \left( W \right)}[R], \varphi^{*} \left( \nu \right) \right) \rightarrow
	\left( \tensf{W}[S], \nu \right)
\end{equation*}
of pre-differential graded Banach coalgebras over $\varphi$.

The coderivation
$\varphi^{*} \left( \nu \right)$ is defined via its corestriction
$\corest{\varphi^{*} \left( \nu \right)} \colon \tensf{\varphi^{*} \left( W \right)}[R] \rightharpoonup
	\varphi^{*} \left( W \right)$ by setting
\begin{align*}
	\varphi^{*} \left( \nu \right)_0 \left( r \right)                                   & \defeq
	\nu_0 \left( \varphi \left( r \right) \right),
	\\
	\varphi^{*} \left( \nu \right)_n \left( w_1 \cotimes_R \dots \cotimes_R w_n \right) & \defeq
	\nu_n \left( w_1 \cotimes_S \dots \cotimes_S w_n \right)
\end{align*}
for $r \in R$ and $w_1, \dots, w_n \in W$.
Note that $\varphi^{*} \left( \nu \right)_1 \colon \varphi^{*} \left( W \right)
	\rightharpoonup \varphi^{*} \left( W \right)$ is a module derivation over $d_R$ as
we have\footnote{Here, given $r \in R$ and $w \in \varphi^{*} \left( W \right)$, we use
	$r \triangleleft w = \varphi \left( r \right) \cdot w$ to denote the action of $R$
	on $\varphi^{*} \left( W \right)$.}
\begin{equation*}
	\begin{aligned}
		\varphi^{*} \left( \nu \right)_1 \left( r \triangleleft w \right) & =
		\varphi^{*} \left( \nu \right)_1 \left( \varphi \left( r \right) \cdot w \right) =
		\nu_1 \left( \varphi(r) \cdot w \right)
		\\
		                                                                  & =
		d_S \left( \varphi(r) \right) \cdot w +
			                                    (-1)^{\braidd{\nu}{\varphi(r)}} \varphi(r) \cdot \nu_1 \left( w \right)
		\\
		                                                                  & = \varphi \left( d_R \left( r \right) \right) \cdot w +
			                                                                                                                        (-1)^{\braidd{\nu}{r}} \varphi(r) \cdot \nu_1 \left( w \right)
		\\
		                                                                  & =
		d_R \left( r \right) \triangleleft w +
			                                   (-1)^{\braidd{\varphi^{*} \left( \nu \right)}{r}} r \triangleleft \nu_1 \left( w \right)
	\end{aligned}
\end{equation*}
and that $\varphi^{*} \left( \nu \right)_k \colon \varphi^{*} \left( W \right)^{\cotimes_R k}
	\rightharpoonup \varphi^{*} \left( W \right)$ for $k \neq 1$ are $R$-linear operators, so the definition is indeed valid.
Note also that the corestriction $\corest{\varphi^{*} \left( \nu \right)}$ of
$\varphi^{*} \left( \nu \right)$, and hence $\varphi^{*} \left( \nu \right)$ itself, is determined uniquely
by the equation
$\resover{\varphi}_1 \circ \corest{\varphi^{*} \left( \nu \right)} = \corest{\nu} \circ \resover{\varphi}$.

\begin{dfn} \label{dfn:pullback-generalized-coderivation}
	The coderivation $\varphi^{*} \left( \nu \right)$ is called the \textbf{pullback} (or
	\textbf{scalar restriction}) of $\nu$ along $\varphi$. The formal tensor $\mathcal{R}$-coalgebra
	$\big( \tensf{\varphi^{*} \left( W \right)}[R], \varphi^{*} \left( \nu \right) \big)$, together
	with the map $\resover{\varphi}$, is called the \textbf{pullback} (or \textbf{scalar restriction}) of
	$\big( \tensf{W}[S], \nu \big)$ along $\varphi$.
\end{dfn}

\begin{figure}
	\begin{tikzcd}
		R &&& S \\
		& {\tensf{\varphi^{*} \left( W \right)}[R]} & {\tensf{W}[S]} \\
		& {\tensf{\varphi^{*} \left( W \right)}[R]} & {\tensf{W}[S]} \\
		R &&& S
		\arrow["\nu", harpoon, from=2-3, to=3-3]
		\arrow["{\resover{\varphi}}", from=2-2, to=2-3]
		\arrow[from=2-3, to=1-4]
		\arrow["{d_S}", harpoon, from=1-4, to=4-4]
		\arrow["\varphi", from=1-1, to=1-4]
		\arrow[from=2-2, to=1-1]
		\arrow["{\exists ! \, \varphi^{*} \left( \nu \right)}"', dotted, harpoon, from=2-2, to=3-2]
		\arrow["{\resover{\varphi}}"', from=3-2, to=3-3]
		\arrow["{d_R}"', harpoon, from=1-1, to=4-1]
		\arrow[from=3-2, to=4-1]
		\arrow["\varphi"', from=4-1, to=4-4]
		\arrow[from=3-3, to=4-4]
	\end{tikzcd}
	\caption{The pullback of a generalized coderivation.}
	\small\textsuperscript{The unnamed maps are the counit maps.}
	\label{fig:pullback-of-coderivation}
\end{figure}

The universal property of the pullback extends naturally to formal tensor coalgebras equipped
with generalized coderivations.
Namely, if $\big( \tensf{V}[R], \mu \big)$ is a formal tensor $\mathcal{R}$-coalgebra and
$\big( \tensf{W}[S], \nu \big)$ is a formal tensor $\mathcal{S}$-coalgebra then
any morphism
$f \colon \big( \tensf{V}[R], \mu \big) \rightarrow \big( \tensf{W}[S], \nu \big)$ of
pre-differential graded Banach coalgebras over $\varphi \colon \mathcal{R} \rightarrow \mathcal{S}$ factors
uniquely as $f = \resover{\varphi} \circ \rescoho{f}$ where
$\rescoho{f} \colon \big( \tensf{V}[R], \mu \big) \rightarrow
	\big( \tensf{\varphi^{*} \left( W \right)}[R], \varphi^{*} \left( \nu \right) \big)$
is a morphism of pre-differential graded Banach $\mathcal{R}$-coalgebras
(see \cref{fig:scalar-rest-formal-tensor-coalg-over-pre-cdga}). The map $\rescoho{f}$ is the
same map as in \cref{fig:scalar-rest-formal-tensor-coalg}.
To verify that $\rescoho{f}$ commutes with the coderivations, it is enough to check that
$\corest{\rescoho{f}} \circ \mu = \corest{\varphi^{*} \left( \nu \right)} \circ \rescoho{f}$
(see \cref{cor:f-circ-mu-nu-circ-f-corest}). This follows from the chain of equalities
\begin{equation*}
	\resover{\varphi}_1 \circ \corest{\rescoho{f}} \circ \mu =
	\corest{f} \circ \mu = \corest{\nu} \circ f =
	\corest{\nu} \circ \resover{\varphi} \circ \rescoho{f} =
	\resover{\varphi}_1 \circ \corest{\varphi^{*} \left( \nu \right)} \circ \rescoho{f}.
\end{equation*}

\begin{figure}
	\begin{tikzcd}
		{\left( \tensf{V}[R], \mu \right)} &
		{\left( \tensf{\varphi^{*} \left( W \right)}[R], \varphi^{*} \left( \nu \right) \right)} &
		{\left( \tensf{W}[S], \nu \right)} \\
		(R, d_R) & (R, d_R) & (S, d_S)
		\arrow[from=1-1, to=2-1] 
		\arrow[from=1-2, to=2-2] 
		\arrow[from=1-3, to=2-3] 
		\arrow["{\exists ! \rescoho{f}}", dashed, from=1-1, to=1-2]
		\arrow["{\resover{\varphi}}", from=1-2, to=1-3]
		\arrow["{\id}", from=2-1, to=2-2]
		\arrow["\varphi", from=2-2, to=2-3]
		\arrow["f", bend left, from=1-1, to=1-3]
	\end{tikzcd}
	\caption{Universal property of the pullback of formal tensor coalgebras.}
	\small\textsuperscript{The unnamed vertical maps are the counit maps.}
	\label{fig:scalar-rest-formal-tensor-coalg-over-pre-cdga}
\end{figure}

Given a morphism $\varphi \colon \mathcal{R} \rightarrow \mathcal{S}$ of pre-differential
graded Banach $\mathbbm{k}$-algebras, the pullback construction
$\varphi^{*}$ defines a covariant functor from the category of formal tensor $\mathcal{S}$-coalgebras
to the category of formal tensor $\mathcal{R}$-coalgebras. That is,
given a morphism $f \colon \big( \tensf{W_1}[S], \nu^1 \big) \rightarrow
	\big( \tensf{W_2}[S], \nu^2 \big)$ of formal tensor $\mathcal{S}$-coalgebras, i.e.,
$f$ satisfies $f \circ \nu^1 = \nu^2 \circ f$, the pullback morphism
$\varphi^{*} \left( f \right) \colon \tensf{\varphi^{*} \left( W_1 \right)}[R]
	\rightarrow \tensf{\varphi^{*} \left( W_2 \right)}[R]$ (see \cref{dfn:pullback-morphism-tensor-coalgebras})
also satisfies $\varphi^{*} \left( \nu^2 \right) \circ \varphi^{*} \left( f \right)
	= \varphi^{*} \left( f \right) \circ \varphi^{*} \left( \nu^1 \right)$. See
\cref{fig:scalar-rest-morphism-formal-tensor-coalg-with-coderivations}.

\begin{figure}
	\begin{tikzcd}
		{\left( \tensf{W_1}[S], \nu^1 \right)} &&
		{\left( \tensf{W_2}[S], \nu^2 \right)} \\
		& {(S,d_S)} \\
		& {(R,d_R)} \\
		{\left( \tensf{\varphi^{*} \left( W_1 \right)}[R], \varphi^{*} \left( \nu^1 \right) \right)} &&
		{\left( \tensf{\varphi^{*} \left( W_2 \right)}[R], \varphi^{*} \left( \nu^2 \right) \right)}
		\arrow["f", from=1-1, to=1-3]
		\arrow["\varphi", from=3-2, to=2-2]
		\arrow[from=1-1, to=2-2] 
		\arrow[from=1-3, to=2-2] 
		\arrow["{\resover{\varphi}^{W_1}}", from=4-1, to=1-1]
		\arrow["{\resover{\varphi}^{W_2}}"', from=4-3, to=1-3]
		\arrow[from=4-1, to=3-2] 
		\arrow[from=4-3, to=3-2] 
		\arrow["{\exists ! \varphi^{*} \left( f \right)}"', dashed, from=4-1, to=4-3]
	\end{tikzcd}
	\caption{Pullback of a morphism between formal tensor $\mathcal{S}$-coalgebras.}
	\small\textsuperscript{The unnamed maps are the counit maps.}
	\label{fig:scalar-rest-morphism-formal-tensor-coalg-with-coderivations}
\end{figure}

In addition, the pullback is functorial with respect to morphisms of the underlying pre-differential
graded Banach $\mathbbm{k}$-algebras. More precisely, given a morphism $\psi \colon \mathcal{P} \rightarrow \mathcal{R}$
of pre-differential graded Banach $\mathbbm{k}$-algebras, we have the identity
$\psi^{*} \circ \varphi^{*} = \left( \varphi \circ \psi \right)^{*}$ as functors from
the category of formal tensor $\mathcal{S}$-coalgebras to the category of formal tensor $\mathcal{P}$-coalgebras.

\begin{rem} \label{rem:pullback-terminology-tensor-coalgebras}
	Our notion of pullback is motivated by the bar construction, in which the structure of
	a graded algebra is encoded by a coderivation on a tensor coalgebra.
	Given a graded-commutative $\mathbbm{k}$-algebra $S$ and a graded $S$-module $A$, endowing
	$A$ with the structure of a graded (non-unital) $S$-algebra is equivalent under the bar construction
	to providing a coderivation $\nu$ on the tensor coalgebra $\tens{A[1]}[S]$ which satisfies $\nu^2 = 0$
	and whose only non-zero component is $\nu_2$. Now, given a morphism $\varphi \colon R \rightarrow S$,
	we know we can endow the pullback $\varphi^{*} \left( A \right)$ with the structure of a graded
	$R$-algebra. Under the bar construction, this corresponds to a coderivation on the
	tensor coalgebra $\tens{\varphi^{*} \left( A \right)[1]}[R] \cong
		\tens{\varphi^{*} \left( A[1] \right)}[R]$. This coderivation is precisely
	the pullback coderivation $\varphi^{*} \left( \nu \right)$ on the pullback tensor coalgebra
	$\tens{\varphi^{*} \left( A[1] \right)}[R]$. More generally, we use our notion of pullback
	in \cref{dfn:scalar-rest-ext-a-infinity} to define the pullback of a Banach $\Ainf$-algebra.
\end{rem}

More generally, we can define the pullback of a morphism
$f \colon \big( \tensf{W}[S], \nu \big) \rightarrow
	\big( \tensf{X}[Q], \xi \big)$
of pre-differential graded Banach coalgebras over a morphism $\alpha \colon (S,d_S) \rightarrow (Q,d_Q)$ of pre-differential graded Banach $\mathbbm{k}$-algebras
as follows. Given a diagram of morphisms of pre-differential graded-commutative Banach $\mathbbm{k}$-algebras
of the form
\begin{equation}
	\label{eq:two-morphism-algebras}
	\begin{tikzcd}
		{\left( S, d_S \right)} && {\left( Q, d_Q \right)} \\
		{\left( R, d_R \right)} && {\left( P, d_P \right)}
		\arrow["\alpha", from=1-1, to=1-3]
		\arrow["\varphi", from=2-1, to=1-1]
		\arrow["\psi"', from=2-3, to=1-3]
		\arrow["{\beta }"', from=2-1, to=2-3]
	\end{tikzcd}
\end{equation}
we will think of the pair $(\varphi, \psi)$ as a morphism between $\alpha$ and $\beta$.

\begin{dfn} \label{dfn:pullback-f-along-diagram}
	Given a morphism $f \colon \big( \tensf{W}[S], \nu \big) \rightarrow
		\big( \tensf{X}[Q], \xi \big)$ over $\alpha$,
	the \textbf{pullback} of $f$ along diagram \eqref{eq:two-morphism-algebras} is the unique morphism
	\begin{equation*}
		f^{\varphi}_{\psi} \colon
		\left( \tensf{\varphi^{*} \left( W \right)}[R], \varphi^{*} \left( \nu \right) \right)
		\rightarrow
		\left( \tensf{\psi^{*} \left( X \right)}[P], \psi^{*} \left( \xi \right) \right)
	\end{equation*}
	over $\beta$ which makes the following diagram commute:
	\begin{figure}[H]
		\begin{tikzcd}
			{\left( \tensf{W}[S], \nu \right)} &&&& {\left( \tensf{X}[Q], \xi \right)} \\
			& {\left( S, d_S \right)} && {\left( Q, d_Q \right)} \\
			& {\left( R, d_R \right)} && {\left( P, d_P \right)} & {} \\
			{\left( \tensf{\varphi^{*} \left( W \right)}[R], \varphi^{*} \left( \nu \right) \right)} &&&&
			{\left( \tensf{\psi^{*} \left( X \right)}[P], \psi^{*} \left( \xi \right) \right)}
			\arrow["f", from=1-1, to=1-5]
			\arrow[from=1-1, to=2-2] 
			\arrow[from=1-5, to=2-4] 
			\arrow["\alpha", from=2-2, to=2-4]
			\arrow["\varphi", from=3-2, to=2-2]
			\arrow["\psi"', from=3-4, to=2-4]
			\arrow["\beta"', from=3-2, to=3-4]
			\arrow[from=4-1, to=3-2] 
			\arrow["{\resover{\varphi}^{W}}", from=4-1, to=1-1]
			\arrow[from=4-5, to=3-4] 
			\arrow["{\resover{\psi}^{X}}"', from=4-5, to=1-5]
			\arrow["{\exists ! \, f^{\varphi}_{\psi}}"', dashed, from=4-1, to=4-5]
		\end{tikzcd}
		\caption{Pullback of a morphism between formal tensor coalgebras with coderivations
			over different pre-differential graded ground algebras.}
		\small\textsuperscript{The unnamed maps are the counit maps.}
		\label{fig:scalar-rest-morphism-formal-tensor-coalg-with-coderivations-different-ground-algebras}
	\end{figure}
\end{dfn}

The existence of $f^{\varphi}_{\psi}$ follows from the following diagram:
\begin{equation*} \adjustbox{scale=0.90,center}{
		\begin{tikzcd}
			{\left( \tensf{W}[S], \nu \right)}
			&&
			{\left( \tensf{\alpha^{*} \left( X \right)}[S], \alpha^{*} \left( \xi \right) \right)} &&
			{\left( \tensf{X}[Q], \xi \right)}
			\\
			& {\left( S,d_S \right)} && {\left( Q,d_Q \right)}
			\\
			& {\left( R, d_R \right)} && {\left( P, d_P \right)}
			\\
			{\left( \tensf{\varphi^{*} \left( W \right)}[R], \varphi^{*} \left( \nu \right) \right)} &&
			{\substack{
					\left( \tensf{\varphi^{*} \left( \alpha^{*} \left( X \right) \right)}[R],
					\varphi^{*} \left( \alpha^{*} \left( \xi \right) \right) \right) \\
					\shortparallel \\
					\left( \tensf{\beta^{*} \left( \psi^{*} \left( X \right) \right)}[R],
					\beta^{*} \left( \psi^{*} \left( \xi \right) \right) \right)
				}
			}
			&&
			{\left( \tensf{\psi^{*} \left( X \right)}[P], \psi^{*} \left( \xi \right) \right)}
			\arrow["{\varphi^{*} \left( \rescoho{f} \right)}"', dashed, from=4-1, to=4-3]
			\arrow["\varepsilon_R", from=4-1, to=3-2]
			\arrow["\varepsilon_R"', from=4-3, to=3-2]
			\arrow["\varepsilon_S"', from=1-1, to=2-2]
			\arrow["{\resover{\varphi}^W}", from=4-1, to=1-1]
			\arrow["{\rescoho{f}}", dashed, from=1-1, to=1-3]
			\arrow["\varepsilon_S", from=1-3, to=2-2]
			\arrow["\varphi", from=3-2, to=2-2]
			\arrow["\quad\alpha", from=2-2, to=2-4]
			\arrow["\psi"', from=3-4, to=2-4]
			\arrow["{\quad \beta}"', from=3-2, to=3-4]
			\arrow["{\resover{\alpha}^X}", from=1-3, to=1-5]
			\arrow["{\resover{\psi}^X}"', from=4-5, to=1-5]
			\arrow["\varepsilon_Q", from=1-5, to=2-4]
			\arrow["\varepsilon_P"', from=4-5, to=3-4]
			\arrow["{\resover{\beta}^{\psi^{*} \left( X \right)}}"', from=4-3, to=4-5]
			\arrow["f", bend left=10, from=1-1, to=1-5] 
			\arrow["f^{\varphi}_{\psi}"',bend right=10, dashed, from=4-1, to=4-5]
			\arrow["{\resover{\varphi}^{\alpha^{*} \left( X \right)}}"', from=4-3, to=1-3]
		\end{tikzcd}}
\end{equation*}
When $(S,d_S) = (Q,d_Q), \alpha = \id_S$ and $(R,d_R) = (P,d_P), \beta = \id_R$, then we must have
$\varphi = \psi$ and $f^{\varphi}_{\psi}$ recovers our previous notion of pullback $\varphi^{*} \left( f \right)$.

The pullback construction allows us to reduce identities for maps over different base algebras to the
corresponding identities over the same algebra. As an example of an application, let us show the following
generalization of \cref{cor:f-circ-mu-nu-circ-f-corest}:

\begin{lm} \label{lm:f-circ-mu-nu-circ-f-corest-different-ground-algebras}
	Let $\mathcal{R} = \left( R, d_R \right)$ and $\mathcal{S} = \left( S, d_S \right)$ be
	two pre-differential graded-commutative Banach $\mathbbm{k}$-algebras and let
	$\varphi \colon \mathcal{R} \rightarrow \mathcal{S}$ be a morphism of pre-differential graded Banach
	$\mathbbm{k}$-algebras.
	Let $V$ be a graded Banach $R$-module and let $W$ be a graded Banach $S$-module, let
	$\mu \colon \tensf{V}[R] \rightharpoonup \tensf{V}[R]$ be a generalized coderivation over $d_R$
	and let $\nu \colon \tensf{W}[S] \rightharpoonup \tensf{W}[S]$ be a generalized coderivation over $d_S$.
	Finally, let $f \colon \tensf{V}[R] \rightarrow \tensf{W}[S]$ be a morphism of graded Banach coalgebras
	over $\varphi$.
	Then $f \circ \mu = \nu \circ f$ if and only if $\corest{f} \circ \mu = \corest{\nu} \circ f$.
\end{lm}
\begin{proof}
	Using the universal property of the pullback, decompose $f$ as
	$f = \resover{\varphi}^W \circ \rescoho{f}$.
	Then
	\begin{equation*}
		\corest{f} \circ \mu = \corest{\left( \resover{\varphi}^W \circ \rescoho{f} \right)} \circ \mu =
		\corest{ \resover{\varphi}^W } \circ \rescoho{f} \circ \mu =
		\resover{\varphi}^W_{1} \circ \corest{\rescoho{f}} \circ \mu
	\end{equation*}
	and
	\begin{equation*}
		\corest{\nu} \circ f = \corest{\nu} \circ \resover{\varphi}^W \circ \rescoho{f} =
		\corest{ \resover{\varphi}^W } \circ \varphi^{*} \left( \nu \right) \circ \rescoho{f} =
		\resover{\varphi}^W_{1} \circ \corest{\varphi^{*} \left( \nu \right)} \circ \rescoho{f}
	\end{equation*}
	(see \cref{fig:proof-f-mu-nu-f-corest-general}). Hence, if $\corest{f} \circ \mu = \corest{\nu} \circ f$
	then $\corest{\rescoho{f}} \circ \mu = \corest{\varphi^{*} \left( \nu \right)} \circ \rescoho{f}$.
	The morphism $\rescoho{f}$ is over the same base algebra $R$ and both $\mu$ and $\varphi^{*} \left( \nu \right)$
	are coderivations over the same derivation $d_R$. Hence, by \cref{cor:f-circ-mu-nu-circ-f-corest}, we
	have $\rescoho{f} \circ \mu = \varphi^{*} \left( \nu \right) \circ \rescoho{f}$ and so
	\begin{equation*}
		f \circ \mu = \resover{\varphi}^W \circ \rescoho{f} \circ \mu =
		\resover{\varphi}^W \circ \varphi^{*} \left( \nu \right) \circ \rescoho{f} =
		\nu \circ \resover{\varphi}^W \circ \rescoho{f} = \nu \circ f.
	\end{equation*}
\end{proof}

\begin{figure}
	\begin{tikzcd}
		{\tensf{V}[R]} && {\tensf{\varphi^{*} \left( W \right)}[R]} && {\tensf{W}[S]} \\
		&& {\tensf{\varphi^{*} \left( W \right)}[R]} && {\tensf{W}[S]} \\
		&& {\varphi^{*} \left( W \right)} && W
		\arrow["{\rescoho{f}}", from=1-1, to=1-3]
		\arrow["f", curve={height=-24pt}, from=1-1, to=1-5]
		\arrow["{\resover{\varphi}^W}", from=1-3, to=1-5]
		\arrow["{\varphi^{*} \left( \nu \right)}", from=1-3, to=2-3]
		\arrow["{\corest{\varphi^{*} \left( \nu \right)}}"', curve={height=40pt}, from=1-3, to=3-3]
		\arrow["\nu", from=1-5, to=2-5]
		\arrow["{\corest{\nu}}", curve={height=-40pt}, from=1-5, to=3-5]
		\arrow["{\resover{\varphi}^W}", from=2-3, to=2-5]
		\arrow["{\pi_1^{\varphi^{*} \left( W \right)}}", from=2-3, to=3-3]
		\arrow["{\pi_1^W}", from=2-5, to=3-5]
		\arrow["{\resover{\varphi}^W_1}"', from=3-3, to=3-5]
	\end{tikzcd}
	\caption{Proof of \cref{lm:f-circ-mu-nu-circ-f-corest-different-ground-algebras}.}
	\label{fig:proof-f-mu-nu-f-corest-general}
\end{figure}

\subsubsection{Scalar Extension for Formal Tensor Coalgebras}
\label{sec:scalar-extension-formal-tensor-coalgebras}
Since the scalar extension functor $\varphi_{!} \colon \GBMod[R] \rightarrow \GBMod[S]$
is strong monoidal, the scalar extension of
a graded Banach $R$-coalgebra has a natural structure of a graded Banach $S$-coalgebra.
Identifying the scalar extension of a formal tensor coalgebra with the formal tensor coalgebra
on the scalar extension module, we obtain a notion of scalar extension for
formal tensor coalgebras which we now describe.

Let $V$ be a graded Banach $R$-module and let $C = \tensf{V}[R]$ be the formal tensor coalgebra.
The scalar extension
$\varphi_{!} \left( C \right) = S \cotimes_R \tensf{V}[R]$ has a natural structure
of a graded Banach $S$-coalgebra (see \cref{subsec:coalg-in-monoidal-cat}) and the tensor constraints
of $\varphi_{!}$ induce
a natural isomorphism
$\Theta \colon \tensf{\varphi_{!} \left( V \right)}[S] \rightarrow
	\varphi_{!} \left( \tensf{V}[R] \right)$ of graded Banach $S$-coalgebras given explicitly by
\begin{equation*}
	\Theta \left(
	\left( s_1 \cotimes_R v_1 \right) \cotimes_S \dots \cotimes_S
	\left( s_k \cotimes_R v_k \right)
	\right)
	=
	(-1)^{\varepsilon}
	\left( s_1 \cdots s_k \right) \cotimes_R
	\left( v_1 \cotimes_R \dots \cotimes_R v_k \right)
\end{equation*}
with
\begin{equation} \label{eq:sign-theta-scalar-ext}
	\varepsilon =
	\sum_{i=2}^k \sum_{j=1}^{i-1} \braidd{s_i}{v_j}.
\end{equation}
Hence, we can (and will) identify the scalar extension $\varphi_{!} \left( \tensf{V}[R] \right)$
of the formal tensor coalgebra
$\tensf{V}[R]$ with the formal tensor coalgebra of
$\varphi_{!} \left( V \right)$ over $S$
and call $\tensf{\varphi_{!} \left( V \right)}[S]$ the \textbf{scalar extension}
of $\tensf{V}[R]$ along $\varphi$.

The scalar extension $\tensf{\varphi_{!} \left( V \right)}[S]$ comes equipped with a
universal Banach coalgebra morphism
$\resunder{\varphi} = \resunder{\varphi}^{V} \colon \tensf{V}[R] \rightarrow
	\tensf{\varphi_{!} \left( V \right)}[S]$ over $\varphi$ given by
\begin{equation} \label{eq:canonical-extension-map}
	\begin{aligned}
		\resunder{\varphi} \left( r \right)                                   & \defeq \varphi(r),
		\\
		\resunder{\varphi} \left( v_1 \cotimes_R \dots \cotimes_R v_n \right) & \defeq
		\left( 1_S \cotimes_R v_1 \right) \cotimes_S \dots \cotimes_S
		\left( 1_S \cotimes_R v_n \right).
	\end{aligned}
\end{equation}
Note that the only non-zero component of the corestriction $\corest{\resunder{\varphi}}$ is given by
$\resunder{\varphi}_1 \colon V \rightarrow S \cotimes_R V$ which is the canonical map
over $\varphi$ given by $v \mapsto 1_S \cotimes_R v$.

The pair $\left( \tensf{\varphi_{!} \left( V \right)}[S], \resunder{\varphi} \right)$
satisfies the following universal property:
Any Banach coalgebra morphism  $f \colon \tensf{V}[R] \rightarrow \tensf{W}[S]$ over
$\varphi$ factors uniquely as $f = \resunder{f} \circ \resunder{\varphi}$ where
$\resunder{f} \colon \tensf{\varphi_{!} \left( V \right)}[S] \rightarrow \tensf{W}[S]$ is a morphism
of graded Banach $S$-coalgebras (see \cref{fig:scalar-ext-formal-tensor-coalg}).
The corestriction $\corest{\resunder{f}}$ of the morphism $\resunder{f}$,
and hence $\resunder{f}$ itself,
is determined uniquely by the equation $\corest{\resunder{f}} \circ \resunder{\varphi} = \corest{f}$.

\begin{figure}
	\begin{tikzcd}
		{\tensf{V}[R]} & {\tensf{\varphi_{!} \left( V \right)}[S]} & {\tensf{W}[S]} \\
		R & S & S
		\arrow[from=1-1, to=2-1]
		\arrow[from=1-2, to=2-2]
		\arrow[from=1-3, to=2-3]
		\arrow["{\exists ! \resunder{f}}", dashed, from=1-2, to=1-3]
		\arrow["{\resunder{\varphi}}", from=1-1, to=1-2]
		\arrow["{\varphi}", from=2-1, to=2-2]
		\arrow["{\id}", from=2-2, to=2-3]
		\arrow["f", bend left, from=1-1, to=1-3]
	\end{tikzcd}
	\caption{Universal property of scalar extension for formal tensor coalgebras.}
	\small\textsuperscript{The unnamed vertical maps are the counit maps.}
	\label{fig:scalar-ext-formal-tensor-coalg}
\end{figure}

The scalar extension construction for formal tensor coalgebras extends also to morphisms as follows.
Given a morphism $f \colon \tensf{V_1}[R] \rightarrow \tensf{V_2}[R]$ of graded Banach $R$-coalgebras,
there exists a unique morphism
$\varphi_{!} \left( f \right) \colon \tensf{\varphi_{!} \left( V_1 \right)}[S]
	\rightarrow \tensf{\varphi_{!} \left( V_2 \right)}[S]$ of graded Banach $S$-coalgebras
such that
$\resunder{\varphi}^{V_2} \circ f = \varphi_{!} \left( f \right) \circ \resunder{\varphi}^{V_1}$
(see \cref{fig:scalar-ext-morphism-formal-tensor-coalg}). The morphism $\varphi_{!} \left( f \right)$
corresponds under the isomorphism $\Theta$ to $\id_S \cotimes_R f$ and its corestriction
$\corest{\varphi_{!} \left( f \right)} \colon \tensf{\varphi_{!} \left( V_1 \right)}[S] \rightarrow
	\varphi_{!} \left( V_2 \right)$ is given by
\begin{align*}
	\varphi_{!} \left( f \right)_{0} \left( s \right) & \defeq
	s \cotimes_R f_0 \left( 1_R \right),
	\\
	\varphi_{!} \left( f \right)_{k} \left(
	\left( s_1 \cotimes_R v_1 \right) \cotimes_S \dots \cotimes_S \left( s_k \cotimes_R v_k \right)
	\right)                                           & \defeq
	                                                    (-1)^{\varepsilon} \left( s_1 \cdots s_k \right) \cotimes_R
	f_k \left( v_1 \cotimes_R \dots \cotimes_R v_k \right),
\end{align*}
where the sign $\varepsilon$ is given by \cref{eq:sign-theta-scalar-ext}. Note that the corestriction
$\corest{\varphi_{!} \left( f \right)}$ of $\varphi_{!} \left( f \right)$, and hence
$\varphi_{!} \left( f \right)$ itself, is determined uniquely
by the equation
$\resunder{\varphi}_1 \circ \corest{f} =
	\corest{\varphi_{!} \left( f \right)} \circ \resunder{\varphi}$.

\begin{figure}
	\begin{tikzcd}
		{\tensf{\varphi_{!} \left( V_1 \right)}[S]} &&
		{\tensf{\varphi_{!} \left( V_2 \right)}[S]} \\
		& {S} \\
		& {R} \\
		{\tensf{V_1}[R]} &&
		{\tensf{V_2}[R]}
		\arrow["f", from=4-1, to=4-3]
		\arrow["\varphi", from=3-2, to=2-2]
		\arrow[from=1-1, to=2-2] 
		\arrow[from=1-3, to=2-2] 
		\arrow["{\resunder{\varphi}^{V_1}}", from=4-1, to=1-1]
		\arrow["{\resunder{\varphi}^{V_2}}"', from=4-3, to=1-3]
		\arrow[from=4-1, to=3-2] 
		\arrow[from=4-3, to=3-2] 
		\arrow["{\exists ! \varphi_{!} \left( f \right)}"', dashed, from=1-1, to=1-3]
	\end{tikzcd}
	\caption{Scalar extension of a morphism between formal tensor $R$-coalgebras.}
	\small\textsuperscript{The unnamed maps are the counit maps.}
	\label{fig:scalar-ext-morphism-formal-tensor-coalg}
\end{figure}

\begin{dfn} \label{dfn:ext-morphism-tensor-coalgebras}
	The morphism $\varphi_{!} \left( f \right) \colon \tensf{\varphi_{!} \left( V_1 \right)}[S]
		\rightarrow \tensf{\varphi_{!} \left( V_2 \right)}[S]$ is called the \textbf{scalar extension}
	of $f \colon \tensf{V_1}[R] \rightarrow \tensf{V_2}[R]$ along $\varphi$.
\end{dfn}

The scalar extension construction $\varphi_{!}$ defines a covariant functor from the category of
formal tensor $R$-coalgebras to the category of formal tensor $S$-coalgebras.
In addition, the scalar extension is functorial \textit{up to a natural isomorphism} with respect to morphisms of the underlying graded Banach $\mathbbm{k}$-algebras.
More precisely, given a morphism $\psi \colon P \rightarrow R$ of graded Banach $\mathbbm{k}$-algebras,
and a graded Banach $P$-module $U$, we have a natural isomorphism
\begin{equation} \label{eq:extension-scalars-iso}
	\varphi_{!} \left( \psi_{!} \left( U \right) \right) =
	S \cotimes_R \left( R \cotimes_P U \right) \cong S \cotimes_P U =
	\left( \varphi \circ \psi \right)_{!} \left( U \right)
\end{equation}
which induces an isomorphism $\tensf{\varphi_{!} \left( \psi_{!} \left( U \right) \right)}[S] \cong
	\tensf{\left( \varphi \circ \psi \right)_{!} \left( U \right)}[S]$
of graded Banach $S$-coalgebras,
so that $\varphi_{!} \circ \psi_{!} \cong \left( \varphi \circ \psi \right)_{!}$ as functors from
the category of formal tensor $P$-coalgebras to the category of formal tensor $S$-coalgebras.

\begin{figure}
	\begin{tikzcd}
		R &&& S \\
		& {\tensf{V}[R]} & {\tensf{\varphi_{!} \left( V \right)}[S]} \\
		& {\tensf{V}[R]} & {\tensf{\varphi_{!} \left( V \right)}[S]} \\
		R &&& S
		\arrow["\mu", harpoon, from=2-2, to=3-2]
		\arrow["{\resunder{\varphi}}", from=2-2, to=2-3]
		\arrow[from=2-3, to=1-4]
		\arrow["{d_S}", harpoon, from=1-4, to=4-4]
		\arrow["\varphi", from=1-1, to=1-4]
		\arrow[from=2-2, to=1-1]
		\arrow["{\exists ! \, \varphi_{!} \left( \mu \right)}"', dotted, harpoon, from=2-3, to=3-3]
		\arrow["{\resunder{\varphi}}"', from=3-2, to=3-3]
		\arrow["{d_R}"', harpoon, from=1-1, to=4-1]
		\arrow[from=3-2, to=4-1]
		\arrow["\varphi"', from=4-1, to=4-4]
		\arrow[from=3-3, to=4-4]
	\end{tikzcd}
	\caption{The scalar extension of a generalized coderivation.}
	\small\textsuperscript{The unnamed maps are the counit maps.}
	\label{fig:scalar-ext-of-coderivation}
\end{figure}
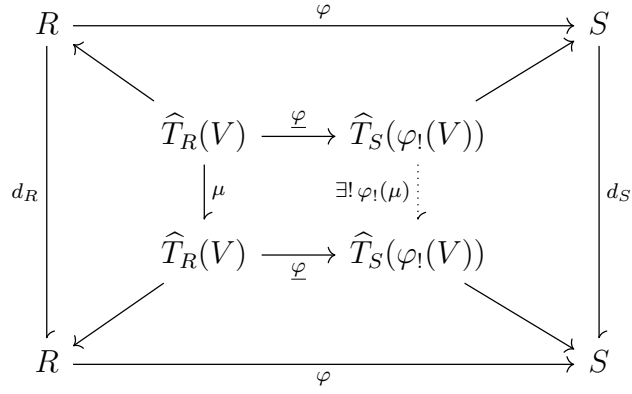

Next, we can extend the scalar extension construction to formal tensor coalgebras equipped with a generalized coderivation. Let $\mathcal{R} = \left( R, d_R \right)$ and
$\mathcal{S} = \left( S, d_S \right)$ be
two pre-differential graded-commutative Banach $\mathbbm{k}$-algebras
and let $\varphi \colon \mathcal{R} \rightarrow \mathcal{S}$ be a morphism of pre-differential graded Banach $\mathbbm{k}$-algebras. Let $\big( \tensf{V}[R], \mu \big)$ be a formal tensor
$\mathcal{R}$-coalgebra. Then there exists a unique coderivation
$\varphi_{!} \left( \mu \right) \colon \tensf{\varphi_{!} \left( V \right)}[S] \rightharpoonup
	\tensf{\varphi_{!} \left( V \right)}[S]$
over $d_S$ with $\degb{\varphi_{!} \left( \mu \right)} = \degb{\mu}$ such that
$\varphi_{!} \left( \mu \right) \circ \resunder{\varphi} = \resunder{\varphi} \circ \mu$
(see \cref{fig:scalar-ext-of-coderivation}). In other words, the coderivation
$\varphi_{!} \left( \mu \right)$ is the unique coderivation with respect to which the morphism
$\resunder{\varphi}$ becomes a morphism
\begin{equation*}
	\resunder{\varphi} \colon
	\left( \tensf{V}[R], \mu \right) \rightarrow
	\left( \tensf{\varphi_{!} \left( V \right)}[S], \varphi_{!} \left( \mu \right) \right)
\end{equation*}
of pre-differential graded Banach coalgebras over $\varphi$.

The coderivation $\varphi_{!} \left( \mu \right)$
corresponds under the isomorphism $\Theta$ to $d_S \cotimes_R \id_{\tensf{V}[R]} + \id_S \cotimes_R \mu$
and its corestriction
$\corest{\varphi_{!} \left( \mu \right)} \colon \tensf{\varphi_{!} \left( V \right)}[S] \rightharpoonup
	\varphi_{!} \left( V \right)$ is given by
\begin{align*}
	\varphi_{!} \left( \mu \right)_0 \left( s \right)
	 & =
	(-1)^{\braidd{\mu}{s}} s \cotimes_R \mu_0 \left( 1 \right),
	\\
	\varphi_{!} \left( \mu \right)_1 \left( s \cotimes_R v \right)
	 & =
	d_S \left( s \right) \cotimes_R v + (-1)^{\braidd{\mu}{s}} s \cotimes_R \mu_1 \left( v \right),
\end{align*}
and
\begin{equation*}
	\begin{aligned}
		\varphi_{!} \left( \mu \right)_{k} \left(
		\left( s_1 \cotimes_R v_1 \right) \cotimes_S \dots \cotimes_S \left( s_k \cotimes_R v_k \right)
		\right)
		={} &
		(-1)^{\varepsilon + \sum_{i=1}^k \braidd{\mu}{s_i}}
		\\
		    & \qquad
		\left( s_1 \cdots s_k \right) \cotimes_R \mu_k \left( v_1 \cotimes_R \dots \cotimes_R v_k \right)
	\end{aligned}
\end{equation*}
for $k \geq 2$. Note that the formula for $\varphi_{!} \left( \mu \right)_1 \colon \varphi_{!} \left( V \right)
	\rightharpoonup \varphi_{!} \left( V \right)$ makes sense by
\cref{rem:tensor-product-d-operators-abuse} and that $\varphi_{!} \left( \mu \right)_1$
is indeed a module derivation over $d_S$ as we have
\begin{equation*}
	\begin{aligned}
		\varphi_{!} \left( \mu \right)_1 \left(
		s' \cdot \left( s \cotimes_R v \right)
		\right)
		={} &
		\varphi_{!} \left( \mu \right)_1 \left(
		\left( s' \cdot s \right) \cotimes_R v
		\right)
		\\
		={} &
		d_S \left( s' \cdot s \right) \cotimes_R v +
			                                         (-1)^{\braid{\degb{\mu}}{\degb{s'} + \degb{s}}} \left( s' \cdot s \right) \cotimes_R
		\mu_1 \left( v \right)
		\\
		={} &
		d_S \left( s' \right) \cdot \left( s \cotimes_R v \right) +
		                            (-1)^{\braidd{s'}{d_S}} \left( s' \cdot d_S \left( s \right) \right) \cotimes_R v
		\\
		    &
		+
		(-1)^{\braid{\degb{\mu}}{\degb{s'} + \degb{s}}} \left( s' \cdot s \right) \cotimes_R
		\mu_1 \left( v \right)
		\\
		={} &
		d_S \left( s' \right) \cdot \left( s \cotimes_R v \right)
		\\
		    &
		+
		(-1)^{\braidd{\mu}{s'}} s' \cdot \left(
		d_S \left( s \right) \cotimes_R v +
			                                (-1)^{\braidd{\mu}{s}} s \cotimes_R \mu_1 \left( v \right)
		\right)
		\\
		={} &
		d_S \left( s' \right) \cdot \left( s \cotimes_R v \right)
		                            +
		                            (-1)^{\braidd{\mu}{s'}} s' \cdot \varphi_{!} \left( \mu \right)_1 \left( s \cotimes_R v \right)
		.
	\end{aligned}
\end{equation*}

\begin{dfn} \label{dfn:scalar-ext-generalized-coderivation}
	The coderivation $\varphi_{!} \left( \mu \right)$ is called the \textbf{scalar extension}
	of $\mu$ along $\varphi$. The formal tensor $\mathcal{S}$-coalgebra
	$\big( \tensf{\varphi_{!} \left( V \right)}[S], \varphi_{!} \left( \mu \right) \big)$, together
	with the map $\resunder{\varphi}$, is called the \textbf{scalar extension} of
	$\big( \tensf{V}[R], \mu \big)$ along $\varphi$.
\end{dfn}

\begin{rem} \label{rem:scalar-extension-depends-differentials}
	Note that the scalar extension $\varphi_{!} \left( \mu \right)$ depends on $\varphi$ as
	a morphism of \textit{pre-differential} graded Banach $\mathbbm{k}$-algebras. Namely,
	given two morphisms $\varphi_i \colon \left( R, d_R \right) \rightarrow \left( S, d_i \right)$
	for $i = 1,2$ which coincide as morphisms of graded Banach $\mathbbm{k}$-algebras, i.e.,
	$\varphi_1 \left( r \right) = \varphi_2 \left( r \right)$ for all $r \in R$, then
	$\left( \varphi_1 \right)_{!} \left( V \right) = \left( \varphi_2 \right)_{!} \left( V \right)$,
	but we still might have
	$\left( \varphi_1 \right)_{!} \left( \mu \right) \neq \left( \varphi_2 \right)_{!} \left( \mu \right)$.
\end{rem}

\begin{rem}
	Although the derivation $d_S$ appears only in the $k = 1$ component of the corestriction
	$\corest{\varphi_{!} \left( \mu \right)}$, it actually plays a role when considering
	the full action of $\varphi_{!} \left( \mu \right)$ on elements of arbitrary weight
	(see \cref{eq:generalized-coder-coextension}). For example, we have
	\begin{equation*}
		\begin{aligned}
			\varphi_{!} \left( \mu \right) \left( s \right) =
			d_S \left( s \right) + (-1)^{\braidd{\mu}{s}} s \cotimes_R \mu_0 \left( 1 \right),
		\end{aligned}
	\end{equation*}
	and
	\begin{equation*}
		\begin{aligned}
			\MoveEqLeft
			\varphi_{!} \left( \mu \right) \left( \left( s_1 \cotimes_R v_1 \right) \cotimes_S \left(  s_2 \cotimes_R v_2 \right) \right)
			\\
			={} &
			\left( 1_S \cotimes_R \mu_0 \left( 1 \right) \right) \cotimes_S
			\left( s_1 \cotimes_R v_1 \right) \cotimes_S
			\left( s_2 \cotimes_R v_2 \right)
			\\
			    & +
			      (-1)^{\braid{\degb{\mu}}{\degb{s_1} + \degb{v_1}}}
			\left( s_1 \cotimes_R v_1 \right) \cotimes_S
			\left( 1_S \cotimes_R \mu_0 \left( 1 \right) \right) \cotimes_S
			\left( s_2 \cotimes_R v_2 \right)
			\\
			    & +
			      (-1)^{\braid{\degb{\mu}}{\degb{s_1} + \degb{v_1} + \degb{s_2} + \degb{v_2}}}
			\left( s_1 \cotimes_R v_1 \right) \cotimes_S
			\left( s_2 \cotimes_R v_2 \right) \cotimes_S
			\left( 1_S \cotimes_R \mu_0 \left( 1 \right) \right)
			\\
			    & +
			\left( d_S \left( s_1 \right) \cotimes_R v_1 \right) \cotimes_S
			\left( s_2 \cotimes_R v_2 \right)
			\\
			    & +
			      (-1)^{\braidd{\mu}{s_1}}
			\left( s_1 \cotimes_R \mu_1 \left( v_1 \right) \right) \cotimes_S
			\left( s_2 \cotimes_R v_2 \right)
			\\
			    & +
			      (-1)^{\braid{\degb{\mu}}{\degb{s_1} + \degb{v_1}}}
			\left( s_1 \cotimes_R v_1 \right) \cotimes_S
			\left( d_S \left( s_2 \right) \cotimes_R v_2 \right)
			\\
			    & +
			      (-1)^{\braid{\degb{\mu}}{\degb{s_1} + \degb{v_1} + \degb{s_2}}}
			\left( s_1 \cotimes_R v_1 \right) \cotimes_S
			\left( s_2 \cotimes_R \mu_1 \left( v_2 \right) \right)
			\\
			    & +
			      (-1)^{\braidd{s_2}{v_1} + \braid{\degb{\mu}}{\degb{s_1} + \degb{s_2}}}
			\left( s_1 \cdot s_2 \right) \cotimes_R
			\mu_2 \left( v_1, v_2 \right).
		\end{aligned}
	\end{equation*}
\end{rem}

The universal property of the scalar extension extends naturally to formal tensor coalgebras equipped
with generalized coderivations.
Namely, if $\big( \tensf{V}[R], \mu \big)$ is a formal tensor $\mathcal{R}$-coalgebra and
$\big( \tensf{W}[S], \nu \big)$ is a formal tensor $\mathcal{S}$-coalgebra then
any morphism
$f \colon \big( \tensf{V}[R], \mu \big) \rightarrow \big( \tensf{W}[S], \nu \big)$ of
pre-differential graded Banach coalgebras over $\varphi \colon \mathcal{R} \rightarrow \mathcal{S}$ factors
uniquely as $f = \resunder{f} \circ \resunder{\varphi}$ where
$\resunder{f} \colon \big( \tensf{\varphi_{!} \left( V \right)}[S], \varphi_{!} \left( \mu \right)
	\big) \rightarrow \big( \tensf{W}[S], \nu \big)$
is a morphism of pre-differential graded Banach $\mathcal{S}$-coalgebras
(see \cref{fig:scalar-ext-formal-tensor-coalg-over-pre-cdga}).
The map $\resunder{f}$ is the
same map as in \cref{fig:scalar-ext-formal-tensor-coalg}.
To verify that $\resunder{f}$ commutes with the coderivations, it is enough to check that
$\corest{\nu} \circ \resunder{f} = \corest{\resunder{f}} \circ \varphi_{!} \left( \mu \right)$
(see \cref{cor:f-circ-mu-nu-circ-f-corest}). This follows from the chain of equalities
\begin{equation*}
	\corest{\resunder{f}} \circ \varphi_{!} \left( \mu \right) \circ \resunder{\varphi} =
	\corest{\resunder{f}} \circ \resunder{\varphi} \circ \mu = \corest{f} \circ \mu =
	\corest{\nu} \circ f = \corest{\nu} \circ \resunder{f} \circ \resunder{\varphi}.
\end{equation*}

\begin{figure}
	\begin{tikzcd}
		{\left( \tensf{V}[R], \mu \right)} &
		{\left( \tensf{\varphi_{!} \left( V \right)}[S], \varphi_{!} \left( \mu \right) \right)} &
		{\left( \tensf{W}[S], \nu \right)} \\
		(R, d_R) & (S, d_S) & (S, d_S)
		\arrow[from=1-1, to=2-1] 
		\arrow[from=1-2, to=2-2] 
		\arrow[from=1-3, to=2-3] 
		\arrow["{\exists ! \resunder{f}}", dashed, from=1-2, to=1-3]
		\arrow["{\resunder{\varphi}}", from=1-1, to=1-2]
		\arrow["{\varphi}", from=2-1, to=2-2]
		\arrow["{\id}", from=2-2, to=2-3]
		\arrow["f", bend left, from=1-1, to=1-3]
	\end{tikzcd}
	\caption{Universal property of the scalar extension of formal tensor coalgebras.}
	\small\textsuperscript{The unnamed vertical maps are the counit maps.}
	\label{fig:scalar-ext-formal-tensor-coalg-over-pre-cdga}
\end{figure}

Given a morphism $\varphi \colon \mathcal{R} \rightarrow \mathcal{S}$ of pre-differential
graded Banach $\mathbbm{k}$-algebras, the scalar extension construction
$\varphi_{!}$ defines a covariant functor from the category of formal tensor $\mathcal{R}$-coalgebras
to the category of formal tensor $\mathcal{S}$-coalgebras. That is,
given a morphism
\begin{equation*}
	f \colon \left( \tensf{V_1}[R], \mu^1 \right) \rightarrow  \left( \tensf{V_2}[R], \mu^2 \right)
\end{equation*}
of formal tensor $\mathcal{R}$-coalgebras, i.e.,
$f$ satisfies $f \circ \mu^1 = \mu^2 \circ f$, the scalar extension
$\varphi_{!} \left( f \right)$ of $f$ given by \cref{dfn:ext-morphism-tensor-coalgebras} also
satisfies
$\varphi_{!} \left( \mu^2 \right) \circ \varphi_{!} \left( f \right) = \varphi_{!} \left( f \right) \circ \varphi_{!} \left( \mu^1 \right)$
(see \cref{fig:scalar-ext-morphism-formal-tensor-coalg-with-coderivations}).

\begin{figure}
	\begin{tikzcd}
		{\left( \tensf{\varphi_{!} \left( V_1 \right)}[S], \varphi_{!} \left( \mu^1 \right) \right)}
		&&
		{\left( \tensf{\varphi_{!} \left( V_2 \right)}[S], \varphi_{!} \left( \mu^2 \right) \right)}
		\\
		& {(S, d_S)} \\
		& {(R, d_R)} \\
		{\left( \tensf{V_1}[R], \mu^1 \right)} &&
		{\left( \tensf{V_2}[R], \mu^2 \right)}
		\arrow["\varphi", from=3-2, to=2-2]
		\arrow[from=1-1, to=2-2] 
		\arrow[from=1-3, to=2-2] 
		\arrow["{\resunder{\varphi}^{V_1}}", from=4-1, to=1-1]
		\arrow["{\resunder{\varphi}^{V_2}}"', from=4-3, to=1-3]
		\arrow[from=4-1, to=3-2] 
		\arrow[from=4-3, to=3-2] 
		\arrow["f", from=4-1, to=4-3]
		\arrow["{\exists ! \varphi_{!} \left( f \right)}"', dashed, from=1-1, to=1-3]
	\end{tikzcd}
	\caption{Scalar extension of a morphism between formal tensor $\mathcal{R}$-coalgebras.}
	\small\textsuperscript{The unnamed maps are the counit maps.}
	\label{fig:scalar-ext-morphism-formal-tensor-coalg-with-coderivations}
\end{figure}

In addition, the scalar extension is functorial \textit{up to a natural isomorphism} with respect to
morphisms of the underlying pre-differential graded Banach $\mathbbm{k}$-algebras.
More precisely, given a morphism $\psi \colon \mathcal{P} \rightarrow \mathcal{R}$
of pre-differential graded Banach $\mathbbm{k}$-algebras, we have
$\varphi_{!} \circ \psi_{!} \cong \left( \varphi \circ \psi \right)_{!}$ as functors from
the category of formal tensor $\mathcal{P}$-coalgebras to the category of formal tensor
$\mathcal{S}$-coalgebras, where the natural isomorphism is induced by \eqref{eq:extension-scalars-iso}.

\begin{rem} \label{rem:scalar-extension-terminology-tensor-coalgebras}
	Our notion of scalar extension is motivated by the bar construction, in which the structure of
	a graded algebra is encoded by a coderivation on a tensor coalgebra.
	Given a graded-commutative $\mathbbm{k}$-algebra $R$ and a graded $R$-module $A$, endowing
	$A$ with the structure of a graded (non-unital) $R$-algebra is equivalent under
	the bar construction to providing a coderivation $\mu$ on the tensor coalgebra $\tens{A[1]}[R]$
	which satisfies $\mu^2 = 0$ and whose only non-zero component is $\mu_2$.
	Now, given a morphism $\varphi \colon R \rightarrow S$,
	we know we can endow the scalar extension module $\varphi_{!} \left( A \right)$
	with the structure of a graded $S$-algebra.
	Under the bar construction, this corresponds to a coderivation on the
	tensor coalgebra $\tens{\varphi_{!} \left( A \right)[1]}[S] \cong
		\tens{\varphi_{!} \left( A[1] \right)}[S]$. This coderivation is precisely
	the scalar extension coderivation $\varphi_{!} \left( \mu \right)$ on the
	scalar extension tensor coalgebra
	$\tens{\varphi_{!} \left( A[1] \right)}[S]$.
	More generally, we use our notion of scalar extension in \cref{dfn:scalar-rest-ext-a-infinity}
	to define the scalar extension of a Banach $\Ainf$-algebra.
\end{rem}

\section{Banach \texorpdfstring{$\Ainf$-}{A-infinity }algebras}
\label{sec:a-inf-algebras}

In this section, we introduce the notion of a (non-Archimedean) Banach $\Ainf$-algebra
$\mathcal{A} = \left( A, \mu \right)$ over a differential graded-commutative
Banach algebra $\mathcal{R}$ and give basic definitions used for the remainder of the text.
Our notion allows $\mathcal{A}$ to be curved, as long
as the curvature $\mu_0 \left( 1 \right)$ satisfies $\nnorm[\mu_0 \left( 1 \right)] < 1$.
We define morphisms $f \colon \mathcal{A} \rightarrow \mathcal{B}$ between Banach $\Ainf$-algebras,
possibly over different ground algebras.
Morphisms $f$ are allowed to have a non-zero change of connection element, as long as
$\nnorm[f_0 \left( 1 \right)] < 1$.
In \cref{subsec:pseudoisotopy-a-inf}, we introduce the notion of pseudoisotopy between
two $\Ainf$-algebras.  In \cref{subsec:bounding-cochains-gauge-equivalence},
we define weak and strong bounding cochains as solutions to the Maurer--Cartan equation,
and show that bounding cochains can be pushed forward along $\Ainf$-morphisms.
Finally, we use our notion of pseudoisotopy to define gauge equivalence between bounding cochains.

The basic idea of working systematically with complete filtered $\Ainf$-algebras appears
in \cite{Fukaya2009}. We work with norms instead of valuations, and our definitions are an
adaptation of the corresponding definitions from
\cite{Solomon2016,Solomon2016a}, appearing in the context of symplectic geometry,
to a more abstract algebraic setting, appropriate to this work.

Let $\mathbbm{k}$ be a commutative ground ring which will be fixed for the duration of the section.
We endow $\mathbbm{k}$ with the trivial norm and consider $\mathbbm{k}$ as a commutative
Banach ring. In what follows, we work with the Koszul grading datum
(see \cref{subsec:grading-data}) in which the grading group is given by $\mathbbm{G} = \ZZ$,
the parity form is given by $\braid{a}{b} = a \cdot b \pmod{2}$, and $\go = 1$.
Thus, objects will be $\ZZ$-graded and the (pre-)differentials on objects
\textit{raise} degree of elements by one.

\subsection{Basic Definitions}

\begin{dfn} \label{def:a-inf-Banach-algebra}
	Let $\mathcal{R} = (R,d)$ be a differential graded-commutative Banach $\mathbbm{k}$-algebra.
	A (\textbf{shifted, weakly curved}) \textbf{Banach} $\Ainf$-\textbf{algebra over} $\mathcal{R}$ is a pair
	$\mathcal{A} = (A,\mu)$ where $A$ is a graded Banach $R$-module
	and $\mu \colon \tensf{A}[R] \rightharpoonup \tensf{A}[R]$ is a degree one coderivation over $d$
	which satisfies the conditions
	\begin{equation} \label{eq:banach-a-inf-mu-conditions}
		(1) \quad \nnorm[\mu] \leq 1, \qquad (2) \quad \nnorm[\mu_0(1)] < 1, \qquad (3) \quad \mu^2 = 0.
	\end{equation}
	The coderivation $\mu$ is called an $\Ainf$\textbf{-structure} on $A$ over $\mathcal{R}$.
\end{dfn}

Let us write more explicitly what \cref{def:a-inf-Banach-algebra} means. By
\cref{prop:classification-generalized-coderivations-formal-tensor-coalgebra}, a coderivation
$\mu \colon \tensf{A}[R] \rightharpoonup \tensf{A}[R]$ over $d$ is uniquely determined by its corestriction
$\corest{\mu} \colon \tensf{A}[R] \rightharpoonup A$. The corestriction $\corest{\mu}$ is
uniquely determined by the sequence of contractive degree one maps
$\mu_k \colon A^{\cotimes k} \rightharpoonup A$,
which in turn, are determined uniquely by the associated contractive, degree one, multilinear maps
$\mu_k \colon A^{\times k} \rightharpoonup A$. Hence, an $\Ainf$-structure on
$A$ over $\mathcal{R}$ is equivalently given by:
\begin{enumerate}
	\item An element $\mu_0(1) \in A^1$ called the \textbf{curvature} which satisfies
	      $\nnorm[\mu_0 \left( 1 \right)] < 1$.\footnote{By \textit{weakly curved}, we mean that we allow
		      non-zero curvature, as long as it satisfies $\nnorm[\mu_0 \left( 1 \right)] < 1$.}
	\item A $\mathbbm{k}$-linear map $\mu_1 \colon A \rightharpoonup A$ of degree one which satisfies
	      \begin{equation*}
		      \mu_1 \left( r \cdot a \right) = dr \cdot a + (-1)^{\degb{r}} r \cdot \mu_1 \left( a \right),
		      \qquad \nnorm[\mu_1 \left( a \right)] \leq \nnorm[a]
	      \end{equation*}
	      for all $a \in A$ and $r \in R$.
	\item A sequence of degree one maps $\mu_k \colon A^{\times k} \rightharpoonup A$ for $k \geq 2$ which are
	      $R$-multilinear in the sense that
	      \begin{equation*}
		      \mu_k \left( a_1, \dots, a_{i-1}, r \cdot a_i, \dots, a_k \right) =
		      (-1)^{\degb{r} \cdot \left( 1 + \degb{a_1} + \dots + \degb{a_{i-1}} \right)}
		      r \cdot \mu_k \left( a_1, \dots, a_k \right)
	      \end{equation*}
	      and satisfy
	      \begin{equation*}
		      \nnorm[\mu_k \left( a_1, \dots, a_k \right)] \leq \nnorm[a_1] \cdots \nnorm[a_k]
	      \end{equation*}
	      for all $a_1, \dots, a_k \in A$ and $r \in R$.
\end{enumerate}
The identity $\mu^2 = 0$, i.e., $\mu$ is a differential,
is equivalent to the identity $\corest{\mu} \circ \mu = 0$, which,
written explicitly, states that
\begin{equation*}
	\corest{\mu} \left( \mu \left( l \right) \right) = (-1)^{\degb{l_{(1)}}} \corest{\mu}
	\left( l_{(1)} \cotimes \corest{\mu} \left( l_{(2)} \right) \cotimes l_{(3)} \right) = 0
\end{equation*}
for all $l \in \tensf{A}[R]$. In terms of the associated multilinear maps,
the identity above is equivalent to a sequence of $\Ainf$\textbf{-identities} given by
\begin{align}\label{eq:ainf_for_mu_k_explicit}
	\sum_{k_1 + k_2 + k_3 = k} & (-1)^{\degb{a_1} + \dots + \degb{a_{k_1}}}
	\\
	                           & \mu_{k_1 + 1 + k_3}
	\left( a_1, \dots, a_{k_1}, \mu_{k_2} \left( a_{k_1+1}, \dots, a_{k_1+k_2} \right),
	a_{k_1+k_2+1}, \dots, a_k \right) = 0 \nonumber
\end{align}
for all $k \geq 0$ and $a_1, \dots, a_k \in A$.

\begin{rem}
	Let $\mathcal{R} = (R,d)$ be a differential graded-commutative Banach $\mathbbm{k}$-algebra
	and let $A$ be a graded Banach $R$-module. A (\textbf{standard, weakly curved}) $\Ainf$\textbf{-structure}
	on $A$ over $\mathcal{R}$ is given by:
	\begin{enumerate}
		\item A graded $R$-linear map $m_0 \colon R \rightharpoonup A$ of degree two which is
		      identified with the element $m_0(1) \in A^2$, called the \textbf{curvature},
		      such that $\nnorm[m_0(1)] < 1$.
		\item A graded $\mathbbm{k}$-linear map $m_1 \colon A \rightharpoonup A$ of degree one
		      which is a module derivation over $d$, i.e.,
		      \begin{equation*}
			      m_1 \left( r \cdot a \right) = dr \cdot a + (-1)^{\degb{r}} r \cdot m_1 \left( a \right),
			      \qquad \nnorm[m_1 \left( a \right)] \leq \nnorm[a]
		      \end{equation*}
		      for all $a \in A$ and $r \in R$.
		\item A sequence $m_k \colon A^{\times k} \rightharpoonup A$ for $k \geq 2$
		      of graded multilinear maps of degree $2 - k$,
		      which are $R$-multilinear in the sense that
		      \begin{equation*}
			      m_k \left( a_1, \dots, a_{i-1}, r \cdot a_i, \dots, a_k \right) =
			      (-1)^{\degb{r} \cdot \left( 2 - k + \degb{a_1} + \dots + \degb{a_{i-1}} \right)}
			      r \cdot m_k \left( a_1, \dots, a_k \right)
		      \end{equation*}
		      and satisfy
		      \begin{equation*}
			      \nnorm[m_k \left( a_1, \dots, a_k \right)] \leq \nnorm[a_1] \cdots \nnorm[a_k]
		      \end{equation*}
		      for all $a_1, \dots, a_k \in A$ and $r \in R$.
	\end{enumerate}
	The maps $\left( m_k \right)_{k \geq 0}$ are required to satisfy the $\Ainf$-identities, given by
	\begin{align}\label{eq:ainf_for_m_k_explicit}
		\sum_{k_1 + k_2 + k_3 = k} & (-1)^{k_1 + k_2 \cdot k_3 +
			                             (2 - k_2) \cdot \left( \degb{a_1} + \dots + \degb{a_{k_1}} \right)}
		\\
		                           & m_{k_1 + 1 + k_3}
		\left( a_1, \dots, a_{k_1}, m_{k_2} \left( a_{k_1+1}, \dots, a_{k_1+k_2} \right),
		a_{k_1+k_2+1}, \dots, a_k \right) = 0, \nonumber
	\end{align}
	for all $k \geq 0$ and $a_1, \dots, a_k \in A$.

	There is a bijective correspondence between standard $\Ainf$-structures on $A$ and
	shifted $\Ainf$-structures on the shifted module $A[1]$. The multilinear maps
	$m_k \colon A^{\times k} \rightharpoonup A$ correspond bijectively
	to multilinear maps $\mu_k \colon A[1]^{\times k} \rightharpoonup A[1]$ via the transformations
	\begin{equation}
		\begin{aligned}
			\mu_k \left( x_1, \dots, x_k \right) & = -(-1)^{\frac{k(k-1)}{2} + \sum_{i=1}^k (k-i) \degb{x_i}}
			\s m_k \left( \s^{-1} x_1, \dots, \s^{-1} x_k \right),
			\\
			m_k \left( a_1, \dots, a_k \right)   & = -(-1)^{\sum_{i=1}^k (k-i) \degb{a_i}}
			\s^{-1} \mu_k \left( \s a_1, \dots, \s a_k \right),
		\end{aligned}
		\label{eq:transformations-shifted-standard-A-inf-algebra}
	\end{equation}
	for $x_1, \dots, x_k \in A[1]$ and $a_1, \dots, a_k \in A$. Using the transformations,
	one can verify that the $\Ainf$-identities \eqref{eq:ainf_for_mu_k_explicit} for $\mu_k$
	are equivalent to the $\Ainf$-identities \eqref{eq:ainf_for_m_k_explicit} for $m_k$.
	For an explanation of the signs and other possible sign conventions, we refer to
	\cref{sub:a-inf-sign-conventions}.

	In order to minimize signs and avoid proliferation of suspensions, we will mostly work with
	shifted $\Ainf$-algebras and refer to them simply as $\Ainf$-algebras.
\end{rem}

\begin{rem}
	Let us explain the relation between our notion of a Banach $\Ainf$-algebra and other similar notions
	appearing in the literature.
	\begin{enumerate}
		\item Let $A$ be a graded $\mathbbm{k}$-module. By taking $R = \mathbbm{k}$ with $d_R = 0$ and
		      endowing $A$ with the trivial norm, with respect to which $A$ becomes a graded Banach
		      $\mathbbm{k}$-module, our notion of a (shifted) $\Ainf$-structure on $A[1]$ reduces
		      to the standard notion of a \textit{non-curved} $\Ainf$-algebra over $\mathbbm{k}$,
		      as discussed for example in \cite{Keller2001}.
		\item In \cite{Kleijn2021}, a \textbf{complete shifted} $\Ainf$-algebra over $\mathbbm{k}$
		      is given by:
		      \begin{enumerate}
			      \item A graded $\mathbbm{k}$-module $A$ endowed with a descending exhaustive filtration
			            \begin{equation*}
				            A = F^1 A \supseteq F^2 A \supseteq \dots
			            \end{equation*}
			            with respect to which $A$ is complete.
			      \item A degree one coderivation $Q \colon \tens{A}[\mathbbm{k}] \rightharpoonup \tens{A}[\mathbbm{k}]$ on the
			            \textit{tensor coalgebra} which preserves the filtration
			            \begin{equation*}
				            Q_n \left( F^{i_1} A \otimes \dots \otimes F^{i_n} A \right) \subseteq
				            F^{i_1 + \dots + i_n} A
			            \end{equation*}
			            and satisfies $Q^2 = 0$.
		      \end{enumerate}
		      Given a complete shifted $\Ainf$-algebra $A$ over $\mathbbm{k}$, one can endow $A$
		      with a non-Archimedean norm by setting
		      \begin{equation*}
			      \nnorm[a] \defeq \frac{1}{2^{\sup \Set{k \in \NZ}[a \in F^k A]}}.
		      \end{equation*}
		      The condition that $A$ is complete is equivalent to the condition that $A$, endowed
		      with the norm above, is a graded Banach $\mathbbm{k}$-module. The condition
		      that $Q_n$ preserves the filtration translates into the condition that each
		      $Q_n \colon A^{\otimes n} \rightharpoonup A$ is contractive and hence
		      $\nnorm[Q] \leq 1$ so we can complete $Q$ and obtain a coderivation
		      $\mu = \widehat{Q} \colon \tensf{A}[\mathbbm{k}] \rightharpoonup \tensf{A}[\mathbbm{k}]$ of the formal
		      tensor coalgebra which satisfies $\mu^2 = 0$. This way we obtain that
		      $\left( A, \mu \right)$ is a Banach $\Ainf$-algebra
		      over $\left( \mathbbm{k}, 0 \right)$.
		\item In \cite[Definition 1.1]{Solomon2016}, the authors introduce the notion of an $n$-dimensional
		      curved cyclic unital $\Ainf$-structure over a differential graded algebra. If one ignores
		      the cyclic structure and rephrases everything in terms of norms instead of valuations, then up
		      to some sign conversions, the resulting notion is equivalent to our definition
		      of a Banach $\Ainf$-algebra over $\mathcal{R}$. For more details, see
		      \cref{subsec:converting-a-inf-jake-to-banach} of \cref{appendix:sign-conversions-jake}.
	\end{enumerate}
\end{rem}

\begin{rem}
	Recall that the derivation $d \colon R \rightharpoonup R$ of $\mathcal{R}$ can be recovered from
	the coderivation $\mu$ by the identity $d = \pi_0 \circ \mu \circ i_0$
	where $i_0 \colon R \rightarrow \tensf{A}[R]$ is the natural inclusion and
	$\pi_0 \colon \tensf{A}[R] \rightarrow R$ is the natural projection (see Remark
	\ref{item:coderivation-underlying-derivation}). 
	Since we require in \cref{def:a-inf-Banach-algebra} that $\nnorm[\mu] \leq 1$,
	we see that we must also have $\nnorm[d] \leq 1$. That is, whenever we talk about a Banach $\Ainf$-algebra
	over a differential graded-commutative Banach algebra $\mathcal{R} = (R,d)$, we assume that the differential
	$d$ on $R$ also satisfies $\nnorm[d] \leq 1$.
\end{rem}

We note that given a Banach $\Ainf$-algebra $\mathcal{A} = \left( A, \mu \right)$,
the pair $\big( \tensf{A}[R], \mu \big)$ is a differential graded Banach $\mathcal{R}$-coalgebra.

\begin{dfn} \label{def:a-inf-Banach-algebra-mor}
	Let $\mathcal{R} = (R,d_R)$ and $\mathcal{S} = (S,d_S)$ be differential
	graded-commutative Banach $\mathbbm{k}$-algebras. Let $\mathcal{A} = (A,\mu)$ be a Banach
	$\Ainf$-algebra over $\mathcal{R}$, and let $\mathcal{B} = (B,\nu)$ be a Banach $\Ainf$-algebra
	over $\mathcal{S}$. A \textbf{morphism} $f \colon \mathcal{A} \rightarrow \mathcal{B}$ of Banach
	$\Ainf$-algebras is defined to be a morphism $f \colon \tensf{A}[R] \rightarrow \tensf{B}[S]$ of
	Banach coalgebras such that $\nnorm[f_0(1)]_B < 1$ and $f \circ \mu = \nu \circ f$.
	The morphism $f$ is called \textbf{strict} if $f_n = 0$ for all $n \neq 1$.
\end{dfn}

Explicitly, by \cref{prop:classification-morphisms-tensor-coalgebra-different-ground-algebras},
the data of a morphism $f \colon \mathcal{A} \rightarrow \mathcal{B}$ is given equivalently
by:
\begin{enumerate}
	\item A morphism $\varphi = \base{f} \colon \mathcal{R} \rightarrow \mathcal{S}$ of differential
	      graded Banach $\mathbbm{k}$-algebras. That is, a degree zero map $\varphi \colon R \rightarrow S$
	      with $\nnorm[\varphi] \leq 1$ which satisfies
	      $\varphi \left( r_1 \cdot r_2 \right) = \varphi(r_1) \cdot \varphi(r_2)$ for all
	      $r_1,r_2 \in R$, $\varphi \left( 1_R \right) = 1_S$, and $\varphi \circ d_R = d_S \circ \varphi$.
	\item An element $f_0(1) \in B^0$ called the \textbf{change of connection element} which
	      satisfies $\nnorm[f_0(1)] < 1$.\footnote{Even though a general Banach coalgebra morphism
		      $f \colon \tensf{A}[R] \rightarrow \tensf{B}[S]$ satisfies only that
		      $f_0 \left( 1 \right)$ is topologically nilpotent, we will work only with
		      morphisms which satisfy in addition that $\nnorm[f_0 \left( 1 \right)] < 1$.}
	\item A sequence of degree zero maps $f_k \colon A^{\times k} \rightarrow B$ for $k \geq 1$ which
	      are $R$-multilinear in the sense that
	      \begin{equation*}
		      f_k \left( a_1, \dots, a_{i-1}, r \cdot a_i, \dots, a_k \right) =
		      (-1)^{\degb{r} \cdot \left( \degb{a_1} + \dots + \degb{a_{i-1}} \right)}
		      \varphi(r) \cdot f_k \left( a_1, \dots, a_k \right)
	      \end{equation*}
	      and satisfy
	      \begin{equation*}
		      \nnorm[f_k \left( a_1, \dots, a_k \right)]  \leq \nnorm[a_1] \cdots \nnorm[a_k]
	      \end{equation*}
	      for all $a_1, \dots, a_k \in A$ and $r \in R$.
\end{enumerate}

The identity $f \circ \mu = \nu \circ f$ is equivalent to the identity
$\corest{f} \circ \mu = \corest{\nu} \circ f$ (see \cref{lm:f-circ-mu-nu-circ-f-corest-different-ground-algebras}), which, written explicitly, states that
\begin{equation*}
	(-1)^{\degb{l_{(1)}}}
	\corest{f} \left( l_{(1)} \otimes_R \corest{\mu} \left( l_{(2)} \right) \otimes_R l_{(3)} \right)
	= \sum_{n = 0}^{\infty} \corest{\nu} \left(
	\corest{f} \left( l_{(1)} \right) \otimes_S \dots \otimes_S \corest{f} \left( l_{(n)} \right) \right)
\end{equation*}
for all $l \in \tensf{A}[R]$. In terms of the associated multilinear maps,
the identity above is equivalent to the sequence of identities given by
\begin{align}\label{eq:ainf_morphism_explicit}
	\sum_{k_1 + k_2 + k_3 = k} (-1)^{\degb{a_1} + \dots + \degb{a_{k_1}}}
	 & f_{k_1 + 1 + k_3} \left( a_1, \dots, a_{k_1},
	\mu_{k_2} \left( a_{k_1 + 1}, \dots, a_{k_1 + k_2} \right),
	a_{k_1 + k_2 + 1}, \dots, a_k \right) = \nonumber
	\\
	\sum_{\substack{n \geq 0 \\ k_1 + \dots + k_n = k}}
	 & \nu_n \left( f_{k_1} \left( a_1, \dots, a_{k_1} \right), \dots,
	f_{k_n} \left( a_{k_1 + \dots + k_{n-1} + 1}, \dots, a_k \right) \right)
\end{align}
for all $k \geq 0$ and $a_1, \dots, a_k \in A$. The identities above involve infinite
summation which makes sense by our assumptions. For example, the identity
for $k = 0$ reads
\begin{equation*}
	f_1 \left( \mu_0 \left( 1 \right) \right) = \sum_{n = 0}^{\infty}
	\nu_n \left( f_0 \left( 1 \right), \dots, f_0 \left( 1 \right) \right)
	= \corest{\nu} \left( \Exp{f_0(1)} \right)
\end{equation*}
and the right-hand side indeed converges as
\begin{equation*}
	\nnorm[\nu_n \left( f_0 \left( 1 \right), \dots, f_0 \left( 1 \right) \right)] \leq
	\nnorm[\nu_n] \cdot \nnorm[f_0(1)]^n \leq \nnorm[f_0(1)]^n \to 0
\end{equation*}
and $B$ is Banach.

\begin{dfn} \label{dfn:a-infinity-unit}
	Let $\mathcal{A} = (A,\mu)$ be a Banach $\Ainf$-algebra over $\mathcal{R} = (R,d_R)$. An element
	$e \in A^{-1}$ is called a \textbf{(strong) unit} for $\mathcal{A}$ if the following conditions hold:
	\begin{enumerate}
		\item $\nnorm[e] \leq 1$.
		\item $\mu_1(e) = 0$.
		\item $\mu_2(e,a) = (-1)^{\degb{a} + 1} \mu_2 \left( a, e \right) = a$ for all $a \in A$.\footnote{We note
			      that with our choice of signs, a unit $e \in A[1]^{-1} = A^0$ for a shifted $\Ainf$-algebra structure on the
			      shifted module $A[1]$ is not a unit for the product $m_2 \colon A \otimes A \rightarrow A$ defined
			      via the transformations \eqref{eq:transformations-shifted-standard-A-inf-algebra}. We do have
			      $m_2 \left( -e, a \right) = m_2 \left( a, -e \right) = a$ for $a \in A$. See also \cref{tab:a-infinity-conversions-mu_k-m_k}.}
		\item $\mu_k(a_1,\dots,a_k) = 0$ whenever $k > 2$ and $a_1,\dots,a_k \in A$ with $a_i = e$
		      for some $1 \leq i \leq k$.
	\end{enumerate}
	When $e$ is a unit for $\mathcal{A}$, the triple $\mathcal{A} = (A,\mu,e)$ is called a
	\textbf{unital Banach} $\Ainf$\textbf{-algebra}.
\end{dfn}

\begin{dfn}
	Let $\mathcal{A} = (A,\mu,e_A)$ be a unital Banach $\Ainf$-algebra over $\mathcal{R}$ and let
	$\mathcal{B} = (B,\nu,e_B)$ be a unital Banach $\Ainf$-algebra over $\mathcal{S}$. A morphism
	$f \colon \mathcal{A} \rightarrow \mathcal{B}$ of Banach $\Ainf$-algebras is called \textbf{unital}
	if $f_1(e_A) = e_B$ and $f_k \left( a_1,\dots,a_k \right) = 0$ whenever $k \geq 2$ and $a_1,\dots,a_k \in A$
	with $a_i = e_A$ for some $1 \leq i \leq k$.
\end{dfn}

\begin{ex} \label{ex:jake-a-inf-algebra-vs-ours}
	Let $\left( X, \omega \right)$ be a symplectic manifold and let $L$ be a connected Lagrangian submanifold of $X$.
	A major family of examples which fall under \cref{def:a-inf-Banach-algebra} arise in symplectic geometry
	as $\Ainf$-structures coming from operations on the space of differential forms on a Lagrangian
	submanifold $L$. In \cite{Solomon2016}, under appropriate conditions on $X$ and $L$, the authors construct a family of
	cyclic unital $\Ainf$-structures $\mathfrak{m}^{J, \gamma}$ on the graded $R$-module $\cdiff{L}[][][R]$ of
	$R$-valued differential forms on $L$, where $R$
	is an extension of the Novikov ring by some formal variables. The family $\mathfrak{m}^{J, \gamma}$ is parametrized
	by almost complex structures $J$ and closed elements $\gamma$ belonging to an ideal of a differential graded-commutative algebra.
	Let us describe briefly the simplest family of such examples, constructed explicitly in
	\cite{Solomon2016}, to which we refer for more details.

	Fix an odd $n$ and consider the pair $(X,L) = (\cp, \rp)$ as a symplectic
	manifold and a Lagrangian submanifold, where the symplectic form $\omega_{\mathrm{FS}}$ is the Fubini-Study form.
	Denote by $\beta \in \homo{\cp}[][\rp][\ZZ] \cong \ZZ$ the element with minimal positive symplectic area
	$\omega_{\mathrm{FS}} \left( \beta \right)$.
	In \cite{Solomon2016}, the authors construct a sequence of $\RR$-linear operations
	$\mathfrak{q}_{k}^{d} \colon \cdiff{L}[]^{\otimes k} \rightharpoonup \cdiff{L}[]$ for $k \geq 0$ and $d \geq 0$,
	having degree $2 - k - d \cdot \left( n + 1 \right)$. The operators are constructed using pullback and pushforward
	operators on the moduli space of holomorphic discs with $k + 1$ marked points representing the relative homology class
	$d \cdot \beta$.

	Let $\mathbbm{k} = \RR$ and let $R = \pows{\RR}[T]$ be the
	graded ring of formal power series in the \textit{even} formal variable $T$, where $\deg{T} = n + 1$ and
	$\nnorm[T] = 2^{-\omega_{\mathrm{FS}} \left( \beta \right)}$.
	Endow the graded ring of real-valued differential forms $\cdiff{L}[]$ with the
	trivial norm and consider the graded Banach $R$-module $\cdiff{L}[][][R] \defeq R \cotimes \cdiff{L}[]$
	of $R$-valued differential forms. The operators $\mathfrak{q}_k^{d}$ extend to contractive operators
	$\overline{\mathfrak{q}}_{k}^{d} \colon \cdiff{L}[][][R]^{\cotimes k} \rightharpoonup \cdiff{L}[][][R]$.
	Summing the operators with coefficients, one obtains operators
	$\mathfrak{m}_k \colon \cdiff{L}[][][R]^{\cotimes k} \rightharpoonup \cdiff{L}[][][R]$ of
	degree $2 - k$ acting by\footnote{In \cite{Solomon2016}, the operators $\mathfrak{q}_k^d$ are denoted by
		$\mathfrak{q}_{k,0}^{d \cdot \beta}$ instead of $\mathfrak{q}_{k}^{d}$, and they give rise to the $\Ainf$-structure
		$\mathfrak{m}^{\gamma,J}$, where $\gamma = 0$ and $J$ is the standard holomorphic structure on $\cp$.}
	\begin{equation*}
		\mathfrak{m}_k \left( c_1, \dots, c_k \right) =
		\sum_{d \geq 0} \overline{\mathfrak{q}}_{k}^{d} \left( c_1, \dots, c_k \right) T^d.
	\end{equation*}
	The formal sum is well-defined since $\nnorm[T^d] \to 0$, and since $\overline{\mathfrak{q}}_{0}^{0} = 0$, we have
	$\nnorm[\mathfrak{m}_0 \left( 1 \right)] < 1$. By shifting and adding signs (see \cref{appendix:sign-conversions-jake}),
	one obtains a unital Banach $\Ainf$-structure on ${\cdiff{L}[][][R]}[1]$ in the sense
	of \cref{def:a-inf-Banach-algebra,dfn:a-infinity-unit}. The unit of the $\Ainf$-algebra is $\s \left( 1_C \right)$, where
	$1_C = 1_R \cotimes 1_{\cdiff{L}[]}$ is the unit of the algebra of $R$-valued differential forms on $L$.
\end{ex}

We end this section with the definitions of scalar restriction and extension for Banach $\Ainf$-algebras.
Let $\mathcal{R}$ and $\mathcal{S}$ be two differential graded-commutative
Banach $\mathbbm{k}$-algebras, and let $\varphi \colon \mathcal{R} \rightarrow \mathcal{S}$ be a
morphism of differential graded Banach $\mathbbm{k}$-algebras.
In \cref{sec:scalar-restriction-formal-tensor-coalgebras}, we discussed scalar restriction
and extension for formal tensor coalgebras, possibly equipped with a generalized coderivation.
When the coderivation satisfies the conditions \eqref{eq:banach-a-inf-mu-conditions}, the
scalar restriction (resp.\ extension) of the coderivation is readily verified to satisfy
the same conditions. Thus, the following definition makes sense:
\begin{dfn} (Scalar Restriction and Extension of Banach $\Ainf$-algebras) \label{dfn:scalar-rest-ext-a-infinity}
	\begin{enumerate}
		\item Let $\mathcal{B} = \left( B, \nu \right)$ be a Banach $\Ainf$-algebra over $\mathcal{S}$.
		      The \textbf{scalar restriction} (or \textbf{pullback}) \textbf{of} $\mathcal{B}$ \textbf{along} $\varphi$
		      is the Banach $\Ainf$-algebra
		      $\varphi^{*} \left( \mathcal{B} \right) \defeq \left( \varphi^{*} \left( B \right), \varphi^{*} \left( \nu \right) \right)$
		      over $\mathcal{R}$, where the coderivation $\varphi^{*} \left( \nu \right)$ is given
		      by \cref{dfn:pullback-generalized-coderivation}.
		\item Let $\mathcal{A} = \left( A, \mu \right)$ be a Banach $\Ainf$-algebra over $\mathcal{R}$.
		      The \textbf{scalar extension} (or \textbf{pushforward}) \textbf{of} $\mathcal{A}$ \textbf{along} $\varphi$
		      is the Banach $\Ainf$-algebra
		      $\varphi_{!} \left( \mathcal{A} \right) \defeq \left( \varphi_{!} \left( A \right), \varphi_{!} \left( \mu \right) \right)$
		      over $\mathcal{S}$, where the coderivation $\varphi_{!} \left( \mu \right)$ is given
		      by \cref{dfn:scalar-ext-generalized-coderivation}.
	\end{enumerate}
\end{dfn}

With the definition above, the canonical morphism
$\resover{\varphi}$ (resp.\ $\resunder{\varphi}$) over $\varphi$, given
by \cref{eq:canonical-pullback-map} (resp.\ \cref{eq:canonical-extension-map}), can
be interpreted as a \textit{strict} morphism
$\resover{\varphi} \colon \varphi^{*} \left( \mathcal{B} \right) \rightarrow \mathcal{B}$
(resp.\ $\resunder{\varphi} \colon \mathcal{A} \rightarrow \varphi_{!} \left( \mathcal{A} \right)$) of
Banach $\Ainf$-algebras. The universal property of the scalar restriction (resp.\ extension)
can then be stated as: Any Banach $\Ainf$-morphism $f \colon \mathcal{A} \rightarrow \mathcal{B}$
over $\varphi$ factors uniquely as $f = \resover{\varphi} \circ \resover{f}$
(resp.\ $f = \resunder{f} \circ \resunder{\varphi}$) for some
morphism $\resover{f} \colon \mathcal{A} \rightarrow \varphi^{*} \left( \mathcal{B} \right)$
(resp.\ $\resunder{f} \colon \varphi_{!} \left( \mathcal{A} \right) \rightarrow \mathcal{B}$)
of Banach $\Ainf$-algebras over $\mathcal{R}$ (resp.\ $\mathcal{S}$).

\subsection{Pseudoisotopy of \texorpdfstring{$\Ainf$-}{A-infinity }algebras}
\label{subsec:pseudoisotopy-a-inf}

\begin{dfn} \label{dfn:pseudo-isotopy}
	Let $\mathcal{S} = \left( S, d_S \right)$ be a differential graded-commutative Banach $\mathbbm{k}$-algebra.
	Let $\mathcal{A}_0 = \left( A_0, \mu^0 \right)$ and $\mathcal{A}_1 = \left( A_1, \mu^1 \right)$ be two
	$\Ainf$-algebras over $\mathcal{S}$.  A \textbf{pseudoisotopy} between $\mathcal{A}_0$
	and $\mathcal{A}_1$ is given by the following data:
	\begin{enumerate}
		\item A differential graded-commutative Banach $\mathbbm{k}$-algebra $\mathfrak{R} = \left( R, d_R \right)$.
		\item An $\Ainf$-algebra $\mathfrak{A} = \left( A, \mu \right)$ over
		      $\mathfrak{R}$ together with two $\Ainf$-morphisms
		      $\evalmf^i \colon \mathfrak{A} \rightarrow \mathcal{A}_i$
		      such that the underlying DGA morphisms
		      $\evalm^i \defeq \base{\evalmf}^i \colon \mathfrak{R} \rightarrow \mathcal{S}$
		      are homotopic as maps of differential graded $\mathbbm{k}$-modules.
	\end{enumerate}
	We will often succinctly denote the data of a pseudoisotopy between
	$\mathcal{A}_0$ and $\mathcal{A}_1$ simply by $\mathfrak{A}$, the notation
	for the $\Ainf$-algebra itself. In this case, we say that $\mathfrak{A}$ is a pseudoisotopy between
	$\mathcal{A}_0$ and $\mathcal{A}_1$, leaving the base differential graded algebras and the
	morphisms implicit.

	When the $\Ainf$-morphisms $\evalmf^i$ of $\mathfrak{A}$ are strict, we will call $\mathfrak{A}$ a
	\textbf{strict pseudoisotopy}. When
	the $\Ainf$-algebras $\mathcal{A}_0, \mathcal{A}_1$ are unital, we also require
	the $\Ainf$-algebra $\mathfrak{A}$ and the morphisms $\evalmf^i$ to be unital
	and say that the pseudoisotopy is \textbf{unital}.
\end{dfn}

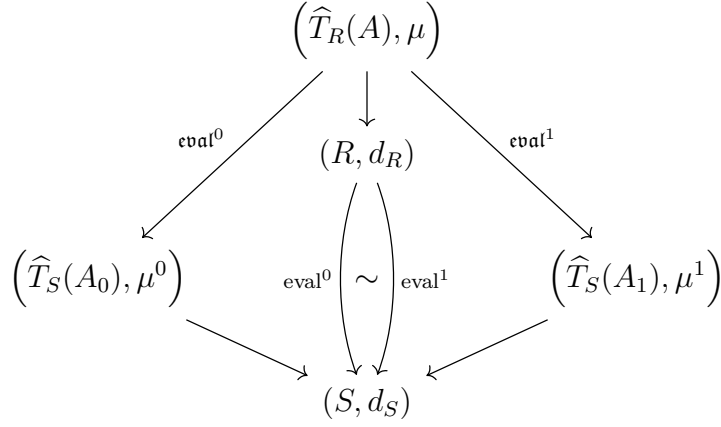
\begin{figure}
	\begin{tikzcd}
		& {\left( \tensf{A}[R], \mu \right)} \\
		& {\left( R, d_R \right)} \\
		{\left( \tensf{A_0}[S], \mu^0 \right)} & \sim & {\left( \tensf{A_1}[S], \mu^1 \right)} \\
		& {\left( S, d_S \right)}
		\arrow[from=1-2, to=2-2]
		\arrow["\evalmf^0"', from=1-2, to=3-1]
		\arrow["\evalmf^1", from=1-2, to=3-3]
		\arrow["\evalm^0"', curve={height=12pt}, from=2-2, to=4-2]
		\arrow["\evalm^1", curve={height=-12pt}, from=2-2, to=4-2]
		\arrow[from=3-1, to=4-2]
		\arrow[from=3-3, to=4-2]
	\end{tikzcd}
	\caption{The data of pseudoisotopy.}
	\label{fig:pseudoisotopy}
\end{figure}

In terms of the associated formal tensor coalgebras, the data of a pseudoisotopy between $\mathcal{A}_0$
and $\mathcal{A}_1$ is depicted in \cref{fig:pseudoisotopy}.

\begin{rem}
	When the pseudoisotopy $\mathfrak{A}$ is strict, the fact that $\evalmf^i$ are $\Ainf$-morphisms
	translates into the following relations:
	\begin{align*}
		\evalmf^i \left( \mu_0(1) \right)                             & = \mu^i_0(1),
		\\
		\evalmf^i \left( \mu_k \left( a_1, \dots, a_k \right) \right) & =
		\mu^i_k \left( \evalmf^{i} \left( a_1 \right), \dots, \evalmf^{i} \left( a_k \right) \right),
		\,\,\,                                                        & \forall a_1, \dots, a_k \in A,
		\\
		\evalmf^i \left( r \cdot a \right)                            & =
		\evalm^i \left( r \right) \cdot \evalmf^i \left( a \right),
		\,\,\,                                                        & \forall r \in R, \, \forall a \in A.
	\end{align*}
\end{rem}

\begin{rem} \label{rem:abstract-pseudoisotopy-vs-jake}
	Our notion of pseudoisotopy is an adaptation and generalization of the notion appearing in \cite[][]{Solomon2016,Solomon2016a}
	to an abstract algebraic setting more appropriate to this work. Let us denote by $I = [0,1]$ the unit interval.
	In \cite[][]{Solomon2016,Solomon2016a}, ignoring the cyclic structure, a pseudoisotopy between two
	$\Ainf$-structures $\left( \mathfrak{m}, \be \right)$ and $\left( \mathfrak{m}', \be' \right)$
	on the \textit{same} underlying $R$-module $\cdiff{L}[][][R]$ is
	an $\Ainf$-structure $\left( \widetilde{\mathfrak{m}}, \widetilde{\be} \right)$ on the module
	$\cdiff{I \times L}[][][R]$ over the differential graded-commutative $\RR$-algebra $\cdiff{I}[][][R]$.
	The module structure of $\cdiff{I \times L}[][][R]$ over $\cdiff{I}[][][R]$ is induced
	by the map $\pi_I^{*}$, where $\pi_I \colon I \times L \rightarrow I$ is the projection.
	The module $\cdiff{I \times L}[][][R]$ comes equipped with two maps
	$j_0^{*}, j_1^{*} \colon \cdiff{I \times L}[][][R] \rightarrow \cdiff{L}[][][R]$,
	where $j_i \colon L \rightarrow I \times L$ are the inclusions $j_i \left( p \right) = \left( i, p \right)$ for $i = 0,1$,
	and the requirement is that
	\begin{equation} \label{eq:pseudoisotopy-jake-style}
		\begin{aligned}
			j_0^{*} \left( \widetilde{\mathfrak{m}}_k \left( \widetilde{c}_1, \dots, \widetilde{c}_k \right) \right) & =
			\mathfrak{m}_k \left( j_0^{*} \left( \widetilde{c}_1 \right), \dots, j_0^{*} \left( \widetilde{c}_k \right) \right),
			\qquad
			j_0^{*} \left( \widetilde{\be} \right) = \be,
			\\
			j_1^{*} \left( \widetilde{\mathfrak{m}}_k \left( \widetilde{c}_1, \dots, \widetilde{c}_k \right) \right) & =
			\mathfrak{m}'_k \left( j_1^{*} \left( \widetilde{c}_1 \right), \dots, j_1^{*} \left( \widetilde{c}_k \right) \right),
			\qquad
			j_1^{*} \left( \widetilde{\be} \right) = \be',
		\end{aligned}
	\end{equation}
	(see \cite[Definition 1.4]{Solomon2016}).

	As discussed in \cref{appendix:sign-conversions-jake}, by rephrasing everything in terms of norms, shifting and adding signs,
	one can convert $\left( \cdiff{L}[][][R], \mathfrak{m}, \be \right)$
	(resp.\ $\left( \cdiff{L}[][][R], \mathfrak{m}', \be' \right)$)
	to a unital Banach $\Ainf$-algebra
	$\mathcal{A}_0 = \left( {\cdiff{L}[][][R]}[1], \mu^0, e_0 \right)$
	(resp.\ $\mathcal{A}_1 = \left( {\cdiff{L}[][][R]}[1], \mu^1, e_1 \right)$)
	over $R$ in the sense of \cref{def:a-inf-Banach-algebra,dfn:a-infinity-unit}.
	Similarly, one can convert the pseudoisotopy
	$\left( \cdiff{I \times L}[][][R], \widetilde{\mathfrak{m}}, \widetilde{\be} \right)$
	to a Banach $\Ainf$-algebra $\mathfrak{A} = \left( {\cdiff{I \times L}[][][R]}[1], \mu, e \right)$ over
	$\mathfrak{R} \defeq \left( \cdiff{I}[][][R], d_{\cdiff{I}[][][R]} \right)$. Requirements \eqref{eq:pseudoisotopy-jake-style}
	then state that the maps $j_0^{*}, j_1^{*}$ induce \textit{strict} unital $\Ainf$-morphisms
	$\evalmf^0, \evalmf^1$. The morphisms $\evalmf^i$ satisfy
	\begin{equation*}
		\evalmf^i \left( \alpha \cdot m \right) = \evalm^i \left( \alpha \right) \cdot \evalmf^i \left( m \right)
	\end{equation*}
	for $\alpha \in \cdiff{I}[][][R]$ and $m \in {\cdiff{I \times L}[][][R]}[1]$. Here, the maps
	$\evalm^i \colon \cdiff{I}[][][R] \rightarrow R$ are morphisms of
	differential graded-commutative Banach $\mathbbm{k}$-algebras, given by
	the evaluation maps at $t = i$ for $i = 0, 1$.
	In other words, the morphisms $\evalmf^i$ are morphisms over the morphisms $\evalm^i$.
	Finally, the morphisms $\evalm^0,\evalm^1$ are homotopic via the integration homotopy
	$h \colon \cdiff{I}[][][R] \rightharpoonup R$ given by $h \left( \alpha \right) = \int_I \alpha$.

	To summarize, a pseudoisotopy between $\left( \mathfrak{m}, \be \right)$ and $\left( \mathfrak{m}', \be' \right)$
	in the sense of \cite[Definition 1.4]{Solomon2016} gives rise to a \textit{strict} unital pseudoisotopy between
	the corresponding shifted $\Ainf$-algebras $\mathcal{A}_0$ and $\mathcal{A}_1$ in the sense of \cref{dfn:pseudo-isotopy}.
\end{rem}

In order to handle several similar constructions later in the text, we introduce
here the notion of a strong pseudoisotopy as follows. Let $\mathcal{F}$ be a
\textbf{cohomological invariant} of Banach $\Ainf$-algebras; that is, a
construction assigning to each Banach $\Ainf$-algebra $\mathcal{A}$ over
a differential graded-commutative Banach $\mathbbm{k}$-algebra $\mathcal{R}$
a graded $\cohom{\mathcal{R}}[]$-module $\mathcal{F}(\mathcal{A})$,
and to each Banach $\Ainf$-morphism $f \colon \mathcal{A} \rightarrow \mathcal{B}$
over $\psi \colon \mathcal{R} \rightarrow \mathcal{S}$ a
morphism of graded modules
$\mathcal{F}(f) \colon \mathcal{F}(\mathcal{A}) \rightarrow \mathcal{F}(\mathcal{B})$
over the induced map $\cohom{\psi}[] \colon \cohom{\mathcal{R}}[] \rightarrow \cohom{\mathcal{S}}[]$,
such that $\mathcal{F}$ respects composition and identities.
\begin{dfn} \label{dfn:F-strong-pseudoisotopy}
	Let $\mathcal{A}_0$ and $\mathcal{A}_1$ be two Banach $\Ainf$-algebras over the same differential graded-commutative Banach $\mathbbm{k}$-algebra $\mathcal{S}$. A pseudoisotopy $\mathfrak{A}$ over $\mathfrak{R}$ between $\mathcal{A}_0$ and $\mathcal{A}_1$ is called $\mathcal{F}$-\textbf{strong} if the following two conditions hold:
	\begin{enumerate}
		\item The underlying homotopic morphisms $\evalm^i = \base{\evalmf}^i \colon \mathfrak{R} \rightarrow \mathcal{S}$, for $i = 0, 1$, induce an \textit{isomorphism} of graded $\mathbbm{k}$-algebras
		      \begin{equation*}
			      \varphi \defeq \cohom{\evalm^0}[] = \cohom{\evalm^1}[] \colon \cohom{\mathfrak{R}}[] \overset{\sim}{\rightarrow} \cohom{\mathcal{S}}[]
		      \end{equation*}
		      between the cohomologies.
		\item The $\Ainf$-morphisms $\evalmf^i \colon \mathfrak{A} \rightarrow \mathcal{A}_i$ over $\evalm^i$, for $i = 0, 1$, induce \textit{isomorphisms} of graded $\mathbbm{k}$-modules
		      \begin{equation*}
			      \mathcal{F}(\evalmf^i) \colon \mathcal{F}(\mathfrak{A}) \overset{\sim}{\rightarrow} \mathcal{F}(\mathcal{A}_i).
		      \end{equation*}
	\end{enumerate}
\end{dfn}

Note that since we assume $\varphi$ is an isomorphism, the second condition is equivalent to the requirement that the induced maps
$\varphi_{!} \left( \mathcal{F}(\mathfrak{A}) \right) \rightarrow \mathcal{F}(\mathcal{A}_i)$
are isomorphisms of graded $\cohom{\mathcal{S}}[]$-modules, or, equivalently, that
the induced maps
$\mathcal{F}(\mathfrak{A}) \rightarrow \varphi^{*} \left( \mathcal{F}(\mathcal{A}_i) \right)$
are isomorphisms of graded $\cohom{\mathcal{R}}[]$-modules.
Thus, for any cohomological invariant $\mathcal{F}$, an $\mathcal{F}$-strong
pseudoisotopy $\mathfrak{A}$ yields a canonical isomorphism of graded $\cohom{\mathcal{S}}[]$-modules
$\mathfrak{a} \colon \mathcal{F}(\mathcal{A}_0) \overset{\sim}{\rightarrow} \mathcal{F}(\mathcal{A}_1)$
given by the composition
\begin{equation*}
	\mathfrak{a} \colon \mathcal{F}(\mathcal{A}_0) \overset{\sim}{\rightarrow} \varphi_{!} \left( \mathcal{F}(\mathfrak{A}) \right) \overset{\sim}{\rightarrow} \mathcal{F}(\mathcal{A}_1).
\end{equation*}

\subsection{Bounding Cochains and Gauge Equivalence}
\label{subsec:bounding-cochains-gauge-equivalence}

\begin{dfn} (Bounding Cochains) \label{dfn:bounding-chains}
	Let $\mathcal{R} = (R,d)$ be a differential graded-commutative Banach $\mathbbm{k}$-algebra.
	\begin{enumerate}
		\item Let $\mathcal{A} = \left( A, \mu \right)$ be a Banach $\Ainf$-algebra over $\mathcal{R}$.
		      An element $b \in A^0$ with $\nnorm[b] < 1$ is called a
		      \textbf{strong bounding cochain} if it satisfies the Maurer--Cartan equation
		      \begin{equation}
			      \corest{\mu} \left( \Exp{b} \right) = \sum_{i=0}^{\infty} \mu_i \left( b^{\cotimes i} \right) = 0.
		      \end{equation}
		      We denote the set of all strong bounding cochains in $A$ by $\mc{\mathcal{A}} \subseteq A^0$.
		\item Let $\mathcal{A} = \left( A, \mu, e \right)$ be a unital Banach $\Ainf$-algebra over $\mathcal{R}$
		      and let $c \in R^2$ with $dc = 0$ and $\nnorm[c] < 1$. An element $b \in A^0$ with $\nnorm[b] < 1$
		      is called a \textbf{weak bounding cochain} if it satisfies the
		      inhomogeneous Maurer--Cartan equation
		      \begin{equation} \label{eq:maurer-cartan-ce}
			      \corest{\mu} \left( \Exp{b} \right) = \sum_{i=0}^{\infty} \mu_i \left( b^{\cotimes i} \right) = c \cdot e.
		      \end{equation}
		      We will denote the set of all weak bounding cochains in $A$ with a fixed $c$ by
		      $\mc{\mathcal{A}}[c] \subseteq A^0$.
	\end{enumerate}
\end{dfn}
Note that the infinite sums appearing in \cref{dfn:bounding-chains} converge
by our assumptions that $A$ is Banach and $\nnorm[b] < 1$. The Maurer--Cartan equation
makes sense in some cases even when $\nnorm[b] \geq 1$. This happens for example when $A$ corresponds
to a differential graded algebra, or, more generally, when $\mu_k = 0$ for $k \gg 0$.
However, we will work only with solutions which satisfy $\nnorm[b] < 1$.

Given a morphism of Banach $\Ainf$-algebras, we can use the pushforward map, defined in \cref{sub:grouplike-exp},
to pushforward bounding cochains:
\begin{lm}[Functoriality of Bounding Cochains Under $\Ainf$-morphisms] \hfill \leavevmode \linebreak
	Let $\mathcal{R} = \left( R, d_R \right), \mathcal{S} = \left( S, d_S \right)$ be two
	differential graded-commutative Banach $\mathbbm{k}$-algebras and let
	$\varphi \colon \mathcal{R} \rightarrow \mathcal{S}$ be a morphism of differential graded Banach
	$\mathbbm{k}$-algebras. Let $\mathcal{A} = \left( A, \mu \right)$ be a Banach $\Ainf$-algebra over
	$\mathcal{R}$, and let $\mathcal{B} = \left( B, \nu \right)$ be a Banach $\Ainf$-algebra over $\mathcal{S}$.
	Let $f \colon \mathcal{A} \rightarrow \mathcal{B}$ be a morphism of Banach $\Ainf$-algebras over $\varphi$
	and let $b \in A^0$.
	\begin{enumerate}
		\item If $b$ is a strong bounding cochain then $\mcfunc{f} \left( b \right) \in B^0$
		      is a strong bounding cochain. Hence, we have an induced pushforward map
		      $\mcfunc{f} \colon \mc{\mathcal{A}} \rightarrow \mc{\mathcal{B}}$.
		\item Assume $\mathcal{A}, \mathcal{B}$ and $f$ are unital. If $b$ is a weak
		      bounding cochain with $\corest{\mu} \left( \Exp{b} \right) = c \cdot e_A$, then
		      $\mcfunc{f} \left( b \right) \in B^0$ is a weak bounding cochain with
		      $\corest{\nu} \left( \Exp{\mcfunc{f} \left( b \right)} \right) = \varphi \left( c \right) \cdot e_B$.
		      Hence, we have an induced pushforward map
		      $\mcfunc{f} \colon \mc{\mathcal{A}}[c] \rightarrow \mc{\mathcal{B}}[\varphi(c)]$.
	\end{enumerate}
\end{lm}
\begin{proof}
	Let $b \in A^0$ with $\nnorm[b] < 1$. Then
	\begin{equation*}
		\nnorm[\mcfunc{f} \left( b \right)] = \nnorm[\sum_{i=0}^{\infty} f_i \left( b^{\cotimes i} \right)]
		\leq \max_{i \geq 0} \, \nnorm[f_i] \cdot \nnorm[b]^i < 1
	\end{equation*}
	since $\nnorm[f_0] < 1$ and $\nnorm[b] < 1$.
	As $\Delta^3 \left( \Exp{b} \right) = \Exp{b} \cootimes \Exp{b} \cootimes \Exp{b}$,
	we have
	\begin{equation} \label{eq:corest-nu-exp-mcfunc-b}
		\corest{\nu} \left( \Exp{\mcfunc{f} \left( b \right)} \right)
		\stackrel{\eqref{eq:exp-mcfunc}}{=}
		\corest{\nu} \left( f \left( \Exp{b} \right) \right) =
		\corest{f} \left( \mu \left( \Exp{b} \right) \right)
		\stackrel{\eqref{eq:generalized-coder-coextension}}{=}
		\corest{f} \left( \Exp{b} \cotimes \corest{\mu} \left( \Exp{b} \right) \cotimes \Exp{b} \right).
	\end{equation}

	When $\corest{\mu} \left( \Exp{b} \right) = 0$, we see that
	$\corest{\nu} \left( \Exp{\mcfunc{f} \left( b \right)} \right) = 0$ and hence
	$\mcfunc{f} \left( b \right)$ is a strong bounding cochain.
	When $\mathcal{A},\mathcal{B}$ and $f$ are unital, and $b$ is a weak bounding cochain with
	$\corest{\mu} \left( \Exp{b} \right) = c \cdot e_A$, then
	\begin{equation*}
		\corest{\nu} \left( \Exp{\mcfunc{f} \left( b \right)} \right)
		\stackrel{\eqref{eq:corest-nu-exp-mcfunc-b}}{=}
		\corest{f} \left( \Exp{b} \cotimes \left( c \cdot e_A \right) \cotimes \Exp{b} \right) =
		f_1 \left( c \cdot e_A \right) = \varphi \left( c \right) \cdot e_B
	\end{equation*}
	by the unitality of $f$. Note that $c' \defeq \varphi \left( c \right)$ satisfies
	$d_S \left( c' \right) = 0$ and $\nnorm[c'] < 1$ because $\varphi$ commutes with the differentials
	and is contractive. Hence, $\mcfunc{f} \left( b \right)$ is a weak bounding cochain belonging to
	$\mc{\mathcal{B}}[\varphi(c)]$.
\end{proof}

\begin{rem} \label{rem:bounding-chains-dga}
	Let $b \in A^0$ be a weak bounding cochain with $\corest{\mu} \left( \Exp{b} \right) = c \cdot e$ for
	some $c \in R^2$. When the map $R \mapsto A$ given by $r \mapsto r \cdot e$ is injective, the element $c$
	is determined uniquely by the weak bounding cochain $b$. In this case the condition $dc = 0$
	which appears in \cref{dfn:bounding-chains}, part $(2)$, is satisfied automatically. To see this,
	apply the $\Ainf$ identity $\corest{\mu} \circ \mu = 0$ to $\Exp{b}$ to get
	\begin{equation*}
		\begin{aligned}
			0 & = \corest{\mu} \left( \mu \left( \Exp{b} \right) \right) =
			\corest{\mu} \left( \Exp{b} \cotimes \corest{\mu} \left( \Exp{b} \right) \cotimes \Exp{b} \right) =
			\corest{\mu} \left( \Exp{b} \cotimes \left( c \cdot e \right) \cotimes \Exp{b} \right)
			\\
			  & =
			\sum_{n,m=0}^{\infty} \mu_{n+m+1} \left( b^{\cotimes n} \cotimes \left( c \cdot e \right) \cotimes
			b^{\cotimes m} \right) =
			\mu_1(c \cdot e) + \mu_2 \left( c \cdot e, b \right) + \mu_2 \left( b, c \cdot e \right)
			\\
			  & = dc \cdot e + c \cdot \mu_1(e) + c \cdot \mu_2 \left(e, b \right) + c \cdot \mu_2 \left( b, e \right)
			= dc \cdot e + c \cdot b - c \cdot b = dc \cdot e.
		\end{aligned}
	\end{equation*}
	Similarly, when the map $r \mapsto r \cdot e$ is an isometry, the condition $\nnorm[c] < 1$ is
	satisfied automatically (since $\nnorm[b], \nnorm[\mu_0(1)] < 1$ and $\nnorm[\mu] \leq 1$).
\end{rem}

\begin{dfn} \label{def:gauge-equivalence-bounding-chains}
	Let $\mathcal{S}$ be a differential graded-commutative Banach $\mathbbm{k}$-algebra.
	\begin{enumerate}
		\item Let $\mathcal{A}_0$ and $\mathcal{A}_1$ be two Banach
		      $\Ainf$-algebras over $\mathcal{S}$ and assume we have a fixed pseudoisotopy $\mathfrak{A}$ between
		      $\mathcal{A}_0$ and $\mathcal{A}_1$. Let $b_0 \in \mc{\mathcal{A}_0}$ and $b_1 \in \mc{\mathcal{A}_1}$
		      be two strong bounding cochains. We will say that $b_0$ and $b_1$ are
		      $\mathfrak{A}$-\textbf{gauge-equivalent} if there exists a $b \in \mc{\mathfrak{A}}$ such that we have
		      \begin{equation} \label{eq:gauge-equivalence-bounding-chains}
			      \mcfunc{\evalmf}^0 \left( b \right) = b_0, \,\,\,
			      \mcfunc{\evalmf}^1 \left( b \right) = b_1.
		      \end{equation}
		\item  Let $\mathcal{A}_0$ and $\mathcal{A}_1$ be two unital Banach
		      $\Ainf$-algebras over $\mathcal{S}$ and assume we have a fixed unital pseudoisotopy $\mathfrak{A}$
		      over $\mathfrak{R} = \left( R, d_R \right)$ between $\mathcal{A}_0$ and $\mathcal{A}_1$.
		      Let $b_0 \in \mc{\mathcal{A}_0}[c_0]$ and $b_1 \in \mc{\mathcal{A}_1}[c_1]$
		      be two weak bounding cochains. We will say that $b_0$ and $b_1$ are
		      $\mathfrak{A}$-\textbf{gauge-equivalent} if there exists an element $c \in R^2$ with
		      $d_R \left( c \right) = 0$ and $\nnorm[c] < 1$ such that
		      \begin{equation*}
			      \evalm^0 \left( c \right) = c_0, \,\,\,
			      \evalm^1 \left( c \right) = c_1
		      \end{equation*}
		      and an element $b \in \mc{\mathfrak{A}}[c]$ such that
		      \begin{equation*}
			      \mcfunc{\evalmf}^0 \left( b \right) = b_0, \,\,\,
			      \mcfunc{\evalmf}^1 \left( b \right) = b_1.
		      \end{equation*}
	\end{enumerate}
\end{dfn}

\begin{rem}
	The condition $\mcfunc{\evalmf}^i \left( b \right) = b_i$ in \cref{def:gauge-equivalence-bounding-chains}
	implies by functoriality that $b_i \in \mc{\mathcal{A}_i}[\evalm^i \left( c \right)]$,
	in addition to $b_i \in \mc{\mathcal{A}_i}[c_i]$. Hence,
	when the map $R \mapsto A_i$ given by $r \mapsto r \cdot e_{A_i}$ is injective, the condition
	$\evalm^i \left( c \right) = c_i$ is redundant and already implied by the condition
	$\mcfunc{\evalmf}^i \left( b \right) = b_i$ (see \cref{rem:bounding-chains-dga}).
\end{rem}

\begin{rem}
	When the $\Ainf$-morphisms $\evalmf^i$ are strict, we have
	$\mcfunc{\evalmf}^i = \evalmf^i$ and so the definition of
	$\mathfrak{A}$-gauge-equivalence reduces to the more familiar relation
	\begin{equation} \label{eq:strict-gauge-equivalence-rel}
		\evalmf^0(b) = b_0, \,\,\, \evalmf^1(b) = b_1.
	\end{equation}
\end{rem}

\begin{rem} \label{rem:abstract-gauge-equivalence-vs-jake}
	Similarly to the case of pseudoisotopy, our notion of $\mathfrak{A}$-gauge-equivalence is an adaptation and
	generalization of the notion appearing in \cite{Solomon2016a} to an abstract algebraic
	setting more appropriate to this work.
	In \cite[Definition 1.1]{Solomon2016a}, the notion of a bounding pair $\left( \gamma, b \right)$
	with respect to an almost complex structure $J$ is introduced and in \cite[Definition 3.12]{Solomon2016a},
	it is defined when a bounding pair $\left( \gamma, b \right)$ with respect to $J$ and
	$\left( \gamma', b' \right)$ with respect to $J'$ are gauge-equivalent using the notion of pseudoisotopy.

	As discussed in \cref{rem:abstract-pseudoisotopy-vs-jake,subsec:converting-a-inf-jake-to-banach},
	in this situation one obtains unital Banach $\Ainf$-algebras
	$\mathcal{A}_0 = \left( {\cdiff{L}[][][R]}[1], \mu^{0}, e_0 \right)$
	from $\mathfrak{m}_k^{\gamma, J}$ and
	$\mathcal{A}_1 = \left( {\cdiff{L}[][][R]}[1], \mu^{1}, e_1 \right)$
	from $\mathfrak{m}_k^{\gamma',J'}$. In addition, one obtains a strict pseudoisotopy
	$\mathfrak{A} = \left( {\cdiff{I \times L}[][][R]}[1], \mu, e \right)$ from
	$\widetilde{\mathfrak{m}}_k^{\widetilde{\gamma},J}$. The element $b_0 \defeq b$ (resp.\ $b_1 \defeq b'$) is a
	weak bounding cochain of $\mathcal{A}_0$ (resp.\ $\mathcal{A}_1$). In addition,
	according to \cite[Definition 3.12]{Solomon2016a}, there exists
	a weak bounding cochain $\mathfrak{b} \defeq \widetilde{b}$ of $\mathfrak{A}$ which satisfies the relations
	of \eqref{eq:strict-gauge-equivalence-rel}. Hence, $b_0$ and $b_1$ are $\mathfrak{A}$-gauge-equivalent according
	to \cref{def:gauge-equivalence-bounding-chains}.
\end{rem}

\begin{lm}
	Let $\mathcal{S} = \left( S, d_S \right)$ be a differential graded-commutative Banach $\mathbbm{k}$-algebra,
	and let $\mathcal{A}_0, \mathcal{A}_1$ be two unital Banach $\Ainf$-algebras over $\mathcal{S}$.
	Let $\mathfrak{A}$ be a unital pseudoisotopy over $\mathfrak{R} = \left( R, d_R \right)$
	between $\mathcal{A}_0$ and $\mathcal{A}_1$,
	and let $b_0 \in \mc{\mathcal{A}_0}[c_0]$ and $b_1 \in \mc{\mathcal{A}_1}[c_1]$ be two weak
	$\mathfrak{A}$-gauge-equivalent bounding cochains.
	Then $\eqcl{c_0} = \eqcl{c_1}$ in $\cohom{\mathcal{S}}[]$.
\end{lm}
\begin{proof}
	Choose a homotopy $h \colon R^{*} \rightharpoonup S^{*-1}$ between $\evalm^0$ and $\evalm^1$, so that
	\begin{equation*}
		d_S \circ h + h \circ d_R = \evalm^1 - \evalm^0.
	\end{equation*}
	Since $b_0$ and $b_1$ are $\mathfrak{A}$-gauge-equivalent, there exists an element $c \in R^2$
	with $d_R \left( c \right) = 0$ and $\evalm^i \left( c \right) = c_i$ for $i = 0,1$. But then
	\begin{equation*}
		c_1 - c_0 =
		\evalm^1 \left( c \right) - \evalm^0 \left( c \right) =
		d_S \left( h \left(c \right) \right) + h \left( d_R \left( c \right) \right) =
		d_S \left( h \left( c \right) \right).
	\end{equation*}
\end{proof}

\section{Cyclization on the Tensor Coalgebra} \label{sec:cyc-tensor-coalgebra}

In this section, we develop the algebraic machinery necessary to show that the actions of generalized coderivations
and coalgebra morphisms on the tensor coalgebra induce corresponding ``cyclic versions'' which behave functorially.
The cyclic versions are introduced by explicit formulas which involve the cyclic rotation
operator $\t$, and their functorial properties are proven via explicit calculations.

In \cref{sec:cyclization-generalized-coderivation}, we introduce the \textbf{cyclization}
$\cycl{\mu} \colon \tensr{V} \rightharpoonup \tensr{V}$ of a generalized coderivation
$\mu \colon \tens{V} \rightharpoonup \tens{V}$ and show in \cref{lm:cycl-comm-bracket}
that the process of cyclization commutes with the Lie bracket.
When $\mu$ is an odd coderivation which is a \textit{differential}, encoding an $\Ainf$-structure,
this implies that $\cycl{\mu}$ is also a differential, known as the \textbf{Hochschild differential} of the corresponding $\Ainf$-algebra.
In \cref{lm:coder-descends-quotient}, we prove an identity relating $\mu, \cycl{\mu}$, and the
rotation operator $\t$, which implies that $\cycl{\mu}$ descends to the cyclic quotient $\tensr{V} / \Im \left( \idd - \t \right)$.

\Cref{sec:cyclization-coalgebra-morphisms} introduces the \textbf{cyclization}
$\cycl{f} \colon \tensr{V} \rightarrow \tensr{W}$ of a tensor coalgebra morphism $f \colon \tens{V} \rightarrow \tens{W}$, possibly
with a non-zero change of connection term. Here, we show in \cref{lm:func-cycl-morphism,lm:func-cycl-morphism-coderivation}
that the process of cyclization commutes with composition of morphisms and is compatible with cyclization of coderivations.
We then prove an identity in \cref{lm:cycl-f-1-t} relating $f, \cycl{f}$, and the rotation operator $\t$, which implies that $\cycl{f}$ also
descends to the cyclic quotient. In the context of $\Ainf$-algebras, the cyclization of an $\Ainf$-morphism is the induced
map on the Hochschild and cyclic complexes.

The cyclization constructions are extended from the reduced tensor module
$\tensr{V}$ to the full tensor module $\tens{V}$ in \cref{sec:extension-cycl-full-tensor-module}.
The extension of $\cycl{\mu}$ involves the curvature term $\mu_0(1)$ while the extension of
$\cycl{f}$ involves the \textbf{cyclic exponential} $\cexp{f_0(1)}$ of the change of connection term $f_0(1)$.
We study the properties of the cyclic exponential and show that
our extensions continue to satisfy the functorial properties described earlier up to $\Im \left( \idd - \t \right)$,
which implies that they descend and behave well on the cyclic quotient $\tens{V} / \Im \left( \idd - \t \right)$.

Finally, in \cref{sec:cyclization-naturality}, we establish the naturality of the cyclization operations with respect to the pullback
of generalized coderivations and morphisms, ensuring that our definitions are robust under base change.

Let us fix a field $\mathbbm{k}$ of characteristic zero, endowed with the trivial norm,
and a grading datum $\left( \GG, \braidop \right)$ (see \cref{subsec:grading-data}).
Fix also a graded-commutative Banach $\mathbbm{k}$-algebra $R$.
Starting from this section, we will always work in the graded Banach framework.
This means that our objects will belong to
various categories of graded Banach objects (such as graded Banach $R$-modules,
graded Banach $R$-coalgebras, etc.), and morphisms will be taken in the corresponding categories.
All the notions we use are to be interpreted in the sense discussed in \cref{sec:non-archimedean-graded-setting}.

Since we consistently work in the non-Archimedean graded Banach framework, from this section onward,
we will simplify our notation by omitting the completion mark from various
symbols. For example, we denote by $\otimes = \otimes_R$ the complete tensor product $\cotimes = \cotimes_R$ and by $\oplus$ the complete direct sum $\coplus$. Furthermore, given a graded
Banach $R$-module $V$, we will denote
by $\tens{V} = \tens{V}[R]$ the formal tensor module/coalgebra $\tensf{V}[R]$ as discussed in
\cref{subsec:formal-tensor-coalgebra} and omit the adjective ``formal'' from the text.
This should not cause any confusion as we will keep referring to our objects as Banach objects and
the presence of infinite sums will suggest that we must work with completions for things to make sense.

\subsection{Cyclization of Generalized Coderivations} \label{sec:cyclization-generalized-coderivation}
Let $V$ be a graded Banach $R$-module.
Given a graded $\mathbbm{k}$-algebra derivation
$d \colon R \rightharpoonup R$, recall that a generalized coderivation on $\tens{V}$ over $d$
is a graded bounded $\mathbbm{k}$-linear map $\mu \colon \tens{V} \rightharpoonup \tens{V}$
of the same degree as $d$ which satisfies
\begin{equation*}
	\mu \left( r \cdot l \right) = dr \cdot l + (-1)^{\braidd{r}{\mu}} r \cdot \mu \left( l \right)
\end{equation*}
for all $l \in \tens{V}$ and $r \in R$ and
\begin{equation*}
	\Delta \circ \mu =  \left( \mu \otimes_R \idd + \idd \otimes_R \mu \right) \circ \Delta.
\end{equation*}

Since the tensor coalgebra $\tens{V}$ is coaugmented, one can recover the underlying algebra derivation
$d \colon R \rightharpoonup R$ from the map $\mu \colon \tens{V} \rightharpoonup \tens{V}$
by the formula $d_{\mu} = d = \varepsilon \circ \mu \circ i_R$ where
$\varepsilon \colon \tens{V} \rightarrow R$ is the counit of $\tens{V}$ and $i_R \colon R \rightarrow \tens{V}$
is the coaugmentation (see \cref{item:coderivation-underlying-derivation} of \cref{sec:pre-differential-graded-coalgebras}).
This means that the data of a generalized coderivation is encoded entirely in the map $\mu$.
We henceforth call $\mu$ a generalized coderivation on $\tens{V}$ without specifying the underlying
derivation $d$.

Let us denote by $\CoDer{\tens{V}} = \CoDer{\tens{V}}[R]$ the graded Banach $R$-module consisting of all generalized coderivations
on $\tens{V}$ and by $\Coder{\tens{V}} = \Coder{\tens{V}}[R]$ the graded Banach $R$-submodule of $\CoDer{\tens{V}}$ consisting of
$R$-linear coderivations. The $R$-action on $\CoDer{\tens{V}}$ is given by the outer action
\eqref{eq:left-R-action-on-hom}\footnote{Note that if $\mu$ is a generalized coderivation on $\tens{V}$ over $d$ and $r \in R$,
	then $r \cdot \mu$ is a generalized coderivation on $\tens{V}$ over the algebra derivation $r \cdot d$.}
and the norm on $\CoDer{\tens{V}}$ is given by the operator norm. The graded Banach $R$-module
$\CoDer{\tens{V}}$ has a natural structure of a graded Banach Lie algebra with respect to the graded
commutator, and the graded Banach $R$-submodule $\Coder{\tens{V}}$ forms a Lie ideal of $\CoDer{\tens{V}}$.

By \cref{prop:classification-generalized-coderivations-formal-tensor-coalgebra},
a generalized coderivation is determined uniquely by the underlying derivation $d \colon R \rightharpoonup R$
and the corestriction $\corest{\mu} \colon \tens{V} \rightharpoonup V$ of $\mu$ via the formula
\begin{equation*}
	\begin{aligned}
		\mu(l) & = (-1)^{\braidd{\corest{\mu}}{l_{(1)}}} l_{(1)} \otimes
		\corest{\mu} \left( l_{(2)} \right) \otimes l_{(3)},
		\qquad l \in V^{\otimes k}, k \geq 1,
		\\
		\mu(r) & = dr + \mu_0(r) = dr + (-1)^{\braidd{\mu}{r}} r \cdot \mu_0(1),
		\qquad r \in R.
	\end{aligned}
\end{equation*}
Given a generalized coderivation $\mu \colon \tens{V} \rightharpoonup \tens{V}$, the
\textbf{cyclization} $\cycl{\mu} \colon \tensr{V} \rightharpoonup \tensr{V}$ of $\mu$
is defined by the formula
\begin{equation} \label{def:cyclization-coder-short}
	\begin{aligned}
		\cycl{\mu}(x \otimes l) \defeq{} &
		(-1)^{\braid{\degb{\mu}}{\degb{x}}} x \otimes \mu \left( l \right)
		\\
		                                 & +
		                                   (-1)^{\braid{\degb{l_{(3)}}}{\degb{x} + \degb{l_{(1)}} + \degb{l_{(2)}}}}
		\corest{\mu} \left( l_{(3)} \otimes x \otimes l_{(1)} \right) \otimes
		l_{(2)}
	\end{aligned}
\end{equation}
for $x \in V$ and $l \in \tens{V}$.
Equivalently, written entirely in terms of the corestriction, we have the formula
\begin{equation} \label{def:cyclization-coder}
	\begin{aligned}
		\cycl{\mu}(x \otimes l) \defeq{} &
		(-1)^{\braid{\degb{\corest{\mu}}}{\degb{x} + \degb{l_{(1)}}}}
		x \otimes l_{(1)} \otimes \corest{\mu} \left( l_{(2)} \right) \otimes l_{(3)}
		\\
		                                 & +
		                                   (-1)^{\braid{\degb{l_{(3)}}}{\degb{x} + \degb{l_{(1)}} + \degb{l_{(2)}}}}
		\corest{\mu} \left( l_{(3)} \otimes x \otimes l_{(1)} \right) \otimes
		l_{(2)}.
	\end{aligned}
\end{equation}
The map $\cycl{\mu}$ is a graded bounded $\mathbbm{k}$-linear map, with $\degb{\cycl{\mu}} = \degb{\mu}$
and $\nnorm[{\cycl{\mu}}] \leq \nnorm[\mu]$. Just like $\mu$, the map $\cycl{\mu}$ is a module derivation
on $\tensr{V}$ over $d$.
Compared with the action of the coderivation $\mu$, the action of the cyclization $\cycl{\mu}$
includes terms which are obtained by cyclically rotating elements from the end of the list to
the beginning of the list and applying the corestriction $\corest{\mu}$ to the resulting sublist.
For example, we have
\begin{equation*}
	\begin{aligned}
		\mu \left( v_1 \otimes v_2 \right) ={} &
		\mu_2 \left( v_1 \otimes v_2 \right) +
		\mu_1 \left( v_1 \right) \otimes v_2 +
			                                 (-1)^{\braidd{\mu}{v_1}} v_1 \otimes \mu_1 \left( v_2 \right)
		\\
		                                       & +
		\mu_0(1) \otimes v_1 \otimes v_2 +
			                             (-1)^{\braidd{\mu}{v_1}} v_1 \otimes \mu_0(1) \otimes v_2
		\\
		                                       & +
		                                         (-1)^{\braid{\degb{\mu}}{\degb{v_1} + \degb{v_2}}} v_1 \otimes v_2 \otimes \mu_0(1)
	\end{aligned}
\end{equation*}
while
\begin{equation*}
	\begin{aligned}
		\cycl{\mu} \left( v_1 \otimes v_2 \right) ={} &
		\mu_2 \left( v_1 \otimes v_2 \right) + (-1)^{\braidd{v_1}{v_2}} \mu_2 \left( v_2 \otimes v_1 \right)
		\\
		                                              & + \mu_1 \left( v_1 \right) \otimes v_2 +
			                                                                                   (-1)^{\braidd{\mu}{v_1}} v_1 \otimes \mu_1 \left( v_2 \right)
		\\
		                                              & +
		                                                (-1)^{\braidd{\mu}{v_1}} v_1 \otimes \mu_0(1) \otimes v_2 +
			                                                                                                      (-1)^{\braid{\degb{\mu}}{\degb{v_1} + \degb{v_2}}} v_1 \otimes v_2 \otimes \mu_0(1).
	\end{aligned}
\end{equation*}
Note that while $\mu ( l )$ contains a summand of the form $\mu_0(1) \otimes l$ where $\mu_0(1)$ appears
in the beginning, such a summand is absent in $\cycl{\mu} \left( l \right)$.

\begin{rem}
	We note that the formula \eqref{def:cyclization-coder} for $\cycl{\mu}$ in terms
	of the corestriction $\corest{\mu}$ involves some abuse of notation, which we now explain
	(see also \cref{rem:coderivation-formula-d-operator-abuse}).
	Since $\corest{\mu}$ is not $R$-linear, each of the terms
	\begin{equation*}
		(-1)^{\braid{\degb{l_{(3)}}}{\degb{x} + \degb{l_{(1)}} + \degb{l_{(2)}}}}
		\corest{\mu} \left( l_{(3)} \otimes x \otimes l_{(1)} \right) \otimes l_{(2)},
		\qquad
		(-1)^{\braid{\degb{\mu}}{\degb{x} + \degb{l_{(1)}}}}
		x \otimes l_{(1)} \otimes \corest{\mu} \left( l_{(2)} \right) \otimes l_{(3)}
	\end{equation*}
	appearing in \cref{def:cyclization-coder} is actually ill-defined.
	However, the only component of $\corest{\mu}$ which is not $R$-linear is the component
	$\mu_1$, and if we write the contribution of $\mu_1$ to the \textit{sum} of the two terms
	above, we see that it has the form
	\begin{equation*}
		\mu_1 \left( x \right) \otimes v_1 \otimes \dots \otimes v_k +
		\sum_{i=1}^k (-1)^{\braid{\degb{\mu}}{\degb{x}+\degb{v_1}+\dots+\degb{v_{i-1}}}}
		x \otimes v_1 \otimes \dots \otimes v_{i-1} \otimes \mu_1 \left( v_i \right)
		\otimes v_{i+1} \otimes \dots \otimes v_k
	\end{equation*}
	for $x \in V$ and $l = v_1 \otimes \dots \otimes v_k \in V^{\otimes k}$. This implies that the sum is
	actually well-defined on $V^{\otimes (k+1)}$ even when $\mu$ is not $R$-linear and
	hence \eqref{def:cyclization-coder} makes sense
	not only for $R$-linear coderivations but also for generalized coderivations. We see also
	that if $\mu$ is a generalized coderivation whose only non-zero component is $\mu_1$, then
	$\cycl{\mu} = \mu$ (i.e., when $\mu$ encodes a single operation of arity one, the cyclization
	has no effect).

	More generally, whenever we use formulas \eqref{def:cyclization-coder-short} or
	\eqref{def:cyclization-coder} and perform certain calculations, certain summands
	appearing in the middle of the calculation might be ill-defined, but one can verify
	that the entire sum will always give
	a well-defined operator on $\tensr{V}$.
\end{rem}

\begin{lm} \label{lm:cycl-comm-bracket}
	The process of cyclization of coderivations commutes with the Lie bracket. Given two
	generalized coderivations $\mu, \nu$ on $\tens{V}$, we have the identity
	$\cycl{[\mu, \nu]} = [\cycl{\mu},\cycl{\nu}]$ on $\tensr{V}$. In other words, $\tensr{V}$ is a
	representation of the Lie algebra $\CoDer{\tens{V}}$
	with respect to the action map $\left( \mu, l \right) \mapsto \cycl{\mu} \left( l \right)$.
\end{lm}
\begin{proof}
	The proof is a straightforward calculation where the only difficulty is to keep track of all the
	signs involved. We have
	\begin{equation*}
		\begin{aligned}
			\left( \cycl{\mu} \circ \cycl{\nu} \right) \left( x \otimes l \right)
			\stackrel{\eqref{def:cyclization-coder-short}}{=}{} &
			\underbrace{\cycl{\mu} \left( (-1)^{\braidd{\nu}{x}} x \otimes \nu(l) \right)}_{I}
			\\
			                                                    & +
			\underbrace{\cycl{\mu} \left(
				(-1)^{\braid{\degb{l_{(3)}}}{\degb{x} + \degb{l_{(1)}} + \degb{l_{(2)}}}}
				\corest{\nu} \left( l_{(3)} \otimes x \otimes l_{(1)} \right) \otimes l_{(2)} \right)}_{II}.
		\end{aligned}
	\end{equation*}
	Let us start by expanding the first term $I$. We have:
	\begin{equation*}
		\begin{aligned}
			I \stackrel{\eqref{def:cyclization-coder-short}}{=}{} &
			\underbrace{
				(-1)^{\braid{\degb{\nu} + \degb{\mu}}{\degb{x}}}
				x \otimes \mu \left( \nu \left( l \right) \right)
			}_{(1)}
			\\
			                                                      & +
			\underbrace{
				(-1)^{\braid{\degb{\nu(l)_{(3)}}}{\degb{x} + \degb{\nu(l)_{(1)}} +
						\degb{\nu(l)_{(2)}}} + \braidd{\nu}{x}}
				\corest{\mu} \left( \nu(l)_{(3)} \otimes x \otimes \nu(l)_{(1)} \right)
				\otimes \nu(l)_{(2)}
			}_{(2)}.
		\end{aligned}
	\end{equation*}
	Since $\nu$ is a coderivation, the generalized co-Leibniz rule tells us that
	\begin{equation}
		\begin{aligned}
			\nu(l)_{(1)} \ootimes \nu(l)_{(2)} \ootimes \nu(l)_{(3)} ={} &
			\nu \left( l_{(1)} \right) \ootimes l_{(2)} \ootimes l_{(3)}
			\\
			                                                             & +
			                                                               (-1)^{\braidd{\nu}{l_{(1)}}}
			l_{(1)} \ootimes \nu \left( l_{(2)} \right) \ootimes l_{(3)}
			\\
			                                                             & +
			                                                               (-1)^{\braid{\degb{\nu}}{\degb{l_{(1)}} + \degb{l_{(2)}}}}
			l_{(1)} \ootimes l_{(2)} \ootimes \nu \left( l_{(3)} \right),
		\end{aligned} \label{eq:generalized-co-leibniz-rule}
	\end{equation}
	and so we can rewrite $(2)$ above as
	\begin{equation*}
		\begin{aligned}
			(2) ={} & (-1)^{\braid{\degb{l_{(3)}}}
				          {\degb{x} + \degb{\nu \left( l_{(1)} \right)} + \degb{l_{(2)}}} +
				          \braidd{\nu}{x}}
			\corest{\mu} \left( l_{(3)} \otimes x \otimes \nu \left( l_{(1)} \right) \right)
			\otimes l_{(2)}
			\\
			        & +
			          (-1)^{\braid{\degb{l_{(3)}}}
				          {\degb{x} + \degb{l_{(1)}} + \degb{\nu \left( l_{(2)}\right)}} +
				          \braidd{\nu}{l_{(1)}} + \braidd{\nu}{x}}
			\corest{\mu} \left( l_{(3)} \otimes x \otimes l_{(1)} \right) \otimes
			\nu \left( l_{(2)} \right)
			\\
			        & +
			          (-1)^{\braid{\degb{\nu \left( l_{(3)} \right)}}
				          {\degb{x} + \degb{l_{(1)}} + \degb{l_{(2)}}} +
				          \braid{\degb{\nu}}{\degb{l_{(1)}} + \degb{l_{(2)}}} + \braidd{\nu}{x}}
			\corest{\mu} \left( \nu \left( l_{(3)} \right) \otimes x \otimes l_{(1)} \right)
			\otimes l_{(2)}
			\\
			={}     &
			(-1)^{\braid{\degb{l_{(3)}}}
				{\degb{x} + \degb{l_{(1)}} + \degb{l_{(2)}}} +
				\braid{\degb{\nu}}{\ul{\degb{l_{(3)}} + \degb{x}}}}
			\corest{\mu} \left( \ul{l_{(3)} \otimes x} \otimes \nu \left( l_{(1)} \right) \right)
			\otimes l_{(2)}
			\\
			        & + (-1)^{\braid{\degb{l_{(3)}}}
				          {\degb{x} + \degb{l_{(1)}} + \degb{l_{(2)}}} +
				          \braid{\degb{\nu}}
				          {\ul{\degb{l_{(3)}} + \degb{x} + \degb{l_{(1)}}}}
			          }
			\corest{\mu} \left( \ul{l_{(3)} \otimes x \otimes l_{(1)}} \right)
			\otimes \nu \left( l_{(2)} \right)
			\\
			        & + (-1)^{\braid{\degb{l_{(3)}}}
				          {\degb{x} + \degb{l_{(1)}} + \degb{l_{(2)}}}
			          }
			\corest{\mu} \left( \nu \left( l_{(3)} \right) \otimes x \otimes l_{(1)}
			\right) \otimes l_{(2)}.
		\end{aligned}
	\end{equation*}
	Next we expand the term $II$:
	\begin{align*}
		II \stackrel{\eqref{def:cyclization-coder-short}}{=}{}        &
		(-1)^{\braid{\degb{l_{(3)}}}
			{\degb{x} + \degb{l_{(1)}} + \degb{l_{(2)}}} +
			\braid{\degb{\mu}}
			{\degb{\corest{\nu}(l_{(3)} \otimes x \otimes l_{(1)})}}
		}
		\corest{\nu} \left( l_{(3)} \otimes x \otimes l_{(1)} \right) \otimes \mu \left( l_{(2)} \right)
		\\
		                                                              & +
		                                                                (-1)^{\braid{\degb{l_{(3)}}}
			                                                                {\degb{x} + \degb{l_{(1)}} + \degb{l_{(21)}} + \degb{l_{(22)}} +
				                                                                \degb{l_{(23)}}} +
			                                                                \braid{\degb{l_{(23)}}}
			                                                                {\degb{\corest{\nu}(l_{(3)} \otimes x \otimes l_{(1)})} + \degb{l_{(21)}} +
				                                                                \degb{l_{(22)}}}}
		\\
		                                                              & \qquad\quad
		\corest{\mu} \left(
		l_{(23)} \otimes \corest{\nu} \left( l_{(3)} \otimes x \otimes l_{(1)}
		\right) \otimes l_{(21)} \right) \otimes l_{(22)}
		\\
		\stackrel{\phantom{\eqref{def:cyclization-coder-short}}}{=}{} &
		(-1)^{\braid{\degb{l_{(3)}}}
			{\degb{x} + \degb{l_{(1)}} + \degb{l_{(2)}}} +
			\braid{\degb{\mu}}
			{\ul{\degb{l_{(3)}} + \degb{x} + \degb{l_{(1)}}}} +
			\braidd{\mu}{\nu}
		}
		\corest{\nu} \left( \ul{l_{(3)} \otimes x \otimes l_{(1)}} \right)
		\otimes \mu \left( l_{(2)} \right)
		\\
		                                                              & +(-1)^{\braid{\degb{l_{(23)}} + \degb{l_{(3)}}}
			                                                                {\degb{x} + \degb{l_{(1)}} + \degb{l_{(21)}} + \degb{l_{(22)}}}
			                                                                + \braid{\degb{\nu}}{\ul{\degb{l_{(23)}}}}}
		\\
		                                                              & \qquad\quad
		\corest{\mu} \left( \ul{l_{(23)}} \otimes \corest{\nu}
		\left( l_{(3)} \otimes x \otimes l_{(1)} \right) \otimes l_{(21)} \right)
		\otimes l_{(22)}.
	\end{align*}
	Combining all the expressions and rearranging, we get
	\begin{equation*}
		\begin{aligned}
			\left( \cycl{\mu} \circ \cycl{\nu} \right) \left( x \otimes l \right) ={} &
			(-1)^{\braid{\degb{\nu}+\degb{\mu}}{\degb{x}}}
			x \otimes \mu \left( \nu \left( l \right) \right)
			\\
			                                                                          & \left.
			\begin{aligned}
				 & + (-1)^{\braid{\degb{l_{(3)}}}{\degb{x} + \degb{l_{(1)}} + \degb{l_{(2)}}}}
				\corest{\mu} \left( \nu \left( l_{(3)} \right) \otimes x \otimes l_{(1)} \right)
				\otimes l_{(2)}
				\\
				 & + (-1)^{\braid{\degb{l_{(23)}} + \degb{l_{(3)}}}
					   {\degb{x} + \degb{l_{(1)}} + \degb{l_{(21)}} + \degb{l_{(22)}}} +
					   \braid{\degb{\nu}}{\ul{\degb{l_{(23)}}}}}
				\\
				 & \qquad\quad
				\corest{\mu} \left( \ul{l_{(23)}} \otimes \corest{\nu}(l_{(3)} \otimes x \otimes
				l_{(1)}) \otimes l_{(21)} \right) \otimes l_{(22)}
				\\
				 & +  (-1)^{\braid{\degb{l_{(3)}}}
					   {\degb{x} + \degb{l_{(1)}} + \degb{l_{(2)}}} +
					   \braid{\degb{\nu}}
					   {\ul{\degb{l_{(3)}} + \degb{x}}}}
				\corest{\mu} \left( \ul{l_{(3)} \otimes x} \otimes \nu \left( l_{(1)}
				\right) \right) \otimes l_{(2)}
			\end{aligned}
			\right\} (*)
			\\
			                                                                          & +(-1)^{\braid{\degb{l_{(3)}}}
				                                                                            {\degb{x} + \degb{l_{(1)}} + \degb{l_{(2)}}} +
				                                                                            \braid{\degb{\nu}}
				                                                                            {\ul{\degb{l_{(3)}} + \degb{x} + \degb{l_{(1)}}}}}
			\\
			                                                                          & \qquad\quad
			\corest{\mu} \left( \ul{l_{(3)} \otimes x \otimes l_{(1)}} \right)
			\otimes \nu \left( l_{(2)} \right)
			\\
			                                                                          & +(-1)^{\braid{\degb{l_{(3)}}}{\degb{x} + \degb{l_{(1)}} + \degb{l_{(2)}}} +
				                                                                            \braid{\degb{\mu}}{\ul{\degb{l_{(3)}} + \degb{x} + \degb{l_{(1)}}}} +
				                                                                            \braidd{\mu}{\nu}}
			\\
			                                                                          & \qquad\quad
			\corest{\nu} \left( \ul{l_{(3)} \otimes x \otimes l_{(1)}} \right) \otimes
			\mu \left( l_{(2)} \right).
		\end{aligned}
	\end{equation*}
	Given $l,s \in \tens{V}$ and $x \in V$, we have the identity
	\begin{equation}
		\begin{aligned}
			\nu \left( l \otimes x \otimes s \right) ={} & \nu \left( l \right) \otimes x \otimes s
			\\
			                                             & +
			                                               (-1)^{\braidd{l_{(1)}}{\nu}}
			l_{(1)} \otimes \corest{\nu} \left( l_{(2)} \otimes x \otimes s_{(1)} \right) \otimes s_{(2)}
			\\
			                                             & +
			                                               (-1)^{\braid{\degb{l} + \degb{x}}{\degb{\nu}}}
			l \otimes x \otimes \nu \left( s \right).
		\end{aligned} \label{eq:coder-l-x-s-identity}
	\end{equation}
	Using the identity above, we can combine all three terms of $(*)$ and obtain the formula
	\begin{align}\label{eq:composition-cyc-coderivations}
		\left( \cycl{\mu} \circ \cycl{\nu} \right) \left( x \otimes l \right) ={} &
		(-1)^{\braid{\degb{\nu}+\degb{\mu}}{\degb{x}}}
		x \otimes \mu \left( \nu \left( l \right) \right)
		\nonumber
		\\
		                                                                          & +(-1)^{\braid{\degb{l_{(3)}}}{\degb{x} + \degb{l_{(1)}} + \degb{l_{(2)}}}}
		\corest{\mu} \left( \nu \left( l_{(3)} \otimes x \otimes l_{(1)} \right)
		\right) \otimes l_{(2)}
		\nonumber
		\\
		                                                                          & + (-1)^{\braid{\degb{l_{(3)}}}
			                                                                            {\degb{x} + \degb{l_{(1)}} + \degb{l_{(2)}}} +
			                                                                            \braid{\degb{\nu}}
			                                                                            {\ul{\degb{l_{(3)}} + \degb{x} + \degb{l_{(1)}}}}
		                                                                            }
		\nonumber
		\\
		                                                                          & \qquad\quad
		\corest{\mu} \left( \ul{l_{(3)} \otimes x \otimes l_{(1)}} \right)
		\otimes \nu \left( l_{(2)} \right)
		\\
		                                                                          & +(-1)^{\braid{\degb{l_{(3)}}}
			                                                                            {\degb{x} + \degb{l_{(1)}} + \degb{l_{(2)}}} +
			                                                                            \braid{\degb{\mu}}
			                                                                            {\ul{\degb{l_{(3)}} + \degb{x} + \degb{l_{(1)}}}} +
			                                                                            \braidd{\mu}{\nu}
		                                                                            }
		\nonumber
		\\
		                                                                          & \qquad\quad
		\corest{\nu} \left( \ul{l_{(3)} \otimes x \otimes l_{(1)}} \right) \otimes
		\mu \left( l_{(2)} \right).
		\nonumber
	\end{align}
	Finally, we get that
	\begin{equation*}
		\begin{aligned}
			\left[ \cycl{\mu}, \cycl{\nu} \right] \left( x \otimes l \right) ={} & \left(
			                                                                       \cycl{\mu} \circ \cycl{\nu} - (-1)^{\braidd{\cycl{\mu}}{\cycl{\nu}}} \cycl{\nu} \circ
			\cycl{\mu} \right) \left( x \otimes l \right)
			\\
			={}                                                                  &
			\left(
			\cycl{\mu} \circ \cycl{\nu} - (-1)^{\braidd{\mu}{\nu}} \cycl{\nu} \circ
			\cycl{\mu} \right) \left( x \otimes l \right)
			\\
			={}                                                                  & (-1)^{\braid{\degb{\nu}+\degb{\mu}}{\degb{x}}} \left(
			x \otimes \mu \left( \nu \left( l \right) \right) -
			                                          (-1)^{\braidd{\mu}{\nu}} x \otimes \nu \left( \mu \left( l
			\right) \right) \right)
			\\
			                                                                     & + (-1)^{\braid{\degb{l_{(3)}}}{\degb{x} + \degb{l_{(1)}} + \degb{l_{(2)}}}}
			\corest{\mu} \left( \nu \left( l_{(3)} \otimes x \otimes l_{(1)} \right)
			\right) \otimes l_{(2)}
			\\
			                                                                     & -
			                                                                       (-1)^{\braidd{\mu}{\nu} + \braid{\degb{l_{(3)}}}{\degb{x} + \degb{l_{(1)}} + \degb{l_{(2)}}}}
			\corest{\nu} \left( \mu \left( l_{(3)} \otimes x \otimes l_{(1)} \right)
			\right) \otimes l_{(2)}
			\\
			={}                                                                  & (-1)^{\braid{\degb{[\mu,\nu]}}{\degb{x}}} x \otimes [\mu,\nu](l)
			\\
			                                                                     & +
			                                                                       (-1)^{\braid{\degb{l_{(3)}}}{\degb{x} + \degb{l_{(1)}} + \degb{l_{(2)}}}}
			\corest{[\mu, \nu]} \left( l_{(3)} \otimes x \otimes l_{(1)} \right) \otimes
			l_{(2)}
			\\
			={}                                                                  & \cycl{[\mu,\nu]} \left( x \otimes l \right).
		\end{aligned}
	\end{equation*}
\end{proof}

We will say that a coderivation $\mu$ is \textbf{odd} if
$\braidd{\mu}{\mu} \equiv 1 \mod 2$.
\begin{cor} \label{cor:odd-coder-cyc-differential}
	Let $\mu$ be an odd coderivation on $\tens{V}$ which satisfies $\mu^2 = 0$. Then
	$\cycl{\mu}^2 = 0$ on $\tensr{V}$.
\end{cor}
\begin{proof}
	For an odd coderivation $\mu$, the condition $\mu^2 = 0$, i.e., $\mu$ is a differential, is
	equivalent to the condition $\left[ \mu, \mu \right] = 2 \cdot \mu^2 = 0$.
	Thus, \cref{lm:cycl-comm-bracket} implies that
	\begin{equation*}
		0 = \cycl{[\mu,\mu]} = \left[ \cycl{\mu}, \cycl{\mu} \right] =
		2 \cdot \cycl{\mu}^2
	\end{equation*}
	and hence $\cycl{\mu}$ is a differential on $\tensr{V}$.\footnote{\cref{cor:odd-coder-cyc-differential}
		holds even without assuming that $2 \neq 0$, as is evident from \cref{eq:composition-cyc-coderivations}.}
\end{proof}

\begin{rem}
	Let $\mathcal{R} = \left( R, d \right)$ be a differential graded-commutative Banach $\mathbbm{k}$-algebra
	and let $\mathcal{A} = \left( A, \mu, e \right)$ be a unital Banach $\Ainf$-algebra over $\mathcal{R}$.
	Then $\cycl{\mu}$ is a differential on $\tensr{A}$, called the \textbf{Hochschild differential},
	and the complex $\left( \tensr{A}, \cycl{\mu} \right)$ is called the \textbf{Hochschild complex}
	of $\mathcal{A}$. When $\mathcal{A}$ corresponds to a DGA, the Hochschild complex
	coincides with the standard one up to a shift and an identification
	(see \cref{appendix:sign-conversions}).
\end{rem}

Let $\t \colon \tensr{V} \rightarrow \tensr{V}$ be the $R$-linear map defined by
\begin{equation}
	\t \left( v_1 \otimes \dots \otimes v_n \right) \defeq
	(-1)^{\braid{\degb{v_n}}{\degb{v_1} + \dots + \degb{v_{n-1}}}}
	v_n \otimes v_1 \otimes \dots \otimes v_{n-1}.
	\label{eq:def-t-rotation}
\end{equation}
The map $\t$ is called the \textbf{rotation operator} and is homogeneous of degree zero with respect to both
the $G$-grading and the weight grading on $\tensr{V}$. The map $\t$ acts on each
weight-homogeneous part $V^{\otimes n}$ as a
generator of the cyclic group action $\ZZ / n \ZZ$.
Denote also by $\N \colon \tensr{V} \rightarrow \tensr{V}$ the $R$-linear map defined by
\begin{equation}
	\rest{\N}{V^{\otimes n}} \defeq \sum_{i=0}^{n-1} {\t}^i.
	\label{eq:def-norm-map}
\end{equation}
The map $\N$ is called the \textbf{norm operator} and is homogeneous of degree zero with respect to both
the $G$-grading and the weight grading on $\tensr{V}$.

Let
$\tensrcyc{V} \defeq \tensr{V} / \Im \left( \idd - \t \right)$ be the quotient of $\tensr{V}$ by the
image of $\idd - \t$. We will call $\tensrcyc{V}$ the \textbf{reduced cyclic tensor module} on $V$.
The module $\tensrcyc{V}$ is naturally $G$-graded and weight-graded,
and elements of weight $n$ are precisely equivalence classes
of elements of $V^{\otimes n}$ which are invariant under the action of $\ZZ /n \ZZ$.
We will not use a special notation for the equivalence classes and
continue to write elements of $\tensrcyc{V}$ as usual. Note that
since we work in the characteristic zero context, we have
$\Im \left( \idd - \t \right) = \ker \left( \N \right)$ (see \cref{subsec:cyc-bicomplex})
and hence $\Im \left( \idd - \t \right)$ is closed, and the quotient $\tensrcyc{V}$ is a graded Banach $R$-module.

Any generalized coderivation $\mu \colon \tens{V} \rightharpoonup \tens{V}$
preserves the reduced tensor module $\tensr{V}$ and hence induces a map
$\mu \colon \tensr{V} \rightharpoonup \tensr{V}$. While
the map $\mu \colon \tensr{V} \rightharpoonup \tensr{V}$ does not descend to a well-defined map
$\mu \colon \tensrcyc{V} \rightharpoonup \tensrcyc{V}$ on $\tensrcyc{V} = \tensr{V} / \Im \left( \idd - \t \right)$,
the cyclization $\cycl{\mu} \colon \tensr{V} \rightharpoonup \tensr{V}$ does descend to a
map $\cycl{\mu} \colon \tensrcyc{V} \rightharpoonup \tensrcyc{V}$, as a consequence of the following lemma:

\begin{lm} \label{lm:coder-descends-quotient}
	We have the identity
	\begin{equation} \label{eq:cycl-mu-idd-t}
		\cycl{\mu} \circ (\idd - \t) = (\idd - \t) \circ \mu
	\end{equation}
	on $\tensr{V}$.
\end{lm}
\begin{proof}
	The proof is a lengthy but straightforward calculation.  We will first verify the
	identity for elementary tensors of weight greater than or equal to two. Such
	tensors can be written in our notation as $x \otimes l \otimes z$ where $x,z
		\in V$ and $l = v_1 \otimes \dots \otimes v_k$ for $k \geq 0$ (when $k = 0$ this
	means that $l = 1$ and we are working with $x \otimes 1 \otimes z = x \otimes
		z$). When we consider all possible splittings of $l \otimes z$ into three
	consecutive lists, the element $z$ can belong to either the first, the second,
	or the third list. This gives us the identity
	\begin{equation}
		\begin{aligned}
			(l \otimes z)_{(1)} \ootimes (l \otimes z)_{(2)} \ootimes (l \otimes
			                                                 z)_{(3)} ={} &
			(l \otimes z) \ootimes 1 \ootimes 1 +
			l_{(1)} \ootimes \left( l_{(2)} \otimes z \right) \ootimes 1
			\\
			                                                                     & + l_{(1)} \ootimes l_{(2)} \ootimes \left( l_{(3)} \otimes z \right).
		\end{aligned} \label{eq:lzsplit}
	\end{equation}
	Thus, we have
	\begin{equation*} \begin{aligned}
			\cycl{\mu} \left(x \otimes l \otimes z \right)
			\stackrel{\eqref{def:cyclization-coder-short}}{=}{} &
			(-1)^{\braidd{\mu}{x}} x \otimes \mu \left( l \otimes z \right)
			\\
			                                                    & + (-1)^{\braid{\degb{(l \otimes z)_{(3)}}}
				                                                      {\degb{x} + \degb{(l \otimes z)_{(1)}} + \degb{(l \otimes
						                                                      z)_{(2)}}}}
			\\
			                                                    & \qquad\quad
			\corest{\mu} \left( \left( l \otimes z \right)_{(3)} \otimes x \otimes \left( l
			                                                                       \otimes z \right)_{(1)} \right) \otimes \left( l \otimes z \right)_{(2)}
			\\
			\stackrel{\eqref{eq:lzsplit}}{=}{}                  &
			(-1)^{\braidd{\mu}{x}}
			x \otimes \mu \left( l \otimes z \right) +
			\corest{\mu} \left( x \otimes l \otimes z \right) +
			\corest{\mu} \left( x \otimes l_{(1)} \right) \otimes l_{(2)} \otimes z
			\\
			                                                    & + (-1)^{\braid{\degb{l_{(3)}} + \degb{z}}{\degb{x} + \degb{l_{(1)}} +
					                                                      \degb{l_{(2)}}}} \corest{\mu} \left( l_{(3)} \otimes z \otimes x \otimes l_{(1)}
			\right) \otimes l_{(2)}.
		\end{aligned} \end{equation*}
	Similarly, we have
	\begin{equation} \label{eq:xlsplit}
		\begin{aligned}
			(x \otimes l)_{(1)} \ootimes (x \otimes l)_{(2)} \ootimes (x \otimes l)_{(3)}
			={} &
			(x \otimes l_{(1)}) \ootimes l_{(2)} \ootimes l_{(3)}
			\\
			    & +
			1 \ootimes (x \otimes l_{(1)}) \ootimes l_{(2)}
			\\
			    & +
			1 \ootimes 1 \ootimes (x \otimes l),
		\end{aligned}
	\end{equation}
	and so
	\begin{equation*} \begin{aligned}
			\left( \cycl{\mu} \circ \t \right) \left( x \otimes l \otimes z \right)
			\eqwithref[eq:def-t-rotation]           &
			(-1)^{\braid{\degb{z}}{\degb{x} + \degb{l}}}
			\cycl{\mu} \left( z \otimes x \otimes l \right)
			\\
			\eqwithref[def:cyclization-coder-short] &
			(-1)^{\braid{\degb{z}}{\degb{x} + \degb{l}} + \braidd{\mu}{z}}
			z \otimes \mu \left( x \otimes l \right)
			\\
			                                        & + (-1)^{\braid{\degb{z}}{\degb{x} + \degb{l}} +
				                                          \braid{\degb{(x \otimes l)_{(3)}}}
				                                          {\degb{z} + \degb{(x \otimes l)_{(1)}} +
					                                          \degb{(x \otimes l)_{(2)}}}
			                                          }
			\\
			                                        & \qquad\quad
			\corest{\mu} \left( (x \otimes l)_{(3)} \otimes z \otimes (x \otimes l)_{(1)}
			\right) \otimes (x \otimes l)_{(2)}
			\\
			\eqwithref[eq:xlsplit]                  &
			(-1)^{\braid{\degb{z}}{\degb{\mu} + \degb{x} + \degb{l}}}
			z \otimes \mu \left( x \otimes l \right)
			\\
			                                        & + (-1)^{\braid{\degb{z}}
				                                          {\degb{x} + \degb{l_{(1)}} + \degb{l_{(2)}} + \degb{l_{(3)}}} +
				                                          \braid{\degb{l_{(3)}}}
				                                          {\degb{z} + \degb{x} + \degb{l_{(1)}} + \degb{l_{(2)}}}
			                                          }
			\\
			                                        & \qquad\quad
			\corest{\mu} \left( l_{(3)} \otimes z \otimes x \otimes l_{(1)} \right) \otimes l_{(2)}
			\\
			                                        & + (-1)^{\braid{\degb{z}}{\degb{x} + \degb{l_{(1)}} + \degb{l_{(2)}}} +
				                                          \braid{\degb{l_{(2)}}}
				                                          {\degb{z} + \degb{x} + \degb{l_{(1)}}}
			                                          }
			\corest{\mu} \left( l_{(2)} \otimes z \right) \otimes x \otimes l_{(1)}
			\\
			                                        & + (-1)^{\braid{\degb{z}}{\degb{x} + \degb{l}} +
				                                          \braid{\degb{x} + \degb{l}}{\degb{z}}
			                                          }
			\corest{\mu} \left( x \otimes l \otimes z \right)
			\\
			\eqwithref                              &
			(-1)^{\braid{\degb{z}}{\degb{\mu} + \degb{x} + \degb{l}}}
			z \otimes \mu \left( x \otimes l \right)
			\\
			                                        & + (-1)^{\braid{\degb{l_{(3)}} + \degb{z}}
				                                          {\degb{x} + \degb{l_{(1)}} + \degb{l_{(2)}}}
			                                          }
			\corest{\mu} \left( l_{(3)} \otimes z \otimes x \otimes l_{(1)} \right) \otimes
			l_{(2)}
			\\
			                                        & + (-1)^{\braid{\degb{l_{(2)}} + \degb{z}}
				                                          {\degb{x} + \degb{l_{(1)}}}
			                                          }
			\corest{\mu} \left( l_{(2)} \otimes z \right) \otimes x \otimes l_{(1)}
			\\
			                                        & + \corest{\mu} \left( x \otimes l \otimes z \right).
		\end{aligned} \end{equation*}
	Subtracting and cancelling the identical terms, we get
	\begin{equation*} \begin{aligned}
			\left( \cycl{\mu} \circ (\idd - \t) \right) \left( x \otimes l \otimes z
			\right) ={} &
			(-1)^{\braidd{\mu}{x}} x \otimes \mu \left( l \otimes z \right)
			+ \corest{\mu} \left( x \otimes l_{(1)} \right) \otimes l_{(2)} \otimes z
			\\
			            & -
			              (-1)^{\braid{\degb{z}}{\degb{\mu} + \degb{x} + \degb{l}}}
			z \otimes \mu \left( x \otimes l \right)
			\\
			            & -
			              (-1)^{\braid{\degb{l_{(2)}} + \degb{z}} {\degb{x} + \degb{l_{(1)}}}}
			\corest{\mu} \left( l_{(2)} \otimes z \right) \otimes x \otimes l_{(1)}.
		\end{aligned} \end{equation*}
	Next, we have
	\begin{equation*}
		\begin{aligned}
			\mu \left( x \otimes l \otimes z \right) ={} & \mu \left( x
			\otimes l \right) \otimes z +
				                          (-1)^{\braid{\degb{\mu}}{\degb{x} + \degb{l_{(1)}}}} x \otimes l_{(1)}
			\otimes \corest{\mu} \left( l_{(2)} \otimes z \right)                                                                   \\
			                                             & + \corest{\mu} \left( x \otimes l \otimes z \right) +
			                                                              (-1)^{\braid{\degb{\mu}}{\degb{x} + \degb{l} + \degb{z}}}
			x \otimes l \otimes z \otimes \mu_{0}(1)
			\\
			={}                                          & (-1)^{\braidd{\mu}{x}} x \otimes \mu \left( l \otimes z \right) +
			\corest{\mu} \left( x \otimes l_{(1)} \right) \otimes l_{(2)} \otimes z                                                 \\
			                                             & + \corest{\mu} \left( x \otimes l \otimes z \right) +
			\mu_{0}(1) \otimes x \otimes l \otimes z.
		\end{aligned}
	\end{equation*}
	Hence,
	\begin{equation*}
		\begin{aligned}
			\left( \t \circ \mu \right) \left( x \otimes l
			\otimes z \right) ={} &
			(-1)^{\braidd{z}{\mu \left( x \otimes l \right)}}
			z \otimes \mu \left( x \otimes l \right)                                               \\
			                      & + (-1)^{\braid{\degb{\mu}}{\degb{x} + \degb{l_{(1)}}} +
				                        \braid{\degb{\mu \left( l_{(2)} \otimes z \right)}}{\degb{x} +
					                        \degb{l_{(1)}}}
			                        }
			\corest{\mu} \left( l_{(2)} \otimes z \right) \otimes x \otimes l_{(1)}
			\\
			                      & + \corest{\mu} \left( x \otimes l \otimes z \right)
			\\
			                      & +
			                        (-1)^{\braid{\degb{\mu}}{\degb{x} + \degb{l} + \degb{z}} +
				                        \braid{\degb{\mu_0(1)}}{\degb{x} + \degb{l} + \degb{z}}
			                        }
			\mu_{0}(1) \otimes x \otimes l \otimes z
			\\
			={}                   & (-1)^{\braid{\degb{z}}{\degb{\mu} + \degb{x} + \degb{l}}}
			z \otimes \mu \left( x \otimes l \right)
			\\
			                      & + (-1)^{\braid{\degb{l_{(2)}} + \degb{z}}
				                        {\degb{x} + \degb{l_{(1)}}}
			                        }
			\corest{\mu} \left( l_{(2)} \otimes z \right) \otimes x \otimes l_{(1)}
			\\
			                      & + \corest{\mu} \left( x \otimes l \otimes z \right)
			\\
			                      & + \mu_{0}(1) \otimes x \otimes l \otimes z
		\end{aligned}
	\end{equation*}
	and so we also have
	\begin{equation*}
		\begin{aligned}
			\left( \left( \idd - \t \right) \circ \mu \right) \left( x \otimes l
			\otimes z \right) ={} &
			(-1)^{\braidd{\mu}{x}} x \otimes \mu \left( l \otimes z \right) +
			\corest{\mu} \left( x \otimes l_{(1)} \right) \otimes l_{(2)} \otimes z
			\\
			                      & -
			                        (-1)^{\braid{\degb{z}}{\degb{\mu} + \degb{x} + \degb{l}}}
			z \otimes \mu \left( x \otimes l \right)
			\\
			                      & -
			                        (-1)^{\braid{\degb{l_{(2)}} + \degb{z}} {\degb{x} + \degb{l_{(1)}}}}
			\corest{\mu} \left( l_{(2)} \otimes z \right) \otimes x \otimes l_{(1)}.
		\end{aligned}
	\end{equation*}
	This shows the identity for elementary tensors of weight greater than or equal
	to two. For elements of weight one, we note that if $x \in V$ then $(\idd - \t)(x) = 0$,
	and so $\left( \cycl{\mu} \circ (\idd - \t) \right)(x) = 0$. Meanwhile,
	\begin{align*}
		\mu(x) = \mu_1(x) + \mu_0(1) \otimes x + (-1)^{\braidd{\mu}{x}} x \otimes \mu_0(1)
	\end{align*}
	which implies that
	\begin{align*}
		\left( \t \circ \mu \right)(x) & = \mu_1(x) +
		                                   (-1)^{\braidd{x}{\mu_0(1)}} x \otimes \mu_0(1) +
		                                                                         (-1)^{\braidd{\mu}{x} + \braidd{\mu_0(1)}{x}}
		\mu_0(1) \otimes x                                                                                                     \\
		                               & = \mu_1(x) + (-1)^{\braidd{x}{\mu}} x \otimes \mu_0(1) + \mu_0(1) \otimes x.
	\end{align*}
	Subtracting both expressions, we get that
	\begin{equation*}
		\left( \cycl{\mu} \circ (\idd - \t) \right)(x) = \left( \left( \idd - \t \right)
		\circ \mu \right)(x) = 0.
	\end{equation*}
\end{proof}

\begin{lm} \label{lm:coder-cyc-norm-map}
	We have the identity
	\begin{equation} \label{eq:mu-N-N-cycl-mu}
		\mu \circ \N = \N \circ \cycl{\mu}
	\end{equation}
	on $\tensr{V}$.
\end{lm}
\begin{proof}
	The proof is similar to the proof of \cref{lm:coder-descends-quotient} and is left to the reader.
\end{proof}

\begin{rem}
	The identities $\cycl{\mu} (\idd - \t) = (\idd - \t) \mu$ and
	$\mu \N = \N \cycl{\mu}$ are generalizations of the identities
	$b \left( \idd - \t \right) = \left( \idd - \t \right) b'$ and $b' \N = \N b$
	for the Hochschild differential $b$ and the bar differential $b'$ of
	an associative algebra (see \cite[Lemma 2.1.1]{Loday1998}).
\end{rem}

\subsection{Cyclization of Tensor Coalgebra Morphisms} \label{sec:cyclization-coalgebra-morphisms}
Let $V$ and $W$ be graded Banach $R$-modules and let
$f \colon \tens{V} \rightarrow \tens{W}$ be a morphism of graded Banach $R$-coalgebras.
By \cref{prop:classification-morphisms-tensor-coalgebra}, the morphism
$f$ is determined uniquely from the corestriction $\corest{f} \colon \tens{V} \rightarrow W$ by the formulas
\begin{align*}
	f(l) & = \sum_{k = 1}^{\infty} \corest{f}(l_{(1)}) \otimes \dots \otimes \corest{f}(l_{(k)}),
	\qquad \left( l \in \tensr{V} \right)
	\\
	f(1) & = \sum_{k=0}^{\infty} f_{0}(1)^{\otimes k} = \Exp{f_0(1)}.
\end{align*}
The \textbf{cyclization} $\cycl{f} \colon \tensr{V} \rightarrow \tensr{W}$ of $f$ is defined by the formula
\begin{equation} \label{def:cyclization-morphism}
	\cycl{f} \left( x \otimes l \right) \defeq
	(-1)^{\braid{\degb{l_{(3)}}}{\degb{x} + \degb{l_{(1)}} + \degb{l_{(2)}}}}
	\corest{f} \left( l_{(3)} \otimes x \otimes l_{(1)} \right) \otimes f \left( l_{(2)} \right)
\end{equation}
for $x \in V$ and $l \in \tens{V}$.
The map $\cycl{f}$ is a graded bounded $R$-linear map of degree zero with $\nnorm[{\cycl{f}}] \leq \nnorm[f]$.
Compared with the action of the morphism $f$, the action of the cyclization $\cycl{f}$
includes terms which are obtained by cyclically rotating elements from the end of the list to
the beginning of the list and applying the corestriction $\corest{f}$ to the resulting sublist.
For example, we have
\begin{equation*}
	f \left( v_1 \otimes v_2 \right) =
	\Exp{f_0(1)} \otimes f_1 \left( v_1 \right) \otimes \Exp{f_0(1)} \otimes f_1 \left( v_2 \right) \otimes \Exp{f_0(1)} +
	\Exp{f_0(1)} \otimes f_2 \left( v_1 \otimes v_2 \right) \otimes \Exp{f_0(1)}
\end{equation*}
while
\begin{equation*}
	\begin{aligned}
		\cycl{f} \left( v_1 \otimes v_2 \right) ={} &
		f_1 \left( v_1 \right) \otimes \Exp{f_0(1)} \otimes f_1 \left( v_2 \right) \otimes \Exp{f_0(1)} +
		f_2 \left( v_1 \otimes v_2 \right) \otimes \Exp{f_0(1)}
		\\
		                                            & +
		                                              (-1)^{\braidd{v_2}{v_1}} f_2 \left( v_2 \otimes v_1 \right) \otimes \Exp{f_0(1)}.
	\end{aligned}
\end{equation*}
Note that while $f \left( x \otimes l \right)$ contains summands of the form $\Exp{f_0(1)} \otimes s$
where $\Exp{f_0(1)}$ appears at the beginning, such summands are absent from $\cycl{f} \left( x \otimes l \right)$.
Note also that when $f$ has a single non-zero component $f_1$, we have
$\cycl{f} = f$, i.e., the cyclization has no effect. In particular, if $V = W$ and $f = \idd$ then $\cycl{f} = \idd$.

\begin{lm} \label{lm:func-cycl-morphism}
	Let $U,V,W$ be graded Banach $R$-modules and let $f \colon \tens{V} \rightarrow \tens{W}$
	and $g \colon \tens{U} \rightarrow \tens{V}$ be two morphisms of graded Banach $R$-coalgebras.
	The process of cyclization of morphisms commutes with composition of morphisms, i.e., we have
	\begin{equation} \label{eq:func-cycl-morphism}
		\cycl{\left( f \circ g \right)} = \cycl{f} \circ \cycl{g}.
	\end{equation}
	In particular, $\tensr{V}$ is a representation of the monoid
	of Banach $R$-coalgebra endomorphisms of $\tens{V}$ with respect to the action map
	$\left( f, l \right) \mapsto \cycl{f} \left( l \right)$.
\end{lm}
\begin{proof}
	Since $g$ is a coalgebra morphism, we have $g^{\otimes 3} \circ \Delta^3_{\tens{U}} = \Delta^3_{\tens{V}} \circ g$,
	which, written in our notation, reads
	\begin{equation} \label{eq:g-commutes-with-delta-3}
		g \left( l_{(1)} \right) \ootimes g \left( l_{(2)} \right) \ootimes g \left( l_{(3)} \right) =
		g \left( l \right)_{(1)} \ootimes g \left( l \right)_{(2)} \ootimes g \left( l \right)_{(3)}
	\end{equation}
	(compare to \cref{eq:generalized-co-leibniz-rule}).
	Furthermore, using the explicit formula for $g$ one can verify that $g$ satisfies the identity
	\begin{equation} \label{eq:cohomo-l-x-s-identity}
		g \left( l \otimes x \otimes s \right) = g \left( l_{(1)} \right) \otimes
		\corest{g} \left( l_{(2)} \otimes x \otimes s_{(1)} \right) \otimes g \left( s_{(2)} \right)
	\end{equation}
	for $x \in U$ and $l, s \in \tens{U}$ (compare to \cref{eq:coder-l-x-s-identity}).
	Hence, we have
	\begin{equation*}
		\begin{aligned}
			\left( \cycl{f} \circ \cycl{g} \right) \left( x \otimes l \right)
			\stackrel{\eqref{def:cyclization-morphism}}{=}{}             &
			(-1)^{\braid{\degb{l_{(3)}}}{\degb{x} + \degb{l_{(1)}} + \degb{l_{(2)}}}}
			\cycl{f} \left(
			\corest{g} \left( l_{(3)} \otimes x \otimes l_{(1)} \right) \otimes
			g \left( l_{(2)} \right)
			\right)
			\\
			\stackrel{\eqref{def:cyclization-morphism}}{=}{}             &
			(-1)^{\braid{\degb{l_{(3)}}}{\degb{x} + \degb{l_{(1)}} + \degb{l_{(2)}}} +
				\braid{\degb{g \left( l_{(2)} \right)_{(3)}}}{
					\degb{\corest{g} \left( l_{(3)} \otimes x \otimes l_{(1)} \right)} +
					\degb{g \left( l_{(2)} \right)_{(1)}} + \degb{g \left( l_{(2)} \right)_{(2)}}}}
			\\
			                                                             & \quad
			\corest{f} \left( g \left( l_{(2)} \right)_{(3)} \otimes
			\corest{g} \left( l_{(3)} \otimes x \otimes l_{(1)} \right) \otimes
			g \left( l_{(2)} \right)_{(1)} \right) \otimes
			f \left( g \left( l_{(2)} \right)_{(2)} \right)
			\\
			\stackrel{\eqref{eq:g-commutes-with-delta-3}}{=}{}           &
			(-1)^{\braid{\degb{l_{(3)}}}{\degb{x} + \degb{l_{(1)}} + \degb{l_{(21)}} +
					\degb{l_{(22)}} + \degb{l_{(23)}}} +
				\braid{\degb{l_{(23)}}}{
					\degb{l_{(3)}} + \degb{x} + \degb{l_{(1)}} + \degb{l_{(21)}} + \degb{l_{(22)}}}}
			\\
			                                                             & \quad
			\corest{f} \left( g \left( l_{(23)} \right) \otimes
			\corest{g} \left( l_{(3)} \otimes x \otimes l_{(1)} \right) \otimes
			g \left( l_{(21)} \right) \right) \otimes
			f \left( g \left( l_{(22)} \right) \right)
			\\
			\stackrel{\phantom{\eqref{eq:g-commutes-with-delta-3}}}{=}{} &
			(-1)^{\braid{\degb{l_{(23)}} + \degb{l_{(3)}}}
				{\degb{x} + \degb{l_{(1)}} + \degb{l_{(21)}} + \degb{l_{(22)}}}}
			\\
			                                                             & \quad
			\corest{f} \left( g \left( l_{(23)} \right) \otimes
			\corest{g} \left( l_{(3)} \otimes x \otimes l_{(1)} \right) \otimes
			g \left( l_{(21)} \right) \right) \otimes
			f \left( g \left( l_{(22)} \right) \right)
		\end{aligned}
	\end{equation*}
	while
	\begin{equation*}
		\begin{aligned}
			\cycl{ \left( f \circ g \right)} \left( x \otimes l \right)
			\stackrel{\eqref{def:cyclization-morphism}}{=}{}           &
			(-1)^{\braid{\degb{l_{(3)}}}{\degb{x} + \degb{l_{(1)}} + \degb{l_{(2)}}}}
			\corest{{\left( f \circ g \right)}} \left( l_{(3)} \otimes x \otimes l_{(1)} \right) \otimes
			\left( f \circ g \right) \left( l_{(2)} \right)
			\\
			\stackrel{\phantom{\eqref{def:cyclization-morphism}}}{=}{} &
			(-1)^{\braid{\degb{l_{(3)}}}{\degb{x} + \degb{l_{(1)}} + \degb{l_{(2)}}}}
			\corest{f} \left( g \left( l_{(3)} \otimes x \otimes l_{(1)} \right) \right) \otimes
			f \left( g \left( l_{(2)} \right) \right)
			\\
			\stackrel{\eqref{eq:cohomo-l-x-s-identity}}{=}{}           &
			(-1)^{\braid{\degb{l_{(31)}} + \degb{l_{(32)}}}
				{\degb{x} + \degb{l_{(11)}} + \degb{l_{(12)}} + \degb{l_{(2)}}}}
			\\
			                                                           & \quad
			\corest{f} \left(
			g \left( l_{(31)} \right) \otimes
			\corest{g} \left( l_{(32)} \otimes x \otimes l_{(11)} \right) \otimes
			g \left( l_{(12)} \right)
			\right)
			\otimes
			f \left( g \left( l_{(2)} \right) \right).
		\end{aligned}
	\end{equation*}
	The coassociativity of the coproduct then implies that
	$\cycl{\left( f \circ g \right)} = \cycl{f} \circ \cycl{g}$.
\end{proof}

\begin{lm} \label{lm:func-cycl-morphism-coderivation}
	Let $f \colon \tens{V} \rightarrow \tens{W}$ be a morphism of graded Banach $R$-coalgebras.
	Let $\mu$ be a generalized coderivation on $\tens{V}$ and $\nu$ be a generalized
	coderivation on $\tens{W}$. Assume that we have $f \circ \mu =
		\nu \circ f$. Then $\cycl{f} \circ \cycl{\mu} = \cycl{\nu} \circ
		\cycl{f}$ on $\tensr{V}$.
\end{lm}
\begin{proof}
	Let $x \in V$ and $l \in \tens{V}$. We have
	\begin{equation*}
		\begin{aligned}
			\left( \cycl{f} \circ \cycl{\mu} \right) \left( x \otimes l \right)
			\stackrel{\eqref{def:cyclization-coder-short}}{=}{} &
			\underbrace{
				(-1)^{\braid{\degb{\mu}}{\degb{x}}}
				\cycl{f} \left( x \otimes \mu \left( l \right) \right)
			}_{(1)}
			\\
			                                                    & +
			\underbrace{
				(-1)^{\braid{\degb{l_{(3)}}}{\degb{x} + \degb{l_{(1)}} + \degb{l_{(2)}}}}
				\cycl{f} \left( \corest{\mu} \left( l_{(3)} \otimes x \otimes l_{(1)} \right) \otimes l_{(2)} \right)
			}_{(2)}.
		\end{aligned}
	\end{equation*}
	Using the generalized co-Leibniz rule (\cref{eq:generalized-co-leibniz-rule}), we can write $(1)$ as
	\begin{align}
		(1)
		\stackrel{\eqref{def:cyclization-morphism}}{=}{}                 &
		(-1)^{\braidd{\mu}{x} +
			\braid{\degb{\mu \left( l \right)_{(3)}}}
			{\degb{x} + \degb{\mu \left( l \right)_{(1)}} + \degb{\mu \left( l \right)_{(2)}}}
		}
		\corest{f} \left( \mu \left( l \right)_{(3)} \otimes x \otimes \mu \left( l \right)_{(1)} \right)
		\otimes f \left( \mu \left( l \right)_{(2)} \right)
		\nonumber
		\\
		\stackrel{\eqref{eq:generalized-co-leibniz-rule}}{=}{}           &
		(-1)^{\braidd{\mu}{x} +
			\braid{\degb{l_{(3)}}}
			{\degb{x} + \degb{\mu \left( l_{(1)} \right)} + \degb{ l_{(2)}}}}
		\corest{f} \left( l_{(3)} \otimes x \otimes \mu \left( l_{(1)} \right) \right)
		\otimes f \left( l_{(2)} \right)
		\nonumber
		\\
		                                                                 & +
		                                                                   (-1)^{\braid{\degb{l_{(3)}}}
			                                                                   {\degb{x} + \degb{l_{(1)}} + \degb{\mu \left( l_{(2)} \right)}} +
			                                                                   \braid{\degb{\mu}}{\degb{x} + \degb{l_{(1)}}}
		                                                                   }
		\corest{f} \left( l_{(3)} \otimes x \otimes l_{(1)} \right)
		\otimes f \left( \mu \left( l_{(2)} \right)\right)
		\nonumber
		\\
		                                                                 & +
		                                                                   (-1)^{\braid{\degb{\mu \left( l_{(3)} \right)}}
			                                                                   {\degb{x} + \degb{l_{(1)}} + \degb{l_{(2)}}} +
			                                                                   \braid{\degb{\mu}}{\degb{x} + \degb{l_{(1)}} + \degb{l_{(2)}}}
		                                                                   }
		\corest{f} \left( \mu \left( l_{(3)} \right) \otimes x \otimes l_{(1)} \right)
		\otimes f \left( l_{(2)} \right)
		\nonumber
		\\
		\stackrel{\phantom{\eqref{eq:generalized-co-leibniz-rule}}}{=}{} &
		(-1)^{\braid{\degb{l_{(3)}}}{\degb{x} + \degb{\mu \left( l_{(1)} \right)} + \degb{ l_{(2)}}} +
			\braid{\degb{\mu}}{\degb{x} + \degb{l_{(3)}}}
		}
		\corest{f} \left( l_{(3)} \otimes x \otimes \mu \left( l_{(1)} \right) \right)
		\otimes f \left( l_{(2)} \right)
		\label{eq:mu-in-first-list}
		\\
		                                                                 & +
		                                                                   (-1)^{\braid{\degb{\mu}}{\degb{l_{(3)}} + \degb{x} + \degb{l_{(1)}}} +
			                                                                   \braid{\degb{l_{(3)}}}{\degb{x} + \degb{l_{(1)}} + \degb{l_{(2)}}}
		                                                                   }
		\corest{f} \left( l_{(3)} \otimes x \otimes l_{(1)} \right)
		\otimes f \left( \mu \left( l_{(2)} \right)\right)
		\nonumber
		\\
		                                                                 & +
		                                                                   (-1)^{\braid{\degb{l_{(3)}}}{\degb{x} + \degb{l_{(1)}} + \degb{l_{(2)}}}}
		\corest{f} \left( \mu \left( l_{(3)} \right) \otimes x \otimes l_{(1)} \right)
		\otimes f \left( l_{(2)} \right).
		\label{eq:mu-in-third-list}
	\end{align}
	The term $(2)$ is given by
	\begin{equation} \begin{aligned}
			(2) ={} &
			(-1)^{\braid{\degb{l_{(23)}}}
				{\degb{\corest{\mu}( l_{(3)} \otimes x \otimes l_{(1)} )} + \degb{l_{(21)}} + \degb{l_{(22)}}} +
				\braid{\degb{l_{(3)}}}{\degb{x} + \degb{l_{(1)}} + \degb{l_{(21)}} + \degb{l_{(22)}} + \degb{l_{(23)}}}
			}
			\\
			        & \qquad\quad
			\corest{f} \left( l_{(23)} \otimes \corest{\mu} \left( l_{(3)} \otimes x \otimes l_{(1)}
			\right) \otimes l_{(21)} \right) \otimes f \left( l_{(22)} \right)
			\\
			={}     &
			(-1)^{\braid{\degb{l_{(31)}}}
				{\degb{\mu} + \degb{l_{(32)}} + \degb{x} + \degb{l_{(11)}} + \degb{l_{(12)}} + \degb{l_{(2)}}} +
				\braid{\degb{l_{(32)}}}{\degb{x} + \degb{l_{(11)}} + \degb{l_{(12)}} + \degb{l_{(2)}} + \degb{l_{(31)}}}
			}
			\\
			        & \qquad\quad
			\corest{f} \left( l_{(31)} \otimes \corest{\mu} \left( l_{(32)} \otimes x \otimes l_{(11)}
			\right) \otimes l_{(12)} \right) \otimes f \left( l_{(2)} \right)
			\\
			={}     &
			(-1)^{\braid{\degb{l_{(31)}} + \degb{l_{(32)}}}
				{\degb{x} + \degb{l_{(11)}} + \degb{l_{(12)}} + \degb{l_{(2)}}} +
				\braidd{\mu}{l_{(31)}}
			}
			\\
			        & \qquad\quad
			\corest{f} \left( l_{(31)} \otimes \corest{\mu} \left( l_{(32)} \otimes x \otimes l_{(11)}
			\right) \otimes l_{(12)} \right) \otimes f \left( l_{(2)} \right).
			\label{eq:mu-in-f}
		\end{aligned}  \end{equation}
	Hence, we get the identity
	\begin{align}
		 & \left( \cycl{f} \circ \cycl{\mu} \right) \left( x \otimes l \right) = (1) + (2)
		\stackrel{\eqref{eq:coder-l-x-s-identity}}{=}{}
		\nonumber
		\\
		 & \qquad
		   (-1)^{\braid{\degb{\mu}}{\degb{l_{(3)}} + \degb{x} + \degb{l_{(1)}}} +
			   \braid{\degb{l_{(3)}}}{\degb{x} + \degb{l_{(1)}} + \degb{l_{(2)}}}}
		\corest{f} \left( l_{(3)} \otimes x \otimes l_{(1)} \right) \otimes \left( f
		\circ \mu \right) \left( l_{(2)} \right)
		\nonumber
		\\
		 & \qquad+
		   (-1)^{\braid{\degb{l_{(3)}}}{\degb{x} + \degb{l_{(1)}} + \degb{l_{(2)}}}}
		\left( \corest{f} \circ \mu \right) \left( l_{(3)} \otimes x \otimes l_{(1)} \right)
		\otimes f \left( l_{(2)} \right) \label{eq:cycl-f-cycl-mu-second}
	\end{align}
	where we used \cref{eq:coder-l-x-s-identity} to combine the terms \eqref{eq:mu-in-third-list},
	\eqref{eq:mu-in-f} and \eqref{eq:mu-in-first-list} and turn them into \eqref{eq:cycl-f-cycl-mu-second}.

	We turn to compute $\cycl{\nu} \circ \cycl{f}$. We have
	\begin{equation*}
		\begin{aligned}
			\left( \cycl{\nu} \circ \cycl{f} \right) \left( x \otimes l \right)
			\stackrel{\eqref{def:cyclization-morphism}}{=}{}    &
			(-1)^{\braid{\degb{l_{(3)}}}{\degb{x} + \degb{l_{(1)}} + \degb{l_{(2)}}}}
			\cycl{\nu} \left( \corest{f} \left( l_{(3)} \otimes x \otimes l_{(1)} \right)
			\otimes f \left( l_{(2)} \right) \right)
			\\
			\stackrel{\eqref{def:cyclization-coder-short}}{=}{} &
			(-1)^{\braid{\degb{l_{(3)}}}
				{\degb{x} + \degb{l_{(1)}} + \degb{l_{(2)}}} +
				\braidd{\nu}{\corest{f}( l_{(3)} \otimes x \otimes l_{(1)} )}}
			\\
			                                                    & \qquad\quad
			\corest{f} \left( l_{(3)} \otimes
			x \otimes l_{(1)} \right) \otimes \left( \nu \circ f \right)
			\left( l_{(2)} \right)
			\\
			                                                    & + (-1)^{\braid{\degb{l_{(3)}}}
				                                                      {\degb{x} + \degb{l_{(1)}} + \degb{l_{(2)}}} +
				                                                      \braid{\degb{f \left( l_{(2)} \right)_{(3)}}}
				                                                      {\degb{\corest{f} \left( l_{(3)} \otimes x \otimes l_{(1)} \right)} +
					                                                      \degb{f \left( l_{(2)} \right)_{(1)}} + \degb{f \left( l_{(2)} \right)_{(2)}}}}
			\\
			                                                    & \qquad\quad
			\corest{\nu} \left( f \left( l_{(2)} \right)_{(3)} \otimes
			\corest{f} \left( l_{(3)} \otimes x \otimes l_{(1)} \right) \otimes
			f \left( l_{(2)} \right)_{(1)} \right) \otimes f \left( l_{(2)} \right)_{(2)}
			\\
			\stackrel{\eqref{eq:g-commutes-with-delta-3}}{=}{}  &
			(-1)^{\braid{\degb{l_{(3)}}}
				{\degb{x} + \degb{l_{(1)}} + \degb{l_{(2)}}} +
				\braid{\degb{\nu}}{\degb{l_{(3)}} + \degb{x} + \degb{l_{(1)}}}}
			\\
			                                                    & \qquad\quad
			\corest{f} \left( l_{(3)} \otimes x \otimes l_{(1)} \right) \otimes \left(
			\nu \circ f \right) \left( l_{(2)} \right)
			\\
			                                                    & +
			                                                      (-1)^{\braid{\degb{l_{(3)}}}
				                                                      {\degb{x} + \degb{l_{(1)}} + \degb{l_{(21)}} + \degb{l_{(22)}} + \degb{l_{(23)}}} +
				                                                      \braid{\degb{l_{(23)}}}
				                                                      {\degb{l_{(3)}} + \degb{x} + \degb{l_{(1)}} +
					                                                      \degb{l_{(21)}} + \degb{l_{(22)}}}}
			\\
			                                                    & \qquad\quad
			\corest{\nu} \left( f \left( l_{(23)} \right) \otimes
			\corest{f} \left( l_{(3)} \otimes x \otimes l_{(1)} \right) \otimes
			f \left( l_{(21)} \right) \right) \otimes f \left( l_{(22)} \right)
			\\
			\stackrel{\,\,\textrm{coasso}\,\,}{=}{}             &
			(-1)^{\braid{\degb{l_{(3)}}}
				{\degb{x} + \degb{l_{(1)}} + \degb{l_{(2)}}} +
				\braid{\degb{\nu}}{\degb{l_{(3)}} + \degb{x} + \degb{l_{(1)}}}}
			\\
			                                                    & \qquad\quad
			\corest{f} \left( l_{(3)} \otimes x \otimes l_{(1)} \right) \otimes \left(
			\nu \circ f \right) \left( l_{(2)} \right)
			\\
			                                                    & +
			                                                      (-1)^{\braid{\degb{l_{(32)}}}
				                                                      {\degb{x} + \degb{l_{(11)}} + \degb{l_{(12)}} + \degb{l_{(2)}} + \degb{l_{(31)}}} +
				                                                      \braid{\degb{l_{(31)}}}
				                                                      {\degb{l_{(32)}} + \degb{x} + \degb{l_{(11)}} +
					                                                      \degb{l_{(12)}} + \degb{l_{(2)}}}}
			\\
			                                                    & \qquad\quad
			\corest{\nu} \left( f \left( l_{(31)} \right) \otimes
			\corest{f} \left( l_{(32)} \otimes x \otimes l_{(11)} \right) \otimes
			f \left( l_{(12)} \right) \right) \otimes f \left( l_{(2)} \right)
			\\
			\stackrel{\eqref{eq:cohomo-l-x-s-identity}}{=}{}    &
			(-1)^{\braid{\degb{l_{(3)}}}
				{\degb{x} + \degb{l_{(1)}} + \degb{l_{(2)}}} +
				\braid{\degb{\nu}}{\degb{l_{(3)}} + \degb{x} + \degb{l_{(1)}}}}
			\\
			                                                    & \qquad\quad
			\corest{f} \left( l_{(3)} \otimes x \otimes l_{(1)} \right) \otimes \left(
			\nu \circ f \right) \left( l_{(2)} \right)
			\\
			                                                    & + (-1)^{\braid{\degb{l_{(3)}}}{\degb{x} + \degb{l_{(1)}} +
					                                                      \degb{l_{(2)}}}}
			\\
			                                                    & \qquad\quad
			\left( \corest{\nu} \circ f \right) \left( l_{(3)} \otimes x \otimes
			l_{(1)} \right) \otimes f \left( l_{(2)} \right).
		\end{aligned}
	\end{equation*}
	The condition $f \circ \mu = \nu \circ f$ implies that (and
	in fact is equivalent to) $\corest{f} \circ \mu = \corest{\nu} \circ f$ and hence we
	obtain the lemma.
\end{proof}

Since our morphisms are counital, any Banach $R$-coalgebra morphism $f \colon \tens{V} \rightarrow \tens{W}$
preserves the reduced tensor modules and induces a map $f \colon \tensr{V} \rightarrow \tensr{W}$. While
the map $f \colon \tensr{V} \rightarrow \tensr{W}$ does not descend to a well-defined map
$f \colon \tensrcyc{V} \rightarrow \tensrcyc{W}$ between the reduced cyclic tensor modules,
we will see that the cyclization $\cycl{f} \colon \tensr{V} \rightarrow \tensr{W}$ does induce a well-defined
map between the reduced cyclic tensor modules.
In order to do that, we will need the following auxiliary lemma:

\begin{lm}
	Assume that $b \in V^{0}$ is topologically nilpotent, $u,v \in V$ and $l \in \tens{V}$. Then we have the
	following identities:
	\begin{align}
		\left( \idd - \t \right) \left( \Exp{b} \otimes u \otimes \Exp{b} \right) ={}                     & 0, \label{eq:tebueb}
		\\
		\left( \idd - \t \right) \left( \Exp{b} \otimes u \otimes l \otimes v \otimes \Exp{b} \right) ={} &
		u \otimes l \otimes v \otimes \Exp{b} -
		                              (-1)^{\braid{\degb{v}}{\degb{u} + \degb{l}}} v \otimes \Exp{b} \otimes u \otimes l. \label{eq:tebulveb}
	\end{align}
\end{lm}
\begin{proof}
	Since $b$ has degree zero, we have
	\begin{equation*}
		\begin{aligned}
			\t \left( \Exp{b} \otimes u \otimes \Exp{b} \right) & =
			\sum_{i \geq 0, j \geq 1} \t \left( b^{\otimes i} \otimes u \otimes
			b^{\otimes j} \right) + \sum_{i \geq 0} \t \left( b^{\otimes i} \otimes u
			\right)                                                                                                                                                                    \\
			                                                    & = \sum_{i \geq 0, j \geq 1} b^{\otimes (i + 1)} \otimes u \otimes b^{\otimes
					                                                                                                                        (j-1)} + \sum_{i \geq 0} u \otimes b^{\otimes i} \\
			                                                    & = \sum_{i \geq 1, j \geq 0} b^{\otimes i} \otimes u \otimes b^{\otimes j}
			+ \sum_{j \geq 0} u \otimes b^{\otimes j}                                                                                                                                  \\
			                                                    & = \Exp{b} \otimes u \otimes \Exp{b}
		\end{aligned}
	\end{equation*}
	which shows \cref{eq:tebueb}. Finally, we have
	\begin{equation*}
		\begin{aligned}
			\t \left( \Exp{b} \otimes u \otimes l \otimes v \otimes \Exp{b} \right)
			={} &
			\sum_{i \geq 0, j \geq 1}
			\t \left( b^{\otimes i} \otimes u \otimes l \otimes v \otimes b^{\otimes j} \right) +
			\sum_{i \geq 0}
			\t \left( b^{\otimes i} \otimes u \otimes l \otimes v \right)
			\\
			={} & \sum_{i \geq 0, j \geq 1} b^{\otimes (i+1)} \otimes u \otimes l
			\otimes v \otimes b^{\otimes (j-1)}                                           \\
			    & + \sum_{i \geq 0} (-1)^{\braid{\degb{v}}{\degb{b^{\otimes i}} + \degb{u}
					        + \degb{l}}} v \otimes b^{\otimes i} \otimes u \otimes l \\
			={} & \sum_{i, j \geq 0} b^{\otimes i} \otimes u \otimes l \otimes v \otimes
			b^{\otimes j} - \sum_{j \geq 0} u \otimes l \otimes v \otimes b^{\otimes j}
			\\
			    & + \sum_{i \geq 0} (-1)^{\braid{\degb{v}}{\degb{u} + \degb{l}}}
			v \otimes b^{\otimes i} \otimes u \otimes l                                   \\
			={} & \Exp{b} \otimes u \otimes l \otimes v \otimes \Exp{b} - u \otimes l
			\otimes v \otimes \Exp{b}
			\\
			    & + \,  (-1)^{\braid{\degb{v}}{\degb{u} + \degb{l}}}
			v \otimes \Exp{b} \otimes u \otimes l
		\end{aligned}
	\end{equation*}
	which shows \cref{eq:tebulveb}.
\end{proof}

\begin{lm} \label{lm:cycl-f-1-t}
	We have the identity $\cycl{f} \circ (\idd - \t) = (\idd - \t) \circ f$ on
	$\tensr{V}$.
\end{lm}
\begin{proof}
	We will first verify the identity for elementary tensors of weight greater
	than or equal to two. Such tensors can be written as $x \otimes l \otimes z$
	with $x, z \in V$ and $l = v_1 \otimes \dots \otimes v_k$ with $k \geq 0$. Let
	us set $b = f_{0}(1) \in W^0$ so that $f(1) = \Exp{b}$. Then using the fact
	that $f$ is a coalgebra morphism, we have the identities
	\begin{align}
		f \left( l \otimes z \right) ={} & f \left( l_{(1)} \right) \otimes
		\corest{f} \left( l_{(2)} \otimes z \right) \otimes \Exp{b}, \label{eq:hatflz}
		\\
		f \left( x \otimes l \right) ={} & \Exp{b} \otimes \corest{f} \left( x \otimes l_{(1)}
		\right) \otimes f \left( l_{(2)} \right), \label{eq:hatfxl}
		\\
		\begin{split}
			f \left( x \otimes l \otimes z \right) ={} &
			\Exp{b} \otimes \corest{f} \left( x \otimes l_{(1)} \right) \otimes f \left(
			l_{(2)} \right) \otimes \corest{f} \left( l_{(3)} \otimes z \right) \otimes
			\Exp{b}
			\\
			                                           & + \Exp{b} \otimes \corest{f} \left( x \otimes l \otimes z \right) \otimes \Exp{b}.
		\end{split} \label{eq:hatfxlz}
	\end{align}
	Plugging the expression for $f \left( x \otimes l \otimes z \right)$
	into $\idd - \t$ and using the fact that $f$ has degree zero we get that
	\begin{equation*}
		\begin{aligned}
			\left( (\idd - \t) \circ f \right) \left( x \otimes l \otimes z \right)
			\stackrel{\eqref{eq:hatfxlz}}{=}{}                     &
			\left( \idd - \t \right) \left(
			\Exp{b} \otimes \corest{f} \left( x \otimes l_{(1)} \right) \otimes f
			\left( l_{(2)} \right) \otimes \corest{f} \left( l_{(3)} \otimes z \right) \otimes \Exp{b}
			\right)
			\\
			                                                       & +
			\left( \idd - \t \right) \left(
			\Exp{b} \otimes \corest{f} \left( x \otimes l \otimes z \right) \otimes \Exp{b}
			\right)
			\\
			\stackrel[\eqref{eq:tebulveb}]{\eqref{eq:tebueb}}{=}{} &
			\corest{f} \left( x \otimes l_{(1)} \right) \otimes f \left( l_{(2)}
			\right) \otimes \corest{f} \left( l_{(3)} \otimes z \right) \otimes \Exp{b}
			\\
			                                                       & - (-1)^{\braid{\degb{l_{(3)}} + \degb{z}}{\degb{x} + \degb{l_{(1)}} + \degb{l_{(2)}}}}
			\\
			                                                       & \qquad
			\corest{f} \left( l_{(3)} \otimes z \right) \otimes \Exp{b} \otimes \corest{f} \left( x
			\otimes l_{(1)} \right) \otimes f \left( l_{(2)} \right)
			\\
			\stackrel[\eqref{eq:hatfxl}]{\eqref{eq:hatflz}}{=}{}   &
			\corest{f} \left( x \otimes l_{(1)} \right) \otimes f \left( l_{(2)}
			\otimes z \right)
			\\
			                                                       & -
			                                                         (-1)^{\braid{\degb{l_{(2)}} + \degb{z}}{\degb{x} + \degb{l_{(1)}}}}
			\corest{f} \left( l_{(2)} \otimes z \right) \otimes f \left( x \otimes l_{(1)} \right).
		\end{aligned}
	\end{equation*}

	Next, we turn to compute
	$\left( \cycl{f} \circ \left( \idd - \t \right) \right) \left( x \otimes l \otimes z \right)$.
	According to the definition of $\cycl{f}$, we have
	\begin{equation*}
		\begin{aligned}
			\cycl{f} \left( x \otimes l \otimes z \right)
			\stackrel{\eqref{def:cyclization-morphism}}{=}{} &
			(-1)^{\braid{\degb{ \left( l \otimes z \right)_{(3)}}}
				{\degb{x} + \degb{ \left( l \otimes z \right)_{(1)}} +
					\degb{\left( l \otimes z \right)_{(2)}}}}
			\\
			                                                 & \qquad\quad
			\corest{f} \left( \left( l \otimes z \right)_{(3)} \otimes x \otimes \left( l
			                                                                     \otimes z \right)_{(1)} \right) \otimes f \left( \left( l \otimes z
			                                                                                                               \right)_{(2)} \right).
		\end{aligned}
	\end{equation*}
	Using \cref{eq:lzsplit} we can write the above expression as
	\begin{equation*}
		\begin{aligned}
			\cycl{f} \left( x \otimes l \otimes z \right)
			\stackrel{\eqref{eq:lzsplit}}{=}{} &
			\corest{f} \left( x \otimes l \otimes z \right) \otimes \Exp{b}
			+
			\corest{f} \left( x \otimes l_{(1)} \right) \otimes f \left( l_{(2)} \otimes
			z \right)
			\\
			                                   & + (-1)^{\braid{\degb{l_{(3)}} + \degb{z}}
				                                     {\degb{x} + \degb{l_{(1)}} + \degb{l_{(2)}}}}
			\corest{f} \left( l_{(3)} \otimes z \otimes x \otimes l_{(1)} \right) \otimes
			f \left( l_{(2)} \right).
		\end{aligned}
	\end{equation*}
	Similarly, using \cref{eq:xlsplit} we have
	\begin{equation*}
		\begin{aligned}
			\left( \cycl{f} \circ \t \right) \left( x \otimes l \otimes z \right)
			\stackrel{\eqref{eq:def-t-rotation}}{=} {}        &
			(-1)^{\braid{\degb{z}}{\degb{x} + \degb{l}}}
			\cycl{f} \left( z \otimes x \otimes l \right)
			\\
			\stackrel{\eqref{def:cyclization-morphism}}{=} {} & (-1)^{\braid{\degb{z}}{\degb{x} + \degb{l}} +
				                                                    \braid{\degb{\left( x \otimes l \right)_{(3)}}}
				                                                    {\degb{z} + \degb{\left( x \otimes l \right)_{(1)}} +
					                                                    \degb{\left( x \otimes l \right)_{(2)}}}}
			\\
			                                                  & \qquad
			\corest{f} \left( \left( x \otimes l \right)_{(3)} \otimes z \otimes \left( x
			                                                                     \otimes l \right)_{(1)} \right) \otimes f \left( \left( x \otimes l
			                                                                                                               \right)_{(2)} \right)
			\\
			\stackrel{\eqref{eq:xlsplit}}{=} {}               &
			(-1)^{\braid{\degb{z}}
				{\degb{x} + \degb{l_{(1)}} + \degb{l_{(2)}} +
					\degb{l_{(3)}}} +
				\braid{\degb{l_{(3)}}}
				{\degb{z} + \degb{x} + \degb{l_{(1)}} + \degb{l_{(2)}}}}
			\\
			                                                  & \qquad
			\corest{f} \left( l_{(3)} \otimes z \otimes x \otimes l_{(1)} \right) \otimes
			f \left( l_{(2)} \right)
			\\
			                                                  & + (-1)^{\braid{\degb{z}}{\degb{x} + \degb{l_{(1)}} + \degb{l_{(2)}}} +
				                                                    \braid{\degb{l_{(2)}}}
				                                                    {\degb{z} + \degb{x} + \degb{l_{(1)}}}}
			\\
			                                                  & \qquad
			\corest{f} \left( l_{(2)} \otimes z \right) \otimes f \left( x \otimes
			l_{(1)} \right)
			\\
			                                                  & + (-1)^{\braid{\degb{z}}{\degb{x} + \degb{l}} +
				                                                    \braid{\degb{x} + \degb{l}}{\degb{z}}}
			\\
			                                                  & \qquad
			\corest{f} \left( x \otimes l \otimes z \right) \otimes \Exp{b}
			\\
			\stackrel{\phantom{\eqref{eq:xlsplit}}}{=} {}     & (-1)^{\braid{\degb{l_{(3)}} + \degb{z}}
				                                                    {\degb{x} + \degb{l_{(1)}} + \degb{l_{(2)}}}}
			\corest{f} \left( l_{(3)} \otimes z \otimes x \otimes l_{(1)} \right) \otimes
			f \left( l_{(2)} \right)
			\\
			                                                  & + (-1)^{\braid{\degb{l_{(2)}} + \degb{z}}{\degb{x} + \degb{l_{(1)}}}}
			\corest{f} \left( l_{(2)} \otimes z \right) \otimes f \left( x \otimes
			l_{(1)} \right)
			\\
			                                                  & + \corest{f} \left( x \otimes l \otimes z \right) \otimes \Exp{b}.
		\end{aligned}
	\end{equation*}
	Canceling the identical terms, we obtain
	\begin{equation*}
		\begin{aligned}
			\left( \cycl{f} \circ (\idd - \t) \right) \left( x \otimes l \otimes z
			\right) ={} & \corest{f} \left( x \otimes l_{(1)} \right) \otimes f \left( l_{(2)} \otimes z \right)
			\\
			            & - (-1)^{\braid{\degb{l_{(2)}} + \degb{z}}{\degb{x} + \degb{l_{(1)}}}}
			\corest{f} \left( l_{(2)} \otimes z \right) \otimes f \left( x \otimes
			l_{(1)} \right)
		\end{aligned}
	\end{equation*}
	which is the same expression we got for $(\idd - \t) \circ f$.

	This shows the identity for elementary tensors of weight greater than or
	equal to two. Finally, if $x \in V$, we have
	\begin{equation*}
		\left( \cycl{f} \circ (\idd - \t) \right)(x) = \cycl{f}(0) = 0
	\end{equation*}
	and also
	\begin{equation*}
		\left( (\idd - \t) \circ f \right)(x) = \left( \idd - \t \right) \left(
		\Exp{b} \otimes f_1(x) \otimes \Exp{b} \right) \stackrel{\eqref{eq:tebueb}}{=} 0.
	\end{equation*}
\end{proof}

\Cref{lm:cycl-f-1-t} implies that $\cycl{f}$ descends to the quotient
$\tensrcyc{V}$. We will continue to denote by $\cycl{f} \colon \tensrcyc{V} \rightarrow
	\tensrcyc{W}$ the induced map between the quotients.

\subsection{Extension To Full Tensor Module} \label{sec:extension-cycl-full-tensor-module}
Let $V$ be a graded Banach $R$-module and extend the definition of the rotation operation
$\t \colon \tensr{V} \rightarrow \tensr{V}$, given by \cref{eq:def-t-rotation},
to an operator $\t \colon \tens{V} \rightarrow \tens{V}$ by setting $\t \left( 1 \right) \defeq 1$. We will call
\begin{equation*}
	\tenscyc{V} \defeq \tens{V} / \Im \left( \idd - \t \right) \cong R \oplus \tensrcyc{V}
\end{equation*}
the \textbf{cyclic tensor module} on $V$. In this subsection, we extend the definitions
of the cyclization of coderivations and morphisms from
the reduced tensor module $\tensr{V}$ to the full tensor module $\tens{V}$ while trying
to preserve as much as possible the ``functorial'' properties of the construction
(\cref{lm:cycl-comm-bracket,lm:func-cycl-morphism,lm:func-cycl-morphism-coderivation}).
We will see that our extension is not functorial on $\tens{V}$ but will be functorial
on the quotient $\tenscyc{V}$.

Let $\mu \colon \tens{V} \rightharpoonup \tens{V}$ be a generalized coderivation. We can naturally extend
the definition of $\cycl{\mu} \colon \tensr{V} \rightharpoonup \tensr{V}$, given by \cref{def:cyclization-coder},
from $\tensr{V}$ to $\tens{V}$ by defining
\begin{equation} \label{def:cycl-mu-1}
	\cycl{\mu}(r) \defeq \mu \left( r \right) = d_{\mu} \left( r \right) + (-1)^{\braidd{\mu}{r}} r \cdot \mu_0(1).
\end{equation}
The resulting map $\cycl{\mu} \colon \tens{V} \rightharpoonup \tens{V}$ is
a module derivation on $\tens{V}$ over $d_{\mu}$ with $\nnorm[\cycl{\mu}] \leq \nnorm[\mu]$,
and induces a map $\cycl{\mu} \colon \tenscyc{V} \rightharpoonup \tenscyc{V}$ on
the cyclic tensor module.

Note that no matter how we define $\cycl{\mu} \left( r \right)$, the identity
$\cycl{\mu} \circ \left( \idd - \t \right) = \left( \idd - \t \right) \circ \mu$
of \Cref{lm:coder-descends-quotient} holds on $\tens{V}$ and not only on
$\tensr{V}$.
With the definition above, \cref{lm:cycl-comm-bracket} does not continue to hold on $\tens{V}$, but it does hold
modulo $\Im(\idd - \t)$:
\begin{lm} \label{lm:cycl-zero-degree-piece-identity}
	Let $\mu, \nu \in \CoDer{\tens{V}}$ and let $r \in R$. Then
	\begin{equation*}
		\cycl{[\mu,\nu]}(r) \equiv \left[ \cycl{\mu}, \cycl{\nu} \right](r) \mod \Im \left( \idd - \t \right).
	\end{equation*}
\end{lm}
\begin{proof}
	We have
	\begin{equation*}
		\begin{aligned}
			\left[ \cycl{\mu}, \cycl{\nu} \right](1) & = \cycl{\mu} \left( \nu_0(1)
			\right) - (-1)^{\braidd{\mu}{\nu}} \cycl{\nu} \left( \mu_0(1) \right)
			\\
			                                         & = \mu_1 \left( \nu_0(1) \right) + (-1)^{\braidd{\mu}{\nu}} \nu_0(1) \otimes
			\mu_0(1)
			\\
			                                         & - (-1)^{\braidd{\mu}{\nu}} \nu_1 \left( \mu_0(1) \right)
			                                                                                            -(-1)^{\braidd{\mu}{\nu} + \braidd{\nu}{\mu}} \mu_0(1) \otimes \nu_0(1)
			\\
			                                         & = \mu_1 \left( \nu(1) \right) - (-1)^{\braidd{\mu}{\nu}} \nu_1 \left(
			\mu(1) \right) - \left(\idd - \t \right) \left( \mu_0(1) \otimes \nu_0(1)
			\right)
			\\
			                                         & = \cycl{[\mu,\nu]}(1) - \left(\idd - \t \right) \left( \mu_0(1) \otimes
			\nu_0(1) \right)
		\end{aligned}
	\end{equation*}
	which shows the lemma for $r = 1$. Since both
	$\cycl{[\mu,\nu]}$ and $\left[ \cycl{\mu}, \cycl{\nu} \right]$ are module
	derivations on $\tens{V}$ over the same derivation $d_{[\mu, \nu]} = \left[ d_{\mu}, d_{\nu} \right]$,
	their difference is $R$-linear and hence
	$\cycl{[\mu,\nu]}(r) \equiv \left[ \cycl{\mu}, \cycl{\nu} \right](r) \mod \Im \left( \idd - \t \right)$
	for all $r \in R$.
\end{proof}

\Cref{lm:cycl-comm-bracket} and \cref{lm:cycl-zero-degree-piece-identity} together
imply the following:

\begin{cor} \label{cor:cycl-comm-bracket-quotient}
	Let $\mu, \nu \in \CoDer{\tens{V}}$ be two generalized coderivations on $\tens{V}$.
	The induced maps $\cycl{\mu}, \cycl{\nu} \colon
		\tenscyc{V} \rightharpoonup \tenscyc{V}$ satisfy $\left[ \cycl{\mu}, \cycl{\nu}
			\right] = \cycl{[\mu,\nu]}$ on the quotient $\tenscyc{V}$.
	Hence, the cyclic tensor module $\tenscyc{V}$ is also a representation of the Lie algebra
	$\CoDer{\tens{V}}$ with respect to the action map $\left( \mu, l \right) \mapsto \cycl{\mu} \left( l \right)$.
	\qed
\end{cor}

We turn next to the cyclization of morphisms. Let $W$ be a graded Banach $R$-module and let
$f \colon \tens{V} \rightarrow \tens{W}$ be a morphism of graded Banach $R$-coalgebras. It will
be useful to introduce the following notation:

\begin{dfn} \label{dfn:cyclic-exponential}
	Let $b \in V^0$ be topologically nilpotent. The \textbf{cyclic exponential} of $b$ is defined
	to be
	\begin{equation} \label{def:exp-cyc}
		\cexp{b} \defeq 1 + b + \frac{b \otimes b}{2} + \frac{b \otimes b \otimes b}{3} +
		\dots \in \tens{V}.
	\end{equation}
\end{dfn}

Let us extend the definition of $\cycl{f} \colon \tensr{V} \rightarrow \tensr{W}$ to
a map $\cycl{f} \colon \tens{V} \rightarrow \tens{W}$ by defining
\begin{equation} \label{def:cycl-f-1}
	\cycl{f}(1) \defeq \cexp{f_0(1)} = 1 + f_0(1) + \frac{f_0(1) \otimes f_0(1)}{2} +
	\dots
\end{equation}
Thus, we have $f(1) = \Exp{f_0(1)}$ while $\cycl{f}(1) = \cexp{f_0(1)}$.
The extended map $\cycl{f}$ is a morphism of graded Banach modules with $\nnorm[\cycl{f}] \leq \nnorm[f]$,
and induces a map $\cycl{f} \colon \tenscyc{V} \rightarrow \tenscyc{W}$ between
the cyclic tensor modules.

Note that no matter how we define $\cycl{f}(1)$, the identity
$\cycl{f} \circ (\idd - \t) = (\idd - \t) \circ f$ of \Cref{lm:cycl-f-1-t}
holds on $\tens{V}$ and not only on $\tensr{V}$,
a consequence of the identity $\t \left( \Exp{b} \right) = \Exp{b}$,
applied to $b = f_0(1)$.

Given a topologically nilpotent $b \in V^0$, we have the identity
\begin{equation*}
	f \left( \Exp{b} \right) = \Exp{ \mcfunc{f} \left( b \right) }
\end{equation*}
for the action of the morphism $f$ on the ``regular'' exponential.
We have an analogous identity for the action of $\cycl{f}$ on the cyclic exponential
(see \cref{fig:exp-cyc-mcfunc} and compare to \cref{fig:def-pushforward-mc}):

\begin{lm} \label{lem:func-mc-cyclic}
	Let $f \colon \tens{V} \rightarrow \tens{W}$ be a morphism of Banach $R$-coalgebras and let
	$b \in V^{0}$ be topologically nilpotent. Then we have the identity
	\begin{equation} \label{eq:func-cyclic-exponent}
		\cycl{f} \left( \cexp{b} \right) \equiv \cexp{ \mcfunc{f} \left( b \right) }
		\mod \Im \left( \idd - \t \right).
	\end{equation}
\end{lm}
\begin{proof}
	Let $b \in V^0$ be topologically nilpotent. We have
	\begin{equation*}
		\begin{aligned}
			\cycl{f} \left( \cexp{b} \right)
			\stackrel{\eqref{def:exp-cyc}}{=}{}                             &
			\cycl{f} \left( 1 + \sum_{n=1}^{\infty} \frac{b^{\otimes n}}{n} \right)
			\\
			\stackrel{\eqref{def:cyclization-morphism}}{=}{}                &
			\cycl{f}(1) +
			\sum_{n=1}^{\infty} \frac{1}{n} \left(
			\sum_{\substack{k_1 + k_2 + k_3 + 1 = n,                                                                    \\ k_1 \geq 0, k_2 \geq 0, k_3 \geq 0}}
			\corest{f} \left( b^{\otimes \left( k_3 + 1 + k_1 \right)} \right) \otimes
			f \left( b^{\otimes k_2} \right) \right)
			\\
			\stackrel{\phantom{\eqref{def:cyclization-morphism}}}{=}{}      &
			\cycl{f}(1) +
			\sum_{n=1}^{\infty}
			\sum_{\substack{l_1 + l_2 = n,                                                                              \\ l_1 \geq 1, l_2 \geq 0}}
			\frac{l_1}{n} f_{l_1} \left( b^{\otimes l_1} \right) \otimes
			f \left( b^{\otimes l_2} \right)
			\\
			\stackrel{\phantom{\eqref{def:cyclization-morphism}}}{=}{}      &
			\cycl{f}(1) + \sum_{n,r=1}^{\infty}
			\sum_{\substack{n_1 + \dots + n_r = n,                                                                      \\ n_1, \dots, n_r \geq 0}}
			\frac{n_1}{n} f_{n_1} \left( b^{\otimes n_1} \right) \otimes \dots \otimes
			f_{n_r} \left( b^{\otimes n_r} \right)
			\\
			\stackrel{\phantom{\eqref{def:cyclization-morphism}}}{\equiv}{} &
			\cycl{f}(1) + \sum_{n,r=1}^{\infty} \frac{1}{r} \left(
			\sum_{\substack{n_1 + \dots + n_r = n,                                                                      \\ n_1, \dots, n_r \geq 0}}
			f_{n_1} \left( b^{\otimes n_1} \right) \otimes \dots \otimes
			f_{n_r} \left( b^{\otimes n_r} \right) \right)
			\\
			                                                                & \qquad \mod \Im \left( \idd - \t \right).
		\end{aligned}
	\end{equation*}
	On the other hand,
	\begin{equation*}
		\begin{aligned}
			\cexp{ \mcfunc{f} \left( b \right)}
			\stackrel{\eqref{def:exp-cyc}}{=}{}           &
			1 + \sum_{r=1}^{\infty} \frac{{\mcfunc{f} \left( b \right)}^{\otimes r}}{r}
			\stackrel{\eqref{eq:mcfunc-explicit}}{=}{}
			1 + \sum_{r=1}^{\infty} \frac{1}{r} \left( \sum_{n=0}^{\infty}
			f_n \left( b^{\otimes n} \right) \right)^{\otimes r}
			\\
			\stackrel{\phantom{\eqref{def:exp-cyc}}}{=}{} &
			1 + \sum_{r = 1}^{\infty} \frac{1}{r} \left(
			                                      f_0 \left( 1 \right)^{\otimes r} +
			\sum_{n=1}^{\infty} \sum_{\substack{n_1 + \dots + n_r = n \\ n_1, \dots, n_r \geq 0}}
			f_{n_1} \left( b^{\otimes n_1} \right) \otimes \dots \otimes
			f_{n_r} \left( b^{\otimes n_r} \right)
			\right)
			\\
			\stackrel{\eqref{def:cycl-f-1}}{=}{}          &
			\cycl{f}(1) +
			\sum_{n,r=1}^{\infty} \frac{1}{r} \left(
			\sum_{\substack{n_1 + \dots + n_r = n                     \\ n_1, \dots, n_r \geq 0}}
			f_{n_1} \left( b^{\otimes n_1} \right) \otimes \dots \otimes
			f_{n_r} \left( b^{\otimes n_r} \right) \right)
		\end{aligned}
	\end{equation*}
	which shows the lemma.
\end{proof}

\begin{figure}
	\begin{tikzcd}
		\tc{V} && {\tenscyc{V}} \\
		\tc{W} && {\tenscyc{W}}
		\arrow["\cexp{-}", from=1-1, to=1-3]
		\arrow["{\mcfunc{f}}"', from=1-1, to=2-1]
		\arrow["\cycl{f}", from=1-3, to=2-3]
		\arrow["\cexp{-}", from=2-1, to=2-3]
	\end{tikzcd}
	\caption{The relation between the pushforward map and the cyclic exponential.}
	\label{fig:exp-cyc-mcfunc}
\end{figure}

The following lemma shows that our choice of defining $\cycl{f} \left( 1 \right)$ using \cref{def:cycl-f-1}
allows us to retain functoriality modulo $\Im \left( \idd - \t \right)$, and hence obtain true functoriality
on the cyclic tensor modules.
\begin{lm} \label{lm:func-cycl-morphism-extended-modulo}
	Let $U,V,W$ be graded Banach $R$-modules and let $f \colon \tens{V} \rightarrow \tens{W}$
	and $g \colon \tens{U} \rightarrow \tens{V}$ be two morphisms of graded Banach $R$-coalgebras.
	Then we have
	\begin{equation*}
		\cycl{\left( f \circ g \right)} \left( 1 \right) \equiv
		\cycl{f} \left( \cycl{g} \left( 1 \right) \right) \mod \Im \left( \idd - \t \right).
	\end{equation*}
\end{lm}
\begin{proof}
	We have
	\begin{equation*}
		\begin{aligned}
			\cycl{\left( f \circ g \right)} \left( 1 \right)
			\stackrel{\eqref{def:cycl-f-1}}{=}{}       &
			\cexp{ \left( f \circ g \right)_0 \left( 1 \right)} =
			\cexp{\left( \corest{f} \circ g \right) \left( 1 \right)} =
			\cexp{\corest{f} \left( \Exp{g_0(1)} \right)}
			\\
			\stackrel{\eqref{eq:mcfunc-explicit}}{=}{} &
			\cexp{\mcfunc{f} \left( g_0(1) \right)}
		\end{aligned}
	\end{equation*}
	while
	\begin{equation*}
		\cycl{f} \left( \cycl{g} \left( 1 \right) \right)
		\stackrel{\eqref{def:cycl-f-1}}{=}{}
		\cycl{f} \left( \cexp{g_0(1)} \right)
		\stackrel{\eqref{eq:func-cyclic-exponent}}{\equiv}{}
		\cexp{\mcfunc{f} \left( g_0(1) \right)} \mod \Im \left( \idd - \t \right).
	\end{equation*}
\end{proof}

\Cref{lm:func-cycl-morphism} and \cref{lm:func-cycl-morphism-extended-modulo} together
imply that

\begin{cor} \label{cor:func-cycl-morphism-extended}
	Let $U,V,W$ be graded Banach $R$-modules and let $f \colon \tens{V} \rightarrow \tens{W}$
	and $g \colon \tens{U} \rightarrow \tens{V}$ be two morphisms of graded Banach $R$-coalgebras.
	The induced maps $\cycl{f} \colon \tenscyc{V} \rightarrow \tenscyc{W}$ and
	$\cycl{g} \colon \tenscyc{U} \rightarrow \tenscyc{V}$ satisfy
	$\cycl{\left( f \circ g \right)} = \cycl{f} \circ \cycl{g}$ on the \textbf{quotient} $\tenscyc{U}$.
	In particular, the cyclic tensor module $\tenscyc{V}$ is a representation of the monoid
	of Banach $R$-coalgebra endomorphisms of $\tens{V}$ with respect to the action map
	$\left( f, l \right) \mapsto \cycl{f} \left( l \right)$. \qed
\end{cor}

Next, we show that \cref{lm:func-cycl-morphism-coderivation} continues to hold on
$\tens{V}$ modulo $\Im \left( \idd - \t \right)$. In order to do that, we will need to study
the interaction between the cyclization of a coderivation and the cyclic exponential.
Given a generalized coderivation $\nu$ on $\tens{V}$ and a topologically nilpotent $b \in V^0$,
we have the following identity for the ``regular'' exponential:
\begin{equation*}
	\nu \left( \Exp{b} \right) = \Exp{b} \otimes \corest{\nu} \left( \Exp{b} \right) \otimes \Exp{b}.
\end{equation*}
Analogously, we have the following identity for the cyclic exponential:
\begin{lm} \label{lm:coder-cycl-exp-identity}
	Let $\nu \colon \tens{V} \rightharpoonup \tens{V}$ be a generalized coderivation on $\tens{V}$
	and let $b \in V^{0}$ be topologically nilpotent. Then:
	\begin{equation} \label{eq:cycl-exp-identity}
		\cycl{\nu} \left( \cexp{b} \right) \equiv
		\corest{\nu} \left( \Exp{b} \right) \otimes \Exp{b} \mod \Im \left( \idd - \t \right).
	\end{equation}
\end{lm}
\begin{proof}
	Given $k \geq 0$, we have
	\begin{equation*}
		\begin{aligned}
			\cycl{\nu} \left( b^{\otimes (k + 1)} \right)
			\stackrel{\eqref{def:cyclization-coder}}{=}{}                &
			\sum_{k_1 + k_2 + k_3 = k}
			b \otimes b^{\otimes k_1} \otimes \nu_{k_2} \left( b^{\otimes k_2} \right) \otimes b^{\otimes k_3} +
			\nu_{k_3 + 1 + k_1} \left( b^{\otimes \left( k_3 + 1 + k_1 \right)} \right) \otimes b^{\otimes k_2}
			\\
			\stackrel{\phantom{\eqref{def:cyclization-coder}}}{\equiv}{} &
			\underbrace{
				\sum_{k_1 + k_2 + k_3 = k}
				\nu_{k_2} \left( b^{\otimes k_2} \right) \otimes b^{\otimes \left( k_3 + 1 + k_1 \right)} +
				\nu_{k_3 + 1 + k_1} \left( b^{\otimes \left( k_3 + 1 + k_1 \right)} \right) \otimes b^{\otimes k_2}
			}_{\bigstar}
			\\
			                                                             & \qquad
			\mod \Im \left( \idd - \t \right),
		\end{aligned}
	\end{equation*}
	and
	\begin{equation*}
		\begin{aligned}
			\bigstar & =
			\sum_{\substack{i + j = k + 1 \\ i \geq 0, j \geq 1}}
			j \cdot \left(
			\nu_{i} \left( b^{\otimes i} \right) \otimes b^{\otimes j} +
			\nu_j \left( b^{\otimes j} \right) \otimes b^{\otimes i}
			\right)
			\\
			         & =
			\sum_{\substack{i + j = k + 1 \\ i, j \geq 0}}
			j \cdot \left(
			\nu_{i} \left( b^{\otimes i} \right) \otimes b^{\otimes j} +
			\nu_j \left( b^{\otimes j} \right) \otimes b^{\otimes i}
			\right)
			\\
			         & =
			\left( k + 1 \right) \cdot
			\sum_{\substack{i + j = k + 1 \\ i, j \geq 0}}
			\nu_{i} \left( b^{\otimes i} \right) \otimes b^{\otimes j}.
		\end{aligned}
	\end{equation*}
	Hence,
	\begin{equation*}
		\begin{aligned}
			\cycl{\nu} \left( \cexp{b} \right) & = \nu_0(1) +
			\sum_{k=0}^{\infty} \frac{ \cycl{\nu} \left( b^{\otimes \left( k + 1 \right)} \right)}{k+1}
			\equiv
			\nu_0(1) +
			\sum_{k=0}^{\infty} \sum_{\substack{i + j = k + 1 \\ i, j \geq 0}}
			\nu_{i} \left( b^{\otimes i} \right) \otimes b^{\otimes j}
			\\
			                                   & =
			\corest{\nu} \left( \Exp{b} \right) \otimes \Exp{b} \mod \Im \left( \idd - \t \right).
		\end{aligned}
	\end{equation*}
\end{proof}

\begin{lm} \label{lm:func-cycl-morphism-coderivation-extended-modulo}
	Let $f \colon \tens{V} \rightarrow \tens{W}$ be a morphism of graded Banach $R$-coalgebras.
	Let $\mu$ be a generalized coderivation on $\tens{V}$ and $\nu$ be a generalized
	coderivation on $\tens{W}$. Assume that we have
	$f \circ \mu = \nu \circ f$. Then
	\begin{equation*}
		\left( \cycl{f} \circ \cycl{\mu} \right) (r) \equiv \left( \cycl{\nu} \circ
		\cycl{f} \right) (r) \mod \Im \left( \idd - \t \right)
	\end{equation*}
	for all $r \in R$.
\end{lm}
\begin{proof}
	Set $b = f_0(1)$. By assumption, we have
	\begin{equation*}
		f_1 \left( \mu_0(1) \right) =
		\left( \corest{f} \circ \mu \right)(1) = \left( \corest{\nu} \circ f \right)(1) =
		\corest{\nu} \left( \Exp{b} \right).
	\end{equation*}
	Hence
	\begin{equation*}
		\left( \cycl{f} \circ \cycl{\mu} \right)(1)
		\stackrel{\eqref{def:cycl-mu-1}}{=}{}
		\cycl{f} \left( \mu_0 \left( 1 \right) \right)
		\stackrel{\eqref{def:cyclization-morphism}}{=}{}
		f_1 \left( \mu_0 \left( 1 \right) \right) \otimes \Exp{b}
		=
		\corest{\nu} \left( \Exp{b} \right) \otimes \Exp{b}
	\end{equation*}
	and
	\begin{equation*}
		\left( \cycl{\nu} \circ \cycl{f} \right)(1)
		\stackrel{\eqref{def:cycl-f-1}}{=}{}
		\cycl{\nu} \left( \cexp{b} \right)
		\stackrel{\eqref{eq:cycl-exp-identity}}{\equiv}{}
		\corest{\nu} \left( \Exp{b} \right) \otimes \Exp{b} \mod \Im \left( \idd - \t \right)
	\end{equation*}
	which shows the lemma for $r = 1$. This is enough since the identity
	$f \circ \mu = \nu \circ f$ implies that $d_{\mu} = d_{\nu}$
	and hence the difference $\cycl{f} \circ \cycl{\mu} - \cycl{\nu} \circ \cycl{f}$ is $R$-linear.
\end{proof}

\Cref{lm:func-cycl-morphism-coderivation} and \cref{lm:func-cycl-morphism-coderivation-extended-modulo}
together imply that
\begin{cor} \label{cor:func-cycl-morphism-coderivation-extended}
	Let $f \colon \tens{V} \rightarrow \tens{W}$ be a morphism of graded Banach $R$-coalgebras.
	Let $\mu$ be a generalized coderivation on $\tens{V}$ and $\nu$ be a generalized
	coderivation on $\tens{W}$. Assume that we have
	$f \circ \mu = \nu \circ f$.
	Then the induced maps
	\begin{equation*}
		\cycl{f} \colon \tenscyc{V} \rightarrow \tenscyc{W}, \quad
		\cycl{\mu} \colon \tenscyc{V} \rightharpoonup \tenscyc{V}, \quad
		\cycl{\nu} \colon \tenscyc{W} \rightharpoonup \tenscyc{W}
	\end{equation*}
	satisfy the identity $\cycl{f} \circ \cycl{\mu} = \cycl{\nu} \circ \cycl{f}$ on
	\textbf{the quotient} $\tenscyc{V}$.
\end{cor}

\subsection{Naturality with Respect to Pullback} \label{sec:cyclization-naturality}
We end the section by noting that the cyclization of generalized coderivations and morphisms
is natural and commutes with pullback in an appropriate sense, and discussing how the cyclization
of tensor coalgebra morphisms extends naturally to morphisms of tensor coalgebras over different ground
algebras.

Let $\mathcal{R} = \left( R, d_R \right)$ and $\mathcal{S} = \left( S, d_S \right)$ be
two pre-differential graded-commutative Banach $\mathbbm{k}$-algebras
and let $\varphi \colon \mathcal{R} \rightarrow \mathcal{S}$ be a morphism of pre-differential graded Banach
$\mathbbm{k}$-algebras.
Given a graded Banach $S$-module $W$ and a generalized coderivation
$\nu \colon \tens{W}[S] \rightharpoonup \tens{W}[S]$ over $d_S$, we can pull $\nu$ along $\varphi$
and obtain a generalized coderivation
\begin{equation*}
	\varphi^{*} \left( \nu \right) \colon
	\tens{\varphi^{*} \left( W \right)}[R] \rightharpoonup \tens{\varphi^{*} \left( W \right)}[R]
\end{equation*}
over $d_R$
(see \cref{dfn:pullback-generalized-coderivation}).
Recall that we have a canonical coalgebra morphism
\begin{equation*}
	\resover{\varphi} = \resover{\varphi}^{W} \colon \tens{\varphi^{*} \left( W \right)}[R] \rightarrow \tens{W}[S]
\end{equation*}
over $\varphi$ given by \cref{eq:canonical-pullback-map}.
One can verify directly from the definitions
that the cyclization of $\nu$ and the cyclization of the pullback
$\varphi^{*} \left( \nu \right)$ are related by the following commutative diagram:
\begin{figure}[h]
	\centering
	\begin{tikzcd}
		{\tens{W}[S]} && {\tens{W}[S]} \\
		{\tens{\varphi^{*} \left( W \right)}[R]} && {\tens{\varphi^{*} \left( W \right)}[R]}
		\arrow["{\cycl{\nu}}", harpoon, from=1-1, to=1-3]
		\arrow["{\cycl{\left( \varphi^{*} \left( \nu \right) \right)}}"', harpoon, from=2-1, to=2-3]
		\arrow["\resover{\varphi}^{W}", from=2-1, to=1-1]
		\arrow["\resover{\varphi}^{W}", from=2-3, to=1-3]
	\end{tikzcd}
	\caption{Naturality of cyclization of coderivations.}
	\label{fig:naturality-cyclization-coderivations}
\end{figure}

Similarly, let $R$ and $S$ be two graded-commutative Banach $\mathbbm{k}$-algebras
and let $\varphi \colon R \rightarrow S$ be a morphism of graded Banach $\mathbbm{k}$-algebras.
Given two graded Banach $S$-modules $W_1, W_2$ and a morphism
$f \colon \tens{W_1}[S] \rightarrow \tens{W_2}[S]$ of Banach $S$-coalgebras, we can pull $f$ along $\varphi$
and obtain a morphism $\varphi^{*} \left( f \right) \colon
	\tens{\varphi^{*} \left( W_1 \right)}[R] \rightarrow \tens{\varphi^{*} \left( W_2 \right)}[R]$
of Banach $R$-coalgebras (see \cref{dfn:pullback-morphism-tensor-coalgebras}).
One can verify directly from the definitions
that the cyclization of $f$ and the cyclization of the pullback
$\varphi^{*} \left( f \right)$ are related by the following commutative diagram:

\begin{figure}[h]
	\centering
	\begin{tikzcd}
		{\tens{W_1}[S]} && {\tens{W_2}[S]} \\
		{\tens{\varphi^{*} \left( W_1 \right)}[R]} && {\tens{\varphi^{*} \left( W_2 \right)}[R]}
		\arrow["{\cycl{f}}", from=1-1, to=1-3]
		\arrow["{\cycl{\left( \varphi^{*} \left( f \right) \right)}}"', from=2-1, to=2-3]
		\arrow["\resover{\varphi}^{W_1}", from=2-1, to=1-1]
		\arrow["\resover{\varphi}^{W_2}", from=2-3, to=1-3]
	\end{tikzcd}
	\caption{Naturality of cyclization of morphisms.}
	\label{fig:naturality-cyclization-morphisms}
\end{figure}

Now, let $V$ be a graded Banach $R$-module and let $W$ be a graded Banach $S$-module.
Given a morphism $f \colon \tens{V}[R] \rightarrow \tens{W}[S]$ of Banach coalgebras over $\varphi$,
we can factor $f$ uniquely as $f = \resover{\varphi} \circ \rescoho{f}$ where
$\rescoho{f} \colon \tens{V}[R] \rightarrow \tens{\varphi^{*} \left( W \right)}[R]$ is
a morphism of Banach $R$-coalgebras and
$\resover{\varphi} \colon \tens{\varphi^{*} \left( W \right)}[R] \rightarrow \tens{W}[S]$ is the canonical
morphism over $\varphi$.
The \textbf{cyclization} $\cycl{f} \colon \tens{V}[R] \rightarrow \tens{W}[S]$ of $f$ is a morphism
of graded Banach modules over $\varphi$ defined by
\begin{equation} \label{def:cyclization-morphism-different-ground-algebras}
	\cycl{f} \defeq \resover{\varphi} \circ \cycl{\rescoho{f}}.
\end{equation}
Explicitly, we have
\begin{align*}
	\cycl{f} \left( x \otimes_R l \right) & =
	(-1)^{\braid{\degb{l_{(3)}}}{\degb{x} + \degb{l_{(1)}} + \degb{l_{(2)}}}}
	\corest{f} \left( l_{(3)} \otimes_R x \otimes_R l_{(1)} \right) \otimes_S f \left(
	l_{(2)} \right),
	\\
	\cycl{f} \left( 1 \right)             & =
	\cexp{f_0(1)} = 1 + f_0(1) + \frac{f_0(1) \otimes_S f_0(1)}{2} + \dots
\end{align*}
for $x \in V$ and $l \in \tens{V}[R]$.
These are the same formulas as \eqref{def:cyclization-morphism} and \eqref{def:cycl-f-1},
with $\otimes = \otimes_R$ replaced by $\otimes_S$ in the appropriate places.

With the definition above, all the results of
\cref{sec:cyclization-coalgebra-morphisms,sec:extension-cycl-full-tensor-module} continue to hold for
morphisms of tensor coalgebras over different ground algebras with obvious modifications.
For example, let us verify that \cref{lm:func-cycl-morphism} continues to hold for morphisms
over different ground algebras. Let $Q, R, S$ be graded-commutative Banach $\mathbbm{k}$-algebras
and let $\psi \colon Q \rightarrow R$ and $\varphi \colon R \rightarrow S$ be morphisms
of graded Banach $\mathbbm{k}$-algebras. Let $U$ be a graded Banach $Q$-module, $V$ be
a graded Banach $R$-module and $W$ be a graded Banach $S$-module. Finally,
let $f \colon \tens{V}[R] \rightarrow \tens{W}[S]$ be a Banach coalgebra morphism over $\varphi$
and $g \colon \tens{U}[Q] \rightarrow \tens{V}[R]$ be a Banach coalgebra morphism over $\psi$.
Then we have the following commutative diagram:

\begin{figure}[h]
	\begin{tikzcd}
		& {\tens{V}[R]} & {\tens{\varphi^{*} \left( W \right)}[R]} & {\tens{W}[S]}
		\\
		{\tens{U}[Q]} & {\tens{\psi^{*} \left( V \right)}[Q]} &
		{\tens{\psi^{*} \left( \varphi^{*} \left( W \right) \right)}[Q]}
		\arrow["\rescoho{f}", from=1-2, to=1-3]
		\arrow["f", curve={height=-24pt}, from=1-2, to=1-4]
		\arrow["\resover{\varphi}^{W}", from=1-3, to=1-4]
		\arrow["g", from=2-1, to=1-2]
		\arrow["\rescoho{g}"', from=2-1, to=2-2]
		\arrow["\rescoho{\left( f \circ g \right)}"', curve={height=24pt}, from=2-1, to=2-3]
		\arrow["\resover{\psi}^{V}"', from=2-2, to=1-2]
		\arrow["\psi^{*} ( \rescoho{f} )"', from=2-2, to=2-3]
		\arrow["\resover{\psi}^{\varphi^{*} \left( W \right)}"', from=2-3, to=1-3]
		\arrow["\resover{\left( \varphi \circ \psi \right)}^{W}"', from=2-3, to=1-4]
	\end{tikzcd}
	\caption{Composition of morphisms over different ground algebras and restrictions.}
	\label{fig:composition-morphisms-different-ground}
\end{figure}

Hence, we see that
\begin{equation*}
	\begin{aligned}
		\cycl{ \left( f \circ g \right)}
		\stackrel{\eqref{def:cyclization-morphism-different-ground-algebras}}{=}{} &
		\resover{\left( \varphi \circ \psi \right)}^{W} \circ
		\cycl{ \left( \rescoho{\left( f \circ g \right)} \right)}
		\stackrel{\textrm{Figure } \ref{fig:composition-morphisms-different-ground}}{=}{}
		\resover{\varphi}^{W} \circ \resover{\psi}^{\varphi^{*} \left( W \right)} \circ
		\cycl{ \left( \psi^{*} ( \rescoho{f} ) \circ \rescoho{g} \right) }
		\\
		\stackrel{\eqref{eq:func-cycl-morphism}}{=}{}                              &
		\resover{\varphi}^{W} \circ \resover{\psi}^{\varphi^{*} \left( W \right)} \circ
		\cycl{ \left( \psi^{*} ( \rescoho{f} ) \right)} \circ
		\cycl{ \rescoho{g} }
		\stackrel{\textrm{Figure } \ref{fig:naturality-cyclization-morphisms}}{=}{}
		\resover{\varphi}^{W} \circ \cycl{ \rescoho{f} } \circ \resover{\psi}^{V} \circ \cycl{ \rescoho{g} }
		\\
		\stackrel{\eqref{def:cyclization-morphism-different-ground-algebras}}{=}{} &
		\cycl{f} \circ \cycl{g}.
	\end{aligned}
\end{equation*}

\begin{rem}
	Note that if we start with knowing that $\cycl{ \left( f \circ g \right)} = \cycl{f} \circ \cycl{g}$
	holds for morphisms over different base algebras, then the commutativity of the diagram in
	\cref{fig:naturality-cyclization-morphisms} becomes a consequence of the identity
	$f \circ \resover{\varphi}^{W_1} = \resover{\varphi}^{W_2} \circ \varphi^{*} \left( f \right)$
	defining the pullback $\varphi^{*} \left( f \right)$ (see \cref{fig:scalar-rest-morphism-formal-tensor-coalg}).
	Similarly, if we know that \cref{lm:func-cycl-morphism-coderivation} holds for morphisms of different base
	algebras, then the commutativity of the diagram in \cref{fig:naturality-cyclization-coderivations} becomes
	a consequence of the identity
	$\mu \circ \resover{\varphi}^W = \resover{\varphi}^W \circ \varphi^{*} \left( \mu \right)$ defining
	the pullback $\varphi^{*} \left( \mu \right)$ (see \cref{fig:pullback-of-coderivation}).
\end{rem}

\begin{rem}
	Note that the canonical map $\resover{\varphi}$ over $\varphi$ commutes with the rotation operator
	in the sense that
	$\resover{\varphi} \circ {\t}^{{\varphi^{*} \left( W \right)}} = {\t}^W \circ \resover{\varphi}$,
	and hence $\resover{\varphi}$ descends to an operator over $\varphi$ between the cyclic quotients.
	The naturality depicted in \cref{fig:naturality-cyclization-coderivations,fig:naturality-cyclization-morphisms}
	continues to hold if one replaces the tensor modules with the cyclic tensor modules.
\end{rem}

\section{Noncommutative Codifferential Forms and Their Calculus} \label{sec:noncomm-diff-calc}

This section introduces noncommutative (cyclic) codifferential forms and their Cartan calculus.
We begin in \cref{subsec:codifferential-forms} by introducing the bigraded module of noncommutative codifferential
forms $\ndf{V}[][]$, constructed as the tensor module on the shifted sum $V \oplus \ul{V}$.
We also introduce the module of cyclic codifferential forms $\ncdf{V}[][]$ by taking the quotient with respect
to the cyclic rotation operator.
In \cref{subsec:de-rham-differential}, we equip these spaces with the de Rham differential $\qdr$,
a degree $(-1,0)$ coderivation that strips the suspension line and satisfies $\qdr^2 = 0$.

\Cref{subsec:lie-derivative,sec:nc-d-calc-cont} systematically extend algebraic structures
from $\tens{V}$ to $\ndf{V}[][]$.
We prove that any generalized coderivation $\mu$ extends uniquely to a Lie derivative $\lie{\mu}$
on $\ndf{V}[][]$ characterized by the fact that it commutes with the de Rham differential ($[\qdr, \lie{\mu}] = 0$).
Furthermore, any linear coderivation induces a contraction operator $\cont{\mu}$,
uniquely determined by the Cartan relation $[\qdr, \cont{\mu}] = \lie{\mu}$.
In \cref{subsec:differential-calculus} we show that the operators satisfy a set of commutation relations, giving
us a Cartan calculus. The operators above have cyclic counterparts on the quotient $\ncdf{V}[][]$,
and we describe them in parallel, giving explicit formulas for their actions.

The remainder of the section explores the homological and categorical properties of this calculus.
\Cref{sec:formal-poincare} establishes the Formal Poincar\'{e} Lemma,
proving that the reduced noncommutative (cyclic) codifferential forms are acyclic with respect to $\qdr$
by constructing an explicit contracting homotopy using the Euler coderivation.
In \cref{sec:functoriality-ndf-and-ncdf}, we show that a morphism $f$ of tensor coalgebras extends
uniquely to an induced morphism $\indmap{f}$ between codifferential forms which commutes with the
de Rham differential $\qdr$. We show that the constructions of $\ndf{V}[][], \lie{\mu}, \cont{\mu}$
and their cyclic counterparts are natural and functorial.
Finally, in \cref{subsec:ndf-compatibility-ground-algebras}, we discuss the compatibility of our constructions
with change of ground algebras.

Our presentation is similar to the one given in \cite{Herscovich2023}, which is dual to the noncommutative
cyclic differential forms introduced by Kontsevich and Soibelman~\cite{Kontsevich:vh}.
In view of our applications, we allow the coderivations we work with to have a curvature term
and the morphisms to have a change of connection term, working with the full (and not reduced) tensor module.
Another difference compared to \cite{Herscovich2023} is that in passing to cyclic codifferential forms,
we use coinvariants instead of invariants as it appears more natural in geometric applications.
Thus, we work with the cyclic tensor module, which is a quotient instead of a submodule, and our operations
on cyclic codifferential forms are defined in terms of cyclization instead of restriction.

In what follows, we fix a field $\mathbbm{k}$ of characteristic zero, endowed with the trivial norm.
Fix also a graded-commutative Banach $\mathbbm{k}$-algebra $R$.

\subsection{Noncommutative Codifferential Forms} \label{subsec:codifferential-forms}

In this section, we will work with the grading group $\GG = \ZZ^2$
so that all our objects will be bigraded.
In our context, there will
be two natural parity forms on $\GG$ given by
\begin{align}
  \braid{(a_1,a_2)}{(b_1,b_2)}_1 & \defeq a_1 \cdot b_1 + a_2 \cdot b_2 \pmod{2},
  \label{eq:parity-inner-product}                                                 \\
  \braid{(a_1,a_2)}{(b_1,b_2)}_2 & \defeq (a_2 - a_1) \cdot (b_2 - b_1) \pmod{2}.
  \label{eq:parity-total-degree}
\end{align}
We will call $\braidop_1$ the \textbf{inner product parity form} and
$\braidop_2$ the \textbf{total degree parity form} for reasons to be
made clear shortly. The parity forms $\braidop_1$ and $\braidop_2$ are
the only possible ones which satisfy
\begin{equation} \label{eq:parity-extends-koszul}
  \braid{(j,0)}{(k,0)} \equiv j \cdot k \pmod{2}, \qquad
  \braid{(0,l)}{(0,m)} \equiv l \cdot m \pmod{2},
\end{equation}
i.e., extend the standard Koszul parity form on each factor.
In what follows, we will always assume
that $\braidop$ is given by either \eqref{eq:parity-inner-product} or
\eqref{eq:parity-total-degree}, and, except when writing explicit formulas, the only
property of $\braidop$ we use is \cref{eq:parity-extends-koszul}.

The parity forms $\braidop_1$ and $\braidop_2$ are
in fact equivalent (see \cref{appendix:parity-forms-equiv}) so we could have chosen to work with
either form. However, the resulting explicit formulas for various
operations are different. We find it convenient to work with both parity forms simultaneously
and state the various formulas and relations for both parity forms.

Let $R$ be a $\ZZ$-graded, graded-commutative, Banach $\mathbbm{k}$-algebra and let $V$ be a
$\ZZ$-graded Banach $R$-module.
We will identify $R$ with the
$\ZZ^2$-graded, graded-commutative, Banach $\mathbbm{k}$-algebra whose $\ZZ^2$-grading is given by
\begin{equation*}
  R^{(k,l)} \defeq
  \begin{cases}
    R^l & k = 0,    \\
    0   & k \neq 0.
  \end{cases}
\end{equation*}
Similarly, we will identify $V$ with the $\ZZ^2$-graded Banach $R$-module whose $\ZZ^2$-grading is given by
\begin{equation*}
  V^{(k,l)} \defeq
  \begin{cases}
    V^l & k = 0,    \\
    0   & k \neq 0.
  \end{cases}
\end{equation*}

Let us denote by $\ul{V}$ the shift $\ul{V} = V[(-1,0)]$. Recall that the $R$-module structure on $\ul{V}$
is defined in a way that makes the suspension map $\s_{(-1,0)} \colon V \rightharpoonup \ul{V}$
an $R$-linear map of degree $(1,0)$. Given $v \in V^d$, we will use the
notation $\ul{v} \defeq \s_{(-1,0)}(v) \in \ul{V}^{(1,d)}$. Using this notation, the $R$-module structure
on $\ul{V}$ is given by
\begin{equation} \label{eq:R-action-ul-v}
  r \ul{v} = (-1)^{\braid{(1,0)}{(0,\degb{r})}} \ul{rv} =
  \begin{cases}
    \ul{rv}                 & \braidop = \braidop_1, \\
    (-1)^{\degb{r}} \ul{rv} & \braidop = \braidop_2.
  \end{cases}
\end{equation}

Let $\ndf{V} \defeq \tens{V \oplus \ul{V}}[R]^{(*,*)}$ be the formal tensor module on $V \oplus \ul{V}$.
Elements of $\ndf{V}[][]$ will be called \textbf{(noncommutative) codifferential forms}.
An element $x \in \ndf{V}[n][m]$ which is an elementary tensor has the form
\begin{equation} \label{eq:elementary-codifferential-form}
  x = l^0 \otimes \ul{v_1} \otimes l^1 \otimes \dots \otimes
  \ul{v_n} \otimes l^{n}
\end{equation}
where $l^0,\dots,l^{n} \in \tens{V}$ are elementary tensors, $v_1, \dots, v_n \in V$ and
\begin{equation*}
  \sum_{i=0}^{n} \degb{l^i} + \sum_{i=1}^n \degb{v_i} = m.
\end{equation*}
We will call $n$ the \textbf{line degree} of $x$, since it corresponds to the number of underlined elements in $x$,
and $m$ the \textbf{cohomological degree} of $x$. Such elementary tensors will be called
\textbf{elementary codifferential forms}, and a general codifferential form is a possibly infinite
sum of elementary codifferential forms.
The reduced tensor module
on $V \oplus \ul{V}$ will be denoted by $\ndfr{V} \defeq \tensr{V \oplus \ul{V}}[R]^{(*,*)}$,
so that $\ndfr{V}[0][] = \ndf{V}[0][] / R = \tensr{V}[R]$ while $\ndfr{V}[n][] = \ndf{V}[n][]$ for $n > 0$.

Replacing the tensor module with the cyclic tensor module, we set
\begin{equation*}
  \ncdf{V} \defeq \tenscyc{V \oplus \ul{V}}[(*,*)] = \ndf{V} / \Im \left( \idd - \t \right).
\end{equation*}
Elements of $\ncdf{V}[n][m]$ will be called \textbf{(noncommutative) cyclic codifferential forms} of
line degree $n$ and cohomological degree $m$.
We will not use an equivalence class notation for the elements of $\ncdf{V}[][]$ and
continue to denote them as usual. Note that by applying the rotation operator $\t$ repeatedly and hiding
the sign inside the tensor product, we can assume that the equivalence class of an elementary codifferential form
$x \in \ncdf{V}[n][]$, i.e., an \textbf{elementary cyclic codifferential form}, has the form
\begin{equation}
  x = \ul{v_1} \otimes l^1 \otimes \ul{v_2} \otimes l^2 \otimes \dots \otimes
  \ul{v_n} \otimes l^{n} \label{eq:representative-starts-with-underline}
\end{equation}
when $n \geq 1$. That is, elementary cyclic codifferential forms in $\ncdf{V}[n][]$ can be written so that
they always start with an underlined element $\ul{v_1}$ and end with either an underlined element or an
element from $\tens{V}$ (compare to \cref{eq:elementary-codifferential-form}).
We will always assume that such elements in $\ncdf{V}[n][]$ are written this way.
The reduced cyclic tensor module on $V \oplus \ul{V}$ will be denoted by
$\ncdfr{V} \defeq \tensrcyc{V \oplus \ul{V}}^{(*,*)}$,
so that $\ncdfr{V}[0][] = \ncdf{V}[0][] / R = \tensrcyc{V}$ while $\ncdfr{V}[i][] = \ncdf{V}[i][]$ for $i > 0$.

\begin{rem}
  The bigraded Banach $R$-modules $\ndf{V}[][]$ and $\ncdf{V}[][]$ depend on the choice of the parity form
  we work with. The dependence enters through the definition of the $R$-action on the shifted
  module $\ul{V}$, through the tensor product $\otimes_R$ and through the action of the rotation operator
  $\t$ on $\tens{V \oplus \ul{V}}$, whose accompanying sign depends on the parity form.
  See \cref{sec:dependence-ndf-braidop} for full details.
\end{rem}

\begin{rem}
  Noncommutative (cyclic) codifferential forms have non-negative line degree so that the
  ``effective'' grading monoid we work with is $\NZ \times \ZZ$. One can generalize the results of this
  section by allowing $V$ to be $\GG$-graded and equipped with a parity form $\braidop_{\GG}$. Then
  codifferential forms will be $\NZ \times \GG$ graded and the symmetry we work with will be given by
  \begin{equation*}
    \braid{(i,g)}{(j,h)} \defeq i \cdot j + \braid{g}{h}_{\GG} \pmod{2},
  \end{equation*}
  for $i,j \in \NZ$ and $g,h \in \GG$, which is the natural generalization of the inner product parity form $\braidop_1$.
\end{rem}

\subsection{The de Rham Differential} \label{subsec:de-rham-differential}
Let us define a coderivation $\qdr \colon \tens{V \oplus \ul{V}} \rightharpoonup \tens{V \oplus \ul{V}}$
by the following commutative diagram:
\begin{figure}[H]
  \begin{tikzcd}
    \tens{V \oplus \ul{V}} \arrow[harpoon]{r}{\qdr}
    \arrow{d}[swap]{\pi_1} & \tens{V \oplus \ul{V}} \arrow[d, "\pi_1"] \\
    V \oplus \ul{V} \arrow[harpoon]{r}{v \mapsto 0}[swap]{\ul{v} \mapsto v} & V \oplus \ul{V}
  \end{tikzcd}
  \caption{Definition of the de Rham differential.}
  \label{fig:def-q-de-rham}
\end{figure}
Namely, $\qdr$ is a coderivation whose corestriction
$\corest{\qdr} = \pi_1 \circ \qdr \colon \tens{V \oplus \ul{V}} \rightharpoonup V \oplus \ul{V}$
has only one component $\qdr_1 \colon V \oplus \ul{V} \rightharpoonup V \oplus \ul{V}$ given by
$\qdr_1(\ul{v}) = v$ and $\qdr_1(v) = 0$ for $v \in V$. In particular, we have
\begin{equation} \label{eq:corest-qdr}
  \corest{\qdr} = \pi_1 \circ \qdr = \qdr_1 \circ \pi_1.
\end{equation}
The coderivation $\qdr$ is $R$-linear of degree $(-1,0)$, satisfies $\nnorm[\qdr] \leq 1$, and preserves
the weight of elements in $\tens{V \oplus \ul{V}}$.

\begin{lm} \label{lm:qdr-squared-zero}
  We have $\qdr \circ \qdr = 0$.
\end{lm}
\begin{proof}
  Since $\braid{(-1,0)}{(-1,0)} \equiv 1 \pmod{2}$ by \eqref{eq:parity-extends-koszul}, the coderivation $\qdr$ is odd and hence
  \begin{equation*}
    \left[ \qdr, \qdr \right] = \qdr \circ \qdr - (-1)^{\braid{(-1,0)}{(-1,0)}} \qdr \circ \qdr =
    2 \left( \qdr \circ \qdr \right).
  \end{equation*}
  In particular, this means that $\qdr \circ \qdr$ is also a coderivation
  and so it is enough to show that its corestriction
  $\corest{ \left( \qdr \circ \qdr \right)} = \corest{\qdr} \circ \qdr$
  is zero. Since $\qdr_1 \circ \qdr_1 = 0$, we have
  \begin{equation*}
    \corest{\left( \qdr \circ \qdr \right)}(x) = \left( \corest{\qdr} \circ \qdr \right)(x)
    \stackrel{\eqref{eq:corest-qdr}}{=}
    \qdr_1 \left( \corest{\qdr} \left( x \right) \right)
    \stackrel{\eqref{eq:corest-qdr}}{=}
    \left( \qdr_1 \circ \qdr_1 \right) \left( \pi_1 \left( x \right) \right) = 0.
  \end{equation*}
\end{proof}
The operator $\qdr$ is called the \textbf{de Rham differential}.
Explicitly, using \cref{eq:generalized-coder-coextension}, we see that the action of $\qdr$
on $\ndf{V}[n][]$ is given by
\begin{equation} \label{eq:qdr-formula}
  \begin{gathered} \qdr \left( l^0 \otimes \ul{v_1} \otimes l^1 \otimes
    \dots \otimes \ul{v_n} \otimes l^{n} \right) = \\
    \sum_{i=1}^{n} (-1)^{\alpha_i} \, l^0 \otimes \ul{v_1} \otimes l^1 \otimes \dots \otimes
    \ul{v_{i-1}} \otimes l^{i-1} \otimes v_i \otimes l^{i} \otimes \ul{v_{i+1}}
    \otimes \dots \otimes \ul{v_n} \otimes l^{n},
  \end{gathered}
\end{equation}
where the sign $\alpha_i$ depends on the parity form and is given by
\begin{equation} \label{eq:qdr-formula-sign}
  \begin{aligned}
    \alpha_i & =
    \braid{(-1,0)}
    {\left( i - 1, \degb{l^0 \otimes v_1 \otimes l^1 \otimes \dots \otimes v_{i-1} \otimes l^{i-1}} \right)}
    \\
             & =
    \begin{cases}
      i - 1                                                                      & \braidop = \braidop_1, \\
      i - 1 + \degb{l^0} + \sum_{j=1}^{i-1} \left(\degb{v_j} + \degb{l^j}\right) &
      \braidop = \braidop_2.
    \end{cases}
  \end{aligned}
\end{equation}
Since the only non-zero component of $\qdr$ is of weight one,
we have $\cycl{\qdr} = \qdr$ and so
$\qdr \colon \ndf{V}[][] \rightharpoonup \ndf{V}[][]$ descends to an operator
$\qdr \colon \ncdf{V}[][] \rightharpoonup \ncdf{V}[][]$ on cyclic codifferential forms.

\subsection{Lie Derivative} \label{subsec:lie-derivative}
\begin{lm} \label{lm:existence-uniqueness-lie}
  Let $\mu \colon \tens{V} \rightharpoonup \tens{V}$ be a generalized coderivation over a graded
  Banach $\mathbbm{k}$-algebra derivation $d \colon R \rightharpoonup R$. Then
  there exists a unique generalized coderivation
  $\lie{\mu} \colon \tens{V \oplus \ul{V}} \rightharpoonup \tens{V \oplus \ul{V}}$
  of degree $(0, \degb{\mu})$ over $d$, called the \textbf{Lie derivative} along $\mu$,
  which satisfies
  \begin{equation} \label{eq:lie-derivative-extends-mu}
    \rest{\lie{\mu}}{\tens{V}} = \mu
  \end{equation}
  and
  \begin{equation} \label{eq:qliecomm}
    \left[ \qdr, \lie{\mu} \right] = 0.
  \end{equation}
  The corestriction of the coderivation $\lie{\mu}$ is given by
  \begin{equation} \label{eq:lie-der-corest}
    \corest{\left( \lie{\mu} \right)} \left( x \right) \defeq
    \begin{cases}
      \corest{\mu} \left( x \right)                          & x \in \ndf{V}[0][],                 \\
      (-1)^{\braid{(-1,0)}{(0, \degb{\mu})}}
      \ul{\corest{\mu} \left( \qdr \left( x \right) \right)} & x \in \ndf{V}[1][],                 \\
      0                                                      & x \in \ndf{V}[n][], \quad n \geq 2.
    \end{cases}
  \end{equation}
\end{lm}
\begin{proof}
  Assume such a coderivation $\lie{\mu}$ exists. Since $\lie{\mu}$ has line degree zero, $\lie{\mu}$ induces
  well-defined maps $\rest{\lie{\mu}}{\ndf{V}[n][]} \colon \ndf{V}[n][*] \rightharpoonup \ndf{V}[n][* + \degb{\mu}]$
  for all $n \geq 0$.
  Let $x_n \in \ndf{V}[n][]$.
  \begin{enumerate}
    \item \Cref{eq:lie-derivative-extends-mu} forces that
          $\corest{ \left( \lie{\mu} \right)} \left( x_0 \right) = \corest{\mu} \left( x_0 \right)$.
    \item \Cref{eq:qliecomm} forces that
          \begin{equation*}
            0 = \corest{\left[ \qdr, \lie{\mu} \right]} \left( x_1 \right) =
            \left( \corest{\qdr} \circ \lie{\mu} -
            (-1)^{\braid{(-1,0)}{(0,\degb{\mu})}} \corest{ \left( \lie{\mu} \right)} \circ \qdr \right)
            \left( x_1 \right)
          \end{equation*}
          and hence
          \begin{equation*}
            \begin{aligned}
              \left( \qdr_1 \circ \corest{ \left( \lie{\mu} \right)} \right) \left( x_1 \right)
              \stackrel{\eqref{eq:corest-qdr}}{=}{}           &
              \left( \corest{\qdr} \circ \lie{\mu} \right) \left( x_1 \right) =
              (-1)^{\braid{(-1,0)}{(0,\degb{\mu})}} \left( \corest{ \left( \lie{\mu} \right)} \circ \qdr \right)
              \left( x_1 \right)
              \\
              \stackrel{\phantom{\eqref{eq:corest-qdr}}}{=}{} &
              (-1)^{\braid{(-1,0)}{(0,\degb{\mu})}} \corest{\mu} \left( \qdr \left( x_1 \right) \right).
            \end{aligned}
          \end{equation*}
          Since $\corest{ \left( \lie{\mu} \right)} \left( x_1 \right)$ has line degree one,
          i.e., $\corest{ \left( \lie{\mu} \right)} \left( x_1 \right) \in \ul{V}$, we have
          \begin{equation} \label{eq:qdr-commute-lie-mu-one-iff}
            \left[ \qdr, \lie{\mu} \right] \left( x_1 \right) = 0 \iff
            \corest{ \left( \lie{\mu} \right)} \left( x_1 \right) =
            (-1)^{\braid{(-1,0)}{(0,\degb{\mu})}}
            \ul{ \corest{\mu} \left( \qdr \left( x_1 \right) \right) }.
          \end{equation}
    \item By degree considerations, any coderivation
          $\eta \colon \tens{V \oplus \ul{V}} \rightharpoonup \tens{V \oplus \ul{V}}$
          of degree $(0,*)$ must satisfy
          $\corest{ \eta } \left( x_n \right) = 0$ for $n \geq 2$. In particular,
          $\corest{ \left( \lie{\mu} \right)} \left( x_n \right) = 0$ for all $n \geq 2$.
  \end{enumerate}
  This shows that \cref{eq:lie-derivative-extends-mu,eq:qliecomm} forces the corestriction of $\lie{\mu}$
  to be given by \cref{eq:lie-der-corest}. Since $\lie{\mu}$ is determined uniquely by its corestriction,
  we obtain uniqueness.

  To show existence, we can use \cref{prop:classification-generalized-coderivations-formal-tensor-coalgebra}
  to define a generalized coderivation
  $\lie{\mu} \colon \tens{V \oplus \ul{V}} \rightharpoonup \tens{V \oplus \ul{V}}$
  by specifying its corestriction to be given by \cref{eq:lie-der-corest}.
  To show this is a legitimate definition, we need to verify that the corestriction
  $\corest{\left( \lie{\mu} \right)}$ satisfies \cref{eq:corest-over-derivation-projection}, i.e.,
  $\left( \lie{\mu} \right)_1$ must be a derivation over $d$ while the maps $\left( \lie{\mu} \right)_k$
  for $k \neq 1$ must be $R$-linear. Since $\qdr$ preserves weight,
  $\left( \lie{\mu} \right)_k$ depend only on $\mu_k$ and hence the maps $\left( \lie{\mu} \right)_k$
  are $R$-linear when $k \neq 1$. For $k = 1$, we have
  \begin{equation*}
    \begin{aligned}
      \left( \lie{\mu} \right)_1 \left( r \cdot v \right)
      \stackrel{\phantom{\eqref{eq:parity-extends-koszul}}}{=}{} &
      \mu_1 \left( r \cdot v \right)
      \\
      \stackrel{\phantom{\eqref{eq:parity-extends-koszul}}}{=}{} &
      dr \cdot v + (-1)^{\degb{\mu} \cdot \degb{r}} r \cdot \mu_1 \left( v \right)
      \\
      \stackrel{\eqref{eq:parity-extends-koszul}}{=}{}           &
      dr \cdot v +
                 (-1)^{\braid{\left( 0, \degb{\mu} \right)}{\left( 0, \degb{r} \right)}} r \cdot
      \left( \lie{\mu} \right)_{1} \left( v \right)
      \\
      \stackrel{\phantom{\eqref{eq:parity-extends-koszul}}}{=}{} &
      dr \cdot v +
                 (-1)^{\braidd{\lie{\mu}}{r}} r \cdot
      \left( \lie{\mu} \right)_{1} \left( v \right)
    \end{aligned}
  \end{equation*}
  and
  \begin{equation*}
    \begin{aligned}
      \left( \lie{\mu} \right)_{1} \left( r \cdot \ul{v} \right)
      \stackrel{\phantom{\eqref{eq:parity-extends-koszul}}}{=}{} &
      \left( \lie{\mu} \right)_{1} \left( (-1)^{\braid{(1,0)}{(0,\degb{r})}} \ul{rv} \right)
      \\
      \stackrel{\phantom{\eqref{eq:parity-extends-koszul}}}{=}{} &
      (-1)^{\braid{(1,0)}{(0,\degb{r})} + \braid{(-1,0)}{(0, \degb{\mu})}}
      \ul{\mu_1 \left( r \cdot v \right)}
      \\
      \stackrel{\phantom{\eqref{eq:parity-extends-koszul}}}{=}{} &
      (-1)^{\braid{(1,0)}{(0,\degb{r})} + \braid{(-1,0)}{(0, \degb{\mu})}}
      \left(
      \ul{dr \cdot v + (-1)^{\degb{\mu} \cdot \degb{r}} r \cdot \mu_1 \left( v \right)}
      \right)
      \\
      \stackrel{\phantom{\eqref{eq:parity-extends-koszul}}}{=}{} &
      (-1)^{\braid{(1,0)}{(0,\degb{r})} + \braid{(-1,0)}{(0, \degb{\mu})} +
        \braid{(1,0)}{(0, \degb{dr})}}
      dr \cdot \ul{v}
      \\
      \quad                                                      & +
                                                                   (-1)^{\braid{(1,0)}{(0,\degb{r})} + \braid{(-1,0)}{(0, \degb{\mu})} +
                                                                     \degb{\mu} \cdot \degb{r} + \braid{(1,0)}{(0, \degb{r})}}
      r \cdot \ul{\mu_1 \left( v \right)}
      \\
      \stackrel{\phantom{\eqref{eq:parity-extends-koszul}}}{=}{} &
      (-1)^{\braid{(-1,0)}{(0, \degb{\mu})} + \braid{(1,0)}{(0, \degb{\mu})}}
      dr \cdot \ul{v} +
               (-1)^{\braid{(-1,0)}{(0, \degb{\mu})} + \degb{\mu} \cdot \degb{r}}
      r \cdot \ul{\mu_1 \left( v \right)}
      \\
      \stackrel{\phantom{\eqref{eq:parity-extends-koszul}}}{=}{} &
      dr \cdot \ul{v} +
               (-1)^{\degb{\mu} \cdot \degb{r}}
      r \cdot \left( \lie{\mu} \right)_1 \left( \ul{v} \right)
      \\
      \stackrel{\eqref{eq:parity-extends-koszul}}{=}{}           &
      dr \cdot \ul{v} +
               (-1)^{\braid{(0, \degb{\mu})}{(0, \degb{r})}}
      r \cdot \left( \lie{\mu} \right)_1 \left( \ul{v} \right)
      \\
      \stackrel{\phantom{\eqref{eq:parity-extends-koszul}}}{=}{} &
      dr \cdot \ul{v} +
               (-1)^{\braidd{\lie{\mu}}{r}}
      r \cdot \left( \lie{\mu} \right)_{1} \left( \ul{v} \right),
    \end{aligned}
  \end{equation*}
  for $r \in R$ and $v \in V$.

  Finally, let us verify that $\lie{\mu}$ as defined above satisfies \cref{eq:lie-derivative-extends-mu,eq:qliecomm}.
  It is clear that $\lie{\mu}$ has degree $(0,\degb{\mu})$ and
  extends $\mu$ on $\tens{V}$. To show that $\left[ \qdr, \lie{\mu} \right] = 0$, it is enough to
  show that $\corest{ \left[ \qdr, \lie{\mu} \right] } = 0$. Since $\left[ \qdr, \lie{\mu} \right]$ has
  degree $(-1, \degb{\mu})$, it is clear that $\corest{ \left[ \qdr, \lie{\mu} \right] }$ vanishes
  on $\ndf{V}[n][]$ for $n \neq 1, 2$. For $x_1 \in \ndf{V}[1][]$ we have
  $\left[ \qdr, \lie{\mu} \right] \left( x_1 \right) = 0$ by construction
  (see \cref{eq:qdr-commute-lie-mu-one-iff}) and for $x_2 \in \ndf{V}[2][]$ we have
  \begin{equation*}
    \begin{aligned}
      \corest{\left[ \qdr, \lie{\mu} \right]} \left( x_2 \right)
      \stackrel{\phantom{\eqref{eq:lie-der-corest}}}{=}{} &
      \left( \corest{\qdr} \circ \lie{\mu} -
      (-1)^{\braid{(-1,0)}{(0,\degb{\mu})}} \corest{ \left( \lie{\mu} \right)} \circ \qdr \right)
      \left( x_2 \right)
      \\
      \stackrel{\eqref{eq:corest-qdr}}{=}{}               &
      \left( \qdr_1 \circ \corest{ \left( \lie{\mu} \right)} -
                          (-1)^{\braid{(-1,0)}{(0,\degb{\mu})}} \corest{ \left( \lie{\mu} \right)} \circ \qdr \right)
      \left( x_2 \right)
      \\
      \stackrel{\eqref{eq:lie-der-corest}}{=}{}           &
      - \ul{ \corest{\mu} \left( \qdr^2 \left( x_2 \right) \right)}
      \stackrel{\cref{lm:qdr-squared-zero}}{=}{} 0.
    \end{aligned}
  \end{equation*}
\end{proof}

The proof of \cref{lm:existence-uniqueness-lie} gives us an explicit formula for the corestriction of $\lie{\mu}$.
Given $x = l^0 \otimes \ul{v_1} \otimes l^1 \otimes \dots \otimes \ul{v_n} \otimes l^{n} \in \ndf{V}[n][]$,
by \cref{eq:lie-der-corest}, we have
\begin{equation}
  \corest{\left( \lie{\mu} \right)} \left( x \right) =
  \begin{cases}
    \corest{\mu} \left( l^0 \right)                              & n = 0,    \\
    (-1)^{\braid{(-1,0)}{(0, \degb{\mu} + \degb{l^0})}}
    \ul{\corest{\mu} \left( l^0 \otimes v_1 \otimes l^1 \right)} & n = 1,    \\
    0                                                            & n \geq 2.
  \end{cases}
  \label{eq:lie-derivation-corestriction}
\end{equation}
The explicit formula for $\corest{\left( \lie{\mu} \right)}$
tells us that $\lie{\mu}$ acts on any consecutive sublist of elements
which does not include an underlined element (including the empty list) by applying $\corest{\mu}$
and acts on any consecutive sublist of elements which includes \textit{only one} underlined
element by ``forgetting'' the line, applying $\corest{\mu}$ and underlining the result, with some signs.
It does not act on sublists which include two underlined elements or more (so it doesn't act ``across lines'').
The signs are determined by the usual braiding conventions. Explicitly, we have
  { \small
    \begin{align}
      \MoveEqLeft
      \lie{\mu} \left( l^0 \otimes \ul{v_1} \otimes l^1 \otimes \dots \otimes
      \ul{v_n} \otimes l^{n} \right) =
      \notag
      \\
      ={} & \sum_{i=0}^{n} (-1)^{\gamma_i} \, l^0 \otimes \ul{v_1}
      \otimes l^1 \otimes \dots \otimes \ul{v_{i}} \otimes l^{i}_{(1)} \otimes
      \corest{\mu} \left( l^{i}_{(2)} \right) \otimes l^{i}_{(3)} \otimes
      \ul{v_{i+1}} \otimes \dots \otimes \ul{v_n} \otimes l^{n} {}+{}
      \label{eq:lie-mu-full-formula}
      \\
          &
      \sum_{i=1}^{n} (-1)^{\delta_i} \, l^0 \otimes \ul{v_1} \otimes l^1 \otimes \dots \otimes
      \ul{v_{i-1}} \otimes l^{i-1}_{(1)} \otimes
      \ul{ \corest{\mu} \left( l^{i-1}_{(2)} \otimes v_{i} \otimes l^{i}_{(1)} \right)} \otimes
      l^{i}_{(2)} \otimes \ul{v_{i+1}} \otimes
      \dots \otimes \ul{v_n} \otimes l^{n} \notag
    \end{align}
  }
where\footnote{There is some abuse of notation involved here, as the sign $\gamma_i$
  depends not only on $i$, but also on the splitting of $l^i$ into three lists as evident from the formula
  for $\gamma_i$. We hope the intention is clear. The same remark applies with the obvious modifications
  to all the other signs we write.}
\begin{equation*}
  \begin{aligned}
    \gamma_i & =
    \braid{\left( 0,\degb{\mu} \right)}
    {\left( i, \degb{ l^0 \otimes v_1
        \otimes l^1 \otimes \dots \otimes v_{i} \otimes l^{i}_{(1)}} \right)}
    \\
             & =
    \begin{cases}
      \degb{\mu} \cdot \left( \sum_{j=0}^{i-1} \left( \degb{l^j} + \degb{v_{j+1}} \right) +
      \degb{l^{i}_{(1)}} \right)
       & \braidop = \braidop_1, \\
      \degb{\mu} \cdot \left( \sum_{j=0}^{i-1} \left( \degb{l^j} + \degb{v_{j+1}} \right) +
      \degb{l^{i}_{(1)}} - i \right)
       & \braidop = \braidop_2
    \end{cases}
  \end{aligned}
\end{equation*}
and
\begin{equation*}
  \begin{aligned}
    \delta_i & =
    \braid{\left( 0,\degb{\mu} \right)}
    {\left( i - 1, \degb{ l^0 \otimes v_1
        \otimes l^1 \otimes \dots \otimes v_{i-1} \otimes l^{i-1}_{(1)}} \right)} +
    \braid{(-1,0)}{(0,\degb{\mu} + \degb{l^{i-1}_{(2)}})}
    \\
             & =
    \begin{cases}
      \degb{\mu} \cdot \left( \sum_{j=0}^{i-2} \left( \degb{l^j} + \degb{v_{j+1}} \right) +
      \degb{l^{i-1}_{(1)}} \right)
                                                              & \braidop = \braidop_1, \\
      \degb{\mu} \cdot \left( \sum_{j=0}^{i-2} \left( \degb{l^j} + \degb{v_{j+1}} \right) +
      \degb{l^{i-1}_{(1)}} + i \right)	+ \degb{l^{i-1}_{(2)}} &
      \braidop = \braidop_2.
    \end{cases}
  \end{aligned}
\end{equation*}
From the explicit formula \eqref{eq:lie-derivation-corestriction} for $\corest{\left( \lie{\mu} \right)}$,
it is also clear that $\lie{\mu}$ is bounded
with $\nnorm[\lie{\mu}] \leq \nnorm[\mu]$, and since $\lie{\mu}$ extends $\mu$ on $\tens{V}$, we have
$\nnorm[\lie{\mu}] = \nnorm[\mu]$.

Let us denote the cyclization of $\lie{\mu}$ by $\clie{\mu} \defeq \cycl{\left( \lie{\mu} \right)}$.
The operator
\begin{equation*}
  \clie{\mu} \colon \ndf{V}[n][*] \rightharpoonup \ndf{V}[n][* + \degb{\mu}]
\end{equation*}
is called the \textbf{cyclic Lie derivative} along $\mu$, and is a bounded module derivation over $d$ with
norm $\nnorm[\clie{\mu}] \leq \nnorm[\mu]$.
By \cref{lm:coder-descends-quotient},
$\clie{\mu}$ descends to a well-defined operator
\begin{equation*}
  \clie{\mu} \colon \ncdf{V}[n][*] \rightharpoonup \ncdf{V}[n][* + \degb{\mu}]
\end{equation*}
on cyclic codifferential forms.
We can describe explicitly the action of
$\clie{\mu}$ on an elementary codifferential form $x \in \ndf{V}[n][]$ as follows.
As usual, the difference between $\clie{\mu}$ and $\lie{\mu}$ is that
$\clie{\mu}$ allows elements from the end of the list to be rotated to the
beginning of the list. When $x \in \ndf{V}[0][]$ then
\begin{equation}
  \clie{\mu} \left( x \right) = \cycl{\mu} \left( x \right),
\end{equation}
i.e., $\clie{\mu}$ extends the action of $\cycl{\mu}$ on $\ndf{V}[0][]$.
When $x \in \ndf{V}[1][]$ is an elementary codifferential form which starts with
an element from $\ul{V}$, i.e., $x = \ul{v} \otimes l$, we have
\begin{equation} \label{eq:cyc-lie-formula-n-eq-1}
  \clie{\mu} \left( \ul{v} \otimes l \right) =
  (-1)^{\gamma_1}
  \ul{v} \otimes l_{(1)} \otimes \corest{\mu} \left( l_{(2)} \right) \otimes l_{(3)}
    +
    (-1)^{\varepsilon}
  \ul{\corest{\mu} \left( l_{(3)} \otimes v \otimes l_{(1)} \right)} \otimes l_{(2)},
\end{equation}
where $\gamma_1$ is as before, and the sign $\varepsilon$ involves the sign from permuting
the factors using $\t$ and is given explicitly by
\begin{equation*}
  \begin{aligned}
    \varepsilon ={} & \braid{\left(0, \degb{l_{(3)}} \right)}{ \left( 1, \degb{v} + \degb{l_{(1)}} + \degb{l_{(2)}} \right)}
    +
    \braid{(-1,0)}{\left( 0, \degb{\mu} + \degb{l_{(3)}} \right)}
    \\
    ={}             &
    \begin{cases}
      \degb{l_{(3)}} \cdot \left( \degb{v} + \degb{l_{(1)}} + \degb{l_{(2)}} \right)              & \braidop = \braidop_1,
      \\
      \degb{l_{(3)}} \cdot \left( \degb{v} + \degb{l_{(1)}} + \degb{l_{(2)}} \right) + \degb{\mu} & \braidop = \braidop_2.
    \end{cases}
  \end{aligned}
\end{equation*}
When $x \in \ndf{V}[n][]$ is an elementary codifferential form which starts with an element from $\ul{V}$ and $n \geq 2$,
i.e., $x = \ul{v_1} \otimes l^1 \otimes \dots \otimes \ul{v_n} \otimes l^n$,
then
\begin{gather}
  \clie{\mu} \left( \ul{v_1} \otimes l^1 \otimes \dots \otimes \ul{v_n} \otimes l^{n} \right) =
  \notag
  \\
  \sum_{i=1}^{n} (-1)^{\gamma_i} \, \ul{v_1}
  \otimes l^1 \otimes \dots \otimes \ul{v_{i}} \otimes l^{i}_{(1)} \otimes
  \corest{\mu} \left( l^{i}_{(2)} \right) \otimes l^{i}_{(3)} \otimes
  \ul{v_{i+1}} \otimes \dots \otimes \ul{v_n} \otimes l^{n} +
  \notag
  \\
  \sum_{i=2}^{n} (-1)^{\delta_i} \, \ul{v_1} \otimes l^1 \otimes \dots \otimes
  \ul{v_{i-1}} \otimes l^{i-1}_{(1)} \otimes
  \ul{ \corest{\mu} \left( l^{i-1}_{(2)} \otimes v_{i} \otimes l^{i}_{(1)} \right)} \otimes
  l^{i}_{(2)} \otimes \ul{v_{i+1}} \otimes \dots \otimes \ul{v_n} \otimes l^{n} +
  \notag
  \\
  (-1)^{\varepsilon}
  \ul{ \corest{\mu} \left( l^n_{(2)} \otimes v_1 \otimes l^1_{(1)} \right)} \otimes
  l^1_{(2)} \otimes \dots \otimes l^{n-1} \otimes \ul{v_n} \otimes
  l^n_{(1)} \label{eq:cyc-lie-formula}
\end{gather}
where the signs $\gamma_i,\delta_i$ are as before (with the understanding that $l^0 = 1$), and the
sign $\varepsilon$ involves the sign from permuting the factors using $\t$ and is given explicitly by
\begin{equation*}
  \begin{aligned}
    \varepsilon & = \braid{\left( 0, \degb{l^n_{(2)}} \right)}
                    {\left( n, \degb{ v_1 \otimes l^1 \otimes \dots \otimes v_{n} \otimes l^n_{(1)} } \right)}
    +
    \braid{\left(-1,0\right)}{\left(0,\degb{\mu} + \degb{l^n_{(2)}}\right)}
    \\
                & =
    \begin{cases}
      \degb{l^n_{(2)}} \cdot \left( \sum_{j=1}^{n-1} \left( \degb{v_j} + \degb{l^j} \right) +
      \degb{v_n} + \degb{l^n_{(1)}} \right)                                     & \braidop = \braidop_1,
      \\
      \degb{l^n_{(2)}} \cdot \left( \sum_{j=1}^{n-1} \left( \degb{v_j} + \degb{l^j} \right) +
      \degb{v_n} + \degb{l^n_{(1)}} - \left( n - 1 \right) \right) + \degb{\mu} &
      \braidop = \braidop_2.
    \end{cases}
  \end{aligned}
\end{equation*}
Finally, when $x \in  \ndf{V}[n][]$ is an elementary codifferential form which
starts with an element from $V$ and $n \geq 1$, i.e.,
$x = v \otimes l^0 \otimes \ul{v_1} \otimes l^1 \otimes \dots \otimes \ul{v_n} \otimes l^n$,
we have
\begin{equation*}
  \begin{gathered}
    \clie{\mu} \left( v \otimes l^0 \otimes \ul{v_1} \otimes l^1 \otimes \dots \otimes
    l^{n-1} \otimes \ul{v_n} \otimes l^n \right) =
    \\
    \sum_{i=0}^{n} \pm \, v \otimes l^0 \otimes \ul{v_1} \otimes l^1 \otimes \dots \otimes
    \ul{v_{i}} \otimes l^{i}_{(1)} \otimes
    \corest{\mu} \left( l^{i}_{(2)} \right) \otimes l^{i}_{(3)} \otimes
    \ul{v_{i+1}} \otimes \dots \otimes \ul{v_n} \otimes l^{n} +
    \\
    \sum_{i=1}^{n} \pm \, v \otimes l^0 \otimes \ul{v_1} \otimes l^1 \otimes \dots \otimes
    \ul{v_{i-1}} \otimes l^{i-1}_{(1)} \otimes
    \ul{ \corest{\mu} \left( l^{i-1}_{(2)} \otimes v_{i} \otimes l^{i}_{(1)} \right)} \otimes
    l^{i}_{(2)} \otimes \ul{v_{i+1}} \otimes \dots \otimes \ul{v_n} \otimes l^{n} \pm
    \\
    \corest{\mu} \left( l^n_{(2)} \otimes v \otimes l^0_{(1)} \right) \otimes
    l^0_{(2)} \otimes \ul{v_1} \otimes l^1 \otimes \dots \otimes l^{n-1} \otimes
    \ul{v_n} \otimes l^n_{(1)} \pm
    \\
    \ul{ \corest{\mu} \left( l^n_{(2)} \otimes v \otimes l^0 \otimes v_1 \otimes l^1_{(1)} \right) } \otimes
    l^1_{(2)} \otimes \ul{v_2} \otimes \dots \otimes l^{n-1} \otimes \ul{v_n} \otimes l^n_{(1)} \pm
    \\
    \ul{ \corest{\mu} \left( l^{n-1}_{(2)} \otimes v_n \otimes l^n \otimes v \otimes l^0_{(1)} \right)} \otimes
    l^0_{(2)} \otimes \ul{v_1} \otimes \dots \otimes \ul{v_{n-1}} \otimes l^{n-1}_{(1)}.
  \end{gathered}
\end{equation*}
Since we will work with $\clie{\mu}$ only on the quotient $\ncdf{V}[][]$, we won't
use the equation above and don't write down the signs explicitly. Finally, note that we have
\begin{equation} \label{eq:lie-clie-1}
  \lie{\mu} \left( 1 \right) = \clie{\mu} \left( 1 \right) = \mu_0(1).
\end{equation}

\subsection{Contraction} \label{sec:nc-d-calc-cont}
\begin{lm} \label{lm:existence-uniqueness-contraction}
  Let $\mu \colon \tens{V} \rightharpoonup \tens{V}$ be an $R$-linear coderivation. Then
  there exists a unique $R$-linear coderivation
  $\cont{\mu} \colon \tens{V \oplus \ul{V}} \rightharpoonup \tens{V \oplus \ul{V}}$
  of degree $(1, \degb{\mu})$, called the \textbf{contraction} along $\mu$,
  which satisfies
  \begin{equation} \label{eq:lie-derivative-commutator-contraction}
    \left[ \qdr, \cont{\mu} \right] = \lie{\mu}.
  \end{equation}
  The corestriction of the coderivation $\cont{\mu}$ is given by
  \begin{equation} \label{eq:corest-cont-mu}
    \corest{\left( \cont{\mu} \right)} \left( x \right) =
    \begin{cases}
      \ul{ \corest{\mu} \left( x \right) } & x \in \ndf{V}[0][],                 \\
      0                                    & x \in \ndf{V}[n][], \quad n \geq 1.
    \end{cases}
  \end{equation}
\end{lm}
\begin{proof}
  Assume such a coderivation $\cont{\mu}$ exists and let $x_n \in \ndf{V}[n][]$.
  Since $\cont{\mu}$ has line degree one, we must have
  $\corest{ \left( \cont{\mu} \right) } \left( x_n \right) = 0$ for all $n \geq 1$. When $n = 0$,
  \cref{eq:lie-derivative-commutator-contraction} forces that
  \begin{equation*}
    \begin{aligned}
      \corest{\mu} \left( x_0 \right)
      \stackrel{\eqref{eq:lie-derivation-corestriction}}{=} &
      \corest{ \left( \lie{\mu} \right)} \left( x_0 \right)
      \stackrel{\eqref{eq:lie-derivative-commutator-contraction}}{=}
      \corest{\left[ \qdr, \cont{\mu} \right]} \left( x_0 \right) =
      \left( \corest{\qdr} \circ \cont{\mu} -
      (-1)^{\braid{(-1,0)}{(1,\degb{\mu})}} \corest{ \left( \cont{\mu} \right)} \circ \qdr \right)
      \left( x_0 \right)
      \\
      \stackrel{\phantom{\eqref{eq:corest-qdr}}}{=}         &
      \left( \corest{\qdr} \circ \cont{\mu} \right) \left( x_0 \right)
      \stackrel{\eqref{eq:corest-qdr}}{=}
      \left( \qdr_1 \circ \corest{ \left( \cont{\mu} \right) } \right) \left( x_0 \right).
    \end{aligned}
  \end{equation*}
  Since $\cont{\mu}$ has line degree one, we must have
  $\corest{ \left( \cont{\mu} \right)} \left( x_0 \right) \in \ul{V}$ and hence
  \begin{equation} \label{eq:commutator-qdr-cont-l-mu-zero-iff}
    \left[ \qdr, \cont{\mu} \right] \left( x_0 \right) = \lie{\mu} \left( x_0 \right) \iff
    \corest{ \left( \cont{\mu} \right)} \left( x_0 \right) = \ul{ \corest{\mu} \left( x_0 \right) }.
  \end{equation}
  This shows that \cref{eq:lie-derivative-commutator-contraction} forces the corestriction of $\cont{\mu}$
  to be given by \cref{eq:corest-cont-mu}. Since $\cont{\mu}$ is determined uniquely by its corestriction,
  we obtain uniqueness.

  To show existence, we can use \cref{prop:classification-generalized-coderivations-formal-tensor-coalgebra}
  and define an $R$-linear coderivation
  $\cont{\mu} \colon \tens{V \oplus \ul{V}} \rightharpoonup \tens{V \oplus \ul{V}}$
  of degree $(1, \degb{\mu})$ by specifying its corestriction to be given by \cref{eq:corest-cont-mu}.
  To verify the resulting $\cont{\mu}$ satisfies \cref{eq:lie-derivative-commutator-contraction},
  it is enough to show that $\corest{ \left[ \qdr, \cont{\mu} \right] } = \corest{ \left( \lie{\mu} \right) }$.
  For $x_0 \in \ndf{V}[0][]$ we have
  $\left[ \qdr, \cont{\mu} \right] \left( x_0 \right) = \lie{\mu} \left( x_0 \right)$ by construction
  (see \cref{eq:commutator-qdr-cont-l-mu-zero-iff}) and for $x_1 \in \ndf{V}[1][]$ we have
  \begin{equation*}
    \begin{aligned}
      \corest{\left[ \qdr, \cont{\mu} \right]} \left( x_1 \right)
      \stackrel{\phantom{\eqref{eq:lie-der-corest}}}{=}{} &
      \left( \corest{\qdr} \circ \cont{\mu} -
      (-1)^{\braid{(-1,0)}{(1,\degb{\mu})}} \corest{ \left( \cont{\mu} \right)} \circ \qdr \right)
      \left( x_1 \right)
      \\
      \stackrel{\eqref{eq:corest-qdr}}{=}{}               &
      \left( \qdr_1 \circ \corest{ \left( \cont{\mu} \right)} -
                          (-1)^{\braid{(-1,0)}{(1,\degb{\mu})}} \corest{ \left( \cont{\mu} \right)} \circ \qdr \right)
      \left( x_1 \right)
      \\
      \stackrel{\eqref{eq:corest-cont-mu}}{=}{}           &
      (-1)^{\braid{(-1,0)}{(1,\degb{\mu})} + 1}
      \ul{ \corest{\mu} \left( \qdr \left( x_1 \right) \right)}
      \\
      \stackrel{\eqref{eq:parity-extends-koszul}}{=}{}    &
      (-1)^{\braid{(-1,0)}{(0, \degb{\mu})}}
      \ul{ \corest{\mu} \left( \qdr \left( x_1 \right) \right)}
      \\
      \stackrel{\eqref{eq:lie-der-corest}}{=}{}           &
      \corest{ \left( \lie{\mu} \right)} \left( x_1 \right).
    \end{aligned}
  \end{equation*}
  Finally, since both $\left[ \qdr, \cont{\mu} \right]$ and $\lie{\mu}$ are of degree $(0, \degb{\mu})$, it is
  clear their corestrictions vanish, and in particular coincide, on $\ndf{V}[n][]$ for $n \geq 2$.
\end{proof}

The explicit formula for $\corest{\left( \cont{\mu} \right)}$
tells us that $\cont{\mu}$ acts on any consecutive list of elements
which does not include an underlined element (including the empty list)
by applying $\corest{\mu}$ and underlining the result, with some signs.
It doesn't act on a list of elements which contains one or more underlined
elements so we can think of the underlined elements as ``barriers''.
The signs are determined by the usual braiding convention. Explicitly, we have

\begin{gather} \label{eq:cont-formula}
  \cont{\mu} \left( l^0 \otimes \ul{v_1} \otimes l^1 \otimes \dots \otimes \ul{v_n} \otimes l^{n} \right) =
  \\
  \sum_{i=0}^{n} (-1)^{\beta_i} \,
  l^0 \otimes \ul{v_1} \otimes l^1 \otimes \dots \otimes \ul{v_{i}} \otimes l^{i}_{(1)} \otimes
  \ul{ \corest{\mu} \left( l^{i}_{(2)} \right)} \otimes l^{i}_{(3)} \otimes
  \ul{v_{i+1}} \otimes l^{i+1} \otimes \dots \otimes \ul{v_n} \otimes l^{n} \notag
\end{gather}
where the sign $\beta_i$ depends on the parity form via \cref{eq:generalized-coder-coextension} and is given by
\begin{equation} \label{eq:cont-formula-signs}
  \begin{aligned}
    \beta_i & =
    \braid{\left( 1,\degb{\mu} \right)}
    {\left( i, \degb{ l^0 \otimes v_1
        \otimes l^1 \otimes \dots \otimes v_{i} \otimes l^{i}_{(1)}} \right)}
    \\
            & =
    \begin{cases}
      i + \degb{\mu} \cdot
      \left( \sum_{j=0}^{i-1} \left(\degb{l^j} + \degb{v_{j+1}}\right) + \degb{l^i_{(1)}} \right)
       & \braidop = \braidop_1,
      \\
      \left( \degb{\mu} - 1 \right) \cdot
      \left( \sum_{j=0}^{i-1} \left(\degb{l^j} + \degb{v_{j+1}}\right) + \degb{l^i_{(1)}} - i
      \right)
       & \braidop = \braidop_2.
    \end{cases}
  \end{aligned}
\end{equation}
From the explicit formula \eqref{eq:corest-cont-mu} for $\corest{\left( \cont{\mu} \right)}$,
and the fact that $v \mapsto \ul{v}$ is an isometry, it is clear that
$\cont{\mu}$ is bounded and $\nnorm[\cont{\mu}] = \nnorm[\mu]$.

Let us denote the cyclization of the coderivation $\cont{\mu}$ by
$\ccont{\mu} \defeq \cycl{ \left( \cont{\mu} \right)}$.
The operator
\begin{equation*}
  \ccont{\mu} \colon \ndf{V}[n][*] \rightharpoonup \ndf{V}[n+1][* + \degb{\mu}]
\end{equation*}
is called
the \textbf{cyclic contraction} along $\mu$, and is a bounded $R$-linear operator with norm
$\nnorm[\ccont{\mu}] \leq \nnorm[\mu]$. By \cref{lm:coder-descends-quotient},
$\ccont{\mu}$ descends to a well-defined operator
\begin{equation*}
  \ccont{\mu} \colon \ncdf{V}[n][*] \rightharpoonup \ncdf{V}[n+1][* + \degb{\mu}]
\end{equation*}
on cyclic codifferential forms.
Verbally, the cyclic contraction $\ccont{\mu}$ acts on any
consecutive list of elements (including the empty list) which do
not include an underlined element by applying $\corest{\mu}$ and underlining the result, allowing
elements from the end of the list to be rotated to the beginning of the list.
More explicitly, let $x \in \ndf{V}[n][]$ be an elementary tensor.
When $x \in \ndf{V}[0][]$ has the form
$x = v \otimes l$ for $v \in V$ and $l \in \tens{V}$, we have
\begin{equation} \label{eq:ccont-mu-on-ncdf-0}
  \begin{aligned}
    \ccont{\mu} \left( v \otimes l \right) ={} &
    (-1)^{\braid{(1, \degb{\mu})}{(0, \degb{v} + \degb{l_{(1)}})}}
    v \otimes l_{(1)} \otimes \ul{ \corest{\mu} \left( l_{(2)} \right) } \otimes l_{(3)}
    \\
                                               & +
                                                 (-1)^{\braid{(0, \degb{l_{(3)}})}{(0, \degb{v} + \degb{l_{(1)}} + \degb{l_{(2)}})}}
    \ul{ \corest{\mu} \left( l_{(3)} \otimes v \otimes l_{(1)} \right) } \otimes
    l_{(2)}.
  \end{aligned}
\end{equation}
When $x$ has the form
$x = \ul{v_1} \otimes l^1 \otimes \dots \otimes \ul{v_n} \otimes l^{n}$,
that is, $x \in \ndf{V}[n][]$ for $n \geq 1$ and starts with an element from $\ul{V}$,
the cyclization $\ccont{\mu} \left( x \right)$ doesn't contain any terms which involve rotation and the sole
difference between $\ccont{\mu}$ and $\cont{\mu}$ is an extra term involving $\ul{ \mu_0 \left( 1 \right) }$:
\begin{align}\label{eq:ccont-mu-on-ncdf-geq-1}
  \ccont{\mu} \left( x \right)
   & =
  \cont{\mu} \left( x \right) - \ul{ \mu_0 \left( 1 \right)} \otimes x
  \\
   & =
  \sum_{i=1}^{n} (-1)^{\beta_i} \,
  \ul{v_1} \otimes l^1 \otimes \dots \otimes \ul{v_{i}} \otimes l^{i}_{(1)} \otimes
  \ul{ \corest{\mu} \left( l^{i}_{(2)} \right)} \otimes l^{i}_{(3)} \otimes
  \ul{v_{i+1}} \otimes l^{i+1} \otimes \dots \otimes \ul{v_n} \otimes l^{n}. \notag
\end{align}
On the other hand, when $x$ has the form $x = v_0 \otimes l^0 \otimes \ul{v_1}
  \otimes l^1 \otimes \dots \otimes \ul{v_n} \otimes l^{n}$, that is, $x \in \ndf{V}[n][]$ for $n \geq 1$ and
starts with an element from $V$, then
\begin{gather}
  \ccont{\mu} \left( v_0 \otimes l^0 \otimes \ul{v_1} \otimes l^1 \otimes
  \dots \otimes \ul{v_n} \otimes l^{n} \right) =
  \notag
  \\
  \sum_{i=0}^{n} \pm \,
  v_0 \otimes l^0 \otimes \ul{v_1} \otimes l^1 \otimes \dots \otimes \ul{v_{i}} \otimes l^{i}_{(1)} \otimes
  \ul{ \corest{\mu} \left( l^{i}_{(2)} \right)} \otimes l^{i}_{(3)} \otimes
  \ul{v_{i+1}} \otimes l^{i+1} \otimes \dots \otimes \ul{v_n} \otimes l^{n}
  \notag
  \\
  \pm \ul{ \corest{\mu} \left( l^{n}_{(2)} \otimes v_0 \otimes l^0_{(1)} \right)} \otimes
  l^0_{(2)} \otimes \ul{v_1} \otimes l^1 \otimes \dots \otimes \ul{v_n} \otimes l^{n}_{(1)}.
  \label{eq:cyc-cont-formula}
\end{gather}
To demonstrate the difference between $\ccont{\mu}$ and
$\cont{\mu}$, we provide an example. Consider $x = v_1 \otimes
  \ul{v_2} \otimes v_3 \in \ndf{V}[1][]$. Then
\begin{equation*}
  \begin{gathered}
    \cont{\mu} \left( v_1 \otimes \ul{v_2} \otimes v_3 \right) =
    \\
    \ul{\mu_0(1)} \otimes v_1 \otimes \ul{v_2} \otimes v_3 \pm
    v_1 \otimes \ul{\mu_0(1)} \otimes \ul{v_2} \otimes v_3 \pm
    v_1 \otimes \ul{v_2} \otimes \ul{\mu_0(1)} \otimes v_3 \pm
    v_1 \otimes \ul{v_2} \otimes v_3 \otimes \ul{\mu_0(1)}
    \\
    \pm \ul{\mu_1(v_1)} \otimes \ul{v_2} \otimes v_3 \pm
    v_1 \otimes \ul{v_2} \otimes \ul{\mu_1(v_3)}
  \end{gathered}
\end{equation*}
while
\begin{equation*}
  \begin{gathered}
    \ccont{\mu} \left( v_1 \otimes \ul{v_2} \otimes v_3 \right) =
    \\
    \pm
    v_1 \otimes \ul{\mu_0(1)} \otimes \ul{v_2} \otimes v_3 \pm
    v_1 \otimes \ul{v_2} \otimes \ul{\mu_0(1)} \otimes v_3 \pm
    v_1 \otimes \ul{v_2} \otimes v_3 \otimes \ul{\mu_0(1)}
    \\
    \pm
    \ul{\mu_1(v_1)} \otimes \ul{v_2} \otimes v_3 \pm
    v_1 \otimes \ul{v_2} \otimes \ul{\mu_1(v_3)}
    \\
    \pm
    \ul{\mu_2(v_3 \otimes v_1)} \otimes \ul{v_2}.
  \end{gathered}
\end{equation*}
Since we will work with $\ccont{\mu}$ only on the quotient $\ncdf{V}[][]$, we won't use
\cref{eq:cyc-cont-formula} and don't write down the signs explicitly.

Finally, note that we have
\begin{equation*}
  \cont{\mu}(1) = \ccont{\mu}(1) = \ul{\mu_0(1)}.
\end{equation*}

\begin{rem} \label{rem:cont-not-defined-over-derivation}
  We have defined the contraction $\cont{\mu}$ only for an $R$-\textit{linear} coderivation
  $\mu \colon \tens{V} \rightharpoonup \tens{V}$. When $R$ is equipped with a non-zero derivation
  $d \colon R \rightharpoonup R$ and $\mu$ is a coderivation over $d$, the construction
  of $\cont{\mu}$ does not make sense. There does not exist a coderivation $\cont{\mu}$,
  generalized or $R$-linear, which satisfies $\left[ \qdr, \cont{\mu} \right] = \lie{\mu}$,
  since the commutator $\left[ \qdr, \cont{\mu} \right]$ is $R$-linear while $\lie{\mu}$
  is a coderivation over $d$.

  Trying to define $\cont{\mu}$ via its corestriction by the same formula \eqref{eq:corest-cont-mu},
  we would have
  \begin{equation*}
    \begin{aligned}
      \left( \cont{\mu} \right)_{1} \left( r \cdot v \right) & =
      \ul{\mu_1 \left( r \cdot v \right)} =
      \ul{dr \cdot v + (-1)^{\braidd{\mu}{r}} r \cdot \mu_1 \left( v \right)}
      \\
                                                             & = (-1)^{\braid{(1,0)}{(0,\degb{d} + \degb{r})}} dr \cdot \ul{v} +
                                                                                                                        (-1)^{\braid{\left( 1, \degb{\mu} \right)}{\left( 0, \degb{r} \right)}}
      r \cdot \left( \cont{\mu} \right)_1 \left( v \right)
    \end{aligned}
  \end{equation*}
  which means that $\left( \cont{\mu} \right)_1$ is neither $R$-linear nor a $d$-operator and
  $\corest{\left( \cont{\mu} \right)} \colon \tens{V \oplus \ul{V}} \rightharpoonup V \oplus \ul{V}$
  cannot be extended to a well-defined operator
  $\cont{\mu} \colon \tens{V \oplus \ul{V}} \rightharpoonup \tens{V \oplus \ul{V}}$.
\end{rem}

\subsection{Differential Calculus} \label{subsec:differential-calculus}

The operators $\qdr, \cont{\mu}$ and $\lie{\mu}$ satisfy
various commutation relations which form a ``differential calculus'':

\begin{lm}
  Let $\mu, \nu \colon \tens{V} \rightharpoonup \tens{V}$ be two generalized
  coderivations. The following commutation relations hold on $\ndf{V}[][]$
  whenever the operators involved are defined:\footnote{That is, whenever the contraction $\cont{\mu}$ appears,
    we assume that $\mu$ is an $R$-linear coderivation.}
  \begin{align}
    \left[ \qdr, \lie{\mu} \right]        & =
    \qdr \circ \lie{\mu} -
               (-1)^{\braid{(-1,0)}{(0,\degb{\mu})}} \lie{\mu} \circ \qdr = 0,
    \nonumber
    \\
    \left[ \qdr, \cont{\mu} \right]       & =
    \qdr \circ \cont{\mu} - (-1)^{\braid{(-1,0)}{(1,\degb{\mu})}} \cont{\mu} \circ \qdr =
    \lie{\mu},
    \nonumber
    \\
    \left[ \lie{\mu}, \lie{\nu} \right]   & =
    \lie{\mu} \circ \lie{\nu} -
                    (-1)^{\braid{(0,\degb{\mu})}{(0,\degb{\nu})}}
    \lie{\nu} \circ \lie{\mu} =
    \lie{\left[ \mu, \nu \right]},
    \label{eq:lieliecomm}
    \\
    \left[ \cont{\nu}, \lie{\mu} \right]  & =
    \cont{\nu} \circ \lie{\mu} -
                     (-1)^{\braid{(1,\degb{\nu})}{(0,\degb{\mu})}} \lie{\mu} \circ \cont{\nu} =
    \cont{\left[ \nu, \mu \right]},
    \label{eq:liecontcomm}
    \\
    \left[ \cont{\mu}, \cont{\nu} \right] & = \cont{\mu} \circ
    \cont{\nu} -
    (-1)^{\braid{(1,\degb{\mu})}{(1,\degb{\nu})}}
    \cont{\nu} \circ
    \cont{\mu} = 0.
    \label{eq:contcontcomm}
  \end{align}
\end{lm}
\begin{proof}
  The first two relations are part of the defining relations for the operators $\lie{\mu}$ and $\cont{\mu}$
  respectively. To prove \cref{eq:lieliecomm}, note that by the graded Jacobi identity
  and the graded antisymmetry of the Lie bracket, we have
  \begin{equation*}
    \pm \left[ \qdr, \left[ \lie{\mu}, \lie{\nu} \right] \right] \pm
    \cancel{ \left[ \lie{\nu}, \left[ \qdr, \lie{\mu} \right] \right] } \pm
    \cancel{ \left[ \lie{\mu}, \left[ \lie{\nu}, \qdr \right] \right] } =
    \pm \left[ \qdr, \left[ \lie{\mu}, \lie{\nu} \right] \right] = 0.
  \end{equation*}
  Hence $\left[ \lie{\mu}, \lie{\nu} \right]$ commutes with $\qdr$, and since
  \begin{equation*}
    \rest{\left[ \lie{\mu}, \lie{\nu} \right]}{\tens{V}} =
    \left[ \rest{\lie{\mu}}{\tens{V}}, \rest{\lie{\nu}}{\tens{V}} \right] =
    \left[ \mu, \nu \right],
  \end{equation*}
  we must have $\left[ \lie{\mu}, \lie{\nu} \right] = \lie{\left[ \mu, \nu \right]}$ by the uniqueness part
  of \cref{lm:existence-uniqueness-lie}.

  To prove \cref{eq:liecontcomm}, we again use the Jacobi identity to see that
  \begin{equation*}
    (-1)^{\braid{(-1,0)}{(0,\degb{\mu})}}
    \left[ \qdr, \left[ \cont{\nu}, \lie{\mu} \right] \right] +
    (-1)^{\braid{(0,\degb{\mu})}{(1,\degb{\nu})}}
    \left[ \lie{\mu}, \left[ \qdr, \cont{\nu} \right] \right] +
    (-1)^{\braid{(1,\degb{\nu})}{(-1,0)}}
    \cancel{\left[ \cont{\nu}, \left[ \lie{\mu}, \qdr \right] \right]}
    = 0.
  \end{equation*}
  Hence
  \begin{equation*}
    \begin{aligned}
      \left[ \qdr, \left[ \cont{\nu}, \lie{\mu} \right] \right]
      \stackrel{\phantom{\eqref{eq:lieliecomm}}}{=}{}                  &
      (-1)^{\braid{(-1,0)}{(0,\degb{\mu})} + \braid{(0,\degb{\mu})}{(1,\degb{\nu})} + 1}
      \left[ \lie{\mu}, \left[ \qdr, \cont{\nu} \right] \right]
      \\
      \stackrel{\eqref{eq:lie-derivative-commutator-contraction}}{=}{} &
      (-1)^{\braid{(0, \degb{\mu})}{(0,\degb{\nu})} + 1} \left[ \lie{\mu}, \lie{\nu} \right]
      =
      \left[ \lie{\nu}, \lie{\mu} \right]
      \\
      \stackrel{\eqref{eq:lieliecomm}}{=}{}                            &
      \lie{\left[ \nu, \mu \right]}.
    \end{aligned}
  \end{equation*}
  Since $\left[ \qdr, \cont{\left[ \nu, \mu \right]} \right] = \lie{\left[ \nu, \mu \right]}$ by definition,
  we must have
  $\left[ \cont{\nu}, \lie{\mu} \right] = \cont{\left[ \nu, \mu \right]}$
  by the uniqueness part of \cref{lm:existence-uniqueness-contraction}.

  Finally, to prove \cref{eq:contcontcomm}, note that the commutator $\left[ \cont{\mu}, \cont{\nu} \right]$
  has degree $(2, \degb{\mu} + \degb{\nu})$. This means that the corestriction
  $\corest{ \left[ \cont{\mu}, \cont{\nu} \right] } \colon \tens{V \oplus \ul{V}} \rightharpoonup V \oplus \ul{V}$
  has line degree two and hence must vanish by degree considerations. Since
  $\corest{ \left[ \cont{\mu}, \cont{\nu} \right] } = 0$, we also have $\left[ \cont{\mu}, \cont{\nu} \right] = 0$.
\end{proof}

Applying \cref{cor:cycl-comm-bracket-quotient}, we see that all the commutation relations satisfied
by the operators also hold for cyclizations of the operators \textit{on the quotient} $\ncdf{V}[][]$:

\begin{lm} \label{lm:cyclic-commutation-relations}
  Let $\mu, \nu \colon \tens{V} \rightharpoonup \tens{V}$ be two generalized
  coderivations. The following commutation relations hold on $\ncdf{V}[][]$
  whenever the operators involved are defined:
  \begin{align}
    \left[ \qdr, \clie{\mu} \right]         & =
    \qdr \circ \clie{\mu} - (-1)^{\braid{(-1,0)}{(0,\degb{\mu})}} \clie{\mu} \circ \qdr = 0,
    \label{eq:qliecommcyc}
    \\
    \left[ \qdr, \ccont{\mu} \right]        & =
    \qdr \circ \ccont{\mu} - (-1)^{\braid{(-1,0)}{(1,\degb{\mu})}} \ccont{\mu} \circ \qdr =
    \clie{\mu},
    \label{eq:qcontliecyc}
    \\
    \left[ \clie{\mu}, \clie{\nu} \right]   & =
    \clie{\mu} \circ \clie{\nu} - (-1)^{\braid{(0,\degb{\mu})}{(0,\degb{\nu})}}
    \clie{\nu} \circ \clie{\mu} =
    \clie{\left[ \mu, \nu \right]},
    \label{eq:lieliecommcyc}
    \\
    \left[ \ccont{\nu}, \clie{\mu} \right]  & =
    \ccont{\nu} \circ \clie{\mu} - (-1)^{\braid{(1,\degb{\nu})}{(0,\degb{\mu})}}
    \clie{\mu} \circ \ccont{\nu} =
    \ccont{\left[ \nu, \mu \right]},
    \label{eq:liecontcommcyc}
    \\
    \left[ \ccont{\mu}, \ccont{\nu} \right] & = \ccont{\mu} \circ \ccont{\nu} -
                                                                  (-1)^{\braid{(1,\degb{\mu})}{(1,\degb{\nu})}} \ccont{\nu} \circ \ccont{\mu} = 0.
    \label{eq:contcontcommcyc}
  \end{align}
  \qed
\end{lm}

\subsection{The Formal Poincar\'{e} Lemma} \label{sec:formal-poincare}
Recall that we denote by $\ndfr{V}[][]$ the reduced tensor module on $V \oplus \ul{V}$.
Since $\qdr$ is a coderivation, it maps $\ndfr{V}[][]$ to itself,
and, fixing $j \in \ZZ$, we obtain a homological complex
\begin{equation} \label{eq:q-row-complex-ndfr}
  \left( \ndfr{V}[*][j], \qdr \right) \colon \qquad
  \cdots \leftarrow 0 \leftarrow \ndfr{V}[0][j] \xleftarrow{\qdr} \ndfr{V}[1][j] \xleftarrow{\qdr}
  \ndfr{V}[2][j] \xleftarrow{\qdr} \cdots
\end{equation}
of $\mathbbm{k}$-modules. The following lemma shows that the resulting complex is contractible:
\begin{lm}[Formal Poincar\'{e} Lemma for ${\ndfr{V}[][]}$] \label{lm:formal-poincare-ndfr}
  \sloppy
  Let $X$ be the Euler coderivation on $\tens{V}$ defined by the formula
  \begin{equation*}
    X \left( v_1 \otimes \dots \otimes v_n \right) \defeq n \cdot \left( v_1 \otimes \dots \otimes v_n \right)
  \end{equation*}
  for all $v_1, \dots, v_n \in V$ and $n \geq 0$, and let $h_{\dr} \colon \ndfr{V} \rightharpoonup \ndfr{V}[*+1]$
  be the map which acts on $\left( V \oplus \ul{V} \right)^{\otimes n} \subset \ndfr{V}$
  by $\frac{1}{n} \cdot \cont{X}$.
  Then $h_{\dr}$ is a bounded $R$-linear map of degree $(1,0)$ with $\nnorm[h_{\dr}] \leq 1$ which satisfies
  \begin{equation*}
    h_{\dr} \qdr + \qdr h_{\dr} = \idd.
  \end{equation*}
  The map $h_{\dr}$ satisfies in addition the following properties:
  \begin{enumerate}
    \item The map $h_{\dr}$ preserves the weight of elements.
    \item $h_{\dr} \, \t = \t \, h_{\dr}$.
    \item $h_{\dr}^2 = 0$.
    \item $\left[ h_{\dr}, \lie{\nu} \right] = 0$ for all coderivations $\nu$ on $\tens{V}$ which
          satisfy $\nu_k = 0$ for all $k \neq 1$.
  \end{enumerate}
\end{lm}
\begin{proof}
  Note that any coderivation $\nu$ on the tensor coalgebra whose corestriction
  satisfies $\nu_k = 0$ for $k \neq 1$ preserves weight and commutes with $\t$.
  Given a coderivation $\nu$ on $\tens{V}$
  which satisfies $\nu_k = 0$ for $k \neq 1$, \cref{eq:lie-der-corest,eq:corest-cont-mu}
  show that the coderivations $\lie{\nu}$ and $\cont{\nu}$ (defined if $\nu$ is $R$-linear)
  on $\tens{V \oplus \ul{V}}$ also satisfy
  $\left( \lie{\nu} \right)_k = \left( \cont{\nu} \right)_k = 0$ for $k \neq 1$ and
  hence preserve the weight of elements on $\tens{V \oplus \ul{V}}$ and commute with $\t$.

  In terms of the corestrictions, the Euler coderivation
  $X$ is the degree zero coderivation which satisfies $X_1 = \idd_{V}$ and $X_k = 0$ for $k \neq 1$.
  By \cref{eq:lie-der-corest} we have $\left( \lie{X} \right)_1 = \idd_{V \oplus \ul{V}}$ and
  $\left( \lie{X} \right)_k = 0$ for $k \neq 1$
  and so $\lie{X}$ is the Euler coderivation on $\tens{V \oplus \ul{V}}$.
  Since $X_k = 0$ for $k \neq 1$, the contraction $\cont{X}$ preserves weight and commutes with $\t$
  and hence the same holds for $h_{\dr}$.

  Let $x \in \ndfr{V}[][]$ be an element of weight $n$.
  Since $\qdr$ preserves weight, we have
  \begin{equation*}
    \left( h_{\dr} \qdr + \qdr h_{\dr} \right) \left( x \right) =
    \frac{1}{n} \cdot \left( \qdr \cont{X} + \cont{X} \qdr \right) \left( x \right)
    \stackrel{\eqref{eq:parity-extends-koszul}}{=}
    \frac{1}{n} \left[ \qdr, \cont{X} \right] \left( x \right)
    \stackrel{\eqref{eq:lie-derivative-commutator-contraction}}{=}
    \frac{1}{n} \cdot \lie{X} \left( x \right) = x.
  \end{equation*}
  We also have
  \begin{equation*}
    h_{\dr}^2 \left( x \right) = \frac{1}{n^2} \cdot \cont{X}^2 \left( x \right)
    \stackrel{\eqref{eq:parity-extends-koszul}}{=}
    \frac{1}{2n^2} \cdot \left[ \cont{X}, \cont{X} \right] \left( x \right)
    \stackrel{\eqref{eq:contcontcomm}}{=} 0.
  \end{equation*}
  Finally, any coderivation $\nu$ on $\tens{V}$ with $\nu_k = 0$ for $k \neq 1$ commutes with $X$.
  Since $\lie{\nu}$ also preserves weight, we have
  \begin{equation*}
    \left[ h_{\dr}, \lie{\nu} \right] =
    \frac{1}{n} \cdot \left[ \cont{X}, \lie{\nu} \right] \left( x \right)
    \stackrel{\eqref{eq:liecontcomm}}{=}
    \frac{1}{n} \cdot \cont{[X, \nu]} \left( x \right) = 0.
  \end{equation*}
\end{proof}

We can also consider the action of $\qdr$ on the reduced cyclic tensor module
$\ncdfr{V}[][]$ and, fixing $j \in \ZZ$, obtain a homological complex
\begin{equation} \label{eq:q-row-complex-ncdfr}
  \left( \ncdfr{V}[*][j], \qdr \right) \colon \qquad
  \cdots \leftarrow 0 \leftarrow \ncdfr{V}[0][j] \xleftarrow{\qdr} \ncdfr{V}[1][j] \xleftarrow{\qdr}
  \ncdfr{V}[2][j] \xleftarrow{\qdr} \cdots
\end{equation}
of $\mathbbm{k}$-modules which is a quotient of the complex $\left( \ndfr{V}[*][j], \qdr \right)$ by
$\Im \left( \idd - \t \right)$. We have a Formal Poincar\'{e} lemma for the resulting complex of cyclic
codifferential forms:

\begin{lm}[Formal Poincar\'{e} Lemma for ${\ncdfr{V}[][]}$] \label{lm:formal-poincare-ncdfr}
  \sloppy
  The map $h_{\dr}$ of \cref{lm:formal-poincare-ndfr} descends to a contraction
  $h_{\dr} \colon \ncdfr{V} \rightharpoonup \ncdfr{V}[*+1]$, i.e., we have
  \begin{equation}
    h_{\dr} \qdr + \qdr h_{\dr} = \idd_{\ncdfr{V}[][]} \label{eq:h-dr-q-contraction}
  \end{equation}
  on $\ncdfr{V}[][]$.
  The induced map $h_{\dr}$ satisfies the following properties:
  \begin{enumerate}
    \item The map $h_{\dr}$ preserves the weight of elements.
    \item $h_{\dr}^2 = 0$. \label{item:hdr-square-zero}
    \item $\left[ h_{\dr}, \clie{\nu} \right] = 0$ for all coderivations $\nu$ on $\tens{V}$ which
          satisfy $\nu_k = 0$ for all $k \neq 1$.
  \end{enumerate}
\end{lm}
\begin{proof}
  Since the contraction $h_{\dr}$ constructed in \cref{lm:formal-poincare-ndfr} commutes with $\t$,
  it descends to a contraction $h_{\dr} \colon \ncdfr{V} \rightharpoonup \ncdfr{V}[*+1]$ on the quotient
  which is given by the same formula and satisfies the same properties as in \cref{lm:formal-poincare-ndfr}.
  The property $\left[ h_{\dr}, \clie{\nu} \right] = 0$ follows from the corresponding property
  of \cref{lm:formal-poincare-ndfr} since $\clie{\nu} = \lie{\nu}$ when $\nu_k = 0$ for all $k \neq 1$.
\end{proof}

\subsection{Functoriality and Naturality} \label{sec:functoriality-ndf-and-ncdf}
Let $f \colon \tens{V}[R] \rightarrow \tens{W}[S]$ be a graded Banach coalgebra
morphism. We show that there exists a unique way to extend $f$ to a graded Banach coalgebra morphism
\begin{equation*}
  \indmap{f} = \indmap{F} \left( f \right) \colon \tens{V \oplus \ul{V}}[R] \rightarrow \tens{W \oplus \ul{W}}[S]
\end{equation*}
which commutes with the de Rham differential $\qdr$:

\begin{lm} \label{lm:existence-uniqueness-morphism-extension}
  Let $f \colon \tens{V}[R] \rightarrow \tens{W}[S]$ be a graded Banach coalgebra
  morphism over a graded Banach $\mathbbm{k}$-algebra morphism $\varphi \colon R \rightarrow S$.
  Then there exists a unique morphism
  $\indmap{f} = \indmap{F} \left( f \right) \colon
    \tens{V \oplus \ul{V}}[R] \rightarrow \tens{W \oplus \ul{W}}[S]$
  of graded Banach coalgebra over $\varphi$, called the \textbf{induced morphism},
  which satisfies
  \begin{equation} \label{eq:ind-map-extends-f}
    \rest{\indmap{f}}{\tens{V}} = f
  \end{equation}
  and
  \begin{equation} \label{eq:func-qdr}
    \left[ \qdr, \indmap{f} \right] = 0.
  \end{equation}
  The corestriction of the morphism $\indmap{F} \left( f \right)$ is given by
  \begin{empheq}[left={\corest{\indmap{f}} \left( x \right) \defeq \empheqlbrace}]{align}
    &\corest{f} \left( x \right) && x \in \ndf{V}[0][],
    \label{eq:def-ind-map-0}
    \\
    &\ul{ \corest{f} \left( \qdr \left( x \right) \right)} && x \in \ndf{V}[1][],
    \label{eq:def-ind-map-1}
    \\
    &0 && x \in \ndf{V}[n][], \quad n \geq 2.
    \label{eq:def-ind-map-twoplus}
  \end{empheq}
\end{lm}
\begin{proof}
  The proof is entirely parallel to the proof of \cref{lm:existence-uniqueness-lie}.
  Assume such a morphism $\indmap{f}$ exists. Since $\indmap{f}$ is a coalgebra morphism, it has degree $(0,0)$ and
  induces well-defined maps $\rest{\indmap{f}}{\ndf{V}[n][]} \colon \ndf{V}[n][*] \rightarrow \ndf{W}[n][*]$
  over $\varphi$ for all $n \geq 0$. Let $x_n \in \ndf{V}[n][]$.
  \begin{enumerate}
    \item \Cref{eq:ind-map-extends-f} forces that
          $\corest{ \indmap{f} } \left( x_0 \right) = \corest{f} \left( x_0 \right)$.
    \item \Cref{eq:func-qdr} forces that
          \begin{equation*}
            0 = \corest{\left[ \qdr, \indmap{f} \right]} \left( x_1 \right) =
            \left( \corest{\qdr} \circ \indmap{f} - \corest{ \indmap{f} } \circ \qdr \right) \left( x_1 \right)
          \end{equation*}
          and hence
          \begin{equation*}
            \left( \qdr_1 \circ \corest{ \indmap{f} } \right) \left( x_1 \right)
            \stackrel{\eqref{eq:corest-qdr}}{=}
            \left( \corest{\qdr} \circ \indmap{f} \right) \left( x_1 \right) =
            \left( \corest{ \indmap{f} } \circ \qdr \right) \left( x_1 \right) =
            \corest{f} \left( \qdr \left( x_1 \right) \right).
          \end{equation*}
          Since $\corest{ \indmap{f} } \left( x_1 \right)$ has line degree one,
          i.e., $\corest{ \indmap{f} } \left( x_1 \right) \in \ul{W}$, we have
          \begin{equation} \label{eq:qdr-commute-f-one-iff}
            \left[ \qdr, \indmap{f} \right] \left( x_1 \right) = 0 \iff
            \corest{ \indmap{f} } \left( x_1 \right) =
            \ul{ \corest{f} \left( \qdr \left( x_1 \right) \right) }.
          \end{equation}
    \item By degree considerations, any morphism
          $h \colon \tens{V \oplus \ul{V}}[R] \rightarrow \tens{W \oplus \ul{W}}[S]$
          must satisfy $\corest{h} \left( x_n \right) = 0$ for $n \geq 2$. In particular,
          $\corest{ \indmap{f} } \left( x_n \right) = 0$ for all $n \geq 2$.
  \end{enumerate}
  This shows that \cref{eq:ind-map-extends-f,eq:func-qdr} force the corestriction of $\indmap{f}$
  to be given by \cref{eq:def-ind-map-0,eq:def-ind-map-1,eq:def-ind-map-twoplus}. Since $\indmap{f}$ is
  determined uniquely by its corestriction, we obtain uniqueness.

  To show existence, we can use \cref{prop:classification-morphisms-tensor-coalgebra-different-ground-algebras}
  to define a Banach coalgebra morphism
  $\indmap{f} \colon \tens{V \oplus \ul{V}}[R] \rightarrow \tens{W \oplus \ul{W}}[S]$
  over $\varphi$ by specifying its corestriction to be given by
  \cref{eq:def-ind-map-0,eq:def-ind-map-1,eq:def-ind-map-twoplus}.
  It is clear that \cref{eq:def-ind-map-0} implies the resulting morphism $\indmap{f}$ extends $f$.
  Since $\left[ \qdr, \indmap{f} \right]$ has degree $(-1, 0)$, it is clear that
  $\corest{ \left[ \qdr, \indmap{f} \right] }$ vanishes on $\ndf{V}[n][]$ for $n \neq 1, 2$.
  For $x_1 \in \ndf{V}[1][]$ we have
  $\left[ \qdr, \indmap{f} \right] \left( x_1 \right) = 0$ by construction
  (see \cref{eq:qdr-commute-f-one-iff}) and for $x_2 \in \ndf{V}[2][]$ we have
  \begin{equation*}
    \begin{aligned}
      \corest{\left[ \qdr, \indmap{f} \right]} \left( x_2 \right)
      \stackrel{\phantom{\eqref{eq:lie-der-corest}}}{=}{}                      &
      \left( \corest{\qdr} \circ \indmap{f} - \corest{\indmap{f}} \circ \qdr \right) \left( x_2 \right)
      \\
      \stackrel{\eqref{eq:corest-qdr}}{=}{}                                    &
      \left( \qdr_1 \circ \corest{ \indmap{f} } - \corest{\indmap{f}} \circ \qdr \right) \left( x_2 \right)
      \\
      \stackrel[\eqref{eq:def-ind-map-twoplus}]{\eqref{eq:def-ind-map-1}}{=}{} &
      - \ul{ \corest{f} \left( \qdr^2 \left( x_2 \right) \right)}
      \stackrel{\cref{lm:qdr-squared-zero}}{=}{} 0.
    \end{aligned}
  \end{equation*}
\end{proof}
\Cref{eq:def-ind-map-0,eq:def-ind-map-1,eq:def-ind-map-twoplus} for $\corest{\indmap{f}}$ imply
that the induced morphism $\indmap{f}$ is bounded with $\nnorm[\indmap{f}] \leq \nnorm[f]$, and,
since $\indmap{f}$ extends $f$ on $\tens{V}$, we have $\nnorm[\indmap{f}] = \nnorm[f]$.

Note that when $f = \idd \colon \tens{V}[R] \rightarrow \tens{V}[R]$ then
$\indmap{F} \left( f \right) = \idd \colon \tens{V \oplus \ul{V}}[R] \rightarrow \tens{V \oplus \ul{V}}[R]$.
Moreover, the construction of the induced morphism is functorial:

\begin{lm}[Functoriality] \label{lm:ind-map-functoriality}
  Let $f \colon \tens{V}[R] \rightarrow \tens{W}[S]$ and $g \colon \tens{U}[Q] \rightarrow \tens{V}[R]$
  be two graded Banach coalgebra morphisms.
  Then
  $\indmap{F} \left( f \circ g \right) = \indmap{F} \left( f \right) \circ \indmap{F} \left( g \right)$.
\end{lm}
\begin{proof}
  Since $\indmap{F} \left( f \right), \indmap{F} \left( g \right)$ both commute with $\qdr$, their composition
  also commutes with $\qdr$ and extends $f \circ g$ on $\tens{U}[Q]$. By the uniqueness part
  of \cref{lm:existence-uniqueness-morphism-extension}, we must have
  $\indmap{F} \left( f \circ g \right) = \indmap{F} \left( f \right) \circ \indmap{F} \left( g \right)$.
\end{proof}

\begin{rem} \label{rem:func-tv-ndf}
  \cref{lm:ind-map-functoriality,lm:qdr-squared-zero,lm:existence-uniqueness-morphism-extension}
  show that the construction
  \begin{equation*}
    \begin{aligned}
      \tens{V}                               & \mapsto \left( \ndf{V}[][], \qdr \right),
      \\
      f \colon \tens{V} \rightarrow \tens{W} & \mapsto
      \indmap{f} \colon \left( \ndf{V}[][], \qdr \right) \rightarrow \left( \ndf{W}[][], \qdr \right)
    \end{aligned}
  \end{equation*}
  is a functor from the category of $\ZZ$-graded tensor Banach coalgebras to the category of
  $\ZZ^2$-differential-graded tensor Banach coalgebras.
\end{rem}

\begin{lm}[Naturality of $\lie{\left(-\right)}$] \label{lm:nat-lie}
  Let $\varphi \colon \left( R, d_R \right) \rightarrow \left( S, d_S \right)$ be a morphism
  of pre-differential graded Banach $\mathbbm{k}$-algebras and let
  $f \colon \tens{V}[R] \rightarrow \tens{W}[S]$ be a graded Banach coalgebra morphism over $\varphi$.
  Let $\mu \colon \tens{V}[R] \rightharpoonup \tens{V}[R]$ be a generalized coderivation over $d_R$
  and let $\nu \colon \tens{W}[S] \rightharpoonup \tens{W}[S]$ be a generalized coderivation over $d_S$.
  Assume that $f \circ \mu = \nu \circ f$. Then
  \begin{equation}
    \indmap{f} \circ \lie{\mu} = \lie{\nu} \circ \indmap{f}. \label{eq:func-lie}
  \end{equation}
\end{lm}
\begin{proof}
  By \cref{lm:f-circ-mu-nu-circ-f-corest-different-ground-algebras}, it is enough to show the
  corestrictions of $\indmap{f} \circ \lie{\mu}$ and $\lie{\nu} \circ \indmap{f}$ coincide.
  Let $x \in \ndf{V}[][]$ be a homogeneous element with respect to line degree. We split into cases:
  \begin{enumerate}
    \item If $x \in \ndf{V}[0][]$ then
          \begin{equation*}
            \begin{aligned}
              \corest{\indmap{f}} \left( \lie{\mu} \left( x \right) \right)
               & \stackrel{\eqref{eq:lie-derivative-extends-mu}}{=}
              \corest{\indmap{f}} \left( \mu \left( x \right) \right)
              \stackrel{\eqref{eq:def-ind-map-0}}{=}
              \corest{f} \left( \mu \left( x \right) \right)
              =
              \corest{\nu} \left( f \left( x \right) \right)
              \\
               & \stackrel{\eqref{eq:lie-derivative-extends-mu}}{=}
              \corest{\left( \lie{\nu} \right)} \left( f \left( x \right) \right)
              \stackrel{\eqref{eq:def-ind-map-0}}{=}
              \corest{\left( \lie{\nu} \right)} \left( \indmap{f} \left( x \right) \right).
            \end{aligned}
          \end{equation*}
    \item If $x \in \ndf{V}[1][]$, then
          \begin{equation*}
            \begin{aligned}
              \corest{\indmap{f}} \left( \lie{\mu} \left( x \right) \right)
               & \stackrel{\eqref{eq:def-ind-map-1}}{=}
              \ul{ \corest{f} \left( \qdr \left( \lie{\mu} \left( x \right) \right) \right)}
              \stackrel{\eqref{eq:qliecomm}}{=}
              (-1)^{\braid{(-1,0)}{(0,\degb{\mu})}}
              \ul{ \corest{f} \left( \lie{\mu} \left( \qdr \left( x \right) \right) \right)}
              \\
               & \stackrel{\eqref{eq:lie-derivative-extends-mu}}{=}
                 (-1)^{\braid{(-1,0)}{(0,\degb{\mu})}}
              \ul{ \corest{f} \left( \mu \left( \qdr \left( x \right) \right) \right)}
              =
              (-1)^{\braid{(-1,0)}{(0,\degb{\nu})}}
              \ul{ \corest{\nu} \left( f \left( \qdr \left( x \right) \right) \right)}
              \\
               & \stackrel{\eqref{eq:func-qdr}}{=}
                 (-1)^{\braid{(-1,0)}{(0,\degb{\nu})}}
              \ul{ \corest{\nu} \left( \qdr \left( \indmap{f} \left( x \right) \right) \right) }
              \stackrel{\eqref{eq:lie-der-corest}}{=}
              \corest{\left( \lie{\nu} \right)} \left( \indmap{f} \left( x \right) \right).
            \end{aligned}
          \end{equation*}
    \item If $x \in \ndf{V}[n][]$ for $n \geq 2$ then
          $\corest{\indmap{f}} \left( \lie{\mu} \left( x \right) \right) =
            \corest{\left( \lie{\nu} \right)} \left( \indmap{f} \left( x \right) \right) = 0$
          by degree considerations.
  \end{enumerate}
\end{proof}

\begin{rem}
  \Cref{lm:nat-lie}, together with \cref{rem:func-tv-ndf} and \cref{eq:qliecomm},
  show that the construction
  \begin{equation*}
    \begin{aligned}
      \left( \tens{V}, \mu \right)                                                   & \mapsto \left( \ndf{V}[][], \qdr, \lie{\mu} \right),
      \\
      f \colon \left( \tens{V}, \mu \right) \rightarrow \left( \tens{W}, \nu \right) & \mapsto
      \indmap{f} \colon \left( \ndf{V}[][], \qdr, \lie{\mu} \right) \rightarrow
      \left( \ndf{W}[][], \qdr, \lie{\nu} \right)
    \end{aligned}
  \end{equation*}
  is a functor from the category of pre-differential $\ZZ$-graded Banach tensor coalgebras
  to the category of $\ZZ^2$-graded Banach tensor coalgebras equipped with two commuting coderivations,
  where the morphisms are required to commute with both coderivations.
\end{rem}

\begin{lm}[Naturality of $\cont{\left(-\right)}$]
  Let $f \colon \tens{V}[R] \rightarrow \tens{W}[S]$ be a graded Banach coalgebra morphism.
  Let $\mu \colon \tens{V}[R] \rightharpoonup \tens{V}[R]$ be an $R$-linear coderivation and let
  $\nu \colon \tens{W}[S] \rightharpoonup \tens{W}[S]$ be an $S$-linear coderivation. Assume that
  $f \circ \mu = \nu \circ f$. Then
  \begin{equation}
    \indmap{f} \circ \cont{\mu} = \cont{\nu} \circ \indmap{f}. \label{eq:func-cont}
  \end{equation}
\end{lm}
\begin{proof}
  By \cref{lm:f-circ-mu-nu-circ-f-corest-different-ground-algebras}, it is enough to show the
  corestrictions of $\indmap{f} \circ \cont{\mu}$ and $\cont{\nu} \circ \indmap{f}$ coincide.
  Let $x \in \ndf{V}[][]$ be a homogeneous element with respect to line degree. We split into cases:
  \begin{enumerate}
    \item If $x \in \ndf{V}[0][]$ then
          \begin{equation*}
            \corest{\indmap{f}} \left( \cont{\mu} \left( x \right) \right)
            \stackrel{\eqref{eq:def-ind-map-1}}{=}
            \ul{ \corest{f} \left( \qdr \left( \cont{\mu} \left( x \right) \right) \right)}
            \stackrel{\eqref{eq:lie-derivative-commutator-contraction}}{=}
            \ul{ \corest{f} \left( \lie{\mu} \left( x \right) \right)}
            \stackrel{\eqref{eq:lie-derivative-extends-mu}}{=}
            \ul{ \corest{f} \left( \mu \left( x \right) \right)}
          \end{equation*}
          while
          \begin{equation*}
            \corest{\left( \cont{\nu} \right)} \left( \indmap{f} \left( x \right) \right)
            \stackrel{\eqref{eq:ind-map-extends-f}}{=}
            \corest{\left( \cont{\nu} \right)} \left( f \left( x \right) \right)
            \stackrel{\eqref{eq:corest-cont-mu}}{=}
            \ul{ \corest{\nu} \left( f \left( x \right) \right)}.
          \end{equation*}
          Since $f \circ \mu = \nu \circ f$ implies $\corest{f} \circ \mu = \corest{\nu} \circ f$,
          the two corestrictions coincide.
    \item If $x \in \ndf{V}[n][]$ for $n \geq 1$ then
          $\corest{\indmap{f}} \left( \cont{\mu} \left( x \right) \right) =
            \corest{\left( \cont{\nu} \right)} \left( \indmap{f} \left( x \right) \right) = 0$
          by degree considerations.
  \end{enumerate}
\end{proof}

We shall denote the cyclization of the morphism $\indmap{F} \left( f \right) = \indmap{f}$ by
$\cindmap{F} \left( f \right) = \cindmap{f}$. The map $\cindmap{f}$ is a bounded morphism of
graded Banach modules with $\nnorm[\cindmap{f}] \leq \nnorm[f]$.
We have the cyclic version of \cref{lm:ind-map-functoriality}:
\begin{lm}[Cyclic Functoriality] \label{lm:cind-map-functoriality}
  Let $f \colon \tens{V}[R] \rightarrow \tens{W}[S]$ and $g \colon \tens{U}[Q] \rightarrow \tens{V}[R]$
  be two graded Banach coalgebra morphisms.
  Then we have
  \begin{equation} \label{eq:cind-map-functoriality}
    \cindmap{F} \left( f \circ g \right) = \cindmap{F} \left( f \right) \circ \cindmap{F} \left( g \right)
    \colon \ncdf{U}[][] \rightarrow \ncdf{W}[][]
  \end{equation}
  on $\ncdf{U}[][]$.
\end{lm}
\begin{proof}
  Using \cref{lm:ind-map-functoriality} and \cref{cor:func-cycl-morphism-extended}, we have
  \begin{equation*}
    \cindmap{F} \left( f \circ g \right) = \cycl{ \indmap{F} \left( f \circ g \right) } =
    \cycl{ \indmap{F} \left( f \right) \circ \indmap{F} \left( g \right) } =
    \cycl{ \indmap{F} \left( f \right) } \circ \cycl{ \indmap{F} \left( g \right) } =
    \cindmap{F} \left( f \right) \circ \cindmap{F} \left( g \right).
  \end{equation*}
\end{proof}

Applying \cref{cor:func-cycl-morphism-coderivation-extended}, we also have:
\begin{cor} \label{cor:homo-cyc-func}
  Let $\mu \in \CoDer{\tens{V}}[R], \nu \in \CoDer{\tens{W}}[S]$ be two generalized coderivations
  and let $f \colon \tens{V}[R] \rightarrow \tens{W}[S]$ be a graded Banach coalgebra morphism such
  that $f \circ \mu = \nu \circ f$.
  Then the induced morphism $\cindmap{f} \colon \ncdf{V}[][] \rightarrow \ncdf{W}[][]$ satisfies
  \begin{align}
    \cindmap{f} \circ \qdr = \qdr \circ \cindmap{f},
    \label{eq:func-qdr-cyc}
    \\
    \cindmap{f} \circ \ccont{\mu} = \ccont{\nu} \circ \cindmap{f},
    \label{eq:func-cont-cyc}
    \\
    \cindmap{f} \circ \clie{\mu} = \clie{\nu} \circ \cindmap{f},
    \label{eq:func-lie-cyc}
  \end{align}
  whenever the operators involved are defined.\footnote{That is, the identity \eqref{eq:func-cont-cyc} holds
    only when $\mu$ is $R$-linear and $\nu$ is $S$-linear.}
\end{cor}

\begin{rem} \label{rem:func-bi-ncdf}
  \cref{lm:cind-map-functoriality,lm:qdr-squared-zero,cor:homo-cyc-func}
  show that the construction
  \begin{equation*}
    \begin{aligned}
      \left( \tens{V}, \mu \right)                                                   & \mapsto \left( \ncdf{V}[][], \qdr, \clie{\mu} \right),
      \\
      f \colon \left( \tens{V}, \mu \right) \rightarrow \left( \tens{W}, \nu \right) & \mapsto
      \cindmap{f} \colon \left( \ncdf{V}[][], \qdr, \clie{\mu} \right) \rightarrow
      \left( \ncdf{W}[][], \qdr, \clie{\nu} \right)
    \end{aligned}
  \end{equation*}
  is a functor from the category of pre-differential $\ZZ$-graded Banach tensor coalgebras
  to the category of $\ZZ^2$-graded Banach modules equipped with two commuting pre-differentials,
  where the morphisms are required to commute with both pre-differentials. When $\mu$ is odd and
  $\mu^2 = 0$, we obtain a functor from the category of differential $\ZZ$-graded Banach tensor coalgebras
  to the category of Banach bicomplexes.
\end{rem}

Finally, let us present explicit formulas for the action of $\indmap{f}$ and $\cindmap{f}$
on $\ndf{V}[n][]$ in terms of $f$. When $x \in \ndf{V}[0][]$ then
$\indmap{f} \left( x \right) = f \left( x \right)$. When $n \geq 1$ and
$x$ has the form
$x = l^0 \otimes \ul{v_1} \otimes l^1 \otimes \dots \otimes \ul{v_n} \otimes l^{n}$,
then
\begin{equation*}
  \begin{aligned}
    \indmap{f} \left( x \right) ={} &
    f \left( l^0_{(1)} \right) \otimes
    \ul{\corest{f} \left( \qdr \left( l^0_{(2)} \otimes \ul{v_1} \otimes l^1_{(1)} \right) \right)} \otimes
    f \left( l^1_{(2)} \right) \otimes
    \ul{\corest{f} \left( \qdr \left( l^1_{(3)} \otimes \ul{v_2} \otimes l^2_{(1)} \right) \right)} \otimes
    f \left( l^2_{(2)} \right) \otimes
    \\
                                    & \qquad \cdots \otimes
    f \left( l^{n-1}_{(2)} \right) \otimes
    \ul{\corest{f} \left( \qdr \left( l^{n-1}_{(3)} \otimes \ul{v_n} \otimes l^n_{(1)} \right) \right)} \otimes
    f \left( l^n_{(2)} \right)
    \\
    ={}                             & (-1)^{\delta} \,
    f \left( l^0_{(1)} \right) \otimes
    \ul{\corest{f} \left( l^0_{(2)} \otimes v_1 \otimes l^1_{(1)} \right)} \otimes
    f \left( l^1_{(2)} \right) \otimes
    \ul{\corest{f} \left( l^1_{(3)} \otimes v_2 \otimes l^2_{(1)} \right)} \otimes
    f \left( l^2_{(2)} \right) \otimes
    \\
                                    & \qquad \cdots \otimes
    f \left( l^{n-1}_{(2)} \right) \otimes
    \ul{\corest{f} \left( l^{n-1}_{(3)} \otimes v_n \otimes l^n_{(1)} \right)} \otimes
    f \left( l^n_{(2)} \right)
  \end{aligned}
\end{equation*}
where the sign $\delta$ comes from applying $\qdr$ and is given explicitly by
\begin{equation*}
  \delta = \braid{\left( -1, 0 \right)}{\left( 0, \degb{l^0_{(2)}} + \sum_{i=1}^{n-1} \degb{l^{i}_{(3)}} \right)}
  =
  \begin{cases}
    0                                                      & \braidop = \braidop_1, \\
    \degb{l^0_{(2)}} + \sum_{i=1}^{n-1} \degb{l^{i}_{(3)}} & \braidop = \braidop_2.
  \end{cases}
\end{equation*}
The action of $\cycl{\indmap{f}}$ coincides with the action of $\cycl{f}$ on $\ndf{V}[0][] = \tens{V}$.
When $n \geq 1$ and $x$ has the form
$x = \ul{v_1} \otimes l^1 \otimes \dots \otimes \ul{v_n} \otimes l^{n}$,
i.e., $x$ starts with an element from $\ul{V}$, then
\begin{align}\label{eq:cindmap-f-explicit-action}
  \cindmap{f} \left( x \right) ={} & (-1)^{\varepsilon'} \,
  \corest{\indmap{f}} \left( l^n_{(2)} \otimes \ul{v_1} \otimes l^1_{(1)} \right) \otimes
  \indmap{f} \left(
  l^1_{(2)} \otimes \ul{v_2} \otimes l^2 \otimes \dots \otimes \ul{v_{n-1}} \otimes l^{n-1} \otimes
  \ul{v_n} \otimes l^n_{(1)} \right)
  \notag
  \\
  ={}                              & (-1)^{\varepsilon} \,
  \ul{\corest{f} \left( l^n_{(3)} \otimes v_1 \otimes l^1_{(1)} \right)} \otimes
  f \left( l^1_{(2)} \right) \otimes
  \ul{\corest{f} \left( l^1_{(3)} \otimes v_2 \otimes l^2_{(1)} \right)} \otimes
  f \left( l^2_{(2)} \right) \otimes
  \notag
  \\
                                   & \qquad \cdots \otimes
  f \left( l^{n-1}_{(2)} \right) \otimes
  \ul{\corest{f} \left( l^{n-1}_{(3)} \otimes v_n \otimes l^n_{(1)} \right)} \otimes
  f \left( l^n_{(2)} \right)
\end{align}
where the sign $\varepsilon$, which comes from permuting the factors using $\t$ and applying $\qdr$,
is given explicitly by
\begin{equation*}
  \begin{aligned}
    \varepsilon & =
    \underbrace{
      \braid{\left( 0, \degb{l^n_{(3)}} \right)}
      {\left( n, \degb{ v_1 \otimes l^1 \otimes \dots \otimes v_{n} \otimes l^n_{(1)} \otimes l^n_{(2)} } \right)}
    }_{\varepsilon'}
    +
    \braid{\left( -1, 0 \right)}{\left( 0, \sum_{i=1}^{n} \degb{l^{i}_{(3)}} \right)}
    \\
                & =
    \begin{cases}
      \degb{l^n_{(3)}} \cdot \left( \sum_{j=1}^{n-1} \left( \degb{v_j} + \degb{l^j} \right) +
      \degb{v_n} + \degb{l^n_{(1)}} + \degb{l^n_{(2)}} \right) &
      \braidop = \braidop_1,
      \\
      \degb{l^n_{(3)}} \cdot \left( \sum_{j=1}^{n-1} \left( \degb{v_j} + \degb{l^j} \right) +
      \degb{v_n} + \degb{l^n_{(1)}} + \degb{l^n_{(2)}} - n \right) +
      \sum_{i=1}^n \degb{l^i_{(3)}}                            &
      \braidop = \braidop_2.
    \end{cases}
  \end{aligned}
\end{equation*}

\subsection{Compatibility with Change of Ground Algebras} \label{subsec:ndf-compatibility-ground-algebras}
Finally, we discuss briefly the compatibility of all our constructions with change of ground algebras.
Let $\varphi \colon R \rightarrow S$ be a morphism of graded Banach $\mathbbm{k}$-algebras
and let $W$ be a graded Banach $S$-module. Consider the canonical morphism
$\resover{\varphi}^W \colon \tens{\varphi^{*} \left( W \right)}[R] \rightarrow \tens{W}[S]$
over $\varphi$ given by \cref{eq:canonical-pullback-map}. Then one can verify directly that
the induced morphism
\begin{equation*}
  \indmap{F} \left( \resover{\varphi}^{W} \right) \colon
  \tens{\varphi^{*} \left( W \right) \oplus \ul{ \varphi^{*} \left( W \right)}}[R] \rightarrow
  \tens{W \oplus \ul{W}}[S]
\end{equation*}
coincides with the canonical morphism
$\resover{\varphi}^{W \oplus \ul{W}} \colon \tens{\varphi^{*} \left( W \oplus \ul{W} \right)}[R] \rightarrow
  \tens{W \oplus \ul{W}}[S]$.\footnote{Here, we implicitly use the natural identifications
  $\varphi^{*} \left( W \oplus \ul{W} \right) \cong \varphi^{*} \left( W \right) \oplus \varphi^{*} \left( \ul{W} \right) \cong
    \varphi^{*} \left( W \right) \oplus \ul{\varphi^{*} \left( W \right)}$.}
We can write this identity succinctly as
\begin{equation} \label{eq:ind-map-commutes-overline}
  \indmap{F} \left( \resover{\varphi} \right) = \resover{\varphi},
\end{equation}
omitting the explicit dependence on the underlying $S$-modules.

\begin{lm}
  Let $f \colon \tens{V}[R] \rightarrow \tens{W}[S]$ be a morphism of graded Banach coalgebras over $\varphi$.
  Then we have
  \begin{equation*}
    \indmap{F} \left( \rescoho{f} \right) = \rescoho{ \indmap{F} \left( f \right) }.
  \end{equation*}
\end{lm}
\begin{proof}
  Decompose $f$ uniquely as $f = \resover{\varphi} \circ \rescoho{f}$ and apply
  functoriality (\cref{lm:ind-map-functoriality}) to obtain the decomposition
  \begin{equation*}
    \indmap{F} \left( f \right) =
    \indmap{F} \left( \resover{\varphi} \right) \circ \indmap{F} \left( \rescoho{f} \right)
    \stackrel{\eqref{eq:ind-map-commutes-overline}}{=}
    \resover{\varphi} \circ \indmap{F} \left( \rescoho{f} \right).
  \end{equation*}
  Since $\rescoho{ \indmap{F} \left( f \right) }$ is the unique morphism for which
  $\indmap{F} \left( f \right) = \resover{\varphi} \circ \rescoho{ \indmap{F} \left( f \right) }$,
  the identity follows by uniqueness.
\end{proof}

\begin{lm}
  Let $f \colon \tens{W_1}[S] \rightarrow \tens{W_2}[S]$ be a morphism of graded Banach $S$-coalgebras. Then we have
  \begin{equation*}
    \resover{\varphi} \circ \indmap{F} \left( \varphi^{*} \left( f \right) \right) =
    \indmap{F} \left( f \right) \circ \resover{\varphi},
    \qquad
    \indmap{F} \left( \varphi^{*} \left( f \right) \right) =
    \varphi^{*} \left( \indmap{F} \left( f \right) \right).
  \end{equation*}
\end{lm}
\begin{proof}
  Applying $\indmap{F}$ to the identity
  $\resover{\varphi} \circ \varphi^{*} \left( f \right) = f \circ \resover{\varphi}$ we obtain
  \begin{equation*}
    \resover{\varphi} \circ \indmap{F} \left( \varphi^{*} \left( f \right) \right)
    \stackrel{\eqref{eq:ind-map-commutes-overline}}{=}
    \indmap{F} \left( \resover{\varphi} \right) \circ \indmap{F} \left( \varphi^{*} \left( f \right) \right)
    =
    \indmap{F} \left( f \right) \circ \indmap{F} \left( \resover{\varphi} \right)
    \stackrel{\eqref{eq:ind-map-commutes-overline}}{=}
    \indmap{F} \left( f \right) \circ \resover{\varphi}.
  \end{equation*}
  Since $\varphi^{*} \left( \indmap{F} \left( f \right) \right)$ is the unique morphism
  for which
  $\resover{\varphi} \circ \varphi^{*} \left( \indmap{F} \left( f \right) \right) =
    \indmap{F} \left( f \right) \circ \resover{\varphi}$, the identity follows by uniqueness.
\end{proof}

\begin{lm}
  Denote by $\qdr^W \colon \tens{W \oplus \ul{W}}[S] \rightharpoonup \tens{W \oplus \ul{W}}[S]$ the
  de Rham differential on $\tens{W \oplus \ul{W}}[S]$ and by
  $\qdr^{\varphi^{*} \left( W \right)} \colon
    \tens{\varphi^{*} \left( W \right) \oplus \ul{ \varphi^{*} \left( W \right)}}[R] \rightharpoonup
    \tens{\varphi^{*} \left( W \right) \oplus \ul{ \varphi^{*} \left( W \right)}}[R]$ the de Rham differential
  on $\tens{\varphi^{*} \left( W \right) \oplus \ul{\varphi^{*} \left( W \right)}}[R]$.
  Then we have
  \begin{equation*}
    \resover{\varphi} \circ \qdr^{\varphi^{*} \left( W \right)} = \qdr^W \circ \resover{\varphi},
    \qquad
    \varphi^{*} \left( \qdr^W \right) = \qdr^{\varphi^{*} \left( W \right)}.
  \end{equation*}
\end{lm}
\begin{proof}
  We have
  \begin{equation*}
    \resover{\varphi} \circ \qdr^{\varphi^{*} \left( W \right)}
    \stackrel{\eqref{eq:ind-map-commutes-overline}}{=}
    \indmap{F} \left( \resover{\varphi} \right) \circ \qdr^{\varphi^{*} \left( W \right)}
    \stackrel{\eqref{eq:func-qdr}}{=}
    \qdr^W \circ \indmap{F} \left( \resover{\varphi} \right)
    \stackrel{\eqref{eq:ind-map-commutes-overline}}{=}
    \qdr^W \circ \resover{\varphi}.
  \end{equation*}
  Since $\varphi^{*} \left( \qdr^W \right)$ is the unique coderivation which satisfies
  $\resover{\varphi} \circ \varphi^{*} \left( \qdr^W \right) = \qdr^W \circ \resover{\varphi}$,
  the identity follows by uniqueness.
\end{proof}

\begin{lm} \label{lm:pullback-cont}
  Let $\nu \colon \tens{W}[S] \rightharpoonup \tens{W}[S]$ be an $S$-linear coderivation.
  Then we have
  \begin{equation*}
    \resover{\varphi} \circ \cont{\varphi^{*} \left( \nu \right)} = \cont{\nu} \circ \resover{\varphi},
    \qquad
    \varphi^{*} \left( \cont{\nu} \right) = \cont{\varphi^{*} \left( \nu \right)}.
  \end{equation*}
\end{lm}
\begin{proof}
  The coderivation $\varphi^{*} \left( \nu \right)$ is the unique coderivation which satisfies
  $\resover{\varphi} \circ \varphi^{*} \left( \nu \right) = \nu \circ \resover{\varphi}$. Hence,
  \begin{equation*}
    \resover{\varphi} \circ \cont{\varphi^{*} \left( \nu \right)}
    \stackrel{\eqref{eq:ind-map-commutes-overline}}{=}
    \indmap{F} \left( \resover{\varphi} \right) \circ \cont{\varphi^{*} \left( \nu \right)}
    \stackrel{\eqref{eq:func-cont}}{=}
    \cont{\nu} \circ \indmap{F} \left( \resover{\varphi} \right)
    \stackrel{\eqref{eq:ind-map-commutes-overline}}{=}
    \cont{\nu} \circ \resover{\varphi}.
  \end{equation*}
  Since $\varphi^{*} \left( \cont{\nu} \right)$ is the unique coderivation which satisfies
  $\resover{\varphi} \circ \varphi^{*} \left( \cont{\nu} \right) = \cont{\nu} \circ \resover{\varphi}$,
  the identity follows by uniqueness.
\end{proof}

\begin{lm}
  Let $d_R \colon R \rightharpoonup R, d_S \colon S \rightharpoonup S$ be two derivations
  such that $\varphi \circ d_R = d_S \circ \varphi$. Let
  $\nu \colon \tens{W}[S] \rightharpoonup \tens{W}[S]$ be a generalized coderivation over $d_S$.
  Then we have
  \begin{equation*}
    \resover{\varphi} \circ \lie{\varphi^{*} \left( \nu \right)} = \lie{\nu} \circ \resover{\varphi},
    \qquad
    \varphi^{*} \left( \lie{\nu} \right) = \lie{\varphi^{*} \left( \nu \right)}.
  \end{equation*}
\end{lm}
\begin{proof}
  The coderivation $\varphi^{*} \left( \nu \right)$ is the unique coderivation which satisfies
  $\resover{\varphi} \circ \varphi^{*} \left( \nu \right) = \nu \circ \resover{\varphi}$. Hence,
  \begin{equation*}
    \resover{\varphi} \circ \lie{\varphi^{*} \left( \nu \right)}
    \stackrel{\eqref{eq:ind-map-commutes-overline}}{=}
    \indmap{F} \left( \resover{\varphi} \right) \circ \lie{\varphi^{*} \left( \nu \right)}
    \stackrel{\eqref{eq:func-lie}}{=}
    \lie{\nu} \circ \indmap{F} \left( \resover{\varphi} \right)
    \stackrel{\eqref{eq:ind-map-commutes-overline}}{=}
    \lie{\nu} \circ \resover{\varphi}.
  \end{equation*}
  Since $\varphi^{*} \left( \lie{\nu} \right)$ is the unique coderivation which satisfies
  $\resover{\varphi} \circ \varphi^{*} \left( \lie{\nu} \right) = \lie{\nu} \circ \resover{\varphi}$,
  the identity follows by uniqueness.
\end{proof}

\begin{rem}
  Since the cyclization construction is compatible with the change of ground algebras
  (see \cref{fig:naturality-cyclization-coderivations,fig:naturality-cyclization-morphisms} of \cref{sec:cyclization-naturality}),
  one can obtain analogous results for the cyclic versions of the various operators by applying cyclization and using the fact
  that $\cycl{\resover{\varphi}} = \resover{\varphi}$. More explicitly, we have the identity
  \begin{equation*}
    \resover{\varphi} \circ \cindmap{F} \left( \varphi^{*} \left( f \right) \right) =
    \cindmap{F} \left( f \right) \circ \resover{\varphi}
  \end{equation*}
  which holds on $\ndf{\varphi^{*} \left( W_1 \right)}[][]$ and the identities
  \begin{equation*}
    \qquad
    \resover{\varphi} \circ \ccont{\varphi^{*} \left( \nu \right)} = \ccont{\nu} \circ \resover{\varphi},
    \qquad
    \resover{\varphi} \circ \clie{\varphi^{*} \left( \nu \right)} = \clie{\nu} \circ \resover{\varphi},
  \end{equation*}
  which hold on $\ndf{\varphi^{*} \left( W \right)}[][]$. The identities above hold on $\ndf{\cdot}[][]$, even before descending
  to the cyclic quotients $\ncdf{\cdot}[][]$.
\end{rem}

\begin{rem}
  The results above describe how noncommutative codifferential forms behave with respect to
  scalar restriction. There are also analogous results for the process of scalar extension.
  Let $V$ be a graded Banach $R$-module and let
  $\resunder{\varphi} \colon \tens{V}[R] \rightarrow \tens{\varphi_{!} \left( V \right)}[S]$
  be the canonical morphism over $\varphi$ (see \cref{eq:canonical-extension-map}).
  We have a natural isomorphism
  \begin{equation*}
    \gamma \colon \varphi_{!} \left( V \oplus \ul{V} \right) \rightarrow \varphi_{!} \left( V \right) \oplus \ul{\varphi_{!} \left( V \right)}
  \end{equation*}
  given by
  \begin{equation*}
    \gamma \left( s \otimes_R v \right) = s \otimes_R v, \qquad
    \gamma \left( s \otimes_R \ul{v} \right) = (-1)^{\braid{(1,0)}{(0,\degb{s})}} \ul{s \otimes_R v}
  \end{equation*}
  which induces an isomorphism
  \begin{equation*}
    \Gamma = \tens{\gamma}[S] \colon \tens{\varphi_{!} \left( V \oplus \ul{V} \right)}[S] \rightarrow
    \tens{\varphi_{!} \left( V \right) \oplus \ul{\varphi_{!} \left( V \right)}}[S]
  \end{equation*}
  between the corresponding tensor algebras.

  Given an $R$-linear coderivation $\mu$ on $\tens{V}[R]$,
  the analogous identity to $\varphi^{*} \left( \cont{\nu} \right) = \cont{\varphi^{*} \left( \nu \right)}$
  of \cref{lm:pullback-cont} is
  $\Gamma \circ \varphi_{!} \left( \cont{\mu} \right) = \cont{\varphi_{!} \left( \mu \right)} \circ \Gamma$,
  where $\Gamma$ is inserted so that the domains and codomains of both sides coincide.
  We leave the interested reader to formulate and prove all other analogous results which might be of interest.
\end{rem}

\section{Complexes for Cyclic Homology} \label{sec:cyclic-homology-models}

In this section, we define the cyclic homology of a curved Banach $\Ainf$-algebra using
several different explicit complexes and establish their equivalence.
We begin in \cref{sec:bar-hoch-cyc-complexes} by defining the bar complex, Hochschild complex and Connes' cyclic
complex of a curved $\Ainf$-algebra $\mathcal{A}$ over a differential graded-commutative Banach $\mathbbm{k}$-algebra $\mathcal{R}$,
and the corresponding homology theories.
We show in \cref{sec:bar-homology-unital-curved-a-inf-alg} that the bar homology of a \textit{unital}
Banach $\Ainf$-algebra is trivial even in the presence of curvature,
provided we work with the complete direct sum complex.
In \cref{subsec:cyc-bicomplex}, we introduce the cyclic bicomplex $\cycbi[A]$
and prove that its complete direct sum totalization is homotopy equivalent to Connes' complex $\cconnes{\mathcal{A}}[]$.
The proof provides an explicit homotopy inverse to the natural projection $\cycbi[A] \rightarrow \cconnes{\mathcal{A}}[]$,
constructed from the contractions of the rows of $\cycbi[A]$. Since we work with infinite sums, this
involves verifying an appropriate convergence condition, naturally stated in terms of
an operator $S$, called the periodicity operator.

The operator $S$ is a \textit{chain-level map} defined on the Hochschild complex and its properties
are studied in greater depth in \cref{subsec:periodicity-operator}. There, we prove various identities
that show that $S$ decomposes as a composition of two boundary operators and descends to Connes' cyclic complex.
We also prove that $S$ is pointwise topologically nilpotent, which implies the convergence condition used
in \cref{subsec:cyc-bicomplex}.
In \cref{sec:connes-exact-sequence}, we derive Connes' long exact sequence for a unital $\mathcal{A}$
from the cyclic bicomplex and show that it is equivalent to the $SBI$ long exact sequence,
expressed in terms of the homology of the Hochschild and Connes' cyclic complexes. We verify
that the map induced by $S$ can be identified, on the level of homology, with the periodicity map $\sigma$
coming from the cyclic bicomplex.
In \cref{sec:three-more-models}, we introduce three more models for cyclic homology
and show they are all homotopy equivalent with explicit chain maps which induce the equivalences. We also
describe how the periodicity operator and Connes' exact sequence are realized in each of the models.

In \cref{sec:bicomplex-models-cyclic-homology}, we use the formalism of cyclic codifferential forms
from \cref{sec:noncomm-diff-calc} to introduce two more models
$\totcomp{\mathcal{A}}[1][]$ and $\totcomp{\mathcal{A}}[2][]$
for cyclic homology.
The model $\totcomp{\mathcal{A}}[1][]$ (resp.\ $\totcomp{\mathcal{A}}[2][]$) is obtained
by taking the (complete direct sum) total complex of
cyclic codifferential forms of degree greater than or equal to one (resp.\ greater than or equal to two).
These total complexes are shown to be homotopy equivalent to Connes' cyclic complex, and we provide
explicit chain maps realizing the equivalences. We also show that the periodicity operator
acting between the models $\totcomp{\mathcal{A}}[1][]$ and $\totcomp{\mathcal{A}}[2][]$
can be realized simply as the projection.

Finally, we conclude the section by defining extended and reduced versions of the cyclic complexes.
The extended version is useful in the presence of curvature while the reduced version gives a smaller complex in the presence of a unit.

For the rest of this section, we fix a field $\mathbbm{k}$ of characteristic zero, endowed with the trivial norm.
Throughout this section, in accordance with \cref{sec:conv-and-not}, unadorned tensor products and direct sums are
understood in the graded Banach sense (and hence completed); completion marks are displayed only for emphasis or
when contrasting with an explicitly algebraic construction.

\subsection{The Bar, Hochschild and Connes Complexes} \label{sec:bar-hoch-cyc-complexes}

Let $\mathcal{R} = (R,d)$ be a differential graded-commutative Banach $\mathbbm{k}$-algebra. Let
$\mathcal{A} = (A, \mu)$ be a Banach $\Ainf$-algebra over $\mathcal{R}$.
To make our notation consistent with the standard notation appearing in the literature, set
\begin{equation*}
	b' \defeq \rest{\mu}{\tensr{A}} \colon \tensr{A} \rightharpoonup \tensr{A}, \quad
	b \defeq \rest{\cycl{\mu}}{\tensr{A}} \colon \tensr{A} \rightharpoonup \tensr{A}.
\end{equation*}

We have $b'^2 = 0$ by the $\Ainf$-relations and $b^2 = 0$ by \cref{cor:odd-coder-cyc-differential}
so that $b'$ and $b$ are differentials. The differential $b'$ is called the \textbf{bar differential} and the differential $b$ is called the \textbf{Hochschild differential}.
In addition, we have the following relations:
\begin{align}
	(\idd - \t) \N & = \N (\idd - \t) = 0, \label{eq:1-t-N=0}                                   \\
	b(\idd - \t)   & \stackrel{\eqref{eq:cycl-mu-idd-t}}{=} (\idd - \t)b', \label{eq:b-1-t-rel} \\
	b' \N          & \stackrel{\eqref{eq:mu-N-N-cycl-mu}}{=} \N b \label{eq:b'-N-rel}.
\end{align}

Using the differentials $b, b'$ we can immediately form the following three complexes:
\begin{dfn}
	The complex $\left( \tensr{A}, b' \right)$ is called the \textbf{bar complex} of
	$\mathcal{A}$ and its cohomology is called the \textbf{bar homology} of $\mathcal{A}$.
\end{dfn}

\begin{dfn}
	Assume that $\mathcal{A}$ is unital.\footnote{While the definition of the Hochschild complex does not require a unit,
		the Hochschild homology for non-unital algebras is typically defined differently. See \cite[Section 1.4]{Loday1998}.}
	The complex $\choch{\mathcal{A}}[] \defeq \left( \tensr{A}, b \right)$ is called the \textbf{Hochschild complex}
	of $\mathcal{A}$ and its cohomology $\hhoch{\mathcal{A}} \defeq \cohom{\tensr{A}}[*][b]$
	is a graded $\cohom{\mathcal{R}}$-module called the \textbf{Hochschild homology} of $\mathcal{A}$.
\end{dfn}

\begin{dfn}
	The complex $\cconnes{\mathcal{A}}[] \defeq \left( \tensrcyc{A}, b \right) =
		\left( \tensr{A} / \Im \left( \idd - \t \right), b \right)$
	is called the \textbf{Connes complex} of $\mathcal{A}$ and its cohomology
	$\hcyc{\mathcal{A}} \defeq \cohom{{\cconnes{\mathcal{A}}[]}}$ is a graded $\cohom{\mathcal{R}}$-module
	called the \textbf{cyclic homology} of $\mathcal{A}$.
\end{dfn}

The underlying graded $R$-modules of $\choch{\mathcal{A}}[]$ and $\cconnes{\mathcal{A}}[]$ will be denoted by
$\choch{A}[]$ and $\cconnes{A}[]$ respectively, emphasizing that they depend only on the underlying graded $R$-module $A$ of
$\mathcal{A}$ and not on the $\Ainf$-structure $\mu$.

\begin{rem}
	Since our $\Ainf$-algebras are of cohomological type, that is, the operator $\mu_1$ has degree $1$ and not
	$-1$, the Hochschild differential $b$ has degree $1$.
	In this setup, it is natural to allow the Hochschild and cyclic
	\textit{homology} to have \textit{cohomological} grading like we do here.
\end{rem}

When $\mathcal{A}$ corresponds to a DGA, our definitions of the Hochschild and cyclic complexes
and homology theories coincide with the standard definitions up to a shift and an identification.
See \cref{appendix:sign-conversions} for details.

\subsection{The Bar Homology of a Unital Curved Banach \texorpdfstring{$\Ainf$}{A-infinity}-algebra}
\label{sec:bar-homology-unital-curved-a-inf-alg}
Let $\mathcal{A} = ( A, \mu, e )$ be a \textit{unital} Banach $\Ainf$-algebra over $\mathcal{R}$.
The purpose of this subsection is to show that the cohomology of the bar complex $( \tensr{A}, b' )$
is trivial even in the presence of curvature. In order for the result to hold, it is necessary
to work with the complete tensor coalgebra in order to guarantee that various infinite sums are well-defined.

In what follows, it will be useful to split the differential $b'$ as a sum $b' = b'_0 + b'_{>0}$
where the operator $b'_0$ involves only the curvature term $\mu_0(1)$ and $b'_{>0}$ involves
all the operators $\mu_k$ for $k \geq 1$.\footnote{Note that even though $b' = \mu$, the notation $b'_0$ in what
	follows \textit{does not} denote the corestriction $\mu_0 \colon R \rightharpoonup A$.}
More precisely, let $b'_0 \colon \tensr{A} \rightharpoonup \tensr{A}$
be given by
\begin{equation*}
	b'_0 \left( a_1 \otimes \dots \otimes a_k \right) \defeq
	\sum_{i=0}^{k} (-1)^{\degb{a_1} + \dots + \degb{a_i}}
	a_1 \otimes \dots \otimes a_{i} \otimes \mu_0(1) \otimes \dots \otimes a_k
\end{equation*}
and let $b'_{>0} \defeq b' - b'_0$.
The map $b'_0$ is an $R$-linear operator of degree $1$ and satisfies $b'_{0} \circ b'_{0} = 0$.
In contrast, $b'_{>0}$ is a derivation over $d$ of degree $1$ which, in general, does not satisfy $b'_{>0} \circ b'_{>0} = 0$.
Instead, we have $b'_{> 0} \circ b'_{> 0} + b'_{> 0} \circ b'_0 + b'_0 \circ b'_{>0} = 0$.

When $\mathcal{A}$ is a unital associative algebra, it is well known that the bar complex is contractible
(see \cite[Page 12]{Loday1998}). The standard argument generalizes verbatim to unital uncurved $\Ainf$-algebras.
More precisely, we have:

\begin{lm}
	Let $H_0 \colon \tensr{A} \rightharpoonup \tensr{A}$ be the degree $-1$ map given by
	\begin{equation}
		H_0 \left( a_1 \otimes \dots \otimes a_k \right) \defeq e \otimes a_1 \otimes \dots \otimes a_k.
		\label{eq:contraction-bar-no-mu-0}
	\end{equation}
	Then we have the identity $b'_{>0} \circ H_0 + H_0 \circ b'_{>0} = \idd$.
\end{lm}
\begin{proof}
	Let $a_1, \dots, a_k \in A$. Since $\mathcal{A}$ is unital, we have
	\begin{equation*}
		\begin{aligned}
			\left( b'_{>0} \circ H_0 \right) \left( a_1 \otimes \dots \otimes a_k \right) ={} &
			b'_{>0} \left( e \otimes a_1 \otimes \dots \otimes a_k \right)
			\\
			={}                                                                               &
			\sum_{i=0}^k \mu_{i+1} \left( e \otimes a_1 \otimes \dots \otimes a_i \right) \otimes a_{i+1} \otimes
			\dots \otimes a_k
			\\
			                                                                                  & + (-1)^{\degb{e}} e \otimes b'_{>0} \left( a_1 \otimes \dots \otimes a_k \right)
			\\
			={}                                                                               &
			\mu_2 \left( e \otimes a_1 \right) \otimes a_2 \otimes \dots \otimes a_k -
			e \otimes b'_{>0} \left( a_1 \otimes \dots \otimes a_k \right)
			\\
			={}                                                                               & a_1 \otimes \dots \otimes a_k - e \otimes b'_{>0} \left( a_1 \otimes \dots \otimes a_k \right)
			\\
			={}                                                                               & a_1 \otimes \dots \otimes a_k -
			\left( H_0 \circ b'_{>0} \right) \left( a_1 \otimes \dots \otimes a_k \right).
		\end{aligned}
	\end{equation*}
\end{proof}

\begin{cor}
	Let $\mathcal{A} = \left( A, \mu, e \right)$ be a unital Banach $\Ainf$-algebra over $\mathcal{R} = (R,d)$
	with $\mu_0(1) = 0$. Then the bar complex $\left( \tensr{A}, b' \right)$ is contractible with
	contracting homotopy $H_0$ given by \cref{eq:contraction-bar-no-mu-0}.
\end{cor}
\begin{proof}
	When $\mathcal{A}$ is uncurved, we have $b' = b'_{>0}$ and hence $b' H_0 + H_0 b' = \idd$.
\end{proof}

When $\mathcal{A}$ is curved, the map defined by \cref{eq:contraction-bar-no-mu-0} fails to be a contraction.
Nonetheless, as already noted in \cite[Lemma 3.2]{Cho2012},
it turns out that even in the curved case, the bar complex is contractible.
To see this, we will think of $\left( \tensr{A}, b' \right) = \left( \tensr{A}, b'_{>0} + b'_{0} \right)$
as a perturbation of the pre-complex $\left( \tensr{A}, b'_{>0} \right)$ and perturb $H_0$ to a
homotopy $H$:

\begin{lm} \label{lm:bar-homology-curved-algebra}
	Let $\mathcal{A} = \left( A, \mu, e \right)$ be a unital Banach $\Ainf$-algebra over $\mathcal{R} = (R,d)$.
	Then the bar complex $\left( \tensr{A}, b' \right)$ is contractible.
	The operator $H \colon \tensr{A} \rightharpoonup \tensr{A}$ given by
	\begin{equation}
		H \left( a_1 \otimes \dots \otimes a_k \right) \defeq
		\sum_{n=0}^{\infty} (-1)^n
		e \otimes \left( \mu_0(1) \otimes e \right)^{\otimes n} \otimes a_1 \otimes \dots \otimes a_k
		\label{eq:contraction-bar-with-mu-0}
	\end{equation}
	satisfies
	\begin{equation}
		b' H + H b' = \idd, \label{eq:b'H-homotopy}
	\end{equation}
	providing a contraction of the bar complex.
\end{lm}
\begin{proof}
	Split the differential $b'$ as $b' = b'_{0} + b'_{>0}$. We want to apply \cref{lm:trivial-perturbation-lemma}
	with $b = b'_{>0}, \delta = b'_{0}$ and $h = H_0$ where $H_0$ is given by \cref{eq:contraction-bar-no-mu-0}.
	Given $a_1, \dots, a_k \in A$ with
	\begin{equation*}
		x = a_1 \otimes \dots \otimes a_k,
	\end{equation*}
	we have
	\begin{equation*}
		\begin{aligned}
			\left( b'_{0} \circ H_0 \right) \left( x \right) ={} &
			b'_{0} \left( e \otimes a_1 \otimes \dots \otimes a_k \right)
			\\
			={}                                                  &
			\mu_0(1) \otimes e \otimes a_1 \otimes \dots \otimes a_k
			\\
			                                                     & +
			\sum_{i=0}^k (-1)^{\degb{e} + \degb{a_1} + \dots + \degb{a_i}}
			e \otimes a_1 \otimes \dots \otimes a_i \otimes \mu_0(1) \otimes a_{i+1} \otimes \dots \otimes a_k
			\\
			={}                                                  &
			\mu_0(1) \otimes e \otimes a_1 \otimes \dots \otimes a_k
			\\
			                                                     & -
			\sum_{i=0}^k (-1)^{\degb{a_1} + \dots + \degb{a_i}}
			e \otimes a_1 \otimes \dots \otimes a_i \otimes \mu_0(1) \otimes a_{i+1} \otimes \dots \otimes a_k,
		\end{aligned}
	\end{equation*}
	and
	\begin{equation*}
		\left( H_0 \circ b'_{0} \right) \left( x \right) =
		\sum_{i=0}^k (-1)^{\degb{a_1} + \dots + \degb{a_i}}
		e \otimes a_1 \otimes \dots \otimes a_i \otimes \mu_0(1) \otimes a_{i+1} \otimes \dots \otimes a_k.
	\end{equation*}
	Hence,
	\begin{equation*}
		\left( b'_{0} \circ H_0 + H_0 \circ b'_{0} \right) \left( x \right) = \mu_0(1) \otimes e \otimes x.
	\end{equation*}

	Since $\nnorm[\mu_0(1)] < 1$ and $\nnorm[e] \leq 1$, we have $\nnorm[ b'_{0} H_0 + H_0 b'_{0}] < 1$ and so
	by \cref{lm:inverse-geometric-series} the map $\idd + b'_{0} H_0 + H_0 b'_{0}$ is invertible with inverse
	given explicitly by
	\begin{equation*}
		\left( \idd + b'_{0} H_0 + H_0 b'_{0} \right)^{-1} \left( x \right) =
		\left( \sum_{n=0}^{\infty} (-1)^n \left( b'_{0} H_0 + H_0 b'_{0} \right)^n \right) \left( x \right) =
		\sum_{n=0}^{\infty} (-1)^n \left( \mu_0(1) \otimes e \right)^{\otimes n} \otimes x.
	\end{equation*}
	Hence, by \cref{lm:trivial-perturbation-lemma} we get that $\left( \tensr{A}, b' \right)$ is contractible
	with contracting homotopy
	$H = H_0 \left( \idd +  b'_{0} H_0 + H_0 b'_{0} \right)^{-1}$ given explicitly by \cref{eq:contraction-bar-with-mu-0}.
\end{proof}

\begin{rem}
	The contraction in \cref{lm:bar-homology-curved-algebra} involves an infinite sum which makes
	sense since we implicitly work with the \textit{completed bar complex}
	$\tensr{A} = \coplus_{i=1}^{\infty} A^{\cotimes i}$. In fact, in the presence of curvature,
	\cref{lm:bar-homology-curved-algebra} generally fails for the standard bar complex
	$\oplus_{i=1}^{\infty} A^{\otimes i}$. See \cref{appendix:bar-complex-not-necessarily-contractible}.
\end{rem}

Now consider the full complex $\left( \tens{A}, b' \right)$ which we will call
the \textbf{extended bar complex}. When $\mu_0(1) = 0$, the
inclusion $R \hookrightarrow \tens{A}$ is a chain map which induces an isomorphism on cohomology. However, when
$\mu_0(1) \neq 0$, the inclusion is not a chain map. Nonetheless, we can show that
$\cohom{\tens{A}}[*][b'] \cong \cohom{R}[*][d]$.

\begin{lm}
	Let $\mathcal{A} = \left( A, \mu, e \right)$ be a unital Banach $\Ainf$-algebra over $\mathcal{R} = (R,d)$.
	Let $p \colon \tens{A} \rightarrow R$ be the natural projection map and let
	$i \colon R \rightarrow \tens{A}$ be the unique degree zero $R$-linear map which satisfies
	\begin{equation*}
		i \left( 1 \right) =
		\left( \sum_{k=0}^{\infty} (-1)^k \left( e \otimes \mu_0(1) \right)^{\otimes k} \right) =
		1 - e \otimes \mu_0(1) + e \otimes \mu_0(1) \otimes e \otimes \mu_0(1) - \dots .
	\end{equation*}
	Then $p$ is a homotopy equivalence with a homotopy inverse $i$. In particular, we have
	$\cohom{\tens{A}}[*][b'] \cong \cohom{R}[*][d]$.
\end{lm}
\begin{proof}
	Let $f \colon R[-1] \rightarrow \tensr{A}$ be the degree zero morphism of DG modules given by
	$f \left( s_{-1} \left( r \right) \right) = (-1)^{\degb{r}} r \cdot \mu_0(1)$. Then we have
	$\Cone{f} = R[-1][1] \oplus \tensr{A} = \tens{A}$ and the differential on $\Cone{f}$ coincides with
	$b'$. Applying \cref{cor:projection-from-mapping-cone-equivalence} and taking into account
	formula \eqref{eq:contraction-bar-with-mu-0} for $H$, we see that the projection $p$ is
	a homotopy equivalence with a homotopy inverse given by $i$.
\end{proof}

\subsection{The Cyclic Bicomplex} \label{subsec:cyc-bicomplex}
Another commonly used model for defining cyclic homology is given by the total complex of a $2$-periodic right half plane bicomplex $\cycbi[A]$ called the \textbf{cyclic bicomplex} shown in \cref{fig:standard-cyclic-bicomplex}.
\begin{figure}[H]
	\centering
	\begin{tikzcd}
		& \vdots & \vdots & \vdots & \vdots \\
		\cdots & 0 & {\tensr{A}^1} & {\tensr{A}^1} & {\tensr{A}^1} & \cdots \\
		\cdots & 0 & {\tensr{A}^0} & {\tensr{A}^0} & {\tensr{A}^0} & \cdots \\
		\cdots & 0 & {\tensr{A}^{-1}} & {\tensr{A}^{-1}} & {\tensr{A}^{-1}} & \cdots \\
		& \vdots & \vdots & \vdots & \vdots
		\arrow[from=3-3, to=3-2]
		\arrow["{\idd-\t}"', from=3-4, to=3-3]
		\arrow["\N"', from=3-5, to=3-4]
		\arrow["b", from=3-3, to=2-3]
		\arrow["{b'}", from=3-4, to=2-4]
		\arrow["b", from=3-5, to=2-5]
		\arrow["{\idd - \t}"', from=2-4, to=2-3]
		\arrow["\N"', from=2-5, to=2-4]
		\arrow[from=2-3, to=2-2]
		\arrow[from=2-2, to=2-1]
		\arrow[from=3-2, to=3-1]
		\arrow[from=4-2, to=4-1]
		\arrow[from=3-2, to=2-2]
		\arrow[from=4-2, to=3-2]
		\arrow["b", from=4-3, to=3-3]
		\arrow["{b'}", from=4-4, to=3-4]
		\arrow["b", from=4-5, to=3-5]
		\arrow["{\idd - \t}"', from=4-4, to=4-3]
		\arrow["\N"', from=4-5, to=4-4]
		\arrow[from=4-3, to=4-2]
		\arrow["{ }"{description}, from=2-2, to=1-2]
		\arrow[from=2-3, to=1-3]
		\arrow[from=2-4, to=1-4]
		\arrow[from=2-5, to=1-5]
		\arrow["{\idd - \t}"', from=2-6, to=2-5]
		\arrow["{\idd - \t}"', from=3-6, to=3-5]
		\arrow["{\idd - \t}"', from=4-6, to=4-5]
		\arrow[from=5-3, to=4-3]
		\arrow[from=5-2, to=4-2]
		\arrow[from=5-4, to=4-4]
		\arrow[from=5-5, to=4-5]
	\end{tikzcd}
	\caption{The Cyclic Bicomplex $\cycbi[A]$.}
	\label{fig:standard-cyclic-bicomplex}
\end{figure}

More precisely, let $\cycbi[A]$ be the double complex defined by
\begin{equation*}
	\cycbi \left( A \right)_{i}^{j} \defeq
	\begin{cases}
		0           & i < 0,    \\
		\tensr{A}^j & i \geq 0.
	\end{cases}
\end{equation*}
The vertical and horizontal differentials
\begin{equation*}
	\delta_{\ver} \colon \cycbi \left( A \right)_{*}^{*} \rightharpoonup \cycbi \left( A \right)_{*}^{*+1},
	\qquad
	\delta_{\hor} \colon \cycbi \left( A \right)_{*}^{*} \rightharpoonup \cycbi \left( A \right)_{*-1}^{*}
\end{equation*}
are given by
\begin{equation*}
	\left( \delta_{\ver} \right)_i^j \defeq
	\begin{cases}
		b  & i \geq 0 \textrm{ even}, \\
		b' & i \geq 0 \textrm{ odd}.
	\end{cases},
	\qquad
	\left( \delta_{\hor} \right)_i^j \defeq
	\begin{cases}
		\N        & i \geq 2 \textrm{ even}, \\
		\idd - \t & i \geq 1 \textrm{ odd}.
	\end{cases}
\end{equation*}
Each column $\cycbi[A]_{*}$ of $\cycbi[A]$ is a differential graded $\mathcal{R}$-module whose differential
$\delta_{\ver} \colon \cycbi[A]_{*} \rightharpoonup \cycbi[A]_{*}$ is a degree one derivation over $d$.
The horizontal map $\delta_{\hor}$ satisfies $\delta_{\hor} \left( r c \right) = r \cdot \delta_{\hor} \left( c \right)$
for all $r \in R$ and $c \in \cycbi[A]$ and is a differential by \cref{eq:1-t-N=0}. By \cref{eq:b-1-t-rel,eq:b'-N-rel} the differentials
$\delta_{\ver}$ and $\delta_{\hor}$ commute so we indeed have a Banach bicomplex (see \cref{subsec:bicomplex-commuting-differentials}).

In the characteristic zero setting which we assume, the rows of the cyclic bicomplex
$\cycbi[A]$ are exact except at the $0$-th column. Explicitly, let us define maps $h, h' \colon \tensr{A} \rightarrow \tensr{A}$ by
\begin{equation} \label{eq:def-h-h'}
	\rest{h}{{A^{\otimes \left( n + 1 \right)}}} \defeq -\frac{1}{(n+1)} \cdot \sum_{i=0}^n i {\t}^i,
	\qquad
	\rest{h'}{{A^{\otimes \left( n + 1 \right)}}} \defeq \frac{1}{(n+1)} \cdot \idd.
\end{equation}
Then we have the identities
\begin{align}
	h' \N + (\idd - \t)h               & = \idd, \label{eq:h'N-homotopy} \\
	\N h' + h \left( \idd - \t \right) & = \idd \label{eq:Nh'-homotopy}
\end{align}
which imply that
\begin{align}
	\Im \left( \idd - \t \right) & = \ker \left( \N \right), \label{eq:im-1-t-ker-N}        \\
	\Im \left( \N \right)        & = \ker \left( \idd - \t \right). \label{eq:im-N-ker-1-t}
\end{align}
Hence, we have $\homo{{\cycbi[A]}^j}[i][\delta_{\hor}] = 0$ when $i > 0$ while
\begin{equation*}
	\homo{{\cycbi[A]}^j}[0][\delta_{\hor}] = \tensr{A}^j / \Im \left( \idd - \t \right) = \cconnes{\mathcal{A}}[j]
\end{equation*}
is the cokernel of the first two columns.

When $A$ corresponds to an associative algebra, or, more generally, $A$ is non-positively graded\footnote{Note
	that a non-positively graded $\Ainf$-algebra $A$ cannot be curved as $\mu_0(1) \in A^1$.},
the bicomplex $\cycbi[A]$ is concentrated in the fourth quadrant and hence all possible totalizations
of $\cycbi[A]$ (direct sum, direct product, complete direct sum, etc.) coincide. In this case,
a standard staircase argument shows that the natural surjection
$p \colon \totc{\cycbi[A]}[][] \rightarrow \cconnes{\mathcal{A}}[]$
between $\totc{\cycbi[A]}[][]$ and Connes' complex $\cconnes{\mathcal{A}}[]$
given by the quotient map $\pi \colon \tensr{A} \rightarrow \tensr{A} / \Im \left( \idd - \t \right)$ on the
zeroth column and $0$ on all other columns is a homotopy equivalence.

In the general case, we will work with the complete direct sum totalization
$\totc{\cycbi[A]}[][][\coplus]$ and argue that the natural surjection $p$ is still
a homotopy equivalence.

\begin{thm} \label{lm:cyc-bicomplex-homotopy-equivalent-connes}
	The natural map $p \colon \totc{\cycbi[A]}[][][\coplus] \rightarrow \cconnes{\mathcal{A}}[]$
	is a homotopy equivalence. In particular,
	$p \colon \cohom{{\totc{\cycbi[A]}[][][\coplus]}} \rightarrow \hcyc{\mathcal{A}}$ is an isomorphism.
\end{thm}

To prove \cref{lm:cyc-bicomplex-homotopy-equivalent-connes}, we will need to analyze the behaviour
of an important operator $S$ called the periodicity operator.
\begin{dfn} \label{dfn:periodicity-operator}
	The \textbf{periodicity operator} $S \colon \tensr{A}^{*} \rightharpoonup \tensr{A}^{*+2}$ is defined by the formula
	\begin{equation*}
		S \defeq - \left( h' b' h b + b h' b' h \right) + b h' h b. 
	\end{equation*}
\end{dfn}

We will justify the name ``periodicity operator'' and investigate its properties more thoroughly in
the next subsection. There, we will show that the operator $S$ is  ``almost nilpotent'' in the sense
that $S^n \left( x \right) \rightarrow 0$ for all $x \in \tensr{A}^{*}$ (see \cref{lm:S-almost-nilpotent}).
Assuming for a moment \cref{lm:S-almost-nilpotent}, let us prove \cref{lm:cyc-bicomplex-homotopy-equivalent-connes}.

\begin{figure}[htb]
	\centering
	\begin{tikzcd}
		& \vdots & \vdots & \vdots & \vdots & \vdots \\
		\cdots & 0 & \cconnes{A}[1] & {\tensr{A}^1} & {\tensr{A}^1} & {\tensr{A}^1} &
		{\phantom{\tensr{A}^1}} \\
		\cdots & 0 & \cconnes{A}[0] & {\tensr{A}^0} & {\tensr{A}^0} & {\tensr{A}^0} &
		{\phantom{\tensr{A}^0}} \\
		\cdots & 0 & \cconnes{A}[-1] & {\tensr{A}^{-1}} & {\tensr{A}^{-1}} & {\tensr{A}^{-1}} &
		{\phantom{\tensr{A}^{-1}}} \\
		& \vdots & {\vdots } & \vdots & \vdots & \vdots
		\arrow[from=2-2, to=2-1]
		\arrow[from=3-2, to=3-1]
		\arrow[from=4-2, to=4-1]
		\arrow[from=5-2, to=4-2]
		\arrow[from=4-2, to=3-2]
		\arrow[from=3-2, to=2-2]
		\arrow[from=2-2, to=1-2]
		\arrow["b", from=4-3, to=3-3]
		\arrow["b", from=3-3, to=2-3]
		\arrow[from=2-3, to=1-3]
		\arrow[from=5-3, to=4-3]
		\arrow[from=4-3, to=4-2]
		\arrow[from=3-3, to=3-2]
		\arrow[from=2-3, to=2-2]
		\arrow["\pi"', from=3-4, to=3-3]
		\arrow["{\idd - \t}"', from=3-5, to=3-4]
		\arrow["\N"', from=3-6, to=3-5]
		\arrow["\pi"', from=4-4, to=4-3]
		\arrow["{\idd - \t}"', from=4-5, to=4-4]
		\arrow["\N"', from=4-6, to=4-5]
		\arrow["\pi"', from=2-4, to=2-3]
		\arrow["{\idd - \t}"', from=2-5, to=2-4]
		\arrow["\N"', from=2-6, to=2-5]
		\arrow["{\idd - \t}"', from=2-7, to=2-6]
		\arrow["{\idd - \t}"', from=3-7, to=3-6]
		\arrow["{\idd - \t}"', from=4-7, to=4-6]
		\arrow[from=2-4, to=1-4]
		\arrow[from=2-5, to=1-5]
		\arrow[from=2-6, to=1-6]
		\arrow["b", from=3-4, to=2-4]
		\arrow["{b'}", from=3-5, to=2-5]
		\arrow["b", from=3-6, to=2-6]
		\arrow["b", from=4-4, to=3-4]
		\arrow["{b'}", from=4-5, to=3-5]
		\arrow["b", from=4-6, to=3-6]
		\arrow[from=5-4, to=4-4]
		\arrow[from=5-5, to=4-5]
		\arrow[from=5-6, to=4-6]
		\arrow["{h' \N}", curve={height=-12pt}, from=3-3, to=3-4]
		\arrow["{h' \N}", curve={height=-12pt}, from=4-3, to=4-4]
		\arrow["{h' \N}", curve={height=-12pt}, from=2-3, to=2-4]
		\arrow["h", curve={height=-12pt}, from=2-4, to=2-5]
		\arrow["{h'}", curve={height=-12pt}, from=2-5, to=2-6]
		\arrow["h", curve={height=-12pt}, from=3-4, to=3-5]
		\arrow["{h'}", curve={height=-12pt}, from=3-5, to=3-6]
		\arrow["h", curve={height=-12pt}, from=4-4, to=4-5]
		\arrow["{h'}", curve={height=-12pt}, from=4-5, to=4-6]
		\arrow["h", curve={height=-12pt}, from=2-6, to=2-7]
		\arrow["h", curve={height=-12pt}, from=3-6, to=3-7]
		\arrow["h", curve={height=-12pt}, from=4-6, to=4-7]
	\end{tikzcd}
	\caption{The Augmented Cyclic Bicomplex $\tilde{\mathcal{CC}} \left( \mathcal{A} \right)$ with Row Contractions.}
	\label{fig:augmented-cyclic-bicomplex}
\end{figure}
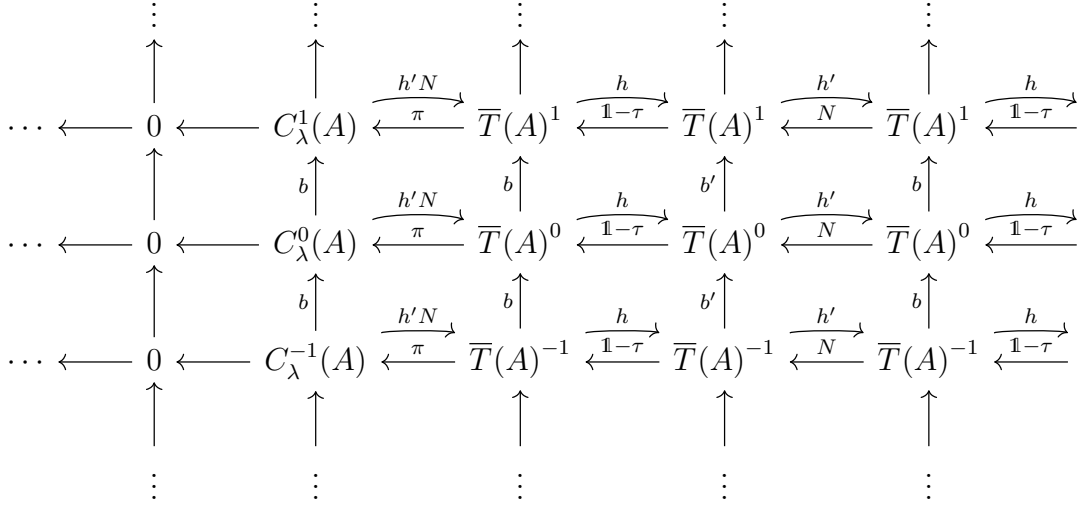

\begin{proof}[Proof of \cref{lm:cyc-bicomplex-homotopy-equivalent-connes}]
	Let us append Connes' complex to the standard cyclic bicomplex at the $(-1)$-th column and denote
	the resulting bicomplex by $\cycbiaug[A]$. We will call $\cycbiaug[A]$ the \textbf{augmented cyclic bicomplex}
	(see \cref{fig:augmented-cyclic-bicomplex}).
	Note that $\totc{\cycbiaug[A]}[][]$ is naturally identified with the mapping cone of the map
	$p[-1] \colon {\totc{\cycbi[A]}[][]}[-1] \rightarrow {\cconnes{\mathcal{A}}[]}[-1]$
	and so, by \cref{lm:contractible-cone-homotopy-equivalence}, it is enough to show that $\totc{\cycbiaug[A]}[][]$ is contractible.

	The rows of the bicomplex $\cycbiaug[A]$ are contractible. Explicitly, let us define a map
	$h_{\cycbiaug} \colon \cycbiaug \left( A \right)_{*}^{*} \rightharpoonup \cycbiaug \left( A \right)_{*+1}^{*}$
	by
	\begin{equation}
		\left( h_{\cycbiaug} \right)_i^j \defeq
		\begin{cases}
			h     & i \geq 0 \textrm { even}, \\
			h'    & i \geq 0 \textrm{ odd},   \\
			h' \N & i = -1,                   \\
			0     & i < -1.
		\end{cases}
	\end{equation}
	The induced map $h_{\cycbiaug}^{\totl} \colon \totc{\cycbiaug[A]} \rightharpoonup \totc{\cycbiaug[A]}[*-1]$
	on $\totc{\cycbiaug[A]}[][]$ is an $R$-linear map of degree $-1$ which satisfies
	\begin{equation}
		\delta_{\hor}^{\totl} h_{\cycbiaug}^{\totl} + h_{\cycbiaug}^{\totl} \delta_{\hor}^{\totl} = \idd,
	\end{equation}
	a consequence of the identities \eqref{eq:h'N-homotopy} and \eqref{eq:Nh'-homotopy}. Hence, to show that
	$\totc{\cycbiaug[A]}[][]$ is contractible, it is enough to verify the convergence condition of
	\cref{lm:total-bicomplex-contractible-R-linear-contraction-1}, namely
	that
	\begin{equation} \label{eq:conv-cond-aug-cyc-bicomplex}
		\left( \delta_{\ver}^{\totl} h_{\cycbiaug}^{\totl} + h_{\cycbiaug}^{\totl} \delta_{\ver}^{\totl} \right)^n
		\left( \s_i \left( c \right) \right) \xrightarrow[n \to \infty]{} 0
	\end{equation}
	for all $c \in \cycbiaug[A]_i^j$.

	Let us set $\hat{T} = \delta_{\ver}^{\totl} h_{\cycbiaug}^{\totl} + h_{\cycbiaug}^{\totl} \delta_{\ver}^{\totl}$.
	First, note that since $\nnorm[\hat{T}] \leq 1$, it is enough to show that
	\begin{equation} \label{eq:conv-cond-aug-cyc-bicomplex-even-powers}
		\hat{T}^{2n} \left( \s_{2k} \left( c \right) \right) \xrightarrow[n \to \infty]{} 0
	\end{equation}
	whenever $c \in \cycbiaug[A]_{2k}^j$, i.e., whenever $c$ belongs to an even column of $\cycbiaug[A]$.
	Given $c \in \cycbiaug[A]_{2k}^j = \tensr{A}^j$ with $k \geq 0$, we have
	\begin{equation*}
		\begin{aligned}
			\hat{T}^2 \left( \s_{2k} \left( c \right) \right)
			 & = \hat{T} \left( \s_{2k+1} \left( hb - b'h \right) \left( c \right) \right)
			= \s_{2k + 2} \left( \left( -h' b' + bh' \right) \left( hb - b' h \right) \left( c \right) \right)
			\\
			 & = \s_{2k+2} \left( \left( -h' b' h b + b h' h b - b h' b' h \right) \left( c \right) \right)
			\\
			 & = \s_{2k + 2} \left( S \left( c \right) \right),
		\end{aligned}
	\end{equation*}
	and so we see that $\hat{T}^2$ acts on the even columns via the periodicity operator $S$. Hence,
	\cref{lm:S-almost-nilpotent} implies that
	\begin{equation*}
		\nnorm[\hat{T}^{2n} \left( \s_{2k} \left( c \right) \right)]
		= \nnorm[\s_{2k + 2n} \left( S^n \left( c \right) \right)] =
		\nnorm[S^n \left( c \right)] \xrightarrow[n \to \infty]{} 0
	\end{equation*}
	showing \cref{eq:conv-cond-aug-cyc-bicomplex-even-powers}.
\end{proof}

\subsection{The Periodicity Operator} \label{subsec:periodicity-operator}
In this subsection we will show various properties of the periodicity operator
$S \colon \tensr{A} \rightharpoonup \tensr{A}$ defined in \cref{dfn:periodicity-operator}. To prove several
identities for $S$, it will be convenient to represent $S$ as a composition $S = T' T$
for some operators $T,T'$ and study their properties.

\begin{figure}[htb]
	\begin{tikzcd}
		{\left( \tensr{A}, b \right)} & {\left( \tensr{A}, b' \right)} &
		{\left( \tensr{A}, b \right)} & {\left( \tensr{A}, b' \right)} &
		{\cdots\phantom{\left( \tensr{A}, b \right)}}
		\arrow["\N", shift left=2, from=1-3, to=1-2]
		\arrow["H", from=1-2, to=1-2, harpoon, loop below]
		\arrow["{\idd - \t}", shift left=2, from=1-4, to=1-3]
		\arrow["h", shift left=2, from=1-1, to=1-2]
		\arrow["{\idd - \t}", shift left=2, from=1-2, to=1-1]
		\arrow["{h'}", shift left=2, from=1-2, to=1-3]
		\arrow["T"', harpoon, curve={height=18pt}, from=1-1, to=1-2]
		\arrow["{T'}"', harpoon, curve={height=18pt}, from=1-2, to=1-3]
		\arrow["\beta", harpoon, curve={height=-48pt}, from=1-1, to=1-3]
		\arrow["S", harpoon, curve={height=-24pt}, from=1-1, to=1-3]
		\arrow["{S'}", harpoon, curve={height=-24pt}, from=1-2, to=1-4]
		\arrow["h", shift left=2, from=1-3, to=1-4]
		\arrow[shift left=2, from=1-5, to=1-4]
		\arrow[shift left=2, from=1-4, to=1-5]
		\arrow[shift right=2, curve={height=-24pt}, from=1-3, to=1-5]
		\arrow["T"', harpoon, curve={height=18pt}, from=1-3, to=1-4]
		\arrow["H", harpoon, from=1-4, to=1-4, loop below]
	\end{tikzcd}
	\caption{The natural domains and codomains of the operators \newline
		$\idd - \t, \N, h, h', T, T', S, S', H, \beta$.
	}
	\label{fig:T-T'-S-S'}
\end{figure}
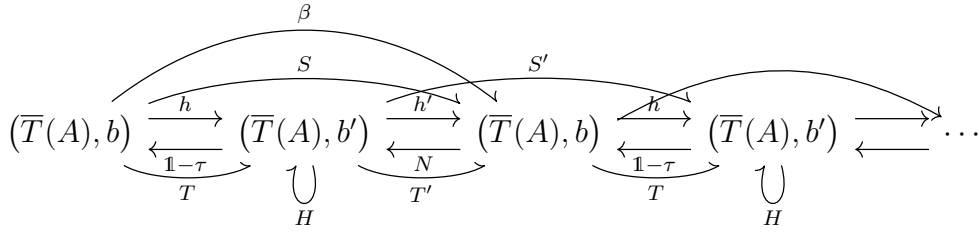

In what follows, we will consider the operators $\idd - \t, h, \N, h'$ as degree zero maps between
chain complexes with domains and codomains as indicated in \cref{fig:T-T'-S-S'}. Working with
the convention above, we have
\begin{align}
	\partial \left( \idd - \t \right) & = b \left( \idd - \t \right) - \left( \idd - \t \right) b'
	\stackrel{\eqref{eq:b-1-t-rel}}{=} 0, \label{eq:partial-1-t}
	\\
	\partial \left( \N \right)        & = b'\N - \N b \stackrel{\eqref{eq:b'-N-rel}}{=} 0. \label{eq:partial-N}
\end{align}

Define operators $T,T' \colon \tensr{A}^{*} \rightharpoonup \tensr{A}^{*+1}$ by
\begin{align}
	T  & \defeq hb - b'h = - \partial \left( h \right), \label{eq:T-def}   \\
	T' & \defeq bh' - h' b' = \partial \left( h' \right) \label{eq:T'-def}
\end{align}
where $h,h'$ are given by \cref{eq:def-h-h'}.

\begin{lm}
	The operators $T,T'$ are $R$-linear of degree $1$ and satisfy the identities
	\begin{align}
		Tb                          & = -b'T, \label{eq:T-b}    \\
		T'b'                        & = -bT', \label{eq:T'-b'}  \\
		T \left( \idd - \t \right)  & = \N T', \label{eq:T-1-t} \\
		\left ( \idd - \t \right) T & = T' \N. \label{eq:1-t-T}
	\end{align}
\end{lm}
\begin{proof}
	Since $T = \partial \left( -h \right)$ and $T' = \partial \left( h' \right)$ are boundaries, they are
	degree one chain maps which shows \cref{eq:T-b,eq:T'-b'}.
	To show \cref{eq:T-1-t}, we compute
	\begin{equation*}
		\begin{aligned}
			T \left( \idd - \t \right)
			\stackrel{\eqref{eq:T-def}}{=}  &
			-\partial \left( h \right) \left( \idd - \t \right)
			\stackrel{\eqref{eq:partial-1-t}}{=}
			-\partial \left( h \left( \idd - \t \right) \right)
			\stackrel{\eqref{eq:Nh'-homotopy}}{=}
			-\partial \left( \idd - \N h' \right)
			\stackrel{\eqref{eq:partial-N}}{=}
			\N \partial \left( h' \right)
			\\
			\stackrel{\eqref{eq:T'-def}}{=} &
			\N T'.
		\end{aligned}
	\end{equation*}
	Similarly, to show \cref{eq:1-t-T}, we compute
	\begin{equation*}
		\begin{aligned}
			\left( \idd - \t \right) T
			\stackrel{\eqref{eq:T-def}}{=}  &
			-\left( \idd - \t \right) \partial \left( h \right)
			\stackrel{\eqref{eq:partial-1-t}}{=}
			- \partial \left( \left( \idd - \t \right) h \right)
			\stackrel{\eqref{eq:h'N-homotopy}}{=}
			- \partial \left( \idd - h' \N \right)
			\stackrel{\eqref{eq:partial-N}}{=}
			\partial \left( h' \right) \N
			\\
			\stackrel{\eqref{eq:T'-def}}{=} &
			T' \N.
		\end{aligned}
	\end{equation*}
\end{proof}

Next, we define an operator $S' \colon \tensr{A} \rightharpoonup \tensr{A}$ analogous to the periodicity
operator $S$ of \cref{dfn:periodicity-operator} by
\begin{equation}
	S' \defeq - \left( h b h' b' + b' h b h' \right) + b' h h' b'. \label{eq:S'-def}
\end{equation}

\begin{lm}
	The operators $S,S' \colon \tensr{A} \rightharpoonup \tensr{A}$ are $R$-linear of degree $2$
	and satisfy the identities
	\begin{align}
		S                          & = T' T, \label{eq:S-T'-T}                      \\
		S'                         & = T T', \label{eq:S'-T-T'}                     \\
		bS                         & = Sb, \label{eq:bS=Sb}                         \\
		b'S'                       & = S'b', \label{eq:b'S'=S'b'}                   \\
		S \left( \idd - \t \right) & = \left( \idd - \t \right) S', \label{eq:S1-t} \\
		S' \N                      & = \N S. \label{eq:S'-N}
	\end{align}
\end{lm}
\begin{proof}
	We have
	\begin{equation*}
		T' T = \left( bh' - h' b' \right) \left( hb - b'h \right) =
		b h' h b - b h' b' h - h' b' h b = S
	\end{equation*}
	and
	\begin{equation*}
		T T' = \left( hb - b'h \right) \left( bh' - h' b' \right) =
		- h b h' b' - b' h b h' + b' h h' b' = S'
	\end{equation*}
	which shows \cref{eq:S-T'-T,eq:S'-T-T'}.
	In particular, $S$ and $S'$ are $R$-linear of degree two as the composition of
	two $R$-linear operators of degree one.\footnote{Note that from the expression
		$S = - \left( h' b' h b + b h' b' h \right) + b h' h b$ it is not immediate that
		$S$ is $R$-linear as $b,b'$ are derivations over $d$; the same applies to $S'$.}
	Identities \eqref{eq:T-b} and \eqref{eq:T'-b'} show that $T,T'$ are chain maps and hence
	the compositions $S = T' T$ and $S' = T T'$ are chain maps which shows \cref{eq:bS=Sb,eq:b'S'=S'b'}.
	Finally,
	\begin{equation*}
		S \left( \idd - \t \right)
		\stackrel{\eqref{eq:S-T'-T}}{=}
		T' T \left( \idd - \t \right)
		\stackrel{\eqref{eq:T-1-t}}{=}
		T' \N T'
		\stackrel{\eqref{eq:1-t-T}}{=}
		\left( \idd - \t \right) T T'
		\stackrel{\eqref{eq:S'-T-T'}}{=}
		\left( \idd - \t \right) S'
	\end{equation*}
	which shows \cref{eq:S1-t}. Identity \eqref{eq:S'-N} is proven analogously.
\end{proof}

\begin{rem}
	The operators $T$ and $T'$ and their relation to $S$ come up naturally in the proof of
	\cref{lm:cyc-bicomplex-homotopy-equivalent-connes}.
	Using the notation of the proof of \cref{lm:cyc-bicomplex-homotopy-equivalent-connes},
	the operator
	$\hat{T} \colon \totc{\cycbiaug[A]} \rightarrow \totc{\cycbiaug[A]}$ satisfies
	\begin{equation*}
		\hat{T} \left( s_{2k} \left( c \right) \right) = s_{2k+1} \left( T \left( c \right) \right),
		\qquad
		\hat{T} \left( s_{2k + 1} \left( c \right) \right) = s_{2k+2} \left( T' \left( c \right) \right)
		\qquad
		\left( k \geq 0 \right).
	\end{equation*}
	During the proof of \cref{lm:cyc-bicomplex-homotopy-equivalent-connes}, we have seen that $\hat{T}^2$ acts on
	even columns via $S$ which explains the identity $S = T' T$.
\end{rem}

In particular, we see that $S \colon \left( \tensr{A}, b \right) \rightharpoonup \left( \tensr{A}, b \right)$
is a chain map of degree two and by \cref{eq:S1-t}, $S$ descends to the quotient
$\tensr{A} / \Im \left( \idd - \t \right) = \cconnes{A}[]$ and induces a map
$S \colon \hcyc{\mathcal{A}} \rightharpoonup \hcyc{\mathcal{A}}[*+2]$ on cyclic homology.

\begin{ex}
	The curvature $\mu_0(1)$ defines a cyclic homology class
	$\eqcl{\mu_0(1)} \in \hcyc{\mathcal{A}}[1]$ since we have
	\begin{equation*}
		b \left( \mu_0 \left( 1 \right) \right) = - \mu_0 \left( 1 \right) \otimes \mu_0 \left( 1 \right) = 0
		\mod \left( \idd - \t \right).
	\end{equation*}
	A direct calculation shows that
	\begin{align*}
		S \left( \mu_0 \left( 1 \right) \right) ={}   &
		\frac{1}{2} \mu_2 \left( \mu_0 \left( 1 \right), \mu_0 \left( 1 \right) \right) +
		\frac{1}{6} \mu_0(1)^{\otimes 3} \in \cconnes{A}[3],
		\\
		S^2 \left( \mu_0 \left( 1 \right) \right) ={} &
		\frac{1}{12} \left(
		\mu_2 \left(
		\mu_0 \left( 1 \right), \mu_2 \left( \mu_0 \left( 1 \right), \mu_0 \left( 1 \right) \right)
		\right) +
		\mu_2 \left(
		\mu_2 \left( \mu_0 \left( 1 \right), \mu_0 \left( 1 \right) \right), \mu_0 \left( 1 \right)
		\right) +
		\mu_4 \left( \mu_0 \left( 1 \right)^{\otimes 4} \right)
		\right)
		\\
		                                              & + \frac{1}{12}
		\left(
		\mu_0 \left( 1 \right)^{\otimes 2} \otimes
		\mu_2 \left( \mu_0 \left( 1 \right), \mu_0 \left( 1 \right) \right)
		\right) +
		\frac{1}{60} \mu_0 \left( 1 \right)^{\otimes 5}
		\in \cconnes{A}[5].
	\end{align*}
\end{ex}

The periodicity operator also induces a map $S \colon \hhoch{\mathcal{A}} \rightharpoonup \hhoch{\mathcal{A}}[*+2]$
on Hochschild homology, but this map is trivial since $S$ is null-homotopic
on $\choch{\mathcal{A}}[]$. Explicitly, let us define
the operator $\beta \colon \tensr{A}^{*} \rightharpoonup \tensr{A}^{*+1}$ by
\begin{equation}
	\beta \defeq T'h = \partial \left( h' \right) h = \left( b h' - h' b' \right) h. \label{eq:def-beta}
\end{equation}
Then
\begin{equation}
	\partial \beta
	=
	\partial \left( \partial \left( h' \right) h \right)
	=
	- \partial \left( h' \right) \partial \left( h \right) = T' T
	\stackrel{\eqref{eq:S-T'-T}}{=} S. \label{eq:partial-beta-S}
\end{equation}

Finally, we will prove that the periodicity operator $S$ is ``almost nilpotent'' in the sense that
\begin{lm} \label{lm:S-almost-nilpotent}
	The operator $S \colon \tensr{A} \rightharpoonup \tensr{A}$ satisfies
	$S^n \left( x \right) \to 0$ for all $x \in \tensr{A}$.
\end{lm}

To prove \cref{lm:S-almost-nilpotent}, we will analyze the action of $S$ on elements based on weight.
Let us denote by $\mathcal{F}_k \defeq \oplus_{i=1}^k A^{\otimes i}$ the weight filtration on $\tensr{A}$,
where we set $\mathcal{F}_k \defeq 0$ if $k \leq 0$.
Denote by ${\mathcal{F}}^d_{k} \defeq \left( \oplus_{i=1}^k A^{\otimes i} \right)^{d}$ the collection of
elements of (cohomological) degree $d$ and weight less than or equal to $k$.
By the definition of $S$, it is clear that we have
$S \left( \mathcal{F}_k \right) \subseteq \mathcal{F}_{k+2}$. The following lemma shows
that if we ignore all terms in $S$ which involve $\mu_0$, the resulting operator actually decreases
the weight of elements by two.

\begin{lm} \label{lm:S-split-based-on-mu_0}
	The operator $S \colon \tensr{A}^{*} \rightharpoonup \tensr{A}^{*+2}$ is a sum $S = S_{-} + S_{\mu_0}$ of
	two operators $S_{-}, S_{\mu_0} \colon \tensr{A}^{*} \rightharpoonup \tensr{A}^{*+2}$ such that:
	\begin{enumerate}
		\item $S_{-} \left( \mathcal{F}_k \right) \subseteq \mathcal{F}_{k-2}$ and $\nnorm[S_{-}] \leq 1$.
		\item $S_{\mu_0} \left( \mathcal{F}_k \right) \subseteq \mathcal{F}_{k+2}$ and
		      $\nnorm[S_{\mu_0}] \leq \nnorm[\mu_0(1)]$.
	\end{enumerate}
\end{lm}
\begin{proof}
	Given $k \geq 0$, let us denote by $b_{> k}$ the part of the operator $b$ which applies only
	the operations $\mu_r$ for $r > k$ and by $b_k$ the part of the operator $b$ which applies only the
	operation $\mu_k$. So for example we have $b = b_0 + b_1 + b_{> 1}$ where
	\begin{equation*}
		b_0 \left( \mathcal{F}_k \right) \subseteq \mathcal{F}_{k+1}, \qquad
		b_1 \left( \mathcal{F}_k \right) \subseteq \mathcal{F}_k, \qquad
		b_{> 1} \left( \mathcal{F}_k \right) \subseteq \mathcal{F}_{k-1}.
	\end{equation*}
	We will use the same notation to describe the splitting of the operator $b'$.\footnote{In the proof of the lemma,
		the notation $b'_k \colon \tensr{A} \rightharpoonup \tensr{A}$ will be used to denote the part of the operator $b'$
		in which we apply only the operation $\mu_k$ and \textit{not} the corestriction
		$b'_k = \mu_k \colon A^{\otimes k} \rightharpoonup A$.}
	Note that by definition we have $b'_1 = b_1$. Note also that the operator $b_1$ preserves the weight of
	elements and commutes with $\t$ and hence commutes with both $h$ and $h'$.

	We have $S = T' T$. Let us split the operators $T,T'$ as
	\begin{align*}
		T  & = hb - b'h = h \left( b_0 + b_1 + b_{> 1} \right) - \left( b'_0 + b'_1 + b'_{>1} \right) h =
		\underbrace{\left( h b_0 - b'_0 h \right)}_{T_0} +
		\underbrace{\left( h b_{>1} - b'_{>1} h \right)}_{T_{>1}},
		\\
		T' & = bh' - h'b' = \left( b_0 + b_1 + b_{> 1} \right)h' - h' \left( b'_0 + b'_1 + b'_{>1} \right) =
		\underbrace{\left( b_0 h' - h' b'_0 \right)}_{T'_0} +
		\underbrace{\left( b_{>1} h' - h' b'_{>1} \right)}_{T'_{>1}},
	\end{align*}
	where the operators $T_0,T_{>1},T'_0,T'_{>1}$ satisfy
	\begin{enumerate}
		\item $T_0 \left( \mathcal{F}_k \right) \subseteq \mathcal{F}_{k+1}$ and $\nnorm[T_0] \leq \nnorm[\mu_0(1)]$.
		\item $T'_0 \left( \mathcal{F}_k \right) \subseteq \mathcal{F}_{k+1}$ and $\nnorm[T'_0] \leq \nnorm[\mu_0(1)]$.
		\item $T_{>1} \left( \mathcal{F}_k \right) \subseteq \mathcal{F}_{k-1}$ and $\nnorm[T_{>1}] \leq 1$.
		\item $T'_{>1} \left( \mathcal{F}_k \right) \subseteq \mathcal{F}_{k-1}$ and $\nnorm[T'_{>1}] \leq 1$.
	\end{enumerate}
	Hence, the decomposition
	\begin{equation*}
		S = T' T = \left( T'_0 + T'_{>1} \right) \left( T_0 + T_{>1} \right) =
		\underbrace{T'_0 T_0 + T'_0 T_{>1} + T'_{>1} T_0}_{S_{\mu_0}} + \underbrace{T'_{>1} T_{>1}}_{S_{-}}
	\end{equation*}
	satisfies the conditions of \cref{lm:S-split-based-on-mu_0}.
\end{proof}

\begin{lm} \label{lm:S-norm-estimate}
	Given $x \in \mathcal{F}_k$ and $N \geq 0$, we have
	$\nnorm[S^{\left \lceil \frac{k}{2} \right \rceil + 2N} \left( x \right)] \leq
		\nnorm[\mu_0(1)]^{N+1} \cdot \nnorm[x]$.
\end{lm}
\begin{proof}
	Let $l \geq 0$. Using \cref{lm:S-split-based-on-mu_0} we can write
	\begin{equation*}
		S^l \left( x \right) = \left( S_{-} + S_{\mu_0} \right)^l \left( x \right) = \sum_{r = 0}^l y_r
	\end{equation*}
	where $y_r \in \mathcal{F}_{k - 2l + 4r}$ and $\nnorm[y_r] \leq \nnorm[\mu_0(1)]^r \cdot \nnorm[x]$ for all
	$0 \leq r \leq l$. The term $y_r$ in the decomposition above is the sum of all terms in which $S_{\mu_0}$
	is applied exactly $r$ times (and hence $S_{-}$ is applied $l - r$ times). In particular,
	if $l = \left \lceil \frac{k}{2} \right \rceil + 2N$ we see that the terms $y_0, \dots, y_N$ all belong to
	$\mathcal{F}_0$ and hence vanish.
\end{proof}

\begin{proof}[Proof of \cref{lm:S-almost-nilpotent}]
	Let $d \in \ZZ$ and $x \in \tensr{A}^d$. Since $\tensr{A}^d = \oplus_{k \geq 1} \left( A^{\otimes k} \right)^d$,
	by \cref{lm:direct-sum-iteration-zero-limit}, we can assume that
	$x \in \left( A^{\otimes k} \right)^d \subseteq \mathcal{F}_k^d$. \Cref{lm:S-norm-estimate}
	together with the fact that $\nnorm[\mu_0(1)] < 1$ and $\nnorm[S] \leq 1$ then imply
	that $S^n \left( x \right) \to 0$.
\end{proof}

\begin{rem}
	Note that by the definition of the complete direct sum totalization, the periodicity operator
	$\sigma \colon \totc{\cycbi[A]} \rightharpoonup \totc{\cycbi[A]}[*+2]$ given by
	\cref{def:sigma} satisfies $\sigma^n \left( x \right) \to 0$ for all $x \in \totc{\cycbi[A]}[]$.
	\Cref{lm:S-almost-nilpotent} shows that the periodicity operator $S$ on
	Connes' complex also satisfies this property.
\end{rem}

\subsection{Connes' Exact Sequence} \label{sec:connes-exact-sequence}
Cyclic homology and Hochschild homology are related via a long exact sequence called
\textbf{Connes' exact sequence} which involves the periodicity operator.
The explicit form of the sequence depends on the models we use for cyclic and Hochschild homology.
We will describe first what the sequence looks like when working with the cyclic bicomplex and then
give a description using Connes' complex and show the relation between both descriptions.

Since the cyclic bicomplex $\cycbi[A]$ is $2$-periodic, we have a natural short exact sequence
of bicomplexes
\begin{equation*}
	0 \rightarrow \cycbi[A]^{\{2\}} \rightarrow \cycbi[A] \rightarrow \cycbi[A][2,0] \rightarrow 0.
\end{equation*}
Here, $\cycbi[A]^{\{2\}}$ denotes the bicomplex which consists of the first two columns of $\cycbi[A]$
and $\cycbi[A][2,0]$ denotes the bicomplex in which we shift the columns of $\cycbi[A]$ to the right by two (so
that $\cycbi[A][2,0]_i^j = \cycbi[A]_{i-2}^j$ for all $i,j$). Totalizing the bicomplexes,
we obtain the short exact sequence
\begin{equation}
	0 \rightarrow \totc{\cycbi[A]^{\{2\}}}[] \xrightarrow{i} \totc{\cycbi[A]}[] \xrightarrow{\sigma}
	{\totc{\cycbi[A]}[]}[2] \rightarrow 0. \label{eq:tot-cyclic-bicomplex-exact-sequence}
\end{equation}
The map $i$ is induced by the inclusion of the first two columns into $\cycbi[A]$ while the map
$\sigma$ is induced from the natural ``periodicity operator''
$\sigma \colon \totc{\cycbi[A]} \rightharpoonup \totc{\cycbi[A]}[*+2]$
given by
\begin{equation}
	\sigma \left( \sum_{k \geq 0} \s_k \left( c_k \right) \right)
	\defeq \sum_{k \geq 0} \s_k \left( c_{k+2} \right). \label{def:sigma}
\end{equation}
The long exact sequence on cohomology associated to \eqref{eq:tot-cyclic-bicomplex-exact-sequence} has the form
\begin{equation} \label{eq:SBI-long-exact-sequence-cyclic-bicomplex}
	\adjustbox{scale=0.93}{
		\begin{tikzcd}
			{\cohom{{\totc{\cycbi[A]^{\{2\}}}[]}}} \rar{i} & {\cohom{{\totc{\cycbi[A]}[]}}} \rar[harpoon]{\sigma} &
			{\cohom{{\totc{\cycbi[A]}[]}}[*+2]}
			\ar[out=-30, in=150, start anchor=real east, end anchor=real west, overlay, swap, harpoon]{dll}{\delta} \\
			{\cohom{{\totc{\cycbi[A]^{\{2\}}}[]}}[*+1]} \rar{i} & \cdots \phantom{\totc{\cycbi[A]}[]} \, & {}
		\end{tikzcd}
	}
\end{equation}

We have already seen in \cref{lm:cyc-bicomplex-homotopy-equivalent-connes} that $\totc{\cycbi[A]}[]$ is
homotopy equivalent to $\cconnes{\mathcal{A}}[]$ via the map $p \colon \totc{\cycbi[A]}[] \rightarrow \cconnes{\mathcal{A}}[]$.
As for the complex
\begin{equation*}
	\totc{\cycbi[A]^{\{2\}}}[] = \Cone{\idd - \t},
\end{equation*}
we have a natural inclusion $j \colon \choch{\mathcal{A}}[] \rightarrow \totc{\cycbi[A]^{\{2\}}}[]$
obtained by identifying $\choch{\mathcal{A}}[]$ with the first column of $\cycbi[A]^{\{2\}}$.
When $\mathcal{A}$ is unital, the second column $\cycbi[A]^{\{2\}}$ is the bar complex which is contractible by
\cref{lm:bar-homology-curved-algebra} and hence by \cref{cor:inclusion-into-mapping-cone-equivalence} the map $j$ is a homotopy equivalence.\footnote{When $\mathcal{A}$ is not unital, the definition of Hochschild chain complex is taken to
	be $\totc{\cycbi[A]^{\{2\}}}[]$ so that everything that follows continues to hold.}
Hence, when $\mathcal{A}$ is unital, both $j$ and $p$ induce isomorphisms on cohomology, and we can
use $j$ and $p$ to obtain an equivalent long exact sequence
\begin{equation} \label{eq:SBI-long-exact-sequence-connes-complex}
	\cdots \rightarrow {\hhoch{\mathcal{A}}} \xrightarrow{I} {\hcyc{\mathcal{A}}} \xrightharpoonup{S}
	{\hcyc{\mathcal{A}}[*+2]} \xrightharpoonup{B} {\hhoch{\mathcal{A}}[*+1]} \rightarrow \cdots
\end{equation}
connecting the cohomology of Connes' complex and the Hochschild complex. The relation between
the sequences \eqref{eq:SBI-long-exact-sequence-cyclic-bicomplex} and \eqref{eq:SBI-long-exact-sequence-connes-complex} is given by the diagram in \cref{fig:SBI-sequence-bicomplex-connes}.
\begin{figure}[htb]
	\centering
	\adjustbox{scale=0.80,center}{
		\begin{tikzcd}[column sep=small]
			\cdots & {\cohom{{\totc{\cycbi[A]^{\{2\}}}[]}}} & {\cohom{{\totc{\cycbi[A]}[]}}} &
			{\cohom{{\totc{\cycbi[A]}[]}}[*+2]} & {\cohom{{\totc{\cycbi[A]^{\{2\}}}[]}}[*+1]} & \cdots \\
			\cdots & {\hhoch{\mathcal{A}}} & {\hcyc{\mathcal{A}}} & {\hcyc{\mathcal{A}}[*+2]} &
			{\hhoch{\mathcal{A}}[*+1]} & \cdots
			\arrow["i", from=1-2, to=1-3]
			\arrow["\sigma", harpoon, from=1-3, to=1-4]
			\arrow["\delta", harpoon, from=1-4, to=1-5]
			\arrow["I", from=2-2, to=2-3]
			\arrow["S", harpoon, from=2-3, to=2-4]
			\arrow["B", harpoon, from=2-4, to=2-5]
			\arrow["j", from=2-2, to=1-2]
			\arrow["p", from=1-3, to=2-3]
			\arrow["p", from=1-4, to=2-4]
			\arrow["j", from=2-5, to=1-5]
			\arrow[from=1-5, to=1-6]
			\arrow[from=2-5, to=2-6]
			\arrow[from=1-1, to=1-2]
			\arrow[from=2-1, to=2-2]
		\end{tikzcd}
	}
	\caption{Connes' exact sequence via the cyclic bicomplex and Connes' complex.}
	\label{fig:SBI-sequence-bicomplex-connes}
\end{figure}

We will describe the operators appearing in \eqref{eq:SBI-long-exact-sequence-connes-complex} explicitly
and show that the diagram in \cref{fig:SBI-sequence-bicomplex-connes} indeed commutes.
The operators $I,S,B$ are chain maps defined on the level of chain complexes and induce operators
(denoted by the same name) on cohomology. The operator
$I \colon \hhoch{\mathcal{A}}[] \rightarrow \hcyc{\mathcal{A}}[]$ is induced
by the natural projection $\choch{\mathcal{A}}[] \rightarrow \cconnes{\mathcal{A}}[]$.
The operator $S$ is the periodicity operator introduced in \cref{dfn:periodicity-operator}. The operator
$B \colon \tensr{A} \rightharpoonup \tensr{A}$ is a degree $-1$ chain map called \textbf{Connes boundary map}
defined by
\begin{equation}
	B \defeq \left( \idd - \t \right) H \N \label{eq:B-def}
\end{equation}
where $H$ is the contraction of the bar complex given by \cref{eq:contraction-bar-with-mu-0}.

\begin{lm} \label{lm:B-properties}
	The operator $B \colon \tensr{A} \rightharpoonup \tensr{A}$ is $R$-linear of degree $-1$. We
	have the identities:
	\begin{align}
		\left[ b, B \right]        & = bB + Bb = 0, \label{eq:bB-Bb} \\
		B \left( \idd - \t \right) & = \N B = 0, \label{eq:B-1-t}    \\
		B^2 = 0.	 \label{eq:B^2=0}
	\end{align}
\end{lm}
\begin{proof}
	To see \cref{eq:bB-Bb}, note that
	\begin{equation*}
		\begin{aligned}
			bB & \stackrel{\eqref{eq:B-def}}{=} b \left( \idd - \t \right) H \N
			\stackrel{\eqref{eq:b-1-t-rel}}{=} \left( \idd - \t \right) b' H \N
			\stackrel{\eqref{eq:b'H-homotopy}}{=} \left( \idd - \t \right) \left( \idd - H b' \right) \N
			\\
			   & \stackrel{\eqref{eq:b'-N-rel}}{=} \left( \idd - \t \right) \N - \left( \idd - \t \right) H \N b
			\stackrel{\eqref{eq:1-t-N=0}}{=} -Bb.
		\end{aligned}
	\end{equation*}
	To see \cref{eq:B-1-t}, note that
	\begin{align*}
		B \left( \idd - \t \right) & \stackrel{\eqref{eq:B-def}}{=} \left( \idd - \t \right) H \N \left( \idd - \t \right)
		\stackrel{\eqref{eq:1-t-N=0}}{=} 0,
		\\
		\N B                       & \stackrel{\eqref{eq:B-def}}{=} \N \left( \idd - \t \right) H \N \stackrel{\eqref{eq:1-t-N=0}}{=} 0.
	\end{align*}
	Finally, we have
	\begin{equation*}
		B^2 \stackrel{\eqref{eq:B-def}}{=} B \left( \idd - \t \right) H \N \stackrel{\eqref{eq:B-1-t}}{=} 0
	\end{equation*}
	which shows \cref{eq:B^2=0}.
\end{proof}

\Cref{lm:B-properties} implies that the map $B$ induces maps $\hhoch{\mathcal{A}} \rightharpoonup \hhoch{\mathcal{A}}[*-1]$
and $\hcyc{\mathcal{A}} \rightharpoonup \hhoch{\mathcal{A}}[*-1]$ on cohomology, denoted by the same name when no confusion
can arise. The relation between $B$ and the connecting morphism $\delta$ is given by the following lemma:

\begin{lm} \label{lm:B-connecting-morphism}
	The following diagram commutes:
	\begin{equation*}
		\begin{tikzcd}
			{\cohom{{\totc{\cycbi[A]}[]}}} & {\cohom{{\totc{\cycbi[A]^{\{2\}}}[]}}[*-1]} \\
			{\hcyc{\mathcal{A}}} & {\hhoch{\mathcal{A}}[*-1].}
			\arrow["\delta", harpoon, from=1-1, to=1-2]
			\arrow["p"', from=1-1, to=2-1]
			\arrow["j", from=2-2, to=1-2]
			\arrow["B"', harpoon, from=2-1, to=2-2]
		\end{tikzcd}
	\end{equation*}
\end{lm}
\begin{proof}
	Let $x = \sum_{k \geq 0} \s_k \left( c_k \right) \in \totc{\cycbi[A]}[]$ be a closed element. In particular,
	we have
	\begin{equation} \label{eq:x-closed-zero-level}
		b c_0 + \left( \idd - \t \right) c_1 = 0.
	\end{equation}
	Diagram chasing shows that $\delta \left( x \right)$ is represented by
	$\s_1 \left( \N c_0 \right)$ while $\left( jBp \right) \left( x \right) = B c_0$. Let us denote
	by $D$ the differential on $\totc{\cycbi[A]^{\{2\}}}[]$ so that
	\begin{equation*}
		D \left( c_0 + \s_1 c_1 \right) = \left( bc_0 + \left( \idd - \t \right) c_1 \right) - \s_1 b' c_1.
	\end{equation*}
	Then
	\begin{equation*}
		\begin{aligned}
			D \left( \s_1 H \N c_0 \right) & \stackrel{\eqref{eq:B-def}}{=}  Bc_0 - \s_1 \left( b' H \N c_0 \right)
			\stackrel{\eqref{eq:b'H-homotopy}}{=} Bc_0 - \s_1 \left( \N c_0 - H b' \N c_0 \right)
			\\
			                               & \stackrel{\eqref{eq:b'-N-rel}}{=} Bc_0 - \s_1 \left( \N c_0 \right) + \s_1 \left( H \N b c_0 \right)
			\stackrel{\eqref{eq:x-closed-zero-level}}{=} Bc_0 - \s_1 \left( \N c_0 \right)
		\end{aligned}
	\end{equation*}
	and hence $jBp $ coincides on cohomology with $\delta$.
\end{proof}
The relation between the map $S$ and the map $\sigma$ is given by the following lemma:
\begin{lm}
	Let $p \colon \totc{\cycbi[A]}[][] \rightarrow \cconnes{\mathcal{A}}[]$
	be the chain map from \cref{lm:cyc-bicomplex-homotopy-equivalent-connes}. Then the
	diagram
	\begin{equation*}
		\begin{tikzcd}
			{\totc{\cycbi[A]}} & {\totc{\cycbi[A]}[*+2]} \\
			\cconnes{\mathcal{A}} & \cconnes{\mathcal{A}}[*+2]
			\arrow["\sigma", harpoon, dashed, from=1-1, to=1-2]
			\arrow["p"', dashed, from=1-1, to=2-1]
			\arrow["p", dashed, from=1-2, to=2-2]
			\arrow["S"', harpoon, dashed, from=2-1, to=2-2]
		\end{tikzcd}
	\end{equation*}
	commutes up to homotopy. In particular, the induced diagram on cohomology
	\begin{equation*}
		\begin{tikzcd}
			{\cohom{{\totc{\cycbi[A]}[][]}}} & {\cohom{{\totc{\cycbi[A]}[][]}}[*+2]} \\
			\hcyc{\mathcal{A}} & \hcyc{\mathcal{A}}[*+2]
			\arrow["\sigma", harpoon, from=1-1, to=1-2]
			\arrow["p"', from=1-1, to=2-1]
			\arrow["p", from=1-2, to=2-2]
			\arrow["S"', harpoon, from=2-1, to=2-2]
		\end{tikzcd}
	\end{equation*}
	strictly commutes.
\end{lm}
\begin{proof}
	Define a degree one map $F \colon \totc{\cycbi[A]} \rightharpoonup \cconnes{\mathcal{A}}[*+1]$ by the formula
	\begin{equation*}
		F \left( \sum_{k \geq 0} \s_k \left( c_k \right) \right) \defeq
		\cyccl{h' \left( c_1 \right) - \beta \left( c_0 \right)}
	\end{equation*}
	where we denote by $\cyccl{x}$ the equivalence class of $x \in \tensr{A}$ in
	$\tensr{A} / \Im \left( \idd - \t \right)$. Then
	\begin{equation*}
		\begin{aligned}
			\left( \partial F \right) \left( \sum_{k \geq 0} \s_k \left( c_k \right) \right)
			 & \stackrel{\phantom{\eqref{eq:h'N-homotopy}}}{=}
			\left( bF + FD \right) \left( \sum_{k \geq 0} \s_k \left( c_k \right) \right)
			\\
			 & \stackrel{\phantom{\eqref{eq:h'N-homotopy}}}{=}
			\cyccl{b h' \left( c_1 \right) - b \beta \left( c_0 \right) +
				h' \left( -b' \left( c_1 \right) + \N \left( c_2 \right) \right) -
				\beta \left( b c_0 + \left( \idd - \t \right) c_1 \right)}
			\\
			 & \stackrel{\eqref{eq:h'N-homotopy}}{=}
			\cyccl{c_2 - \left( \partial \beta \right) \left( c_0 \right) +
				\left( b h' - h' b' - \beta \left( \idd - \t \right)\right) \left( c_1 \right)}
			\\
			 & \stackrel{\eqref{eq:partial-beta-S}, \eqref{eq:T'-def}, \eqref{eq:def-beta}}{=}
			\cyccl{c_2 - S \left( c_0 \right) +
				T' \left( \idd -  h \left( \idd - \t \right) \right) \left( c_1 \right)}
			\\
			 & \stackrel{\eqref{eq:Nh'-homotopy}}{=}
			\cyccl{c_2 - S \left( c_0 \right) + \left( T' \N h' \right) \left( c_1 \right)}
			\\
			 & \stackrel{\eqref{eq:1-t-T}}{=}
			\cyccl{c_2 - S \left( c_0 \right) + \left( \idd - \t \right) T h' \left( c_1 \right)}
			\\
			 & \stackrel{\phantom{\eqref{eq:h'N-homotopy}}}{=}
			\cyccl{c_2 - S \left( c_0 \right)} =
			\left( p \sigma - S p \right) \left( \sum_{k \geq 0} \s_k \left( c_k \right) \right).
		\end{aligned}
	\end{equation*}
\end{proof}

\begin{rem}
	We have chosen to define the periodicity operator on Connes' complex by the formula
	\begin{equation*}
		S = - \left( h' b' h b + b h' b' h \right) + b h' h b.
	\end{equation*}
	This has several advantages:
	\begin{enumerate}
		\item The operator $S$ is defined on the chain level and becomes a chain map on $\left( \tensr{A}, b \right)$.
		\item As a consequence of the identity $S \left( \idd - \t \right) = \left( \idd - \t \right) S'$, the operator
		      $S$ descends to a well-defined chain level operator on Connes' complex
		      $\left( \tensr{A} / \Im \left( \idd - \t \right), b \right)$.
		\item The operator $S$ is $R$-linear of degree $2$, even when $b$ and $b'$ are derivations over $d$.
	\end{enumerate}
	In \cite[Section 2.2]{Loday1998} the periodicity operator on Connes' complex for
	an associative algebra has two descriptions:
	\begin{enumerate}
		\item On the level of cyclic homology, it is given by $-h'b'hb$. This clearly coincides with our definition
		      up to two extra boundary terms. However, the expression $-h'b'hb$ is not a chain map on
		      $\left( \tensr{A}, b \right)$, nor does it give a well-defined map on the quotient
		      $\left( \tensr{A} / \Im \left( \idd - \t \right), b \right)$. It only becomes well-defined
		      on cyclic homology classes.
		\item Another description is given in terms of a \textit{chain map}
		      \begin{equation*}
			      b^{[2]} \colon \left( \tensr{A}^{*}, b \right) \rightharpoonup \left( \tensr{A}^{*+2}, b \right)
		      \end{equation*}
		      which has the form $b^{[2]} = \left[ b, \beta \right]$ for some
		      degree one map $\beta \colon \tensr{A}^{*} \rightharpoonup \tensr{A}^{*+1}$. The periodicity
		      operator then acts on a cyclic homology class $x \in \hcyc{\mathcal{A}}[d]$ as
		      $c_d \cdot b^{[2]} \left( x \right)$	for some constant $c_d \in \mathbb{Q}$.

		      This is somewhat analogous to our description of $S$ as $S = [b, \beta]$ but our $\beta$
		      (defined by \cref{eq:def-beta}) does not coincide with Loday's $\beta$. While the operators
		      $\beta$ and $b^{[2]}$ of Loday have simple descriptions in terms of the face maps,
		      the operator $b^{[2]}$ does not descend to a well-defined operator on Connes' complex
		      (see \cite[Exercise E.2.2.6]{Loday1998}) and is only well-defined on homology classes,
		      unlike our $S$. The method of constructing $b^{[2]}$ depends heavily on the
		      simplicial structure one has when working with an associative algebra and does not seem
		      to generalize easily to the case of a general $\Ainf$-algebra.
	\end{enumerate}
\end{rem}

\subsection{Three More Complexes} \label{sec:three-more-models}
The $2$-periodic long exact sequence
\begin{equation*}
	\cdots \leftarrow \left( \tensr{A}, b \right) \xleftarrow{\idd - \t} \left( \tensr{A}, b' \right)
	\xleftarrow{\N} \left( \tensr{A}, b \right) \leftarrow \cdots
\end{equation*}
comes associated with two natural short exact sequences:
\begin{align}
	0 \rightarrow {\left( \Im \left( \idd - \t \right), b \right)} \rightarrow{} & {\left( \tensr{A}, b \right)}
	\rightarrow {\left( \tensr{A} / \Im \left( \idd - \t \right), b \right)} \rightarrow 0
	\label{eq:im-1-t-exact-sequence}
	\\
	0 \rightarrow {\left( \Im \left( \N \right), b' \right)} \rightarrow{}       & {\left( \tensr{A}, b' \right)}
	\rightarrow {\left( \tensr{A} / \Im \left( \N \right), b' \right)} \rightarrow 0.
	\label{eq:im-N-exact-sequence}
\end{align}
Since we have $\Im \left( \idd - \t \right) = \ker \left( \N \right)$, the chain map
$\N \colon \left( \tensr{A}, b \right) \rightarrow \left( \tensr{A}, b' \right)$ induces an isomorphism
$\overline{\N} \colon \left( \tensr{A} / \Im \left( \idd - \t \right), b \right) \rightarrow
	\left( \Im \left( \N \right), b' \right)$. The inverse of $\overline{\N}$ is
the map $\overline{h'} \colon \left( \Im \left( \N \right), b' \right) \rightarrow
	\left( \tensr{A} / \Im \left( \idd - \t \right), b \right)$ induced by the map $h'$ by precomposing
$h'$ with the inclusion $\Im \left( \N \right) \hookrightarrow \tensr{A}$ and projecting the result onto
$\tensr{A} / \Im \left( \idd - \t \right)$. Note that even though
$h' \colon \left( \tensr{A}, b' \right) \rightarrow \left( \tensr{A}, b \right)$ is not a chain map, the induced map
$\overline{h'}$ becomes a chain map, a consequence of \cref{eq:b-1-t-rel,eq:b'-N-rel,eq:h'N-homotopy}.
Similarly, since $\Im \left( \N \right) = \ker \left( \idd - \t \right)$, the chain map
$\idd - \t \colon \left( \tensr{A}, b' \right) \rightarrow \left( \tensr{A}, b \right)$
induces an isomorphism $\overline{\idd - \t} \colon \left( \tensr{A} / \Im \left( \N \right), b' \right) \rightarrow
	\left( \Im \left( \idd - \t \right), b \right)$ whose inverse
$\overline{h}$ is a chain map induced by the map $h$.

\begin{figure}[htb]
	\begin{tikzcd}
		0 & {\left( \Im \left( \idd - \t \right), b \right)} & {\left( \tensr{A}, b \right)} &
		{\left( \tensr{A} / \Im \left( \idd - \t \right), b \right)} & 0 \\
		\\
		0 & {\left( \Im \left( \N \right), b' \right)} & {\left( \tensr{A}, b' \right)} &
		{\left( \tensr{A} / \Im \left( \N \right), b' \right)} & 0
		\arrow[from=1-1, to=1-2]
		\arrow[from=1-2, to=1-3]
		\arrow[from=1-3, to=1-4]
		\arrow[from=1-4, to=1-5]
		\arrow[from=3-1, to=3-2]
		\arrow[from=3-2, to=3-3]
		\arrow[from=3-3, to=3-4]
		\arrow[from=3-4, to=3-5]
		\arrow["{\overline{\N}}"'{pos=0.9}, from=1-4, to=3-2]
		\arrow["{\overline{\idd - \t}}"{pos=0.9}, from=3-4, to=1-2]
		\arrow["{\overline{h'}}"'{pos=0.9}, curve={height=6pt}, from=3-2, to=1-4]
		\arrow["{\overline{h}}"{pos=0.9}, curve={height=-6pt}, from=1-2, to=3-4]
	\end{tikzcd}
	\caption{Two short exact sequences associated to a 2-periodic long exact sequence.}
	\label{fig:two-short-exact-sequences-2-periodic}
\end{figure}

The relations between the two short exact sequences \eqref{eq:im-1-t-exact-sequence} and \eqref{eq:im-N-exact-sequence}
are depicted in \cref{fig:two-short-exact-sequences-2-periodic}. The six complexes appearing in
\cref{fig:two-short-exact-sequences-2-periodic} are:
\begin{enumerate}
	\item The Hochschild complex $\left( \tensr{A}, b \right)$ and the bar complex $\left( \tensr{A}, b' \right)$.
	\item The pair of complexes $\left( \tensr{A} / \Im \left( \idd - \t \right), b \right)$ (Connes' complex)
	      and $\left( \Im \left( \N \right), b' \right)$ which are isomorphic via
	      $\overline{\N}$ and $\overline{h'}$.
	\item The pair of complexes $\left( \tensr{A} / \Im \left( \N \right), b' \right)$ and
	      $\left( \Im \left( \idd - \t \right), b \right)$ which are isomorphic via
	      $\overline{\idd - \t}$ and $\overline{h}$.
\end{enumerate}

When $\mathcal{A}$ is unital, the bar complex $\left( \tensr{A}, b' \right)$ is contractible. By considering
the long exact sequence in cohomology associated to \eqref{eq:im-N-exact-sequence},
we see that the connecting morphism $\delta$ gives us isomorphisms $\cohom{\tensr{A} / \Im \left( \N \right)}[*][b'] \cong \cohom{\Im \left( \N \right)}[*+1][b']$ and hence we have isomorphisms
\begin{equation*}
	\cohom{\Im \left( \idd - \t \right)}[*][b] \cong
	\cohom{\tensr{A} / \Im \left( \N \right)}[*][b'] \cong
	\cohom{\Im \left( \N \right)}[*+1][b'] \cong
	\cohom{\tensr{A} / \Im \left( \idd - \t \right)}[*+1][b].
\end{equation*}

In what follows, we will describe an explicit chain map $\left( \tensr{A} / \Im \left( \N \right), b' \right)
	\rightarrow \left( \Im \left( \N \right), b' \right)[1]$ which induces the connecting morphism on cohomology
and show that in the unital case, all four complexes
\begin{equation}
	\left( \Im \left( \idd - \t \right), b \right), \
	\left( \tensr{A} / \Im \left( \N \right), b' \right), \
	\left( \Im \left( \N \right), b' \right)[1], \
	\left( \tensr{A} / \Im \left( \idd - \t \right), b \right)[1]
	\label{eq:four-models-cyclic-homology}
\end{equation}
are homotopy equivalent with explicit chain maps which induce the equivalences. The homotopy equivalences
between the four models are described by the following lemma:

\begin{lm} \label{lm:T-T'-h-h'-H-commute-up-to-homotopy}
	Assume that $\mathcal{A}$ is unital and let $H$ be a contraction of the bar complex such as the one given by
	\cref{eq:contraction-bar-with-mu-0}.
	The maps $T, T', h,  h', H \colon \tensr{A} \rightharpoonup \tensr{A}$ induce \textbf{chain maps} between the complexes
	described in the following diagram:
	\begin{figure}[H]
		\begin{tikzcd}
			{\left( \Im \left( \idd - \t \right), b \right)} && {\left( \Im \left( \N \right), b' \right)} &&
			{\left( T = -\partial \left( h \right) \right)} \\
			{\left( \tensr{A} / \Im \left( \N \right), b' \right)} &&
			{\left( \tensr{A} / \Im \left( \idd - \t \right), b \right)} &&
			{\left( T' = \partial \left( h' \right) \right)}
			\arrow["{\overline{T}}", dashed, harpoon, from=1-1, to=1-3] 
			\arrow["\overline{h}"', dashed, from=1-1, to=2-1]
			\arrow["{\overline{h'}}", dashed, from=1-3, to=2-3]
			\arrow["{\overline{T'}}"', harpoon, dashed, from=2-1, to=2-3] 
			\arrow["{-\overline{H}}"', harpoon, dashed, from=1-3, to=2-1]
		\end{tikzcd}
		\caption{Chain maps induced by $T,T',h,h',H$.}
		\label{fig:chain-maps-T-T'-h-h'-H}
	\end{figure}
	The diagram commutes up to homotopy and all the maps in the diagram are homotopy equivalences.
\end{lm}
\begin{proof}
	The identities \eqref{eq:T-b} and \eqref{eq:T-1-t} show that $T \colon \left( \tensr{A}, b \right) \rightharpoonup
		\left( \tensr{A}, b' \right)$ is a degree one chain map which maps $\Im \left( \idd - \t \right)$ into
	$\Im \left( \N \right)$ and hence induces $\overline{T}$.
	Similarly, identities \eqref{eq:T'-b'} and \eqref{eq:1-t-T} show
	that $T' \colon \left( \tensr{A}, b' \right) \rightharpoonup \left( \tensr{A}, b \right)$ is a degree one chain map
	which maps $\Im \left( \N \right)$ into $\Im \left( \idd - \t \right)$ and hence induces a chain map
	$\overline{T'}$ on the quotients. The identity $\left( b' H + H b' \right) \N = \N$ shows that
	if we restrict $H$ to $\Im \left( \N \right)$ and project the result onto $\tensr{A} / \Im \left( \N \right)$,
	the map $H$ becomes a degree $-1$ chain map $\overline{H}$.

	The outer square of \cref{fig:chain-maps-T-T'-h-h'-H} commutes up to homotopy because of the identity
	\begin{equation*}
		\partial \left( h' h \right) = \partial \left( h' \right) h + h' \partial \left( h \right) =
		T' h - h' T.
	\end{equation*}
	The upper triangle of \cref{fig:chain-maps-T-T'-h-h'-H} commutes up to homotopy because of the identity
	\begin{equation*}
		\partial \left( H h \right) = h - H \partial \left( h \right) = h + HT
	\end{equation*}
	while the lower triangle of \cref{fig:chain-maps-T-T'-h-h'-H} commutes up to homotopy because of the identity
	\begin{equation*}
		\partial \left( h' H \right) = \partial \left( h' \right) H + h' = T' H + h'.
	\end{equation*}

	Next, let us show that the map $\overline{H}$ is a homotopy equivalence.
	We have $-\overline{H} \, \overline{T} \sim \overline{h}$ and by composing with $\overline{\idd - \t}$ on the
	right we see that
	\begin{equation*}
		-\overline{H} \, \overline{T} \left( \overline{\idd - \t} \right) \sim
		\overline{h} \left( \overline{\idd - \t} \right) = \idd_{\tensr{A} / \Im \left( \N \right)}.
	\end{equation*}
	Similarly, we have $-\overline{T'} \, \overline{H} \sim \overline{h'}$ and by composing with $\overline{\N}$
	on the left, we see that
	\begin{equation*}
		- \overline{T} \left( \overline{\idd - \t} \right) \overline{H}
		\stackrel{\eqref{eq:T-1-t}}{=}
		- \overline{\N} \, \overline{T'} \, \overline{H} \sim \overline{\N} \, \overline{h'} =
		\idd_{\Im \left( \N \right)}.
	\end{equation*}
	Hence, $\overline{H}$ is a homotopy equivalence with homotopy inverse
	$- \overline{T} \left( \overline{\idd - \t} \right)$. Finally, since the maps
	$\overline{h}, \overline{h'}$ are isomorphisms and $\overline{H}$ is a homotopy equivalence, the
	identities $-\overline{H} \, \overline{T} \sim \overline{h}$ and
	$-\overline{T'} \, \overline{H} \sim \overline{h'}$ imply that $\overline{T}$ and $\overline{T'}$ are also
	homotopy equivalences.
\end{proof}

\begin{lm} \label{lm:-T-1-t-induces-connecting-morphism}
	The degree one chain map $-\overline{T} \left( \overline{\idd - \t} \right) \colon
		\left( \tensr{A} / \Im \left( \N \right), b' \right) \rightharpoonup \left( \Im \left( \N \right), b' \right)$
	induces a map $\cohom{\tensr{A} / \Im \left( \N \right)}[*][b'] \rightarrow
		\cohom{\Im \left( \N \right)}[*+1][b']$ on cohomology which coincides with the connecting morphism $\delta$
	of the long exact sequence in cohomology associated to \eqref{eq:im-N-exact-sequence}. When
	$\mathcal{A}$ is unital, the map $-\overline{T} \left( \overline{\idd - \t} \right)$ is a homotopy equivalence
	with homotopy inverse given by $\overline{H}$.
\end{lm}
\begin{proof}
	Let us denote by $\eqcl{x}_{\N}$ the projection of $x \in \tensr{A}$ onto
	$\tensr{A} / \Im \left( \N \right)$. Then, if $\eqcl{x}_{\N}$ is a closed element of
	$\left( \tensr{A} / \Im \left( \N \right), b' \right)$, then $b'x = \N y$ for some $y \in \tensr{A}$
	and $\delta \left( \eqcl{x}_{\N} \right)$ is represented by $\N y$. Now,
	\begin{equation*}
		\begin{aligned}
			-T \left( \idd - \t \right) x
			 & \stackrel{\eqref{eq:T-1-t}}{=} -\N T' x
			\stackrel{\eqref{eq:T'-def}}{=} -\N \left( bh' - h'b' \right)x
			\stackrel{\eqref{eq:b'-N-rel}}{=}
			\N h' \N y - b' \N h' x
			\\
			 & \stackrel{\eqref{eq:Nh'-homotopy}}{=}
			\N y - h\left( \idd - \t \right) \N y - b' \N h' x
			\stackrel{\eqref{eq:1-t-N=0}}{=} \N y - b' \N h' x
		\end{aligned}
	\end{equation*}
	which shows that $-T \left( \idd - \t \right) x$ coincides with $\N y$ up to a boundary term.
	The last part of \cref{lm:-T-1-t-induces-connecting-morphism} was already shown in the proof of
	\cref{lm:T-T'-h-h'-H-commute-up-to-homotopy}.
\end{proof}

\begin{cor} \label{cor:beta-homotopy-equivalence}
	The map $\beta$ given by \cref{eq:def-beta} induces a degree one chain map
	\begin{equation*}
		\overline{\beta} \colon \left( \Im \left( \idd - \t \right), b \right) \rightharpoonup
		\left( \tensr{A} / \Im \left( \idd - \t \right), b \right)
	\end{equation*}
	which is a homotopy equivalence when $\mathcal{A}$ is unital.
	A homotopy inverse for $\overline{\beta}$ is given by the map $-\overline{B}$
	where
	\begin{equation*}
		\overline{B} \colon \left( \tensr{A} / \Im \left( \idd - \t \right), b \right) \rightharpoonup
		\left( \Im \left( \idd - \t \right), b \right)
	\end{equation*}
	is induced by the map $B$ given by \cref{eq:B-def}.
\end{cor}
\begin{proof}
	We have $\overline{\beta} = \overline{T'} \, \overline{h}$ and hence by
	\cref{lm:T-T'-h-h'-H-commute-up-to-homotopy}, $\overline{\beta}$ is a
	homotopy equivalence when $\mathcal{A}$ is unital.
	We have
	\begin{equation*}
		-\overline{\beta} \, \overline{B} = -\overline{T'} \, \overline{h} \left( \overline{\idd - \t} \right)
		\overline{H} \, \overline{\N} = -\overline{T'} \, \overline{H} \, \overline{\N}
		\sim
		\overline{h'} \, \overline{\N} = \idd_{\tensr{A} / \Im \left( \idd - \t \right)}
	\end{equation*}
	and hence $-\overline{B}$ is a right homotopy inverse of $\overline{\beta}$ and, since $\overline{\beta}$ is an equivalence, also
	a left homotopy inverse of $\overline{\beta}$.
\end{proof}

\begin{rem}
	One can also prove \cref{cor:beta-homotopy-equivalence} directly using the identities
	\begin{align*}
		-\beta B                          & = \idd - \left[ b, h' H \N \right]
		+ \left( \idd - \t \right) h \left( \left[ b, h' H \N \right] - \idd \right),
		\\
		-B \beta \left( \idd - \t \right) & = \left( \idd - \left[ b, \left( \idd - \t \right) H h \right] \right)
		\left( \idd - \t \right),
	\end{align*}
	which hold on $\tensr{A}$ and are proven using straightforward computations.
\end{rem}

Finally, since the four complexes in \eqref{eq:four-models-cyclic-homology} are homotopy equivalent (up to
a shift) and compute cyclic homology, we can ask how the periodicity operator and Connes' exact sequence
are realized in each of the models. To state the relations, note that the Connes boundary map
$B \colon \tensr{A} \rightharpoonup \tensr{A}$ factors as
\begin{equation*}
	\begin{tikzcd}
		{\left( \tensr{A}, b \right)} &&& {\left( \tensr{A}, b \right)} \\
		{\left( \tensr{A} / \Im \left( \idd - \t \right), b \right)} &
		{\left( \Im \left( \N \right), b' \right)} & {\left( \tensr{A} / \Im \left( \N \right), b' \right)} &
		{\left( \Im \left( \idd - \t \right), b \right)}
		\arrow["\pi"', from=1-1, to=2-1]
		\arrow["{\overline{\N}}", from=2-1, to=2-2]
		\arrow["{\overline{H}}", harpoon, from=2-2, to=2-3]
		\arrow["{\overline{\idd - \t}}", from=2-3, to=2-4]
		\arrow["k"', from=2-4, to=1-4]
		\arrow["B", harpoon, from=1-1, to=1-4]
		\arrow["{\overline{B}}"', harpoon, curve={height=18pt}, from=2-1, to=2-4]
	\end{tikzcd}
\end{equation*}
Then we have the following commutative diagram:
\begin{figure}[H]
	\centering
	\adjustbox{scale=0.75, center}
	{
		\begin{tikzcd}[column sep=small]
			\cdots & {\cohom{{\totc{\cycbi[A]^{\{2\}}}[]}}} & {\cohom{{\totc{\cycbi[A]}[]}}} &
			{\cohom{{\totc{\cycbi[A]}[]}}[*+2]} & {\cohom{{\totc{\cycbi[A]^{\{2\}}}[]}}[*+1]} & \cdots \\
			\cdots & {\hhoch{\mathcal{A}}} & {\hcyc{\mathcal{A}}} & {\hcyc{\mathcal{A}}[*+2]} &
			{\hhoch{\mathcal{A}}[*+1]} & \cdots \\
			\cdots & {\hhoch{\mathcal{A}}} & {\cohom{\Im \left( \N \right)}[*][b']} &
			{\cohom{\Im \left( \N \right)}[*+2][b']} & {\hhoch{\mathcal{A}}[*+1]} & \cdots \\
			\cdots & {\hhoch{\mathcal{A}}} & {\cohom{{\tensr{A} / \Im \left( \N \right)}}[*-1][b']} &
			{\cohom{{\tensr{A} / \Im \left( \N \right)}}[*+1][b']} & {\hhoch{\mathcal{A}}[*+1]} & \cdots \\
			\cdots & {\hhoch{\mathcal{A}}} & {\cohom{\Im \left( \idd - \t \right)}[*-1][b]} &
			{\cohom{\Im \left( \idd - \t \right)}[*+1][b]} & {\hhoch{\mathcal{A}}[*+1]} & \cdots \\
			\cdots & {\hhoch{\mathcal{A}}} & {\hcyc{\mathcal{A}}} & {\cohom{\Im \left( \idd - \t \right)}[*+1][b]} &
			{\hhoch{\mathcal{A}}[*+1]} & \cdots
			\arrow["i", from=1-2, to=1-3]
			\arrow["\sigma", harpoon, from=1-3, to=1-4]
			\arrow[draw=none, from=1-4, to=1-5]
			\arrow["p", from=1-3, to=2-3]
			\arrow["p", from=1-4, to=2-4]
			\arrow["S", harpoon, from=2-3, to=2-4]
			\arrow["{k \overline{B}}", harpoon, from=2-4, to=2-5]
			\arrow["j", from=2-2, to=1-2]
			\arrow["I", from=2-2, to=2-3]
			\arrow["{\overline{\N}}", from=2-3, to=3-3]
			\arrow["{\overline{H}}", harpoon, from=3-3, to=4-3]
			\arrow["{\overline{H} \, \overline{\N} I}", harpoon, from=4-2, to=4-3]
			\arrow["{\overline{\N} I}", from=3-2, to=3-3]
			\arrow["\delta", harpoon, from=1-4, to=1-5]
			\arrow["j", from=2-5, to=1-5]
			\arrow["{\overline{\N}}", from=2-4, to=3-4]
			\arrow["{S'}", harpoon, from=3-3, to=3-4]
			\arrow["{k \left( \overline{\idd - \t} \right) \overline{H}}", harpoon, from=3-4, to=3-5]
			\arrow["{\overline{H}}", harpoon, from=3-4, to=4-4]
			\arrow["{S'}", harpoon, from=4-3, to=4-4]
			\arrow["{k \left( \overline{\idd - \t} \right)}", from=4-4, to=4-5]
			\arrow["{\overline{\idd - \t}}", from=4-3, to=5-3]
			\arrow["{\overline{\idd - \t}}", from=4-4, to=5-4]
			\arrow["S", harpoon, from=5-3, to=5-4]
			\arrow["{k}", from=5-4, to=5-5]
			\arrow["-\beta", harpoon, from=5-3, to=6-3]
			\arrow["{\overline{B}}"', harpoon, curve={height=85pt}, from=2-3, to=5-3]
			\arrow["{\overline{B}}", harpoon, curve={height=-85pt}, from=2-4, to=5-4]
			\arrow["\delta", harpoon, from=6-3, to=6-4]
			\arrow[equal, from=5-4, to=6-4]
			\arrow["k", from=6-4, to=6-5]
			\arrow["I", from=6-2, to=6-3]
			\arrow["{\overline{B} I}", harpoon, from=5-2, to=5-3]
			\arrow[equal, from=6-5, to=5-5]
			\arrow[equal, from=5-5, to=4-5]
			\arrow[equal, from=4-5, to=3-5]
			\arrow[equal, from=3-5, to=2-5]
			\arrow[equal, from=2-2, to=3-2]
			\arrow[equal, from=3-2, to=4-2]
			\arrow[equal, from=4-2, to=5-2]
			\arrow[equal, from=5-2, to=6-2]
			\arrow[from=1-1, to=1-2]
			\arrow[from=2-1, to=2-2]
			\arrow[from=3-1, to=3-2]
			\arrow[from=4-1, to=4-2]
			\arrow[from=5-1, to=5-2]
			\arrow[from=6-1, to=6-2]
			\arrow[from=1-5, to=1-6]
			\arrow[from=2-5, to=2-6]
			\arrow[from=3-5, to=3-6]
			\arrow[from=4-5, to=4-6]
			\arrow[from=5-5, to=5-6]
			\arrow[from=6-5, to=6-6]
		\end{tikzcd}
	}
	\caption{Different versions of Connes' long exact sequence.}
\end{figure}

Hence, we see that the periodicity operator on
$\left( \Im \left( \idd - \t \right), b \right), \left( \tensr{A} / \Im \left( \idd - \t \right), b \right)$
is induced by $S$ while the periodicity operator on
$\left( \Im \left( \N \right), b' \right), \left( \tensr{A} / \Im \left( \N \right), b' \right)$ is induced by
$S'$. We also see that the short exact sequence \eqref{eq:im-1-t-exact-sequence} can be identified with
Connes' exact sequence where now the connecting morphism $\delta$ plays the role of the periodicity operator
$S$ while the inclusion $j$ plays the role of Connes' boundary map.

\begin{rem}
	We offer another perspective on the operators $T,T'$ and how they appear naturally.
	The two short exact sequences \eqref{eq:im-1-t-exact-sequence} and \eqref{eq:im-N-exact-sequence} are
	termwise split exact sequences, where the splittings are induced naturally from the homotopies of the rows,
	i.e., they come from the identities \eqref{eq:h'N-homotopy} and \eqref{eq:Nh'-homotopy}.
	The situation is depicted in \cref{fig:two-short-exact-sequences-2-periodic-splitting}.
	In the figure, the arrows with tails are the termwise splittings
	and \textit{are not} chain maps. All other arrows are chain maps.
	The dashed arrows are isomorphisms, the squiggly arrows are homotopy equivalences,
	the arrows which end in a harpoon correspond to degree one maps, while all other arrows correspond to
	maps of degree zero.

	Since the first exact sequence is split exact, the connecting homomorphism of the sequence is induced by the chain map
	\begin{equation*}
		\begin{aligned}
			\overline{\left( \idd - \t \right) h} b \overline{h' \N}
			\stackrel{\eqref{eq:h'N-homotopy}}{=}{} &
			\overline{ \left( \idd - h' \N \right) b h' \N}
			\stackrel{\eqref{eq:b'-N-rel}}{=}{}
			\overline{bh' \N - h' b' \N h' \N}
			\\
			\stackrel{\eqref{eq:Nh'-homotopy}}{=}{} &
			\overline{bh' \N - h' b' \left( \idd - h \left( \idd - \t \right) \right) \N}
			\stackrel{\eqref{eq:1-t-N=0}}{=}
			\overline{bh' \N - h' b' \N}
			\\
			\stackrel{\eqref{eq:T'-def}}{=}{}       &
			\imsub{\overline{T'}} \, \overline{\N}.
		\end{aligned}
	\end{equation*}
	Hence, up to the isomorphism $\overline{\N}$, the map
	$\imsub{\overline{T'}} \colon \left( \Im \left( \N \right), b' \right) \rightharpoonup \left( \Im \left( \idd - \t \right), b \right)$
	is a chain map which induces the connecting homomorphism and,
	thinking of $\imsub{\overline{T'}}$ as a degree zero chain map
	$\left( \Im \left( \N \right), b' \right)[-1] \rightarrow \left( \Im \left( \idd - \t \right), b \right)$,
	we have
	\begin{equation*}
		\Cone{\imsub{\overline{T'}}} \cong \Cone{\imsub{\overline{T'}} \, \overline{\N}} \cong \left( \tensr{A}, b \right)
	\end{equation*}
	so the mapping cone of $\imsub{\overline{T'}}$ computes the Hochschild homology.

	Similarly, the connecting homomorphism of the second exact sequence is induced by the chain map
	\begin{equation*}
		\begin{aligned}
			\overline{\N h'} b' \overline{h \left( \idd - \t \right)}
			\stackrel{\eqref{eq:Nh'-homotopy}}{=} &
			\overline{ \left( \idd - h \left( \idd - \t \right) \right) b' h \left( \idd - \t \right)}
			\stackrel{\eqref{eq:b-1-t-rel}}{=}
			\overline{b'h \left( \idd - \t \right) - h b \left( \idd - \t \right) h \left( \idd - \t \right)}
			\\
			\stackrel{\eqref{eq:h'N-homotopy}}{=} &
			\overline{b'h \left( \idd - \t \right) - h b \left( \idd - h' \N \right) \left( \idd - \t \right)}
			\stackrel{\eqref{eq:1-t-N=0}}{=}
			\overline{b'h \left( \idd - \t \right) - h b \left( \idd - \t \right)}
			\\
			\stackrel{\eqref{eq:T-def}}{=}        &
			- \imsub{\overline{T}} \left( \overline{ \idd - \t } \right),
		\end{aligned}
	\end{equation*}
	as we have already seen in \cref{lm:-T-1-t-induces-connecting-morphism}.
	Hence, up to the isomorphism $\overline{ \left( \idd - \t \right)}$, the map
	$\imsub{\overline{T}} \colon \left( \Im \left( \idd - \t \right), b \right) \rightharpoonup \left( \Im \left( \N \right), b' \right)$
	is a chain map which induces the connecting homomorphism of the second sequence and,
	thinking of $\imsub{\overline{T}}$ as a degree zero chain map
	$\left( \Im \left( \idd - \t \right), b \right)[-1] \rightarrow \left( \Im \left( \N \right), b' \right)$,
	we have $\Cone{\imsub{\overline{T}}} \cong \left( \tensr{A}, b' \right)$.
\end{rem}

\begin{figure}[htb]
	\begin{tikzcd}[column sep = large]
		0 & {\left( \Im \left( \idd - \t \right), b \right)} & {\left( \tensr{A}, b \right)} &
		{\left( \tensr{A} / \Im \left( \idd - \t \right), b \right)} & 0 \\
		\\
		0 & {\left( \Im \left( \N \right), b' \right)} & {\left( \tensr{A}, b' \right)} &
		{\left( \tensr{A} / \Im \left( \N \right), b' \right)} & 0
		\arrow[from=1-1, to=1-2]
		\arrow[from=1-2, to=1-3]
		\arrow[start anchor=north, end anchor=north east, "{\overline{\left( \idd - \t \right)h}}"',
			from=1-3, to=1-2, curve={height=15pt}, tail]
		\arrow[from=1-3, to=1-4]
		\arrow[start anchor=north, end anchor=north east, "{\overline{h' \N}}"',
			from=1-4, to=1-3, curve={height=15pt}, tail]
		\arrow[from=1-4, to=1-5]
		\arrow[from=3-1, to=3-2]
		\arrow[from=1-2, to=3-2, "\imsub{\overline{T}}", harpoon, squiggly, curve={height=-10pt}]
		\arrow[from=3-2, to=1-2, "\imsub{\overline{T'}}", harpoon, curve={height=-10pt}]
		\arrow[from=3-2, to=3-3]
		\arrow[start anchor=south, end anchor=south east, "{\overline{\N h'}}",
			from=3-3, to=3-2, curve={height=-15pt}, tail]
		\arrow[from=3-3, to=3-4]
		\arrow[start anchor=south, end anchor=south east, "{\overline{h \left( \idd - \t \right)}}",
			from=3-4, to=3-3, curve={height=-15pt}, tail]
		\arrow[from=3-4, to=3-5]
		\arrow[from=3-4, to=1-4, "\quotsub{\overline{T'}}", harpoon, squiggly, curve={height=-10pt}]
		\arrow[from=1-4, to=3-4, "\quotsub{\overline{T}}", harpoon, curve={height=-10pt}]
		\arrow["{\overline{\N}}"'{pos=0.9}, dashed, from=1-4, to=3-2]
		\arrow["{\overline{\idd - \t}}"{pos=0.9}, dashed, from=3-4, to=1-2]
		\arrow["{\overline{h'}}"'{pos=0.9}, dashed, curve={height=6pt}, from=3-2, to=1-4]
		\arrow["{\overline{h}}"{pos=0.9}, dashed, curve={height=-6pt}, from=1-2, to=3-4]
	\end{tikzcd}
	\caption{The two short exact sequences associated to a 2-periodic long exact sequence, together with
		their splittings and maps between their components.}
	\label{fig:two-short-exact-sequences-2-periodic-splitting}
\end{figure}

\subsection{Total Complexes of Cyclic Codifferential Forms}
\label{sec:bicomplex-models-cyclic-homology}
Let $\mathcal{R} = (R,d)$ be a differential graded-commutative Banach $\mathbbm{k}$-algebra and let
$\mathcal{A} = \left( A, \mu\right)$ be a Banach $\Ainf$-algebra.
Consider the bigraded Banach $R$-module $\ncdfr{A} = \tensrcyc{A \oplus \ul{A}}[(*,*)]$
together with the operators $\qdr \colon \ncdfr{A} \rightharpoonup \ncdfr{A}[* - 1][*]$
and $\clie{\mu} \colon \ncdfr{A} \rightharpoonup \ncdfr{A}[*][* + 1]$.
Note that $\ncdfr{A}[0][] = \ncdf{A}[0][] / R$ while $\ncdfr{A}[i][] = \ncdf{A}[i][]$ for $i > 0$.
In what follows, we will work with the inner product parity form given by \cref{eq:parity-inner-product}.
As a consequence of our definitions and \cref{lm:cyclic-commutation-relations}, we
get the following relations on $\ncdfr{A}$:
\begin{align}
	\qdr^2                          & = 0,
	\label{eq:qdr-differential}
	\\
	\clie{\mu}^2                    & = \frac{1}{2} \left[ \clie{\mu}, \clie{\mu} \right] = \frac{1}{2} \clie{[\mu,\mu]} = 0,
	\label{eq:clie-differential}
	\\
	\left[ \qdr, \clie{\mu} \right] & = \qdr \circ \clie{\mu} - \clie{\mu} \circ \qdr = 0.
	\label{eq:qdr-clie-commute}
\end{align}
In particular, we see that both $\qdr$ and $\clie{\mu}$ are differentials of
degree $(-1,0)$ and $(0,1)$ respectively which commute with each other.
Hence, we obtain a right half-plane Banach bicomplex $\left( \ncdfr{A}, \qdr, \clie{\mu} \right)$ (see \cref{fig:bicomplex-ncdf}).
We will denote the bicomplex by $\ncdfr{\mathcal{A}}[][]$, suppressing
the differentials from the notation. Similarly, we will denote by $\ncdfr{\mathcal{A}}[0][]$
the zeroth column of the bicomplex and by $\ncdf{\mathcal{A}}[i][]$ the $i$-th column of the bicomplex
when $i > 0$, endowed with the differential $\clie{\mu}$.
\begin{figure}[htb]
	\centering
	\begin{tikzcd}
		& \vdots & {\vdots} & {\vdots } & {\vdots } \\
		{\cdots} & 0 & \ncdfr{A}[0][1] & \ncdf{A}[1][1] & \ncdf{A}[2][1] & {\cdots } \\
		{\cdots} & 0 & \ncdfr{A}[0][0] & \ncdf{A}[1][0] & \ncdf{A}[2][0] & {\cdots } \\
		{\cdots} & 0 & \ncdfr{A}[0][-1] & \ncdf{A}[1][-1] & \ncdf{A}[2][-1] & {\cdots } \\
		& {\vdots} & {\vdots } & {\vdots } & {\vdots } \\
		\arrow[from=2-2, to=2-1]
		\arrow[from=2-3, to=2-2]
		\arrow["{\qdr}"', from=2-4, to=2-3]
		\arrow["{\qdr}"', from=2-5, to=2-4]
		\arrow[from=2-6, to=2-5]
		\arrow[from=3-2, to=3-1]
		\arrow[from=3-3, to=3-2]
		\arrow["{\qdr}"', from=3-4, to=3-3]
		\arrow["{\qdr}"', from=3-5, to=3-4]
		\arrow[from=3-6, to=3-5]
		\arrow[from=4-2, to=4-1]
		\arrow[from=4-3, to=4-2]
		\arrow["{\qdr}"', from=4-4, to=4-3]
		\arrow["{\qdr}"', from=4-5, to=4-4]
		\arrow[from=4-6, to=4-5]
		\arrow[from=2-2, to=1-2]
		\arrow[from=3-2, to=2-2]
		\arrow[from=4-2, to=3-2]
		\arrow[from=5-2, to=4-2]
		\arrow["{\clie{\mu}}", from=2-3, to=1-3]
		\arrow["{\clie{\mu}}", from=3-3, to=2-3]
		\arrow["{\clie{\mu}}", from=4-3, to=3-3]
		\arrow["{\clie{\mu}}", from=5-3, to=4-3]
		\arrow["{\clie{\mu}}", from=2-4, to=1-4]
		\arrow["{\clie{\mu}}", from=3-4, to=2-4]
		\arrow["{\clie{\mu}}", from=4-4, to=3-4]
		\arrow["{\clie{\mu}}", from=5-4, to=4-4]
		\arrow["{\clie{\mu}}", from=2-5, to=1-5]
		\arrow["{\clie{\mu}}", from=3-5, to=2-5]
		\arrow["{\clie{\mu}}", from=4-5, to=3-5]
		\arrow["{\clie{\mu}}", from=5-5, to=4-5]
	\end{tikzcd}
	\caption{The Bicomplex $\ncdfr{\mathcal{A}}[][]$.}
	\label{fig:bicomplex-ncdf}
\end{figure}

The first two columns of $\ncdfr{\mathcal{A}}[][]$ we already introduced before. The zeroth column
$\ncdfr{\mathcal{A}}[0][]$ coincides with Connes' complex
$\left( \tensr{A} / \Im \left( \idd - \t \right), b \right)$
while the first column $\ncdf{\mathcal{A}}[1][]$ can be naturally identified with the Hochschild complex
$\left( \tensr{A}, b \right)$. More precisely, let
$\varphi \colon \left( \ndf{A}[1][], \clie{\mu} \right) \rightarrow \left( \tensr{A}, b \right)$
be the map given by
\begin{equation*}
	\varphi \left( l \otimes \ul{x} \otimes s \right) =
	(-1)^{\degb{l} \cdot \left( \degb{x} + \degb{s} \right)} x \otimes s \otimes l
\end{equation*}
where $x \in A$ and $l,s \in \tens{A}$. The map $\varphi$ is a chain map which satisfies
$\varphi \circ \left( \idd - \t \right) = 0$ and hence induces a chain map
\begin{equation} \label{eq:varphi-proj-hoch}
	\varphi \colon \left( \ncdf{A}[1][], \clie{\mu} \right) \rightarrow \left( \tensr{A}, b \right),
\end{equation}
denoted by the same name, on the quotient
$\ncdf{\mathcal{A}}[1][] = \ndf{\mathcal{A}}[1][] / \Im \left( \idd - \t \right)$. The induced
map $\varphi$ is an isomorphism identifying
$\ncdf{\mathcal{A}}[1][]$ with $\choch{\mathcal{A}}[] = \left( \tensr{A}, b \right)$.
Since we adopt the convention of
writing elements of $\ncdf{\mathcal{A}}[1][]$ in the form $\ul{x} \otimes l$ where
$x \in A$ and $l \in \tens{A}$ (see page~\pageref{eq:representative-starts-with-underline}),
the identification $\varphi$ becomes $\ul{x} \otimes l \mapsto x \otimes l$ with no signs involved.
Identifying $\ncdf{\mathcal{A}}[1][]$ with $\choch{\mathcal{A}}[]$ via $\varphi$, the map
$\qdr^1 \colon \ncdf{\mathcal{A}}[1][] \rightarrow \ncdfr{\mathcal{A}}[0][]$ becomes
the natural projection map $\choch{\mathcal{A}}[] \rightarrow \cconnes{\mathcal{A}}[]$.

\begin{lm} \label{lm:tot-ncdfr-contractible}
	The total complex $\totc{\ncdfr{\mathcal{A}}[][]}[][][\coplus]$ is contractible.
\end{lm}
\begin{proof}
	By the formal Poincar\'{e} \cref{lm:formal-poincare-ncdfr}, the rows of $\ncdfr{\mathcal{A}}[][]$
	are contractible with a contraction compatible with the $R$-action and so it is enough to verify the convergence
	condition  \eqref{eq:bicomplex-convergence-cond-1} of
	\cref{lm:total-bicomplex-contractible-R-linear-contraction-1}, namely that
	\begin{equation} \label{eq:conv-cond-ncdfr-bicomplex}
		{\underbrace{\left( h_{\dr} \clie{\mu} - \clie{\mu} h_{\dr} \right)}_{T}}^n \left( x \right) \to 0
	\end{equation}
	for all $x \in \ncdfr{A}[i][j]$.
	To prove that \cref{eq:conv-cond-ncdfr-bicomplex} holds,
	we will analyze the action of $T$ on $x$ based on weight.
	Let us denote by $\prescript{}{i}{\mho}^{j}_{\textrm{cyc},k} \left( A \right)$ the collection
	of elements of line degree $i$, cohomological degree $j$ and weight $k$. Then
	$\prescript{}{i}{\mho}^{j}_{\textrm{cyc},k} \left( A \right) = 0$ if $k < i$ and
	$\ncdfr{A}[i][j] = \oplus_{k=i}^{\infty} \prescript{}{i}{\mho}^{j}_{\textrm{cyc},k} \left( A \right)$.
	Finally, denote by $\prescript{}{i}{\mathcal{F}}_l^j \defeq
		\oplus_{k \leq l} \prescript{}{i}{\mho}^{j}_{\textrm{cyc},k} \left( A \right)$
	the weight filtration on $\ncdfr{A}[i][j]$.

	By \cref{lm:direct-sum-iteration-zero-limit},
	it is enough to show that $T^n \left( x \right) \to 0$ whenever
	$x \in \prescript{}{i}{\mho}^{j}_{\textrm{cyc},k} \left( A \right)$ for $k \geq i$.
	Split the coderivation $\mu$ into three parts $\mu = \mu^0 + \mu^1 + \mu^{>1}$ where
	\begin{align*}
		\mu^0_0    & = \mu_0, \hspace{70pt} \mu^0_k = 0 \textrm{ for } k \neq 0,                          \\
		\mu^1_1    & = \mu_1, \hspace{70pt} \mu^1_k = 0 \textrm{ for } k \neq 1,                          \\
		\mu^{>1}_k & = \mu_k \textrm{ for } k > 1, \hspace{20pt} \mu^{>1}_k = 0 \textrm { for } k \leq 1.
	\end{align*}
	Set $\clie{\mu}^0 = \clie{{\mu^0}}, \clie{\mu}^1 = \clie{{\mu^1}}$, and $\clie{\mu}^{>1} = \clie{{\mu^{>1}}}$,
	so that $\clie{\mu} = \clie{\mu}^0 + \clie{\mu}^1 + \clie{\mu}^{>1}$
	and we have
	\begin{equation*}
		\clie{\mu}^0  \left( \prescript{}{i}{\mathcal{F}}_l^j \right) \subseteq
		\prescript{}{i}{\mathcal{F}}_{l+1}^{j+1}, \qquad
		\clie{\mu}^1 \left( \prescript{}{i}{\mathcal{F}}_l^j \right) \subseteq
		\prescript{}{i}{\mathcal{F}}_l^{j+1}, \qquad
		\clie{\mu}^{>1} \left( \prescript{}{i}{\mathcal{F}}_l^j \right) \subseteq
		\prescript{}{i}{\mathcal{F}}_{l-1}^{j+1}.
	\end{equation*}
	Since $h_{\dr}$ commutes with $\clie{\mu}^1$ by \cref{lm:formal-poincare-ncdfr}, we have
	\begin{equation*}
		h_{\dr} \clie{\mu} - \clie{\mu} h_{\dr} =
		\underbrace{\left( h_{\dr} \clie{\mu}^0 - \clie{\mu}^0 h_{\dr} \right)}_{T_{\mu_0}} +
		\underbrace{\left( h_{\dr} \clie{\mu}^{>1} - \clie{\mu}^{>1} h_{\dr} \right)}_{T_{-}},
	\end{equation*}
	where
	\begin{enumerate}
		\item $T_{-} \left( \prescript{}{i}{\mathcal{F}}_l^j \right) \subseteq
			      \prescript{}{i+1}{\mathcal{F}}_{l-1}^{j+1}$ and $\nnorm[T_{-}] \leq 1$.
		\item $T_{\mu_0} \left( \prescript{}{i}{\mathcal{F}}_l^j \right) \subseteq
			      \prescript{}{i+1}{\mathcal{F}}_{l+1}^{j+1}$ and $\nnorm[T_{\mu_0}] \leq \nnorm[\mu_0(1)]$.
	\end{enumerate}
	Hence, when $x \in \prescript{}{i}{\mho}^{j}_{\textrm{cyc},k} \left( A \right)$, we can write
	\begin{equation*}
		T^n \left( x \right) = \left( T_{-} + T_{\mu_0} \right)^n \left( x \right) = \sum_{r=0}^n y_r
	\end{equation*}
	where $y_r \in \prescript{}{i+n}{\mathcal{F}}_{k-n+2r}^{j+n}$ and
	$\nnorm[y_r] \leq \nnorm[\mu_0(1)]^r \cdot \nnorm[x]$ for all $0 \leq r \leq n$.
	The term $y_r$ in the decomposition above is the sum of all terms in which $T_{\mu_0}$
	is applied exactly $r$ times (and hence $T_{-}$ is applied $\left( n - r \right)$ times).
	Therefore, if we fix $r \geq 0$ and take $n > \frac{k-i}{2} + r$ then $y_0 = \dots = y_r = 0$
	and so $\nnorm[T^n \left( x \right)] \leq \nnorm[\mu_0(1)]^{r+1} \cdot \nnorm[x]$. Since
	$\nnorm[\mu_0(1)] < 1$, we deduce that $\nnorm[T^n \left( x \right)] \to 0$.
\end{proof}

Since  $\totc{\ncdfr{A}[][]}[][][\coplus]$ is contractible, by erasing the zeroth column we can obtain another
model for cyclic homology. Denote by
$\totcomp{\mathcal{A}}[1] \defeq \totc{\ncdfr{\mathcal{A}}[][]}[*][\geq 1][\coplus]$
the complete direct sum totalization of the bicomplex obtained from $\ncdfr{\mathcal{A}}[][]$
by removing the zeroth column.
Let $p_1 \colon \totcomp{\mathcal{A}}[1] \rightharpoonup \ncdfr{\mathcal{A}}[0][*+1]$
be the natural surjection given by
\begin{equation} \label{eq:def-p}
	p_1 \left( \sum_{i \geq 1} s_i \left( x_i \right) \right) \defeq \qdr^1 \left( x_1 \right).
\end{equation}

\begin{thm} \label{lm:p_1-homotopy-equivalence}
	The degree one chain map
	$p_1 \colon \totcomp{\mathcal{A}}[1] \rightharpoonup \ncdfr{\mathcal{A}}[0][*+1]$
	is a homotopy equivalence. A homotopy inverse
	$i_1 \colon \ncdfr{\mathcal{A}}[0] \rightharpoonup \totcomp{\mathcal{A}}[1][*-1]$
	for $p_1$ is given by the formula
	\begin{equation} \label{eq:def-p_1'}
		i_1 \left( x \right) = \sum_{k \geq 0} (-1)^{\frac{k(k-1)}{2}}
		s_{k+1} \left( \left( h_{\dr} \clie{\mu} \right)^k \left( h_{\dr} (x) \right) \right).
	\end{equation}
	We have $p_1 \circ i_1 = \idd$ and $i_1 \circ p_1 = \idd - \partial \left( H_{+} \right)$
	for a degree $-1$ map $H_{+} \colon \totcomp{A}[1][] \rightharpoonup \totcomp{A}[1][]$
	which satisfies $H_{+} \circ i_1 = 0$.
\end{thm}
\begin{proof}
	The theorem follows from \cref{lm:tot-ncdfr-contractible} by applying \cref{cor:total-bicomplex-projection-equivalence-1}.
	The formula \eqref{eq:def-p_1'} for $i_1$
	follows from \cref{eq:total-bicomplex-projection-equivalence-inverse-1} since $h_{\dr}^2 = 0$.
\end{proof}

\begin{ex}
	Let us assume that $\left( \mathcal{A}, \mu \right)$ corresponds to an associative algebra, so that
	$\mu_k = 0$ whenever $k \neq 2$. In this special case, $A$ is concentrated at cohomological degree $-1$
	and hence the cohomological grading on each column $\ncdf{A}[i]$ corresponds to the weight of elements in
	$\tens{A \oplus \ul{A}}$. The resulting bicomplex
	\begin{equation*}
		\left( \ncdf{A}[k][l], \qdr, \clie{\mu} \right)_{k \geq 1}^{l \in \ZZ}
	\end{equation*}
	is depicted in \cref{fig:alt-bicomplex-assoc-alg}, and is somewhat similar to the $(b,B)$-bicomplex
	computing the cyclic homology of an associative algebra (see \cite[Page 57]{Loday1998}).
	The de Rham differential $\qdr$ takes the role of $B$, and instead of having shifted copies of
	the Hochschild complex $\ncdf{\mathcal{A}}[1][]$, we have the complexes $\ncdf{\mathcal{A}}[k][]$ for $k \geq 1$.
\end{ex}

\begin{figure}[htb]
	\centering
	\adjustbox{scale=0.86,center}{
		\begin{tikzcd}
			0 \\
			{\ul{A}} & 0 \\
			{\ul{A} \otimes A} & {\ul{A}\otimes \ul{A}} & 0 \\
			{\ul{A} \otimes A^{\otimes 2}} & {\ul{A} \otimes  A
				\otimes \ul{A}} & {\ul{A} \otimes \ul{A} \otimes
				\ul{A}} \\
			{\ul{A} \otimes A^{\otimes 3}} & {\left( \ul{A} \otimes
				A^{\otimes 2} \otimes \ul{A} \right) \oplus \left(
				\ul{A} \otimes A \otimes \ul{A} \otimes A
				\right)_{\mathbb{Z}_2} } & {\ul{A} \otimes A\otimes \ul{A}\otimes \ul{A}} \\
			{\ul{A} \otimes A^{\otimes 4}} & {\left( \ul{A} \otimes
				A^{\otimes 3} \otimes \ul{A} \right) \oplus \left( \ul{A}
				\otimes A^{\otimes 2} \otimes \ul{A} \otimes A \right) } &
			{\left( \ul{A} \otimes A^{\otimes 2} \otimes
				\ul{A}^{\otimes 2} \right) \oplus \left( \left( \ul{A}
				\otimes A \right)^{\otimes 2} \otimes \ul{A} \right)} \\
			{\vdots } & {\vdots } & {\vdots } \\
			\ncdf{\mathcal{A}}[1] & \ncdf{\mathcal{A}}[2] & \ncdf{\mathcal{A}}[3]
			\arrow["{\clie{\mu}}", from=2-1, to=1-1]
			\arrow["{\clie{\mu}}", from=3-1, to=2-1]
			\arrow["{\clie{\mu}}", from=4-1, to=3-1]
			\arrow["{\clie{\mu}}", from=5-1, to=4-1]
			\arrow["{\clie{\mu}}", from=6-1, to=5-1]
			\arrow["{\clie{\mu}}", from=3-2, to=2-2]
			\arrow["{\clie{\mu}}", from=4-2, to=3-2]
			\arrow["{\clie{\mu}}", from=5-2, to=4-2]
			\arrow["{\clie{\mu}}", from=6-2, to=5-2]
			\arrow["{\clie{\mu}}", from=4-3, to=3-3]
			\arrow["{\clie{\mu}}", from=5-3, to=4-3]
			\arrow["{\clie{\mu}}", from=6-3, to=5-3]
			\arrow["{\qdr^2}"', from=2-2, to=2-1]
			\arrow["{\qdr^2}"', from=3-2, to=3-1]
			\arrow["{\qdr^2}"', from=4-2, to=4-1]
			\arrow["{\qdr^2}"', from=5-2, to=5-1]
			\arrow["{\qdr^2}"', from=6-2, to=6-1]
			\arrow["{\qdr^3}"', from=3-3, to=3-2]
			\arrow["{\qdr^3}"', from=4-3, to=4-2]
			\arrow["{\qdr^3}"', from=5-3, to=5-2]
			\arrow["{\qdr^3}"', from=6-3, to=6-2]
			\arrow[from=7-1, to=6-1]
			\arrow[from=7-2, to=6-2]
			\arrow[from=7-3, to=6-3]
		\end{tikzcd}
	}
	\caption{The bicomplex $\left( \ncdf{A}[k][l], \qdr, \clie{\mu} \right)_{k \geq 1}^{l \in \ZZ}$
		for an associative algebra.}
	\label{fig:alt-bicomplex-assoc-alg}
\end{figure}

Somewhat surprisingly, when $\mathcal{A}$ is unital,
we can get yet another model for cyclic homology by erasing the first
two columns of $\ncdfr{\mathcal{A}}[][]$. To see that, we will need to show
that the cokernel of $\qdr^3 \colon \ncdf{A}[3][] \rightarrow \ncdf{A}[2][]$ is
homotopy equivalent to Connes' complex. Denote by
$\ncdf{\mathcal{A}}[2][] / \Im \left( \qdr^3 \right)$ the complex
$\big( {{\ncdf{A}[2][]}/ \Im \left( \qdr^3 \right)}, \clie{\mu} \big)$,
suppressing the differential $\clie{\mu}$ from the notation.

Since $\qdr^2 \circ \qdr^3 = 0$, and $\qdr$ commutes with $\clie{\mu}$,
we have an induced chain map
\begin{equation*}
	\overline{\qdr}^2 \colon \big( {\ncdf{A}[2][]}/ \Im \left( \qdr^3 \right), \clie{\mu} \big)
	\rightarrow
	\left( \Im \left( \qdr^2 \right), \clie{\mu} \right).
\end{equation*}
Restricting
the isomorphism $\varphi \colon \big( \ncdf{A}[1][], \clie{\mu} \big) \rightarrow \left( \tensr{A}, b \right)$
of \eqref{eq:varphi-proj-hoch} to $\Im \left( \qdr^2 \right) = \ker \left( \qdr^1 \right)$, we obtain
an isomorphism
\begin{equation*}
	\overline{\varphi} \colon \left( \Im \left( \qdr^2 \right), \clie{\mu} \right) \rightarrow
	\left( \Im \left( \idd - \t \right), b \right).
\end{equation*}

\begin{dfn} \label{dfn:psi-ncdf-2-0}
	The degree one chain map $\psi \colon {\ncdf{\mathcal{A}}[2][]}/ \Im \left( \qdr^3 \right) \rightharpoonup
		\ncdfr{\mathcal{A}}[0][]$ is defined by the composition
	\begin{equation} \label{eq:psi-def}
		\psi \defeq \overline{\beta} \circ \overline{\varphi} \circ \overline{\qdr}^2,
	\end{equation}
	where $\overline{\beta}$ is the degree one chain map from \cref{cor:beta-homotopy-equivalence}.
\end{dfn}

\begin{thm}[Homotopy Equivalence of ${\ncdf{\mathcal{A}}[2][] / \Im \left( \qdr^3 \right)}$ and ${{\ncdf{\mathcal{A}}[0][]}[1]}$]
	\label{lm:ncdf-0-ncdf-2-mod-q3-equiv}
	Let $\mathcal{A}$ be a unital Banach $\Ainf$-algebra.
	Then the degree one chain map $\psi$ of \eqref{eq:psi-def} is a homotopy equivalence.
	A homotopy inverse for $\psi$ is given by the degree minus one chain map
	$\psi' \colon \ncdfr{\mathcal{A}}[0][] \rightharpoonup
		{\ncdf{\mathcal{A}}[2][]}/ \Im \left( \qdr^3 \right)$ defined by
	\begin{equation} \label{eq:psi'-def}
		\psi' \defeq - \overline{h}_{\dr} \circ \overline{\varphi}^{-1} \circ \overline{B},
	\end{equation}
	where $\overline{B}$ is the map from \cref{cor:beta-homotopy-equivalence}.
\end{thm}
\begin{proof}
	By the formal Poincar\'{e} \cref{lm:formal-poincare-ncdfr}, we have
	$\Im \left( \qdr^3 \right) = \ker \left( \qdr^2 \right)$, and hence the induced chain map
	\begin{equation*}
		\overline{\qdr}^2 \colon \big( {\ncdf{A}[2][]}/ \Im \left( \qdr^3 \right), \clie{\mu} \big)
		\rightarrow
		\left( \Im \left( \qdr^2 \right), \clie{\mu} \right)
	\end{equation*}
	is an isomorphism. The inverse of $\overline{\qdr}^2$ is the map
	$\overline{h}_{\dr} \colon \left( \Im \left( \qdr^2 \right), \clie{\mu} \right) \rightarrow
		\big( {\ncdf{A}[2][]}/ \Im \left( \qdr^3 \right), \clie{\mu} \big)$
	induced by $h_{\dr} \colon \ncdf{A}[1][] \rightarrow \ncdf{A}[2][]$.

	When $\mathcal{A}$ is unital, \cref{cor:beta-homotopy-equivalence} shows that $\overline{\beta}$ is
	a homotopy equivalence with homotopy inverse $-\overline{B}$.
	Since $\psi = \overline{\beta} \circ \overline{\varphi} \circ \overline{\qdr}^2$ is the composition
	of two isomorphisms and a homotopy equivalence, it follows that $\psi$ is a homotopy equivalence
	with homotopy inverse given by $\psi'$.
\end{proof}

Set
$\totcomp{\mathcal{A}}[2] \defeq \totc{\ncdfr{\mathcal{A}}[][]}[*][\geq 2][\coplus]$ and let
$p_2 \colon \totcomp{\mathcal{A}}[2] \rightharpoonup
	\left( {\ncdf{\mathcal{A}}[2][]}/ \Im \left( \qdr^3 \right) \right)^{*+2}$ be the natural surjection given by
\begin{equation}
	p_2 \left( \sum_{k \geq 2} s_k \left( x_k \right) \right) = \eqcl{x_2}_{\qdr^3},
\end{equation}
where we denote by $\eqcl{x_2}_{\qdr^3}$ the image of $x_2$ by the projection
$\ncdf{A}[2][] \rightarrow \ncdf{A}[2][] / \Im \left( \qdr^3 \right)$.

\begin{thm} \label{lm:p_2-homotopy-equivalence}
	The degree $2$ chain map
	$p_2 \colon \totcomp{\mathcal{A}}[2][] \rightharpoonup {{\ncdf{\mathcal{A}}[2][]}/
			\Im \left( \qdr^3 \right)}$ is a homotopy equivalence. A homotopy inverse
	$i_2 \colon {\ncdf{\mathcal{A}}[2][]} / \Im \left( \qdr^3 \right) \rightharpoonup
		\totcomp{\mathcal{A}}[2][]$ for $p_2$ is given by the degree $-2$ chain map
	\begin{equation}
		i_2 \left( \eqcl{x}_{\qdr^3} \right) = \sum_{k \geq 0} (-1)^{\frac{k(k+1)}{2}}
		s_{k+2} \left( \left( h_{\dr} \clie{\mu} \right)^k \left( \left( h_{\dr} \qdr \right) (x) \right) \right).
	\end{equation}
\end{thm}
\begin{proof}
	Let $\left( p_2 \right)^{[-1]}_{[1]} \colon {\totcomp{\mathcal{A}}[2][]}[-1] \rightarrow
		\left({\ncdf{\mathcal{A}}[2][]}/ \Im \left( \qdr^3 \right) \right)[1]$ be defined by the
	diagram
	\begin{equation*}
		\begin{tikzcd}
			{{\totcomp{\mathcal{A}}[2][]}[-1]} &
			{\totcomp{\mathcal{A}}[2][]} &
			{{\ncdf{\mathcal{A}}[2][]}/ \Im \left( \qdr^3 \right)} &
			{\left({\ncdf{\mathcal{A}}[2][]}/ \Im \left( \qdr^3 \right) \right)[1],}
			\arrow[harpoon, from=1-1, to=1-2]
			\arrow["\left( p_2 \right)^{[-1]}_{[1]}",curve={height=-25pt}, from=1-1, to=1-4]
			\arrow["p_2", harpoon, from=1-2, to=1-3]
			\arrow[harpoon, from=1-3, to=1-4]
		\end{tikzcd}
	\end{equation*}
	where the unnamed maps are the natural suspension and desuspension maps.
	Then $\left( p_2 \right)^{[-1]}_{[1]}$ is a degree zero chain map and the mapping cone
	of $\left( p_2 \right)^{[-1]}_{[1]}$
	is isomorphic to the total complex of the bicomplex
	\begin{equation}
		\begin{tikzcd}[column sep=0.58cm]
			\cdots & 0 & {\ncdf{\mathcal{A}}[2][]}/ \Im \left( \qdr^3 \right) & \ncdf{\mathcal{A}}[2][] &
			\ncdf{\mathcal{A}}[3][] & \cdots\phantom{\ncdf{\mathcal{A}}[3][]}
			\arrow["\qdr^4"', from=1-6, to=1-5]
			\arrow["\qdr^3"', from=1-5, to=1-4]
			\arrow[from=1-4, to=1-3]
			\arrow[from=1-3, to=1-2]
			\arrow[from=1-2, to=1-1]
			\arrow["h_{\dr} \circ \qdr^2", curve={height=-12pt}, from=1-3, to=1-4]
			\arrow["h_{\dr}", curve={height=-12pt}, from=1-4, to=1-5]
			\arrow["h_{\dr}", curve={height=-12pt}, from=1-5, to=1-6]
		\end{tikzcd}
		\label{eq:tot-complex-geq-2-augmented}
	\end{equation}
	where ${\ncdf{\mathcal{A}}[2][]}/ \Im \left( \qdr^3 \right)$ is placed at the first column. The bicomplex
	\eqref{eq:tot-complex-geq-2-augmented} has contractible rows and satisfies the
	convergence condition of \cref{lm:total-bicomplex-contractible-R-linear-contraction-1} (this follows
	from the proof of \cref{lm:tot-ncdfr-contractible}) and hence the total complex is contractible. By
	\cref{lm:contractible-cone-homotopy-equivalence}, the map
	$\left( p_2 \right)^{[-1]}_{[1]}$, and hence $p_2$, is a homotopy equivalence.
	The formula for $i_2$ follows as in \cref{eq:total-bicomplex-projection-equivalence-inverse-1} since
	$h_{\dr}^2 = 0$.
\end{proof}

Finally, we give a simple description of the periodicity operator $S$ acting between
the models $\totcomp{\mathcal{A}}[1][]$ and $\totcomp{\mathcal{A}}[2][]$ for cyclic homology.
Let $\pi_{\geq 2} \colon \totcomp{\mathcal{A}}[1][] \rightarrow \totcomp{\mathcal{A}}[2][]$ denote the projection
\begin{equation*}
	\pi_{\geq 2} \left( \sum_{k \geq 1} s_k \left( x_k \right) \right) \defeq
	\sum_{k \geq 2} s_k \left( x_k \right).
\end{equation*}

\begin{lm} \label{lm:periodicity-operator-as-projection}
	For any Banach $\Ainf$-algebra $\mathcal{A}$, the diagram
	\begin{equation} \label{eq:diag-S-p-psi}
		\begin{tikzcd}
			\totcomp{\mathcal{A}}[2][] && \totcomp{\mathcal{A}}[1][] \\
			{{\ncdf{\mathcal{A}}[2][]}/ \Im \left( \qdr^3 \right)} \\
			\ncdfr{\mathcal{A}}[0][] && \ncdfr{\mathcal{A}}[0][]
			\arrow["p_2"', "{[2]}", harpoon, from=1-1, to=2-1]
			\arrow["\psi"', "{[1]}", harpoon, from=2-1, to=3-1]
			\arrow["S", "{[2]}"', harpoon, from=3-3, to=3-1]
			\arrow["i_1"', "{[-1]}", harpoon, from=3-3, to=1-3]
			\arrow["\pi_{\geq 2}"', "{[0]}", from=1-3, to=1-1]
		\end{tikzcd}
	\end{equation}
	strictly commutes.
\end{lm}
\begin{proof}
	During the proof, given $x \in \tensr{A}$, we will denote the cyclic equivalence class of
	$x$ in $\ncdfr{A}[0][] = \tensr{A} / \Im \left( \idd - \t \right)$ by $\cyccl{x}$.
	We have
	\begin{equation*}
		\begin{aligned}
			\left( \overline{\qdr}^2 \circ p_2 \circ \pi_{\geq 2} \circ i_1 \right) \left( \cyccl{x} \right)
			\stackrel{\eqref{eq:def-p_1'}}{=}{}                     &
			\left( \qdr^2 \circ h_{\dr} \circ \clie{\mu} \circ h_{\dr} \right) \left( x \right)
			\\
			\stackrel{\eqref{eq:h-dr-q-contraction}}{=}{}           &
			\left( \left( \idd - h_{\dr} \circ \qdr^1 \right) \circ \clie{\mu} \circ h_{\dr} \right)
			\left( x \right)
			\\
			\stackrel{\eqref{eq:qdr-clie-commute}}{=}{}             &
			\left( \clie{\mu} \circ h_{\dr} - h_{\dr} \circ \clie{\mu} \circ \qdr^1 \circ h_{\dr} \right)
			\left( x \right)
			\\
			\stackrel{\eqref{eq:h-dr-q-contraction}}{=}{}           &
			\left( \clie{\mu} \circ h_{\dr} - h_{\dr} \circ \clie{\mu} \right) \left( x \right)
			\\
			\stackrel{\phantom{\eqref{eq:h-dr-q-contraction}}}{=}{} &
			\left( \clie{\mu} \circ h_{\dr} \right) \left( x \right) -
			\left( h_{\dr} \circ b \right) \left( x \right).
		\end{aligned}
	\end{equation*}
	The explicit formula for $h_{\dr}$ shows that we have $\varphi \circ h_{\dr} = h' \N$ and since
	$\varphi \circ \clie{\mu} = b \circ \varphi$ we get that
	\begin{equation*}
		\begin{aligned}
			\left( \psi \circ p_2 \circ \pi_{\geq 2} \circ i_1 \right) \left( \cyccl{x} \right)
			\stackrel{\eqref{eq:psi-def}}{=}{}           &
			\left( \overline{\beta} \circ \overline{\varphi} \circ \overline{\qdr}^2 \circ p_2
			\circ \pi_{\geq 2} \circ i_1 \right) \left( \cyccl{x} \right)
			\\
			\stackrel{\phantom{\eqref{eq:psi-def}}}{=}{} &
			\cyccl{ \left( \beta \left( b h' \N - h' \N b \right) \right) \left( x \right)}.
		\end{aligned}
	\end{equation*}
	We have
	\begin{equation*}
		\begin{aligned}
			\beta \left( b h' \N - h' \N b \right)
			\stackrel{\phantom{\eqref{eq:partial-1-t}}}{=}{} &
			\beta \partial \left( h' \N \right)
			\stackrel{\eqref{eq:h'N-homotopy}}{=}
			\beta \partial \left( \idd - \left( \idd - \t \right)h \right)
			\stackrel{\eqref{eq:partial-1-t}}{=}
			- \beta \left( \idd - \t \right) \partial \left( h \right)
			\\
			\stackrel{\eqref{eq:T-def}}{=}{}                 &
			\beta \left( \idd - \t \right) T
			\stackrel{\eqref{eq:def-beta}}{=} T' h \left( \idd - \t \right) T
			\stackrel{\eqref{eq:Nh'-homotopy}}{=}
			T' \left( \idd - \N h' \right) T
			\\
			\stackrel{\eqref{eq:1-t-T}}{=}{}                 &
			T' T - \left( \idd - \t \right) T h' T
			\stackrel{\eqref{eq:S-T'-T}}{=}
			S - \left( \idd - \t \right) T h' T
		\end{aligned}
	\end{equation*}
	and hence
	\begin{equation*}
		\left( \psi \circ p_2 \circ \pi_{\geq 2} \circ i_1 \right) \left( \cyccl{x} \right) =
		\cyccl{Sx}
	\end{equation*}
	as required.
\end{proof}

When $\mathcal{A}$ is unital, the vertical maps in diagram \eqref{eq:diag-S-p-psi}
are homotopy equivalences, and thus all the complexes in the diagram compute
the cyclic homology up to a shift.

\subsection{Extended and Reduced Cyclic Complexes} \label{sec:extended-reduced-cyclic-complexes}

Let $\mathcal{R} = (R,d)$ be a differential graded-commutative Banach $\mathbbm{k}$-algebra and let
$\mathcal{A} = \left( A, \mu \right)$ be a Banach $\Ainf$-algebra over $\mathcal{R}$.
In the previous sections, we worked with Connes' cyclic complex
\begin{equation*}
	\cconnes{\mathcal{A}}[] = \left( \tensrcyc{A}, \cycl{\mu} \right) = \left( \ncdfr{A}[0][], \clie{\mu} \right)
\end{equation*}
whose underlying $R$-module is the \textit{reduced} cyclic tensor module $\tensrcyc{A} = \tensr{A} / \Im \left( \idd - \t \right)$.
In \cref{sec:extension-cycl-full-tensor-module}, we showed that it also makes sense to work with the ``full''
cyclic tensor module $\ncdf{A}[0][] = \tenscyc{A} = \tens{A} / \Im \left( \idd - \t \right)$ in which we have
the extra element $1 \in R$ which gets mapped by $\cycl{\mu} = \clie{\mu}$ to the closed
curvature element $\mu_0 \left( 1 \right)$.

\begin{dfn}
	The complex
	\begin{equation*}
		\ncdf{\mathcal{A}}[0][] =
		\left( \ncdf{A}[0][], \clie{\mu} \right) = \left( \tenscyc{A}, \cycl{\mu} \right)
	\end{equation*}
	is called the \textbf{extended Connes complex} of $\mathcal{A}$ and its cohomology
	\begin{equation*}
		\hcyce{\mathcal{A}} \defeq \cohom{{\ncdf{\mathcal{A}}[0][]}}
	\end{equation*}
	is a graded $\cohom{\mathcal{R}}$-module called the \textbf{extended cyclic homology} of
	$\mathcal{A}$.
\end{dfn}

\begin{rem} \label{rem:extended-to-standard-les}
	The extended Connes complex of $\mathcal{A}$ is the cone
	of the degree zero morphism $f \colon \mathcal{R}[-1] \rightarrow \ncdfr{\mathcal{A}}[0][]$ of
	differential graded $\mathcal{R}$-modules given by
	$f \left( s_{-1} \left( r \right) \right) = (-1)^{\degb{r}} r \cdot \mu_0(1)$. Thus, the extended cyclic
	homology and the standard cyclic homology are related to each other via the long exact sequence
	\begin{equation*}
		\dots \xrightarrow{\delta} \cohom{\mathcal{R}}[*-1] \xrightharpoonup{f} \hcyc{\mathcal{A}} \xrightarrow{i} \hcyce{\mathcal{A}} \xrightarrow{\delta}
		\cohom{\mathcal{R}} \xrightharpoonup{f} \dots
	\end{equation*}
	When $\mu_0 \left( 1 \right) = 0$, or, more generally, when $f$ is null-homotopic, i.e.,
	$\mu_0 \left( 1 \right)$ is exact in $\ncdfr{\mathcal{A}}[0][]$,
	we have $\hcyce{\mathcal{A}} \cong \cohom{\mathcal{R}} \oplus \hcyc{\mathcal{A}}$.
\end{rem}

In the next sections, we will also need reduced versions of the cyclic complexes, defined in the case $\mathcal{A}$ is unital.
Let $\mathcal{A} = \left( A, \mu, e \right)$ be a unital $\Ainf$-algebra over $\mathcal{R}$.
We have the following expression for the action of $\clie{\mu}$ on elements of $\ndfr{A}[][]$ which start with the unit $e$:

\begin{lm} \label{lm:clie-mu-on-elem-starts-with-e}
	Let $y \in \ndf{A}[][]$. Then
	\begin{equation} \label{eq:clie-mu-on-elem-starts-with-e}
		\clie{\mu} \left( e \otimes y \right) = -e \otimes \lie{\mu} \left( y \right) + \left( \idd - \t \right) \left( y \right).
	\end{equation}
\end{lm}
\begin{proof}
	Recall that we have $\clie{\mu} = \cycl{ \left( \lie{\mu} \right) }$. Then by \cref{def:cyclization-coder-short}, we have
	\begin{align}
		\clie{\mu} \left( e \otimes y \right)
		\stackrel{\phantom{\eqref{eq:parity-extends-koszul}}}{=}{} &
		(-1)^{\braid{(0,-1)}{(0,1)}} e \otimes \lie{\mu} \left( y \right)
		\notag
		\\
		                                                           & +
		                                                             (-1)^{\braid{\degb{y_{(3)}}}{(0,-1) + \degb{y_{(1)}} + \degb{y_{(2)}}}}
		\corest{\left( \lie{\mu} \right)} \left( y_{(3)} \otimes e \otimes y_{(1)} \right) \otimes y_{(2)}
		\label{eq:clie-mu-e-y-open-def}
		\\
		\stackrel{\eqref{eq:parity-extends-koszul}}{=}{}           &
		- e \otimes \lie{\mu} \left( y \right) +
		                      (-1)^{\braid{\degb{y_{(3)}}}{(0,-1) + \degb{y_{(1)}} + \degb{y_{(2)}}}}
		\corest{\left( \lie{\mu} \right)} \left( y_{(3)} \otimes e \otimes y_{(1)} \right) \otimes y_{(2)}.
		\notag
	\end{align}

	Let $y_1, \dots, y_k \in A \oplus \ul{A}$. The explicit formula \eqref{eq:lie-derivation-corestriction} for
	$\corest{\left( \lie{\mu} \right)}$ and the unit property show that
	\begin{equation} \label{eq:corest-lie-mu-with-e-k-neq-2}
		\left( \lie{\mu} \right)_k \left( y_1, \dots, y_k \right) = 0
	\end{equation}
	for $k \neq 2$ whenever $y_i = e$ for some $1 \leq i \leq k$. For the case $k = 2$, we have
	\begin{equation} \label{eq:corest-lie-mu-with-unit}
		\begin{aligned}
			\left( \lie{\mu} \right)_2 \left( e \otimes a \right)      & = \mu_2 \left( e, a \right) = a,
			\\
			\left( \lie{\mu} \right)_2 \left( e \otimes \ul{a} \right) & = (-1)^{\braid{(-1,0)}{(0,1-1)}} \ul{\mu_2 \left( e, a \right)} = \ul{a},
			\\
			\left( \lie{\mu} \right)_2 \left( a \otimes e \right)      & = \mu_2 \left( a, e \right) = (-1)^{\degb{a} + 1} a,
			\\
			\left( \lie{\mu} \right)_2 \left( \ul{a} \otimes e \right) & =
			(-1)^{\braid{(-1,0)}{(0,1)}} \ul{\mu_2 \left( a, e \right)} = (-1)^{\braid{(-1,0)}{(0,1)} + \degb{a} + 1} \ul{a}
		\end{aligned}
	\end{equation}
	for $a \in A$. Since we assume condition \eqref{eq:parity-extends-koszul} on the parity form, \cref{eq:corest-lie-mu-with-unit}
	can be rewritten uniformly as
	\begin{equation} \label{eq:corest-lie-mu-with-e-k-eq-2}
		\begin{aligned}
			\left( \lie{\mu} \right)_{2} \left( e, y \right) & = y,
			\\
			\left( \lie{\mu} \right)_{2} \left( y, e \right) & = (-1)^{\braid{\degb{y}}{(0,1)} + 1} y
		\end{aligned}
	\end{equation}
	whenever $y \in A \oplus \ul{A}$.\footnote{\Cref{eq:corest-lie-mu-with-e-k-neq-2,eq:corest-lie-mu-with-e-k-eq-2} state
		that the element $e$ acts as ``unit'' for the degree $(0,1)$ coderivation $\lie{\mu}$ in an appropriate sense.}

	When $y = 1$, \cref{eq:clie-mu-on-elem-starts-with-e} clearly holds. Assume that $y$ has the form
	$y = y_1 \otimes \dots \otimes y_r$ for $r \geq 1$, where $y_i \in A \oplus \ul{A}$ for $1 \leq i \leq r$.
	Then
	\begin{equation*}
		\begin{aligned}
			\clie{\mu} \left( e \otimes y \right)
			\stackrel[\eqref{eq:corest-lie-mu-with-e-k-neq-2}]{\eqref{eq:clie-mu-e-y-open-def}}{=}{} &
			- e \otimes \lie{\mu} \left( y \right) +
			\left( \lie{\mu} \right)_2 \left( e, y_1 \right) \otimes y_2 \otimes \dots \otimes y_r
			\\
			                                                                                         & +
			                                                                                           (-1)^{\braid{\degb{y_r}}{(0,-1) + \degb{y_1} + \dots + \degb{y_{r-1}}}}
			\left( \lie{\mu} \right)_2 \left( y_r, e \right) \otimes y_1 \otimes \dots \otimes y_{r-1}
			\\
			\stackrel{\eqref{eq:corest-lie-mu-with-e-k-eq-2}}{=}{}                                   &
			- e \otimes \lie{\mu} \left( y \right) + \left( \idd - \t \right) \left( y \right),
		\end{aligned}
	\end{equation*}
	as required.
\end{proof}

Consider the elements of $\ndfr{A}[][] = \tensr{A \oplus \ul{A}}[R]$ of the form
\begin{equation} \label{eq:elements-ndfr-free-unit}
	a_0^1 \otimes \dots \otimes a_0^{r_0} \otimes \ul{a_1} \otimes a_1^1 \otimes \dots \otimes a_1^{r_1} \otimes \ul{a_2} \otimes
	\dots \otimes \ul{a_k} \otimes a_k^1 \otimes \dots \otimes a_k^{r_k}
\end{equation}
where $k \geq 0$ and $r_0, \dots, r_k \geq 0$ and $a_i, a_i^j \in A$, such that $a_i^j = e$ for some $0 \leq i \leq k$ and
$1 \leq j \leq r_i$.
Those are the elementary codifferential forms on $A$ in which the unit $e$ of $A$ appears at least once without an underline.
We will denote elements of the form \eqref{eq:elements-ndfr-free-unit} succinctly by $*e*$.
Let $\left< *e* \right> \subseteq \ndfr{A}[][]$ be the bigraded Banach $R$-submodule of
$\ndfr{A}[][] = \tensr{A \oplus \ul{A}}[R]$ generated by elements of the form \eqref{eq:elements-ndfr-free-unit}.
Similarly, denote by $\left< *e*, \Im \left( \idd - \t \right) \right>$ the bigraded Banach $R$-submodule of
$\ndfr{A}[][]$ generated by elements of the form \eqref{eq:elements-ndfr-free-unit} and the image of $\idd - \t$.

\begin{cor} \label{cor:star-e-star-1-t-invariant}
	The graded $R$-module $\left< *e*, \Im \left( \idd - \t \right) \right>$ is invariant under $\clie{\mu}$.
	Hence, $\left< *e*, \Im \left( \idd - \t \right) \right>$ is a subcomplex of $\ndfr{\mathcal{A}}[][]$.
\end{cor}
\begin{proof}
	Denote by $e*$ elements of the form \eqref{eq:elements-ndfr-free-unit} in which $a_0^1 = e$, i.e., start with the unit
	$e$. Note that we have $\left< *e*, \Im \left( \idd - \t \right) \right> = \left< e*, \Im \left( \idd - \t \right) \right>$
	as any element of the form $*e*$ can be rotated using $\t$ to begin with $e$.
	\Cref{lm:clie-mu-on-elem-starts-with-e} shows that $\clie{\mu} \left( \left< e* \right> \right) \subseteq \left< *e*, \Im \left( \idd - \t \right) \right>$
	and the identity $\clie{\mu} \circ \left( \idd - \t \right) = \left( \idd - \t \right) \circ \lie{\mu}$
	(see \cref{lm:coder-descends-quotient})
	shows that $\Im \left( \idd - \t \right)$ is invariant under $\clie{\mu}$.
\end{proof}

Let $\degen[A][][]$ be the closure of the image of $\left< *e* \right>$ under the projection
$\ndfr{A}[][] \rightarrow \ncdfr{A}[][]$.
Given $k \geq 0$, we denote by $\degen[A][k][]$ the line degree $k$ component of $\degen[A][][]$.
The graded Banach $R$-submodule $\degen[A][k][]$ is generated by the \textit{cyclic} equivalence classes of elements
of the form \eqref{eq:elements-ndfr-free-unit} with fixed $k$.
By \cref{cor:star-e-star-1-t-invariant}, $\degen[A][k][]$ is invariant under $\clie{\mu}$ for each $k \geq 0$.
\begin{dfn}
	The quotient complex
	\begin{equation*}
		\ncdfrred{\mathcal{A}}[0][] \defeq \ncdfr{\mathcal{A}}[0][] / \degen[A][0][] \cong
		\left( \ndfr{A}[0][] / \left< *e*, \Im \left( \idd - \t \right) \right>, \clie{\mu} \right)
	\end{equation*}
	is called the \textbf{reduced Connes complex} of $\mathcal{A}$ and its cohomology
	\begin{equation*}
		\hcycred{\mathcal{A}} \defeq \cohom{{\ncdfrred{\mathcal{A}}[0][]}}
	\end{equation*}
	is a graded $\cohom{\mathcal{R}}$-module called the \textbf{reduced cyclic homology} of $\mathcal{A}$.
\end{dfn}

\begin{dfn}
	The quotient complex
	\begin{equation*}
		\ncdfred{\mathcal{A}}[0][] \defeq \ncdf{\mathcal{A}}[0][] / \degen[A][0][] \cong
		\left( \ndf{A}[0][] / \left< *e*, \Im \left( \idd - \t \right) \right>, \clie{\mu} \right)
	\end{equation*}
	is called the \textbf{extended reduced Connes complex} of $\mathcal{A}$ and its cohomology
	\begin{equation*}
		\hcycered{\mathcal{A}} \defeq \cohom{{\ncdfred{\mathcal{A}}[0][]}}
	\end{equation*}
	is a graded $\cohom{\mathcal{R}}$-module called the \textbf{extended reduced cyclic homology} of $\mathcal{A}$.
\end{dfn}

Here, $\left< *e*, \Im \left( \idd - \t \right) \right>$ denotes the graded \textit{Banach} $R$-submodule
generated by the Banach $R$-submodules $\left< *e* \right>$ and $\Im \left( \idd - \t \right)$.

\begin{rem}
	When $\mathcal{A}$ corresponds to an associative algebra, i.e., $\mu_k = 0$ for $k \neq 2$, the
	complex $\ncdfrred{\mathcal{A}}[0][]$ coincides up to an identification with
	the reduced Connes complex \cite[Section 2.2.13]{Loday1998} appearing in the literature.
\end{rem}

\begin{rem} \label{rem:cindmap-maps-degen-to-degen}
	Let $\mathcal{B}$ be a Banach $\Ainf$-algebra over a differential graded-commutative Banach $\mathbbm{k}$-algebra $\mathcal{S}$.
	Given a unital morphism $f \colon \mathcal{A} \rightarrow \mathcal{B}$ of Banach $\Ainf$-algebras, and $k \geq 0$,
	the explicit formulas for the induced morphism $\cindmap{f} \colon \ncdf{\mathcal{A}}[k][] \rightarrow \ncdf{\mathcal{B}}[k][]$
	(\cref{def:cyclization-morphism,eq:cindmap-f-explicit-action}) show that $\cindmap{f}$ maps
	$\degen[A][k][]$ to $\degen[B][k][]$. In particular, the constructions
	$\mathcal{A} \mapsto \ncdfrred{\mathcal{A}}[0][], \mathcal{A} \mapsto \ncdfred{\mathcal{A}}[0][]$
	enjoy the same functorial properties
	as the non-reduced versions $\mathcal{A} \mapsto \ncdfr{\mathcal{A}}[0][], \mathcal{A} \mapsto \ncdf{\mathcal{A}}[0][]$.
\end{rem}

\section{The \texorpdfstring{$\infty$}{Infinity}-modulus Function and its Invariance}
\label{sec:generalized-trace}

In this section, we define the $\infty$-modulus function and show that the $\infty$-modulus
of a bounding cochain gives a gauge-invariant cohomology class.
We begin in \cref{subsec:canonical-elements-ncdf-0} by considering
the cyclic exponential $\G{b}[0]$ of a topologically nilpotent element $b$.
The construction $b \mapsto \G{b}[0]$ is natural and gives us a canonical chain
in the extended Connes complex. When $b$ is a strong (resp.\ weak) bounding cochain,
we show that $\G{b}[0]$ is a cycle and thus
defines a homology class in the extended
(resp.\ extended reduced)
cyclic homology. We also show that the cyclic homology class is invariant under
$\mathfrak{A}$-gauge-equivalence, when $\mathfrak{A}$ is a strong pseudoisotopy.

\Cref{sec:pre-infinity-traces} introduces the notion of a pre-$\infty$-trace
on a Banach $\Ainf$-algebra as a closed linear functional on Connes' cyclic complex.
We explicitly spell out the relations a pre-$\infty$-trace must satisfy in terms
of its components, and also discuss a unital version.

In \cref{sec:traces-modulus-invariance}, we extend the concept of a pre-$\infty$-trace
by defining an $\infty$-trace as a closed linear functional on the \textit{extended} Connes complex,
and discuss its properties.
The $\infty$-modulus of a topologically nilpotent element $b$ is defined by evaluating the
$\infty$-trace on the cyclic exponential $\G{b}[0]$ associated to $b$.
In \cref{thm:invariance-extended-infty-modulus}, we prove that when $b$ is a bounding cochain,
the $\infty$-modulus of $b$ gives a cohomology class which is invariant under
$\mathfrak{A}$-gauge-equivalence.
The theorem is proven for both strong bounding cochains in the non-unital version,
and weak bounding cochains in the unital version.

Let $\mathbbm{k}$ be a field of characteristic zero, endowed with the trivial norm, and fixed for the duration
of the section. Let $\mathcal{R} = \left( R, d \right)$ be a differential graded-commutative Banach $\mathbbm{k}$-algebra and
let $\mathcal{A} = \left( A, \mu \right)$ be a Banach $\Ainf$-algebra over $\mathcal{R}$.

\subsection{Canonical Cycle in Connes' Cyclic Complex} \label{subsec:canonical-elements-ncdf-0}
Given a topologically nilpotent element $b \in \tc{A}$, consider the
chains
\begin{align}
	\G{b}[0] = \G{b}[0][\mathcal{A}]   & \defeq \cexp{b} = 1 + \sum_{n = 1}^{\infty} \frac{b^{\otimes n}}{n}
	\in \ncdf{\mathcal{A}}[0][0],
	\label{eq:G_0-def}
	\\
	\Ho{b}[0] = \Ho{b}[0][\mathcal{A}] & \defeq \corest{\mu} \left( \Exp{b} \right) \otimes \Exp{b}
	\in \ncdfr{\mathcal{A}}[0][1] \subset \ncdf{\mathcal{A}}[0][1].
	\label{eq:Ho_0-def}
\end{align}
The chain $\G{b}[0]$ is the cyclic exponential of $b$ (see \cref{dfn:cyclic-exponential}).
The cyclic exponential $\G{b}[0]$ includes the term $1 \in R$ and hence belongs to the
extended Connes complex $\ncdf{A}[0][] = \tenscyc{A}[]$,
while the chain $\Ho{b}[0]$ lives in the standard Connes
complex $\ncdfr{A}[0][] = \tensrcyc{A}[]$.
In \cref{sec:extension-cycl-full-tensor-module} we showed how to extend the definition
of $\cycl{\mu} = \clie{\mu}$ from $\ncdfr{A}[0][]$ to $\ncdf{A}[0][]$ and how, given
a morphism, to extend the action of $\cindmap{f} = \cycl{f}$ from $\ncdfr{A}[0][]$ to $\ncdf{A}[0][]$.
With the extensions in place, we have the following properties of $\G{b}[0]$ and
$\Ho{b}[0]$:

\begin{lm} \label{lm:G-0-properties}
	Let $\mathcal{A} = \left( A, \mu \right)$ be a Banach $\Ainf$-algebra over a differential
	graded-commutative Banach $\mathbbm{k}$-algebra $\mathcal{R} = \left( R, d_R \right)$
	and let $\mathcal{B} = \left( B, \nu \right)$ be a Banach $\Ainf$-algebra over a
	differential graded-commutative Banach $\mathbbm{k}$-algebra $\mathcal{S} = \left( S, d_S \right)$.
	Let $f \colon \mathcal{A} \rightarrow \mathcal{B}$ be a Banach $\Ainf$-morphism.
	Given a topologically nilpotent element $b \in \tc{A}$, we have the following identities
	between elements of $\ncdf{\cdot}[0][]$:
	\begin{enumerate}
		\item{(Naturality)}
		      \begin{equation} \label{eq:cycl-f-G0}
			      \cycl{f} \left( \G{b}[0][\mathcal{A}] \right) = \G{\mcfunc{f} \left( b \right)}[0][\mathcal{B}].
		      \end{equation}
		\item{(Differential)}
		      \begin{equation} \label{eq:clie-mu-G0}
			      \clie{\mu} \left( \G{b}[0][\mathcal{A}] \right) = \Ho{b}[0][\mathcal{A}].
		      \end{equation}
	\end{enumerate}
\end{lm}
\begin{proof}
	\Cref{eq:cycl-f-G0} was proven in \cref{lem:func-mc-cyclic}, while
	\cref{eq:clie-mu-G0} was proven in \cref{lm:coder-cycl-exp-identity}.
\end{proof}

When $b \in \tc{A}$ is a strong (resp.\ weak) bounding cochain, \cref{lm:G-0-properties} implies
that $\G{b}[0]$ is a cycle in the extended (resp.\ extended reduced) Connes complex.
Thus, we have a way of associating a canonical cycle and homology class to
any strong (resp.\ weak) bounding cochain:

\begin{lm}[Canonical Cycle and Homology Class for a Bounding Cochain] \hfill
	\label{lm:G_0-bounding-chain-closed}
	\begin{enumerate}
		\item Let $\mathcal{A} = \left( A, \mu \right)$ be a Banach $\Ainf$-algebra over
		      $\mathcal{R}$ and let $b \in \mc{\mathcal{A}}$ be a strong bounding cochain.
		      Then $\G{b}[0]$ is a cycle in $\ncdf{\mathcal{A}}[0][0]$ which defines
		      a homology class $\eqcl{ \G{b}[0] } \in \hcyce{\mathcal{A}}[0]$ in
		      the extended cyclic homology of $\mathcal{A}$.
		\item Let $\mathcal{A} = \left( A, \mu, e \right)$ be a unital Banach $\Ainf$-algebra over
		      $\mathcal{R}$ and let $b \in \mc{\mathcal{A}, c}$ be a weak bounding cochain.
		      Then $\G{b}[0]$ is a cycle in $\ncdfred{\mathcal{A}}[0][0]$
		      which defines a homology class $\eqcl{ \G{b}[0] } \in  \hcycered{\mathcal{A}}[0]$
		      in the extended reduced cyclic homology of $\mathcal{A}$.
	\end{enumerate}
\end{lm}
\begin{proof}
	By \cref{lm:G-0-properties}, we have
	\begin{equation*}
		\clie{\mu} \left( \G{b}[0] \right) = \Ho{b}[0] = \corest{\mu} \left( \Exp{b} \right) \otimes \Exp{b}.
	\end{equation*}
	When $b \in \mc{\mathcal{A}}$ is a strong bounding cochain,
	we have $\corest{\mu} \left( \Exp{b} \right) = 0$ and hence $\G{b}[0]$ is a cycle.
	When $\mathcal{A}$ is unital and $b \in \mc{\mathcal{A}, c}$ is a weak bounding cochain,
	then $\corest{\mu} \left( \Exp{b} \right) = c \cdot e$,
	so
	\begin{equation*}
		\clie{\mu} \left( \G{b}[0] \right) =
		c \cdot e \otimes \Exp{b} = c \cdot \left( \sum_{n = 0}^{\infty} e \otimes b^{\otimes n} \right).
	\end{equation*}
	Since the differential of $\G{b}[0]$ contains the unit $e$, it belongs to
	$\degen[A][0][]$ and hence $\G{b}[0]$ is a cycle in the \textit{quotient}
	$\ncdfred{\mathcal{A}}[0][] = \ncdf{\mathcal{A}}[0][] / \degen[A][0][]$, which is the extended
	reduced Connes complex of $\mathcal{A}$.
\end{proof}

Let $\mathcal{A}_0$ and $\mathcal{A}_1$ be two Banach
$\Ainf$-algebras over the same differential graded-commutative $\mathbbm{k}$-algebra
$\mathcal{S}$.
Applying \cref{dfn:F-strong-pseudoisotopy} to the extended cyclic homology functor
$\mathcal{F} = \hcyce{}[]$, we obtain the notion of an $\hcyce{}[]$-\textbf{strong}
pseudoisotopy $\mathfrak{A}$ between $\mathcal{A}_0$ and $\mathcal{A}_1$.
Such a pseudoisotopy yields a canonical isomorphism
\begin{equation*}
	\mathfrak{a} \colon  \hcyce{\mathcal{A}_0 / \mathcal{S}}[] \overset{\sim}{\rightarrow} \hcyce{\mathcal{A}_1 / \mathcal{S}}[]
\end{equation*}
of graded $\cohom{\mathcal{S}}[]$-modules.

Similarly, when $\mathcal{A}_0,\mathcal{A}_1$ and $\mathfrak{A}$ are unital, applying
\cref{dfn:F-strong-pseudoisotopy} to the extended reduced cyclic homology functor
$\mathcal{F} = \hcycered{}[]$ yields the notion of an $\hcycered{}[]$-\textbf{strong}
pseudoisotopy $\mathfrak{A}$ between $\mathcal{A}_0$ and $\mathcal{A}_1$,
which induces a canonical isomorphism
\begin{equation*}
	\mathfrak{a} \colon \hcycered{\mathcal{A}_0 / \mathcal{S}}[] \overset{\sim}{\rightarrow} \hcycered{\mathcal{A}_1 / \mathcal{S}}[].
\end{equation*}
Then we have:
\begin{thm}[Invariance of the Canonical Cycle under Strong Pseudoisotopy] \hfill
	\label{lm:invariance-Gb0-pseudoisotopy}
	Let $\mathcal{S}$ be a differential graded-commutative Banach $\mathbbm{k}$-algebra, and let
	$\mathcal{A}_0$ and $\mathcal{A}_1$ be two non-unital (resp.\ unital) Banach $\Ainf$-algebras over
	$\mathcal{S}$. Let $\mathfrak{A}$ be a $\hcyce{}[]$-strong
	(resp.\ $\hcycered{}[]$-strong) pseudoisotopy between $\mathcal{A}_0$
	and $\mathcal{A}_1$.
	Let $b_0 \in \tc{A_0}$ and $b_1 \in \tc{A_1}$ be two strong (resp.\ weak) $\mathfrak{A}$-gauge-equivalent
	bounding cochains. Then
	\begin{equation*}
		\mathfrak{a}(\eqcl{ \G{b_0}[0][\mathcal{A}_0] }) = \eqcl{ \G{b_1}[0][\mathcal{A}_1] }
		\textrm{ in } \hcyce{\mathcal{A}_1}[0]
		\qquad \left(
		\textrm{resp.\ in } \hcycered{\mathcal{A}_1}[0]
		\right).
	\end{equation*}
\end{thm}
\begin{proof}
	Assume we are in the non-unital case.
	Let $b_0 \in \mc{\mathcal{A}_0}$ and $b_1 \in \mc{\mathcal{A}_1}$ be two $\mathfrak{A}$-gauge-equivalent bounding cochains.
	Choose $b \in \mc{\mathfrak{A}}$ with $\mcfunc{\evalmf^i} \left( b \right) = b_i$ for $i = 0, 1$.
	Then, by \cref{lm:G_0-bounding-chain-closed}, $\G{b}[0][\mathfrak{A}]$ is a cycle defining a homology class
	$\eqcl{ \G{b}[0][\mathfrak{A}] }$ in $\hcyce{\mathfrak{A}}[0]$, and
	\begin{equation*}
		\begin{aligned}
			\hcyce{\evalmf^i}[] \left( \eqcl{ \G{b}[0][\mathfrak{A}] } \right)
			\stackrel{\phantom{\eqref{eq:cycl-f-G0}}}{=}{} &
			\cohom{\cycl{{\evalmf}}^i}[] \left( \eqcl{ \G{b}[0][\mathfrak{A}] } \right) =
			\eqcl{ \cycl{{\evalmf}}^i \left( \G{b}[0][\mathfrak{A}] \right) }
			\\
			\stackrel{\eqref{eq:cycl-f-G0}}{=}{}           &
			\eqcl{ \G{\mcfunc{\evalmf^i} \left( b \right)}[0][\mathcal{A}_i] } =
			\eqcl{ \G{b_i}[0][\mathcal{A}_i] }
		\end{aligned}
	\end{equation*}
	for $i = 0, 1$. Hence, $\mathfrak{a}(\eqcl{ \G{b_0}[0][\mathcal{A}_0] }) = \eqcl{ \G{b_1}[0][\mathcal{A}_1] }$.
	The proof for the unital case proceeds the same way, using the fact
	that $\G{b}[0][\mathfrak{A}]$ is a cycle in the extended reduced cyclic complex.
\end{proof}

\begin{rem}
	In \cite{Fukaya2009}, the authors construct several different models for ``the'' mapping cylinder of
	a filtered $\Ainf$-algebra $\mathcal{A}$ and use them to define the notion of homotopy between
	$\Ainf$-morphisms (see also \cite{Hicks2019}). Each such model gives a specific pseudoisotopy $\mathfrak{A}$
	between $\mathcal{A}$ and itself, in which we expect that the induced maps
	$\hcyce{\evalmf^i}[] \colon \hcyce{\mathfrak{A}}[] \overset{\sim}{\rightarrow} \hcyce{\mathcal{A}}[]$
	are not only isomorphisms, but also coincide for $i = 0, 1$. In this case, the isomorphism
	$\mathfrak{a}$ is the identity, and two $\mathfrak{A}$-gauge-equivalent bounding
	cochains $b_0,b_1 \in \tc{A}$ define the same canonical homology class
	$\eqcl{ \G{b_0}[0][\mathcal{A}] } = \eqcl{ \G{b_1}[0][\mathcal{A}] } \in \hcyce{\mathcal{A}}[0]$.
\end{rem}

We finish this section by considering a variant of the cyclic exponential $\G{b}[0]$ defined
in the standard Connes complex $\ncdfr{\mathcal{A}}[0][]$, and study its properties.
Consider the truncated cyclic exponential
\begin{equation} \label{eq:Go_0-def}
	\Go{b}[0] = \Go{b}[0][\mathcal{A}] \defeq \sum_{n = 0}^{\infty} \frac{b^{\otimes \left( n + 1 \right)}}{n + 1} =
	\G{b}[0] - 1 \in \ncdfr{\mathcal{A}}[0][0],
\end{equation}
which doesn't involve the term $1$ and lives in $\ncdfr{\mathcal{A}}[0][]$.
In the standard Connes complex, we have a ``relative'' version of \cref{lm:G-0-properties},
describing the properties of $\Go{b}[0][\mathcal{A}]$:
\begin{lm} \label{lm:Go_0-and-Ho_0-properties}
	Let $\mathcal{A} = \left( A, \mu \right)$ be a Banach $\Ainf$-algebra over a differential
	graded-commutative Banach $\mathbbm{k}$-algebra $\mathcal{R} = \left( R, d_R \right)$
	and let $\mathcal{B} = \left( B, \nu \right)$ be a Banach $\Ainf$-algebra over a
	differential graded-commutative Banach $\mathbbm{k}$-algebra $\mathcal{S} = \left( S, d_S \right)$.
	Let $f \colon \mathcal{A} \rightarrow \mathcal{B}$ be a Banach $\Ainf$-morphism.
	Given a topologically nilpotent element $b \in \tc{A}$, we have the following identities between
	elements of $\ncdfr{\cdot}[0][]$:
	\begin{align}
		\clie{\mu} \left( \Ho{b}[0][\mathcal{A}] \right) & = 0,
		\label{eq:clie-mu-Ho-0}
		\\
		\clie{\mu} \left( \Go{b}[0][\mathcal{A}] \right) & = \Ho{b}[0][\mathcal{A}] - \Ho{0}[0][\mathcal{A}],
		\label{eq:clie-mu-Go-0}
		\\
		\cycl{f} \left( \Ho{b}[0][\mathcal{A}] \right)   & = \Ho{\mcfunc{f} \left( b \right)}[0][\mathcal{B}],
		\label{eq:cycl-f-Ho-0}
		\\
		\cycl{f} \left( \Go{b}[0][\mathcal{A}] \right)   & = \Go{\mcfunc{f} \left( b \right)}[0][\mathcal{B}] -
		\Go{\mcfunc{f} \left( 0 \right)}[0][\mathcal{B}].
		\label{eq:cycl-f-Go-0}
	\end{align}
\end{lm}
\begin{proof}
	We shall perform the calculations in $\ncdf{A}[0][]$, but final results will give us identities
	in $\ncdfr{A}[0][]$. We have
	\begin{equation*}
		\begin{aligned}
			\clie{\mu} \left( \Go{b}[0][\mathcal{A}] \right)
			\stackrel{\eqref{eq:Go_0-def}}{=}{}                          &
			\clie{\mu} \left( \G{b}[0] - 1 \right)
			=
			\clie{\mu} \left( \G{b}[0] \right) - \clie{\mu} \left( 1 \right)
			\\
			\stackrel[\eqref{eq:lie-clie-1}]{\eqref{eq:clie-mu-G0}}{=}{} &
			\Ho{b}[0][\mathcal{A}] - \mu_0 \left( 1 \right)
			\\
			\stackrel{\eqref{eq:Ho_0-def}}{=}{}                          &
			\Ho{b}[0][\mathcal{A}] - \Ho{0}[0][\mathcal{A}]
		\end{aligned}
	\end{equation*}
	which shows \cref{eq:clie-mu-Go-0}.
	Moreover,
	\begin{equation*}
		0 = \clie{\mu}^2 \left( \Go{b}[0][\mathcal{A}] \right) =
		\clie{\mu} \left( \Ho{b}[0][\mathcal{A}] \right) - \clie{\mu} \left( \mu_0 \left( 1 \right) \right) =
		\clie{\mu} \left( \Ho{b}[0][\mathcal{A}] \right) - \clie{\mu}^2 \left( 1 \right) =
		\clie{\mu} \left( \Ho{b}[0][\mathcal{A}] \right)
	\end{equation*}
	which shows \cref{eq:clie-mu-Ho-0}.
	Similarly,
	\begin{equation*}
		\begin{aligned}
			\cycl{f} \left( \Go{b}[0][\mathcal{A}] \right)
			\stackrel{\eqref{eq:Go_0-def}}{=}{}                        &
			\cycl{f} \left( \G{b}[0] - 1 \right) = \cycl{f} \left( \G{b}[0] \right) - \cycl{f} \left( 1 \right)
			\\
			\stackrel[\eqref{def:cycl-f-1}]{\eqref{eq:cycl-f-G0}}{=}{} &
			\G{ \mcfunc{f} \left( b \right) }[0] - \G{ \mcfunc{f} \left( 0 \right) }[0]
			\\
			\stackrel{\eqref{eq:Go_0-def}}{=}{}                        &
			\Go{\mcfunc{f} \left( b \right)}[0][\mathcal{B}] -
			\Go{\mcfunc{f} \left( 0 \right)}[0][\mathcal{B}]
		\end{aligned}
	\end{equation*}
	which shows \cref{eq:cycl-f-Go-0}. Moreover,
	\begin{equation*}
		\begin{aligned}
			\cycl{f} \left( \Ho{b}[0][\mathcal{A}] \right)
			\stackrel{\eqref{eq:clie-mu-G0}}{=}{}                                &
			\cycl{f} \left( \clie{\mu} \left( \G{b}[0][\mathcal{A}] \right) \right)
			\\
			\stackrel{\eqref{cor:func-cycl-morphism-coderivation-extended}}{=}{} &
			\clie{\nu} \left( \cycl{f} \left( \G{b}[0][\mathcal{A}] \right) \right)
			\\
			\stackrel{\eqref{eq:cycl-f-G0}}{=}{}                                 &
			\clie{\nu} \left( \G{ \mcfunc{f} \left( b \right) }[0][\mathcal{B}] \right)
			\\
			\stackrel{\eqref{eq:clie-mu-G0}}{=}{}                                &
			\Ho{ \mcfunc{f} \left( b \right) }[0][\mathcal{B}]
		\end{aligned}
	\end{equation*}
	which shows \cref{eq:cycl-f-Ho-0}.
\end{proof}

\begin{rem} \label{rem:Go-0-closed-no-curvature}
	Assume that $b \in \mc{\mathcal{A}}$ is a strong bounding cochain. Then in general,
	$\Go{b}[0]$ is not a cycle of $\ncdfr{\mathcal{A}}[0][]$. Instead,
	we have $\clie{\mu} \left( \Go{b}[0] \right) = -\mu_0 \left( 1 \right)$ by
	\cref{eq:clie-mu-Go-0}. Working in the extended cyclic complex $\ncdf{\mathcal{A}}[0][]$
	allows us to cancel the contribution of $-\mu_0 \left( 1 \right)$ by adding $1$ to
	$\Go{b}[0]$. Similarly, \cref{eq:cycl-f-Go-0} shows that $\Go{b}[0]$ is not natural
	in general, but is natural for morphisms with $f_0 \left( 1 \right) = 0$.

	If we restrict our
	attention to uncurved Banach $\Ainf$-algebras and morphisms without a change of connection term,
	then the truncated cyclic exponential $\Go{b}[0]$ gives a canonical cycle in the
	standard Connes complex
	satisfying \cref{lm:G_0-bounding-chain-closed,lm:invariance-Gb0-pseudoisotopy}, without
	needing to work with the extended complex or extended cyclic homology.
\end{rem}

\subsection{Pre-\texorpdfstring{$\infty$}{Infinity}-Traces} \label{sec:pre-infinity-traces}

\begin{dfn} \label{dfn:pre-infty-trace}
	An $n$-\textbf{dimensional pre}-$\infty$-\textbf{trace} on $\mathcal{A}$ is a morphism
	\begin{equation*}
		\theta \colon \ncdfr{\mathcal{A}}[0][] \rightarrow \mathcal{R}[1-n]
	\end{equation*}
	of differential graded Banach $\mathcal{R}$-modules between the cyclic complex
	$\ncdfr{\mathcal{A}}[0][]$ and $\mathcal{R}[1-n]$, i.e., a contractive degree zero $R$-linear map which satisfies
	$d_{\mathcal{R}[1-n]} \circ \theta = \theta \circ \clie{\mu}$.
\end{dfn}

In what follows, the dimension $n$ will be fixed, and we often omit it for brevity.
Let us write more explicitly what \cref{dfn:pre-infty-trace} means. A
pre-$\infty$-trace $\theta$ is determined uniquely by the associated sequence
$\left( \theta_k \colon A^{\times k} \rightharpoonup R \right)_{k \geq 1}$ of contractive,
$R$-multilinear maps of degree $1 - n$, called the \textbf{components} of $\theta$, related
to $\theta$ by
\begin{equation*}
	\s_{1 - n} \left( \theta_k \left( a_1, \dots, a_k \right) \right) =
	\theta \left( a_1 \otimes \dots \otimes a_k \right)
\end{equation*}
for $k \geq 1$ and $a_1, \dots, a_k \in A$.
The components $\theta_k$ of $\theta$ are cyclically invariant in the sense that
\begin{equation} \label{eq:theta-k-cyc-invariant}
	\theta_k \left( a_1, \dots, a_k \right) =
	(-1)^{\degb{a_k} \left( \degb{a_1} + \dots + \degb{a_{k-1}} \right)}
	\theta_k \left( a_k, a_1, \dots, a_{k-1} \right)
\end{equation}
and are required to satisfy the identities
\begin{gather} \label{eq:infty-trace-k-geq-1-rel}
	(-1)^{1 - n} d \left( \theta_{k+1} \left( a_0, \dots, a_k \right) \right) =
	\\
	\sum_{k_1 + k_2 + k_3 = k}
	(-1)^{\varepsilon_1}
	\theta_{1 + k_2} \left(
	\mu_{k_3 + 1 + k_1} \left( a_{k_1 + k_2 + 1}, \dots, a_k, a_0, a_1, \dots, a_{k_1} \right),
	a_{k_1+1}, \dots, a_{k_1 + k_2} \right) +
	\notag
	\\
	\sum_{k_1 + k_2 + k_3 = k}
	(-1)^{\varepsilon_2}
	\theta_{1 + k_1 + 1 + k_3} \left(
	a_0, \dots, a_{k_1}, \mu_{k_2} \left( a_{k_1 + 1}, \dots, a_{k_1 + k_2} \right),
	a_{k_1 + k_2 + 1}, \dots, a_k
	\right)
	\notag
\end{gather}
for $k \geq 0$, where
\begin{align*}
	\varepsilon_1 & \defeq
	\left( \degb{a_{k_1 + k_2 + 1}} + \dots + \degb{a_k} \right)
	\left( \degb{a_0} + \dots + \degb{a_{k_1 + k_2}} \right),
	\\
	\varepsilon_2 & \defeq \degb{a_0} + \dots + \degb{a_{k_1}}.
\end{align*}

A pre-$\infty$-trace $\theta$ is called \textbf{strict} if $\theta_k = 0$ for $k \neq 1$. A
strict pre-$\infty$-trace is determined by a single $R$-linear map $\theta_1 \colon A \rightharpoonup R$ of
degree $1 - n$, called the \textbf{trace}, which is required to satisfy
\begin{equation} \label{eq:infty-trace-strict-rel}
	\sum_{\textrm{cyc}} (-1)^{\varepsilon}
	\theta_1 \left( \mu_{k} \left( a_{i+1}, \dots, a_k, a_1, \dots, a_i \right) \right) =
	\begin{cases}
		0                                                         & k \neq 1, \\
		(-1)^{1 - n} d \left( \theta_1 \left( a_1 \right) \right) & k = 1,
	\end{cases}
\end{equation}
for $k \geq 1$ and $a_1, \dots, a_k \in A$,
where
\begin{equation*}
	\varepsilon = \left( \degb{a_{i+1}} + \dots + \degb{a_k} \right) \left( \degb{a_1} + \dots + \degb{a_i} \right).
\end{equation*}
When $d = 0$ and $\mu_k = 0$ for $k \neq 1,2$, we see that the relations
\eqref{eq:infty-trace-strict-rel} satisfied by the trace
reduce to the relations satisfied by the standard notion of trace
on a differential graded algebra corresponding to $\mathcal{A}$.

\begin{dfn} \label{dfn:pre-infty-trace-unital}
	Let $\mathcal{A} = \left( A, \mu, e \right)$ be a unital Banach $\Ainf$-algebra over $\mathcal{R}$.
	A pre-$\infty$-trace $\theta$ on $\mathcal{A}$ is called \textbf{unital} if
	\begin{equation} \label{eq:pre-infty-trace-unital-condition}
		\theta_k \left( a_1, \dots, a_k \right) = 0
	\end{equation}
	whenever $k \geq 1$ and $a_i = e$ for some $1 \leq i \leq k$.
	Equivalently, a unital pre-$\infty$-trace $\theta$
	is a morphism $\theta \colon \ncdfrred{\mathcal{A}}[0][] \rightarrow \mathcal{R}[1-n]$
	between the \textit{reduced} cyclic complex $\ncdfrred{\mathcal{A}}[0][]$ and
	$\mathcal{R}[1-n]$.
\end{dfn}

\subsection{\texorpdfstring{$\infty$}{Infinity}-Traces, the \texorpdfstring{$\infty$}{infinity}-modulus and its Invariance}
\label{sec:traces-modulus-invariance}

\begin{dfn} \label{dfn:infty-trace}
	An $n$-\textbf{dimensional} $\infty$-\textbf{trace} on $\mathcal{A}$ is a morphism
	\begin{equation*}
		\theta \colon \ncdf{\mathcal{A}}[0][] \rightarrow \mathcal{R}[1-n]
	\end{equation*}
	of differential graded Banach $\mathcal{R}$-modules between the \textit{extended} cyclic complex
	$\ncdf{\mathcal{A}}[0][]$ and $\mathcal{R}[1-n]$, i.e., a contractive degree zero $R$-linear map which satisfies
	$d_{\mathcal{R}[1-n]} \circ \theta = \theta \circ \clie{\mu}$.

	An $n$-\textbf{dimensional} $\infty$-\textbf{trace Banach} $\Ainf$-\textbf{algebra} over
	$\mathcal{R}$ is a triple $\mathcal{A} = \left( A, \mu, \theta \right)$ where $\left( A, \mu \right)$ is a
	Banach $\Ainf$-algebra over $\mathcal{R}$ and $\theta$ is an $n$-dimensional $\infty$-trace on
	$\mathcal{A}$.
\end{dfn}

In what follows, the dimension $n$ will be fixed, and we often omit it for brevity.
The data of an $\infty$-trace $\theta$ on $\mathcal{A}$ is equivalently given by
a pre-$\infty$-trace $\rest{\theta}{\ncdfr{\mathcal{A}}[0][]}$, called
the \textbf{associated pre}-$\infty$-\textbf{trace}, together
with an element $\theta_0 = \theta \left( 1 \right) \in R^{1 - n}$ which satisfies $\nnorm[\theta_0] \leq 1$
and is related to the associated pre-$\infty$-trace via the relation \eqref{eq:infty-trace-0-rel} below.
Thinking of $\theta_0$ as an operator of arity $0$, we see that an $\infty$-trace is determined uniquely by
the associated sequence $\left( \theta_k \colon A^{\times k} \rightharpoonup R \right)_{k \geq 0}$ of contractive,
$R$-multilinear maps of degree $1 - n$, called the \textbf{components} of $\theta$, related
to $\theta$ by $\s_{1 - n} \left( \theta_0 \right) = \theta \left( 1 \right)$, and
\begin{equation*}
	\s_{1 - n} \left( \theta_k \left( a_1, \dots, a_k \right) \right) =
	\theta \left( a_1 \otimes \dots \otimes a_k \right)
\end{equation*}
for $k \geq 1$ and $a_1, \dots, a_k \in A$.
The maps $\theta_k$ for
$k \geq 1$ are cyclically invariant in the sense of \eqref{eq:theta-k-cyc-invariant} and satisfy
the identities \eqref{eq:infty-trace-k-geq-1-rel}.

In addition, we have an extra identity
\begin{equation} \label{eq:infty-trace-0-rel}
	(-1)^{1-n} d \left( \theta_0 \right) = \theta_1 \left( \mu_0 \left( 1 \right) \right)
\end{equation}
which comes from plugging in the element $1$ of the extended cyclic complex into both sides of the identity
$d_{\mathcal{R}[1-n]} \circ \theta = \theta \circ \clie{\mu}$.
We emphasize that our notion of $\infty$-trace is ``extended'' in the sense that $\theta$ is defined on the
\textit{extended} cyclic complex and hence comes equipped with the extra element $\theta_0$,
considered as an operator of arity $0$, which satisfies the identity
\eqref{eq:infty-trace-0-rel}. For more on the role of $\theta_0$, see \cref{rem:extended-vs-standard-trace}.

An $\infty$-trace $\theta$ is called \textbf{strict} if $\theta_k = 0$ for $k \neq 1$. A
strict $\infty$-trace is determined by a single $R$-linear map $\theta_1 \colon A \rightharpoonup R$ of
degree $1 - n$, called the \textbf{trace}, which is required to satisfy
the identities \eqref{eq:infty-trace-strict-rel} for $k \geq 0$, and not only for $k \geq 1$.
Note that a strict $\infty$-trace $\theta_1$ is the same thing as a strict pre-$\infty$-trace $\theta_1$ which satisfies
$\theta_1 \left( \mu_0 \left( 1 \right) \right) = 0$.

\begin{dfn} \label{dfn:infty-trace-unital}
	Let $\mathcal{A} = \left( A, \mu, e \right)$ be a unital Banach $\Ainf$-algebra over $\mathcal{R}$.
	An $\infty$-trace $\theta$ on $\mathcal{A}$ is called \textbf{unital} if
	\begin{equation} \label{eq:infty-trace-unital-condition}
		\theta_k \left( a_1, \dots, a_k \right) = 0
	\end{equation}
	whenever $k \geq 1$ and $a_i = e$ for some $1 \leq i \leq k$.
	Equivalently, a unital $\infty$-trace $\theta$
	is a morphism $\theta \colon \ncdfred{\mathcal{A}}[0][] \rightarrow \mathcal{R}[1-n]$
	between the \textit{extended reduced} cyclic complex $\ncdfred{\mathcal{A}}[0][]$ and
	$\mathcal{R}[1-n]$.\footnote{Note that the unitality condition involves $\theta_k$ for $k \geq 1$, i.e.,
		$\theta$ is unital if and only if the associated pre-$\infty$-trace $\rest{\theta}{\ncdfr{\mathcal{A}}[0][]}$ is unital.}

	An $n$-\textbf{dimensional unital} $\infty$-\textbf{trace Banach} $\Ainf$-\textbf{algebra}
	over $\mathcal{R}$ is a quadruple $\mathcal{A} = \left( A, \mu, e, \theta \right)$ where
	$\left( A, \mu, e \right)$ is a unital Banach $\Ainf$-algebra over $\mathcal{R}$ and $\theta$ is an
	$n$-dimensional unital $\infty$-trace on $\mathcal{A}$.
\end{dfn}
Since the constructions $\mathcal{A} \mapsto \ncdf{\mathcal{A}}[0][], \ncdfred{\mathcal{A}}[0][]$ are functorial,
we have a natural notion of a morphism between $\Ainf$-algebras equipped with $\infty$-traces:
\begin{dfn} \label{dfn:morphism-inf-trace}
	Given an $\infty$-trace Banach $\Ainf$-algebra
	$\mathcal{A} = \left( A, \mu, \theta_A \right)$ over $\mathcal{R} = \left( R, d_R \right)$ and
	an $\infty$-trace Banach $\Ainf$-algebra $\mathcal{B} = \left( B, \nu, \theta_B \right)$ over $\mathcal{S} = \left( S, d_S \right)$,
	a \textbf{morphism of} $\infty$-\textbf{trace Banach} $\Ainf$-\textbf{algebras} is a
	morphism
	$f \colon \left( A, \mu \right) \rightarrow \left( B, \nu \right)$ of Banach $\Ainf$-algebras such that
	\begin{equation} \label{eq:morphism-extended-infty-trace}
		\theta_B \circ \cycl{f} = \base{f}[1-n] \circ \theta_A,
	\end{equation}
	see \cref{fig:morphism-extended-infty-trace-algebras}. When $\mathcal{A}$ and $\mathcal{B}$ are unital,
	we require that the morphism $f$ is also unital.
\end{dfn}

\begin{figure}
	\begin{tikzcd}
		{\ncdf{\mathcal{A} / \mathcal{R}}[0][]}
		\arrow{r}{\theta_A}
		\arrow{d}[swap]{\cycl{f}}
		&
		{\mathcal{R}[1 - n]}
		\arrow{d}{{\base{f}}[1 - n]}
		\\
		{\ncdf{\mathcal{B} / \mathcal{S}}[0][]}
		\arrow{r}{\theta_B}
		&
		{\mathcal{S}[1 - n]}
	\end{tikzcd}
	\caption{Morphism of $\infty$-trace Banach $\Ainf$-algebras.}
	\label{fig:morphism-extended-infty-trace-algebras}
\end{figure}

Given a topologically nilpotent element $b \in \tc{A}$, consider the canonical chain
\begin{equation*}
	\G{b}[0][] = \G{b}[0][\mathcal{A}] = \cexp{b} = 1 + \sum_{k = 1}^{\infty} \frac{b^{\otimes k}}{k}
	\in \ncdf{\mathcal{A}}[0][0]
\end{equation*}
of the extended cyclic complex, discussed in \cref{subsec:canonical-elements-ncdf-0}.

\begin{dfn} \label{dfn:extended-infty-modulus}
	Let $\mathcal{A} = \left( A, \mu, \theta \right)$ be an $\infty$-trace Banach $\Ainf$-algebra
	and let $b \in \tc{A}$ be a topologically nilpotent element. The $\infty$-\textbf{modulus}
	function $\TM{} \colon \tc{A} \rightarrow R^{1-n}$ is defined by
	\begin{equation} \label{eq:dfn-extended-infty-modulus}
		\TM{b} = \TM{b}[\mathcal{A}] \defeq \theta \left( \G{b}[0][\mathcal{A}] \right) =
		\theta_0 + \sum_{k = 1}^{\infty} \frac{1}{k} \theta_k \left( b^{\otimes k} \right)
	\end{equation}
	and the element $\TM{b}$ is called the $\infty$-\textbf{modulus} of $b$.
\end{dfn}

\begin{lm}[Naturality of the $\infty$-modulus] \label{lm:functoriality-infty-modulus}
	Let $\mathcal{A} = \left( A, \mu, \theta_A \right)$ be an $n$-dimensional $\infty$-trace
	Banach $\Ainf$-algebra over a differential graded-commutative Banach $\mathbbm{k}$-algebra
	$\mathcal{R} = \left( R, d_R \right)$ and let $\mathcal{B} = \left( B, \nu, \theta_B \right)$
	be an $n$-dimensional $\infty$-trace Banach $\Ainf$-algebra
	over a differential graded-commutative Banach $\mathbbm{k}$-algebra $\mathcal{S} = \left( S, d_S \right)$.
	Given a morphism $f \colon \mathcal{A} \rightarrow \mathcal{B}$ of $\infty$-trace Banach
	$\Ainf$-algebras and $b \in \tc{A}$, we have
	\begin{equation} \label{eq:infty-modulus-functoriality}
		\base{f} \left( \TM{b}[\mathcal{A}] \right) =
		\TM{\mcfunc{f} \left( b \right)}[\mathcal{B}].
	\end{equation}
\end{lm}
\begin{proof}
	By \cref{lm:G-0-properties}, we have
	\begin{equation*}
		\begin{aligned}
			\base{f}[1-n] \left( \TM{b}[\mathcal{A}] \right)
			\stackrel{\eqref{eq:dfn-extended-infty-modulus}}{=}{}    &
			\base{f}[1-n] \left( \theta_A \left( \G{b}[0][\mathcal{A}] \right) \right)
			\\
			\stackrel{\eqref{eq:morphism-extended-infty-trace}}{=}{} &
			\theta_B \left( \cycl{f} \left( \G{b}[0][\mathcal{A}] \right) \right)
			\\
			\stackrel{\eqref{eq:cycl-f-G0}}{=}{}                     &
			\theta_B \left( \G{\mcfunc{f} \left( b \right)}[0][\mathcal{B}] \right)
			\\
			\stackrel{\eqref{eq:dfn-extended-infty-modulus}}{=}{}    &
			\TM{\mcfunc{f} \left( b \right)}[\mathcal{B}].
		\end{aligned}
	\end{equation*}
\end{proof}

\begin{lm}[$\infty$-modulus of a Bounding Cochain is a Cocycle]
	\label{lm:infty-modulus-bounding-chain-closed}
	Let $\mathcal{R}$ be a differential graded-commutative Banach $\mathbbm{k}$-algebra,
	and let $\mathcal{A}$ be a non-unital (resp.\ unital) $n$-dimensional $\infty$-trace
	Banach $\Ainf$-algebra over $\mathcal{R}$. Let $b \in \tc{A}$ be a strong (resp.\ weak) bounding cochain.
	Then $\TM{b}$ is a cocycle which defines a cohomology class
	$\eqcl{ \TM{b} } \in \cohom{\mathcal{R}}[1-n]$.
\end{lm}
\begin{proof}
	When $\mathcal{A}$ is non-unital (resp.\ unital) $\infty$-trace Banach $\Ainf$-algebra,
	and $b$ is a strong (resp.\ weak) bounding cochain, \cref{lm:G_0-bounding-chain-closed} shows that
	the chain $\G{b}[0]$ is a cycle of $\ncdf{\mathcal{A}}[0][0]$ (resp.\ $\ncdfred{\mathcal{A}}[0][0]$).
	Since $\TM{b} = \theta(\G{b}[0])$ is the application of a chain map of degree $1 - n$
	to a cycle of degree $0$, the result follows.
\end{proof}

\Cref{dfn:pseudo-isotopy} of pseudoisotopy between two Banach $\Ainf$-algebras
extends naturally to the setting of $\infty$-trace Banach $\Ainf$-algebras as follows.

\begin{dfn}[Pseudoisotopy of $\infty$-trace Banach $\Ainf$-algebras]
	\label{dfn:pseudo-isotopy-inf-trace-algebras}
	Let $\mathcal{A}_0$ and $\mathcal{A}_1$ be two Banach $\Ainf$-algebras over $\mathcal{S}$
	and let $\mathfrak{A}$ be a pseudoisotopy over $\mathfrak{R}$
	between $\mathcal{A}_0$ and $\mathcal{A}_1$.

	Assume that $\mathcal{A}_i$ are equipped with $\infty$-traces $\theta^i \colon \ncdf{\mathcal{A}_i}[0][] \rightarrow \mathcal{S}[1-n]$
	for $i=0,1$, and that $\mathfrak{A}$ is also equipped with an $\infty$-trace
	$\theta \colon \ncdf{\mathfrak{A}}[0][] \rightarrow \mathfrak{R}[1-n]$,
	such that the pseudoisotopy maps $\evalmf^i \colon \mathfrak{A} \rightarrow \mathcal{A}_i$
	become morphisms of $\infty$-trace Banach $\Ainf$-algebras (see \cref{fig:pseudoisotopy-infty-trace-algebras}).
	In this case, we say that $\mathfrak{A}$ is a
	\textbf{pseudoisotopy of} $\infty$-\textbf{trace Banach} $\Ainf$-\textbf{algebras}.

	When $\mathcal{A}_0, \mathcal{A}_1$ and $\mathfrak{A}$ are unital, and the traces
	$\theta^0,\theta^1$ and $\theta$ are also unital, we say that $\mathfrak{A}$ is a
	\textbf{unital pseudoisotopy of} (\textbf{unital}) $\infty$-\textbf{trace Banach} $\Ainf$-\textbf{algebras}.
\end{dfn}

\begin{figure}
	\begin{tikzcd}
		{\ncdf{\mathfrak{A} / \mathfrak{R}}[0][]}
		\arrow{r}{\theta} \arrow{d}[swap]{\cycl{\evalmf^i}}
		&
		{\mathfrak{R}[1 - n]} \arrow{d}{{\evalm^i}[1-n]}
		\\
		{\ncdf{\mathcal{A}_i / \mathcal{S}}[0][]}
		\arrow{r}{\theta^i}
		&
		{\mathcal{S}[1-n]}
	\end{tikzcd}
	\caption{Pseudoisotopy of $\infty$-trace Banach $\Ainf$-algebras.}
	\label{fig:pseudoisotopy-infty-trace-algebras}
\end{figure}

\begin{thm}[Invariance of the $\infty$-modulus under Gauge Equivalence]
	\label{thm:invariance-extended-infty-modulus}
	Let $\mathcal{A}_0$ and $\mathcal{A}_1$ be two non-unital (resp.\ unital) $n$-dimensional
	$\infty$-trace Banach $\Ainf$-algebras over $\mathcal{S} = \left( S, d_S \right)$
	and let $\mathfrak{A}$ be a non-unital (resp.\ unital) pseudoisotopy
	of $\infty$-trace Banach $\Ainf$-algebras between $\mathcal{A}_0$ and $\mathcal{A}_1$,
	defined over $\mathfrak{R} = \left( R, d_R \right)$.
	Let $b_0 \in \tc{A_0}$ and $b_1 \in \tc{A_1}$ be two
	$\mathfrak{A}$-gauge-equivalent strong (resp.\ weak) bounding cochains.
	Then
	\begin{equation*}
		\eqcl{ \TM{b_0}[\mathcal{A}_0] } = \eqcl{ \TM{b_1}[\mathcal{A}_1] }
		\in \cohom{\mathcal{S}}[1-n].
	\end{equation*}
\end{thm}
\begin{proof}
	Let $b \in \tc{A}$ with $\mcfunc{\evalmf}^i \left( b \right) = b_i$ for $i = 0, 1$
	and let $\theta \colon \ncdf{\mathfrak{A}}[0][] \rightarrow \mathfrak{R}[1-n]$ be such that
	the diagram of \cref{fig:pseudoisotopy-infty-trace-algebras} commutes for $i = 0, 1$.
	Choose a homotopy $h \colon R^{*} \rightharpoonup S^{*-1}$ between
	$\evalm^0 = \base{\evalmf}^0$ and $\evalm^1 = \base{\evalmf}^1$ with
	\begin{equation*}
		d_S \circ h + h \circ d_R = \evalm^1 - \evalm^0.
	\end{equation*}
	Then
	\begin{equation*}
		\begin{aligned}
			\TM{b_1}[\mathcal{A}_1] - \TM{b_0}[\mathcal{A}_0]
			\stackrel{\eqref{eq:infty-modulus-functoriality}}{=}{}           &
			\evalm^1 \left( \TM{b}[\mathfrak{A}] \right) -
			\evalm^0 \left( \TM{b}[\mathfrak{A}] \right)
			\\
			\stackrel{\phantom{\eqref{eq:infty-modulus-functoriality}}}{=}{} &
			\left( d_S \circ h + h \circ d_R \right) \left( \TM{b}[\mathfrak{A}] \right)
			\\
			\stackrel{\phantom{\eqref{eq:infty-modulus-functoriality}}}{=}{} &
			d_S \left( h \left( \TM{b}[\mathfrak{A}] \right) \right)
		\end{aligned}
	\end{equation*}
	where we used the fact that $d_R \left( \TM{b}[\mathfrak{A}] \right) = 0$ by
	\cref{lm:infty-modulus-bounding-chain-closed}.
\end{proof}

Assembling the results of this section, we obtain a complete proof of
\cref{thm:infinity-modulus-properties} from the introduction:
\begin{proof}[Proof of \cref{thm:infinity-modulus-properties}]
	\Cref{item:infinity-modulus-1} of \cref{thm:infinity-modulus-properties}
	was proved in \cref{lm:functoriality-infty-modulus},
	\cref{item:infinity-modulus-2} was proved in \cref{lm:infty-modulus-bounding-chain-closed},
	and  \cref{item:infinity-modulus-3} was proved in \cref{thm:invariance-extended-infty-modulus}.
\end{proof}

\begin{ex} \label{ex:integral-invariant-gauge-equivalence}
	As discussed in \cref{rem:abstract-pseudoisotopy-vs-jake,rem:abstract-gauge-equivalence-vs-jake}, our notions of
	pseudoisotopy and gauge equivalence generalize the notions appearing in \cite{Solomon2016,Solomon2016a}.
	In \cite[Lemma 3.16]{Solomon2016a} the authors show that given a bounding pair $\left( \gamma, b \right)$
	with respect to an almost complex structure $J$ and a bounding pair $\left( \gamma', b' \right)$ with respect to an
	almost complex structure $J'$, the integral of the bounding cochain is invariant under gauge equivalence, i.e., we have
	$\int_{L} b = \int_{L} b'$. This also follows from \cref{thm:invariance-extended-infty-modulus} as the integral
	defines a strict unital $\infty$-trace on the $\Ainf$-algebras involved.

	More precisely, given the premises of \cite[Lemma 3.16]{Solomon2016a}, as described in \cref{rem:abstract-gauge-equivalence-vs-jake},
	let us set $A = {\cdiff{L}[][][R]}[1]$ and $\widetilde{A} = {\cdiff{I \times L}[][][R]}[1]$.
	Then we have:
	\begin{enumerate}
		\item A Banach $\Ainf$-algebra $\mathcal{A}_0 = \left( A, \mu^0, e_0 \right)$ over $\left( R, 0 \right)$
		      associated to $\mathfrak{m}^{\gamma,J}$, together with a weak bounding cochain $b_0$ of $\mathcal{A}_0$.
		\item A Banach $\Ainf$-algebra $\mathcal{A}_1 = \left( A, \mu^1, e_1 \right)$ over $\left( R, 0 \right)$
		      associated to $\mathfrak{m}^{\gamma',J'}$, together with a weak bounding cochain $b_1$ of $\mathcal{A}_1$.
		\item A strict pseudoisotopy $\mathfrak{A} = \left( \widetilde{A}, \mu, e \right)$ over $\left( \cdiff{I}[][][R], d \right)$
		      between $\mathcal{A}_0$ and $\mathcal{A}_1$,	associated to $\widetilde{\mathfrak{m}}^{\widetilde{\gamma},\widetilde{J}}$,
		      with respect to which $b_0$ and $b_1$ are $\mathfrak{A}$-gauge-equivalent.
	\end{enumerate}
	Note that algebras $\mathcal{A}_0,\mathcal{A}_1$ have the same underlying module $A$ and even the same unit, but are
	equipped with different $\Ainf$-structures $\mu^0$ (resp.\ $\mu^1$) because of the dependence on $\left( \gamma, J \right)$
	(resp.\ $\left( \gamma', J' \right)$).

	Let $t \colon \ncdf{A}[0][] \rightarrow R[1-n]$ be determined by\footnote{The sign factor is immaterial and is there to
		guarantee that $t_1$ is $R$-linear of degree $1 - n$ with the conventions used in \cite{Solomon2016,Solomon2016a}.}
	\begin{equation*}
		t_1 \left( a \right) \defeq (-1)^{\left( -n \right) \left( \degb{a} + 1 \right)} \int_L \s^{-1} a, \qquad t_k = 0, \quad k \neq 1.
	\end{equation*}
	Then $\theta^0 = \theta^1 = t$ determine a strict $\infty$-trace on $\mathcal{A}_i$ for $i = 0, 1$, i.e., $t$
	satisfies the identities of \eqref{eq:infty-trace-strict-rel} for $k \geq 0$ with respect to both $\mu^0$ and $\mu^1$.
	This follows from the properties proved in \cite{Solomon2016}, see also \cref{sec:converting-jake-to-cyclic-structure}.
	Similarly, let $\theta \colon \ncdf{\widetilde{A}}[0][] \rightarrow {\cdiff{I}[][][R]}[1 - n]$ be determined by
	\begin{equation*}
		\theta_1 \left( \widetilde{a} \right) \defeq
		(-1)^{\left( -n \right) \left( \degb{\widetilde{a}} + 1 \right)} \left( \pi_{I} \right)_{*} \left( \s^{-1} \left( \widetilde{a} \right) \right),
		\qquad
		\theta_k = 0, \quad k \neq 1,
	\end{equation*}
	where $\left( \pi_{I} \right)_{*}$ is the pushforward (fiberwise integration) along the projection
	$\pi_I \colon I \times L \rightarrow I$. Then $\theta$ determines a strict $\infty$-trace on $\mathfrak{A}$
	with respect to which $\mathfrak{A}$ becomes a pseudoisotopy of $\infty$-trace Banach $\Ainf$-algebras.
	When $\dim L > 0$, the $\infty$-traces are unital and thus the invariance of the integral
	under gauge equivalence follows from \cref{thm:invariance-extended-infty-modulus}.
\end{ex}

\begin{rem}
	A version of the $\infty$-modulus function adapted to an $\infty$-trace that is not strictly unital is used in~\cite[Theorem 11]{Sela2024} to classify bounding cochains in the endomorphism algebra of a normed matrix factorization.
\end{rem}

\begin{rem} \label{rem:extended-vs-standard-trace}
	Our notion of $\infty$-trace and $\infty$-modulus involves an inhomogeneous
	term $\theta_0$, coming from the fact that $\theta$ is defined on the \textit{extended}
	cyclic complex. Let us discuss briefly the role of the term $\theta_0$.

	Given a pre-$\infty$-trace $\theta \colon \ncdfr{\mathcal{A}}[0][] \rightarrow \mathcal{R}[1-n]$
	on $\mathcal{A}$, one can define the \textbf{pre}-$\infty$-\textbf{modulus} function
	$\TM{}[> 0] \colon \tc{A} \rightarrow R^{1-n}$ by
	\begin{equation} \label{eq:pre-infinity-modulus-function}
		\TM{b}[> 0] \defeq \theta \left( \Go{b}[0] \right) =
		\sum_{k = 1}^{\infty} \frac{1}{k} \theta_k \left( b^{\otimes k} \right) \in R^{1-n}.
	\end{equation}
	In general, the pre-$\infty$-modulus of a bounding cochain is not a cocycle and the
	pre-$\infty$-modulus is not invariant under gauge equivalence.
	However, one can deduce from \cref{lm:Go_0-and-Ho_0-properties} that $\TM{b}[>0]$ does give an invariant cocycle,
	if one restricts attention to $\Ainf$-algebras which are not curved and to $\Ainf$-morphisms
	without a change of connection term. See also \cref{rem:Go-0-closed-no-curvature}.

	Going back to the general case, an $\infty$-trace $\theta$ on $\mathcal{A}$ is given by the same data as a pre-$\infty$-trace
	$\rest{\theta}{\ncdfr{\mathcal{A}}[0][]}$, together
	with an additional $0$-th component
	$\theta_0 = \theta \left( 1 \right) \in R^{1-n}$, satisfying the extra identity
	$(-1)^{1-n} d \left( \theta_0 \right) = \theta_1 \left( \mu_0 \left( 1 \right) \right)$
	of  \cref{eq:infty-trace-0-rel}.
	Note that even when $\theta_0 = 0$, being an $\infty$-trace imposes the extra condition
	$\theta_1 \left( \mu_0 \left( 1 \right) \right) = 0$, compared to a pre-$\infty$-trace.

	Any pre-$\infty$-trace $\theta$ which satisfies $\theta_1 \left( \mu_0 \left( 1 \right) \right) = 0$
	defines an $\infty$-trace $\theta$ by setting $\theta_0 = 0$. When working with such
	pre-$\infty$-traces and pseudoisotopies which involve no change of connection term,
	our $\infty$-modulus $\TM{b}$ reduces to the pre-$\infty$-modulus $\TM{b}[>0]$,
	so $\TM{b}[>0]$ gives a gauge invariant cocycle.
	However, in the general case, only the \textit{sum} $\TM{b} = \theta_0 + \TM{b}[> 0]$
	is natural and gives an invariant cocycle.
	This is similar to what happens in \cite{Fukaya2011,Joyce2008,Solomon2016a}
	when one adds an inhomogeneous term $\mathfrak{m}_{-1}$ to the superpotential to make it invariant.
	See also \cref{rem:pre-total-vs-total-inner-product}.
\end{rem}

\section{Total Inner Products and the Superpotential} \label{sec:generalized-superpotential}

In this section,
we define the notion of a total inner product, which is
a chain map defined on a total complex constructed from the bicomplex of
cyclic codifferential forms. We construct a cyclic Chern--Simons form
in the total complex and use it to define the superpotential function, which
is shown to be natural up to exact terms
and gauge invariant on (weak) bounding cochains. We start in
\cref{subsec:generalized-superpotential-intro} with a summary of previous work on inner products on $\Ainf$-algebras and associated definitions of the superpotential, and describe
our objectives. \Cref{sec:sp-technical-overview} gives a technical overview which explains the organization of the rest of the section.

In what follows, we fix a field $\mathbbm{k}$ of characteristic zero, endowed with the trivial norm.

\subsection{Motivation} \label{subsec:generalized-superpotential-intro}
Let $\mathcal{A} = \left( A, \mu \right)$ be a Banach $\Ainf$-algebra over a
differential graded-commutative Banach $\mathbbm{k}$-algebra $\mathcal{R} = (R,d)$.
Recall that a cyclic structure on $\mathcal{A}$ is given by:

\begin{dfn}[{\cite[Definition 1.1]{Solomon2016}}] \label{dfn:cyclic-structure}
	An $n$-\textbf{dimensional cyclic structure} on $\mathcal{A}$ is an
	$R$-linear map $\inncur \colon A \otimes A \rightarrow R[2 - n]$ which satisfies:
	\begin{enumerate}
		\item (Antisymmetry)
		      \begin{equation} \label{eq:antisymmetry}
			      \inncur[a][b] = (-1)^{\degb{a} \cdot \degb{b} + 1} \inncur[b][a].
		      \end{equation}
		\item (Cyclic Pairing)
		      \begin{align}
			      \inncur[\mu_k \left( a_1, \dots, a_k \right)][a_{k+1}] ={} &
			      (-1)^{\degb{a_{k+1}} \cdot \left( \degb{a_1} + \dots + \degb{a_k} \right)}
			      \inncur[\mu_k \left( a_{k+1}, a_1, \dots, a_{k-1} \right)][a_k] \notag
			      \\
			                                                                 & +
			      \delta_{1,k} \cdot d_{\mathcal{R}[2-n]} \left( \inncur[a_1][a_2] \right).
			      \label{eq:cyclic-pairing}
		      \end{align}
		\item (Contractive)
		      \begin{equation}
			      \nnorm[{\inncur[a][b]}] \leq \nnorm[a] \cdot \nnorm[b].
		      \end{equation}
	\end{enumerate}
\end{dfn}
Given a cyclic structure $\inncur$ on $\mathcal{A}$, the superpotential function
$\SP[][>0] \colon \tc{A} \rightarrow R^{3-n}$ associated to $\mathcal{A}$ and $\inncur$ is given by
\begin{equation} \label{eq:sp-gt-0-cyclic-structure}
	\SP[b][>0] \defeq \sum_{k=0}^{\infty} \frac{1}{k+1} \inncur[\mu_k \left( b, \dots, b \right)][b].
\end{equation}
When $\mathcal{A}$ is curved, the superpotential $\SP[][>0]$ is often modified to include
an inhomogeneous term, denoted by $\mathfrak{m}_{-1} \in R^{3-n}$ in~\cite{Fukaya2011,Joyce2008,Solomon2016a,Solomon2024},
usually coming naturally from the geometric construction of $\mathcal{A}$, so that the superpotential
takes the form
\begin{equation} \label{eq:sp-cyclic-structure}
	\SP[b] \defeq \mathfrak{m}_{-1} + \sum_{k=0}^{\infty} \frac{1}{k+1} \inncur[\mu_k \left( b, \dots, b \right)][b].
\end{equation}
The term $\mathfrak{m}_{-1}$ is related to the cyclic structure via an identity~\cite[Proposition 4.20]{Solomon2016} of the form
\begin{equation} \label{eq:d-m-minus-1-identity}
	d \mathfrak{m}_{-1} \equiv -\frac{1}{2} \inncur[\mu_0 \left( 1 \right)][\mu_0 \left( 1 \right)],
\end{equation}
possibly modulo terms which are later removed from the superpotential to keep it invariant~\cite{Solomon2016a}.

There is a close relation between cyclic structures and the complex
$\ncdf{\mathcal{A}}[2][]$ of cyclic codifferential 2-forms, which we now describe.
In what follows, we will work with the inner product pairing $\braidop_1$
given by \cref{eq:parity-inner-product}. Given an $n$-dimensional
cyclic structure $\inncur$ on $\mathcal{A}$, define an $R$-linear contractive map
$\phi_2 \colon \ncdf{A}[2][] \rightarrow R[2-n]$ by
\begin{empheq}[left={ \phi_2 \left( \ul{a} \otimes a_1 \otimes \dots \otimes a_i \otimes \ul{b} \otimes b_1 \otimes \dots \otimes b_j \right) \defeq \empheqlbrace}]{align}
	&\inncur[a][b] && i = j = 0,
	\label{eq:cyclic-structure-to-phi-2-0-0}
	\\
	& \,\, 0 && i > 0 \textrm{ or } j > 0.
	\label{eq:cyclic-structure-to-phi-2-0-gt-0}
\end{empheq}
The map $\phi_2$ is well-defined by the antisymmetry of $\inncur$, and the cyclic pairing property
implies that $\phi_2$ is a \textit{chain map}, i.e., $\phi_2 \circ \clie{\mu} = d_{\mathcal{R}[2-n]} \circ \phi_2$.
Conversely, let us say that a map $\phi_2 \colon \ncdf{A}[2][] \rightarrow R[2-n]$
is \textbf{strict} if it satisfies \eqref{eq:cyclic-structure-to-phi-2-0-gt-0}, i.e.,
only the values $\phi_2 \left( \ul{a} \otimes \ul{b} \right)$ may be non-zero.
Then a strict contractive chain map $\phi_2$ defines a cyclic structure via $\inncur[a][b] \defeq \phi_2 \left( \ul{a} \otimes \ul{b} \right)$.
For details, see \cref{appendix:cyclic-structures}.

The condition of being a strict map is quite rigid: it is not a cohomological condition,
nor preserved under general $\Ainf$-morphisms. More precisely,
if $\phi_2$ is a strict chain map and
$\psi \colon \ncdf{A}[2][] \rightarrow R[1-n]$ is an arbitrary map,
then $\phi_2 + \partial \psi  \colon \ncdf{\mathcal{A}}[2][] \rightarrow \mathcal{R}[2-n]$
is a chain map but not necessarily strict.
In addition, given an $\Ainf$-morphism
$f \colon \mathcal{B} \rightarrow \mathcal{A}$, the composition
$\phi_2 \circ \cindmap{f} \colon \ncdf{\mathcal{B}}[2][] \rightarrow \mathcal{R}[2-n]$
is a chain map, but won't be strict in general unless $f$ is also strict, and will not give us a
cyclic structure on $\mathcal{B}$.

Note however that a strict map $\phi_2$ also trivially satisfies $\phi_2 \circ \qdr^3 = 0$,
and the condition $\phi_2 \circ \qdr^3 = 0$ is a cohomological condition which is preserved by
naturality under general $\Ainf$-morphisms.
Here, $\qdr^3 \colon \ncdf{A}[3][] \rightarrow \ncdf{A}[2][]$ is the de Rham differential.
By relaxing the condition of strictness and replacing it with the condition
$\phi_2 \circ \qdr^3 = 0$, we obtain the notion of a pre-homotopy inner product:

\begin{dfn} \label{dfn:pre-homotopy-inner-product}
	An $n$\textbf{-dimensional pre-homotopy inner product} on $\mathcal{A}$ is a contractive
	$R$-linear map $\phi_2 \colon \ncdf{A}[2][] \rightarrow R[2-n]$ which satisfies
	$d_{\mathcal{R}[2-n]} \circ \phi_2 = \phi_2 \circ \clie{\mu}$ and $\phi_2 \circ \qdr^3 = 0$.
	Equivalently, an $n$-dimensional pre-homotopy inner product is a morphism
	$\phi_{2} \colon \ncdf{\mathcal{A}}[2][] / \Im \left( \qdr^3 \right) \rightarrow \mathcal{R}[2-n]$
	of differential graded Banach $\mathcal{R}$-modules.
\end{dfn}

The condition $\phi_2 \circ \qdr^3 = 0$ appears under the name \textit{closedness} in~\cite{Cho2008,cho-homotopy-superpotential}.
In the setting of \cite{cho-homotopy-superpotential}, working with uncurved $\Ainf$-algebras
and morphisms without change of connection elements,
a pre-homotopy inner product satisfying a non-degeneracy condition is called
a \textit{strong homotopy inner product}. Given a strong homotopy inner product, Cho and Lee
introduce in \cite[Definition 3.1]{cho-homotopy-superpotential} a generalized superpotential function via the explicit formula
\begin{equation} \label{eq:sp-cho-lee}
	\SP[b][>0] \defeq \sum_{k, l, m = 0}^{\infty} \frac{1}{k + l + 1 + m} \phi_2 \left(
	\ul{ \mu_{k} \left( b^{\otimes k} \right)} \otimes b^{\otimes l} \otimes
	\ul{b} \otimes b^{\otimes m} \right),
\end{equation}
calculate its derivative and show its invariance in an appropriate sense. When $\phi_2$ corresponds
to a cyclic structure, the superpotential \eqref{eq:sp-cho-lee} reduces to the standard
superpotential \eqref{eq:sp-gt-0-cyclic-structure}.

In the following sections, we introduce the notion of a \textbf{total inner product}, which is a chain
map $\phi \colon \totcompe{\mathcal{A}}[2][] \rightarrow \mathcal{R}[4-n]$ between an
extension of the total complex $\totcomp{\mathcal{A}}[2][]$ of cyclic codifferential forms
of degree greater than or equal to two and $\mathcal{R}[4-n]$. A total inner product
is determined uniquely by its components, which are a sequence $\left( \phi_k \right)_{k=2}^{\infty}$
of $R$-linear contractive maps $\phi_k \colon \ncdf{A}[k][] \rightharpoonup R$
of degree $4 - n - k$, together with an extra element $\phi_{\ul{1}} \in R^{3-n}$, required to satisfy
the identities
\begin{align}
	\label{eq:d-phi-ul-1-intro}
	(-1)^{4 - n} d \left( \phi_{\ul{1}} \right)     & =
	\sum_{k=2}^{\infty} \frac{(-1)^{\frac{k \left( k + 1 \right)}{2}}}{k}
	\phi_k \left( \ulz{\mu}^{\otimes k} \right),
	\\
	\label{eq:partial-phi-2}
	\partial \left( \phi_2 \right)                  & = 0,
	\\
	\label{eq:partial-phi-k-plus-1}
	(-1)^{4 - n} \partial \left( \phi_{k+1} \right) & =
	\phi_k \circ \qdr^{k+1}, \qquad k \geq 2,
\end{align}
where $\partial$ denotes the differentials on $\InnHom{\ncdf{\mathcal{A}}[k][]}{\mathcal{R}}$.
The notion of a total inner product replaces the closedness condition $\phi_2 \circ \qdr^3 = 0$
with the infinite family of conditions \eqref{eq:partial-phi-k-plus-1}
requiring that the maps $\phi_k \circ \qdr^{k+1}$ are exact
in $\InnHom{\ncdf{\mathcal{A}}[k][]}{\mathcal{R}}$. The component $\phi_{\ul{1}}$ (or, more precisely, $(-1)^{4-n} \phi_{\ul{1}}$)
	plays the role of the term $\mathfrak{m}_{-1}$, with the identity \eqref{eq:d-phi-ul-1-intro}
	generalizing the identity \eqref{eq:d-m-minus-1-identity}.

	Given a total inner product, we define the (generalized) superpotential function by
	\begin{align}
		\SP[b] ={} & \phi_{\ul{1}} +
		\sum_{\substack{k=2 \\ i_1,\dots,i_{k-1} = 0 \\ j_1,\dots,j_k=0}}^{\infty}
		\frac{(-1)^{\frac{(k - 2) \cdot (k - 1)}{2}}}{1 + \sum_{r=1}^{k-1} i_r + \sum_{r = 1}^k j_r}
		\label{eq:sp-general-braidop-1}
		\\
		           & \qquad\qquad\qquad
		\phi_k \left(
		\ul{ \mu_{i_1} \left( b^{\otimes i_1} \right)} \otimes b^{\otimes j_1} \otimes \dots
		\otimes \ul{ \mu_{i_{k-1}} \left( b^{\otimes i_{k-1}} \right) } \otimes
		b^{\otimes j_{k-1}} \otimes \ul{b} \otimes b^{\otimes j_k}
		\right),
		\notag
	\end{align}
	which generalizes both \eqref{eq:sp-cyclic-structure} and \eqref{eq:sp-cho-lee}.
	The formula for the superpotential comes from applying $\phi$ to a canonical form
	in the extended total complex $\totcompe{\mathcal{A}}[2][]$ which we call
	the cyclic Chern--Simons form.
	The cyclic Chern--Simons form is closed when $b$ is a bounding cochain, giving a cohomological interpretation of the superpotential.
	We show that the superpotential is invariant under gauge equivalence of bounding
	cochains and calculate its derivative. Our results hold for both strong bounding cochains in the non-unital case,
	and weak bounding cochains in the unital case, where modifications are needed to
	account for the presence of the unit.

	\subsection{Technical Overview} \label{sec:sp-technical-overview}
	In \cref{sec:construction-special-elements}, we work with the complex
	\begin{equation*}
		\totcomp{\mathcal{A}}[1][] = \totcomp{\mathcal{A}, \braidop_2}[1][] =
		\totc{\ncdf{A, \braidop_2}, -\qdr^{\braidop_2},
			\clie{\mu}^{\braidop_2}}[][\geq 1][\braidop_2],
	\end{equation*}
	the total complex of cyclic codifferential forms with degree greater than or equal to one,
	constructed using the total degree parity form $\braidop_2$. The complex $\totcomp{\mathcal{A}}[1][]$
	is homotopy equivalent to the shift ${\ncdfr{\mathcal{A}}[0][]}[1]$ of Connes' cyclic complex
$\ncdfr{\mathcal{A}}[0][]$ via the natural degree one projection
$p_1 \colon \totcomp{\mathcal{A}}[1][] \rightharpoonup \ncdfr{\mathcal{A}}[0][]$.

	Given a topologically nilpotent element $b$,
	we construct chains $\Go{b}[\geq 1]$ and $\Ho{b}[\geq 1]$ in $\totcomp{\mathcal{A}}[1][]$
	which satisfy properties analogous to the ones satisfied by the chains $\Go{b}[0]$ and
$\Ho{b}[0]$ of \cref{subsec:canonical-elements-ncdf-0}. The only difference is
	that the ``relative naturality'' of $\Go{b}[\geq 1]$ holds only up to an exact term
	(see \cref{lm:Go_0-and-Ho_0-properties} vs.\ \cref{lm:Go-and-Ho-geq1-properties-general}).
	The chain $\Go{b}[\geq 1]$ (resp.\ $\Ho{b}[\geq 1]$) projects onto the chain
$\Go{b}[0]$ (resp.\ $-\Ho{b}[0]$) under $p_1$.
	The construction of $\Go{b}[\geq 1]$ and $\Ho{b}[\geq 1]$ is carried out in \cref{sec:lift-construction-R-linear}
	for a Banach $\Ainf$-algebra $\mathcal{A} = \left( A, \mu \right)$ over a graded-commutative
	Banach $\mathbbm{k}$-algebra $R$, i.e., when $\mu$ is $R$-linear, using the exponential of the
	cyclic contraction $\exp \left( \ccont{\mu} \right)$. The construction is then generalized
	in \cref{sec:lift-general-construction} for a Banach $\Ainf$-algebra over a
	\textit{differential} graded-commutative Banach $\mathbbm{k}$-algebra $\mathcal{R}$.

	In \cref{sec:extended-tot-comp-geq1-braidop-2}, we define an extended total complex
$\totcompe{\mathcal{A}}[1][]$ by adjoining a shifted copy of $R$, generated by the symbol $\ul{1}$,
	to $\totcomp{\mathcal{A}}[1][]$. The differential of $\ul{1}$ is defined in such a way that the
	projection $p_1$ extends naturally to a homotopy equivalence
$p_1^{+} \colon \totcompe{\mathcal{A}}[1][] \rightharpoonup \ncdf{\mathcal{A}}[0][]$ between
$\totcompe{\mathcal{A}}[1][]$ and the \textit{extended Connes complex} $\ncdf{\mathcal{A}}[0][]$,
	obtaining another model for the \textit{extended} cyclic homology.
	We show that a Banach $\Ainf$-morphism $f \colon \mathcal{A} \rightarrow \mathcal{B}$ induces a
	chain map $\cindmape{f} \colon \totcompe{\mathcal{A}}[1][] \rightarrow \totcompe{\mathcal{B}}[1][]$
	between the extended total complexes, and that this construction is functorial on the level of cohomology.

	With the extended total complex in hand, we define the chain
$\G{b}[\geq 1] = \ul{1} + \Go{b}[\geq 1]$ of $\totcompe{\mathcal{A}}[1][]$,
	analogous to $\G{b}[0] = 1 + \Go{b}[0]$ in the extended Connes complex $\ncdf{\mathcal{A}}[0][]$,
	and show that $\G{b}[\geq 1]$ satisfies properties analogous to the properties
	satisfied by $\G{b}[0]$ (see \cref{lm:G-0-properties} vs.\ \cref{lm:G-geq-1-properties}).

	In \cref{sec:extended-tot-comp-geq2-braidop-2},
	we consider the extended total complex $\totcompe{\mathcal{A}}[2][]$, defined
	by taking the quotient of $\totcompe{\mathcal{A}}[1][]$ by the first column $\ncdf{\mathcal{A}}[1][]$.
	The cyclic Chern--Simons form $\G{b}[\geq 2] \in\totcompe{\mathcal{A}}[2][]$ is defined to be the projection of $\G{b}[\geq 1] \in \totcompe{\mathcal{A}}[1][]$ to $\totcompe{\mathcal{A}}[2][]$.
	The cyclic Chern--Simons form satisfies properties similar to those of $\G{b}[\geq 1]$
	(see \cref{lm:G-geq-2-properties}). When $\mathcal{A}$ is unital, reduced
	versions $\totcompered{\mathcal{A}}[2][], \totcompesred{\mathcal{A}}[2][]$
	of $\totcompe{\mathcal{A}}[2][]$ are introduced in \cref{sec:extended-reduced-total-complexes},
	and are shown to enjoy the same functorial properties as $\totcompe{\mathcal{A}}[2][]$.

	In \cref{sec:canonical-chains-totcompe-geq-2}, we show that when $b$ is a strong (resp.\ weak)
	bounding cochain, the cyclic Chern--Simons form $\G{b}[\geq 2]$ is closed and defines a class
	in $\cohom{{\totcompe{\mathcal{A}}[2][]}}[-1]$ (resp.\ $\cohom{{\totcompered{\mathcal{A}}[2][]}}[-1]$).
	We introduce the notion of a $\cohom{}[] \totcompe{}[2]$-strong (resp.\ $\cohom{}[] \totcompered{}[2]$-strong)
	pseudoisotopy $\mathfrak{A}$, and deduce that the
	class of $\G{b}[\geq 2]$ is invariant under $\mathfrak{A}$-gauge-equivalence,
	when $\mathfrak{A}$ is $\cohom{}[] \totcompe{}[2]$-strong (resp.\ $\cohom{}[] \totcompered{}[2]$-strong).

	\Cref{sec:generalized-inner-product-superpotential} is devoted to total inner products and the
	superpotential. We define an $n$-dimensional total inner product to be a
	morphism $\phi \colon \totcompe{\mathcal{A}}[2][] \rightarrow \mathcal{R}[4-n]$ of
	differential graded Banach $\mathcal{R}$-modules, and we define the superpotential function
	by $\SP[b] = \phi \left( \G{b}[\geq 2] \right)$.
	We introduce morphisms between total inner product Banach $\Ainf$-algebras, i.e.,
	Banach $\Ainf$-algebras equipped with total inner products, and show
	that the superpotential is natural up to an exact term. When $b$ is a bounding cochain, the
	superpotential gives us a natural cohomology class $\eqcl{ \SP[b] } \in \cohom{\mathcal{R}}[3-n]$.
	The notion of a pseudoisotopy
$\mathfrak{A}$ extends naturally to Banach $\Ainf$-algebras equipped with total inner products, and
	we prove that the class $\eqcl{ \SP[b] }$ is invariant under
$\mathfrak{A}$-gauge-equivalence.
	We show the results both for strong bounding cochains, and, in the unital case, for weak bounding cochains,
	where we work with the extended reduced total complex $\totcompered{\mathcal{A}}[2][]$ instead of $\totcompe{\mathcal{A}}[2][]$.

	In \cref{sec:formal-derivative-sp}, we give an explicit formula
	for the formal derivative of the superpotential. We discuss how to perform
	scalar extension for total inner product Banach $\Ainf$-algebras, use it
	to adjoin a formal even variable $t$ to the algebra,
	and use previous results to calculate
	the formal derivative $\partial_t \left( \SP[b] \right)$ of the superpotential,
	where $b$ depends on $t$.

	In \cref{sec:description-using-braid-op-1}, we demonstrate how to obtain
	equivalent notions of the total inner product and the superpotential, starting
	with the total complex
	\begin{equation*}
		\totcomp{\mathcal{A}, \braidop_1}[2][] = \totc{\ncdf{A, \braidop_1}, \qdr^{\braidop_1},
			\clie{\mu}^{\braidop_1}}[][\geq 2][\braidop_1],
	\end{equation*}
	constructed using the inner product parity form $\braidop_1$, instead of
$\totcomp{\mathcal{A}, \braidop_2}[2][]$. The passage between
	the two equivalent descriptions involves signs which we make explicit.
	This is used to deduce \Crefrange{thm:superpotential-properties}{thm:Gb-geq-2-properties}
	of the introduction from the corresponding results we proved working with
	the total degree parity form $\braidop_2$.

	Finally, in \cref{sec:total-inner-product-explicit-relations}, we define the components of a total inner product, which encode the information
	contained in the total inner product completely and describe the relations the components
	must satisfy. We give a formula for the superpotential in terms of its components,
	and relate our notions to those appearing in
	the introduction and \cref{subsec:generalized-superpotential-intro}.

	\subsection{Construction of Canonical Chains in the Total Complex} \label{sec:construction-special-elements}
	Let $\mathcal{R} = (R,d)$ be a differential graded-commutative Banach $\mathbbm{k}$-algebra and let
$\mathcal{A} = \left( A, \mu \right)$ be a Banach $\Ainf$-algebra over $\mathcal{R}$.
	In this section and the following ones, up to \cref{sec:description-using-braid-op-1},
	we will work with the noncommutative differential calculus
	described in \cref{sec:noncomm-diff-calc} using the total degree parity form $\braidop_2$ given by
	\cref{eq:parity-total-degree}.

	Consider the bigraded Banach $R$-module $\ncdf{A} = \tenscyc{A \oplus \ul{A}}[(*,*)]$
	together with the operators $\qdr \colon \ncdf{A} \rightharpoonup \ncdf{A}[* - 1][*]$
	and $\clie{\mu} \colon \ncdf{A} \rightharpoonup \ncdf{A}[*][* + 1]$.
	As a consequence of our definitions and \cref{lm:cyclic-commutation-relations}, we
	have the following relations on $\ncdf{A}$:
	\begin{align}
		\qdr^2                          & = 0,
		\label{eq:qdr-differential-total-degree}
		\\
		\clie{\mu}^2                    & = \frac{1}{2} \left[ \clie{\mu}, \clie{\mu} \right] = \frac{1}{2} \clie{[\mu,\mu]} = 0,
		\label{eq:clie-differential-total-degree}
		\\
		\left[ \qdr, \clie{\mu} \right] & = \qdr \circ \clie{\mu} + \clie{\mu} \circ \qdr = 0.
		\label{eq:qdr-clie-anticommute}
	\end{align}
	Thus, $\qdr$ and $\clie{\mu}$ are \textit{anticommuting} differentials of
	degree $(-1,0)$ and $(0,1)$ respectively, and we obtain a right half-plane Banach bicomplex
$\left( \ncdf{A}, -\qdr, \clie{\mu} \right)$, which we denote by $\ncdf{\mathcal{A}}[][]$, omitting
	the total differential $D = \clie{\mu} - \qdr$ from the notation.

	Consider also the bigraded Banach $R$-submodule $\ncdfr{A} = \tensrcyc{A\oplus \ul{A}}[(*,*)]$ of
$\ncdf{A}$. Recall that $\ncdfr{A}[0][] = \ncdf{A}[0][] / R$, while $\ncdfr{A}[i][] = \ncdf{A}[i][]$ for $i > 0$.
	The differentials $\qdr$ and $\clie{\mu}$ leave $\ncdfr{A}[][]$ invariant, so we get a double
	Banach subcomplex $\left( \ncdfr{A}, -\qdr, \clie{\mu} \right)$, which we denote by
$\ncdfr{\mathcal{A}}[][]$, omitting the total differential from the notation.
	Similarly, given $i \geq 0$, we denote by $\ncdf{\mathcal{A}}[i][]$ (resp.\ $\ncdfr{\mathcal{A}}[i][]$) the
$i$-th column of $\ncdf{\mathcal{A}}[][]$ (resp.\ $\ncdfr{\mathcal{A}}[][]$), endowed with the differential
$\clie{\mu}$.

	Given $k \geq 1$, let us denote by
$\totcomp{\mathcal{A}}[k][] \defeq \totc{\ncdf{\mathcal{A}}[][]}[][\geq k][]$ the total complex
	of the bicomplex obtained from $\ncdf{\mathcal{A}}[][]$ (or $\ncdfr{\mathcal{A}}[][]$)
	by erasing the columns $0, \dots, k - 1$. Since
	the horizontal and vertical differentials of $\ncdf{\mathcal{A}}[][]$ anticommute,
	the totalization is defined as in \cref{subsec:bicomplex-anticommuting-differentials},
	and we denote the total differential on
$\totcomp{\mathcal{A}}[k][]$ by
	\begin{equation*}
		D = D_{\mu} = D_{\mathcal{A}} = D_{\mu}^{\geq k} = D_{\mathcal{A}}^{\geq k} \defeq \clie{\mu} - \qdr,
	\end{equation*}
	depending on whether we want to emphasize the dependence of $D$ on $\mathcal{A}$ (via $\mu$)
	and the underlying graded module $\totcomp{A}[k][]$.

	\begin{rem}
	The bicomplex $\ncdfr{\mathcal{A}}[][]$, and the total complexes
$\totcomp{\mathcal{A}}[1][]$ and $\totcomp{\mathcal{A}}[2][]$ were studied in
	\cref{sec:bicomplex-models-cyclic-homology}, where we used the inner product parity form $\braidop_1$
	and worked with the bicomplex $\left( \ncdfr{A}, \qdr, \clie{\mu} \right)$,
	while in this section we work with $\braidop_2$ and use $-\qdr$ as the horizontal differential.
	The advantage of working with $\braidop_2$ in our setting is that
	the differentials $\qdr$ and $\clie{\mu}$ anticommute, and the totalization
	is defined without any twisting or signs involved
	(see \cref{subsec:bicomplex-anticommuting-differentials}). This helps to minimize
	signs as much as possible,
	but the difference is immaterial and one can always convert the formulas and
	constructions from one convention to the other. We describe how to do that in
	\cref{sec:description-using-braid-op-1} and \cref{appendix:parity-forms-equiv}.
	\end{rem}

	Let $p_1 = p_1^{A} \colon \totcomp{\mathcal{A}}[1] \rightharpoonup \ncdfr{\mathcal{A}}[0][*+1]$
	be the natural degree one projection given by
	\begin{equation} \label{eq:def-p1-anticommute}
		p_1 \left( \sum_{i \geq 1} x_i \right) \defeq \qdr^1 \left( x_1 \right).
	\end{equation}
	Since the totalization of $\ncdfr{\mathcal{A}}[][]$
	is contractible,
	the map $p_1$ is a homotopy equivalence.\footnote{See \cref{lm:tot-ncdfr-contractible,lm:p_1-homotopy-equivalence}
	for the statement and proof when we work with $\braidop_1$ and the horizontal differential $\qdr$ instead of
$-\qdr$. One can either repeat the arguments presented there, or use the isomorphism
$\Psi$ of \cref{sec:dependence-ndf-braidop} to translate everything to our current context. See also \cref{sec:description-using-braid-op-1}.}
	Hence, we have two different chain complex models $\ncdfr{\mathcal{A}}[0][]$ and $\totcomp{\mathcal{A}}[1][]$
	for the cyclic homology of $\mathcal{A}$. Both models are functorial, and in a compatible way.
	More precisely, recall that given a Banach $\Ainf$-morphism $f \colon \mathcal{A} \rightarrow \mathcal{B}$,
	we have an induced morphism $\cindmap{f} \colon \ncdf{A} \rightarrow \ncdf{B}$, which
	commutes on the nose with both $\clie{\mu}$ and $\qdr$ (see \cref{cor:homo-cyc-func,rem:func-bi-ncdf}),
	and coincides with $\cycl{f}$ on $\ncdf{A}[0][]$.
	We have a functorial chain map
$\cycl{f} \colon \ncdfr{\mathcal{A}}[0][] \rightarrow \ncdfr{\mathcal{B}}[0][]$,
	and taking the totalization of $\cindmap{f}$, i.e., letting $\cindmap{f}$ act
	on $\totcomp{\mathcal{A}}[1][]$ by acting on each line degree separately, we get a chain
	map $\cindmap{f} \colon \totcomp{\mathcal{A}}[1][] \rightarrow \totcomp{\mathcal{B}}[1][]$
	between the total complexes. The map $p_1$ is compatible with the induced chain maps
	in the sense that $p_1^{B} \circ \cindmap{f} = \cycl{f} \circ p_1^{A}$
	(see \cref{fig:induced-maps-two-models-cyclic-complexes}).

	\begin{figure}[!htb]
	\begin{tikzcd}
		&& {\totcomp{\mathcal{A}}[1][]} && {\totcomp{\mathcal{B}}[1][]} \\
		{f \colon \mathcal{A} \rightarrow \mathcal{B}} & \longmapsto \\
		&& {\ncdfr{\mathcal{A}}[0][]} && {\ncdfr{\mathcal{B}}[0][]}
		\arrow["\cindmap{f}", from=1-3, to=1-5]
		\arrow["p_1^{A}"', harpoon, from=1-3, to=3-3]
		\arrow["p_1^{B}", harpoon, from=1-5, to=3-5]
		\arrow["\cycl{f} = \cindmap{f}"', from=3-3, to=3-5]
	\end{tikzcd}
	\caption{Compatible chain maps induced by $f$ between the two models of cyclic homology.}
	\label{fig:induced-maps-two-models-cyclic-complexes}
	\end{figure}

	In what follows, given a topologically nilpotent element $b \in \tc{A}$, we construct chains
$\Go{b}[\geq 1][\mathcal{A}] \in \totcomp{A}[1][-1]$ and $\Ho{b}[\geq 1][\mathcal{A}] \in \totcomp{A}[1][0]$,
	which are lifts of the chains $\Go{b}[0][\mathcal{A}]$ and $\Ho{b}[0][\mathcal{A}]$
	from \cref{sec:generalized-trace} in the sense that
	\begin{equation*}
		p_1 \left( \Go{b}[\geq 1][\mathcal{A}] \right) = \Go{b}[0][\mathcal{A}],
		\qquad
		p_1 \left( \Ho{b}[\geq 1][\mathcal{A}] \right) = -\Ho{b}[0][\mathcal{A}],
	\end{equation*}
	and which satisfy properties almost identical to the ones in \cref{lm:Go_0-and-Ho_0-properties}.
	We first discuss the construction in the case where the $\Ainf$-algebra $\mathcal{A} = \left( A, \mu \right)$
	is a Banach $\Ainf$-algebra over a differential graded-commutative Banach $\mathbbm{k}$-algebra
$\mathcal{R} = (R,0)$ with a \textit{zero differential}, i.e., we assume that the coderivation $\mu$ is
$R$-linear. This will allow us to work with the operator $\ccont{\mu}$, which is not defined when $\mu$
	is a generalized coderivation over a non-zero differential. The general case will be treated in
	\cref{sec:lift-general-construction}.

	\subsubsection{\texorpdfstring{The operator $\exp \left( \ccont{\mu} \right)$}{The Exponential of the Cyclic Contraction}}
	We start with a simple commutator relation calculation:
	\begin{lm} \label{lm:qdr-exp-cont-rel}
	For all $n \geq 0$, the following relation on $\ncdf{A}[][]$ holds:
	\begin{equation}
		\qdr \circ \frac{\ccont{\mu}^{n+1}}{(n+1)!} =
		\clie{\mu} \circ \frac{\ccont{\mu}^n}{n!} +
		\frac{\ccont{\mu}^{n+1}}{(n+1)!} \circ \qdr.
	\end{equation}
	\end{lm}
	\begin{proof}
	Given a graded module $W$, graded maps $f, g \in \InnEnd{W}$, and $k \geq 0$,
	we have the following identity for the graded commutator of $f$ and $g^{k+1}$:
	\begin{equation}
		\label{eq:commutator-power-identity}
		\left[ f, g^{k+1} \right] =
		\sum_{i=0}^k (-1)^{\braid{\degb{f}}{i \cdot \degb{g}}} g^i \circ \left[ f, g \right] \circ g^{k-i}.
	\end{equation}
	Hence,
	\begin{equation*}
		\begin{aligned}
			\left[ \qdr, \frac{\ccont{\mu}^{n+1}}{(n+1)!} \right]
			\stackrel{\eqref{eq:commutator-power-identity}}{=}{} &
			\frac{1}{(n+1)!} \sum_{i=0}^n (-1)^{\braid{(-1,0)}{(i,i)}_2}
			\ccont{\mu}^i \circ \left[ \qdr, \ccont{\mu} \right] \circ \ccont{\mu}^{n-i}
			\\
			\stackrel{\eqref{eq:qcontliecyc}}{=}{}               &
			\frac{1}{(n+1)!} \sum_{i=0}^n \ccont{\mu}^i \circ \clie{\mu} \circ \ccont{\mu}^{n-i}
			\\
			\stackrel{\eqref{eq:liecontcommcyc}}{=}{}            &
			\frac{n+1}{(n+1)!} \, \clie{\mu} \circ \ccont{\mu}^n = \clie{\mu} \circ \frac{\ccont{\mu}^n}{n!}.
		\end{aligned}
	\end{equation*}
	\end{proof}

	Consider the operator $\exp \left( \ccont{\mu} \right) \colon \totc{\ncdf{A}}[][][] \rightarrow \totc{\ncdf{A}}[][][]$
	defined by\footnote{There is a potential conflict of notation with the exponential
	of the formal tensor coalgebra, defined in \cref{eq:exp-tensor-coalgebra}, without any factorials. Our only use
	of the ``regular'' exponential is in defining $\exp \left( \ccont{\mu} \right)$, so this shouldn't cause any
	confusion.}
	\begin{equation} \label{eq:def-exp-ccont-mu}
		\exp \left( \ccont{\mu} \right) \left( x \right) \defeq
		\sum_{n = 0}^{\infty} \frac{\ccont{\mu}^n \left( x \right)}{n!}.
	\end{equation}
	The operator $\ccont{\mu}$ acts on elements of the form $x = \ul{a_1} \otimes \dots \otimes \ul{a_n}
\in \ncdf{A}[n][]$ having
	line degree equal to weight by adding terms of the form $\ul{\mu_0 \left( 1 \right)}$:
	\begin{equation*}
		\ccont{\mu} \left( x \right) = \sum_{i=1}^n
		\ul{a_1} \otimes \dots \otimes \ul{a_i} \otimes \ul{\mu_0 \left( 1 \right)} \otimes
		\ul{a_{i+1}} \otimes \dots \otimes \ul{a_n}.
	\end{equation*}
	Hence, for such elements we have
$\nnorm[\ccont{\mu} \left( x \right)] \leq \nnorm[\mu_0 \left( 1 \right)] \cdot \nnorm[x]$.
	Since repeated applications of $\ccont{\mu}$ increase the line degree of elements and
	either preserve the weight filtration or add terms of the form $\ul{\mu_0 \left( 1 \right)}$,
	an argument similar to that in \cref{lm:S-almost-nilpotent,lm:S-split-based-on-mu_0,lm:S-norm-estimate}
	shows that the infinite sum in \cref{eq:def-exp-ccont-mu}
	converges. As $\ccont{\mu}$ has bi-degree $(1,1)$, the operators $\ccont{\mu}^n$ have total degree zero
	on the totalization, so $\exp \left( \ccont{\mu} \right)$ has degree zero and satisfies
$\nnorm[\exp \left( \ccont{\mu} \right)] \leq \max \Set{1, \nnorm[\mu]} \leq 1$.

	\begin{lm}
	Let $x \in \totc{\ncdf{A}[][]}[][][]$. Then we have
	\begin{equation} \label{eq:basic-identity}
		\left( \left( \clie{\mu} - \qdr \right) \circ \exp (\ccont{\mu}) \right) \left( x \right) =
		\left( \exp \left( \ccont{\mu} \right) \circ \left( -\qdr \right) \right) \left( x \right).
	\end{equation}
	Hence, on the level of the total complexes $\totcomp{\mathcal{A}}[k][] = \totc{\ncdf{\mathcal{A}}[][]}[][\geq k][]$,
	we have the identity
	\begin{equation} \label{eq:basic-identity-tot-comp-k}
		D_{\mu}^{\geq k} \left( \left( \pi_{\geq k} \circ \exp(\ccont{\mu}) \right) (x) \right) =
		-\left( \pi_{\geq k} \circ \exp \left( \ccont{\mu} \right) \right) \left( \qdr(x) \right),
	\end{equation}
	where $\pi_{\geq k} \colon \totc{\ncdf{A}[][]}[][][] \rightarrow \totc{\ncdf{A}[][]}[][\geq k][]$
	denotes the natural projection.
	\end{lm}
	\begin{proof}
	By \cref{lm:qdr-exp-cont-rel}, we have
	\begin{equation*}
		\begin{aligned}
			\clie{\mu} \circ \exp(\ccont{\mu}) & =
			\sum_{n=0}^{\infty} \clie{\mu} \circ \frac{\ccont{\mu}^n}{n!}
			\\
			                                   & =
			\sum_{n=0}^{\infty} \qdr \circ \frac{\ccont{\mu}^{n+1}}{(n+1)!} -
			\sum_{n=0}^{\infty} \frac{\ccont{\mu}^{n+1}}{(n+1)!} \circ \qdr
			\\
			                                   & =
			\qdr \circ \exp(\ccont{\mu}) - \qdr - \left( \exp \left( \ccont{\mu} \right) \circ \qdr - \qdr \right)
			\\
			                                   & = \qdr \circ \exp(\ccont{\mu}) +  \exp \left( \ccont{\mu} \right) \circ \left( -\qdr \right).
		\end{aligned}
	\end{equation*}
	The calculation is formal but makes sense when applied to $x \in \totc{\ncdf{A}[][]}[][][]$ since all the
	infinite sums converge, and gives us \cref{eq:basic-identity}. Applying the projection
$\pi_{\geq k}$ and using the fact that $\pi_{\geq k}$ is a chain map, i.e.,
$\pi_{\geq k} \circ \left( \clie{\mu} - \qdr \right) = D_{\mu}^{\geq k} \circ \pi_{\geq k}$,
	we obtain \cref{eq:basic-identity-tot-comp-k}.
	\end{proof}

	Here, and in what follows, given $k \geq 1$, we will generically denote by
$\pi_{\geq k}$
	the natural projection onto the truncated total complex
$\totcomp{\mathcal{A}}[k][] = \tot_{\geq k}(\ncdf{A}[][])$ from \textit{any}
	total complex containing it
	(e.g., from the full total complex
$\tot(\ncdf{A}[][])$ or from $\totcomp{\mathcal{A}}[j][]$ for
$j < k$).\footnote{Note that when $j < k$, $\totcomp{\mathcal{A}}[k][]$ is not
	a \textit{subcomplex} of $\totcomp{\mathcal{A}}[j][]$, but rather a quotient.}

	\subsubsection{Construction and Properties in the \texorpdfstring{$R$}{R}-linear Case}
	\label{sec:lift-construction-R-linear}

	Let $b \in \tc{A}$ be a topologically nilpotent element. Consider the chain
	\begin{equation} \label{eq:G_1-def}
		\G{b}[1] = \G{b}[1][\mathcal{A}] \defeq \sum_{i=0}^{\infty} \frac{1}{i+1} \ul{b} \otimes b^{\otimes i}
		\in \ncdf{A}[1][0]
	\end{equation}
	which is a lift to $\ncdf{A}[1][0]$ of the cyclic exponential $\Go{b}[0] \in \ncdfr{A}[0][0]$ given by \cref{eq:Go_0-def},
	in the sense that
	\begin{equation} \label{eq:qdr-G-1-is-Go-0}
		\qdr \left( \G{b}[1] \right) = \Go{b}[0].
	\end{equation}
	Let us define
	\begin{align}
		\Go{b}[\geq 1] = \Go{b}[\geq 1][\mathcal{A}]
		 & \defeq
		\exp \left( \ccont{\mu} \right) \left( \G{b}[1] \right) \in \totcomp{A}[1][-1].
		\label{eq:Go-geq-1-def}
	\end{align}
	The chain $\Go{b}[\geq 1]$ is a lift of the truncated cyclic exponential $\Go{b}[0]$ in the sense that
$p_1 \left(  \Go{b}[\geq 1] \right) = \Go{b}[0]$. We can represent $\Go{b}[\geq 1]$ as a sum
$\Go{b}[\geq 1] = \sum_{n=1}^{\infty} \G{b}[n]$, where each cyclic codifferential form
	\begin{equation*}
		\G{b}[n+1] = \G{b}[n+1][\mathcal{A}] =
		\frac{1}{n!} \ccont{\mu}^{n} \left( \G{b}[1] \right) \in \ncdf{A}[n + 1][n]
	\end{equation*}
	is of line degree $n + 1$. Explicitly, we have:
	\begin{lm} \label{lm:G_n-explicit-formula}
	The form $\G{b}[n] \in \ncdf{A}[n][n - 1]$ is given by
	\begin{equation} \label{eq:G_n-explicit-formula}
		\begin{aligned}
			\G{b}[n] = \sum_{\substack{j_1, \dots, j_n = 0 \\ i_2, \dots, i_n = 0}}^{\infty}
			 & \frac{1}{1 + \sum_{r = 1}^n j_r + \sum_{r=2}^{n} i_r}
			\\
			 & \quad
			\ul{b} \otimes b^{\otimes j_1} \otimes \ul{ \mu_{i_2} \left( b^{\otimes i_2} \right)}
			\otimes b^{\otimes j_2} \otimes \dots \otimes \ul{ \mu_{i_n} \left( b^{\otimes i_n} \right) }
			\otimes b^{\otimes j_n}.
		\end{aligned}
	\end{equation}
	\end{lm}
	\begin{proof}
	The form $\G{b}[2] \in \ncdf{A}[2][1]$ is given by
	\begin{equation*}
		\begin{aligned}
			\G{b}[2]
			\stackrel{\phantom{\eqref{eq:ccont-mu-on-ncdf-geq-1}}}{=}{} &
			\ccont{\mu} (\G{b}[1])
			=
			\ccont{\mu} \left( \sum_{i=0}^{\infty} \frac{1}{i+1} \ul{b} \otimes b^{\otimes i} \right)
			\\
			\stackrel{\eqref{eq:ccont-mu-on-ncdf-geq-1}}{=}{}           &
			\sum_{i_1,i_2,i_3=0}^{\infty} \frac{1}{i_1 + i_2 + i_3 + 1} \,
			\ul{b} \otimes b^{\otimes i_1} \otimes
			\ul{ \mu_{i_2} \left( b^{\otimes i_2} \right)} \otimes b^{\otimes i_3}.
		\end{aligned}
	\end{equation*}
	Similarly, the form $\G{b}[3] \in \ncdf{A}[3][2]$ is given by
	\begin{align*}
		\MoveEqLeft
		\G{b}[3]
		\stackrel{\phantom{\eqref{eq:ccont-mu-on-ncdf-geq-1}}}{=}{}
		\frac{1}{2!} \ccont{\mu}^2 \left(\G{b}[1] \right) =
		\frac{1}{2} \ccont{\mu} \left( \G{b}[2] \right) \notag
		\\
		\stackrel{\phantom{\eqref{eq:ccont-mu-on-ncdf-geq-1}}}{=}{} &
		\frac{1}{2} \ccont{\mu} \left(
		\sum_{i_1,i_2,i_3=0}^{\infty} \frac{1}{i_1 + i_2 + i_3 + 1} \,
		\ul{b} \otimes b^{\otimes i_1} \otimes
		\ul{ \mu_{i_2} \left( b^{\otimes i_2} \right)} \otimes b^{\otimes i_3}
		\right)
		\\
		\stackrel{\eqref{eq:ccont-mu-on-ncdf-geq-1}}{=}{}           &
		\frac{1}{2}  \sum_{\substack{j_1,j_2,j_3=0 \\ i_2,i_3=0}}^{\infty}
		\frac{1}{j_1+j_2+j_3+i_2+i_3+1} \,
		\ul{b} \otimes b^{\otimes j_1} \otimes \ul{\mu_{j_2} \left( b^{\otimes j_2} \right)}
		\otimes b^{\otimes j_3} \otimes
		\ul{ \mu_{i_2} \left( b^{\otimes i_2} \right)} \otimes b^{\otimes i_3}
		\notag
		\\
		                                                            & +
		\frac{1}{2} \sum_{\substack{i_1,i_2=0 \\j_1,j_2,j_3=0}}^{\infty}
		\frac{1}{i_1+i_2+j_1+j_2+j_3+1} \,
		\ul{b} \otimes b^{\otimes i_1} \otimes \ul{ \mu_{i_2} \left( b^{\otimes i_2} \right)} \otimes
		b^{\otimes j_1} \otimes \ul{\mu_{j_2} \left( b^{\otimes j_2} \right)} \otimes b^{\otimes j_3}
		\notag
		\\
		\stackrel{\phantom{\eqref{eq:ccont-mu-on-ncdf-geq-1}}}{=}{} &
		\sum_{\substack{j_1,j_2,j_3=0\\i_2,i_3=0}}^{\infty}
		\frac{1}{j_1+j_2+j_3+i_2+i_3+1} \,
		\ul{b} \otimes b^{\otimes j_1} \otimes \ul{\mu_{i_2} \left( b^{\otimes i_2} \right)}
		\otimes b^{\otimes j_2} \otimes
		\ul{ \mu_{i_3} \left( b^{\otimes i_3} \right)} \otimes b^{\otimes j_3}.
		\notag
	\end{align*}
	More generally, an induction argument using
$\G{b}[n+1] = \frac{1}{n} \ccont{\mu} \left( \G{b}[n] \right)$
	and \cref{eq:ccont-mu-on-ncdf-geq-1} shows that \cref{eq:G_n-explicit-formula} holds for
	all $n \geq 1$.
	\end{proof}

	Consider also the chain
	\begin{equation} \label{eq:Ho-geq-1-def}
		\Ho{b}[\geq 1] = \Ho{b}[\geq 1][\mathcal{A}] \defeq
		- \left( \pi_{\geq 1} \circ \exp \left( \ccont{\mu} \right) \right) \left( \G{b}[0] \right)
		\in \totcomp{A}[1][0].
	\end{equation}
	Then we have:
	\begin{lm} \label{lm:Go-geq-1-differential}
	The differential of $\Go{b}[\geq 1]$ is given by
	\begin{equation} \label{eq:D-Go-geq-1-linear}
		D_{\mathcal{A}} \left( \Go{b}[\geq 1][\mathcal{A}] \right) =
		\Ho{b}[\geq 1][\mathcal{A}] - \Ho{0}[\geq 1][\mathcal{A}].
	\end{equation}
	\end{lm}
	\begin{proof}
	Using \cref{eq:basic-identity-tot-comp-k}, and the fact that $\G{b}[0] = 1 + \Go{b}[0]$, with $\Go{0}[0] = 0$,
	we see that
	\begin{equation*}
		\begin{aligned}
			D_{\mathcal{A}} \left( \Go{b}[\geq 1][\mathcal{A}] \right) \!
			\eqwithref[eq:Go-geq-1-def]              &
			D_{\mu}^{\geq 1} \left(
			\left( \pi_{\geq 1} \circ \exp \left( \ccont{\mu} \right) \right) \left( \G{b}[1] \right)
			\right)
			\\
			\eqwithref[eq:basic-identity-tot-comp-k] &
			- \left( \pi_{\geq 1} \circ \exp \left( \ccont{\mu} \right) \right) \left(
			\qdr \left( \G{b}[1] \right)
			\right)
			\\
			\eqwithref[eq:qdr-G-1-is-Go-0]           &
			- \left( \pi_{\geq 1} \circ \exp \left( \ccont{\mu} \right) \right) \left(
			\Go{b}[0]
			\right)
			\\
			\eqwithref                               &
			- \left( \pi_{\geq 1} \circ \exp \left( \ccont{\mu} \right) \right) \left(
			\G{b}[0] - 1
			\right)
			\\
			\eqwithref                               &
			- \left( \pi_{\geq 1} \circ \exp \left( \ccont{\mu} \right) \right) \left(
			\G{b}[0] \right)
			-
			\left( -
			\left( \pi_{\geq 1} \circ \exp \left( \ccont{\mu} \right) \right) \left(
			\G{0}[0]
			\right)
			\right)
			\\
			\eqwithref[eq:Ho-geq-1-def]              &
			\Ho{b}[\geq 1][\mathcal{A}] - \Ho{0}[\geq 1][\mathcal{A}].
		\end{aligned}
	\end{equation*}
	\end{proof}

	We can represent $\Ho{b}[\geq 1]$ as a sum $\Ho{b}[\geq 1] = \sum_{n = 1}^{\infty} \H{b}[n]$, where
	each cyclic codifferential form
	\begin{equation*}
		\H{b}[n] = \H{b}[n][\mathcal{A}] =
		-\frac{1}{n!} \ccont{\mu}^n \left( \G{b}[0] \right) \in \ncdf{A}[n][n]
	\end{equation*}
	is of line degree $n$. Then we have the following explicit formula for $\H{b}[n]$:
	\begin{lm} \label{lm:H_n-formula}
	The form $\H{b}[n]$ is given by
	\begin{equation} \label{eq:H_n(b)-formula}
		\begin{aligned}
			\H{b}[n] & =
			-\frac{1}{n} \left( \ul{\corest{\mu} \left( \Exp{b} \right)} \otimes \Exp{b} \right)^{\otimes n}
			\\
			         & =
			-\frac{1}{n} \left(
			\sum_{\substack{i_1, \dots, i_n = 0 \\ j_1, \dots, j_n = 0}}^{\infty}
			\ul{ \mu_{i_1} \left( b^{\otimes i_1} \right) } \otimes b^{\otimes j_1} \otimes \dots \otimes
			\ul{ \mu_{i_n} \left( b^{\otimes i_n} \right) } \otimes b^{\otimes j_n}
			\right).
		\end{aligned}
	\end{equation}
	\end{lm}
	\begin{proof}
	Given $k \geq 0$, we can make a calculation similar to the one appearing in
	the proof of \cref{lm:coder-cycl-exp-identity}, with $\ccont{\mu}$ replacing $\cycl{\mu}$:
	\begin{equation*}
		\begin{aligned}
			\ccont{\mu} \left( b^{\otimes (k + 1)} \right)
			\stackrel{\eqref{eq:ccont-mu-on-ncdf-0}}{=}{}           &
			\sum_{k_1 + k_2 + k_3 = k}
			b \otimes b^{\otimes k_1} \otimes \ul{ \mu_{k_2} \left( b^{\otimes k_2} \right) }
			\otimes b^{\otimes k_3} +
			\ul{ \mu_{k_3 + 1 + k_1} \left( b^{\otimes \left( k_3 + 1 + k_1 \right)} \right) } \otimes
			b^{\otimes k_2}
			\\
			\stackrel{\phantom{\eqref{eq:ccont-mu-on-ncdf-0}}}{=}{} &
			\sum_{k_1 + k_2 + k_3 = k}
			\ul{ \mu_{k_2} \left( b^{\otimes k_2} \right) } \otimes
			b^{\otimes \left( k_3 + 1 + k_1 \right)}
			+
			\ul{ \mu_{k_3 + 1 + k_1} \left( b^{\otimes \left( k_3 + 1 + k_1 \right)} \right) }
			\otimes b^{\otimes k_2}
			\\
			\stackrel{\phantom{\eqref{eq:ccont-mu-on-ncdf-0}}}{=}{} &
			\sum_{\substack{i + j = k + 1                           \\ i \geq 0, j \geq 1}}
			j \cdot \left(
			\ul{ \mu_{i} \left( b^{\otimes i} \right) } \otimes b^{\otimes j} +
			\ul{ \mu_j \left( b^{\otimes j} \right) } \otimes b^{\otimes i}
			\right)
			\\
			\stackrel{\phantom{\eqref{eq:ccont-mu-on-ncdf-0}}}{=}{} &
			\sum_{\substack{i + j = k + 1                           \\ i, j \geq 0}}
			j \cdot \left(
			\ul{ \mu_{i} \left( b^{\otimes i} \right) } \otimes b^{\otimes j} +
			\ul{ \mu_j \left( b^{\otimes j} \right) } \otimes b^{\otimes i}
			\right)
			\\
			\stackrel{\phantom{\eqref{eq:ccont-mu-on-ncdf-0}}}{=}{} &
			\left( k + 1 \right) \cdot
			\sum_{\substack{i + j = k + 1                           \\ i, j \geq 0}}
			\ul{ \mu_{i} \left( b^{\otimes i} \right) } \otimes b^{\otimes j}
		\end{aligned}
	\end{equation*}
	and hence
	\begin{equation*}
		\begin{aligned}
			- \H{b}[1] & =
			\ccont{\mu} \left( \G{b}[0] \right)
			=
			\ccont{\mu} \left( 1 \right) +
			\sum_{k = 0}^{\infty} \ccont{\mu} \left( \frac{b^{\otimes (k + 1)}}{k + 1} \right)
			=
			\ul{	\mu_0 \left( 1 \right)} + \sum_{\substack{i + j = k + 1 \\ i, j, k \geq 0}}
			\ul{ \mu_{i} \left( b^{\otimes i} \right) } \otimes b^{\otimes j}
			\\
			           & =
			\sum_{\substack{i + j = k                                   \\ i, j, k \geq 0}}
			\ul{ \mu_{i} \left( b^{\otimes i} \right) } \otimes b^{\otimes j}
			= \ul{\corest{\mu} \left( \Exp{b} \right)} \otimes \Exp{b}
		\end{aligned}
	\end{equation*}
	which shows the lemma for $n = 1$. Assuming the lemma for $n$, and using
	the explicit formula for the action of $\ccont{\mu}$ given in
	\cref{sec:nc-d-calc-cont}, we have
	\begin{equation*}
		\begin{aligned}
			\MoveEqLeft
			\ccont{\mu} \left( \H{b}[n] \right)
			\stackrel{\phantom{\eqref{eq:ccont-mu-on-ncdf-geq-1}}}{=}{}
			-\frac{1}{n} \ccont{\mu} \left( \left(
			\ul{\corest{\mu} \left( \Exp{b} \right)} \otimes \Exp{b} \right)^{\otimes n} \right)
			\\
			\stackrel{\eqref{eq:ccont-mu-on-ncdf-geq-1}}{=}{}           &
			- \frac{1}{n}
			\sum_{i=1}^n
			\left( \ul{\corest{\mu} \left( \Exp{b} \right)} \otimes \Exp{b} \right)^{\otimes
				                                                                \left( i - 1 \right)}
			\otimes
			\left( \ul{\corest{\mu} \left( \Exp{b} \right)} \otimes \Exp{b} \otimes
			\ul{\corest{\mu} \left( \Exp{b} \right)} \otimes \Exp{b} \right)
			\otimes
			\left( \ul{\corest{\mu} \left( \Exp{b} \right)} \otimes \Exp{b} \right)^{\otimes
				                                                                \left( n - i \right)}
			\\
			\stackrel{\phantom{\eqref{eq:ccont-mu-on-ncdf-geq-1}}}{=}{} &
			-
			\left( \ul{\corest{\mu} \left( \Exp{b} \right)} \otimes \Exp{b} \right)^{\otimes \left( n + 1 \right)}
		\end{aligned}
	\end{equation*}
	and hence
	\begin{equation*}
		\H{b}[n+1] = \frac{1}{n + 1} \cdot \ccont{\mu} \left( \H{b}[n] \right) =
		- \frac{1}{n+1} \cdot
		\left( \ul{\corest{\mu} \left( \Exp{b} \right)} \otimes \Exp{b} \right)^{\otimes \left( n + 1 \right)}.
	\end{equation*}
	\end{proof}

	\begin{rem}
	Note that the first term $\H{b}[1] = - \ul{\corest{\mu} \left( \Exp{b} \right)} \otimes \Exp{b}$
	of $\Ho{b}[\geq 1]$ is a lift to $\ncdf{A}[1][1]$ of
$- \Ho{b}[0] = - \corest{\mu} \left( \Exp{b} \right) \otimes \Exp{b} \in \ncdfr{A}[0][1]$,
	so that we have $p_1 \left( \Ho{b}[\geq 1] \right) = - \Ho{b}[0]$.
	\end{rem}

	\begin{lm}
	Given $b \in \tc{A}$, we have
	\begin{equation} \label{eq:D-Ho-geq-1-linear}
		D_{\mathcal{A}} \left( \Ho{b}[\geq 1][\mathcal{A}] \right) = 0.
	\end{equation}
	\end{lm}
	\begin{proof}
		We have
		\begin{equation*}
			\begin{aligned}
				D_{\mathcal{A}} \left( \Ho{b}[\geq 1][\mathcal{A}] \right) \!
				\eqwithref[eq:Ho-geq-1-def]              &
				D_{\mu}^{\geq 1} \left( -
				\left( \pi_{\geq 1} \circ \exp \left( \ccont{\mu} \right) \right) \left( \G{b}[0] \right)
				\right)
				\\
				\eqwithref[eq:basic-identity-tot-comp-k] &
				\left( \pi_{\geq 1} \circ \exp \left( \ccont{\mu} \right) \right) \left(
				\qdr \left( \G{b}[0] \right)
				\right)
				=
				\left( \pi_{\geq 1} \circ \exp \left( \ccont{\mu} \right) \right) \left( 0 \right) = 0.
			\end{aligned}
		\end{equation*}
	\end{proof}

	Next, we discuss the behaviour of $\Go{b}[\geq 1]$ and $\Ho{b}[\geq 1]$ under $\Ainf$-morphisms.
	Let $\mathcal{B} = \left( B,\nu \right)$ be another Banach $\Ainf$-algebra over
	a differential graded-commutative Banach $\mathbbm{k}$-algebra $\mathcal{S} = \left( S, 0 \right)$,
	again with a zero differential.

	\begin{lm} \label{lm:Go-geq-1-pseudo-func}
	Given $b \in \tc{A}$, we have the identity
	\begin{equation} \label{eq:cindmap-f-Go-geq-1-R-linear}
		\cindmap{f} \left( \Go{b}[\geq 1][\mathcal{A}] \right) =
		\Go{ \mcfunc{f} \left( b \right) }[\geq 1][\mathcal{B}] -
		\Go{ \mcfunc{f} \left( 0 \right) }[\geq 1][\mathcal{B}] +
		D_{\mathcal{B}} \left( R \left( b; f \right) \right)
	\end{equation}
	for some $R \left( b; f \right) \in \totcomp{B}[1][-2]$.
	When $b = 0$ or when $f$ is strict, one can choose $R \left( b; f \right) = 0$.
	\end{lm}
	\begin{proof}
	Consider first the situation for the form $\G{b}[1][\mathcal{A}]$. We have
	\begin{equation} \label{eq:Gb-1-func-modulo-q}
		\begin{aligned}
			\qdr \left( \cindmap{f} \left( \G{b}[1][\mathcal{A}] \right) \right)
			\stackrel{\eqref{eq:func-qdr-cyc}}{=}{}    &
			\cindmap{f} \left( \qdr \left( \G{b}[1][\mathcal{A}] \right) \right)
			\\
			\stackrel{\eqref{eq:qdr-G-1-is-Go-0}}{=}{} &
			\cindmap{f} \left( \Go{b}[0][\mathcal{A}] \right)
			\\
			\stackrel{\eqref{eq:cycl-f-Go-0}}{=}{}     &
			\Go{\mcfunc{f} \left( b \right)}[0][\mathcal{B}] -
			\Go{\mcfunc{f} \left( 0 \right)}[0][\mathcal{B}]
			\\
			\stackrel{\eqref{eq:qdr-G-1-is-Go-0}}{=}{} &
			\qdr \left(
			\G{\mcfunc{f} \left( b \right)}[1][\mathcal{B}] - \G{\mcfunc{f} \left( 0 \right)}[1][\mathcal{B}]
			\right).
		\end{aligned}
	\end{equation}
	Let us set
$x \defeq \left(
\G{ \mcfunc{f} \left( b \right) }[1][\mathcal{B}] - \G{ \mcfunc{f} \left( 0 \right) }[1][\mathcal{B}]
\right) - \cindmap{f} \left( \G{b}[1][\mathcal{A}] \right) \in \ncdf{B}[1][0]$.
	By \cref{eq:Gb-1-func-modulo-q}, we have $\qdr \left( x \right) = 0$,
	so by the formal Poincar\'{e} \cref{lm:formal-poincare-ncdfr}, there exists some
$R_2 \left( b; f \right) \in \ncdf{B}[2][0]$ such that $\qdr \left( R_2 \left( b; f \right) \right) = x$,
	and thus
	\begin{equation} \label{eq:cindmap-f-G-b-1}
		\cindmap{f} \left( \G{b}[1][\mathcal{A}] \right) =
		\G{\mcfunc{f} \left( b \right)}[1][\mathcal{B}] - \G{\mcfunc{f} \left( 0 \right)}[1][\mathcal{B}] -
		\qdr \left( R_2 \left( b; f \right)  \right).
	\end{equation}
	As a particular choice of $R_2 \left( b; f \right)$, we can take
$R_2 \left( b; f \right) = h_{\dr} \left( x \right)$.
	Note that we have $h_{\dr} \left( \Go{b}[0] \right) = \G{b}[1]$, and since $h_{\dr}^2 = 0$
	(\cref{item:hdr-square-zero} of \cref{lm:formal-poincare-ncdfr}), we have
	\begin{equation} \label{eq:R_2-b-f}
		R_2 \left( b; f \right) = - h_{\dr} \left( \cindmap{f} \left( \G{b}[1] \right) \right).
	\end{equation}
	When $b = 0$, then $\G{b}[1] = 0$, and so $x = 0$ and $R_2 \left( b; f \right) = 0$. When $f$ is strict,
	we actually have $\cindmap{f} \left( \G{b}[1] \right) = \G{ \mcfunc{f} \left( b \right) }[1]$, so
	in this case we also have $x = 0$ and $R_2 \left( b; f \right) = 0$.

	Now consider the chain $\Go{b}[\geq 1][\mathcal{A}]$. We have
	\begin{equation*} 
		\begin{aligned}
			\cindmap{f} \left( \Go{b}[\geq 1][\mathcal{A}] \right) \!
			\eqwithref[eq:Go-geq-1-def]              &
			\cindmap{f} \left( \exp \left( \ccont{\mu} \right) \left( \G{b}[1][\mathcal{A}] \right) \right)
			\\
			\eqwithref[eq:func-cont-cyc]             &
			\exp \left( \ccont{\nu} \right) \left( \cindmap{f} \left( \G{b}[1][\mathcal{A}] \right) \right)
			\\
			\eqwithref[eq:cindmap-f-G-b-1]           &
			\exp \left( \ccont{\nu} \right) \left(
			\G{ \mcfunc{f} \left( b \right) }[1][\mathcal{B}] -
			\G{ \mcfunc{f} \left( 0 \right) }[1][\mathcal{B}] -
			\qdr \left( R_2 \left( b; f \right) \right)
			\right)
			\\
			\eqwithref[eq:Go-geq-1-def]              &
			\Go{ \mcfunc{f} \left( b \right) }[\geq 1][\mathcal{B}] -
			\Go{ \mcfunc{f} \left( 0 \right) }[\geq 1][\mathcal{B}] +
			\exp \left( \ccont{\nu} \right) \left( -\qdr \left( R_2 \left( b; f \right) \right) \right)
			\\
			\eqwithref[eq:basic-identity-tot-comp-k] &
			\Go{ \mcfunc{f} \left( b \right) }[\geq 1][\mathcal{B}] -
			\Go{ \mcfunc{f} \left( 0 \right) }[\geq 1][\mathcal{B}] +
			D_{\nu} \left( \exp \left( \ccont{\nu} \right) \left( R_2 \left( b; f \right) \right) \right)
			\\
			\eqwithref                               &
			\Go{ \mcfunc{f} \left( b \right) }[\geq 1][\mathcal{B}] -
			\Go{ \mcfunc{f} \left( 0 \right) }[\geq 1][\mathcal{B}] +
			D_{\mathcal{B}} \left( R \left( b; f \right) \right),
		\end{aligned}
	\end{equation*}
	when $R \left( b; f \right) \defeq \exp \left( \ccont{\nu} \right) \left( R_2 \left( b; f \right) \right)$,
	as required.
	\end{proof}

	\begin{rem}
	Given $b \in \tc{A}$, \cref{lem:func-mc-cyclic} implies the identity
	\begin{equation*}
		\cycl{f} \left( \Go{b}[0][\mathcal{A}] \right) =
		\Go{ \mcfunc{f} \left( b \right) }[0][\mathcal{B}] -
		\Go{ \mcfunc{f} \left( 0 \right) }[0][\mathcal{B}]
		\mod \Im \left( \idd - \t \right)
	\end{equation*}
	in the tensor module $\ndfr{B}[0][] = \tensr{B}$. From the proof of \cref{lem:func-mc-cyclic}, it is clear that
	if $f$ is strict then we actually have the identity
	\begin{equation*}
		\cycl{f} \left( \Go{b}[0][\mathcal{A}] \right) =
		\Go{ \mcfunc{f} \left( b \right) }[0][\mathcal{B}],
	\end{equation*}
	and not only modulo $\Im \left( \idd - \t \right)$.	Similarly, the proof of \cref{lm:Go-geq-1-pseudo-func}
	shows that we have the identity
	\begin{equation*}
		\cindmap{f} \left( \G{b}[1][\mathcal{A}] \right) = \G{ \mcfunc{f} \left( b \right) }[1][\mathcal{B}]
		- \G{ \mcfunc{f} \left( 0 \right) }[1][\mathcal{B}] \mod \Im \left( \qdr^2 \right)
	\end{equation*}
	in $\ncdf{B}[1][] \cong \tensr{B}$, and, if $f$ is strict, we have the identities
	\begin{equation*}
		\cindmap{f} \left( \G{b}[n][\mathcal{A}] \right) =
		\G{ \mcfunc{f} \left( b \right) }[n][\mathcal{B}]
	\end{equation*}
	in $\ncdf{B}[n][]$ for all $n \geq 1$.
	\end{rem}

	\begin{lm} \label{lm:Ho-geq-1-func}
	Given $b \in \tc{A}$, we have the identity
	\begin{equation} \label{eq:cindmap-f-Ho-geq-1-R-linear}
		\cindmap{f} \left( \Ho{b}[\geq 1][\mathcal{A}] \right) =
		\Ho{ \mcfunc{f} \left( b \right) }[\geq 1][\mathcal{B}].
	\end{equation}
	\end{lm}
	\begin{proof}
		We have
		\begin{equation*}
			\begin{aligned}
				\cindmap{f} \left( \Ho{b}[\geq 1][\mathcal{A}] \right) \!
				\eqwithref[eq:Ho-geq-1-def]  &
				- \cindmap{f} \left( \left( \pi_{\geq 1} \circ \exp \left( \ccont{\mu} \right) \right) \left( \G{b}[0] \right) \right)
				\\
				\eqwithref                   &
				- \pi_{\geq 1} \left( \left( \cindmap{f} \circ \exp \left( \ccont{\mu} \right) \right) \left( \G{b}[0] \right) \right)
				\\
				\eqwithref[eq:func-cont-cyc] &
				- \left( \pi_{\geq 1} \circ \exp \left( \ccont{\nu} \right) \right) \left( \cindmap{f} \left( \G{b}[0] \right) \right)
				\\
				\eqwithref[eq:cycl-f-G0]     &
				- \left( \pi_{\geq 1} \circ \exp \left( \ccont{\nu} \right) \right) \left( \G{ \mcfunc{f} \left( b \right) }[0][\mathcal{B}] \right)
				\\
				\eqwithref[eq:Ho-geq-1-def]  &
				\Ho{ \mcfunc{f} \left( b \right) }[\geq 1][\mathcal{B}].
			\end{aligned}
		\end{equation*}
	\end{proof}

	\subsubsection{General Construction} \label{sec:lift-general-construction}
	We wish to extend the definitions of $\Go{\cdot}[\geq 1][\mathcal{A}]$ and $\Ho{\cdot}[\geq 1][\mathcal{A}]$
	from \cref{sec:lift-construction-R-linear}
	so that they will be defined for a Banach $\Ainf$-algebra $\mathcal{A} = \left( A,\mu \right)$
	over an arbitrary differential graded-commutative Banach $\mathbbm{k}$-algebra
$\mathcal{R} = \left( R,d_R \right)$, \textit{and} so that
	all the results from \cref{sec:lift-construction-R-linear} will continue to hold in the extended setting.
	However, when $\mu$ is a derivation over $d_R$ and not $R$-linear, the operator $\ccont{\mu}$ which appears
	in the definitions of $\Go{\cdot}[\geq 1][\mathcal{A}]$ and $\Ho{\cdot}[\geq 1][\mathcal{A}]$
	is not defined on $\ncdf{A/R}[][]$.

	Note that the operator $\ccont{\mu}$ appears only in the definitions of
$\Go{\cdot}[\geq 1][\mathcal{A}]$ and $\Ho{\cdot}[\geq 1][\mathcal{A}]$, and in the proofs of the results,
	but not in their statements. Instead of defining the chains $\Go{\cdot}[\geq 1][\mathcal{A}]$ and
$\Ho{\cdot}[\geq 1][\mathcal{A}]$ using $\ccont{\mu}$, we can use the explicit formulas
	\eqref{eq:G_n-explicit-formula} and \eqref{eq:H_n(b)-formula} as definitions, and then
	repeat the calculations of \cref{sec:lift-construction-R-linear} using only the operators
$\clie{\mu}$ and $\qdr$ which make sense in the general setting. Alternatively, and this is
	the path we choose, we can circumvent the problem using restriction of scalars and consider $\mathcal{A}$
	as an $\Ainf$-algebra over a smaller subalgebra such that
	the restriction of $\mu$ will become linear.
	The natural choice for such a subalgebra is the subalgebra
$\clsub{R} \defeq \Set{r \in R}[d_R \left( r \right) = 0]$ of cocycles of $\mathcal{R}$.
	Consider the following restriction of scalars diagram
	\begin{equation*}
		\begin{tikzcd}
			{\left( \tens{i^{*} \left( A \right)}[\clsub{R}], i^{*} \left( \mu \right) \right)} &&
			{\left( \tens{A}[R], \mu \right)} \\
			{\left( \clsub{R}, 0 \right)} && {\left( R, d_R \right),}
			\arrow["{\resover{i}}", from=1-1, to=1-3]
			\arrow["\varepsilon_{\clsub{R}}", from=1-1, to=2-1]
			\arrow["\varepsilon_R"', from=1-3, to=2-3]
			\arrow["i", from=2-1, to=2-3]
		\end{tikzcd}
	\end{equation*}
	where $i \colon \left( \clsub{R}, 0 \right) \rightarrow \left( R, d_R \right)$ is the natural inclusion,
$\resover{i}$ is the canonical morphism over $i$ (see \cref{sec:scalar-restriction-formal-tensor-coalgebras}),
	and the maps $\varepsilon$ are the counit maps.
	We will denote the restriction $i^{*} \left( \mu \right)$ by $\clsub{\mu}$ and the pair
$i^{*} \left( \mathcal{A} \right) = \left( i^{*} \left( A \right), \clsub{\mu} \right)$,
	which is an $\Ainf$-algebra over $\clsub{\mathcal{R}} = \left( \clsub{R}, 0 \right)$, by
$\clsub{\mathcal{A}}$.

	\begin{dfn}
	Let $\mathcal{A} = \left( A, \mu \right)$ be a Banach $\Ainf$-algebra over
	a differential graded-commutative Banach $\mathbbm{k}$-algebra $\mathcal{R} = \left( R,d_R \right)$.
	Given a topologically nilpotent element $b \in \tc{A}$, define\footnote{Here, we use the standard
	abuse of notation and think of $b$ both as an element of $A$ and of $i^{*} \left( A \right)$.}
	\begin{align}
		\Go{b}[\geq 1][\mathcal{A}] & \defeq
		\cindmap{\resover{i}} \left( \Go{b}[\geq 1][\clsub{\mathcal{A}}] \right) \in
		\totcomp{\mathcal{A}}[1][-1],
		\label{eq:Go-geq-1-def-general}
		\\
		\Ho{b}[\geq 1][\mathcal{A}] & \defeq
		\cindmap{\resover{i}} \left( \Ho{b}[\geq 1][\clsub{\mathcal{A}}] \right) \in
		\totcomp{\mathcal{A}}[1][0].
		\label{eq:Ho-geq-1-def-general}
	\end{align}
	\end{dfn}
	Note that when $d_R = 0$, we have $\clsub{\mathcal{R}} = \mathcal{R}$, $\clsub{\mathcal{A}} = \mathcal{A}$, and
$\cindmap{\resover{i}} \colon \totcomp{\clsub{\mathcal{A}}}[1][] \rightarrow \totcomp{\mathcal{A}}[1][]$
	is the identity map, so this indeed extends the previous definitions. Note also that since $\cindmap{\resover{i}}$
	is given by
	\begin{gather*}
		\cindmap{\resover{i}} \left(
		\ul{a_1} \otimes_{\clsub{R}} a_1^1 \otimes_{\clsub{R}} \dots \otimes_{\clsub{R}}
		a_1^{d_1} \otimes_{\clsub{R}} \dots \otimes_{\clsub{R}}
		\ul{a_k} \otimes_{\clsub{R}} a_k^1 \otimes_{\clsub{R}} \dots \otimes_{\clsub{R}}
		a_k^{d_k}
		\right) =
		\\
		\ul{a_1} \otimes_{R} a_1^1 \otimes_{R} \dots \otimes_{R} a_1^{d_1} \otimes_{R} \dots \otimes_{R}
		\ul{a_k} \otimes_{R} a_k^1 \otimes_{R} \dots \otimes_{R} a_k^{d_k},
	\end{gather*}
	the chain $\Go{b}[\geq 1][\mathcal{A}]$ (resp.\ $\Ho{b}[\geq 1][\mathcal{A}]$) is given by
	the same formula \eqref{eq:G_n-explicit-formula} (resp.\ \eqref{eq:H_n(b)-formula}) as in the $R$-linear case.

	\begin{lm} \label{lm:Go-and-Ho-geq1-properties-general}
	Let $\mathcal{A} = \left( A, \mu \right)$ be a Banach $\Ainf$-algebra over a differential
	graded-commutative Banach $\mathbbm{k}$-algebra $\mathcal{R} = \left( R, d_R \right)$
	and let $\mathcal{B} = \left( B, \nu \right)$ be a Banach $\Ainf$-algebra over a
	differential graded-commutative Banach $\mathbbm{k}$-algebra $\mathcal{S} = \left( S, d_S \right)$.
	Let $f \colon \mathcal{A} \rightarrow \mathcal{B}$ be a Banach $\Ainf$-morphism.
	Given a topologically nilpotent element $b \in \tc{A}$, we have the following identities between chains
	of $\totcomp{\cdot}[1][]$:
	\begin{align}
		D_{\mathcal{A}} \left( \Ho{b}[\geq 1][\mathcal{A}] \right) & = 0,
		\label{eq:D-Ho-geq-1}
		\\
		D_{\mathcal{A}} \left( \Go{b}[\geq 1][\mathcal{A}] \right) & =
		\Ho{b}[\geq 1][\mathcal{A}] - \Ho{0}[\geq 1][\mathcal{A}],
		\label{eq:D-Go-geq-1}
		\\
		\cindmap{f} \left( \Ho{b}[\geq 1][\mathcal{A}] \right)     & =
		\Ho{\mcfunc{f} \left( b \right)}[\geq 1][\mathcal{B}].
		\label{eq:cindmap-f-Ho-geq-1}
	\end{align}
	We also have the identity
	\begin{equation} \label{eq:cindmap-f-Go-geq-1}
		\cindmap{f} \left( \Go{b}[\geq 1][\mathcal{A}] \right) =
		\Go{\mcfunc{f} \left( b \right)}[\geq 1][\mathcal{B}] -
		\Go{\mcfunc{f} \left( 0 \right)}[\geq 1][\mathcal{B}] +
		D_{\mathcal{B}} \left( R \left( b; f \right) \right)
	\end{equation}
	for some chain $R \left( b; f \right) \in \totcomp{B}[1][-2]$.
	When $b = 0$ or when $f$ is strict, one can choose $R \left( b; f \right) = 0$. \qed
	\end{lm}
	\begin{proof}
	The identities \eqref{eq:D-Ho-geq-1} and \eqref{eq:D-Go-geq-1} follow from the $R$-linear versions
	using the fact that
$\cindmap{\resover{i}} \colon \totcomp{\clsub{\mathcal{A}}/\clsub{\mathcal{R}}}[1][] \rightarrow
\totcomp{\mathcal{A}/\mathcal{R}}[1][]$
	is a chain map. We have
	\begin{equation*}
		D_{\mathcal{A}} \left( \Ho{b}[\geq 1][\mathcal{A}] \right)
		\stackrel{\eqref{eq:Ho-geq-1-def-general}}{=}
		\left( D_{\mathcal{A}} \circ \cindmap{\resover{i}} \right) \left( \Ho{b}[\geq 1][\clsub{\mathcal{A}}] \right)
		=
		\left( \cindmap{\resover{i}} \circ D_{\clsub{\mathcal{A}}} \right) \left( \Ho{b}[\geq 1][\clsub{\mathcal{A}}] \right)
		\stackrel{\eqref{eq:D-Ho-geq-1-linear}}{=}
		\cindmap{\resover{i}} \left( 0 \right) = 0,
	\end{equation*}
	and
	\begin{equation*}
		\begin{aligned}
			D_{\mathcal{A}} \left( \Go{b}[\geq 1][\mathcal{A}] \right)
			\stackrel{\eqref{eq:Go-geq-1-def-general}}{=} &
			\left( D_{\mathcal{A}} \circ \cindmap{\resover{i}} \right) \left(
			\Go{b}[\geq 1][\clsub{\mathcal{A}}]
			\right)
			=
			\left( \cindmap{\resover{i}} \circ D_{\clsub{\mathcal{A}}} \right) \left(
			\Go{b}[\geq 1][\clsub{\mathcal{A}}]
			\right)
			\\
			\stackrel{\eqref{eq:D-Go-geq-1-linear}}{=}    &
			\left( \cindmap{\resover{i}} \right) \left(
			\Ho{b}[\geq 1][\clsub{\mathcal{A}}] - \Ho{0}[\geq 1][\clsub{\mathcal{A}}]
			\right)
			\\
			\stackrel{\eqref{eq:Ho-geq-1-def-general}}{=} &
			\Ho{b}[\geq 1][\mathcal{A}] - \Ho{0}[\geq 1][\mathcal{A}].
		\end{aligned}
	\end{equation*}
	To show \cref{eq:cindmap-f-Ho-geq-1,eq:cindmap-f-Go-geq-1}, consider the following pullback diagram,
	constructed as in \cref{dfn:pullback-f-along-diagram}:
	\begin{equation*}
		\begin{tikzcd}
			{\left( \tens{A}[R], \mu \right)} &&&& {\left( \tens{B}[S], \nu \right)} \\
			& {\left( R, d_R \right)} && {\left( S, d_S \right)} \\
			& {\left( \clsub{R}, 0 \right)} && {\left( \clsub{S}, 0 \right)} \\
			{\left( \tens{i^{*} \left( A \right)}[\clsub{R}], \clsub{\mu} \right)} &&&&
			{\left( \tens{j^{*} \left( B \right)}[\clsub{S}], \clsub{\nu} \right)}
			\arrow["{\base{f}}", from=2-2, to=2-4]
			\arrow["{i}", from=3-2, to=2-2]
			\arrow["{j}"', from=3-4, to=2-4]
			\arrow["f_{\textrm{base}}^{\textrm{cl}}"', from=3-2, to=3-4] 
			\arrow["f", from=1-1, to=1-5]
			\arrow["{\resover{i}}", from=4-1, to=1-1]
			\arrow["{\resover{j}}"', from=4-5, to=1-5]
			\arrow["{\clsup{f}}"', from=4-1, to=4-5]
			\arrow["{\varepsilon_{\clsub{R}}}", from=4-1, to=3-2]
			\arrow["{\varepsilon_R}"', from=1-1, to=2-2]
			\arrow["{\varepsilon_S}", from=1-5, to=2-4]
			\arrow["{\varepsilon_{\clsub{S}}}"', from=4-5, to=3-4]
		\end{tikzcd}
	\end{equation*}
	In the diagram above, the morphism $f_{\textrm{base}}^{\textrm{cl}}$ coincides with
	the restriction of the morphism $\base{f}$, the underlying morphism of $f$ to $\clsub{R}$.
	The morphism $\clsup{f}$ is the pullback of the morphism $f$ along the middle square diagram
	and its corestriction satisfies
	\begin{equation*}
		f_k^{\textrm{cl}} \left( a_1 \otimes_{\clsub{R}} \dots \otimes_{\clsub{R}} a_k \right) =
		f_k \left( a_1 \otimes_R \dots \otimes_R a_k \right)
	\end{equation*}
	for all $k \geq 0$.	Hence, we have $\mcfunc{\clsup{f}} \left( b \right) = \mcfunc{f} \left( b \right)$ and so
	\begin{equation*}
		\begin{aligned}
			\cindmap{f} \left( \Go{b}[\geq 1][\mathcal{A}] \right) \!
			\eqwithref[eq:Go-geq-1-def-general]          &
			\cindmap{f} \left( \cindmap{\resover{i}} \left( \Go{b}[\geq 1][\clsub{\mathcal{A}}] \right) \right)
			\\
			\eqwithref[eq:cind-map-functoriality]        &
			\cindmap{\left( f \circ \resover{i} \right)} \left( \Go{b}[\geq 1][\clsub{\mathcal{A}}] \right) =
			\cindmap{\left( \resover{j} \circ \clsup{f} \right)} \left(
			\Go{b}[\geq 1][\clsub{\mathcal{A}}]
			\right)
			\\
			\eqwithref[eq:cind-map-functoriality]        &
			\cindmap{\resover{j}} \left( \clsup{\cindmap{f}} \left(
			\Go{b}[\geq 1][\clsub{\mathcal{A}}] \right) \right)
			\\
			\eqwithref[eq:cindmap-f-Go-geq-1-R-linear]   &
			\cindmap{\resover{j}} \left(
			\Go{\mcfunc{\clsup{f}} \left( b \right)}[\geq 1][\clsub{\mathcal{B}}] -
			\Go{\mcfunc{\clsup{f}} \left( 0 \right)}[\geq 1][\clsub{\mathcal{B}}] +
			D_{\clsub{\mathcal{B}}} \left( R \left( b; \clsup{f} \right) \right)
			\right)
			\\
			\eqwithref[eq:func-qdr-cyc][eq:func-lie-cyc] &
			\cindmap{\resover{j}} \left(
			\Go{\mcfunc{f} \left( b \right)}[\geq 1][\clsub{\mathcal{B}}] -
			\Go{\mcfunc{f} \left( 0 \right)}[\geq 1][\clsub{\mathcal{B}}]
			\right) +
			D_{\mathcal{B}} \left( \cindmap{\resover{j}} \left( R \left( b; \clsup{f} \right) \right) \right)
			\\
			\eqwithref[eq:Go-geq-1-def-general]          &
			\Go{\mcfunc{f} \left( b \right)}[\geq 1][\mathcal{B}] -
			\Go{\mcfunc{f} \left( 0 \right)}[\geq 1][\mathcal{B}] +
			D_{\mathcal{B}} \left( \cindmap{\resover{j}} \left( R \left( b; \clsup{f} \right) \right) \right),
		\end{aligned}
	\end{equation*}
	which shows \cref{eq:cindmap-f-Go-geq-1}. The proof of \cref{eq:cindmap-f-Ho-geq-1} is similar.
	\end{proof}

	\begin{rem} \label{rem:Go-Ho-geq1-vs-0}
	The properties of the chains $\Go{b}[\geq 1][\mathcal{A}]$ and $\Ho{b}[\geq 1][\mathcal{A}]$
	belonging to $\totcomp{\mathcal{A}}[1][]$ stated in
	\cref{lm:Go-and-Ho-geq1-properties-general}
	are almost identical to the properties of the chains $\Go{b}[0][\mathcal{A}]$ and $\Ho{b}[0][\mathcal{A}]$
	belonging to $\ncdfr{\mathcal{A}}[0][]$ stated in \cref{lm:Go_0-and-Ho_0-properties}, except
	that \cref{eq:cycl-f-Go-0} holds only up to a boundary term, so it is replaced by
	\cref{eq:cindmap-f-Go-geq-1}.
	\end{rem}

	\subsection{The Extended Total Complex \texorpdfstring{$\totcompe{\mathcal{A}}[1][]$}{with Non-Negative Line Degree}}
	\label{sec:extended-tot-comp-geq1-braidop-2}

	Recall that the extended cyclic complex $\ncdf{\mathcal{A}}[0][]$ is obtained from the
	standard cyclic complex $\ncdfr{\mathcal{A}}[0][]$ by adjoining the element $1$ and setting
$\clie{\mu} \left( 1 \right) = \mu_0 \left( 1 \right)$, i.e., adding a primitive to the curvature
$\mu_0 \left( 1 \right)$. This works as the curvature $\mu_0 \left( 1 \right) \in \ncdfr{A}[0][1]$
	is a cycle in the cyclic complex. Moreover, given a Banach $\Ainf$-morphism
$f \colon \mathcal{A} \rightarrow \mathcal{B}$
	between $\mathcal{A} = \left( A, \mu \right)$ and $\mathcal{B} = \left( B, \nu \right)$,
	one can extend the induced chain map
$\cycl{f} \colon \ncdfr{\mathcal{A}}[0][] \rightarrow \ncdfr{\mathcal{B}}[0][]$
	to a map $\cycl{f} \colon \ncdf{\mathcal{A}}[0][] \rightarrow \ncdf{\mathcal{B}}[0][]$
	by setting
	\begin{equation*}
		\cycl{f} \left( 1 \right) =
		\cexp{f_0 \left( 1 \right)} =
		\G{f_0 \left( 1 \right)}[0][\mathcal{B}] =
		1 + \Go{f_0 \left( 1 \right)}[0][\mathcal{B}]
	\end{equation*}
	and obtain a functorial induced chain map.

	We will do something similar and extend the total complex $\totcomp{\mathcal{A}}[1][]$ to a complex
$\totcompe{\mathcal{A}}[1][]$, and the projection
$p_1 \colon \totcomp{\mathcal{A}}[1][*] \rightharpoonup \ncdfr{\mathcal{A}}[0][*+1]$, given by
	\eqref{eq:def-p1-anticommute},
	to a projection $p_1^{+} \colon \totcompe{\mathcal{A}}[1][*] \rightharpoonup \ncdf{\mathcal{A}}[0][*+1]$
	which will remain a homotopy equivalence. This will give us another
	chain complex model $\totcompe{\mathcal{A}}[1][]$ for the \textit{extended cyclic homology}.

	Consider the chain $\Ho{0}[\geq 1][\mathcal{A}] \in \totcomp{\mathcal{A}}[1][0]$ given by
	\begin{equation} \label{eq:Ho_geq_1-at-0-formula}
		\Ho{0}[\geq 1][\mathcal{A}] = \sum_{n = 1}^{\infty} \H{0}[n][\mathcal{A}]
		\stackrel{\eqref{eq:H_n(b)-formula}}{=}
		- \sum_{n=1}^{\infty} \frac{1}{n} \, \ul{\mu_0 \left( 1 \right)}^{\otimes n}.
	\end{equation}
	By \cref{eq:D-Ho-geq-1}, we have $D_{\mathcal{A}} \left( \Ho{0}[\geq 1][\mathcal{A}] \right) = 0$,
	so $\Ho{0}[\geq 1][\mathcal{A}]$ is a cycle lifting the cycle $-\mu_0 \left( 1 \right)$ under the map $p_1$.
	Let $f \colon R \rightarrow \totcomp{\mathcal{A}}[1][]$ be the degree zero $R$-linear map determined by
$f \left( 1 \right) = \Ho{0}[\geq 1][\mathcal{A}]$. Since $\Ho{0}[\geq 1][\mathcal{A}]$ is a cycle,
	the map $f$ is a chain map. Let us define $\totcompe{\mathcal{A}}[1][] \defeq \Cone{f}$ and call
$\totcompe{\mathcal{A}}[1][]$ the \textbf{extended total complex}.
	The underlying graded module of
$\totcompe{\mathcal{A}}[1][]$ is given by  $\totcompe{A}[1][] = R[1] \oplus \totcomp{A}[1][]$
	and we continue to denote the
	differential on $\totcompe{\mathcal{A}}[1][]$ by $D_{\mathcal{A}}$, which is consistent with the fact
	that $\totcomp{\mathcal{A}}[1][]$ is a subcomplex of $\totcompe{\mathcal{A}}[1][]$. We denote the degree $-1$
	element $\s \left( 1_R \right) \in \totcompe{A}[1][]$ by $\ul{1}$ so that we have
$r \cdot \ul{1} = (-1)^{\degb{r}} \ul{r}$, and
	\begin{equation} \label{eq:D-geq-1-action-on-ul-1}
		\begin{aligned}
			D_{\mathcal{A}} \left( \ul{1} \right) ={} &
			D_{\mathcal{A}}^{\geq 1, +} \left( \ul{1} \right) =
			\Ho{0}[\geq 1][\mathcal{A}] =
			- \sum_{n=1}^{\infty} \frac{1}{n} \ul{\mu_0(1)}^{\otimes n},
			\\
			D_{\mathcal{A}} \left( \ul{r} \right) ={} &
			D_{\mathcal{A}}^{\geq 1, +} \left( \ul{r} \right) =
			(-1)^{\degb{r}} D_{\mathcal{A}} \left( r \cdot \ul{1} \right) =
			(-1)^{\degb{r}} dr \cdot \ul{1} + r \cdot D_{\mathcal{A}} \left( \ul{1} \right).
		\end{aligned}
	\end{equation}
	Note that the notation $\ul{1}$ is
	consistent with thinking of $1$ as an element of cohomological degree $0$ so that $\ul{1}$ has ``bidegree'' $(1,0)$
	and total degree $-1$.
	Finally, we extend the projection $p_1$ to a projection
$p_1^{+} \colon \totcompe{\mathcal{A}}[1][*] \rightharpoonup \ncdf{\mathcal{A}}[0][*+1]$
	by defining $p_1^{+}(\ul{1}) \defeq 1 \in \ncdf{A}[0][0]$. Then
	\begin{lm}
	The extended projection $p_1^{+} \colon  \totcompe{A}[1][*] \rightharpoonup \ncdf{A}[0][*+1]$
	is a chain map which is a homotopy equivalence.
	\end{lm}
	\begin{proof}
	Let $g \colon R \rightarrow {\ncdfr{A}[0][]}[1]$ be the chain map given by
$g \left( r \right) = -r \cdot \s \left( \mu_0 \left( 1 \right) \right)$.
	The map $\s \circ p_1 \colon \totcomp{A}[1][] \rightarrow {\ncdfr{A}[0][]}[1]$ is a
	degree zero homotopy equivalence, and we have $\left( \s \circ p_1 \right) \circ f = g$, so by
	\cref{cor:triangle-homotopy-equivalence-mapping-cone}, the map
$\left( \s \circ p_1 \right)^{+} \colon \Cone{f} \rightarrow \Cone{g}$ given by
$\left( \s r, x \right) \mapsto \left( \s r, \s p_1 \left( x \right) \right)$ is a homotopy equivalence.
	Note that $\Cone{g}$ is isomorphic to ${\ncdf{A}[0][]}[1]$ via
$\left( \s r, \s x_0 \right) \mapsto \s \left( r + x_0 \right)$ so the composition
$\totcompe{A}[1][] \rightarrow {\ncdf{A}[0][]}[1]$ given by
$\left( \s r, x \right) \mapsto \s \left( r + p_1 \left( x \right) \right)$ is a homotopy equivalence
	which is precisely $\s \circ p_1^{+}$. Hence, $p_1^{+}$ is a degree one homotopy equivalence.
	\end{proof}

	Given $b \in \tc{A}$ a topologically nilpotent element, let us define the chain
	\begin{equation} \label{eq:G-geq-1-def}
		\G{b}[\geq 1] = \G{b}[\geq 1][\mathcal{A}] \defeq
		\ul{1} + \Go{b}[\geq 1][\mathcal{A}] \in \totcompe{A}[1][-1]
	\end{equation}
	which belongs to the extended total complex.
	Given a Banach $\Ainf$-morphism
$f \colon \mathcal{A} \rightarrow \mathcal{B}$ between $\mathcal{A} = \left( A, \mu \right)$ and
$\mathcal{B} = \left( B, \nu \right)$, we extend the induced map
$\cindmap{f} \colon \totcomp{\mathcal{A}}[1][] \rightarrow \totcomp{\mathcal{B}}[1][]$
	to a map $\cindmape{f} \colon \totcompe{\mathcal{A}}[1][] \rightarrow \totcompe{\mathcal{B}}[1][]$
	between the extended total complexes by defining
	\begin{equation} \label{def:cycl-f-ul-1}
		\begin{aligned}
			\cindmape{f} \left( \ul{1} \right) \defeq{} &
			\G{ \mcfunc{f} \left( 0 \right) }[\geq 1][\mathcal{B}] =
			\G{ f_0 \left( 1 \right) }[\geq 1][\mathcal{B}] =
			\ul{1} + \Go{ f_0 \left( 1 \right) }[\geq 1][\mathcal{B}]
			\\
			={}                                         &
			\ul{1} + \sum_{i=0}^{\infty} \frac{1}{i+1} \, \ul{f_0(1)} \otimes f_0(1)^{\otimes i} +
			\G{f_0 \left( 1 \right)}[2][\mathcal{B}] + \dots
		\end{aligned}
	\end{equation}
	With the definitions above, we have the following properties of the chain $\G{b}[\geq 1][\mathcal{A}]$:
	\begin{lm} \label{lm:G-geq-1-properties}
	Let $\mathcal{A} = \left( A, \mu \right)$ be a Banach $\Ainf$-algebra over a differential
	graded-commutative Banach $\mathbbm{k}$-algebra $\mathcal{R} = \left( R, d_R \right)$
	and let $\mathcal{B} = \left( B, \nu \right)$ be a Banach $\Ainf$-algebra over a
	differential graded-commutative Banach $\mathbbm{k}$-algebra $\mathcal{S} = \left( S, d_S \right)$.
	Let $f \colon \mathcal{A} \rightarrow \mathcal{B}$ be a Banach $\Ainf$-morphism.
	Given a topologically nilpotent element $b \in \tc{A}$, we have the following identities
	between chains of $\totcompe{\cdot}[1][]$:
	\begin{enumerate}
	\item{(Differential)}
	\begin{equation} \label{eq:D-G-geq-1}
		D_{\mathcal{A}} \left( \G{b}[\geq 1][\mathcal{A}] \right) = \Ho{b}[\geq 1][\mathcal{A}].
	\end{equation}
	\item{(Weak Naturality)}
	\begin{equation} \label{eq:cindmap-f-G-geq-1}
		\cindmape{f} \left( \G{b}[\geq 1][\mathcal{A}] \right) =
		\G{\mcfunc{f} \left( b \right)}[\geq 1][\mathcal{B}] +
		D_{\mathcal{B}} \left( R \left( b; f \right) \right)
	\end{equation}
	for some chain $R \left( b; f \right) \in \totcomp{B}[1][-2]$.
	When $b = 0$ or when $f$ is strict, one can choose $R \left( b; f \right) = 0$.
	\end{enumerate}
	\end{lm}
	\begin{proof}
		We have
		\begin{equation*}
			D_{\mathcal{A}} \left( \G{b}[\geq 1][\mathcal{A}] \right)
			\stackrel{\eqref{eq:G-geq-1-def}}{=}
			D_{\mathcal{A}} \left( \ul{1} + \Go{b}[\geq 1][\mathcal{A}] \right)
			\stackrel[\eqref{eq:D-Go-geq-1}]{\eqref{eq:D-geq-1-action-on-ul-1}}{=}
			\Ho{0}[\geq 1][\mathcal{A}] +
			\left( \Ho{b}[\geq 1][\mathcal{A}] - \Ho{0}[\geq 1][\mathcal{A}] \right)
			=
			\Ho{b}[\geq 1][\mathcal{A}]
		\end{equation*}
		which shows \cref{eq:D-G-geq-1}. Similarly,
		\begin{equation*}
			\begin{aligned}
				\cindmape{f} \left( \G{b}[\geq 1][\mathcal{A}] \right) \!
				\eqwithref[eq:G-geq-1-def]                         &
				\cindmape{f} \left( \ul{1} + \Go{b}[\geq 1][\mathcal{A}] \right)
				=
				\cindmape{f} \left( \ul{1} \right) +
				\cindmap{f} \left( \Go{b}[\geq 1][\mathcal{A}] \right)
				\\
				\eqwithref[eq:cindmap-f-Go-geq-1][def:cycl-f-ul-1] &
				\G{ \mcfunc{f} \left( 0 \right) }[\geq 1][\mathcal{B}] +
				\left(
				\Go{\mcfunc{f} \left( b \right)}[\geq 1][\mathcal{B}] -
				\Go{\mcfunc{f} \left( 0 \right)}[\geq 1][\mathcal{B}] +
				D_{\mathcal{B}} \left( R \left( b; f \right) \right)
				\right)
				\\
				\eqwithref[eq:G-geq-1-def]                         &
				\ul{1} + \Go{ \mcfunc{f} \left( 0 \right) }[\geq 1][\mathcal{B}] +
				\left(
				\Go{\mcfunc{f} \left( b \right)}[\geq 1][\mathcal{B}] -
				\Go{\mcfunc{f} \left( 0 \right)}[\geq 1][\mathcal{B}] +
				D_{\mathcal{B}} \left( R \left( b; f \right) \right)
				\right)
				\\
				\eqwithref[eq:G-geq-1-def]                         &
				\G{\mcfunc{f} \left( b \right)}[\geq 1][\mathcal{B}] +
				D_{\mathcal{B}} \left( R \left( b; f \right) \right)
			\end{aligned}
		\end{equation*}
		which shows \cref{eq:cindmap-f-G-geq-1}.
	\end{proof}

	\begin{lm} \label{lm:cindmape-f-geq-1-chain-map}
	The map $\cindmape{f} \colon \totcompe{A}[1][] \rightarrow \totcompe{B}[1][]$ between the extended
	total complexes is a chain map, i.e., we have
	\begin{equation} \label{eq:D-cindmap-f-ul-1}
		D_{\mathcal{B}} \left( \cindmape{f} \left( \ul{1} \right) \right) =
		\cindmape{f} \left( D_{\mathcal{A}} \left( \ul{1} \right) \right).
	\end{equation}
	\end{lm}
	\begin{proof}
		We have
		\begin{equation*}
			\begin{aligned}
				D_{\mathcal{B}} \left( \cindmape{f} \left( \ul{1} \right) \right)
				\stackrel{\eqref{def:cycl-f-ul-1}}{=}{}           &
				D_{\mathcal{B}} \left( \G{ \mcfunc{f} \left( 0 \right) }[\geq 1][\mathcal{B}] \right)
				\stackrel{\eqref{eq:D-G-geq-1}}{=}{}
				\Ho{ \mcfunc{f} \left( 0 \right) }[\geq 1][\mathcal{B}]
				\stackrel{\eqref{eq:cindmap-f-Ho-geq-1}}{=}{}
				\cindmap{f} \left( \Ho{0}[\geq 1][\mathcal{A}] \right)
				\\
				\stackrel{\eqref{eq:D-geq-1-action-on-ul-1}}{=}{} &
				\cindmap{f} \left( D_{\mathcal{A}} \left( \ul{1} \right) \right)
				=
				\cindmape{f} \left( D_{\mathcal{A}} \left( \ul{1} \right) \right).
			\end{aligned}
		\end{equation*}
	\end{proof}

	To summarize, given a Banach $\Ainf$-algebra $\mathcal{A} = \left( A, \mu \right)$
	over a graded-commutative Banach $\mathbbm{k}$-algebra $\mathcal{R}$,
	we have constructed a chain complex $\totcompe{\mathcal{A}}[1][]$ over $\mathcal{R}$ which is homotopy equivalent to
$\ncdf{\mathcal{A}}[0][]$. Also, given another Banach $\Ainf$-algebra $\mathcal{B} = \left( B, \nu \right)$
	and a Banach $\Ainf$-morphism $f \colon \mathcal{A} \rightarrow \mathcal{B}$,
	possibly with a change of connection term $f_0 \left( 1 \right)$,
	we have constructed an induced \textit{chain map}
$\cindmape{f} \colon \totcompe{\mathcal{A}}[1][] \rightarrow \totcompe{\mathcal{B}}[1][]$,
	extending the chain map
$\cindmap{f} \colon \totcomp{\mathcal{A}}[1][] \rightarrow \totcomp{\mathcal{B}}[1][]$.

	\begin{rem} \label{rem:totcompe-geq-1-not-functorial-nose}
	Unlike the constructions
	\begin{align*}
		                                              & \mathcal{A}                             &  & \mapsto &   & \totcomp{\mathcal{A}}[1][],
		\\
		f \colon \mathcal{A}                          & \rightarrow \mathcal{B}                 &  & \mapsto &
		\cindmap{f} \colon \totcomp{\mathcal{A}}[1][] & \rightarrow \totcomp{\mathcal{B}}[1][],
	\end{align*}
	and
	\begin{align*}
		                                        & \mathcal{A}                          &  & \mapsto &   & \ncdf{\mathcal{A}}[0][],
		\\
		f \colon \mathcal{A}                    & \rightarrow \mathcal{B}              &  & \mapsto &
		\cycl{f} \colon \ncdf{\mathcal{A}}[0][] & \rightarrow \ncdf{\mathcal{B}}[0][],
	\end{align*}
	the extended total complex construction $\totcompe{\cdot}[1][]$ is not functorial on the nose: Given another Banach
$\Ainf$-morphism $g \colon \mathcal{B} \rightarrow \mathcal{C}$,
	the weak naturality property \eqref{eq:cindmap-f-G-geq-1} implies that we have
	\begin{equation*}
		\begin{aligned}
			\left( \cindmape{g} \circ \cindmape{f} \right) \left( \ul{1} \right) & =
			\cindmape{g} \left( \G{ \mcfunc{f} \left( 0 \right) }[\geq 1][\mathcal{B}] \right) =
			\G{ \mcfunc{g} \left( \mcfunc{f} \left( 0 \right) \right) }[\geq 1][\mathcal{C}] +
			D_{\mathcal{C}} \left( R \left( \mcfunc{f} \left( 0 \right); g \right) \right)
			\\
			                                                                     & =
			\G{ \mcfunc{ \left( g \circ f \right) } \left( 0 \right) }[\geq 1][\mathcal{C}] +
			D_{\mathcal{C}} \left( R \left( \mcfunc{f} \left( 0 \right); g \right) \right)
			\\
			                                                                     & =
			\cindmape{ \left( g \circ f \right) } \left( \ul{1} \right) +
			D_{\mathcal{C}} \left( R \left( f, g \right) \right)
		\end{aligned}
	\end{equation*}
	for some
$R \left( f, g \right) = R \left( \mcfunc{f} \left( 0 \right); g \right) \in
\totcomp{C}[1][-2] \subset \totcompe{C}[1][-2]$ which depends on $f$,
	via $f_0 \left( 1 \right)$, and $g$. Hence, the chain maps
	\begin{equation*}
		\cindmape{g} \circ \cindmape{f}, \cindmape{ \left( g \circ f \right) } \colon
		\totcompe{\mathcal{A}}[1][] \rightarrow
		\totcompe{\mathcal{C}}[1][]
	\end{equation*}
	agree only \textit{up to an exact term}. However, if we denote by
	\begin{gather*}
		\cohom{\cindmape{f}}[] \colon \cohom{{\totcompe{\mathcal{A}}[1][]}}[] \rightarrow
		\cohom{{\totcompe{\mathcal{B}}[1][]}}[],
		\qquad
		\cohom{\cindmape{g}}[] \colon \cohom{{\totcompe{\mathcal{B}}[1][]}}[] \rightarrow
		\cohom{{\totcompe{\mathcal{C}}[1][]}}[],
		\\
		\cohom{\cindmape{g} \circ \cindmape{f}}[] \colon \cohom{{\totcompe{\mathcal{A}}[1][]}}[] \rightarrow
		\cohom{{\totcompe{\mathcal{C}}[1][]}}[],
	\end{gather*}
	the induced maps between the cohomologies, we do have
$\cohom{\cindmape{g} \circ \cindmape{f}}[] = \cohom{\cindmape{g}}[] \circ \cohom{\cindmape{f}}[]$.
	Note also that if either $f_0 \left( 1 \right) = 0$ or if $g$ is strict, then the weak naturality
	property \eqref{eq:cindmap-f-G-geq-1}
	shows that we can take $R \left( f, g \right) = 0$, and in this case we do have
$\cindmape{ \left( g \circ f \right) } = \cindmape{g} \circ \cindmape{f}$ on $\totcompe{\cdot}[1][]$.
	In particular, we see that the construction is functorial on the level of chain complexes if we
	restrict ourselves to $\Ainf$-morphisms without change of connection terms.
	\end{rem}

	\begin{rem}
	In \cref{sec:cyclization-generalized-coderivation} (resp.\ \cref{sec:cyclization-coalgebra-morphisms}),
	we defined the cyclization of coderivations $\mu \colon \tens{A} \rightharpoonup \tens{A}$
	(resp.\ morphisms $f \colon \tens{A} \rightarrow \tens{B}$), resulting in operators
	defined on the \textit{reduced} tensor module $\tensr{A}$. Since in our setting, coderivations
	and coalgebra morphisms have $0$-arity terms and actually live on the \textit{full} tensor module $\tens{A}$,
	it was natural to see if we can extend the construction to $\tens{A}$, which is the ``extended''
	version of $\tensr{A}$.

	To do that, we introduced the cyclic exponential $\cexp{b} = \G{b}[0] \in \tens{A}$, and used it
	to extend the action of $\cycl{f}$ by setting
$\cycl{f} \left( 1 \right) = \cexp{f_0 \left( 1 \right)}$
	(cf.\ $f \left( 1 \right) = \Exp{f_0 \left( 1 \right)}$),
	and $\cycl{\mu} \left( 1 \right) = \mu_0 \left( 1 \right)$.
	We then studied the properties of $\G{b}[0]$
	and used them to show our extension is functorial on the cyclic quotient $\tenscyc{A}$
	(see \cref{sec:extension-cycl-full-tensor-module}). Finally, in \cref{lm:Go_0-and-Ho_0-properties},
	we deduced the properties of $\Go{b}[0] = \cexp{b} - 1$
	in $\tensrcyc{A} = \ncdfr{\mathcal{A}}[0][]$
	from the corresponding properties of $\G{b}[0]$ in $\tenscyc{A} = \ncdf{\mathcal{A}}[0][]$.

	With $\totcompe{\mathcal{A}}[1][]$, we went the other way around. We first
	constructed a chain $\Go{b}[\geq 1]$ in $\totcomp{\mathcal{A}}[1][]$
	and studied its properties in \cref{sec:construction-special-elements}.
	We then adjoined $\ul{1}$ to $\totcomp{\mathcal{A}}[1][]$,
	defined $\G{b}[\geq 1] = \ul{1} + \Go{b}[\geq 1]$, extended the action
	of $\cindmap{f}$ and $D_{\mu}$ from $\totcomp{\mathcal{A}}[1][]$ to $\totcompe{\mathcal{A}}[1][]$, and finally
	deduced the properties of $\G{b}[\geq 1]$ in $\totcompe{\mathcal{A}}[1][]$ from the corresponding
	properties of $\Go{b}[\geq 1]$ in $\totcomp{\mathcal{A}}[1][]$.
	\end{rem}

	\subsection{The Extended Total Complex \texorpdfstring{$\totcompe{\mathcal{A}}[2][]$}{with Line Degree Greater Than One}}
	\label{sec:extended-tot-comp-geq2-braidop-2}
	Recall that we have defined the complex
$\totcomp{\mathcal{A}}[2][] \defeq \totc{\ncdf{\mathcal{A}}[][]}[][\geq 2][]$
	as the total complex of the bicomplex obtained from $\ncdf{\mathcal{A}}[][]$
	by erasing the first two columns. The construction $\totcomp{\cdot}[2][]$ is functorial:
	Given a Banach $\Ainf$-morphism
$f \colon \mathcal{A} \rightarrow \mathcal{B}$, the induced functorial morphism
$\cindmap{f} \colon \ncdf{\mathcal{A}}[][] \rightarrow \ncdf{\mathcal{B}}[][]$
	acts on each column separately and commutes with both $\qdr$ and $\clie{\mu}$, so
	the totalization of its restriction to columns $\geq 2$ gives us a functorial chain map
$\cindmap{f} \colon \totcomp{\mathcal{A}}[2][] \rightarrow \totcomp{\mathcal{B}}[2][]$.

	Let us denote by $\totc{\ncdf{\mathcal{A}}[][]}[][{[1]}][]$ the total complex of the bicomplex
	obtained from $\ncdf{\mathcal{A}}[][]$ by erasing all columns except column one. By our conventions,
	we have $\totc{\ncdf{\mathcal{A}}[][]}[][{[1]}][*] = \ncdf{\mathcal{A}}[1][*+1]$, endowed
	with the differential $\clie{\mu}$  (see \cref{eq:tot-com-def-parity-2}).
	On the level of graded $\mathbbm{k}$-modules,
$\totc{\ncdf{A}[][]}[][{[1]}][]$ coincides with ${\ncdf{A}[1][]}[1]$,
	except that we don't twist the $R$-action on $\totc{\ncdf{A}[][]}[][{[1]}][]$ and
	don't twist the differential by a sign.
	We then have a short exact sequence
	\begin{equation*}
		0 \rightarrow \totc{\ncdf{\mathcal{A}}[][]}[][{[1]}][] \rightarrow
		\totcomp{\mathcal{A}}[1][] \xrightarrow{\pi_{\geq 2}} \totcomp{\mathcal{A}}[2][] \rightarrow 0,
	\end{equation*}
	which can be taken as an alternative definition of $\totcomp{\mathcal{A}}[2][]$.

	Similarly, we
	define the \textbf{extended total complex} $\totcompe{\mathcal{A}}[2][]$ as the quotient
	of the extended total complex $\totcompe{\mathcal{A}}[1][]$ by $\totc{\ncdf{\mathcal{A}}[][]}[][{[1]}][]$.
	We have a short exact sequence
	\begin{equation*}
		0 \rightarrow \totc{\ncdf{\mathcal{A}}[][]}[][{[1]}][] \rightarrow
		\totcompe{\mathcal{A}}[1][] \xrightarrow{\pi_{\geq 2}^{+}} \totcompe{\mathcal{A}}[2][] \rightarrow 0,
	\end{equation*}
	where the projection chain map $\pi_{\geq 2}^{+}$ maps elements
$\sum_{i \geq 1} x_i$ with $x_i \in \ncdf{A}[i][]$ to $\sum_{i \geq 2} x_i$, and $\ul{1}$ to $\ul{1}$,
	or, more precisely, $\ul{1}_{\geq 1}$ to $\ul{1}_{\geq 2}$.
	The extended total complex $\totcompe{\mathcal{A}}[2][]$
	contains $\totcomp{\mathcal{A}}[2][]$ as a subcomplex, and can be thought of as obtained
	from $\totcomp{\mathcal{A}}[2][]$ by adjoining $\ul{1} = \ul{1}_{\geq 2}$
	and defining
	\begin{equation} \label{eq:action-D-geq-2-on-ul-1}
		D_{\mathcal{A}} \left( \ul{1} \right) = D_{\mathcal{A}}^{\geq 2, +} \left( \ul{1} \right)
		= \pi_{\geq 2}^{+} \left( \Ho{0}[\geq 1][\mathcal{A}] \right)
		= \sum_{n = 2}^{\infty} \H{0}[n][\mathcal{A}] =
		-\sum_{n=2}^{\infty} \frac{1}{n} \ul{\mu_0(1)}^{\otimes n}
	\end{equation}
	(compare with \cref{eq:D-geq-1-action-on-ul-1}).

	Given $b \in \tc{A}$, let us set
	\begin{align}
		\Go{b}[\geq 2][\mathcal{A}] & \defeq \pi_{\geq 2} \left( \Go{b}[\geq 1][\mathcal{A}] \right) =
		\sum_{n = 2}^{\infty} \G{b}[n] \in \totcomp{A}[2][-1],
		\label{eq:Go-geq2-def}
		\\
		\G{b}[\geq 2][\mathcal{A}]  & \defeq \pi_{\geq 2}^{+} \left( \G{b}[\geq 1][\mathcal{A}] \right) =
		\ul{1} + \Go{b}[\geq 2] \in \totcompe{A}[2][-1],
		\label{eq:G-geq-2-def}
		\\
		\Ho{b}[\geq 2][\mathcal{A}] & \defeq \pi_{\geq 2} \left( \Ho{b}[\geq 1][\mathcal{A}] \right) =
		\sum_{n=2}^{\infty} \H{b}[n] \in \totcomp{A}[2][0] \subset \totcompe{A}[2][0].
		\label{eq:Ho-geq-2-def}
	\end{align}

	Given a Banach $\Ainf$-morphism
$f \colon \mathcal{A} \rightarrow \mathcal{B}$ between Banach $\Ainf$-algebras
$\mathcal{A} = \left( A, \mu \right)$ and
$\mathcal{B} = \left( B, \nu \right)$, we can extend the induced map
$\cindmap{f} \colon \totcomp{\mathcal{A}}[2][] \rightarrow \totcomp{\mathcal{B}}[2][]$
	to a map $\cindmape{f} \colon \totcompe{\mathcal{A}}[2][] \rightarrow \totcompe{\mathcal{B}}[2][]$
	between the extended total complexes by defining
	\begin{align}
		\cindmape{f} \left( \ul{1} \right) \defeq{} &
		\G{ \mcfunc{f} \left( 0 \right) }[\geq 2][\mathcal{B}] =
		\G{ f_0 \left( 1 \right) }[\geq 2][\mathcal{B}]
		\label{def:cycl-f-ul-1-geq2}
		\\
		={}                                         &
		\ul{1} +
		\sum_{i_1,i_2,i_3=0}^{\infty} \frac{1}{i_1 + i_2 + i_3 + 1} \,
		\ul{f_0 \left( 1 \right)} \otimes {f_0 \left( 1 \right)}^{\otimes i_1} \otimes
		\ul{ \nu_{i_2} \left( {f_0 \left( 1 \right)}^{\otimes i_2} \right)} \otimes
		                                                                    {f_0 \left( 1 \right)}^{\otimes i_3} + \dots
		\notag
	\end{align}
	(compare with \cref{def:cycl-f-ul-1}).

	To be more explicit, let us temporarily denote the induced maps on $\totcomp{\mathcal{A}}[k][]$
	(resp.\ $\totcompe{\mathcal{A}}[k][]$)
	by $\cindmap{f}[k] \colon \totcomp{\mathcal{A}}[k][] \rightarrow \totcomp{\mathcal{B}}[k][]$
	(resp.\ $\cindmape{f}[k] \colon \totcompe{\mathcal{A}}[k][] \rightarrow \totcompe{\mathcal{B}}[k][]$)
	for $k = 1, 2$.
	Since $\cindmap{f}$ has bi-degree $(0,0)$, we have
$\cindmap{f}[2] \circ \pi_{\geq 2} = \pi_{\geq 2} \circ \cindmap{f}[1]$. Definition
	\eqref{def:cycl-f-ul-1-geq2} then guarantees that we also have
$\cindmape{f}[2] \circ \pi_{\geq 2}^{+} = \pi_{\geq 2}^{+} \circ \cindmape{f}[1]$.

	\begin{lm} \label{lm:cindmape-f-geq-2-chain-map}
	The map $\cindmape{f} \colon \totcompe{A}[2][] \rightarrow \totcompe{B}[2][]$ between the extended
	total complexes is a chain map.
	\end{lm}
	\begin{proof}
	Using the fact that $\cindmape{f}[\cdot]$ commutes with $\pi_{\geq 2}^{+}$, and that
	both $\pi_{\geq 2}^{+}$ and $\cindmape{f}[1]$ are chain maps (\cref{lm:cindmape-f-geq-1-chain-map}),
	we have
	\begin{equation*}
		\begin{aligned}
			\left( D_{\mathcal{B}}^{\geq 2, +} \circ \cindmape{f}[2] \right) \circ \pi_{\geq 2}^{+}
			 & =
			D_{\mathcal{B}}^{\geq 2, +} \circ  \pi_{\geq 2}^{+} \circ \cindmape{f}[1]
			=
			\pi_{\geq 2}^{+} \circ \left( D_{\mathcal{B}}^{\geq 1, +} \circ \cindmape{f}[1] \right)
			\\
			 & =
			\pi_{\geq 2}^{+} \circ \left( \cindmape{f}[1] \circ D_{\mathcal{A}}^{\geq 1, +} \right) =
			\cindmape{f}[2] \circ \pi_{\geq 2}^{+} \circ D_{\mathcal{A}}^{\geq 1, +}
			\\
			 & =
			\left( \cindmape{f}[2] \circ D_{\mathcal{A}}^{\geq 2, +} \right) \circ \pi_{\geq 2}^{+}.
		\end{aligned}
	\end{equation*}
	Since $\pi_{\geq 2}^{+}$ is surjective, we conclude that
$D_{\mathcal{B}}^{\geq 2, +} \circ \cindmape{f}[2] = \cindmape{f}[2] \circ D_{\mathcal{A}}^{\geq 2, +}$.
	\end{proof}

	\begin{dfn} \label{dfn:cyclic-chern-simons-form-braidop-2}
	Given $b \in \tc{A}$, the chain
	\begin{equation*}
		\begin{aligned}
			\G{b}[\geq 2] ={} & \ul{1} +
			\sum_{\substack{n = 2 \\ j_1, \dots, j_n = 0 \\ i_2, \dots, i_n = 0}}^{\infty}
			\frac{1}{1 + \sum_{r = 1}^n j_r + \sum_{r=2}^{n} i_r}
			\\
			                  & \qquad \qquad \qquad
			\ul{b} \otimes b^{\otimes j_1} \otimes \ul{ \mu_{i_2} \left( b^{\otimes i_2} \right)}
			\otimes b^{\otimes j_2} \otimes \dots \otimes \ul{ \mu_{i_n} \left( b^{\otimes i_n} \right) }
			\otimes b^{\otimes j_n} \in \totcompe{A}[2][-1]
		\end{aligned}
	\end{equation*}
	of the extended total complex $\totcompe{A}[2][]$ is called the
	\textbf{cyclic Chern--Simons form} associated to $b$.\footnote{For a relation to the classical Chern--Simons $3$-form and action,
	see \cref{sec:recovering-chern-simons}.}
	\end{dfn}

	We also have the analogous version of \cref{lm:G-geq-1-properties} describing the
	properties of the cyclic Chern--Simons form in $\totcompe{\cdot}[2][]$:
	\begin{lm} \label{lm:G-geq-2-properties}
	Let $\mathcal{A} = \left( A, \mu \right)$ be a Banach $\Ainf$-algebra over a differential
	graded-commutative Banach $\mathbbm{k}$-algebra $\mathcal{R} = \left( R, d_R \right)$
	and let $\mathcal{B} = \left( B, \nu \right)$ be a Banach $\Ainf$-algebra over a
	differential graded-commutative Banach $\mathbbm{k}$-algebra $\mathcal{S} = \left( S, d_S \right)$.
	Let $f \colon \mathcal{A} \rightarrow \mathcal{B}$ be a Banach $\Ainf$-morphism.
	Given a topologically nilpotent element $b \in \tc{A}$, we have the following identities
	between elements of $\totcompe{\cdot}[2][]$:
	\begin{enumerate}
	\item{(Differential)}
	\begin{equation} \label{eq:D-G-geq-2}
		D_{\mathcal{A}} \left( \G{b}[\geq 2][\mathcal{A}] \right) = \Ho{b}[\geq 2][\mathcal{A}].
	\end{equation}
	\item{(Weak Naturality)}
	\begin{equation} \label{eq:cindmap-f-G-geq-2}
		\cindmape{f} \left( \G{b}[\geq 2][\mathcal{A}] \right) =
		\G{\mcfunc{f} \left( b \right)}[\geq 2][\mathcal{B}] +
		D_{\mathcal{B}} \left( R \left( b; f \right) \right)
	\end{equation}
	for some element $R \left( b; f \right) \in \totcomp{B}[2][-2]$.
	When $b = 0$ or when $f$ is strict, one can choose $R \left( b; f \right) = 0$.
	\end{enumerate}
	\end{lm}
	\begin{proof}
	The lemma follows immediately from \cref{lm:G-geq-1-properties} using the facts
	that $\cindmape{f}[\cdot]$ commutes with $\pi_{\geq 2}^{+}$, and that
$\pi_{\geq 2}^{+}$ is a chain map. For example, we have
	\begin{equation*}
		\begin{aligned}
			D_{\mathcal{A}}^{\geq 2, +} \left( \G{b}[\geq 2][\mathcal{A}] \right)
			\stackrel{\eqref{eq:G-geq-2-def}}{=} &
			\left( D_{\mathcal{A}}^{\geq 2, +} \circ \pi_{\geq 2}^{+} \right) \left(
			\G{b}[\geq 1][\mathcal{A}]
			\right)
			=
			\left( \pi_{\geq 2}^{+} \circ D_{\mathcal{A}}^{\geq 1, +} \right) \left(
			\G{b}[\geq 1][\mathcal{A}]
			\right)
			\\
			\stackrel{\eqref{eq:D-G-geq-1}}{=}   &
			\pi_{\geq 2}^{+} \left( \Ho{b}[\geq 1][\mathcal{A}] \right)
			\stackrel{\eqref{eq:Ho-geq-2-def}}{=}
			\Ho{b}[\geq 2][\mathcal{A}]
		\end{aligned}
	\end{equation*}
	which shows \cref{eq:D-G-geq-2}. The proof of \cref{eq:cindmap-f-G-geq-2} is similar.
	\end{proof}

	To summarize, given a Banach $\Ainf$-algebra $\mathcal{A} = \left( A, \mu \right)$
	over a graded-commutative Banach $\mathbbm{k}$-algebra $\mathcal{R}$,
	we have constructed a chain complex $\totcompe{\mathcal{A}}[2][]$ over $\mathcal{R}$.
	Also, given another Banach $\Ainf$-algebra $\mathcal{B} = \left( B, \nu \right)$
	and a Banach $\Ainf$-morphism $f \colon \mathcal{A} \rightarrow \mathcal{B}$,
	possibly with a change of connection term $f_0 \left( 1 \right)$,
	we have constructed an induced \textit{chain map}
$\cindmape{f} \colon \totcompe{\mathcal{A}}[2][] \rightarrow \totcompe{\mathcal{B}}[2][]$,
	extending the chain map
$\cindmap{f} \colon \totcomp{\mathcal{A}}[2][] \rightarrow \totcomp{\mathcal{B}}[2][]$.

	\begin{rem} \label{rem:totcompe-geq-2-not-functorial-nose}
	Similarly to the situation with $\totcompe{\cdot}[1][]$ described in \cref{rem:totcompe-geq-1-not-functorial-nose},
	given another Banach $\Ainf$-morphism $g \colon \mathcal{B} \rightarrow \mathcal{C}$,
	the weak naturality property \eqref{eq:cindmap-f-G-geq-2} implies that we have
	\begin{equation} \label{eq:totcompe-geq2-func-basic}
		\left( \cindmape{g} \circ \cindmape{f} \right) \left( x \right) =
		\left( \cindmap{g} \circ \cindmap{f} \right) \left( x \right) =
		\cindmap{ \left( g \circ f \right) } \left( x \right) =
		\cindmape{ \left( g \circ f \right) } \left( x \right)
	\end{equation}
	when $x \in \totcomp{\mathcal{A}}[2][] \subseteq \totcompe{\mathcal{A}}[2][]$ and
	\begin{equation} \label{eq:totcompe-geq2-func-ul1}
		\left( \cindmape{g} \circ \cindmape{f} \right) \left( \ul{1} \right) =
		\cindmape{ \left( g \circ f \right) } \left( \ul{1} \right) +
		D_{\mathcal{C}} \left( R \left( f, g \right) \right)
	\end{equation}
	for some
$R \left( f, g \right) = R \left( \mcfunc{f} \left( 0 \right); g \right) \in
\totcomp{C}[2][-2] \subset \totcompe{C}[2][-2]$, which depends on
$f$, via $f_0 \left( 1 \right)$, and $g$.
	Hence, the construction
	\begin{align*}
		                                                & \mathcal{A}                              &  & \mapsto &   & \totcompe{\mathcal{A}}[2][],
		\\
		f \colon \mathcal{A}                            & \rightarrow \mathcal{B}                  &  & \mapsto &
		\cindmape{f} \colon \totcompe{\mathcal{A}}[2][] & \rightarrow \totcompe{\mathcal{B}}[2][],
	\end{align*}
	is not functorial on the level of chain complexes, unless one restricts attention to Banach
$\Ainf$-morphisms without change of connection terms. However, the construction is functorial
	\textit{on the level of cohomology}, i.e., we have
$H \left( \cindmape{g} \circ \cindmape{f} \right) = H \left( \cindmape{g} \right) \circ H \left( \cindmape{f} \right)$
	as maps from $\cohom{{\totcompe{\mathcal{A}}[2][]}}[]$ to $\cohom{{\totcompe{\mathcal{C}}[2][]}}[]$.
	\end{rem}

	\begin{rem}
	The total complex $\totcomp{\mathcal{A}}[1][]$ is homotopy equivalent to the
	shift ${\ncdfr{\mathcal{A}}[0][]}[1]$ of Connes' cyclic complex $\ncdfr{\mathcal{A}}[0][]$,
	and the extended total complex $\totcompe{\mathcal{A}}[1][]$ is homotopy equivalent to a
	shift ${\ncdf{\mathcal{A}}[0][]}[1]$ of the extended Connes complex $\ncdf{\mathcal{A}}[0][]$.

	When $\mathcal{A}$ is unital, the total complex $\totcomp{\mathcal{A}}[2][]$ is
	also homotopy equivalent to a different shift
${\ncdfr{\mathcal{A}}[0][]}[3]$ of Connes' cyclic complex.\footnote{This is a consequence of
	\cref{lm:ncdf-0-ncdf-2-mod-q3-equiv,lm:p_2-homotopy-equivalence} of \cref{sec:bicomplex-models-cyclic-homology}.
	Note that in \cref{sec:bicomplex-models-cyclic-homology}, we worked with
$\totcomp{\mathcal{A}}[2][] = \totc{\ncdf{A}, \qdr, \clie{\mu}}[][\geq 2][]$ and used the inner product parity form,
	but both constructions are explicitly isomorphic, see \cref{sec:description-using-braid-op-1}
	and \cref{sec:dependence-ndf-braidop}.}
	However, the extended total complex $\totcompe{\mathcal{A}}[2][]$ \textit{is not} homotopy equivalent
	to the shift ${\ncdf{\mathcal{A}}[0][]}[3]$. Instead, it is homotopy equivalent to the shift
$( \ncdfr{\mathcal{A}}[0][] )^{+}[3]$ of a different extension $( \ncdfr{\mathcal{A}}[0][] )^{+}$
	of $\ncdfr{\mathcal{A}}[0][]$. While the extension $\ncdf{\mathcal{A}}[0][]$ is obtained from
$\ncdfr{\mathcal{A}}[0][]$ by adjoining a generator $1$ of degree $0$ and extending the differential by
$1 \mapsto \mu_0 \left( 1 \right)$, the extension $( \ncdfr{\mathcal{A}}[0][] )^{+}$
	is obtained by adjoining a generator $S(1)$ of degree two and extending the differential by
$S(1) \mapsto S \left( \mu_0 \left( 1 \right) \right)$, where $S \left( \mu_0 \left( 1 \right) \right)$
	is the image of $\mu_0 \left( 1 \right)$ under the periodicity operator $S$,
	a consequence of the diagram of \cref{lm:periodicity-operator-as-projection}.
	\end{rem}

	For completeness, we also state the analogous version of \cref{lm:Go-and-Ho-geq1-properties-general},
	describing the properties of the chains $\Go{b}[\geq 2]$ and $\Ho{b}[\geq 2]$ in
	the \textit{non-extended} total complex $\totcomp{\cdot}[2][]$:
	\begin{lm} \label{lm:Go-and-Ho-geq2-properties-general}
	Let $\mathcal{A} = \left( A, \mu \right)$ be a Banach $\Ainf$-algebra over a differential
	graded-commutative Banach $\mathbbm{k}$-algebra $\mathcal{R} = \left( R, d_R \right)$
	and let $\mathcal{B} = \left( B, \nu \right)$ be a Banach $\Ainf$-algebra over a
	differential graded-commutative Banach $\mathbbm{k}$-algebra $\mathcal{S} = \left( S, d_S \right)$.
	Let $f \colon \mathcal{A} \rightarrow \mathcal{B}$ be a Banach $\Ainf$-morphism.
	Given a topologically nilpotent element $b \in \tc{A}$, we have the following identities between chains
	of $\totcomp{\cdot}[2][]$:
	\begin{align}
		D_{\mathcal{A}} \left( \Ho{b}[\geq 2][\mathcal{A}] \right) & = 0,
		\label{eq:D-Ho-geq-2}
		\\
		D_{\mathcal{A}} \left( \Go{b}[\geq 2][\mathcal{A}] \right) & =
		\Ho{b}[\geq 2][\mathcal{A}] - \Ho{0}[\geq 2][\mathcal{A}],
		\label{eq:D-Go-geq-2}
		\\
		\cindmap{f} \left( \Ho{b}[\geq 2][\mathcal{A}] \right)     & =
		\Ho{\mcfunc{f} \left( b \right)}[\geq 2][\mathcal{B}].
		\label{eq:cindmap-f-Ho-geq-2}
	\end{align}
	We also have the identity
	\begin{equation} \label{eq:cindmap-f-Go-geq-2}
		\cindmap{f} \left( \Go{b}[\geq 2][\mathcal{A}] \right) =
		\Go{\mcfunc{f} \left( b \right)}[\geq 2][\mathcal{B}] -
		\Go{\mcfunc{f} \left( 0 \right)}[\geq 2][\mathcal{B}] +
		D_{\mathcal{B}} \left( R \left( b; f \right) \right)
	\end{equation}
	for some chain $R \left( b; f \right) \in \totcomp{B}[2][-2]$.
	When $b = 0$ or when $f$ is strict, one can choose $R \left( b; f \right) = 0$. \qed
	\end{lm}
	\begin{proof}
	The lemma follows immediately from \cref{lm:Go-and-Ho-geq1-properties-general} using
	the facts that $\cindmap{f}[\cdot]$ commutes with $\pi_{\geq 2}$, and that
$\pi_{\geq 2}$ is a chain map, as in the proof of \cref{lm:G-geq-2-properties}.
	\end{proof}

	\begin{rem} \label{rem:Go-Ho-geq2-vs-0}
	As in \cref{rem:Go-Ho-geq1-vs-0},
	the properties of the chains $\Go{b}[\geq 2][\mathcal{A}]$ and $\Ho{b}[\geq 2][\mathcal{A}]$
	belonging to $\totcomp{\mathcal{A}}[2][]$ stated in
	\cref{lm:Go-and-Ho-geq2-properties-general}
	are almost identical to the properties of the chains $\Go{b}[0][\mathcal{A}]$ and $\Ho{b}[0][\mathcal{A}]$
	belonging to $\ncdfr{\mathcal{A}}[0][]$ stated in \cref{lm:Go_0-and-Ho_0-properties}, except
	that \cref{eq:cycl-f-Go-0} holds only up to an exact term, so it is replaced by
	\cref{eq:cindmap-f-Go-geq-2}.
	\end{rem}

	\subsection{Reduced Versions of \texorpdfstring{$\totcompe{\mathcal{A}}[2][]$}{the Extended Total Complex with Line Degree Greater Than One}} \label{sec:extended-reduced-total-complexes}
	In what follows, we will need reduced versions of the extended total complex
$\totcompe{\mathcal{A}}[2][]$, defined when $\mathcal{A}$ is unital. Let us start
	with some preliminaries.
	Let $V$ be a graded $R$-module. Given an element $x \in V$, we can define the $R$-linear coderivation
$\nu_x$ on $\tens{V}$ of degree $\degb{x}$ by
	\begin{equation*}
		\nu_x \left( v_1 \otimes \dots \otimes v_k \right) \defeq
		\sum_{i=1}^{k+1}
		(-1)^{\braid{\degb{x}}{\degb{v_1} + \dots + \degb{v_{i-1}}}}
		v_1 \otimes \dots \otimes v_{i-1} \otimes x \otimes v_i \otimes \dots \otimes v_k.
	\end{equation*}
	In terms of the corestriction $\corest{ \left( \nu_x \right)}$, we have $\left( \nu_x \right)_0(1) = x$, while
$\left( \nu_x \right)_k = 0$ for $k \geq 1$.
	Since $\nu_x$ increases the weight degree, we have $\corest{\left( \nu_x \right)} \circ \nu_x = 0$.
	When $x$ is odd in the sense that $\braid{x}{x} = 1$, this implies that
$\nu_x^2 = \frac{1}{2} \left[ \nu_x, \nu_x \right] = 0$.

	Now let $\mathcal{A} = (A,\mu,e)$ be a unital Banach $\Ainf$-algebra over the differential graded-commutative
	Banach $\mathbbm{k}$-algebra $\mathcal{R} = (R,d)$. Then the degree $-1$ coderivation $\nu_e$ on $\tens{A}$ satisfies
$\nu_e^2 = 0$.

	\begin{lm}
	The coderivation $\nu_e$ commutes with $\mu$:
$\left[ \nu_e, \mu \right] = \nu_e \circ \mu + \mu \circ \nu_e = 0$.
	\end{lm}
	\begin{proof}
	It is enough to verify that
$\corest{\left[ \nu_e, \mu \right]} = \corest{ \left( \nu_e \right)} \circ \mu + \corest{\mu} \circ \nu_e = 0$.
	We have
	\begin{equation*}
		\left[ \nu_e, \mu \right]_k \left( a_1, \dots, a_k \right) =
		\sum_{i=1}^{k+1} (-1)^{\degb{a_1} + \dots + \degb{a_{i-1}}}
		\mu_{k+1} \left( a_1, \dots, a_{i-1}, e, a_i, \dots, a_k \right) = 0
	\end{equation*}
	for all $k \geq 0$, by the fact that $e$ is a unit for $\mu$. Note that when $k \neq 1$, not only the sum
	above is zero but also each of the summands.
	\end{proof}

	\begin{lm} \label{lm:f-commutes-nu-e}
	Let $\mathcal{A} = \left( A, \mu, e_A \right)$ be a unital Banach $\Ainf$-algebra over $\mathcal{R} = \left( R, d_R \right)$
	and let $\mathcal{B} = \left( B, \nu, e_B \right)$ be a unital Banach $\Ainf$-algebra over $\mathcal{S} = \left( S, d_S \right)$.
	Given a unital $\Ainf$-morphism $f \colon \mathcal{A} \rightarrow \mathcal{B}$, we have
$f \circ \nu_{e_A} = \nu_{e_B} \circ f$.
	\end{lm}
	\begin{proof}
	It is enough to verify that $\corest{f} \circ \nu_{e_A} = \corest{ \left( \nu_{e_B} \right) } \circ f$.
	Indeed, we have
	\begin{equation*}
		\begin{aligned}
			\left( \corest{f} \circ \nu_{e_A} \right) \left( a_1, \dots, a_k \right) & =
			\sum_{i=1}^{k+1} (-1)^{\degb{a_1} + \dots + \degb{a_{i-1}}} f_{k+1} \left( a_1, \dots, a_{i-1}, e_A, a_i, \dots, a_k \right)
			\\
			                                                                         & = \delta_{k,0} f_1 \left( e_A \right) = \delta_{k,0} \, e_B
		\end{aligned}
	\end{equation*}
	when $k \geq 0$, by the unitality of $f$. Note again that when $k \neq 0$, not only the sum above is zero but also each of the summands.
	We also have
	\begin{equation*}
		\left( \corest{ \left( \nu_{e_B} \right) } \circ f \right) \left( a_1, \dots, a_k \right) =
		\delta_{k,0} \corest{ \left( \nu_{e_B} \right) } \left( f \left( 1 \right) \right) =
		\delta_{k,0} \corest{ \left( \nu_{e_B} \right) } \left( 1 \right) = \delta_{k,0} \, e_B,
	\end{equation*}
	as required.
	\end{proof}

	In what follows, we will slightly abuse notation and denote the coderivation $\nu_e$ on $\tens{A}$ simply by $e$,
	relying on context to differentiate between the unit element $e$ and the coderivation $\nu_e$ defined by $e$.
	Using our notation, we have $\left[ e, e \right] = 0$ and $\left[ \mu, e \right] = 0$.
	Since we work with the total degree parity form $\braidop_2$, we have:

	\begin{lm} \label{lm:coder-e-op-comm-relations}
	The following commutation relations hold on $\ncdf{A}[][]$:
	\begin{align}
		\left[ \ccont{e}, \clie{\mu} \right] & =
		\ccont{e} \circ \clie{\mu} - \clie{\mu} \circ \ccont{e} =
		0,
		\label{eq:lie-mu-cont-e-commcyc}
		\\
		\left[ \ccont{e}, \clie{e} \right]   & =
		\ccont{e} \circ \clie{e} - \clie{e} \circ \ccont{e} =
		0,
		\label{eq:lie-e-cont-e-commcyc}
		\\
		\left[ \qdr, \ccont{e}^{k+1} \right] & =
		\qdr \circ \ccont{e}^{k+1} - \ccont{e}^{k+1} \circ \qdr =
		(k+1) \left( \ccont{e}^k \circ \clie{e} \right).
		\label{eq:q-cont-e-commcyc}
	\end{align}
	\end{lm}
	\begin{proof}
	Since $e$ has degree $-1$, the operator $\ccont{e}$ of bidegree $(1,-1)$ has even total degree, hence
	the explicit form of the commutation relations.
	By the relation \eqref{eq:liecontcommcyc}, we have $\left[ \ccont{e}, \clie{\mu} \right] = \ccont{\left[ e, \mu \right]} = 0$,
	and $\left[ \ccont{e}, \clie{e} \right] = \ccont{\left[ e, e \right]} = 0$.
	Finally, by \cref{eq:commutator-power-identity,eq:lie-e-cont-e-commcyc}, we have
	\begin{equation*}
		\left[ \qdr, \ccont{e}^{k+1} \right] =
		\sum_{i=0}^k (-1)^{\braid{(-1,0)}{(i,-i)}}
		\ccont{e}^{i} \circ \left[ \qdr, \ccont{e} \right] \circ \ccont{e}^{k - i} =
		\sum_{i=0}^k \ccont{e}^i \circ \clie{e} \circ \ccont{e}^{k-i} =
		(k+1) \left( \ccont{e}^k \circ \clie{e} \right).
	\end{equation*}
	\end{proof}

	Given $k \geq 1$,
	denote by $\degenu[A][k][]$ the closure of the image of the map
	\begin{equation*}
		\rest{\ccont{e}^k}{\ncdf{A}[0][]} \colon \ncdf{A}[0][] \rightharpoonup \ncdf{A}[k][].
	\end{equation*}
	Cyclic codifferential forms in the image of $\rest{\ccont{e}^k}{\ncdf{A}[0][]}$
	are specific linear combinations of forms having the shape
	\begin{equation*}
		\ul{e} \otimes l^1 \otimes \dots \otimes \ul{e} \otimes l^k,
	\end{equation*}
	in which the only underlined element which appears is the unit $\ul{e}$.
	Note that in particular, we have $\ul{e}^{\otimes k} = \frac{1}{(k-1)!} \ccont{e}^k \left( 1 \right) \in \degenu[A][k][]$.
	For $x_k \in \degenu[A][k][]$ of the form $x_k = \ccont{e}^k \left( a_0 \right)$
	with $a_0 \in \ncdf{A}[0][]$, we have
	\begin{align*}
		\clie{\mu} \left( x_k \right) \stackrel{\eqref{eq:lie-mu-cont-e-commcyc}}{=}{} &
		\ccont{e}^k \left( \clie{\mu} \left( a_0 \right) \right) \in \degenu[A][k][],
		\\
		\qdr \left( x_k \right) \stackrel{\eqref{eq:q-cont-e-commcyc}}{=}{}            &
		k \cdot \ccont{e}^{k-1} \left( \clie{e} \left( a_0 \right) \right) \in
		\degenu[A][k-1][].
	\end{align*}
	Since $\clie{\mu}$ and $\qdr$ are contractive, and in particular continuous, the identities
	imply that we have
	\begin{equation} \label{eq:degenu-clie-qdr-invariance}
		\clie{\mu} \left( \degenu[A][k][] \right) \subseteq \degenu[A][k][],
		\qquad
		\qdr \left( \degenu[A][k][] \right) \subseteq \degenu[A][k-1][]
	\end{equation}
	for all $k \geq 1$.

	Now consider the graded Banach $R$-submodule of $\totcomp{A}[2][]$ given by
	\begin{equation} \label{eq:degenu-geq-2}
		\degenu[A][][\geq 2] \defeq \bigoplus_{k \geq 2} \degenu[A][k][].
	\end{equation}
	Given an element $x = \sum_{k \geq 2} x_k \in \degenu[A][][\geq 2]$ with $x_k \in \degenu[A][k][]$,
	we have $D_{\mathcal{A}} \left( x \right) = \sum_{k \geq 2} y_k$
	with $y_k = \clie{\mu} \left( x_k \right) - \qdr \left( x_{k+1} \right)$.
	By \eqref{eq:degenu-clie-qdr-invariance}, we have $y_k \in \degenu[A][k][]$ for all $k \geq 2$,
	so $\degenu[A][][\geq 2]$ is a subcomplex of $\totcomp{\mathcal{A}}[2][]$.

	\begin{dfn} \label{dfn:extended-reduced-total-complex}
	Let $\mathcal{A} = \left( A, \mu, e \right)$ be a unital Banach $\Ainf$-algebra over the differential graded-commutative
	Banach $\mathbbm{k}$-algebra $\mathcal{R} = \left( R, d \right)$. The quotient complex
	\begin{equation*}
		\totcompred{\mathcal{A}}[2][] \defeq \totcomp{\mathcal{A}}[2][] / \degenu[A][][\geq 2]
	\end{equation*}
	is called the \textbf{reduced total complex} associated to $\mathcal{A}$. Similarly, the quotient complex
	\begin{equation*}
		\totcompered{\mathcal{A}}[2][] \defeq \totcompe{\mathcal{A}}[2][] / \degenu[A][][\geq 2]
	\end{equation*}
	is called the \textbf{extended reduced total complex} associated to $\mathcal{A}$.
	\end{dfn}

	\begin{rem}	\label{rem:extended-reduced-functorial}
	Given a unital morphism $f \colon \mathcal{A} \rightarrow \mathcal{B}$ between two
	unital Banach $\Ainf$-algebras $\mathcal{A}$ and $\mathcal{B}$,
	\cref{lm:f-commutes-nu-e,cor:homo-cyc-func} show that the induced morphism
$\cindmap{f} \colon \ncdf{\mathcal{A}}[k][] \rightarrow \ncdf{\mathcal{B}}[k][]$
	maps $\degenu[A][k][]$ to $\degenu[B][k][]$. In particular, the induced map $\cindmap{f}$
	between the total complexes
	maps $\degenu[A][][\geq 2]$ to $\degenu[B][][\geq 2]$.
	Hence, the constructions $\mathcal{A} \mapsto \totcompred{\mathcal{A}}[2][]$
	(resp.\ $\mathcal{A} \mapsto \totcompered{\mathcal{A}}[2][]$)
	enjoy the same functorial properties
	as the non-reduced versions $\mathcal{A} \mapsto \totcomp{\mathcal{A}}[2][]$
	(resp.\ $\mathcal{A} \mapsto \totcompe{\mathcal{A}}[2][]$).
	\end{rem}

	Finally, we introduce a stronger version of the reduced total complex. Given $k \geq 2$, recall from
	\cref{sec:extended-reduced-cyclic-complexes} the graded Banach $R$-module
$\degen[A][k][] \subseteq \ncdf{\mathcal{A}}[k][]$ generated by cyclic codifferential forms of line degree $k$
	which contain the unit $e$ as long as it is not underlined.
	By \cref{cor:star-e-star-1-t-invariant}, we have $\clie{\mu} \left( \degen[A][k][] \right) \subseteq \degen[A][k][]$,
	and we also have $\qdr \left( \degen[A][k][] \right) \subseteq \degen[A][k-1][]$.
	Thus, the graded Banach $R$-module
	\begin{equation}
		\degen[A][][\geq 2] \defeq \bigoplus_{k = 2}^{\infty} \degen[A][k][] \subseteq \totcomp{\mathcal{A}}[2][]
	\end{equation}
	is a subcomplex of the total complex $\totcomp{\mathcal{A}}[2][]$. Let
$\left< \degen[A][][\geq 2], \degenu[A][][\geq 2] \right>$ be the graded
	\textit{Banach} $R$-submodule of $\totcomp{\mathcal{A}}[2][]$ generated by the elements of
$\degen[A][][\geq 2]$ and $\degenu[A][][\geq 2]$.

	\begin{dfn} \label{dfn:strongly-extended-reduced-total-complex}
	Let $\mathcal{A} = \left( A, \mu, e \right)$ be a unital Banach $\Ainf$-algebra over the differential graded-commutative
	Banach $\mathbbm{k}$-algebra $\mathcal{R} = \left( R, d \right)$. The quotient complex
	\begin{equation*}
		\totcompsred{\mathcal{A}}[2][] \defeq \totcomp{\mathcal{A}}[2][] / \left< \degen[A][][\geq 2], \degenu[A][][\geq 2] \right>
	\end{equation*}
	is called the \textbf{strongly reduced total complex} associated to $\mathcal{A}$. Similarly, the quotient complex
	\begin{equation*}
		\totcompesred{\mathcal{A}}[2][] \defeq \totcompe{\mathcal{A}}[2][] / \left< \degen[A][][\geq 2], \degenu[A][][\geq 2] \right>
	\end{equation*}
	is called the \textbf{extended strongly reduced total complex} associated to $\mathcal{A}$.
	\end{dfn}

	\begin{rem} \label{rem:extended-strongly-reduced-functorial}
	Given a unital morphism $f \colon \mathcal{A} \rightarrow \mathcal{B}$ between two unital
	Banach $\Ainf$-algebras,
	\cref{rem:cindmap-maps-degen-to-degen} shows that the induced morphism
$\cindmap{f} \colon \ncdf{\mathcal{A}}[k][] \rightarrow \ncdf{\mathcal{B}}[k][]$
	maps $\degen[A][k][]$ to $\degen[B][k][]$. In particular, the induced map $\cindmap{f}$
	between the total complexes
	maps $\degen[A][][\geq 2]$ to $\degen[B][][\geq 2]$. Taking into account \cref{rem:extended-reduced-functorial}, we see
	that the constructions $\mathcal{A} \mapsto \totcompsred{\mathcal{A}}[2][]$ (resp.\ $\mathcal{A} \mapsto \totcompesred{\mathcal{A}}[2][]$)
	enjoy the same functorial properties
	as the non-reduced versions $\mathcal{A} \mapsto \totcomp{\mathcal{A}}[2][]$ (resp.\ $\mathcal{A} \mapsto \totcompe{\mathcal{A}}[2][]$).
	\end{rem}

	\subsection{The Cyclic Chern--Simons Form Associated to a Bounding Cochain} \label{sec:canonical-chains-totcompe-geq-2}
	Similar to the situation with the cyclic exponential $\G{b}[0]$ described in
	\cref{subsec:canonical-elements-ncdf-0}, the cyclic Chern--Simons form
$\G{b}[\geq 2]$ of \cref{dfn:cyclic-chern-simons-form-braidop-2} is closed when $b$ is a bounding cochain:
	\begin{lm}[Cyclic Chern--Simons Form is Closed for a Bounding Cochain]  \hfill
	\label{lm:G-geq-2-bounding-chain-closed}
	\begin{enumerate}
	\item Let $\mathcal{A} = \left( A, \mu \right)$ be a Banach $\Ainf$-algebra over
$\mathcal{R}$ and let $b \in \mc{\mathcal{A}}$ be a strong bounding cochain.
	Then $\G{b}[\geq 2]$ is closed in the extended total complex $\totcompe{\mathcal{A}}[2][]$.
	\item Let $\mathcal{A} = \left( A, \mu, e \right)$ be a unital Banach $\Ainf$-algebra over
$\mathcal{R}$ and let $b \in \mc{\mathcal{A}, c}$ be a weak bounding cochain.
	Then $\G{b}[\geq 2]$ is closed in the extended reduced total complex $\totcompered{\mathcal{A}}[2][]$.
	\label{item::G-geq-2-weak-bounding-chain-closed}
	\end{enumerate}
	\end{lm}
	\begin{proof}
	By \cref{lm:G-geq-2-properties}, we have
	\begin{equation*}
		D \left( \G{b}[\geq 2] \right) = \Ho{b}[\geq 2]
		\stackrel{\eqref{eq:Ho-geq-2-def}}{=}
		\sum_{n \geq 2} \H{b}[n]
		\stackrel{\eqref{eq:H_n(b)-formula}}{=}
		- \sum_{n=2}^{\infty}
		\frac{1}{n} \left( \ul{\corest{\mu} \left( \Exp{b} \right)} \otimes \Exp{b} \right)^{\otimes n}.
	\end{equation*}
	When $b \in \mc{\mathcal{A}}$ is a strong bounding cochain, we have
$\corest{\mu} \left( \Exp{b} \right) = 0$ and hence $\G{b}[\geq 2]$ is closed in
$\totcompe{\mathcal{A}}[2][-1]$.
	When $\mathcal{A}$ is unital and $b \in \mc{\mathcal{A}, c}$ is a weak bounding cochain,
	then $\corest{\mu} \left( \Exp{b} \right) = c \cdot e$,
	so
	\begin{equation*}
		\begin{aligned}
			D \left( \G{b}[\geq 2] \right) & =
			- \sum_{n=2}^{\infty}
			\frac{1}{n} \left( \ul{c \cdot e} \otimes \Exp{b} \right)^{\otimes n}
			=
			- \sum_{n=2}^{\infty}
			\frac{c^n}{n} \left( \ul{e} \otimes \Exp{b} \right)^{\otimes n}
			\\
			                               & =
			- \sum_{n=2}^{\infty} \frac{c^n}{n!} \ccont{e}^n \left( \G{b}[0] \right)
			\in \degenu[A][][\geq 2],
		\end{aligned}
	\end{equation*}
	which shows that $\G{b}[\geq 2]$ is closed, when considered as an element
	of the extended reduced total complex $\totcompered{\mathcal{A}}[2][]$.
	\end{proof}

	Let $\mathcal{A}_0$ and $\mathcal{A}_1$ be two Banach
$\Ainf$-algebras over the same differential graded-commutative $\mathbbm{k}$-algebra
$\mathcal{S}$.
	Applying \cref{dfn:F-strong-pseudoisotopy} to the homology of the extended total
	complex functor $\mathcal{F} = \cohom{}[] \totcompe{}[2]$, we obtain the notion of
	a $\cohom{}[] \totcompe{}[2]$-\textbf{strong} pseudoisotopy $\mathfrak{A}$
	between $\mathcal{A}_0$ and $\mathcal{A}_1$.
	Such a pseudoisotopy yields a canonical isomorphism
	\begin{equation*}
		\mathfrak{a} \colon \cohom{{\totcompe{\mathcal{A}_0 / \mathcal{S}}[2][]}}[] \overset{\sim}{\rightarrow} \cohom{{\totcompe{\mathcal{A}_1 / \mathcal{S}}[2][]}}[]
	\end{equation*}
	of graded $\cohom{\mathcal{S}}[]$-modules.

	Similarly, when $\mathcal{A}_0,\mathcal{A}_1$ and $\mathfrak{A}$ are unital,
	applying \cref{dfn:F-strong-pseudoisotopy} to the homology of the extended
	reduced total complex functor $\mathcal{F} = \cohom{}[] \totcompered{}[2]$ yields
	the notion of a $\cohom{}[] \totcompered{}[2]$-\textbf{strong} pseudoisotopy
$\mathfrak{A}$ between $\mathcal{A}_0$ and $\mathcal{A}_1$,
	which induces a canonical isomorphism
	\begin{equation*}
		\mathfrak{a} \colon \cohom{{\totcompered{\mathcal{A}_0 / \mathcal{S}}[2][]}}[] \overset{\sim}{\rightarrow} \cohom{{\totcompered{\mathcal{A}_1 / \mathcal{S}}[2][]}}[].
	\end{equation*}
	Then we have:
	\begin{thm}[Cyclic Chern--Simons Form Invariance under Strong Pseudoisotopy]
	\label{lm:invariance-Gb-geq-2-pseudoisotopy}
	Let $\mathcal{S}$ be a differential graded-commutative Banach $\mathbbm{k}$-algebra, and let
$\mathcal{A}_0$ and $\mathcal{A}_1$ be two non-unital (resp.\ unital) Banach $\Ainf$-algebras over
$\mathcal{S}$. Let $\mathfrak{A}$ be a $\cohom{}[] \totcompe{}[2]$-strong
	(resp.\ $\cohom{}[] \totcompered{}[2]$-strong) pseudoisotopy between $\mathcal{A}_0$
	and $\mathcal{A}_1$.
	Let $b_0 \in \tc{A_0}$ and $b_1 \in \tc{A_1}$ be
$\mathfrak{A}$-gauge-equivalent strong (resp.\ weak) bounding cochains.
	Then
	\begin{equation*}
		\mathfrak{a}(\eqcl{ \G{b_0}[\geq 2][\mathcal{A}_0] }) = \eqcl{ \G{b_1}[\geq 2][\mathcal{A}_1] }
		\textrm{ in } \cohom{{\totcompe{\mathcal{A}_1}[2][]}}[-1]
		\qquad \left(
		\textrm{resp.\ in } \cohom{{\totcompered{\mathcal{A}_1}[2][]}}[-1]
		\right).
	\end{equation*}
	\end{thm}
	\begin{proof}
	The proof is the same as for \cref{lm:invariance-Gb0-pseudoisotopy}, except we use
	the weak naturality property \eqref{eq:cindmap-f-G-geq-2} of $\G{\cdot}[\geq 2]$.

	Assume we are in the non-unital case.
	Let $b_0 \in \mc{\mathcal{A}_0}$ and $b_1 \in \mc{\mathcal{A}_1}$ be two $\mathfrak{A}$-gauge-equivalent bounding cochains.
	Choose $b \in \mc{\mathfrak{A}}$ with $\mcfunc{\evalmf}^i \left( b \right) = b_i$ for $i = 0, 1$.
	By \cref{lm:G-geq-2-bounding-chain-closed}, $\G{b}[\geq 2][\mathfrak{A}]$ is closed and
	defines a class $\eqcl{ \G{b}[\geq 2][\mathfrak{A}] } \in \cohom{{\totcompe{\mathfrak{A}}[2][]}}[-1]$.
	We have
	\begin{equation*}
		\begin{aligned}
			\cohom{}[] \totcompe{{\evalmf^i}}[2][] \left( \eqcl{ \G{b}[\geq 2][\mathfrak{A}] } \right)
			\stackrel{\phantom{\eqref{eq:cindmap-f-G-geq-2}}}{=}{} &
			\cohom{{\cindmape{F} \left( \evalmf^i \right)}}[]  \left( \eqcl{ \G{b}[\geq 2][\mathfrak{A}] } \right)
			=
			\eqcl{ \cindmape{F} \left( \evalmf^i \right) \left( \G{b}[\geq 2][\mathfrak{A}] \right) }
			\\
			\stackrel{\eqref{eq:cindmap-f-G-geq-2}}{=}{}           &
			\eqcl{ \G{\mcfunc{\evalmf}^i \left( b \right)}[\geq 2][\mathcal{A}_i] + D_{\mathcal{A}_i} \left( R \left( b; \evalmf^i \right) \right)}
			\\
			\stackrel{\phantom{\eqref{eq:cindmap-f-G-geq-2}}}{=}{} &
			\eqcl{ \G{\mcfunc{\evalmf}^i \left( b \right)}[\geq 2][\mathcal{A}_i] } =
			\eqcl{ \G{b_i}[\geq 2][\mathcal{A}_i] }
		\end{aligned}
	\end{equation*}
	for $i = 0, 1$, and hence,
$\mathfrak{a}(\eqcl{ \G{b_0}[\geq 2][\mathcal{A}_0] }) = \eqcl{ \G{b_1}[\geq 2][\mathcal{A}_1] }$.
	The proof for the unital case proceeds the same way, using the fact
	that $\G{b}[\geq 2][\mathfrak{A}]$ is a cycle in $\totcompered{\mathfrak{A}}[2][]$ by
	part \eqref{item::G-geq-2-weak-bounding-chain-closed} of \cref{lm:G-geq-2-bounding-chain-closed}.
	\end{proof}

	\begin{rem} \label{rem:Go-geq-2-closed-no-curvature}
	Assume that $b \in \mc{\mathcal{A}}$ is a strong bounding cochain. Then in general,
$\Go{b}[\geq 2]$ is not a cycle of $\totcomp{\mathcal{A}}[2][]$. Instead,
	we have
	\begin{equation*}
		\begin{aligned}
			D \left( \Go{b}[\geq 2] \right)
			\stackrel{\eqref{eq:G-geq-2-def}}{=}
			D \left( \G{b}[\geq 2] - \ul{1} \right)
			= -D \left( \ul{1} \right)
			\stackrel{\eqref{eq:action-D-geq-2-on-ul-1}}{=}
			\sum_{n=2}^{\infty} \frac{1}{n} \ul{\mu_0(1)}^{\otimes n}.
		\end{aligned}
	\end{equation*}
	Working in the extended total complex $\totcompe{\mathcal{A}}[2][]$
	allows us to cancel the contributions of $\mu_0 \left( 1 \right)$
	by adding $\ul{1}$ to $\Go{b}[\geq 2]$.
	Similarly, \cref{eq:cindmap-f-Go-geq-2} shows that $\Go{b}[\geq 2]$ is not natural in general,
	but is natural up to an exact term for morphisms with $f_0 \left( 1 \right) = 0$.

	If we restrict our
	attention to uncurved Banach $\Ainf$-algebras and morphisms without a change of connection term,
	\cref{lm:Go-and-Ho-geq2-properties-general} implies that $\Go{b}[\geq 2]$ gives us a
	version of the cyclic Chern--Simons form in the total complex $\totcomp{\mathcal{A}}[2][]$
	satisfying \cref{lm:G-geq-2-bounding-chain-closed,lm:invariance-Gb-geq-2-pseudoisotopy},
	without needing to work in the extended total complex $\totcompe{\mathcal{A}}[2][]$.
	Compare to \cref{rem:Go-0-closed-no-curvature}.
	\end{rem}

	\subsection{Total Inner Products and the Superpotential}
	\label{sec:generalized-inner-product-superpotential}
	Let $\mathcal{R} = (R,d)$ be a differential graded-commutative Banach $\mathbbm{k}$-algebra, and let
$\mathcal{A} = \left( A, \mu \right)$ be a Banach $\Ainf$-algebra over $\mathcal{R}$.

	\begin{dfn} \label{dfn:generalized-inner-product}
	An $n$-\textbf{dimensional total inner product} on $\mathcal{A}$ is a
	morphism
	\begin{equation*}
		\phi \colon \totcompe{\mathcal{A}}[2][] \rightarrow \mathcal{R}[4-n]
	\end{equation*}
	of differential graded Banach $\mathcal{R}$-modules between the \textit{extended} total complex
$\totcompe{\mathcal{A}}[2][]$ and $\mathcal{R}[4-n]$,
	i.e., a contractive degree zero $R$-linear map
$\phi \colon \totcompe{A}[2][] \rightarrow R[4-n]$ which satisfies
$d_{R[4-n]} \circ \phi = \phi \circ D_{\mathcal{A}}$. An
$n$-\textbf{dimensional pre-total inner product} on $\mathcal{A}$ is defined
	similarly but using the total complex $\totcomp{\mathcal{A}}[2][]$ instead of the
	extended total complex $\totcompe{\mathcal{A}}[2][].$
	An $n$-\textbf{dimensional total inner product Banach} $\Ainf$-\textbf{algebra over} $\mathcal{R}$
	is a triple $\mathcal{A} = \left( A, \mu, \phi \right)$ where $\left( A, \mu \right)$ is a
	Banach $\Ainf$-algebra over $\mathcal{R}$ and $\phi$ is an $n$-dimensional total inner product
	on $\left( A, \mu \right)$.
	\end{dfn}

	In what follows, the dimension $n$ will be fixed, and we often omit it for brevity.

	\begin{dfn} \label{dfn:total-inner-product-unital}
	Let $\mathcal{A} = \left( A, \mu, e \right)$ be a unital Banach $\Ainf$-algebra over $\mathcal{R}$.
	A total inner product $\phi$ on $\mathcal{A}$ is called \textbf{unital} (or $\degenu$-\textbf{unital})
	if
	\begin{equation} \label{eq:total-inner-product-i-e-cond}
		\phi \circ \rest{\ccont{e}^k}{\ncdf{A}[0][]} = 0
	\end{equation}
	for all $k \geq 2$, where $\ccont{e}$ is the cyclic contraction with the unit, described in \cref{sec:extended-reduced-total-complexes}.
	A total inner product $\phi$ on $\mathcal{A}$ is called \textbf{strongly unital} if it is unital, and in addition,
	\begin{equation}	\label{eq:total-inner-product-free-unit-cond}
		\phi \left(
		\ul{a_1} \otimes a_1^1 \otimes \dots \otimes a_1^{r_1} \otimes \ul{a_2} \otimes
		\dots \otimes \ul{a_k} \otimes a_k^1 \otimes \dots \otimes a_k^{r_k}
		\right) = 0
	\end{equation}
	whenever $k \geq 2$ and $r_1, \dots, r_k \geq 0$ and $a_i, a_i^j \in A$, such that $a_i^j = e$ for some $1 \leq i \leq k$ and $j > 0$.

	Equivalently, a unital (resp.\ strongly unital) total inner product $\phi$
	is a morphism $\phi \colon \totcompered{\mathcal{A}}[2][] \rightarrow \mathcal{R}[4-n]$
	(resp.\ $\phi \colon \totcompesred{\mathcal{A}}[2][] \rightarrow \mathcal{R}[4-n]$)
	between the \textit{extended reduced} (resp.\ \textit{extended strongly reduced}) total complex
$\totcompered{\mathcal{A}}[2][]$ (resp.\ $\totcompesred{\mathcal{A}}[2][]$) and
$\mathcal{R}[4-n]$. See \cref{sec:extended-reduced-total-complexes}.

	A \textbf{(strongly) unital total inner product Banach} $\Ainf$-\textbf{algebra}
	over $\mathcal{R}$ is a quadruple $\mathcal{A} = \left( A, \mu, e, \phi \right)$ where
$\left( A, \mu, e \right)$ is a unital Banach $\Ainf$-algebra over $\mathcal{R}$ and $\phi$ is a (strongly) unital
	total inner product on $\mathcal{A}$.
	\end{dfn}

	Next, we introduce a notion of morphisms for total inner product Banach $\Ainf$-algebras. Given
	a differential graded-commutative Banach $\mathbbm{k}$-algebra $\mathcal{R} = \left( R, d \right)$, let us denote
	by $\pi_R \colon R \rightarrow R / dR$ the projection onto the equivalence classes of cochains modulo coboundaries.
	Also, given a morphism $\varphi \colon \left( R, d_R \right) \rightarrow \left( S, d_S \right)$ between two
	differential graded-commutative Banach $\mathbbm{k}$-algebras, denote by
$\varphi^{\diamond} \colon R / d_R R \rightarrow S / d_S S$ the induced map. Given a total inner product
$\phi \colon \totcompe{\mathcal{A}}[2][] \rightarrow \mathcal{R}[4-n]$
	on a Banach $\Ainf$-algebra $\mathcal{A}$ over $\mathcal{R}$, set
	\begin{equation} \label{eq:phi-diamond-def}
		\phi^{\diamond} \defeq \pi_R[4 - n] \circ \phi \colon \totcompe{\mathcal{A}}[2][] \rightarrow \left( R / dR \right)[4-n].
	\end{equation}

	\begin{dfn} \label{dfn:morphism-total-inner-products}
	Given two total inner product Banach $\Ainf$-algebras $\mathcal{A} = \left( A, \mu_A, \phi_A \right)$ over
$\mathcal{R} = \left( R, d_R \right)$ and $\mathcal{B} = \left( B, \mu_B, \phi_B \right)$ over
$\mathcal{S} = \left( S, d_S \right)$,
	a \textbf{morphism of total inner product Banach} $\Ainf$-\textbf{algebras} is a
	morphism $f \colon \left( A, \mu_A \right) \rightarrow \left( B, \mu_B \right)$ of Banach $\Ainf$-algebras such that
	\begin{equation} \label{eq:morphism-total-inner-products}
		\phi_B^{\diamond} \circ \cindmape{f} =
		\base{f}^{\diamond}[4 - n] \circ \phi_A^{\diamond}
	\end{equation}
	(see \cref{fig:morphism-total-inner-product}).

	When $\mathcal{A}$ and $\mathcal{B}$ are (strongly) unital, we require that the morphism $f$ is also unital.
	\end{dfn}

	\begin{figure}
	\begin{tikzcd}
		{\totcompe{\mathcal{A} / \mathcal{R}}[2][]} && {\mathcal{R}[4-n]} && {\left( R / dR \right)[4-n]}
		\\
		{\totcompe{\mathcal{B} / \mathcal{S}}[2][]} && {\mathcal{S}[4-n]} && {\left( S / dS \right)[4-n]}
		\arrow["{\phi_A}", from=1-1, to=1-3]
		\arrow["{\phi_A^{\diamond}}", curve={height=-20pt}, from=1-1, to=1-5]
		\arrow["\cindmape{f}"', from=1-1, to=2-1]
		\arrow["{\pi_R[4-n]}", from=1-3, to=1-5]
		\arrow["{{\base{f}^{\diamond}[4-n]}}", from=1-5, to=2-5]
		\arrow["{\phi_B}"', from=2-1, to=2-3]
		\arrow["{\phi_B^{\diamond}}"', curve={height=20pt}, from=2-1, to=2-5]
		\arrow["{\pi_S[4-n]}"', from=2-3, to=2-5]
	\end{tikzcd}
	\caption{Morphism of total inner product Banach $\Ainf$-algebras.}
	\label{fig:morphism-total-inner-product}
	\end{figure}

	\begin{rem}
	While the more natural condition is $\phi_B \circ \cindmape{f} = \base{f}[4 - n] \circ \phi_A$, this condition
	is not preserved under composition because the construction of the \textit{extended} total complex
$\mathcal{A} \mapsto \totcompe{\mathcal{A}}[2][]$ is not functorial (see \cref{rem:totcompe-geq-2-not-functorial-nose}).
	The weaker condition of \cref{dfn:morphism-total-inner-products}, requires only that $\phi_B \circ \cindmape{f}$
	coincides with $\base{f}[4 - n] \circ \phi_A$ up to exact terms and guarantees that the composition of morphisms is a morphism.

	Namely, let $\mathcal{C} = \left( C, \mu_C, \phi_C \right)$ be a total inner product Banach $\Ainf$-algebra
	over the differential graded-commutative Banach $\mathbbm{k}$-algebra $\mathcal{T} = \left( T, d_T \right)$
	and let $g \colon \mathcal{B} \rightarrow \mathcal{C}$ be a morphism
	of total inner product Banach $\Ainf$-algebras. Let us verify the identity
	\begin{equation} \label{eq:proj-composition-total-inner-product-morphisms}
		\left( \phi_C^{\diamond} \circ \cindmape{g} \circ \cindmape{f} \right) \left( x \right) =
		\left( \phi_C^{\diamond} \circ \cindmape{\left( g \circ f \right)} \right) \left( x \right)
	\end{equation}
	for $x \in \totcompe{\mathcal{A}}[2][]$.

	When $x \in \totcomp{\mathcal{A}}[2][] \subseteq \totcompe{\mathcal{A}}[2][]$, by \cref{eq:totcompe-geq2-func-basic}, we have
$\left(  \cindmape{g} \circ \cindmape{f} \right) \left( x \right) = \cindmape{ \left( g \circ f \right)} \left( x \right)$
	which in particular implies \cref{eq:proj-composition-total-inner-product-morphisms}.
	When $x = \ul{1}$, by \cref{eq:totcompe-geq2-func-ul1}, and using the fact that $\phi_C$ is a chain map, we have
	\begin{equation*}
		\begin{aligned}
			\left( \phi_C^{\diamond} \circ \cindmape{g} \circ \cindmape{f} \right) \left( \ul{1} \right) & =
			\phi_C^{\diamond} \left(
			\cindmape{ \left( g \circ f \right) } \left( \ul{1} \right) +
			D_{\mathcal{C}} \left( R \left( f, g \right) \right)
			\right)
			\\
			                                                                                             & =
			\left( \phi_C^{\diamond} \circ \cindmape{\left( g \circ f \right)} \right) \left( \ul{1} \right)
			+
			\left( \pi_T[4-n] \circ \phi_C \right) \left( D_{\mathcal{C}} \left( R \left( f, g \right) \right) \right)
			\\
			                                                                                             & =
			\left( \phi_C^{\diamond} \circ \cindmape{\left( g \circ f \right)} \right) \left( \ul{1} \right)
			+
			\left( \pi_T[4-n] \circ d_{T[4-n]} \circ \phi_C \right) \left( R \left( f, g \right) \right)
			\\
			                                                                                             & =
			\left( \phi_C^{\diamond} \circ \cindmape{\left( g \circ f \right)} \right) \left( \ul{1} \right).
		\end{aligned}
	\end{equation*}

	Hence, we have
	\begin{equation*}
		\begin{aligned}
			\base{\left( g \circ f \right)}^{\diamond}[4 - n] \circ \phi_A^{\diamond} & =
			\base{g}^{\diamond}[4 - n] \circ \base{f}^{\diamond}[4 - n] \circ \phi_A^{\diamond}
			\\
			                                                                          & =
			\base{g}^{\diamond}[4 - n] \circ \phi_B^{\diamond} \circ \cindmape{f}
			\\
			                                                                          & =
			\phi_C^{\diamond} \circ \cindmape{g} \circ \cindmape{f}
			\\
			                                                                          & =
			\phi_C^{\diamond} \circ \cindmape{\left( g \circ f \right)}
		\end{aligned}
	\end{equation*}
	which shows that $g \circ f$ is a morphism of total inner product Banach $\Ainf$-algebras.
	\end{rem}

	\begin{dfn} \label{dfn:superpotential}
	Let $\mathcal{A} = \left( A, \mu, \phi \right)$ be an
$n$-dimensional total inner product Banach $\Ainf$-algebra and let $b \in \tc{A}$ be
	a topologically nilpotent element. The \textbf{superpotential} function
$\SP \colon \tc{A} \rightarrow R^{3-n}$ associated to $\mathcal{A}$ is defined by
	\begin{equation} \label{eq:def-superpotential}
		\begin{aligned}
			\SP[b] & = \SP[b][\mathcal{A}] \defeq \phi \left( \G{b}[\geq 2][\mathcal{A}] \right)
			\\
			       & =
			\phi \left( \ul{1} \right) +
			\sum_{\substack{k = 2 \\ j_1, \dots, j_k = 0 \\ i_2, \dots, i_k = 0}}^{\infty}
			\frac{1}{1 + \sum_{r = 1}^k j_r + \sum_{r=2}^{k} i_r}
			\\
			       & \phantom{=
				         \phi \left( \ul{1} \right) +
				         \sum_{\substack{k = 2 \\ j_1, \dots, j_k = 0 \\ i_2, \dots, i_k = 0}}^{\infty}
			         }
			\phi \left(
			\ul{b} \otimes b^{\otimes j_1} \otimes \ul{ \mu_{i_2} \left( b^{\otimes i_2} \right)}
			\otimes b^{\otimes j_2} \otimes \dots \otimes \ul{ \mu_{i_k} \left( b^{\otimes i_k} \right) }
			\otimes b^{\otimes j_k}
			\right).
		\end{aligned}
	\end{equation}
	\end{dfn}

	\begin{lm}[Naturality of the Superpotential up to Exact Terms] \label{lm:functoriality-superpotential}
	Let $\mathcal{A} = \left( A, \mu, \phi_A \right)$ be an $n$-dimensional total inner product
	Banach $\Ainf$-algebra over a differential graded-commutative Banach $\mathbbm{k}$-algebra
$\mathcal{R} = \left( R, d_R \right)$ and let $\mathcal{B} = \left( B, \nu, \phi_B \right)$
	be an $n$-dimensional total inner product Banach $\Ainf$-algebra
	over a differential graded-commutative Banach $\mathbbm{k}$-algebra $\mathcal{S} = \left( S, d_S \right)$.
	Given a morphism $f \colon \mathcal{A} \rightarrow \mathcal{B}$ of total inner product Banach $\Ainf$-algebras
	and $b \in \tc{A}$, we have
	\begin{equation} \label{eq:sp-pseudo-invariance}
		\base{f} \left( \SP[b][\mathcal{A}] \right) =
		\SP[\mcfunc{f} \left( b \right)][\mathcal{B}] + d_S \left( x \right)
	\end{equation}
	for some $x$ in $S^{2-n}$.
	\end{lm}
	\begin{proof}
	Using the fact that $\phi$ is a chain map, we have
	\begin{equation*}
		\begin{aligned}
			\eqcl{ \base{f} \left( \SP[b][\mathcal{A}] \right) }_{S/dS}
			\stackrel{\phantom{\eqref{eq:morphism-total-inner-products}}}{=}{}        &
			\base{f}^{\diamond} \left( \eqcl{ \SP[b][\mathcal{A}] }_{R/dR} \right)
			\\
			\stackrel{\eqref{eq:def-superpotential}}{=}{}                             &
			\base{f}^{\diamond} \left( \eqcl{ \phi_A \left( \G{b}[\geq 2][\mathcal{A}] \right)}_{R/dR} \right)
			\\
			\stackrel{\eqref{eq:phi-diamond-def}}{=}{}                                &
			\left( \base{f}^{\diamond}[4-n] \circ \phi_A^{\diamond} \right) \left( \G{b}[\geq 2][\mathcal{A}] \right)
			\\
			\stackrel{\eqref{eq:morphism-total-inner-products}}{=}{}                  &
			\phi_B^{\diamond} \left( \cindmape{f} \left( \G{b}[\geq 2][\mathcal{A}] \right) \right)
			\\
			\stackrel{\eqref{eq:cindmap-f-G-geq-2}}{=}{}                              &
			\phi_B^{\diamond} \left(
			\G{\mcfunc{f} \left( b \right)}[\geq 2][\mathcal{B}] +
			D_{\mathcal{B}} \left( R \left( b; f \right) \right)
			\right)
			\\
			\stackrel[\eqref{eq:phi-diamond-def}]{\eqref{eq:def-superpotential}}{=}{} &
			\eqcl{
				\SP[\mcfunc{f} \left( b \right)][\mathcal{B}]  +
				d_{S[4-n]} \left( \phi_B \left( R \left( b; f \right) \right) \right)
			}_{S/dS}
			=
			\eqcl{ \SP[\mcfunc{f} \left( b \right)][\mathcal{B}] }_{S/dS}
		\end{aligned}
	\end{equation*}
	which proves \cref{eq:sp-pseudo-invariance}.
	\end{proof}

	\begin{lm}[Superpotential of a Bounding Cochain is a Cocycle] \label{lm:superpotential-bounding-chain-closed}
	Let $\mathcal{R}$ be a differential graded-commutative Banach $\mathbbm{k}$-algebra,
	and let $\mathcal{A}$ be a non-unital (resp.\ unital) $n$-dimensional total inner product
	Banach $\Ainf$-algebra over $\mathcal{R}$. Let $b \in \tc{A}$ be a strong (resp.\ weak) bounding cochain.
	Then $\SP[b]$ is a cocycle which defines a cohomology class
$\eqcl{ \SP[b] } \in \cohom{\mathcal{R}}[3-n]$.
	\end{lm}
	\begin{proof}
	When $\mathcal{A}$ is non-unital (resp.\ unital) total inner product Banach $\Ainf$-algebra,
	and $b$ is a strong (resp.\ weak) bounding cochain, \cref{lm:G-geq-2-bounding-chain-closed} shows that
	the cyclic Chern--Simons form $\G{b}[\geq 2]$ is closed in $\totcompe{\mathcal{A}}[2][-1]$ (resp.\ $\totcompered{\mathcal{A}}[2][-1]$).
	Since $\SP[b] = \phi(\G{b}[\geq 2])$ is the application of a chain map of degree $4 - n$
	to a closed element of degree $-1$, the result follows.
	\end{proof}

	\Cref{dfn:pseudo-isotopy} of pseudoisotopy between two Banach $\Ainf$-algebras
	extends naturally to the setting of total inner product Banach $\Ainf$-algebras as follows.

	\begin{dfn}[Pseudoisotopy of Total Inner Product Banach $\Ainf$-algebras]
	\label{dfn:pseudo-isotopy-total-inner-product-algebras}
	Let $\mathcal{A}_0$ and $\mathcal{A}_1$ be two Banach $\Ainf$-algebras over $\mathcal{S}$
	and let $\mathfrak{A}$ be a pseudoisotopy over $\mathfrak{R}$
	between $\mathcal{A}_0$ and $\mathcal{A}_1$.
	Assume that $\mathcal{A}_i$ are equipped with total inner products
$\phi^i \colon \totcompe{\mathcal{A}_i}[2][] \rightarrow \mathcal{S}[4-n]$
	for $i=0,1$, and that $\mathfrak{A}$ is also equipped with a total inner product
$\phi \colon \totcompe{\mathfrak{A}}[2][] \rightarrow \mathfrak{R}[4-n]$ on $\mathfrak{A}$
	such that the pseudoisotopy maps $\evalmf^i \colon \mathfrak{A} \rightarrow \mathcal{A}_i$
	become morphisms of total inner product Banach $\Ainf$-algebras (see \cref{fig:pseudoisotopy-total-inner-product-algebras}).
	In this case, we say that $\mathfrak{A}$ is a
	\textbf{pseudoisotopy of total inner product Banach} $\Ainf$-\textbf{algebras}.

	When $\mathcal{A}_0, \mathcal{A}_1$ and $\mathfrak{A}$ are unital, and the total inner products
$\phi^0,\phi^1$ and $\phi$ are also (strongly) unital, we say that $\mathfrak{A}$ is a
	\textbf{(strongly) unital pseudoisotopy of (unital) total inner product Banach}
$\Ainf$-\textbf{algebras}.
	\end{dfn}

	\begin{figure}
	\begin{tikzcd}
		{\totcompe{\mathfrak{A} / \mathfrak{R}}[2][]} && {\mathfrak{R}[4 - n]} && {\left( R / dR \right)[4-n]}
		\\
		{\totcompe{\mathcal{A}_i / \mathcal{S}}[2][]} && {\mathcal{S}[4-n]} && {\left( S / dS \right)[4-n]}
		\arrow["{\phi}", from=1-1, to=1-3]
		\arrow["{\phi^{\diamond}}", curve={height=-20pt}, from=1-1, to=1-5]
		\arrow["\cindmape{\evalmf^i}"', from=1-1, to=2-1]
		\arrow["{\pi_R[4-n]}", from=1-3, to=1-5]
		\arrow["{{\evalm^{i,\diamond}}[4-n]}", from=1-5, to=2-5]
		\arrow["{\phi^i}"', from=2-1, to=2-3]
		\arrow["{\phi^{i, \diamond}}"', curve={height=20pt}, from=2-1, to=2-5]
		\arrow["{\pi_S[4-n]}"', from=2-3, to=2-5]
	\end{tikzcd}
	\caption{Pseudoisotopy of total inner product Banach $\Ainf$-algebras.}
	\label{fig:pseudoisotopy-total-inner-product-algebras}
	\end{figure}

	\begin{thm}[Invariance of the Superpotential under Gauge Equivalence]
	\label{thm:invariance-superpotential}
	Let $\mathcal{A}_0$ and $\mathcal{A}_1$ be two non-unital (resp.\ unital) $n$-dimensional
	total inner product Banach $\Ainf$-algebras over $\mathcal{S} = \left( S, d_S \right)$
	and let $\mathfrak{A}$ be a non-unital (resp.\ unital) pseudoisotopy
	of total inner product Banach $\Ainf$-algebras
	between $\mathcal{A}_0$ and $\mathcal{A}_1$, defined over $\mathfrak{R} = \left( R, d_R \right)$.
	Let $b_0 \in \tc{A_0}$ and $b_1 \in \tc{A_1}$ be two
$\mathfrak{A}$-gauge-equivalent strong (resp.\ weak) bounding cochains.
	Then
	\begin{equation*}
		\eqcl{ \SP[b_0][\mathcal{A}_0] } = \eqcl{ \SP[b_1][\mathcal{A}_1] }
		\textrm{ in } \cohom{\mathcal{S}}[3-n].
	\end{equation*}
	\end{thm}
	\begin{proof}
	The proof is the same as the proof of \cref{thm:invariance-extended-infty-modulus},
	except we use the naturality up to exact terms established in \cref{lm:functoriality-superpotential},
	which is enough for our purposes.

	Let $b \in \tc{A}$ with $\mcfunc{\evalmf}^i \left( b \right) = b_i$ for $i = 0, 1$
	and let $\phi \colon \totcompe{\mathfrak{A} / \mathfrak{R}}[2][] \rightarrow \mathfrak{R}[4-n]$
	be such that the diagram of \cref{fig:pseudoisotopy-total-inner-product-algebras} commutes for $i = 0, 1$.
	Choose a homotopy $h \colon R^{*} \rightharpoonup S^{*-1}$ between
$\evalm^0 = \base{\evalmf}^0$ and $\evalm^1 = \base{\evalmf}^1$ with
	\begin{equation*}
		d_S \circ h + h \circ d_R = \evalm^1 - \evalm^0.
	\end{equation*}
	By \cref{lm:functoriality-superpotential}, we have
	\begin{equation*}
		\SP[b_i][\mathcal{A}_i]
		\stackrel{\eqref{eq:gauge-equivalence-bounding-chains}}{=}
		\SP[ \mcfunc{\evalmf}^i \left( b \right) ][\mathcal{A}_i]
		\stackrel{\eqref{eq:sp-pseudo-invariance}}{=}
		\evalm^i \left( \SP[b][\mathfrak{A}] \right) - d_S \left( x_i \right)
	\end{equation*}
	for some elements $x_i \in S^{2 - n}$ and $i = 0, 1$. Then
	\begin{equation*}
		\begin{aligned}
			\SP[b_1][\mathcal{A}_1] - \SP[b_0][\mathcal{A}_0]
			={} &
			\left(
			\evalm^1 \left( \SP[b][\mathfrak{A}] \right) -
			\evalm^0 \left( \SP[b][\mathfrak{A}] \right)
			\right) + d_S \left( x_0 - x_1 \right)
			\\
			={} &
			\left( d_S \circ h + h \circ d_R \right) \left( \SP[b][\mathfrak{A}] \right) +
			d_S \left( x_0 - x_1 \right)
			\\
			={} &
			d_S \left( x_0 - x_1 + h \left( \SP[b][\mathfrak{A}] \right) \right),
		\end{aligned}
	\end{equation*}
	where we used the fact that $d_R \left( \SP[b][\mathfrak{A}] \right) = 0$ by
	\cref{lm:superpotential-bounding-chain-closed}.
	\end{proof}

	\subsection{Formal Derivative of the Superpotential}
	\label{sec:formal-derivative-sp}
	Let $S$ be a graded-commutative Banach $\mathbbm{k}$-algebra and let $\mathcal{B} = \left( B, \nu, \phi \right)$
	be an $n$-dimensional total inner product Banach $\Ainf$-algebra over $S$.
	Consider the superpotential $\SP \colon \tc{B} \rightarrow S^{3-n}$ given by
$\SP[b] = \phi \left( \G{b}[\geq 2] \right)$, and assume
	that $b = b(t) \in \tc{B}$ depends on some formal parameter $t$, and that
$\phi$ and $\nu$ do not depend on $t$. In this section we show that the formal
	derivative $\partial_t \left( \SP[b] \right)$ of the superpotential is given by
	\begin{equation*}
		\partial_t \left( \SP[b] \right) =
		\phi \left(
		\sum_{k=2}^{\infty}
		\left( \ul{ \partial_t \left( b \right) } \otimes \Exp{b} \right) \otimes
		\left( \ul{ \corest{\nu} \left( \Exp{b} \right) } \otimes \Exp{b} \right)^{\otimes \left( k - 1 \right)}
		\right).
	\end{equation*}

	While the calculation of the formal derivative can be done directly (see \cref{rem:partial-t-Gb-explicit}),
	in order to reduce calculations, we adopt the strategy of translating the calculation of the derivative
	into the problem of calculating the differential of $\SP$ in an appropriate scalar extension and
	using previous results.

	\subsubsection{Scalar Extension for Total Inner Products} \label{sec:scalar-extension-total-inner-products}
	Before going into the details of the calculation, we need to discuss the notion of scalar extension
	for total inner products. Let $\mathcal{R} = \left( R, d_R \right)$
	and $\mathcal{S} = \left( S, d_S \right)$ be differential graded-commutative Banach $\mathbbm{k}$-algebras
	and let $i \colon \mathcal{R} \rightarrow \mathcal{S}$ be a morphism of differential
	graded Banach $\mathbbm{k}$-algebras.
	Given a Banach $\Ainf$-algebra $\mathcal{A} = \left( A, \mu \right)$ over $\mathcal{R}$, let
$i_{!} \left( \mathcal{A} \right) = \left( i_{!} \left( A \right), i_{!} \left( \mu \right) \right)$
	be the scalar extension of $\mathcal{A}$ along $i$, which is a Banach $\Ainf$-algebra over $\mathcal{S}$
	(see \cref{dfn:scalar-rest-ext-a-infinity}).
	Given a graded $R$-linear map $\phi \colon \totcompe{A/R}[2][] \rightharpoonup R$, there exists
	a unique graded $S$-linear map $i_{!} \left( \phi \right) \colon
\totcompe{i_{!} \left( A \right) / S}[2][] \rightharpoonup S$ with
$\degb{i_{!} \left( \phi \right)} = \degb{\phi}$, called the \textbf{scalar extension} of $\phi$, such that
$i_{!} \left( \phi \right) \circ \cindmape{F} \left( \resunder{i} \right) = i \circ \phi$,
	where $\resunder{i} \colon \mathcal{A} \rightarrow i_{!} \left( \mathcal{A} \right)$ is the canonical
	map over $i$ and $\cindmape{F} \left( \resunder{i} \right)$ is the induced map
	by $\resunder{i}$ between the extended total complexes. More precisely, we have a concrete isomorphism
	\begin{align}
		\totcompe{i_{!} \left( A \right) / S}[2][] & =
		S[1] \oplus \totc{\ncdf{i_{!} \left( A \right) / S}[][]}[][\geq 2][] \cong
		\left( S \otimes_R R[1] \right) \oplus \totc{S \otimes_R \ncdf{A / R}[][]}[][\geq 2][]
		\notag
		\\
		                                           & \cong
		\left( S \otimes_R R[1] \right) \oplus \left( S \otimes_R \totc{\ncdf{A / R}[][]}[][\geq 2][] \right)
		\label{eq:totcompe-commutes-scalar-extension}
		\\
		                                           & \cong
		S \otimes_R \left( R[1] \oplus \totc{\ncdf{A / R}[][]}[][\geq 2][] \right) =
		S \otimes_R \totcompe{A / R}[2][],
		\notag
	\end{align}
	which is a \textit{chain map}, and under this isomorphism, the scalar extension $i_{!} \left( \phi \right)$ corresponds to
	standard scalar extension $\id_{S} \otimes_R \phi$. See \cref{fig:scalar-extension-phi}.

	Since the horizontal maps of \cref{fig:scalar-extension-phi} are chain
	maps, uniqueness implies that we have
$\partial \left( {i_{!} \left( \phi \right)} \right) = i_{!} \left( \partial \phi \right)$,
	where $\partial$ denotes the differentials on the inner hom complexes.
	In particular, if $\phi$ corresponds to an $n$-dimensional total inner product on $\mathcal{A}$,
	then $i_{!} \left( \phi \right)$ corresponds to an $n$-dimensional total inner product on
$i_{!} \left( \mathcal{A} \right)$,\footnote{By corresponding, we mean that we identify an $n$-dimensional total inner product, i.e.,
	a contractive degree zero $R$-linear chain map $\phi \colon \totcompe{A/R}[2][] \rightarrow R[4-n]$, with
	the contractive $R$-linear chain map $\totcompe{A/R}[2][] \rightharpoonup R$ of degree $4 - n$
	obtained by composing $\phi$ with the canonical chain map $R[4-n] \rightharpoonup R$, and similarly for $i_{!} \left( \phi \right)$.}
	i.e., we have
	\begin{equation*}
		d_R \circ \phi = (-1)^{4-n} \phi \circ D_{\mu} \implies
		d_S \circ i_{!} \left( \phi \right) =
		(-1)^{4-n} i_{!} \left( \phi \right) \circ D_{i_{!} \left( \mu \right)}.
	\end{equation*}

	When $\mathcal{A} = \left( A, \mu, \phi \right)$ is a total inner
	product Banach $\Ainf$-algebra, the total inner product Banach
$\Ainf$-algebra
$i_{!} \left( \mathcal{A} \right) = \left( i_{!} \left( A \right),
i_{!} \left( \mu \right), i_{!} \left( \phi \right) \right)$
	is called the \textbf{scalar extension of} $\mathcal{A}$ \textbf{along} $i$,
	and the canonical $\Ainf$-morphism
$\resunder{i} \colon \mathcal{A} \rightarrow i_{!} \left( \mathcal{A} \right)$
	over $i$ associated to the scalar extension becomes
	a morphism of total inner product Banach $\Ainf$-algebras.

	The process of scalar extension for total inner products respects compositions, i.e.,
	if we have a third differential graded-commutative Banach $\mathbbm{k}$-algebra
$\mathcal{T}$ and a morphism $j \colon \mathcal{S} \rightarrow \mathcal{T}$,
	the natural strict isomorphism
$j_{!} \left( i_{!} \left( \mathcal{A} \right) \right) \cong \left( j \circ i \right)_{!} \left( \mathcal{A} \right)$
	of Banach $\Ainf$-algebras over $\mathcal{T}$ is also a morphism of total inner
	product Banach $\Ainf$-algebras.

	\begin{figure}
	\adjustbox{scale=0.85,center}{
		\begin{tikzcd}
			{\left( \totcompe{A/R}[2][], D_{\mu} \right)} &&
			{\left( \totcompe{i_{!} \left( A \right) / S}[2][], D_{i_{!} \left( \mu \right)} \right)}
			&&
			{\left( S \otimes_R \totcompe{A/R}[2][], d_S \otimes_R \id + \id \otimes_R D_{\mu} \right)}
			\\
			&&  &&
			{\left( S \otimes_R R, d_S \otimes_R \id + \id \otimes_R d_R \right)}
			\\
			\left( R, d_R \right) && \left( S, d_S \right) && \left( S, d_S \right)
			\arrow["{\cindmape{F} \left( \resunder{i} \right)}", from=1-1, to=1-3]
			\arrow["{\cong}", from=1-3, to=1-5]
			\arrow["{x \mapsto 1_S \otimes_R x}", bend left=10, from=1-1, to=1-5,
				start anchor=north, end anchor=north]
			\arrow["\phi"', harpoon, from=1-1, to=3-1]
			\arrow["{!\exists i_{!} \left( \phi \right)}", dotted, harpoon, from=1-3, to=3-3]
			\arrow["{\cong}", from=2-5, to=3-5]
			\arrow["{\id \otimes_R \phi}", harpoon, from=1-5, to=2-5]
			\arrow["i"', from=3-1, to=3-3]
			\arrow["\id"', from=3-3, to=3-5]
		\end{tikzcd}
	}
	\caption{Scalar extension for maps $\totcompe{A/R}[2][] \rightharpoonup R$.}
	\label{fig:scalar-extension-phi}
	\end{figure}

	\subsubsection{Derivative Calculation} \label{sec:derivative-calc}
	Let $R$ be a graded-commutative Banach $\mathbbm{k}$-algebra and let
$\mathcal{A} = \left( A, \mu, \phi_A \right)$ be a Banach $\Ainf$-algebra over $R$ endowed
	with an $n$-dimensional total inner product $\phi_A$.
	Let $t$ be a formal variable of even degree $\degb{t}$ with $\nnorm[t] = \alpha < 1$
	and consider the graded-commutative Banach $\mathbbm{k}$-algebra $S \defeq \pows{R}[t]$. A
	degree $m$ element of $S$ is a formal power series $p = p(t) = \sum_{k=0}^{\infty} a_k \cdot t^k$
	with $a_k \in R^{m - k \cdot \degb{t}}$ such that $\nnorm[a_k] \cdot \alpha^k \to 0$ and
	we have $\nnorm[p] = \max_{k} \Set{\nnorm[a_k] \cdot \alpha^k}$. The \textbf{formal derivative}
	of a power series $p(t) = \sum_{n=0}^{\infty} a_n \cdot t^n \in \pows{R}[t]^{m}$
	is given by
	\begin{equation*}
		\partial_t \left( p \right) = \partial_t \left( p(t) \right) =
		\sum_{n=0}^{\infty} \left( \left( n + 1 \right) \cdot a_{n+1} \right) \cdot t^{n} \in
		\pows{R}[t]^{m - \degb{t}}
	\end{equation*}
	so that $\partial_t \colon S \rightharpoonup S$ is an $R$-linear derivation of degree
$-\degb{t}$.

	Let $i \colon R \rightarrow S$ be the inclusion and consider the scalar extension
$i_{!} \left( \mathcal{A} \right)$ along $i$. Set $B = i_{!} \left( A \right), \nu = i_{!} \left( \mu \right)$
	and $\phi_B = i_{!} \left( \phi_A \right)$ so that
$\mathcal{B} = \left( B, \nu, \phi_B \right) = i_{!} \left( \mathcal{A} \right)$
	is an $n$-dimensional total inner product $\Ainf$-algebra over $S$.

	\begin{thm} \label{thm:derivative-sp}
	Let $b \in \tc{B}$. Then we have\footnote{We use standard abuse of notation to denote
	by $\partial_t$ both the derivative operator $\partial_t \colon S \rightharpoonup S$ and
	the induced module derivative $\partial_t \otimes_R \id \colon B \rightharpoonup B$
	over $\partial_t$ on $B = S \otimes_R A$.}
	\begin{equation} \label{eq:derivative-sp-braid-2}
		\partial_t \left( \SP[b][\mathcal{B}] \right) =
		\partial_t \left( \phi_B \left( \G{b}[\geq 2][\nu] \right) \right) =
		\phi_B \left(
		\sum_{k=2}^{\infty}
		\left( \ul{ \partial_t \left( b \right) } \otimes \Exp{b} \right) \otimes
		\left( \ul{ \corest{\nu} \left( \Exp{b} \right) } \otimes \Exp{b} \right)^{\otimes \left( k - 1 \right)}
		\right).
	\end{equation}
	\end{thm}
	\begin{proof}
	Let $dt$ be a formal variable of degree $\degb{t} + 1$ with $\nnorm[dt] = \nnorm[t]$
	and consider the differential graded-commutative Banach $\mathbbm{k}$-algebra
	\begin{equation*}
		W \defeq \underbrace{\pows{R}[t]}_{W_0} \oplus \underbrace{\pows{R}[t] \left< dt \right>}_{W_1}
	\end{equation*}
	of K\"ahler differentials on the formal line. An element of $W$ of (total) degree $m$ has the form
$p(t) + q(t) \, dt$ with $p(t) \in \pows{R}[t]^{m}$ and
$q(t) \in \pows{R}[t]^{m - \degb{t} - 1}$. In terms of the formal derivative,
	the contractive degree one differential
$d \colon W \rightharpoonup W$ on $W$ is given by
	\begin{equation*}
		\begin{aligned}
			d \left( p(t) + q(t) \, dt \right) = dt \cdot \partial_t \left( p \left( t \right) \right) & =
			(-1)^{\degb{dt} \cdot \degb{\partial_t \left( p \left( t \right) \right)}}
			\partial_t \left( p \left( t \right) \right) \cdot dt
			\\
			                                                                                           & =
			(-1)^{\left( \degb{t} + 1 \right) \cdot \left( \degb{p} - \degb{t} \right)}
			\partial_t \left( p \left( t \right) \right) \, dt.
		\end{aligned}
	\end{equation*}
	Let $b = b(t) \in \tc{B}$ be a topologically nilpotent element of $B = \pows{R}[t] \otimes_R A$
	which depends on the formal variable $t$. Instead of calculating the formal derivative
$\partial_t \left( \SP[b][\mathcal{B}] \right)$,
	we can take the differential of $\SP[b][\mathcal{B}]$ and extract the coefficient of $dt$.

	Consider the commutative diagram of differential graded Banach $\mathbbm{k}$-algebras
	and morphism given in the upper part of \cref{fig:base-algebras-morphisms-der-calc}.
	The morphisms $i,j,k,l$ are
	the natural inclusions, while $p \colon \left( W, d \right) \rightarrow \left( W_0, 0 \right)$
	is the projection onto $W_0$ given by $f(t) + g(t) \, dt \mapsto f(t)$. Note that
$j = l$ as morphisms of the underlying graded Banach $\mathbbm{k}$-algebras, and that
$p \circ k = \id$, but the identity map $\left( W, 0 \right) \rightarrow \left( W, d \right)$
	does not respect the differentials and hence we do not put the corresponding arrow in
	\cref{fig:base-algebras-morphisms-der-calc}.

	The upper part of
	\cref{fig:base-algebras-morphisms-der-calc} induces a commutative diagram
	between the corresponding scalar extensions given in the lower part of
	\cref{fig:base-algebras-morphisms-der-calc}. The objects are total inner product Banach $\Ainf$-algebras
	and the morphisms are \textit{strict} and respect all the structures.
	The morphisms $\resunder{i},\resunder{j},\resunder{l}$
	are the canonical morphisms over $i,j,l$ respectively, while $\resunder{p}$ (resp.\ $\resunder{k}$) is given
	by $\resunder{p}_{1} = p \otimes_R \id_A$ (resp.\ $\resunder{k}_{1} = k \otimes_R \id_A$).\footnote{Even
	though $\resunder{p}$ is not the canonical morphism $j_{!} \left( A \right) \rightarrow
p_{!} \left( j_{!} \left( A \right) \right)$, it can be identified with the canonical
	morphism using the natural isomorphism $p_{!} \left( j_{!} \left( A \right) \right) \cong
\left( p \circ j \right)_{!} \left( A \right) = i_{!} \left( A \right)$ so
	we continue to use the same notation. The same applies to $\resunder{k}$.}
	Note that $l_{!} \left( A \right) = j_{!} \left( A \right)$ and that $\resunder{j} = \resunder{l}$
	and $l_{!} \left( \phi \right) = j_{!} \left( \phi \right)$, so the underlying graded $W$-modules
	and the total inner products are the same, but $l_{!} \left( \mu \right) \neq j_{!} \left( \mu \right)$
	so the $\Ainf$-structures are different (see \cref{rem:scalar-extension-depends-differentials}).
	The coderivation $l_{!} \left( \mu \right)$ is $W$-linear, while $j_{!} \left( \mu \right)$ is a
	coderivation over $d$, so the identity map does not respect the coderivations and hence we do not
	put the corresponding arrow in \cref{fig:base-algebras-morphisms-der-calc}.

	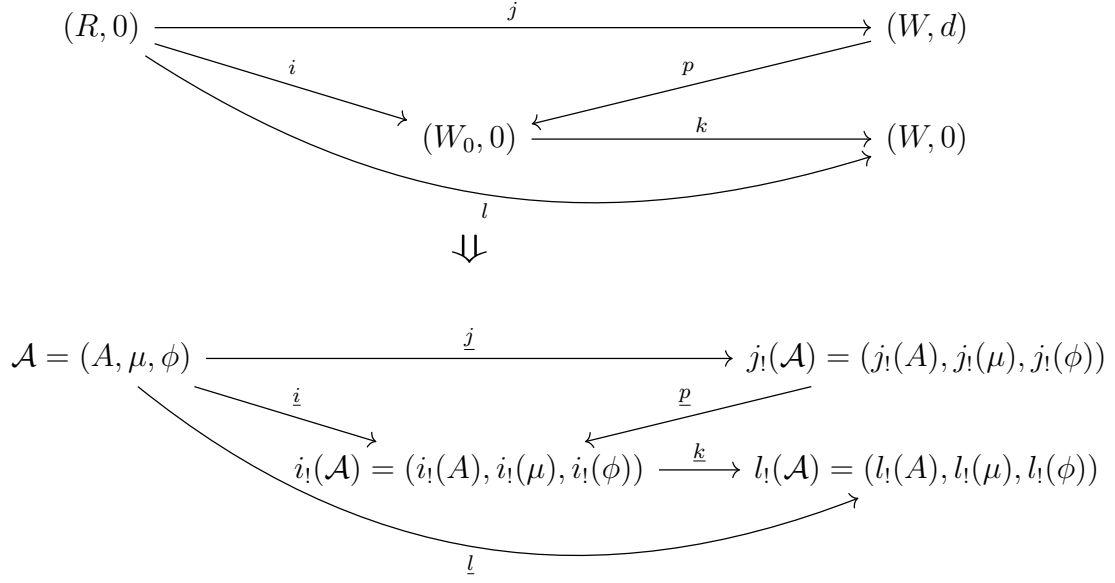
\begin{figure}
	\begin{tikzcd}
	{\left( R,0 \right)} && {\left( W, d \right)} \\
	& {\left( W_0, 0 \right)} & {\left( W ,0 \right)} \\
	& {\mathrel{\scalebox{2}[1]{$\Downarrow$}}}
	\\
	{\mathcal{A} = \left( A, \mu, \phi \right)} &&
	{j_{!} \left( \mathcal{A} \right) = \left( j_{!} \left( A \right), j_{!} \left( \mu \right),
		j_{!} \left( \phi \right) \right)}
	\\
	&
	{i_{!} \left( \mathcal{A} \right) = \left( i_{!} \left( A \right), i_{!} \left( \mu \right),
		i_{!} \left( \phi \right) \right)}
	&
	{l_{!} \left( \mathcal{A} \right) = \left(
		l_{!} \left( A \right), l_{!} \left( \mu \right), l_{!} \left( \phi \right)
		\right)}
	\arrow["{\resunder{j}}", from=4-1, to=4-3]
	\arrow["{\resunder{i}}", from=4-1, to=5-2]
	\arrow["{\resunder{p}}"', from=4-3, to=5-2]
	\arrow["{\resunder{k}}", from=5-2, to=5-3]
	\arrow["{\resunder{l}}"', from=4-1, to=5-3, bend left=-30]
	\arrow["j", from=1-1, to=1-3]
	\arrow["i", from=1-1, to=2-2]
	\arrow["p"', from=1-3, to=2-2]
	\arrow["k", from=2-2, to=2-3]
	\arrow["l"', from=1-1, to=2-3, bend left=-25]
	\end{tikzcd}
	\caption{Morphisms between various scalar extensions.}
	\label{fig:base-algebras-morphisms-der-calc}
	\end{figure}

	Let us set $\mathfrak{b} \defeq \mcfunc{\resunder{k}} \left( b \right) = \resunder{k}_{1} \left( b \right)
\in \tc{l_{!} \left( A \right)} = \tc{j_{!} \left( A \right)}$,
	so that
	\begin{equation*}
		\mcfunc{\resunder{p}} \left( \mathfrak{b} \right) = \resunder{p}_{1} \left( \mathfrak{b} \right) =
		\left( \resunder{p}_1 \circ \resunder{k}_1 \right) \left( b \right) =
		\left( \left( p \otimes_R \id \right) \circ \left( k \otimes_R \id \right) \left( b \right) \right) =
		\left( \left( p \circ k \right) \otimes_R \id \right) \left( b \right) = b.
	\end{equation*}
	Denote by $\cont{\partial_t} \colon W \rightharpoonup W$ the degree $-\degb{t} - 1$ operator extracting
	the coefficient of $dt$, i.e.,
$\cont{\partial_t} \left( p \left( t \right) + dt \, q \left( t \right) \right) = q \left( t \right)$.
	Then $\partial_t =  p \circ \cont{\partial_t} \circ d \circ k$ and
$d \circ k \circ p = d$, and hence we have
	\begin{equation*}
		\begin{aligned}
			\partial_t \left( \SP[b][\mathcal{B}] \right) & =
			\left( p \circ \cont{\partial_t} \circ d \circ k \right) \left(
			\SP[\mcfunc{\resunder{p}} \left( \mathfrak{b} \right)][\mathcal{B}]
			\right)
			\\
			                                              & =
			\left( p \circ \cont{\partial_t} \circ \left( d \circ k \circ p \right) \right)
			\left(
			\SP[\mathfrak{b}][j_{!} \left( \mathcal{A} \right)]
			\right)
			\\
			                                              & =
			\left( p \circ \cont{\partial_t} \circ d \right)
			\left(
			\SP[\mathfrak{b}][j_{!} \left( \mathcal{A} \right)]
			\right),
		\end{aligned}
	\end{equation*}
	so it is enough to calculate $d \left( \SP[\mathfrak{b}][j_{!} \left( \mathcal{A} \right)] \right)$.
	Unwinding the definitions and using previous results, we have:
	\begin{equation*}
		\begin{aligned}
			d \left( \SP[\mathfrak{b}][j_{!} \left( \mathcal{A} \right)] \right)
			\stackrel{\phantom{\eqref{eq:D-G-geq-2}}}{=} &
			\left( d \circ j_{!} \left( \phi \right) \right) \left(
			\G{\mathfrak{b}}[\geq 2][j_{!} \left( \mu \right)]
			\right)
			\\
			\stackrel{\phantom{\eqref{eq:D-G-geq-2}}}{=} &
			(-1)^{\degb{\phi}}
			\left( j_{!} \left( \phi \right) \circ D_{j_{!} \left( \mu \right)} \right) \left(
			\G{\mathfrak{b}}[\geq 2][j_{!} \left( \mu \right)]
			\right)
			\\
			\stackrel{\eqref{eq:D-G-geq-2}}{=}           &
			(-1)^{\degb{\phi}} j_{!} \left( \phi \right) \left(
			\Ho{\mathfrak{b}}[\geq 2][j_{!} \left( \mu \right)]
			\right),
		\end{aligned}
	\end{equation*}
	where
	\begin{equation*}
		\Ho{\mathfrak{b}}[\geq 2][j_{!} \left( \mu \right)]
		\stackrel[\eqref{eq:Ho-geq-2-def}]{\eqref{eq:H_n(b)-formula}}{=}
		- \sum_{r = 2}^{\infty} \frac{1}{r}
		\left(
		\ul{ \corest{j_{!} \left( \mu \right)} \left( \Exp{\mathfrak{b}} \right)} \otimes \Exp{\mathfrak{b}}
		\right)^{\otimes r}.
	\end{equation*}
	The $\Ainf$-structure $j_{!} \left( \mu \right)$ is a sum of two $\Ainf$-structures
$j_{!} \left( \mu \right) = l_{!} \left( \mu \right) + \hat{d}$, where
	the $\Ainf$-structure $l_{!} \left( \mu \right)$ is $W$-linear, and $\hat{d}$ is an $\Ainf$-structure over
$d$ whose only non-zero component $\hat{d}_1 \colon W \otimes_R A \rightharpoonup W \otimes_R A$
	is $\hat{d}_1 = d \otimes_R \id$. Since we have
	\begin{equation*}
		\corest{\hat{d}} \left( \Exp{\mathfrak{b}} \right) = \hat{d}_1 \left( \mathfrak{b} \right)
		= dt \cdot \resunder{k}_1 \left( \partial_t \left( b \right) \right),
	\end{equation*}
	and $\left( dt \right)^2 = 0$, we see that we can decompose $\Ho{\mathfrak{b}}[\geq 2][j_{!} \left( \mu \right)]$
	as a sum
	\begin{equation*}
		\begin{aligned}
			\Ho{\mathfrak{b}}[\geq 2][j_{!} \left( \mu \right)] ={} &
			- \sum_{r = 2}^{\infty} \frac{1}{r}
			\left(
			\ul{ \corest{\left( l_{!} \left( \mu \right) + \hat{d} \right)} \left( \Exp{\mathfrak{b}} \right)}
			\otimes
			\Exp{\mathfrak{b}}
			\right)^{\otimes r}
			\\
			={}                                                     &
			- \sum_{r = 2}^{\infty} \frac{1}{r}
			\left(
			\ul{ \corest{l_{!} \left( \mu \right)} \left( \Exp{\mathfrak{b}} \right)}
			\otimes
			\Exp{\mathfrak{b}}
			\right)^{\otimes r}
			-
			\sum_{r = 2}^{\infty}
			\left(
			\ul{ \hat{d}_1 \left( \mathfrak{b} \right) } \otimes \Exp{\mathfrak{b}}
			\right) \otimes
			\left(
			\ul{ \corest{l_{!} \left( \mu \right)} \left( \Exp{\mathfrak{b}} \right)} \otimes
			\Exp{\mathfrak{b}}
			\right)^{\otimes \left( r - 1 \right)}
			\\
			={}                                                     &
			\Ho{\mathfrak{b}}[\geq 2][l_{!} \left( \mu \right)]
			+
			(-1)^{\degb{dt} + 1}
			dt \, \left(
			\sum_{r = 2}^{\infty}
			\left(
			\ul{ \resunder{k}_1 \left( \partial_t \left( b \right) \right) } \otimes \Exp{\mathfrak{b}}
			\right) \otimes
			\left(
			\ul{ \corest{l_{!} \left( \mu \right)} \left( \Exp{\mathfrak{b}} \right)} \otimes
			\Exp{\mathfrak{b}}
			\right)^{\otimes \left( r - 1 \right)}
			\right)
			\\
			={}                                                     &
			\Ho{\mathfrak{b}}[\geq 2][l_{!} \left( \mu \right)]
			+
			dt \cdot \cindmap{F} \left( \resunder{k} \right) \left(
			\sum_{r = 2}^{\infty}
			\left( \ul{\partial_t \left( b \right)} \otimes \Exp{b} \right) \otimes
			\left(
			\ul{ \corest{i_{!} \left( \mu \right)} \left( \Exp{b} \right)} \otimes \Exp{b}
			\right)^{\otimes \left( r - 1 \right)}
			\right).
		\end{aligned}
	\end{equation*}
	Then
	\begin{equation*}
		\begin{aligned}
			\MoveEqLeft
			d \left( \SP[\mathfrak{b}][j_{!} \left( \mathcal{A} \right)] \right) ={}
			                                                                      (-1)^{\degb{\phi}} j_{!} \left( \phi \right) \left(
			\Ho{\mathfrak{b}}[\geq 2][j_{!} \left( \mu \right)]
			\right)
			\\
			={} &
			(-1)^{\degb{\phi}} j_{!} \left( \phi \right) \left(
			\Ho{\mathfrak{b}}[\geq 2][l_{!} \left( \mu \right)]
			\right)
			\\
			    & +
			      (-1)^{\degb{\phi}} j_{!} \left( \phi \right) \left(
			dt \cdot \cindmap{F} \left( \resunder{k} \right)
			\left(
			\sum_{r = 2}^{\infty}
			\left( \ul{\partial_t \left( b \right)} \otimes \Exp{b} \right) \otimes
			\left(
			\ul{ \corest{i_{!} \left( \mu \right)} \left( \Exp{b} \right)} \otimes \Exp{b}
			\right)^{\otimes \left( r - 1 \right)}
			\right)
			\right)
			\\
			={} &
			(-1)^{\degb{\phi}}
			\cancel{ {\left( l_{!} \left( \phi \right) \circ D_{l_{!} \left( \mu \right)} \right)} }
			\left(
			\G{\mathfrak{b}}[\geq 2][l_{!} \left( \mu \right)]
			\right)
			\\
			    & +
			      (-1)^{\degb{\phi} \cdot \left( 1 + \degb{dt} \right)}
			dt \,
			\left(
			j_{!} \left( \phi \right) \circ \cindmap{F} \left( \resunder{k} \right)
			\right)
			\left(
			\sum_{r = 2}^{\infty}
			\left( \ul{\partial_t \left( b \right)} \otimes \Exp{b} \right) \otimes
			\left(
			\ul{ \corest{i_{!} \left( \mu \right)} \left( \Exp{b} \right)} \otimes \Exp{b}
			\right)^{\otimes \left( r - 1 \right)}
			\right)
			\\
			={} &
			dt \, \left( k \circ i_{!} \left( \phi \right) \right)
			\left(
			\sum_{r = 2}^{\infty}
			\left( \ul{\partial_t \left( b \right)} \otimes \Exp{b} \right) \otimes
			\left(
			\ul{ \corest{i_{!} \left( \mu \right)} \left( \Exp{b} \right)} \otimes \Exp{b}
			\right)^{\otimes \left( r - 1 \right)}
			\right).
		\end{aligned}
	\end{equation*}
	Finally
	\begin{equation*}
		\begin{aligned}
			\partial_t \left( \SP[b][\mathcal{B}] \right) & =
			\left( p \circ \cont{\partial_t} \circ d \right)
			\left(
			\SP[\mathfrak{b}][j_{!} \left( \mathcal{A} \right)]
			\right)
			\\
			                                              & =
			\left( p \circ k \circ i_{!} \left( \phi \right) \right)
			\left(
			\sum_{r = 2}^{\infty}
			\left( \ul{\partial_t \left( b \right)} \otimes \Exp{b} \right) \otimes
			\left(
			\ul{ \corest{i_{!} \left( \mu \right)} \left( \Exp{b} \right)} \otimes \Exp{b}
			\right)^{\otimes \left( r - 1 \right)}
			\right)
			\\
			                                              & =
			\phi_B \left(
			\sum_{r=2}^{\infty}
			\left( \ul{ \partial_t \left( b \right) } \otimes \Exp{b} \right) \otimes
			\left( \ul{ \corest{\nu} \left( \Exp{b} \right) } \otimes \Exp{b} \right)^{\otimes \left( r - 1 \right)}
			\right)
		\end{aligned}
	\end{equation*}
	as required.
	\end{proof}

	\begin{rem} \label{rem:lie-partial-t-G-geq-2}
	Let $\eta$ be a coderivation on $\tens{B}[S]$ such that $\left[ \eta, \nu \right] = 0$.
	Then one can show that
	\begin{equation*}
		\begin{aligned}
			\clie{\eta} \left( \Go{b}[\geq 2][\nu] \right) ={} &
			(-1)^{\degb{\eta}}
			\sum_{k=2}^{\infty}
			\left( \ul{ \corest{\eta} \left( \Exp{b} \right) } \otimes \Exp{b} \right) \otimes
			\left( \ul{ \corest{\nu} \left( \Exp{b} \right) } \otimes \Exp{b} \right)^{\otimes (k - 1)}
			\\
			                                                   & +
			                                                     (-1)^{\degb{\eta} - 1}
			\sum_{k=2}^{\infty} \ulz{\eta} \otimes \ulz{\nu}^{\otimes (k - 1)}
			+
			D_{\nu} \left( z \right)
		\end{aligned}
	\end{equation*}
	for some $z \in \totcomp{B}[2][\degb{\eta} - 2]$. In particular, if we take
$\eta = \partial_t$, we see that
	\begin{equation} \label{eq:clie-partial-t-Go-b}
		\clie{\partial_t} \left( \Go{b}[\geq 2][\nu] \right) =
		\underbrace{
			\sum_{k=2}^{\infty}
			\left( \ul{ \partial_t \left( b \right) } \otimes \Exp{b} \right) \otimes
			\left( \ul{ \corest{\nu} \left( \Exp{b} \right) } \otimes \Exp{b} \right)^{\otimes \left( k - 1 \right)}
		}_{\star}
		+
		D_{\nu} \left( z \right)
	\end{equation}
	for some $z \in \totcomp{B}[2][-\degb{t} - 2]$. Thus, it is not true
	that the formal derivative of $\Go{b}[\geq 2][\nu]$ is given by $\star$, this holds
	\textit{only up to an exact element} of the total complex.
	This gives a different proof of \cref{thm:derivative-sp}, as
	\begin{equation*}
		\begin{aligned}
			\partial_t \left( \SP[b][\mathcal{B}] \right) & =
			\partial_t \left( \phi_B \left( \Go{b}[\geq 2][\nu] \right) \right) =
			\phi_B \left( \clie{\partial_t} \left( \Go{b}[\geq 2][\nu] \right) \right) =
			\phi_B \left( \star + D_{\nu} \left( z \right) \right)
			\\
			                                              & =
			\phi_B \left(
			\sum_{k=2}^{\infty}
			\left( \ul{ \partial_t \left( b \right) } \otimes \Exp{b} \right) \otimes
			\left( \ul{ \corest{\nu} \left( \Exp{b} \right) } \otimes \Exp{b} \right)^{\otimes \left( k - 1 \right)}
			\right).
		\end{aligned}
	\end{equation*}
	\end{rem}

	\begin{rem} \label{rem:partial-t-Gb-explicit}
	We briefly describe how to compute the formal derivative of $\G{b}[2]$ directly.
	To simplify notation, we denote the term
	\begin{equation*}
		\G{b}[2] = \sum_{i,j,k=0}^{\infty} \frac{1}{i + j + k + 1} \,
		\ul{b} \otimes b^{\otimes i} \otimes
		\ul{ \nu_j \left( b^{\otimes j} \right)} \otimes b^{\otimes k}
	\end{equation*}
	by $( \ul{b} \pstar \ul{\nu (\pstar)} \pstar )_{w}$, where $\pstar$ denotes
	arbitrary non-negative powers of $b$, the brackets $\left( \quad \right)$ indicate we perform
	a summation over all powers of $b$ appearing in the terms $\pstar$,
	and the subscript $w$ indicates that we divide each term of the summation
	by the sum of all the powers appearing in that term.
	Using our notation, we have
	\begin{equation*}
		\begin{aligned}
			\partial_t \left( \G{b}[2] \right) ={} &
			\left( \ul{\dot{b}}  \pstar  \ul{\nu \left( \pstar \right)} \pstar \right)_{w}
			+
			\left( \ul{b} \pstar \dot{b} \pstar \ul{\nu \left( \pstar \right)} \pstar \right)_{w}
			+
			\left( \ul{b} \pstar \ul{\nu \left( \pstar \dot{b} \pstar \right)} \pstar \right)_{w}
			\\
			                                       & +
			\left( \ul{b} \pstar \ul{\nu \left( \pstar \right)} \pstar \dot{b} \pstar \right)_{w},
		\end{aligned}
	\end{equation*}
	where $\dot{b} = \partial_t \left( b \right)$.
	Working with the formulas for $\braidop_2$, and using the fact that
$b, \pstar, \dot{b}$ are all of even degree, we have
	\begin{align*}
		\qdr \left( \underbrace{ \left(
			            \ul{b} \pstar \ul{\dot{b}} \pstar \ul{\nu \left( \pstar \right)} \pstar
			            \right)_{w}
		            }_{A} \right) \stackrel{\eqref{eq:qdr-formula}}{=} &
		\left( \ul{\dot{b}} \pstar \ul{\nu \left( \pstar \right)} \hpstar \right)_{w}
		-
		\left( \ul{b} \pstar \dot{b} \pstar \ul{\nu \left( \pstar \right)} \pstar \right)_{w}
		+
		\left( \ul{b} \pstar \ul{\dot{b}} \pstar \nu \left( \pstar \right) \pstar \right)_{w},
		\\
		\qdr \left( \underbrace{ \left(
			            \ul{b} \pstar \ul{\nu \left( \pstar \right)} \pstar \ul{\dot{b}} \pstar
			            \right)_{w}
		            }_{B} \right) \stackrel{\eqref{eq:qdr-formula}}{=} &
		\left( \ul{\nu \left( \pstar \right)} \pstar \ul{\dot{b}} \hpstar \right)_{w}
		-
		\left( \ul{b} \pstar \nu \left( \pstar \right) \pstar \ul{\dot{b}} \pstar \right)_{w}
		-
		\left( \ul{b} \pstar \ul{\nu \left( \pstar \right)} \pstar \dot{b} \pstar \right)_{w},
	\end{align*}
	and
	\begin{equation*}
		\begin{aligned}
			\clie{\nu} \left( \underbrace{ \left(
			           \ul{b} \pstar \ul{\dot{b}} \pstar
			           \right)_{w}
			}_{C} \right)
			\stackrel{\eqref{eq:cyc-lie-formula}}{=} &
			-
			\left( \ul{ \nu \left( \hpstar \right) } \pstar \ul{\dot{b}} \pstar \right)_{w}
			-
			\left( \ul{b} \pstar \nu \left( \pstar \right) \pstar \ul{\dot{b}} \pstar \right)_{w}
			+
			\left( \ul{b} \pstar \ul{ \nu \left( \pstar \dot{b} \pstar \right) } \pstar \right)_{w}
			\\
			                                         &
			+
			\left( \ul{b} \pstar \ul{\dot{b}} \pstar \nu \left( \pstar \right) \pstar \right)_{w},
		\end{aligned}
	\end{equation*}
	where $\hpstar$ denotes terms of the form $b^{\otimes i} \otimes b \otimes b^{\otimes j}$,
	or equivalently, terms of the form $k \cdot b^{\otimes k}$, where $k$ is \textit{positive}.
	Set $C_3 \defeq A + B \in \ncdf{B}[3][1 - \degb{t}]$ and $C_2 = C \in \ncdf{B}[2][-\degb{t}]$.
	Then we have
	\begin{equation*}
		\begin{aligned}
			\partial_t \left( \G{b}[2] \right) ={} &
			\left( \ul{\dot{b}}  \pstar  \ul{\nu \left( \pstar \right)} \pstar \right)_{w} +
			\left(
			{\ul{\dot{b}} \pstar \ul{\nu \left( \pstar \right)} \hpstar} +
			{\ul{\dot{b}} \hpstar \ul{\nu \left( \pstar \right)} \pstar}
			\right)_{w}
			\\
			                                       & +
			\left(
			{\ul{b} \pstar \ul{\dot{b}} \pstar \nu \left( \pstar \right) \pstar} -
			{\ul{b} \pstar \nu \left( \pstar \right) \pstar \ul{\dot{b}} \pstar} +
			{\ul{b} \pstar \ul{\nu \left( \pstar \dot{b} \pstar \right)} \pstar}
			\right)_{w}
			\\
			                                       & - \qdr \left( A + B \right)
			\\
			={}                                    &
			\left( {\ul{\dot{b}}  \pstar  \ul{\nu \left( \pstar \right)} \pstar} +
			{\ul{\dot{b}} \hpstar \ul{\nu \left( \pstar \right)} \pstar} +
			{\ul{\dot{b}} \pstar \ul{\nu \left( \hpstar \right)} \pstar} +
			{\ul{\dot{b}} \pstar \ul{\nu \left( \pstar \right)} \hpstar}
			\right)_{w}
			\\
			                                       & + \clie{\nu} \left( C \right) - \qdr \left( A + B \right)
			\\
			={}                                    &
			\left( \ul{\dot{b}} \pstar \ul{\nu \left( \pstar \right)} \pstar \right) + \clie{\nu} \left( C_2 \right)
			- \qdr \left( C_3 \right).
		\end{aligned}
	\end{equation*}

	This calculation is equivalent to the one done by Cho and Lee in
	\cite[Lemma 3.4]{cho-homotopy-superpotential}.
	In particular, if we have a total inner product $\phi$ with $\phi_k = 0$ for $k \geq 3$,
	then $\phi_2 \circ \qdr = 0$ and $\phi_2 \circ \clie{\mu} = 0$, and we have
	\begin{equation*}
		\begin{aligned}
			\partial_t \left( \SP[b] \right) & =
			\partial_t \left( \phi_2 \left( \G{b}[2] \right) \right) =
			\phi_2 \left( \partial_t \left( \G{b}[2] \right) \right) =
			\phi_2 \left( \left( \ul{\dot{b}} \pstar \ul{\nu \left( \pstar \right)} \pstar \right) \right)
			\\
			                                 & =
			\sum_{i,j,k} \phi_2 \left(
			\ul{\dot{b}} \otimes b^{\otimes i} \otimes \ul{\nu_j \left( b^{\otimes j} \right)} \otimes
			b^{\otimes k}
			\right),
		\end{aligned}
	\end{equation*}
	which is consistent with \cref{thm:derivative-sp}.

	One can proceed similarly and show that $\partial_t \left( \G{b}[3] \right)$ has the form
	\begin{equation*}
		\partial_t \left( \G{b}[3] \right) = \left( \ul{\dot{b}} \pstar \ul{\nu \left( \pstar \right)} \pstar \ul{\nu \left( \pstar \right)} \pstar \right)
		+ \clie{\nu} \left( C_3 \right) - \qdr \left( C_4 \right),
	\end{equation*}
	for some explicit $C_4 \in \ncdf{B}[4][2 - \degb{t}]$, and so on, obtaining
	the formula \eqref{eq:clie-partial-t-Go-b} of \cref{rem:lie-partial-t-G-geq-2}.
	\end{rem}

	\subsection{Description Using \texorpdfstring{$\braidop_1$}{the Inner Product Parity Form}}
	\label{sec:description-using-braid-op-1}
	In this section, we give an equivalent description of our constructions
	using the inner product pairing $\braidop_1$ instead of the total degree pairing $\braidop_2$.
	We use the passage between the two descriptions in \cref{sec:sp-proof-intro-theorems} to deduce
	\Crefrange{thm:superpotential-properties}{thm:Gb-geq-2-properties} of the introduction
	which were stated working with $\braidop_1$,
	from the corresponding results we proved working with $\braidop_2$.

	In \Crefrange{sec:construction-special-elements}{sec:formal-derivative-sp}, we worked with the complexes
	\begin{equation*}
		\totcomp{\mathcal{A}, \braidop_2}[k][] = \totc{\ncdf{A, \braidop_2}, -\qdr^{\braidop_2},
			\clie{\mu}^{\braidop_2}}[][\geq k][\braidop_2]
	\end{equation*}
	for $k = 1, 2$ and their extended versions, constructed using the total degree parity form $\braidop_2$.
	Let us set
	\begin{equation*}
		\totcomp{\mathcal{A}, \braidop_1}[k][] \defeq \totc{\ncdf{A, \braidop_1}, \qdr^{\braidop_1},
			\clie{\mu}^{\braidop_1}}[][\geq k][\braidop_1].
	\end{equation*}
	By \cref{sec:dependence-ndf-braidop}, we have a natural isomorphism
	\begin{equation*}
		\Psi \colon \totc{\ncdfr{A, \braidop_2}}[][][\braidop_2] \rightarrow
		\totc{\ncdfr{A, \braidop_1}}[][][\braidop_1]
	\end{equation*}
	of graded Banach $R$-modules, given by
	\begin{equation} \label{eq:Psi-explicit-formula-in-generalized-sp}
		\begin{aligned}
			\Psi \left(
			\underbrace{
				\ul{a^1}_2 \otimes_2 l^1 \otimes_2 \dots \otimes_2
				\ul{a^k}_2 \otimes_2 l^k
			}_{x_k}
			\right) ={} &
			(-1)^{\varepsilon \left( x_k \right) + \left( k - 1 \right)}
			\\
			            & \qquad \s_k \left(
			\ul{a^1}_1 \otimes_1 l^1 \otimes_1 \dots \otimes_1
			\ul{a^k}_1 \otimes_1 l^k
			\right),
		\end{aligned}
	\end{equation}
	where
	\begin{equation*}
		\varepsilon \left( x_k \right) =
		\sum_{i=1}^{k} \left( k - i \right) \cdot \left( \degb{a^i} + \degb{l^i} \right),
	\end{equation*}
	and we use the notation $\otimes_i$ (resp.\ $\ul{a}_i$)
	to emphasize that the tensor product (resp.\ shift) is taken using
$\braidop_i$ for $i = 1, 2$.
	The isomorphism $\Psi$ commutes with the horizontal and vertical differentials,
	i.e., we have
	\begin{equation*}
		\Psi \circ \totc{-\qdr^{\braidop_2}}[][][\braidop_2] =
		\totc{\qdr^{\braidop_1}}[][][\braidop_1] \circ \Psi,
		\qquad
		\Psi \circ \totc{\clie{\mu}^{\braidop_2}}[][][\braidop_2] =
		\totc{\clie{\mu}^{\braidop_1}}[][][\braidop_1] \circ \Psi,
	\end{equation*}
	and hence $\Psi$ induces a chain map
	\begin{equation*}
		\Psi \colon \totc{\ncdf{A, \braidop_2}, -\qdr^{\braidop_2}, \clie{\mu}^{\braidop_2}}[][][\braidop_2]
		\rightarrow
		\totc{\ncdf{A, \braidop_1}, \qdr^{\braidop_1}, \clie{\mu}^{\braidop_1}}[][][\braidop_1]
	\end{equation*}
	between the total complexes.
	Truncating $\Psi$ induces an isomorphism
	\begin{equation*}
		\Psi_{\geq k} = \Psi_{\geq k}^{A}\colon \totcomp{\mathcal{A}, \braidop_2}[k][] \rightarrow \totcomp{\mathcal{A}, \braidop_1}[k][]
	\end{equation*}
	of differential graded Banach $\mathcal{R}$-modules acting by the same formula.

	The isomorphism $\Psi_{\geq k}$ respects the induced actions by
	Banach $\Ainf$-morphisms.
	Given a Banach $\Ainf$-morphism $f \colon \mathcal{A} \rightarrow \mathcal{B}$,
	and $i \in \Set{1,2}$, we have an induced morphism
$\cindmap{f}^{\braidop_i} \colon \ncdf{A, \braidop_i}[][] \rightarrow \ncdf{B, \braidop_i}[][]$
	which commutes with both $\qdr^{\braidop_i}$ and $\clie{\mu}^{\braidop_i}$.
	Hence, we get an induced chain map
	\begin{equation*}
		\totc{\cindmap{f}^{\braidop_i}}[][\geq k][\braidop_i] \colon
		\totcomp{\mathcal{A}, \braidop_i}[k][] \rightarrow \totcomp{\mathcal{B}, \braidop_i}[k][]
	\end{equation*}
	between the truncated total complexes, which for simplicity of notation, we continue
	to denote by $\cindmap{f}^{\braidop_i}$.\footnote{When working with $\braidop_2$, we denoted
	the map simply by $\cindmap{f}$, leaving the explicit dependence on the parity
	form $\braidop_2$ and $k$ to be determined from the context.}
	The isomorphism
$\Psi_{\geq k}$ is natural with respect to Banach $\Ainf$-morphisms, i.e., we have
$\Psi_{\geq k}^{B} \circ \cindmap{f}^{\braidop_2} = \cindmap{f}^{\braidop_1} \circ \Psi_{\geq k}^{A}$.

	Next, we construct an extended total complex for $\braidop_1$. We
	define the extended total complex $\totcompe{\mathcal{A}, \braidop_1}[2][]$ by
	adjoining a shifted copy of $R$ whose elements are denoted by $\ul{r}_1$, with
$\degb{\ul{r}_1} = \degb{r} - 1$, and the $R$-action
	given by $r \cdot \ul{r'}_1 = \ul{r \cdot r'}_1$.\footnote{The action of $R$ on the shifted copy
	was chosen to be consistent with the conventions
	for $\braidop_1$, see \cref{eq:R-action-ul-v}.}
	The differential of the generator $\ul{1}_1$ is defined to be
	\begin{equation*}
		\begin{aligned}
			D_{\mathcal{A}, \braidop_1} \left( \ul{1}_{1} \right)
			\defeq{} &
			\Psi_{\geq 2} \left( D_{\mathcal{A}, \braidop_2} \left( \ul{1}_{2} \right) \right)
			\stackrel{\eqref{eq:action-D-geq-2-on-ul-1}}{=}
			\Psi_{\geq 2} \left(
			-\sum_{k=2}^{\infty} \frac{1}{k} {\ulz{\mu}_2}^{\otimes_2 k}
			\right)
			\\
			={}      &
			\sum_{k = 2}^{\infty}
			(-1)^{\frac{(k-1)k}{2} + k} \frac{1}{k}
			\s_k {\ulz{\mu}_1}^{\otimes_1 k}
			\\
			={}      &
			\sum_{k = 2}^{\infty}
			(-1)^{\frac{k \left( k + 1 \right)}{2}} \frac{1}{k}
			\s_k {\ulz{\mu}_1}^{\otimes_1 k}.
		\end{aligned}
	\end{equation*}
	Then $\Psi_{\geq 2}$ extends to an isomorphism
$\Psi_{\geq 2}^{+} = \Psi_{\geq 2}^{+, A} \colon \totcompe{\mathcal{A}, \braidop_2}[2][] \rightarrow \totcompe{\mathcal{A}, \braidop_1}[2][]$
	between the extended total complexes, with $\Psi_{\geq 2}^{+} \left( \ul{1}_2 \right) = \ul{1}_1$.

	The extended isomorphism allows us to transport the cyclic Chern--Simons form.
	Given $b \in \tc{A}$,
	let $\G{b}[\geq 2][\mathcal{A}, \braidop_1] \in \totcompe{\mathcal{A}, \braidop_1}[2][-1]$
	be the chain corresponding to the cyclic Chern--Simons form of
	\cref{dfn:cyclic-chern-simons-form-braidop-2}
	under the isomorphism $\Psi_{\geq 2}^{+}$.
	Explicitly, we have
	\begin{align}
		\G{b}[\geq 2][\mathcal{A}, \braidop_1] \defeq{} &
		\Psi_{\geq 2}^{+} \left( \G{b}[\geq 2][\mathcal{A}, \braidop_2] \right) =
		\ul{1}_1 + \Psi_{\geq 2} \left( \Go{b}[\geq 2][\mathcal{A}, \braidop_2] \right) =
		\notag
		\\
		\stackrel{\phantom{\star}}{=}{}                 &
		\ul{1}_1 +
		\sum_{\substack{k=2 \\ i_2,\dots,i_k = 0 \\ j_1,\dots,j_k=0}}^{\infty}
		\frac{(-1)^{k - 1 + \frac{(k - 2) \cdot (k - 1)}{2}}}{1 + \sum_{r = 1}^k j_r + \sum_{r=2}^{k} i_r}
		\notag
		\\
		                                                & \quad\quad
		\s_k \left(
		\ul{b}_1 \otimes_1 b^{\otimes_1 \, j_1} \otimes_1 \ul{ \mu_{i_2} \left( b^{\otimes i_2} \right)}_1
		\otimes_1 b^{\otimes_1 \, j_2} \otimes_1 \dots \otimes_1
		\ul{ \mu_{i_k} \left( b^{\otimes i_k} \right) }_1 \otimes_1 b^{\otimes_1 \, j_k}
		\right)
		\notag
		\\
		\stackrel{\star}{=}{}                           &
		\ul{1}_1 +
		\sum_{\substack{k=2 \\ i_1,\dots,i_{k-1} = 0 \\ j_1,\dots,j_k=0}}^{\infty}
		\frac{(-1)^{\frac{(k - 2) \cdot (k - 1)}{2}}}{1 + \sum_{r=1}^{k-1} i_r + \sum_{r = 1}^k j_r}
		\label{eq:G2-braid-op-1-explicit-formula}
		\\
		\notag
		                                                & \quad\quad
		\s_k \left(
		\ul{ \mu_{i_1} \left( b^{\otimes i_1} \right)}_{1} \otimes_1 b^{\otimes_1 \, j_1} \otimes_1 \dots
		\otimes \ul{ \mu_{i_{k-1}} \left( b^{\otimes i_{k-1}} \right) }_{1} \otimes_1 b^{\otimes_1 \, j_{k-1}}
		\otimes_1 \ul{b}_1 \otimes_1 b^{\otimes_1 \, j_k}
		\right),
	\end{align}
	where in $\star$ we used the rotation operator suitable for $\braidop_1$ to rewrite the terms
	so that they start with $\corest{\mu}$, and adjusted the sign accordingly.

	We can now state the formal definitions of the corresponding
	notions from previous sections, working with $\braidop_1$:
	\begin{dfn} \label{dfn:cyclic-chern-simons-form-braidop-1}
	The chain $\G{b}[\geq 2][\mathcal{A}, \braidop_1]$ is called
	the \textbf{cyclic Chern--Simons form} in the total complex $\totcompe{\mathcal{A}, \braidop_1}[2][]$.
	\end{dfn}
	\begin{dfn} \label{dfn:generalized-inner-product-braidop-1}
	An $n$\textbf{-dimensional total inner product} on $\mathcal{A}$ is a
	morphism
	\begin{equation*}
		\phi^{\braidop_1} \colon \totcompe{\mathcal{A}, \braidop_1}[2][] \rightarrow \mathcal{R}[4-n]
	\end{equation*}
	of differential graded Banach $\mathcal{R}$-modules.
	An $n$-\textbf{dimensional total inner product Banach} $\Ainf$-\textbf{algebra over} $\mathcal{R}$
	is a triple $\mathcal{A} = \left( A, \mu, \phi^{\braidop_1} \right)$ where $\left( A, \mu \right)$ is a
	Banach $\Ainf$-algebra over $\mathcal{R}$ and $\phi^{\braidop_1}$ is an $n$-dimensional total inner product
	on $\left( A, \mu \right)$.
	\end{dfn}

	\begin{dfn} \label{dfn:superpotential-braidop-1}
	The \textbf{superpotential function}
$\SP[][\mathcal{A}] \colon \tc{A} \rightarrow R^{3-n}$
	associated to the $n$-dimensional total inner product Banach $\Ainf$-algebra
$\mathcal{A} = \left( A, \mu, \phi^{\braidop_1} \right)$
	is defined by
	\begin{equation*}
		\SP[b][\mathcal{A}] \defeq
		\phi^{\braidop_1} \left( \G{b}[\geq 2][\mathcal{A}, \braidop_1] \right).
	\end{equation*}
	\end{dfn}

	\Crefrange{dfn:cyclic-chern-simons-form-braidop-1}{dfn:superpotential-braidop-1} were
	used in \cref{sec:statement-results}, where we gave an overview of our results.
	For the explicit relations a total
	inner product $\phi^{\braidop_1}$ must satisfy, see \cref{sec:total-inner-product-explicit-relations}.

	We can relate total inner products across the two parity forms.
	Given an $R$-linear map $\phi^{\braidop_1} \colon \totcompe{A, \braidop_1}[2][] \rightarrow R[4-n]$,
	we define
	\begin{equation} \label{eq:phi-braidop-2-from-braidop-1}
		\phi^{\braidop_2} \defeq \phi^{\braidop_1} \circ \Psi_{\geq 2}^{+} \colon \totcompe{A, \braidop_2}[2][] \rightarrow R[4-n].
	\end{equation}
	Then $\phi^{\braidop_1}$ is a total inner product in the sense of \cref{dfn:generalized-inner-product-braidop-1}
	if and only if
$\phi^{\braidop_2}$ is a total inner product in the sense of \cref{dfn:generalized-inner-product},
	and we have a bijection between the two notions.
	To differentiate between the two different notions of total inner products, we will always
	decorate them with the superscript $\braidop_i$, with $i \in \Set{1, 2}$, indicating
	the parity form we work with.
	When $\phi^{\braidop_1}$ is a total
	inner product, we will call $\phi^{\braidop_2}$ the \textbf{corresponding} total inner product
	to $\phi^{\braidop_1}$.
	Given a total inner product Banach $\Ainf$-algebra $\mathcal{A} = \left( A, \mu, \phi^{\braidop_1} \right)$
	in the sense of \cref{dfn:generalized-inner-product-braidop-1},
	we will denote by $\mathcal{A}^{\sharp} = \left( A, \mu, \phi^{\braidop_2} \right)$ the
	corresponding total inner product Banach $\Ainf$-algebra.
	With the definitions above, we have
	\begin{equation} \label{eq:sp-A-vs-A-sharp}
		\begin{aligned}
			\SP[b][\mathcal{A}^{\sharp}] & =
			\phi^{\braidop_2} \left( \G{b}[\geq 2][\mathcal{A}, \braidop_2] \right) =
			\left( \phi^{\braidop_1} \circ \Psi_{\geq 2}^{+} \right) \left(
			\G{b}[\geq 2][\mathcal{A}, \braidop_2]
			\right)
			\\
			                             & =
			\phi^{\braidop_1} \left( \G{b}[\geq 2][\mathcal{A}, \braidop_1] \right) =
			\SP[b][\mathcal{A}],
		\end{aligned}
	\end{equation}
	so the superpotential function in the sense of \cref{dfn:superpotential} for
$\mathcal{A}^{\sharp}$ coincides
	with the superpotential function in the sense of \cref{dfn:superpotential-braidop-1} for $\mathcal{A}$.

	We now define morphisms of total inner products directly in the $\braidop_1$ setting.
	Given a Banach $\Ainf$-morphism $f \colon \mathcal{A} \rightarrow \mathcal{B}$,
	there is a unique way to extend the induced chain map
$\cindmap{f}^{\braidop_1} \colon \totcomp{\mathcal{A}, \braidop_1}[2][] \rightarrow \totcomp{\mathcal{B}, \braidop_1}[2][]$
	to a chain map
	\begin{equation*}
		\cindmap{f}^{+,\braidop_1} \colon \totcompe{\mathcal{A}, \braidop_1}[2][] \rightarrow
		\totcompe{\mathcal{B}, \braidop_1}[2][]
	\end{equation*}
	which satisfies\footnote{As noted in \cref{rem:totcompe-geq-2-not-functorial-nose},
	the construction $\mathcal{A} \mapsto \totcompe{\mathcal{A}, \braidop_i}[2][]$ is not a functor,
	so we can't say that $\Psi_{\geq 2}^{+}$ is a natural transformation, but it becomes a natural transformation
	if one takes the cohomology of $\totcompe{\mathcal{A}, \braidop_i}[2][]$.}
	\begin{equation} \label{eq:quasi-naturality-psi-2-plus}
		\Psi_{\geq 2}^{+, B} \circ \cindmap{f}^{+, \braidop_2} =
		\cindmap{f}^{+, \braidop_1} \circ \Psi_{\geq 2}^{+, A}.
	\end{equation}
	The extension is obtained by setting
	\begin{equation*}
		\cindmap{f}^{+,\braidop_1} \left( \ul{1}_1 \right) \defeq
		\Psi_{\geq 2}^{+, B} \left( \cindmap{f}^{+, \braidop_2} \left( \ul{1}_2 \right) \right)
		\stackrel{\eqref{def:cycl-f-ul-1-geq2}}{=}
		\Psi_{\geq 2}^{+, B} \left( \G{f_0 \left( 1 \right)}[\geq 2][\mathcal{B}, \braidop_2] \right) =
		\G{f_0 \left( 1 \right)}[\geq 2][\mathcal{B}, \braidop_1],
	\end{equation*}
	which is the same formula as \eqref{def:cycl-f-ul-1-geq2}, working
	with $\G{f_0 \left( 1 \right)}[\geq 2][\mathcal{B}, \braidop_1]$ instead of
$\G{f_0 \left( 1 \right)}[\geq 2][\mathcal{B}, \braidop_2]$.
	\begin{dfn} \label{dfn:morphism-total-inner-products-braidop-1}
	Given two total inner product Banach $\Ainf$-algebras
$\mathcal{A} = \left( A, \mu_A, \phi_A^{\braidop_1} \right)$ over $\mathcal{R}$
	and $\mathcal{B} = \left( B, \mu_B, \phi_B^{\braidop_1} \right)$ over $\mathcal{S}$,
	a \textbf{morphism of total inner product Banach} $\Ainf$-\textbf{algebras},
	denoted $f \colon \mathcal{A} \rightarrow \mathcal{B}$,
	is a morphism $f \colon \left( A, \mu_A \right) \rightarrow \left( B, \mu_B \right)$
	of Banach $\Ainf$-algebras which satisfies
	\begin{equation} \label{eq:morphism-total-inner-products-braidop-1}
		\left( \phi_B^{\braidop_1} \right)^{\diamond} \circ \cindmap{f}^{+,\braidop_1} =
		\base{f}^{\diamond}[4 - n] \circ \left( \phi_A^{\braidop_1} \right)^{\diamond}.
	\end{equation}
	\end{dfn}
	This is the same as \cref{dfn:morphism-total-inner-products}, replacing
$\cindmap{f}^{+, \braidop_2}$ with $\cindmap{f}^{+, \braidop_1}$.
	It follows from \eqref{eq:quasi-naturality-psi-2-plus} that a morphism
$f \colon \left( A, \mu_A \right) \rightarrow \left( B, \mu_B \right)$ of Banach $\Ainf$-algebras
	is a morphism
$f \colon \mathcal{A} \rightarrow \mathcal{B}$ in the sense of \cref{dfn:morphism-total-inner-products-braidop-1},
	if and only if
$f \colon \mathcal{A}^{\sharp} \rightarrow \mathcal{B}^{\sharp}$ is a morphism
	in the sense of \cref{dfn:morphism-total-inner-products}.

	Working with unital $\Ainf$-algebras, we can define a reduced version
$\totcompered{\mathcal{A}}[2][]$ in the $\braidop_1$ setting.
	Assume $\mathcal{A} = \left( A, \mu, e \right)$ is unital. Given $k \geq 2$,
	denote by $\degenu[A][k][]^{\braidop_1}$ the closure of the image of the map
	\begin{equation*}
		\rest{\left( \ccont{e}^{\braidop_1} \right)^{k}}{\ncdf{A}[0][]} \colon \ncdf{A}[0][] \rightharpoonup \ncdf{A, \braidop_1}[k][],
	\end{equation*}
	where we think of $e$ as the coderivation $\nu_e$ of \cref{sec:extended-reduced-total-complexes}.
	Consider the graded Banach $R$-submodule of $\totcomp{A, \braidop_1}[2][]$ given by
	\begin{equation}
		\degenu[A][][\geq 2]^{\braidop_1} \defeq \bigoplus_{k \geq 2} \s_k \left( \degenu[A][k][]^{\braidop_1} \right).
	\end{equation}
	Since the commutation relations of \cref{lm:coder-e-op-comm-relations} hold for both
	parity forms, where only the interpretation of the graded commutator is different,
	we have that $\degenu[A][][\geq 2]^{\braidop_1}$ is a subcomplex of the total complex $\totcomp{\mathcal{A}, \braidop_1}[2][]$.
	One can verify that given $x_0 \in \ncdf{A}[0][]$, we have
	\begin{equation*}
		\Psi_{\geq 2} \left( \left( \ccont{e}^{\braidop_2} \right)^k \left( x_0 \right) \right)
		=
		(-1)^{\frac{k(k+1)}{2}} \s_k \, \left( \ccont{e}^{\braidop_1} \right)^k \left( x_0 \right),
	\end{equation*}
	and hence $\Psi_{\geq 2}$ maps $\degenu[A][][\geq 2]^{\braidop_2}$ (given by \eqref{eq:degenu-geq-2})
	onto $\degenu[A][][\geq 2]^{\braidop_1}$.
	The \textbf{extended reduced} total complex associated to $\mathcal{A}$ via $\braidop_1$ is defined to be
	\begin{equation*}
		\totcompered{\mathcal{A}, \braidop_1}[2][] \defeq \totcompe{\mathcal{A}, \braidop_1}[2][] / \degenu[A][][\geq 2]^{\braidop_1},
	\end{equation*}
	and the map $\Psi_{\geq 2}^{+}$ induces an isomorphism
$\Psi_{\geq 2}^{+} \colon \totcompered{\mathcal{A}, \braidop_2}[2][] \rightarrow \totcompered{\mathcal{A}, \braidop_1}[2][]$.

	A total inner product $\phi^{\braidop_1}$ on $\mathcal{A}$ is called \textbf{unital} if
	\begin{equation} \label{eq:total-inner-product-i-e-cond-braidop-1}
		\phi^{\braidop_1} \circ \s_k \circ \rest{\left( \ccont{e}^{\braidop_1} \right)^k}{\ncdf{A}[0][]} = 0
	\end{equation}
	for all $k \geq 2$, which is equivalent to the condition that $\phi^{\braidop_1}$
	descends to $\totcompered{\mathcal{A}, \braidop_1}[2][]$. Note that
	this is the same condition \eqref{eq:total-inner-product-i-e-cond} as for $\phi^{\braidop_2}$,
	except we use $\ccont{e}^{\braidop_1}$ instead of $\ccont{e}^{\braidop_2}$ and add
	suspensions because of our conventions. The total inner product $\phi^{\braidop_1}$
	is unital in the sense of \cref{eq:total-inner-product-i-e-cond-braidop-1} if and only if
	the corresponding total inner product $\phi^{\braidop_2}$ is unital in the sense of
	\cref{eq:total-inner-product-i-e-cond}.

	Finally, we describe the modifications needed when working
	with pseudoisotopies in the $\braidop_1$ setting.
	We start with the notion of
	a $\cohom{}[] \totcompe{}[2][]$-strong pseudoisotopy.
	Consider the two functors
	\begin{equation*}
		\mathcal{F}_1 \colon \mathcal{A} \mapsto \cohom{}[] \totcompe{\mathcal{A}, \braidop_1}[2][],
		\qquad
		\mathcal{F}_2 \colon \mathcal{A} \mapsto \cohom{}[] \totcompe{\mathcal{A}, \braidop_2}[2][],
	\end{equation*}
	and the corresponding notions of $\mathcal{F}_i$-strong pseudoisotopies in
	the sense of \cref{dfn:F-strong-pseudoisotopy}, for $i = 1, 2$.
	When working with $\braidop_2$, we used the functor $\mathcal{F}_2$, while
	if we work with $\braidop_1$, we can use the functor $\mathcal{F}_1$.
	Given a pseudoisotopy $\mathfrak{A}$ between $\mathcal{A}_0$ and $\mathcal{A}_1$,
	\cref{eq:quasi-naturality-psi-2-plus} implies that we have
	\begin{equation*}
		\cohom{\Psi_{\geq 2}^{+,A_j}}[] \circ \mathcal{F}_2 \left( \evalmf^j \right) =
		\mathcal{F}_1 \left( \evalmf^j \right) \circ \cohom{\Psi_{\geq 2}^{+,\mathfrak{A}}}[]
	\end{equation*}
	for $j = 0, 1$, and hence $\mathfrak{A}$ is $\mathcal{F}_1$-strong if and only if
	it is $\mathcal{F}_2$-strong, so both notions are equivalent, and we can
	call them $\cohom{}[] \totcompe{}[2][]$-\textbf{strong} without any ambiguity.
	A $\mathcal{F}_i$-strong pseudoisotopy induces
	an isomorphism
$\mathfrak{a}^{\braidop_i} \colon \cohom{}[] \totcompe{\mathcal{A}_0, \braidop_i}[2][]
\rightarrow \cohom{}[] \totcompe{\mathcal{A}_1, \braidop_i}[2][]$ for $i = 1, 2$,
	and we have
	\begin{equation} \label{eq:mathfrak-a-psi}
		\mathfrak{a}^{\braidop_1} \circ \cohom{\Psi_{\geq 2}^{+,A_0}}[] =
		\cohom{\Psi_{\geq 2}^{+,A_1}}[] \circ \mathfrak{a}^{\braidop_2}.
	\end{equation}
	When $\mathcal{A}_0, \mathcal{A}_1$ and $\mathfrak{A}$ are unital, we work
	with the functors
$\mathcal{A} \xmapsto{\mathcal{F}_i} \cohom{}[] \totcompered{\mathcal{A}, \braidop_i}[2][]$,
	and the statements above hold for $\cohom{}[] \totcompered{}[2][]$-strong pseudoisotopies
	with obvious modifications.

	Now assume that $\mathcal{A}_0$ and $\mathcal{A}_1$ are
	two total inner product Banach $\Ainf$-algebras over $\mathcal{S}$. A
	\textbf{pseudoisotopy} $\mathfrak{A}$ between $\mathcal{A}_0$ and $\mathcal{A}_1$
	is defined as in \cref{dfn:pseudo-isotopy-total-inner-product-algebras},
	by requiring the pseudoisotopy maps $\evalmf^i$ to be
	morphisms of total inner product Banach $\Ainf$-algebras in the sense of
	\cref{dfn:morphism-total-inner-products-braidop-1}. As a consequence,
$\mathfrak{A}$ is a pseudoisotopy between $\mathcal{A}_0$ and $\mathcal{A}_1$
	if and only if $\mathfrak{A}^{\sharp}$ is a pseudoisotopy between
$\mathcal{A}_0^{\sharp}$ and $\mathcal{A}_1^{\sharp}$ in the sense
	of \cref{dfn:pseudo-isotopy-total-inner-product-algebras}.
	When $\mathcal{A}_0, \mathcal{A}_1$ are unital, then $\mathfrak{A}$
	is a unital pseudoisotopy if and only if $\mathfrak{A}^{\sharp}$
	is a unital pseudoisotopy.

	\subsubsection{Proofs of the Main Theorems} \label{sec:sp-proof-intro-theorems}

	\begin{proof}[Proof of \cref{thm:superpotential-properties}]
	\Cref{item:superpotential-1} follows from the corresponding \cref{lm:functoriality-superpotential}
	by applying it to $\mathcal{A}^{\sharp}$ and $\mathcal{B}^{\sharp}$ and using \eqref{eq:sp-A-vs-A-sharp}.
	Similarly, \cref{item:superpotential-2} follows from
	the corresponding \cref{lm:superpotential-bounding-chain-closed}
	by applying it to $\mathcal{A}^{\sharp}$ and using \eqref{eq:sp-A-vs-A-sharp}.
	Finally, \cref{item:superpotential-3} follows from
	the corresponding \cref{thm:invariance-superpotential}
	by applying it to $\mathcal{A}_0^{\sharp}, \mathcal{A}_1^{\sharp}$ and $\mathfrak{A}^{\sharp}$,
	and using \eqref{eq:sp-A-vs-A-sharp}.
	\end{proof}

	\begin{proof}[Proof of \cref{thm:superpotential-derivative}]
	\phantomsection\label{proof:thm:superpotential-derivative}
	Let $\mathcal{A} = \left( A, \mu, \phi_A^{\braidop_1} \right)$ be a total inner product
	Banach $\Ainf$-algebra over a differential graded-commutative $\mathbbm{k}$-algebra
$\mathcal{R}$. Given a morphism $i \colon \mathcal{R} \rightarrow \mathcal{S}$ of
	differential graded-commutative Banach $\mathbbm{k}$-algebras,
	there is a natural notion of scalar extension
$i_{!} \left( \phi_A^{\braidop_1} \right)$ for total inner products defined using $\braidop_1$.
	We just follow the isomorphisms as in \eqref{eq:totcompe-commutes-scalar-extension},
	working with $\totcompe{\cdot, \braidop_1}[2][]$ instead of $\totcompe{\cdot, \braidop_2}[2][]$.
	The isomorphism $\Psi_{\geq 2}^{+}$ is natural with respect to both notions of scalar
	extension, i.e., we have
	\begin{equation} \label{eq:phi-Psi-scalar-extension}
		\left( i_{!} \left( \phi_A^{\braidop_1} \right) \right)^{\braidop_2} =
		i_{!} \left( \phi_A^{\braidop_2} \right).
	\end{equation}
	If we denote by $i_{!} \left( \mathcal{A} \right)$ (resp.\ $i_{!} \left( \mathcal{A}^{\sharp} \right)$)
	the scalar extension of $\mathcal{A}$ (resp.\ $\mathcal{A}^{\sharp}$), \eqref{eq:phi-Psi-scalar-extension}
	can be written succinctly as
	\begin{equation} \label{eq:scalar-extension-sharp}
		i_{!} \left( \mathcal{A}^{\sharp} \right) = i_{!} \left( \mathcal{A} \right)^{\sharp}.
	\end{equation}

	Now assume $i \colon R \rightarrow \pows{R}[t]$ where $t$ is an even formal variable,
	as described in the beginning of \cref{sec:derivative-calc}, and write
	\begin{equation*}
		\mathcal{B} = i_{!} \left( \mathcal{A} \right) =
		\left( i_{!} \left( A \right), i_{!} \left( \mu \right), i_{!} \left( \phi_A^{\braidop_1} \right) \right)
		= \left( B, \nu, \phi_B^{\braidop_1} \right),
	\end{equation*}
	as in \cref{thm:derivative-sp}.
	Given $\mathcal{A}$, take $\mathcal{A}^{\sharp}$ and apply \cref{thm:derivative-sp}
	with $\mathcal{A}^{\sharp}$ and $i_{!} \left( \mathcal{A}^{\sharp} \right)$
	to obtain a formula for the derivative of the superpotential on
$i_{!} \left( \mathcal{A}^{\sharp} \right)$.
	We have
	\begin{equation*}
		\begin{aligned}
			\partial_t \left( \SP[b][i_{!} \left( \mathcal{A} \right)] \right)
			\eqwithref[eq:sp-A-vs-A-sharp]                        &
			\partial_t \left( \SP[b][i_{!} \left( \mathcal{A} \right)^{\sharp}] \right)
			\stackrel{\eqref{eq:scalar-extension-sharp}}{=}
			\partial_t \left( \SP[b][i_{!} \left( \mathcal{A}^{\sharp} \right)] \right)
			\\
			\eqwithref[eq:derivative-sp-braid-2]                  &
			i_{!} \left( \phi_A^{\braidop_2} \right) \left(
			\sum_{k=2}^{\infty}
			\left( \ul{ \partial_t \left( b \right) } \otimes \Exp{b} \right) \otimes
			\left( \ul{ \corest{\nu} \left( \Exp{b} \right) } \otimes \Exp{b} \right)^{\otimes \left( k - 1 \right)}
			\right)
			\\
			\eqwithref[eq:phi-Psi-scalar-extension]               &
			\left( i_{!} \left( \phi_A^{\braidop_1} \right) \right)^{\braidop_2}
			\left(
			\sum_{k=2}^{\infty}
			\left( \ul{ \partial_t \left( b \right) } \otimes \Exp{b} \right) \otimes
			\left( \ul{ \corest{\nu} \left( \Exp{b} \right) } \otimes \Exp{b} \right)^{\otimes \left( k - 1 \right)}
			\right)
			\\
			\eqwithref[eq:phi-braidop-2-from-braidop-1]           &
			\phi_B^{\braidop_1} \left(
			\Psi_{\geq 2}^{+} \left(
			\sum_{k=2}^{\infty}
			\left( \ul{ \partial_t \left( b \right) } \otimes \Exp{b} \right) \otimes
			\left( \ul{ \corest{\nu} \left( \Exp{b} \right) } \otimes \Exp{b} \right)^{\otimes \left( k - 1 \right)}
			\right)
			\right)
			\\
			\eqwithref[eq:Psi-explicit-formula-in-generalized-sp] &
			\phi_B^{\braidop_1} \left(
			                    \sum_{k=2}^{\infty} (-1)^{\frac{(k-2)(k-1)}{2} + k - 1}
			\s_k \left(
			\left( \ul{ \partial_t \left( b \right) } \otimes \Exp{b} \right) \otimes
			\left( \ul{ \corest{\nu} \left( \Exp{b} \right) } \otimes \Exp{b} \right)^{\otimes \left( k - 1 \right)}
			\right)
			\right)
			\\
			\eqwithref                                            &
			\phi_B^{\braidop_1} \left(
			                    \sum_{k=2}^{\infty} (-1)^{\frac{(k-2)(k-1)}{2}}
			\s_k \left(
			\left( \ul{ \corest{\nu} \left( \Exp{b} \right) } \otimes \Exp{b} \right)^{\otimes \left( k - 1 \right)}
			\otimes
			\ul{ \partial_t \left( b \right) } \otimes \Exp{b}
			\right)
			\right).
		\end{aligned}
	\end{equation*}
	The formula above is precisely formula \eqref{eq:derivative-sp-braid-1} we gave for the derivative in
	\cref{thm:superpotential-derivative}, up to replacing $\mu$ with $\nu$ and $\phi_B^{\braidop_1}$ with $\phi$,
	as in the statement of the theorem.
	\end{proof}

	\begin{proof}[Proof of \cref{thm:Gb-geq-2-properties}]
	\Cref{item:Gb-geq-2-properties-1} follows from the corresponding \cref{eq:cindmap-f-G-geq-2} of \cref{lm:G-geq-2-properties},
	using the fact that $\Psi_{\geq 2}^{+}$ is a chain map which commutes with the induced
	morphisms between the extended complexes. More precisely, we have
	\begin{equation*}
		\begin{aligned}
			\cindmap{f}^{+,\braidop_1} \left( \G{b}[\geq 2][\mathcal{A}, \braidop_1] \right)
			\eqwithref[eq:G2-braid-op-1-explicit-formula] &
			\cindmap{f}^{+,\braidop_1} \left( \Psi_{\geq 2}^{+} \left( \G{b}[\geq 2][\mathcal{A}, \braidop_2] \right) \right)
			\\
			\eqwithref[eq:quasi-naturality-psi-2-plus]    &
			\Psi_{\geq 2}^{+} \left( \cindmap{f}^{+,\braidop_2} \left( \G{b}[\geq 2][\mathcal{A}, \braidop_2] \right) \right)
			\\
			\eqwithref[eq:cindmap-f-G-geq-2]              &
			\Psi_{\geq 2}^{+} \left(
			\G{\mcfunc{f} \left( b \right)}[\geq 2][\mathcal{B}, \braidop_2] +
			D_{\mathcal{B}, \braidop_2} \left( R \left( b; f \right) \right)
			\right)
			\\
			\eqwithref[eq:G2-braid-op-1-explicit-formula] &
			\G{\mcfunc{f} \left( b \right)}[\geq 2][\mathcal{B}, \braidop_1] +
			D_{\mathcal{B}, \braidop_1} \left( \Psi_{\geq 2}^{+} \left( R \left( b; f \right) \right) \right).
		\end{aligned}
	\end{equation*}

	\Cref{item:Gb-geq-2-properties-2} immediately follows from the corresponding
	\cref{lm:G-geq-2-bounding-chain-closed} using the fact that $\Psi_{\geq 2}^{+}$
	is a chain map.

	Finally, \cref{item:Gb-geq-2-properties-3} follows from the corresponding
	\cref{lm:invariance-Gb-geq-2-pseudoisotopy}	by applying \cref{eq:mathfrak-a-psi}
	and using \cref{dfn:cyclic-chern-simons-form-braidop-1} of the
	cyclic Chern--Simons form in $\totcompe{\mathcal{A}, \braidop_1}[2][]$.
	\end{proof}

	\subsection{Total Inner Product Components and Relations} \label{sec:total-inner-product-explicit-relations}
	In this section, we unwind our definitions and give an explicit description of the data contained
	in the notion of a total inner product in terms of its components and the relations
	satisfied by the components. We also write down explicit formula for the superpotential
	in terms of the components of the total inner product.
	To compare the relations satisfied by total inner products to
	the relations satisfied by other notions appearing in the
	introduction and the literature, we will work with the inner product parity form $\braidop_1$.

	Given an $n$-dimensional total inner product
	\begin{equation*}
		\phi = \phi^{\braidop_1} \colon \totcompe{\mathcal{A}, \braidop_1}[2][] \rightarrow \mathcal{R}[4-n]
	\end{equation*}
	and $k \geq 2$, we define the \textbf{components} $\phi_k \colon \ncdf{A, \braidop_1}[k][] \rightharpoonup R$
	of $\phi$, which are $R$-linear contractive maps of degree $4 - n - k$, by the composition
	\begin{equation} \label{eq:components-phi_k}
		\begin{tikzcd}[column sep=small]
			{\ncdf{A, \braidop_1}[k][]} & {\ncdf{A, \braidop_1}[k][][k]} &
			{\totcompe{A, \braidop_1}[2][]} & {R[4 - n]} & R.
			\arrow[harpoon, from=1-1, to=1-2, "\s_k \,\,", "{[-k]}"']
			\arrow["\phi_k", bend left=20, from=1-1, to=1-5, start anchor=north, end anchor=north,
				harpoon, "{[4 - n - k]}"']
			\arrow[from=1-2, to=1-3, hookrightarrow]
			\arrow[from=1-3, to=1-4, "\phi"]
			\arrow[harpoon, from=1-4, to=1-5, "{[4 - n]}"']
		\end{tikzcd}
	\end{equation}
	Define also the $\ul{1}$-\textbf{component} $\phi_{\ul{1}} \in R^{3 - n}$ by the equation
	\begin{equation*}
		\phi \left( \ul{1} \right)  \defeq \s_{4 - n} \phi_{\ul{1}} \in R[4 - n]^{-1} = R^{3-n}.
	\end{equation*}
	The map $\phi$ is completely determined by the sequence of components
$\left( \phi_{\ul{1}}, \phi_2, \phi_3, \dots \right)$, and the requirement that $\phi$ is a chain map,
	i.e., $d_{\mathcal{R}[4-n]} \circ \phi = \phi \circ D_{\mathcal{A}}$, translates into the following relations in terms of the components of $\phi$:
	\begin{align}
		 & \begin{aligned}\label{eq:d-phi-ul-1}
		                  (-1)^{4 - n} d \left( \phi_{\ul{1}} \right) \, & =
		                  \sum_{k=2}^{\infty} \frac{(-1)^{\frac{k \left( k + 1 \right)}{2}}}{k}
		                  \phi_k \left( \ulz{\mu}^{\otimes k} \right),
		   \end{aligned}
		\\
		 & \begin{aligned} \label{eq:d-phi-2}
		                   (-1)^{4 - n} d \circ \phi_2 & = \phi_2 \circ \clie{\mu},
		   \end{aligned}
		\\
		 & \begin{aligned} \label{eq:d-phi-3}
		                   (-1)^{4 - n} d \circ \phi_3 & = - \phi_3 \circ \clie{\mu} + \phi_2 \circ \qdr,
		   \end{aligned}
		\\
		 & \begin{aligned} \label{eq:d-phi-k}
		                   (-1)^{4 - n} d \circ \phi_k & = (-1)^k \phi_k \circ \clie{\mu} + \phi_{k-1} \circ \qdr,
		                   \qquad k \geq 4.
		   \end{aligned}
	\end{align}
	Let us write down the relations of \cref{eq:d-phi-2,eq:d-phi-3} explicitly.
	The relation \eqref{eq:d-phi-2} reads
	\begin{align}
		(-1)^{4-n} d \left( \phi_2 \left( \ul{x}, a, \ul{y}, b \right) \right)
		\stackrel{\eqref{eq:cyc-lie-formula}}{=}{} &
		(-1)^{\degb{x} + \degb{a_{(1)}}}
		\phi_2 \left( \ul{x}, a_{(1)}, \corest{\mu} \left( a_{(2)} \right), a_{(3)}, \ul{y}, b \right)
		\notag
		\\
		                                           & +
		                                             (-1)^{\degb{x} + \degb{a_{(1)}}}
		\phi_2 \left( \ul{x}, a_{(1)}, \ul{\corest{\mu} \left( a_{(2)}, y, b_{(1)} \right)}, b_{(2)} \right)
		\label{eq:d-phi-2-explicit}
		\\
		                                           & +
		                                             (-1)^{\degb{x} + \degb{a} + \degb{y} + \degb{b_{(1)}}}
		\phi_2 \left( \ul{x}, a, \ul{y}, b_{(1)}, \corest{\mu} \left( b_{(2)} \right), b_{(3)} \right)
		\notag
		\\
		                                           & +
		                                             (-1)^{\degb{b_{(2)}} \cdot \left( \degb{x} + \degb{a} + \degb{y} + \degb{b_{(1)}} \right)}
		\phi_2 \left(
		\ul{\corest{\mu} \left( b_{(2)}, x, a_{(1)} \right)}, a_{(2)}, \ul{y}, b_{(1)}
		\right)
		\notag
	\end{align}
	for $x,y \in A$ and $a, b \in \tens{A}$. Similarly, relation \eqref{eq:d-phi-3} reads
	\begin{equation} \label{eq:d-phi-3-explicit}
		\begin{aligned}
			\MoveEqLeft
			(-1)^{4-n} d \left( \phi_3 \left( \ul{x}, a, \ul{y}, b, \ul{z}, c \right) \right)
			\stackrel{\eqref{eq:cyc-lie-formula}}{=}{}
			\\
			={} &
			(-1)^{\degb{x} + \degb{a_{(1)}} + 1}
			\phi_3 \left(
			\ul{x}, a_{(1)}, \corest{\mu} \left( a_{(2)} \right), a_{(3)}, \ul{y}, b, \ul{z}, c
			\right)
			\\
			    & +
			      (-1)^{\degb{x} + \degb{a_{(1)}} + 1}
			\phi_3 \left(
			\ul{x}, a_{(1)}, \ul{\corest{\mu} \left( a_{(2)}, y, b_{(1)} \right)}, b_{(2)}, \ul{z}, c
			\right)
			\\
			    & +
			      (-1)^{\degb{x} + \degb{a} + \degb{y} + \degb{b_{(1)}} + 1}
			\phi_3 \left(
			\ul{x}, a, \ul{y}, b_{(1)}, \corest{\mu} \left( b_{(2)} \right), b_{(3)}, \ul{z}, c
			\right)
			\\
			    & +
			      (-1)^{\degb{x} + \degb{a} + \degb{y} + \degb{b_{(1)}} + 1}
			\phi_3 \left(
			\ul{x}, a, \ul{y}, b_{(1)}, \ul{\corest{\mu} \left( b_{(2)}, z, c_{(1)} \right)}, c_{(2)}
			\right)
			\\
			    & +
			      (-1)^{\degb{x} + \degb{a} + \degb{y} + \degb{b} + \degb{z} + \degb{c_{(1)}} + 1}
			\phi_3 \left(
			\ul{x}, a, \ul{y}, b, \ul{z}, c_{(1)}, \corest{\mu} \left( c_{(2)} \right), c_{(3)}
			\right)
			\\
			    & +
			      (-1)^{\degb{c_{(2)}} \cdot \left( \degb{x} + \degb{a} + \degb{y} + \degb{b} + \degb{c_{(1)}} \right) + 1}
			\phi_3 \left(
			\ul{\corest{\mu} \left( c_{(2)}, x, a_{(1)} \right)}, a_{(2)}, \ul{y}, b, \ul{z}, c_{(1)}
			\right)
			\\
			    & +
			      (-1)^{\left( \degb{x} + \degb{a} \right) \cdot \left( \degb{y} + \degb{b} + \degb{z} + \degb{c} \right)}
			\phi_2 \left( \ul{y}, b, \ul{z}, c, x, a
			\right)
			\\
			    & -
			\phi_2 \left( \ul{x}, a, y, b, \ul{z}, c \right)
			\\
			    & +
			\phi_2 \left( \ul{x}, a, \ul{y}, b, z, c \right)
		\end{aligned}
	\end{equation}
	for $x,y,z \in A$ and $a,b,c \in \tens{A}$.

	We can go down one more level and, given $k \geq 2$ and $r_1, \dots, r_k \geq 0$,
	define the $R$-linear maps
$\phi_k^{r_1,\dots,r_k} \colon
\ul{A} \otimes A^{\otimes r_1} \otimes \dots \otimes \ul{A} \otimes A^{\otimes r_k} \rightharpoonup
R$
	of degree $4 - n - k$, by the composition
	\begin{equation} \label{eq:components-phi_k-r_1-r_k}
		\begin{tikzcd}[column sep=2em]
			{\ul{A} \otimes A^{\otimes r_1} \otimes \dots \otimes \ul{A} \otimes A^{\otimes r_k}} &
			{\ndf{A, \braidop_1}[k][]} & {\ncdf{A, \braidop_1}[k][]} & R.
			\arrow[hookrightarrow, from=1-1, to=1-2]
			\arrow["\phi_k^{r_1,\dots,r_k}", bend left=20, from=1-1, to=1-4, start anchor=north, end anchor=north,
				harpoon, "{[4 - n - k]}"']
			\arrow[two heads, from=1-2, to=1-3]
			\arrow[harpoon, "\phi_k", from=1-3, to=1-4, "{[4 - n - k]}"']
		\end{tikzcd}
	\end{equation}
	Here, the first two horizontal maps are the natural inclusion and the natural projection respectively,
	and we abuse notation and think of $\ul{A}$ both as a $\ZZ$-graded $R$-module which is a copy of $A$ and
	as a bigraded module consisting of a copy of $A$ in line degree one. Since we work with
$\braidop_1$, the $R$-module structure on $\ul{A}$ is the same as the $R$-module structure
	on $A$, and the tensor product $\otimes_1$ is compatible with $\otimes_R$, so
	this works out (see the discussion in \cref{sec:tot-inner-product-components-both-parity-forms}).
	We will call the maps $\phi_k^{r_1,\dots,r_k}$ the \textbf{components} of $\phi_k$.
	Each $\phi_k$ is completely determined by the sequence of components
$\left( \phi_k^{r_1, \dots, r_k} \right)_{r_1, \dots, r_k \geq 0}$, so the components contain
	all the information of $\phi_k$, with some redundancy because of the cyclic symmetry.
	In terms of the components $\phi_{\ul{1}}$ and $\phi_k^{r_1, \dots, r_k}$, we have the following relations:
	\begin{enumerate}
	\item{(Cyclic Symmetry)} For $k \geq 2$ and $r_1, \dots, r_k \geq 0$, we have
	\begin{equation*}
		\phi_k^{r_1,\dots,r_k} \left( \ul{a_1}, l^1, \dots, \ul{a_k}, l^k \right)
		= (-1)^{\varepsilon}
		\phi_k^{r_k,r_1,\dots,r_{k-1}} \left(
		\ul{a_k}, l^k, \ul{a_1}, l^1, \dots, \ul{a_{k-1}}, l^{k-1}
		\right)
	\end{equation*}
	where $a_1, \dots, a_k \in A$, $l^i \in A^{\otimes r_i}$ for $i = 1, \dots, k$ and
	\begin{equation*}
		\begin{aligned}
			\varepsilon & =
			\braid{\left( 1, \degb{a_k} + \degb{l^k} \right)}
			{\left( k - 1, \degb{a_1} + \degb{l^1} + \dots + \degb{a_{k-1}} + \degb{l^{k-1}} \right)}_{1}
			\\
			            & =
			\left( \degb{a_k} + \degb{l^k} \right) \cdot
			\left( \degb{a_1} + \degb{l^1} + \dots + \degb{a_{k-1}} + \degb{l^{k-1}} \right) + (k-1).
		\end{aligned}
	\end{equation*}
	This is a consequence of the fact that the maps $\phi_k$ are defined on $\ncdf{A}[k][]$
	which is a quotient of $\ndf{A}[k][]$ by the image of $\idd - \t$. In particular, note that the map
$\phi_2^{0,0} \colon \ul{A}^{\otimes 2} \rightharpoonup R$ is graded-antisymmetric while
$\phi_3^{0,0,0} \colon \ul{A}^{\otimes 3} \rightharpoonup R$ is cyclically invariant.
	\item{($\phi_{\ul{1}}$ relation)} The relation \eqref{eq:d-phi-ul-1}, coming from
$d_{\mathcal{R}[4-n]} \left( \phi \left( \ul{1} \right) \right) =
\phi \left( D_{\mathcal{A}} \left( \ul{1} \right) \right)$,
	is given explicitly by:
	\begin{equation*}
		\begin{aligned}
			(-1)^{4 - n} d \left( \phi_{\ul{1}} \right)
			 & =
			\sum_{k=2}^{\infty} \frac{(-1)^{\frac{k \left( k + 1 \right)}{2}}}{k}
			\phi_k^{0,\dots,0} \left( \ulz{\mu}, \dots, \ulz{\mu} \right)
			\\
			 & =
			-\frac{1}{2} \phi_2^{0,0} \left( \ulz{\mu}, \ulz{\mu} \right) +
			\frac{1}{3} \phi_3^{0,0,0} \left( \ulz{\mu}, \ulz{\mu}, \ulz{\mu} \right) + \dots
		\end{aligned}
	\end{equation*}
	\item{($\phi_2$ relation)} The relation $d_{\mathcal{R}[4-n]} \circ \phi_2 = \phi_2 \circ \clie{\mu}$
	of \eqref{eq:d-phi-2} is given explicitly by:
	\begin{gather*}
		(-1)^{4-n} d \left( \phi_2^{k,l} \left( \ul{x}, a_1, \dots, a_k,
		\ul{y}, b_1, \dots, b_l \right) \right) =
		\\
		\sum_{\substack{k_1 + k_2 = k, \\ l_1 + l_2 = l}} \pm
		\phi_2^{k_2,l_1} \left(
		\ul{\mu_{l_2 + 1 + k_1} \left( b_{l_1+1}, \dots, b_l, x, a_1, \dots, a_{k_1} \right)},
		a_{k_1+1}, \dots, a_k, \ul{y}, b_1, \dots, b_{l_1} \right) +
		\\
		\sum_{k_1 + k_2 + k_3 = k} \pm
		\phi_2^{k_1 + 1 + k_3,l} \left(
		\ul{x}, a_1, \dots, a_{k_1}, \mu_{k_2} \left( a_{k_1+1}, \dots, a_{k_1 + k_2} \right),
		a_{k_1 + k_2 + 1}, \dots, a_{k}, \ul{y}, b_1, \dots, b_l \right) +
		\\
		\sum_{\substack{k_1 + k_2 = k, \\ l_1 + l_2 = l}} \pm
		\phi_2^{k_1,l_2} \left(
		\ul{x}, a_1, \dots, a_{k_1}, \ul{\mu_{k_2 + 1 + l_1} \left( a_{k_1+1}, \dots, a_k, y, b_1, \dots,
			b_{l_1} \right)}, b_{l_1+1}, \dots, b_{l} \right) +
		\\
		\sum_{l_1 + l_2 + l_3 = l} \pm
		\phi_2^{k, l_1 + 1 + l_3} \left(
		\ul{x}, a_1, \dots, a_k, \ul{y}, b_1, \dots, b_{l_1}, \mu_{l_2} \left( b_{l_1 + 1}, \dots,
		b_{l_1 + l_2} \right), b_{l_1 + l_2 + 1}, \dots, b_l \right),
	\end{gather*}
	where the signs are as in \cref{eq:d-phi-2-explicit}.
	\item{($\phi_3$ relation)} The relation
$d_{\mathcal{R}[4-n]} \circ \phi_3 = -\phi_3 \circ \clie{\mu} + \phi_{2} \circ \qdr$
	of \eqref{eq:d-phi-3} is given explicitly by:
	{ \small
		\begin{gather*}
			(-1)^{4-n} d \left(
			\phi_3^{k,l,m} \left(
			\ul{x}, \underbrace{a_1, \dots, a_k}_{a},
			\ul{y}, \underbrace{b_1, \dots, b_l}_{b},
			\ul{z}, \underbrace{c_1, \dots, c_m}_{c}
			\right)
			\right) =
			\\
			\sum_{\substack{k_1 + k_2 = k \\ m_1 + m_2 = m}} \pm
			\phi_3^{k_2,l,m_1} \left(
			\ul{ \mu_{m_2 + 1 + k_1} \left( c_{m_1+1}, \dots, c_m, x, a_1, \dots, a_{k_1} \right) },
			a_{k_1+1}, \dots, a_k, \ul{y}, b, \ul{z}, c_1, \dots, c_{m_1} \right) +
			\\
			\sum_{k_1 + k_2 + k_3 = k} \pm
			\phi_3^{k_1 + 1 + k_3, l, m} \left(
			\ul{x}, a_1, \dots, a_{k_1}, \mu_{k_2} \left( a_{k_1+1}, \dots, a_{k_1 + k_2} \right), a_{k_1 + k_2 + 1},
			\dots, a_k, \ul{y},b, \ul{z}, c \right) +
			\\
			\sum_{\substack{k_1 + k_2 = k \\ l_1 + l_2 = l}} \pm
			\phi_3^{k_1, l_2, m} \left(
			\ul{x}, a_1, \dots, a_{k_1}, \ul{ \mu_{k_2 + 1 + l_1} \left( a_{k_1+1}, \dots, a_k, y, b_1, \dots,
				b_{l_1} \right)}, b_{l_1+1}, \dots, b_l, \ul{z}, c \right) +
			\\
			\sum_{l_1 + l_2 + l_3 = l} \pm
			\phi_3^{k, l_1 + 1 + l_3, m} \left(
			\ul{x}, a, \ul{y}, b_1, \dots, b_{l_1}, \mu_{l_2} \left( b_{l_1+1}, \dots, b_{l_1+l_2} \right),
			b_{l_1+l_2+1}, \dots, b_l, \ul{z}, c \right) +
			\\
			\sum_{\substack{l_1 + l_2 = l \\ m_1 + m_2 = m}} \pm
			\phi_3^{k,l_1,m_2} \left(
			\ul{x}, a, \ul{y}, b_1, \dots, b_{l_1}, \ul{ \mu_{l_2 + 1 + m_1} \left( b_{l_1+1}, \dots, b_l,
				z, c_1, \dots, c_{m_1} \right)}, c_{m_1+1}, \dots, c_m \right) +
			\\
			\sum_{m_1 + m_2 + m_3 = m} \pm
			\phi_3^{k,l,m_1 + 1 + m_3} \left(
			\ul{x}, a, \ul{y}, b, \ul{z}, c_1, \dots, c_{m_1},
			\mu_{m_2} \left( c_{m_1+1}, \dots, c_{m_1 + m_2} \right), c_{m_1 + m_2 + 1}, \dots, c_m \right) +
			\\
			\pm
			\phi_2^{l, m + 1 + k} \left(
			\ul{y}, b_1, \dots, b_l, \ul{z}, c_1, \dots, c_m, x, a_1, \dots, a_k
			\right) -
			\\
			\phi_2^{k + 1 + l, m} \left(
			\ul{x}, a_1, \dots, a_k, y, b_1, \dots, b_l, \ul{z}, c_1, \dots, c_m
			\right) +
			\\
			\phi_2^{k, l + 1 + m} \left(
			\ul{x}, a_1, \dots, a_k, \ul{y}, b_1, \dots, b_l, z, c_1, \dots, c_m
			\right),
		\end{gather*}
	}
	where the signs are as in \cref{eq:d-phi-3-explicit}.
	\end{enumerate}
	It is clear that if we want, we can write down explicit formulas
	for the relations of \eqref{eq:d-phi-k} involving $\phi_k$ for $k \geq 4$, but we stop here.

	Given $b \in \tc{A}$, the superpotential
	\begin{equation*}
		\SP[b] = \SP[b][\mathcal{A}, \braidop_1] = \phi \left( \G{b}[\geq 2][\mathcal{A}, \braidop_1] \right)
	\end{equation*}
	is given in terms of the components $\phi_{\ul{1}}$ and $\phi_k^{r_1, \dots, r_k}$ by
	\begin{align}
		\SP[b] ={} & \phi_{\ul{1}} +
		\sum_{\substack{k = 2 \\ i_1,\dots,i_{k-1} = 0 \\ j_1,\dots,j_k=0}}^{\infty}
		\frac{(-1)^{\frac{(k - 2) \cdot (k - 1)}{2}}}{1 + \sum_{r=1}^{k-1} i_r + \sum_{r = 1}^k j_r}
		\label{eq:sp-explicit-formula}
		\\
		           & \qquad\qquad\qquad
		\phi_k^{j_1,\dots,j_k}
		\left(
		\ul{ \mu_{i_1} \left( b^{\otimes i_1} \right)} \otimes b^{\otimes j_1} \otimes \dots
		\otimes \ul{ \mu_{i_{k-1}} \left( b^{\otimes i_{k-1}} \right) } \otimes
		b^{\otimes j_{k-1}} \otimes \ul{b} \otimes b^{\otimes j_k}
		\right)
		\notag
		\\
		={}        &
		\phi_{\ul{1}} +
		\sum_{i_1,j_1,j_2=0}^{\infty} \frac{1}{i_1 + j_1 + 1 + j_2} \phi_2^{j_1,j_2} \left(
		\ul{ \mu_{i_1} \left( b^{\otimes i_1} \right)} \otimes b^{\otimes j_1} \otimes
		\ul{b} \otimes b^{\otimes j_2}
		\right) - \cdots
		\notag
	\end{align}

	When $\phi_3 = \phi_4 = \dots = 0$, the data of a total inner product is encoded
	solely in the component $\phi_2 \colon \ncdf{A}[2][] \rightharpoonup R$ of degree $2 - n$
	and in the constant $\phi_{\ul{1}} \in R^{3-n}$.
	Such total inner products are called \textbf{homotopy inner products}, and
	the map $\phi_2$, without the constant term $\phi_{\ul{1}}$, gives us
	a pre-homotopy inner product (\cref{dfn:pre-homotopy-inner-product}).
	When, in addition, $\phi_{\ul{1}} = 0$, the superpotential \eqref{eq:sp-explicit-formula} reduces
	to Cho's potential function \eqref{eq:sp-cho-lee},
	introduced in \cite[Definition 3.1]{cho-homotopy-superpotential}.
	Finally, when in addition $\phi_2$ is strict, we are left with a cyclic structure, and the
	superpotential \eqref{eq:sp-explicit-formula} reduces to the classic superpotential
	\eqref{eq:sp-gt-0-cyclic-structure}.

	Homotopy and pre-homotopy inner products are further discussed in \cref{sec:homotopy-inner-products}.

	\begin{rem} \label{rem:pre-total-vs-total-inner-product}
	Our notion of a total inner product and the superpotential involves the inhomogeneous
	term $\phi_{\ul{1}}$, coming from the fact that $\phi$ is defined on the
	\textit{extended} total complex $\totcompe{\mathcal{A}}[2][]$.
	We discuss briefly the meaning of the inhomogeneous term.

	Given a pre-total inner product $\phi \colon \totcomp{\mathcal{A}}[2][] \rightarrow \mathcal{R}[4-n]$
	on $\mathcal{A}$, one can define the components $\phi_k$ of $\phi$
	as in \eqref{eq:components-phi_k}, replacing the extended total complex $\totcompe{\mathcal{A}}[2][]$
	with the total complex $\totcomp{\mathcal{A}}[2][]$. A pre-total inner product doesn't come
	equipped with a $\phi_{\ul{1}}$ component, and satisfies only the relations of
	\eqref{eq:d-phi-2}--\eqref{eq:d-phi-k}.
	The \textbf{pre-superpotential} function
$\SP[][> 0] \colon \tc{A} \rightarrow R^{3-n}$
	associated to a pre-total inner product $\phi$ is defined by
	\begin{align}
		\SP[b][> 0] \defeq{} &
		\phi \left( \Go{b}[\geq 2][\mathcal{A}] \right) =
		\sum_{\substack{k = 2 \\ i_1, \dots, i_{k-1} = 0 \\ j_1, \dots, j_k = 0}}^{\infty}
		\frac{(-1)^{\frac{(k - 2) \cdot (k - 1)}{2}}}{1 + \sum_{r=1}^{k-1} i_r + \sum_{r = 1}^k j_r}
		\label{eq:pre-sp-phi-braidop-1}
		\\
		                     & \qquad
		\phi_k^{j_1,\dots,j_k}
		\left(
		\ul{ \mu_{i_1} \left( b^{\otimes i_1} \right)} \otimes b^{\otimes j_1} \otimes \dots
		\otimes \ul{ \mu_{i_{k-1}} \left( b^{\otimes i_{k-1}} \right) } \otimes
		b^{\otimes j_{k-1}} \otimes \ul{b} \otimes b^{\otimes j_k}
		\right).
		\notag
	\end{align}

	In general, the pre-superpotential of a bounding cochain is not a cocycle.
	There are two cases in which we can construct a gauge-invariant cocycle
	from a bounding cochain using the pre-superpotential:
	\begin{enumerate}
	\item Assume we restrict our attention to $\Ainf$-algebras which are not curved,
	and work with $\Ainf$-morphisms without a change of connection element.
	In this case, \cref{lm:Go-and-Ho-geq2-properties-general} shows
	that the chain $\Go{b}[\geq 2]$ of the total complex satisfies the same
	properties as the cyclic Chern--Simons form $\G{b}[\geq 2]$ of the extended total complex.
	Repeating the arguments of \cref{sec:generalized-inner-product-superpotential}
	with pre-total inner products and the pre-superpotential replacing
	total inner products and the superpotential, we can deduce that $\SP[][>0]$
	gives us a gauge-invariant cohomology class, when restricted to bounding cochains.
	See also \cref{rem:Go-geq-2-closed-no-curvature}.
	\item Assume we allow curvature but work with
	pre-total inner products which satisfy the extra condition
	\begin{equation} \label{eq:phi-vanishes-mu-0-lift}
		\sum_{k=2}^{\infty} \frac{(-1)^{\frac{k \left( k + 1 \right)}{2}}}{k}
		\phi_k \left( \ulz{\mu}^{\otimes k} \right) = 0.
	\end{equation}
	Such pre-total inner products can be upgraded to total inner products
	by taking $\phi_{\ul{1}} = 0$, in which case, the identity \eqref{eq:d-phi-ul-1}
	reduces to \eqref{eq:phi-vanishes-mu-0-lift}. For total inner products with
$\phi_{\ul{1}} = 0$, the superpotential $\SP$ reduces to the
	pre-superpotential $\SP[][>0]$. Morphisms of total inner products
	without a change of connection term respect the condition $\phi_{\ul{1}} = 0$,
	so if we restrict our attention to such morphisms, the results of
	\cref{sec:generalized-inner-product-superpotential}
	show that we also get a gauge-invariant cohomology class from $\SP[][>0]$.
	\end{enumerate}
	However, in the general case, only the \textit{sum}
$\SP[] = \phi_{\ul{1}} + \SP[][> 0]$
	gives a gauge-invariant cohomology class on bounding cochains.

	This phenomenon was observed in the context of the Fukaya $\Ainf$-algebra associated
	to a compact Lagrangian. The Fukaya $\Ainf$-algebra is curved and comes equipped
	with a geometrically meaningful inhomogeneous term, often denoted by $\mathfrak{m}_{-1}$,
	which comes from counting disks with zero marked boundary points.
	The term $\mathfrak{m}_{-1}$ is related to the cyclic structure via an identity of
	the form \eqref{eq:d-m-minus-1-identity}, of which \eqref{eq:d-phi-ul-1} is a generalization,
	and it is added to the classical superpotential \eqref{eq:sp-gt-0-cyclic-structure}, i.e.,
	our pre-superpotential, to guarantee that it remains invariant (see \cite{Fukaya2011,Joyce2008,Solomon2016a}).
	\end{rem}

	\begin{rem}
	We have defined the components of total inner product
	\begin{equation*}
		\phi^{\braidop_1} \colon \totcompe{\mathcal{A}, \braidop_1}[2][] \rightarrow \mathcal{R}[4-n],
	\end{equation*}
	working with the inner product parity form $\braidop_1$. The components
	of a total inner product
$\phi^{\braidop_2} \colon \totcompe{\mathcal{A}, \braidop_2}[2][] \rightarrow \mathcal{R}[4-n]$
	can be defined by passing to the corresponding
$\phi^{\braidop_1}$, as described in \cref{sec:description-using-braid-op-1},
or directly, and they involve extra signs. See \cref{sec:tot-inner-product-components-both-parity-forms}.
\end{rem}

\section{Pre-Homotopy and Homotopy Inner Products} \label{sec:homotopy-inner-products}

Let $\mathbbm{k}$ be a field of characteristic zero, endowed with the trivial norm,
let $\mathcal{R} = (R,d)$ be a differential graded-commutative Banach $\mathbbm{k}$-algebra,
and let $\mathcal{A} = \left( A, \mu \right)$ be a Banach $\Ainf$-algebra over $\mathcal{R}$.

In this section, we study \textit{homotopy inner products}, which are total inner
products $\phi$ with components $\phi_k = 0$ for $k \geq 3$.
Homotopy inner products are determined by a constant component $\phi_{\ul{1}}$ and the
component
$\phi_2 \colon \ncdf{\mathcal{A}}[2][] / \Im \left( \qdr^3 \right) \rightharpoonup \mathcal{R}$,
which is a pre-homotopy inner product, related to $\phi_{\ul{1}}$ via an extra identity.
(Pre)-homotopy inner products and their associated (pre)-superpotential functions
are discussed in \Crefrange{sec:pre-homotopy-inner-products}{subsec:homotopy-inner-products}.
Focusing on pre-homotopy inner products, in \cref{sec:rel-pre-homotopy-pre-trace}
we use the homotopy equivalence between Connes' cyclic complex
$\ncdfr{\mathcal{A}}[0][]$ and the complex
$\ncdf{\mathcal{A}}[2][] / \Im \left( \qdr^3 \right)$ to provide explicit formulas
translating pre-$\infty$-traces to pre-homotopy inner products and vice versa.
We show that the formulas respect unitality conditions for both notions.
When $\mathcal{A}$ corresponds to a differential graded algebra,
this recovers the correspondence between traces and cyclic structures on $\mathcal{A}$.
We also study the conditions
under which a pre-$\infty$-trace (resp.\ pre-homotopy inner product) corresponds
to a homotopy inner product (resp.\ $\infty$-trace) with a vanishing constant component.

In \cref{sec:recovering-chern-simons}, we show
how to recover the classical Chern--Simons $3$-form and action from
the cyclic Chern--Simons form and the superpotential, working in
the differential graded algebra of $\End{E}$-valued differential forms on
a closed oriented manifold $M$, where $E$ is a flat vector bundle over $M$.

In \cref{sec:derivative-and-trace}, we consider the derivative of the superpotential
associated to a strongly unital homotopy inner product $\phi$. We show that the
derivative of the superpotential on weak bounding cochains coincides
with a multiple of the derivative of the pre-$\infty$-modulus function, constructed using
the pre-$\infty$-trace associated to $\phi_2$.

Finally, in \cref{sec:sp-and-periodicity}, we establish a direct relation
between the pre-superpotential function,
the pre-$\infty$-trace, and the periodicity
operator on cyclic homology.

In this section, when writing tensor products, we sometimes omit the tensor symbol
for brevity, and write expressions such as $a_1 \otimes \dots \otimes a_n$
in the form $a_1 \, \dots \, a_n$, with the understanding that concatenation
represents the tensor product. To facilitate comparison with the literature, we will
work in this section with the inner product parity form $\braidop_1$.

\subsection{Pre-Homotopy Inner Products} \label{sec:pre-homotopy-inner-products}

Fix $n \in \ZZ$. Recall from \cref{dfn:pre-homotopy-inner-product} that an $n$-dimensional
pre-homotopy inner product is a morphism
$\phi_{2} \colon \ncdf{\mathcal{A}}[2][] / \Im \left( \qdr^3 \right) \rightarrow
	\mathcal{R}[2-n]$
of differential graded Banach $\mathcal{R}$-modules.
When convenient, we will think equivalently of $\phi_2$ as
\begin{enumerate}
	\item A graded contractive $R$-linear \textit{chain} map $\ncdf{\mathcal{A}}[2][] / \Im \left( \qdr^3 \right) \rightharpoonup \mathcal{R}$
	      of degree $2 - n$, obtained by composing $\phi_2$ with the canonical map $\mathcal{R}[2-n] \rightharpoonup \mathcal{R}$.
	\item A graded contractive $R$-linear \textit{chain} map $\ncdf{\mathcal{A}}[2][] \rightharpoonup \mathcal{R}$ of degree $2 - n$
	      which vanishes on the image of $\qdr^3$, obtained by composing the map from the previous item with the projection
	      $\ncdf{\mathcal{A}}[2][] \twoheadrightarrow \ncdf{\mathcal{A}}[2][] / \Im \left( \qdr^3 \right)$.
\end{enumerate}
In what follows, the dimension $n$ of $\phi_2$ will be fixed, and
we often omit it for brevity.

Given a pre-homotopy inner product $\phi_2$, and $k,l \geq 0$, consider the $R$-linear maps
\begin{equation*}
	\phi_2^{k,l} \colon \ul{A} \otimes A^{\otimes k} \otimes \ul{A} \otimes A^{\otimes l} \rightharpoonup R
\end{equation*}
of degree $2 - n$, called the \textbf{components} of $\phi_2$, defined by the composition
\begin{equation*}
	\adjustbox{scale=0.85,center}
	{
		\begin{tikzcd}
			{\ul{A} \otimes A^{\otimes k} \otimes \ul{A} \otimes A^{\otimes l}} &
			{\ndf{A}[2][]} & {\ncdf{A}[2][]} & {\ncdf{\mathcal{A}}[2][] / \Im \left( \qdr^3 \right)} & R[2-n] & R.
			\arrow[hookrightarrow, from=1-1, to=1-2]
			\arrow["\phi_{2}^{k,l}", bend left=10, from=1-1, to=1-6, start anchor=north, end anchor=north,
				harpoon, "{[2 - n]}"']
			\arrow[two heads, from=1-2, to=1-3]
			\arrow[two heads, from=1-3, to=1-4]
			\arrow["\phi_2", from=1-4, to=1-5]
			\arrow[harpoon, from=1-5, to=1-6, "{[2 - n]}"']
		\end{tikzcd}
	}
\end{equation*}
In terms of components,
\cref{dfn:pre-homotopy-inner-product} implies the following relations:
\begin{enumerate}
	\item{(Antisymmetry)} For $k,l \geq 0$ we have
	      \begin{equation*}
		      \phi_2^{k,l} \left( \ul{x}, r, \ul{y}, s \right)
		      = (-1)^{\left( \degb{x} + \degb{r} \right) \left( \degb{y} + \degb{s} \right) + 1}
		      \phi_2^{l,k} \left( \ul{y}, s, \ul{x}, r \right)
	      \end{equation*}
	      where $x, y \in A, r \in A^{\otimes k}$ and $s \in A^{\otimes l}$.
	\item{($\qdr$-closed)} For $k,l,m \geq 0$ we have
	      \begin{equation*}
		      \begin{aligned}
			      0 ={} &
			      (-1)^{\left( \degb{y} + \degb{b} + \degb{z} + \degb{c} \right) \left( \degb{x} + \degb{a} \right)}
			      \phi_2^{l, m + 1 + k} \left(
			      \ul{y}, b, \ul{z}, c, x, a
			      \right)
			      \\
			            & -
			      \phi_2^{k + 1 + l, m} \left(
			      \ul{x}, a, y, b, \ul{z}, c
			      \right)
			      \\
			            & +
			      \phi_2^{k, l + 1 + m} \left(
			      \ul{x}, a, \ul{y}, b, z, c
			      \right).
		      \end{aligned}
	      \end{equation*}
	      where $x,y,z \in A$ and $a \in A^{\otimes k}, b \in A^{\otimes l}, c \in A^{\otimes m}$.
	\item{($\clie{\mu}$-closed)} For $k,l \geq 0$, we have
	      \begin{gather*}
		      (-1)^{2-n} d \left( \phi_2^{k,l} \left( \ul{x}, a_1, \dots, a_k,
		      \ul{y}, b_1, \dots, b_l \right) \right) =
		      \\
		      \sum_{\substack{k_1 + k_2 = k, \\ l_1 + l_2 = l}} \pm
		      \phi_2^{k_2,l_1} \left(
		      \ul{\mu_{l_2 + 1 + k_1} \left( b_{l_1+1}, \dots, b_l, x, a_1, \dots, a_{k_1} \right)},
		      a_{k_1+1}, \dots, a_k, \ul{y}, b_1, \dots, b_{l_1} \right) +
		      \\
		      \sum_{k_1 + k_2 + k_3 = k} \pm
		      \phi_2^{k_1 + 1 + k_3,l} \left(
		      \ul{x}, a_1, \dots, a_{k_1}, \mu_{k_2} \left( a_{k_1+1}, \dots, a_{k_1 + k_2} \right),
		      a_{k_1 + k_2 + 1}, \dots, a_{k}, \ul{y}, b_1, \dots, b_l \right) +
		      \\
		      \sum_{\substack{k_1 + k_2 = k, \\ l_1 + l_2 = l}} \pm
		      \phi_2^{k_1,l_2} \left(
		      \ul{x}, a_1, \dots, a_{k_1}, \ul{\mu_{k_2 + 1 + l_1} \left( a_{k_1+1}, \dots, a_k, y, b_1, \dots,
			      b_{l_1} \right)}, b_{l_1+1}, \dots, b_{l} \right) +
		      \\
		      \sum_{l_1 + l_2 + l_3 = l} \pm
		      \phi_2^{k, l_1 + 1 + l_3} \left(
		      \ul{x}, a_1, \dots, a_k, \ul{y}, b_1, \dots, b_{l_1}, \mu_{l_2} \left( b_{l_1 + 1}, \dots,
		      b_{l_1 + l_2} \right), b_{l_1 + l_2 + 1}, \dots, b_l \right)
	      \end{gather*}
	      where $x,y \in A$ and $a_1, \dots, a_k, b_1, \dots, b_l \in A$. The signs are as in \cref{eq:d-phi-2-explicit}.
\end{enumerate}
A pre-homotopy inner product is called \textbf{strict} if $\phi_2^{i,j} = 0$ when $i > 0$ or $j > 0$.
As discussed in \cref{subsec:generalized-superpotential-intro} and proven in \cref{appendix:cyclic-structures},
a strict pre-homotopy inner product is the same thing as a cyclic structure on $\mathcal{A}$.

Next, we discuss unitality conditions for pre-homotopy inner products.

\begin{dfn} \label{dfn:pre-homotopy-inner-product-unitality-conditions}
	Let $\mathcal{A} = \left( A, \mu, e \right)$ be a unital Banach $\Ainf$-algebra over $\mathcal{R}$
	and let $\phi_2$ be a pre-homotopy inner product on $\mathcal{A}$.
	\begin{enumerate}
		\item The pre-homotopy inner product $\phi_2$ is called $\degen$\textbf{-unital} if
		      \begin{equation}
			      \phi_2^{k,l} \left( \ul{x}, a_1, \dots, a_k, \ul{y}, b_1, \dots, b_l \right) = 0
		      \end{equation}
		      whenever
		      $a_i = e$ for some $1 \leq i \leq k$ or $b_j = e$ for some $1 \leq j \leq l$.
		      Equivalently, $\phi_2$ is $\degen$-unital if $\phi_2$ descends to the quotient complex
		      $\ncdf{\mathcal{A}}[2][] / \left< \Im \left( \qdr^3 \right), \degen[A][2][] \right>$.
		\item The pre-homotopy inner product $\phi_2$ is called $\degenu$\textbf{-unital} if
		      \begin{equation}
			      \phi_2 \circ \rest{\ccont{e}^2}{\ncdf{A}[0][]} = 0.
		      \end{equation}
		      Equivalently, $\phi_2$ is $\degenu$-unital if $\phi_2$ descends to the quotient complex
		      \begin{equation*}
			      \ncdf{\mathcal{A}}[2][] / \left< \Im \left( \qdr^3 \right), \degenu[A][2][] \right>.
		      \end{equation*}
		\item The pre-homotopy inner product $\phi_2$ is called \textbf{strongly unital} if
		      $\phi_2$ is both $\degen$-unital and $\degenu$-unital.
		      Equivalently, $\phi_2$ is strongly unital if $\phi_2$ descends to the quotient complex
		      \begin{equation*}
			      \ncdf{\mathcal{A}}[2][] / \left< \Im \left( \qdr^3 \right), \degen[A][2][], \degenu[A][2][] \right>.
		      \end{equation*}
	\end{enumerate}
\end{dfn}

We can describe the condition of being strongly unital more explicitly as follows.
Consider the elements of $\ncdf{\mathcal{A}}[2][]$ which have the form
\begin{equation} \label{eq:elem-two-underlined-units}
	\ul{e} \, a_1 \, \dots \, a_k \, \ul{e} \, b_1 \, \dots \, b_l,
\end{equation}
with $k,l \geq 0$ and $a_i, b_j \in A$. Those are the elementary tensors which contain two underlined units and
we denote them succinctly by $\ul{e} {*} \ul{e} {*}$.
Let us also use the notation $*e*$ to denote
elementary tensors of $\ncdf{\mathcal{A}}[2][]$ in which the unit $e$ of $\mathcal{A}$ appears at least once without an underline,
i.e., elements of $\degen[A][2]$.

\begin{lm} \label{lm:desc-ncdf-2-strongly-unital-subcomplex}
	We have
	\begin{equation} \label{eq:strongly-unital-subcomplex-diff-forms}
		\left< \Im \left( \qdr^3 \right), \degen[A][2][], \degenu[A][2][] \right> =
		\left< \Im \left( \qdr^3 \right), *e*, \ul{e} {*} \ul{e} {*} \right> =
		\left< \Im \left( \qdr^3 \right), *e*, \ul{e} \, \ul{e} \right>.
	\end{equation}
\end{lm}
\begin{proof}
	Recall that $\degenu[A][2][]$ is generated as a graded Banach $R$-module by elements of
	$\Im{\rest{\ccont{e}^2}{\ncdf{A}[0][]}}$. All such elements are of the form $\ul{e} {*} \ul{e} {*}$.
	Note also that the complex $\left< \Im \left( \qdr^3 \right), \degen[A][2][] \right>$ already
	contains all elements of the form \eqref{eq:elem-two-underlined-units}
	whenever $k > 0$ or $l > 0$. In particular, this implies that $\left< \Im \left( \qdr^3 \right), \degen[A][2][] \right>$
	contains all elements of the form $\ccont{e}^2 \left( x \right)$ for $x \in \ncdfr{A}[0][]$, but
	it does not contain $\ccont{e}^2 \left( 1 \right) = \ul{e} \, \ul{e}$. Hence, we have
	\cref{eq:strongly-unital-subcomplex-diff-forms}.
\end{proof}

\begin{cor}
	Let $\mathcal{A} = \left( A, \mu, e \right)$ be a unital Banach $\Ainf$-algebra over $\mathcal{R}$
	and let $\phi_2$ be a pre-homotopy inner product on $\mathcal{A}$. Then $\phi_2$ is strongly unital if
	and only if it satisfies:
	\begin{enumerate}
		\item $\phi_2^{k,l} \left( \ul{x}, a_1, \dots, a_k, \ul{y}, b_1, \dots, b_l \right) = 0$
		      whenever $a_i = e$ or $b_j = e$.
		\item $\phi_2^{k,l} \left( \ul{e}, a_1, \dots, a_k, \ul{e}, b_1, \dots, b_l \right) = 0$
		      whenever $k,l \geq 0$.
	\end{enumerate}
	\qed
\end{cor}

\begin{lm} \label{lm:ud-unital-plus-cond-implies-strongly-unital}
	Let $\mathcal{A} = \left( A, \mu, e \right)$ be a unital Banach $\Ainf$-algebra over $\mathcal{R}$
	and let $\phi_2$ be a $\degen$-unital pre-homotopy inner product on $\mathcal{A}$.
	If $\phi_2 \left( \ul{e}, \ul{e} \right) = 0$ then $\phi_2$ is also $\degenu$-unital, hence, strongly unital.
	\qed
\end{lm}

\subsection{Homotopy Inner Products and the Superpotential} \label{subsec:homotopy-inner-products}
\begin{dfn} \label{dfn:homotopy-inner-product}
	An $n$-dimensional total inner product $\phi \colon \totcompe{\mathcal{A}}[2][] \rightarrow \mathcal{R}[4-n]$
	for which $\phi_k = 0$ for $k \geq 3$ will be called an $n$\textbf{-dimensional homotopy inner product}.
\end{dfn}

More explicitly, by \cref{sec:total-inner-product-explicit-relations}, the data of a homotopy inner product $\phi$
is equivalently described by the pair $\left( \phi_{\ul{1}}, \phi_2 \right)$ of its non-zero components,
where $\phi_{\ul{1}} \in R^{3-n}$ with $\nnorm[\phi_{\ul{1}}] \leq 1$, and
$\phi_2 \colon \ncdf{A}[2][] \rightharpoonup R$ is an
$R$-linear contractive map of degree $2 - n$.
The pair $\left( \phi_{\ul{1}}, \phi_2 \right)$ is required to satisfy
the identities \eqref{eq:d-phi-ul-1} to \eqref{eq:d-phi-3}, which take the form:
\begin{align}
	(-1)^{2-n} d \left( \phi_{\ul{1}} \right) & = -\frac{1}{2} \phi_2 \left( \ulz{\mu}, \ulz{\mu} \right),
	\label{eq:d-phi-ul-1-id}
	\\
	(-1)^{2-n} d \circ \phi_2                 & = \phi_2 \circ \clie{\mu}, \label{eq:d-phi-2-id}
	\\
	\phi_2 \circ \qdr^3                       & = 0. \label{eq:phi-2-q-id}
\end{align}
Identity \eqref{eq:phi-2-q-id} implies that $\phi_2$ descends to a map
$\phi_{2} \colon \ncdf{A}[2][] / \Im \left( \qdr^3 \right) \rightharpoonup R$,
denoted by the same name, which, by \cref{eq:d-phi-2-id}, is a chain map.
Hence, equivalently, the data of a homotopy inner product is given by:
\begin{enumerate}
	\item A pre-homotopy inner product $\phi_{2}$ in the sense of \cref{dfn:pre-homotopy-inner-product}, called
	      the \textbf{associated pre-homotopy inner product}.
	\item An element $\phi_{\ul{1}} \in R^{3-n}$ which satisfies $\nnorm[\phi_{\ul{1}}] \leq 1$.
\end{enumerate}
The map $\phi_2$ and the element $\phi_{\ul{1}}$ are required to be related by the identity \eqref{eq:d-phi-ul-1-id}.

Thus, a homotopy inner product $\phi$ is given by the same data as a
pre-homotopy inner product $\phi_2$,
except that it comes with an additional ``$0$-th component'' $\phi_{\ul{1}}$,
and we have the extra identity \eqref{eq:d-phi-ul-1-id}.
Note that even when $\phi_{\ul{1}} = 0$, being a homotopy inner product imposes the
extra condition
\begin{equation} \label{eq:cond-pre-homotopy-upgrade-homotopy-inner-product}
	\phi_2 \left( \ulz{\mu}, \ulz{\mu} \right) = 0,
\end{equation}
which a pre-homotopy inner product is not required a priori to satisfy.
Any pre-homotopy inner product $\phi_2$ which satisfies \cref{eq:cond-pre-homotopy-upgrade-homotopy-inner-product}
can be upgraded to a homotopy inner product by setting $\phi_{\ul{1}} = 0$.

\begin{rem} \label{rem:ulz-otimes-2-mod-q}
	When thinking of $\phi_2$ as defined on $\ncdf{A}[2][] / \Im \left( \qdr^3 \right)$,
	the element $-\frac{1}{2} \ulz{\mu}^{2}$ appearing in \cref{eq:d-phi-ul-1-id} is also considered
	as an element of the quotient $\ncdf{A}[2][] / \Im \left( \qdr^3 \right)$.
	Note that we have
	\begin{equation*}
		\begin{aligned}
			\clie{\mu} \left( \frac{1}{2} \ulz{\mu}^{2} \right) & =
			\frac{1}{2} \left(
			- \ulz{\mu} \, \mu_0 \left( 1 \right) \, \ulz{\mu} + \ulz{\mu} \, \ulz{\mu} \, \mu_0 \left( 1 \right)
			\right)
			\\
			                                                    & =
			\ulz{\mu} \, \ulz{\mu} \, \mu_0 \left( 1 \right) = \frac{1}{3} \qdr^3 \left( \ulz{\mu}^{3} \right),
		\end{aligned}
	\end{equation*}
	so $\ulz{\mu}^{2}$ is in fact a \textit{closed} element of $\ncdf{A}[2][] / \Im \left( \qdr^3 \right)$, which is
	consistent with \cref{eq:d-phi-ul-1-id}.
\end{rem}

Since a homotopy inner product is a particular instance of a total inner product, all the definitions and results from
\cref{sec:generalized-inner-product-superpotential} apply to homotopy inner products. When $\mathcal{A}$ is unital,
the homotopy inner product $\phi$ is $\degenu$-unital (resp.\ strongly unital) in the sense of \cref{dfn:total-inner-product-unital},
if and only if the associated pre-homotopy inner product $\phi_2$ is $\degenu$-unital (resp.\ strongly unital)
in the sense of \cref{dfn:pre-homotopy-inner-product-unitality-conditions}.

Given a topologically nilpotent element $b \in \tc{A}$, the superpotential
function $\SP \colon \tc{A} \rightarrow R^{3-n}$
from \cref{sec:generalized-inner-product-superpotential} associated to $\phi$,
given by \eqref{eq:sp-explicit-formula}, takes the form
\begin{equation} \label{eq:sp-explicit-formula-homotopy-inner-product-braidop-1}
	\SP[b] =
	\phi_{\ul{1}} +
	\sum_{i,j,k=0}^{\infty} \frac{1}{i + j + 1 + k} \phi_2^{j,k} \left(
	\ul{ \mu_{i} \left( b^{i} \right)}, b^{j}, \ul{b}, b^{k}
	\right).
\end{equation}
By \cref{lm:superpotential-bounding-chain-closed}, $\SP[b]$ is a closed element of $R$ when
$b$ is a strong bounding cochain, as well as when $b$ is a weak bounding cochain provided that
the Banach $\Ainf$-algebra $\mathcal{A}$ is unital and $\phi_2$ is
$\degenu$-unital. In both cases, \cref{thm:invariance-superpotential} shows that the cohomology
class of $\SP[b]$ is invariant under gauge equivalence.

The value $\SP[b]$ of the superpotential is the sum of two terms: A constant term $\phi_{\ul{1}}$
and a term depending on $b$ and $\phi_2$, which comes from applying $\phi_2$ to
the cyclic codifferential form $\G{b}[2] \in \ncdf{A}[2][1]$ given by
\begin{equation} \label{eq:G-b-2-braidop-1}
	\G{b}[2] = \sum_{i,j,k=0}^{\infty} \frac{1}{i + j + 1 + k}
	\ul{ \mu_{i} \left( b^{i} \right)} \, b^{j} \, \ul{b} \, b^{k}.
\end{equation}
The form $\G{b}[2]$ is precisely the line degree two component of the
cyclic Chern--Simons form \eqref{eq:G2-braid-op-1-explicit-formula}.
Given a pre-homotopy inner product $\phi_2$ on $\mathcal{A}$, we define the
\textbf{pre-superpotential} function
$\SP[][> 0] \colon \tc{A} \rightarrow R^{3-n}$ associated to $\phi_2$ by
\begin{equation} \label{eq:pre-sp-phi-2-braidop-1}
	\SP[b][>0] \defeq \phi_2 \left( \G{b}[2] \right) =
	\sum_{i,j,k=0}^{\infty} \frac{1}{i + j + 1 + k} \phi_2^{j,k} \left(
	\ul{ \mu_{i} \left( b^{i} \right)}, b^{j}, \ul{b}, b^{k}
	\right).
\end{equation}

\begin{rem}
	The pre-superpotential \eqref{eq:pre-sp-phi-2-braidop-1} is precisely
	the pre-superpotential \eqref{eq:pre-sp-phi-braidop-1} of
	\cref{rem:pre-total-vs-total-inner-product}, thinking of $\phi_2$ as a pre-total
	inner product whose only non-zero component is $\phi_2$. As discussed
	in \cref{rem:pre-total-vs-total-inner-product}, in general, only
	the sum $\SP = \phi_{\ul{1}} + \SP[][>0]$ gives a gauge-invariant cohomology class
	on bounding cochains. However, if we restrict our attention to
	pre-homotopy inner products which satisfy \cref{eq:cond-pre-homotopy-upgrade-homotopy-inner-product},
	and also to $\Ainf$-morphisms without a change of connection term,
	the pre-superpotential $\SP[][>0]$ alone also gives us a gauge-invariant
	cohomology class on bounding cochains.
	This makes sense for example if we work with $\Ainf$-algebras with no curvature.
\end{rem}

\begin{rem}
	The de Rham model of the Fukaya $\Ainf$-algebra for a Lagrangian submanifold
	constructed in \cite{Solomon2016} comes with a cyclic structure
	and an $\mathfrak{m}_{-1}$ term, which give rise to a homotopy inner product with
	a caveat. The associated superpotential is then the classical superpotential
	of \cite{Solomon2016a}. See \cref{appendix:sign-conversions-jake}.
\end{rem}

\subsection{The Relation Between Pre-Homotopy Inner Products and Pre-\texorpdfstring{$\infty$}{Infinity}-Traces} \label{sec:rel-pre-homotopy-pre-trace}

Recall from \cref{sec:pre-infinity-traces} that an $n$-dimensional pre-$\infty$-trace
on $\mathcal{A}$
is a chain map
\begin{equation*}
	\theta \colon \ncdfr{\mathcal{A}}[0][] \rightarrow \mathcal{R}[1-n],
\end{equation*}
while an $n$-dimensional pre-homotopy inner product is a chain map
\begin{equation*}
	\phi_2 \colon \ncdf{\mathcal{A}}[2][] / \Im \left( \qdr^3 \right) \rightarrow \mathcal{R}[2-n].
\end{equation*}
In \cref{lm:ncdf-0-ncdf-2-mod-q3-equiv}, we have shown that when $\mathcal{A}$ is unital, the complexes
$\ncdf{\mathcal{A}}[2][] / \Im \left( \qdr^3 \right)$ and ${\ncdfr{\mathcal{A}}[0][]}[1]$ are homotopy equivalent and
constructed explicit chain maps $\psi, \psi'$ inducing the equivalence. We can use the maps $\psi,\psi'$
to convert pre-$\infty$-traces to pre-homotopy inner products and vice versa.
We give explicit formulas converting $\theta$ to $\phi_2$ and vice versa,
show that they respect appropriate unitality conditions, and study their
behaviour when applied to strict pre-$\infty$-traces and pre-homotopy inner products,
i.e., traces and cyclic structures.

As described in \cref{sec:traces-modulus-invariance}, an $n$-dimensional
$\infty$-trace is given by an $n$-dimensional pre-$\infty$-trace
$\theta$, together with an extra component $\theta_0 \in R^{1-n}$,
related to $\theta$ via the identity \eqref{eq:infty-trace-0-rel}.
Similarly, an $n$-dimensional homotopy
inner product is given by an $n$-dimensional pre-homotopy inner
product $\phi_2$ together with an extra component $\phi_{\ul{1}} \in R^{3-n}$,
related to $\phi_2$ via the identity \eqref{eq:d-phi-ul-1-id}.

The correspondence between pre-$\infty$-traces and pre-homotopy inner products
does not extend to a correspondence between
$\infty$-traces and homotopy inner products, as we don't have a way to relate the
extra elements $\theta_0$ and $\phi_{\ul{1}}$ of $R$, but we study the conditions
under which a pre-$\infty$-trace (resp.\ pre-homotopy inner product) corresponds
to a homotopy inner product (resp.\ $\infty$-trace) with a vanishing extra component.

In what follows, we use freely the notation and maps from \cref{sec:cyclic-homology-models}.
We also slightly abuse notation and often do not differentiate between
elements or submodules of $\ndfr{A}[0][]$ (resp.\ $\ncdf{\mathcal{A}}[2][]$) and their images under the projection
$\pi_0 \colon \ndfr{A}[0][] \twoheadrightarrow \ncdfr{A}[0][]$ (resp.\
$\pi_2 \colon \ncdf{\mathcal{A}}[2][] \twoheadrightarrow \ncdf{\mathcal{A}}[2][] / \Im \left( \qdr^3 \right)$),
relying on context to understand whether we work with elements or equivalence classes.

\subsubsection{Pre-Homotopy Inner Products to Pre-\texorpdfstring{$\infty$}{Infinity}-Traces}
In what follows we assume that $\mathcal{A}$ is unital with unit $e$.
Given a pre-homotopy inner product $\phi_2$ on $\mathcal{A}$, we can use the map
$\psi'$, given by \cref{eq:psi'-def} of \cref{lm:ncdf-0-ncdf-2-mod-q3-equiv}, to
define a pre-$\infty$-trace $\theta$,
called the \textbf{corresponding pre}-$\infty$-\textbf{trace}, by
\begin{equation} \label{eq:naive-trace-from-hip}
	\theta \defeq \phi_2 \circ \psi'.\footnote{Here, we think of $\theta, \phi_2, \psi'$
		as graded $R$-linear maps with codomain $R$ of degrees $1-n,2-n,-1$, respectively.}
\end{equation}

To make the relation between $\theta$ and $\phi_2$ more explicit, let us write down several convenient formulas
for the map $\psi'$.

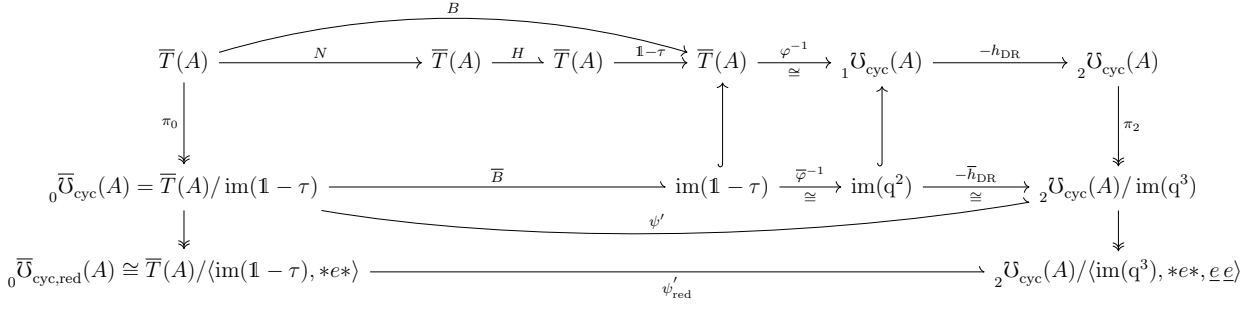
\begin{figure}[htb]
	\centering
	\adjustbox{scale=0.7,center}
	{
		\begin{tikzcd}
			{\tensr{A}} & {\tensr{A}} & {\tensr{A}} & {\tensr{A}} & {\ncdf{A}[1][]} & {\ncdf{A}[2][]}
			\\
			\\
			{{\ncdfr{A}[0][] = \tensr{A} / \Im \left( \idd - \t \right)}} &&& {\Im \left( \idd - \t \right)} &
			{\Im \left( \qdr^2 \right)} & {\ncdf{A}[2][] / \Im \left( \qdr^3 \right)}
			\\
			{\ncdfrred{A}[0][] \cong \tensr{A} / \left< \Im \left( \idd - \t \right), *e* \right>}
			&&&&&
			{\ncdf{A}[2][] / \left< \Im \left( \qdr^3 \right), *e*, \ul{e} \, \ul{e} \right>}
			\arrow["B", from=1-1, to=1-4, curve={height=-30pt}]
			\arrow["\N", from=1-1, to=1-2]
			\arrow["{\pi_0}"', two heads, from=1-1, to=3-1]
			\arrow["H", harpoon, from=1-2, to=1-3]
			\arrow["{\idd - \t}", from=1-3, to=1-4]
			\arrow["{\varphi^{-1}}", "\cong"', from=1-4, to=1-5]
			\arrow["-h_{\dr}", from=1-5, to=1-6]
			\arrow["\pi_2", two heads, from=1-6, to=3-6]
			\arrow["{\overline{B}}", harpoon, from=3-1, to=3-4]
			\arrow["{\psi'}", harpoon, curve={height=30pt}, from=3-1, to=3-6]
			\arrow[two heads, from=3-1, to=4-1]
			\arrow[hook, from=3-4, to=1-4]
			\arrow["{\overline{\varphi}^{-1}}", "\cong"', from=3-4, to=3-5]
			\arrow[hook, from=3-5, to=1-5]
			\arrow["{-\overline{h}_{\dr}}", "\cong"', from=3-5, to=3-6]
			\arrow[two heads, from=3-6, to=4-6]
			\arrow["{\psi_{\textrm{red}}^{'}}"', harpoon, from=4-1, to=4-6]
		\end{tikzcd}
	}
	\caption{Maps involved in the definition of $\psi'$.}
	\label{fig:psi'-definition-and-reduced}
\end{figure}

\begin{lm} \label{lm:psi-tag-explicit-formula}
	The map $\psi' \colon \ncdfr{\mathcal{A}}[0][] \rightharpoonup \ncdf{\mathcal{A}}[2][] / \Im \left( \qdr^3 \right)$
	has the explicit form
	\begin{equation*}
		\psi' \left( a_1 \, \dots \, a_m \right) =
		\sum_{i=1}^{m} \sum_{k=0}^{\infty}
		(-1)^{\degb{a_1} + \dots + \degb{a_{i}} + k + 1}
		a_1 \, \dots \, \ul{a_{i}} \, \ul{e} \, \left( \mu_0 \left( 1 \right) \, e \right)^k \, a_{i+1} \,
		\dots \, a_m.
	\end{equation*}
\end{lm}
\begin{proof}
	As discussed in \cref{lm:ncdf-0-ncdf-2-mod-q3-equiv}, the map $\psi'$ is induced by the composition of the maps
	\begin{equation*}
		- h_{\dr} \circ \varphi^{-1} \circ B =
		- h_{\dr} \circ \varphi^{-1} \circ \left( \idd - \t \right) \circ H \circ \N,
	\end{equation*}
	where by induced we mean that
	we let it act on elements of the quotient $\tensr{A} / \Im \left( \idd - \t \right)$ and project the result
	from $\ncdf{A}[2][]$ to the quotient $\ncdf{A}[2][] / \Im \left( \qdr^3 \right)$. See the
	first two rows of \cref{fig:psi'-definition-and-reduced}.

	Given $a_1, \dots, a_m \in A$, we have
	\begin{equation*}
		\begin{aligned}
			\left( H \circ \N \right) \left( a_1, \dots, a_m \right) ={} &
			H \left(
			  \sum_{i=0}^{m-1}
			  (-1)^{
				  \left( \degb{a_1} + \dots + \degb{a_i} \right) \cdot
				  \left( \degb{a_{i+1}} + \dots + \degb{a_m} \right)
			  }
			a_{i+1} \, \dots \, a_m \, a_1 \, \dots \, a_{i}
			\right)
			\\
			={}                                                          &
			\sum_{i=0}^{m-1} \sum_{k=0}^{\infty}
			                 (-1)^{k +
				                 \left( \degb{a_1} + \dots + \degb{a_i} \right) \cdot
				                 \left( \degb{a_{i+1}} + \dots + \degb{a_m} \right)
			                 }
			\\
			                                                             & \qquad\qquad\quad\quad
			e \, \left( \mu_0 \left( 1 \right) \, e \right)^k \, a_{i+1} \, \dots \, a_m \, a_1 \, \dots \, a_{i}.
		\end{aligned}
	\end{equation*}
	We also have, for $m \geq 2$, the identity
	\begin{equation*}
		\begin{aligned}
			\left( \varphi^{-1} \circ \left( \idd - \t \right) \right) \left(
			a_1 \, \dots \, a_m \right) & =
			\ul{a_1} \, \dots \, a_m -
				                       (-1)^{\degb{a_m} \cdot \left( \degb{a_1} + \dots + \degb{a_{m-1}} \right)}
			\ul{a_m} \, a_1 \, \dots \, a_{m-1}
			\\
			                            & =
			\ul{a_1} \, \dots \, a_m - a_1 \, \dots \, \ul{a_m}
			\\
			                            & =
			- \qdr^2 \left( \ul{a_1} \, \dots \, \ul{a_m} \right),
		\end{aligned}
	\end{equation*}
	and hence
	\begin{equation*}
		\begin{aligned}
			\left( -h_{\dr} \circ \varphi^{-1} \circ \left( \idd - \t \right) \right) \left(
			a_1 \, \dots \, a_m
			\right) & =
			\left( h_{\dr} \circ \qdr^2 \right) \left( \ul{a_1} \, \dots \, \ul{a_m} \right)
			\\
			        & =
			\ul{a_1} \, \dots \, \ul{a_m} -
			\left( \qdr^3 \circ h_{\dr} \right) \left( \ul{a_1} \, \dots \, \ul{a_m} \right).
		\end{aligned}
	\end{equation*}
	Our choice of notation is such that we always write explicitly \textit{all} the underlined
	elements of a form, so that all the terms between the ellipses are always without lines.
	Thus,
	\begin{equation*}
		\begin{aligned}
			\psi' \left( a_1 \, \dots \, a_m \right) ={} &
			- \left( \pi_2 \circ h_{\dr} \circ \varphi^{-1} \circ \left( \idd - \t \right) \circ H \circ \N \right) \left( a_1 \, \dots \, a_m \right)
			\\
			={}                                          &
			\sum_{k=0}^{\infty}
			(-1)^{k}
			\ul{e} \, \left( \mu_0 \left( 1 \right) \, e \right)^k \, a_{1} \, \dots \, a_{m-1} \, \ul{a_m}
			\\
			                                             & +
			\sum_{i=1}^{m-1} \sum_{k=0}^{\infty}
			                 (-1)^{k +
				                 \left( \degb{a_1} + \dots + \degb{a_{i}} \right) \cdot
				                 \left( \degb{a_{i+1}} + \dots + \degb{a_m} \right)
			                 }
			\ul{e} \, \left( \mu_0 \left( 1 \right) \, e \right)^k \, a_{i+1} \, \dots \, a_m \, a_1 \, \dots \,
			\ul{a_{i}}
			\\
			={}                                          &
			\sum_{i=1}^{m} \sum_{k=0}^{\infty}
			               (-1)^{\degb{a_1} + \dots + \degb{a_{i}} + k + 1}
			a_1 \, \dots \, \ul{a_{i}} \, \ul{e} \, \left( \mu_0 \left( 1 \right) \, e \right)^k \, a_{i+1} \,
			\dots \, a_m.
		\end{aligned}
	\end{equation*}
\end{proof}

\begin{cor} \label{cor:psi-tag-i-underlined}
	The map $\psi' \colon \ncdfr{\mathcal{A}}[0][] \rightharpoonup \ncdf{\mathcal{A}}[2][] / \Im \left( \qdr^3 \right)$
	has the explicit form
	\begin{equation} \label{eq:psi-tag-i-underlined}
		\psi' \left( a_1 \, \dots \, a_m \right) =
		\sum_{i=1}^{m}
		(-1)^{\degb{a_1} + \dots + \degb{a_{i}} + 1}
		a_1 \, \dots \, \ul{a_{i}} \, \ul{e} \, a_{i+1} \, \dots \, a_m
		\mod \left< *e* \right>.
	\end{equation}
	\qed
\end{cor}

\begin{cor} \label{cor:psi-tag-conv-form}
	The map $\psi' \colon \ncdfr{\mathcal{A}}[0][] \rightharpoonup \ncdf{\mathcal{A}}[2][] / \Im \left( \qdr^3 \right)$
	has the explicit form
	\begin{equation} \label{eq:psi-tag-conv-form}
		\psi' \left( a \, l \right) =
		(-1)^{\degb{a} + \degb{l_{(1)}} + 1} \ul{a} \, l_{(1)} \, \ul{e} \, l_{(2)}
		\mod \left< *e* \right>
	\end{equation}
	for $a \in A$ and $l \in \tens{A}$.
\end{cor}
\begin{proof}
	In the expression \eqref{eq:psi-tag-i-underlined} for $\psi'$, the underlined elements are $a_i$ and $e$.
	Since we are working modulo $\Im \left( \qdr^3 \right)$ and $\left< *e* \right> = \degen[A][2]$, we can use the relations
	\begin{equation*}
		\begin{aligned}
			\qdr^3 \left(
			\ul{a_1} \, \dots \, \ul{a_{i}} \, \ul{e} \, 	a_{i+1} \, \dots \, a_m
			\right)
			={} &
			a_1  \, \dots \, \ul{a_{i}} \, \ul{e} \,	a_{i+1} \, \dots \, a_m
			\\
			    & -
			\ul{a_1} \, \dots \, a_i \, \ul{e} \, a_{i+1} \, \dots \, a_m
			\\
			    & +
			\ul{a_1} \, \dots \, \ul{a_{i}} \, e \, a_{i+1} \, \dots \, a_m
		\end{aligned}
	\end{equation*}
	for $2 \leq i \leq m$ to replace \cref{eq:psi-tag-i-underlined} with
	\begin{equation*}
		\psi' \left( a_1 \, \dots \, a_m \right) =
		\sum_{i=1}^{m}
		(-1)^{\degb{a_1} + \dots + \degb{a_{i}} + 1}
		\ul{a_1} \, \dots \, a_{i} \, \ul{e} \, a_{i+1} \, \dots \, a_m,
	\end{equation*}
	which is the more explicit version of \cref{eq:psi-tag-conv-form}.
\end{proof}

We can use the formulas we wrote for $\psi'$ to describe the relation between $\phi_2$ and the corresponding
pre-$\infty$-trace $\theta$, related via \cref{eq:naive-trace-from-hip}.

\begin{cor}
	Let $\phi_2$ be a $\degen$-unital pre-homotopy inner product on $\mathcal{A}$.
	Then the corresponding pre-$\infty$-trace $\theta$ is given by the formula
	\begin{equation} \label{eq:theta-phi-unital-rel}
		\theta \left( a \, l \right) = (-1)^{\degb{a} + \degb{l_{(1)}} + 1} \phi_2 \left( \ul{a} \, l_{(1)} \, \ul{e} \, l_{(2)} \right).
	\end{equation}
\end{cor}
\begin{proof}
	Immediate from \cref{cor:psi-tag-conv-form}.
\end{proof}

\begin{cor} \label{cor:psi-tag-descends-reduced}
	The map $\psi'$ satisfies
	\begin{equation*}
		\psi' \left( \degen[\mathcal{A}][0][] \right) \subseteq
		\left< \Im \left( \qdr^3 \right), \degen[A][2][], \degenu[A][2][] \right> / \Im \left( \qdr^3 \right)
	\end{equation*}
	and hence $\psi'$ descends to a well-defined chain map
	\begin{equation*}
		\psi'_{\textrm{red}} \colon \ncdfrred{A}[0][] \rightharpoonup
		\ncdf{\mathcal{A}}[2][] / \left< \Im \left( \qdr^3 \right), \degen[A][2][], \degenu[A][2][] \right>.
	\end{equation*}
	See the bottom two rows of \cref{fig:psi'-definition-and-reduced}.
\end{cor}
\begin{proof}
	By cyclic symmetry, the degenerate subcomplex $\degen[\mathcal{A}][0][] \subset \ncdfr{A}[0][]$
	is generated by elements of the form $e*$.
	The explicit formula of \cref{cor:psi-tag-conv-form} shows that $\psi'$ maps $e*$ to
	a linear combination of elements of the form $\ul{e}{*}\ul{e}{*}$ modulo
	$\left< \Im \left( \qdr^3 \right), *e* \right> / \Im \left( \qdr^3 \right)$,
	so the image of $\psi'$ lies in $\left<  \Im \left( \qdr^3 \right), *e*, \ul{e}{*}\ul{e}{*} \right> / \Im \left( \qdr^3 \right)$.
	By \cref{lm:desc-ncdf-2-strongly-unital-subcomplex}, this target space coincides with
	$\left< \Im \left( \qdr^3 \right), \degen[A][2][], \degenu[A][2][] \right> / \Im \left( \qdr^3 \right)$,
	completing the proof.
\end{proof}

\begin{cor}
	Let $\phi_2$ be a strongly unital pre-homotopy inner product on $\mathcal{A}$.
	Then the corresponding pre-$\infty$-trace $\theta$ is unital in the sense of \cref{dfn:pre-infty-trace-unital}.
\end{cor}
\begin{proof}
	Immediate from \cref{cor:psi-tag-descends-reduced}.
\end{proof}

\begin{cor} \label{cor:cyclic-structure-gives-pre-infinity-trace}
	Let $\phi_2$ be a strict pre-homotopy inner product on $\mathcal{A}$, i.e.,
	a cyclic structure. Then the corresponding pre-$\infty$-trace
	$\theta$ is also strict, i.e., $\theta_k = 0$ for $k > 1$, and we have
	\begin{equation*}
		\theta \left( a \right) = \theta_1 \left( a \right) = (-1)^{\degb{a} + 1} \phi_2^{0,0} \left( \ul{a}, \ul{e} \right) =
		\phi_2^{0,0} \left( \ul{e}, \ul{a} \right) = \phi_2 \left( \ul{e}, \ul{a} \right).
	\end{equation*}
\end{cor}
\begin{proof}
	Immediate from \cref{lm:psi-tag-explicit-formula}.
\end{proof}

\begin{rem} \label{rem:pre-homotopy-inn-prod-to-trace-zero-term}
	Let $\phi_2$ be a $\degen$-unital pre-homotopy inner product on $\mathcal{A}$.
	Then, by \cref{cor:psi-tag-conv-form}, the corresponding pre-$\infty$-trace $\theta$ satisfies
	\begin{equation*}
		\theta_1 \left( \mu_0 \left( 1 \right) \right) = \phi_2 \left( \psi' \left( \mu_0 \left( 1 \right) \right) \right) =
		\phi_2 \left( \ulz{\mu}, \ul{e} \right).
	\end{equation*}
	Thus, if
	\begin{equation} \label{eq:phi-2-mu-e-cond}
		\phi_2 \left( \ulz{\mu}, \ul{e} \right) = 0,
	\end{equation}
	the corresponding pre-$\infty$-trace $\theta$
	satisfies $\theta_1 \left( \mu_0 \left( 1 \right) \right) = 0$, i.e.,
	condition \eqref{eq:infty-trace-0-rel}. In this case, $\theta$ can be upgraded to a full
	$\infty$-trace by setting $\theta_0 = 0$.

	Condition \eqref{eq:phi-2-mu-e-cond} appears in the literature, for example
	in \cite[Condition 9, Definition 1.1]{Solomon2016}, where it is part of
	the definition of a cyclic unital structure.
	See also \cref{rem:extended-vs-standard-trace}.
\end{rem}

\subsubsection{Pre-\texorpdfstring{$\infty$}{Infinity}-Traces to Pre-Homotopy Inner Products}
In what follows, unless explicitly stated, we do not assume that $\mathcal{A}$ is unital.
Given a pre-$\infty$-trace $\theta$ on $\mathcal{A}$, we can use the map $\psi$,
given by \cref{eq:psi-def} of \cref{dfn:psi-ncdf-2-0}, to
define a pre-homotopy inner product $\phi_2$, called the \textbf{corresponding pre-homotopy inner product}, by
\begin{equation} \label{eq:hip-from-naive-trace}
	\phi_2 \defeq \theta \circ \psi.\footnote{Here, we think of $\phi_2, \theta, \psi$
		as graded $R$-linear maps with codomain $R$ of degrees $2-n,1-n,1$, respectively.}
\end{equation}
Note that the map $\psi$ is defined for every Banach $\Ainf$-algebra $\mathcal{A}$
and allows us to convert pre-$\infty$-traces to pre-homotopy inner products.
Unitality is needed only for the inverse map $\psi'$ and for the resulting
homotopy-equivalence statement.

To make the relation between $\phi_2$ and $\theta$ more explicit, let us write down several convenient formulas for
the map $\psi$.

\begin{figure}[htb]
	\centering
	\adjustbox{scale=0.75,center}
	{
		\begin{tikzcd}
			{\ncdf{A}[2][]} & {\ncdf{A}[1][]} & {\tensr{A}} & {\tensr{A}} & {\tensr{A}} \\
			{\ncdf{A}[2][] / \Im \left( \qdr^3 \right)} & {\Im \left( \qdr^2 \right)} & {\Im \left( \idd - \t \right)} &&
			{\tensr{A} / \Im \left( \idd - \t \right) = \ncdfr{A}[0][]} \\
			{\ncdf{A}[2][] / \left< \Im \left( \qdr^3 \right), *e*, \ul{e} \ul{e} \right>} &&&&
			{\tensr{A} / \left< \Im \left( \idd - \t \right), *e* \right> = \ncdfrred{A}[0][]}
			\arrow["\qdr^2", from=1-1, to=1-2]
			\arrow[two heads, from=1-1, to=2-1]
			\arrow["\varphi", "\cong"', from=1-2, to=1-3]
			\arrow["h", from=1-3, to=1-4]
			\arrow[hook, from=2-2, to=1-2]
			\arrow[hook, from=2-3, to=1-3]
			\arrow["\beta", harpoon, curve={height=-25pt}, from=1-3, to=1-5]
			\arrow["{T'}", harpoon, from=1-4, to=1-5]
			\arrow[two heads, from=1-5, to=2-5]
			\arrow["{\overline{\qdr}^2}", "{\cong}"', from=2-1, to=2-2]
			\arrow["\psi", harpoon, curve={height=25pt}, from=2-1, to=2-5]
			\arrow[two heads, from=2-1, to=3-1]
			\arrow["{\overline{\varphi}}", "{\cong}"', from=2-2, to=2-3]
			\arrow["{\overline{\beta}}", harpoon, from=2-3, to=2-5]
			\arrow[two heads, from=2-5, to=3-5]
			\arrow["{\psi_{\textrm{red}}}"', harpoon, from=3-1, to=3-5]
		\end{tikzcd}
	}
	\caption{Maps involved in the definition of $\psi$.}
	\label{fig:psi-definition-and-reduced}
\end{figure}

\begin{lm} \label{lm:psi-T'-rot-identity}
	We have the identity
	\begin{equation} \label{eq:lm:psi-T'-rot-identity}
		\begin{aligned}
			\psi \left( \ul{x} \, \underbrace{u_1 \, \dots \, u_i}_{u} \, \ul{y} \, \underbrace{v_1 \, \dots \, v_j}_{v} \right) ={} &
			\left( h'b' - bh' \right) \left( \idd + \t + \dots + {\t}^j \right) \left( x \, u \, y \, v \right).
		\end{aligned}
	\end{equation}
\end{lm}
\begin{proof}
	As discussed in \cref{lm:ncdf-0-ncdf-2-mod-q3-equiv}, the map $\psi$ is induced by the composition of the maps
	\begin{equation*}
		\beta \circ \varphi \circ \qdr^2 =
		T' \circ h \circ \varphi \circ \qdr^2,
	\end{equation*}
	where by induced we mean that
	we let it act on elements of the quotient $\ncdf{A}[2][] / \Im \left( \qdr^3 \right)$
	and project the result from $\tensr{A}[]$ to the quotient $\tensrcyc{A}[] = \ncdfr{A}[0][]$
	(see the first two rows of \cref{fig:psi-definition-and-reduced}).

	We have
	\begin{equation*}
		\qdr^2 \left( \ul{x} \, u \, \ul{y} \, v \right) =
		x \, u \, \ul{y} \, v - \ul{x} \, u \, y \, v =
			(-1)^{\left( \degb{y} + \degb{v} \right) \left( \degb{x} + \degb{u} \right)} \ul{y} \, v \, x \, u - \ul{x} \, u \, y \, v
	\end{equation*}
	and hence
	\begin{equation*}
		\begin{aligned}
			\left( \varphi \circ \qdr^2 \right) \left( \ul{x} \, u \, \ul{y} \, v \right) & =
			(-1)^{\left( \degb{y} + \degb{v} \right) \left( \degb{x} + \degb{u} \right)} y \, v \, x \, u - x \, u \, y \, v
			\\
			                                                                              & =
			- \left( \idd - {\t}^{j+1} \right) \left( x \, u \, y \, v \right)
			\\
			                                                                              & =
			- \left( \idd - \t \right) \left( \idd + \t + \dots + {\t}^j \right) \left( x \, u \, y \, v \right).
		\end{aligned}
	\end{equation*}
	Then using the fact that $\beta = T' h$ with $T' = bh' - h'b'$ and $T'\N = \left( \idd - \t \right)T$, we have
	\begin{equation*}
		\begin{aligned}
			\left( \beta \, \varphi \qdr^2 \right) \left( \ul{x} \, u \, \ul{y} \, v \right)
			\eqwithref[eq:def-beta]     &
			- \beta \left( \idd - \t \right) \left( \idd + \t + \dots + {\t}^j \right) \left( x \, u \, y \, v \right)
			\\
			\eqwithref                  &
			- T' h \left( \idd - \t \right) \left( \idd + \t + \dots + {\t}^j \right) \left( x \, u \, y \, v \right)
			\\
			\eqwithref[eq:Nh'-homotopy] &
			T' \left( \N h' - \idd \right) \left( \idd + \t + \dots + {\t}^j \right) \left( x \, u \, y \, v \right)
			\\
			\eqwithref[eq:1-t-T]        &
			\left( \idd - \t \right) T h' \left( \idd + \t + \dots + {\t}^j \right) \left( x \, u \, y \, v \right)
			\\
			                            &
			- T' \left( \idd + \t + \dots + {\t}^j \right) \left( x \, u \, y \, v \right)
			\\
			\eqwithref[eq:T'-def]       &
			\left( h'b' - bh' \right) \left( \idd + \t + \dots + {\t}^j \right) \left( x \, u \, y \, v \right)
			\mod \Im \left( \idd - \t \right).
		\end{aligned}
	\end{equation*}
\end{proof}

Since the operator $h'$ divides an element by its weight, writing a more explicit formula for $\psi$ than the one given in
\cref{lm:psi-T'-rot-identity} is quite cumbersome. However, we can write a simple formula for $\psi$
when the $\Ainf$-structure $\mu$ encodes a single operation:

\begin{cor}
	Assume that the $\Ainf$-structure $\mu$ satisfies $\mu_k = 0$ for $k \neq l$. Then
	\begin{align}\label{eq:psi-single-operation}
		\psi \left( \ul{x} \, \underbrace{u_1 \, \dots \, u_i}_{u} \, \ul{y} \, \underbrace{v_1 \, \dots \, v_j}_{v} \right) ={} &
		\frac{\left( l - 1 \right) \left( j + 1 \right)}{\left(i + j + 3 - l \right) \left( i + j + 2 \right)}
		b \left( x \, u \, y \, v \right)
		\\
		                                                                                                                         & -
		\frac{1}{i + j + 3 - l} \left( b - b' \right) \left( \idd + \t + \dots + {\t}^j \right) \left( x \, u \, y \, v \right).
		\notag
	\end{align}
\end{cor}
\begin{proof}
	When $\mu_k = 0$ for $k \neq l$, we have
	\begin{equation*}
		\begin{aligned}
			\left( h'b' - bh' \right) \left( x \, u \, y \, v \right) & =
			\frac{1}{i + j + 2 + \left( 1 - l \right)} b' \left( x \, u \, y \, v \right) -
			\frac{1}{i+j+2} b \left( x \, u \, y \, v \right)
			\\
			                                                          & =
			\frac{1}{i + j + 3 - l} \left( b + \left( b' - b \right) \right) \left( x \, u \, y \, v \right) -
			\frac{1}{i+j+2} b \left( x \, u \, y \, v \right)
			\\
			                                                          & =
			\frac{l - 1}{\left(i + j + 3 - l \right) \left( i + j + 2 \right)} b \left( x \, u \, y \, v \right)
			- \frac{1}{i + j + 3 - l} \left( b - b' \right) \left( x \, u \, y \, v \right).
		\end{aligned}
	\end{equation*}
	Using \cref{lm:psi-T'-rot-identity} and the identity $b \t = b \mod \left( \idd - \t \right)$, we have
	\begin{equation*}
		\begin{aligned}
			\psi \left( \ul{x} \, u \, \ul{y} \, v \right) \equiv{} &
			\left( h'b' - bh' \right) \left( \idd + \t + \dots + {\t}^j \right) \left( x \, u \, y \, v \right)
			\\
			\equiv{}                                                &
			\frac{\left( l - 1 \right) \left( j + 1 \right)}{\left(i + j + 3 - l \right) \left( i + j + 2 \right)}
			b \left( x \, u \, y \, v \right)
			\\
			                                                        & -
			\frac{1}{i + j + 3 - l} \left( b - b' \right) \left( \idd + \t + \dots + {\t}^j \right) \left( x \, u \, y \, v \right)
			\mod \Im \left( \idd - \t \right).
		\end{aligned}
	\end{equation*}
\end{proof}

\begin{cor} \label{cor:psi-dga}
	Assume that $\mu_k = 0$ for $k \neq 1, 2$, i.e., $\mathcal{A}$ corresponds to a differential graded algebra.
	Then $\psi$ has the form
	\begin{equation*}
		\begin{gathered}
			\psi \left( \ul{x_0} \, x_1 \, \dots \, x_i \, \ul{x_{i+1}} \, x_{i+2} \, \dots \, x_{i+j+1} \right) =
			\\
			\frac{j+1}{\left( i + j + 1 \right) \left( i + j + 2 \right)}
			\sum_{r=0}^i (-1)^{\degb{x_0} + \dots + \degb{x_{r-1}}}
			x_0 \, \dots \, x_{r-1} \, \mu_2 \left( x_r, x_{r+1} \right) \, x_{r+2} \, \dots \, x_{i+j+1} +
			\\
			\frac{- \left( i+1 \right)}{\left( i + j + 1 \right) \left( i + j + 2 \right)}
			\sum_{r=i+1}^{i+j} (-1)^{\degb{x_0} + \dots + \degb{x_{r-1}}}
			x_0 \, \dots \, x_{r-1} \, \mu_2 \left( x_r, x_{r+1} \right) \, x_{r+2} \, \dots \, x_{i+j+1} +
			\\
			\frac{- \left( i+1 \right)}{\left( i + j + 1 \right) \left( i + j + 2 \right)}
			(-1)^{\degb{x_{i+j+1}} \left( \degb{x_0} + \dots + \degb{x_{i+j}} \right)}
			\mu_2 \left( x_{i+j+1}, x_0 \right) \, x_1 \, \dots \, x_{i+j}.
		\end{gathered}
	\end{equation*}
\end{cor}
\begin{proof}
	We can split $\mu$ as a sum $\mu = \mu^1 + \mu^2$ where $\mu^i_j = \delta_{i,j} \cdot \mu_j$, i.e., $\mu^i$ involves
	only the contributions from the operation $\mu_i$. Since the formula for $\psi = \psi_{\mu}$ is linear in $\mu$, we can consider
	the contribution of each $\mu^i$ separately.

	When $\mu = \mu^1$, we have $b = b'$ and \cref{eq:psi-single-operation} shows that $\psi_{\mu^1} = 0$. Thus, we
	can assume that $\mu_k = 0$ for $k \neq 2$. Then \cref{eq:psi-single-operation} gives us the formula
	\begin{gather*}
		\psi \left( \ul{x_0} \, x_1 \, \dots \, x_i \, \ul{x_{i+1}} \, x_{i+2} \, \dots \, x_{i+j+1} \right) =
		\frac{j + 1}{\left(i + j + 1 \right) \left( i + j + 2 \right)}
		b \left( x_0 \, \dots \, x_{i+j+1} \right) +
		\\
		\frac{-1}{i + j + 1} \left( b - b' \right) \left( \idd + \t + \dots + {\t}^j \right) \left( x_0 \, \dots \, x_{i+j+1} \right).
	\end{gather*}
	Writing the above formula explicitly and identifying cyclically equivalent terms gives us the formula of \cref{cor:psi-dga}.
\end{proof}

\begin{ex}
	Let $a,b,c \in A$. Using \cref{eq:lm:psi-T'-rot-identity}, we have
	\begin{equation} \label{eq:psi-ula-ulb}
		\psi \left( \ul{a} \, \ul{b} \right) =
		\frac{1}{2} \left( \mu_2 \left( a, b \right) - (-1)^{\degb{a} \degb{b}} \mu_2 \left( b, a \right) \right) +
		\frac{1}{6} \left( \mu_0 \left( 1 \right) \, a \, b - (-1)^{\degb{a}} a \, \mu_0 \left( 1 \right) \, b  \right),
	\end{equation}
	which, when $\mu_0 \left( 1 \right) = 0$, coincides with the formula
	\begin{equation} \label{eq:psi-ula-ulb-dga}
		\psi \left( \ul{a} \, \ul{b} \right) = \frac{1}{2} \mu_2 \left( a, b \right) - \frac{1}{2} (-1)^{\degb{a} \degb{b}} \mu_2 \left( b, a \right)
	\end{equation}
	from \cref{cor:psi-dga}. When $\mu_k = 0$ for $k \neq 1,2$, \cref{cor:psi-dga} yields the formulas
	\begin{align}
		\psi \left( \ul{a} \, b \, \ul{c} \right) & =
		\frac{1}{6} \left( \mu_2 \left( a, b \right) \, c + (-1)^{\degb{a}} a \, \mu_2 \left( b, c \right) \right) -
		\frac{1}{3} (-1)^{\degb{c} \left( \degb{a} + \degb{b} \right)} \mu_2 \left( c, a \right) \, b,
		\label{eq:psi-ula-b-ulc-dga}
		\\
		\psi \left( \ul{a} \, \ul{b} \, c \right) & =
		\frac{1}{3} \mu_2 \left( a, b \right) \, c -
		\frac{1}{6} \left( (-1)^{\degb{a}} a \, \mu_2 \left( b, c \right) +
		                                              (-1)^{\degb{c} \left( \degb{a} + \degb{b} \right)} \mu_2 \left( c, a \right) \, b
		\right).
		\label{eq:psi-ula-ulb-c-dga}
	\end{align}

	Note that \cref{eq:psi-ula-ulb} is consistent with the symmetry
	$\psi \left( \ul{a} \, \ul{b} \right) = (-1)^{\degb{a} \degb{b} + 1} \psi \left( \ul{b}, \ul{a} \right)$,
	and that \cref{eq:psi-ula-b-ulc-dga,eq:psi-ula-ulb-c-dga} are consistent with the symmetry
	\begin{equation*}
		\psi \left( \ul{a} \, b \, \ul{c} \right) = (-1)^{\degb{c} \left( \degb{a} + \degb{b} \right) + 1} \psi \left( \ul{c} \, \ul{a} \, b \right).
	\end{equation*}
\end{ex}

\begin{lm} \label{lm:psi-respects-units}
	Assume $\mathcal{A}$ is unital with unit $e$. Then the map $\psi$ satisfies
	\begin{equation*}
		\psi \left( \left< \Im \left( \qdr^3 \right), \degen[\mathcal{A}][2][], \degenu[\mathcal{A}][2][] \right> / \Im \left( \qdr^3 \right) \right) \subseteq \degen[\mathcal{A}][0][]
	\end{equation*}
	and hence $\psi$ descends to a well-defined chain map
	\begin{equation*}
		\psi_{\textrm{red}} \colon \ncdf{\mathcal{A}}[2][] / \left< \Im \left( \qdr^3 \right), \degen[A][2][], \degenu[A][2][] \right>
		\rightharpoonup \ncdfrred{A}[0][],
	\end{equation*}
	see the bottom two rows of \cref{fig:psi-definition-and-reduced}.
\end{lm}
\begin{proof}
	By \cref{lm:desc-ncdf-2-strongly-unital-subcomplex}, we have
	\begin{equation*}
		\left< \Im \left( \qdr^3 \right), \degen[\mathcal{A}][2][], \degenu[\mathcal{A}][2][] \right> / \Im \left( \qdr^3 \right) =
		\left< \Im \left( \qdr^3 \right), \ul{e} \, \ul{e}, *e* \right> / \Im \left( \qdr^3 \right).
	\end{equation*}
	Note that we have $b \left( a{*}e{*}b \right), b' \left( a{*}e{*}b \right) \in \degen[\mathcal{A}][0][]$ for $a,b \in A$.
	That is, both $b'$ and $b$ map elementary tensors in which the unit $e$ appears somewhere \textit{as long as it is not} the
	first or last element, to a sum of
	elements which contain the unit $e$ somewhere. Now let $\alpha \in \ncdf{A}[2][] / \Im \left( \qdr^3 \right)$
	be an equivalence class generated by an element of the form $*e*$, so that
	\begin{equation*}
		\alpha \equiv \ul{x} \, \underbrace{u_1 \, \dots \, u_i}_{u} \, \ul{y} \, \underbrace{v_1 \, \dots \, v_j}_{v} \mod \Im \left( \qdr^3 \right),
	\end{equation*}
	where we have either $u_r = e$ for some $1 \leq r \leq i$, or $v_r = e$ for some $1 \leq r \leq j$.
	Modulo $\Im \left( \qdr^3 \right)$, we can use the relations to cyclically shift the underlines and hence we can assume that the unit $e$ appears in $u$. From \cref{eq:lm:psi-T'-rot-identity},
	we see that $\psi \left( \alpha \right)$ is a sum of terms of the form
	$\left( h'b' - bh' \right) \left( {\t}^r \right) \left( x \, u \, y \, v \right)$ for $0 \leq r \leq j$. Since we
	rotate at most $j$ times, the unit $e$ never appears as the first or last element of ${\t}^r \left( x \, u \, y \, v \right)$
	and so $\psi \left( \alpha \right) \in \degen[\mathcal{A}][0][]$.
	Finally, \cref{eq:psi-ula-ulb} shows that $\psi \left( \ul{e} \, \ul{e} \right) = e + \frac{1}{3} \mu_0 \left( 1 \right) \, e \, e$
	also belongs to $\degen[\mathcal{A}][0][]$.

	Thus, we have shown that
	$\psi \left( \left< \Im \left( \qdr^3 \right), \degen[\mathcal{A}][2][], \degenu[\mathcal{A}][2][] \right> / \Im \left( \qdr^3 \right) \right) \subseteq \degen[\mathcal{A}][0][]$.
	Since we have $\ncdfrred{\mathcal{A}}[0][] = \ncdfr{\mathcal{A}}[0][] / \degen[\mathcal{A}][0][]$, the conclusion follows.
\end{proof}

Given a pre-$\infty$-trace $\theta$, we can use the calculations we did to draw some conclusions about the corresponding
pre-homotopy inner product $\phi_2 = \theta \circ \psi$.

\begin{cor}
	Assume $\mathcal{A}$ is unital and let $\theta$ be a unital pre-$\infty$-trace on $\mathcal{A}$. Then the corresponding
	pre-homotopy inner product $\phi_2$ is strongly unital.
\end{cor}
\begin{proof}
	Immediate from \cref{lm:psi-respects-units}.
\end{proof}

When $\theta$ is strict, the corresponding pre-homotopy inner product $\phi_2$ isn't necessarily strict. However, we do have:
\begin{cor} \label{cor:strict-theta-cdga-to-strict-phi}
	Assume that $\mu_k = 0$ for $k > 2$, i.e., $\mathcal{A}$ corresponds to a curved differential graded algebra,
	and let $\theta$ be a strict pre-$\infty$-trace on $\mathcal{A}$. Then
	the corresponding pre-homotopy inner product $\phi_2$ is also strict and we have
	\begin{equation} \label{eq:phi_2-from-strict-trace-cdga}
		\phi_2 \left( \ul{a}, \ul{b} \right) = \theta \left( \mu_2 \left( a, b \right) \right).
	\end{equation}
\end{cor}
\begin{proof}
	When $\mu_k = 0$ for $k > 2$, the formula \eqref{eq:lm:psi-T'-rot-identity} of \cref{lm:psi-T'-rot-identity} implies
	that $\psi$ maps elements of weight $r$ to a sum of elements which have weight at least $r - 1$.
	Hence, when $\theta$ is strict, we have $\phi_2^{i,j} = 0$ for $(i,j) \neq 0$. When $i = j = 0$, we have
	\begin{equation*}
		\begin{aligned}
			\phi_2 \left( \ul{a}, \ul{b} \right)
			\stackrel{\phantom{\eqref{eq:psi-ula-ulb}}}{=} &
			\theta \left( \psi \left( \ul{a}, \ul{b} \right) \right)
			\\
			\stackrel{\eqref{eq:psi-ula-ulb}}{=}           &
			\frac{1}{2} \cdot \theta_1 \left( \mu_2 \left( a, b \right) - (-1)^{\degb{a} \cdot \degb{b}} \mu_2 \left( b, a \right) \right)
			\\
			                                               & +
			\frac{1}{6} \cdot \left(
			\theta_3 \left( \mu_0 \left( 1 \right), a, b \right) - (-1)^{\degb{a}} \theta_3 \left( a, \mu_0 \left( 1 \right), b \right)
			\right)
			\\
			\stackrel{\phantom{\eqref{eq:psi-ula-ulb}}}{=} &
			\frac{1}{2} \cdot \theta \left( \mu_2 \left( a, b \right) - (-1)^{\degb{a} \cdot \degb{b}} \mu_2 \left( b, a \right) \right).
		\end{aligned}
	\end{equation*}
	Since $\theta$ is a strict pre-$\infty$-trace, it satisfies the relations
	of \eqref{eq:infty-trace-strict-rel}, for which the $k = 2$ case reads
	\begin{equation*}
		\theta \left( \mu_2 \left( a, b \right) \right) =
		(-1)^{\degb{a}\degb{b} + 1} \theta \left( \mu_2 \left( b, a \right) \right),
	\end{equation*}
	which implies \cref{eq:phi_2-from-strict-trace-cdga}.
\end{proof}

\begin{cor}
	Assume that $\mu_k = 0$ for $k \neq 1, 2$, i.e., $\mathcal{A}$ corresponds to a
	differential graded algebra,
	and let $\theta$ be a strict pre-$\infty$-trace on $\mathcal{A}$.
	Let $\phi_2$ be the corresponding pre-homotopy inner product to $\theta$.
	Then the pre-superpotential function associated to $\phi_2$ takes the form
	\begin{equation} \label{eq:pre-sp-strict-trace-dga}
		\SP[b][>0] = \theta \left(
		\frac{1}{2} \mu_2 \left( \mu_1 \left( b \right), b \right) +
		\frac{1}{3} \mu_2 \left( \mu_2 \left( b, b \right), b \right)
		\right).
	\end{equation}
\end{cor}
\begin{proof}
	Let $b \in \tc{A}$ and consider the element $\G{b}[2] \in \ncdf{A}[2][1]$,
	given by \eqref{eq:G-b-2-braidop-1}, so that
	\begin{equation*}
		\SP[b][>0] \stackrel{\eqref{eq:pre-sp-phi-2-braidop-1}}{=}
		\phi_2 \left( \G{b}[2] \right) \stackrel{\eqref{eq:hip-from-naive-trace}}{=}
		\theta \left( \psi \left( \G{b}[2] \right) \right).
	\end{equation*}
	By \cref{cor:strict-theta-cdga-to-strict-phi}, $\phi_2$ is strict
	so only the summands of weight two in $\G{b}[2]$ survive in $\SP[b][>0]$,
	and we get
	\begin{equation*}
		\SP[b][>0] =
		\frac{1}{2} \phi_2 \left( \ul{\mu_1 \left( b \right)}, \ul{b} \right) +
		\frac{1}{3} \phi_2 \left( \ul{\mu_2 \left( b, b \right)}, \ul{b} \right)
		\stackrel{\eqref{eq:phi_2-from-strict-trace-cdga}}{=}
		\theta \left(
		\frac{1}{2} \mu_2 \left( \mu_1 \left( b \right), b \right) +
		\frac{1}{3} \mu_2 \left( \mu_2 \left( b, b \right), b \right)
		\right).
	\end{equation*}
\end{proof}

\begin{rem}
	Let $\theta$ be a pre-$\infty$-trace on $\mathcal{A}$. Then, by \cref{eq:psi-ula-ulb},
	the corresponding pre-homotopy inner product $\phi_2$ satisfies
	\begin{equation*}
		\begin{aligned}
			-\frac{1}{2} \phi_2 \left( \ulz{\mu}, \ulz{\mu} \right) & =
			-\frac{1}{2} \theta \left( \psi \left( \ulz{\mu}, \ulz{\mu} \right) \right)
			\\
			                                                        & =
			-\frac{1}{2} \theta_1 \left(  \mu_2 \left( \mu_0 \left( 1 \right), \mu_0 \left( 1 \right) \right) \right) -
			\frac{1}{6} \theta_3 \left( \mu_0 \left( 1 \right), \mu_0 \left( 1 \right), \mu_0 \left( 1 \right) \right).
		\end{aligned}
	\end{equation*}
	Thus, if $\theta$ satisfies
	\begin{equation*}
		\theta_1 \left(  \mu_2 \left( \mu_0 \left( 1 \right), \mu_0 \left( 1 \right) \right) \right) = -
		\frac{1}{3} \theta_3 \left( \mu_0 \left( 1 \right), \mu_0 \left( 1 \right), \mu_0 \left( 1 \right) \right),
	\end{equation*}
	then the corresponding pre-homotopy inner product $\phi_2$ satisfies
	$\phi_2 \left( \ulz{\mu}, \ulz{\mu} \right) = 0$, i.e., condition
	\eqref{eq:cond-pre-homotopy-upgrade-homotopy-inner-product}.
	In this case, $\phi_2$ can be upgraded to a full homotopy inner product
	by setting $\phi_{\ul{1}} = 0$.
	This happens for example if $\mu_0 \left( 1 \right) = 0$, or, if $\mathcal{A}$ and $\theta$
	are unital and $\mu_0 \left( 1 \right) = c \cdot e$ for some $c \in R^2$.
	See also \cref{rem:pre-homotopy-inn-prod-to-trace-zero-term}.
\end{rem}

\subsection{Recovering the Classical Chern--Simons Form and Action}
\label{sec:recovering-chern-simons}
Assume that $\mathcal{A} = \left( A, \mu \right)$ is a Banach $\Ainf$-algebra
which corresponds to a differential graded algebra. In this case,
the line degree two component $\G{b}[2]$ of the cyclic Chern--Simons form
is given by
\begin{equation*}
	\G{b}[2]                           =
	{\underbrace{ \frac{1}{2} \ul{\mur[\mu][b]{1}} \, \ul{b} + \frac{1}{3} \ul{\mur[\mu][b]{2}} \, \ul{b}}_{\eqcl{ \G{b}[2] }_{(2)}}}
	+ {\text{terms of weight} \geq 3}
\end{equation*}
for $b \in \tc{A}$.
When $\mathcal{A}$ is equipped with a strict pre-$\infty$-trace $\theta$, the
associated pre-superpotential function
$\SP[b][>0] = \theta \left( \psi \left( \G{b}[2] \right) \right)$ reduces to the
finite sum \eqref{eq:pre-sp-strict-trace-dga} and
makes sense for any $b \in A^0$, without assuming that $b$ is topologically nilpotent.
Although $\G{b}[2]$ is not defined for arbitrary $b \in A^0$, the pre-superpotential
$\SP[b][>0]$ depends only on the weight-two component $\eqcl{ \G{b}[2] }_{(2)}$
of $\G{b}[2]$, which is well-defined.
Since $\psi$ is weight homogeneous of degree minus one, we have
$\psi ( \eqcl{ \G{b}[2] }_{(2)} ) = \eqcl{ \psi \left( \G{b}[2] \right) }_{(1)}$,
and we can write
\begin{equation*}
	\SP[b][>0] =
	\left( \theta \circ \psi \right) \left( \eqcl{ \G{b}[2] }_{(2)} \right) =
	\theta \left( \eqcl{ \psi \left( \G{b}[2] \right) }_{(1)} \right),
\end{equation*}
avoiding the infinite sum of $\G{b}[2]$. Explicitly,
\cref{eq:psi-ula-ulb-dga} gives the formula
\begin{equation} \label{eq:psi-G-b-2-w-1}
	\begin{aligned}
		\eqcl{ \psi \left( \G{b}[2] \right) }_{(1)}
		={} &
		\frac{1}{2} \left(
		\frac{1}{2} \left(
		\mu_2 \left( \mur[\mu][b]{1}, b \right) -
		\mu_2 \left( b, \mur[\mu][b]{1} \right)
		\right)
		\right)
		\\
		    & +
		\frac{1}{2} \left(
		\frac{1}{3} \left(
		\mu_2 \left( \mur[\mu][b]{2}, b \right) -
		\mu_2 \left( b, \mur[\mu][b]{2} \right)
		\right)
		\right)
	\end{aligned}
\end{equation}
for $\eqcl{ \psi \left( \G{b}[2] \right) }_{(1)}$.
In what follows, we will
relate $\eqcl{ \psi \left( \G{b}[2] \right) }_{(1)}$ and $\SP[b][>0]$
to the classical Chern--Simons form and action.

Let $M$ be a closed oriented $n$-dimensional manifold, and let $E$ be a vector bundle
over $M$ endowed with a flat connection $\nabla$.
Let $B = \df{M}[\End{E}]$ be the differential graded algebra
of $\End{E}$-valued differential forms on $M$. The
product $m_2 \left( \omega_1, \omega_2 \right) = \omega_1 \wedge \omega_2$
on $B$ combines both the product on differential forms
and the composition of $\End{E}$, and there is also a module
action $\wedge \colon \df{M} \otimes_{\RR} \df{M}[\End{E}] \rightarrow \df{M}[\End{E}]$
compatible with the product on $B$.
The differential
$m_1 \left( \omega \right) = D \left( \omega \right)$ on $B$
is the exterior covariant
derivative, constructed from the exterior derivative $d$ on $\df{M}$
and the induced connection from $\nabla$ on $\End{E}$.
The operator $D$ is a graded module derivation over $d$ of degree one, and since
$\nabla$ is flat, we have $D^2 = 0$. Thus, we can think of
$\mathcal{B} = \left( B, m_1, m_2 \right)$ as
a differential graded algebra over both $\RR$ and
$\mathcal{S} = \left( \df{M}, d \right)$.

Let $\tr \colon \df{M}[\End{E}] \rightarrow \df{M}$ be the
pointwise trace.
Denote by $b^1$ and $b^2$ the two (commuting) Hochschild differentials
associated to $D$ and $\wedge$ respectively, working in the unshifted framework
(see \cref{sec:signs-hoch-unshifted}).
The wedge product on $\End{E}$-valued forms is not graded-commutative, but
since the standard trace on matrices satisfies $\tr[AB] = \tr[BA]$, and
the wedge product on forms is graded-commutative,
we have
\begin{equation*}
	\left( \tr \circ b^2 \right) \left( \omega_1, \omega_2 \right) =
	\tr[\omega_1 \wedge \omega_2] - (-1)^{\degb{\omega_1} \cdot \degb{\omega_2}}
	\tr[\omega_2 \wedge \omega_1] = 0.
\end{equation*}
Since taking trace is natural, we have
\begin{equation*}
	\left( \tr \circ b^1 \right) \left( \omega \right) =
	\tr[D\omega] = d\tr[\omega] = \left( d \circ \tr \right) \left( \omega \right).
\end{equation*}
Thus, $\tr \colon \mathcal{B} \rightarrow \mathcal{S}$ is a trace on $\mathcal{B}$
over $\mathcal{S}$.

Let $\theta \colon \df{M}[\End{E}] \rightharpoonup \RR$ be defined by
$\theta \left( \omega \right) = \int_M \tr[\omega]$.
By Stokes' theorem, the integration map $\int_M$ is a chain map,
and hence the composition $\theta = \int_M \circ \tr$ also satisfies
$\theta \circ b^1 = 0$. The condition $\theta \circ b^2 = 0$ follows
from $\tr \circ b^2 = 0$, so $\theta \colon \mathcal{B} \rightharpoonup \RR$
is also a trace on $\mathcal{B}$, working over $\RR$.

Now let $\mathcal{A} = \left( A, \mu_1, \mu_2 \right)$, with $A = B[1]$,
be the shifted $\Ainf$-algebra corresponding to $\left( B, D, \wedge \right)$
with the conventions of \eqref{eq:ainf_for_m_k_explicit}.
Given $\omega \in B^i$,
we use the notation $\bm{\omega}$ or $\s \omega$ to denote the corresponding
(same) element of $B[1]^{i-1}$.
Let $b \in B^1 = \df{M}[\End{E}][][1]$ and
let $\bm{b} \in A^0$ be the corresponding element in the shifted $\Ainf$-algebra.
Then we have
\begin{align*}
	\eqcl{\psi \left( \G{\bm{b}}[2] \right)}_{(1)}
	\stackrel{\eqref{eq:psi-G-b-2-w-1}}{=}{}           &
	\frac{1}{2} \left(
	\frac{1}{2} \left(
	\mu_2 \left( \mur[\mu][\bm{b}]{1}, \bm{b} \right) -
	\mu_2 \left( \bm{b}, \mur[\mu][\bm{b}]{1} \right)
	\right)
	\right)
	\\
	                                                   & +
	\frac{1}{2} \left(
	\frac{1}{3} \left(
	\mu_2 \left( \mur[\mu][\bm{b}]{2}, \bm{b} \right) -
	\mu_2 \left( \bm{b}, \mur[\mu][\bm{b}]{2} \right)
	\right)
	\right)
	\\
	\stackrel{\phantom{\eqref{eq:psi-G-b-2-w-1}}}{=}{} &
	\frac{1}{2} \s \left(
	\frac{1}{2} \left( Db \wedge b + b \wedge Db \right) -
	\frac{2}{3} b \wedge b \wedge b
	\right).
\end{align*}

Define $\trb \colon \ncdfr{\mathcal{A} / \mathcal{S}}[0][] \rightarrow S$
by $\trb[\bm{\omega}] = \tr[\omega]$.
Then $\trb$ is a zero-dimensional strict pre-$\infty$-trace
on $\mathcal{A}$ over $\mathcal{S}$ in the sense of \eqref{eq:infty-trace-strict-rel}.
By \cref{cor:strict-theta-cdga-to-strict-phi},
the corresponding strict pre-homotopy inner product is given by
\begin{equation*}
	\phi_2 \left( \ul{\bm{\omega_1}}, \ul{\bm{\omega_2}} \right) =
	\trb \left( \mu_2 \left( \bm{\omega_1}, \bm{\omega_2} \right) \right) =
	(-1)^{\degb{\omega_1} + 1} \tr[\omega_1 \wedge \omega_2],
\end{equation*}
and by \cref{eq:pre-sp-strict-trace-dga}, the pre-superpotential is given by
\begin{equation} \label{eq:sp-cs-form}
	\begin{aligned}
		{\SP}^{\mathcal{A}/S}_{>0} \left( -\bm{b} \right) ={} &
		\trb \left( \eqcl{\psi \left( \G{-\bm{b}}[2] \right)}_{(1)} \right)
		=
		\trb \left(
		\frac{1}{2} \mu_2 \left( \mu_1 \left( \bm{b} \right), \bm{b} \right) -
		\frac{1}{3} \mu_2 \left( \mu_2 \left( \bm{b}, \bm{b} \right), \bm{b} \right)
		\right)
		\\
		={}                                                   &
		\frac{1}{2} \tr[Db \wedge b + \frac{2}{3} b \wedge b \wedge b]
		\in \df{M}[][][3].
	\end{aligned}
\end{equation}

Similarly, we can define $\bm{\theta} \colon \ncdfr{\mathcal{A} / \RR}[0][] \rightharpoonup \RR$
by $\bm{\theta} \left( \bm{\omega} \right) = \theta \left( \omega \right)$.
Then $\bm{\theta}$ is an $n$-dimensional strict pre-$\infty$-trace on
$\mathcal{A}$ over $\RR$ in the sense of \eqref{eq:infty-trace-strict-rel},
whose corresponding pre-superpotential is given by
\begin{equation} \label{eq:sp-cs-action}
	{\SP}^{\mathcal{A}/\RR}_{>0} \left( -\bm{b} \right) =
	\int_{M} {\SP}^{\mathcal{A}/S}_{>0} \left( -\bm{b} \right) =
	\frac{1}{2}\int_M \tr[Db \wedge b + \frac{2}{3} b \wedge b \wedge b] \in \RR.
\end{equation}
When $\dim M = n = 3$, \cref{eq:sp-cs-form} gives the Chern--Simons form on $M$
and \cref{eq:sp-cs-action} gives the Chern--Simons action up to a factor of $\frac{1}{2}$.
The term $\eqcl{\psi \left( \G{-\bm{b}}[2] \right)}_{(1)}$ gives $\frac{1}{2}$ of the $\End{E}$-valued
three form $\frac{1}{2} \left( Db \wedge b + b \wedge Db \right) + \frac{2}{3} b \wedge b \wedge b$,
from which the classical Chern--Simons form
$\tr[Db \wedge b + \frac{2}{3} b \wedge b \wedge b]$ is obtained by applying the trace.

Note that in this setting with no curvature, everything above could have been
stated just as well for the superpotential, dropping all the ``pre'' adjectives.

\subsection{Derivative of the Superpotential and the \texorpdfstring{$\infty$}{Infinity}-Modulus} \label{sec:derivative-and-trace}

Let $R$ be a graded-commutative Banach $\mathbbm{k}$-algebra and let
$\mathcal{A}$ be a Banach $\Ainf$-algebra over $R$ equipped with
an $n$-dimensional homotopy inner product.
Assume that $\mathcal{B} = \left( B, \mu, \phi \right)$ is obtained from $\mathcal{A}$ by
scalar extension along $R \rightarrow \pows{R}[t]$, where $t$ is an even formal variable.
\Cref{thm:superpotential-derivative}, proven in \cref{sec:derivative-calc}, shows that the
formal derivative of the superpotential
\begin{equation*}
	\SP[b] = \phi_{\ul{1}} +
	\sum_{i,j,k=0}^{\infty} \frac{1}{i + j + 1 + k} \phi_2^{j,k} \left(
	\ul{ \mu_{i} \left( b^{i} \right)}, b^{j}, \ul{b}, b^{k}
	\right)
\end{equation*}
associated to the homotopy inner product $\phi$ is given by
\begin{equation} \label{eq:derivative-sp-homotopy-inner-product}
	\partial_t \, \SP[b] =
	\sum_{i,j,k=0}^{\infty} \phi_2^{j,k} \left(
	\ul{\mu_i \left( b^{i} \right)}, b^{j}, \ul{\dot{b}}, b^{k}
	\right),
\end{equation}
where $b = b \left( t \right) \in \tc{B}$ and $\dot{b} = \partial_t \, b \left( t \right)$.
When $b \in \mc{\mathcal{B}}$ is a strong
bounding cochain, \cref{eq:derivative-sp-homotopy-inner-product} implies that
$\partial_t \, \SP[b] = 0$ so that
strong bounding cochains are critical points of the superpotential function.

Now assume that $\mathcal{A}$ is unital. In what follows, we give a relation between
the derivative of the superpotential function at \textit{weak bounding cochains} and the derivative
of the pre-$\infty$-modulus function \eqref{eq:pre-infinity-modulus-function},
introduced in \cref{sec:traces-modulus-invariance}.

\begin{thm}[\cref{thm:derivative-sp-multiple-modulus} of the Introduction] \label{thm:derivative-sp-homotopy-inner-product-trace}
	Let $R$ be a graded-commutative Banach $\mathbbm{k}$-algebra and let
	$\mathcal{A}$ be a unital Banach $\Ainf$-algebra over $R$ equipped
	with a strongly unital homotopy inner product.
	Assume that $\mathcal{B} = \left( B, \mu, e, \phi \right)$ is obtained from $\mathcal{A}$ by
	scalar extension along $R \rightarrow \pows{R}[t]$, where $t$ is an even formal variable
	with $\nnorm[t] < 1$.

	Let $\theta = \phi_2 \circ \psi'$ be the pre-$\infty$-trace corresponding to the pre-homotopy inner product $\phi_2$,
	and let $b = b(t) \in \tc{B}$.
	Then
	\begin{enumerate}
		\item The pre-$\infty$-modulus function $\TM{}[>0]$ is expressed in terms of $\phi_2$ by
		      \begin{equation} \label{eq:tm-b-in-terms-of-phi-2}
			      \TM{b}[> 0] = \sum_{j,k=0}^{\infty} \frac{1}{j+1+k} \phi_2^{j,k} \left( \ul{e}, b^j, \ul{b}, b^k \right).
		      \end{equation}
		\item We have
		      \begin{equation} \label{eq:partial-tm-b-in-terms-of-phi-2}
			      \partial_t \, \TM{b}[> 0] = \sum_{j,k=0}^{\infty} \phi_2^{j,k} \left( \ul{e}, b^j, \ul{\dot{b}}, b^k \right).
		      \end{equation}
		\item When $b \in \mc{\mathcal{B}}[c]$ is a weak bounding cochain, we have
		      \begin{equation} \label{eq:partial-sp-b-c-tm-b}
			      \partial_t \, \SP[b] = c \cdot \partial_t \, \TM{b}[>0].
		      \end{equation}
	\end{enumerate}
\end{thm}
\begin{proof}
	The scalar extension of a unital Banach $\Ainf$-algebra is unital, and the scalar
	extension of total inner products described in \cref{sec:scalar-extension-total-inner-products}
	is readily seen to respect the unitality conditions, so that
	$\phi$ is a strongly unital homotopy inner product on $\mathcal{B}$.

	Given $a \in B$ and $b \in \tc{B}$, \cref{eq:psi-tag-conv-form} for $\psi'$ implies the identities
	\begin{gather}
		\psi' \left( \sum_{i=0}^{\infty} \frac{1}{i+1} b^{i+1} \right) =
		\sum_{j,k=0}^{\infty} \frac{1}{j + 1 + k} \ul{e} \, b^j \, \ul{b} \, b^k
		\mod \left< *e* \right>,
		\label{eq:psi-sum-b-k-identity}
		\\
		\psi' \left( \sum_{i=0}^{\infty} a \, b^i \right) \equiv (-1)^{\degb{a} + 1}
		\sum_{j,k = 0}^{\infty} \ul{a} \, b^j \, \ul{e} \, b^k
		=
		\sum_{j,k = 0}^{\infty} \ul{e} \, b^j \, \ul{a} \, b^k
		\mod \left< *e* \right>.
		\label{eq:psi'-sum-a-b-k-identity}
	\end{gather}
	Since $\phi_2$ is strongly unital, in particular $\degen$-unital, we have
	\begin{equation*}
		\begin{aligned}
			\TM{b}[> 0]
			\stackrel{\eqref{eq:pre-infinity-modulus-function}}{=} &
			\sum_{i = 1}^{\infty} \frac{1}{i} \theta_i \left( b^i \right) =
			\theta \left( \sum_{i = 0}^{\infty} \frac{1}{i+1} b^{i+1} \right)
			\stackrel{\eqref{eq:naive-trace-from-hip}}{=}
			\left( \phi_2 \circ \psi' \right)  \left(  \sum_{i = 0}^{\infty} \frac{1}{i+1} b^{i+1} \right)
			\\
			\stackrel{\eqref{eq:psi-sum-b-k-identity}}{=}          &
			\sum_{j,k=0}^{\infty} \frac{1}{j+1+k} \phi_2^{j,k} \left( \ul{e}, b^j, \ul{b}, b^k \right),
		\end{aligned}
	\end{equation*}
	which shows \cref{eq:tm-b-in-terms-of-phi-2}. Similarly, we have
	\begin{equation*}
		\begin{aligned}
			\partial_t \, \TM{b}[> 0]
			\stackrel{\eqref{eq:pre-infinity-modulus-function}}{=}{} &
			\partial_t \left( \sum_{i = 1}^{\infty} \frac{1}{i} \theta_i \left( b^{i} \right) \right)
			=
			\theta \left( \sum_{i = 0}^{\infty} \dot{b} \, b^i \right)
			\stackrel{\eqref{eq:naive-trace-from-hip}}{=}
			\left( \phi_2 \circ \psi' \right) \left( \sum_{i = 0}^{\infty} \dot{b} \, b^i \right)
			\\
			\stackrel{\eqref{eq:psi'-sum-a-b-k-identity}}{=}{}       &
			\sum_{j,k = 0}^{\infty} \phi_2^{j,k} \left( \ul{e}, b^j, \ul{\dot{b}}, b^k \right),
		\end{aligned}
	\end{equation*}
	which shows \cref{eq:partial-tm-b-in-terms-of-phi-2}.

	Finally, when $b \in \mc{\mathcal{B}}[c]$ is a weak bounding cochain, the formula \eqref{eq:derivative-sp-homotopy-inner-product}
	for the derivative $\partial_t \, \SP[b]$ simplifies to
	\begin{equation*}
		\partial_t \, \SP[b]
		=
		\sum_{j,k = 0}^{\infty} \phi_2^{j,k} \left( \ul{c \cdot e}, b^{j}, \ul{\dot{b}}, b^{k} \right)
		=
		c \cdot \sum_{j,k = 0}^{\infty} \phi_2^{j,k} \left( \ul{e}, b^{j}, \ul{\dot{b}}, b^{k} \right)
		\stackrel{\eqref{eq:partial-tm-b-in-terms-of-phi-2}}{=}
		c \cdot \partial_t \, \TM{b}[>0],
	\end{equation*}
	as $c$ is even, which shows \cref{eq:partial-sp-b-c-tm-b}.
\end{proof}

\begin{rem}
	The proof of \cref{thm:derivative-sp-homotopy-inner-product-trace} shows that given a $\degen$-unital pre-homotopy inner product
	$\phi_2$ and a pre-$\infty$-trace $\theta$, which are related via $\theta = \phi_2 \circ \psi'$, and given a weak bounding cochain
	$b \in \mc{\mathcal{B}}[c]$, we have the identity
	\begin{equation*}
		\partial_t \, \SP[b][> 0] = c \cdot \partial_t \, \TM{b}[> 0]
	\end{equation*}
	between the associated pre-superpotential and pre-$\infty$-modulus functions.
\end{rem}

\begin{rem}
	\Cref{eq:tm-b-in-terms-of-phi-2} shows that our pre-$\infty$-modulus function $\TM{b}[> 0]$, when
	expressed in terms of $\phi_2$, coincides with the potential $\Psi$, introduced by
	Cho and Lee in \cite[Definition 1.4]{cho-homotopy-superpotential}.
\end{rem}

\subsection{The Superpotential and the Periodicity Operator} \label{sec:sp-and-periodicity}
As a final application, we present a relation between the truncated cyclic exponential
\begin{equation*}
	\Go{b}[0] = \sum_{k \geq 1} \frac{1}{k} b^{k} \in \ncdfr{A}[0][0]
\end{equation*}
of \cref{subsec:canonical-elements-ncdf-0}, the line degree two component
\eqref{eq:G-b-2-braidop-1} of the cyclic Chern--Simons form
\begin{equation*}
	\G{b}[2] = \sum_{i,j,k=0}^{\infty} \frac{1}{i + j + 1 + k}
	\ul{ \mu_{i} \left( b^{i} \right)} \, b^{j} \, \ul{b} \, b^{k} \in
	\ncdf{A}[2][] / \Im \left( \qdr^3 \right),
\end{equation*}
considered as an element in the quotient
$\ncdf{A}[2][] / \Im \left( \qdr^3 \right)$,
and the periodicity operator $S$ on Connes' complex. We deduce a connection
between the pre-superpotential function, the pre-$\infty$-trace and the periodicity
operator.

\begin{thm}
	Let $\mathcal{A} = \left( A, \mu \right)$ be a Banach $\Ainf$-algebra over
	a differential graded-commutative Banach $\mathbbm{k}$-algebra $\mathcal{R} = \left( R, d \right)$.
	Let $\theta \colon \ncdfr{\mathcal{A}}[0][] \rightarrow \mathcal{R}[1-n]$ be a pre-$\infty$-trace on $\mathcal{A}$
	and let $\phi_2 = \theta \circ \psi \colon \ncdf{\mathcal{A}}[2][] / \Im \left( \qdr^3 \right) \rightarrow \mathcal{R}[2-n]$
	be the corresponding pre-homotopy inner product.
	Given a strong bounding cochain $b \in \mc{\mathcal{A}}$, we have the relations
	\begin{align}
		\psi \left( \G{b}[2] \right)
		           & \equiv S \left( \Go{b}[0] \right) \mod
		\clie{\mu} \left( \ncdfr{A}[0][] \right),
		\label{eq:sp-psi-Gb-minus-S}
		\\
		\SP[b][>0] & \equiv \theta \left( S \left( \Go{b}[0] \right) \right) \mod
		dR.
		\label{eq:sp-b-minus-theta-S-Go}
	\end{align}
	In particular, when $\mu_0 \left( 1 \right) = 0$, we get
	\begin{align*}
		\eqcl{ \psi \left( \G{b}[2] \right) } & = \eqcl{ S \left( \Go{b}[0] \right) }
		\in \hcyc{\mathcal{A}}[2],
		\\
		\eqcl{ \SP[b][>0] }                   & = \eqcl{ \theta \left( S \left( \Go{b}[0] \right) \right) }
		\in \cohom{\mathcal{R}}[3-n].
	\end{align*}
\end{thm}
\begin{proof}
	Given $b \in \tc{A}$, consider the chain
	$\Go{b}[\geq 1] = \Go{b}[\geq 1][\braidop_1] \in \totcomp{\mathcal{A}}[1][-1]$ given by
	\begin{equation*}
		\Go{b}[\geq 1] =
		\sum_{\substack{k=1 \\ i_1,\dots,i_{k-1} = 0 \\ j_1,\dots,j_k=0}}^{\infty}
		\frac{(-1)^{\frac{(k - 2) \cdot (k - 1)}{2}}}{1 + \sum_{r=1}^{k-1} i_r + \sum_{r = 1}^k j_r}
		\s_k \left(
		\ul{ \mu_{i_1} \left( b^{i_1} \right)} \, b^{j_1} \, \dots
		\, \ul{ \mu_{i_{k-1}} \left( b^{i_{k-1}} \right) } \, b^{j_{k-1}}
		\, \ul{b} \, b^{j_k}
		\right).
	\end{equation*}
	The chain $\Go{b}[\geq 1]$ is the image of the chain
	$\Go{b}[\geq 1][\braidop_2]$, defined in \cref{eq:Go-geq-1-def},
	under the isomorphism $\Psi_{\geq 1}$, converting
	constructions done using $\braidop_2$ to $\braidop_1$.
	(see \cref{sec:description-using-braid-op-1}).

	We have the formula
	\begin{equation*}
		\begin{aligned}
			D_{\mathcal{A}} \left( \Go{b}[\geq 1][\braidop_1] \right)
			\eqwithref                                                               &
			\Psi_{\geq 1} \left( D_{\mathcal{A}} \left( \Go{b}[\geq 1][\braidop_2] \right) \right)
			\\
			\eqwithref[eq:D-Go-geq-1-linear]                                         &
			\Psi_{\geq 1} \left( \Ho{b}[\geq 1][\braidop_2] - \Ho{0}[\geq 1][\braidop_2] \right)
			\\
			\eqwithref[eq:H_n(b)-formula][eq:Psi-explicit-formula-in-generalized-sp] &
			\sum_{k \geq 1} (-1)^{\frac{k(k+1)}{2}}
			\frac{1}{k} \s_k \left( \left( \ul{\corest{\mu} \left( \Exp{b} \right)} \, \Exp{b} \right)^{k} - \ulz{\mu}^{k} \right).
		\end{aligned}
	\end{equation*}
	In particular, when $b \in \mc{\mathcal{A}}$ is a strong bounding cochain, we have
	\begin{equation}
		D_{\mathcal{A}} \left( \Go{b}[\geq 1] \right) = \sum_{k \geq 1} (-1)^{\frac{k(k+1)}{2} + 1}
		\frac{1}{k} \s_k \left( \ulz{\mu}^{k} \right).
	\end{equation}

	By \cref{lm:p_1-homotopy-equivalence}, the maps $p_1 \colon \totcomp{\mathcal{A}}[1] \rightharpoonup
		\ncdfr{\mathcal{A}}[0][*+1]$ and $i_1 \colon \ncdfr{\mathcal{A}}[0] \rightharpoonup \totcomp{\mathcal{A}}[1][*-1]$
	induce a homotopy equivalence between $\totcomp{\mathcal{A}}[1][]$ and ${\ncdfr{\mathcal{A}}[0][]}[1]$,
	and there is a map $H_{+} \colon \totcomp{\mathcal{A}}[1][] \rightharpoonup \totcomp{\mathcal{A}}[1][]$ such that
	$i_1 \circ p_1 = \idd - \partial \left( H_{+} \right)$ with $H_{+} \circ i_1 = 0$. Note that
	by the explicit formula \eqref{eq:def-p_1'} for $i_1$, we have
	\begin{equation*}
		\begin{aligned}
			i_1 \left( \mu_0 \left( 1 \right) \right)
			 & =
			\sum_{k \geq 0} (-1)^{\frac{k(k-1)}{2}} \s_{k+1} \left( \left( h_{\dr} \clie{\mu} \right)^k \ulz{\mu} \right)
			=
			\sum_{k \geq 0} (-1)^{\frac{k(k-1)}{2}} \frac{1}{k + 1} \s_{k+1} \left( \ulz{\mu}^{k+1} \right)
			\\
			 & =
			\sum_{k \geq 1} (-1)^{\frac{(k-1)(k-2)}{2}} \frac{1}{k} \s_k \left( \ulz{\mu}^{k} \right)
			=
			\sum_{k \geq 1} (-1)^{\frac{k(k+1)}{2} + 1} \frac{1}{k} \s_k \left( \ulz{\mu}^{k} \right)
			=
			D_{\mathcal{A}} \left( \Go{b}[\geq 1] \right),
		\end{aligned}
	\end{equation*}
	which means that, when $b$ is a strong bounding cochain, the differential of $\Go{b}[\geq 1]$ lies in the image of the map
	$i_1$.
	Then
	\begin{equation} \label{eq:go-b-geq-1-rel-p-bar}
		\begin{aligned}
			\Go{b}[\geq 1] & = \left( i_1 \circ p_1 + \partial \left( H_{+} \right) \right) \left( \Go{b}[\geq 1] \right)
			\\
			               & = i_1 \left( \Go{b}[0] \right) + D_{\mathcal{A}} \left( H_{+} \left( \Go{b}[\geq 1] \right) \right) +
			\left( H_{+} \circ D_{\mathcal{A}} \right) \left( \Go{b}[\geq 1] \right)
			\\
			               & = i_1 \left( \Go{b}[0] \right) + D_{\mathcal{A}} \left( H_{+} \left( \Go{b}[\geq 1] \right) \right)
			+ \left( H_{+} \circ i_1 \right) \left( \mu_0 \left( 1 \right) \right)
			\\
			               & =
			i_1 \left( \Go{b}[0] \right) - D_{\mathcal{A}} \left( x \right)
		\end{aligned}
	\end{equation}
	for $x = -H_{+} \left( \Go{b}[\geq 1] \right) \in \totcomp{\mathcal{A}}[1][-2]$.

	We now do some diagram chasing based on the diagram
	\eqref{eq:diag-S-p-psi} of \cref{lm:periodicity-operator-as-projection}.
	We have
	\begin{equation*}
		\begin{aligned}
			\psi \left( \G{b}[2] \right)
			\stackrel{\phantom{\eqref{eq:go-b-geq-1-rel-p-bar}}}{=}{} &
			\left( \psi \circ p_2 \right) \left( \Go{b}[\geq 2] \right)
			\\
			\stackrel{\phantom{\eqref{eq:go-b-geq-1-rel-p-bar}}}{=}{} &
			\left( \psi \circ p_2 \circ \pi_{\geq 2} \right) \left( \Go{b}[\geq 1] \right)
			\\
			\stackrel{\eqref{eq:go-b-geq-1-rel-p-bar}}{=}{}           &
			\left( \psi \circ p_2 \circ \pi_{\geq 2} \circ i_1 \right) \left( \Go{b}[0] \right) -
			\left( \psi \circ p_2 \circ \pi_{\geq 2} \right) \left( D_{\mathcal{A}} \left( x \right) \right)
			\\
			\stackrel{\phantom{\eqref{eq:go-b-geq-1-rel-p-bar}}}{=}{} &
			S \left( \Go{b}[0] \right) +
			\clie{\mu} \left( \left( \psi \circ p_2 \circ \pi_{\geq 2} \right) \left( x \right) \right)
			\\
			\stackrel{\phantom{\eqref{eq:go-b-geq-1-rel-p-bar}}}{=}{} &
			S \left( \Go{b}[0] \right) +
			\clie{\mu} \left( y \right)
		\end{aligned}
	\end{equation*}
	for $y = \left( \psi \circ p_2 \circ \pi_{\geq 2} \right) \left( x \right) \in \ncdfr{\mathcal{A}}[0][1]$,
	which shows \cref{eq:sp-psi-Gb-minus-S}.
	Hence,
	\begin{equation*}
		\begin{aligned}
			\SP[b][>0] & = \phi_2 \left( \G{b}[2] \right) = \left( \theta \circ \psi \right) \left( \G{b}[2] \right) =
			\theta \left( S \left( \Go{b}[0] \right) \right) + \theta \left( \clie{\mu} \left( y \right) \right)
			\\
			           & =
			\theta \left( S \left( \Go{b}[0] \right) \right) + d_{\mathcal{R}} \left( (-1)^{1-n} \theta \left( y \right) \right),
		\end{aligned}
	\end{equation*}
	which shows \cref{eq:sp-b-minus-theta-S-Go}.

	When $\mu_0 \left( 1 \right) = 0$, both $\Go{b}[0]$ and $\G{b}[2]$ are closed,
	and $\psi$ and $S$ are chain maps, so we have equality between the cyclic homology classes
	in $\hcyc{\mathcal{A}}[2]$, and also between the cohomology classes in $\cohom{\mathcal{R}}[3-n]$.
\end{proof}

\clearpage
\appendix

\section{Category Theory Background} \label{appendix:cat-theory-background}

\subsection{Algebra Objects in a Monoidal Category} \label{subsec:alg-in-monoidal-cat}
Let $\mathcal{C} = \left( C, \otimes, \monunit \right)$ be a monoidal category.
A (\textbf{unital}) \textbf{algebra object} in $\mathcal{C}$ is an object $A$ of $C$ together with a multiplication morphism
$m \colon A \otimes A \rightarrow A$ and a unit morphism $u \colon \monunit \rightarrow A$
such that the following diagrams commute:

\begin{figure}[H]
	\centering
	\begin{tikzcd}
		{\left( A \otimes A \right) \otimes A} && {A \otimes \left( A \otimes A \right)} & {\monunit \otimes A} &
		{A \otimes A } & {A \otimes \monunit} \\
		{A \otimes A } && {A \otimes A } && A \\
		& A
		\arrow["m"', from=2-1, to=3-2]
		\arrow["m", from=2-3, to=3-2]
		\arrow["{m \otimes \id}"', from=1-1, to=2-1]
		\arrow["{\id \otimes m}", from=1-3, to=2-3]
		\arrow["{\alpha_{A,A,A}}","\cong"', from=1-1, to=1-3]
		\arrow["{\lambda_A}"', "\cong", from=1-4, to=2-5]
		\arrow["{u \otimes \id}", from=1-4, to=1-5]
		\arrow["{\id \otimes u}"', from=1-6, to=1-5]
		\arrow["{\rho_A}", "\cong"', from=1-6, to=2-5]
		\arrow["m", from=1-5, to=2-5]
	\end{tikzcd}
	\caption{Associativity and unitality of the multiplication $m$.}
	\label{fig:alg-assoc-unit}
\end{figure}

The maps $\alpha_{A,A,A}, \lambda_A$ and $\rho_A$ in \cref{fig:alg-assoc-unit} are the
associator, left unitor, and right unitor of the monoidal structure of $\mathcal{C}$.\footnote{
	The associator and unitor maps constitute part of the data of a monoidal category
	$\mathcal{C}$ but are usually suppressed in our notation.}
We note that the monoidal unit $\monunit$ admits a canonical structure of an algebra
with multiplication $m_{\monunit} = \lambda_{\monunit} = \rho_{\monunit}$ and unit $u_{\monunit} = \id$.
A (\textbf{unital}) \textbf{morphism} between two algebras $(A,m_A,u_A)$ and $(B,m_B,u_B)$ is a morphism $f \colon A \rightarrow B$ in $C$ which makes the following diagrams commute:

\begin{figure}[H]
	\begin{tikzcd}
		{A \otimes A} && A && A && B \\
		{B \otimes B} && B &&& \monunit
		\arrow["{m_A}", from=1-1, to=1-3]
		\arrow["f", from=1-3, to=2-3]
		\arrow["{f\otimes f}"', from=1-1, to=2-1]
		\arrow["{m_B}"', from=2-1, to=2-3]
		\arrow["f", from=1-5, to=1-7]
		\arrow["{u_B}"', from=2-6, to=1-7]
		\arrow["{u_A}", from=2-6, to=1-5]
	\end{tikzcd}
	\caption{Unital morphism of algebra objects.}
	\label{fig:alg-mor}
\end{figure}

\begin{rem}
	An algebra object in a monoidal category is more commonly called a  monoid object.
	Since we use these notions only in linear monoidal categories,
	we have decided to use the terminology of algebra objects as in \cite{Brandenburg2014}.
\end{rem}

Now assume the category $C$ is equipped with a symmetry $\sigma_{M,N} \colon M \otimes N \rightarrow N \otimes M$,
which makes $\mathcal{C} = \left( C, \otimes, \monunit, \sigma \right)$ into a symmetric monoidal category.
An algebra $A$ in a symmetric monoidal category $\left( C, \otimes, \monunit, \sigma \right)$
is said to be \textbf{commutative} (or $\sigma$-\textbf{commutative}, if we want to emphasize the role of
the symmetry $\sigma$) if the following diagram commutes:

\begin{figure}[H]
	\begin{tikzcd}
		{A \otimes A} && {A \otimes A } \\
		& A
		\arrow["{\sigma_{A,A}}", from=1-1, to=1-3]
		\arrow["m_A"', from=1-1, to=2-2]
		\arrow["m_A", from=1-3, to=2-2]
	\end{tikzcd}
	\caption{Commutative algebra object in a symmetric monoidal category.}
	\label{fig:alg-comm}
\end{figure}

We will denote the category of algebra objects in $\mathcal{C}$ by $\Alg[\mathcal{C}]$ and the full
subcategory of commutative algebras in $\mathcal{C}$ by $\CAlg[\mathcal{C}]$. Given two algebra objects
$A,B \in \Alg[\mathcal{C}]$, their \textbf{tensor product} is the object $A \otimes B$ of $C$ endowed
with the unit morphism $u_{A \otimes B}$ given by the diagram
\begin{figure}[H]
	\begin{tikzcd}[column sep=large]
		\monunit && {\monunit \otimes \monunit} && {A \otimes B}.
		\arrow["{\lambda_{\monunit} = \rho_{\monunit}}", from=1-1, to=1-3]
		\arrow["{u_A \otimes u_B}", from=1-3, to=1-5]
		\arrow["{u_{A \otimes B}}", curve={height=-25pt}, from=1-1, to=1-5]
	\end{tikzcd}
	\caption{Unit morphism of the tensor product of algebras.}
	\label{fig:alg-tensor-product-unit}
\end{figure}
The multiplication morphism $m_{A \otimes B}$ of $A \otimes B$ is given by the diagram
\begin{figure}[H]
	\begin{tikzcd}
		{\left( A \otimes B \right) \otimes \left(A \otimes B \right)} &
		{\left( A \otimes \left( B \otimes A \right) \right) \otimes B} &&&
		{\left( A \otimes \left( A \otimes B \right) \right) \otimes B} \\
		{A \otimes B} &&&& {\left( A \otimes A \right) \otimes \left(B \otimes B \right).}
		\arrow["\cong", from=1-1, to=1-2]
		\arrow["{\left( \id_A \otimes \sigma_{B,A} \right) \otimes \id_B}", from=1-2, to=1-5]
		\arrow["{m_{A \otimes B}}"', from=1-1, to=2-1]
		\arrow["\cong", from=1-5, to=2-5]
		\arrow["{m_A \otimes m_B}"', from=2-5, to=2-1]
	\end{tikzcd}
	\caption{Multiplication morphism of the tensor product of algebras.}
	\label{fig:alg-tensor-product-mult}
\end{figure}
With the definitions above, the symmetry morphism $\sigma \colon A \otimes B \rightarrow B \otimes A$ becomes a morphism of algebras, the tensor product of commutative algebras is commutative and hence the categories
$\Alg[\mathcal{C}]$ and $\CAlg[\mathcal{C}]$ become symmetric monoidal.

Let $\mathcal{C}_1 = \left( C_1, \otimes_1, \monunit_1 \right)$ and
$\mathcal{C}_2 = \left( C_2, \otimes_2, \monunit_2 \right)$ be two monoidal categories.
Let $\mathcal{F} \colon \mathcal{C}_1 \rightarrow \mathcal{C}_2$ be a lax monoidal functor with a unit
constraint
$\mu_0 \colon \monunit_2 \rightarrow \mathcal{F} \left( \monunit_1 \right)$ and tensor constraints
$\mu_{M,N} \colon \mathcal{F} \left( M \right) \otimes_2 \mathcal{F} \left( N \right) \rightarrow
	\mathcal{F} \left( M \otimes_1 N \right)$. A general fact about lax monoidal functors is that they
send algebra objects to algebra objects.  More precisely, given an algebra object
$\left( A, m_A, u_A \right)$ in $\mathcal{C}_1$, we can endow $\mathcal{F} \left( A \right)$
with an algebra structure by setting
$m_{\mathcal{F}(A)} \defeq \mathcal{F} \left( m_A \right) \circ \mu_{A,A}$ and
$u_{\mathcal{F}(A)} \defeq \mathcal{F} \left( u_A \right) \circ \mu_0$. Then
$\left( \mathcal{F} \left( A \right), m_{\mathcal{F} \left( A \right)}, u_{\mathcal{F} \left( A \right)} \right)$
is an algebra object in $\mathcal{C}_2$. The fact that
$\left( \mathcal{F} \left( A \right), m_{\mathcal{F} \left( A \right)}, u_{\mathcal{F} \left( A \right)} \right)$
satisfies the axioms of an algebra follows from the definition of a monoidal functor.

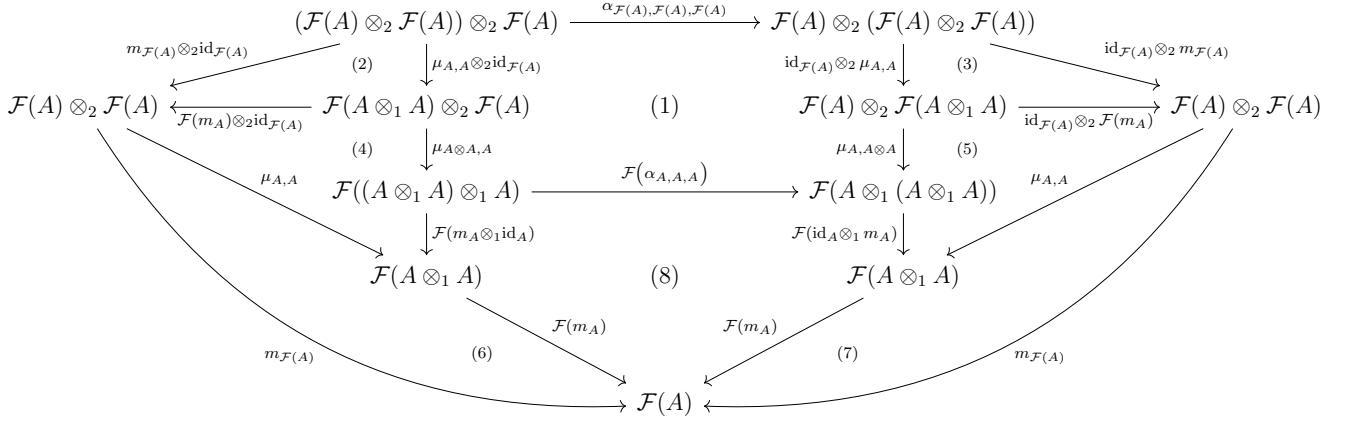
\begin{figure}[htb]
	\adjustbox{scale=0.76,center}{
		\begin{tikzcd}
			&& {\left( \mathcal{F} \left( A \right) \otimes_2 \mathcal{F} \left( A \right) \right) \otimes_2 \mathcal{F}
				\left( A \right)} &&
			{\mathcal{F} \left( A \right) \otimes_2 \left( \mathcal{F} \left( A \right) \otimes_2 \mathcal{F} \left( A
				\right) \right)}
			\\
			{\mathcal{F} \left( A \right) \otimes_2 \mathcal{F} \left( A \right)} && {\mathcal{F} \left( A \otimes_1 A \right) \otimes_2 \mathcal{F} \left( A \right)} & {(1)} & {\mathcal{F} \left( A \right) \otimes_2 \mathcal{F} \left( A \otimes_1 A \right)} && {\mathcal{F} \left( A \right) \otimes_2 \mathcal{F} \left( A \right)} \\
			&& {\mathcal{F} \left( \left( A \otimes_1 A \right) \otimes_1 A \right)} && {\mathcal{F} \left( A \otimes_1 \left( A \otimes_1 A \right) \right)} \\
			&& {\mathcal{F} \left( A \otimes_1 A \right)} & {(8)} & {\mathcal{F} \left( A \otimes_1 A \right)} \\
			\\
			&&& {\mathcal{F} \left( A \right)}
			\arrow["{\mathcal{F} \left( m_A \otimes_1 \id_A \right)}", from=3-3, to=4-3]
			\arrow["{\mathcal{F} \left( m_A \right)}", "(6)\qquad"', from=4-3, to=6-4]
			\arrow["{\mathcal{F} \left( m_A \right)}"', "\qquad(7)", from=4-5, to=6-4]
			\arrow["{\mathcal{F} \left( \id_A \otimes_1 \, m_A \right)}"', from=3-5, to=4-5]
			\arrow["{\mathcal{F} \left( \alpha_{A,A,A} \right)}", from=3-3, to=3-5]
			\arrow["{\mu_{A,A}}", from=2-1, to=4-3]
			\arrow["{m_{\mathcal{F} \left( A \right)}}"', bend right, from=2-1, to=6-4]
			\arrow["{\mu_{A,A}}"', from=2-7, to=4-5]
			\arrow["{m_{\mathcal{F} \left( A \right)}}", bend left, from=2-7, to=6-4]
			\arrow["{\mu_{A \otimes A, A}}", "(4)\qquad"', from=2-3, to=3-3]
			\arrow["{\mu_{A,A\otimes A}}"', "\qquad(5)", from=2-5, to=3-5]
			\arrow["{\mu_{A,A} \otimes_2 \id_{\mathcal{F} \left( A \right)}}", "(2)\qquad"', from=1-3, to=2-3]
			\arrow["{\id_{\mathcal{F} \left( A \right)} \otimes_2 \, \mu_{A,A} }"', "\qquad(3)",
				from=1-5, to=2-5]
			\arrow["{\alpha_{\mathcal{F} \left( A \right),\mathcal{F} \left( A \right),\mathcal{F} \left( A \right)}}", from=1-3, to=1-5]
			\arrow["{m_{\mathcal{F} \left( A \right)} \otimes_2 \id_{\mathcal{F} \left( A \right)}}"', from=1-3, to=2-1]
			\arrow["{\quad\id_{\mathcal{F} \left( A \right)} \otimes_2 \, m_{\mathcal{F} \left( A \right)}}", from=1-5, to=2-7]
			\arrow["{\mathcal{F} \left( m_A \right) \otimes_2 \id_{\mathcal{F} \left( A \right)}}", from=2-3, to=2-1]
			\arrow["{\id_{\mathcal{F} \left( A \right)} \otimes_2 \, \mathcal{F} \left( m_A \right)}"', from=2-5, to=2-7]
		\end{tikzcd}
	}
	\caption{Associativity of the induced multiplication on $\mathcal{F}(A)$.}
	\label{fig:assoc-induced-mult-lax-func}
\end{figure}

For example, to verify the associativity axiom for $\mathcal{F} \left( A \right)$, consider
the diagram in \cref{fig:assoc-induced-mult-lax-func}. Diagram $(1)$ commutes by definition
of a monoidal functor. Diagrams $(2)$ and $(3)$ commute by definition of
$m_{\mathcal{F}} \left( A \right)$ and the functoriality of $\otimes_2$.
Diagrams $(4)$ and $(5)$ commute by naturality of the transformations $\mu_{M,N}$,
which is part of the definition of a monoidal functor. Diagrams $(6)$ and $(7)$ commute by
definition of $m_{\mathcal{F} \left( A \right)}$. Finally, diagram $(8)$ is obtained by
applying $\mathcal{F}$ to the commutative diagram expressing the associativity of $m_A$ in
$\mathcal{C}_1$. The unitality of $\mathcal{F} \left( A \right)$ is proved similarly.

Given a morphism $f \colon A \rightarrow B$ between algebra objects of $\mathcal{C}_1$,
$\mathcal{F} \left( f \right) \colon \mathcal{F} \left( A \right) \rightarrow \mathcal{F} \left( B \right)$,
which is a priori a morphism in $C_2$, becomes a morphism between algebra objects of $\mathcal{C}_2$,
i.e., it respects the multiplications and units.
Hence, the construction described above gives us an induced functor
$\mathcal{F}_{\Alg} \colon \Alg[\mathcal{C}_1] \rightarrow \Alg[\mathcal{C}_2]$.

In addition, if the categories $\mathcal{C}_i$ are endowed with symmetries $\sigma_i$ which make them into
symmetric monoidal categories, and the functor $\mathcal{F}$ is symmetric lax monoidal, then the construction
above sends commutative algebras in $\mathcal{C}_1$ to commutative algebras in $\mathcal{C}_2$,
and we also obtain an induced functor
$\mathcal{F}_{\CAlg} \colon \CAlg[\mathcal{C}_1] \rightarrow \CAlg[\mathcal{C}_2]$. In this case,
both functors $\mathcal{F}_{\Alg}$ and $\mathcal{F}_{\CAlg}$ become symmetric lax monoidal
with the same coherence maps of $\mathcal{F}$.
When the functor $\mathcal{F}$ is a strong symmetric monoidal equivalence between $\mathcal{C}_1$ and
$\mathcal{C}_2$, the induced functors $\mathcal{F}_{\Alg}, \mathcal{F}_{\CAlg}$ are also strong symmetric monoidal equivalences.

\subsection{Modules over Algebra Objects}
\label{subsec:modules-in-monoidal-cat}
Let $\mathcal{C} = \left( C, \otimes, \monunit \right)$ be a monoidal category and let
$A \in \Alg[\mathcal{C}]$ be an algebra object. A \textbf{(unital left) module} over $A$ is an object
$M$ of $C$ together with a left action morphism $l_M \colon A \otimes M \rightarrow M$ such that the following
diagrams commute:
\begin{figure}[H]
	\centering
	\begin{tikzcd}
		{\left( A \otimes A \right) \otimes M} && {A \otimes \left( A \otimes M \right)} & {\monunit \otimes M} &
		{A \otimes M} \\
		{A \otimes M } && {A \otimes M} && M \\
		& M
		\arrow["{l_M}"', from=2-1, to=3-2]
		\arrow["{l_M}", from=2-3, to=3-2]
		\arrow["{m_A \otimes \id_M}"', from=1-1, to=2-1]
		\arrow["{\id_{A} \otimes \, l_M}", from=1-3, to=2-3]
		\arrow["{\alpha_{A,A,M}}", "\cong"', from=1-1, to=1-3]
		\arrow["{\lambda_M}"', from=1-4, to=2-5]
		\arrow["{u \otimes \id_M}", from=1-4, to=1-5]
		\arrow["{l_M}", from=1-5, to=2-5]
	\end{tikzcd}
	\caption{Associativity and unitality of the module action $l_M$.}
	\label{fig:module-assoc-unit}
\end{figure}

A morphism between two left $A$-modules $\left( M, l_M \right)$ and $\left( N, l_N \right)$ is a morphism
$f \colon M \rightarrow N$ in $C$ such that $f \circ l_M = l_N \circ \left( \id_A \otimes f \right)$.
We will denote the category of left $A$-modules by $\Mod[A]$, or by $\Mod[A][\mathcal{C}]$, when
we want to emphasize that $A$ is an algebra object of the monoidal category $\mathcal{C}$ and
that modules are objects of $C$.

Now assume $C$ is equipped with a symmetry $\sigma_{M,N} \colon M \otimes N \rightarrow N \otimes M$,
which makes the category $\mathcal{C} = \left( C, \otimes, \monunit, \sigma \right)$
into a symmetric monoidally cocomplete category.\footnote{This means that $\mathcal{C}$ is symmetric monoidal,
	the underlying category $C$ is cocomplete, and, for each $A \in C$, the endofunctor
	$A \otimes -$, hence also $- \otimes A$, is cocontinuous. This is the case for example if $\mathcal{C}$ is
	closed symmetric monoidal and $C$ is cocomplete.}
In this case, given a commutative algebra object $A \in \CAlg[\mathcal{C}]$, the category $\Mod[A][\mathcal{C}]$
can be endowed with a natural symmetric monoidal structure which makes $\Mod[A][\mathcal{C}]$
into a symmetric monoidally cocomplete category. The monoidal product of two $A$-modules
$M$ and $N$, denoted by $M \otimes_A N$, is given by
the coequalizer of the two action morphisms $M \otimes A \otimes N \rightrightarrows M \otimes N$,
endowed with an action of $A$, and the monoidal unit is given by $A$, considered
as a module object over itself (see \cite[Chapter 4]{Brandenburg2014} for details).

Let $\mathcal{C}_1 = \left( C_1, \otimes_1, \monunit_1 \right)$ and
$\mathcal{C}_2 = \left( C_2, \otimes_2, \monunit_2 \right)$ be two monoidal categories. Let
$\mathcal{F} \colon \mathcal{C}_1 \rightarrow \mathcal{C}_2$ be a lax monoidal functor with a unit constraint
$\mu_0 \colon \monunit_2 \rightarrow \mathcal{F} \left( \monunit_1 \right)$ and tensor constraints
$\mu_{M,N} \colon \mathcal{F} \left( M \right) \otimes_2 \mathcal{F} \left( N \right)
	\rightarrow \mathcal{F} \left( M \otimes_1 N \right)$.
Given an algebra object $A_1 \in \Alg[\mathcal{C}_1]$,
denote the algebra object $\mathcal{F}_{\Alg} \left( A_1 \right) \in \Alg[\mathcal{C}_2]$ by $A_2$.
Similarly to the situation with algebra objects, a lax monoidal functor $\mathcal{F}$ will also
send module objects over $A_1$ to module objects over $A_2$. More precisely,
given an $A_1$-module $\left( M, l_M \right)$ in $\mathcal{C}_1$,
we can endow $\mathcal{F} \left( M \right)$ with a $A_2$-module structure by setting
$l_{\mathcal{F}(M)} \defeq \mathcal{F} \left( l_M \right) \circ \mu_{A_1,M}$. One
can verify that $\mathcal{F} \left( M \right)$ then becomes a module object over $A_2$, and
that $\mathcal{F}$ maps morphisms of $A_1$-modules to morphisms of $A_2$-modules.
This construction gives us an induced functor $\mathcal{F}_{\Mod} \colon \Mod[A_1][\mathcal{C}_1] \rightarrow \Mod[A_2][\mathcal{C}_2]$.

When $\mathcal{C}_1$ and $\mathcal{C}_2$ are symmetric monoidally cocomplete categories,
the algebra $A_1$ is commutative, and $\mathcal{F} \colon \mathcal{C}_1 \rightarrow \mathcal{C}_2$
is a symmetric lax monoidal functor,
the induced functor $\mathcal{F}_{\Mod} \colon \Mod[A_1][\mathcal{C}_1] \rightarrow \Mod[A_2][\mathcal{C}_2]$
can be enhanced into a symmetric lax monoidal functor. The tensor constraints
\begin{equation*}
	\mu_{M,N}^{\mathcal{F}_{\Mod}} \colon
	\mathcal{F}_{\Mod} \left( M \right) \otimes_{A_2} \mathcal{F}_{\Mod} \left( N \right) \rightarrow
	\mathcal{F}_{\Mod} \left( M \otimes_{A_1} N \right)
\end{equation*}
of $\mathcal{F}_{\Mod}$ are induced from the tensor constraints
$\mu_{M,N}^{\mathcal{F}} \colon \mathcal{F} \left( M \right) \otimes_2 \mathcal{F} \left( N \right) \rightarrow \mathcal{F} \left( M \otimes_1 N \right)$ of $\mathcal{F}$,
and the unit constraint
$\mu_0^{\mathcal{F}_{\Mod}} \colon A_2 \rightarrow \mathcal{F}_{\Mod} \left( A_1 \right) = A_2$ is just the identity map.
When the functor $\mathcal{F}$ is a strong symmetric monoidal equivalence between
$\mathcal{C}_1$ and $\mathcal{C}_2$, then the induced functor
$\mathcal{F}_{\Mod} \colon \Mod[A_1][\mathcal{C}_1] \rightarrow
	\Mod[A_2][\mathcal{C}_2]$ is also a strong symmetric monoidal equivalence between the module categories.

\subsection{Coalgebra Objects in a Monoidal Category} \label{subsec:coalg-in-monoidal-cat}
Let $\mathcal{C} = \left( C, \otimes, \monunit \right)$ be a monoidal category. By reversing the arrows
in the diagrams of \cref{subsec:alg-in-monoidal-cat}, one can define a coalgebra object in a monoidal category
and the notion of a morphism of coalgebras.

A (\textbf{counital}) \textbf{coalgebra object} in $\mathcal{C}$ is an object $D$ of $C$ together with a comultiplication morphism
$\Delta \colon D \rightarrow D \otimes D$ and a counit morphism $\varepsilon \colon D \rightarrow \monunit$ such that the following diagrams commute:

\begin{figure}[H]
	\centering
	\begin{tikzcd}
		{\left( D \otimes D \right) \otimes D} && {D \otimes \left( D \otimes D \right)} & {\monunit \otimes D} &
		{D \otimes D} & {D \otimes \monunit} \\
		{D \otimes D } && {D \otimes D } && D  \\
		& D
		\arrow["\Delta", from=3-2, to=2-1]
		\arrow["\Delta"', from=3-2, to=2-3]
		\arrow["{\Delta \otimes \id}", from=2-1, to=1-1]
		\arrow["{\id \otimes \Delta}"', from=2-3, to=1-3]
		\arrow["{\alpha^{-1}_{D,D,D}}"', "\cong", from=1-3, to=1-1]
		\arrow["{\lambda_D^{-1}}", "\cong"', from=2-5, to=1-4]
		\arrow["{\varepsilon \otimes \id}"', from=1-5, to=1-4]
		\arrow["{\id \otimes \varepsilon}", from=1-5, to=1-6]
		\arrow["{\rho_D^{-1}}"',"\cong", from=2-5, to=1-6]
		\arrow["\Delta", from=2-5, to=1-5]
	\end{tikzcd}
	\caption{Coassociativity and counitality of the comultiplication $\Delta$.}
	\label{fig:coalg-coassoc-counit}
\end{figure}
We note that the monoidal unit $\monunit$ admits a canonical structure of a coalgebra
with comultiplication $\Delta_{\monunit} = \lambda_{\monunit}^{-1} = \rho_{\monunit}^{-1}$
and counit $\varepsilon_{\monunit} = \id$.
A (\textbf{counital}) \textbf{morphism} between two coalgebras $(D,\Delta_D,\varepsilon_D)$ and
$(E,\Delta_E,\varepsilon_E)$ is a morphism $f \colon D \rightarrow E$ in $C$ which makes the following diagrams commute:

\begin{figure}[H]
	\begin{tikzcd}
		{D \otimes D} && D && D && E \\
		{E \otimes E} && E &&& \monunit
		\arrow["{\Delta_D}"', from=1-3, to=1-1]
		\arrow["f", from=1-3, to=2-3]
		\arrow["{f\otimes f}"', from=1-1, to=2-1]
		\arrow["{\Delta_E}", from=2-3, to=2-1]
		\arrow["f", from=1-5, to=1-7]
		\arrow["{\varepsilon_E}", from=1-7, to=2-6]
		\arrow["{\varepsilon_D}"', from=1-5, to=2-6]
	\end{tikzcd}
	\caption{Counital morphism of coalgebra objects.}
	\label{fig:coalg-mor}
\end{figure}

We will denote the category of coalgebra objects in $\mathcal{C}$ by $\CoAlg[\mathcal{C}]$.
Given a counital coalgebra object $\left( D, \Delta, \varepsilon \right)$, we define the sequence
of iterated coproducts $\Delta^n \colon D \rightarrow D^{\otimes n}$ for $n \geq 0$ as follows:
\begin{enumerate}
	\item The map $\Delta^0 \colon D \rightarrow \monunit$ is the counit map $\varepsilon$.
	\item The map $\Delta^1 \colon D \rightarrow D$ is the identity map.
	\item The map $\Delta^2 \colon D \rightarrow D \otimes D$ is the coproduct $\Delta^2 \defeq \Delta$.
	\item The maps $\Delta^n \colon D \rightarrow D^{\otimes n}$ are defined inductively by
	      $\Delta^n \defeq \left( \Delta^{n-1} \otimes \id \right) \circ \Delta^2$ when $n > 2$.
\end{enumerate}
The coassociativity condition on $\Delta^2$ (left-hand side of \cref{fig:coalg-coassoc-counit})
says that
\begin{equation*}
	\left( \Delta^2 \otimes \id \right) \circ \Delta^2 = \left( \id \otimes \Delta^2 \right) \circ \Delta^2,
\end{equation*}
where we suppress the natural associativity isomorphism $\alpha_{D,D,D}$ used to identify both sides.
More generally,
for $l \geq 0$ and $k_1, \dots, k_l \geq 0$ we have the identity
\begin{equation}
	\left( \Delta^{k_1} \otimes \dots \otimes \Delta^{k_l} \right) \circ \Delta^l
	= \Delta^{k_1 + \dots + k_l},
	\label{eq:iterated-coproducts-coalgebra-identity}
\end{equation}
where we suppress the appropriate natural isomorphism, induced from the monoidal structure, used
to identify both sides. In particular, we obtain the identity:
\begin{equation*}
	\left( \Delta^n \otimes \Delta^m \right) \circ \Delta = \Delta^{n+m}.
	\qquad
\end{equation*}
Given a counital coalgebra morphism $f \colon D \rightarrow E$, one can show that $f$ commutes with
all the iterated coproducts in the sense that, for all $n \geq 0$, we have
\begin{equation}
	f^{\otimes n} \circ \Delta^n_D = \Delta^n_{E} \circ f.
	\label{eq:morphism-commutes-iterated-coproducts}
\end{equation}

Let $\mathcal{C}_1 = \left( C_1, \otimes_1, \monunit_1 \right)$ and
$\mathcal{C}_2 = \left( C_2, \otimes_2, \monunit_2 \right)$ be two monoidal categories. Let
$\mathcal{F} \colon \mathcal{C}_1 \rightarrow \mathcal{C}_2$ be an oplax monoidal functor
with a unit constraint $\mu_0 \colon \mathcal{F} \left( \monunit_1 \right) \rightarrow \monunit_2$ and
tensor constraints $\mu_{M,N} \colon \mathcal{F} \left( M \otimes_1 N \right) \rightarrow
	\mathcal{F} \left( M \right) \otimes_2 \mathcal{F} \left( N \right)$.
Dually to what we have described in
\cref{subsec:alg-in-monoidal-cat}, an oplax functor sends coalgebra objects in $\mathcal{C}_1$
to coalgebra objects in $\mathcal{C}_2$. More precisely, given a coalgebra object
$(D,\Delta_D,\varepsilon_D)$ in $\mathcal{C}_1$, we can endow $\mathcal{F}(D)$ with a coalgebra structure by setting
$\Delta_{\mathcal{F}(D)} \defeq \mu_{D,D} \circ \mathcal{F} \left( \Delta_D \right)$ and
$\varepsilon_{\mathcal{F}(D)} \defeq \mu_0 \circ \mathcal{F} \left( \varepsilon_D \right)$.
Then $\left( \mathcal{F}(D), \Delta_{\mathcal{F}(D)}, \varepsilon_{\mathcal{F}(D)} \right)$ is a coalgebra
object in $\mathcal{C}_2$.

Given a morphism $f \colon D \rightarrow E$ between coalgebra objects of $\mathcal{C}_1$,
$\mathcal{F} \left( f \right) \colon \mathcal{F} \left( D \right) \rightarrow \mathcal{F} \left( E \right)$,
which is a priori a morphism in $C_2$, becomes a morphism between coalgebra objects of
$\mathcal{C}_2$, i.e., it respects the comultiplication and counit maps.
Hence, the construction described above gives us an induced functor
$\mathcal{F}_{\CoAlg} \colon \CoAlg[\mathcal{C}_1] \rightarrow \CoAlg[\mathcal{C}_2]$.

\subsection{Comodules over Coalgebra Objects}
\label{subsec:comodules-in-monoidal-cat}
Let $\mathcal{C} = \left( C, \otimes, \monunit \right)$ be a monoidal category.
By reversing the arrows in the diagrams of \cref{subsec:modules-in-monoidal-cat}, one can define
the notion of a comodule over a coalgebra object in a monoidal category.

Let $D \in \CoAlg[\mathcal{C}]$ be a coalgebra object. A \textbf{(counital) left comodule} over $D$ is an object
$M$ of $C$ together with a left coaction morphism $\Delta^{1|0} \colon M \rightarrow D \otimes M$ such that the following diagrams commute:
\begin{figure}[H]
	\centering
	\begin{tikzcd}
		{\left( D \otimes D \right) \otimes M} && {D \otimes \left( D \otimes M \right)} & {\monunit \otimes M} &
		{D \otimes M} \\
		{D \otimes M } && {D \otimes M} && M \\
		& M
		\arrow["{\Delta_{M}^{1|0}}", from=3-2, to=2-1]
		\arrow["{\Delta_{M}^{1|0}}"', from=3-2, to=2-3]
		\arrow["{\Delta_D \otimes \id_M}", from=2-1, to=1-1]
		\arrow["{\id_{D} \otimes \Delta_{M}^{1|0}}"', from=2-3, to=1-3]
		\arrow["{\alpha^{-1}_{D,D,M}}"', "\cong", from=1-3, to=1-1]
		\arrow["{\lambda_M^{-1}}", from=2-5, to=1-4]
		\arrow["{\varepsilon \otimes \id_M}"', from=1-5, to=1-4]
		\arrow["{\Delta_{M}^{1|0}}"', from=2-5, to=1-5]
	\end{tikzcd}
	\caption{Coassociativity and counitality of the comodule coaction $\Delta^{1|0}$.}
	\label{fig:comodule-coassoc-counit}
\end{figure}

A morphism between two left $D$-comodules $( M, \Delta_{M}^{1|0} )$ and
$( N, \Delta_{N}^{1|0} )$ is a morphism
$f \colon M \rightarrow N$ in $C$ such that
$\Delta_{N}^{1|0} \circ f = \left( \id_D \otimes f \right) \circ \Delta_{M}^{1|0}$.
A comodule will be taken by default to mean counital left comodule unless stated otherwise.

Similarly, one can define a \textbf{right} $D$\textbf{-comodule}, which is an object $M$
of $C$, together with a right coaction morphism $\Delta^{0|1} \colon M \rightarrow M \otimes D$,
satisfying identities analogous to those in \cref{fig:comodule-coassoc-counit}. Given two
coalgebra objects $D$ and $E$, a $\left( D, E \right)$\textbf{-bicomodule} is an object
$M$ of $C$ which is simultaneously a left $D$-comodule and a right $E$-comodule
such that the following diagram commutes:
\begin{figure}[H]
	\centering
	\begin{tikzcd}
		{\left( D \otimes M \right) \otimes E} && {D \otimes \left( M \otimes E \right)}  \\
		{M \otimes E} && {D \otimes M} \\
		& M
		\arrow["{\Delta^{0|1}}", from=3-2, to=2-1]
		\arrow["{\Delta^{1|0}}"', from=3-2, to=2-3]
		\arrow["{\Delta^{1|0} \otimes \id_E}", from=2-1, to=1-1]
		\arrow["{\id_{D} \otimes \Delta^{0|1}}"', from=2-3, to=1-3]
		\arrow["{\alpha_{D,M,E}}", "\cong"', from=1-1, to=1-3]
	\end{tikzcd}
	\caption{Compatibility of left and right coactions in a bicomodule.}
	\label{fig:bicomodule-compatibility-condition}
\end{figure}

We note that any object $M$ of $C$ has a canonical structure of a $\left( \monunit, \monunit \right)$-bicomodule,
where $\monunit$ is endowed with the canonical coalgebra structure. The left (resp.\ right) coaction map on $M$
is given by the inverse of the left unitor $\lambda_{M}^{-1}$ (resp.\ inverse of the right unitor $\rho_{M}^{-1}$).
Similarly, any left $D$-comodule can be thought of as a $\left( D, \monunit \right)$-bicomodule,
and any right $D$-comodule can be thought of as a $\left( \monunit, D \right)$-bicomodule.

Given a $(D,E)$-bicomodule $M$ and $n,m \geq 0$, we can use the bicoaction maps $\Delta^{1|0}, \Delta^{0|1}$,
the coalgebra coproducts $\Delta_D, \Delta_E$, the counit maps $\varepsilon_D, \varepsilon_E$,
and the identity maps $\id_D$,$\id_M$ and $\id_E$, to construct higher coaction maps $\Delta^{m|n}$ of the form
$\Delta^{m|n} \colon M \rightarrow D^{\otimes m} \otimes M \otimes E^{\otimes n}$ with
$\Delta^{0|0} = \id_M$. This can be done in various different ways. Let us give a few examples:
\begin{enumerate}
	\item Consider the two maps
	      \begin{align*}
		       & M \xrightarrow{\Delta^{1|0}} D \otimes M \xrightarrow {\Delta_D \otimes \id_M}
		      \left( D \otimes D \right) \otimes M,
		      \\
		       & M \xrightarrow{\Delta^{1|0}} D \otimes M \xrightarrow {\id_D \otimes \Delta^{1|0}}
		      D \otimes \left( D \otimes M \right).
	      \end{align*}
	      The coassociativity condition for the left coaction map (see left-hand side of
	      \cref{fig:comodule-coassoc-counit}) tells us that both maps are the same when
	      we identify their codomains using the associativity isomorphisms of the monoidal structure.
	      Hence, we get a unique map $\Delta^{2|0} \colon M \rightarrow D^{\otimes 2} \otimes M$.
	\item Consider the two maps
	      \begin{align*}
		       & M \xrightarrow{\Delta^{0|0} = \id_M} M,
		      \\
		       & M \xrightarrow{\Delta^{1|0}} D \otimes M \xrightarrow{\varepsilon_D \otimes \id_M}
		      \monunit \otimes M.
	      \end{align*}
	      The counitality condition for the left coaction map (see right-hand side of
	      \cref{fig:comodule-coassoc-counit}) tells us that both maps are the same when
	      we identify their codomains using the left unitor isomorphisms of the monoidal structure.
	\item Consider the two maps
	      \begin{align*}
		       & M \xrightarrow{\Delta^{0|1}} M \otimes E \xrightarrow{\Delta^{1|0} \otimes \id_E}
		      \left( D \otimes M \right) \otimes E,                                                \\
		       & M \xrightarrow{\Delta^{1|0}} D \otimes M \xrightarrow{\id_D \otimes \Delta^{0|1}}
		      D \otimes \left( M \otimes E \right).
	      \end{align*}
	      The bicomodule compatibility condition (see \cref{fig:bicomodule-compatibility-condition}) tells
	      us that both maps are the same when we identify their codomains using the associativity isomorphism
	      of the monoidal structure. Hence, we get a unique map
	      $\Delta^{1|1} \colon M \rightarrow D \otimes M \otimes E$, called the \textbf{bicoaction map} on $M$.
	\item Consider the five maps
	      \begin{align*}
		       & M \xrightarrow{\Delta^{0|1}} M \otimes E \xrightarrow{\Delta^{1|0} \otimes \id_E}
		      \left( D \otimes M \right) \otimes E
		      \xrightarrow{\left( \Delta_D \otimes \id_M \right) \otimes \id_E}
		      \left( \left( D \otimes D \right) \otimes M \right) \otimes E
		      \\
		       & M \xrightarrow{\Delta^{0|1}} M \otimes E \xrightarrow{\Delta^{1|0} \otimes \id_E}
		      \left( D \otimes M \right) \otimes E
		      \xrightarrow{\left( \id_D \otimes \Delta^{1|0} \right) \otimes \id_E}
		      \left( D \otimes \left( D \otimes M \right) \right) \otimes E,
		      \\
		       & M \xrightarrow{\Delta^{1|0}} D \otimes M \xrightarrow{\id_D \otimes \Delta^{0|1}}
		      D \otimes \left( M \otimes E \right)
		      \xrightarrow{\Delta_D \otimes \left( \id_M \otimes \id_E \right)}
		      \left( D \otimes D \right) \otimes \left( M \otimes E \right),
		      \\
		       & M \xrightarrow{\Delta^{1|0}} D \otimes M \xrightarrow{\id_D \otimes \Delta^{0|1}}
		      D \otimes \left( M \otimes E \right)
		      \xrightarrow{\id_D \otimes \left( \Delta^{1|0} \otimes \id_E \right)}
		      D \otimes \left( \left( D \otimes M \right) \otimes E \right),
		      \\
		       & M \xrightarrow{\Delta^{1|0}} D \otimes M \xrightarrow{\Delta_D \otimes \id_M}
		      \left( D \otimes D \right) \otimes M
		      \xrightarrow{\left( \id_D \otimes \id_D \right) \otimes \Delta^{0|1}}
		      \left( D \otimes D \right) \otimes \left( M \otimes E \right).
	      \end{align*}
	      The maps have the same codomain up to parenthesization. Since in a monoidal category there
	      exists a \textit{unique} isomorphism induced by the monoidal structure between any
	      two choices of parenthesizations, we can identify the codomains of all the maps using the induced
	      isomorphisms and obtain maps with the same domain and codomain.
	      As a consequence of the axioms for coalgebras, comodules and the monoidal structure, one
	      can verify that all the maps above become the same map when identifying their codomains.
	      Hence, we get a unique map $\Delta^{2|1} \colon M \rightarrow D^{\otimes 2} \otimes M \otimes E$.
\end{enumerate}

More generally, it turns out that as a consequence of the coassociativity and counitality of
the comodules and coalgebras, and the axioms of a monoidal category, all different ways of
iteratively constructing the maps $\Delta^{m|n}$ from the basic building blocks will result
in the same map, up to an identification using the natural isomorphisms induced by the
monoidal structure. The identifications are usually implicitly understood and suppressed
from the notation. In particular, given a $\left( C, D \right)$-bicomodule $M$, we have the
identity
\begin{equation}
	\left( \Delta_{C}^{k_1} \otimes \dots \otimes \Delta_{C}^{k_m} \otimes \Delta_{M}^{r|s} \otimes
	\Delta_{D}^{l_1} \otimes \dots \otimes \Delta_{D}^{l_n} \right) \circ \Delta_{M}^{m|n} =
	\Delta_{M}^{(k_1 + \dots + k_m + r)|(s + l_1 + \dots + l_n)},
	\label{eq:higher-coproducts-bicomodule-identity}
\end{equation}
and there are analogous identities for left (resp.\ right) comodules obtained by thinking
of a left (resp.\ right) $D$-comodule as a $\left( D, \monunit \right)$-bicomodule (resp.\
$\left( \monunit, D \right)$-bicomodule).
The identity \eqref{eq:higher-coproducts-bicomodule-identity} is the analogue of
\cref{eq:iterated-coproducts-coalgebra-identity} for bicomodules
and reduces to \cref{eq:iterated-coproducts-coalgebra-identity} if one considers a coalgebra $D$ as
a $\left( D, D \right)$-bicomodule over itself in the natural way.

Given a morphism $f \colon M \rightarrow N$ of left $D$-comodules, one can show
that $f$ commutes with all the higher coaction maps $\Delta^{m|0}$ in the sense that,
for all $m \geq 0$, we have
\begin{equation}
	\left( \left( \id_D \right)^{\otimes m} \otimes f \right) \circ \Delta_M^{m|0} = \Delta_{N}^{m|0} \circ f.
	\label{eq:morphism-commutes-higher-coactions}
\end{equation}

\begin{rem}
	The structure of a $\left( D, E \right)$-bicomodule is encoded entirely in the single bicoaction map
	$\Delta^{1|1} \colon M \rightarrow D \otimes M \otimes E$.
	Given the bicoaction $\Delta^{1|1}$, the left coaction $\Delta^{1|0} \colon M \rightarrow D \otimes M$
	can be recovered using the counit of $E$ as
	\begin{equation*}
		\Delta^{1|0} = \left( \id_D \otimes \id_M \otimes \varepsilon_E \right) \circ \Delta^{1|1}
	\end{equation*}
	and similarly for the right coaction.
\end{rem}

Let $\mathcal{C}_1 = \left( C_1, \otimes_1, \monunit_1 \right)$ and
$\mathcal{C}_2 = \left( C_2, \otimes_2, \monunit_2 \right)$ be two monoidal categories.
Let $\mathcal{F} \colon \mathcal{C}_1 \rightarrow \mathcal{C}_2$ be an oplax monoidal functor
with a unit constraint $\mu_0 \colon \mathcal{F} \left( \monunit_1 \right) \rightarrow \monunit_2$ and
tensor constraints $\mu_{M,N} \colon \mathcal{F} \left( M \otimes_1 N \right) \rightarrow
	\mathcal{F} \left( M \right) \otimes_2 \mathcal{F} \left( N \right)$.
Given a coalgebra object $D_1 \in \CoAlg[\mathcal{C}_1]$,
denote the coalgebra object $\mathcal{F} \left( D_1 \right) \in \CoAlg[\mathcal{C}_2]$ by $D_2$.
Dually to what we have described in
\cref{subsec:modules-in-monoidal-cat}, an oplax monoidal functor $\mathcal{F}$ will send
comodules over $D_1$ to comodules over $D_2$.
More precisely, given a $D_1$-comodule $( M, \Delta_M^{1|0} )$ in $\mathcal{C}_1$,
we can endow $\mathcal{F} \left( M \right)$ with a $D_2$-comodule structure by setting
$\Delta_{\mathcal{F}(M)}^{1|0} \defeq \mu_{D_1,M} \circ \mathcal{F} ( \Delta_M^{1|0} )$.
One can verify that $\mathcal{F} \left( M \right)$ then becomes a comodule object over $D_2$,
and that $\mathcal{F}$ maps morphisms of $D_1$-comodules to morphisms of $D_2$-comodules.

\section{Homological Algebra} \label{appendix:homological-algebra}

\subsection{A Trivial Perturbation Lemma}

\begin{lm} \label{lm:trivial-perturbation-lemma}
	Let $\mathcal{R} = (R,d)$ be a differential graded-commutative Banach $\mathbbm{k}$-algebra and
	let $C$ be a graded Banach $R$-module equipped with two degree one maps
	$b,\delta \colon C \rightharpoonup C$ such that $b + \delta$ is a derivation over $d$ and
	$\left( b + \delta \right)^2 = 0$. \footnote{We do not require that $b^2 = 0$ or $\delta^2 = 0$. Note also
		that $b$ can be $R$-linear and $\delta$ a derivation over $d$ or $\delta$ can be $R$-linear and then
		$b$ is a derivation over $d$.}

	Assume we have a degree $-1$ operator $h \colon C \rightharpoonup C$ with $bh + hb = \idd$
	such that $\idd + \delta h + h \delta$ is invertible. Then the differentiable graded $\mathcal{R}$-module
	$\mathcal{C} = \left( C, b + \delta \right)$ is
	contractible. Two possible contracting homotopies $H, H' \colon C \rightharpoonup C$ of
	$\mathcal{C}$ are given by
	\begin{align}
		H  & \defeq h \left( \idd + \delta h + h \delta \right)^{-1}, \label{eq:trivial-perturbation-lemma-H}
		\\
		H' & \defeq \left( \idd + \delta h + h \delta \right)^{-1} h. \label{eq:trivial-perturbation-lemma-H-tag}
	\end{align}
\end{lm}
\begin{proof}
	Denote by $\partial = \partial_{b + \delta}$ the differential on 	$\InnEnd{\mathcal{C}}$. We have
	\begin{equation*}
		\partial \left( h \right) = \left( b + \delta \right) h + h \left( b + \delta \right)
		= \idd + \delta h + h \delta
	\end{equation*}
	and hence $\idd + \delta h + h \delta$ is a chain map, i.e.,
	$\partial \left( \idd + \delta h + h \delta \right) = \partial^2 \left( h \right) = 0$.
	Since $\idd + \delta h + h \delta$ is invertible, the inverse map
	$\left( \idd + \delta h + h \delta \right)^{-1}$ is also a chain map and the product
	rule implies that
	\begin{equation*}
		\begin{aligned}
			\partial \left( H \right) & =
			\partial \left( h \left( \idd + \delta h + h \delta \right)^{-1} \right)
			= \left( \partial h \right) \left( \idd + \delta h + h \delta \right)^{-1}
			- h \partial \left( \left( \idd + \delta h + h \delta \right)^{-1} \right)
			\\
			                          & = \left( \idd + \delta h + h \delta \right) \left( \idd + \delta h + h \delta \right)^{-1} = \idd.
		\end{aligned}
	\end{equation*}
	A similar calculation shows that we also have $\partial \left( H' \right) = \idd$.
\end{proof}

\begin{rem}
	\Cref{lm:trivial-perturbation-lemma} can be viewed as a simple version of the homological perturbation lemma.
	Namely, when $b^2 = 0$, the assumptions of \cref{lm:trivial-perturbation-lemma} imply that $\left( C, b \right)$
	is contractible. Considering $\delta$ as a ``small perturbation'' of $b$,
	(in the sense that $\idd + \delta h + h \delta$ is invertible), the conclusion is that the perturbed
	complex $\left( C, b + \delta \right)$ is also contractible.
\end{rem}

\begin{rem} 
	In general, the explicit contracting homotopies $H'$ and $H$ of
	\cref{lm:trivial-perturbation-lemma} are different, but they coincide if $h^2 = 0$.
\end{rem}

\subsection{Mapping Cones and Homotopy Equivalence}

Let $\mathcal{R} = (R,d)$ be a differential graded $\mathbbm{k}$-algebra and let $\mathcal{A} = (A,d_A)$ and
$\mathcal{B} = (B,d_B)$ be two differential graded $\mathcal{R}$-modules. Given a morphism
$f \colon \mathcal{A} \rightarrow \mathcal{B}$ of DG-modules, the mapping cone of $f$ is given as a graded $R$-module by $\Cone{f} \defeq A[1] \oplus B$.
We denote by $\s \colon A \rightharpoonup A[1]$ and $\sigma = \s^{-1} \colon A[1] \rightharpoonup A$ the natural shift maps.
Elements of degree $i$ in $\Cone{f}$ will be written as column vectors
\begin{equation*}
	\begin{pmatrix}
		\s a \\ b
	\end{pmatrix}
\end{equation*}
where $a \in A[1]^i = A^{i+1}$ and $b \in B^i$. The differential on $\Cone{f}$ is given in matrix notation
by
\begin{equation*}
	d_{\Cone{f}} =
	\begin{pmatrix}
		d_{A[1]} & 0   \\
		f \sigma & d_B
	\end{pmatrix}.
\end{equation*}
More explicitly, we have
\begin{equation*}
	d_{\Cone{f}} \begin{pmatrix} \s a \\ b \end{pmatrix} =
	\begin{pmatrix}
		d_{A[1]} \left( \s a \right) \\
		f \sigma \left( \s a \right) + d_B \left( b \right)
	\end{pmatrix} =
	\begin{pmatrix}
		-\s \left( d_A \left( a \right) \right) \\
		f \left( a \right) + d_B \left( b \right)
	\end{pmatrix}.
\end{equation*}

Now let $\mathcal{C} = (C,d_C)$ and $\mathcal{D} = (D,d_D)$ be another pair of differential graded $\mathcal{R}$-modules and let $g \colon \mathcal{C} \rightarrow \mathcal{D}$ be a morphism of DG-modules.
Given a homogeneous map $H \colon \Cone{f} \rightharpoonup \Cone{g}$ of graded modules,
we can write $H$ in matrix notation as
\begin{equation*}
	H = \begin{pmatrix} H^{A[1]}_{C[1]} & H^{B}_{C[1]} \\ H^{A[1]}_D & H^B_D \end{pmatrix}.
\end{equation*}
The differential of $H$ is then given by
\begin{equation} \label{eq:differential-map-on-mapping-cones}
	\begin{aligned}
		\partial H & = d_{\Cone{g}} \circ H - (-1)^{\degb{H}} H \circ d_{\Cone{f}}
		\\
		           & =
		\begin{pmatrix}
			\partial \left( H^{A[1]}_{C[1]} \right) - (-1)^{\degb{H}} H^B_{C[1]} f \sigma                  &
			\partial \left( H^B_{C[1]} \right)                                                               \\
			\partial \left( H^{A[1]}_D \right) + g \sigma H^{A[1]}_{C[1]} - (-1)^{\degb{H}} H^B_D f \sigma &
			\partial \left( H^B_D \right) + g \sigma H^B_{C[1]}
		\end{pmatrix}.
	\end{aligned}
\end{equation}
In particular, if $H$ is lower triangular we see that
\begin{equation} \label{eq:differential-lower-triangular-map-on-mapping-cones}
	\partial \begin{pmatrix} H^{A[1]}_{C[1]} & 0 \\ H^{A[1]}_D & H^B_D \end{pmatrix} =
	\begin{pmatrix}
		\partial \left( H^{A[1]}_{C[1]} \right)                                                        & 0 \\
		\partial \left( H^{A[1]}_D \right) + g \sigma H^{A[1]}_{C[1]} - (-1)^{\degb{H}} H^B_D f \sigma &
		\partial \left( H^B_D \right)
	\end{pmatrix}.
\end{equation}

The following basic observation is sometimes called Vogt's lemma. It says that given a homotopy equivalence,
one can choose the homotopies in such a way that they satisfy a higher coherence condition.
\begin{lm}[Vogt's lemma] \label{lm:vogts-lemma}
	Let $\varphi \colon A \rightarrow B$ be a homotopy equivalence with homotopy inverse
	$\varphi' \colon B \rightarrow A$. Given a homotopy $h^B \colon B \rightharpoonup B$
	with $\partial h^B = \idd_B - \varphi \varphi'$, there exists a homotopy $h^A \colon A \rightharpoonup A$
	with $\partial h^A = \idd_A - \varphi' \varphi$ and ``homotopies''
	\begin{equation*}
		h^A_B \colon A \rightharpoonup B, \quad h^B_A \colon B \rightharpoonup A
	\end{equation*}
	with $\degb{h^A_B} = \degb{h^B_A} = -2$ such that
	\begin{equation*}
		\partial h^A_B = \varphi h^A - h^B \varphi, \quad \partial h^B_A = \varphi' h^B - h^A \varphi'.
	\end{equation*}
\end{lm}
\begin{proof}
	By assumption, there exists a homotopy $\hat{h}^A \colon A \rightharpoonup A$ such that
	$\partial \hat{h}^A = \idd_A - \varphi' \varphi$. Note that $\varphi \hat{h}^A$ and $h^B \varphi$
	are both homotopies between $\varphi \varphi' \varphi$ and $\varphi$ and hence we have
	$\partial \left( h^B \varphi - \varphi \hat{h}^A \right) = 0$. We will modify the homotopy
	$\hat{h}^A$ by setting $h^A \defeq \hat{h}^A + \varphi' \chi$ where
	$\chi \colon A \rightharpoonup B$ is a degree $-1$ map with $\partial \chi = 0$ to be determined shortly. Then
	\begin{equation*}
		\partial h^A = \partial \hat{h}^A + \partial \left( \varphi' \chi \right) =
		\partial \hat{h}^A = \idd_A - \varphi' \varphi
	\end{equation*}
	so that $h^A$ is also a homotopy between $\idd_A$ and $\varphi' \varphi$. With respect to the new homotopy
	$h^A$, we have
	\begin{equation*}
		\begin{aligned}
			\varphi h^A - h^B \varphi & = \varphi \left( \hat{h}^A + \varphi' \chi \right) - h^B \varphi =
			\left( \varphi \varphi' - \idd_B \right) \chi +
			\left( \chi - \left( h^B \varphi - \varphi \hat{h}^A \right) \right)                           \\
			                          & =
			\partial \left( -h^B \chi \right) +
			\left( \chi - \left( h^B \varphi - \varphi \hat{h}^A \right) \right)
		\end{aligned}
	\end{equation*}
	Hence, we can take $\chi \defeq h^B \varphi - \varphi \hat{h}^A$ and
	$h^A_B \defeq -h^B \chi = h^B \varphi \hat{h}^A - h^B h^B \varphi$ to obtain
	$\partial h^A_B = \varphi h^A - h^B \varphi$. Finally,
	\begin{equation*}
		\begin{aligned}
			\varphi' h^B - h^A \varphi' & = \varphi' h^B - \left( \hat{h}^A + \varphi' \chi \right) \varphi'
			=
			\varphi' h^B - \left( \hat{h}^A + \varphi' h^B \varphi - \varphi' \varphi \hat{h}^A \right)
			\varphi'
			\\
			                            & = \varphi' h^B \left( \idd_B - \varphi \varphi' \right)
			- \left( \idd_A - \varphi' \varphi \right) \hat{h}^A \varphi'
			\\
			                            & = - \partial \left( \varphi' h^B h^B \right) +
			\varphi' \partial \left( h^B \right) h^B -
			\partial \left( \hat{h}^A \hat{h}^A \varphi' \right) -
			\hat{h}^A \partial \left( \hat{h}^A \right) \varphi'
			\\
			                            & = -\partial \left( \varphi' h^B h^B + \hat{h}^A \hat{h}^A \varphi' \right)
			+ \varphi' \left( \idd_B - \varphi \varphi' \right) h^B -
			\hat{h}^A \left( \idd_A - \varphi' \varphi \right) \varphi'
			\\
			                            & = -\partial \left( \varphi' h^B h^B + \hat{h}^A \hat{h}^A \varphi' \right)
			+ \left( \idd_A - \varphi' \varphi \right) \varphi' h^B -
			\hat{h}^A \varphi' \left( \idd_B  - \varphi \varphi' \right)
			\\
			                            & = \partial \left( \hat{h}^A \varphi' h^B - \varphi' h^B h^B -
			\hat{h}^A \hat{h}^A \varphi' \right)
		\end{aligned}
	\end{equation*}
	and so we can take
	$h^B_A \defeq \hat{h}^A \varphi' h^B - \varphi' h^B h^B - \hat{h}^A \hat{h}^A \varphi'$
	to obtain $\partial h^B_A = \varphi' h^B - h^A \varphi'$.
\end{proof}

\begin{lm} \label{lm:contractible-cone-homotopy-equivalence}
	Let $\varphi \colon \mathcal{A} \rightarrow \mathcal{B}$ be a morphism of DG-modules.
	Then $\varphi$ is a homotopy equivalence if and only if $\Cone{\varphi}$ is contractible.
\end{lm}
\begin{proof}
	Let $H \colon \Cone{\varphi} \rightharpoonup \Cone{\varphi}$ be a degree $-1$ map. Using
	\cref{eq:differential-map-on-mapping-cones}, the condition for $H$ to be a contraction can be written
	explicitly as
	\begin{equation} \label{eq:contraction-on-cone-matrix-equation}
		\partial H = \begin{pmatrix}
			\partial \left( H^{A[1]}_{A[1]} \right) + H^B_{A[1]} \varphi \sigma                        &
			\partial \left( H^B_{A[1]} \right)                                                           \\
			\partial \left( H^{A[1]}_B \right) + \varphi \sigma H^{A[1]}_{A[1]} + H^B_B \varphi \sigma &
			\partial \left( H^B_B \right) + \varphi \sigma H^B_{A[1]}
		\end{pmatrix}
		=
		\begin{pmatrix}
			\idd_{A[1]} & 0 \\ 0 & \idd_B
		\end{pmatrix}.
	\end{equation}
	Let us write $H^{A[1]}_{A[1]} = -\s h^A \sigma, H^B_{A[1]} = \s \varphi', H^{A[1]}_B = h^A_B \sigma$ and $H^B_B = h^B$.
	Then $\varphi' \colon B \rightarrow A$ is a degree $0$ map,
	$h^A \colon A \rightharpoonup A$ and $h^B \colon B \rightharpoonup B$ are degree $-1$ maps and
	$h^A_B \colon A \rightharpoonup B$ is a degree $-2$ map. In terms
	of the maps $\varphi,h^A,h^B,h^A_B$, the equation \eqref{eq:contraction-on-cone-matrix-equation} can be written
	as
	\begin{align}
		\partial h^A      & = \idd_A - \varphi' \varphi, \\
		\partial \varphi' & = 0,                         \\
		\partial h^A_B    & = \varphi h^A - h^B \varphi, \\
		\partial h^B      & = \idd_B - \varphi \varphi'.
	\end{align}
	Hence, the data of a contraction $H$ on $\Cone{\varphi}$ consists of:
	\begin{enumerate}
		\item A DG morphism $\varphi' \colon \mathcal{B} \rightarrow \mathcal{A}$ which is a homotopy inverse
		      of $\varphi$.
		\item Two homotopies $h^A \colon A \rightharpoonup A$ and $h^B \colon B \rightharpoonup B$ which satisfy
		      $\partial h^A = \idd_A - \varphi' \varphi$ and $\partial h^B = \idd_B - \varphi \varphi'$.
		\item A ``second order'' homotopy $h^A_B \colon A \rightharpoonup B$ which satisfies
		      $\partial h^A_B = \varphi h^A - h^B \varphi$.
	\end{enumerate}
	Hence, if $\Cone{\varphi}$ is contractible, we see that $\varphi$ is a homotopy equivalence. Conversely,
	if $\varphi$ is a homotopy equivalence, we can choose a homotopy inverse
	$\varphi' \colon \mathcal{B} \rightarrow \mathcal{A}$ and a homotopy $h^B \colon B \rightharpoonup B$
	which satisfies $\partial h^B = \idd_B - \varphi \varphi'$. Then Vogt's lemma (\cref{lm:vogts-lemma})
	guarantees that we can find a homotopy $h^A \colon A \rightharpoonup A$ with
	$\partial h^A = \idd_A - \varphi' \varphi$, \textit{and} a second order homotopy $h^A_B \colon A \rightharpoonup B$
	which satisfies $\partial h^A_B = \varphi h^A - h^B \varphi$, and use all the
	data to construct a contraction of $\Cone{\varphi}$.
\end{proof}

\begin{figure}[htb]
	\begin{tikzcd}
		{\mathcal{A}} & {\mathcal{B}} \\
		{\mathcal{C}} & {\mathcal{D}}
		\arrow["f", from=1-1, to=1-2]
		\arrow["\psi", from=1-2, to=2-2]
		\arrow["\varphi"', from=1-1, to=2-1]
		\arrow["g"', from=2-1, to=2-2]
		\arrow["{h^A_D}", squiggly, harpoon, from=1-1, to=2-2]
	\end{tikzcd}
	\caption{A square commuting up to homotopy, columns are homotopy equivalences.}
	\label{fig:square-commutes-up-to-homotopy}
\end{figure}

\begin{lm} \label{lm:square-commute-up-to-homotopy-induced-map-cones}
	Let $f \colon \mathcal{A} \rightarrow \mathcal{B}, g \colon \mathcal{C} \rightarrow \mathcal{D}$
	and $\varphi \colon \mathcal{A} \rightarrow \mathcal{C}, \psi \colon \mathcal{B} \rightarrow \mathcal{D}$
	be morphisms of DG modules. Assume that $\varphi$ and $\psi$ are homotopy
	equivalences and that there exists a homotopy $h^A_D \colon A \rightharpoonup D$ with
	$\partial h^A_D = \psi f - g \varphi$ (see \cref{fig:square-commutes-up-to-homotopy}).
	Then the map $\Theta \colon \Cone{f} \rightarrow \Cone{g}$ given by
	\begin{equation} \label{eq:Theta-chain-map-between-two-mapping-cones}
		\Theta \defeq \begin{pmatrix} \s \varphi \sigma & 0 \\ h^A_D \sigma & \psi \end{pmatrix}
	\end{equation}
	is a homotopy equivalence.
\end{lm}
\begin{proof}
	Choose homotopy inverses $\varphi' \colon C \rightarrow A$ for $\varphi$ and $\psi' \colon D \rightarrow B$
	for $\psi$. Then choose homotopies $h^C \colon C \rightharpoonup C$ and $h^B \colon B \rightharpoonup B$
	such that $\partial h^C = \idd_C - \varphi \varphi'$ and $\partial h^B = \idd_B - \psi' \psi$.
	Set $h^C_B \defeq \psi' g h^C - \psi' h^A_D \varphi' - h^B f \varphi'$. Then
	\begin{equation*}
		\begin{aligned}
			\partial h^C_B & = \psi' g \left( \idd_C - \varphi \varphi' \right) -
			\psi' \left( \psi f - g\varphi \right) \varphi' - \left( \idd_B - \psi' \psi \right) f \varphi'
			\\
			               & = \psi' g - \psi' g \varphi \varphi'  - \psi' \psi f \varphi' + \psi' g \varphi \varphi' -
			f \varphi' + \psi' \psi f \varphi' = \psi' g - f \varphi'.
		\end{aligned}
	\end{equation*}
	Define $\Theta' \colon \Cone{g} \rightarrow \Cone{f}$ by
	\begin{equation} \label{eq:Theta'-chain-map-inverse-between-two-mapping-cones}
		\Theta' \defeq \begin{pmatrix} \s \varphi' \sigma & 0 \\ h^C_B \sigma & \psi' \end{pmatrix}.
	\end{equation}
	We claim that both $\Theta$ and $\Theta'$ are morphisms of DG-modules and that
	$\Theta'$ is a homotopy inverse of $\Theta$.
	Using \cref{eq:differential-lower-triangular-map-on-mapping-cones}, we see that
	\begin{equation*}
		\partial \Theta =
		\begin{pmatrix}
			\partial \left( \s \varphi \sigma \right)                                & 0 \\
			\left( \partial \left( h^A_D \right) + g \varphi - \psi f \right) \sigma &
			\partial \left( \psi \right)
		\end{pmatrix} = 0
	\end{equation*}
	since $\partial h^A_D = \psi f - g \varphi$. Hence $\Theta$ is a morphism of DG-modules. A similar calculation
	shows that $\Theta'$ is also a morphism of DG-modules. Now, by Vogt's lemma (\cref{lm:vogts-lemma})
	we can choose homotopies $h^A \colon A \rightharpoonup A$ and $h^A_C \colon A \rightharpoonup C$ such that
	$\partial h^A = \idd_A - \varphi' \varphi$ and $\partial h^A_C = \varphi h^A - h^C \varphi$.
	Similarly, choose homotopies $h^D \colon D \rightharpoonup D$ and $h^B_D \colon B \rightharpoonup D$ such that
	$\partial h^D = \idd_D - \psi \psi'$ and $\partial h^B_D = \psi h^B - h^D \psi$. Then set
	\begin{equation*}
		\begin{aligned}
			H^{\Cone{f}} & \defeq
			\begin{pmatrix}
				-\s h^A \sigma                                                    & 0   \\
				\left( h^B f h^A + \psi' g h^A_C + \psi' h^A_D h^A \right) \sigma & h^B
			\end{pmatrix}, \\
			H^{\Cone{g}} & \defeq
			\begin{pmatrix}
				-\s h^C \sigma                                                          & 0   \\
				\left( h^D g h^C + h^B_D f \varphi' - h^D h^A_D \varphi' \right) \sigma & h^D
			\end{pmatrix}.
		\end{aligned}
	\end{equation*}
	Using \cref{eq:differential-lower-triangular-map-on-mapping-cones}, we calculate
	\begin{equation*}
		\begin{aligned}
			\partial H^{\Cone{f}} & =
			\begin{pmatrix}
				\s \partial \left( h^A \right) \sigma & 0            \\
				\left( \partial \left( h^B f h^A + \psi' g h^A_C + \psi' h^A_D h^A \right)
				+ h^B f - f h^A \right) \sigma
				                                      & \partial h^B
			\end{pmatrix}.
		\end{aligned}
	\end{equation*}
	Let us calculate the boundary term appearing in the lower left corner of $\partial H^{\Cone{f}}$.
	We have
	\begin{equation*}
		\begin{aligned}
			\partial \left( h^B f h^A + \psi' g h^A_C + \psi' h^A_D h^A \right) & =
			\left( \idd_B - \psi' \psi \right) f h^A - h^B f \left( 1 - \varphi' \varphi \right)
			+ \psi' g \left( \varphi h^A - h^C \varphi \right)
			\\
			                                                                    & \qquad
			+ \psi' \left( \psi f - g \varphi \right) h^A - \psi' h^A_D \left( \idd_A - \varphi' \varphi \right)
			\\
			                                                                    & = fh^A - \psi' \psi f h^A - h^B f + h^B f \varphi' \varphi + \psi' g \varphi h^A - \psi' g h^C \varphi
			\\
			                                                                    & \qquad
			+ \psi' \psi f h^A - \psi' g \varphi h^A - \psi' h^A_D + \psi' h^A_D \varphi' \varphi
			\\
			                                                                    & = f h^A - h^B f +
			\left( \psi' h^A_D \varphi' + h^B f \varphi' - \psi' g h^C \right) \varphi - \psi' h^A_D
			\\
			                                                                    & = f h^A - h^B f - \left( h^C_B \varphi + \psi' h^A_D \right).
		\end{aligned}
	\end{equation*}
	Hence,
	\begin{equation*}
		\begin{aligned}
			\partial H^{\Cone{f}} & =
			\begin{pmatrix}
				\s \left( \idd_A - \varphi' \varphi \right) \sigma  & 0                   \\
				- \left( h^C_B \varphi + \psi' h^A_D \right) \sigma & \idd_B - \psi' \psi
			\end{pmatrix}
			= \idd_{\Cone{f}} - \Theta' \Theta.
		\end{aligned}
	\end{equation*}
	Similarly, using \cref{eq:differential-lower-triangular-map-on-mapping-cones}, we calculate
	\begin{equation*}
		\begin{aligned}
			\partial H^{\Cone{g}} & =
			\begin{pmatrix}
				\s \partial \left( h^C \right) \sigma & 0            \\
				\left( \partial \left( h^D g h^C + h^B_D f \varphi' - h^D h^A_D \varphi' \right)
				+ h^D g - g h^C \right) \sigma
				                                      & \partial h^D
			\end{pmatrix}.
		\end{aligned}
	\end{equation*}
	Concentrating on the boundary term appearing in the lower left corner of $\partial H^{\Cone{g}}$,
	we see that
	\begin{equation*}
		\begin{aligned}
			\partial \left( h^D g h^C + h^B_D f \varphi' - h^D h^A_D \varphi' \right) & =
			\left( \idd_D - \psi \psi' \right) g h^C - h^D g \left( \idd_C - \varphi \varphi' \right)
			+ \left( \psi h^B - h^D \psi \right) f \varphi'
			\\
			                                                                          & \qquad
			- \left( \idd_D - \psi \psi' \right) h^A_D \varphi' + h^D \left( \psi f - g \varphi \right) \varphi'
			\\
			                                                                          & = gh^C - \psi \psi' g h^C - h^D g + h^D g \varphi \varphi' + \psi h^B f \varphi' - h^D \psi f \varphi'
			\\
			                                                                          & \qquad
			- h^A_D \varphi' + \psi \psi' h^A_D \varphi' + h^D \psi f \varphi' - h^D g \varphi \varphi'
			\\
			                                                                          & = g h^C - h^D g
			- \psi \left( \psi' g h^C - h^B f \varphi' - \psi' h^A_D \varphi' \right) - h^A_D \varphi'
			\\
			                                                                          & = g h^C - h^D g - \left( h^A_D \varphi' + \psi h^C_B \right).
		\end{aligned}
	\end{equation*}
	Hence,
	\begin{equation*}
		\begin{aligned}
			\partial H^{\Cone{g}} & =
			\begin{pmatrix}
				\s \left( \idd_C - \varphi \varphi' \right) \sigma  & 0                   \\
				- \left( h^A_D \varphi' + \psi h^C_B \right) \sigma & \idd_D - \psi \psi'
			\end{pmatrix}
			= \idd_{\Cone{g}} - \Theta \Theta'.
		\end{aligned}
	\end{equation*}
\end{proof}

\begin{cor} \label{cor:triangle-homotopy-equivalence-mapping-cone}
	Let $\mathcal{A},\mathcal{B},\mathcal{D}$ be differential graded $\mathcal{R}$-modules. Let
	$f \colon \mathcal{A} \rightarrow \mathcal{B}$, $g \colon \mathcal{A} \rightarrow \mathcal{D}$
	and $\psi \colon \mathcal{B} \rightarrow \mathcal{D}$ be morphisms such that $\psi f = g$
	and assume that $\psi$ is a homotopy equivalence.
	Then the map $\psi^{+} \colon \Cone{f} \rightarrow \Cone{g}$ given by
	\begin{equation} \label{eq:homotopy-equivalence-mapping-cones-formula}
		\begin{pmatrix} \s a \\ b \end{pmatrix} \mapsto \begin{pmatrix} sa \\ \psi(b) \end{pmatrix}
	\end{equation}
	is a homotopy equivalence. If $\psi' \colon \mathcal{D} \rightarrow \mathcal{B}$ is a homotopy
	inverse for $\psi$ and $h \colon B \rightharpoonup B$ is such that
	$\partial h = \idd - \psi' \psi$ then the map $\Theta' \colon \Cone{g} \rightarrow \Cone{f}$ given by
	\begin{equation} \label{eq:homotopy-inverse-equivalence-mapping-cones-formula}
		\begin{pmatrix} \s a \\ d \end{pmatrix} \mapsto
		\begin{pmatrix} \s a \\ -h \left( f \left( a \right) \right) + \psi' \left( d \right) \end{pmatrix}
	\end{equation}
	is a homotopy inverse for $\Theta$.
\end{cor}
\begin{proof}
	Let $\mathcal{C} = \mathcal{A}$ and $\varphi \colon \mathcal{A} \rightarrow \mathcal{C}$ be the identity map.
	Then $\varphi$ and $\psi$ are homotopy equivalences and $\psi f = g \varphi$ so we can take $h^A_D = 0$. Applying
	\cref{lm:square-commute-up-to-homotopy-induced-map-cones} we see that the map
	$\Theta \colon \Cone{f} \rightarrow \Cone{g}$, which coincides with $\psi^{+}$, is
	a homotopy equivalence.

	In the proof of \cref{lm:square-commute-up-to-homotopy-induced-map-cones},
	we can take $\varphi' = \idd, h^C = 0$ and $h^B = h$ and then $h^C_B = -hf$. Then the
	resulting formula \eqref{eq:Theta'-chain-map-inverse-between-two-mapping-cones} for a homotopy inverse $\Theta'$
	coincides with \cref{eq:homotopy-inverse-equivalence-mapping-cones-formula}.
\end{proof}

\begin{cor} \label{cor:projection-from-mapping-cone-equivalence}
	Let $\mathcal{A}, \mathcal{B}$ be differential graded $\mathcal{R}$-modules and let
	$f \colon \mathcal{A} \rightarrow \mathcal{B}$ be a morphism. Assume that $\mathcal{B}$ is
	contractible with contracting homotopy $h$. Then the natural projection map
	$p \colon \Cone{f} \rightarrow \mathcal{A}[1]$ is a homotopy equivalence with homotopy inverse
	given by the map
	\begin{equation} \label{eq:homotopy-inverse-projection-of-cone}
		sa \mapsto
		\begin{pmatrix}
			sa \\ -h \left( f \left( a \right) \right)
		\end{pmatrix}.
	\end{equation}
\end{cor}
\begin{proof}
	Apply \cref{cor:triangle-homotopy-equivalence-mapping-cone} with $\mathcal{D} = 0, g = 0$ and $\psi = 0$.
\end{proof}

\begin{cor} \label{cor:inclusion-into-mapping-cone-equivalence}
	Let $\mathcal{C}, \mathcal{D}$ be differential graded $\mathcal{R}$-modules and let
	$g \colon \mathcal{C} \rightarrow \mathcal{D}$ be a morphism. Assume that $\mathcal{C}$ is
	contractible with contracting homotopy $h$. Then the natural inclusion map
	$i \colon \mathcal{D} \rightarrow \Cone{g}$ is a homotopy equivalence with homotopy inverse
	given by the map
	\begin{equation} \label{eq:homotopy-inverse-inclusion-into-cone}
		\begin{pmatrix} sc \\ d \end{pmatrix}
		\mapsto
		g \left( h \left( c \right) \right) + d.
	\end{equation}
\end{cor}
\begin{proof}
	Applying \cref{lm:square-commute-up-to-homotopy-induced-map-cones} with $\mathcal{A} = 0, \mathcal{B} = \mathcal{D}$
	and $f = 0, \varphi = 0, \psi = \id$, we see that $i \colon \Cone{f} = \mathcal{D} \rightarrow \Cone{g}$
	is a homotopy equivalence.

	In the proof of \cref{lm:square-commute-up-to-homotopy-induced-map-cones},
	we can take $\varphi' = 0, \psi' = \idd, h^C = h, h^B = 0$ and then $h^C_B = gh$. Then the
	resulting formula \eqref{eq:Theta'-chain-map-inverse-between-two-mapping-cones} for a homotopy inverse $\Theta'$
	coincides with \cref{eq:homotopy-inverse-inclusion-into-cone}.
\end{proof}

\section{Non-Archimedean Groups, Rings and Modules}  \label{appendix:non-archimedean-groups-rings-modules}

\subsection{Non-Archimedean Abelian Groups} \label{sub:non-arch-abelian-groups}

\begin{dfn}
	Let $G$ be an abelian group. A \textbf{non-Archimedean group seminorm} on $G$ is a function
	$\nnorm \colon G \rightarrow \RPL$ which satisfies the following properties:
	\begin{enumerate}
		\item $\nnorm[0] = 0$.
		\item $\nnorm[g - h] \leq \max \Set{\nnorm[g], \nnorm[h]}$ for all $g,h \in G$.
	\end{enumerate}
	If, in addition, $\nnorm[g] = 0 \iff g = 0$, we call $\nnorm$ a \textbf{non-Archimedean group norm}.
\end{dfn}

A pair $\left( G, \nnorm \right)$ where $G$ is an abelian group and $\nnorm$ is a (semi)norm on $G$ is called
a \textbf{(semi)normed group}. When there is no possibility of confusion, we will often suppress the seminorm
from our notation. The seminorm $\nnorm$ defines a pseudometric $d(g,h) \defeq \nnorm[g-h]$ which turns $G$ into a topological group. If the seminorm is a norm, the pseudometric is an honest metric. A normed group
$\left( G, \nnorm \right)$ is called a \textbf{Banach group} if $G$ with the induced metric is a
complete metric space. We emphasize that for us, a Banach group is both normed and complete.

A group morphism $f \colon \left( G, \nnorm_G \right) \rightarrow \left( H, \nnorm_H \right)$
between seminormed groups is said to be \textbf{bounded} if there exists $C > 0$ such that
$\nnorm[f(g)]_H \leq C \cdot \nnorm[g]_G$ for all $g \in G$. Bounded morphisms are continuous with respect
to the induced topologies but in general continuous group morphisms need not be bounded.
Assuming $f$ is bounded, the \textbf{operator seminorm} of $f$ is defined by
\begin{equation} \label{def:norm-of-morphism}
	\nnorm[f] \defeq \inf \Set{C \geq 0}[{\nnorm[f(g)]_H \leq C \cdot \nnorm[g]_G \,\,\, \forall g \in G}].
\end{equation}
The morphism $f$ is called \textbf{contractive}
if $\nnorm[f(g)]_H \leq \nnorm[g]_G$ for all $g \in G$ and an \textbf{isometry} if $\nnorm[f(g)]_H =
	\nnorm[g]_G$ for all $g \in G$. Note that an isometry need not be injective if the domain
$\left( G, \nnorm_G \right)$ is not normed.

The definitions of bounded maps and the operator seminorm extend naturally to multilinear maps. A multilinear
map $B \colon G_1 \times \dots \times G_n \rightarrow H$ between seminormed groups is said to be \textbf{bounded} if there exists $C > 0$ such that
\begin{equation*}
	\nnorm[B(g_1,\dots,g_n)]_H \leq C \cdot \nnorm[g_1]_{G_1} \cdots \nnorm[g_n]_{G_n}
\end{equation*}
for all $g_1 \in G_1, \dots, g_n \in G_n$ and the
\textbf{operator seminorm} of a bounded multilinear map $B$ is defined by
\begin{equation} \label{def:norm-of-multilinear-map}
	\nnorm[B] \defeq \inf \Set{C \geq 0}
	[{\nnorm[B(g_1,\dots,g_n)]_H \leq C \cdot \nnorm[g_1]_{G_1} \cdots \nnorm[g_n]_{G_n} \,\,\, \forall g_i \in G_i}].
\end{equation}
The multilinear map $B$ is said to be \textbf{contractive} if $\nnorm[B(g_1,\dots,g_n)]_H \leq \nnorm[g_1]_{G_1} \cdots \nnorm[g_n]_{G_n}$ for all $g_1 \in G_1, \dots, g_n \in G_n$ (i.e., if $\nnorm[B] \leq 1$).

Given a seminormed group $\left( G, \nnorm_G \right)$ and a subgroup $H \leq G$, we can endow $H$ with the restriction seminorm $h \mapsto \nnorm[h]_G$ which we call the \textbf{induced seminorm}. Note that if
$\nnorm_G$ is a norm then $H$ is normed and if $G$ is complete and $H$ is closed then $H$ is also complete. Similarly, we can endow the quotient group $G / H$ with the \textbf{quotient seminorm} defined by
\begin{equation*}
	\nnorm[\eqcl{g}]_{G / H} \defeq \inf_{h \in H} \nnorm[g + h]_G.
\end{equation*}
Note that if $\nnorm_G$ is a norm, the quotient seminorm is a norm if and only if $H$ is closed and if $G$ is complete and $H$ is closed, then $G / H$ is complete.

When working with seminormed or Banach groups, there are at least two natural categories one can work
in. We can choose the morphisms in the category to be either bounded maps or contractive maps. We will
work with contractive maps as morphisms since the resulting category behaves more nicely (it admits infinite
coproducts and products). Let us denote the category of seminormed groups with contractive group morphisms by
$\SNAb$ and by $\BAb$ the full subcategory of $\SNAb$ whose objects are Banach groups.
\begin{rem}
	Even though the morphisms in our categories are contractive, we will see that bounded morphisms appear naturally
	as the ``internal hom'' object of our categories.
	Since we often work with both bounded and contractive morphisms,
	we introduce the notation $f \colon G \rightharpoonup H$ for bounded group morphisms while reserving the
	notation $f \colon G \rightarrow H$ for contractive group morphisms (the actual morphisms in our category)
	and use similar notation to differentiate visually between contractive and bounded multilinear maps.
\end{rem}
Given two seminormed groups $G$ and $H$, we will denote the set of morphisms between $G$ and $H$ by
\begin{equation}
	\Hom{G}{H}[] \defeq
	\Set{f \colon G \rightarrow H}
	[{f \textrm{ is a group homomorphism with } \nnorm[f] \leq 1}]
	\label{eq:hom-seminormed-abelian-groups}
\end{equation}
and the seminormed group of all bounded group homomorphisms between $G$ and $H$ by
\begin{equation}
	\InnHom{G}{H}[][] \defeq \Set{f \colon G \rightharpoonup H}[f \textrm{ is a bounded group homomorphism}].
	\label{eq:inn-hom-seminormed-abelian-groups}
\end{equation}
where we equip $\InnHom{G}{H}[][]$ with the operator seminorm given by \cref{def:norm-of-morphism}.

While our interest lies mainly in $\BAb$, some constructions such as the algebraic tensor product of
Banach groups result in a seminormed group which needs to be completed (see \cref{rem:tensor-product-normed-groups-only-seminormed}).
Hence, we find it useful to discuss both categories $\SNAb$ and $\BAb$ and the relations between them.

\begin{rem}
	Note that any abelian group $G$ can be endowed with the \textbf{trivial norm}
	\begin{equation}
		\trivnorm[g] \defeq \begin{cases}
			1 & g \neq 0, \\
			0 & g = 0
		\end{cases}
		\label{eq:def-triv-norm}
	\end{equation}
	which makes $\left( G, \trivnorm \right)$ a Banach group with the discrete topology. Any group morphism
	$f \colon G \rightarrow H$ becomes a contractive group morphism
	$f \colon \left( G, \trivnorm \right) \rightarrow \left( H, \trivnorm \right)$,
	and so we can think of the category of $\Ab$ of abelian groups as a full subcategory of $\BAb$.
	The inclusion $\Ab \hookrightarrow \BAb$ is left adjoint to the ``underlying set'' functor which sends a
	Banach group $\left( G, \nnorm_G \right)$ to its closed unit ball
	$\clball{G}{1} \defeq \Set{g \in G}[{\nnorm[g]_G \leq 1}]$.

	When thinking of a group without a predefined norm as a Banach group, we will always assume
	that it is endowed with the trivial norm.
\end{rem}

\subsubsection{The category \texorpdfstring{$\SNAb$}{SNAb}} \label{sub:cat-SNAb}

First, note that since we are in the non-Archimedean setting, the hom set
$\Hom{G}{H}[]$ given by \cref{eq:hom-seminormed-abelian-groups}
is an abelian group\footnote{That is, the sum of two contractive group homomorphisms $f,g \colon G \rightarrow H$ is still contractive as we have $\nnorm[f + g] \leq \max \Set{\nnorm[f],\nnorm[g]} \leq 1$.} and composition of morphisms is bilinear so $\SNAb$ is preadditive.
Endowing $\Hom{G}{H}[]$ with the operator seminorm, it even becomes an object of $\SNAb$, but this is not
``the'' internal hom object of $\SNAb$ which will be discussed shortly.

The category $\SNAb$ is bicomplete. Specific limits and colimits of interest are constructed as follows:
\begin{enumerate}
	\item The categorical product in $\SNAb$ of a family $\left( G_i \right)_{i \in I}$ of seminormed groups
	      is given by the \textbf{bounded product}
	      \begin{equation}
		      \prod_{i \in I}^{\B} G_i \defeq
		      \Set{(g_i)_{i \in I} \in \prod_{i \in I} G_i}[\sup_{i \in I} {\nnorm[g_i]} < \infty]
		      \label{eq:direct-product-seminormed-groups}
	      \end{equation}
	      equipped with the seminorm $\nnorm[(g_i)_{i \in I}] \defeq \sup_{i \in I} \nnorm[g_i]$.
	\item The categorical coproduct $\oplus_{i \in I}^{\SNAb} G_i$ in $\SNAb$
	      of a family $\left( G_i \right)_{i \in I}$ of seminormed groups
	      is given by the algebraic direct sum $\oplus_{i \in I} G_i$ of the underlying
	      groups equipped with the seminorm
	      \begin{equation}
		      \nnorm[\sum_{i \in I} g_i] = \max_{i \in I} \nnorm[g_i]. \label{eq:direct-sum-seminorm}
	      \end{equation}
	      Since the coproduct in $\SNAb$ is the algebraic direct sum equipped with a seminorm, we will
	      often denote it simply by $\oplus_{i \in I} G_i$, leaving the seminorms implicit.
	\item Given a morphism $f \colon G \rightarrow H$ in $\SNAb$ between two seminormed groups:
	      \begin{enumerate}
		      \item The categorical kernel of $f$ is the algebraic kernel $\ker(f)$ endowed with the induced
		            seminorm. The morphism $f$ is a monomorphism if and only if $f$ is injective.
		      \item The categorical cokernel of $f$ is the
		            algebraic cokernel $\coker{f} = H / \Im(f)$ endowed with the quotient norm.
		            The morphism $f$ is an epimorphism if and only if $f$ is surjective.
		      \item The map $f$ is a categorical isomorphism if and only if $f$ is an isometric isomorphism
		            of groups.
	      \end{enumerate}
\end{enumerate}
In particular, we see that $\SNAb$ has finite biproducts, kernels and cokernels and hence $\SNAb$ is
pre-abelian.
\begin{rem} \label{rem:SNAb-not-abelian}
	Note that the categorical image $\ker \left( \coker{f} \right)$ is
	the algebraic image $\Im(f)$ endowed with the induced norm (from $H$) while
	the categorical coimage $\coker{\ker \left( f \right)}$ is given by the algebraic coimage
	$G / \ker \left( f \right)$ endowed with the quotient seminorm (from $G$). In general,
	the canonical map $\coker{\ker \left( f \right)} = G /  \ker \left( f \right) \rightarrow
		\Im(f) = \ker \left( \coker{f} \right)$ need not be an isometry and hence
	the category $\SNAb$ is only pre-abelian and not abelian.
\end{rem}

Given two seminormed groups,
$G$ and $H$, we can endow the algebraic tensor product $G \otimes H$ with a seminorm by setting
\begin{equation}
	\nnorm[x]_{G \otimes H} \defeq \inf \Set{\max_{i \in I} \, \nnorm[g_i]_G \cdot \nnorm[h_i]_H}
	[x = \sum_{i \in I} g_i \otimes h_i, \, |I| < \infty] \label{eq:projective-tensor-seminorm}
\end{equation}
for $x \in G \otimes H$. The seminorm $\nnorm_{G \otimes H}$ is sometimes called the
\textbf{non-Archimedean projective tensor seminorm} and the pair $\left( G \otimes H, \nnorm_{G \otimes H} \right)$
is called the \textbf{seminormed tensor product} of $G$ and $H$.

Similar to the algebraic tensor product, the seminormed tensor product can be characterized by a universal property involving
\textit{bounded} bilinear maps. The seminormed tensor product $G \otimes H$ comes equipped with a canonical
contractive bilinear map $\otimes \colon G \times H \rightarrow G \otimes H$ characterized by
the following universal property: Given a seminormed group $L$ and a bounded bilinear map
$B \colon G \times H \rightharpoonup L$ there exists a unique bounded morphism
$\varphi_B \colon G \otimes H \rightharpoonup L$ with $\nnorm[\varphi_B] = \nnorm[B]$ such that
$\varphi_B \left( g \otimes h \right) = B(g,h)$ for all $g \in G$ and $h \in H$ (see
\cref{fig:seminormed-tensor-product-groups-universal-property}).
\begin{figure}[htb]
	\centering
	\begin{tikzcd}
		{G \times H} && {G \otimes H} \\
		&& L
		\arrow["\otimes", from=1-1, to=1-3]
		\arrow["{\substack{\exists! \, \varphi_B \\ \textrm{bounded}}}", dashed, harpoon, from=1-3, to=2-3]
		\arrow["\substack{B  \textrm{ bounded} \\ \textrm{bilinear}}"', harpoon, from=1-1, to=2-3]
	\end{tikzcd}
	\caption{Universal property of the seminormed tensor product.}
	\label{fig:seminormed-tensor-product-groups-universal-property}
\end{figure}

Given two bounded morphisms $\varphi \colon G \rightharpoonup G'$ and $\psi \colon H \rightharpoonup H'$ between
seminormed groups, their tensor product $\varphi \otimes \psi \colon G \otimes H \rightharpoonup G' \otimes H'$
is also bounded with the bound $\nnorm[\varphi \otimes \psi] \leq \nnorm[\varphi] \cdot \nnorm[\psi]$.
Restricting our attention to contractive morphisms yields a bifunctor $\otimes \colon \SNAb \times \SNAb \rightarrow \SNAb$ which,
together with the standard symmetries, associators and unitors, endows the category
$\SNAb$ with the structure of a symmetric monoidal category whose unit is the
seminormed group $\left( \ZZ, \trivnorm \right)$.\footnote{This means that the standard symmetry maps $G \otimes H \rightarrow H \otimes G$ given by
	$g \otimes h \mapsto h \otimes g$ are \textbf{isometries} with respect to the projective tensor seminorm. Similarly,
	the standard associators and unitors also become isometries.}
It follows from the universal property of the seminormed tensor product that the symmetric monoidal category
$\SNAb$ is closed with the internal hom object given by the seminormed group
$\InnHom{G}{H}[][]$ of all \textit{bounded} group homomorphisms (see \cref{eq:inn-hom-seminormed-abelian-groups}).

\begin{rem}
	Sometimes it is useful to characterize the seminormed tensor product $G \otimes H$ up to an isometric
	isomorphism using a weaker universal property than the one we described above. Assume we have a group $T$
	equipped with a contractive bilinear map $c \colon G \times H \rightarrow T$. We can
	write two universal properties:
	\begin{enumerate}[label=Property \Alph*., ref=Property \Alph*, itemindent=*]
		\item Given a seminormed group $L$ and a \textit{contractive} bilinear map
		      $B \colon G \times H \rightarrow L$, there exists a unique contractive morphism
		      $\varphi^c_B \colon T \rightarrow L$ such that
		      \begin{equation*}
			      \varphi^c_B \left( c \left( g, h \right) \right) = B \left( g, h \right)
		      \end{equation*}
		      for all $g \in G$ and $h \in H$. \label{prop:A}
		\item Given a seminormed group $L$ and a bounded bilinear map
		      $B \colon G \times H \rightharpoonup L$, there exists a unique bounded morphism
		      $\varphi^c_B \colon T \rightharpoonup L$ such that $\nnorm[\varphi^c_B] = \nnorm[B]$ and
		      \begin{equation*}
			      \varphi^c_B \left( c \left( g, h \right) \right) = B \left( g, h \right)
		      \end{equation*}
		      for all $g \in G$ and $h \in H$. \label{prop:B}
	\end{enumerate}
	Both \labelcref{prop:A} and \labelcref{prop:B} characterize $G \otimes H$ uniquely up to an isometric isomorphism.
	Note that \labelcref{prop:A} is stated only with respect to contractive maps, and we do not require that $\nnorm[\varphi^c_B] = \nnorm[B]$.
	\labelcref{prop:A} comes up naturally by thinking of $\cdot \otimes H$ as the left adjoint of the functor
	$\Hom{H}{\cdot}[][\SNAb]$ in the category $\SNAb$ whose morphisms are contractive maps leading to the natural isomorphisms
	of \textit{sets}
	\begin{equation*}
		\Hom{G \otimes H}{L}[][\SNAb] \cong \Hom{G}{\Hom{H}{L}[][\SNAb]}[][\SNAb].
	\end{equation*}
	By abstract nonsense, the adjunction extends naturally to the internal homs, giving us natural isomorphisms of
	\textit{seminormed groups}
	\begin{equation*}
		\InnHom{G \otimes H}{L} \cong \InnHom{G}{\InnHom{H}{L}}.
	\end{equation*}
	This implies in particular that we can work with bounded and not only contractive maps and that
	the correspondence $\varphi_B \leftrightarrow B$ is an isometry.
\end{rem}

\subsubsection{The category \texorpdfstring{$\BAb$}{BAb}} \label{sub:cat-BAb}
Before discussing the properties of $\BAb$, we describe the notion of a (separated) completion
of a seminormed group. We start with a basic extension lemma whose proof is standard when the groups are normed but
works just as well in the seminormed case:
\begin{lm} \label{lm:extending-multilinear-maps-by-density}
	Let $L$ be a Banach group and let $\eta_i \colon G_i \rightarrow H_i$ for $i = 1, \dots, n$ be
	isometries of seminormed groups with dense image. Then, given a bounded multilinear map
	$B \colon G_1 \times \dots \times G_n \rightharpoonup L$ there exists a unique bounded
	multilinear map $\tilde{B} \colon H_1 \times \dots \times H_n \rightharpoonup L$
	such that $\tilde{B} \left( \eta_1 \left( g_1 \right), \dots, \eta_n \left( g_n \right) \right) =
		B \left( g_1, \dots, g_n \right)$ for all $g_1 \in G_1, \dots, g_n \in G_n$. In addition we have
	$\lVert \tilde{B} \rVert  = \nnorm[B]$. \qed
	\begin{figure}[htb]
		\centering
		\begin{tikzcd}
			G_1 \times \dots \times G_n && {H_1 \times \dots \times H_n} \\
			&& L
			\arrow["{\left( \eta_1, \dots, \eta_n \right)}", from=1-1, to=1-3]
			\arrow["B"', harpoon, from=1-1, to=2-3]
			\arrow["{\exists! \, \tilde{B}}", dotted, harpoon, from=1-3, to=2-3]
		\end{tikzcd}
		\captionof{figure}{Extending multilinear maps by density.}
		\label{fig:extension-multilinear-maps}
	\end{figure}
\end{lm}

A (separated) \textbf{completion} of a seminormed group $G$ is defined to be a Banach group $\widehat{G}$
together with an isometry $\eta = \eta_G \colon G \rightarrow \widehat{G}$ whose image is dense in
$\widehat{G}$. By \cref{lm:extending-multilinear-maps-by-density}, any completion
$\eta_G \colon G \rightarrow \widehat{G}$ satisfies the following universal property: Given
a Banach group $H$ and a bounded morphism $f \colon G \rightharpoonup H$, there exists a unique
bounded extension $\tilde{f} \colon \widehat{G} \rightharpoonup H$ (the \textbf{adjunct} of $f$)
with $\lVert \tilde{f} \rVert = \nnorm[f]$ such that $\tilde{f} \circ \eta_G = f$ (see \cref{fig:completion-adjunction}). In addition, if $f$ happens to be an isometry then so is $\tilde{f}$.

\begin{rem}
	Let us clarify in what sense the map $\tilde{f}$ ``extends'' $f$. Note that the kernel of the map
	$\eta \colon G \rightarrow \widehat{G}$ is precisely
	the subgroup $\clball{G}{0} \defeq \Set{g \in G}[{\nnorm[g] = 0}]$. When $G$ is normed, $\eta$ is injective
	and by identifying $G$ isometrically with its image $\eta \left( G \right) \subseteq \widehat{G}$,
	we can really think of $\tilde{f}$ as extending $f$ from the dense subgroup
	$\eta \left( G \right) \subseteq \widehat{G}$ to the completion $\widehat{G}$. In the general case,
	a bounded morphism $f \colon G \rightarrow H$ still induces
	a well-defined bounded morphism $f' \colon \eta \left( G \right) \rightarrow H$ by setting
	$f' \left( \eta \left( g \right) \right) = f \left( g \right)$ and the map
	$\tilde{f} \colon \widehat{G} \rightarrow H$ is an extension of $f'$ from the dense subgroup
	$\eta \left( G \right) \subseteq \widehat{G}$ to the completion $\widehat{G}$.
	To see that $f'$ is well-defined, note that if $g,g' \in G$ with $\eta(g) = \eta(g')$ then
	\begin{equation*}
		\nnorm[g - g'] = \nnorm[\eta \left( g - g' \right)] = 0 \implies
		\nnorm[f(g) - f(g')] = \nnorm[f \left( g - g' \right)] \leq \nnorm[f] \cdot \nnorm[g - g'] = 0
	\end{equation*}
	which implies that $f(g) = f(g')$ since we assumed that $H$ is Banach, and, in particular, normed.
\end{rem}

\phantomsection
\label{par:completion-of-Banach-identity}
One can show that every seminormed group has a completion (see \cite[Proposition 5, Section 1.1.7]{Bosch1984}
for a construction) and the universal property
satisfied by a completion shows that it is unique up to a unique \textbf{isometric} isomorphism.
Thus, we may speak of \textit{the} completion $\widehat{G}$ of a seminormed group $G$
(the map $\eta_G$ being implicit in the background).
When $G$ is Banach, we will assume that $\widehat{G} = G$ and $\eta_G = \id_{G}$.

Given a bounded morphism $f \colon G \rightharpoonup H$, we will denote by
$\widehat{f} \colon \widehat{G} \rightharpoonup \widehat{H}$ the unique bounded morphism such that
$\widehat{f} \circ \eta_G = \eta_H \circ f$ (see \cref{fig:completion-functoriality}).
We have $\lVert \widehat{f} \rVert = \nnorm[f]$ and if $f$ is an isometry then so is $\widehat{f}$. In addition,
the universal property implies that the completion of morphisms is functorial in the sense that
$\widehat{f \circ g} = \widehat{f} \circ \widehat{g}$ for all bounded morphisms $f \colon G \rightharpoonup H$ and
$g \colon F \rightharpoonup G$ and $\widehat{\id_G} = \id_{\widehat{G}}$.

\begin{figure}[htb]
	\centering
	\begin{minipage}[b]{.5\textwidth}
		\centering
		\begin{tikzcd}
			G && {\widehat{G}} \\
			&& H
			\arrow["{\eta_G}", from=1-1, to=1-3]
			\arrow["f"', harpoon, from=1-1, to=2-3]
			\arrow["{\exists ! \tilde{f}}", harpoon, dotted, from=1-3, to=2-3]
		\end{tikzcd}
		\captionof{figure}{Universal property of the completion.}
		\label{fig:completion-adjunction}
	\end{minipage}%
	\begin{minipage}[b]{.5\textwidth}
		\centering
		\begin{tikzcd}
			G && {\widehat{G}} \\
			H && {\widehat{H}}
			\arrow["{\eta_G}", from=1-1, to=1-3]
			\arrow["f"', harpoon, from=1-1, to=2-1]
			\arrow["{\exists ! \widehat{f}}", harpoon, dotted, from=1-3, to=2-3]
			\arrow["{\eta_H}", from=2-1, to=2-3]
		\end{tikzcd}
		\captionof{figure}{\raggedright Completion of a bounded group morphism.}
		\label{fig:completion-functoriality}
	\end{minipage}
\end{figure}

Since
the completion preserves the norm of morphisms, we obtain a completion functor
$\wedge \colon \SNAb \rightarrow \BAb$ which is left adjoint to the forgetful functor
$U \colon \BAb \rightarrow \SNAb$ with
the map $\eta_G \colon G \rightarrow \widehat{G}$ being the unit of the adjunction.
Hence, we see that
$\BAb$ is a reflective replete full subcategory of $\SNAb$ with reflector $\wedge$.

Since $\SNAb$ is bicomplete, the category $\BAb$ is also bicomplete. Limits in $\BAb$ are the same as limits in $\SNAb$ (and the limit
is automatically Banach) while colimits are computed by applying the completion functor. In particular,
we have the following:
\begin{enumerate}
	\item The categorical product of a family $\Set{G_i}_{i \in I}$
	      of Banach groups is given by the bounded product $\prod_{i \in I}^{\B} G_i$
	      of \eqref{eq:direct-product-seminormed-groups}.
	\item The categorical coproduct in $\BAb$ of a family $\Set{G_i}_{i \in I}$ of Banach groups
	      is given by the completion of the algebraic direct sum $\oplus_{i \in I} G_i$ with respect
	      to the norm given by \cref{eq:direct-sum-seminorm}. Equivalently, it is given
	      by the \textbf{complete direct sum}
	      \begin{equation}
		      \cbigoplus_{i \in I} G_i \defeq
		      \Set{(g_i)_{i \in I} \in \prod_{i \in I} G_i}[g_i \to 0] \subseteq
		      \prod_{i \in I}^{\B} G_i,
		      \label{eq:completed-direct-sum}
	      \end{equation}
	      where by $g_i \to 0$ we mean that for any $\varepsilon > 0$ we have $\nnorm[g_i] < \varepsilon$
	      for all but finitely many $i \in I$. The norm
	      on $\coplus_{i \in I} G_i$ is given by
	      $\nnorm[\left( g_i \right)_{i \in I}] \defeq \max_{i \in I} \nnorm[g_i]$.
	\item Given a morphism $f \colon G \rightarrow H$ in $\BAb$ between two Banach groups:
	      \begin{enumerate}
		      \item The categorical kernel of $f$ is the algebraic kernel $\ker(f)$ endowed with the induced
		            norm. The kernel is closed and hence Banach.
		            The morphism $f$ is a monomorphism if and only if $f$ is injective.
		      \item The categorical cokernel of $f$ is given by $H / \overline{\Im(f)}$, which is the
		            separated completion of the algebraic cokernel $\coker{f} = H / \Im(f)$
		            endowed with the quotient norm. The morphism $f$ is an epimorphism if and only if $f$
		            has dense set-theoretic image.
		      \item The map $f$ is a categorical isomorphism if and only if $f$ is an isometric isomorphism
		            of groups.
	      \end{enumerate}
\end{enumerate}
In particular, we see that $\BAb$ has finite biproducts, kernels and cokernels and hence $\BAb$ is
pre-abelian.
\begin{rem} \label{rem:BAb-not-abelian}
	The canonical map
	$\coker{\ker \left( f \right)} = G / \ker \left( f \right) \rightarrow
		\overline{\Im(f)} = \ker \left( \coker{f} \right)$ is injective but not necessarily surjective, and
	even when it is bijective, it is not necessarily an isometry
	so $\BAb$ (like $\SNAb$) is only pre-abelian and not abelian.
\end{rem}

Next, we discuss how to endow $\BAb$ with the structure of a symmetric monoidal category.
Given two seminormed groups $G$ and $H$, the \textbf{complete tensor product}\footnote{A possibly more appropriate
	name would be \textbf{Banach tensor product} since we take the separated completion of $G \otimes H$
	and not just complete it, but we follow standard terminology.} $G \cotimes H$ is defined
to be the completion $\widehat{G \otimes H}$ of the algebraic tensor product $G \otimes H$
with respect to the projective tensor seminorm given by \cref{eq:projective-tensor-seminorm}.
Given $g \in G$ and $h \in H$, we will denote by $g \cotimes h$ the image of $g \otimes h$
under the completion map $\eta \colon G \otimes H \rightarrow \widehat{G \otimes H}$ and call
such elements \textbf{elementary tensors}. Note that unlike in the algebraic or seminormed case,
the complete tensor product $G \cotimes H$ is not spanned by $\ZZ$-linear combinations of
elementary tensors $g \cotimes h$. Instead, the group generated by elementary tensors is dense in $G \cotimes H$.
In other words, $G \cotimes H$ is generated by elementary tensors as a \textit{Banach group}.

An important property of the complete tensor product is that, given two
seminormed groups $G$ and $H$, the canonical map
\begin{equation}
	G \cotimes H = \widehat{G \otimes H} \xrightarrow[\cong]{\widehat{\eta_G \otimes \eta_H}}
	\widehat{\widehat{G} \otimes \widehat{H}} = \widehat{G} \cotimes \widehat{H}
	\label{eq:complete-tensor-product-strong-monoidal}
\end{equation}
is an isometric isomorphism of Banach groups. To construct the inverse map in the other direction,
one can use \cref{lm:extending-multilinear-maps-by-density} and the universal properties of the completion
and the seminormed tensor product.

\begin{rem} \label{rem:tensor-product-normed-groups-only-seminormed}
	Note that even when $G$ and $H$ are normed, the algebraic tensor product $G \otimes H$
	endowed with the projective tensor seminorm given by \cref{eq:projective-tensor-seminorm}
	need not be normed. Hence, we really need to take the \textit{separated completion}
	of $G \otimes H$ to obtain a Banach group.

	For a simple example, fix a prime $p$ and let $G = \ZZ$ be endowed with the $p$-\textbf{adic norm}
	given by $\nnorm[p^n \cdot a]_{p} = p^{-n}$ where $p \nmid a$ and $\nnorm[0]_{p} = 0$.
	Let $H = \QQ$ be endowed with the trivial norm $\trivnorm$. Then both $G$ and $H$ are normed
	and $G \otimes H \cong \QQ$ but the projective seminorm of every $x \in G \otimes H$ is zero. To see this,
	write $x = 1 \otimes q$ for some $q \in \QQ$ and note that we have
	\begin{equation*}
		\nnorm[x] = \nnorm[1 \otimes q] = \nnorm[p^n \otimes \left( p^{-n} \cdot q \right)] \leq
		\nnorm[p^n]_p \cdot \trivnorm[p^{-n} \cdot q] \leq \nnorm[p^n]_p = p^{-n}
	\end{equation*}
	for all $n \in \NN$ which implies that $\nnorm[x] = 0$. Thus, $G \cotimes H = 0$.

	Although $G$ is not Banach, one can replace $G$ with its completion $\widehat{G} = \ZZ_p$, the $p$-adic integers, and
	obtain an example where both $\widehat{G},H$ are Banach and $\widehat{G} \otimes H \cong \QQ_p$ is non-zero
	but the norm of every element of $\widehat{G} \otimes H \cong \QQ_p$ is zero so $\widehat{G} \cotimes H = 0$.
\end{rem}

\begin{figure}[htb]
	\centering
	\begin{tikzcd}
		{G \times H} && {G \otimes H} &&
		{G \cotimes H = \extrawidehat{G \otimes H}} \\
		&& L
		\arrow["\cotimes", curve={height=-25pt}, from=1-1, to=1-5]
		\arrow["\otimes", from=1-1, to=1-3]
		\arrow["\eta", from=1-3, to=1-5]
		\arrow["{\exists!}", dashed, harpoon, from=1-3, to=2-3]
		\arrow["{\exists! \, \varphi_B}", dashed, harpoon, from=1-5, to=2-3]
		\arrow["B"', harpoon, from=1-1, to=2-3]
	\end{tikzcd}
	\caption{Universal property of the complete tensor product of seminormed groups.}
	\label{fig:complete-tensor-product-universal-property}
\end{figure}
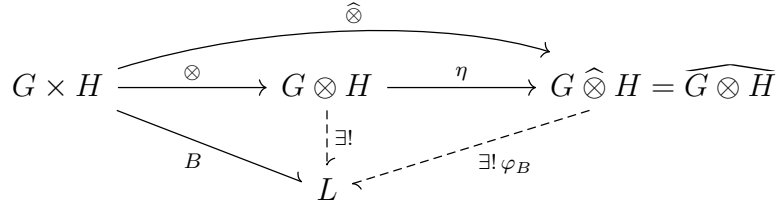

Similar to the seminormed tensor product, the complete tensor product $G \cotimes H$ is also
characterized by a universal property involving bounded bilinear maps into \textit{Banach} groups.
The Banach group $G \cotimes H$ comes equipped with a canonical contractive bilinear map
$\cotimes \colon G \times H \rightarrow G \cotimes H$ characterized by the
following universal property: Given a Banach group $L$ and a bounded bilinear map
$B \colon G \times H \rightharpoonup L$ there exists a unique bounded morphism
$\varphi_B \colon G \cotimes H \rightharpoonup L$ with $\nnorm[\varphi_B] = \nnorm[B]$ such that
$\varphi_B \left( g \cotimes h \right) = B(g,h)$ for all $g \in G$ and $h \in H$ (see
\cref{fig:complete-tensor-product-universal-property}).

Given two bounded morphisms $\varphi \colon G \rightharpoonup G'$ and $\psi \colon H \rightharpoonup H'$ between
seminormed groups, their complete tensor product
$\varphi \cotimes \psi \colon G \cotimes H \rightharpoonup G' \cotimes H'$
is defined by $\varphi \cotimes \psi \defeq \widehat{\varphi \otimes \psi}$ and we have
$\nnorm[\varphi \cotimes \psi] \leq \nnorm[\varphi] \cdot \nnorm[\psi]$. In particular, we see that the
complete tensor product construction gives us a bifunctor $\cotimes \colon \SNAb \times \SNAb \rightarrow \BAb$.

Even though we have defined the complete tensor product for seminormed groups, we now restrict our
attention to the complete tensor product of Banach groups and obtain a bifunctor
$\cotimes \colon \BAb \times \BAb \rightarrow \BAb$. Given three Banach groups $G,H,L$,
we can use the isomorphism of \cref{eq:complete-tensor-product-strong-monoidal} to obtain natural
isometric associativity isomorphisms
\begin{equation*}
	\begin{aligned}
		\left( G \cotimes H \right) \cotimes L & = \extrawidehat{\widehat{G \otimes H} \otimes L} =
		\extrawidehat{\widehat{G \otimes H} \otimes \widehat{L}} \cong
		\extrawidehat{\left( G \otimes H \right) \otimes L}
		\\
		                                       & \cong
		\extrawidehat{G \otimes \left( H \otimes L \right)} \cong
		\extrawidehat{\widehat{G} \otimes \widehat{H \otimes L}} =
		\extrawidehat{G \otimes \widehat{H \otimes L}} =
		G \cotimes \left( H \cotimes L \right)
	\end{aligned}
\end{equation*}
which act on elementary tensors by the expected formula $\left( g \cotimes h \right) \cotimes l \mapsto
	g \cotimes \left( h \cotimes l \right)$. In addition, we also have the natural isometric isomorphisms
\begin{align*}
	 & \ZZ \cotimes G \cong G \qquad           & k \cotimes g & \mapsto k \cdot g,    \\
	 & G \cotimes \ZZ \cong G \qquad           & g \cotimes k & \mapsto k \cdot g,    \\
	 & G \cotimes H \cong H \cotimes G, \qquad & g \cotimes h & \mapsto h \cotimes g,
\end{align*}
where by $\ZZ$ we mean the group $\ZZ$ endowed with the trivial norm \cref{eq:def-triv-norm}.
The bifunctor $\cotimes \colon \BAb \times \BAb \rightarrow \BAb$, together with the
associativity isomorphisms, the unitors and the symmetry maps described above, endow the category $\BAb$
with the structure of a symmetric monoidal category whose unit is the Banach group $\left( \ZZ, \trivnorm \right)$
(the same unit of $\SNAb$).

Given two seminormed groups $G$ and $H$, the seminormed group $\InnHom{G}{H}[][]$
of all bounded group homomorphisms is Banach once $H$ is Banach. In particular, $\InnHom{G}{H}[][]$
is an object of $\BAb$ when both $G$ and $H$ are objects of $\BAb$.
It follows from the universal property of the complete tensor product that the symmetric monoidal
category $\BAb$ is closed with the internal hom object given by
$\InnHom{G}{H}[]$ (the same internal hom object as in $\SNAb$).

With respect to the monoidal structures defined on $\SNAb$ and $\BAb$, the adjunction
\begin{equation*}
	\wedge \colon \SNAb \stackrel[]{\dashv}{\rightleftarrows} \BAb \colon U
\end{equation*}
is naturally enhanced into a monoidal adjunction in which the left adjoint functor $\wedge$ is strong symmetric
monoidal via the coherence isomorphism of \cref{eq:complete-tensor-product-strong-monoidal}
and the right adjoint functor $U$ is lax symmetric monoidal via the canonical map $G \otimes H \rightarrow G \cotimes H$.

\begin{rem}
	From a categorical perspective, the construction of the monoidal structure $\cotimes$ as the completion
	of $\otimes$ is an example of a general construction in which a reflective full subcategory inherits a monoidal
	structure from the ambient category. See \cite{Day1972}.
\end{rem}

\subsubsection{Unordered Infinite Sums in non-Archimedean Banach Groups}
We discuss briefly unordered infinite sums in non-Archimedean Banach groups.

\begin{dfn}
	Let $G$ be a seminormed group, let $I$ be a set and let $(g_{i})_{i \in I}$ be an $I$-indexed family of
	elements of $G$. We say that the \textbf{unordered sum} $\sum_{i \in I} g_i$ \textbf{converges} to an
	element $g \in G$ if for every $\varepsilon > 0$ there exists a finite subset $J \subseteq I$
	such that for any finite subset $J \subseteq K \subseteq I$ we have
	$\nnorm[\sum_{k \in K} g_k - g] < \varepsilon$. An unordered sum $\sum_{i \in I} g_i$ which converges
	to some $g \in G$ is called \textbf{convergent}.
\end{dfn}

When $G$ is normed, it is readily seen that an unordered sum $\sum_{i \in I} g_i$ can converge to
at most one element $g \in G$, in which case we write $\sum_{i \in I} g_i = g$.
Given a bounded homomorphism $f \colon G \rightharpoonup H$ of normed groups and a
convergent unordered sum $\sum_{i \in I} g_i$ in $G$, the unordered sum
$\sum_{i \in I} f \left( g_i \right)$ also converges in $H$
and we have $f \left( \sum_{i \in I} g_i \right) = \sum_{i \in I} f(g_i)$.
In what follows, we restrict our attention to the case where the group $G$ is Banach.
In the non-Archimedean Banach setting, we have the following basic criteria:

\begin{lm} \label{lm:convergence-un-ordered-sums}
	Let $G$ be a \textbf{non-Archimedean Banach} group.
	\begin{enumerate}
		\item Let $(g_n)_{n=0}^{\infty}$ be a sequence of elements of $G$.
		      The ordered sum
		      \begin{equation*}
			      \sum_{n=0}^{\infty} g_n = \lim_{N \to \infty} \sum_{n=0}^N g_n
		      \end{equation*}
		      converges if and only if $\nnorm[g_n] \to 0$.
		\item Let $(g_i)_{i \in I}$ be an $I$-indexed family of elements of $G$.
		      The unordered sum $\sum_{i \in I} g_i$ converges if and only if $g_i \to 0$ in the
		      sense that for any $\varepsilon > 0$ we have $\nnorm[g_i] < \varepsilon$ for all but
		      finitely many $i \in I$.
	\end{enumerate} \qed
\end{lm}

Given an unordered sum $\sum_{i \in I} g_i$, let us set $I_{\neq 0} = \Set{i \in I}[g_i \neq 0]$.
The unordered sum $\sum_{i \in I} g_i$ converges if and only if $\sum_{i \in I_{\neq 0}} g_i$ converges.
\Cref{lm:convergence-un-ordered-sums} implies in particular that for $\sum_{i \in I} g_i$ to be convergent, the
set $I_{\neq 0}$ must be countable. Hence, when discussing convergent unordered sums, we can assume
that $I$ is countable. When $I$ is finite, we have a finite sum. When $I$ is (infinite) countable,
we can choose a bijection $\sigma \colon \NZ \rightarrow I$ and consider
the ordered sum $\sum_{n = 0}^{\infty} g_{\sigma(n)}$. The ordered sum $\sum_{n = 0}^{\infty} g_{\sigma(n)}$
converges if and only if the unordered sum $\sum_{i \in I} g_i$ converges and they converge to the same value.
In particular, an arbitrary rearrangement of a convergent ordered sum remains convergent and converges to the
same value (see \cite[Corollary 4, Section 1.1.8]{Bosch1984}).

Next, we discuss rearrangement of unordered sums:
\begin{lm}
	Let $(g_i)_{i \in I}$ be an $I$-indexed family of elements of $G$
	such that $\sum_{i \in I} g_i$ converges. Let $I = \cup_{\alpha \in A} I_{\alpha}$
	be a decomposition of $I$ into pairwise disjoint subsets. Then for each $\alpha \in A$ the
	sum $\sum_{i \in I_{\alpha}} g_i$ converges and we have
	$\sum_{i \in I} g_i = \sum_{\alpha \in A} \left( \sum_{i \in I_{\alpha}} g_i \right)$. \qed
\end{lm}
In particular, convergent double infinite sums can be rearranged at will. More precisely, let
$(g_k^n)$ be an $\NZ \times \NZ$-indexed family of elements of $G$. When $g_k^n \to 0$ (in the sense
of \cref{lm:convergence-un-ordered-sums}, part $(2)$), we have
\begin{equation*}
	\sum_{(k,n) \in \NZ \times \NZ} g_k^n = \sum_{k=0}^{\infty} \left( \sum_{n=0}^{\infty} g_k^n \right)
	= \sum_{n=0}^{\infty} \left( \sum_{k=0}^{\infty} g_k^n \right) = \lim_{N \to \infty} \sum_{k,n=0}^N g_k^n.
\end{equation*}

Let us comment on the relation between unordered infinite sums and complete direct sums of Banach groups.
Given a family of Banach groups $\left( G_i \right)_{i \in I}$,
let $G = \coplus_{i \in I} G_i$. We can think of an element $g \in G$ in two equivalent ways:
\begin{enumerate}
	\item As an $I$-indexed family $g = (g_i)_{i \in I}$ where $g_i \in G_i$ and $g_i \to 0$
	      as in \eqref{eq:completed-direct-sum}.
	\item As a convergent unordered sum $g = \sum_{i \in I} g_i$ in $G$
	      where each $g_i \in G_i \subseteq G$.
\end{enumerate}
Using the language of unordered sums, an element $g \in G$ has a unique representation
as a convergent sum $g = \sum_{i \in I} g_i$ where $g_i \in G_i$. With
respect to the unique representation $g = \sum_{i \in I} g_i$, the norm of $g$ is given by
$\nnorm[g] = \max_{i \in I} \nnorm[g_i]$. When $I$ is countable, one can choose some enumeration
$I = \left( i_n \right)_{n \in \NZ}$ and work with ordered sums instead. Then any element
$g \in G$ has a unique representation as a convergent ordered sum $g = \sum_{n=0}^{\infty} g_n$ where
each $g_n \in G_{i_n}$ and $g_n \to 0$.

We end this subsection with a simple but useful criterion for checking the convergence of iterations of operators defined on a direct sum.

\begin{lm} \label{lm:direct-sum-iteration-zero-limit}
	Let $G = \coplus_{i \in I} G_i$ and let $\varphi \colon G \rightarrow G$ be a contractive group morphism.
	The following conditions are equivalent:
	\begin{enumerate}
		\item $\varphi^n(g) \to 0$ for each $i \in I$ and $g \in G_i$.
		\item $\varphi^n(g) \to 0$ for each $g \in G$.
	\end{enumerate}
\end{lm}
\begin{proof}
	It is clear that $(2)$ implies $(1)$, so let us show that $(1)$ implies $(2)$. Let $\varepsilon > 0$ and
	let $g = \sum_{i \in I} g_i \in G$ where each $g_i \in G_i$. Let $I_0 \subseteq I$
	be a finite subset such that whenever $i \notin I_0$ we have $\nnorm[g_i] < \varepsilon$.
	Then for each $i \notin I_0$ we have
	\begin{equation*}
		\nnorm[\varphi^n(g_i)] \leq \nnorm[\varphi]^n \nnorm[g_i] < \varepsilon
	\end{equation*}
	for all $n \geq 0$.
	For each $i \in I_0$ choose $N_i$ such that if $n > N_i$ we have $\nnorm[\varphi^n(g_i)] < \varepsilon$.
	Then if $n > \max_{i \in I_0} N_i$ we have
	\begin{equation*}
		\nnorm[\varphi^n(g)] = \nnorm[\sum_{i \in I} \varphi^n \left( g_i \right)]
		\leq \sup_{i \in I} \nnorm[\varphi^n(g_i)] \leq \varepsilon.
	\end{equation*}
\end{proof}

\subsubsection{Relation Between Norms and Filtrations}
Let us comment briefly on the relation between working with filtrations and working with seminorms.
We discuss separately the case of $\NZ$-indexed and $\RR_{\geq 0}$-indexed filtrations.

\begin{enumerate}
	\item Let us denote by $\mathbf{FiltAb}_{\NZ}$ the category of abelian groups $G$, equipped
	      with a descending filtration of subgroups
	      \begin{equation*}
		      G = F_0 G \supseteq F_1 G \supseteq \dots
	      \end{equation*}
	      with morphisms being group homomorphisms $f \colon G \rightarrow H$ compatible with
	      the filtrations in the sense that $f \left( F_{n} G \right) \subseteq F_{n} H$ for all $n \geq 0$.
	      Denote also by $\clball{\SNAb}{1}$ the full subcategory of $\SNAb$ whose objects
	      are seminormed groups $A$ for which $\nnorm[a] \leq 1$ for all $a \in A$.

	      We have an adjunction
	      \begin{equation}
		      L \colon \mathbf{FiltAb}_{\NZ} \stackrel[]{\dashv}{\rightleftarrows} \clball{\SNAb}{1} \colon R
		      \label{eq:adjunction-filt-ab-nz-snab}
	      \end{equation}
	      in which the left adjoint $L$ sends a filtered group $\left( F_n G \right)_{n \geq 0}$ (with
	      $F_0 = G$) to the group $G$ endowed with the seminorm
	      \begin{equation*}
		      \nnorm[g]_{\textrm{filt}} \defeq 2^{-\sup \Set{n \in \NZ}[g \in F_n G]}
		      = \inf \Set{2^{-n}}[g \in F_n G].
	      \end{equation*}
	      The right adjoint $R$ sends a seminormed group $A$ to the group $A$ endowed with
	      filtration $F_n A \defeq \clball{A}{2^{-n}}$ by closed balls. We have
	      $RL = \id$ as functors and the unit map of the adjunction is the identity map (an isomorphism).
	      Hence $L$ is full and faithful and the category $\mathbf{FiltAb}_{\NZ}$ is embedded
	      in $\clball{\SNAb}{1}$.
	      On the other hand, the counit $\varepsilon \colon (LR)(A) \rightarrow A$
	      is the identity map, but it is not an isomorphism as the norm on $(LR)(A) = A$
	      is given by $\nnorm[a]_{\textrm{filt}} = \inf \Set{2^{-n}}[{\nnorm[a] \leq 2^{-n}}]$.\footnote{Note
		      however that the norm $\nnorm_{\textrm{filt}}$ is equivalent to $\nnorm$ in the sense that
		      $\nnorm[a] \leq \nnorm[a]_{\textrm{filt}} \leq 2 \cdot \nnorm[a]$ for all $a \in A$.}
	      Restricting to the full subcategory of $\clball{\SNAb}{1}$ on which $\varepsilon$ is an isomorphism,
	      one obtains an isomorphism of categories between $\mathbf{FiltAb}_{\NZ}$ and the category
	      of seminormed groups whose seminorms take values in $\Set{1,2^{-1},2^{-2},\dots,0}$.

	      Taking into account the monoidal structure on the category of filtered groups
	      (see for example \cite[Chapter 1]{Kleijn2021}),
	      one can verify that the adjunction \eqref{eq:adjunction-filt-ab-nz-snab}
	      is monoidal with $R$ being lax monoidal and $L$ strong monoidal.
	      The tensor constraints for both $R$ and $L$ are the identity maps, and we have
	      $L \left( G \otimes_{\textrm{filt}} H \right) =
		      L \left( G \right) \otimes_{\textrm{SN}} L \left( H \right)$,
	      so that constructions done using the filtered tensor product of filtered groups can
	      be interpreted via $L$ in terms of the seminormed tensor product of seminormed groups.

	      Finally, the algebraic notion of filtered completion corresponds nicely under $L$
	      to the completion of a seminormed group. Recall that a filtered group
	      $\left( F_n G \right)_{n \geq 0}$ is complete if the canonical map
	      $G \rightarrow \varprojlim G / F_n G$ is an isomorphism. The group
	      $\widetilde{G} \defeq \varprojlim G / F_n G$ is called the
	      \textbf{filtered completion} of $G$ and has
	      a natural filtration given by $F_k \left( \widetilde{G} \right) \defeq
		      \ker \left( \widetilde{G} \rightarrow G / F_k G \right)$.
	      In terms of the embedding $L$, a filtered group
	      $G$ is complete if and only if $L \left( G \right)$
	      is a Banach group and the image of the canonical (filtered) morphism
	      $G \rightarrow \widetilde{G}$ under $L$ forms a separated completion of the
	      seminormed group $L(G)$, i.e., $L ( \widetilde{G} )$ and
	      $\widehat{L \left( G \right)}$ are naturally isomorphic as Banach groups.
	\item Let us denote by $\mathbf{FiltAb}_{\RR_{\geq 0}}$ the category of abelian groups $G$, equipped
	      with a descending filtration $\left( F_{\lambda} G \right)_{\lambda \geq 0}$, indexed
	      by non-negative real numbers, such that $F_0 G = G$. The morphisms of $\mathbf{FiltAb}_{\RR_{\geq 0}}$
	      are group homomorphisms $f \colon G \rightarrow H$ compatible with the filtrations in
	      the sense that $f \left( F_{\lambda} G \right) \subseteq F_{\lambda} H$ for all $\lambda \geq 0$.
	      We have an adjunction
	      \begin{equation}
		      L \colon \mathbf{FiltAb}_{\RR_{\geq 0}}
		      \stackrel[]{\dashv}{\rightleftarrows} \clball{\SNAb}{1} \colon R
		      \label{eq:adjunction-filt-ab-rr-geq-snab}
	      \end{equation}
	      where the left adjoint and the right adjoint are given by the same formulas as in the previous
	      item.
	      In this case, we have $LR = \id$ as functors and
	      the counit of the adjunction is the identity map (an isomorphism) so that $R$
	      is full and faithful and $\clball{\SNAb}{1}$ is a reflective subcategory of
	      $\mathbf{FiltAb}_{\RR_{\geq 0}}$.
	      However, the unit of the adjunction is not an isomorphism but instead the identity map
	      $\eta \colon G \rightarrow (RL)(G)$ where the filtration on $(RL)(G)$ is given by
	      \begin{equation*}
		      F_{\lambda} \left( (RL)(G) \right) = \bigcap_{\alpha < \lambda} F_{\alpha} G.
	      \end{equation*}
	      A descending filtration $\left( F_{\lambda} G \right)_{\lambda \geq 0}$
	      is called \textbf{left continuous} if
	      $F_{\lambda} G = \bigcap_{\alpha < \lambda} F_{\alpha} G$ for all $\lambda \geq 0$. Note
	      that the descending filtration by closed balls is always left continuous and $\eta$ is an isomorphism
	      if and only if the filtration of $G$ is left continuous. Hence, by restricting to the full
	      subcategory of $\mathbf{FiltAb}_{\RR_{\geq 0}}$ for which $\eta$ is an isomorphism,
	      we obtain an isomorphism	of categories between descending,
	      left continuous, $\RR_{\geq 0}$-filtered groups and
	      $\clball{\SNAb}{1}$.

	      Taking into account the monoidal structure on the category of filtered groups
	      (see for example \cite[Section 2.2]{DeDeken2018}), one can verify that the adjunction
	      \eqref{eq:adjunction-filt-ab-rr-geq-snab}
	      is monoidal with $R$ being lax monoidal and $L$ strong monoidal.
	      The tensor constraints for both $R$ and $L$ are the identity maps, and we have
	      the identity
	      \begin{equation*}
		      L \left( R \left( A \right) \otimes_{\textrm{filt}} R \left( B \right) \right) =
		      A \otimes_{\textrm{SN}} B.
	      \end{equation*}
	      This means that the seminormed tensor product is obtained from the filtered tensor
	      product via the reflector $L$ and constructions done using the seminormed tensor product
	      can be interpreted in terms of the filtered tensor product of filtered groups.

	      Finally, the completion of a seminormed group corresponds nicely under $R$
	      to the algebraic notion of filtered completion.
	      A filtered group
	      $\left( F_{\lambda} G \right)_{\lambda \geq 0}$ is \textbf{complete} if the canonical map
	      $G \rightarrow \varprojlim G / F_{\lambda} G$ is an isomorphism. The group
	      $\widetilde{G} \defeq \varprojlim G / F_{\lambda} G$ is called the
	      \textbf{filtered completion} of $G$ and has
	      a natural filtration given by $F_{\mu} \left( \widetilde{G} \right) \defeq
		      \ker \left( \widetilde{G} \rightarrow G / F_{\mu} G \right)$.\footnote{Note
		      that we have $\varprojlim_{\lambda \geq 0} G / F_{\lambda} G \cong
			      \varprojlim_{n \in \NZ} G / F_{n} G$ as $\NZ$ is cofinal in $\RR_{\geq 0}$ so that
		      the filtered completion of a $\RR_{\geq 0}$-indexed filtered group, as a group,
		      can be defined in the same way as for a $\NZ$-indexed filtered group.}
	      A seminormed group $A$ is Banach if and only if $R \left( A \right)$
	      is complete and the groups $R ( \widehat{A} )$ and $\widetilde{R \left( A \right)}$
	      are naturally isomorphic as filtered groups.
\end{enumerate}

Although we discussed the relations between the seminormed and filtered framework in the context of abelian
groups, the discussion can be generalized to cover ring, modules, graded modules, etc.

\subsection{Non-Archimedean Rings}
Since the category $\SNAb$ is monoidal, one can talk about algebra objects in $\SNAb$
(see \cref{subsec:alg-in-monoidal-cat}). A \textbf{(non-Archimedean) seminormed ring} is just an
algebra object of $\SNAb$. Unwinding the definition, we see that a seminormed ring
is a pair $\left( R, \nnorm \right)$ where $R$ is a ring and $\nnorm$ is a non-Archimedean seminorm on $R$
which satisfies certain compatibility conditions:
\begin{dfn} \label{def:seminorm-ring}
	Let $R$ be a ring. A \textbf{non-Archimedean ring (semi)norm}\footnote{
		One usually uses the terminology of ``absolute value'' and the notation $| \cdot |$, reserving the
		notation $\nnorm$ for norms on modules. However, since we work with graded objects and use
		$| \cdot |$ to denote the degree of an element, we chose to use $\nnorm$ for ring norms
		and module norms.}
	on $R$ is a function $\nnorm \colon R \rightarrow \RPL$, such that:
	\begin{enumerate}
		\item $\nnorm$ is a non-Archimedean group (semi)norm on $(R,+)$.
		\item $\nnorm[r \cdot s] \leq \nnorm[r] \cdot \nnorm[s]$ for all $r,s \in R$.
		\item $\nnorm[1_R] \leq 1$.\footnote{Note that conditions $(3)$ and $(2)$ imply
			      that in fact $\nnorm[1_R] = 1$ or $\nnorm[1_R] = 0$. If $\nnorm[1_R] = 0$ then $\nnorm[r] = 0$
			      for all $r \in R$. When $R$ is normed, this is possible only when $R = 0$ is the zero ring.}
	\end{enumerate}
\end{dfn}
Hence, a (semi)normed ring is a pair $\left( R, \nnorm \right)$ where $R$ is a ring and $\nnorm$ is a
non-Archimedean ring (semi)norm. The topology on $R$ induced by the associated pseudometric turns $R$ into
a topological ring. A normed ring $\left( R, \nnorm_R \right)$ is
called a \textbf{Banach ring} if $R$ with the induced metric is a complete metric space.

\begin{rem} \label{rem:banach-ring-as-algebra-object}
	Note that a Banach ring is the same thing as an algebra object in the category $\BAb$. A priori
	an algebra object of $\BAb$ is a Banach group $R$ equipped with a ``multiplication'' of the form
	$m \colon R \cotimes R \rightarrow R$ which satisfies associativity conditions stated using $\cotimes$.
	By the universal property of the complete tensor product, such a multiplication corresponds
	bijectively to a contractive bilinear map $\cdot \colon R \times R \rightarrow R$ which endows $R$ with an
	associative multiplication.

	Equivalently, but rephrased differently, the forgetful functor $\BAb \rightarrow \SNAb$ is lax monoidal
	and hence sends an algebra object $R$ of $\BAb$ to an algebra object of $\SNAb$. The multiplication
	$R \otimes R \rightarrow R$ is obtained by precomposing the multiplication $R \cotimes R \rightarrow R$
	with the canonical map $R \otimes R \rightarrow R \cotimes R$.
\end{rem}

Given a seminormed ring $R$, the completion $\widehat{R}$, which is a priori an abelian Banach group,
has a unique structure of a Banach ring with respect to which the canonical morphism $\eta_R \colon R \rightarrow \widehat{R}$ becomes an isometric morphism of rings.
The multiplication $\cdot_{\widehat{R}}$ on $\widehat{R}$ is obtained by extending the
multiplication $\cdot_R$ on $R$ using \cref{lm:extending-multilinear-maps-by-density} (see
\cref{fig:extending-ring-multiplication-to-completion}) and the unit of $\widehat{R}$ is given
by $\eta_R \left( 1_R \right)$.

\begin{figure}[htb]
	\centering
	\begin{tikzcd}
		R \times R && {\widehat{R} \times \widehat{R}} \\
		R && {\widehat{R}}
		\arrow["{\left( \eta_R, \eta_R \right)}", from=1-1, to=1-3]
		\arrow["\cdot_R"', from=1-1, to=2-1]
		\arrow["{\exists ! \,\, \cdot_{\widehat{R}}}", dotted, from=1-3, to=2-3]
		\arrow["{\eta_R}"', from=2-1, to=2-3]
	\end{tikzcd}
	\captionof{figure}{Extending ring multiplication to the completion.}
	\label{fig:extending-ring-multiplication-to-completion}
\end{figure}

Given a bounded ring morphism $f \colon R \rightharpoonup S$ between seminormed rings, the completion
$\widehat{f} \colon \widehat{R} \rightharpoonup \widehat{S}$ becomes a bounded morphism of Banach rings with
respect to the previously described ring structures on $\widehat{R}$ and $\widehat{S}$.

\begin{ex} \label{ex:formal-power-laurent-series-Banach-ring}
	Let $k$ be a ring, let $N \in \ZZ$ and let $p \left( X \right) = \sum_{n \geq N} a_n X^n$ be a formal Laurent
	series with coefficients in $k$. Recall that the \textbf{order} of $p \left( X \right)$ is given by
	\begin{equation*}
		\ord  p \left( X \right) \defeq
		\begin{cases}
			\min \Set{n \geq 0}[a_n \neq 0] & p \left( X \right) \neq 0, \\
			+\infty                         & p \left( X \right) = 0.
		\end{cases}
	\end{equation*}
	In particular, the notion of order makes sense also for polynomials and formal power series where we identify
	polynomials and formal power series with formal Laurent series in the obvious way.
	Given $0 < \alpha < 1$, we defined the \textbf{order norm}
	of a polynomial or a formal power/Laurent series by the formula
	\begin{equation}	 \label{eq:order-norm}
		\nnorm[p \left( X \right)] = \nnorm[p \left( X \right)] _{\alpha} \defeq \alpha^{\ord p \left( X \right)}.
	\end{equation}

	Endowed with order norm, the polynomial ring $k[X]$ becomes a normed ring which is not complete.
	Its completion is given by the formal power series ring $\pows{k}[X]$, also endowed with the order norm.
	Hence, $\pows{k}[X]$ is a Banach ring and inside it, the formal power series
	$\sum_{n \geq 0} a_n X^n$ actually ``converges''
	in the sense that we have
	\begin{equation*}
		\lim_{N \to \infty} \sum_{n = 0}^N a_n X^n = \sum_{n \geq 0} a_n X^n.
	\end{equation*}

	Similarly, the ring $k \left( \left(  X \right) \right)$ of formal Laurent series endowed with order
	norm is also a Banach ring. When $k$ is a field, then
	$k \left( \left( X \right) \right)$ is an example of a Banach field, the field of fractions of
	$\pows{k}[X]$.
\end{ex}

Let us denote by $\SNRing$ the category of seminormed rings with contractive ring morphisms as morphisms
and by $\BRing$ the full subcategory of $\SNRing$ whose objects are Banach rings. Since completion
preserves the norm of morphisms, by restricting our attention to contractive ring morphisms, we obtain
a completion functor for rings $\wedge \colon \SNRing \rightarrow \BRing$ which is
left adjoint to the forgetful functor $\BRing \rightarrow \SNRing$.

\begin{rem}
	Since the completion functor $\wedge \colon \SNAb \rightarrow \BAb$ is lax (even strong) monoidal, it
	extends by abstract nonsense to a completion functor $\wedge \colon \SNRing \rightarrow \BRing$
	(see \cref{subsec:alg-in-monoidal-cat}).
	The ring structure obtained on $\widehat{R}$ via the categorical construction coincides
	with the ring structure described in \cref{fig:extending-ring-multiplication-to-completion}.
\end{rem}

\begin{rem}
	Given a ring $R$, we can endow it with the trivial norm \eqref{eq:def-triv-norm}
	which is a ring norm on $R$ with respect to which $R$ is complete. Hence, we can
	think of the category $\Ring$ of rings as a full subcategory of $\BRing$.
	When thinking of a ring without a predefined norm as a Banach ring, we will always assume
	that it is endowed with the trivial norm.
\end{rem}

\subsection{Non-Archimedean Modules over Non-Archimedean Rings}
Since the category $\SNAb$ is monoidal, one can talk about module objects over algebra objects of $\SNAb$
(see \cref{subsec:modules-in-monoidal-cat}).
Let $\left( R, \nnorm_R \right)$ be a seminormed ring, i.e., an algebra object of $\SNAb$.
A \textbf{(non-Archimedean) seminormed} $R$-\textbf{module} is just a module object
over $\left( R, \nnorm_R \right)$ in $\SNAb$. Unwinding the definition, we see that a
seminormed $R$-module is a pair $\left( M, \nnorm_M \right)$ where $M$ is an $R$-module
and $\nnorm_M$ is a non-Archimedean seminorm on $M$ which satisfies certain compatibility conditions:

\begin{dfn} \label{def:seminorm-module}
	Let $\left( R, \nnorm_R \right)$ be a (semi)normed non-Archimedean ring and let $M$ be a left $R$-module.
	A \textbf{non-Archimedean module (semi)norm} on $M$ over $\left( R, \nnorm_R \right)$
	is a function $\nnorm_M \colon M \rightarrow \RPL$, such that:
	\begin{enumerate}
		\item $\nnorm_M$ is a non-Archimedean group (semi)norm on $(M,+)$.
		\item $\nnorm[r \cdot m]_M \leq \nnorm[r]_R \cdot \nnorm[m]_M$ for all $r \in R$ and $m \in M$.
	\end{enumerate}
\end{dfn}
Hence, a seminormed $R$-module is a pair $\left( M, \nnorm_M \right)$ where $M$ is an $R$-module
and $\nnorm_M$ is a non-Archimedean module (semi)norm over $\left( R, \nnorm_R \right)$.
The topology on $M$ induced by the associated pseudometric turns $M$
into a topological module over $R$. A normed module $\left( M, \nnorm_M \right)$ is called
a \textbf{Banach module} if $M$ with the induced metric is a complete metric space.

\begin{rem} \label{rem:banach-module-as-module-object}
	Note that for the same reasons as in \cref{rem:banach-ring-as-algebra-object},
	a Banach module over a Banach ring is the same thing as a module object over an algebra object
	in the category $\BAb$.
\end{rem}

Given a seminormed ring $R$, let us denote the category of seminormed $R$-modules with contractive
$R$-linear module homomorphisms as morphisms by $\SNMod[R]$ (or by $\SNMod[R, \nnorm_R]$ when
we want to emphasize the role of the seminorm on $R$).
Denote also by $\BMod[R]$ the full subcategory of $\SNMod[R]$ whose objects are Banach $R$-modules.
In parallel to \cref{sub:non-arch-abelian-groups}, we now discuss the properties of the categories
$\SNMod[R]$ and $\BMod[R]$. Note that we have $\SNAb = \SNMod[\ZZ, \trivnorm]$ and $\BAb = \BMod[\ZZ, \trivnorm]$
so that the categories $\SNMod[R]$ and $\BMod[R]$ (for an arbitrary Banach ring $R$) generalize
the categories $\SNAb$ and $\BAb$. Since the properties of $\SNMod[R]$ and $\BMod[R]$ are analogous
to the properties of $\SNAb$ and $\BAb$, we choose to be more succinct and leave out some details.
Although much of the following can be generalized by discussing
bimodules, we will assume that the ground ring $R$ is commutative.

\subsubsection{The category \texorpdfstring{$\SNMod[R]$}{SNMod(R)}}
The category $\SNMod[R]$ is preadditive and bicomplete. Limits (resp.\ colimits) are given by
limits (resp.\ colimits) of the underlying seminormed abelian groups endowed with the natural
$R$-action. We will use the notation we have set up for limits and colimits of seminormed groups
(see \cref{sub:cat-SNAb}) also for limits and colimits of seminormed $R$-modules.
In particular, the category $\SNMod[R]$ has finite biproducts, kernels and cokernels
and hence is pre-abelian but not abelian for the same reasons given in \cref{rem:SNAb-not-abelian}.

Given two seminormed $R$-modules $M$ and $N$, we can endow the algebraic tensor product $M \otimes_R N$,
considered as an $R$-module, with a module seminorm over $\left( R, \nnorm_R \right)$ by setting
\begin{equation}
	\nnorm[x]_{M \otimes_R N} \defeq \inf \Set{\max_{i \in I} \, \nnorm[m_i]_M \cdot \nnorm[n_i]_N}
	[x = \sum_{i \in I} m_i \otimes_R n_i, \, |I| < \infty] \label{eq:projective-tensor-seminorm-modules}
\end{equation}
for $x \in M \otimes_R N$. This is the same formula as for seminormed groups (see \cref{eq:projective-tensor-seminorm})
with $\otimes$ replaced by $\otimes_R$. Note that the seminorm on $R$ does not play a direct
role in the definition of $\nnorm_{M \otimes_R N}$. It plays an indirect role as the seminorms on $M$ and $N$ satisfy
a compatibility condition with $\nnorm_R$ and this guarantees that the seminorm on $M \otimes_R N$
given by \cref{eq:projective-tensor-seminorm-modules} is indeed a module seminorm compatible with $\nnorm_R$
with respect to the natural structure of $M \otimes_R N$ as an $R$-module.
The seminorm $\nnorm_{M \otimes_R N}$ is called the \textbf{non-Archimedean projective tensor seminorm} and the pair
$\left( M \otimes_R N, \nnorm_{M \otimes_R N} \right)$
is called the \textbf{seminormed tensor product} of $M$ and $N$ over $R$.

Similar to the situation for seminormed groups (see \cref{sub:cat-SNAb}), the seminormed tensor product of seminormed $R$-modules
can be characterized by a universal property involving bounded
$R$-bilinear maps. The seminormed tensor product $M \otimes_R N$ comes equipped with a canonical
contractive $R$-bilinear map $\otimes \colon M \times N \rightarrow M \otimes_R N$ characterized by
the following universal property: Given a seminormed $R$-module $L$ and a bounded $R$-bilinear map
$B \colon M \times N \rightharpoonup L$ there exists a unique bounded $R$-linear map
$\varphi_B \colon M \otimes_R N \rightharpoonup L$ with $\nnorm[\varphi_B] = \nnorm[B]$ such that
$\varphi_B \left( m \otimes_R n \right) = B(m,n)$ for all $m \in M$ and $n \in N$ (see
\cref{fig:seminormed-tensor-product-modules-universal-property}).
\begin{figure}[htb]
	\centering
	\begin{tikzcd}
		{M \times N} && {M \otimes_R N} \\
		&& L
		\arrow["\otimes_R", from=1-1, to=1-3]
		\arrow["\substack{\exists! \, \varphi_B \\ \textrm{bounded } R\textrm{-linear}}",
			dashed, harpoon, from=1-3, to=2-3]
		\arrow["\substack{B \textrm{ bounded} \\ R\textrm{-bilinear}}"', harpoon, from=1-1, to=2-3]
	\end{tikzcd}
	\caption{Universal property of the seminormed tensor product for seminormed $R$-modules.}
	\label{fig:seminormed-tensor-product-modules-universal-property}
\end{figure}

Given two $R$-linear maps $\varphi \colon M \rightharpoonup M'$ and $\psi \colon N \rightharpoonup N'$ between
seminormed $R$-modules, their tensor product
$\varphi \otimes_R \psi \colon M \otimes_R N \rightharpoonup M' \otimes_R N'$
is also bounded with $\nnorm[\varphi \otimes_R \psi] \leq \nnorm[\varphi] \cdot \nnorm[\psi]$.
The seminormed tensor product, together with the standard symmetries, associators and unitors, endow the category
$\SNMod[R]$ with the structure of a symmetric monoidal category whose unit is the seminormed ground ring
$\left( R, \nnorm_R \right)$, considered as a seminormed $\left( R, \nnorm_R \right)$-module.
It follows from the universal property of the seminormed tensor product that the symmetric monoidal category
$\SNMod[R]$ is closed with the internal hom object given by the $R$-module
\begin{equation}
	\InnHom{M}{N}[][R] \defeq
	\Set{f \colon M \rightharpoonup N}[f \textrm{ is a bounded } R\textrm{-linear map}]
	\label{eq:inn-hom-seminormed-R-modules}
\end{equation}
endowed with the operator seminorm (\cref{def:norm-of-morphism}) which is compatible with
$\nnorm_R$.

\begin{rem} \label{rem:tensor-product-of-seminormed-modules-internally}
	Since the category $\SNMod[R]$ is the category of module objects over a commutative algebra object
	$R$ of $\SNAb$ and the category $\SNAb$ is a cocomplete tensor category, it naturally inherits a
	monoidal structure by taking $M \otimes_R N$ to be the coequalizer in $\SNAb$ of the two obvious
	action morphisms $M \otimes R \otimes N \rightrightarrows M \otimes N$
	(see \cite[Proposition 4.1.10]{Brandenburg2014}). Since cokernels in $\SNAb$ are computed just like
	cokernels in $\Ab$, the resulting object is the algebraic
	tensor product $M \otimes_R N$, endowed with a quotient seminorm. One can verify that the resulting
	seminorm is precisely the one given by \cref{eq:projective-tensor-seminorm-modules}. The fact
	that the resulting symmetric monoidal structure $\otimes_R$ is closed with internal hom given by
	the \textbf{bounded} $R$-linear maps also follows by abstract nonsense (see
	\cite[Remark 4.1.18]{Brandenburg2014}).
\end{rem}

\subsubsection{The category \texorpdfstring{$\BMod[R]$}{BMod(R)}} \label{sec:banach-modules-over-banach-rings}
Before discussing the properties of $\BMod[R]$, we discuss the completion procedure for seminormed $R$-modules.
We have two variants:
\begin{enumerate}
	\item Given a seminormed ring $R$ and a seminormed $R$-module $M$, the completion $\widehat{M}$,
	      which is a priori an abelian Banach group, has a unique structure of Banach $R$-module with
	      respect to which the canonical morphism $\eta_M \colon M \rightarrow \widehat{M}$ becomes an
	      isometric $R$-linear map of $R$-modules. The action of $R$ on $\widehat{M}$ is obtained
	      by extending the action $R \times M \rightarrow M$ using
	      \cref{lm:extending-multilinear-maps-by-density} (see \cref{fig:extending-module-action-to-completion}).
	      Given a bounded $R$-linear map $f \colon M \rightharpoonup N$ from a seminormed $R$-module into
	      a Banach $R$-module, the adjunct $\tilde{f} \colon \widehat{M} \rightharpoonup N$ becomes
	      a bounded $R$-linear map of modules. In particular, this implies that given a bounded $R$-linear map
	      $f \colon M_1 \rightharpoonup M_2$ between seminormed $R$-modules,
	      the completion $\widehat{f} \colon \widehat{M_1} \rightharpoonup \widehat{M_2}$ also becomes a bounded
	      $R$-linear map of Banach modules.
	\item Given a seminormed ring $R$ and a Banach $R$-module $N$, we can extend the $R$-action
	      on $N$ to an $\widehat{R}$-action and endow $N$ with the structure of a Banach $\widehat{R}$-module.
	      The action of $\widehat{R}$ on $N$ is obtained by extending the action $R \times N \rightarrow N$ using
	      \cref{lm:extending-multilinear-maps-by-density} (see
	      \cref{fig:extending-module-action-on-Banach-module-to-base-ring-completion}). Given
	      a bounded $R$-linear map $f \colon N_1 \rightharpoonup N_2$ between Banach $R$-modules,
	      the same map is also $\widehat{R}$-linear.
\end{enumerate}

\begin{figure}[htb]
	\centering
	\begin{minipage}[b]{.5\textwidth}
		\centering
		\begin{tikzcd}
			{R \times M} && {R \times \widehat{M}} \\
			M && {\widehat{M}}
			\arrow["{\left( \id_R, \eta_M \right)}", from=1-1, to=1-3]
			\arrow[from=1-1, to=2-1]
			\arrow["{\eta_M}", from=2-1, to=2-3]
			\arrow["{\exists !}"', dashed, from=1-3, to=2-3]
		\end{tikzcd}
		\caption{\raggedright Extending $R$-action from $M$ to $\widehat{M}$.}
		\label{fig:extending-module-action-to-completion}
	\end{minipage}%
	\begin{minipage}[b]{.5\textwidth}
		\centering
		\begin{tikzcd}
			{R \times N} && {\widehat{R} \times N \phantom{\widehat{M}}} \\
			N && \phantom{\widehat{M}}
			\arrow["{\left( \eta_R, \id_N \right)}", from=1-1, to=1-3]
			\arrow[from=1-1, to=2-1]
			\arrow["{\phantom{\eta_M}}", phantom, from=2-1, to=2-3]
			\arrow["{\exists !}", dashed, from=1-3, to=2-1]
		\end{tikzcd}
		\caption{\raggedright Extending $R$-action to an $\widehat{R}$-action.}
		\label{fig:extending-module-action-on-Banach-module-to-base-ring-completion}
	\end{minipage}
\end{figure}

The two variants give us two completion functors
$\wedge \colon \SNMod[R] \rightarrow \BMod[R]$ and
$\wedge_R^{\hat{R}} \colon \BMod[R] \rightarrow \BMod (\widehat{R})$ whose composition
$\SNMod[R] \rightarrow \BMod (\widehat{R})$ is a ``relative'' completion functor.
The functor $\wedge_{\widehat{R}/R} \colon \BMod[R] \rightarrow \BMod (\widehat{R})$ is not very interesting
as it is an isomorphism of categories with inverse given by treating a Banach $\widehat{R}$-module $N$
as an $R$-module via restriction of scalars along the canonical morphism $\eta_R \colon R \rightarrow \widehat{R}$.
When working with Banach $R$-modules, we will always assume that the ground ring $R$ is Banach
and always work with the completion functor $\wedge \colon \SNMod[R] \rightarrow \BMod[R]$ which is
left adjoint to the natural forgetful functor $U \colon \BMod[R] \rightarrow \SNMod[R]$ with the map
$\eta_M \colon M \rightarrow \widehat{M}$ being the unit of the adjunction. Hence, we see that
$\BMod[R]$ is a reflective replete full subcategory of $\SNMod[R]$ with reflector $\wedge$.

The category $\BMod[R]$ is preadditive and bicomplete. Limits (resp.\ colimits) are given by
limits (resp.\ colimits) of the underlying abelian Banach groups endowed with the natural
$R$-action. Alternatively, since $\BMod[R]$ is a reflective replete full subcategory of $\SNMod[R]$,
limits in $\BMod[R]$ are the same as limits in $\SNMod[R]$ (and are automatically Banach) while
colimits are obtained by computing colimits in $\SNMod[R]$ and completing them.
We will use the notation we have set up for limits and colimits of Banach groups
(see \cref{sub:cat-BAb}) also for limits and colimits of Banach $R$-modules.
In particular, the category $\BMod[R]$ has finite biproducts, kernels and cokernels
and hence is pre-abelian but not abelian for the same reasons given in \cref{rem:BAb-not-abelian}.

Next, we discuss how to endow $\BMod[R]$ with the structure of a symmetric monoidal category.
Given two seminormed $R$-modules $M$ and $N$, the \textbf{complete tensor product} $M \cotimes_R N$ is
defined to be the (separated) completion $\extrawidehat{M \otimes_R N}$ of the algebraic tensor product
$M \otimes_R N$ with respect to the projective tensor seminorm given by
\cref{eq:projective-tensor-seminorm-modules}.
Given $m \in M$ and $n \in N$, we will denote by $m \cotimes_R n$ the image of $m \otimes_R n$
under the completion map $\eta \colon M \otimes_R N \rightarrow \extrawidehat{M \otimes_R N}$ and call
such elements \textbf{elementary tensors}. Note that unlike in the algebraic or seminormed case,
the complete tensor product $M \cotimes_R N$ is not spanned by $R$-linear combinations of
elementary tensors $m \cotimes_R n$. Instead, the normed $R$-module generated by elementary tensors is dense in
$M \cotimes_R N$. In other words, $M \cotimes_R N$ is generated by elementary tensors as a \textbf{Banach}
$R$\textbf{-module}.

\begin{rem}
	Note that even when $M$ and $N$ are normed, the algebraic tensor product $M \otimes_R N$
	endowed with the projective tensor seminorm given by \cref{eq:projective-tensor-seminorm-modules}
	need not be normed. Hence, we really need to take the \textit{separated}
	completion of $M \otimes_R N$ to obtain a Banach $R$-module. We have already seen an example in
	\cref{rem:tensor-product-normed-groups-only-seminormed}. For another example
	having a somewhat different flavor, we refer to \cref{sec:example-pows-several-variables-completion}.
\end{rem}

An important property of the complete tensor product is that, given two
seminormed $R$-modules $M$ and $N$, the canonical map
\begin{equation}
	M \cotimes_R N = \extrawidehat{M \otimes_R N} \xrightarrow[\cong]{\extrawidehat{\eta_M \otimes_R \eta_N}}
	\extrawidehat{\widehat{M} \otimes_R \widehat{N}} = \widehat{M} \cotimes_R \widehat{N}
	\label{eq:complete-tensor-product-modules-strong-monoidal}
\end{equation}
is an isometric isomorphism of Banach $R$-modules.\footnote{To construct the inverse map in the other direction,
	one can use \cref{lm:extending-multilinear-maps-by-density} which also works for $R$-modules and $R$-bilinear maps.}

\begin{figure}[htb]
	\centering
	\begin{tikzcd}
		{M \times N} && {M \otimes_R N} &&
		{M \cotimes_R N = \extrawidehat{M \otimes_R N}} \\
		&& L
		\arrow["\cotimes_R", curve={height=-25pt}, from=1-1, to=1-5]
		\arrow["\otimes_R", from=1-1, to=1-3]
		\arrow["\eta", from=1-3, to=1-5]
		\arrow["{\exists!}", dashed, harpoon, from=1-3, to=2-3]
		\arrow["\substack{\exists! \, \varphi_B \textrm{ bounded} \\ R\textrm{-linear} }",
			dashed, harpoon, from=1-5, to=2-3]
		\arrow["\substack{B \textrm{ bounded} \\ R\textrm{-bilinear}}"', harpoon, from=1-1, to=2-3]
	\end{tikzcd}
	\caption{Universal property of the complete tensor product of seminormed Banach $R$-modules.}
	\label{fig:complete-tensor-product-modules-universal-property}
\end{figure}
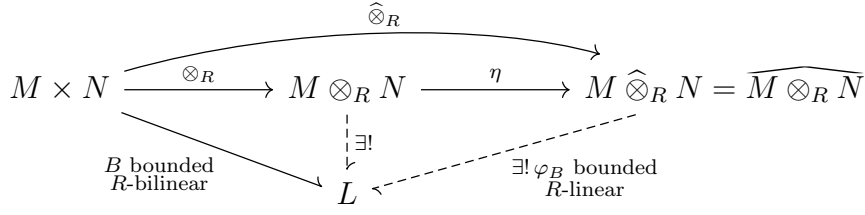

Similar to the seminormed tensor product of abelian groups, the complete tensor product $M \cotimes_R N$ is also
characterized by a universal property involving bounded $R$-bilinear maps into \textit{Banach} $R$-modules.
The Banach $R$-module $M \cotimes_R N$ comes equipped with a canonical contractive bilinear map
$\cotimes \colon M \times N \rightarrow M \cotimes_R N$ characterized by the
following universal property: Given a Banach $R$-module $L$ and a bounded $R$-bilinear map
$B \colon M \times N \rightharpoonup L$ there exists a unique bounded $R$-linear map
$\varphi_B \colon M \cotimes_R N \rightharpoonup L$ with $\nnorm[\varphi_B] = \nnorm[B]$ such that
$\varphi_B \left( m \cotimes_R n \right) = B(m,n)$ for all $m \in M$ and $n \in N$ (see
\cref{fig:complete-tensor-product-modules-universal-property}).

Given two bounded morphisms $\varphi \colon M \rightharpoonup M'$ and $\psi \colon N \rightharpoonup N'$ between
seminormed $R$-modules, their complete tensor product
$\varphi \cotimes_R \psi \colon M \cotimes_R N \rightharpoonup M' \cotimes_R N'$
is defined by $\varphi \cotimes_R \psi \defeq \extrawidehat{\varphi \otimes_R \psi}$ and we have
$\nnorm[\varphi \cotimes_R \psi] \leq \nnorm[\varphi] \cdot \nnorm[\psi]$. In particular, we see that the
complete tensor product construction gives us a bifunctor
\begin{equation*}
	\cotimes_R \colon \SNMod[R] \times \SNMod[R] \rightarrow \BMod[R].
\end{equation*}

Even though we have defined the complete tensor product for seminormed $R$-modules, we now restrict our
attention to the complete tensor product of Banach $R$-modules and obtain a bifunctor
$\cotimes_R \colon \BMod[R] \times \BMod[R] \rightarrow \BMod[R]$. Given three Banach $R$-modules $M,N,L$,
we can use the isomorphisms of \cref{eq:complete-tensor-product-modules-strong-monoidal} to obtain natural
associativity isomorphisms $\left( M \cotimes_R N \right) \cotimes_R L \cong M \cotimes_R \left( N \cotimes_R L \right)$
which act on elementary tensors by the expected formula $\left( m \cotimes_R n \right) \cotimes_R l \mapsto
	m \cotimes_R \left( n \cotimes_R l \right)$.
In addition, we also have natural isometric isomorphisms
\begin{align*}
	 & R \cotimes_R M \cong M \qquad               & r \cotimes_R m & \mapsto r \cdot m,      \\
	 & M \cotimes_R R \cong M \qquad               & m \cotimes_R r & \mapsto r \cdot m,      \\
	 & M \cotimes_R N \cong N \cotimes_R M, \qquad & m \cotimes_R n & \mapsto n \cotimes_R m.
\end{align*}
The bifunctor $\cotimes_R \colon \BMod[R] \times \BMod[R] \rightarrow \BMod[R]$, together with the
associativity isomorphisms, the unitors and the symmetry maps described above, endow the category $\BMod[R]$
with the structure of a symmetric monoidal category whose unit is the Banach ground ring $R$,
considered as a Banach $R$-module over itself.\footnote{Recall
	that we assume that $R$ is Banach, so this is the same unit of the monoidal structure on $\SNMod[R]$.}

Given two seminormed $R$-modules $M$ and $N$, the seminormed $R$-module $\InnHom{M}{N}[][R]$
of all bounded $R$-linear maps is Banach once $N$ is Banach. In particular,
$\InnHom{M}{N}[][R]$ is an object of $\BMod[R]$ when both $M$ and $N$ are objects of
$\BMod[R]$. It follows from the universal property of the complete tensor product $\cotimes_R$
that the symmetric monoidal category $\BMod[R]$ is closed with the internal hom object given by
$\InnHom{M}{N}[][R]$, the same internal hom object as in $\SNMod[R]$.

\begin{rem} \label{rem:tensor-product-of-Banach-modules-internally}
	Since the category $\BMod[R]$ is the category of module objects over a commutative algebra object
	$R$ of $\BAb$ and the category $\BAb$ is a cocomplete tensor category, it naturally inherits a
	monoidal structure by taking $M \cotimes_R N$ to be the coequalizer in $\BAb$ of the two obvious
	action morphisms $M \cotimes R \cotimes N \rightrightarrows M \cotimes N$
	(see \cite[Proposition 4.1.10]{Brandenburg2014}). The resulting object is given by
	\begin{equation*}
		\left( M \cotimes N \right) \Big/ \overline{
			\left< \left( r \cdot m \right) \cotimes n - m \cotimes \left( r \cdot n \right) \, \middle| \,
			r \in R, m \in M, n \in N \right>
		}
	\end{equation*}
	endowed with the quotient norm induced from the norm on $M \cotimes N$. One can verify that the resulting
	object is isomorphically isometric to the completion $\extrawidehat{M \otimes_R N} = M \cotimes_R N$.
	The fact that the resulting symmetric monoidal structure $\cotimes_R$ is closed with internal hom given by
	the bounded $R$-linear maps also follows by abstract nonsense (see \cite[Remark 4.1.18]{Brandenburg2014}).
\end{rem}

With respect to the monoidal structures defined on $\SNMod[R]$ and $\BMod[R]$, the adjunction
\begin{equation*}
	\wedge \colon \SNMod[R] \stackrel[]{\dashv}{\rightleftarrows} \BMod[R] \colon U
\end{equation*}
is naturally enhanced into a monoidal adjunction in which the left adjoint functor $\wedge$ is strong symmetric
monoidal via the coherence isomorphism of \cref{eq:complete-tensor-product-modules-strong-monoidal}
and the right adjoint forgetful functor $U$ is lax symmetric monoidal via the
map $M \otimes_R N \rightarrow M \cotimes_R N$.

\subsection{An Example} \label{sec:example-pows-several-variables-completion}

Given a commutative ring $R$ and formal variables $u,v$, we show that the complete tensor
product $\pows{R}[u] \cotimes_R \pows{R}[v]$ can be identified with $\pows{R}[u,v]$, with
the completion map being the natural multiplication map
$\pows{R}[u] \otimes_R \pows{R}[v] \rightarrow \pows{R}[u,v]$. We present an example of a ring $R$ for which
the natural multiplication map is not injective. In particular, this shows that
it is possible that $\pows{R}[u] \otimes_R \pows{R}[v]$ is only seminormed and not normed, so
that the projective tensor seminorm is only a seminorm, even though $\pows{R}[u],\pows{R}[v]$ are Banach.

Let $R$ be a commutative ring and let $u_1,\dots,u_n$ be formal variables. Recall that
the \textbf{order} of a formal power series
$p(u_1,\dots,u_n) = \sum_{r_1,\dots,r_n} a_{r_1,\dots,r_n} u_1^{r_1} \dots u_n^{r_n}$ is given by
\begin{equation*}
	\ord \left( p(u_1,\dots,u_n) \right) \defeq
	\begin{cases}
		\min \Set{r_1 + \dots + r_n}[a_{r_1,\dots,r_n} \neq 0] & p(u_1,\dots,u_n) \neq 0, \\
		+\infty                                                & p(u_1,\dots,u_n) = 0.
	\end{cases}
\end{equation*}
The ring of formal power series $\pows{R}[u_1,\dots,u_n]$, endowed with the norm $\nnorm[p(u)] \defeq 2^{-\ord p(u)}$,
is complete and is naturally a Banach $\left( R, \trivnorm \right)$-module and a Banach
$\left( R, \trivnorm \right)$-algebra.

\begin{lm} \label{lm:mult-map-is-completion}
	The natural multiplication map $\mu \colon \pows{R}[u] \otimes_R \pows{R}[v] \rightarrow \pows{R}[u,v]$
	is an isometry where the tensor product $\pows{R}[u] \otimes_R \pows{R}[v]$ is endowed with the projective
	tensor seminorm. Since $\pows{R}[u,v]$ is normed, complete and the image of $\mu$ is dense,
	this means that $\pows{R}[u,v] \cong \pows{R}[u] \cotimes_R \pows{R}[v]$ and $\mu$ can be identified with
	the natural isometry mapping $\pows{R}[u] \otimes_R \pows{R}[v]$ into the completion
	$\extrawidehat{\pows{R}[u] \otimes_R \pows{R}[v]}$.
\end{lm}
\begin{proof}
	Note that given $p(u) \in \pows{R}[u]$ and $q(v) \in \pows{R}[v]$ we have
	\begin{equation*}
		\ord_{\pows{R}[u]} p(u) + \ord_{\pows{R}[v]} q(v) \leq \ord_{\pows{R}[u,v]} \left( p(u) \cdot q(v) \right).
	\end{equation*}
	In terms of norms, we get the equivalent inequality
	\begin{equation*}
		\nnorm[p(u) \cdot q(v)]_{\pows{R}[u,v]} \leq \nnorm[p(u)]_{\pows{R}[u]} \cdot \nnorm[q(v)]_{\pows{R}[v]}.
	\end{equation*}
	Let $\xi \in \pows{R}[u] \otimes_R \pows{R}[v]$ and choose some decomposition
	$\xi = \sum_{k=0}^K p_k(u) \otimes q_k(v)$ where $p_k(u) \in \pows{R}[u]$ and $q_k(v) \in \pows{R}[v]$.
	Then
	\begin{equation*}
		\begin{aligned}
			\nnorm[\mu \left( \xi \right)]_{\pows{R}[u,v]} & =
			\nnorm[\sum_{k=0}^K p_k(u) \cdot q_k(v)]_{\pows{R}[u,v]} \leq
			\max_{0 \leq k \leq K} \nnorm[p_k(u) \cdot q_k(v)]_{\pows{R}[u,v]}
			\\
			                                               & \leq
			\max_{0 \leq k \leq K} \nnorm[p_k(u)]_{\pows{R}[u]} \cdot \nnorm[q_k(v)]_{\pows{R}[v]}.
		\end{aligned}
	\end{equation*}
	Since this holds for any decomposition of $\xi$, we see that
	$\nnorm[\mu \left( \xi \right)] \leq \nnorm[\xi]$.
	Next, let us show that if $l \geq 0$ and $\nnorm[\mu \left( \xi \right)] \leq 2^{-l}$ then we also have
	$\nnorm[\xi] \leq 2^{-l}$. This, together with the inequality
	$\nnorm[\mu \left( \xi \right)] \leq \nnorm[\xi]$ will immediately imply that in fact we have
	$\nnorm[\mu \left( \xi \right)] = \nnorm[\xi]$.
	Write
	\begin{equation*}
		p_k(u) = \sum_{i=0}^{\infty} a_i^k u^i, \qquad q_k(v) = \sum_{j=0}^{\infty} b_j^k v^j,
		\qquad \mu(\xi) = \sum_{i,j=0}^{\infty}
		\underbrace{\left( \sum_{k=0}^K a_i^k \cdot b_j^k \right)}_{c_{ij}} u^i v^j.
	\end{equation*}
	Let $l \geq 0$ and assume that $\nnorm[\mu \left( \xi \right)] \leq 2^{-l}$, or, equivalently,
	that $\ord_{\pows{R}[u,v]} \mu \left( \xi \right) \geq l$. This means that $c_{i,j} = 0$ whenever $i + j < l$.
	Write
	\begin{equation*}
		p_k(u) = \underbrace{\sum_{i=0}^l a_i^k u^i}_{p'_k(u)} +
		\underbrace{\sum_{i=l+1}^{\infty} a_i^k u^i}_{p''_k(u)}, \quad
		q_k(v) = \underbrace{\sum_{j=0}^l b_j^k v^j}_{q'_k(v)} +
		\underbrace{\sum_{j=l+1}^{\infty} b_j^k v^j}_{q''_k(v)}
	\end{equation*}
	so that $\nnorm[p''_k(u)]_{\pows{R}[u]}, \nnorm[q''_k(v)]_{\pows{R}[v]} \leq 2^{-(l+1)}$.
	Then we can rewrite $\xi$ as
	\begin{equation*}
		\begin{aligned}
			\xi ={} & \sum_{k=0}^K \left( p'_k(u) + p''_k(u) \right) \otimes \left( q'_k(v) + q''_k(v) \right)
			= \sum_{k = 0}^K \left( \sum_{i=0}^l a_i^k u^i \right) \otimes \left( \sum_{j=0}^l b_j^k v^j \right)
			\\
			        & +
			\sum_{k=0}^K p'_k(u) \otimes q''_k(v) +
			\sum_{k=0}^K p''_k(u) \otimes q'_k(v) + \sum_{k=0}^K p''_k(u) \otimes q''_k(v)
			\\
			={}     &
			\sum_{\substack{0 \leq i,j \leq l \\ i + j \geq l}} \left( c_{ij} u^i \right) \otimes v^j +
			\sum_{k=0}^K p'_k(u) \otimes q''_k(v) + \sum_{k=0}^K p''_k(u) \otimes q'_k(v) +
			\sum_{k=0}^K p''_k(u) \otimes q''_k(v).
		\end{aligned}
	\end{equation*}
	Hence, we have found a decomposition $\xi = \sum_s f_s(u) \otimes g_s(v)$ (possibly different from the initial
	one $\xi = \sum_k p_k(u) \otimes q_k(v)$) with respect to which
	$\nnorm[f_s(u)]_{\pows{R}[u]} \nnorm[g_s(v)]_{\pows{R}[v]} \leq 2^{-l}$ for all $s$.
	This shows that $\nnorm[\xi] \leq 2^{-l}$.
\end{proof}

\begin{rem}
	One can deduce \cref{lm:mult-map-is-completion} more abstractly as follows.
	Let us denote by $R_i$ the ring $R$, considered as an $R$-module and endowed with the norm
	$\nnorm[r]_i = 2^{-i}$ whenever $0 \neq r \in R_i$. Then each $\left( R_i, \nnorm_i \right)$ is an
	$(R,\trivnorm)$-module\footnote{Note that according to our definitions, $R_i$ is not a normed \textit{ring}
		when $i \geq 1$ since if $r,r' \in R_i$ and $r \cdot r' \neq 0$ then $\nnorm[r \cdot r']_i = 2^{-i}$ while
		$\nnorm[r]_i \cdot \nnorm[r']_i = 2^{-2i}$ and hence $\nnorm[r \cdot r']_i > \nnorm[r]_i \cdot \nnorm[r']_i$.
		It is however an $\left( R, \trivnorm \right)$ \textit{module}.}
	and the power series ring $\pows{R}[u]$ is
	isometric to the \textit{complete} direct sum
	$\widehat{\oplus}_{i=0}^{\infty} R_i$ of the modules $R_i$. In addition, one can easily
	verify that the multiplication map $R_i \otimes_R R_j \mapsto R_{i+j}$ is an isometry and hence
	$R_i \cotimes_R R_j = R_i \otimes R_j \cong R_{i+j}$. Hence, we have
	\begin{equation*}
		\pows{R}[u] \cotimes_R \pows{R}[v] \cong \left( \widehat{\oplus}_{i=0}^{\infty} R_i \right) \cotimes_R
		\left( \widehat{\oplus}_{j=0}^{\infty} R_j \right) \cong
		\widehat{\oplus}_{i,j=0}^{\infty} R_i \cotimes_R R_j \cong
		\widehat{\oplus}_{i,j=0}^{\infty} R_{i+j} \cong \pows{R}[u,v].
	\end{equation*}
\end{rem}

\begin{cor} \label{cor:kernel-mult-map-is-ker-seminorm}
	The kernel of the multiplication map $\mu \colon \pows{R}[u] \otimes_R \pows{R}[v] \rightarrow \pows{R}[u,v]$
	is equal to the kernel of the seminorm $\nnorm_{\pows{R}[u] \otimes_R \pows{R}[v]}$:
	\begin{equation*}
		\ker \mu = \Set{\xi \in \pows{R}[u] \otimes_R \pows{R}[v]}[{\nnorm[\xi]= 0}].
	\end{equation*}
	The seminorm $\nnorm_{\pows{R}[u] \otimes_R \pows{R}[v]}$ is a norm if and only if the map $\mu$ is injective.
	\qed
\end{cor}

\Cref{cor:kernel-mult-map-is-ker-seminorm} translates the question of whether the seminorm on
$\pows{R}[u] \otimes_R \pows{R}[v]$ is a norm into an algebraic question of whether the multiplication map $\mu$ is
injective. In general, it is known that the map $\mu$ might not be injective (see \cite{Goodearl1972}).
Our purpose is to provide an explicit example of $R$ and $\xi \neq 0$ which lies in the kernel of $\mu$.

How can we produce an element in the kernel of $\mu$? Let $a_i,b_i,a'_i,b'_i \in R$ for $0 \leq i < \infty$
be elements of $R$ which satisfy
\begin{equation} \label{eq:def-xi-cond}
	a_i \cdot b_j = a'_i \cdot b'_j \textrm{ for all } 0 \leq i,j < \infty.
\end{equation}
Consider the element $\xi \in \pows{R}[u] \otimes_R \pows{R}[v]$ given by
\begin{equation} \label{eq:def-xi-counter-example-power-series}
	\xi \defeq \left( \sum_{i=0}^{\infty} a_i u^i \right) \otimes_R \left( \sum_{j=0}^{\infty} b_j v^j \right) -
	\left( \sum_{i=0}^{\infty} a'_i u^i \right) \otimes_R \left( \sum_{j=0}^{\infty} b'_j v^j \right).
\end{equation}
Clearly we have
\begin{equation*}
	\mu \left( \xi \right) = \sum_{i,j=0}^{\infty} \left( a_i \cdot b_j \right) u^i v^j -
	\sum_{i,j=0}^{\infty} \left( a'_i \cdot b'_j \right) u^i v^j = 0
\end{equation*}
by our assumption and so $\nnorm[\xi] = \nnorm[\mu \left( \xi \right)] = 0$.
Hence, if we manage to choose $a_i,b_i,a'_i,b'_i$ satisfying \eqref{eq:def-xi-cond}
\textbf{in such a way that} $\xi \neq 0$, we will have our desired counterexample. For example, we might
try and take $a'_i = t a_i$ and $b'_j = t^{-1} b_j$ for an invertible $t \in R$, but this clearly won't work
since if $t \in R$ then we can distribute $t$ across $\otimes_R$ and cancel it.

The following example, taken from \cite[Appendix 5.3]{duchamp:hal-02943601} and adapted to our context, shows that we
can make the idea above work. We start with the technical \cref{lm:when-tensor-product-nonzero} and construct
the example in \cref{lm:universal-xi}.

\begin{lm} \label{lm:when-tensor-product-nonzero}
	Let $\mathbbm{k}$ be a commutative ring and let $R$ be a $\mathbbm{k}$-algebra. Let $M$ be a right $R$-module
	and $N$ be a left $R$-module. Let $\xi_R = \sum_{k=1}^K m_k \otimes_R n_k \in M \otimes_R N$ where $m_k \in M$
	and $n_k \in N$. Given a $\mathbbm{k}$-subalgebra $S$ of $R$, set
	$\xi_S \defeq \sum_{k=1}^K m_k \otimes_S n_k \in M \otimes_S N$ where we consider $M$ and $N$ as $S$-modules
	by scalar restriction. Then the following statements are equivalent:
	\begin{enumerate}
		\item $\xi_R \neq 0$.
		\item $\xi_S \neq 0$ for all \textbf{finitely generated} $\mathbbm{k}$-subalgebra $S$ of $R$.
	\end{enumerate}
\end{lm}
\begin{proof}
	Given $\mathbbm{k}$-subalgebras $S_1,S_2$ of $R$ such that $S_1 \subseteq S_2$, let
	$\varphi_{S_1,S_2} \colon M \otimes_{S_1} N \rightarrow M \otimes_{S_2} N$
	be the natural surjective map of $\mathbbm{k}$-modules satisfying
	$\varphi_{S_1,S_2} \left( m \otimes_{S_1} n \right) = m \otimes_{S_2} n$. Then
	$\varphi_{S_1,S_2} \left( \xi_{S_1} \right) = \xi_{S_2}$.

	Assume first that $\xi_S = 0$ for some $\mathbbm{k}$-subalgebra $S$ of $R$. Then
	$\xi_R = \varphi_{S,R} \left( \xi_S \right) = 0$.
	Conversely, assume that $\xi_R = 0$ and consider $\xi_{\mathbbm{k}} \in M \otimes_{\mathbbm{k}} N$.
	We have $0 = \xi_R = \varphi_{\mathbbm{k},R} \left( \xi_{\mathbbm{k}} \right)$.
	The kernel of $\varphi_{\mathbbm{k},R}$ is spanned by elements of
	the form $\left( mr \right) \otimes_{\mathbbm{k}} n - m \otimes_{\mathbbm{k}} \left( rn \right)$
	where $m \in M, n \in N$ and $r \in R$. Hence, we can find \textbf{finitely} many elements
	$m'_l \in M, n'_l \in N, r_l \in R, \lambda_l \in \mathbbm{k}$ such that
	\begin{equation*}
		\xi_{\mathbbm{k}} = \sum_{l=1}^L \lambda_l \left( \left( m'_l r_l \right) \otimes_{\mathbbm{k}} n'_l
		- m'_l \otimes_{\mathbbm{k}} \left( r_l n'_l \right) \right).
	\end{equation*}
	Let $S$ be the $\mathbbm{k}$-subalgebra of $R$ generated by $r_1, \dots, r_L$. Then
	\begin{equation*}
		\xi_S = \varphi_{\mathbbm{k}, S} \left( \xi_{\mathbbm{k}} \right) =
		\sum_{l=1}^L \lambda_l \left( \left( m'_l r_l \right) \otimes_{S} n'_l
		- m'_l \otimes_{S} \left( r_l n'_l \right) \right) = \sum_{l=1}^L \lambda_l \cdot 0 = 0
	\end{equation*}
	since $r_1,\dots,r_L \in S$.
\end{proof}

\begin{lm} \label{lm:universal-xi}
	Let $\mathbbm{k}$ be a field and let $t, (x_i)_{i=0}^{\infty}, (y_j)_{j=0}^{\infty}$ be formal variables.
	Let $R$ be the $\mathbbm{k}$-subalgebra of $\mathbbm{k}(t)[x_i,y_j]$ given by
	\begin{equation} \label{def:R-xi-yj-txi}
		R = \left< x_i, y_j, tx_i, t^{-1}y_j \right>_{\mathbbm{k}-\textrm{alg}} \subseteq \mathbbm{k}(t)[x_i,y_j],
	\end{equation}
	and let $\xi \in \pows{R}[u] \otimes_R \pows{R}[v]$ be the element
	\begin{equation*}
		\xi \defeq
		\left( \sum_{i=0}^{\infty} tx_i u^i \right) \otimes_R \left( \sum_{j=0}^{\infty} t^{-1}y_j v^j \right) -
		\left( \sum_{i=0}^{\infty} x_i u^i \right) \otimes_R \left( \sum_{j=0}^{\infty} y_j v^j \right).
	\end{equation*}
	Then $\xi \neq 0$.
\end{lm}
\begin{proof}
	Given a $\mathbbm{k}$-subalgebra $S$ of $R$, set
	\begin{equation*}
		\xi_S \defeq
		\left( \sum_{i=0}^{\infty} tx_i u^i \right) \otimes_S \left( \sum_{j=0}^{\infty} t^{-1}y_j v^j \right) -
		\left( \sum_{i=0}^{\infty} x_i u^i \right) \otimes_S \left( \sum_{j=0}^{\infty} y_j v^j \right)
		\in \pows{R}[u] \otimes_S \pows{R}[v]
	\end{equation*}
	By \cref{lm:when-tensor-product-nonzero}, it is enough to show that $\xi_S \neq 0$ for any finitely
	generated $\mathbbm{k}$-subalgebra $S$ of $R$.

	Let $S$ be a finitely generated $\mathbbm{k}$-subalgebra of $R$. Since each element of $R$
	is a polynomial in finitely many of the generators $x_i,y_j,tx_k,t^{-1}y_l$ and $S$ is finitely generated,
	one can choose $N \geq 0$ such that $S$ is also a subalgebra of
	$\mathbbm{k}(t)[(x_i)_{i \neq N}, (y_i)_{i \neq N}]$.
	Let
	\begin{equation*}
		\hat{p} \colon \mathbbm{k}(t)[x_i,y_i] \rightarrow \mathbbm{k}(t)[x,y]
	\end{equation*}
	be the $\mathbbm{k}(t)$-algebra
	homomorphism which maps $(x_i)_{i \neq N}, (y_i)_{i \neq N}$ to $0$ and $x_N,y_N$ to $x,y$ respectively.
	The image $\overline{R} \defeq \hat{p} \left( R \right)$ of $R$ under $\hat{p}$ is the
	$\mathbbm{k}$-subalgebra of $\mathbbm{k}(t)[x,y]$ generated by $x,y,tx, t^{-1}y$, the ``single variable''
	version of the algebra \eqref{def:R-xi-yj-txi}. A straightforward induction shows that the product of
	finitely many
	elements from the set $\Set{x,y,tx, t^{-1}y}$ results in a monomial of the form $t^q x^r y^s$ where
	$r,s \geq 0, q \in \ZZ$ and $-s \leq q \leq r$. Conversely, any such monomial can be obtained as a product of
	elements from the set $\Set{x,y,tx,t^{-1}y}$. Hence, $\overline{R}$ is generated as a free $\mathbbm{k}$-module
	by the monomials $\left( t^q x^r y^s \right)_{-s \leq q \leq r}$. In particular, we see that
	$\overline{R} \cap \mathbbm{k}(t) = \mathbbm{k}$. Since
	$S \subseteq \mathbbm{k}(t)[(x_i)_{i \neq N}, (y_i)_{i \neq N}]$ and $S \subseteq R$, we see that
	\begin{equation*}
		\hat{p} \left( S \right) \subseteq \hat{p} \left( \mathbbm{k}(t)[(x_i)_{i \neq N}, (y_i)_{i \neq N}] \right)
		\cap \hat{p} \left( R \right) = \mathbbm{k}(t) \cap \overline{R} = \mathbbm{k} \implies
		\hat{p} \left( S \right) = \mathbbm{k}
	\end{equation*}
	The situation is summarized in the following diagram:
	\[\begin{tikzcd}
			{\mathbbm{k}(t)[(x_i)_{i\neq N},(y_i)_{i \neq N}]} &&& {\mathbbm{k}(t)[x_i,y_i]} \\
			& S & {R } \\
			& \mathbbm{k} & {\overline{R}} \\
			{\mathbbm{k}(t)} &&& {\mathbbm{k}(t)[x,y]}
			\arrow["i", hook, from=2-2, to=2-3]
			\arrow["{\hat{p}}", two heads, from=1-4, to=4-4]
			\arrow["p", two heads, from=2-3, to=3-3]
			\arrow[hook, from=3-3, to=4-4]
			\arrow[hook, from=2-3, to=1-4]
			\arrow["q"', two heads, from=2-2, to=3-2]
			\arrow["j"', hook, from=3-2, to=3-3]
			\arrow[hook, from=2-2, to=1-1]
			\arrow["{\hat{q}}"', two heads, from=1-1, to=4-1]
			\arrow[hook, from=3-2, to=4-1]
			\arrow["{\hat{i}}", hook, from=1-1, to=1-4]
			\arrow["{\hat{j}}"', hook, from=4-1, to=4-4]
		\end{tikzcd}\]
	In the diagram above, hooked arrows are subset inclusions and the maps $p,q,\hat{q}$ are restrictions
	of the map $\hat{p}$ discussed above.

	Let us define a map
	$\varphi \colon \pows{R}[u] \otimes_S \pows{R}[v] \rightarrow \overline{R} \otimes_{\mathbbm{k}} \overline{R}$
	by
	\begin{equation*}
		\varphi \left( \left( \sum_{i=0}^{\infty} r_i u^i \right) \otimes_S
		\left( \sum_{j=0}^{\infty} r'_j v^j \right) \right) \defeq p(r_N) \otimes_{\mathbbm{k}} p(r'_N).
	\end{equation*}
	The fact that $\varphi$ is well-defined follows from the commutativity of the inner square in the diagram
	above. To see this more explicitly, let $s \in S$. Then
	\begin{equation*}
		\begin{aligned}
			\varphi \left( \left( s \cdot \left( \sum_{i=0}^{\infty} r_i u^i \right) \right) \otimes_S
			\left( \sum_{j=0}^{\infty} r'_j v^j \right) \right) & =
			\varphi \left( \left( \sum_{i=0}^{\infty} \left( i(s) \cdot r_i \right) u^i \right) \otimes_S
			\left( \sum_{j=0}^{\infty} r'_j v^j \right) \right)
			\\
			                                                    & =
			p(i \left( s \right) \cdot r_N) \otimes_{\mathbbm{k}} p(r'_N)
			\\
			                                                    & =
			\left( \left( p \circ i \right) \left( s \right) \cdot p(r_N) \right) \otimes_{\mathbbm{k}} p(r'_N)
			\\
			                                                    & =
			\left( \left( j \circ q \right) \left( s \right) \cdot p(r_N) \right) \otimes_{\mathbbm{k}} p(r'_N)
			\\
			                                                    & =
			p(r_N) \otimes_{\mathbbm{k}} \left( \left( j \circ q \right) \left( s \right) \cdot p(r'_N) \right)
			\\
			                                                    & =
			p(r_N) \otimes_{\mathbbm{k}} \left( \left( p \circ i \right) \left( s \right) \cdot p(r'_N) \right)
			\\
			                                                    & =
			p(r_N) \otimes_{\mathbbm{k}}  p \left( i(s) \cdot r'_N \right)
			\\
			                                                    & =
			\varphi \left( \left( \sum_{i=0}^{\infty} r_i u^i \right) \otimes_S
			\left( \sum_{j=0}^{\infty} \left( i(s) \cdot r'_j \right) v^j \right) \right)
			\\
			                                                    & =
			\varphi \left( \left( \sum_{i=0}^{\infty} r_i u^i \right) \otimes_S
			\left( s \cdot \left( \sum_{j=0}^{\infty} r'_j v^j \right) \right) \right).
		\end{aligned}
	\end{equation*}
	Applying $\varphi$ to $\xi_S$, we see that
	\begin{equation*}
		\varphi \left( \xi_S \right) = p \left( tx_N \right) \otimes_{\mathbbm{k}} p \left( t^{-1}y_N \right)
		- p \left( x_N \right) \otimes_{\mathbbm{k}} p \left( y_N \right) =
		\left( tx \right) \otimes_{\mathbbm{k}} \left( t^{-1} y \right) - x \otimes_{\mathbbm{k}} y \neq 0
	\end{equation*}
	since $tx,t^{-1}y,x,y \in \overline{R}$ are $\mathbbm{k}$-linearly independent. This shows that $\xi_S \neq 0$.
\end{proof}

\begin{rem}
	The construction of the ring $R$ and the element $\xi$ of \cref{lm:universal-xi} is natural if one tries
	to find a ``universal'' $\xi$ satisfying \cref{eq:def-xi-cond} and use it to demonstrate the
	obstruction to the injectivity of the multiplication map.
\end{rem}

\begin{rem}
	We can show directly that the element $\xi$ defined by \cref{eq:def-xi-counter-example-power-series} satisfies
	$\nnorm[\xi] = 0$ without appealing to \cref{lm:mult-map-is-completion} as follows:

	Let $N \geq 0$. Since the tensor product distributes over finite sums, we can rewrite the first term of
	\cref{eq:def-xi-counter-example-power-series} as
	\begin{equation*}
		\begin{aligned}
			\left( \sum_{i=0}^{\infty} a_i u^i \right) \otimes_R \left( \sum_{j=0}^{\infty} b_j v^j \right)
			={} &
			\left( \sum_{i=0}^{N} a_i u^i + \sum_{i=N+1}^{\infty} a_i u^i \right) \otimes_R
			\left( \sum_{j=0}^{N} b_j v^j + \sum_{j=N+1}^{\infty} b_j v^j \right)
			\\
			={} &
			\sum_{i,j=0}^{N} \left( a_i \cdot b_j \right) \cdot \left( u^i \otimes_R v^j \right) +
			\sum_{i=0}^N u^i \otimes_R \left( \sum_{j=N+1}^{\infty} \left( a_i \cdot b_j \right) v^j \right)
			\\
			    & +
			\sum_{j=0}^N \left( \sum_{i=N+1}^{\infty} \left( a_i \cdot b_j \right) u^i \right) \otimes_R v^j
			\\
			    & +
			\left( \sum_{i=N+1}^{\infty} a_i u^i \right) \otimes_R \left( \sum_{j=N+1}^{\infty} b_j v^j \right).
		\end{aligned}
	\end{equation*}
	Performing the same decomposition for the second summand of \cref{eq:def-xi-counter-example-power-series}
	and subtracting, we see that $\xi$ can also be written as
	\begin{equation*}
		\xi = \left( \sum_{i=N+1}^{\infty} a_i u^i \right) \otimes_R \left( \sum_{j=N+1}^{\infty} b_j v^j \right) -
		\left( \sum_{i=N+1}^{\infty} a'_i u^i \right) \otimes_R
		\left( \sum_{j=N+1}^{\infty} b'_j v^j \right).
	\end{equation*}
	Hence, we have
	\begin{equation*}
		\begin{aligned}
			\nnorm[\xi] & \leq \max
			\Set{\nnorm[\sum_{i=N+1}^{\infty} a_i u^i]_{\pows{R}[u]} \cdot \nnorm[\sum_{j=N+1}^{\infty} b_j
					v^j]_{\pows{R}[v]},
				\nnorm[\sum_{i=N+1}^{\infty} a'_i u^i]_{\pows{R}[u]} \cdot \nnorm[\sum_{j=N+1}^{\infty} b'_j v^j]_{\pows{R}[v]}}
			\\
			            & \leq \frac{1}{2^{2N+2}}
		\end{aligned}
	\end{equation*}
	for all $N \geq 0$ which shows that $\nnorm[\xi] = 0$.
\end{rem}

\section{Banach Bicomplexes with Contractible Rows} \label{appendix:bicomplexes}

\begin{figure}[H]
	\begin{tikzcd}
		& \vdots & \vdots \\
		\cdots & {C_0^1} & {C_1^1} & \cdots \\
		\cdots & {C_0^0} & {C_1^0} & \cdots \\
		& \vdots & \vdots
		\arrow["b"', from=3-3, to=3-2]
		\arrow["b"', from=2-3, to=2-2]
		\arrow["\delta", from=3-2, to=2-2]
		\arrow["\delta", from=3-3, to=2-3]
		\arrow["b"', from=3-2, to=3-1]
		\arrow["b"', from=3-4, to=3-3]
		\arrow["\delta", from=4-2, to=3-2]
		\arrow["\delta", from=4-3, to=3-3]
		\arrow["\delta", from=2-2, to=1-2]
		\arrow["\delta", from=2-3, to=1-3]
		\arrow["b"', from=2-4, to=2-3]
		\arrow["b"', from=2-2, to=2-1]
	\end{tikzcd}
	\caption{General Bicomplex.}
	\label{fig:general-bicomplex}
\end{figure}

In what follows, we will work with $\ZZ^2$-graded objects $C^{*}_{*}$ where
elements of $C^{y}_{x}$ are of degree $(x,y)$. When drawing bigraded objects, we will place $C_x^y$ at the
$(x,y)$-coordinate of the plane, as in \cref{fig:general-bicomplex}. Objects $C^{*}$
which are $\ZZ$-graded (with an upper index) will be identified with bigraded objects as necessary by placing them at the zeroth column, i.e., by setting $C^{*}_0 = C^{*}$ and $C^{*}_i = 0$ for $i \neq 0$.

Let $\mathcal{R} = \left( R ,d \right)$ be a differential graded-commutative Banach $\mathbbm{k}$-algebra.
A \textbf{Banach bicomplex} $\mathcal{C} = \left( C, b, \delta \right)$
\textbf{over} $\mathcal{R}$ consists of a $\ZZ^2$-graded Banach
$R$-module $C = C^{*}_{*}$, equipped with a horizontal differential
$b \colon C^{*}_{*} \rightharpoonup C^{*}_{* - 1}$ of degree $(-1,0)$
and a vertical differential $\delta \colon C^{*}_{*} \rightharpoonup C^{* + 1}_{*}$ of degree $(0,1)$ which
satisfy $\nnorm[b], \nnorm[\delta] \leq 1$. The horizontal differential $b$ is required to be $R$-linear while the vertical differential $\delta$ is required to be a module derivation over $d$. In addition, we require that
$b^2 = \delta^2 = \left[ b, \delta \right] = 0$.

When working with bigraded objects, there are two equivalent pairings one can work with: The inner
product parity form $\braidop_1$ or the total degree parity form $\braidop_2$ (see
\cref{appendix:parity-forms-equiv}).
The choice dictates whether we work with commuting or anticommuting differentials and has several
other implications. Our purpose is to describe the complete direct sum totalization of a bicomplex
and give a criterion for the contractibility of the total complex when its rows are contractible.

Since we use both conventions in this work, we will first describe the constructions and results
working with $\braidop_2$ and the state the necessary modifications needed when working with
$\braidop_1$.

\subsection{Banach Bicomplexes with Anticommuting Differentials}
\label{subsec:bicomplex-anticommuting-differentials}
Assume we choose to work with the total degree parity form $\braidop_2$ given by \cref{eq:parity-total-degree}.
This choice has several implications on the form of our bicomplexes which is worthwhile to make explicit.
First, since
\begin{equation*}
	\left[ b, \delta \right] = b \delta - (-1)^{\braid{(-1,0)}{(0,1)}_2} \delta b =
	b \delta + \delta b = 0,
\end{equation*}
we see that our differentials in fact anticommute. The vertical differential $\delta$ is a module derivation over $d$ and so satisfies
\begin{equation} \label{eq:vertical-differential-derivation}
	\delta \left( r \cdot c \right) =
	d \left( r \right) \cdot c + (-1)^{\braid{(0,1)}{(0,\degb{r})}_2} r \cdot \delta \left( c \right) =
	d \left( r \right) \cdot c + (-1)^{\degb{r}} r \cdot \delta \left( c \right) \in C_i^{j+l+1}
\end{equation}
for $r \in R^l$ and $c \in C_i^j$. In particular, we see that each column
$\mathcal{C}_i = \left( C_i^{*}, \delta \right)$
of the bicomplex is a differential graded Banach $\mathcal{R}$-module.
The fact that the horizontal differential $b$ is $R$-linear
translates into the identity
\begin{equation} \label{eq:horizontal-differential-linear-pairing-2}
	b \left( r \cdot c \right) = (-1)^{\braid{(-1,0)}{(0,\degb{r})}_2} r \cdot b \left( c \right) =
	(-1)^{\degb{r}} r \cdot b \left( c \right) \in C_{i-1}^{j+l}
\end{equation}
for $r \in R^l$ and $c \in C_i^j$. Note that the rows $C_{*}^j$ are \textit{not} graded $R$-modules:
multiplying an element in $j$-th row by $r \in R^l$ gives an element in the $(j + l)$-th row. However,
each row $\mathcal{C}^j = \left( C_{*}^j, b \right)$ is a homological chain complex of Banach $\mathbbm{k}$-modules.

The \textbf{complete direct sum totalization} of a $\ZZ^2$-graded Banach $R$-module $C$
is the $\ZZ$-graded Banach $R$-module given by
\begin{equation} \label{eq:tot-com-def-parity-2-module}
	\totc{C}[d][][\coplus] = \totc{C}[d][][\coplus, \braidop_2] \defeq
	\mathop{\coplus} \limits_{i \in \ZZ} C_i^{i+d}.
\end{equation}
A general element $x \in \totc{C}[d][][\coplus]$ will be denoted by $x = \sum_{i \in \ZZ} c_i^{i+d}$,
where $c_i^{i+d} \in C_i^{i+d}$ with $\nnorm[c_i^{i+d}] \xrightarrow[\abs{i} \to \infty]{} 0$.
The $R$-action on $\totc{C}[][][\coplus]$ is the natural one, i.e., we have
\begin{equation*}
	r^j \cdot \left( \sum_{i \in \ZZ} c_i^{i+d} \right) =
	\sum_{i \in \ZZ} \left( r^j \cdot c_i^{i+d} \right) \in \totc{C}[d+j][][\coplus]
\end{equation*}
for $r^j \in R^j$.
Given two $\ZZ^2$-graded Banach $R$-modules $C$ and $D$ and a graded bounded $\mathbbm{k}$-linear map
$f \colon C \rightharpoonup D$ of degree $(a,b)$, we define the induced map
\begin{equation*}
	f^{\totl} = \totc{f}[][][] = \totc{f}[][][\braidop_2] \colon \totc{C}[][][\coplus] \rightharpoonup \totc{D}[][][\coplus]
\end{equation*}
by
\begin{equation} \label{eq:tot-graded-maps-braid-op-2}
	\totc{f}[][][] \left( \sum_{i \in \ZZ} c_i^{i+d} \right) \defeq
	\sum_{i \in \ZZ} f \left( c_i^{i+d} \right).
\end{equation}
Since we work with the parity form $\braidop_2$, the definition above involves no signs
and guarantees that if $f$ is $R$-linear then $\totc{f}[][][]$ is also $R$-linear of degree $b - a$,
and that if $C = D$ and $\delta \colon C \rightharpoonup C$ is a module derivation over $d$
then $\totc{\delta}[][][]$ is also a module derivation over $d$.

\begin{rem} \label{rem:identification-tot-2}
	Given an element $c \in C$, which by our conventions is always a homogeneous
	element $c \in C_i^j$ with degree $(i,j)$, we can also think of $c$ as
	an element of $\totc{C}[][][\coplus]$ of degree $j - i$, i.e.,
	an element $c \in \totc{C}[j-i][][\coplus] = \coplus_{r \in \ZZ} C_r^{r + \left( j - i \right)}$.
	Since we work in the graded setting in which the degree of elements plays a role in contributing
	various signs, this could a priori cause various ambiguities but is in fact quite harmless.

	To expand, let us temporarily denote the element $c$, considered as an element of
	$\totc{C}[][][\coplus]$, by $\totc{c}[][][]$. Let $f \colon C \rightharpoonup D$
	be a graded bounded map of degree $(a,b)$. Then we have
	\begin{equation*}
		r \cdot \totc{c}[][][] = \totc{rc}[][][], \qquad
		\totc{f}[][][] \left( \totc{c}[][][] \right) = \totc{f \left( c \right)}[][][]
	\end{equation*}
	as both the $R$-action and the totalization of maps on $\totc{C}[][][\coplus]$ is not twisted by a sign.
	Note also that the condition of being $R$-linear looks the same whether we work on $C$
	or on the totalization $\totc{C}[][][\coplus]$ as
	\begin{gather*}
		f \left( r \cdot c \right) =
		(-1)^{\braid{\left( a,b \right)}{\left( 0,\degb{r} \right)}_2} r \cdot f \left( c \right) =
		(-1)^{\degb{r} \left( b - a \right)} r \cdot f \left( c \right),
		\\
		\totc{f}[][][] \left( r \cdot \totc{c}[][][] \right) = (-1)^{\degb{r} \cdot \degb{\totc{f}[][][]}} r \cdot
		\totc{f}[][][] \left( \totc{c}[][][] \right) = (-1)^{\degb{r} \left( b - a \right)} r \cdot \totc{f}[][][] \left( \totc{c}[][][] \right)
	\end{gather*}
	and the same holds for the condition of being a module derivation over $d$.

	Whenever working with Banach bicomplexes with the parity form $\braidop_2$,
	we will freely identify $c$ with $\totc{c}[][][]$ and $f$ with $\totc{f}[][][]$, relying on
	context to understand whether we work on the original bigraded $R$-module $C$ or on its totalization
	$\totc{C}[][][\coplus]$. See also \cref{rem:identification-tot-1}.
\end{rem}

Given a bicomplex $\mathcal{C} = \left( C_{*}^{*}, b, \delta \right)$, the
\textbf{complete direct sum totalization} of the bicomplex $\mathcal{C}$ is the
differential $\ZZ$-graded Banach $\mathcal{R}$-module given by
\begin{equation} \label{eq:tot-com-def-parity-2}
	\totc{\mathcal{C}}[][][\coplus] \defeq
	\left(
	\totc{C}[][][\coplus],
	D = b + \delta
	\right).
\end{equation}
Note that $D$ is a degree one derivation over $d$ and we have $D^2 = 0$ so
$\totc{\mathcal{C}}[][][\coplus]$ is indeed a differential graded Banach $\mathcal{R}$-module.

The following lemma is a version of the standard staircase argument adapted to our context which states that
under an appropriate convergence condition, the total complex is contractible when each row is contractible:

\begin{lm} \label{lm:total-bicomplex-contractible-R-linear-contraction-2}
	Let $\mathcal{C} = (C_{*}^{*}, b, \delta)$ be a Banach bicomplex and assume that the rows
	$\left( C_{*}^j, b \right)$
	are contractible with contracting $\mathbbm{k}$-linear homotopies
	$h^j \colon C_{*}^j \rightharpoonup C_{*+1}^j$ which satisfy $\nnorm[h^j] \leq 1$ for all $j \in \ZZ$.
	Assume also that the rows homotopies $h^j$ are compatible with
	the $R$-action on $C$, namely, assume that\footnote{The reason
		for the sign in \cref{eq:rows-homotopies-compatability-2} is to guarantee that $h$ defines an $R$-linear map
		of degree $(1,0)$ on $C_{*}^{*}$ and that it induces an $R$-linear map of degree $-1$ on the totalization.}
	\begin{equation} \label{eq:rows-homotopies-compatability-2}
		h^{j + k} \left( r \cdot c \right) = (-1)^{k} r \cdot h^{j} \left( c \right)
	\end{equation}
	whenever $c \in C_i^j$ and $r \in R^k$. Assume also that for each $c \in C_i^j$ we have
	\begin{equation}
		\left( h \delta + \delta h \right)^n (c) \to 0. \label{eq:bicomplex-convergence-cond-2}
	\end{equation}
	Then the infinite sum
	\begin{equation} \label{eq:bicomplex-R-linear-contraction-2}
		H \defeq \sum_{n=0}^{\infty} (-1)^n \left(  h \delta + \delta h \right)^n h
	\end{equation}
	converges pointwise and defines an $R$-linear degree $-1$ map
	$H \colon \totc{C}[*][][\coplus] \rightharpoonup \totc{C}[* - 1][][\coplus]$
	which satisfies $DH + HD = \idd$ and $\nnorm[H] \leq 1$. In particular,
	$\totc{\mathcal{C}}[][][\coplus]$ is contractible as a differential graded module over $\mathcal{R}$.

	When the row homotopies $h^j$ satisfy $\left( h^j \right)^2 = 0$ for all $j \in \ZZ$, the contraction $H$ satisfies $H^2 = 0$.
\end{lm}
\begin{proof}
	First, note that since $\nnorm[h^j] \leq 1$ for all $j \in \ZZ$, the maps $h^j$ extend to a well-defined
	bounded map $h \colon \totc{C}[*][][\coplus] \rightharpoonup \totc{C}[* - 1][][\coplus]$
	on the totalization with	$\nnorm[h] \leq 1$. Condition \eqref{eq:rows-homotopies-compatability-2} guarantees
	that $h$ is $R$-linear of degree $-1$ on the total complex. We have $bh + hb = \idd$ on
	$\totc{C}[][][\coplus]$ and so $\left( \totc{C}[][][\coplus], b \right)$
	is contractible.

	The total complex
	$\totc{\mathcal{C}}[][][\coplus] = \left( \totc{C}[][][\coplus], b + \delta \right)$
	is a perturbation of the contractible complex $\left( \totc{\mathcal{C}}[][][\coplus], b \right)$ by $\delta$,
	so by \cref{lm:trivial-perturbation-lemma}
	it is enough to show that $\idd + h \delta + \delta h$ is invertible. The map
	$h \delta + \delta h \colon \totc{C}[*][][\coplus] \rightarrow \totc{C}[*][][\coplus]$ is
	a bounded $R$-linear map of degree $0$ with $\nnorm[h \delta + \delta h] \leq 1$.
	Given $d \in \ZZ$, \cref{lm:direct-sum-iteration-zero-limit} implies that we also have
	$\left( h \delta + \delta h \right)^n \left( c \right) \to 0$ for any
	$c \in \totc{C}[d][][\coplus] = \coplus_{i \in \ZZ} C_i^{i+d}$. Hence, by
	\cref{lm:inverse-geometric-series}, the map $\idd + h \delta + \delta h$ is invertible with
	an inverse given by the pointwise converging geometric series
	$\sum_{n=0}^{\infty} (-1)^n \left( h \delta + \delta h \right)^n$.

	The formula \eqref{eq:bicomplex-R-linear-contraction-2} for $H$ then follows from
	the formula of the contraction \eqref{eq:trivial-perturbation-lemma-H-tag} of
	\cref{lm:trivial-perturbation-lemma}. When $\left( h^j \right)^2 = 0$,
	i.e., $h^2 = 0$, the formula for $H$ takes the form
	\begin{equation*}
		H = \sum_{n=0}^{\infty} (-1)^n \left( h \delta \right)^n h = \sum_{n=0}^{\infty} (-1)^n h \left( \delta h \right)^n,
	\end{equation*}
	from which it is clear that $H^2 = 0$.
\end{proof}

\begin{rem}
	Let $\mathcal{C} = \left(C^{*}_{*}, b, \delta \right)$ be a bicomplex of $\mathbbm{k}$-modules
	with anticommuting differentials. Endow $\mathbbm{k}$ with the trivial norm and differential and set the norm
	of each nonzero $c \in C_i^j$ to be $2^{-j}$. This turns the bicomplex $\mathcal{C}$ into a Banach bicomplex.
	In this special case, the complete direct sum totalization of $\mathcal{C}$ coincides with the
	\textbf{Laurent totalization}
	\begin{equation*}
		\totc{C}[d][][\coplus] = \mathop{\coplus} \limits_{i \in \ZZ} C_i^{i+d}
		= \left( \bigoplus_{i < 0} C_i^{i+d} \right) \times \left( \prod_{i \geq 0} C_i^{i+d} \right)
		\eqdef \totc{C}[d][][L]
	\end{equation*}
	which is a mix of the direct sum and direct product totalizations. When $\mathcal{C}$ is a right
	half-plane complex, the Laurent totalization coincides with the direct product totalization
	$\totc{\mathcal{C}}[][][\prod]$. When $\mathcal{C}$ is a lower half-plane complex, the Laurent
	totalization coincides with the direct sum totalization $\totc{\mathcal{C}}[][][\oplus]$.

	Now assume that the rows are contractible with contractions $h^j \colon C^j_{*} \rightharpoonup C^j_{*+1}$.
	Then we automatically have $\nnorm[h^j] \leq 1$ and for any $c \in C_i^j$ we have
	$\left( h \delta + \delta h \right)^n \left( c \right) \in C_{i+n}^{j+n}$ and so
	$\nnorm[ \left( h \delta + \delta h \right)^n \left( c \right)] \leq 2^{-j-n} \to 0$. Hence, the conditions of
	\cref{lm:total-bicomplex-contractible-R-linear-contraction-2} are satisfied and we deduce the known
	fact that $\totc{\mathcal{C}}[][][L]$ is contractible.

	In particular, when $\mathcal{C}$ is a right half-plane (resp.\ lower half-plane) complex, we deduce
	that $\totc{\mathcal{C}}[][][\prod]$ (resp.\ $\totc{\mathcal{C}}[][][\oplus]$) is contractible,
	and in particular acyclic (see \cite[Acyclic Assembly Lemma 2.7.3]{Weibel1994}).
\end{rem}

Given a Banach bicomplex $\mathcal{C} = \left( C, b, \delta \right)$ and $k \in \ZZ$, we can truncate
the bicomplex by setting $C_i^j = 0$ whenever $i < k$. We will denote the complete direct sum totalization
of the truncated bicomplex by $\totc{\mathcal{C}}[][\geq k][\coplus]$.

\begin{cor} \label{cor:total-bicomplex-projection-equivalence-2}
	Let $\mathcal{C} = (C_{*}^{*}, b, \delta)$ be a Banach bicomplex concentrated at the right
	half plane, i.e., $C_i = 0$ whenever $i < 0$,
	which satisfies the conditions of \cref{lm:total-bicomplex-contractible-R-linear-contraction-2}.
	Let $p \colon \totc{\mathcal{C}}[][\geq 1][\coplus] \rightharpoonup \mathcal{C}_0$ be the degree one chain map
	given by
	$p \left( \sum_{i \geq 1} c_i^{i+d} \right) \defeq b \left( c_1^{1 + d} \right)$. Then $p$
	is a homotopy equivalence. A homotopy inverse
	$i \colon \mathcal{C}_0 \rightharpoonup \totc{\mathcal{C}}[][\geq 1][\coplus]$
	for $p$ is given explicitly by
	\begin{equation*}
		i = \rest{H}{C_0} = \sum_{n \geq 0} (-1)^n \left( h \delta + \delta h \right)^n \rest{h}{C_0},
	\end{equation*}
	where $H$ is the contraction from \cref{lm:total-bicomplex-contractible-R-linear-contraction-2}.
	We have $pi = \idd$ and $ip = \idd - \partial(H_{+})$ where
	$H_{+} \defeq \rest{H}{\totc{C}[][\geq 1][\coplus]} \colon
		\totc{C}[][\geq 1][\coplus] \rightharpoonup \totc{C}[][\geq 1][\coplus]$. In other words,
	\begin{equation*}
		\left( C_0, \delta \right) \stackbin[i]{p}{\leftrightharpoons} \left(\totc{C}[][\geq 1][\coplus], D \right), H_{+}
	\end{equation*}
	is a deformation retract. We also have $pH_{+} = 0$.
	When the row homotopies $h^j$ satisfy $\left( h^j \right)^2 = 0$, then we also have $H_{+}i = 0$ and $H_{+}^2 = 0$
	so that the deformation retract is in fact a special deformation retract.
\end{cor}
\begin{proof}
	\sloppy
	Let $p' \colon {\totc{C}[][\geq 1][\coplus]}[-1] \rightarrow C_0$ be the composition of
	$p$ with the natural map
	${\totc{C}[][\geq 1][\coplus]}[-1] \rightharpoonup \totc{C}[][\geq 1][\coplus]$.
	Then $p'$ is a degree zero chain map and the mapping cone of $p'$ is given as an $R$-module by
	\begin{equation*}
		\Cone{p'} = \left( {\totc{C}[][\geq 1][\coplus]}[-1] \right)[1] \oplus C_0 =
		\totc{C}[][\geq 1][\coplus] \oplus C_0
	\end{equation*}
	and hence it is isomorphic to $\totc{C}[][\geq 0][\coplus] = \totc{C}[][][\coplus]$.
	Under the isomorphism, the differential on $\Cone{p'}$ coincides with the differential
	of the totalization $\totc{\mathcal{C}}[][][\coplus]$.
	By \cref{lm:total-bicomplex-contractible-R-linear-contraction-2}, $\Cone{p'}$ is contractible
	with an explicit contraction $H$ and so
	\cref{lm:contractible-cone-homotopy-equivalence} implies that $p'$, and hence $p$,
	is a homotopy equivalence. Let us show that $i$ is a homotopy inverse of $p$ satisfying
	the required properties.\footnote{The fact that $i$ is a homotopy inverse
		of $p$ can also be deduced from the proof of \cref{lm:contractible-cone-homotopy-equivalence}.}

	Note that the contraction $H$ maps elements of $C_i$ to $\coplus_{j \geq i + 1} C_j$ and that the
	map $i$ coincides with the contraction $H$ of $\Cone{p'}$, the only difference
	being that it has different domain and codomain.
	Given $c_0 \in C_0$, we have
	\begin{equation*}
		\left( pi \right) \left( c_0 \right) =  \left( bh \right) \left( c_0 \right) = c_0
	\end{equation*}
	so $pi = \idd$. Let us denote the differential on $\totc{C}[][\geq 0][\coplus]$ (resp.\ $\totc{C}[][\geq 1][\coplus]$)
	by $D$ (resp.\ $D_{+}$). Then, given $c_i \in C_i$ with $i > 1$, we have
	\begin{equation*}
		\left( \partial H_{+} \right) \left( c_i \right) =
		\left( D_{+} H_{+} + H_{+} D_{+} \right) \left( c_i \right) = \left( D H + H D \right) \left( c_i \right) = c_i =
		\left( \idd - ip \right) \left( c_i \right)
	\end{equation*}
	since $p \left( c_i \right) = 0$. Finally, given $c_1 \in C_1$, we have
	\begin{equation*}
		\begin{aligned}
			\left( \partial H_{+} \right) \left( c_1 \right) & =
			\left( D_{+} H_{+} + H_{+} D_{+} \right) \left( c_1 \right) =
			DH \left( c_1 \right) + H \left( D - b \right) \left( c_1 \right) =
			\left( \idd - ib \right) \left( c_1 \right)
			\\
			                                                 & =
			\left( \idd - ip \right) \left( c_1 \right).
		\end{aligned}
	\end{equation*}
	It is clear that $pH_{+} = 0$. When $\left( h^j \right)^2 = 0$, we have $H^2 = 0$ and so $H_{+} i = 0$ and $H_{+}^2 = 0$.
\end{proof}

\subsubsection{An Alternative Contraction}
In \cref{lm:total-bicomplex-contractible-R-linear-contraction-2}, we have shown that,
under a suitable convergence condition, the total complex of a Banach bicomplex is contractible
with a contraction given by
\begin{equation*}
	H = \sum_{n=0}^{\infty} (-1)^n \left(  h \delta + \delta h \right)^n h.
\end{equation*}
In the literature one can find an arguably simpler contraction for the total complex, given by
\begin{equation*}
	H' = \sum_{n=0}^{\infty} (-1)^n \left( h \delta \right)^n h
\end{equation*}
(see e.g.\ \cite[Section 3.5]{Crainic2004}). In this subsection, which is not used elsewhere in the work,
we show that the contraction $H'$ also works in our setting but in general, since we work over
a differential graded algebra, the resulting contraction $H'$ is not $R$-linear, unlike $H$.
The advantage of working with the commutator $h \delta + \delta h = \left[ h, \delta \right]$ appearing
in $H$ over $h \delta$ appearing in $H'$, is that it is guaranteed to be $R$-linear even when $\delta$ is a derivation over $d$.

\begin{lm} \label{lm:total-bicomplex-contractible-rows}
	Let $\mathcal{C} = (C_{*}^{*}, b, \delta)$ be a Banach bicomplex and assume that the rows
	$\left( C_{*}^j, b \right)$
	are contractible with contracting $\mathbbm{k}$-linear homotopies
	$h^j \colon C_{*}^j \rightharpoonup C_{*+1}^j$ which satisfy $\nnorm[h^j] \leq 1$ for all $j \in \ZZ$.
	Assume also that for each $c \in C_i^j$ we have
	\begin{equation*}
		\left( h \delta \right)^n (c) \to 0.
	\end{equation*}
	Then the infinite sum
	\begin{equation} \label{eq:bicomplex-k-linear-contraction-2}
		H' \defeq \sum_{n=0}^{\infty} (-1)^n \left( h \delta \right)^n h = \sum_{n=0}^{\infty} (-1)^n h \left( \delta h \right)^n
	\end{equation}
	converges pointwise and defines a $\mathbbm{k}$-linear degree $-1$ map
	$H' \colon \totc{C}[*][][\coplus] \rightharpoonup \totc{C}[* - 1][][\coplus]$
	which satisfies $DH' + H'D = \idd$ and $\nnorm[H'] \leq 1$.
	In particular, $\totc{\mathcal{C}}[][][\coplus]$
	is contractible as a differential graded module over $\mathbbm{k}$.
\end{lm}
\begin{proof}
	First, note that since $\nnorm[h^j] \leq 1$ for all $j \in \ZZ$, the maps $h^j$ extend to a well-defined
	bounded map $h \colon \totc{\mathcal{C}}[*][][\coplus] \rightharpoonup \totc{\mathcal{C}}[* - 1][][\coplus]$
	on the totalization with	$\nnorm[h] \leq 1$. The map
	$h \delta \colon \totc{\mathcal{C}}[*][][\coplus] \rightarrow \totc{\mathcal{C}}[*][][\coplus]$ is also
	bounded with $\nnorm[h \delta] \leq 1$.

	Let $d, i \in \ZZ$ and $c \in C_i^{i+d}$. Since $\left( h \delta \right)^n \left( c \right) \to 0$, we also have
	\begin{equation}	 \label{eq:homotopy-infinite-sum-general-term-tends-to-zero}
		(-1)^n \left( \left( h \delta \right)^n h \right) \left( c \right) \to 0
	\end{equation}
	whenever $c \in C_i^{i+d}$. \Cref{lm:direct-sum-iteration-zero-limit} then implies that
	\cref{eq:homotopy-infinite-sum-general-term-tends-to-zero} continues to hold whenever
	$c \in \totc{\mathcal{C}}[d][][\coplus] = \coplus_{i \in \ZZ} C_i^{i+d}$.
	Hence, the infinite sum $H' = \sum_{n=0}^{\infty} (-1)^n \left( h \delta \right)^n h$
	indeed converges pointwise whenever we evaluate it at any $c \in \totc{\mathcal{C}}[d][][\coplus]$.
	Note that $H'$ is also bounded with $\nnorm[H'] \leq 1$.

	Denote by $\partial$ the differential on
	$\InnHom{\totc{\mathcal{C}}[][][\coplus]}{\totc{\mathcal{C}}[][][\coplus]}$. We have the following identities:
	\begin{align*}
		\partial \left( \delta \right)    & = \left( b + \delta \right) \delta + \delta \left( b + \delta \right) =
		b \delta + \delta b + 2\delta^2 = 0,
		\\
		\partial \left( h \right)         & =
		\left( b + \delta \right)h + h \left( b + \delta \right) = bh + hb + \delta h + h \delta =
		\idd + \delta h + h \delta,
		\\
		\partial \left ( h \delta \right) & = \partial \left( h \right) \delta - h \partial \left( \delta \right) =
		\delta + \delta h \delta = \delta \left( \idd + h \delta \right) = \left( \idd + \delta h \right) \delta.
	\end{align*}
	Induction using the Leibniz identity shows that for $n \geq 1$ we have
	\begin{equation*}
		\partial \left( \left( h \delta \right)^n \right) =
		\delta \left( \left( h \delta \right)^{n-1} + \left( h \delta \right)^n \right).
	\end{equation*}
	Hence, we see that
	\begin{equation*}
		\partial \left( \left( h \delta \right)^n h \right) =
		\begin{cases}
			\left( h \delta \right)^n + \left( h \delta \right)^{n+1} + \left( \delta h \right)^n + \left( \delta
			                                                                                        h \right)^{n+1} & n > 0, \\
			\idd + \delta h + h \delta                                                                            & n = 0.
		\end{cases}
	\end{equation*}
	Thus,
	\begin{equation*}
		\begin{aligned}
			\partial \left( H' \right) & = \partial \left(
			\lim_{N \to \infty} \sum_{n=0}^N (-1)^n \left( h \delta \right)^n h
			\right) =
			\lim_{N \to \infty} \sum_{n=0}^{N} (-1)^n \partial \left( \left( h \delta \right)^n h \right)
			\\
			                           & =
			\lim_{N \to \infty} \left( \idd + \delta h + h \delta +
			\sum_{n = 1}^{N} (-1)^n \left(
			\left( h \delta \right)^n + \left( h \delta \right)^{n+1} + \left( \delta h \right)^n + \left( \delta
			                                                                                        h \right)^{n+1} \right) \right)
			\\
			                           & = \lim_{N \to \infty} \left(
			\idd + (-1)^N \left( \left( \delta h \right)^{N+1} + \left( h \delta \right)^{N+1} \right)
			\right)
			= \idd
		\end{aligned}
	\end{equation*}
	where all limits in the computation above are interpreted pointwise.
\end{proof}

The homotopy $H'$ of \cref{lm:total-bicomplex-contractible-rows} is $\mathbbm{k}$-linear but not
necessarily $R$-linear, even when $h$ is compatible with the $R$-action.
By imposing additional conditions, we can guarantee that $H'$ is in fact $R$-linear.

\begin{lm}
	Let $\mathcal{C} = (C_{*}^{*}, b, \delta)$ be a Banach bicomplex which satisfies the conditions of
	\cref{lm:total-bicomplex-contractible-rows}. Assume that the rows homotopies $h^j$ are compatible with
	the $R$-action on $C$, namely, assume that
	\begin{equation} \label{eq:rows-homotopies-compatability}
		h^{j + k} \left( r \cdot c \right) = (-1)^{k} r \cdot h^{j} \left( c \right)
	\end{equation}
	whenever $c \in C_i^j$ and $r \in R^k$. Assume in addition that $\left( h^j \right)^2 = 0$ for all $j \in \ZZ$.
	Then the homotopy $H'$ of \cref{lm:total-bicomplex-contractible-rows} is $R$-linear.
\end{lm}
\begin{proof}
	Condition \ref{eq:rows-homotopies-compatability} guarantees
	that the induced map $h \colon \totc{\mathcal{C}}[*][][\coplus] \rightharpoonup \totc{\mathcal{C}}[* - 1][][\coplus]$
	is $R$-linear of degree $-1$ on the total complex. The composition $h \delta $ of an $R$-linear map
	$h$ and a derivation $\delta$ is in general not $R$-linear. Instead, we have
	\begin{equation*}
		\left( h \delta \right) \left( r c \right) =
		h \left( dr \cdot c + (-1)^{\degb{r}} r \cdot \delta \left( c \right) \right) =
		(-1)^{\degb{r} + 1} dr \cdot h \left( c \right) + r \cdot \left( h \delta \right) \left( c \right).
	\end{equation*}
	Then
	\begin{equation*}
		\begin{aligned}
			\left( h \delta h \right) \left( r c \right) & =
			(-1)^{\degb{r}} \left( h \delta \right) \left( r \cdot h \left( c \right) \right) =
			-dr \cdot h^2 \left( c \right) + (-1)^{\degb{r}} r \cdot \left( h \delta \right) \left( h \left( c \right) \right)
			\\
			                                             & =
			(-1)^{\degb{r}} r \cdot \left( h \delta h \right) \left( c \right)
		\end{aligned}
	\end{equation*}
	so the map $h \delta h$ is $R$-linear if $h^2 = 0$. More generally, one can verify by induction
	that the maps $\left( h \delta \right)^n h$ are all $R$-linear and hence $H'$ is also $R$-linear.
\end{proof}

\begin{rem}
	A contraction $h \colon C_{*} \rightharpoonup C_{* + 1}$ which satisfies $h^2 = 0$ is sometimes called
	a \textbf{special contraction}. Being a special contraction is not really a restrictive condition. Given
	a contraction $h' \colon C_{*} \rightharpoonup C_{* + 1}$ with $bh' + h'b = \idd$, we can set $h = h'bh'$.
	Then
	\begin{equation*}
		bh + hb = bh'bh' + h'bh'b = \left( \idd - h'b \right)bh' + h' b \left( \idd - b h' \right) =
		bh' + h' b = \idd
	\end{equation*}
	and
	\begin{equation*}
		\begin{aligned}
			h^2 & = h'b h' \left( h' b \right) h' = h' b h' \left( \idd - b h' \right) h' =
			h' b \left( h' \right)^2 - h' b \left( h' b \right) \left( h' \right)^2
			\\
			    & =
			h' b \left( h' \right)^2 - h' b \left( \idd - b h' \right) \left( h' \right)^2 = 0
		\end{aligned}
	\end{equation*}
	and we obtain a special contraction.
	Alternatively, if we have row contractions $h$ which don't satisfy $h^2 = 0$ and we do not want to modify them,
	we can use
	the contraction
	\begin{equation*}
		H = \sum_{n=0}^{\infty} (-1)^n \left( h \delta + \delta h \right)^n h
	\end{equation*}
	of \cref{lm:total-bicomplex-contractible-R-linear-contraction-2}. The contraction $H$ coincides
	with $H'$ of \cref{lm:total-bicomplex-contractible-rows}
	when $h^2 = 0$ but since $h \delta + \delta h$ is the commutator of
	an $R$-linear map and a derivation, the map $h \delta + \delta h$ is $R$-linear and so is $H$.
\end{rem}

\subsection{Banach Bicomplexes with Commuting Differentials} \label{subsec:bicomplex-commuting-differentials}
Assume we choose to work with the inner product parity form $\braidop_1$ given by \cref{eq:parity-inner-product}
instead of working with the total degree parity form $\braidop_2$. This choice dictates that
the differentials $b,\delta$ satisfy
\begin{equation*}
	\left[ b, \delta \right] = b \delta - (-1)^{\braid{(-1,0)}{(0,1)}_1} \delta b =
	b \delta - \delta b = 0
\end{equation*}
and hence the horizontal and vertical differentials \textit{commute} instead of \textit{anticommute}.
The vertical differential $\delta$ satisfies the same identity as \cref{eq:vertical-differential-derivation}
so that we still have that each column $\mathcal{C}_i = \left( C_i^{*}, \delta \right)$ of the bicomplex is a
differential graded Banach $\mathcal{R}$-module. The fact that the horizontal differential $b$ is $R$-linear translates now into the identity
\begin{equation} \label{eq:horizontal-differential-linear-pairing-1}
	b \left( r \cdot c \right) = (-1)^{\braid{(-1,0)}{(0,\degb{r})}_1} r \cdot b \left( c \right) =
	r \cdot b \left( c \right) \in C_{i-1}^{j+l}
\end{equation}
for $r \in R^l$ and $c \in C_i^j$, as opposed to \cref{eq:horizontal-differential-linear-pairing-2},
and so
$b \colon \left( C_{*}, \delta \right) \rightarrow \left( C_{*-1}, \delta \right)$ is a
morphism of differential graded Banach $\mathcal{R}$-modules.
Hence, we see that a Banach bicomplex with respect to $\braidop_1$
is the same thing as a homological complex of differential graded Banach $\mathcal{R}$-modules
\begin{equation*}
	\cdots \xleftarrow{b} \left( C_{-1}, \delta \right) \xleftarrow{b} \left( C_{0}, \delta \right)
	\xleftarrow{b} \left( C_{1}, \delta \right) \xleftarrow{b} \cdots
\end{equation*}
i.e., a chain complex in the category $\DGBMod[\mathcal{R}]$.

When working with the parity form $\braidop_1$, the definition of the total complex needs to be modified (see
also \cite[Page 247]{Yekutieli2020}).
As a graded $R$-module, the \textbf{complete direct sum totalization} of a $\ZZ^2$-graded Banach $R$-module $C$
is given by
\begin{equation} \label{eq:tot-com-def-parity-1}
	\totc{C}[d][][\coplus] = \totc{C}[d][][\coplus, \braidop_1] \defeq
	\mathop{\coplus} \limits_{i \in \ZZ} C_i[i]^d.
\end{equation}
The resulting $R$-module is not the same as in \eqref{eq:tot-com-def-parity-2} since
shifting the columns twists the $R$-module structure on each column by a sign.
A general element $x \in \totc{C}[d][][\coplus]$ will be denoted by
$x = \sum_{i \in \ZZ} s_i c_i \in \totc{\mathcal{C}}[d][][\coplus]$,
where $c_i \in C_i^{i+d}$ with $\nnorm[c_i] \xrightarrow[\abs{i} \to \infty]{} 0$.
Given two $\ZZ^2$-graded Banach $R$-modules $C$ and $D$ and a graded bounded $\mathbbm{k}$-linear map
$f \colon C \rightharpoonup D$ of degree $(a,b)$, we define the induced map
\begin{equation*}
	f^{\totl} = \totc{f}[][][] = \totc{f}[][][\braidop_1] \colon \totc{C}[][][\coplus] \rightharpoonup \totc{D}[][][\coplus]
\end{equation*}
by
\begin{equation} \label{eq:tot-graded-maps-braid-op-1}
	f^{\totl}
	\left( \sum_{i \in \ZZ} s_i c_i \right) \defeq
	\sum_{i \in \ZZ} (-1)^{bi} \s_{i+a} f \left( c_i \right).
\end{equation}
The definition is not the same as in \eqref{eq:tot-graded-maps-braid-op-2} and involves a sign. For
an explanation of the appearing sign, see \cref{rem:tot-sign-morphisms-from-tensor-constraints}.

\begin{rem} \label{rem:identification-tot-1}
	When working with the parity form $\braidop_1$, given an element $c \in C$ of degree $(i,j)$, i.e., $c \in C_i^j$,
	when thinking of $c$ as an element of the total complex $\totc{C}[][][\coplus]$ of degree $j - i$, i.e.,
	an element of $\totc{C}[j-i][][\coplus] = \coplus_{r \in \ZZ} C_r[r]^{j - i}$,
	we denote it by $\s_i c_i \in C_i[i]^{j-i}$. This is important because the $R$-action
	inside the total complex on the column $C_i$ is twisted by a sign coming from the shift $C_i[i]$,
	i.e., we have $r \cdot \s_i c_i = (-1)^{\degb{r} \cdot i} \s_i \left( r \cdot c_i \right)$.
	Compare to \cref{rem:identification-tot-2}.
\end{rem}

Definition \eqref{eq:tot-graded-maps-braid-op-1} guarantees that if $f$ is $R$-linear
then $f^{\totl}$ is also $R$-linear of degree $b - a$,
and that if $C = D$ and $\delta \colon C \rightharpoonup C$ is a module derivation over $d$
then $\delta^{\totl}$ is also a module derivation over $d$. Thus, the differential $b^{\totl}$ is $R$-linear
while $\delta^{\totl}$ is a derivation over $d$.
We also have
\begin{equation*}
	\left[ b^{\totl}, \delta^{\totl} \right] = b^{\totl} \delta^{\totl} + \delta^{\totl} b^{\totl} = 0,
\end{equation*}
i.e., the degree one differentials $b^{\totl}, \delta^{\totl}$ \textit{anticommute} on the total complex.

The \textbf{complete direct sum totalization} of the bicomplex $\mathcal{C} = \left( C_{*}^{*}, b, \delta \right)$
is then defined to be the differential $\ZZ$-graded Banach $\mathcal{R}$-module given by
\begin{equation} 
	\totc{\mathcal{C}}[][][\coplus] = \totc{\mathcal{C}}[][][\coplus, \braidop_1] \defeq
	\left(
	\totc{C}[][][\coplus],
	D = b^{\totl} + \delta^{\totl}
	\right).
\end{equation}
Note that $D$ is a degree one derivation over $d$ and we have $D^2 = 0$ so
$\totc{\mathcal{C}}[][][\coplus]$ is indeed a differential graded Banach $\mathcal{R}$-module.
The differential $D$ acts on an element
$\sum_{i \in \ZZ} s_i c_i \in \totc{\mathcal{C}}[][][\coplus]$ by the formula
\begin{equation*}
	D \left( \sum_{i \in \ZZ} s_i c_i \right) = \sum_{i \in \ZZ}
	s_i \left( (-1)^i \delta \left( c_i \right) + b_{i+1} \left( c_{i+1} \right) \right).
\end{equation*}

The modifications above changes the $R$-module structure and twists the signs of the vertical differentials
to guarantee that $\totc{\mathcal{C}}[][][\coplus]$ is indeed an $R$-module with two anticommuting
degree one differentials $b^{\totl}, \delta^{\totl}$. Note also that that the action of $\delta^{\totl}$
on elements in the $i$-th column inside $\totc{\mathcal{C}}[][][\coplus]$ coincides with the action
of the differential $\rest{\delta}{C_i}[i]$ on the shifted column $C_i[i]$.

Then we have the following versions
of \cref{lm:total-bicomplex-contractible-R-linear-contraction-2} and
\cref{cor:total-bicomplex-projection-equivalence-2}:

\begin{lm} \label{lm:total-bicomplex-contractible-R-linear-contraction-1}
	Let $\mathcal{C} = (C_{*}^{*}, b, \delta)$ be a Banach bicomplex and assume that the rows
	$\left( C_{*}^j, b \right)$ are contractible with contracting $\mathbbm{k}$-linear homotopies
	$h^j \colon C_{*}^j \rightharpoonup C_{*+1}^j$ which satisfy $\nnorm[h^j] \leq 1$ for all $j \in \ZZ$.
	Assume also that the rows homotopies $h^j$ are compatible with
	the $R$-action on $C$, namely, assume that
	\begin{equation} \label{eq:rows-homotopies-commuting-compatability-1}
		h^{j + k} \left( r \cdot c \right) = r \cdot h^{j} \left( c \right)
	\end{equation}
	whenever $c \in C_i^j$ and $r \in R^k$. Hence, the maps $h^j$ induce a degree minus one
	$R$-linear map $h^{\totl} \colon \totc{\mathcal{C}}[*][][\coplus] \rightharpoonup \totc{\mathcal{C}}[*-1][][\coplus]$
	with $\nnorm[h^{\totl}] \leq 1$ acting on
	$s_i \left( c_i^{i+d} \right) \in \totc{\mathcal{C}}[d][][\coplus]$
	by $h^{\totl} \left( s_i \left( c_i^{i+d} \right) \right) =
		s_{i+1} \left( h^{i+d} \left( c_i^{i+d} \right) \right)$.
	Assume also that for each $c \in C_i^j$ we have
	\begin{equation}
		\left( h \delta - \delta h \right)^n (c) \to 0. \label{eq:bicomplex-convergence-cond-1}
	\end{equation}
	Then the infinite sum
	\begin{equation} \label{eq:bicomplex-commuting-R-linear-contraction}
		H \defeq
		\sum_{n=0}^{\infty} (-1)^n \left( \delta^{\totl} h^{\totl} + h^{\totl} \delta^{\totl} \right)^n h^{\totl}
	\end{equation}
	converges pointwise and defines an $R$-linear degree $-1$ map
	$H \colon \totc{C}[*][][\coplus] \rightharpoonup \totc{C}[* - 1][][\coplus]$
	which satisfies $DH + HD = \idd$ and $\nnorm[H] \leq 1$. In particular, $\totc{\mathcal{C}}[][][\coplus]$
	is contractible as a differential graded module over $\mathcal{R}$.

	When the row homotopies $h^j$ satisfy $\left( h^j \right)^2 = 0$ for all $j \in \ZZ$, the contraction $H$ satisfies $H^2 = 0$.
\end{lm}
\begin{proof}
	Note that on the level of total complexes, we have
	\begin{equation} \label{eq:tot-delta-tot-h-plus-first-parity}
		\left( \delta^{\totl} h^{\totl} + h^{\totl} \delta^{\totl} \right)^n \left( \s_i c_i \right) =
		(-1)^{ni + \frac{n(n-1)}{2}} \s_{i+n} \left( \left( h \delta - \delta h \right)^n \left( c_i \right) \right)
	\end{equation}
	and then the proof is the same as the proof of \cref{lm:total-bicomplex-contractible-R-linear-contraction-2}.
\end{proof}

\begin{cor} \label{cor:total-bicomplex-projection-equivalence-1}
	Let $\mathcal{C} = (C_{*}^{*}, b, \delta)$ be a Banach bicomplex concentrated at the right
	half plane (i.e., $C_i = 0$ whenever $i < 0$)
	which satisfies the conditions of \cref{lm:total-bicomplex-contractible-R-linear-contraction-1}.
	Let $p \colon \totc{\mathcal{C}}[][\geq 1][\coplus] \rightharpoonup \mathcal{C}_0$ be the degree one
	chain map given by
	$p \left( \sum_{i \geq 1} s_i \left( c_i^{i+d} \right) \right) = b \left( c_1^{1 + d} \right)$.
	Then $p$ is a homotopy equivalence. A homotopy inverse
	$i \colon \mathcal{C}_0 \rightharpoonup \totc{\mathcal{C}}[][\geq 1][\coplus]$
	for $p$ is given explicitly by
	\begin{equation*}
		i = \rest{H}{C_0} = \sum_{n \geq 0} (-1)^n \left( \delta^{\totl} h^{\totl} + h^{\totl} \delta^{\totl} \right)^n \rest{h^{\totl}}{C_0},
	\end{equation*}
	where $H$ is the contraction from \cref{lm:total-bicomplex-contractible-R-linear-contraction-1},
	and we have
	\begin{equation} \label{eq:total-bicomplex-projection-equivalence-inverse-1}
		i \left( c_0 \right) =
		\sum_{n \geq 0} (-1)^{\frac{n(n-1)}{2}}
		s_{n+1} \left( \left( h\delta - \delta h \right)^n h \left( c_0 \right) \right).
	\end{equation}
	We have $pi = \idd$ and $ip = \idd - \partial(H_{+})$ where
	$H_{+} \defeq \rest{H}{\totc{C}[][\geq 1][\coplus]} \colon
		\totc{C}[][\geq 1][\coplus] \rightharpoonup \totc{C}[][\geq 1][\coplus]$. In other words,
	\begin{equation*}
		\left( C_0, \delta \right) \stackbin[i]{p}{\leftrightharpoons} \left(\totc{C}[][\geq 1][\coplus], D \right), H_{+}
	\end{equation*}
	is a deformation retract. We also have $pH_{+} = 0$.
	When the row homotopies $h^j$ satisfy $\left( h^j \right)^2 = 0$, then we also have $H_{+}i = 0$ and $H_{+}^2 = 0$
	so that the deformation retract is in fact a special deformation retract.
\end{cor}
\begin{proof}
	The proof is the same as the proof of \cref{cor:total-bicomplex-projection-equivalence-2}
	and \cref{eq:total-bicomplex-projection-equivalence-inverse-1} follows from
	\cref{eq:tot-delta-tot-h-plus-first-parity}.
\end{proof}

\begin{rem}
	Note that although condition \eqref{eq:rows-homotopies-commuting-compatability-1} of
	\cref{lm:total-bicomplex-contractible-R-linear-contraction-1} looks different from condition
	\eqref{eq:rows-homotopies-compatability-2} of \cref{lm:total-bicomplex-contractible-R-linear-contraction-2},
	both conditions state that $h$ defines an $R$-linear map of degree $(1,0)$ on $C$, only with
	respect to different parity forms $\braidop_i$, for $i = 1,2$.
	Similarly, the convergence conditions \eqref{eq:bicomplex-convergence-cond-1} of
	\cref{lm:total-bicomplex-contractible-R-linear-contraction-1} and
	\eqref{eq:bicomplex-convergence-cond-2} of \cref{lm:total-bicomplex-contractible-R-linear-contraction-2} look different
	but they both can be rephrased uniformly as
	\begin{equation*}
		\left[ h, \delta \right]^n \left( c \right) \to 0,
	\end{equation*}
	where the commutator $\left[ h, \delta \right] = h \delta - (-1)^{\braidd{h}{\delta}_i} \delta h$ involves
	the different parity forms $\braidop_i$ for $i = 1,2$. In addition, both conditions can be written
	as
	\begin{equation*}
		\left( h^{\totl} \delta^{\totl} + \delta^{\totl} h^{\totl} \right)^n \left( x \right)
		=
		\left( \totc{h}[][][\braidop_i] \totc{\delta}[][][\braidop_i] +
		\totc{\delta}[][][\braidop_i] \totc{h}[][][\braidop_i] \right)^n \left( x \right)
		\to 0
	\end{equation*}
	for $x \in \totc{C}[][][\coplus, \braidop_i]$, where the meaning of the
	$\totl$ depends on the parity form.
\end{rem}

\section{Bar Complex Without Completion is not Acyclic}
\label{appendix:bar-complex-not-necessarily-contractible}

\counterwithin{thm}{section}

In this appendix, we show that the bar complex $( \tensr{A}, \mu )$ of an $\Ainf$-algebra $( A, \mu )$
with $\mu_0 \left( 1 \right) \neq 0$ over a field is not contractible. In what follows, $\tens{A}$ (resp.\ $\tensr{A}$) will denote
the standard tensor module (resp.\ reduced tensor module), as defined using $\otimes$ and $\oplus$, without any completions.

\begin{lm} \label{lm:contraction-only-mu-0}
	Let $R$ be a graded-commutative $\mathbbm{k}$-algebra and let
	$\mathcal{A} = \left( A, \mu \right)$ be an $\Ainf$-algebra over $R$ with $\mu_k = 0$ for $k > 0$.
	This is equivalent to the data of a graded $R$-module $A$ together with a choice of a degree one
	element $\mu_0(1) \in A$.

	Assume that there exists a degree $-1$ $R$-linear map $\varphi \colon A \rightharpoonup R$ such that
	$\varphi \left( \mu_0 \left( 1 \right) \right) = 1$. Then the extended bar complex
	$\left( \tens{A}, \mu \right)$ is contractible. A contracting homotopy is given explicitly
	by the formula
	\begin{align*}
		h \left( a_0 \otimes \dots \otimes a_k \right) & =
		\varphi \left( a_0 \right) \cdot \left( a_1 \otimes \dots \otimes a_k \right), \,\,\, k \geq 0
		\\
		h \left( 1 \right)                             & = 0.
	\end{align*}
\end{lm}
\begin{proof}
	Let us verify that $\mu h + h \mu = \idd$.
	We have
	\begin{equation*}
		\left( \mu h + h \mu \right) \left( 1 \right) = h \left( \mu_0 \left( 1 \right) \right) =
		\varphi \left( \mu_0 \left( 1 \right) \right) = 1.
	\end{equation*}
	Next, given $k \geq 0$ and $a_0, \dots, a_k \in A$ with $x = a_0 \otimes \dots \otimes a_k$, we have
	\begin{equation*}
		\begin{aligned}
			\left( \mu h \right) \left( x \right) & =
			\mu \left( \varphi \left( a_0 \right) \cdot \left( a_1 \otimes \dots \otimes a_k \right) \right)
			\\
			                                      & =
			(-1)^{\degb{\varphi \left( a_0 \right)}} \varphi \left( a_0 \right) \cdot
			\mu \left( a_1 \otimes \dots \otimes a_k \right)
			\\
			                                      & =
			(-1)^{\degb{a_0} - 1} \varphi \left( a_0 \right) \left(
			                                                 \sum_{i=0}^{k} (-1)^{\sum_{j=1}^i \degb{a_j}}
			a_1 \otimes \dots \otimes a_i \otimes \mu_0(1) \otimes a_{i+1} \otimes \dots \otimes a_{k}
			\right)
			\\
			                                      & =
			\sum_{i=0}^{k} (-1)^{\sum_{j=0}^i \degb{a_j} - 1}
			\varphi \left( a_0 \right) \cdot
			\left(a_1 \otimes \dots \otimes a_i \otimes \mu_0(1) \otimes a_{i+1} \otimes
			\dots \otimes a_{k} \right)
		\end{aligned}
	\end{equation*}
	while
	\begin{equation*}
		\begin{aligned}
			\left( h \mu \right) \left( x \right) ={} &
			h \left( \mu_0(1) \otimes a_0 \otimes \dots \otimes a_k \right)
			\\
			                                          & +
			\sum_{i=0}^{k} (-1)^{\sum_{j=0}^i \degb{a_j}} h \left(
			a_0 \otimes \dots \otimes a_i \otimes \mu_0(1) \otimes a_{i+1} \otimes \dots \otimes a_{k} \right)
			\\
			={}                                       &
			\varphi \left( \mu_0(1) \right) \cdot \left( a_0 \otimes \dots \otimes a_k \right)
			\\
			                                          & +
			\sum_{i=0}^{k} (-1)^{\sum_{j=0}^i \degb{a_j}}
			\varphi \left( a_0 \right) \cdot
			\left(a_1 \otimes \dots \otimes a_i \otimes \mu_0(1) \otimes a_{i+1} \otimes
			\dots \otimes a_{k} \right)
			\\
			={}                                       &
			a_0 \otimes \dots \otimes a_k
			\\
			                                          & +
			\sum_{i=0}^{k} (-1)^{\sum_{j=0}^i \degb{a_j}}
			\varphi \left( a_0 \right) \cdot
			\left(a_1 \otimes \dots \otimes a_i \otimes \mu_0(1) \otimes a_{i+1} \otimes
			\dots \otimes a_{k} \right).
		\end{aligned}
	\end{equation*}
\end{proof}

\begin{cor}
	Under the assumptions of \cref{lm:contraction-only-mu-0},
	the cohomology of the bar complex $\left( \tensr{A}, \mu \right)$
	is generated as a graded module by $\mu_0(1)$.
\end{cor}
\begin{proof}
	Let $x \in \tensr{A}$ with $\mu \left( x \right) = 0$. The bar complex $\tensr{A}$ is a subcomplex of the
	extended bar complex $\tens{A}$. Since $\tens{A}$ is acyclic by \cref{lm:contraction-only-mu-0},
	we can write
	\begin{equation*}
		x = \mu \left( r + y \right) = (-1)^{\degb{r}} r \cdot \mu_0 \left( 1 \right) + \mu \left( y \right)
	\end{equation*}
	for some $r \in R$ and $y \in \tensr{A}$. Hence $x$ is cohomologous in $\tensr{A}$ to
	$(-1)^{\degb{r}} r \cdot \mu_0 \left( 1 \right)$.
\end{proof}

\begin{lm}
	Let $\mathcal{A} = \left( A, \mu \right)$ be an $\Ainf$-algebra over a field $\mathbbm{k}$ with
	$\mu_0(1) \neq 0$. Then the bar complex $\left( \tensr{A}, \mu \right)$ is not acyclic.
	The curvature $\mu_0(1) \in \tensr{A}$ is closed but not exact.
\end{lm}
\begin{proof}
	Split the differential $\mu$ as a sum $\mu = \delta + b$ where $\delta, b$ are the unique coderivations
	on $\tens{A}$ which satisfy
	\begin{equation*}
		\begin{aligned}
			\delta_0(1) & = \mu_0(1), & \delta_k & = 0     &  & \left( k \geq 1 \right), \\
			b_0(1)      & = 0,        & b_k      & = \mu_k &  & \left( k \geq 1 \right).
		\end{aligned}
	\end{equation*}
	Denote by $F_k = \oplus_{i=1}^k A^{\otimes i}$ the weight filtration on $\tensr{A}$.
	Note that $b \left( F_k \right) \subseteq F_k$ while $\delta$ increases weight by one and satisfies
	$\delta^2 = 0$. We also have
	\begin{equation*}
		0 = \mu^2 = \left( \delta + b \right)^2 = \delta^2 + \delta b + b \delta + b^2 = \delta b + b \delta + b^2.
	\end{equation*}

	Assume by contradiction that $\mu_0(1)$ is exact and let $x \in F_n$ be a primitive for $\mu_0(1)$.
	Let us show that if $n \geq 2$, we can find a different primitive $x' \in F_{n-1}$.

	Write $x = x_1 + \dots + x_n$ where each $x_i \in A^{\otimes i}$ has weight $i$. By comparing
	weights, we see that the equation $( \delta + b ) ( x ) = \mu_0 ( 1 )$
	implies that $\delta ( x_n ) = 0$. By \cref{lm:contraction-only-mu-0},
	$( \tens{A}, \delta )$ is contractible
	and we can choose $y_{n-1} = h \left( x_n \right) \in A^{\otimes \left( n - 1 \right)}$ with
	$\delta ( y_{n-1} ) = x_n$. Set $x' \defeq x_1 + \dots + x_{n-1} - b y_{n-1} \in F_{n-1}$. Then
	\begin{equation*}
		\begin{aligned}
			( \delta + b ) ( x' ) & =
			( \delta + b ) ( x_1 + \dots + x_{n-1} ) -
			( \delta b + b^2 ) ( y_{n-1} )
			\\
			                      & = ( \delta + b ) ( x_1 + \dots + x_{n-1} ) +
			( b \delta ) ( y_{n-1} )
			\\
			                      & = ( \delta + b ) ( x_1 + \dots + x_{n-1} ) +
			b ( x_n )
			\\
			                      & = ( \delta + b ) ( x_1 + \dots + x_{n-1} ) +
			( b + \delta ) ( x_n )
			\\
			                      & = ( \delta + b ) ( x ) = \mu_0(1).
		\end{aligned}
	\end{equation*}

	Repeating the argument inductively, we see that we can find a primitive $x \in F_1 = A$ with
	$( \delta + b ) ( x ) = \mu_0 ( 1 )$. This implies that
	$\delta ( x ) = 0$ and so we can write $x = \delta ( \lambda )$ for some
	$\lambda \in \mathbbm{k}$. But then
	\begin{equation*}
		\mu_0 ( 1 ) = ( \delta + b ) ( x ) =
		( b \delta ) ( \lambda ) =
		- ( \delta b + b^2 ) ( \lambda ) = 0,
	\end{equation*}
	contradicting the assumption that $\mu_0(1) \neq 0$.
\end{proof}

\counterwithin{thm}{subsection}

\section{Sign Conversions} \label{appendix:sign-conversions}

We briefly discuss various sign conventions for $\Ainf$-algebras and
show that our definitions of the Hochschild complex and Connes' cyclic complex coincide,
after an identification, with the standard definitions in the case of a
differential graded algebra.

Let $\mathbbm{k}$ be a fixed field of characteristic zero. Let
$\mathcal{R} = \left( R, d \right)$ be a differential graded-commutative
$\mathbbm{k}$-algebra, and let $A$ be an $R$-module.
In what follows, when applying our definitions, which are stated
in the non-Archimedean graded Banach framework, we implicitly endow
all the objects with the trivial norm.

\subsection{\texorpdfstring{$\Ainf$-}{A-infinity }Algebra Sign Conventions} \label{sub:a-inf-sign-conventions}
We adopt the point of view that a (cohomological) $\Ainf$-structure on $A$
is a degree one coderivation $\mu$ on the shifted tensor coalgebra $\tens{A[1]}$ which satisfies $\mu^2 = 0$.
The coderivation $\mu$ is uniquely determined by the corestrictions $\mu_k \colon A[1]^{\otimes k} \rightharpoonup A[1]$
and the equation $\corest{\mu} \circ \mu = 0$ (which is equivalent to $\mu^2 = 0$) gives us
the sequence of quadratic relations
\begin{equation} \label{eq:ainf_for_mu_k}
	\sum_{k_1 + k_2 + k_3 = k} \mu_{k_1 + 1 + k_3} \circ \left( \idd^{\otimes k_1}
	\otimes \mu_{k_2} \otimes \idd^{\otimes k_3} \right) = 0
\end{equation}
for $k \geq 0$.

For simplicity of notation, denote by $\sigma = \s^{-1} \colon A[1] \rightharpoonup A$ the desuspension map.
Given bijections $\phi_k \colon A^{\otimes k} \rightharpoonup A[1]^{\otimes k}$ of degree $-k$,
we can transfer the operators $\mu_k$ to operators $m_k \colon A^{\otimes k} \rightharpoonup A$ of
degree $2 - k$ by setting $m_k \defeq \sigma \circ \mu_k \circ \phi_k$. The quadratic relations of
\cref{eq:ainf_for_mu_k} will then translate into quadratic relations for the operators $m_k$ which will involve
extra signs. The specific form of the relations for $m_k$ will depend on the bijections $\phi_k$.
Let us discuss briefly three common choices for the bijections $\phi_k$. For a summary of the different
bijections and their implications, we refer the reader to \cref{tab:a-infinity-relations-and-equations} and
\cref{tab:a-infinity-conversions-mu_k-m_k}.

\begin{figure}[h]
	\centering
	\begin{minipage}[b]{.5\textwidth}
		\centering
		\begin{tikzcd}
			{A[1]^{\otimes k}} & {A[1]} \\
			{A^{\otimes k}} & A
			\arrow["{\mu_k}", harpoon, from=1-1, to=1-2]
			\arrow["{\phi_k = -{ \s }^{\otimes k}}", harpoon, from=2-1, to=1-1]
			\arrow["{m_k}"', harpoon, from=2-1, to=2-2]
			\arrow["{\sigma}", harpoon, from=1-2, to=2-2]
		\end{tikzcd}
		\captionof{figure}{$\phi_k = - {\s}^{\otimes k}$.}
		\label{fig:translate-mu_k-m_k-start-mu_k}
	\end{minipage}%
	\begin{minipage}[b]{.5\textwidth}
		\centering
		\begin{tikzcd}
			{A[1]^{\otimes k}} & {A[1]} \\
			{A^{\otimes k}} & A
			\arrow["{\mu_k}", harpoon, from=1-1, to=1-2]
			\arrow["{\phi_k^{-1} = -{ \sigma }^{\otimes k}}"', harpoon, from=1-1, to=2-1]
			\arrow["{m_k}"', harpoon, from=2-1, to=2-2]
			\arrow["s"', harpoon, from=2-2, to=1-2]
		\end{tikzcd}
		\captionof{figure}{$\phi_k^{-1} = - {\sigma}^{\otimes k}$.}
		\label{fig:translate-mu_k-m_k-start-m_k}
	\end{minipage}
	\caption{Two different ways of converting $\mu_k$ to $m_k$.}
	\label{fig:two-ways-converting-mu_k-m_k}
\end{figure}

\begin{enumerate}
	\item In this work, we choose to work with the bijections
	      \begin{equation} \label{eq:phi_m_k_from_mu_k}
		      \begin{aligned}
			      \phi_k \left( a_1 \otimes \dots \otimes a_k \right) \defeq{} &
			      \left( - { \s }^{\otimes k} \right) \left( a_1 \otimes \dots \otimes a_k \right)
			      \\
			      ={}                                                          &
			      - (-1)^{\sum_{i=1}^k (k-i) \degb{a_i}} \s a_1 \otimes \dots \otimes \s a_k.
		      \end{aligned}
	      \end{equation}
	      See \cref{fig:translate-mu_k-m_k-start-mu_k}. The signs appearing in the formula \eqref{eq:phi_m_k_from_mu_k}
	      follow from the Koszul sign convention.
	      A direct calculation shows that the relations of \cref{eq:ainf_for_mu_k} translate into the relations
	      \begin{equation} \label{eq:ainf_for_m_k}
		      \sum_{k_1 + k_2 + k_3 = k} (-1)^{k_1 + k_2 \cdot k_3} m_{k_1 + 1 + k_3} \circ
		      \left( \idd^{\otimes k_1} \otimes m_{k_2} \otimes \idd^{\otimes k_3} \right) = 0
	      \end{equation}
	      for the operators $m_k$. The relations for $m_k$ are the most common form of relations for
	      an $\Ainf$-algebra (see for example \cite[Page 7]{Keller2001}),
	      and generalize the relations satisfied by the multiplication and differential of a DGA.

	      Our choice of $\phi_k$ implies in particular that $\phi_1 = - \s$ and hence
	      $m_1$ and $\mu_1$ are related by $\mu_1 = -\s \circ m_1 \circ \sigma$. When $\mu_0 = 0$,
	      then $\mu_1, m_1$ are differentials and the relation above implies that
	      $\mu_1 \colon A[1] \rightharpoonup A[1]$ is the shift of the differential
	      $m_1 \colon A \rightharpoonup A$ according to our conventions.

	      The minus sign before the suspension map in $\phi_k = - { \s }^{\otimes k}$
	      is not necessary when working over an ungraded ring. In this case, one can also work with the bijections
	      $\phi_k = \s^{\otimes k}$, and the resulting operators $m_k$ will still satisfy the same relations
	      \eqref{eq:ainf_for_m_k}. However, when working over a differential graded-commutative ground algebra, with our conventions, we must have
	      $\phi_1 = -\s$ to guarantee that $m_1$ is a derivation over $d$ (and not over $-d$).
	      An alternative bijection which guarantees that $m_1$ is a derivation over $d$ without changing
	      the form of the $\Ainf$-relations for the $m_k$ is given by
	      $\phi_k = \left( - \s \right)^{\otimes k} = (-1)^k \s^{\otimes k}$.
	\item Another common bijection is given by
	      \begin{equation*}
		      \begin{aligned}
			      \phi_k \left( a_1 \otimes \dots \otimes a_k \right) \defeq{} &
			      - \left( {\sigma}^{\otimes k} \right)^{-1}
			      \left( a_1 \otimes \dots \otimes a_k \right)
			      \\
			      ={}                                                          &
			      - (-1)^{\frac{k(k-1)}{2} + \sum_{i=1}^k (k-i) \degb{a_i}} \s a_1 \otimes \dots \otimes \s a_k.
		      \end{aligned}
	      \end{equation*}
	      This bijection is more natural if one thinks of the operators $m_k$ as the basic objects and then
	      uses \cref{fig:translate-mu_k-m_k-start-m_k} to define the operators $\mu_k$.
	      In this case, the relations of \cref{eq:ainf_for_mu_k} translate into different relations
	      \begin{equation} \label{eq:ainf_for_m_k2}
		      \sum_{k_1 + k_2 + k_3 = k} (-1)^{(k_1 + 1) \cdot (k_2  + 1)} m_{k_1 + 1 + k_3} \circ
		      \left( \idd^{\otimes k_1} \otimes m_{k_2} \otimes \idd^{\otimes k_3} \right) = 0
	      \end{equation}
	      for the operators $m_k$.

	      Note that the two diagrams in \cref{fig:two-ways-converting-mu_k-m_k} are not equivalent since
	      \begin{equation*}
		      \left( \s^{\otimes k} \right)^{-1} =
		      (-1)^{\frac{k(k-1)}{2}} \left( \s^{-1} \right)^{\otimes k} =
		      (-1)^{\frac{k(k-1)}{2}} {\sigma}^{\otimes k}
		      \neq
		      \left( \sigma \right)^{\otimes k},
	      \end{equation*}
	      and hence the resulting relations for $m_k$ are different. The relations in \cref{eq:ainf_for_m_k2}
	      appear for example in \cite[Page 529]{Polishchuk2004}, and also generalize
	      the relations satisfied by the multiplication and differential of a DGA.
	\item Finally, one can choose to work with the bijection
	      \begin{equation*}
		      \phi_k \left( a_1 \otimes \dots \otimes a_k \right)
		      \defeq \s a_1 \otimes \dots \otimes \s a_k
	      \end{equation*}
	      which involves no signs. In this case, the relations of \cref{eq:ainf_for_mu_k} translate into
	      the relations
	      \begin{align}
		      \sum_{k_1 + k_2 + k_3 = k} (-1)^{\sum_{i=1}^{k_1} \degb{a_i} + k_1}
		       & m_{k_1 + 1 + k_3} \left( a_1, \dots, a_{k_1}, \right. \label{eq:ainf_for_mk_fukaya}
		      \\
		       & \qquad m_{k_2} \left( a_{k_1+1}, \dots, a_{k_1 + k_2} \right),
		      a_{k_1 + k_2 + 1}, \dots, a_k \left. \right) = 0. \nonumber
	      \end{align}
	      for $a_1, \dots, a_k \in A$.
	      This convention appears in \cite{Fukaya2009} and \cite{Seidel2008}. Note that the relations
	      in \cref{eq:ainf_for_mk_fukaya} \textit{do not} coincide with the relations satisfied by the multiplication
	      and differential of a DGA without a sign conversion.

	      In our setting, we work over a graded-commutative base algebra $R$. In this case, one cannot work with the
	      ``bijection'' $\phi_k$ since it is actually not well-defined! For example, we have
	      \begin{equation*}
		      \phi_2 \left( (r a_1) \otimes_R a_2 \right) = \s(ra_1) \otimes_R \s a_2 =
			      (-1)^{\degb{r}} r \left( \s a_1 \otimes_R \s a_2 \right),
	      \end{equation*}
	      while
	      \begin{equation*}
		      \phi_2 \left( (-1)^{\degb{r} \cdot \degb{a_1}} a_1 \otimes_R \left( r a_2 \right) \right) =
		      (-1)^{\degb{r} \cdot \degb{a_1}} \s a_1 \otimes_R \s \left( r a_2 \right) =
		      r \left( \s a_1 \otimes_R \s a_2 \right).
	      \end{equation*}
	      Working over a (differential) graded-commutative base algebra, the signs implied by the Koszul ``convention''
	      stop being a convention and are forced on us to guarantee that various maps are well-defined
	      (see also \cref{foot:sign-necessary-tensor-product}).
\end{enumerate}

\begin{table}[h!]
	\centering
	\begin{tabular}{||c c c c||}
		\hline
		$\Ainf$-structure                                                                               & $\Ainf$-relations             & Unit Equation & MC Equation \\ [0.5ex]
		\hline\hline\hline
		$\mu_k \colon A[1]^{\otimes k} \rightharpoonup A[1]$                                            & \eqref{eq:ainf_for_mu_k}      &
		$\subalign{x &= \mu_2 \left( e, x \right) \\ &= (-1)^{\degb{x} + 1} \mu_2 \left( x, e \right)}$ &
		$\sum_{k=0}^{\infty} \mu_k \left( b^{\otimes k} \right) = ce$                                                                                                 \\
		$m_k \colon A^{\otimes k} \rightharpoonup A$                                                    & \eqref{eq:ainf_for_m_k}       &
		$\subalign{m_2 \left( \be, a \right) &= m_2 \left( a, \be \right) = a}$                         &
		$\sum_{k=0}^{\infty} (-1)^{\frac{k(k+1)}{2}} m_k \left( \mathbf{b}^{\otimes k} \right) = c\be$                                                                \\
		$m_k \colon A^{\otimes k} \rightharpoonup A$                                                    & \eqref{eq:ainf_for_m_k2}      &
		$\subalign{m_2 \left( \be, a \right) &= m_2 \left( a, \be \right) = a}$                         &
		$\sum_{k=0}^{\infty} m_k \left( \mathbf{b}^{\otimes k} \right) + c\be = 0$                                                                                    \\
		$m_k \colon A^{\otimes k} \rightharpoonup A$                                                    & \eqref{eq:ainf_for_mk_fukaya} &
		$\subalign{a &= m_2 \left( \be, a \right) \\ &= (-1)^{\degb{a}} m_2 \left( a, \be \right)}$     &
		$\sum_{k=0}^{\infty} m_k \left( \mathbf{b}^{\otimes k} \right) = c\be$                                                                                        \\
		\hline
	\end{tabular}
	\caption{Different versions of $\Ainf$-relations, unit and MC equations.}
	\label{tab:a-infinity-relations-and-equations}
\end{table}

\begin{table}[h!]
	\centering
	\begin{tabular}{||c c c c||}
		\hline
		Relation between $\mu_k$ and $m_k$                              &
		$\substack{ \Ainf\textrm{-relations} \\ \textrm{for } m_k}$     & Unit                                       & Bounding Cochain \\ [0.5ex]
		\hline\hline\hline
		$m_k = - \sigma \circ \mu_k \circ \s^{\otimes k}$               & \eqref{eq:ainf_for_m_k}                    &
		$e \leftrightarrow \be = -\sigma e$                             & $b \leftrightarrow \mathbf{b} = -\sigma b$                    \\
		$m_k = \sigma \circ \mu_k \circ \left( -\s \right)^{\otimes k}$ & \eqref{eq:ainf_for_m_k}                    &
		$e \leftrightarrow \be = \sigma e$                              & $b \leftrightarrow \mathbf{b} = \sigma b$                     \\
		$\mu_k = - \s \circ m_k \circ {\sigma}^{\otimes k}$             & \eqref{eq:ainf_for_m_k2}                   &
		$e \leftrightarrow \be = \sigma e$                              & $b \leftrightarrow \mathbf{b} = \sigma b$                     \\
		$\mu_k = \s \circ m_k \circ \left( -\sigma \right)^{\otimes k}$ & \eqref{eq:ainf_for_m_k2}                   &
		$e \leftrightarrow \be = -\sigma e$                             & $b \leftrightarrow \mathbf{b} = -\sigma b$                    \\
		\hline
	\end{tabular}
	\caption{Different conversions between $\mu_k$ and $m_k$ operators.}
	\label{tab:a-infinity-conversions-mu_k-m_k}
\end{table}

\subsection{The Hochschild Complex} \label{sec:signs-hoch-unshifted}
Let $\left( A, m_1, m_2 \right)$ be a differential graded $\mathbbm{k}$-algebra. We recall the construction of the
Hochschild complex for $A$ as described in \cite[Section 5.3.2]{Loday1998}.\footnote{In
	\cite[Section 5.3.2]{Loday1998} one works with homological DGAs and the Hochschild differential has
	degree $-1$. The relation between our description and the one given in \cite[Section 5.3.2]{Loday1998} is obtained
	by converting a homological DGA $A_{*}$ into a cohomological DGA $A^{*}$ by setting $A^{*} = A_{-*}$.}
Define a degree one map $b^1 \colon \tensr{A} \rightharpoonup \tensr{A}$ by
\begin{equation*}
	b^1 \left( a_1 \otimes \dots \otimes a_k \right) = \sum_{i=1}^k
		(-1)^{\degb{a_1} + \dots + \degb{a_{i-1}}}
	a_1 \otimes \dots \otimes a_{i-1} \otimes m_1 \left( a_i \right) \otimes a_{i+1} \otimes \dots \otimes a_k,
\end{equation*}
and a degree zero map $b^2 \colon \tensr{A} \rightarrow \tensr{A}$ by
\begin{equation*}
	\begin{aligned}
		b^2 \left( a_0 \otimes \dots \otimes a_k \right) ={} &
		\sum_{i=0}^{k-1} (-1)^i a_0 \otimes \dots \otimes a_{i-1} \otimes m_2 \left( a_i, a_{i+1} \right) \otimes
		a_{i+2} \otimes \dots \otimes a_k                                                                                                     \\
		                                                     & +
		                                                       (-1)^{\degb{a_k} \cdot \left( \degb{a_0} + \dots + \degb{a_{k-1}} \right) + k}
		m_2 \left( a_k, a_0 \right) \otimes a_1 \otimes \dots \otimes a_{k-1},
		\\
		b^2 \left( a_0 \right) = 0.
	\end{aligned}
\end{equation*}

Both $b^1$ and $b^2$ are differentials and the differential $b^1$ coincides with the
induced differential on $\left( A, m_1 \right)^{\otimes k}$. As one can readily verify,
we have $b^1 \circ b^2 = b^2 \circ b^1$, and so we can form a right half plane double complex
$\mathcal{E}$ with commuting differentials whose columns are $\left( A, m_1 \right)^{\otimes k}$ as in \cref{fig:hochschild-double-complex-dga}.

\begin{figure}[h]
	\begin{tikzcd}
		&& \vdots & \vdots & \vdots \\
		\cdots & 0 & {A^1} & {\left( A^{\otimes 2} \right)^1} & {\left( A^{\otimes 3} \right)^1} & \cdots \\
		\cdots & 0 & {A^0} & {\left( A^{\otimes 2} \right)^0} & {\left( A^{\otimes 3} \right)^0} & \cdots \\
		\cdots & 0 & {A^{-1}} & {\left( A^{\otimes 2} \right)^{-1}} & {\left( A^{\otimes 3} \right)^{-1}} & \cdots \\
		&& \vdots & \vdots & \vdots
		\arrow[from=2-3, to=2-2]
		\arrow[from=3-3, to=3-2]
		\arrow[from=4-3, to=4-2]
		\arrow[from=2-2, to=2-1]
		\arrow[from=4-2, to=4-1]
		\arrow[from=3-2, to=3-1]
		\arrow["{b^2}"', from=2-6, to=2-5]
		\arrow["{b^2}"', from=3-6, to=3-5]
		\arrow["{b^2}"', from=4-6, to=4-5]
		\arrow["{b^2}"', from=2-5, to=2-4]
		\arrow["{b^2}"', from=2-4, to=2-3]
		\arrow["{b^2}"', from=3-4, to=3-3]
		\arrow["{b^2}"', from=3-5, to=3-4]
		\arrow["{b^2}"', from=4-4, to=4-3]
		\arrow["{b^2}"', from=4-5, to=4-4]
		\arrow["{b^1}", from=2-3, to=1-3]
		\arrow["{b^1}", from=2-4, to=1-4]
		\arrow["{b^1}", from=2-5, to=1-5]
		\arrow["{b^1}", from=5-3, to=4-3]
		\arrow["{b^1}", from=5-4, to=4-4]
		\arrow["{b^1}", from=5-5, to=4-5]
		\arrow["{b^1}", from=4-3, to=3-3]
		\arrow["{b^1}", from=4-4, to=3-4]
		\arrow["{b^1}", from=4-5, to=3-5]
		\arrow["{b^1}", from=3-5, to=2-5]
		\arrow["{b^1}", from=3-4, to=2-4]
		\arrow["{b^1}", from=3-3, to=2-3]
	\end{tikzcd}
	\caption{The Hochschild Double Complex Associated to a DGA.}
	\label{fig:hochschild-double-complex-dga}
\end{figure}
The \textbf{Hochschild complex} is then defined to be the total complex of the bicomplex $\mathcal{E}$.
More explicitly, the Hochschild complex is given by
$\totc{\mathcal{E}}[][][\oplus] = \oplus_{k \geq 0} A^{\otimes (k+1)}[k]$,
endowed with the differential $D = b^1_{\totl} + b^2_{\totl}$ given by
\begin{equation} \label{eq:hochschild-bicomplex-total-differential}
	D \left( \s_k \left( l \right) \right) = (-1)^k \s_k \left( b^1 \left( l \right) \right) +
	\s_{k-1} \left( b^2 \left( l \right) \right)
\end{equation}
(see \cref{subsec:bicomplex-commuting-differentials}).

Let us relate the definition given above to the definition used in this work.
Given an $\Ainf$-algebra $\left( A, m \right)$, which corresponds to the shifted $\Ainf$-algebra
$\left( A[1], \mu \right)$ using the bijections given in \cref{eq:phi_m_k_from_mu_k}, consider the complex
$\left( \tensr{A[1]}, \cycl{\mu} \right)[-1]$. Denote by $C \defeq \oplus_{k \geq 0} A^{\otimes (k+1)}[k]$ and let
$\varphi \colon C \rightarrow {\tensr{A[1]}}[-1]$ be the natural degree zero isomorphism given by
\begin{equation*}
	\begin{aligned}
		\varphi \left( \s_{k} \left( a_0 \otimes \dots \otimes a_k \right) \right) \defeq{} &
		\s_{-1} \left( {\s}^{\otimes \left( k + 1 \right)}
		\left( a_0 \otimes \dots \otimes a_k \right) \right)
		\\
		={}                                                                                 &
		(-1)^{\sum_{i=0}^k (k - i) \cdot \degb{a_i}}
		\s_{-1} \left( \s a_0 \otimes \dots \otimes \s a_k \right).
	\end{aligned}
\end{equation*}
We can transfer the differential on ${\tensr{A[1]}}[-1]$ to $C$ by setting
\begin{equation*}
	b \defeq \varphi^{-1} \circ \left( \cycl{\mu} \right)_{{\tensr{A[1]}}[-1]} \circ \varphi,
\end{equation*}
so that
$\varphi \colon \left( C, b \right) \rightarrow
	\left( {\tensr{A[1]}}[-1], \left( \cycl{\mu} \right)_{{\tensr{A[1]}}[-1]} \right)$ becomes an isomorphism
of differential graded modules.
The expression of $b$ on $C$ in terms of the operators $m_k$ is given by the following lemma:

\begin{lm}
	Let $a_0 \in A$ and let $l = a_1 \otimes \dots \otimes a_k \in A^{\otimes k}$. Then we have
	the formula
	\begin{equation} \label{eq:formula-for-b}
		\begin{aligned}
			b \left( \s_k \left( a_0 \otimes l \right) \right) ={} &
			(-1)^{\varepsilon_1}
			\s_{k_1 + 1 + k_3} \left(
			a_0 \otimes l_{(1)} \otimes m_{k_2} \left( l_{(2)} \right) \otimes l_{(3)}
			\right)
			\\
			                                                       & +
			                                                         (-1)^{\varepsilon_2} \s_{k_2} \left(
			m_{k_3 + 1 + k_1} \left( l_{(3)} \otimes a_0 \otimes l_{(1)} \right) \otimes l_{(2)}
			\right),
		\end{aligned}
	\end{equation}
	where $k_i = \weight{l_{(i)}}$ is the weight of $l_{(i)}$ and the signs $\varepsilon_1$ and $\varepsilon_2$
	are given by
	\begin{align}
		\varepsilon_1 & = (2 - k_2) \cdot \left( \degb{a_0} + \degb{l_{(1)}} \right) + (1 + k_1) + k_2 \cdot k_3,
		\\
		\varepsilon_2 & = \degb{l_{(3)}} \cdot \left( \degb{a_0} + \degb{l_{(1)}} + \degb{l_{(2)}} \right) +
		k_3 \cdot \left( 1 + k_1 + k_2 \right) +
		\left( k_3 + 1 + k_1 \right) \cdot k_2.
	\end{align}
\end{lm}
\begin{proof}
	It is clear from the formula \eqref{def:cyclization-coder} for $\cycl{\mu}$ that the expression
	for $b$ has the form given by \cref{eq:formula-for-b}. It only remains to compute the signs.
	Given an element $l = a_1 \otimes \dots \otimes a_k \in A^{\otimes k}$, denote by
	$\varepsilon \left( l \right) \defeq \sum_{i=1}^k (k - i) \degb{a_i}$ the sign factor appearing in all
	our identifications. Let us also denote
	$\s l \defeq \s a_1 \otimes \dots \otimes \s a_k$, so that using our notation we have
	$\varphi \left( \s_k l \right) = (-1)^{\varepsilon \left( l \right)} \s_{-1} \left( \s l \right)$
	for $l \in A^{\otimes (k+1)}$. Since we work with the bijections given by \cref{eq:phi_m_k_from_mu_k},
	we have the relation
	$\mu_{k} \left( \s l \right) = -(-1)^{\varepsilon \left( l \right)} \s m_k \left( l \right)$ between $\mu_k$
	and $m_k$. Unwinding all the definitions, we see that the sign accompanying the term
	$\s_{k_1 + 1 + k_3} \left( a_0 \otimes l_{(1)} \otimes m_{k_2} \left( l_{(2)} \right) \otimes l_{(3)} \right)$
	of $b$ is given by
	\begin{equation*}
		\varepsilon_1' = \varepsilon \left( a_0 \otimes l \right) + 1 + \degb{\s a_0} + \degb{\s l_{(1)}} +
		\varepsilon \left( l_{(2)} \right) + 1 +
		\varepsilon \left( a_0 \otimes l_{(1)} \otimes m_{k_2} \left( l_{(2)} \right) \otimes l_{(3)} \right).
	\end{equation*}
	Given a specific splitting
	$l = l_{\left< 1 \right>} \otimes \dots \otimes l_{\left< r \right>}$
	of $l$ into $r$ consecutive, possibly empty, lists, we have the identity
	\begin{equation}
		\varepsilon \left( l \right) = \sum_{i=1}^r \varepsilon \left( l_{\left< i \right>} \right) +
		\sum_{j=2}^r k_j \cdot \sum_{i=1}^{j-1} \degb{l_{\left< i \right>}},
		\label{eq:identity-varepsilon-sign}
	\end{equation}
	where $k_j$ is the weight (i.e., length) of $l_{\left< j \right>}$.
	Using the identity above, we have
	\begin{equation*}
		\begin{aligned}
			\varepsilon_1' ={} & \varepsilon \left( l \right) + k \cdot \degb{a_0} + 1 + \left( \degb{a_0} - 1 \right)
			+ \left( \degb{l_{(1)}} - k_1 \right) + \varepsilon \left( l_{(2)} \right) + 1
			\\
			                   & + \varepsilon \left( l_{(1)} \right) + \varepsilon \left( l_{(3)} \right) + k_1 \cdot \degb{a_0}
			+ \left( \degb{a_0} + \degb{l_{(1)}} \right) + k_3 \cdot \left( \degb{a_0} + \degb{l_{(1)}}
			+ \degb{m_{k_2} \left( l_{(2)} \right)} \right)
			\\
			={}                & \varepsilon \left( l \right) + \varepsilon \left( l_{(1)} \right) + \varepsilon \left( l_{(2)} \right)
			+ \varepsilon \left( l_{(3)} \right)
			+ \left( k + k_1 + k_3 + 2 \right) \cdot \degb{a_0}
			+ \left( 2 + k_3 \right) \cdot \degb{l_{(1)}}
			\\
			                   & + k_3 \cdot \degb{l_{(2)}} + k_3 \cdot (2 - k_2) - k_1 + 1
			\\
			\equiv{}           &
			k_2 \cdot \degb{l_{(1)}} + k_3 \cdot \left( \degb{l_{(1)}} + \degb{l_{(2)}} \right) +
			(2 - k_2) \cdot \degb{a_0} + k_3 \cdot \degb{l_{(1)}}
			\\
			                   & + k_3 \cdot \degb{l_{(2)}} + k_2 \cdot k_3 + k_1 + 1
			\\
			\equiv{}           &
			(2 - k_2) \cdot \left( \degb{a_0} + \degb{l_{(1)}} \right) + \left( 1 + k_1 \right) + k_2 \cdot k_3
		\end{aligned}
	\end{equation*}
	which is precisely $\varepsilon_1$.
	The sign accompanying the term
	$\s_{k_2} \left( m_{k_3 + 1 + k_1} \left( l_{(3)} \otimes a_0 \otimes l_{(1)} \right) \otimes l_{(2)} \right)$
	of $b$ is obtained similarly.
\end{proof}

Now, assume that $A$ is a differential graded $\mathbbm{k}$-algebra, so that $m_k = 0$ for $k \neq 1, 2$. In this case, we see
that the differential $b$ given by \cref{eq:formula-for-b} coincides precisely with
the total differential $D$ given by \cref{eq:hochschild-bicomplex-total-differential}.

\begin{rem}
	It is tempting to ignore the suspension maps, change the grading of elements on $\tensr{A}$
	and write the formula \eqref{eq:formula-for-b} for $b$ directly on $\tensr{A}$.
	Let us define a new graded $R$-module $\tensr{A}_{\star}$ by setting
	\begin{equation*}
		\tensr{A}_{\star}^d = \bigoplus_{k \geq 0} \left( A^{\otimes_R (k+1)} \right)^{d + k}.
	\end{equation*}
	This way, an elementary tensor $a_0 \otimes \dots \otimes a_k \in \tensr{A}_{\star}$ has
	(star) degree
	\begin{equation*}
		\degb{a_0 \otimes \dots \otimes a_k}_{\star} = \degb{a_0} + \dots + \degb{a_k} - k.
	\end{equation*}
	Then the map $\tensr{A}_{\star} \rightarrow C$ given by
	\begin{equation*}
		a_0 \otimes \dots \otimes a_k \mapsto \s_k \left( a_0 \otimes \dots \otimes a_k \right)
	\end{equation*}
	appears to be a degree zero isomorphism between $\tensr{A}_{\star}$ and $C$, and we can transport $b$ along
	this isomorphism to a map on $\tensr{A}_{\star}$ given by the same formula \eqref{eq:formula-for-b} as above,
	only without the suspensions.

	This works well when the base ring $R$ is not graded but in our context, the map described above
	is actually not well-defined. The reason is that
	while the $R$-action on $\tensr{A}_{\star}$ is given by the standard $R$-action on each factor
	$A^{\otimes (k + 1)}$, in $C$, the $R$-action on each factor $A^{\otimes (k+1)}[k]$ is twisted by the shift.

	It turns out that the $R$-modules $C$ and $\tensr{A}_{\star}$ are isomorphic, but the map above needs to be
	modified with a sign factor, and this changes the formula for the Hochschild differential on $\tensr{A}_{\star}$.
	We will not pursue this issue further.
\end{rem}

\subsection{The Rotation Operator and Connes' Complex} \label{sec:signs-connes-unshifted}
Let $(A,m_1,m_2)$ be a differential graded $\mathbbm{k}$-algebra. We recall the construction
of Connes' cyclic complex for $A$ as described in \cite[Section 5.3.2]{Loday1998}.
Denote by $\tilde{t} = \t_{\tensr{A}} \colon \tensr{A} \rightarrow \tensr{A}$
the rotation operator on $\tensr{A}$, using the unshifted signs. We have
\begin{equation} \label{eq:formula-for-tilde-t}
	\tilde{t} \left( a_0 \otimes \dots \otimes a_k \right) \defeq
	(-1)^{\degb{a_k} \cdot \left( \sum_{i=0}^{k-1} \degb{a_i} \right) + k}
	a_k \otimes a_0 \otimes \dots \otimes a_{k-1}.
\end{equation}
Define an operator $t \colon C \rightarrow C$ by the same formula as above, taking into
account the shifts. Namely,
\begin{equation*}
	t \left( \s_k \left( a_0 \otimes \dots \otimes a_k \right) \right) \defeq
	(-1)^{\degb{a_k} \cdot \left( \sum_{i=0}^{k-1} \degb{a_i} \right) + k}
	\s_k \left( a_k \otimes a_0 \otimes \dots \otimes a_{k-1} \right).
\end{equation*}
Then $b \left( \Im \left( \idd - t \right) \right) \subseteq \Im \left( \idd - t \right)$
and so, $b$ descends to the quotient $C / \Im \left( \idd - t \right)$.
\textbf{Connes' complex} is then defined to be the complex
$\left( C / \Im \left( \idd - t \right), b \right)$.

The relation between the rotation operator $\t = \t_{\tensr{A[1]}}$ on $\tensr{A[1]}$, which applies a cyclic rotation with
signs dictated by the Koszul convention, and the operator $\tilde{t}$ on $\tensr{A}$ given
by \cref{eq:formula-for-tilde-t}, is given by the following lemma:

\begin{figure}[H]
	\begin{tikzcd}
		{A[1]^{\otimes (k+1)}} & {A[1]^{\otimes (k+1)}} \\
		{A^{\otimes (k+1)}} & {A^{\otimes (k+1)}}
		\arrow["{\s^{\otimes (k+1)}}", harpoon, from=2-1, to=1-1]
		\arrow["{\left( \s^{\otimes (k+1)} \right)^{-1}}", harpoon, from=1-2, to=2-2]
		\arrow["{\tilde{t}}"', from=2-1, to=2-2]
		\arrow["{\t}", from=1-1, to=1-2]
	\end{tikzcd}
	\caption{Translating $\t \colon \tensr{A[1]} \rightarrow \tensr{A[1]}$ to
		$\tilde{t} \colon \tensr{A} \rightarrow \tensr{A}$.}
	\label{fig:translating-tau-to-tilde-t}
\end{figure}

\begin{lm}
	The operators $\t$ and $\tilde{t}$ are related by the diagram in \cref{fig:translating-tau-to-tilde-t}.
\end{lm}
\begin{proof}
	Let $a_0, \dots, a_k \in A$. Clearly, we have
	\begin{equation*}
		\left( \left( \s^{\otimes (k+1)} \right)^{-1} \circ \t \circ \s^{\otimes (k+1)} \right)
		\left( a_0 \otimes \dots \otimes a_k \right) =
		(-1)^{\varepsilon} a_k \otimes a_0 \otimes \dots \otimes a_{k-1},
	\end{equation*}
	so we need to verify the sign $\varepsilon$ coincides with the sign appearing in \cref{eq:formula-for-tilde-t}.
	Indeed, we have
	\begin{equation*}
		\begin{aligned}
			\varepsilon & = \sum_{i=0}^k (k-i) \degb{a_i} +
			\degb{ \s a_k} \left( \sum_{i=0}^{k-1} \degb{ \s a_i} \right) + k \degb{a_k} +
			\sum_{i=0}^{k-1} \left( k - i - 1 \right) \degb{a_i}
			\\
			            & \equiv
			\left( \degb{a_k} - 1 \right) \cdot \left( \sum_{i=0}^{k-1} \degb{a_i} - k \right)
			+ k \degb{a_k} - \sum_{i=0}^{k-1} \degb{a_i}
			\\
			            & \equiv
			\degb{a_k} \cdot \left( \sum_{i=0}^{k-1} \degb{a_i} \right) + k \mod 2.
		\end{aligned}
	\end{equation*}
\end{proof}

\begin{rem}
	Note that the sign $(-1)^{\varepsilon}$ appearing in \cref{eq:formula-for-tilde-t} has two contributions.
	Set $\varepsilon_1 =\degb{a_k} \cdot \left( \sum_{i=0}^{k-1} \degb{a_i} \right)$
	and $\varepsilon_2 = k$, so that $\varepsilon = \varepsilon_1 + \varepsilon_2$.
	Then $(-1)^{\varepsilon_1}$ is a sign coming from applying the Koszul convention to
	$a_0 \otimes \dots \otimes a_k \mapsto a_k \otimes a_0 \otimes \dots \otimes a_{k-1}$,
	while $(-1)^{\varepsilon_2}$ is the sign of the cyclic permutation $(0 1 \cdots k)$.
	The sign $(-1)^{\varepsilon}$ is sometimes called the antisymmetric Koszul sign associated to the permutation
	and the elements.
	Another way to think about the sign $(-1)^{\varepsilon}$ is to note that $\tensr{A}$ is actually bigraded
	with respect to cohomological degree and weight, and $\tilde{t}$ preserves both gradings. If
	we apply the Koszul sign rule with the symmetry induced by the standard inner product pairing
	(given by \cref{eq:parity-inner-product}), we obtain precisely the sign above.
\end{rem}

In our work, we have defined Connes' complex to be
$\left( \tensr{A[1]} / \Im \left( \idd - \t \right), \cycl{\mu} \right)$. The relation between
$\left( \tensr{A[1]} / \Im \left( \idd - \t \right), \cycl{\mu} \right)$ and the complex
$\left( C / \Im \left( \idd - t \right), b \right)$ is given by the following lemma:

\begin{lm}
	Let $\left( A, m \right)$ be an $\Ainf$-algebra corresponding to the
	shifted $\Ainf$-algebra $\left( A[1], \mu \right)$. Then the differential $b \colon C \rightarrow C$
	given by \cref{eq:formula-for-b} descends to the quotient
	$C / \Im \left( \idd - t \right)$. The map
	$\left( C / \Im \left( \idd - t \right), b \right) \rightarrow
		\left( \tensr{A[1]} / \Im \left( \idd - \t \right), \cycl{\mu} \right)[-1]$ induced by $\varphi$,
	given explicitly by
	\begin{equation*}
		\eqcl{\s_k \left( a_0 \otimes \dots \otimes a_k \right)} \mapsto
		(-1)^{\sum_{i=0}^k (k-i) \cdot \degb{a_i}} \s_{-1} \eqcl{ \s a_0 \otimes \dots \otimes \s a_k },
	\end{equation*}
	is a well-defined isomorphism of differential graded modules. \qed
\end{lm}

In particular, when $A$ is a differential graded $\mathbbm{k}$-algebra,
we see that our Connes complex $\left( \tensr{A[1]} / \Im \left( \idd - \t \right), \cycl{\mu} \right)$
coincides with the standard Connes complex $\left( C / \Im \left( \idd - t \right), b \right)$
up to a shift and an identification.

\section{Equivalence of Codifferential Forms Constructed Using the Inner Product and Total Degree Parity Forms}
\label{appendix:parity-forms-equiv}

In \cref{sec:noncomm-diff-calc}, we have described the construction of noncommutative
codifferential forms and their calculus, working with a fixed background parity form $\braidop$
which was either $\braidop_1$ or $\braidop_2$.
Recall their definitions from
\cref{eq:parity-inner-product,eq:parity-total-degree}:
\begin{align*}
	\braid{(a_1,a_2)}{(b_1,b_2)}_1 & \defeq a_1 b_1 + a_2 b_2 \mod 2,
	\\
	\braid{(a_1,a_2)}{(b_1,b_2)}_2 & \defeq (a_2 - a_1) \cdot (b_2 - b_1) \mod 2.
\end{align*}
In what follows, we explicitly describe how
the parity form $\braidop$ plays a role in the construction and specify the precise sense
in which the constructions using $\braidop_i$ for $i \in \Set{1,2}$ are equivalent. For simplicity,
we state everything for graded modules over graded algebras, but everything we do here works verbatim
in the Banach context, by adding the Banach adjective and using complete direct sums and tensor products,
as we implicitly do in \cref{sec:noncomm-diff-calc}.

\subsection{Dependence of Categorical Constructions on \texorpdfstring{$\braidop$}{the Parity Form}} \label{sec:cat-dependence-braidop}
We start by setting up notation and recalling the role the parity form $\braidop$ plays in various
notions. Our categorical setup was described at length in \cref{sec:prelim}, and we refer there for further details.

Consider the category $\GMod[\mathbbm{k}][\ZZ^2]$ of $\ZZ^2$-graded $\mathbbm{k}$-modules
$M = \left( M_i^j \right)$. We use the notation $\degb{m} = (i,j)$ for the degree of an
element $m \in M_i^j$.
The tensor product $M \otimes_{\mathbbm{k}} N$ of two $\ZZ^2$-graded $\mathbbm{k}$-modules $M$ and $N$
is given by
\begin{equation*}
	\left( M \otimes_{\mathbbm{k}} N \right)_i^j =
	\bigoplus_{\substack{i_1 + i_2 = i \\ j_1 + j_2 = j}} M_{i_1}^{j_1} \otimes_{\mathbbm{k}} N_{i_2}^{j_2}
\end{equation*}
and endows $\GMod[\mathbbm{k}][\ZZ^2]$ with a monoidal structure
whose unit is $\mathbbm{k}$, considered as a $\ZZ^2$-graded module concentrated in degree $(0,0)$.
We denote the resulting monoidal category by
\begin{equation*}
	\mathcal{C} \defeq
	\left( \GMod[\mathbbm{k}][\ZZ^2], \otimes_{\mathbbm{k}}, \mathbbm{k} \right)
\end{equation*}
where, as usual, we suppress from the notation the associators and unitors which are part of the
definition of the monoidal structure.

We can endow the same category $\mathcal{C}$ with different symmetries
$\sigma_i \colon M \otimes_{\mathbbm{k}} N \rightarrow N \otimes_\mathbbm{k} M$ given by
\begin{equation*}
	m \otimes_{\mathbbm{k}} n \xmapsto{\sigma_i} (-1)^{\braid{\degb{m}}{\degb{n}}_i} n \otimes_{\mathbbm{k}} m
\end{equation*}
for $i \in \Set{1,2}$. This gives us two distinct closed symmetric monoidal categories
\begin{equation*}
	\mathcal{C}_i \defeq
	\left( \GMod[\mathbbm{k}][\ZZ^2], \otimes_{\mathbbm{k}}, \mathbbm{k}, \sigma_i \right)
\end{equation*}
with the same objects and morphisms, the same tensor product $\otimes_{\mathbbm{k}}$, the same unit
$\mathbbm{k}$ and the same internal hom object $\InnHom{M}{N}[][\mathbbm{k}]$ consisting of graded maps
but with different symmetries. The basic role the different symmetries play in $\mathcal{C}_i$ is in defining
the tensor product of two graded maps $f \colon M \rightharpoonup M'$ and
$g \colon N \rightharpoonup N'$ which is given by
\begin{equation*}
	\left( f \otimes_{i} g \right) \left( m \otimes_{\mathbbm{k}} n \right) =
	(-1)^{\braidd{g}{m}_i} f \left( m \right) \otimes_{\mathbbm{k}} g \left( n \right)
\end{equation*}
(see \cref{eq:tensor-product-graded-maps}). We use the notation $f \otimes_{i} g$ instead of
$f \otimes_{\mathbbm{k}} g$ to emphasize the dependence of the tensor product of graded maps on
the parity form $\braidop_i$. Any notion that is related to the tensor product of graded
maps will take a different form based on $\braidop_i$.

The fact that the categories $\mathcal{C}_i$ are monoidal allows us to talk about
algebra objects in $\mathcal{C}_i$ (see \cref{subsec:alg-in-monoidal-cat}).
Since the monoidal structures of $\mathcal{C}_1$ and $\mathcal{C}_2$ are strictly the same and coincide
with $\mathcal{C}$, there is no difference between algebra objects in $\mathcal{C}_1$ or in $\mathcal{C}_2$,
and we have $\Alg[\mathcal{C}] = \Alg[\mathcal{C}_1] = \Alg[\mathcal{C}_2]$.
Furthermore, $\mathcal{C}_i$ are \textit{symmetric} monoidal, so one can talk about \textit{commutative} algebra
objects of $\mathcal{C}_i$. An algebra object
$A$ of $\mathcal{C}_i$, i.e., a $\ZZ^2$-graded $\mathbbm{k}$-algebra, is commutative if
\begin{equation*}
	a \cdot b =
	m \left( a \otimes b \right) =
	(-1)^{\braidd{a}{b}_i} m \left( b \otimes a \right) =
	(-1)^{\braidd{a}{b}_i} b \cdot a
\end{equation*}
for all $a,b \in A$. The notion of commutativity depends on the parity $\braidop_i$,
so we have $\CAlg[\mathcal{C}_1] \neq \CAlg[\mathcal{C}_2]$.

Given $i \in \Set{1,2}$, let $A_i$ be a commutative algebra object of $\mathcal{C}_i$,
and consider the category $\GMod[A_i][\ZZ^2]$ of $\ZZ^2$-graded left $A_i$-modules,
i.e., module objects over commutative algebra objects of $\mathcal{C}_i$.
The category $\GMod[A_i][\ZZ^2]$
has a natural structure of a closed symmetric monoidal category which we denote by
\begin{equation*}
	\mathcal{C}_i \left( A_i \right) \defeq \left( \GMod[A_i][\ZZ^2], \otimes_i, A_i, \widehat{\sigma}_i \right).
\end{equation*}
We have the following differences between $\mathcal{C}_i \left( A_i \right)$ and $\mathcal{C}_i$:
\begin{enumerate}
	\item The tensor product $\otimes_i = \otimes_{A_i, \braidop_i}$ of two graded $A_i$-modules $M$ and $N$
	      depends both on $A_i$ and $\braidop_i$ as it is constructed via the quotient
	      of $M \otimes_{\mathbbm{k}} N$
	      by the relations
	      \begin{equation*}
		      \left( a \cdot m \right) \otimes_{\mathbbm{k}} n =
			      (-1)^{\braidd{a}{m}_i} m \otimes_{\mathbbm{k}} \left( a \cdot n \right)
	      \end{equation*}
	      (see \cref{eq:algebraic-graded-tensor-product,eq:R-action-R-tensor-product}).
	      The resulting $\ZZ^2$-graded $\mathbbm{k}$-module $M \otimes_{A_i} N$, even ignoring
	      the $A_i$-action, depends on $\braidop_i$.
	\item The unit object $A_i$ of $\mathcal{C}_i$ depends on $i$.
	\item A map $f \colon M \rightharpoonup N$ is graded $A_i$-linear if it satisfies
	      \begin{equation} \label{eq:R_i-linear-map}
		      f \left( a \cdot m \right) = (-1)^{\braidd{f}{a}_i} a \cdot f \left( m \right)
	      \end{equation}
	      so the notion of a graded $A_i$-linear map depends both on $A_i$ and $\braidop_i$.
	\item The internal hom object $\InnHom{M}{N}[][A_i]$ consists of graded $A_i$-linear maps
	      and depends both on $A_i$ and $\braidop_i$.
\end{enumerate}

The symmetries $\widehat{\sigma}_i \colon M \otimes_i N \rightarrow N \otimes_i M$ of $\mathcal{C}_i \left( A_i \right)$
are induced from the symmetries $\sigma_i$ on $\mathcal{C}_i$, and are given by
\begin{equation*}
	m \otimes_{i} n \xmapsto{\widehat{\sigma}_i} (-1)^{\braidd{m}{n}_i} n \otimes_{i} m.
\end{equation*}
Like in $\mathcal{C}_i$, the symmetries
play a role in defining
the tensor product of two graded $A_i$-linear maps $f \colon M \rightharpoonup M'$ and
$g \colon N \rightharpoonup N'$ which is given by
\begin{equation*}
	\left( f \otimes_{i} g \right) \left( m \otimes_{i} n \right) =
	(-1)^{\braidd{g}{m}_i} f \left( m \right) \otimes_{i} g \left( n \right)
\end{equation*}
(see \cref{eq:tensor-product-graded-R-linear-maps}).

Now, it might be the case that we have an algebra object $R$ of $\mathcal{C}$ which is both
$\sigma_1$-commutative and $\sigma_2$-commutative. In this case, we can take $A_1 = A_2 = R$ and
then the categories $\mathcal{C}_i \left( R \right)$ have the same objects and morphisms,
i.e., the underlying categories of $\mathcal{C}_i$ for $i = 1,2$ are both $\GMod[R][\ZZ^2]$,
but even then, they usually have different monoidal products $\otimes_i = \otimes_{R, \braidop_i}$
and different closed structures $\InnHom{M}{N}[][R,\braidop_i]$ because the symmetries are different!
When $R = \mathbbm{k}$, we have $\mathcal{C}_i \left( \mathbbm{k} \right) = \mathcal{C}_i$
and in this case, the monoidal and closed structures are the same, only the symmetries are different.

\subsection{Dependence of the Noncommutative Differential Calculus on \texorpdfstring{$\braidop$}{the Parity Form}}
\label{sec:dependence-ndf-braidop}
Let $R$ be a $\ZZ$-graded, graded-commutative, $\mathbbm{k}$-algebra and let $V$ be a $\ZZ$-graded
$R$-module. We start by going over the constructions in \cref{sec:noncomm-diff-calc} and describing explicitly
the dependence of the constructions on the background parity form. Fix $i \in \Set{1,2}$.
\begin{enumerate}
	\item We think of $R$ as a $\ZZ^2$-graded $\mathbbm{k}$-algebra by placing $R^{*}$
	      in bidegree $(0,*)$, i.e., placing it in line degree zero. Since both parity forms
	      $\braidop_i$ extend the standard Koszul parity form on the second factor, the resulting
	      $\ZZ^2$-graded $\mathbbm{k}$-algebra is $\sigma_i$-commutative with respect to both parity forms.
	\item We think of $V$ as a $\ZZ^2$-graded $\mathbbm{k}$-module by placing $V^{*}$ in bidegree
	      $(0,*)$. Then $V$ becomes a $\ZZ^2$-graded module over the $\ZZ^2$-graded-commutative
	      algebra $R$, i.e., an object of $\mathcal{C}_i \left( R \right)$.
	\item We consider the module $\ul{V}_i = V[(-1,0)]_i$, which is the shift of $V$ concentrated
	      in line degree one. Elements of $\ul{V}_i$ are denoted by $\ul{v}_i$ and, by our conventions
	      regarding shifts, the $R$-action on $\ul{V}_i$ is given by
	      $r \cdot \ul{v}_i = (-1)^{\braid{(0,\degb{r})}{(1,0)}_i} \ul{r \cdot v}_i$
	      (see \cref{eq:R-action-on-suspension,eq:R-action-ul-v}).
	      The resulting $R$-module is an object of $\mathcal{C}_i \left( R \right)$
	      which depends on $\braidop_i$.
	\item We then form the tensor module
	      \begin{equation*}
		      \ndf{V}[][] = \ndf{V, \braidop_i}[][] = \tens{V \oplus \ul{V}_i}[R, \otimes_i] =
		      \bigoplus_{k=0}^{\infty} \left( V \oplus \ul{V}_i \right)^{\otimes_i k}
	      \end{equation*}
	      \textit{using the monoidal structure on} $\mathcal{C}_i \left( R \right)$.
	\item We work with $\ndf{V, \braidop_i}[][] = \tens{V \oplus \ul{V}_i}[R, \otimes_i]$ as a
	      coalgebra object of the category $\mathcal{C}_i \left( R \right)$. In particular,
	      when we discuss coderivations $\eta$ on $\ndf{V, \braidop_i}[][]$, the coderivation
	      equation
	      \begin{equation} \label{eq:braidop_i_coderivation}
		      \Delta \circ \eta =
		      \left( \eta \otimes_i \id + \id \otimes_i \eta \right) \circ \Delta
	      \end{equation}
	      involves the tensor product of graded $R$-linear maps, so there is dependence on $\braidop_i$.
	      In addition, the formula \eqref{eq:generalized-coder-coextension} describing the action of the coderivation
	      $\eta$ in terms of the corestriction $\corest{\eta}$ depends on $\braidop_i$.
	\item The action of the rotation operator $\t_i$ on $\tens{V \oplus \ul{V}_i}[R, \otimes_i]$
	      (see \cref{eq:def-t-rotation}) depends explicitly on $\braidop_i$.
	\item When we form the cyclic tensor module $\ncdf{V, \braidop_i}[][]$ as the quotient
	      \begin{equation*}
		      \ncdf{V, \braidop_i}[][] =
		      \ndf{V, \braidop_i}[][] / \Im \left( \idd - \t_i \right)=
		      \tens{V \oplus \ul{V}_i}[R, \otimes_i] / \Im \left( \idd - \t_i \right),
	      \end{equation*}
	      the dependence of the result on $\braidop_i$ enters both via the $R$-module
	      $\ndf{V, \braidop_i}[][]$, and the $R$-submodule $\Im \left( \idd - \t_i \right)$
	      by which we quotient.
	\item The definition of the cyclization of a coderivation $\eta$ on $\ndf{V, \braidop_i}[][]$
	      (see \cref{def:cyclization-coder-short} of \cref{sec:cyclization-generalized-coderivation})
	      involves signs computed using $\braidop_i$.
	      In particular, the cyclic versions $\clie{\mu}, \ccont{\mu}$, obtained as cyclizations
	      of the coderivations $\lie{\mu}, \cont{\mu}$, depend on $\braidop_i$.
\end{enumerate}

Since the objects and morphisms of $\mathcal{C}_i \left( R \right)$ are the same, it makes sense
to ask whether $\ndf{V, \braidop_i}[][]$ (resp.\ $\ncdf{V, \braidop_i}[][]$) are isomorphic as $\ZZ^2$-graded
$R$-modules, and while \textit{they are} isomorphic,
since the closed monoidal structures of $\mathcal{C}_i \left( R \right)$ are different,
any isomorphism won't respect all the structures involved.
From a categorical perspective, it actually
makes more sense to compare instead their totalizations, which are both objects of the same
monoidal category $\GMod[R][\ZZ]$ of $\ZZ$-graded $R$-modules.
Recall from \cref{appendix:bicomplexes} that the definition of the totalization also depends on the parity form $\braidop_i$.
Then we have the following:
\begin{enumerate}
	\item We have a natural isomorphism
	      \begin{equation*}
		      \Phi_2 \colon \tens{V \oplus V[1]}[R] \rightarrow \totc{\ndf{V, \braidop_2}}[][][\braidop_2]
	      \end{equation*}
	      of $\ZZ$-graded $R$-modules given by
	      \begin{equation*}
		      l^0 \otimes_R \s v_1 \otimes_R \dots \otimes_R \s v_k \otimes_R l^k \xmapsto{\Phi_2}
		      l^0 \otimes_2 \ul{v_1} \otimes_2 \dots \otimes_2 \ul{v_k} \otimes_2 l^k
	      \end{equation*}
	      which involves no signs.\footnote{In order to avoid making the notation more cumbersome than it
		      already is, instead of writing
		      $\left( l^0 \otimes_2 \ul{v_1}_2 \otimes_2 \dots \otimes_2 \ul{v_k}_2 \otimes_2 l^k \right)
			      \xmapsto{\Phi_2}
			      l^0 \otimes_R \s v_1 \otimes_R \dots \otimes_R \s v_k \otimes_R l^k$,
		      we have chosen to leave the dependence of $\ul{v}_i$ on $i$ implicit, relying on the context
		      to understand we use $\ul{v}_2$. The same applies to the rest of the section whenever the notation
		      gets too convoluted.}
	\item We have a natural isomorphism
	      \begin{equation*}
		      \Phi_1 \colon \tens{V \oplus V[1]}[R] \rightarrow \totc{\ndf{V, \braidop_1}}[][][\braidop_1]
	      \end{equation*}
	      of $\ZZ$-graded $R$-modules given by
	      \begin{equation*}
		      l^0 \otimes_R \s v_1 \otimes_R \dots \otimes_R \s v_k \otimes_R l^k
		      \xmapsto{\Phi_1}
		      (-1)^{\varepsilon}
		      \s_k \left(
		      l^0 \otimes_1 \ul{v_1} \otimes_1 \dots \otimes_1 \ul{v_k} \otimes_1 l^k
		      \right)
	      \end{equation*}
	      where the sign factor $\varepsilon$ is given by
	      \begin{equation} \label{eq:sign-factor-Phi-1}
		      \begin{aligned}
			      \varepsilon ={} &
			      \degb{l^0} + \left( \degb{l^0} + \degb{v_1} + \degb{l^1} \right) + \dots
			      \\
			                      & +
			      \left( \degb{l^0} + \degb{v_1} + \degb{l^1} + \dots + \degb{v_{k-1}} + \degb{l^{k-1}} \right).
		      \end{aligned}
	      \end{equation}
	      The sign factor can be understood informally as coming from moving
	      all the $\s$ factors in $l^0 \otimes_R \s v_1 \otimes_R \dots \otimes_R \s v_k \otimes_R l^k$
	      to the beginning and then replacing $v_i$ with $\ul{v_i}$. For an explanation
	      of the appearance of the sign factor as a consequence of non-trivial tensor constraints, we refer to
	      \cref{sec:equiv-symmetries-bi-graded-categorical}.
\end{enumerate}

Let us fix $i \in \Set{1,2}$. Given a graded $\mathbbm{k}$-linear map
\begin{equation*}
	\eta \colon \tens{V \oplus \ul{V}_i}[R, \otimes_i] \rightharpoonup \tens{V \oplus \ul{V}_i}[R, \otimes_i]
\end{equation*}
of degree $(a,b)$, set
\begin{equation*}
	\Phi_i \left[ \eta \right] \defeq \Phi_i^{-1} \circ \totc{\eta}[][][\braidop_i]
	\circ \Phi_i \colon
	\tens{V \oplus V[1]}[R] \rightharpoonup \tens{V \oplus V[1]}[R].
\end{equation*}
Then if $\eta$ is a graded $R$-linear map of degree $(a,b)$ (resp.\ a module derivation over a derivation
$d \colon R \rightharpoonup R$ of degree $(0,\degb{d})$), the map
$\Phi_i \left[ \eta \right]$ is a graded $R$-linear map of degree $b - a$ (resp.\ a module derivation
over $d$ of degree $\degb{d}$) on the tensor module $\tens{V \oplus V[1]}[R]$.
The isomorphism $\Phi_i$ is compatible with the coalgebra structures on both sides
in the following sense: When $\eta$ is a coderivation in the sense of
\cref{eq:braidop_i_coderivation} (resp.\ coalgebra morphism),
then $\Phi_i \left[ \eta \right]$ is also a coderivation
(resp.\ coalgebra morphism).

The isomorphism $\Phi_i$ commutes with the rotation operators on both sides in the sense that
$\Phi_i \left[ \t_i \right] = \t$, where $\t$ is the rotation operator on $\tens{V \oplus V[1]}[R]$.
Hence, $\Phi_i$ also commutes with the cyclization process and we have
$\Phi_i \left[ \cycl{\eta} \right] = \cycl{ \left( \Phi_i \left[ \eta \right] \right)}$.
In addition, $\Phi_i$ descends to a well-defined isomorphism
\begin{equation*}
	\Phi_i \colon \tenscyc{V \oplus V[1]} \rightarrow \totc{\ncdf{V, \braidop_i}}[][][\braidop_i],
\end{equation*}
which we continue to denote by the same name.

Let us denote the isomorphism $\Phi_1 \circ \Phi_2^{-1}$ by
\begin{equation*}
	\Phi \defeq \Phi_1 \circ \Phi_2^{-1} \colon \totc{\ndf{V, \braidop_2}}[][][\braidop_2] \rightarrow
	\totc{\ndf{V, \braidop_1}}[][][\braidop_1],
\end{equation*}
given explicitly by
\begin{equation*}
	l^0 \otimes_2 \ul{v_1} \otimes_2 \dots \otimes_2 \ul{v_k} \otimes_2 l^k
	\xmapsto{\Phi}
	(-1)^{\varepsilon}
	\s_k \left(
	l^0 \otimes_1 \ul{v_1} \otimes_1 \dots \otimes_1 \ul{v_k} \otimes_1 l^k
	\right).
\end{equation*}
Given a coderivation $\mu$ on $\tens{V}[R]$, the isomorphism $\Phi$ identifies the basic
operators of the noncommutative calculus constructed using the different pairings $\braidop_i$
with each other, i.e., we have
\begin{align}
	\Phi \circ \totc{\qdr^{\braidop_2}}[][][\braidop_2]       & =
	\totc{\qdr^{\braidop_1}}[][][\braidop_1] \circ \Phi,
	\label{eq:Phi-qdr-rel}
	\\
	\Phi \circ \totc{\lie{\mu}^{\braidop_2}}[][][\braidop_2]  & =
	\totc{\lie{\mu}^{\braidop_1}}[][][\braidop_1] \circ \Phi,
	\label{eq:Phi-lie-mu-rel}
	\\
	\Phi \circ \totc{\cont{\mu}^{\braidop_2}}[][][\braidop_2] & =
	\totc{\cont{\mu}^{\braidop_1}}[][][\braidop_1] \circ \Phi.
	\label{eq:Phi-cont-mu-rel}
\end{align}
The isomorphism $\Phi$ is also natural with respect to morphisms of tensor coalgebras.
Namely, let $\varphi \colon R \rightarrow S$ be a morphism of graded $\mathbbm{k}$-algebras,
and let $f \colon \tens{V}[R] \rightarrow \tens{W}[S]$ be a morphism of graded
coalgebras over $\varphi$. Then we have
\begin{equation}	\label{eq:Phi-indmap-rel}
	\Phi \circ \totc{\indmap{f}^{\braidop_2}}[][][\braidop_2] =
	\totc{\indmap{f}^{\braidop_1}}[][][\braidop_1] \circ \Phi,
\end{equation}
where $\indmap{f}^{\braidop_i} \colon \ndf{V/R, \braidop_i}[][] \rightarrow \ndf{W/S, \braidop_i}[][]$
is the morphism induced by $f$ on codifferential forms constructed using $\braidop_i$,
for $i = 1, 2$.

The same holds for the cyclic versions obtained via cyclization, i.e., we have
\begin{align}
	\Phi \circ \totc{\clie{\mu}^{\braidop_2}}[][][\braidop_2]  & =
	\totc{\clie{\mu}^{\braidop_1}}[][][\braidop_1] \circ \Phi, \label{eq:Phi-clie-rel}
	\\
	\Phi \circ \totc{\ccont{\mu}^{\braidop_2}}[][][\braidop_2] & =
	\totc{\ccont{\mu}^{\braidop_1}}[][][\braidop_1] \circ \Phi, \label{eq:Phi-ccont-rel}
	\\
	\Phi \circ \totc{\cindmap{f}^{\braidop_2}}[][][\braidop_2] & =
	\totc{\cindmap{f}^{\braidop_1}}[][][\braidop_1] \circ \Phi. \label{eq:Phi-cindmap-rel}
\end{align}
The isomorphism $\Phi$ also induces a well-defined isomorphism
\begin{equation*}
	\Phi \colon \totc{\ncdf{V, \braidop_2}}[][][\braidop_2]
	\rightarrow
	\totc{\ncdf{V, \braidop_1}}[][][\braidop_1]
\end{equation*}
and \cref{eq:Phi-qdr-rel,eq:Phi-clie-rel,eq:Phi-ccont-rel,eq:Phi-cindmap-rel} continue to hold on the totalizations of the quotients.
In particular, if $\mu$ is an $\Ainf$-structure, i.e., $\mu^2 = 0$, then
\begin{equation*}
	\Phi \colon
	\totc{\ncdf{V, \braidop_2}, \qdr^{\braidop_2}, \clie{\mu}^{\braidop_2}}[][][\braidop_2]
	\rightarrow
	\totc{\ncdf{V, \braidop_1}, \qdr^{\braidop_1}, \clie{\mu}^{\braidop_1}}[][][\braidop_1]
\end{equation*}
gives a \textit{natural} isomorphism between the total complexes which is a \textit{chain map}.
Since we always write elementary
cyclic codifferential forms of line degree $k$ as
\begin{equation*}
	x_k = \ul{a_1} \otimes_2 l^1 \otimes_2 \dots \otimes_2 \ul{a_k} \otimes_2 l^k \in \ncdf{A, \braidop_2}[k][],
\end{equation*}
starting with an underlined element of $A$, the isomorphism $\Phi$ acts on $x_k$ via
\begin{equation*}
	\Phi \left( x_k \right) = (-1)^{\varepsilon \left( x_k \right)}
	\s_k \left( \ul{a_1} \otimes_1 l^1 \otimes_1 \dots \otimes_1 \ul{a_k} \otimes_1 l^k \right),
\end{equation*}
where now, the sign factor $\varepsilon = \varepsilon \left( x_k \right)$ of \eqref{eq:sign-factor-Phi-1} takes the form
\begin{equation} \label{eq:sign-factor-Psi-cyc}
	\varepsilon \left( x_k \right) = \sum_{i=1}^{k} \left( k - i \right) \cdot \left( \degb{a_i} + \degb{l^i} \right).
\end{equation}

All the claims above can be verified in a straightforward manner, and also follow from the more abstract
discussion in \cref{sec:equiv-symmetries-bi-graded-categorical}.
To demonstrate how the isomorphism $\Phi$ identifies the different constructions,
we give two examples:

\begin{ex}
	Let $\mu$ be a coderivation on $\tens{V}[R]$. We compare the action of
	$\lie{\mu}^{\braidop_i}$ on $\ndf{V, \braidop_i}[1][]$ for $i \in \Set{1,2}$.
	Let $x_i = l^0 \otimes_i \ul{v}_i \otimes_i l^1 \in \ndf{V, \braidop_i}[1][]$. Then
	by \cref{eq:lie-mu-full-formula}, we have
	\begin{equation} \label{eq:lie-mu-x-1-braidop-1}
		\begin{aligned}
			\lie{\mu}^{\braidop_1} \left( x_1 \right) ={} &
			(-1)^{\degb{l^0_{(1)}} \cdot \degb{\mu}}
			l^0_{(1)} \otimes_1 \corest{\mu} \left( l^0_{(2)} \right) \otimes_1 l^0_{(3)} \otimes_1 \ul{v}_{1}
			\otimes_1 l^1
			\\
			                                              & +
			                                                (-1)^{\degb{l^0_{(1)}} \cdot \degb{\mu}}
			l^0_{(1)} \otimes_1 \ul{ \corest{\mu} \left( l^0_{(2)} \otimes_R v \otimes_R l^1_{(1)} \right) }_{1}
			\otimes_1 l^1_{(2)}
			\\
			                                              & +
			                                                (-1)^{\left( \degb{l^0} + \degb{v} + \degb{l^1_{(1)}} \right) \cdot \degb{\mu}}
			l^0 \otimes_1 \ul{v}_{1} \otimes_1 l^1_{(1)} \otimes_1 \corest{\mu} \left( l^1_{(2)} \right)
			\otimes_1 l^1_{(3)}
		\end{aligned}
	\end{equation}
	while
	\begin{equation} \label{eq:lie-mu-x-1-braidop-2}
		\begin{aligned}
			\lie{\mu}^{\braidop_2} \left( x_2 \right) ={} &
			(-1)^{\degb{l^0_{(1)}} \cdot \degb{\mu}}
			l^0_{(1)} \otimes_2 \corest{\mu} \left( l^0_{(2)} \right) \otimes_2 l^0_{(3)} \otimes_2 \ul{v}_{2}
			\otimes_2 l^1
			\\
			                                              & +
			                                                (-1)^{\degb{l^0_{(1)}} \cdot \degb{\mu} + \degb{\mu} + \degb{l^0_{(2)}}}
			l^0_{(1)} \otimes_2 \ul{ \corest{\mu} \left( l^0_{(2)} \otimes_R v \otimes_R l^1_{(1)} \right) }_{2}
			\otimes_2 l^1_{(2)}
			\\
			                                              & +
			                                                (-1)^{\left( \degb{l^0} + \degb{v} + \degb{l^1_{(1)}} + 1 \right) \cdot \degb{\mu}}
			l^0 \otimes_2 \ul{v}_{2} \otimes_2 l^1_{(1)} \otimes_2 \corest{\mu} \left( l^1_{(2)} \right)
			\otimes_2 l^1_{(3)}.
		\end{aligned}
	\end{equation}
	We have $\Phi \left( x_2 \right) = (-1)^{\degb{l^0}} \s x_1$ and so
	\begin{equation*}
		\begin{aligned}
			\left( \totc{\lie{\mu}^{\braidop_1}}[][] \circ \Phi \right) \left( x_2 \right)
			\stackrel{\phantom{\eqref{eq:tot-graded-maps-braid-op-1}}}{=} &
			(-1)^{\degb{l^0}} \totc{\lie{\mu}^{\braidop_1}}[][] \left( \s x_1 \right)
			\stackrel{\eqref{eq:tot-graded-maps-braid-op-1}}{=}
			(-1)^{\degb{l^0} + \degb{\mu}} \s \left( \lie{\mu}^{\braidop_1} \left( x_1 \right) \right)
			\\
			\stackrel{\eqref{eq:lie-mu-x-1-braidop-1}}{=}                 &
			(-1)^{\degb{l^0_{(1)}} \cdot \degb{\mu} + \degb{l^0} + \degb{\mu}}
			\s \left(
			l^0_{(1)} \otimes_1 \corest{\mu} \left( l^0_{(2)} \right) \otimes_1 l^0_{(3)} \otimes_1 \ul{v}_{1}
			\otimes_1 l^1
			\right)
			\\
			                                                              & +
			                                                                (-1)^{\degb{l^0_{(1)}} \cdot \degb{\mu} + \degb{l^0} + \degb{\mu}}
			\s \left(
			l^0_{(1)} \otimes_1 \ul{ \corest{\mu} \left( l^0_{(2)} \otimes_R v \otimes_R l^1_{(1)} \right) }_{1}
			\otimes_1 l^1_{(2)}
			\right)
			\\
			                                                              & +
			                                                                (-1)^{\left( \degb{l^0} + \degb{v} + \degb{l^1_{(1)}} + 1 \right) \cdot \degb{\mu} + \degb{l^0}}
			\s \left(
			l^0 \otimes_1 \ul{v}_{1} \otimes_1 l^1_{(1)} \otimes_1 \corest{\mu} \left( l^1_{(2)} \right)
			\otimes_1 l^1_{(3)}
			\right),
		\end{aligned}
	\end{equation*}
	where the extra signs, compared to \cref{eq:lie-mu-x-1-braidop-1}, come from both $\Phi$ and the totalization process.
	On the other hand, we have
	\begin{equation*}
		\begin{aligned}
			\MoveEqLeft
			\left( \Phi \circ \totc{\lie{\mu}^{\braidop_2}}[][][\braidop_2] \right) \left( x_2 \right)
			\stackrel{\eqref{eq:tot-graded-maps-braid-op-2}}{=}
			\Phi \left( \lie{\mu}^{\braidop_2} \left( x_2 \right) \right)
			\\
			\stackrel{\eqref{eq:lie-mu-x-1-braidop-2}}{=} &
			(-1)^{\degb{l^0} + \degb{\mu}}
			\s \left(
			   (-1)^{\degb{l^0_{(1)}} \cdot \degb{\mu}}
			l^0_{(1)} \otimes_1 \corest{\mu} \left( l^0_{(2)} \right) \otimes_1 l^0_{(3)} \otimes_1 \ul{v}_{1}
			\otimes_1 l^1
			\right)
			\\
			                                              & +
			                                                (-1)^{\degb{l^0_{(1)}}}
			\s \left(
			   (-1)^{\degb{l^0_{(1)}} \cdot \degb{\mu} + \degb{\mu} + \degb{l^0_{(2)}}}
			l^0_{(1)} \otimes_1 \ul{ \corest{\mu} \left( l^0_{(2)} \otimes_R v \otimes_R l^1_{(1)} \right) }_{1}
			\otimes_1 l^1_{(2)}
			\right)
			\\
			                                              & +
			                                                (-1)^{\degb{l^0}}
			\s \left(
			   (-1)^{\left( \degb{l^0} + \degb{v} + \degb{l^1_{(1)}} + 1 \right) \cdot \degb{\mu}}
			l^0 \otimes_1 \ul{v}_{1} \otimes_1 l^1_{(1)} \otimes_1 \corest{\mu} \left( l^1_{(2)} \right)
			\otimes_1 l^1_{(3)}
			\right).
		\end{aligned}
	\end{equation*}
	Here, the extra signs compared to \cref{eq:lie-mu-x-1-braidop-2} come only from $\Phi$. From
	the expressions above, it is clear that we indeed have
	\begin{equation*}
		\left( \Phi \circ \totc{\lie{\mu}^{\braidop_2}}[][][\braidop_2] \right) \left( x_2 \right) = \left( \totc{\lie{\mu}^{\braidop_1}}[][] \circ \Phi \right) \left( x_2 \right).
	\end{equation*}
\end{ex}

\begin{ex}
	Let $\mu$ be a coderivation on $\tens{V}[R]$.
	We compare the action of the cyclization
	$\clie{\mu}^{\braidop_i}$, which involves the rotation operator,
	on $\ndf{V, \braidop_i}[1][]$ for $i \in \Set{1,2}$.
	Let $x_i = \ul{v}_i \otimes_i l \in \ndf{V, \braidop_i}[1][]$. Then
	by \cref{eq:cyc-lie-formula-n-eq-1}, we have
	\begin{equation*}
		\begin{aligned}
			\clie{\mu}^{\braidop_1} \left( x_1 \right) ={} &
			(-1)^{\left( \degb{v} + \degb{l_{(1)}} \right) \cdot \degb{\mu}}
			\ul{v}_{1} \otimes_1 l_{(1)} \otimes_1 \corest{\mu} \left( l_{(2)} \right) \otimes_1 l_{(3)}
			\\
			                                               & +
			                                                 (-1)^{\degb{l_{(3)}} \cdot \left( \degb{v} + \degb{l_{(1)}} + \degb{l_{(2)}} \right)}
			\ul{ \corest{\mu} \left( l_{(3)} \otimes_R v \otimes_R l_{(1)} \right) }_{1} \otimes_1 l_{(2)}
		\end{aligned}
	\end{equation*}
	while
	\begin{equation} \label{eq:clie-mu-x-1-braidop-2}
		\begin{aligned}
			\clie{\mu}^{\braidop_2} \left( x_2 \right) ={} &
			(-1)^{\left( \degb{v} + \degb{l_{(1)}} \right) \cdot \degb{\mu} + \degb{\mu}}
			\ul{v}_{2} \otimes_2 l_{(1)} \otimes_2 \corest{\mu} \left( l_{(2)} \right) \otimes_2 l_{(3)}
			\\
			                                               & +
			                                                 (-1)^{\degb{l_{(3)}} \cdot \left( \degb{v} + \degb{l_{(1)}} + \degb{l_{(2)}} \right) + \degb{\mu}}
			\ul{ \corest{\mu} \left( l_{(3)} \otimes_R v \otimes_R l_{(1)} \right) }_{2} \otimes_2 l_{(2)},
		\end{aligned}
	\end{equation}
	so the difference in the sign factors is precisely $(-1)^{\degb{\mu}}$. We have
	\begin{equation*}
		\begin{aligned}
			\left( \Phi \circ \totc{\clie{\mu}^{\braidop_2}}[][][\braidop_2] \right) \left( x_2 \right) ={} &
			(-1)^{\left( \degb{v} + \degb{l_{(1)}} \right) \cdot \degb{\mu} + \degb{\mu}}
			\s \left( \ul{v}_{1} \otimes_1 l_{(1)} \otimes_1 \corest{\mu} \left( l_{(2)} \right) \otimes_1 l_{(3)}
			\right)
			\\
			                                                                                                & +
			                                                                                                  (-1)^{\degb{l_{(3)}} \cdot \left( \degb{v} + \degb{l_{(1)}} + \degb{l_{(2)}} \right) + \degb{\mu}}
			\s \left(
			\ul{ \corest{\mu} \left( l_{(3)} \otimes_R v \otimes_R l_{(1)} \right) }_{1} \otimes_1 l_{(2)}
			\right),
		\end{aligned}
	\end{equation*}
	where there are no extra signs compared to \cref{eq:clie-mu-x-1-braidop-2}
	coming from $\Phi$, since all expressions have line degree one and start with an underlined element.
	Since we also have $\Phi \left( x_2 \right) = \s \left( x_1 \right)$, we indeed have
	\begin{equation*} 
		\begin{aligned}
			\left( \totc{\clie{\mu}^{\braidop_1}}[][][\braidop_1] \circ \Phi \right) \left( x_2 \right)
			\stackrel{\phantom{\eqref{eq:tot-graded-maps-braid-op-1}}}{=} &
			\totc{\clie{\mu}^{\braidop_1}}[][][\braidop_1] \left( \s \left( x_1 \right) \right)
			\\
			\stackrel{\eqref{eq:tot-graded-maps-braid-op-1}}{=}           &
			(-1)^{\degb{\mu}} \s \left( \clie{\mu}^{\braidop_1} \left( x_1 \right) \right)
			\\
			\stackrel{\phantom{\eqref{eq:tot-graded-maps-braid-op-1}}}{=} &
			\left( \Phi \circ \totc{\clie{\mu}^{\braidop_2}}[][][\braidop_2] \right) \left( x_2 \right).
		\end{aligned}
	\end{equation*}
\end{ex}

Finally, in \cref{sec:generalized-superpotential}, we work with the total complex
\begin{equation*}
	\totc{\ncdf{A, \braidop_2}, -\qdr^{\braidop_2}, \clie{\mu}^{\braidop_2}}[][][\braidop_2],
\end{equation*}
in which the horizontal differential is $-\qdr^{\braidop_2}$
instead of $\qdr^{\braidop_2}$.
We have a natural isomorphism
\begin{equation*}
	\gamma \colon \left( \ncdf{A, \braidop_2}, -\qdr^{\braidop_2}, \clie{\mu}^{\braidop_2} \right)
	\rightarrow
	\left( \ncdf{A, \braidop_2}, \qdr^{\braidop_2}, \clie{\mu}^{\braidop_2} \right)
\end{equation*}
of bicomplexes, i.e., commuting with both differentials, given by
$\gamma \left( x_k \right) = (-1)^{k - 1} x_k$ for
$x_k \in \ncdf{A, \braidop_2}[k][]$.
Let $\Gamma \defeq	\totc{\gamma}[][][\braidop_2]$. The composition
$\Psi = \Phi \circ \Gamma$ gives us a \textit{natural} isomorphism
\begin{equation*}
	\Psi \colon \totc{\ncdf{A, \braidop_2}, -\qdr^{\braidop_2}, \clie{\mu}^{\braidop_2}}[][][\braidop_2]
	\rightarrow
	\totc{\ncdf{A, \braidop_1}, \qdr^{\braidop_1}, \clie{\mu}^{\braidop_1}}[][][\braidop_1],
\end{equation*}
which is a \textit{chain map}, acting via
\begin{equation} \label{eq:Psi-explicit-formula}
	\Psi \left( \ul{a_1} \otimes_2 l^1 \otimes_2 \dots \otimes_2 \ul{a_k} \otimes_2 l^k \right)
	= (-1)^{\varepsilon \left( x_k \right) + k - 1}
	\s_k \left( \ul{a_1} \otimes_1 l^1 \otimes_1 \dots \otimes_1 \ul{a_k} \otimes_1 l^k \right).
\end{equation}

\subsection{Equivalence of the Inner Product and Total Degree Parity Forms}
\label{sec:equiv-symmetries-bi-graded-categorical}

\subsubsection{Equivalence for \texorpdfstring{$\ZZ^2$-graded $\mathbbm{k}$-modules}{Bigraded k-modules}}
\label{sec:braidop-eq-graded-k-modules}

Recall from \cref{sec:cat-dependence-braidop} the symmetric closed monoidal categories
\begin{equation*}
	\mathcal{C}_i =
	\left( \GMod[\mathbbm{k}][\ZZ^2], \otimes_{\mathbbm{k}}, \mathbbm{k}, \sigma_i \right)
\end{equation*}
of $\ZZ^2$-graded $\mathbbm{k}$-modules endowed with the symmetries
$\sigma_i = \sigma_i^{M,N} \colon M \otimes_{\mathbbm{k}} N \rightarrow N \otimes_{\mathbbm{k}} M$
given by
\begin{equation*}
	m \otimes_{\mathbbm{k}} n \xmapsto{\sigma_i} (-1)^{\braid{\degb{m}}{\degb{n}}_i} n \otimes_{\mathbbm{k}} m,
\end{equation*}
for $i \in \Set{1,2}$.
The categories $\mathcal{C}_i$ have
the same objects and morphisms, the same tensor product $\otimes_{\mathbbm{k}}$, the same unit $\mathbbm{k}$
and the same internal hom object $\InnHom{M}{N}[][\mathbbm{k}]$ consisting of graded maps,
but are endowed with different symmetries $\sigma_i$.

In the context of supergeometry, the symmetry $\sigma_1$ is sometimes called the \textbf{Deligne symmetry}, while the symmetry
$\sigma_2$ is called the \textbf{Bernstein symmetry} (see \cite{SignsInSupergeometry}).
To describe the relation between the pairings, it will be useful to introduce the \textit{non-symmetric} map
$\varphi \colon \ZZ^2 \times \ZZ^2 \rightarrow \ZZ_2$ given by
\begin{equation} \label{eq:varphi-formula}
	\varphi \left( \left( a_1, a_2 \right), \left( b_1, b_2 \right) \right) \defeq a_2 \cdot b_1 \mod 2.
\end{equation}

\begin{lm}
	Given $a,b \in \ZZ^2$, the relation between $\braidop_1$ and $\braidop_2$ in terms
	of $\varphi$ is given by
	\begin{equation} \label{eq:braid-phi-relation}
		\varphi \left( a, b \right) + \braid{a}{b}_1 = \braid{a}{b}_2 + \varphi \left( b, a \right).
	\end{equation}
\end{lm}
\begin{proof}
	Writing $a = \left( a_1, a_2 \right), b = \left( b_1, b_2 \right)$, we have
	\begin{equation*}
		\begin{aligned}
			\braid{a}{b}_2 + \varphi \left( b, a \right) & =
			\left( a_2 - a_1 \right) \cdot \left( b_2 - b_1 \right) + a_1 \cdot b_2
			\\
			                                             & = \left( a_1 \cdot b_1 + a_2 \cdot b_2 \right) - a_2 \cdot b_1
			\\
			                                             & \equiv
			\varphi \left( a, b \right) + \braid{a}{b}_1 \mod 2.
		\end{aligned}
	\end{equation*}
\end{proof}

The relation \eqref{eq:braid-phi-relation} implies that the symmetric monoidal categories $\mathcal{C}_1$ and
$\mathcal{C}_2$ are equivalent. More precisely, we have the following proposition:
\begin{prop}
	The identity functor $\mathcal{F} \colon \mathcal{C}_1 \rightarrow \mathcal{C}_2$, endowed with
	the unit constraint
	$\mu_0 \colon \mathbbm{k} \rightarrow \mathcal{F} \left( \mathbbm{k} \right)$ given by the identity,
	and with the tensor constraints
	\begin{equation*}
		\mu_{M,N} = \mu^{\mathcal{F}}_{M,N} \colon
		\mathcal{F} \left( M \right) \otimes_{\mathbbm{k}} \mathcal{F} \left( N \right)
		\rightarrow \mathcal{F} \left( M \otimes N \right)
	\end{equation*}
	given by
	\begin{equation} \label{eq:tensor-constraints-mathcal-F-k-modules}
		\mu_{M,N} \left( m \otimes_{\mathbbm{k}} n \right) \defeq
		(-1)^{\varphi \left( \degb{m}, \degb{n} \right)} m \otimes_{\mathbbm{k}} n
	\end{equation}
	is a strong symmetric monoidal functor which gives us an equivalence of symmetric monoidal categories
	between $\mathcal{C}_1$ and $\mathcal{C}_2$.
\end{prop}
\begin{proof}
	The proof appears in \cite[Proposition 4.3]{SignsInSupergeometry} in the context of supergeometry. In
	\cite{SignsInSupergeometry}, one works with the category of chain complexes of super vector spaces
	endowed with the two symmetries $\sigma_1$ and $\sigma_2$. The objects are
	$\ZZ_2 \times \ZZ$-graded
	$\mathbbm{k}$-modules $M = \left( M_i^j \right)$ endowed with a differential
	$d \colon M_i^j \rightharpoonup M_i^{j-1}$. Ignoring the differential and working with
	$\ZZ^2$-graded modules
	instead of $\ZZ_2 \times \ZZ$-graded modules, the proof in \cite{SignsInSupergeometry} goes
	through verbatim. The condition that
	$\mathcal{F}$ is \textit{symmetric} is expressed via the following
	commutative diagram
	\begin{equation*}
		\adjustbox{scale=0.85,center}{
			\begin{tikzcd}
				{\mathcal{F} \left( M \right) \otimes_{\mathbbm{k}} \mathcal{F} \left( N \right)} &&
				{\mathcal{F} \left( N \right) \otimes_{\mathbbm{k}} \mathcal{F} \left( M \right)} &&
				{M \otimes_{\mathbbm{k}} N} && {N \otimes_{\mathbbm{k}} M}
				\\
				&&& \iff
				\\
				{\mathcal{F} \left( M \otimes_{\mathbbm{k}} N \right)} &&
				{\mathcal{F} \left( N \otimes_{\mathbbm{k}} M \right)} &&
				{M \otimes_{\mathbbm{k}} N} && {N \otimes_{\mathbbm{k}} M}
				\arrow["{\sigma_2^{\mathcal{F} \left( M \right), \mathcal{F} \left( N \right)}}", from=1-1, to=1-3]
				\arrow["{\mu_{M,N}}"', from=1-1, to=3-1]
				\arrow["{\mu_{N,M}}", from=1-3, to=3-3]
				\arrow["{(-1)^{\braid{\degb{m}}{\degb{n}}_2}}", from=1-5, to=1-7]
				\arrow["{(-1)^{\varphi \left( \degb{m}, \degb{n} \right)}}"', from=1-5, to=3-5]
				\arrow["{(-1)^{\varphi \left( \degb{n}, \degb{m} \right)}}", from=1-7, to=3-7]
				\arrow["{\mathcal{F} \left( \sigma_1^{M,N} \right)}"', from=3-1, to=3-3]
				\arrow["{(-1)^{\braid{\degb{m}}{\degb{n}}_1}}"', from=3-5, to=3-7]
			\end{tikzcd}
		}
	\end{equation*}
	whose commutativity is equivalent to relation \eqref{eq:braid-phi-relation}.
\end{proof}

Consider also the category
\begin{equation*}
	\mathcal{D} \defeq
	\left( \GMod[\mathbbm{k}][\ZZ], \otimes_{\mathbbm{k}}, \mathbbm{k}, \sigma \right)
\end{equation*}
of $\ZZ$-graded $\mathbbm{k}$-modules, endowed with the Koszul symmetries $\sigma$. We have
a totalization functor
\begin{equation*}
	\tot \colon \GMod[\mathbbm{k}][\ZZ^2] \rightarrow \GMod[\mathbbm{k}][\ZZ]
\end{equation*}
given by $\totc{M}[][][k] \defeq \oplus_{i} M_i^{i+k}$ on objects and by
$\totc{f}[][][] \left( \sum_{i} m_i^{i+k} \right) = \sum_{i} f \left( m_i^{i+k} \right)$
on morphisms. To enhance $\tot$ into a \textit{symmetric} strong monoidal functor between
$\mathcal{C}_i$ and $\mathcal{D}$, one needs to check how the symmetry maps $\sigma_i$ relate
to $\sigma$ after totalization. When working with $\braidop_2$, thinking of an element
$m_i^j \otimes n_k^l$ as an element of $M \otimes N$ in $\mathcal{C}_2$, the symmetry $\sigma_2$ sends it
to $(-1)^{(j - i) \cdot (l - k)} n_k^l \otimes_{\mathbbm{k}} m_i^j$, while
thinking of each $m_i^j, n_k^l$ as degree $j - i, l - k$ elements of
$\tot \left( M \right), \tot \left( N \right)$ respectively, the Koszul symmetry $\sigma$ sends
$m \otimes n$ to the same element (see also \cref{rem:identification-tot-2}). This means that we can use trivial tensor constraints for $\braidop_2$ but
non-trivial constraints are necessary for $\braidop_1$.
Let us denote by
$\tot^{\braidop_2} \colon \mathcal{C}_2 \rightarrow \mathcal{D}$ the functor $\tot$ endowed
with the trivial unit constraint and the trivial tensor constraints given by
\begin{equation*}
	m_i^j \otimes_{\mathbbm{k}} n_k^l \xmapsto{\mu_{M,N}^{\tot^{\braidop_2}}}
	m_i^j \otimes_{\mathbbm{k}} n_k^l
\end{equation*}
and denote by
$\tot^{\braidop_1} \colon \mathcal{C}_1 \rightarrow \mathcal{D}$ the functor $\tot$ endowed
with the trivial unit constraint and the tensor constraints
$\totc{M}[][][\braidop_1] \otimes_{\mathbbm{k}} \totc{N}[][][\braidop_1] \rightarrow
	\totc{M \otimes_{\mathbbm{k}} N}[][][\braidop_1]$ constructed using $\varphi$ by
\begin{equation*}
	\s_i m_i^j \otimes_{\mathbbm{k}} \s_k n_k^l \xmapsto{\mu_{M,N}^{\tot^{\braidop_1}}}
	(-1)^{\varphi \left( \degb{m}, \degb{n} \right)}
	\s_{i+k} \left( m_i^j \otimes_{\mathbbm{k}} n_k^l \right) =
	(-1)^{jk} \s_{i + k} \left( m_i^j \otimes_{\mathbbm{k}} n_k^l \right).
\end{equation*}
Here, we adopt the convention that when $M \in \mathcal{C}_1$ and $m = m_i^j \in M$ then
the corresponding (same) element of $\tot^{\braidop_1} \left( M \right)$ having degree $j - i$
is denoted by $\s_i m_i^j$,
and a general element of $\tot^{\braidop_1} \left( M \right)$ of degree $k$ is denoted by
$\sum_{i} \s_i m_i^{i + k}$ (see also \cref{rem:identification-tot-1}). On the contrary, when $M \in \mathcal{C}_2$ and $m = m_i^j \in M$
then the corresponding (same) element of $\tot^{\braidop_2} \left( M \right)$ is still denoted by
$m_i^j$, relying on context to understand whether we think of it as an element of bidegree $(i,j)$
of $M$ or an element of total degree $j - i$ of $\tot^{\braidop_2} \left( M \right)$
(see also \cref{rem:identification-tot-2}).
This allows us to differentiate between $\tot^{\braidop_i}$ and will be extremely convenient for the rest of the discussion.

The three functors $\mathcal{F}, \tot^{\braidop_1}$ and $\tot^{\braidop_2}$ sit in the following commutative diagram
\begin{equation} \label{eq:diag-f-tot-1-tot-2-k-modules}
	\begin{tikzcd}
		{\mathcal{C}_1} &&&&&& {\mathcal{C}_2} \\
		\\
		\\
		&&& {\mathcal{D}}
		\arrow["{\mathcal{F}}", "(-1)^{\varphi \left( \degb{m}, \degb{n} \right)}"', from=1-1, to=1-7]
		\arrow["{\tot^{\braidop_1}}"', "(-1)^{\varphi \left( \degb{m}, \degb{n} \right)}", from=1-1, to=4-4]
		\arrow["{\tot^{\braidop_2}}", "(-1)^{0}"',
			from=1-7, to=4-4]
	\end{tikzcd}
\end{equation}
of strong symmetric monoidal functors, where the composition is in the sense of monoidal functors (i.e.,
the tensor constraints are also composed) and we emphasize which functors have non-trivial tensor constraints
by writing the signs appearing in the tensor constraints in the diagram. The underlying functor
of $\mathcal{F}$ is the identity functor and the underlying functors of $\tot^{\braidop_i}$ coincide with $\tot$
but have different tensor constraints.

\begin{rem}
	The only reason the tensor constraints $\mu^{\mathcal{F}}_{M,N}$ (resp.\ $\mu^{\tot^{\braidop_1}}_{M,N}$)
	are not trivial and involve a sign is to make the identity (resp.\ totalization)
	functor \textit{symmetric}.
\end{rem}

\subsubsection{Equivalence for \texorpdfstring{$\ZZ^2$-graded}{Bigraded} Algebras} \label{sec:equivalence-pairing-algebras}
The fact that the categories $\mathcal{C}_i$ are symmetric monoidal allows us to talk about
commutative algebra objects in $\mathcal{C}_i$. Since the monoidal structure of $\mathcal{C}_1$ and $\mathcal{C}_2$ are strictly the same, there is no difference between algebra objects in $\mathcal{C}_1$ or in $\mathcal{C}_2$.
Unwinding the definitions, one sees that an algebra object $A$ in $\mathcal{C}_i$ is the same thing as a $\ZZ^2$-graded, associative, unital, $\mathbbm{k}$-algebra. An algebra object $A$ in $\mathcal{C}_i$ is commutative,
or $\sigma_i$-commutative if we want to emphasize the symmetry involved,
if
\begin{equation*}
	a \cdot b = m \left( a \otimes b \right) =
	(-1)^{\braidd{a}{b}_i} m \left( b \otimes a \right) =
	(-1)^{\braidd{a}{b}_i} b \cdot a
\end{equation*}
for all $a,b \in A$.

By abstract nonsense (see \cref{subsec:alg-in-monoidal-cat}), the functor $\mathcal{F}$ induces
a strong symmetric monoidal equivalence
$\mathcal{F}_{\CAlg} \colon \CAlg[\mathcal{C}_1] \rightarrow \CAlg[\mathcal{C}_2]$ converting $\sigma_1$-commutative algebras to $\sigma_2$-commutative algebras.
Given $A \in \CAlg[\mathcal{C}_1]$, we will denote the resulting algebra $\mathcal{F}_{\CAlg} \left( A \right)$ by $\wt{A}$. The algebra $\wt{A}$ coincides with $A$ as a $\ZZ^2$-graded $\mathbbm{k}$-module but has a different product coming from the tensor constraints of $\mathcal{F}$.
To differentiate between the product on $A$ and $\wt{A}$, given an element $a \in A$,
we will denote the corresponding (same) element of $\wt{A}$ by $\wt{a}$.
Then the product on $\wt{A}$ is given by
\begin{equation} \label{eq:product-wt-A}
	\wt{a} \cdot \wt{b} = (-1)^{\varphi \left( \degb{a}, \degb{b} \right)} \wt{ a \cdot b},
\end{equation}
i.e., the product is twisted using $\varphi$.
Note that even though the functor $\mathcal{F}$ is the identity functor,
the induced functor $\mathcal{F}_{\CAlg}$ is not the identity functor because the product
on $\wt{A}$ is twisted.

Similarly, one can obtain induced totalization functors
$\tot^{\braidop_i}_{\CAlg} \colon \CAlg[\mathcal{C}_i] \rightarrow \CAlg[\mathcal{D}]$
for $i \in \Set{1,2}$ which endow the totalization of a $\sigma_i$-commutative algebra with the structure of
a $\ZZ$-graded-commutative $\mathbbm{k}$-algebra (using Koszul signs).
The functor
$\tot^{\braidop_2}_{\CAlg}$ does not twist the multiplication by a sign, while
$\tot^{\braidop_1}_{\CAlg}$ does twist the multiplication, i.e., we have
\begin{equation*}
	\s_i a_i^j \cdot \s_k b_k^l =
		(-1)^{\varphi \left( \degb{a}, \degb{b} \right)} \s_{i+k} \left( a_i^j \cdot b_k^l \right) =
	(-1)^{jk} \s_{i+k} \left( a_i^j \cdot b_k^l \right),
\end{equation*}
and the diagram \eqref{eq:diag-f-tot-1-tot-2-k-modules} induces a commutative diagram
\begin{equation} \label{eq:diag-f-tot-1-tot-2-k-algebras}
	\begin{tikzcd}
		{\CAlg[\mathcal{C}_1]} &&&&&& {\CAlg[\mathcal{C}_2]} \\
		\\
		\\
		&&& {\CAlg[\mathcal{D}]}
		\arrow["{\mathcal{F}_{\CAlg}}", "(-1)^{\varphi \left( \degb{m}, \degb{n} \right)}"',
			from=1-1, to=1-7]
		\arrow["{\tot^{\braidop_1}_{\CAlg}}"', "(-1)^{\varphi \left( \degb{m}, \degb{n} \right)}",
			from=1-1, to=4-4]
		\arrow["{\tot^{\braidop_2}_{\CAlg}}", "(-1)^{0}"',
			from=1-7, to=4-4]
	\end{tikzcd}
\end{equation}
of strong symmetric monoidal functors.
Unlike the situation of diagram \eqref{eq:diag-f-tot-1-tot-2-k-modules},
all three underlying functors appearing in diagram \eqref{eq:diag-f-tot-1-tot-2-k-algebras} are
in general different.
Both $\mathcal{F}_{\CAlg}$ and $\tot^{\braidop_1}_{\CAlg}$
twist the algebra structures, while $\tot^{\braidop_2}_{\CAlg}$ does not.

\subsubsection{Equivalence for Graded Modules over Graded-Commutative Algebras} \hfill \\
\label{sec:braidop-eq-graded-modules-over-graded-algebras}
Let $A_1 \in \CAlg[\mathcal{C}_1]$ be a $\sigma_1$-commutative algebra object
and let $\wt{A_1} = A_2$ be the corresponding $\sigma_2$-commutative algebra object discussed in
\cref{sec:equivalence-pairing-algebras}.
Recall that given $i \in \Set{1,2}$, we denoted by $\mathcal{C}_i \left( A_i \right)$ the closed symmetric monoidal category of left $A_i$-modules endowed with the tensor product $\otimes_i = \otimes_{A_i, \braidop_i}$
and the symmetries coming from $\braidop_i$.
By abstract nonsense (see \cref{subsec:modules-in-monoidal-cat}),
the functor $\mathcal{F}$ induces
a strong symmetric monoidal equivalence
$\mathcal{F}_{\Mod} \colon \mathcal{C}_1 \left( A_1 \right) \rightarrow
	\mathcal{C}_2 \left( A_2 \right)$ converting modules over $A_1$ to modules over $A_2$.

Given $M = M_1 \in \mathcal{C}_1 \left( A_1 \right)$, denote the resulting
$\wt{A_1}$-module $\mathcal{F}_{\Mod} \left( M \right) \in \mathcal{C}_2 \left( A_2 \right)$
by $\wt{M} = M_2$. Since the underlying functor of $\mathcal{F}$ is the identity functor,
the module $\wt{M}$ coincides with $M$ as a $\ZZ^2$-graded $\mathbbm{k}$-module and
is endowed with an $A_2$-action using the $A_1$-action on $M$ and the tensor constraints of
$\mathcal{F}$. Given an element $m \in M$, we will denote the corresponding (same) element of
$\wt{M}$ by $\wt{m}$. The tensor constraints
\begin{equation} \label{eq:tensor-constraints-F-mod-signature}
	\mu^{\mathcal{F}_{\Mod}}_{M,N} \colon \wt{M} \otimes_{2} \wt{N} \rightarrow \wt{M \otimes_1 N}
\end{equation}
of $\mathcal{F}_{\Mod}$ are the isomorphisms inherited from the tensor constraints of $\mathcal{F}$,
given by
\begin{equation} \label{eq:tensor-constraints-F-mod}
	\wt{m} \otimes_{2} \wt{n} \xmapsto{\mu^{\mathcal{F}_{\Mod}}_{M,N}}
	(-1)^{\varphi \left( \degb{m}, \degb{n} \right)} \wt{ m \otimes_{1} n},
\end{equation}
and the trivial unit constraint
\begin{equation} \label{eq:unit-constraint-F-mod}
	\monunit_{\mathcal{C}_2 \left( A_2 \right)} = A_2 \xrightarrow{\id = \mu_0^{\mathcal{F}_{\Mod}}} \mathcal{F}_{\Mod} \left( \monunit_{\mathcal{C}_1 \left( A_1 \right)} \right) = \wt{A_1} = A_2,
\end{equation}
also inherited from $\mathcal{F}$.
Thus, the action of $A_2 = \wt{A_1}$ on $\wt{M}$ is given by
\begin{equation} \label{eq:action-wt-A-2-on-wt-M}
	\wt{a} \cdot \wt{m} = (-1)^{\varphi \left( \degb{a}, \degb{m} \right)} \wt{ a \cdot m}
\end{equation}
for $a \in A_1$ and $m \in M$.

Given two $A_1$-modules $M, N \in \mathcal{C}_1 \left( A_1 \right)$, we denote the
action of the functor $\mathcal{F}_{\Mod}$ on morphisms $f \colon M \rightarrow N$ by
$\mathcal{F}_{\Mod} \left( f \right) = \wt{f}$. It is given explicitly by
\begin{equation} \label{eq:action-F-mod-on-morphisms}
	\wt{f} \left( \wt{m} \right) = \wt{ f \left( m \right)}.
\end{equation}
The functor $\mathcal{F}_{\Mod}$ also induces maps converting graded $A_1$-linear maps
to graded $A_2$-linear maps as follows:
Consider the natural isomorphisms
\begin{equation*}
	\eta_{M,N} \colon \wt{\InnHom{M}{N}[][A_1]} \rightarrow \InnHom{\wt{M}}{\wt{N}}[][A_2]
\end{equation*}
coming from the tensor-hom adjunction and the tensor constraints
via
\begin{equation*}
	\begin{aligned}
		\eta_{M,N} & \in \Hom{\wt{\InnHom{M}{N}[][A_1]}}{\InnHom{\wt{M}}{\wt{N}}[][A_2]}[][A_2] \\
		           & \cong
		\Hom{\wt{\InnHom{M}{N}[][A_1]} \otimes_{2} \wt{M}}{\wt{N}}[][A_2]
		\\
		           & \cong
		\Hom{\wt{\InnHom{M}{N}[][A_1] \otimes_{1} M}}{\wt{N}}[][A_2] \ni
		\wt{ \mathrm{ev}_{M,N} },
	\end{aligned}
\end{equation*}
where $\mathrm{ev}_{M,N} \colon \InnHom{M}{N}[][A_1] \otimes_1 M \rightarrow N$ is the natural
evaluation map (the counit of the tensor-hom adjunction in $\mathcal{C}_1 \left( A_1 \right)$).
The isomorphism $\eta$ allows us to convert graded $A_1$-linear maps $f \colon M \rightharpoonup N$
to graded $A_2$-linear maps $\eta ( \wt{f} ) \colon \wt{M} \rightharpoonup \wt{N}$
of the same degree. To simplify the notation, given $f \in \InnHom{M}{N}[][A_1]$, we will denote
the action of $\eta ( \wt{f} )$ on $\wt{m} \in \wt{M}$ simply by $\wt{f} \left( \wt{m} \right)$,
leaving the isomorphism $\eta$ implicit. Thus, when we think of
$\wt{f}$ as acting on $\wt{M}$, we let it act through $\eta$.
Explicitly, using our tilde notation, we have
\begin{equation} \label{eq:action-F-mod-on-graded-morphisms}
	\wt{f} \left ( \wt{m} \right) =
	\left( \eta_{M,N} \left( \wt{f} \right) \right) \left( \wt{m} \right) =
	(-1)^{\varphi \left( \degb{f}, \degb{m} \right)} \wt{f \left( m \right)}
\end{equation}
for $f \in \InnHom{M}{N}[][A_1]$ and $m \in M$. Note that the notation
of \eqref{eq:action-F-mod-on-graded-morphisms} is consistent with both the action of
$\mathcal{F}_{\Mod}$ on morphisms (graded \textit{degree zero} $A_2$-linear maps) given
by \cref{eq:action-F-mod-on-morphisms}, and with
the $A_i$-actions of both sides given by \cref{eq:action-wt-A-2-on-wt-M,eq:left-R-action-on-hom},
i.e., we have
\begin{equation*}
	\begin{aligned}
		\left( \wt{a} \cdot \wt{f} \right) \left( \wt{m} \right)
		 & =
		(-1)^{\varphi \left( \degb{a}, \degb{f} \right)}
		\wt{\left( a \cdot f \right)} \left( \wt{m} \right)
		=
		(-1)^{\varphi \left( \degb{a}, \degb{f} \right) +
			\varphi \left( \degb{a} + \degb{f}, \degb{m} \right)}
		\wt{ \left( a \cdot f \right) \left( m \right)}
		\\
		 & =
		(-1)^{\varphi \left( \degb{a}, \degb{f} \right) +
			\varphi \left( \degb{a} + \degb{f}, \degb{m} \right)}
		\wt{ a \cdot f \left( m \right) }
		\\
		 & =
		(-1)^{\varphi \left( \degb{a}, \degb{f} \right) +
			\varphi \left( \degb{a} + \degb{f}, \degb{m} \right) +
			\varphi \left( \degb{a}, \degb{f} + \degb{m} \right)}
		\wt{a} \cdot \wt{f \left( m \right)}
		\\
		 & =
		(-1)^{\varphi \left( \degb{a}, \degb{f} \right) +
			\varphi \left( \degb{a} + \degb{f}, \degb{m} \right) +
			\varphi \left( \degb{a}, \degb{f} + \degb{m} \right) +
			\varphi \left( \degb{f}, \degb{m} \right)}
		\wt{a} \cdot \left( \wt{f} \left( \wt{m} \right) \right)
		\\
		 & =
		\wt{a} \cdot \left( \wt{f} \left( \wt{m} \right) \right)
	\end{aligned}
\end{equation*}
for $a \in A_1, m \in M$ and $f \in \InnHom{M}{N}[][A_1]$.

With the definitions and our choice of notation, one can easily verify that the following relations hold:
\begin{enumerate}
	\item{(Compatibility with Composition)}
	      Given two graded $A_1$-linear maps
	      $f \colon M \rightharpoonup N$ and $g \colon L \rightharpoonup M$, we have
	      \begin{equation} \label{eq:wt-g-circ-f}
		      \wt{f \circ g} =
		      (-1)^{\varphi \left( \degb{f}, \degb{g} \right)} \wt{f} \circ \wt{g} \colon
		      \wt{L} \rightharpoonup \wt{N}.
	      \end{equation}
	\item{(Compatibility with Graded Commutator)}
	      Given two graded $A_1$-linear maps
	      $f \colon M \rightharpoonup M$ and $g \colon M \rightharpoonup M$, we have
	      \begin{equation} \label{eq:comp-wt-graded-comm}
		      \wt{\left[ f, g \right]_1} = (-1)^{\varphi \left( \degb{f}, \degb{g} \right)}
		      \left[ \wt{f}, \wt{g} \right]_2 \colon \wt{M} \rightharpoonup \wt{M},
	      \end{equation}
	      where $\left[ \cdot, \cdot \right]_i$ denotes the graded commutator taken in
	      $\mathcal{C}_i \left( A_i \right)$ for $i \in \Set{1,2}$, i.e.,
	      \begin{equation*}
		      \left[ u, v \right]_i = u \circ v - (-1)^{\braidd{u}{v}_i} v \circ u.
	      \end{equation*}
	\item{(Compatibility with Tensor Product of Graded Maps)}
	      Given two graded $A_1$-linear maps
	      $f \colon M \rightharpoonup M'$ and $g \colon N \rightharpoonup N'$,
	      the following diagrams commutes
	      \begin{equation*}
		      \begin{tikzcd}
			      {\wt{M \otimes_1 N}} &&& {\wt{M' \otimes_1 N'}} \\
			      {\wt{M} \otimes_{2} \wt{N}} &&&
			      {\wt{M'} \otimes_{2} \wt{N'}}
			      \arrow["{(-1)^{\varphi \left( \degb{f}, \degb{g} \right)} \wt{f \otimes_1 g}}",
				      harpoon, from=1-1, to=1-4]
			      \arrow["{\mu^{\mathcal{F}_{\Mod}}_{M,N}}", from=2-1, to=1-1]
			      \arrow["{\wt{f} \otimes_{2} \wt{g}}"', harpoon, from=2-1, to=2-4]
			      \arrow["{\mu^{\mathcal{F}_{\Mod}}_{M',N'}}"', from=2-4, to=1-4]
		      \end{tikzcd}
	      \end{equation*}
	      (compare the horizontal arrows to the tensor constraints given by \cref{eq:tensor-constraints-F-mod}).
\end{enumerate}

Let us denote by
\begin{equation*}
	A \defeq \tot^{\braidop_1}_{\CAlg} \left( A_1 \right) =
	\tot^{\braidop_2}_{\CAlg} \left( A_2 \right)
\end{equation*}
the totalization of $A_1$, or, equivalently, by the commutativity of
the diagram \eqref{eq:diag-f-tot-1-tot-2-k-algebras}, the totalization of $A_2$. Then $A$
is a $\ZZ$-graded-commutative $\mathbbm{k}$-algebra with respect to the Koszul signs.
Denote by
\begin{equation*}
	\mathcal{D} \left( A \right) \defeq
	\left( \GMod[A][\ZZ], \otimes_{A}, A, \sigma \right)
\end{equation*}
the category of $\ZZ$-graded $A$-modules, endowed with the tensor product $\otimes_A$
and the Koszul symmetries $\sigma$.
Similarly to the discussion above, one can obtain induced totalization functors
\begin{equation*}
	\tot^{\braidop_i}_{\Mod} \colon \mathcal{C}_i \left( A_i \right) \rightarrow
	\mathcal{D} \left( \tot^{\braidop_i}_{\CAlg} \left( A_i \right) \right)
\end{equation*}
for $i \in \Set{1,2}$ which endow the totalization of a
$\ZZ^2$-graded $A_i$-module with the structure of a $\ZZ$-graded module
over the $\ZZ$-graded-commutative algebra $\tot^{\braidop_i}_{\CAlg} \left( A_i \right)$.
The functor
$\tot^{\braidop_2}_{\Mod}$ does not twist the module action by a sign while
$\tot^{\braidop_1}_{\Mod}$ does twist the module action by a sign and
the diagrams \eqref{eq:diag-f-tot-1-tot-2-k-modules} and \eqref{eq:diag-f-tot-1-tot-2-k-algebras}
induce a commutative diagram
\begin{equation} \label{eq:diag-f-tot-1-tot-2-A_2-modules}
	\begin{tikzcd}
		{\mathcal{C}_1 \left( A_1 \right)} &&&&&& {\mathcal{C}_2 \left( A_2 \right)} \\
		\\
		\\
		&&& {\mathcal{D} \left( A \right)}
		\arrow["{\mathcal{F}_{\Mod}}",
			"(-1)^{\varphi \left( \degb{m}, \degb{n} \right)}"', from=1-1, to=1-7]
		\arrow["{\tot^{\braidop_1}_{\Mod}}"', "(-1)^{\varphi \left( \degb{m}, \degb{n} \right)}",
			from=1-1, to=4-4]
		\arrow["{\tot^{\braidop_2}_{\Mod}}", "(-1)^{0}"', from=1-7, to=4-4]
	\end{tikzcd}
\end{equation}
of strong symmetric monoidal functors. Note that the diagram
\eqref{eq:diag-f-tot-1-tot-2-A_2-modules} generalizes
diagram \eqref{eq:diag-f-tot-1-tot-2-k-modules} and reduces to it when
$A_2 = A_1 = A = \mathbbm{k}$. Unlike the situation of diagram \eqref{eq:diag-f-tot-1-tot-2-k-modules},
all three underlying functors of the functors appearing in diagram \eqref{eq:diag-f-tot-1-tot-2-A_2-modules}
are in general different. Both $\mathcal{F}_{\Mod}$ and $\tot^{\braidop_1}_{\Mod}$
twist the module structures, while $\tot^{\braidop_2}_{\Mod}$ does not.

\begin{rem} \label{rem:tot-sign-morphisms-from-tensor-constraints}
	Similar to the discussion for $\mathcal{F}_{\Mod}$ above, the functors
	$\tot^{\braidop_i}_{\Mod}$ also induce maps converting graded
	$A_i$-linear maps of degree $\left( a, b \right)$ to graded $A$-linear maps
	of degree $b - a$ via the tensor-hom adjunction. Since
	$\tot^{\braidop_2}_{\Mod}$ is equipped with trivial tensor constraints, given
	a graded $A_2$-linear map $f \colon M_2 \rightharpoonup M_2$, the induced map
	is given by
	\begin{equation*}
		\tot^{\braidop_2}_{\Mod} \left( f \right) \left( \sum_{i \in \ZZ} m_i^{i+d} \right) =
		\sum_{i \in \ZZ} f \left( m_i^{i+d} \right),
	\end{equation*}
	as in \cref{eq:tot-graded-maps-braid-op-2}, while for a graded $A_1$-linear map
	$f \colon M_1 \rightharpoonup M_1$, the induced map is given by
	\begin{equation*}
		\tot^{\braidop_1}_{\Mod} \left( f \right)
		\left( \sum_{i \in \ZZ} s_i m_i^{i+d} \right) =
		\sum_{i \in \ZZ}
		(-1)^{\varphi \left( \left( a, b \right), \left( i, i + d \right) \right)}
		\s_{i + a} f \left( m_i^{i+d} \right) =
		\sum_{i \in \ZZ} (-1)^{bi} \s_{i+a} f \left( m_i^{i+d} \right),
	\end{equation*}
	as in \cref{eq:tot-graded-maps-braid-op-1}.
\end{rem}

\begin{rem} \label{rem:extension-working-over-pre-cdga}
	The equivalence $\mathcal{F}_{\Mod}$ discussed above extends to the
	setting of Banach bicomplexes. Assume that $A_1$ is equipped with a graded
	algebra derivation $d_1 \colon A_1 \rightharpoonup A_1$ of degree
	$\left( 0, 1 \right)$, which is a differential. Assume also that $M = M_1$ is equipped
	with a horizontal $A_1$-linear differential $b \colon M \rightharpoonup M$
	of degree $\left( -1 ,0 \right)$, and with a vertical differential
	$\delta \colon M \rightharpoonup M$ of degree $\left( 0, 1 \right)$ which
	is required to be a module derivation over $d_1$. The differentials $b$ and
	$\delta$ are assumed to satisfy $\left[ b, \delta \right]_1 = 0$, i.e.,
	$\left( M, b, \delta \right)$ is a bicomplex with \textit{commuting} differentials
	(see \cref{subsec:bicomplex-commuting-differentials}).

	Then $\left( \wt{M}, \wt{b}, \wt{\delta} \right)$ is a bicomplex with \textit{anticommuting}
	differentials over $\left( \wt{A_1}, \wt{d} \right) \defeq \left( A_2, d_2 \right)$, i.e.,
	$\wt{b}$ is an $A_2$-linear differential, $\wt{\delta}$ is a differential which
	is also module derivation over
	$d_2$, and $\left[ \wt{b}, \wt{\delta} \right]_2 = 0$ (see \cref{subsec:bicomplex-anticommuting-differentials}).
	Note that we have
	\begin{align*}
		\wt{b} \left( \wt{m_i^j} \right)      & =
		(-1)^{\varphi \left( \left( -1, 0 \right), \left( i, j \right) \right)}
		\wt{b \left( m_i^j \right)} = \wt{b \left( m_i^j \right)},
		\\
		\wt{\delta} \left( \wt{m_i^j} \right) & =
		(-1)^{\varphi \left( \left( 0, 1 \right), \left( i, j \right) \right)}
		\wt{\delta \left( m_i^j \right)} = (-1)^i \wt{\delta \left( m_i^j \right)},
	\end{align*}
	so that the horizontal differential $b$ stays the same, while the vertical differential
	$\delta$ is twisted by a sign when acting on odd columns.

	The extension of $\mathcal{F}_{\Mod}$ described above is compatible with the
	totalization functors for bicomplexes (see \cref{appendix:bicomplexes}),
	so that the commutative diagram \eqref{eq:diag-f-tot-1-tot-2-A_2-modules} continues
	to hold in the more general setting, replacing
	the categories $\mathcal{C}_i \left( A_i \right)$ with the categories
	of bicomplexes over $A_i$ (working with $\braidop_i$), and the category
	$\mathcal{D} \left( A \right)$ with the category of differential graded modules
	over the differential graded-commutative algebra $\mathcal{A} = \left( A, d \right)$,
	where $d$ is the totalization of $d_1$.

	Further generalizations are possible. One can work with pre-differentials instead
	of differentials or allow $A_1$ to be equipped with two (pre)-differentials.
	The interested reader is invited to come up with the appropriate formulations and
	fill in the details.
\end{rem}

\begin{rem}
	If done carefully, the statements of this section can be rephrased to hold in the
	abstract categorical setting, starting with equivalent closed symmetric monoidal categories.
	For example, \cref{eq:wt-g-circ-f} can be equivalently stated as the commutativity of the
	diagram
	\begin{equation} \label{eq:diag-wt-g-circ-f-abstract}
		\begin{tikzcd}
			{\wt{\InnHom{M}{N}[][A_1] \otimes_{1} \InnHom{L}{M}[][A_1]}} &&
			{\wt{\InnHom{L}{N}[][A_1]}} \\
			{\wt{\InnHom{M}{N}[][A_1]} \otimes_{2} \wt{\InnHom{L}{M}[][A_1]}} \\
			{\InnHom{\wt{M}}{\wt{N}}[][A_2] \otimes_{2} \InnHom{\wt{L}}{\wt{M}}[][A_2]} &&
			{\InnHom{\wt{L}}{\wt{N}}[][A_2]}
			\arrow["{\wt{c_{L,M,N}}}", from=1-1, to=1-3]
			\arrow["{\mu^{\mathcal{F}_{\Mod}}_{\InnHom{M}{N}[][A_1],\InnHom{L}{M}[][A_1]}}",
				from=2-1, to=1-1]
			\arrow["{\eta_{L,N}}"', from=1-3, to=3-3]
			\arrow["\wt{c}_{L,M,N}"', from=2-1, to=1-3]
			\arrow["{\eta_{M,N} \otimes_{2} \eta_{L,M}}"', from=2-1, to=3-1]
			\arrow["{c_{\wt{L},\wt{M},\wt{N}}}"', from=3-1, to=3-3]
		\end{tikzcd}
	\end{equation}
	where the maps $c_{L,M,N}$ (resp.\ $c_{\wt{L},\wt{M},\wt{N}}$) are the natural composition maps
	defined for $\mathcal{C}_1 \left( A_1 \right)$ (resp.\ $\mathcal{C}_2 \left( A_2 \right)$)
	using the tensor-hom adjunctions and the evaluation maps of
	$\mathcal{C}_1 \left( A_1 \right)$ (resp.\ $\mathcal{C}_2 \left( A_2 \right)$).
	Note that the commutativity of \eqref{eq:diag-wt-g-circ-f-abstract} is a completely
	abstract statement which doesn't involve talking about elements
	of $\InnHom{\cdot}{\cdot}[][A_1]$ or thinking about them as graded maps. When applied to
	our specific context, one obtains \cref{eq:wt-g-circ-f}.

	From this perspective, the commutativity of \eqref{eq:diag-wt-g-circ-f-abstract} can
	be interpreted as saying that $\mathcal{F}_{\Mod}$, together with the maps
	$\eta$, defines an \textit{enriched functor}
	\begin{equation*}
		\mathcal{F}_{\Mod} \colon
		\left( \mathcal{F}_{\Mod} \right)_{!} \left( \mathcal{C}_1 \left( A_1 \right) \right)
		\rightarrow
		\mathcal{C}_2 \left( A_2 \right)
	\end{equation*}
	between $\mathcal{C}_2 \left( A_2 \right)$-enriched categories. Here,
	$\mathcal{C}_2 \left( A_2 \right)$ is the self-enriched
	$\mathcal{C}_2 \left( A_2 \right)$-category and
	$\left( \mathcal{F}_{\Mod} \right)_{!} \left( \mathcal{C}_1 \left( A_1 \right) \right)$
	is the enriched $\mathcal{C}_2 \left( A_2 \right)$-category obtained by changing
	the base of enrichment of the self-enriched $\mathcal{C}_1 \left( A_1 \right)$-category
	$\mathcal{C}_1 \left( A_1 \right)$ from $\mathcal{C}_1 \left( A_1 \right)$ to
	$\mathcal{C}_2 \left( A_2 \right)$ using $\mathcal{F}_{\Mod}$
	(see \cite{EnrichedCategoryBaseChange}). As $\wt{c}_{L,M,N}$ is the composition in the category
	$\left( \mathcal{F}_{\Mod} \right)_{!} \left( \mathcal{C}_1 \left( A_1 \right) \right)$,
	diagram \eqref{eq:diag-wt-g-circ-f-abstract} simply depicts the condition that
	the enriched functor respects compositions.
	We will not pursue this point of view further, but for a reference that discusses
	such generalizations, see \cite[Appendix B]{Niles2025}.
\end{rem}

\subsubsection{Equivalence for Tensor Coalgebras} \label{sec:equiv-tensor-coalgebras}
Let $A_1 \in \CAlg[\mathcal{C}_1]$ be a $\sigma_1$-commutative algebra object
and let $\wt{A_1} = A_2$ be the corresponding $\sigma_2$-commutative algebra object discussed in
\cref{sec:equivalence-pairing-algebras}.
By abstract nonsense (see \cref{subsec:coalg-in-monoidal-cat}), the functor $\mathcal{F}_{\Mod}$
induces a strong symmetric monoidal equivalence
$\mathcal{F}_{\CoAlg} \colon \CoAlg[\mathcal{C}_1 \left( A_1 \right)] \rightarrow
	\CoAlg[\mathcal{C}_2 \left( A_2 \right)]$
converting coalgebra objects of $\mathcal{C}_1 \left( A_1 \right)$, i.e., graded $A_1$-coalgebras,
to coalgebra objects of $\mathcal{C}_2 \left( A_2 \right)$, i.e., graded $A_2$-coalgebras.
Given a graded $A_1$-coalgebra $C = C_1 \in \CoAlg[\mathcal{C}_1 \left( A_1 \right)]$, we will denote
the resulting coalgebra $\mathcal{F}_{\CoAlg} \left( C \right)$ by $\wt{C} = C_2$. The notation
is consistent with the previous section as $\wt{C}$ coincides as a graded $A_2$-module with
$\mathcal{F}_{\Mod} \left( C \right)$, which we also denoted by $\wt{C}$. The coproduct
structure on $\wt{C}$ comes from the coproduct structure on $C$, twisted by the tensor constraints.
Explicitly, if we use Sweedler's notation and write the coproduct
$\Delta_C \colon C \rightarrow C \otimes_1 C$ on $C$ as
$\Delta_C \left( c \right) = c_{(1)} \otimes_1 c_{(2)}$, then the coproduct
$\Delta_{\wt{C}} \colon \wt{C} \rightarrow \wt{C} \otimes_2 \wt{C}$ on
$\wt{C}$ is given by
\begin{equation*}
	\Delta_{\wt{C}} \left( \wt{c} \right) =
	(-1)^{\varphi \left( \degb{c_{(1)}}, \degb{c_{(2)}} \right)} \wt{c_{(1)}} \otimes_2 \wt{c_{(2)}}.
\end{equation*}
The construction $C \xmapsto{\mathcal{F}_{\CoAlg}} \wt{C}$ is compatible with coderivations in
the following sense: Given
a graded $A_1$-linear coderivation $\mu \colon C \rightharpoonup C$, i.e., we have
\begin{equation*}
	\Delta_C \circ \mu =
	\left( \mu \otimes_1 \id + \id \otimes_1 \mu \right) \circ \Delta_C,
\end{equation*}
the graded $A_2$-linear map $\wt{\mu} \colon \wt{C} \rightharpoonup \wt{C}$ given by
\cref{eq:action-F-mod-on-graded-morphisms} is a coderivation on the coalgebra $\wt{C}$, i.e.,
we have
\begin{equation*}
	\Delta_{\wt{C}} \circ \wt{\mu} =
	\left( \wt{\mu} \otimes_2 \id + \id \otimes_2 \wt{\mu} \right) \circ
	\Delta_{\wt{C}}.
\end{equation*}

Similarly to the discussion above, one can obtain induced totalization functors
\begin{equation*}
	\tot^{\braidop_i}_{\CoAlg} \colon \CoAlg[\mathcal{C}_i \left( A_i \right)] \rightarrow
	\CoAlg[\mathcal{D} \left( \tot^{\braidop_i}_{\CAlg} \left( A_i \right) \right)]
\end{equation*}
which endow the totalization of a $\ZZ^2$-graded $A_i$-coalgebra with the structure of
a $\ZZ$-graded $\tot^{\braidop_i}_{\CAlg} \left( A_i \right)$-coalgebra.
The functor $\tot^{\braidop_2}_{\CoAlg}$ does not twist the coproduct by a sign while
the functor $\tot^{\braidop_1}_{\CoAlg}$ does twist the coproduct by a sign coming
from the tensor constraints of $\tot^{\braidop_1}_{\Mod}$.
If we denote as before by $A \defeq \tot^{\braidop_1}_{\CAlg} \left( A_1 \right) =
	\tot^{\braidop_2}_{\CAlg} \left( A_2 \right)$
the totalization of $A_1$, or, equivalently, the totalization of $A_2$,
then we obtain a commutative diagram
\begin{equation} \label{eq:diag-f-tot-1-tot-2-A_i-coalgebras}
	\begin{tikzcd}
		{\CoAlg \left( \mathcal{C}_1 \left( A_1 \right) \right)} & &
		{\CoAlg \left( \mathcal{C}_2 \left( A_2 \right) \right)} \\
		\\
		\\
		&
		{\CoAlg \left( \mathcal{D} \left( A \right) \right)}
		\arrow["{\mathcal{F}_{\CoAlg}}",
			"(-1)^{\varphi \left( \degb{m}, \degb{n} \right)}"', from=1-1, to=1-3]
		\arrow["{\tot^{\braidop_1}_{\CoAlg}}"', "(-1)^{\varphi \left( \degb{m}, \degb{n} \right)}",
			from=1-1, to=4-2]
		\arrow["{\tot^{\braidop_2}_{\CoAlg}}", "(-1)^{0}"', from=1-3, to=4-2]
	\end{tikzcd}
\end{equation}
of strong symmetric monoidal functors.

Now let $M \in \mathcal{C}_1 \left( A_1 \right)$ be a graded $A_1$-module and
let $C_1 = \tens{M}[A_1] = \oplus_{k = 0}^{\infty} M^{\otimes_{1} k}$ be the tensor coalgebra on $M$, constructed in $\mathcal{C}_1 \left( A_1 \right)$. By applying the functor $\mathcal{F}_{\CoAlg}$,
we obtain a graded $A_2$-coalgebra $\wt{C_1} = \wt{\tens{M}[A_1]}$ whose coproduct is given
explicitly by
\begin{equation*}
	\begin{aligned}
		\Delta_{\wt{\tens{M}[A_1]}} \left( \wt{ m_1 \otimes_1 \dots \otimes_1 m_k} \right) =
		\sum_{i=0}^k & (-1)^{\varphi \left( \sum_{j=1}^i \degb{m_j}, \sum_{j=i+1}^k \degb{m_j} \right)}
		\\
		             & \qquad
		\wt{ m_1 \otimes_1 \dots \otimes_1 m_i} \otimes_2 \wt{m_{i+1} \otimes_1 \dots \otimes_1 m_k}.
	\end{aligned}
\end{equation*}

\begin{figure}[htbp]
	\centering
	\adjustbox{scale=0.85,center}{
		\begin{tikzcd}
			{\mathcal{C}_1 \left( A_1 \right)} &&&&&& {\mathcal{C}_2 \left( A_2 \right)} \\
			&& {\CoAlg \left( \mathcal{C}_1 \left( A_1 \right) \right)} && {\CoAlg \left( \mathcal{C}_2 \left( A_2 \right) \right)} \\
			&&& \circlearrowright \\
			&&& {\CoAlg \left( \mathcal{D} \left( A \right) \right)} \\
			\\
			&&& {\mathcal{D} \left( A \right)}
			\arrow["{\mathcal{F}_{\Mod}}", from=1-1, to=1-7]
			\arrow["{\tens{}[A_1]}"', from=1-1, to=2-3]
			\arrow["{\tot^{\braidop_1}_{\Mod}}"', curve={height=20pt}, from=1-1, to=6-4]
			\arrow["{\Theta_{\mathcal{F}}}", between={0.1}{0.9}, Rightarrow, from=1-7, to=2-3]
			\arrow["{\tens{}[A_2]}", from=1-7, to=2-5]
			\arrow["{\tot^{\braidop_2}_{\Mod}}", curve={height=-20pt}, from=1-7, to=6-4]
			\arrow["{\mathcal{F}_{\CoAlg}}"', from=2-3, to=2-5]
			\arrow["{\tot^{\braidop_1}_{\CoAlg}}", from=2-3, to=4-4]
			\arrow["{\tot^{\braidop_2}_{\CoAlg}}"', from=2-5, to=4-4]
			\arrow["{\Theta_{\tot^{\braidop_1}}}"{pos=0.5}, between={0.1}{0.9}, Rightarrow, from=6-4, to=2-3]
			\arrow["{\Theta_{\tot^{\braidop_2}}}"'{pos=0.5}, between={0.1}{0.9}, Rightarrow, from=6-4, to=2-5]
			\arrow["{\tens{}[A]}"', from=6-4, to=4-4]
		\end{tikzcd}
	}
	\caption{Functors and Natural Transformations Relating the Tensor Coalgebras for Different Pairings.}
	\label{fig:Theta-func-nat}
\end{figure}
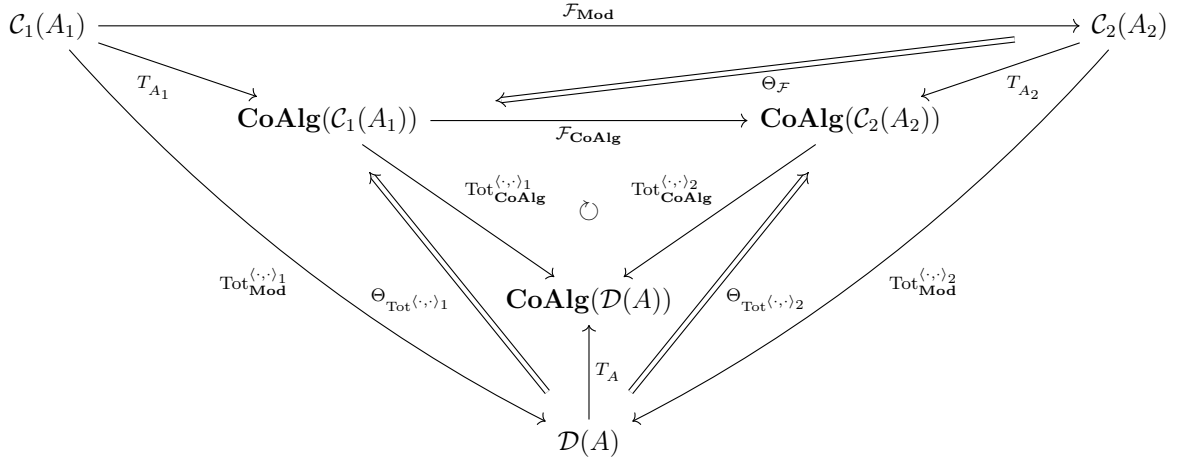

The unit constraint \eqref{eq:unit-constraint-F-mod} and the (iterated versions of the)
tensor constraints \eqref{eq:tensor-constraints-F-mod-signature} of $\mathcal{F}_{\Mod}$
induce a natural \textit{isomorphism}
\begin{equation*}
	\Theta_{\mathcal{F}} \colon \tens{\wt{M}}[A_2] \rightarrow \wt{\tens{M}[A_1]}
\end{equation*}
of graded $A_2$-modules, coming from the series of identifications
\begin{equation*}
	\wt{\tens{M}[A_1]} = \wt{\oplus_{k=0}^{\infty} M^{\otimes_{1} k}} \cong
	\oplus_{k=0}^{\infty} \wt{M^{\otimes_{1} k}} \cong
	\oplus_{k=0}^{\infty} \wt{M}^{\otimes_{2} k} = \tens{\wt{M}}[A_2],
\end{equation*}
given explicitly by
\begin{equation} \label{eq:Theta-explicit-formula}
	\Theta_{\mathcal{F}} \left( \wt{m_1} \otimes_2 \dots \otimes_2 \wt{m_k} \right) =
	(-1)^{\sum_{j=1}^{k-1} \varphi \left( \degb{m_1} + \dots + \degb{m_j}, \degb{m_{j+1}} \right)}
	\wt{ m_1 \otimes_1 \dots \otimes_1 m_k}.
\end{equation}
While a priori only an isomorphism of graded $A_2$-modules, in fact the isomorphism $\Theta$ is
an isomorphism of graded $A_2$-\textit{coalgebras}, identifying the graded $A_2$-coalgebra
$\wt{\tens{M}[A_1]}$ with the tensor coalgebra $\tens{\wt{M}}[A_2]$ as coalgebras.\footnote{In fact, since
	$\mathcal{F}_{\Mod}$ is \textit{strong} monoidal, $\wt{\tens{M}[A_1]}$ has also a natural structure
	of a graded $A_2$-algebra coming from the product of the tensor algebra $\tens{M}[A_1]$ and
	the tensor constraints. The isomorphism $\Theta_{\mathcal{F}}$ is an isomorphism between \textit{both} the coalgebra
	and the algebra structures on both sides.}
Similarly, for $i \in \Set{1,2}$, the unit and tensor constraints of $\tot^{\braidop_i}_{\Mod}$ induce natural
\textit{isomorphisms}
\begin{equation*}
	\Theta_{i} \colon \tens{\tot^{\braidop_i}_{\Mod} \left( M_i \right)}[A]
	\rightarrow
	\tot^{\braidop_i}_{\CoAlg} \left( \tens{M_i}[A_i] \right)
\end{equation*}
of graded $A$-coalgebras, determined by
\begin{equation} \label{eq:Theta-1-explicit-formula}
	\begin{aligned}
		\Theta_{1} \left( s_{i_1} \left( m_{i_1}^{j_1} \right) \otimes_A \dots \otimes_A s_{i_k} \left( m_{i_k}^{j_k} \right) \right) ={} &
		(-1)^{\sum_{r=1}^{k-1} \varphi \left( \left( i_1, j_1 \right) + \dots + \left( i_r, j_r \right), \left( i_{r+1}, j_{r+1} \right) \right)}
		\\
		                                                                                                                                  & \qquad
		\s_{i_1 + \dots + i_k} \left( m_{i_1}^{j_1} \otimes_1 \dots \otimes_1 m_{i_k}^{j_k} \right),
	\end{aligned}
\end{equation}
and
\begin{equation*}
	\Theta_{2} \left( m_{i_1}^{j_1} \otimes_A \dots \otimes_A m_{i_k}^{j_k} \right) =
	m_{i_1}^{j_1} \otimes_2 \dots \otimes_2 m_{i_k}^{j_k}.
\end{equation*}
The relation between the functors $\mathcal{F}, \tot^{\braidop_1}, \tot^{\braidop_2}$ and
the natural transformations $\Theta_{\mathcal{F}}, \Theta_1, \Theta_2$ is summarized in
\cref{fig:Theta-func-nat}, in which both the inner and outer triangles strictly commute
while the squares commute only up to the natural isomorphisms $\Theta_{\bullet}$.

Given an $A_1$-linear coderivation $\mu \colon \tens{M}[A_1] \rightharpoonup \tens{M}[A_1]$
on the tensor coalgebra $\tens{M}[A_1]$, we will denote by
\begin{equation} \label{eq:def-overline-coder}
	\overline{\mu} \defeq \Theta_{\mathcal{F}}^{-1} \circ \wt{\mu} \circ \Theta_{\mathcal{F}} \colon \tens{\wt{M}}[A_2]
	\rightharpoonup \tens{\wt{M}}[A_2]
\end{equation}
the corresponding $A_2$-linear coderivation on the tensor coalgebra $\tens{\wt{M}}[A_2]$.
This gives us a construction sending the pair $\left( \tens{M}[A_1], \mu \right)$
to the pair $\left( \tens{\wt{M}}[A_2], \overline{\mu} \right)$. In terms of
the corestrictions, $\corest{\overline{\mu}} \colon \tens{\wt{M}}[A_2] \rightharpoonup \wt{M}$
is related to $\corest{\mu} \colon \tens{M}[A_1] \rightharpoonup M$ via the following
formula:
\begin{equation} \label{eq:corest-overline-mu}
	\begin{aligned}
		\corest{\overline{\mu}}
		\left( \wt{m_1} \otimes_2 \dots \otimes_2 \wt{m_k} \right)
		={} &
		(-1)^{\sum_{j=1}^{k-1} \varphi \left( \degb{m_1} + \dots + \degb{m_j}, \degb{m_{j+1}} \right)
			+ \varphi \left( \degb{\mu}, \sum_{j=1}^k \degb{m_j} \right)}
		\\
		    & \qquad
		\wt{\corest{\mu} \left( m_1 \otimes_1 \dots \otimes_1 m_k \right)}.
	\end{aligned}
\end{equation}

\subsubsection{Equivalence \texorpdfstring{of $\ndf{\cdot, \braidop_1}[][]$ and $\ndf{\cdot, \braidop_2}[][]$}{for Codifferential Forms}}
The discussion in the previous subsections was quite general, but now we can apply it
to the specific situation we are interested in, described in \cref{sec:dependence-ndf-braidop}.

Let $R$ be a $\ZZ$-graded, graded-commutative, $\mathbbm{k}$-algebra. We think of $R$
as a $\ZZ^2$-graded $\mathbbm{k}$-algebra by placing $R^{*}$
in bidegree $(0,*)$. The resulting $\ZZ^2$-graded algebra is graded-commutative with respect
to both $\sigma_1$ and $\sigma_2$. Since $R$ is concentrated in bidegree $(0,*)$,
the product on $\mathcal{F}_{\Alg} \left( R \right) = \wt{R}$, given by \cref{eq:product-wt-A},
is $\wt{r_1} \cdot \wt{r_2} = \wt{r_1 \cdot r_2}$. Thus, we can identify $\wt{R}$
with $R$ via $r \leftrightarrow \wt{r}$ with no signs and drop the tilde notation.
Similarly, we have natural identifications
$\tot^{\braidop_2}_{\CAlg} \left( R \right) = \tot^{\braidop_1}_{\CAlg} \left( R \right) = R$ with no signs.
In this situation, the commutative diagram \eqref{eq:diag-f-tot-1-tot-2-A_2-modules} takes the form
\begin{equation*} 
	\begin{tikzcd}
		{\mathcal{C}_1 \left( R \right)} &&&&&& {\mathcal{C}_2 \left( R \right)} \\
		\\
		\\
		&&& {\mathcal{D} \left( R \right).}
		\arrow["{\mathcal{F}_{\Mod}}",
			"(-1)^{\varphi \left( \degb{m}, \degb{n} \right)}"', from=1-1, to=1-7]
		\arrow["{\tot^{\braidop_1}_{\Mod}}"', "(-1)^{\varphi \left( \degb{m}, \degb{n} \right)}",
			from=1-1, to=4-4]
		\arrow["{\tot^{\braidop_2}_{\Mod}}", "(-1)^0"', from=1-7, to=4-4]
	\end{tikzcd}
\end{equation*}

Let $V$ be a $\ZZ$-graded $R$-module. We think of $V$ as a $\ZZ^2$-graded $\mathbbm{k}$-module
by placing $V^{*}$ in bidegree $(0,*)$, so that $V$ becomes an object of
$\mathcal{C}_i \left( R \right)$ for both $i = 1, 2$. Since $V$ is concentrated in
bidegree $(0,*)$, we have $\mathcal{F}_{\Mod} \left( V \right) = \wt{V} = V$, without
any sign twisting the action. Now consider the modules $\ul{V}_i = V[(-1,0)]_i$,
which are the horizontal shifts of $V$ concentrated in line degree one,
as objects of $\mathcal{C}_i \left( R \right)$. The modules $\ul{V}_i$ are
identical as $\ZZ^2$-graded $\mathbbm{k}$-modules, but have different $R$-actions given by
\begin{equation} \label{eq:R-action-ul-v-appendix}
	r \ul{v}_i = (-1)^{\braid{(1,0)}{(0,\degb{r})}_i} \ul{r v}_i =
	\begin{cases}
		\ul{rv}_{1}                 & i = 1, \\
		(-1)^{\degb{r}} \ul{rv}_{2} & i = 2.
	\end{cases}
\end{equation}
We can identify $\mathcal{F}_{\Mod} \left( \ul{V}_1 \right) = \wt{\ul{V}_1}$
with $\ul{V}_2$ via $\wt{\ul{v}_1} \leftrightarrow \ul{v}_2$ without any signs,
since we have
\begin{equation*}
	r \cdot \wt{\ul{v}_1} = \wt{r} \cdot \wt{\ul{v}_1}
	\stackrel{\eqref{eq:action-wt-A-2-on-wt-M}}{=}
	(-1)^{\varphi \left( \left( 0, \degb{r} \right), \left( 1, \degb{v} \right) \right)}
	\wt{ \left( r \cdot \ul{v}_1 \right) }
	\stackrel{\eqref{eq:varphi-formula}}{=}
	(-1)^{\degb{r}}
	\wt{ \left( r \cdot \ul{v}_1 \right) }
	\stackrel{\eqref{eq:R-action-ul-v-appendix}}{=}
	(-1)^{\degb{r}} \wt{ \ul{r \cdot v}_1 },
\end{equation*}
which is the same formula \eqref{eq:R-action-ul-v-appendix} as for the action of $r$ on $\ul{v}_2$.

Let $M_i = V \oplus \ul{V}_i \in \mathcal{C}_i \left( R \right)$ and construct
the tensor coalgebra
\begin{equation*}
	\tens{M_i}[R] = \tens{V \oplus \ul{V}_i}[R, \otimes_i] = \ndf{V, \braidop_i}[][]
\end{equation*}
in $\mathcal{C}_i \left( R \right)$ for $i = 1, 2$.
By the previous paragraph, we have a natural identification
$\mathcal{F}_{\Mod} \left( M_1 \right) = \wt{M_1} \cong M_2$, which doesn't involve any signs.
By the definition of the totalization functors, we have
$\tot^{\braidop_i}_{\Mod} \left( M_i \right) = V \oplus V[1]$ for both $i = 1, 2$.

We now apply the results of \cref{sec:equiv-tensor-coalgebras} to our situation described
above, in which $A_1 = A_2 = A = R$ is concentrated in line degree zero, and the modules
$M_i$ are concentrated in line degrees zero and one.
The tensor constraints of the strong
monoidal functor $\tot^{\braidop_2}_{\Mod}$
induce a natural isomorphism
\begin{equation*}
	\tens{V \oplus V[1]}[R] = \tens{\tot^{\braidop_2}_{\Mod} \left( M_2 \right)}[R]
	\xrightarrow{\Theta_2}
	\tot^{\braidop_2}_{\CoAlg} \left( \tens{M_2}[R] \right) = \totc{\ndf{V, \braidop_2}}[][][\braidop_2].
\end{equation*}
Since the tensor constraints are trivial, $\Theta_2$ involves no signs, and we see
that $\Theta_2 = \Phi_2$.

Similarly, the tensor constraints of the strong
monoidal functor $\tot^{\braidop_1}_{\Mod}$
induce a natural isomorphism
\begin{equation*}
	\tens{V \oplus V[1]}[R] = \tens{\tot^{\braidop_1}_{\Mod} \left( M_1 \right)}[R]
	\xrightarrow{\Theta_1}
	\tot^{\braidop_1}_{\CoAlg} \left( \tens{M_1}[R] \right) = \totc{\ndf{V, \braidop_1}}[][][\braidop_1].
\end{equation*}
The isomorphism $\Theta_1$, given by \cref{eq:Theta-1-explicit-formula}, involves a sign coming
from the non-trivial tensor constraints of $\tot^{\braidop_1}_{\Mod}$.
Since $M_1$ is concentrated in line degrees zero and one,
$\tens{\tot^{\braidop_1}_{\Mod} \left( M_1 \right)}[R]$ is generated by tensors of the
form
$x = l^0 \otimes_R \s v_1 \otimes_R \dots \otimes_R \s v_k \otimes_R l^k$,
where $k \geq 0$, $l^i \in \tens{V}[R]$ for $0 \leq i \leq k$, and
$v_i \in V$ for $i = 1, \dots, k$. In this case, using the explicit formula
\eqref{eq:varphi-formula} of $\varphi$, we see that the sign is given by
\begin{equation*}
	\varepsilon = \degb{l^0} \cdot k + \left( \degb{v_1} + \degb{l^1} \right) \cdot (k-1) + \dots + \left( \degb{v_{k-1}} + \degb{l^{k-1}} \right),
\end{equation*}
which is the same as the sign \eqref{eq:sign-factor-Phi-1} appearing in the formula for $\Phi_1$.
Hence, $\Theta_1$ is precisely $\Phi_1$.

Finally, the tensor constraints of the strong monoidal functor $\mathcal{F}_{\Mod}$ induce a natural
coalgebra isomorphism
\begin{equation*}
	\ndf{V, \braidop_2}[][] = \tens{M_2}[R] = \tens{\wt{M_1}}[R]
	\xrightarrow{\Theta_{\mathcal{F}}}
	\wt{\tens{M_1}[R]} = \wt{\ndf{V, \braidop_1}[][]},
\end{equation*}
given explicitly by
\begin{equation*}
	\Theta_{\mathcal{F}} \left(
	l^0 \otimes_{2} \ul{v_1} \otimes_{2} \dots \otimes_{2} \ul{v_k} \otimes_2 l^k
	\right) =
	(-1)^{\varepsilon}
	\wt{l^0 \otimes_{1} \ul{v_1} \otimes_{1} \dots \otimes_{1} \ul{v_k} \otimes_1 l^k},
\end{equation*}
where $\varepsilon$ is the same sign as for $\Theta_1 = \Phi_1$.
By the commutativity of the inner triangle of \cref{fig:Theta-func-nat}, we have
$\tot^{\braidop_2} \left( \Theta_{\mathcal{F}} \right) = \Phi$,
so that
$\Theta_{\mathcal{F}} \colon \ndf{V, \braidop_2}[][] \rightarrow
	\wt{\ndf{V, \braidop_1}[][]}$
is the lift of
$\Phi \colon \totc{\ndf{V, \braidop_2}}[][][\braidop_2] \rightarrow
	\totc{\ndf{V, \braidop_1}}[][][\braidop_1]$,
to the level of $\ZZ^2$-graded coalgebras, where both the $R$-module and coalgebra structures
on the codomain are twisted (see \cref{fig:Theta-trans}).
\begin{figure}[htbp]
	\centering
	\adjustbox{scale=0.93,center}{
		\begin{tikzcd}
			{\totc{\ndf{V, \braidop_2}[][]}[][][\braidop_2]}
			&& {\tot^{\braidop_2} ( \wt{ \ndf{V, \braidop_1}[][]} ) = \totc{\ndf{V, \braidop_1}[][]}[][][\braidop_1]} \\
			& {\tens{V \oplus V[1]}[R]}
			\arrow["{\Phi = \tot^{\braidop_2} \left( \Theta_{\mathcal{F}} \right)}", from=1-1, to=1-3]
			\arrow["{\Phi_2 = \Theta_2}"', from=2-2, to=1-1, end anchor=south]
			\arrow["{\Phi_1 = \Theta_1}"', from=2-2, to=1-3, end anchor=south]
		\end{tikzcd}
	}
	\caption{The Relation Between $\Phi_{\bullet}$ and $\Theta_{\bullet}$.}
	\label{fig:Theta-trans}
\end{figure}

Let us show that the process $\eta \mapsto \overline{\eta}$ identifies the
operators on codifferential forms constructed using $\braidop_i$ for $i = 1, 2$:
\begin{lm} \label{lm:overline-calc-diff-pairings}
	Let $\mu$ be an $R$-linear coderivation on $\tens{V}[R]$. Then we have the identities:
	\begin{align}
		\overline{\qdr^{\braidop_1}}       & = \qdr^{\braidop_2},
		\label{eq:overline-qdr-1-qdr-2}
		\\
		\overline{\lie{\mu}^{\braidop_1}}  & = \lie{\mu}^{\braidop_2},
		\label{eq:overline-lie-mu-1-lie-mu-2}
		\\
		\overline{\cont{\mu}^{\braidop_1}} & = \cont{\mu}^{\braidop_2}.
		\label{eq:overline-cont-mu-1-cont-mu-2}
	\end{align}
\end{lm}
\begin{proof}
	To verify the identities, it is enough to compare the corestriction of both sides.
	First note
	that by \cref{eq:corest-overline-mu}, the corestriction $\corest{\overline{\qdr^{\braidop_1}}}$ vanishes
	on elements of weight greater than one (since this is true for $\qdr^{\braidop_1}$), while
	for elements of weight one, we have
	\begin{align*}
		\overline{\qdr^{\braidop_1}}_1 \left( v \right)
		\stackrel{\eqref{eq:corest-overline-mu}}{=} &
		(-1)^{\varphi \left( (-1,0), (0, \degb{v}) \right)}
		\qdr^{\braidop_1}_{1} \left( v \right)
		\stackrel{\eqref{eq:varphi-formula}}{=}
		\qdr^{\braidop_1}_{1} \left( v \right) = 0 = \qdr^{\braidop_2}_{1} \left( v \right), \\
		\overline{\qdr^{\braidop_1}}_1 \left( \ul{v} \right)
		\stackrel{\eqref{eq:corest-overline-mu}}{=} &
		(-1)^{\varphi \left( (-1,0), (1, \degb{v}) \right)}
		\qdr^{\braidop_1}_{1} \left( \ul{v} \right)
		\stackrel{\eqref{eq:varphi-formula}}{=}
		\qdr^{\braidop_1}_{1} \left( \ul{v} \right) = v = \qdr^{\braidop_2}_{1} \left( \ul{v} \right),
	\end{align*}
	which shows \cref{eq:overline-qdr-1-qdr-2}.

	Next, note that the corestriction of $\lie{\mu}^{\braidop_1}$ on elements of line
	degree zero is given by
	\begin{equation*}
		\begin{aligned}
			\left( \overline{\lie{\mu}^{\braidop_1}} \right)_{k}
			\left( v_1 \otimes_2 \dots \otimes_2 v_k \right)
			\eqwithref[eq:corest-overline-mu] &
			(-1)^{\sum_{j=1}^{k-1} \varphi \left( (0, \sum_{i=1}^j \degb{v_i}), (0, \degb{v_{j+1}}) \right)
				+ \varphi \left( (0, \degb{\mu}), (0, \sum_{j=1}^k \degb{v_j}) \right)}
			\\
			                                  & \qquad
			\left( \lie{\mu}^{\braidop_1} \right)_k \left( v_1 \otimes_1 \dots \otimes_1 v_k \right)
			\\
			\eqwithref[eq:varphi-formula]     &
			\left( \lie{\mu}^{\braidop_1} \right)_k \left( v_1 \otimes_1 \dots \otimes_1 v_k \right)
			\\
			\eqwithref[eq:lie-der-corest]     &
			\mu_k \left( v_1 \otimes_R \dots \otimes_R v_k \right),
		\end{aligned}
	\end{equation*}
	so that $\overline{\lie{\mu}^{\braidop_1}}$ extends the action of $\mu$ on elements
	of line degree zero. We also have
	\begin{equation*}
		\begin{aligned}
			\left[ \qdr^{\braidop_2}, \overline{\lie{\mu}^{\braidop_1}} \right]_{2}
			\stackrel{\eqref{eq:overline-qdr-1-qdr-2}}{=} &
			\left[ \overline{\qdr^{\braidop_1}}, \overline{\lie{\mu}^{\braidop_1}} \right]_{2}
			\stackrel{\eqref{eq:def-overline-coder}}{=}
			\Theta_{\mathcal{F}}^{-1} \circ \left[ \wt{\qdr^{\braidop_1}}, \wt{\lie{\mu}^{\braidop_1}} \right]_{2} \circ \Theta_{\mathcal{F}}
			\\
			\stackrel{\eqref{eq:comp-wt-graded-comm}}{=}  &
			(-1)^{\varphi \left( (-1, 0), (0, \degb{\mu}) \right)}
			\Theta_{\mathcal{F}}^{-1} \circ
			\wt{\left[ \qdr^{\braidop_1}, \lie{\mu}^{\braidop_1} \right]_{1}} \circ \Theta_{\mathcal{F}}
			\stackrel{\eqref{eq:qliecomm}}{=} 0,
		\end{aligned}
	\end{equation*}
	so by the uniqueness part of \cref{lm:existence-uniqueness-lie}, we get
	\cref{eq:overline-lie-mu-1-lie-mu-2}.

	Finally, we have
	\begin{equation*}
		\begin{aligned}
			\left[ \qdr^{\braidop_2}, \overline{\cont{\mu}^{\braidop_1}} \right]_{2}
			\eqwithref[eq:overline-qdr-1-qdr-2]                                     &
			\left[ \overline{\qdr^{\braidop_1}}, \overline{\cont{\mu}^{\braidop_1}} \right]_{2}
			\stackrel{\eqref{eq:def-overline-coder}}{=}
			\Theta_{\mathcal{F}}^{-1} \circ \left[ \wt{\qdr^{\braidop_1}}, \wt{\cont{\mu}^{\braidop_1}} \right]_{2} \circ \Theta_{\mathcal{F}}
			\\
			\eqwithref[eq:comp-wt-graded-comm]                                      &
			(-1)^{\varphi \left( (-1, 0), (1, \degb{\mu}) \right)}
			\Theta_{\mathcal{F}}^{-1} \circ
			\wt{\left[ \qdr^{\braidop_1}, \cont{\mu}^{\braidop_1} \right]_{1}} \circ \Theta_{\mathcal{F}}
			\\
			\eqwithref[eq:varphi-formula][eq:lie-derivative-commutator-contraction] &
			\Theta_{\mathcal{F}}^{-1} \circ \wt{\lie{\mu}^{\braidop_1}} \circ \Theta_{\mathcal{F}}
			\stackrel{\eqref{eq:def-overline-coder}}{=}
			\overline{\lie{\mu}^{\braidop_1}}
			\stackrel{\eqref{eq:overline-lie-mu-1-lie-mu-2}}{=}
			\lie{\mu}^{\braidop_2},
		\end{aligned}
	\end{equation*}
	so by the uniqueness part of \cref{lm:existence-uniqueness-contraction}, we get
	\cref{eq:overline-cont-mu-1-cont-mu-2}.
\end{proof}

Given two coderivations $\eta_i$ on $\ndf{V, \braidop_i}[][]$ for $i = 1, 2$,
the relation $\overline{\eta_1} = \eta_2$ holds if and only if
$\wt{\eta_1} \circ \Theta_{\mathcal{F}} = \Theta_{\mathcal{F}} \circ \eta_2$. In this case,
applying $\tot^{\braidop_2}$ to both sides, we get the relation
\begin{equation*}
	\begin{aligned}
		\tot^{\braidop_1} \left( \eta_1 \right) \circ \Phi & =
		\tot^{\braidop_2} \left( \wt{\eta_1} \right) \circ \tot^{\braidop_2} \left( \Theta_{\mathcal{F}} \right)
		\\
		                                                   & =
		\tot^{\braidop_2} \left( \Theta_{\mathcal{F}} \right) \circ \tot^{\braidop_2}  \left( \eta_2 \right) =
		\Phi \circ \tot^{\braidop_2}  \left( \eta_2 \right).
	\end{aligned}
\end{equation*}
Thus, applying \cref{lm:overline-calc-diff-pairings},
\cref{eq:overline-qdr-1-qdr-2,eq:overline-lie-mu-1-lie-mu-2,eq:overline-cont-mu-1-cont-mu-2}
immediately imply
\cref{eq:Phi-qdr-rel,eq:Phi-lie-mu-rel,eq:Phi-cont-mu-rel} from
\cref{sec:dependence-ndf-braidop}.

Similarly, one can show that
$\widetilde{\indmap{f}^{\braidop_1}} \circ \Theta_{\mathcal{F}} =
	\Theta_{\mathcal{F}} \circ \indmap{f}^{\braidop_2}$ for a
morphism $f$, which implies \cref{eq:Phi-indmap-rel}.

We note that the construction $\mu \mapsto \widetilde{\mu}$ of \eqref{eq:action-F-mod-on-graded-morphisms}
works also for module derivations over algebra derivations, which amounts
to extending the equivalences of
\crefrange{sec:braidop-eq-graded-k-modules}{sec:braidop-eq-graded-modules-over-graded-algebras}
to objects equipped with a pre-differential (see also \cref{rem:extension-working-over-pre-cdga}).
In this case the map $\mu \mapsto \overline{\mu}$ of \eqref{eq:def-overline-coder} will
respect generalized coderivations and
\cref{eq:overline-lie-mu-1-lie-mu-2} will also hold when $\mu$ is a generalized coderivation.
Since all our functors are \textit{symmetric} monoidal, the maps $\Theta_{\bullet}$
interleave the rotation operators on the domain and codomain, from which one
can deduce \cref{eq:Phi-clie-rel,eq:Phi-ccont-rel,eq:Phi-cindmap-rel} (the cyclic versions of
\cref{eq:Phi-lie-mu-rel,eq:Phi-cont-mu-rel,eq:Phi-indmap-rel}). We invite the interested reader to
fill in the details.

\subsection{Total Inner Product Components for Both Parity Forms}
\label{sec:tot-inner-product-components-both-parity-forms}

We discuss briefly how to define the components of a total inner product
$\phi^{\braidop_2} \colon \totcompe{\mathcal{A}, \braidop_2}[2][] \rightarrow \mathcal{R}[4-n]$,
where $\totcompe{\mathcal{A}, \braidop_2}[2][]$ is the extension
of the truncated total complex
\begin{equation*}
	\totcomp{\mathcal{A}, \braidop_2}[2][] = \totc{\ncdfr{A, \braidop_2}, -\qdr^{\braidop_2}, \clie{\mu}^{\braidop_2}}[][\geq 2][\braidop_2]
\end{equation*}
constructed using the parity form $\braidop_2$.

Let us set
\begin{equation*}
	\totcomp{\mathcal{A}, \braidop_1}[2][] = \totc{\ncdf{A, \braidop_1}, \qdr^{\braidop_1},
		\clie{\mu}^{\braidop_1}}[][\geq 2][\braidop_1].
\end{equation*}
As discussed in \cref{sec:description-using-braid-op-1}, the isomorphism
\begin{equation*}
	\Psi \colon \totc{\ncdf{A, \braidop_2}, -\qdr^{\braidop_2}, \clie{\mu}^{\braidop_2}}[][][\braidop_2]
	\rightarrow
	\totc{\ncdf{A, \braidop_1}, \qdr^{\braidop_1}, \clie{\mu}^{\braidop_1}}[][][\braidop_1],
\end{equation*}
of \eqref{eq:Psi-explicit-formula} induces by truncation an isomorphism
$\Psi_{\geq 2} \colon \totcomp{\mathcal{A}, \braidop_2}[2][] \rightarrow \totcomp{\mathcal{A}, \braidop_1}[2][]$,
which extends to an isomorphism
$\Psi_{\geq 2}^{+} \colon \totcompe{\mathcal{A}, \braidop_2}[2][] \rightarrow \totcompe{\mathcal{A}, \braidop_1}[2][]$
between the extended total complexes.

Then we have a bijection between total inner products
$\phi^{\braidop_2} \colon \totcompe{\mathcal{A}, \braidop_2}[2][] \rightarrow \mathcal{R}[4-n]$
and total inner products
$\phi^{\braidop_1} \colon \totcompe{\mathcal{A}, \braidop_1}[2][] \rightarrow \mathcal{R}[4-n]$,
given by
\begin{equation*}
	\phi^{\braidop_1} \mapsto \phi^{\braidop_2} = \phi^{\braidop_1} \circ \Psi_{\geq 2}^{+}, \qquad
	\phi^{\braidop_2} \mapsto \phi^{\braidop_1} = \phi^{\braidop_2} \circ \left( \Psi_{\geq 2}^{+} \right)^{-1}.
\end{equation*}
Thus, given a total inner product $\phi^{\braidop_2}$, we can define the components of
$\phi^{\braidop_2}$ by converting it to $\phi^{\braidop_1} = \phi^{\braidop_2} \circ \left( \Psi_{\geq 2}^{+} \right)^{-1}$
and taking its components as in \cref{sec:total-inner-product-explicit-relations}.

Let us recall how the components of $\phi = \phi^{\braidop_1}$ are defined.
When working with $\braidop_1$, we have the property that the inclusion
${\ncdf{A, \braidop_1}[k][]}[k] \hookrightarrow \totcomp{A, \braidop_1}[2][]$
as the $k$-th column is well-defined, and the vertical differential on
$\totcomp{A, \braidop_1}[2][]$, restricted to the $k$-th column, coincides with
${\clie{\mu}^{\braidop_1}}[k]$.
The components
$\phi_k \colon \ncdf{A, \braidop_1}[k][] \rightharpoonup R$ of $\phi$ are defined
via the composition \eqref{eq:components-phi_k}, giving us the expression
\begin{multline*}
	\phi_k \left(
	\ul{a_1} \otimes_1 b_1^1 \otimes_1 \dots \otimes_1 b_1^{r_1} \otimes_1 \dots \otimes_1
	\ul{a_k} \otimes_1 b_k^1 \otimes_1 \dots \otimes_1 b_k^{r_k}
	\right) =
	\\
	\sigma_{4-n} \phi \left( \s_k \left(
	\ul{a_1} \otimes_1 b_1^1 \otimes_1 \dots \otimes_1 b_1^{r_1} \otimes_1 \dots \otimes_1
	\ul{a_k} \otimes_1 b_k^1 \otimes_1 \dots \otimes_1 b_k^{r_k} \right)
	\right),
\end{multline*}
which does not involve any signs. Here, $\sigma_{4-n} \colon R[4-n] \rightharpoonup R$ is
the natural desuspension map.

We also have the property that the natural maps
\begin{equation*}
	A \otimes_R A^{\otimes_R r_1} \otimes_R \dots \otimes_R A \otimes_R A^{\otimes_R r_k}
	\rightarrow \ndf{A, \braidop_1}[k][]
\end{equation*}
given by
\begin{multline*}
	a_1 \otimes_R b_1^1 \otimes_R \dots \otimes_R b_1^{r_1} \otimes_R \dots \otimes_R
	a_k \otimes_R b_k^1 \otimes_R \dots \otimes_R b_k^{r_k} \mapsto
	\\
	\ul{a_1} \otimes_1 b_1^1 \otimes_1 \dots \otimes_1 b_1^{r_1} \otimes_1 \dots \otimes_1
	\ul{a_k} \otimes_1 b_k^1 \otimes_1 \dots \otimes_1 b_k^{r_k},
\end{multline*}
without any sign, are well-defined. The components
\begin{equation*}
	\phi_k^{r_1,\dots,r_k} \colon
	A \otimes_R A^{\otimes_R r_1} \otimes_R \dots \otimes_R A \otimes_R A^{\otimes_R r_k}
	\rightharpoonup R
\end{equation*}
of $\phi_k$ are defined via the composition \eqref{eq:components-phi_k-r_1-r_k},
so that
\begin{multline*}
	\phi_k^{r_1, \dots, r_k} \left(
	a_1 \otimes_R b_1^1 \otimes_R \dots \otimes_R b_1^{r_1} \otimes_R \dots \otimes_R
	a_k \otimes_R b_k^1 \otimes_R \dots \otimes_R b_k^{r_k}
	\right) =
	\\
	\phi_k \left(
	\ul{a_1} \otimes_1 b_1^1 \otimes_1 \dots \otimes_1 b_1^{r_1} \otimes_1 \dots \otimes_1
	\ul{a_k} \otimes_1 b_k^1 \otimes_1 \dots \otimes_1 b_k^{r_k}
	\right),
\end{multline*}
again without any signs.

Now assume we want to define the components of $\phi^{\braidop_2}$ directly,
without passing through $\phi^{\braidop_1}$. Given a $\ZZ$-graded $R$-module $M$
and $k \in \ZZ$, denote by $M \left< k \right>$ the shifted module in which
$M \left< k \right>^i = M^{i + k}$, but in which the $R$-action is not twisted.
Given an element $m \in M$, the corresponding (same) element of $M \left< k \right>$
is conveniently denoted by $m \s_k$, and we have
$r \left( m \s_k \right) = \left( rm \right) \s_k$. Note that the degree $-k$
``identity'' map $m \mapsto m \s_k$ is not $R$-linear of degree $-k$ because the $R$
action is not twisted. Instead, the suspension map $M \rightharpoonup M \left < k \right>$
is given by $m \mapsto (-1)^{k \cdot \degb{m}} m \s_k$, where the sign is introduced
to make it $R$-linear of degree $-k$. When $M$ is equipped with a (pre)-differential
$\mu$, we endow $M \left< k \right>$ with the (pre)-differential $\mu \left< k \right>$,
acting via $\mu \left< k \right> \left( m \s_k \right) = \mu \left( m \right) \s_k$,
with no sign twisting, so that the suspension map $M \rightharpoonup M \left< k \right>$
becomes a degree $-k$ chain map.

When working with $\braidop_2$, the inclusion
${\ncdf{A, \braidop_2}[k][]} \left< k \right> \hookrightarrow \totcomp{A, \braidop_2}[2][]$
as the $k$-th column, given by $x_k \s_k \mapsto x_k$, is well-defined, and the vertical differential on
$\totcomp{A, \braidop_2}[2][]$, restricted to the $k$-th column, coincides with
$\clie{\mu}^{\braidop_2} \left< k \right>$, since we don't twist the $R$-action
on the columns nor the vertical differential when using $\tot^{\braidop_2}$.
Thus, to make the components $\phi_k \colon \ncdf{A, \braidop_2}[k][] \rightharpoonup R$
of $\phi = \phi^{\braidop_2}$ well-defined
and $R$-linear of the correct degree, it is natural to define them via
\begin{multline*}
	\phi_k \left(
	\underbrace{
		\ul{a_1} \otimes_2 b_1^1 \otimes_2 \dots \otimes_2 b_1^{r_1} \otimes_2 \dots \otimes_2
		\ul{a_k} \otimes_2 b_k^1 \otimes_2 \dots \otimes_2 b_k^{r_k}
	}_{x_k}
	\right) =
	\\
	(-1)^{k \cdot \degb{x_k}}
	\sigma_{4-n} \phi \left(
	\ul{a_1} \otimes_2 b_1^1 \otimes_2 \dots \otimes_2 b_1^{r_1} \otimes_2 \dots \otimes_2
	\ul{a_k} \otimes_2 b_k^1 \otimes_2 \dots \otimes_2 b_k^{r_k}
	\right),
\end{multline*}
where $\degb{x_k}$ is the cohomological degree of $x_k$.

Also, when working with $\braidop_2$, the natural maps
\begin{equation*}
	A \otimes_R A^{\otimes_R r_1} \otimes_R \dots \otimes_R A \otimes_R A^{\otimes_R r_k}
	\rightarrow \ndf{A, \braidop_2}[k][]
\end{equation*}
given by
\begin{multline*}
	a_1 \otimes_R b_1^1 \otimes_R \dots \otimes_R b_1^{r_1} \otimes_R \dots \otimes_R
	a_k \otimes_R b_k^1 \otimes_R \dots \otimes_R b_k^{r_k} \mapsto
	\\
	\ul{a_1} \otimes_2 b_1^1 \otimes_2 \dots \otimes_2 b_1^{r_1} \otimes_2 \dots \otimes_2
	\ul{a_k} \otimes_2 b_k^1 \otimes_2 \dots \otimes_2 b_k^{r_k},
\end{multline*}
without an extra sign factor, are ill-defined. The maps $a \mapsto \ul{a}$ are not
$R$-linear and the tensor product $\otimes_R$ is not compatible with $\otimes_2$
when mixing regular and underlined elements. For example, we have
\begin{equation*}
	\begin{aligned}
		a \otimes_R \left( r \cdot b \right) \otimes_R c \mapsto
		\ul{a} \otimes_2 \left( r \cdot b \right) \otimes_2 \ul{c} & =
		(-1)^{\braid{(0,\degb{r})}{(1,\degb{a})}_{2}} \left( r \cdot \ul{a} \right) \otimes_2 b \otimes_2 \ul{c}
		\\
		                                                           & =
		(-1)^{\degb{r} \cdot \left( \degb{a} - 1 \right)} r \cdot \left( \ul{a} \otimes_2 b \otimes_2 \ul{c} \right),
	\end{aligned}
\end{equation*}
while
\begin{equation*}
	\begin{aligned}
		a \otimes_R b \otimes_R \left( r \cdot c \right) \mapsto
		\ul{a} \otimes_2 b \otimes_2 \ul{r \cdot c} & =
		(-1)^{\degb{r}} \ul{a} \otimes_2 b \otimes_2 \left( r \cdot \ul{c} \right)
		\\
		                                            & =
		(-1)^{\degb{r} + \braid{(0, \degb{r})}{(1, \degb{a} + \degb{b})}_2}
		\left( r \cdot \ul{a} \right) \otimes_2 b \otimes_2 \ul{c}
		\\
		                                            & =
		(-1)^{(\degb{a} + \degb{b}) \cdot \degb{r}} r \cdot \left( \ul{a} \otimes_2 b \otimes_2 \ul{c} \right),
	\end{aligned}
\end{equation*}
so that
\begin{equation*}
	(-1)^{\degb{r} \cdot \degb{b}} a \otimes_R b \otimes_R \left( r \cdot c \right) \mapsto
	(-1)^{\degb{r} \cdot \degb{a}} r \cdot \left( \ul{a} \otimes_2 b \otimes_2 \ul{c} \right) \neq
	(-1)^{\degb{r} \cdot \left( \degb{a} - 1 \right)} r \cdot \left( \ul{a} \otimes_2 b \otimes_2 \ul{c} \right).
\end{equation*}

The ill-definedness of the maps can be fixed by adding a sign factor. We choose
the sign so that the components
$\phi_k^{r_1,\dots,r_k} \colon
	A \otimes_R A^{\otimes_R r_1} \otimes_R \dots \otimes_R A \otimes_R A^{\otimes_R r_k}
	\rightharpoonup R$
of $\phi_k$ will coincide with the components of the corresponding $\phi^{\braidop_1}$,
giving us the expression
\begin{multline*}
	\phi_k^{r_1, \dots, r_k} \left(
	a_1 \otimes_R b_1^1 \otimes_R \dots \otimes_R b_1^{r_1} \otimes_R \dots \otimes_R
	a_k \otimes_R b_k^1 \otimes_R \dots \otimes_R b_k^{r_k}
	\right) =
	\\
	(-1)^{k \cdot \degb{x_k} + \varepsilon \left( x_k \right) + k - 1}
	\phi_k \left(
	\ul{a_1} \otimes_2 b_1^1 \otimes_2 \dots \otimes_2 b_1^{r_1} \otimes_2 \dots \otimes_2
	\ul{a_k} \otimes_2 b_k^1 \otimes_2 \dots \otimes_2 b_k^{r_k}
	\right),
\end{multline*}
where $\varepsilon \left( x_k \right)$ is given by \eqref{eq:sign-factor-Psi-cyc}.
Note that the sign $k \cdot \degb{x_k}$ gets cancelled when describing $\phi_{k}^{r_1, \dots r_k}$
directly in terms of $\phi_k$, i.e., we have
\begin{multline*}
	\phi_k^{r_1, \dots, r_k} \left(
	a_1 \otimes_R b_1^1 \otimes_R \dots \otimes_R b_1^{r_1} \otimes_R \dots \otimes_R
	a_k \otimes_R b_k^1 \otimes_R \dots \otimes_R b_k^{r_k}
	\right) =
	\\
	(-1)^{\varepsilon \left( x_k \right) + k - 1}
	\sigma_{4 - n} \phi \left(
	\ul{a_1} \otimes_2 b_1^1 \otimes_2 \dots \otimes_2 b_1^{r_1} \otimes_2 \dots \otimes_2
	\ul{a_k} \otimes_2 b_k^1 \otimes_2 \dots \otimes_2 b_k^{r_k}
	\right).
\end{multline*}

The situation is summarized in \cref{fig:total-inner-product-components-both-pairings}.
The various compositions in the bottom part of
\cref{fig:total-inner-product-components-both-pairings} are used to define the components
of $\phi^{\braidop_1}$ and do not involve any signs, while the upper part describes
the corresponding maps and signs involved in defining the components of $\phi^{\braidop_2}$.

\begin{landscape}
	\begin{figure}[h]
		\centering
		\adjustbox{scale=0.70,center}{
			\begin{tikzcd}[column sep=normal, row sep=huge]
				&& {\ul{A} \otimes_2 A^{\otimes_2 i_1} \otimes_2 \dots \otimes_2 \ul{A}
					\otimes_2 A^{\otimes_2 i_k}} &
				{\left( \ul{A} \otimes_2 A^{\otimes_2 i_1} \otimes_2
					\dots \otimes_2 \ul{A} \otimes_2 A^{\otimes_2 i_k} \right) \left< k \right>}
				\\
				&& {\ndf{\mathcal{A}, \braidop_2}[k][]} & {{\ndf{\mathcal{A}, \braidop_2}[k][]} \left< k \right>}
				\\
				&& {\ncdf{\mathcal{A}, \braidop_2}[k][]} & {{\ncdf{\mathcal{A}, \braidop_2}[k][]} \left< k \right>} &
				{\totc{\ncdfr{A, \braidop_2}, -\qdr, \clie{\mu}}[][\geq 2][\braidop_2]}
				\\
				{A \otimes_R A^{\otimes_R i_1} \otimes_R \dots \otimes_R A \otimes_R A^{\otimes_R i_k}} &&&&&
				{R[4-n]} & R
				\\
				&& {\ncdf{\mathcal{A}, \braidop_1}[k][]} & {{\ncdf{\mathcal{A}, \braidop_1}[k][]}[k]} &
				{\totc{\ncdfr{A, \braidop_1}, \qdr, \clie{\mu}}[][\geq 2][\braidop_1]}
				\\
				&& {\ndf{\mathcal{A}, \braidop_1}[k][]} & {{\ndf{\mathcal{A}, \braidop_1}[k][]}[k]}
				\\
				&& {\ul{A} \otimes_1 A^{\otimes_1 i_1} \otimes_1 \dots \otimes_1 \ul{A} \otimes_1
					A^{\otimes_1 i_k}} &
				{\left( \ul{A} \otimes_1 A^{\otimes_1 i_1} \otimes_1 \dots \otimes_1
					\ul{A} \otimes_1 A^{\otimes_1 i_k} \right) [k]}
				\arrow["(-1)^{k \cdot \degb{x_k}}", from=1-3, to=1-4]
				\arrow[from=1-3, to=2-3]
				\arrow[from=1-4, to=2-4]
				\arrow["(-1)^{k \cdot \degb{x_k}}", harpoon, from=2-3, to=2-4]
				\arrow[from=2-3, to=3-3]
				\arrow[from=2-4, to=3-4]
				\arrow["(-1)^{k \cdot \degb{x_k}}", harpoon, from=3-3, to=3-4]
				\arrow["(-1)^{k \cdot \degb{x_k} + \varepsilon \left( x_k \right) + (k - 1)}"', from=3-3, to=5-3]
				\arrow[from=3-4, to=3-5]
				\arrow["(-1)^{\varepsilon \left( x_k \right) + (k-1)}"', from=3-4, to=5-4]
				\arrow["{\rest{\phi^{\braidop_2}}{\totcomp{\mathcal{A}, \braidop_2}[2][]}}", from=3-5, to=4-6]
				\arrow["{\Psi_{\geq 2}}", "\sum_{k \geq 2} (-1)^{\varepsilon \left( x_k \right) + (k-1)}"', from=3-5, to=5-5]
				\arrow["(-1)^{k \cdot \degb{x_k} + \varepsilon \left( x_k \right) + (k - 1)}", from=4-1, to=1-3]
				\arrow[from=4-1, to=7-3]
				\arrow[harpoon, from=4-6, to=4-7]
				\arrow[from=5-3, to=5-4, harpoon]
				\arrow[from=5-4, to=5-5]
				\arrow["{\rest{\phi^{\braidop_1}}{\totcomp{\mathcal{A}, \braidop_1}[2][]}}"', from=5-5, to=4-6]
				\arrow[from=6-3, to=5-3]
				\arrow[from=6-3, to=6-4, harpoon]
				\arrow[from=6-4, to=5-4]
				\arrow[from=7-3, to=6-3]
				\arrow[from=7-3, to=7-4, harpoon]
				\arrow[from=7-4, to=6-4]
			\end{tikzcd}
		}
		\caption{Components of Total Inner Products for Both Parity Forms.}
		\label{fig:total-inner-product-components-both-pairings}
	\end{figure}
\end{landscape}

\section{Cyclic Structures and Strict Chain Maps on \texorpdfstring{$\ncdf{A}[2][]$}{Cyclic Codifferential 2-Forms}}
\label{appendix:cyclic-structures}

\counterwithin{thm}{section}
\counterwithin{equation}{section}

Let $\mathcal{A} = \left( A, \mu \right)$ be a Banach $\Ainf$-algebra over a
differential graded-commutative Banach $\mathbbm{k}$-algebra $\mathcal{R} = (R,d)$.
In this appendix, we establish the bijection between strict contractive chain maps
$\ncdf{\mathcal{A}}[2][] \rightarrow \mathcal{R}[2-n]$ and $n$-dimensional cyclic structures on
$\mathcal{A}$ in the sense of \cref{dfn:cyclic-structure}.

In what follows, we will work with the inner product pairing $\braidop_1$ given by \cref{eq:parity-inner-product}
to minimize sign conversions.
Recall from \cref{subsec:codifferential-forms} that the notation $\ul{A}$ is a shorthand for the shift $A[(-1,0)]$,
where we consider $A$ as a $\ZZ^2$-graded Banach module concentrated in the zeroth column.
Since we work with the symmetry $\braidop_1$, our
conventions regarding the shift dictate that the $R$-module structure on $\ul{A}$ is the same as the
$R$-module structure on $A$, because we shift in the line degree direction.
Namely, since the degree of the shift map
$s_{A[(-1,0)]} \colon A \rightharpoonup A[(-1,0)]$ is $(1,0)$, we have
\begin{equation*}
	r \cdot \ul{a} = r \cdot s_{A[(-1,0)]} \left( a \right) =
	(-1)^{\braid{\left( 0, \degb{r} \right)}{\left( 1, 0 \right)}_1} s_{A[(-1,0)]} \left( r \cdot a \right)
	= s_{A[(-1,0)]} \left(r \cdot a \right) = \ul{r \cdot a}.
\end{equation*}
In addition, the $\ZZ^2$-graded tensor product over $R$
\begin{equation*} \begin{aligned}
		\ul{A} \otimes_{R} \ul{A} & = \left( \ul{A} \otimes_{\mathbbm{k}} \ul{A} \right) /
		\gen{\left( r \cdot \ul{a_1} \right) \otimes \ul{a_2} - (-1)^{\braid{(0,\degb{r})}{(1,\degb{a_1})}_1}
			\ul{a_1} \otimes \left( r \cdot \ul{a_2} \right)}
		\\
		                          & =
		\left( \ul{A} \otimes_{\mathbbm{k}} \ul{A} \right) /
		\gen{\left( r \cdot \ul{a_1} \right) \otimes \ul{a_2} - (-1)^{\degb{r} \cdot \degb{a_1}}
			\ul{a_1} \otimes \left( r \cdot \ul{a_2} \right)}
	\end{aligned} \end{equation*}
can be identified, after forgetting the line degree, with the tensor product $A \otimes_{R} A$ of $\ZZ$-graded
modules over $R$. More generally, the tensor product $\ul{A} \otimes A^{\otimes i} \otimes \ul{A} \otimes A^{\otimes j}$,
computed in the category of bigraded $R$-modules can be identified with the tensor product
$A \otimes A^{\otimes i} \otimes A \otimes A^{\otimes j}$ of $\ZZ$-graded $R$-modules.

Given a graded contractive $R$-linear map $\phi \colon \ncdf{A}[2][] \rightharpoonup M$, let us denote by
\begin{equation*}
	\phi^{i,j} \colon A \otimes A^{\otimes i} \otimes A \otimes A^{\otimes j} \rightharpoonup M
\end{equation*}
the \textbf{components} of $\phi$, given by
\begin{equation*}
	\phi^{i,j} \left(
	a \otimes  b_1 \otimes \dots \otimes b_i \otimes c \otimes d_1 \otimes \dots \otimes d_j
	\right) =
	\phi \left(
	\ul{a} \otimes b_1 \otimes \dots \otimes b_i \otimes \ul{c} \otimes d_1 \otimes \dots \otimes d_j
	\right).
\end{equation*}
By the previous paragraph, the components are well-defined, and are contractive and $R$-linear of degree $\degb{\phi}$.
The sequence of components $\left( \phi^{i,j} \right)_{i,j \geq 0}$ determines the map $\phi$
completely but contains some redundancy.
Since $\phi$ is defined on $\ncdf{A}[2][]$, and we work with $\braidop_1$, we have
\begin{equation} \label{eq:components-phi-2-symmetry}
	\begin{aligned}
		\phi^{i,j} \left( a \otimes l \otimes b \otimes s \right) & =
		\phi \left( \ul{a} \otimes l \otimes \ul{b} \otimes s \right)
		\\
		                                                          & =
		(-1)^{\braid{\left( 1, \degb{a} + \degb{l} \right)}{\left( 1, \degb{b} + \degb{s} \right)}_{1}}
		\phi \left( \ul{b} \otimes s \otimes \ul{a} \otimes l \right)
		\\
		                                                          & =
		(-1)^{1 + \left( \degb{a} + \degb{l} \right) \left( \degb{b} + \degb{s} \right)}
		\phi \left( \ul{b} \otimes s \otimes \ul{a} \otimes l \right)
		\\
		                                                          & =
		(-1)^{1 + \left( \degb{a} + \degb{l} \right) \left( \degb{b} + \degb{s} \right)}
		\phi^{j,i} \left( b \otimes s \otimes a \otimes l \right).
	\end{aligned}
\end{equation}
Conversely, given a sequence of graded contractive $R$-linear maps $\left( \phi^{i,j} \right)_{i,j \geq 0}$
which satisfy \eqref{eq:components-phi-2-symmetry},
it determines a unique graded contractive $R$-linear map $\phi \colon \ncdf{A}[2][] \rightharpoonup M$
whose components are given by $\phi^{i,j}$.

Let us say that $\phi$ is \textbf{strict} if $\phi^{i,j} = 0$ for all $(i,j) \neq (0,0)$. Then we have:

\begin{lm} \label{lm:constant-cyclic-structure}
	Let $\inncur \colon A \otimes A \rightarrow R[2 - n]$ be an $n$-dimensional cyclic structure on $\mathcal{A}$.
	Then the map $\phi \colon \ncdf{A}[2][] \rightarrow R[2-n]$ determined by
	\begin{align}
		\phi^{0,0} \left( a \otimes b \right) \defeq \inncur[a][b],
		\label{eq:phi2-0-0-def}
		\\
		\phi^{i,j} \defeq 0, \qquad \forall (i,j) \neq (0,0)
		\label{eq:phi2-constant}
	\end{align}
	is a well-defined contractive map which satisfies $d_{\mathcal{R}[2-n]} \circ \phi = \phi \circ \clie{\mu}$.
	Conversely, given a strict contractive chain map
	$\phi \colon \ncdf{\mathcal{A}}[2][] \rightarrow \mathcal{R}[2-n]$,
	the map $\inncur \colon A \otimes A \rightarrow R[2 - n]$
	given by $\inncur[a][b] \defeq \phi^{0,0} \left( a \otimes b \right)$
	is an $n$-dimensional cyclic structure on $\mathcal{A}$.
\end{lm}
\begin{proof}
	First, note that a simple induction argument shows that the cyclic pairing property
	\eqref{eq:cyclic-pairing} is equivalent to the property
	\begin{equation} \begin{gathered} \label{eq:cyclic-pairing-extended}
			\inncur[\mu_k \left( a_1, \dots, a_k \right)][a_{k+1}] =
			\delta_{1,k} \cdot d_{\mathcal{R}[2-n]} \left( \inncur[a_1][a_2] \right) +
			\\
			(-1)^{\left( \degb{a_{m+1}} + \dots + \degb{a_{k+1}} \right) \cdot
				\left( \degb{a_1} + \dots + \degb{a_m} \right)}
			\inncur[\mu_k \left( a_{m+1}, \dots, a_{k+1}, a_1, \dots, a_{m-1} \right)][a_m]
		\end{gathered} \end{equation}
	for all $1 \leq m \leq k$ and $a_1,\dots,a_{k+1} \in A$.

	Given an $n$-dimensional cyclic structure $\inncur$, the antisymmetry of $\inncur$ guarantees that $\phi^{0,0}$ is
	antisymmetric and so determines a degree zero map $\phi \colon \ncdf{A}[2][] \rightarrow R[2-n]$,
	which is contractive as $\inncur$ is contractive.
	We need to show that $\phi$ is a chain map. Let $x \in \ncdf{A}[2][]$ be of the form
	\begin{equation*}
		x = \ul{a_{k+1}} \otimes {\underbrace{a_1 \otimes \dots \otimes a_{k-1}}_{l}} \otimes \ul{a_k}
	\end{equation*}
	where $k \geq 1$ and $a_1, \dots, a_{k+1} \in A$. Then
	\begin{equation} \label{eq:phi2-clie-calc}
		\begin{aligned}
			\clie{\mu} \left( x \right) \stackrel{\eqref{eq:cyc-lie-formula}}{=}
			 & \ul{\corest{\mu} \left( a_{k+1} \otimes l_{(1)} \right)} \otimes l_{(2)} \otimes \ul{a_k}
			\pm \ul{a_{k+1}} \otimes l_{(1)} \otimes \corest{\mu} \left( l_{(2)} \right) \otimes l_{(3)} \otimes
			\ul{a_{k}}
			\\
			 & \pm
			   (-1)^{\degb{a_{k+1}} + \degb{l_{(1)}}}
			\ul{a_{k+1}} \otimes l_{(1)} \otimes \ul{\corest{\mu} \left( l_{(2)} \otimes a_{k} \right)}
			\pm \ul{a_{k+1}} \otimes l \otimes \ul{a_k} \otimes \mu_0(1)
		\end{aligned}
	\end{equation}
	and so
	\begin{flalign}
		\phi \left( \clie{\mu} \left( x \right) \right)
		 & \eqwithref[eq:phi2-clie-calc][eq:phi2-constant]
		\phi \left( \ul{\mu_k \left( a_{k+1}, a_1, \dots, a_{k-1} \right)} \otimes \ul{a_k} \right) +
		                                                                                    (-1)^{\degb{a_{k+1}}} \phi \left( \ul{a_{k+1}} \otimes \ul{\mu_k \left( a_1, \dots, a_k \right)} \right)
		\notag
		 &                                                                                           \\
		 & \eqwithref[eq:phi2-0-0-def]
		\inncur[\mu_k \left( a_{k+1}, a_1, \dots, a_{k-1} \right)][a_k] +
		(-1)^{\degb{a_{k+1}}} \inncur[a_{k+1}][\mu_k \left( a_1, \dots, a_k \right)]
		\notag
		 &                                                                                           \\
		 & \eqwithref[eq:antisymmetry]
		\inncur[\mu_k \left( a_{k+1}, a_1, \dots, a_{k-1} \right)][a_k] +
		\notag
		 &                                                                                           \\
		 & \qquad\quad\quad
		              (-1)^{\degb{a_{k+1}} \cdot \left( \degb{a_1} + \dots + \degb{a_k} \right) + 1}
		\inncur[\mu_k \left( a_1, \dots, a_k \right)][a_{k+1}]
		\notag
		 &                                                                                           \\
		 & \eqwithref[eq:cyclic-pairing]
		\delta_{1,k} \cdot d_{\mathcal{R}[2-n]}
		\left( (-1)^{\degb{a_2} \cdot \degb{a_1} + 1} \inncur[a_1][a_2] \right)
		\label{eq:phi2-closed-calc}
		 &                                                                                           \\
		 & \eqwithref[eq:antisymmetry]
		\delta_{1,k} \cdot d_{\mathcal{R}[2-n]} \left( \inncur[a_2][a_1] \right)
		\notag
		 &                                                                                           \\
		 & \eqwithref[eq:phi2-0-0-def]
		d_{\mathcal{R}[2-n]} \left( \delta_{1,k} \cdot \phi \left( \ul{a_2} \otimes \ul{a_1} \right) \right)
		\notag
		 &                                                                                           \\
		 & \eqwithref[eq:phi2-constant]
		d_{\mathcal{R}[2-n]} \left( \phi \left( \ul{a_{k+1}} \otimes a_1 \otimes \dots \otimes a_{k-1} \otimes \ul{a_k} \right) \right) =
		d_{\mathcal{R}[2-n]} \left( \phi \left( x \right) \right).
		\notag
	\end{flalign}
	More generally, let $x \in \ncdf{A}[2][]$ be of the form
	\begin{equation*}
		x = \ul{a_{k+1}} \otimes \underbrace{a_1 \otimes \dots \otimes a_{m-1}}_l \otimes
		\ul{a_{m}} \otimes \underbrace{a_{m+1} \otimes \dots \otimes a_k}_s
	\end{equation*}
	where $1 \leq m \leq k$ and $a_1, \dots, a_{k+1} \in A$.
	A calculation similar to \eqref{eq:phi2-closed-calc},
	based on \cref{eq:cyclic-pairing-extended} instead of \cref{eq:cyclic-pairing}, shows that
	$\phi \left( \clie{\mu} \left( x \right) \right) =
		d_{\mathcal{R}[2-n]} \left( \phi \left( x \right) \right)$.

	Conversely, if $\phi$ is a strict contractive chain map, then $\inncur$ is clearly
	antisymmetric and contractive, and the calculation
	\eqref{eq:phi2-closed-calc} shows that $\inncur$ also satisfies the cyclic pairing property
	\eqref{eq:cyclic-pairing}.
\end{proof}

Thus, \cref{lm:constant-cyclic-structure} establishes a bijection
between $n$-dimensional cyclic structures on $\mathcal{A}$ and \textit{strict morphisms}
$\ncdf{\mathcal{A}}[2][] \rightarrow \mathcal{R}[2-n]$ of differential graded
Banach $\mathcal{R}$-modules.

\counterwithin{thm}{subsection}
\counterwithin{equation}{subsection}

\section{Cyclic Structure and Superpotential -- Relation to Previous Works} \label{appendix:sign-conversions-jake}

In this appendix we explain how to relate the notion
of a cyclic unital $\Ainf$-structure appearing in \cite[][Definition 1.1]{Solomon2016},
and the notion of the superpotential appearing in \cite[Page 7]{Solomon2016a}
to our definitions.
This will involve some sign conversions which we make explicit.

To start, let us comment on the difference between the conventions used in our work and the ones used
in \cite[][]{Solomon2016,Solomon2016a}. Let $C$ be a $\ZZ$-graded left module over
a graded $\mathbbm{k}$-algebra $R$. In our work, the shifted module $C[1]$ has the same underlying module
as $C$ with grading $C[1]^{*} = C^{* + 1}$, but we twist the module structure on $C[1]$
(see \cref{sec:suspension-graded-R-module}). This implies that the suspension map
$\s \colon C \rightharpoonup C[1]$, which has the identity as the underlying map, is
$R$-linear of degree $-1$ and the inverse map
$\sigma = \s^{-1} \colon C[1] \rightharpoonup C$ is $R$-linear of degree $1$:
\begin{align}
	\s \left( r \cdot c \right)                                & = (-1)^{\degb{r} \cdot (-1)} r \cdot \s \left( c \right)
	= (-1)^{\degb{r}} r \cdot \s \left( c \right), \quad       & c \in C,
	\label{eq:shift-linear}
	\\
	{\sigma} \left( r \cdot x \right)                          & = (-1)^{\degb{r} \cdot 1} r \cdot {\sigma} \left( x \right)
	= (-1)^{\degb{r}} r \cdot {\sigma} \left( x \right), \quad & x \in C[1].
	\label{eq:shift-inverse-linear}
\end{align}
In contrast, the suspension map does not appear explicitly in \cite[][]{Solomon2016,Solomon2016a}.
Elements $c \in C^{*}$ are considered freely as elements $c \in C[1]^{*-1}$, and the module structure
on $C[1]$ does not appear to be twisted. In addition, it appears that the (complete) tensor product $\otimes$
in \cite[][]{Solomon2016,Solomon2016a} is never taken over the graded algebra $R$, working
implicitly with $\otimes_{\mathbbm{k}}$.
We however assume that $R$ is graded-commutative and work directly with $\otimes_R$.
In what follows, we will freely use notation from \cite[][Section 1]{Solomon2016}.

\subsection{Sign Explanation} \label{sub:sign-explanation-jake-work}
Let $C = \cdiff{L} \otimes R$ and fix $\gamma \in \mathcal{I}_Q R$ with $d \gamma = 0$
and $\degb{\gamma} = 2$. In \cite{Solomon2016}, an $\Ainf$-structure is constructed on
$C$ using the moduli spaces of $J$-holomorphic disks.
Given $k \geq 0$ and $\beta \in \Pi$, define maps
$\mathfrak{m}^{\gamma,\beta}_k \colon C^{\otimes k} \rightharpoonup C$ by
\begin{equation}
	\begin{aligned}
		\mathfrak{m}^{\gamma,\beta}_k \left( c_1, \dots, c_k \right) \defeq
		 & (-1)^{1 + \sum_{j=1}^k j \cdot \left( \degb{c_j} + 1 \right)}
		\\
		 & \qquad \sum\limits_{l \geq 0} \frac{1}{l!} \left( {evb_0^\beta} \right)_{*} \left(
		                                                                               \bigwedge_{j=1}^l \left( evi_j^\beta \right)^{*} \gamma \wedge
		\bigwedge_{j=1}^k \left( evb_j^\beta \right)^{*} c_j
		\right)
		\label{eq:mathfrak-m-gamma-beta-k}
	\end{aligned}
\end{equation}
for $(k,\beta) \neq (1, \beta_0), (0, \beta_0)$ and
\begin{equation*}
	\mathfrak{m}^{\gamma, \beta_0}_1 \left( c \right) \defeq d c,
	\qquad
	\mathfrak{m}^{\gamma, \beta_0}_0 \defeq 0.
\end{equation*}
Define also $\mathfrak{m}_k^{\gamma} \colon C^{\otimes k} \rightharpoonup C$ by
\begin{equation}
	\mathfrak{m}_k^{\gamma} \defeq \sum_{\beta \in \Pi} \mathfrak{m}^{\gamma,\beta}_k T^{\beta}.
	\label{eq:mathfrak-m-gamma-k}
\end{equation}
Then $\Set{\mathfrak{m}_k^{\gamma}}_{k \geq 0}$ is an $\Ainf$-structure on $C$ in the sense of \cite{Solomon2016}.

Let us motivate the appearance of the signs
\begin{equation}
	\varepsilon \left( c_1, \dots, c_k \right) \defeq 1 + \sum_{j=1}^k j \cdot \left( \degb{c_j} + 1 \right)
	\label{eq:sign-factor-jake}
\end{equation}
in \cref{eq:mathfrak-m-gamma-beta-k}. Consider the maps
$\mathfrak{q}_{k,l}^{\gamma,\beta} \colon C^{\times k} \rightharpoonup C$ given by
\begin{equation*}
	\mathfrak{q}_{k,l}^{\gamma,\beta} \left( c_1, \dots, c_k \right) \defeq
	\begin{cases}
		\left( {evb_0^\beta} \right)_{*} \left(
		                                 \bigwedge_{j=1}^k \left( evb_j^\beta \right)^{*} c_j \wedge
		\bigwedge_{j=1}^l \left( evi_j^\beta \right)^{*} \gamma \right)
		T^{\beta}            & (k,l,\beta) \neq \substack{(1,0,\beta_0) \\ (0,0,\beta_0)},
		\\
		d \left( c_1 \right) & (k,l,\beta) = (1,0,\beta_0),
		\\
		0                    & (k,l,\beta) = (0,0,\beta_0),
	\end{cases}
\end{equation*}
without any sign factors.\footnote{The maps
	$\mathfrak{q}_{k,l}^{\gamma,\beta}$ are almost the same as the maps
	$\mathfrak{q}_{k,l}^{\beta} \left( - ; \gamma^{\times l} \right)$
	from \cite[][Page 12]{Solomon2016}, the only difference being that $\mathfrak{q}_{k,l}^{\gamma,\beta}$
	includes the $T^{\beta}$ factor and lacks the sign factor.}
Then in terms of the maps $\mathfrak{q}_{k,l}^{\gamma,\beta}$ we have
\begin{equation*}
	\mathfrak{m}_k^{\gamma} \left( c_1, \dots, c_k \right) =
	\begin{cases}
		\sum\limits_{(l,\beta) \in \NZ \times \Pi} \frac{1}{l!}
		                                           (-1)^{\varepsilon \left( c_1, \dots, c_k \right)}
		\mathfrak{q}_{k,l}^{\gamma,\beta} \left( c_1, \dots, c_k \right)
		 & k \neq 1,
		\\
		d c_1 +
		\sum\limits_{(l,\beta) \in \NZ \times \Pi \setminus \Set{(0,\beta_0)}} \frac{1}{l!}
		                                                                       (-1)^{\varepsilon \left( c_1 \right)}
		\mathfrak{q}_{1,l}^{\gamma,\beta} \left( c_1 \right)
		 & k = 1.
	\end{cases}
\end{equation*}
The maps $\mathfrak{q}_{k,l}^{\gamma,\beta}$ are $\mathbbm{k}$-multilinear of degree $2 - k$.
As a consequence of the conventions chosen in \cite{Solomon2016},\footnote{The conventions
	in \cite{Solomon2016} dictate that the differential $d$ is extended so that it becomes $R$-linear of degree one.
	Given a proper submersion $f \colon M \rightarrow N$, the pullback map
	$f^{*} \colon \cdiff{N}[*] \rightarrow \cdiff{M}[*]$ on $R$-valued differential forms
	is $R$-linear of degree zero, while the pushforward map
	$f_{*} \colon \cdiff{M}[*] \rightharpoonup \cdiff{N}[* - \left( \dim M - \dim N \right)]$
	has degree $\dim N - \dim M$ but is chosen to be $R$-linear of \textit{degree zero}.}
they satisfy
\begin{equation*}
	\mathfrak{q}_{k,l}^{\gamma,\beta} \left( c_1, \dots, r \cdot c_i, \dots, c_k \right) =
	\begin{cases}
		(-1)^{\degb{r} \cdot \left( \sum_{j=1}^{i-1} \degb{c_j} \right)}
		r \cdot \mathfrak{q}_{k,l}^{\gamma,\beta} \left( c_1, \dots, c_k \right)       &
		(k,l,\beta) \neq (1,0,\beta_0),
		\\
		(-1)^{\degb{r}} r \cdot \mathfrak{q}_{1,0}^{\gamma,\beta_0} \left( c_1 \right) &
		(k,l,\beta) = (1,0,\beta_0).
	\end{cases}
\end{equation*}
Hence, the maps $\mathfrak{q}_{k,l}^{\gamma,\beta}$ for $(k,l,\beta) \neq (0,0,\beta_0)$ \textit{are not}
$R$-multilinear of degree $2 - k$. Instead, they look like $R$-multilinear maps of degree \textit{zero}.
By replacing
$\mathfrak{q}_{k,l}^{\gamma, \beta} \left( c_1, \dots, c_k \right)$ with
$(-1)^{\left( 2 - k \right) \cdot \left( \sum_{j=1}^k \degb{c_j} \right)} \mathfrak{q}_{k,l}^{\gamma,\beta}
	\left( c_1, \dots, c_k \right)$, we can obtain honest $R$-multilinear operators of degree $2 - k$.
Thus, we can define $R$-linear operators $m_k \colon C^{\otimes_R k} \rightharpoonup C$ of degree $2 - k$
by
\begin{equation*}
	m_k^{\gamma} \left( c_1, \ldots, c_k \right) \defeq
	\begin{cases}
		\sum\limits_{(l,\beta) \in \NZ \times \Pi} \frac{1}{l!}
		                                           (-1)^{\left( 2 - k \right) \cdot \left( \sum_{j=1}^k \degb{c_j} \right)}
		\mathfrak{q}_{k,l}^{\gamma,\beta} \left( c_1, \dots, c_k \right) &
		k \neq 1,
		\\
		d \left( c_1 \right) +
		\sum\limits_{(l,\beta) \in \NZ \times \Pi \setminus \Set{(0,\beta_0)}} \frac{1}{l!}
		                                                                       (-1)^{\degb{c_1}} \mathfrak{q}_{1,l}^{\gamma,\beta} \left( c_1 \right)
		                                                                 & k = 1.
	\end{cases}
\end{equation*}

The operators $\Set{m_k^{\gamma}}_{k \geq 0}$ satisfy the identities of \cref{eq:ainf_for_m_k2} for
an unshifted $\Ainf$-algebra and give us an $\Ainf$-perturbation of the standard DGA structure on
$C$ in the sense that
\begin{equation*}
	m_1^{\gamma} \left( c \right) = dc + \cdots, \qquad
	m_2^{\gamma} \left( c_1, c_2 \right) = c_1 \wedge c_2 + \cdots.
\end{equation*}
To convert $\Set{m_k^{\gamma}}_{k \geq 0}$ to a shifted $\Ainf$-algebra,
as discussed in \cref{sub:a-inf-sign-conventions},
we can use the transformations of \cref{fig:translate-mu_k-m_k-start-m_k} and set
$\mu_k \defeq -\s \circ m_k^{\gamma} \circ {\sigma}^{\otimes k}$.
Then
\begin{equation*}
	\begin{aligned}
		\sigma \mu_k \left( \s c_1, \dots, \s c_k \right) & = (-1)^{1 + \sum_{j=1}^k \left( k - j \right) \degb{\s c_j}}
		m_k \left( c_1, \dots, c_k \right)
		\\
		                                                  & =
		\begin{cases}
			\sum\limits_{(l,\beta) \in \NZ \times \Pi} \frac{1}{l!}
			                                           (-1)^{1 + \sum_{j=1}^k \left( k - j \right) \degb{\s c_j} +
				                                           \left( 2 - k \right) \cdot \left( \sum_{j=1}^k \degb{c_j} \right)}
			\mathfrak{q}_{k,l}^{\gamma,\beta} \left( c_1, \dots, c_k \right)
			 & k \neq 1,
			\\
			-d \left( c_1 \right) +
			\sum\limits_{(l,\beta) \in \NZ \times \Pi \setminus \Set{(0,\beta_0)}} \frac{1}{l!}
			                                                                       (-1)^{1 + \degb{c_1}} \mathfrak{q}_{1,l}^{\gamma,\beta} \left( c_1 \right)
			 & k = 1.
		\end{cases}
		\\
		                                                  & =
		\begin{cases}
			\sum\limits_{(l,\beta) \in \NZ \times \Pi} \frac{1}{l!}
			                                           (-1)^{k + \varepsilon \left( c_1, \dots, c_k \right)}
			\mathfrak{q}_{k,l}^{\gamma,\beta} \left( c_1, \dots, c_k \right)
			 & k \neq 1,
			\\
			-d \left( c_1 \right) +
			\sum\limits_{(l,\beta) \in \NZ \times \Pi \setminus \Set{(0,\beta_0)}} \frac{1}{l!}
			                                                                       (-1)^{1 + \varepsilon \left( c_1 \right)}
			\mathfrak{q}_{1,l}^{\gamma,\beta} \left( c_1 \right)
			 & k = 1,
		\end{cases}
	\end{aligned}
\end{equation*}
and hence
\begin{equation}
	\mu_k \left( \s c_1, \dots, \s c_k \right) = (-1)^{k} \s \mathfrak{m}_k^{\gamma} \left( c_1, \dots, c_k \right).
	\label{eq:relation-mu_k-mathfrak-m_k-concrete}
\end{equation}
Thus, up to an immaterial factor of $(-1)^k$, we recover the operators $\mathfrak{m}^{\gamma}_k$ together
with the signs as in \cref{eq:mathfrak-m-gamma-beta-k,eq:mathfrak-m-gamma-k}.

\begin{rem}
	Let $f \colon M \rightarrow N$ be a proper submersion.
	By changing the definition of the pushforward map
	$f_{*} \colon \cdiff{M}[*] \rightharpoonup \cdiff{N}[* - \left( \dim M - \dim N \right)]$
	so that it becomes left $R$-linear instead of right $R$-linear, i.e., satisfy
	\begin{equation*}
		f_{*} \left( r \cdot \alpha \right) = f_{*} \left( f^{*} \left( r \right) \wedge \alpha \right)
		= (-1)^{\degb{r} \cdot \left( \dim M - \dim N \right)} r \wedge f_{*} \left( \alpha \right)
		= (-1)^{\degb{r} \cdot \degb{f_{*}}} r \cdot f_{*} \left( \alpha \right)
	\end{equation*}
	for $\alpha \in \cdiff{M}, r \in \cdiff{N}$, we can get rid of
	the signs appearing in $m_k^{\gamma}$ altogether and still obtain an unshifted $\Ainf$-algebra.

	The conclusion is that the $\Ainf$-relations are a consequence of the consistent choice
	of orientations for the moduli spaces used in the definitions of the operators.
	The hard thing about obtaining the $\Ainf$-relations is to consistently orient the moduli spaces.
	The $\varepsilon$ sign factor given by \eqref{eq:sign-factor-jake} has no intrinsic geometric meaning
	and is merely an algebraic factor used to adapt the operations $\mathfrak{q}_{k,l}^{\gamma,\beta}$
	to adhere to certain conventions.
\end{rem}

\subsection{Converting the \texorpdfstring{$\Ainf$-algebra}{A-infinity Algebra} Structure}
\label{subsec:converting-a-inf-jake-to-banach}
Let $\mathcal{R} = (R,d)$ be a differential graded $\mathbbm{k}$-algebra endowed with a valuation
$\varsigma_R$, and let $C$ be a graded module over $R$ with valuation $\varsigma_C$
as in \cite[][Section 1.2]{Solomon2016}. Although the definitions and the precise properties of the valuations
are left unspecified in \cite[][]{Solomon2016}, we follow a reasonable interpretation and assume that
the valuations can be used to define seminorms in the standard way via
$\nnorm[r]_R \defeq e^{-\varsigma_R \left( r \right)}, \nnorm[c]_C \defeq e^{-\varsigma_C \left( c \right)}$.
We assume that the resulting seminorms are compatible with the algebra and module structures and
that $R$ is graded-commutative. We also assume that $\mathcal{R}$ becomes a differential graded
\textit{Banach} $\mathbbm{k}$-algebra, and $\left( C, \nnorm_C \right)$ becomes a graded \textit{Banach} module over
$\left( R, \nnorm_R \right)$.

Let $\left( \{ \mathfrak{m}_k \}_{k \ge 0}, \inncur, \be \right)$
be an $n$-dimensional curved cyclic unital $\Ainf$-structure on $C$ as given by \cite[][Definition 1.1]{Solomon2016}.
We assume that the unit $\be$ satisfies $\varsigma_C \left( \be \right) \geq 0$.\footnote{We note
	that all our assumptions hold for the explicit cyclic unital $\Ainf$-algebra constructed
	in \cite[][]{Solomon2016} and discussed in \cref{sub:sign-explanation-jake-work}.}
Starting with $ \{ \mathfrak{m}_k \}_{k \ge 0}$, define a sequence of operations
$\mu_k \colon C[1]^{\times k} \rightharpoonup C[1]$
by the formula
\begin{equation} \label{def:conv-jake-mk-uk}
	\mu_k \left( a_1, \dots, a_k \right) \defeq
	(-1)^k  \s \left( \mathfrak{m}_k \left( \sigma a_1, \dots, \sigma a_k \right) \right).
\end{equation}
In the case $k = 0$, we have $\mu_0 \left( 1 \right) = \s \mathfrak{m}_0$.

\begin{lm}	\label{lm:convert-jake-mk-uk-linearity}
	The operations $\mu_k$ for $k \neq 1$ are $R$-multilinear of degree one, while $\mu_1$ is
	a derivation over $d$. More explicitly, we have
	\begin{equation*}
		\mu_k \left( a_1, \dots, a_{i-1}, r \cdot a_{i}, \dots, a_k \right) =
		\delta_{1,k} \cdot dr \cdot a_1 +
			(-1)^{\degb{r} \cdot \left( 1 + \sum_{j=1}^{i-1} \degb{a_j} \right)}
		r \cdot \mu_k \left( a_1, \dots, a_k \right)
	\end{equation*}
	for $a_1, \dots, a_k \in C[1]$ and $r \in R$.
\end{lm}
\begin{proof}
	Since $\mathfrak{m}_k$ are of degree $2 - k$, the operators $\mu_k$ are of degree one.
	Using property $(1)$ of \cite[][Definition 1.1]{Solomon2016} we have
	\begin{equation*}
		\begin{aligned}
			\mu_k \left( a_1, \dots, a_{i-1}, r \cdot a_{i}, \dots, a_k \right)
			\eqwithref[def:conv-jake-mk-uk]     &
			(-1)^k \cdot \s \left(
			\mathfrak{m}_k \left( \sigma a_1, \dots, \sigma \left( r \cdot a_i \right),
			\dots, \sigma a_k \right) \right)
			\\
			\eqwithref[eq:shift-inverse-linear] &
			(-1)^k \cdot \s \left(
			                \mathfrak{m}_k \left( \sigma a_1, \dots,
			                (-1)^{\degb{r}} r \cdot \left( \sigma a_i \right), \dots, \sigma a_k \right) \right)
			\\
			\eqwithtext[(P1)]                   &
			(-1)^{k + \degb{r}} \cdot \s \left( \delta_{1,k} \cdot dr \cdot \sigma a_1 \right) +
			\\
			                                    &
			(-1)^{k + \degb{r}} \cdot \s \left(
			                             (-1)^{\degb{r} \cdot \left(i + \sum_{j=1}^{i-1} \degb{\sigma a_j} \right)}
			r \cdot \mathfrak{m}_k \left( \sigma a_1, \dots, \sigma a_k \right) \right)
			\\
			\eqwithref[eq:shift-linear]         &
			\delta_{1,k} \cdot (-1)^{k + \degb{r} + \degb{dr}} \cdot dr \cdot a_1 +
			\\
			                                    &
			(-1)^{\degb{r} \cdot \left( 1 + \sum_{j=1}^{i-1} \degb{a_j} \right)} r \cdot
			(-1)^k \cdot \s \left( \mathfrak{m}_k \left( \sigma a_1, \dots, \sigma a_k \right) \right)
			\\
			\eqwithref[def:conv-jake-mk-uk]     &
			\delta_{1,k} \cdot dr \cdot a_1 +
				                            (-1)^{\degb{r} \cdot \left( 1 + \sum_{j=1}^{i-1} \degb{a_j} \right)}
			r \cdot \mu_k \left( a_1, \dots, a_k \right).
		\end{aligned}
	\end{equation*}
\end{proof}

\begin{lm} \label{lm:convert-jake-mk-uk-continuity}
	The operators $\mu_k$ satisfy
	\begin{equation*}
		\nnorm[\mu_k \left( a_1, \dots, a_k \right)] \leq \nnorm[a_1] \cdots \nnorm[a_k]
	\end{equation*}
	for all $k \geq 1$ and $a_1, \dots, a_k \in C[1]$. We also have $\nnorm[\mu_0 \left( 1 \right)] < 1$.
\end{lm}
\begin{proof}
	Using property $(4)$ of \cite[][Definition 1.1]{Solomon2016} we have
	\begin{equation*}
		\begin{aligned}
			\nnorm[\mu_k \left( a_1, \dots, a_k \right)] &
			\stackrel{\eqref{def:conv-jake-mk-uk}}{=}{}
			\nnorm[\mathfrak{m}_k \left( \sigma a_1, \dots, \sigma a_k \right)] =
			e^{-\varsigma_C \left( \mathfrak{m}_k \left( \sigma a_1, \dots, \sigma a_k \right) \right)}
			\\
			                                             & \stackrel{\,\,\text{(P4)}\,\,}{\leq}{}
			e^{-\sum_{i=1}^k \varsigma_C \left( \sigma a_i \right)} =
			\nnorm[\sigma a_1] \cdots \nnorm[\sigma a_k] = \nnorm[a_1] \cdots \nnorm[a_k],
		\end{aligned}
	\end{equation*}
	and
	\begin{equation*}
		\nnorm[\mu_0] = \nnorm[\s \mathfrak{m}_0] = \nnorm[\mathfrak{m}_0]
		= e^{-\varsigma_C \left( \mathfrak{m}_0 \right)} \stackrel{\,\,\text{(P4)}\,\,}{<} e^0 = 1.
	\end{equation*}
\end{proof}

\Cref{lm:convert-jake-mk-uk-linearity,lm:convert-jake-mk-uk-continuity} imply that the maps
$a_1 \otimes_R \dots \otimes_R a_k \mapsto \mu_k \left( a_1, \dots, a_k \right)$ are well-defined
and so induce degree
one contractive maps $\mu_k \colon C[1]^{\cotimes_R k} \rightharpoonup C[1]$. Let us denote by
$\mu \colon \tens{C[1]}[R] \rightharpoonup \tens{C[1]}[R]$ the degree one coderivation
over $d$ whose components are given by $\mu_k$.

\begin{lm}
	The coderivation $\mu$ satisfies $\mu^2 = 0$.
\end{lm}
\begin{proof}
	We need to show that $\mu_k$ satisfy the $\Ainf$-relations of \cref{eq:ainf_for_mu_k_explicit}.
	This will follow directly from the $\Ainf$-relations for $\mathfrak{m}_k$
	(property $(3)$ of \cite[][Definition 1.1]{Solomon2016}).
	We have
	\begin{equation*}
		\mu_1 \left( \mu_0 \left( 1 \right) \right) = \mu_1 \left( \s \, \mathfrak{m}_0 \right) =
		-\s \, \mathfrak{m}_1 \left( \sigma \left( \s \, \mathfrak{m}_0  \right) \right)
		= -\s \, \mathfrak{m}_1 \left( \mathfrak{m}_0 \right) = 0.
	\end{equation*}
	The $\Ainf$-relations for $k \geq 1$ follow similarly. Given 	$a_1, \dots, a_k \in C[1]$,
	set $\alpha_i \defeq \sigma a_i$ so that we have
	$\degb{a_i} = \degb{\alpha_i} - 1 \equiv \degb{\alpha_i} + 1 \mod 2$ and
	$\mu_k \left( a_1, \dots, a_k \right) =
		(-1)^k \cdot \s \left( \mathfrak{m}_k \left( \alpha_1, \dots, \alpha_k \right) \right)$.
	Then
	\begin{gather*}
		\sum_{\substack{k_1 + k_2 = k + 1 \\ 1 \leq i \leq k_1}} (-1)^{\sum_{j=1}^{i-1} \degb{a_j}}
		\mu_{k_1} \left( a_1, \dots, a_{i-1}, \mu_{k_2} \left( a_i, \dots, a_{i+k_2-1} \right),
		a_{i+k_2}, \dots, a_k \right)
		\stackrel{\eqref{def:conv-jake-mk-uk}}{=}
		\\
		\sum_{\substack{k_1 + k_2 = k + 1 \\ 1 \leq i \leq k_1}} (-1)^{\sum_{j=1}^{i-1} \degb{a_j}}
		\mu_{k_1} \left( a_1, \dots, a_{i-1},
		(-1)^{k_2} \cdot \s \left( \mathfrak{m}_{k_2} \left( \alpha_i, \dots, \alpha_{i+k_2-1} \right) \right),
		a_{i+k_2}, \dots, a_k \right)
		\stackrel{\eqref{def:conv-jake-mk-uk}}{=}
		\\
		\sum_{\substack{k_1 + k_2 = k + 1 \\ 1 \leq i \leq k_1}} (-1)^{\sum_{j=1}^{i-1} \left( \degb{\alpha_j} + 1 \right) + k_1 + k_2}
		\s \, \mathfrak{m}_{k_1} \left( \alpha_1, \dots, \alpha_{i-1},
		\mathfrak{m}_{k_2} \left( \alpha_i, \dots, \alpha_{i+k_2-1} \right),
		\alpha_{i+k_2}, \dots, \alpha_k \right) =
		\\
		(-1)^{k+1} \s \left( \sum_{\substack{k_1 + k_2 = k + 1 \\ 1 \leq i \leq k_1}} (-1)^{\sum_{j=1}^{i-1} \left( \degb{\alpha_j} + 1 \right)}
		\mathfrak{m}_{k_1} \left( \alpha_1, \dots, \alpha_{i-1},
		\mathfrak{m}_{k_2} \left( \alpha_i, \dots, \alpha_{i+k_2-1} \right),
		\alpha_{i+k_2}, \dots, \alpha_k \right) \right),
	\end{gather*}
	which vanishes by the $\Ainf$-relations for $\mathfrak{m}_k$.
\end{proof}

Let us set $A = C[1]$. We have shown that
\begin{cor}
	The pair $\mathcal{A} = \left( A, \mu \right)$ is a (shifted)
	Banach $\Ainf$-algebra over $\mathcal{R} = (R, d)$ in the sense
	of \cref{def:a-inf-Banach-algebra}.
	\qed
\end{cor}

Next, we convert the unit of the $\Ainf$-algebra. Set
\begin{equation} \label{eq:conv-jake-unit}
	e \defeq \s \be \in C[1].
\end{equation}

\begin{lm} \label{lm:conv-jake-unit}
	The element $e \in A^{-1}$ is a unit for $\mathcal{A}$ in the sense of
	\cref{dfn:a-infinity-unit}.
\end{lm}
\begin{proof}
	By our assumption on $\be$, we have
	$\nnorm[e] = \nnorm[\be] = e^{-\varsigma_C \left( \be \right)} \leq e^0 = 1$.
	Using property $(8)$ of \cite[][Definition 1.1]{Solomon2016}, we have
	\begin{equation*}
		\begin{aligned}
			\mu_k \left( a_1, \dots, a_{i-1}, e, a_{i+1}, \dots, a_k \right)
			\eqwithref[def:conv-jake-mk-uk] &
			(-1)^k \s \mathfrak{m}_k \left( \sigma a_1, \dots, \sigma a_{i-1}, \be,
			\sigma a_{i+1}, \dots, \sigma a_k \right)
			\\
			\eqwithtext[(P8)]               &
			0
		\end{aligned}
	\end{equation*}
	for $k \neq 2$.
	By property $(10)$ of \cite[][Definition 1.1]{Solomon2016}, we have
	\begin{align*}
		\mu_2 \left( a, e \right) & \stackrel{\eqref{def:conv-jake-mk-uk}}{=}
		\s \mathfrak{m}_2 \left( \sigma a, \sigma e \right) =
		\s \mathfrak{m}_2 \left( \sigma a, \be \right) \stackrel{\text{(P10)}}{=}
		\s \sigma a = a,
		\\
		\mu_2 \left( e, a \right) & \stackrel{\eqref{def:conv-jake-mk-uk}}{=}
		\s \mathfrak{m}_2 \left( \sigma e, \sigma a \right) =
		\s \mathfrak{m}_2 \left( \be, \sigma a \right) \stackrel{\text{(P10)}}{=}
		\s (-1)^{\degb{\sigma a}} \sigma a = (-1)^{\degb{a} + 1} a.
	\end{align*}
\end{proof}

\begin{rem}
	Our choice of the factor $(-1)^k$ in the translation between $\mathfrak{m}_k$ and $\mu_k$
	of \cref{def:conv-jake-mk-uk} is motivated by the following points:
	\begin{enumerate}
		\item Since we work over a differential graded algebra, our conventions regarding the shifted
		      modules force us to define
		      \begin{equation*}
			      \mu_1 \left( a \right) = - \s \mathfrak{m}_1 \left( \sigma a \right)
		      \end{equation*}
		      with a minus sign to guarantee that $\mu_1$ is a derivation over $d$ and not over $-d$.
		\item We have
		      \begin{equation*}
			      \mu_2 \left( \s c_1, \s c_2 \right) = \s \mathfrak{m}_2 \left( c_1, c_2 \right)
		      \end{equation*}
		      without a sign which guarantees that the unit for $\left( C[1], \mu \right)$ is $\s \be$ and not
		      $-\s \be$.
		\item With the factor $(-1)^k$, the relation between $\mu_k$ and $\mathfrak{m}_k$ is consistent
		      with the discussion in \cref{sub:sign-explanation-jake-work}
		      (see \cref{eq:relation-mu_k-mathfrak-m_k-concrete}).
	\end{enumerate}

	We note that one can also choose to define
	\begin{equation*}
		\mu_k \left( a_1, \dots, a_k \right) \defeq
		- \s \left( \mathfrak{m}_k \left( \sigma a_1, \dots, \sigma a_k \right) \right)
	\end{equation*}
	and still obtain an $\Ainf$-structure over $(R,d)$ with unit $- \s \be$. See also the
	discussion in \cref{sub:a-inf-sign-conventions}.
\end{rem}

\subsection{Converting the Cyclic Structure} \label{sec:converting-jake-to-cyclic-structure}
Define an operator $\varphi \colon A \times A \rightharpoonup R$
by the formula
\begin{equation} \label{eq:def-conv-jake-phi}
	\varphi \left( a_1, a_2 \right) \defeq
	(-1)^{(1-n) \cdot \left( \degb{a_1} + \degb{a_2} \right)} \prec \sigma a_1, \sigma a_2 \succ.
\end{equation}

\begin{lm} \label{lm:convert-jake-phi-linearity}
	The operator $\varphi$ is $R$-multilinear of degree $2-n$. More explicitly, we have
	\begin{equation*}
		\begin{aligned}
			\varphi \left( r \cdot a_1, a_2 \right) & =
			(-1)^{\degb{r} \cdot (2-n)} r \cdot \varphi \left( a_1, a_2 \right),
			\\
			\varphi \left( a_1, r \cdot a_2 \right) & =
			(-1)^{\degb{r} \cdot \left(2 - n + \degb{a_1} \right)} r \cdot \varphi \left( a_1, a_2 \right).
		\end{aligned}
	\end{equation*}
\end{lm}
\begin{proof}
	Using property $(2)$ of \cite[][Definition 1.1]{Solomon2016}, we have
	\begin{equation*}
		\begin{aligned}
			\varphi \left( r \cdot a_1, a_2 \right)
			\eqwithref[eq:def-conv-jake-phi]    &
			(-1)^{(1-n) \cdot \left( \degb{r} + \degb{a_1} + \degb{a_2} \right)}
			\prec \sigma \left( r \cdot a_1 \right), \sigma a_2 \succ
			\\
			\eqwithref[eq:shift-inverse-linear] &
			(-1)^{\degb{r} \cdot (2-n) + (1-n) \cdot \left( \degb{a_1} + \degb{a_2} \right)}
			\prec r \cdot \left( \sigma a_1 \right), \sigma a_2 \succ
			\\
			\eqwithtext[(P2)]                   &
			(-1)^{\degb{r} \cdot (2-n) + (1-n) \cdot \left( \degb{a_1} + \degb{a_2} \right)}
			r \cdot \prec \sigma a_1, \sigma a_2 \succ
			\\
			\eqwithref[eq:def-conv-jake-phi]    &
			(-1)^{\degb{r} \cdot (2-n)} r \cdot \varphi \left( a_1, a_2 \right).
		\end{aligned}
	\end{equation*}
	Similarly, we have
	\begin{equation*}
		\begin{aligned}
			\varphi \left( a_1, r \cdot a_2 \right)
			\eqwithref[eq:def-conv-jake-phi]    &
			(-1)^{(1-n) \cdot \left( \degb{r} + \degb{a_1} + \degb{a_2} \right)}
			\prec \sigma a_1 , \sigma \left( r \cdot a_2 \right) \succ
			\\
			\eqwithref[eq:shift-inverse-linear] &
			(-1)^{\degb{r} \cdot (2-n) + (1-n) \cdot \left( \degb{a_1} + \degb{a_2} \right)}
			\prec \sigma a_1 , r \cdot \left( \sigma a_2 \right) \succ
			\\
			\eqwithtext[(P2)]                   &
			(-1)^{\degb{r} \cdot (2-n) + (1-n) \cdot \left( \degb{a_1} + \degb{a_2} \right) +
				\degb{r} \cdot \left( \degb{\sigma a_1} + 1 \right)}
			r \cdot \prec \sigma a_1, \sigma a_2 \succ
			\\
			\eqwithref[eq:def-conv-jake-phi]    &
			(-1)^{\degb{r} \cdot \left(2 - n + \degb{a_1} \right)} r \cdot \varphi \left( a_1, a_2 \right).
		\end{aligned}
	\end{equation*}
\end{proof}

\begin{lm} \label{lm:convert-jake-phi-continuity}
	The operator $\varphi$ satisfies
	\begin{equation*}
		\nnorm[\varphi \left( a_1, a_2 \right)] \leq \nnorm[a_1] \nnorm[a_2]
	\end{equation*}
	for $a_1, a_2 \in C[1]$.
\end{lm}
\begin{proof}
	Using property $(5)$ of \cite[][Definition 1.1]{Solomon2016}, we have
	\begin{equation*}
		\begin{aligned}
			\nnorm[\varphi \left( a_1, a_2 \right)] & \stackrel{\eqref{eq:def-conv-jake-phi}}{=}
			\nnorm[\prec \sigma a_1, \sigma a_2 \succ] =
			e^{-\varsigma_R \left( \prec \sigma a_1, \sigma a_2 \succ \right)}
			\\
			                                        & \stackrel{\,\,\text{(P5)}\,\,}{\leq}
			e^{- \left( \varsigma_C \left( \sigma a_1 \right) +  \varsigma_C \left( \sigma a_2 \right) \right)} =
			\nnorm[\sigma a_1] \nnorm[\sigma a_2] = \nnorm[a_1] \nnorm[a_2].
		\end{aligned}
	\end{equation*}
\end{proof}

\Cref{lm:convert-jake-phi-linearity,lm:convert-jake-phi-continuity} imply that $\varphi$
induces a well-defined contractive $R$-linear map
\begin{equation*}
	\varphi \colon A \otimes A \rightharpoonup R
\end{equation*}
of degree $2 - n$
given by $a_1 \otimes a_2 \mapsto \varphi \left( a_1, a_2 \right)$.

\begin{lm} \label{lm:convert-jake-phi-antisymmetry}
	The operator $\varphi$ is antisymmetric:
	\begin{equation} \label{eq:convert-jake-phi-antisymmetry}
		\varphi \left( a_1, a_2 \right) = (-1)^{\degb{a_1} \cdot \degb{a_2} + 1}
		\varphi \left( a_2, a_1 \right).
	\end{equation}
\end{lm}
\begin{proof}
	Using property $(6)$ of \cite[][Definition 1.1]{Solomon2016}, we have
	\begin{equation*}
		\begin{aligned}
			\varphi \left( a_1, a_2 \right) & \eqwithref[eq:def-conv-jake-phi]
			                                  (-1)^{(1-n) \cdot \left( \degb{a_1} + \degb{a_2} \right)} \prec \sigma a_1, \sigma a_2 \succ
			\\
			                                & \eqwithtext[(P6)]
			                                  (-1)^{(1-n) \cdot \left( \degb{a_1} + \degb{a_2} \right) +
				                                  \left( \degb{\sigma a_1} + 1 \right) \cdot \left( \degb{\sigma a_2} + 1 \right) + 1}
			\prec \sigma a_2, \sigma a_1 \succ
			\\
			                                & \eqwithref[eq:def-conv-jake-phi]
			                                  (-1)^{\degb{a_1} \cdot \degb{a_2} + 1} \varphi \left( a_2, a_1 \right).
		\end{aligned}
	\end{equation*}
\end{proof}

\begin{lm}
	The operator $\varphi$ satisfies the cyclic pairing property:
	\begin{equation}
		\begin{aligned}\label{eq:jake-phi-cyclic-symmetry}
			\varphi \left( \mu_k \left( a_1, \dots, a_k \right), a_{k+1} \right) ={} &
			(-1)^{\degb{a_{k+1}} \cdot \left( \degb{a_1} + \dots + \degb{a_k} \right)}
			\varphi \left( \mu_k \left( a_{k+1}, a_1, \dots, a_{k-1} \right), a_k \right)
			\\
			                                                                         & +
			\delta_{1,k} \cdot (-1)^{2-n} d \left( \varphi \left( a_1, a_2 \right) \right).
		\end{aligned}
	\end{equation}
\end{lm}
\begin{proof}
	Set
	\begin{equation*} \begin{aligned}
			\varepsilon \defeq & (1-n) \cdot \left( \degb{a_1} + \dots + \degb{a_{k+1}} + 1 \right)                      \\
			={}                & (1-n) \cdot \left( \degb{\mu_k \left( a_1, \dots, a_k \right)} + \degb{a_{k+1}} \right)
			= (1-n) \cdot \left( \degb{\mu_k \left( a_{k+1}, a_1, \dots, a_{k-1} \right)} + \degb{a_k} \right).
		\end{aligned} \end{equation*}
	Then using property (7) of \cite[][Definition 1.1]{Solomon2016}, we have
	\begin{equation*} \begin{aligned}
			\varphi \left( \mu_k \left( a_1, \dots, a_k \right), a_{k+1} \right)
			\eqwithref[eq:def-conv-jake-phi] &
			(-1)^{(1-n) \cdot \left( \degb{a_1} + \dots + \degb{a_{k+1}} + 1 \right)}
			\prec \sigma \mu_k \left( a_1, \dots, a_k \right), \sigma a_{k+1} \succ
			\\
			\eqwithref[def:conv-jake-mk-uk]  &
			(-1)^{\varepsilon + k}
			\prec \mathfrak{m}_k \left( \sigma a_1, \dots, \sigma a_k \right), \sigma a_{k+1} \succ
			\\
			\eqwithtext[(P7)]                &
			(-1)^{\varepsilon + k +
				\degb{a_{k+1}} \cdot \left( \degb{a_1} + \dots + \degb{a_k} \right)}
			\prec \mathfrak{m}_k \left( \sigma a_{k+1}, \sigma a_1, \dots, \sigma a_{k-1} \right), \sigma a_k \succ
			\\
			                                 & \qquad + (-1)^{\varepsilon + k} \delta_{1,k} \cdot d \prec \sigma a_1, \sigma a_2 \succ
			\\
			\eqwithref[def:conv-jake-mk-uk]  &
			(-1)^{\degb{a_{k+1}} \cdot \left( \degb{a_1} + \dots + \degb{a_k} \right) + \varepsilon}
			\prec \sigma \mu_k \left( a_{k+1}, a_1, \dots, a_{k-1} \right), \sigma a_k \succ
			\\
			                                 & \qquad+
			                                   (-1)^{(1-n) \cdot \left( \degb{a_1} + \degb{a_2} \right) + (1 - n) + 1}
			\delta_{1,k} \cdot d \prec \sigma a_1, \sigma a_2 \succ
			\\
			\eqwithref[eq:def-conv-jake-phi] &
			(-1)^{\degb{a_{k+1}} \cdot \left( \degb{a_1} + \dots + \degb{a_k} \right)}
			\varphi \left( \mu_k \left( a_{k+1}, a_1, \dots, a_{k-1} \right), a_k \right)
			\\
			                                 & \qquad+
			                                   (-1)^{2-n} \delta_{1,k} \cdot d \left( \varphi \left( a_1, a_2 \right) \right).
			\\
		\end{aligned} \end{equation*}
\end{proof}

\begin{cor}
	The map $A \otimes A \rightarrow R[2-n]$ given by
	$a_1 \otimes a_2 \mapsto \s_{2-n} \varphi \left( a_1, a_2 \right)$
	is an $n$-dimensional cyclic structure on the $\Ainf$-algebra $\mathcal{A}$
	in the sense of \cref{dfn:cyclic-structure}.
\end{cor}

\begin{lm} \label{lm:jake-phi-m0-unit}
	The operator $\varphi$ satisfies
	\begin{equation} \label{eq:jake-phi-m0-unit}
		\varphi \left( \mu_0 \left( 1 \right), e \right) = 0.
	\end{equation}
\end{lm}
\begin{proof}
	We have
	\begin{equation*}
		\varphi \left( \mu_0 \left( 1 \right), e \right)
		\stackrel{\eqref{eq:def-conv-jake-phi}}{=}
		\prec \sigma \mu_0 \left( 1 \right), \sigma e \succ
		\stackrel[\eqref{eq:conv-jake-unit}]{\eqref{def:conv-jake-mk-uk}}{=}
		\prec \mathfrak{m}_0, \be \succ
		\stackrel{\text{(P9)}}{=}
		0.
	\end{equation*}
\end{proof}

\begin{rem}
	Applying \cref{cor:cyclic-structure-gives-pre-infinity-trace},
	we see that the map $\theta_1 \colon A \rightharpoonup R$
	given by
	\begin{equation*}
		\theta_1 \left( a \right) = \varphi \left( e, a \right) =
		(-1)^{\left( 1 - n \right) \left( \degb{a} - 1 \right)}
		\inncur[\sigma e][\sigma a] =
		(-1)^{\left( 1 - n \right) \left( \degb{a} - 1 \right)}
		\inncur[\be][\sigma a]
	\end{equation*}
	is an $n$-dimensional trace on $\mathcal{A}$ in the sense of
	\cref{eq:infty-trace-strict-rel}.
	By \cref{lm:jake-phi-m0-unit}, the map $\theta \colon \ncdf{A}[0][] \rightharpoonup R[1-n]$
	whose only non-zero component is $\theta_1$, gives a strict
	$n$-dimensional trace on $\mathcal{A}$ in the sense of \cref{dfn:infty-trace}.
	When $n > 0$, the cyclic structure $\inncur$ of \cite{Solomon2016} satisfies
	$\inncur[\be][\be] = 0$, which implies that
	$\theta_1 \left( e \right) = \inncur[\be][\be] = 0$, i.e.,
	$\theta$ is unital in the sense of \cref{dfn:infty-trace-unital}.
\end{rem}

\subsection{Converting the Superpotential}
In what follows, we will work with the inner product pairing $\braidop_1$ given by \cref{eq:parity-inner-product}.
\cref{lm:constant-cyclic-structure} shows that if we define a map
$\phi_2 \colon \ncdf{A}[2][] \rightarrow R[2-n]$
by specifying its components to be
\begin{equation*}
	\phi_2^{0,0} \left( \ul{a_1} \otimes \ul{a_2} \right) = \varphi \left( a_1, a_2 \right),
	\quad \phi_2^{k,l} = 0 \quad \forall (k,l) \neq (0,0),
\end{equation*}
then $\phi_2$ is a chain map,
i.e., an $n$-dimensional (strict) pre-homotopy inner product.
To upgrade $\phi_2$ to a homotopy inner product, we need extra data.
The cyclic $\Ainf$-structures constructed in~\cite{Solomon2016} all come equipped
with a constant term $\mathfrak{m}_{-1}$ of degree $3 - n$,
coming from the moduli space of $J$-holomorphic
disks with zero boundary marked points. The term $\mathfrak{m}_{-1}$ is not
part of the definition of a cyclic $\Ainf$-structure, but it satisfies
$\varsigma_R \left( \mathfrak{m}_{-1} \right) > 0$ (in our language,
$\nnorm[\mathfrak{m}_{-1}] < 1$), and is proven to be related to the cyclic structure
in the case of pseudoisotopies.

Working in the context of~\cite{Solomon2016} with the $\Ainf$-structures
$\mathfrak{m}_k = \mathfrak{m}_k^{\gamma}$ described there, define a map
$\phi \colon \totcompe{A}[2][] \rightarrow R[4-n]$ by specifying
that the only non-zero components of $\phi$ are $\phi_2$ and
$\phi_{\ul{1}} \defeq (-1)^{2-n} \mathfrak{m}_{-1}$.

To verify that $\phi$ is a homotopy inner product,
we need to check the extra identity
\begin{equation*}
	(-1)^{2 - n} d \left( \phi_{\ul{1}} \right) = d \left( \mathfrak{m}_{-1} \right) =
	-\frac{1}{2} \phi_2 \left( \ul{\mu_0 \left( 1 \right)}, \ul{\mu_0 \left( 1 \right)} \right)
\end{equation*}
of \eqref{eq:d-phi-ul-1-id}. We have
\begin{equation*}
	\begin{aligned}
		-\frac{1}{2} \phi_2 \left( \ul{\mu_0 \left( 1 \right)}, \ul{\mu_0 \left( 1 \right)} \right)
		\stackrel{\phantom{\eqref{eq:def-conv-jake-phi}}}{=}{} &
		-\frac{1}{2} \varphi \left( \mu_0 \left( 1 \right), \mu_0 \left( 1 \right) \right)
		\\
		\stackrel{\eqref{eq:def-conv-jake-phi}}{=}{}           &
		-\frac{1}{2} (-1)^{\left( 1 - n \right) \cdot
			 \left( \degb{\mu_0 \left( 1 \right)} + \degb{\mu_0 \left( 1 \right)} \right)}
		\inncur[\sigma \mu_0 \left( 1 \right)][\sigma \mu_0 \left( 1 \right)]
		\\
		\stackrel{\eqref{def:conv-jake-mk-uk}}{=}{}            &
		-\frac{1}{2}
		\inncur[ \mathfrak{m}_0  ][ \mathfrak{m}_0 ].
	\end{aligned}
\end{equation*}

When $C = \cdiff{L} \otimes R$, it can be shown that
$\inncur[ \mathfrak{m}_0  ][ \mathfrak{m}_0 ] = 0$ if $\eqcl{L} = 0 \in \homo{X}[n]$
or if $L$ admits a strong bounding cochain. See~\cite[Proposition~2.15]{Solomon2016} and~\cite[Corollary~4.2]{Solomon2024}. In this case, we work over $R$ so the differential
$d = 0$ is trivial and the extra identity holds. When working with pseudoisotopies,
we have $C = \cdiff{I \times L} \otimes R$, and we work over the differential
graded-commutative algebra $\cdiff{I} \otimes R$ of differential forms on the
interval. In this case, \cite[Proposition 4.20]{Solomon2016} shows that
\begin{equation} \label{eq:d-m-minus-1-GW}
	d \mathfrak{m}_{-1} = -\frac{1}{2}
	\inncur[ \mathfrak{m}_0  ][ \mathfrak{m}_0 ] + \widetilde{GW},
\end{equation}
where $\widetilde{GW}$ are extra terms coming from pseudoholomorphic spheres,
so strictly speaking, $\phi$ is not a homotopy inner product.

Finally, we check that our superpotential $\SP$
associated to $\phi$ coincides with the standard superpotential
$\hat{\Omega}$ of \cite[Page 7]{Solomon2016a}.
Let $\mathbf{b} \in C^1$ be a weak bounding cochain for $\mathfrak{m}$ in the sense that
$\mathbf{b} \in \mathcal{I}_R \left( C \right)$, and
$\sum_{k \geq 0} \mathfrak{m}_k \left( \mathbf{b}^{\otimes k} \right) = c \cdot \be$
for $c \in \mathcal{I}_R$ as in~\cite[Definition 1.1]{Solomon2016a}.
Set $b \defeq -\s \mathbf{b} \in A^0 = C[1]^0$. Then
$\nnorm[b] < 1$ and
\begin{equation*}
	\corest{\mu} \left( \Exp{b} \right) =
	\sum_{k \geq 0} \mu_k \left( b^{\otimes k} \right) =
	\sum_{k \geq 0} (-1)^k \mu_k \left( \left( \s \mathbf{b} \right)^{\otimes k} \right)
	\stackrel{\eqref{def:conv-jake-mk-uk}} =
	\s \left( \sum_{k \geq 0}  \mathfrak{m}_k \left( \mathbf{b}^{\otimes k} \right) \right)
	= \s \left( c \cdot \be \right) = c \cdot e,
\end{equation*}
so $b$ is a weak bounding cochain for $\mu$.
Evaluating $\SP[b]$ in terms of $\mathfrak{m}_k$ and $\mathbf{b}$, we have
\begin{equation*}
	\begin{aligned}
		\SP[b]
		\eqwithref[eq:sp-explicit-formula-homotopy-inner-product-braidop-1] &
		\phi_{\ul{1}} +
		\sum_{k \geq 0} \frac{1}{k+1} \phi_2^{0,0}
		\left( \ul{\mu_k \left( b^{\otimes k} \right)}, \ul{b} \right)
		\\
		\eqwithref                                                          &
		(-1)^{2-n} \cdot \mathfrak{m}_{-1} +
		\sum_{k \geq 0} \frac{1}{k+1}
		\varphi \left( \mu_k \left( b^{\otimes k} \right), b \right)
		\\
		\eqwithref[eq:def-conv-jake-phi]                                    &
		(-1)^{2 - n} \cdot \mathfrak{m}_{-1} +
		\sum_{k \geq 0} (-1)^{(1 - n) \cdot \left( \degb{\mu_k \left( b^{\otimes k} \right)} + \degb{b} \right)}
		\frac{1}{k+1} \inncur[\sigma \mu_k \left( b^{\otimes k} \right)][\sigma b]
		\\
		\eqwithref[def:conv-jake-mk-uk]                                     &
		(-1)^{2 - n} \cdot \mathfrak{m}_{-1} +
		                   (-1)^{1 - n} \sum_{k \geq 0}
		\frac{1}{k+1}
		\inncur[(-1)^k \mathfrak{m}_k \left( \left( \sigma b \right)^{\otimes k} \right)][\sigma b]
		\\
		\eqwithref                                                          &
		(-1)^{2 - n} \cdot \mathfrak{m}_{-1} +
		                   (-1)^{n} \sum_{k \geq 0}
		\frac{1}{k+1} \inncur[\mathfrak{m}_k \left( \mathbf{b}^{\otimes k} \right)][\mathbf{b}]
		\\
		\eqwithref                                                          &
		(-1)^n \left(
		\mathfrak{m}_{-1} +
		\sum_{k \geq 0}
		\frac{1}{k+1} \inncur[\mathfrak{m}_k \left( \mathbf{b}^{\otimes k} \right)][\mathbf{b}]
		\right) = \hat{\Omega} \left( \mathbf{b} \right).
	\end{aligned}
\end{equation*}

\begin{rem}
	Because of the $\widetilde{GW}$ term in~\eqref{eq:d-m-minus-1-GW}, $\phi$ fails to be a homotopy inner product in the case of pseudoisotopies. Consequently, the standard superpotential $\hat{\Omega}$ is not pseudoisotopy invariant unless terms which come from sphere bubbling are removed.
	Following the ideas of \cite{Solomon2024}, one can generalize the
	notion of a total inner product, allowing the total inner product to take values in
	a module instead of the ground algebra. One can then reinterpret $\phi$
	as an honest total inner product, and our superpotential $\SP$ as the
	superpotential $\Omega$ of \cite{Solomon2016a}, which is pseudoisotopy invariant.
\end{rem}

\printbibliography

\end{document}